\documentclass{amsbook}

\newtheorem{thm}{Theorem}[section]
\newtheorem*{thm*}{Theorem}

\newtheorem{claim}[thm]{Claim}
\newtheorem{cond}[thm]{Condition}
\newtheorem{cor}[thm]{Corollary}

\newtheorem{lem}[thm]{Lemma}
\newtheorem*{lem*}{Lemma}
\newtheorem{mainthm}{Theorem}
\newtheorem*{mainthm*}{Theorem}
\newtheorem{maincor}[mainthm]{Corollary}

\newtheorem{mainprop}[mainthm]{Proposition}

\newtheorem{prop}[thm]{Proposition}

\theoremstyle{definition}

\newtheorem*{case*}{Case}

\newtheorem{defn}[thm]{Definition}
\newtheorem*{defn*}{Definition}
\newtheorem{exmp}[thm]{Example}
\newtheorem*{exmp*}{Example}

\newtheorem{mainconj}[mainthm]{Conjecture}
\newtheorem{maindefn}[mainthm]{Definition}

\renewcommand{\thestep}{}

\theoremstyle{remark}

\newtheorem{case}{Case}\renewcommand{\thecase}{}

\numberwithin{subcase}{case}

\numberwithin{subsubcase}{subcase}

\newtheorem{rmk}[thm]{Remark}
\newtheorem*{rmk*}{Remark}

\makeatletter
\def\alphenumi{
  \def\theenumi{\alph{enumi}}
  \def\p@enumi{\theenumi}
  \def\labelenumi{(\@alph\c@enumi)}}
\makeatother

\makeatletter
\def\thecase{\@arabic\c@case}
\makeatother

\makeatletter
\def\thestep{\@arabic\c@step}
\makeatother

\newcommand{\transv}{\mathrel{\text{\tpitchfork}}}
\makeatletter
\newcommand{\tpitchfork}{%
  \vbox{
    \baselineskip\z@skip
    \lineskip-.52ex
    \lineskiplimit\maxdimen
    \m@th
    \ialign{##\crcr\hidewidth\smash{$-$}\hidewidth\crcr$\pitchfork$\crcr}
  }%
}
\makeatother

\DeclareFontFamily{U}{mathx}{\hyphenchar\font45}
\DeclareFontShape{U}{mathx}{m}{n}{
      <5> <6> <7> <8> <9> <10>
      <10.95> <12> <14.4> <17.28> <20.74> <24.88>
      mathx10
      }{}
\DeclareSymbolFont{mathx}{U}{mathx}{m}{n}
\DeclareFontSubstitution{U}{mathx}{m}{n}
\DeclareMathAccent{\widecheck}{0}{mathx}{"71}
\DeclareMathAccent{\wideparen}{0}{mathx}{"75}

\newcount\hh
\newcount\mm
\mm=\time
\hh=\time
\divide\hh by 60
\divide\mm by 60
\multiply\mm by 60
\mm=-\mm
\advance\mm by \time
\def\hhmm{\number\hh:\ifnum\mm<10{}0\fi\number\mm}

\let\oldmarginpar\marginpar
\renewcommand\marginpar[1]{\-\oldmarginpar[\raggedleft\footnotesize #1]%
{\raggedright\footnotesize #1}}

\newcommand\embed{\hookrightarrow}

\newcommand\barz{{\bar z}}

\newcommand\ubarCC{{\underline{\mathbb{C}}}}

\newcommand\ubarRR{{\underline{\mathbb{R}}}}

\newcommand\CC{\mathbb{C}}

\newcommand\EE{\mathbb{E}}
\newcommand\FF{\mathbb{F}}

\newcommand\HH{\mathbb{H}}

\newcommand\KK{\mathbb{K}}

\newcommand\NN{\mathbb{N}}

\newcommand\PP{\mathbb{P}}
\newcommand\QQ{\mathbb{Q}}
\newcommand\RR{\mathbb{R}}

\newcommand\ZZ{\mathbb{Z}}

\newcommand\cE{{\mathcal{E}}}
\newcommand\cF{{\mathcal{F}}}

\newcommand\cK{{\mathcal{K}}}

\newcommand\cM{{\mathcal{M}}}

\newcommand\cS{{\mathcal{S}}}
\newcommand\cT{{\mathcal{T}}}
\newcommand\cU{{\mathcal{U}}}

\newcommand\cV{{\mathcal{V}}}
\newcommand\cW{{\mathcal{W}}}
\newcommand\cX{{\mathcal{X}}}

\newcommand\fg{{\mathfrak{g}}}

\newcommand\fH{{\mathfrak{H}}}

\newcommand\fK{{\mathfrak{K}}}

\newcommand\fM{{\mathfrak{M}}}

\newcommand\fs{{\mathfrak{s}}}
\newcommand\fS{{\mathfrak{S}}}
\newcommand\ft{{\mathfrak{t}}}

\newcommand\fu{{\mathfrak{u}}}

\newcommand\fV{{\mathfrak{V}}}

\newcommand\sA{{\mathscr{A}}}
\newcommand\sB{{\mathscr{B}}}
\newcommand\sC{{\mathscr{C}}}
\newcommand\sD{{\mathscr{D}}}
\newcommand\sE{{\mathscr{E}}}
\newcommand\sF{{\mathscr{F}}}
\newcommand\sG{{\mathscr{G}}}
\newcommand\sH{{\mathscr{H}}}
\newcommand\sI{{\mathscr{I}}}

\newcommand\sK{{\mathscr{K}}}
\newcommand\sL{{\mathscr{L}}}
\newcommand\sM{{\mathscr{M}}}

\newcommand\sO{{\mathscr{O}}}

\newcommand\sR{{\mathscr{R}}}
\newcommand\sS{{\mathscr{S}}}
\newcommand\sT{{\mathscr{T}}}
\newcommand\sU{{\mathscr{U}}}
\newcommand\sV{{\mathscr{V}}}
\newcommand\sW{{\mathscr{W}}}
\newcommand\sX{{\mathscr{X}}}
\newcommand\sY{{\mathscr{Y}}}
\newcommand\sZ{{\mathscr{Z}}}

\newcommand\balpha{{\boldsymbol{\alpha}}}

\newcommand\bga{{\boldsymbol{\gamma}}}
\newcommand\bgamma{{\boldsymbol{\gamma}}}
\newcommand\bchi{{\boldsymbol{\chi}}}

\newcommand\bkappa{{\boldsymbol{\kappa}}}

\newcommand\bmu{{\boldsymbol{\mu}}}

\newcommand\bomega{{\boldsymbol{\omega}}}

\newcommand\bB{{\mathbf{B}}}

\newcommand\bg{{\mathbf{g}}}
\newcommand\bG{{\mathbf{G}}}

\newcommand\bh{{\mathbf{h}}}
\newcommand\bH{{\mathbf{H}}}

\newcommand\bJ{{\mathbf{J}}}

\newcommand\bL{{\mathbf{L}}}

\newcommand\bM{{\mathbf{M}}}

\newcommand\bx{{\mathbf{x}}}
\newcommand\bX{{\mathbf{X}}}

\newcommand\bY{{\mathbf{Y}}}

\newcommand\bZ{{\mathbf{Z}}}

\newcommand{\cov}{\nabla}

\newcommand{\rd}{\partial}

\newcommand\eps{\varepsilon}

\newcommand\la{\lambda}
\newcommand\La{\Lambda}
\newcommand\ka{\kappa}
\newcommand\om{\omega}
\newcommand\Om{\Omega}
\newcommand\si{\sigma}

\let\gimel\relax

\DeclareFontFamily{U}{rcjhbltx}{}
\DeclareFontShape{U}{rcjhbltx}{m}{n}{<->rcjhbltx}{}
\DeclareSymbolFont{hebrewletters}{U}{rcjhbltx}{m}{n}
\DeclareMathSymbol{\alef}{\mathord}{hebrewletters}{39}
\DeclareMathSymbol{\bet}{\mathord}{hebrewletters}{98}
\DeclareMathSymbol{\gimel}{\mathord}{hebrewletters}{103}
\DeclareMathSymbol{\dalet}{\mathord}{hebrewletters}{100}
\DeclareMathSymbol{\heh}{\mathord}{hebrewletters}{104}
\DeclareMathSymbol{\vav}{\mathord}{hebrewletters}{119}
\DeclareMathSymbol{\zayin}{\mathord}{hebrewletters}{122}
\DeclareMathSymbol{\chet}{\mathord}{hebrewletters}{120}
\DeclareMathSymbol{\tet}{\mathord}{hebrewletters}{84}
\DeclareMathSymbol{\yod}{\mathord}{hebrewletters}{121}
\DeclareMathSymbol{\kaf}{\mathord}{hebrewletters}{107}
\DeclareMathSymbol{\lamed}{\mathord}{hebrewletters}{108}
\DeclareMathSymbol{\mem}{\mathord}{hebrewletters}{109}
\DeclareMathSymbol{\nun}{\mathord}{hebrewletters}{110}
\DeclareMathSymbol{\samech}{\mathord}{hebrewletters}{115}
\DeclareMathSymbol{\ayin}{\mathord}{hebrewletters}{96}
\DeclareMathSymbol{\peh}{\mathord}{hebrewletters}{112}
\DeclareMathSymbol{\tzadi}{\mathord}{hebrewletters}{118}
\DeclareMathSymbol{\kof}{\mathord}{hebrewletters}{113}
\DeclareMathSymbol{\resh}{\mathord}{hebrewletters}{114}
\DeclareMathSymbol{\shin}{\mathord}{hebrewletters}{115}
\DeclareMathSymbol{\taw}{\mathord}{hebrewletters}{116}

\newcommand\gl{{\mathfrak{g}\mathfrak{l}}}
\newcommand\fsl{{\mathfrak{s}\mathfrak{l}}}

\newcommand\su{{\mathfrak{s}\mathfrak{u}}}
\newcommand\GL{\operatorname{GL}}
\newcommand\Or{\operatorname{O}}

\newcommand\SL{\operatorname{SL}}
\newcommand\SO{\operatorname{SO}}

\newcommand\SU{\operatorname{SU}}
\newcommand\U{\operatorname{U}}

\newcommand\less{\setminus}

\newcommand\ad{{\operatorname{ad}}}

\DeclareMathOperator{\Ann}{Ann}

\DeclareMathOperator{\Aut}{Aut}

\newcommand\CCl{\operatorname{{\mathbb{C}\ell}}}
\DeclareMathOperator{\Coh}{Coh}

\DeclareMathOperator{\codim}{codim}

\newcommand\Coker{\operatorname{Coker}}

\DeclareMathOperator{\Crit}{Crit}

\newcommand\dist{\operatorname{dist}}

\newcommand\Dom{\operatorname{Dom}}

\DeclareMathOperator{\Epim}{Epim}
\newcommand\End{\operatorname{End}}

\DeclareMathOperator{\Euler}{Euler}

\newcommand\Fred{\operatorname{Fred}}

\newcommand\grad{\operatorname{grad}}
\newcommand\Graph{\operatorname{Graph}}

\newcommand\hess{\operatorname{hess}}
\newcommand\Hess{\operatorname{Hess}}

\DeclareMathOperator{\Hom}{Hom}

\DeclareMathOperator{\Imag}{Im}

\DeclareMathOperator{\Ind}{Index}

\DeclareMathOperator{\Isom}{Isom}

\newcommand\Ker{\operatorname{Ker}}

\newcommand\Length{\operatorname{Length}}

\newcommand\Map{\operatorname{Map}}

\newcommand\Met{\operatorname{Met}}

\newcommand\PD{\operatorname{PD}}

\DeclareMathOperator{\pr}{pr}

\newcommand\Ran{\operatorname{Ran}}
\newcommand\rank{\operatorname{rank}}

\newcommand\Real{\operatorname{Re}}
\newcommand\Red{\operatorname{Red}}

\newcommand\Ric{\operatorname{Ric}}

\newcommand\Stab{\operatorname{Stab}}

\DeclareMathOperator{\SW}{SW}
\DeclareMathOperator{\Sym}{Sym}

\newcommand\tr{\operatorname{tr}}
\newcommand\Tr{\operatorname{Tr}}

\newcommand\vol{\operatorname{vol}}
\newcommand\Vol{\operatorname{Vol}}

\newcommand\can{{\mathrm{can}}}

\newcommand\disc{{\mathrm{disc}}}

\newcommand\ess{{\mathrm{ess}}}

\DeclareMathOperator{\expdim}{\mathrm{exp\ dim}}

\newcommand\id{{\mathrm{id}}}

\newcommand\loc{{\mathrm{loc}}}

\newcommand\mutatis{{\emph{mutatis mutandis }}}

\newcommand\ps{{\mathrm{ps}}}
\newcommand\red{{\mathrm{red}}}

\newcommand\reg{{\mathrm{reg}}}

\newcommand\sss{{\mathrm{ss}}}

\newcommand\spinc{\text{$\mathrm{spin}^c$ }}
\newcommand\spinu{\text{$\mathrm{spin}^u$ }}
\newcommand\Spinc{\text{$\mathrm{Spin}^c$}}

\newcommand\vir{\mathrm{vir}}
\numberwithin{equation}{section}
\numberwithin{section}{chapter}
\numberwithin{figure}{section}

\usepackage{centernot}
\usepackage{graphicx}
\usepackage{hyperref}
\usepackage{mathabx}
\usepackage{mathrsfs}
\usepackage[usenames]{color}
\usepackage{paralist}
\usepackage{scrextend}
\usepackage{slashed}
\usepackage{txfonts}
\usepackage{url}
\usepackage{verbatim}
\usepackage{xr}

\usepackage{amsmath}
\usepackage{amscd}
\usepackage{amssymb}
\usepackage{lmodern}

\newcommand{\seminorm}[1]{\left\lvert\hspace{-1 pt}\left\lvert\hspace{-1 pt}\left\lvert#1\right\lvert\hspace{-1 pt}\right\lvert\hspace{-1 pt}\right\lvert}

\hypersetup{pdftitle={Almost Hermitian structures on moduli spaces of non-Abelian monopoles and applications to the topology of symplectic four-manifolds}}
\hypersetup{pdfauthor={Paul M. N. Feehan and Thomas G. Leness}}

\begin{document}

\frontmatter

\title[Almost Hermitian structures on moduli spaces of non-Abelian monopoles]{Almost Hermitian structures on moduli spaces of non-Abelian monopoles and applications to the topology of symplectic four-manifolds}

\author[Paul M. N. Feehan]{Paul M. N. Feehan}
% \address{Department of Mathematics, Rutgers, The State University of New Jersey, 110 Frelinghuysen Road, Piscataway, NJ 08854-8019, United States}
% \email{feehan@math.rutgers.edu}
% \urladdr{\url{https://math.rutgers.edu/~feehan}}
\author[Thomas G. Leness]{Thomas G. Leness}
% \address{Department of Mathematics, Florida International University, Miami, FL 33199, United States}
% \email{lenesst@fiu.edu}
% \urladdr{\url{http://faculty.fiu.edu/~lenesst}}

\dedicatory{}

%COMMENT Remove hours and minutes for arXiv and journal versions and fix date
%\date{\today{ }\hhmm}
%\date{This version: October 17, 2024}

% AMS 2010 subject classifications (used in AMS journals)


% AMS keywords (used in AMS journals)

\maketitle

% Change page number to 7 if a dedication is present.
\setcounter{page}{7}

\tableofcontents
%\listoffigures

\chapter*{Preface}
\label{chap:Preface}
This work is a sequel to our previous monograph \cite{Feehan_Leness_introduction_virtual_morse_theory_so3_monopoles}, where we initiated our program to prove that the Bogomolov--Miyaoka--Yau inequality holds for closed, symplectic four-manifolds and, more generally, for closed, smooth four-manifolds with a Seiberg--Witten basic class. This inequality was first proved for compact, complex surfaces of general type by Miyaoka (1977) and Yau (1978). Our program to prove these conjectures uses a version of Morse theory for a natural Hamiltonian, the square of the $L^2$ norm of the coupled spinors, for the circle action on the moduli space of non-Abelian monopoles over a closed four-manifold. It has the aim of proving the existence of a projectively anti-self-dual connection on a rank-two Hermitian vector bundle over a blow-up of the four-manifold, where the first Pontrjagin number of the vector bundle is negative and greater than or equal to minus the Euler characteristic of the blown-up four-manifold. Our Morse theory arguments rely on positivity of virtual Morse--Bott indices for critical points of Hamiltonians for circle actions on complex analytic spaces (or real analytic spaces that, locally, are sufficiently well-approximated by complex analytic model spaces), as developed by the first author in \cite{Feehan_analytic_spaces}. In our application to the moduli space of non-Abelian monopoles, the critical points are fixed points of the circle action and thus represented by Seiberg--Witten monopoles.

\chapter*{Acknowledgments}
\label{chap:Acknowledgments}
We are grateful to Ian Agol, Denis Auroux, Fernando Chamizo, Kai Cieliebak, Simon Donaldson, Adam Earnst, Mariano Echeverria, Nicolas Ginoux, Fabio Gironella, Kristen Hendricks, Henryk Iwaniec, Dieter Kotschick, Peter Kronheimer, John Lott, Feng Luo, Vicente Mu\~noz, Pedro Ontaneda, Julie Rowlett, Daniel Ruberman, 
% TL12-18-2025: Is this Nik Saveliev? or someone else?
%PF12-18-2025 Someone else.
Nikhil Savale, Cliff Taubes, Richard Thomas, Richard Wentworth, Chris Woodward, and Zhengyi Zhou for helpful communications during the preparation of this work. We are extremely grateful to Cliff Taubes for kindly alerting us to a serious error in an earlier version of this monograph. We thank the National Science Foundation for its support of our research via the grants DMS-2104865 (Feehan) and DMS-2104871 (Leness). This work is also based in part on research supported by the National Science Foundation under Grant No. 1440140, while one of the authors (Feehan) was in residence at the Simons Laufer Mathematical Sciences Institute in Berkeley, California, during Fall 2022 as a Research Professor in the program \emph{Analytic and Geometric Aspects of Gauge Theory}. We also thank the Institute for Advanced Study, Princeton, for providing a welcoming scientific environment for its neighboring mathematicians.

\bigskip
\bigskip

\leftline{Paul M. N. Feehan}
\leftline{Department of Mathematics}
\leftline{Rutgers, The State University of New Jersey}
\leftline{Piscataway, NJ 08854}
\leftline{United States}
\medskip

\leftline{\texttt{feehan@math.rutgers.edu}}
\leftline{\url{math.rutgers.edu/~feehan}}
\bigskip

\leftline{Thomas G. Leness}
\leftline{Department of Mathematics}
\leftline{Florida International University}
\leftline{Miami, FL 33199}
\leftline{United States}
\medskip

\leftline{\texttt{lenesst@fiu.edu}}
\leftline{\url{fiu.edu/~lenesst}}
\bigskip

%COMMENT Remove hours and minutes for arXiv and journal versions and fix date
\leftline{This version: December 18, 2025}
%\leftline{This version: \date{\today{ }\hhmm}}

\mainmatter

\chapter{Introduction}
\label{chap:Introduction}
% PF9-10-2024 Add in references or reference details
% PF7-14-2025 Change \Ind to \Index = \Ind
%PF7-22-2025 Define almost complex, almost Hermitian, almost K\"ahler, complex, Hermitian, K\"ahler
This work is a sequel to our previous monograph \cite{Feehan_Leness_introduction_virtual_morse_theory_so3_monopoles}, where we initiated our program to prove the \emph{symplectic Bogomolov--Miyaoka--Yau conjecture} --- see Conjecture \ref{mainconj:BMY_symplectic} for a precise statement and Conjecture \ref{mainconj:BMY_Seiberg-Witten} more general assertion that implies Conjecture \ref{mainconj:BMY_symplectic} as a corollary when combined with results due to Taubes \cite{TauSymp, TauSympMore, TauSWGromov}. Our program to prove Conjecture \ref{mainconj:BMY_Seiberg-Witten} relies on gauge theory, employing Morse theory for a natural Hamiltonian for the circle action on an analytic compactification of the moduli space of non-Abelian monopoles over a closed symplectic four-manifold (more generally, a closed four-manifold with a Seiberg--Witten basic class). We use Morse theory to prove the existence of a projectively anti-self-dual connection on a rank-two Hermitian vector bundle over a blow-up of the four-manifold. The first Pontrjagin number of the vector bundle is assumed to be negative and greater than or equal to minus the Euler characteristic of the blown-up four-manifold. Our Morse theory arguments relies on a new definition of a Morse--Bott type index, which we call a \emph{virtual Morse--Bott index} \cite{Feehan_analytic_spaces}, that is appropriate for critical points of Hamiltonians for circle actions on complex analytic spaces (or real analytic spaces that are locally well-approximated by complex analytic spaces), as developed by the first author in \cite{Feehan_analytic_spaces}.

The moduli space of non-Abelian monopoles has a circle action induced by scalar multiplication by $S^1$ on the coupled spinor sections. Fixed points of this circle action are represented by Seiberg--Witten monopoles if the coupled spinor is not identically zero and otherwise are represented by projectively anti-self-dual connections \cite{FL2a}. A natural candidate for a Hamiltonian for the circle action on the moduli space of non-Abelian monopoles is the square of the $L^2$ norm of the coupled spinor, so points represented by anti-self-dual connections are absolute minima of this Hamiltonian. When the four-manifold is \emph{complex K\"ahler}, then the moduli space of non-Abelian monopoles inherits the structure of a complex K\"ahler manifold near regular points and admits circle-equivariant virtual neighborhoods near singular points that are complex K\"ahler manifolds. In particular, the fundamental two-form given by the K\"ahler metric on the moduli space of non-Abelian monopoles is a circle-invariant, non-degenerate two-form. In this situation, Frankel's Theorem (see \cite[Section 3]{Frankel_1959} for symplectic manifolds and Theorem \ref{thm:Frankel_almost_Hermitian} for almost Hermitian manifolds) implies that points in the moduli space of non-Abelian monopoles with non-vanishing coupled spinor are critical points of the Hamiltonian if and only if they are fixed points of the circle action, in which case they are represented by Seiberg--Witten monopoles. One of the main results in our previous monograph, namely 
%TL11-26-2025: Update to \cite[Section 2.3, Corollary 11]{Feehan_Leness_introduction_virtual_morse_theory_so3_monopoles}?
\cite[Section 1.44, Corollary 7]{Feehan_Leness_introduction_virtual_morse_theory_so3_monopoles}, established that the virtual Morse--Bott index of each such critical point is positive and so by \cite{Feehan_analytic_spaces}, it cannot be a local minimum of the Hamiltonian. One thus expects gradient flow for this Hamiltonian to ultimately yield the desired projectively anti-self-connection and prove Conjecture \ref{mainconj:BMY_Seiberg-Witten}.

While our previous results \cite{Feehan_Leness_introduction_virtual_morse_theory_so3_monopoles} provide compelling evidence in support of our program to prove Conjecture \ref{mainconj:BMY_Seiberg-Witten}, they must be extended in two directions in order to prove that conjecture in full. First, one must allow for \emph{energy bubbling} and include circle-equivariant virtual neighborhoods of \emph{ideal points} in the analytic compactification of the moduli space of non-Abelian monopoles. Second, one must allow the four-manifold to be \emph{symplectic} or more generally, a \emph{four-manifold with a Seiberg--Witten basic class}. The construction of circle-equivariant virtual neighborhoods requires \emph{gluing}, albeit a local construction (similar to \cite{FL3}) that we expect to be simpler than the global construction assumed (as a technical hypothesis) in our monograph \cite{FL5} and giving a cobordism formula for Donaldson invariants in terms of Seiberg--Witten invariants, and ultimately a proof in \cite{FL6, FL7, FL8} of Witten's formula \cite{Witten} for the Donaldson series in terms of Seiberg--Witten invariants and basic classes.

Our present work extends our main results in \cite{Feehan_Leness_introduction_virtual_morse_theory_so3_monopoles} by replacing the constraint that the four-manifolds be complex K\"ahler by the hypothesis in Conjecture \ref{mainconj:BMY_symplectic} that $X$ be symplectic.

In parallel work, joint with Wentworth \cite{Feehan_Leness_Wentworth_virtual_morse_theory_stable_pairs_bmy_kaehler}, we use the identification of non-Abelian monopoles over complex K\"ahler surfaces with oriented, polystable pairs of holomorphic vector bundles and sections to give a gauge-theoretic proof of the Bogomolov--Miyaoka--Yau inequality for complex surfaces of general type. While our paradigm in this special case is similar to that described in 
%TL11-26-2025: Does this need updating?
\cite[Section 1.3]{Feehan_Leness_introduction_virtual_morse_theory_so3_monopoles}, we may use a Gieseker compactification of the moduli space of oriented, polystable pairs of holomorphic vector bundles and sections instead of an analytic compactification constructed by gluing.

In Section \ref{sec:Review_BMY_inequality_compact_complex_surfaces}, we review the Bogomolov--Miyaoka--Yau inequality for complex surfaces of general type. Section \ref{sec:Symplectic_BMY_conjecture} introduces Conjecture \ref{mainconj:BMY_Seiberg-Witten} (the Bogomolov--Miyaoka--Yau inequality for four-manifolds with a Seiberg--Witten basic class), and its corollary, Conjecture \ref{mainconj:BMY_symplectic} (the Bogomolov--Miyaoka--Yau inequality for symplectic four-manifolds), along with a brief discussion of previous efforts to prove or disprove Conjecture \ref{mainconj:BMY_symplectic}. In Section \ref{sec:Existence_ASD_connections_and_BMY_inequality}, we describe Conjecture \ref{mainconj:Existence_ASD_connection} (on the existence of projectively anti-self-dual connections with small instanton number), which we expect to ultimately prove using methods of geometric analysis, and explain how it implies Conjecture \ref{mainconj:BMY_Seiberg-Witten}. In Section \ref{sec:Frankel_theorem_circle_actions_almost_Hermitian_manifolds}, we recall our version of Frankel's Theorem for the Hamiltonian function of a circle action on a smooth almost symplectic manifold that we proved in \cite{Feehan_Leness_introduction_virtual_morse_theory_so3_monopoles}. Section \ref{sec:Virtual_Morse-Bott_signature_Hamiltonian_function_circle_action_complex_analytic_space} contains an introduction to the definition of the virtual Morse--Bott signature of a critical point for the Hamiltonian function of a circle action on a complex analytic space. In Section \ref{sec:Perturbed_non-Abelian_monopole_equations_almost_Hermitian_4-manifolds}, we define the system of perturbed non-Abelian monopole equations, over almost Hermitian four-manifolds, that underlies our approach to prove Conjecture \ref{mainconj:BMY_Seiberg-Witten}. Section \ref{sec:Main_results} contains the statements of our main results and Section \ref{sec:Outline} gives an outline of our present work. In Section \ref{sec:Further_work}, we discuss the directions in which we expect to extend our present work and ultimately prove Conjecture \ref{mainconj:BMY_Seiberg-Witten}.

\section{Review of the Bogomolov--Miyaoka--Yau inequality for complex surfaces of general type}
\label{sec:Review_BMY_inequality_compact_complex_surfaces}
We begin by recalling the
  
\begin{thm}[Bogomolov--Miyaoka--Yau inequality for complex surfaces of general type]
\label{thm:BMY}
(See Miyaoka \cite[Theorem 4]{Miyaoka_1977} and Yau \cite[Theorem 4]{YauPNAS}.)
If $X$ is a compact, complex surface of general type, then
\begin{equation}
  \label{eq:BMY}
  c_1(X)^2 \leq 3c_2(X).
\end{equation}
\end{thm}

Here, $c_1(X)$ and $c_2(X)$ are the Chern class and Chern number of the holomorphic tangent bundle, $\sT_X \cong T^{1,0}X$. In \cite{Miyaoka_1977}, Miyaoka proved Theorem \ref{thm:BMY} using methods of algebraic geometry. See Barth, Hulek, Peters, and Van de Ven \cite[Section VII.4]{Barth_Hulek_Peters_Van_de_Ven_compact_complex_surfaces} for a simplification of Miyaoka's proof of Theorem \ref{thm:BMY}. Bogomolov \cite{Bogomolov_1978} proved a weaker version of \eqref{eq:BMY}, namely $c_1(X)^2\leq 4c_2(X)$.

According to the Enriques--Kodaira classification of minimal, compact, complex surfaces (see Barth, Hulek, Peters, and Van de Ven \cite[Chapter VI, Section 1, Theorem 1.1 and Table 10, pp. 243--244]{Barth_Hulek_Peters_Van_de_Ven_compact_complex_surfaces}, every such complex surface obeys the inequality \eqref{eq:BMY} when $b^+(X) > 1$ or $c_2(X) \geq 0$. (See Remarks \ref{rmk:Hypotheses_conjecture_BMY_Seiberg-Witten_c_2(X)=0} and \ref{rmk:Hypotheses_conjecture_BMY_Seiberg-Witten_b^+(X)=1} for a discussion of counterexamples when $b^+(X)=1$ and $c_2(X)<0$.)

Yau proved \eqref{eq:BMY} as a consequence of his proof of the Calabi Conjectures. Let $X$ be a complex manifold of dimension $n\geq 2$ with K\"ahler metric $g$. Let $\omega_g$ be the closed, real $(1,1)$-form and $\Ric(g)$ be the Ricci curvature $2$-form associated to $g$. One calls $g$ a \emph{K\"ahler--Einstein metric} if there exists a $\lambda\in\RR$ such that
\begin{equation}
  \label{eq:KE_metric}
  \Ric(g) = \lambda\omega_g.
\end{equation}
As part of his work \cite{Yau} on the Calabi Conjectures, Yau proved the
%PF9-10-2024 Reference the actual PDEs he solves; maybe mention Weinkove's work. 

\begin{thm}[Yau]
\label{thm:Calabi_conjecture}  
(See Yau \cite[Theorems 1 and 2, pp. 363--364]{Yau} and \cite[Theorem 1]{YauPNAS}, or Barth, Hulek, Peters, and Van de Ven \cite[Theorem 15.2, p. 52]{Barth_Hulek_Peters_Van_de_Ven_compact_complex_surfaces}.)
Let $X$ be a compact, complex manifold of dimension $n\geq 2$ such that $-c_1(X)$ can be represented by a K\"ahler form. Then $X$ admits a K\"ahler--Einstein metric.
\end{thm}

Tosatti, and Weinkove \cite[Section 1, Main Theorem, p. 1187]{Tosatti_Weinkove_2010} extended the main result of Yau in \cite{Yau} to compact, complex Hermitian manifolds, while Chu, Tosatti, and Weinkove \cite[Theorem 1.1, p. 1950]{Chu_Tosatti_Weinkove_2019} extended it further to compact almost Hermitian manifolds.
%PF10-17-2024 We should add explanation for why their results does not lead to BMY.
The K\"ahler--Einstein metric $g$ produced by Theorem \ref{thm:Calabi_conjecture} obeys \eqref{eq:KE_metric} with constant $\lambda<0$.

\begin{lem}[Chern--Weil inequality]
\label{lem:Chern-Weil_inequality}
(See Tosatti \cite[Lemma 2.6]{Tosatti_2017} for an exposition.)  
If $(X, \omega)$ is a compact K\"ahler--Einstein manifold of dimension $n \geq 2$, then
\begin{equation}
  \label{eq:Chern-Weil_inequality}
  \left(\frac{2(n+1)}{n}c_2(X) - c_1(X)^2\right)\cdot[\omega]^{n-2} \geq 0,
\end{equation}
with equality if and only if $\omega$ has constant holomorphic sectional curvature.
\end{lem}

By applying Theorem \ref{thm:Calabi_conjecture} and Lemma \ref{lem:Chern-Weil_inequality}, Yau proved

\begin{thm}[Bogomolov--Miyaoka--Yau inequality and uniformization]
\label{thm:Bogomolov-Miyaoka-Yau_inequality_and_uniformization}  
(See Yau \cite[Theorem 4]{YauPNAS}.)
Let $X$ be a compact, complex K\"ahler surface with ample canonical bundle. Then inequality \eqref{eq:BMY} holds and equality occurs if and only if $X$ is covered biholomorphically by the ball in $\CC^2$.
\end{thm}
%PF11-24-2025 Check our definitions of canonical bundle everywhere
%COMMENT https://en.wikipedia.org/wiki/Canonical_bundle
%COMMENT https://en.wikipedia.org/wiki/Chern_class
%COMMENT https://en.wikipedia.org/wiki/Holomorphic_tangent_bundle

Recall that a holomorphic line bundle $L$ over a complex manifold $X$ is defined to be \emph{ample} if there exists a holomorphic embedding $f : X \to \CC\PP^N$ such that $L^{\otimes k} \cong f^*\sO_{\CC\PP^N}(1)$ for some positive integers $k$ and $N$.

A well-known theorem due to Kodaira (see Angella and Spotti \cite[Theorem 2.5, p. 207]{Angella_Spotti_2017}) asserts that a holomorphic line bundle over a compact complex manifold admits a positive Hermitian metric if and only if it is ample. If $X$ is \emph{Fano}, so $c_1(X)>0$, that is, the first Chern class can be represented by a K\"ahler form, Kodaira's Theorem says that this condition is equivalent to the anti-canonical line bundle being ample. If $c_1(X)<0$, Kodaira's Theorem says that this condition is equivalent to the canonical line bundle being ample and, in particular, $X$ is a \emph{surface of general type}.

Simpson \cite[p. 871]{Simpson_1988} proved \eqref{eq:BMY} as a corollary of his main theorem \cite[p. 870]{Simpson_1988} on existence of a Hermitian--Einstein connection on a stable Higgs bundle of rank $3$ over $X$ and the following

\begin{thm}[Bogomolov--Gieseker inequality]
\label{thm:Bogomolov-Gieseker_inequality}
(See Kobayashi \cite[Theorem 4.4.7]{Kobayashi} or L\"ubke and Teleman \cite[Corollary 2.2.4]{Lubke_Teleman_1995}.)
Let $(E,H)$ be a Hermitian vector bundle over of rank $r$ over a compact, complex K\"ahler manifold of dimension $n \geq 2$. If $(E,H)$ admits a Hermitian--Einstein connection, then
\begin{equation}
  \label{eq:BG_inequality}
  \int_X\left(2rc_2(E) - (r-1)c_1(E)^2\right)\wedge \omega^{n-2} \geq 0.
\end{equation}
\end{thm}

According to Bogomolov \cite{Bogomolov_1978} and Gieseker \cite{Gieseker_1977}, a version of inequality \eqref{eq:BG_inequality} holds for any \emph{slope semi-stable, torsion-free sheaf} over a smooth complex projective surface (see Huybrechts and Lehn \cite[Theorem 3.4.1, p. 80]{Huybrechts_Lehn_geometry_moduli_spaces_sheaves}). 

\section{Symplectic Bogomolov--Miyaoka--Yau conjecture}
\label{sec:Symplectic_BMY_conjecture}
We now state a conjectural version of Theorem \ref{thm:BMY} for a broader class of four-manifolds than compact, complex surfaces of general type. For a closed topological four-manifold $X$, we define
\begin{equation}
\label{eq:Define_c1Squred_c_2}
  c_1(X)^2 := 2e(X)+3\sigma(X)
  \quad\text{and}\quad
  c_2(X) := e(X),
\end{equation}
where $e(X)=2-2b_1(X)+b_2(X)$ and $\sigma(X)=b^+(X)-b^-(X)$ are the \emph{Euler characteristic} and \emph{signature} of $X$, respectively. If $Q_X$ is the intersection form on $H_2(X;\ZZ)$, then $b^\pm(X)$ are the dimensions of the maximal positive and negative subspaces of $Q_X$ on $H_2(X;\RR)$.

\begin{defn}[Seiberg--Witten basic class]
\label{defn:Seiberg-Witten_basic_class}
(See Salamon \cite[Section 7.5, Definition 7.35, p. 261]{SalamonSWBook} or Morgan \cite[Section 6.7]{MorganSWNotes}.)
Let $X$ be a closed, connected, oriented, smooth four-manifold with $b^+(X) > 1$ and odd $b^+(X) - b_1(X)$. We call $K \in H^2(X;\ZZ)$ a \emph{Seiberg--Witten basic class} if
\begin{subequations}
  \label{eq:Seiberg-Witten_basic_class}
  \begin{gather}
    \label{eq:Seiberg-Witten_basic_class_K_squared=c1_squared}
    K^2 = 2e(X)+3\sigma(X),
    \\
    \label{eq:Seiberg-Witten_basic_class_c1spinc=K}
    K = c_1(\fs),
    \\
    \label{eq:Seiberg-Witten_basic_class_SWspinc_neq_0}
    \SW_X(\fs) \neq 0,
  \end{gather}
\end{subequations}
for some \spinc structure $\fs$ on $X$, where $\SW_X(\fs)$ is its Seiberg--Witten invariant.
\qed\end{defn}

%PF9-29-2024 Recheck GS references against latest edition. Numbering has changed.
A closed, orientable, smooth manifold $X$ admits a \emph{spin${}^c$ structure} (see Gompf and Stipsicz \cite[Definition 2.4.15, p. 55]{GompfStipsicz}, Lawson and Michelsohn \cite[Definition D.1, p. 391]{LM}, Morgan \cite[Section 3.1, pp. 24--25]{MorganSWNotes}, or Salamon \cite[Section 5.1, Definition 1.5.2, p. 153]{SalamonSWBook}) it and only if $w_2(X) \in H^2(X;\ZZ/2\ZZ)$ has an integral lift to an element of $H^2(X;\ZZ)$ (see Gompf and Stipsicz \cite[Proposition 2.4.16, p. 56]{GompfStipsicz} or Lawson and Michelsohn \cite[Corollary D.4, p. 391]{LM}). If $\fs = (\rho,W)$ is a \spinc structure on $X$, for some Riemannian metric $g$ on $X$, then $c_1(\fs) := c_1(W^+)$ by the forthcoming definition \eqref{eq:DefineChernClassOfSpinc}.

See Kronheimer and Mrowka \cite[Definitions 1.5.3 and 1.5.4, p. 13]{KMBook}, Morgan \cite[Theorem 6.7.3, p. 100]{MorganSWNotes}, Nicolaescu \cite[Theorem 2.3.5, p. 153]{NicolaescuSWNotes}, Salamon \cite[Equation (7.20), p. 247]{SalamonSWBook}, or Witten \cite[Equation (2.7), p. 774]{Witten} for the definition of \emph{Seiberg--Witten invariants} using generic geometric perturbations and Ruan \cite{RuanSW} for the definition of Seiberg--Witten invariants using virtual fundamental classes. If $b^+(X) - b_1(X)$ is even, then the Seiberg--Witten invariants are defined to be zero (see Morgan \cite[Section 6.7, p. 99]{MorganSWNotes}), so when we write that $\SW_X(\fs)\neq 0$ for some \spinc structure $\fs$ on a four-manifold $X$ with $b^+(X) > 1$, we implicitly assume as well that $b^+(X) - b_1(X)$ is odd.

An \emph{almost complex structure} $J$ on a smooth manifold $X$ is a smooth section of $\End(TX)$ with $J^2 = -\id_{TX}$. Let $X$ be a closed, oriented, smooth four-manifold. By Wu's Theorem (see Gompf and Stipsicz \cite[Theorem 1.4.15, p. 30]{GompfStipsicz}, Salamon \cite[Section 13.1, Proposition 13.1, p. 408]{SalamonSWBook}, or Wu \cite{Wu_1952}), there is a one-to-one correspondence between isomorphism classes of almost complex structures on $X$ which are compatible with the orientation and integral cohomology classes $K \in H^2(X;\ZZ)$ which satisfy Equation \eqref{eq:Seiberg-Witten_basic_class_K_squared=c1_squared}, namely $K^2 = c_1(X)^2$, and 
\begin{equation}
\label{eq:K_integral_lift_w2}
% PF9-29-2024 Salamon says <K,\alpha> = <\alpha,\alpha> for all alpha in H^2(X;\ZZ)? Is that equivalent?
%PF10-14-2024 Please address
%10-14-2024: <K,\alpha> \equiv <\alpha,\alpha> \pmod 2 is equivalent to the below by Poincare duality over \ZZ/2 and GS Prop. 1.4.18.  Not true if \pmod 2 excluded (just multiply \alpha by r...)
  K \equiv w_2(X) \pmod{2\ZZ}.
\end{equation}
One says that a pair $J_1, J_2$ of almost complex structures on a smooth manifold $X$ are \emph{isomorphic} if $(TX,J_1)$ and $(TX,J_2)$ are isomorphic as complex vector bundles.

%PF10-16-2024 I edited. I think this remark should be turned into a proof and a statement of Wu's theorem would be better moved to a short first subsection in section 4. References above to Wu's theorem can cite that.
%TL10-16-2024: I'm hesitant to highlight it much --it's not far removed from an exercise in GompfStipsicz.
%PF10-17-2024 We can leave it alone for now, but a long side remark like this is better tucked into the body of the paper. A long remark in the intro is a highlight!
\begin{rmk}[On Wu's Theorem]
\label{rmk:Wu_theorem}
Because the article \cite{Wu_1952}
by Wu is difficult to find, we give an exposition of the proof of Wu's Theorem here. We suppose as in \cite{Wu_1952} that $X$ is closed, oriented, smooth four-manifold. An almost complex structure $J$ defines a smooth rank two complex vector bundle, $E_J:= (TX,J)$. Conversely, if $E_\RR$ denotes the smooth rank four real vector bundle underlying a smooth rank two complex vector bundle $E$, then an isomorphism $E_\RR\cong TX$ of smooth real vector bundles defines an almost complex structure on $X$. Hence, there is a bijection between almost complex structures on $X$, up to the notion of isomorphism defined above, and isomorphism classes of smooth complex vector bundles $E$ with $E_\RR\cong TX$.

Wu's criterion classifies smooth rank two complex vector bundles $E$ with $E_\RR\cong TX$ up to an isomorphism of smooth complex vector bundles as follows. We claim that a smooth complex vector bundle $E$ satisfies $E_\RR\cong TX$ if and only if $c_2(E)=e(TX)$ and $K = c_1(E)$ satisfy the criteria \eqref{eq:Seiberg-Witten_basic_class_K_squared=c1_squared} and \eqref{eq:K_integral_lift_w2} of Wu's Theorem, namely
\[
  K\equiv w_2(TX)\pmod{2\ZZ} \quad\text{and}\quad K^2=2e(X)+3\sigma(X),
\]
where we write $e(X) = \langle e(TX),[X]\rangle$ and $w_2(TX) = w_2(X)$. To see this, first note that by Gompf and Stipsicz \cite[Theorem 1.4.20b, p. 31]{GompfStipsicz}, a rank two complex vector bundle $E$ has $E_\RR\cong TX$ if and only if the characteristic classes $w_2 \in H^2(X;\ZZ/2\ZZ)$, and $p_1, e \in H^4(X;\ZZ)$ are the same for both vector bundles. The equality $w_2(E_\RR)\equiv c_1(E)\pmod{2\ZZ}$ appears in \cite[Theorem 1.4.9, p. 28]{GompfStipsicz} and $e(E_\RR)=c_2(E)$ is the definition of the top Chern class \cite[Section 14, p. 158]{MilnorStasheff}. It remains to check that these two equalities and $c_1(E)^2 = 3\sigma(X)+2e(X)$ imply that $p_1(E_\RR)=p_1(TX)$. We compute that
\[
  p_1(E_\RR) = -c_2(E\otimes_\RR\CC) = -c_2(E\oplus \bar E)=2c_2(E)-c_1(E)^2,
\]
where the second equality appears in \cite[Lemma 15.4, p. 176]{MilnorStasheff}. Because $\langle p_1(TX),[X]\rangle=3\sigma(X)$ by the Hirzebruch signature theorem, \cite[Theorem 1.4.12, p. 28]{GompfStipsicz}, we see that $p_1(E_\RR)=p_1(TX)$ if and only if $3\sigma(X)=2c_2(E)-c_1(E)^2$. Thus, $E_\RR\cong TX$ if and only if $c_2(E)=e(TX)$ and $c_1(E)$ satisfies the criteria \eqref{eq:Seiberg-Witten_basic_class_K_squared=c1_squared} and \eqref{eq:K_integral_lift_w2} of Wu's theorem.

Finally, smooth rank two complex vector bundles are isomorphic if and only if their Chern classes are equal. Hence, there is a bijection between classes $K\in H^2(X;\ZZ)$ satisfying \eqref{eq:Seiberg-Witten_basic_class_K_squared=c1_squared} and \eqref{eq:K_integral_lift_w2} and isomorphism classes of rank two complex vector bundles $E$ with $E_\RR\cong TX$.
\qed\end{rmk}

If $K \in H^2(X;\ZZ)$ obeys Equations \eqref{eq:Seiberg-Witten_basic_class_K_squared=c1_squared} and
\eqref{eq:Seiberg-Witten_basic_class_c1spinc=K}, for some \spinc structure $\fs$ on $X$, then $K$ also obeys \eqref{eq:K_integral_lift_w2} by Gompf and Stipsicz \cite[Proposition 2.4.16, p. 56]{GompfStipsicz}. If $(X,g)$ is a closed, smooth Riemannian $2n$-manifold that admits a $g$-orthogonal almost complex structure $J$, then 
as described in Salamon \cite[Lemma 4.52, p. 141]{SalamonSWBook}, the four-manifold
$X$ admits a natural \spinc structure $\fs_\can = (\rho_\can,W_\can)$, called the \emph{canonical spin${}^c$ structure}, such that $\det W_\can^+ \cong K_X^*$, where
% TL5-20-2025: Use \Lambda or \wedge here?
%PF8-31-2025 \wedge
$K_X = \Lambda^{0,n}(X) = \wedge^{0,n}(T^*X)$ is the canonical line bundle (see Morgan \cite[Corollary 3.4.5, p. 49]{MorganSWNotes} or Salamon \cite[Section 5.3, Corollary 5.21, p. 168]{SalamonSWBook}).

\begin{mainconj}[Bogomolov--Miyaoka--Yau inequality for four-manifolds with a Seiberg--Witten basic class]
\label{mainconj:BMY_Seiberg-Witten}
Let $X$ be a closed, connected, oriented, smooth, four-dimensional manifold with
$b^+(X) > 1$
%PF3-19-2025 Below is redundant
%TL8-24-2025: The redundancy is somewhat buried in the definition of SW basic class (I think) so leaving it here is ok for emphasis on the limitations of this conjecture.
and odd $b^+(X)-b_1(X)$.
If $X$ has a Seiberg--Witten basic class, then \eqref{eq:BMY} holds:
\[
  c_1(X)^2 \leq 3c_2(X).
\]
\end{mainconj}

\begin{rmk}[On the hypotheses of Conjecture \ref{mainconj:BMY_Seiberg-Witten} when $c_2(X)=0$]
\label{rmk:Hypotheses_conjecture_BMY_Seiberg-Witten_c_2(X)=0} 
According to the Enriques--Kodaira classification of minimal, compact, complex surfaces (see Barth, Hulek, Peters, and Van de Ven \cite[Chapter VI, Section 1, Theorem 1.1 and Table 10, pp. 243--244]{Barth_Hulek_Peters_Van_de_Ven_compact_complex_surfaces}, there are many examples of complex surfaces which obey the inequality \eqref{eq:BMY} with $c_2(X)= 0$ and $c_1(X)^2 \leq 0$:
\begin{inparaenum}[\itshape i\upshape)]
\item minimal rational surfaces,
\item minimal surfaces of class VII,
\item ruled surfaces of genus $g=0,1$,
\item bi-elliptic surfaces,
\item Kodaira surfaces,
\item tori, and
\item minimal properly elliptic surfaces.
\end{inparaenum}  
% PF9-29-2024 Add ref to explanation
Our gauge theory program to prove Conjecture \ref{mainconj:BMY_Seiberg-Witten} requires that $c_2(X) > 0$, but that does not exclude the possibility that Conjecture \ref{mainconj:BMY_Seiberg-Witten} holds when $c_2(X) = 0$. Conjecture \ref{mainconj:BMY_Seiberg-Witten} does not hold when $b^+(X)=1$ and $c_2(X)<0$, as we explain in the forthcoming Remark \ref{rmk:Hypotheses_conjecture_BMY_Seiberg-Witten_b^+(X)=1}.
\qed\end{rmk}  

\begin{rmk}[On the hypotheses of Conjecture \ref{mainconj:BMY_Seiberg-Witten} when $b^+(X)=1$]
\label{rmk:Hypotheses_conjecture_BMY_Seiberg-Witten_b^+(X)=1}
The approach that we describe in our present work and Feehan and Leness 
%TL11-26-2025: I think this reference is still okay but I still wanted to highlight it
\cite[Chapter 1]{Feehan_Leness_introduction_virtual_morse_theory_so3_monopoles} to prove Conjecture \ref{mainconj:BMY_Seiberg-Witten} fails when $Q_X$ is odd and $b^+(X) = 1 = b^-(X)$. This finding is consistent with results from the Kodaira--Enriques classification of compact, connected, minimal complex surfaces in Barth, Hulek, Peters, and Van de Ven \cite[Section VI.1, Table 10, p. 244]{Barth_Hulek_Peters_Van_de_Ven_compact_complex_surfaces}.

Ruled complex surfaces $X$ fibering over a Riemann surface $\Sigma_g$ of genus $g$ violate inequality \eqref{eq:BMY} when $g \geq 2$ since $c_1(X)^2 = 8(1-g)$ and $c_2(X) = 4(1-g)$, so $c_1(X)^2 \leq 3c_2(X)$ if and only if $g=0,1$. Ruled complex surfaces have $b^+(X)=1=b^-(X)$ and $b_1(X)=2g$ and non-zero Seiberg--Witten invariants: see Li and Liu \cite{LiLiu}, Ohta and Ono \cite{Ohta_Ono_1996a, Ohta_Ono_1996b}, and Salamon \cite[Proposition 9.17, p. 309]{SalamonSWBook} for ruled complex surfaces, Salamon \cite[Proposition 9.17, p. 309 and Theorem 12.10, p. 378]{SalamonSWBook} for K\"ahler surfaces with $b^+(X)=1$, and Taubes \cite{TauSymp, TauSympMore}, \cite[Proposition 13.9, p. 412]{SalamonSWBook} for symplectic four-manifolds with $b^+(X)=1$. See Taubes \cite[Theorem 1]{TauSymp}, \cite[Proof of Theorem 3]{TauSympMore} for $\CC\PP^2$, and Kotschick \cite[Section 3.2]{KotschickSW} for comments on Taubes' theorems in the case $b^+(X)=1$. We refer to the forthcoming Remark \ref{rmk:bPlusOne_and_Feasibility} for a further discussion of the case $b^+(X)=1$.
\qed\end{rmk}

% PF9-29-2024 Get precise references for assertion below, including Kotschick, Salamon
Taubes \cite{TauSymp, TauSympMore} proved that symplectic four-manifolds with $b^+(X) > 1$ have Seiberg--Witten simple type with non-zero Seiberg--Witten invariants.
%PF10-3-2024 Did he also do this for b^+ =  1?
According to Gompf and Stipsicz \cite[Theorem 1.4.13, p. 29 and comments, pp. 29--30]{GompfStipsicz} and Salamon \cite[Section 13.1, Remark 13.3, p. 409]{SalamonSWBook}, if $X$ is a closed almost complex four-manifold, then $b^+(X)-b_1(X)$ is odd. Thus, Conjecture \ref{mainconj:BMY_Seiberg-Witten} implies the following well-known claim as a special case:

\begin{mainconj}[Bogomolov--Miyaoka--Yau inequality for symplectic four-manifolds]
\label{mainconj:BMY_symplectic}
(Compare Gompf and Stipsicz \cite[Remark 10.2.16 (c), p. 401]{GompfStipsicz} or Stipsicz \cite[Section 5, Conjecture 5.1, p. 276]{Stipsicz_2000}.)
Let $X$ be a closed, connected, oriented, four-dimensional, smooth manifold with $b^+(X) > 1$. If $X$ is symplectic, then \eqref{eq:BMY} holds:
\[
  c_1(X)^2 \leq 3c_2(X).
\]
\end{mainconj}

Our Conjecture \ref{mainconj:BMY_Seiberg-Witten} is motivated in part by the fact that Szab\'o \cite{SzaboNoSymplectic} proved existence of four-dimensional, \emph{non-symplectic}, smooth manifolds with non-zero Seiberg--Witten invariants. See also the remarks by Stipsicz in \cite[Section 5, second paragraph below Conjecture 5.1, p. 276]{Stipsicz_2000}.

Conjecture \ref{mainconj:BMY_Seiberg-Witten} has inspired constructions by topologists of examples to shed light on inequality \eqref{eq:BMY}, including work of Akhmedov, Hughes, and Park \cite{Akhmedov_Hughes_Park_2013, Akhmedov_Park_2008jggt, Akhmedov_Park_2010mrl}, Baldridge, Kirk, and Li \cite{Baldridge_Kirk_2006, Baldridge_Kirk_2007, Baldridge_Li_2005}, Bryan, Donagi, and Stipsicz \cite{Bryan_Donagi_Stipsicz_2001}, Fintushel and Stern \cite{Fintushel_Stern_1994im}, Gompf and Mrowka \cite{GompfNewSymplectic, GompfMrowka}, Hamenst\"adt \cite{Hamenstaedt_2012arxiv}, Park and Stipsicz \cite{Park_Stipsicz_2015, Stipsicz_1998, Stipsicz_2000}, Smith \cite{Smith_2023-7-26_private}, Torres \cite{Torres_2014}, and others.

Conjecture \ref{mainconj:BMY_Seiberg-Witten} holds for all examples that satisfy the hypotheses, but the Bogomolov--Miyaoka--Yau inequality \eqref{eq:BMY} can fail for four-manifolds with zero Seiberg--Witten invariants, such as a connected sum of two or more copies of $\CC\PP^2$.
%COMMENT Does it fail for connected sums obeying the topological hypotheses
%TL9-27-2024: Take $X_{m,n}=m\CC\PP^2\# n\overline{\CC\PP}^2$, i.e. the connected sum of $m$ copies of $\CC\PP^2$ and $n$ copies of $\overline{\CC\PP^2}$.  Then $e(X_{m,n})=2+m+n$ and $\sigma(X_{m,n})=m-n$.  Thus
%$c_1(X_{m,n})^2-c_2(X_{m,n})=2m-4n-2$ which we can make positive by taking e.g. $m=5$ and $n=2$.

The congruence $e(X)+\sigma(X)\equiv 0 \pmod{4}$ is a \emph{necessary} condition for closed, smooth $4$-manifold $X$ to be almost complex by Hirzebruch \cite[Kommentare, p. 777]{Hirzebruch_gesammelte_abhandlungen_I_and_II}; thus, if $X_1$ and $X_2$ are closed, smooth, almost complex $4$-manifolds, then $X_1\# X_2$ cannot be almost complex (see Albanese and Milivojevic \cite[Section 1, p. 1]{Albanese_Milivojevic_2019} or Goertsches and Konstantis \cite[Section 1, p. 131]{Goertsches_Konstantis_2019}). Note that $e(X)+\sigma(X) = 2(1+b^+(X)-b_1(X))$ when $X$ is a closed, smooth, almost complex $4$-manifold by Salamon \cite[Remark 13.3, p. 409]{SalamonSWBook}
%PF10-19-2024 Added
and so the necessary condition $e(X)+\sigma(X)\equiv 0 \pmod{4}$ for $X$ to be almost complex is equivalent to $1+b^+(X)-b_1(X) \equiv 2 \pmod{2}$, in other words, $b^+(X)-b_1(X)$ must be odd, as assumed in Conjecture \ref{mainconj:BMY_Seiberg-Witten}.

%PF10-19-2024 Added
On the other hand, the existence of an almost complex structure on $X$ is \emph{not} sufficient for \eqref{eq:BMY} to hold. Indeed, if $X = \CC\PP^2\# \CC\PP^2\#\CC\PP^2$ (smooth connected sum), then it admits an almost complex structure by Albanese and Milivojevic \cite[Remark 2.1, p. 3]{Albanese_Milivojevic_2019} (based on Goertsches and Konstantis
%PF10-19-2024 Get exact reference
\cite{Goertsches_Konstantis_2019}) but $c_1^2(X) = 19$ and $c_2(X) = 15$ (see Milivojevi\'c \cite[Section 3]{Milivojevic_AC_but_no_CX_structure_4-manifold}), and so it does not obey \eqref{eq:BMY}. Because this example is a non-trivial connected sum, its Seiberg--Witten invariants are zero (see, for example, Donaldson \cite[Section 5, pp. 61--62]{DonSW}, Kotschick \cite[Theorem 5.4, p. 216]{KotschickSW}, Nicolaescu \cite[Theorem 4.6.1, p. 461]{NicolaescuSWNotes}, or Salamon \cite[Theorem 11.1, p. 353]{SalamonSWBook}) and so it cannot have a Seiberg--Witten basic class in the sense of in Definition \ref{defn:Seiberg-Witten_basic_class}.

LeBrun \cite{LeBrunEinsteinMostow} proved that the Bogomolov--Miyaoka--Yau inequality \eqref{eq:BMY} holds for Einstein four-manifolds with non-zero Seiberg--Witten invariants.

\section{Anti-self-dual connections and applications to a proof of the Bogomolov--Miyaoka--Yau inequality}
\label{sec:Existence_ASD_connections_and_BMY_inequality}
Our approach to proving Conjecture \ref{mainconj:BMY_Seiberg-Witten} is to prove the following assertion and then deduce \eqref{eq:BMY} as an easy consequence:

\begin{mainconj}[Existence of anti-self-dual connections with small instanton number]
\label{mainconj:Existence_ASD_connection}
Assume the hypotheses of Conjecture \ref{mainconj:BMY_Seiberg-Witten} and let $(E,H)$ be a complex rank two, Hermitian vector bundle over $X$ whose associated bundle $\su(E)$ of skew-Hermitian, traceless endomorphisms of $E$ has first Pontrjagin number obeying the \emph{fundamental bounds},
\begin{equation}
  \label{eq:p1_lower_bound}
  0 > p_1(\su(E))[X] \geq -c_2(X).
\end{equation}
Let $g$ be a Riemannian metric on $X$ that is generic in the sense of Freed and Uhlenbeck \cite{DK, FU}. Then there exists a smooth, projectively $g$-anti-self-dual unitary connection $A$ on $E$, so the curvature $F_A \in \Omega^2(\fu(E))$ obeys
\begin{equation}
\label{eq:ASD}
  (F_A^+)_0 = 0 \in \Omega^+(X;\su(E)),
\end{equation}
where ${}^+:\wedge^2(T^*X)\to\wedge^+(T^*X)$ and $(\,\cdot\,)_0:\fu(E)\to\su(E)$ are orthogonal projections. 
\end{mainconj}

Recall that $\su(E)$ is a real vector bundle over $X$ and that its first Pontrjagin class $p_1(\su(E))$ is an element of $H^4(X;\ZZ)$, and its first Pontrjagin number is $p_1(\su(E))[X] = \langle p_1(\su(E)),[X]\rangle \in \ZZ$, where $[X] \in H_0(X;\ZZ)$ denotes the fundamental class of $X$. Taubes \cite{TauIndef}
%PF9-10-2024 Get exact ref
proved existence of solutions to the anti-self-dual equation \eqref{eq:ASD} only when the \emph{instanton number},
\begin{equation}
  \label{eq:DefineKappa}
  \kappa(E) := -\frac{1}{4}p_1(\su(E))[X],
\end{equation}
is sufficiently large. The hypotheses of Conjecture \ref{mainconj:Existence_ASD_connection} imply that $\kappa(E) > 0$, but the difficulty in proving Conjecture \ref{mainconj:Existence_ASD_connection} is apparent from the fact that the bound \eqref{eq:p1_lower_bound} implies
\[
  \kappa(E) \leq \frac{c_2(X)}{4},
\]
and so $\kappa(E)$ could be small, in which case Taubes' gluing method does not apply.

If $X$ is a minimal symplectic four-manifold with $b^+(X)>1$, then $c_1(X)^2 \geq 0$ (see Gompf and Stipsicz \cite[Theorem 10.1.18, p. 392]{GompfStipsicz}, Taubes \cite[Proposition 4.2, p. 231]{TauSWisGRAnnounce}, or Salamon \cite[Corollary 13.20, p. 417]{SalamonSWBook}). A compact symplectic four-manifold $(X,\omega)$ is called \emph{minimal} if it does not contain any symplectically embedded sphere with self-intersection number $-1$. If such a sphere does exist then $X$ decomposes as a connected sum of some symplectic four-manifold $(X', \omega')$ with $\overline{\CC\PP^2}$ (see McDuff and Salamon \cite{McDuffSalamonSympTop3}).
%PF9-10-2024 Get exact ref
By induction, every compact symplectic four-manifold is a connected sum of
a minimal one with finitely many copies of $\CC\PP^2$. In particular, if \eqref{eq:BMY} holds for a minimal symplectic four-manifold with $b^+(X)>1$, then $3c_2(X) \geq c_1(X)^2 \geq 0$. If $X$ is a minimal, compact, connected, complex surface of general type, then $3c_2(X) \geq c_1(X)^2 > 0$.

We aim to prove Conjecture \ref{mainconj:Existence_ASD_connection} via existence of  projectively anti-self-dual connections as absolute minima of a certain Hamiltonian function $f$ for the circle action on the \emph{singular} moduli space of non-Abelian monopoles. We now recall our explanation from Feehan and Leness 
%TL11-26-2025: Update?
\cite[Section 1.3]{Feehan_Leness_introduction_virtual_morse_theory_so3_monopoles} of why Conjecture \ref{mainconj:Existence_ASD_connection} implies Conjecture \ref{mainconj:BMY_Seiberg-Witten} when $c_2(X) > 0$.

For $w\in H^2(X;\ZZ)$ and $4\kappa\in \ZZ$, let $(E,H)$ be a rank-two Hermitian vector bundle over $X$ with $c_1(E)=w$, fixed unitary connection $A_d$ on $\det E$, and
\[
  p_1(\su(E))[X] = -4\kappa,
\]
where by Donaldson and Kronheimer \cite[Equation (2.1.39), p. 42]{DK} one has
\begin{equation}
\label{eq:DK_2-1-39}
  p_1(\su(E))=c_1(E)^2-4c_2(E).
\end{equation}
The moduli space of projectively $g$-anti-self-dual connections on $E$ is
\[
  M_\kappa^w(X,g) := \{A: (F_A^{+,g})_0=0\}/\SU(E),
\]
where $\SU(E) = \Aut(E,H,A_d)$ is the group of determinant-one, unitary automorphisms of $(E,H)$. The expected dimension of $M_\kappa^w(X,g)$ is given by \cite[Equation (4.2.22), p. 137]{DK}
\begin{equation}
  \label{eq:Expected_dimension_moduli_space_ASD_connections}
  \expdim M_\kappa^w(X,g) = -2p_1(\su(E))[X] - \frac{1}{2}(c_2(X) + c_1(X)^2).
\end{equation}
%PF3-24-2025 Updated. @Tom: please check carefully. Note KM's meaning for w (in \ZZ, not \ZZ/2\ZZ)
We now recall some important results that bear on the smoothness of the moduli space $M_\kappa^w(X,g)$. If $b^+(X) > 0$, then for any fixed integers $l > 0$ and $r>0$, there is an open dense subset $\sC_l(X)$ of the Banach space $\sC(X)$ of conformal equivalence classes $[g]$ of $C^r$ Riemannian metrics on $X$ such that the only reducible $g$-anti-self-dual connection on an $\SO(3)$ bundle $V$ over $X$ with $\kappa(V) = -\frac{1}{4}p_1(V) \leq l$ is a
% TL8-24-2025: DK say product connection but they may be working on simply connected manifolds.  I think "flat" is more accurate.
%PF8-31-2025 Your comment is vague. Please recheck this whole discussion around here and make sure it's accurate from start to finish.
flat connection (see Donaldson and Kronheimer \cite[Corollary 4.3.15, p. 148]{DK}).
%TL3-25-2025: Edited to give following
(Equivalence classes of flat connections may be excluded from the Uhlenbeck compactification of the moduli space of $g$-anti-self-dual connections on $V$ by requiring  that $w_2(V)$ pair non-trivially with a homology class represented by an embedded sphere
as noted by Morgan and Mrowka \cite[sentence prior to Corollary 2.2, p. 226]{MorganMrowkaPoly} in their ``blow-up trick''.)
%(Flat connections may be excluded from the Uhlenbeck compactification of the moduli space of $g$-anti-self-dual connections on $V$ by requiring that no integral lift of the class $w_2(V) \in H^2(X;\ZZ/2\ZZ)$ is torsion --- see the discussion following Definition \ref{defn:Good} based on Morgan and Mrowka \cite{MorganMrowkaPoly}.)
If $X$ is \emph{simply connected}, then by Donaldson and Kronheimer \cite[Section 4.3.3, p. 146 and Corollary 4.3.18, p. 149]{DK}, there is a dense (second category) subset $\sC'(X) \subset \sC(X)$ such that for any $[g] \in \sC'(X)$ and $\SO(3)$ bundle $\su(E)$ with $\kappa(E) \leq l$, the moduli subspace $M_\kappa^{w,*}(X,g)$ of irreducible $g$-anti-self-dual connections on $\su(E)$ is regular. More generally, even if $X$ is not necessarily simply connected, by Kronheimer and Mrowka \cite[Section 2 (i), Corollary 2.5, p. 587]{KMStructure} all points in $M_\kappa^w(X,g)$ are regular except for those represented by flat or reducible connections. Thus, Kronheimer and Mrowka obtain a generalization of Donaldson and Kronheimer \cite[Section 4.3.3, Corollary 4.3.19, p. 149]{DK} that removes their assumption that $X$ be simply connected: for a dense (second category) set $\sC_l(X) \cap \sC'(X)$ of conformal equivalence classes $[g]$ of $C^r$ metrics on $X$, the moduli spaces $M_\kappa^w(X,g)$ for $\SO(3)$ bundles $\su(E)$ with $0 < \kappa(E) \leq l$
%PF3-24-2025 We could say $0 \leq \kappa(E) \leq l$ because the case \kappa(E) = 0 would imply (?) that such connections are flat and those are excluded by w_2(\su(E)) not having an integral lift to a torsion class?
contain no flat or reducible connections and are regular.

When $g$ is \emph{generic} in the above sense of \cite{DK, FU}, then $M_\kappa^w(X,g)$ is a smooth (usually non-compact) manifold of the expected dimension if non-empty. If Conjecture \ref{mainconj:Existence_ASD_connection} holds, then $\su(E)$ admits a $g$-anti-self-dual connection when the basic lower bound \eqref{eq:p1_lower_bound} holds and the metric $g$ on $X$ is generic. The moduli space $M_\kappa^w(X,g)$ is thus a non-empty, smooth manifold and so
\[
  \dim M_\kappa^w(X,g) \geq 0.
\]
This yields the Bogomolov--Miyaoka--Yau inequality \eqref{eq:BMY} since
\begin{align*}
  0 &\leq 2\dim M_\kappa^w(X,g)
  \\
    &= -4p_1(\su(E))[X]-c_2(X) - c_1(X)^2 \quad\text{(by \eqref{eq:Expected_dimension_moduli_space_ASD_connections})}
  \\
    &\leq 3c_2(X) - c_1(X)^2  \quad\text{(by \eqref{eq:p1_lower_bound})}.
\end{align*}
Note that if $\kappa(E) = 0$, then an anti-self-dual connection on $\su(E)$ would be flat by \cite[Proposition 2.1.42, p. 43]{DK} and, although $M_0^w(X)$ would be non-empty, the Freed--Uhlenbeck generic metrics theorem would no longer imply that $M_0^w(X)$ is a smooth manifold and have $\dim M_0^w(X) \geq 0$.
% PF9-10-2024 Can we say anything about the case c_2(X) = 0?

The Bogomolov--Miyaoka--Yau inequality \eqref{eq:BMY} also follows from the non-negativity of the expected dimension of the moduli space of projectively anti-self-dual connections on a rank-two Hermitian bundle over the smooth blow-up of $X$. By smooth blow-up of $X$, we mean the connected sum $\widetilde X:=X\#\overline{\CC\PP}^2$ as in Gompf and Stipsicz \cite[Definition 2.2.7, p. 43]{GompfStipsicz}.  The characteristic numbers of $X$ and $\widetilde X$ satisfy
\begin{equation}
\label{eq:BlowUpCharNumbers}
c_1(\widetilde X)^2=c_1(X)-1 \quad\text{and}\quad c_2(\widetilde X)=c_2(X)+1.
\end{equation}
If $E$ is a complex rank-two Hermitian vector bundle on $\widetilde X$ satisfying the fundamental bounds,
\begin{equation}
\label{eq:p1_lower_bound_blowup}
0> p_1(\su(E))[\widetilde X] \ge -c_2(X),
\end{equation}
and if, for $w=c_1(E)$ and $\kappa=\kappa(E)$, the expected dimension of $M_\kappa^w(\widetilde X,g)$ is non-negative, then
%PF9-16-2024 The implication ("If ..., then") is confusing the logic. We first assume p_1 has above property. Period. We then say that if \dim M_\kappa^w \geq 0 (which is not implied by the p_1 property), then the below is a consequence of that and the p_1 lower bound. Please reference the previous argument giving BMY from \dim M_\kappa^w \geq 0 since the below is a repeat for \widetilde X
%then positivity of the expected dimension of the moduli space of
%projectively anti-self-dual connections on $E$
%would imply
%TL8-24-2025: Is current form OK with you?
%PF8-31-2025 No, it's not. Please rewrite it to be parallel to the previous version for X. Why are you writing expected dimension?
\begin{align*}
  0 &\leq 2\expdim M_\kappa^w(\widetilde X,g)
  \\
    &= -4p_1(\su(E))[\widetilde X]-c_2(\widetilde X) - c_1(\widetilde X)^2 \quad\text{(by \eqref{eq:Expected_dimension_moduli_space_ASD_connections})}
  \\
    &\leq 4c_2(X) -c_2(\widetilde X) - c_1(\widetilde X)^2    \quad\text{(by \eqref{eq:p1_lower_bound_blowup})}
  \\
   &= 3c_2(X)-c_1(X)^2\quad\text{(by \eqref{eq:BlowUpCharNumbers}).}
\end{align*}
Thus, to prove the Bogomolov--Miyaoka--Yau inequality \eqref{eq:BMY} holds for $X$, it suffices to prove
the following reformulation of Conjecture \ref{mainconj:Existence_ASD_connection}.

\begin{mainconj}[Existence of anti-self-dual connections with small instanton number on the smooth blow-up]
\label{mainconj:Existence_ASD_connection_BlowUp}
Assume the hypotheses of Conjecture \ref{mainconj:BMY_Seiberg-Witten} for $X$ and let $(E,H)$ be a complex rank two, Hermitian vector bundle over the smooth blow-up $\widetilde X=X\#\overline{\CC\PP}^2$ obeying the fundamental bounds \eqref{eq:p1_lower_bound_blowup}. Let $\tilde g$ be a Riemannian metric on $\widetilde X$ that is generic in the sense of Freed and Uhlenbeck \cite{DK, FU}. Then there exists a smooth, projectively $\tilde g$-anti-self-dual unitary connection $A$ on $E$, so the curvature $F_A \in \Omega^2(\fu(E))$ obeys
\begin{equation}
\label{eq:ASD_blowup}
  (F_A^+)_0 = 0 \in \Omega^+(\widetilde X;\su(E)).
\end{equation}
\end{mainconj}

\section[Frankel's Theorem for circle actions on an almost symplectic manifold]{Frankel's Theorem for the Hamiltonian function of a circle action on a smooth almost symplectic manifold}
\label{sec:Frankel_theorem_circle_actions_almost_Hermitian_manifolds}
The version, Theorem \ref{thm:Frankel_almost_Hermitian}, of Frankel's Theorem \cite[Section 3]{Frankel_1959} that we proved in our monograph, Feehan and Leness 
%TL11-26-2025: Update to \cite[Section 2.2, Theorem 6]{Feehan_Leness_introduction_virtual_morse_theory_so3_monopoles}
\cite[Section 1.1, Theorem 2]{Feehan_Leness_introduction_virtual_morse_theory_so3_monopoles}, is more general than that stated in \cite{Frankel_1959} because we allow for circle actions on closed, smooth manifolds $(M,g,J)$ that are only assumed to be \emph{almost Hermitian}\label{page:Almost_Hermitian_manifold}, so the $g$-orthogonal almost complex structure $J$ need not be integrable and the fundamental two-form
%PF9-9-2024 Switching convention to match almost everywhere else. Any unintended consequences?
\begin{equation}
  \label{eq:Fundamental_two-form}
  \omega = g(J\cdot,\cdot)
\end{equation}
is non-degenerate but not required to be closed, whereas Frankel assumed in \cite[Section 3]{Frankel_1959} that $\omega$ was closed. Recall that\footnote{If $E$ is a smooth vector bundle over a smooth manifold, we let $\Omega^0(E) = C^\infty(E)$ denote the Fr\'echet space of smooth sections of $E$.} $J \in C^\infty(\End(TM))$ is an \emph{almost complex structure} 
\label{page:Almost_complex_structure} on $M$ if $J^2 = -\id_{TM}$ and $(M,J)$ is thus an \label{page:Almost_complex_manifold} \emph{almost complex manifold}. One says that $J$ is \emph{orthogonal with respect to} or \emph{compatible with} a Riemannian metric $g$ on $M$ if
\begin{equation}
  \label{eq:g_compatible_J}
  g(JX,JY) = g(X,Y),
\end{equation}
for all vector fields $X, Y \in C^\infty(TM)$. Recall that a smooth manifold is called \emph{almost symplectic} \label{page:almost_symplectic_manifold} if it admits a non-degenerate two-form and \emph{symplectic} \label{page:symplectic_manifold} if that two-form is closed (see Libermann and Marle \cite[Definition 12.4]{Libermann_Marle_symplectic_geometry_analytical_mechanics}).
% COMMENT https://encyclopediaofmath.org/wiki/Almost-symplectic_structure

A circle action is called \emph{Hamiltonian} with respect to a circle-invariant smooth two-form $\omega$ on $M$ if \footnote{By analogy with the usual meaning \cite[Definition 2.1]{Dwivedi_Herman_Jeffrey_van_den_Hurk} of a Hamiltonian vector field and Hamiltonian function.} there exists a smooth function $f:M\to\RR$ such that
\begin{equation}
\label{eq:MomentMap}
  df = \iota_X\omega,
\end{equation}
where $X \in C^\infty(TM)$ is the vector field generated by the circle action, so $X_p = D_1\rho(1,p) \in T_pM$ for all $p\in M$, with $D_1\rho$ denoting the differential of $\rho$ in directions tangent to $S^1$.

\begin{rmk}[Behavior of Hamiltonian functions under pullback]
\label{rmk:Functorial_properties_Hamiltonian_functions}  
The following property of Hamiltonian functions will be useful in our applications. Suppose that $N$ is a smooth manifold with a circle action and that $Y \in C^\infty(TN)$ is the vector field generated by the circle action. If $F:N \to M$ is a smooth circle-equivariant map, then the vector fields $X$ on $M$ and $Y$ on $N$ are $F$-\emph{related} (see Lee \cite[Chapter 8, p. 182]{Lee_john_smooth_manifolds}) in the sense that
\[
  (dF)_p(Y_p) = X_{F(p)}, \quad\text{for all } p \in N,
\]
since $F(e^{i\theta}\cdot y) = e^{i\theta}\cdot F(y)$ for all $y \in N$. Equation \eqref{eq:MomentMap} implies that
\begin{multline*}
  F^*(df)_p(v) = d(F^*f)_p(v) = F^*(\iota_X\omega)_p(v) = \omega\left(X_p,(dF)_p(v)\right)
  \\
  = \omega\left((dF)_p(Y_p),(dF)_p(v)\right) = (F^*\omega)_p(Y_p,v) = (\iota_YF^*\omega)_p(v),
  \\
  \text{for all } p \in N \text{ and } v \in T_pN.
\end{multline*}
Therefore, $F^*f$ is a Hamiltonian function for the circle action on $N$ with respect to the circle-invariant smooth two-form $F^*\omega$ on $N$. 
\qed\end{rmk}

Adapting Bott \cite[Definition, p. 248]{Bott_1954}, \cite{Bott_1959} and Nicolaescu \cite[Definition 2.41]{Nicolaescu_morse_theory} (see also Feehan \cite[Definition 1.2, p. 3279]{Feehan_lojasiewicz_inequality_all_dimensions}), we make the

%PF5-22-2025 I don't think we need double labels like the below? They seem redundant. Won't a page ref to {defn:Morse-Bott_function} generate the page just as well as {page:Morse-Bott_function}?
\label{page:Morse-Bott_function}
\begin{defn}[Morse--Bott function]
\label{defn:Morse-Bott_function}  
Let $(M,g)$ be a smooth manifold and $f:M\to\RR$ be a smooth function. We let
\[
  \Crit f := \{p\in M: df(p)=0\}
\]
denote the \emph{critical set} of $f$. The function $f$ is \emph{Morse--Bott at $p$} if there exists an open neighborhood $U\subset M$ of $p$ such that $U\cap \Crit f$ is a smooth submanifold with tangent space
\[
  T_p\Crit f = \Ker\hess f(p),
\]
where $\hess f(p) \in \Hom(T_pM,T_p^*M)$ is the Hessian form defined in Feehan and Leness 
%TL11-26-2026: update to \cite[Section 3.1]{Feehan_Leness_introduction_virtual_morse_theory_so3_monopoles}?
\cite[Section 2.1]{Feehan_Leness_introduction_virtual_morse_theory_so3_monopoles}. The function $f$ is \emph{Morse--Bott} if it is Morse--Bott at every point of $\Crit f$.
\qed\end{defn}

\begin{thm}[Frankel's theorem for circle actions on almost Hermitian manifolds]
\label{thm:Frankel_almost_Hermitian}
(See Feehan and Leness 
%TL11-26-2025: Update to \cite[Section 2.2, Theorem 6]{Feehan_Leness_introduction_virtual_morse_theory_so3_monopoles}
\cite[Section 1.1, Theorem 2]{Feehan_Leness_introduction_virtual_morse_theory_so3_monopoles}; compare Frankel \cite[Section 3]{Frankel_1959}.)
Let $M$ be a smooth manifold endowed with a smooth circle action,
\begin{equation}
\label{eq:Circle_action_smooth_manifold}
\rho:S^1\times M\to M,
\end{equation}
and a non-degenerate two-form $\om$ that is circle-invariant. Then
\begin{enumerate}
\item
\label{item:Frankel_almost_Hermitian_FixedPointsAreZerosOfVField}
A point $p\in M$ is a fixed point of the action \eqref{eq:Circle_action_smooth_manifold} if and only if $X_p=0$, where $X \in C^\infty(TM)$ is the vector field generated by the circle action.

\item
\label{item:Frankel_almost_Hermitian_ComponentsOfFixedPointsAreSmoothSubmanifolds}
Each connected component of the fixed-point set of the circle action \eqref{eq:Circle_action_smooth_manifold} is a smooth submanifold of even dimension in $M$.
\end{enumerate}
In addition, let $f:M\to\RR$ be a smooth function that is Hamiltonian in the sense of \eqref{eq:MomentMap}. Then
\begin{enumerate}
\setcounter{enumi}{2}
\item
\label{item:Frankel_almost_Hermitian_FixedPointsAreCriticalPoints}
A point $p\in M$ is a critical point of $f$ if and only if $p$ is a fixed point of the circle action \eqref{eq:Circle_action_smooth_manifold}.
\end{enumerate}

\begin{enumerate}
\setcounter{enumi}{3}
\item
  \label{item:Frankel_almost_Hermitian_f_is_MB}
  The function $f$ is Morse--Bott at a critical point $p$ in the sense of Definition \ref{defn:Morse-Bott_function}.
\end{enumerate}
Furthermore, assume that there is a smooth Riemannian metric on $M$ that is circle-invariant. Then
\begin{enumerate}
\setcounter{enumi}{4}
\item
\label{item:Frankel_almost_Hermitian_ACisS1Invariant}
There are a smooth almost complex structure $J$ on $TM$ and a smooth Riemannian metric $g$ on $M$ such that  $(\omega,g,J)$ is a compatible triple in the sense that it obeys \eqref{eq:Fundamental_two-form} and \eqref{eq:g_compatible_J}, and $g$ and $J$ are circle-invariant.

\item
\label{item:Frankel_almost_Hermitian_WeightsAreEigenvalues}
The eigenvalues of the Hessian operator $\Hess_g f(p) \in \End(T_pM)$ are given by the weights of the circle action on $(T_pM,J)$ if the signs of the weights are chosen to be compatible with $J$ in the sense of Feehan and Leness
%TL11-26-2025: Updated 
\cite[Lemma 4.3.2 and Definition 4.3.3]{Feehan_Leness_introduction_virtual_morse_theory_so3_monopoles}.

\item
  \label{item:Frankel_almost_Hermitian_Sylvester}
  The signature $(\lambda_p^+(f),\lambda_p^0(f),\lambda_p^-(f))$ of the Hessian operator $\Hess_g f(p)$ is independent of the Riemannian metric $g$, where $\lambda_p^\pm(f)$ denotes the number of positive (negative) eigenvalues of $\Hess_g f(p)$ and $\lambda_p^0(f)$ denotes the nullity of $\Hess_g f(p)$.
\end{enumerate}
\end{thm}

If $p\in M$ is a critical point of $f$, then Theorem \ref{thm:Frankel_almost_Hermitian} implies that subspace $T_p^-M \subset T_pM$ on which the Hessian $\Hess_g f(p) \in \End(T_pM)$ is \emph{negative definite} is equal to the subspace of $T_pM$ on which the circle acts with \emph{negative weight}. Hence, the \emph{Morse--Bott index} \label{page:Morse--Bott_index} of $f$ at a critical point $p$,
\[
  \lambda_p^-(f) := \dim_\RR T_p^-M,
\]
is equal to the dimension of the subspace of $T_pM$ on which the circle acts with \emph{negative weight}; one calls $\lambda_p^+(f)$ and $\lambda_p^0(f)$ the \emph{Morse--Bott co-index} and \emph{nullity}, respectively.

\section[Virtual Morse--Bott signature of a critical point in a complex analytic space]{Virtual Morse--Bott signature of a critical point for the Hamiltonian function of a circle action on a complex analytic space}
\label{sec:Virtual_Morse-Bott_signature_Hamiltonian_function_circle_action_complex_analytic_space}
We recall the definition of virtual Morse--Bott signature from Feehan and Leness 
%TL11-26-2025: Update?
\cite[Section 1.3]{Feehan_Leness_introduction_virtual_morse_theory_so3_monopoles} and refer to Feehan \cite{Feehan_analytic_spaces} for a more detailed discussion and proofs.

Suppose that $(X,g,J,\omega)$ is a complex, finite-dimensional, Hermitian manifold with a real analytic circle action such that the almost complex structure $J$ is circle invariant and orthogonal with respect to the Riemannian metric\footnote{Following the convention of McDuff and Salamon \cite[Equation (1.4.1), p. 152]{McDuffSalamonSympTop3}.}, $g = \omega(\cdot,J\cdot)$, where $\omega = g(J\cdot,\cdot)$ is the fundamental two-form as in \eqref{eq:Fundamental_two-form}. We further assume that the circle action is Hamiltonian with real analytic Hamiltonian function $f:X\to\RR$ such that $df = \iota_\xi\omega$, where $\xi$ is the real analytic vector field on $X$ generated by the circle action. We let $Y\subset X$ be a circle-invariant, closed, complex analytic subspace (in the sense of Grauert and Remmert \cite[Chapter 1, Section 1.5, p. 7]{Grauert_Remmert_coherent_analytic_sheaves}), $p \in Y$ be a point, and $F:U\to\CC^m$ be a complex analytic, local defining function for $Y$ on an open neighborhood $U\subset X$ of $p$ in the sense that $Y\cap U = F^{-1}(0)$ and $\sO_{Y\cap U} = (\sO_U/\sI)\restriction Y\cap U$ is the structure sheaf for $Y\cap U$ defined by the ideal $\sI = (f_1,\ldots,f_m) \subset \sO_U$ given by the component functions $f_j$ of $F$. We let the complex linear subspace $\bH_p^2 \subset \CC^m$ be the orthogonal complement of $\Ran dF(p) \subset \CC^m$ and let $\bH_p^1 = \Ker dF(p) \subset T_pX$ denote the (complex linear)
\emph{Zariski tangent space} to $Y$ at $p$
\label{Zariski_Tangent_Space}. We let $S \subset X$ be the circle-invariant, complex, Hermitian submanifold given by $F^{-1}(\bH_p^2)$ and observe that $T_pS = \bH_p^1$. If $p$ is a \emph{critical point}
\label{page:Critical_point} of $f:Y\to\RR$ in the sense that $\bH_p^1 \subseteq \Ker df(p)$, then $p$ is a fixed point of the induced circle action on $S$ by Frankel's Theorem \ref{thm:Frankel_almost_Hermitian}. Because $p\in X$ is also a fixed point of the circle action on $X$, one can show \cite{Feehan_analytic_spaces} that the holomorphic map $F$ determines a circle action on $\CC^m$ and that the subspace $\bH_p^2$ is circle-invariant with respect to this induced action. We let $\bH_p^{k,1} \subset \bH_p^1$ and $\bH_p^{k,2} \subset \bH_p^2$ be the subspaces on which the circle acts with positive, zero, or negative weight depending on whether $k=+$, $0$, or $-$, respectively, so that
\[
  \bH_p^1 = \bH_p^{+,1} \oplus \bH_p^{0,1} \oplus \bH_p^{-,1}
  \quad\text{and}\quad
  \bH_p^2 = \bH_p^{+,2} \oplus \bH_p^{0,2} \oplus \bH_p^{-,2}.
\]  
We define the \emph{virtual Morse--Bott co-index, nullity} and \emph{index}, respectively, for $f$ at $p$ by
\begin{equation}
\label{eq:Virtual_Morse-Bott_signature}  
  \lambda_p^k(f) := \dim_\RR\bH_p^{k,1} - \dim_\RR\bH_p^{k,2}, \quad\text{for } k = +, 0, -.
\end{equation}
By Feehan \cite[Theorem 8]{Feehan_analytic_spaces}, if $\lambda_p^-(f) > 0$ (respectively, $\lambda_p^+(f) > 0$), then $p$ cannot be a local minimum (respectively, maximum) for $f:Y\to\RR$.

The preceding result suffices for the application discussed in Feehan and Leness \cite{Feehan_Leness_introduction_virtual_morse_theory_so3_monopoles} for the moduli space of non-Abelian monopoles over a compact, complex K\"ahler surface, where the open subset of that moduli space represented by pairs with non-vanishing coupled spinor (these are the non-zero-section points) has the structure of a circle-invariant complex analytic space. Each non-zero-section point in the moduli space that is regular has an open neighborhood that is a circle-invariant, complex Hermitian manifold. Each non-zero-section point in the moduli space that is singular has a circle-invariant open neighborhood defined by a circle-equivariant, holomorphic local model that is equivariantly embedded in a circle-invariant, complex Hermitian manifold $S$, as described above.

However, our present work is concerned with moduli spaces of non-Abelian monopoles over closed, symplectic four-manifolds. In this case, it no longer follows that the open subset of non-zero-section, regular points is necessarily a (circle-invariant) complex Hermitian manifold. However, Corollary \ref{maincor:Almost_Hermitian_structure_moduli_space_non-Abelian_monopoles_symplectic_4-manifolds} provides (circle-invariant) local virtual moduli spaces of non-Abelian monopoles that are almost Hermitian. Moreover, we expect that each non-zero-section point in the moduli space that is singular will have a circle-invariant open neighborhood defined by a real analytic local model that is approximately holomorphic in the sense of Chirka \cite{Chirka_1969}, Donaldson \cite{DonSympAlmostCx}, and Auroux \cite{AurouxFamilySymp}, where the approximation improves when $r \to \infty$, where $r \in [1,\infty)$ is the perturbation parameter in our choice of perturbed non-Abelian monopole equations \eqref{eq:SO(3)_monopole_equations_almost_Hermitian_perturbed_intro_regular}, analogous to Taubes' perturbation parameter $r$ is his choice of perturbed Seiberg--Witten monopole equations \cite[Section 1 (d), Equations (1.19) and (1.20)]{TauSWGromov}. We expect to provide an extension of \cite[Theorem 8]{Feehan_analytic_spaces} to the approximately holomorphic setting in a planned revision of \cite{Feehan_analytic_spaces}.

\section{Perturbed non-Abelian monopole equations over almost Hermitian four-manifolds}
\label{sec:Perturbed_non-Abelian_monopole_equations_almost_Hermitian_4-manifolds}
%PF9-12-2024 Add refs to this new section and edit out references to later versions of the same equations
In Feehan and Leness
%TL11-26-2026: Updated to \cite[Lemma 8.3.7, Remark 8.3.13, and Section 8.2.2]{Feehan_Leness_introduction_virtual_morse_theory_so3_monopoles}
\cite[Lemma 8.3.7, Remark 8.3.13, and Section 8.2.2]{Feehan_Leness_introduction_virtual_morse_theory_so3_monopoles}, we derived the form of the (unperturbed) non-Abelian monopole equations over almost Hermitian four-manifolds; see Okonek and Teleman \cite[Proposition 2.6, p. 900]{OTVortex} for the structure of the non-Abelian monopole equations over complex K\"ahler surfaces. In order to prove most of our main results in our present work (see Section \ref{sec:Main_results}), we shall need a perturbation of the non-Abelian monopole equations that extends the perturbation to the Seiberg--Witten monopole equations that was famously introduced by Taubes \cite{TauSWGromov} in his proof of the relation between Seiberg--Witten and Gromov invariants of symplectic four-manifolds. (See also the expositions of Taubes' work due to Donaldson \cite{DonSW} and Kotschick \cite{KotschickSW}.) That such a perturbation should exist for the non-Abelian monopole equations is not obvious.

\subsection{Non-Abelian monopole equations with a singular perturbation over almost Hermitian four-manifolds}
\label{subsec:Perturbed_non-Abelian_monopole_equations_almost_Hermitian_4-manifolds_singular}
Our first choice of perturbation --- see equations \eqref{eq:SO(3)_monopole_equations_almost_Hermitian_perturbed_intro} and \eqref{eq:Definition_wp_intro} --- allows us to achieve our goal of deriving pointwise estimates, similar to those of Taubes \cite{TauSWGromov}, for the components of the non-Abelian monopole coupled spinor. However, it does so at the cost of introducing a singularity in the non-Abelian monopole equations and this means that Sobolev solutions are no longer gauge equivalent to smooth solutions. While the latter technical issue is inconvenient but manageable, a far more serious problem is that the resulting system \eqref{eq:SO(3)_monopole_equations_almost_Hermitian_perturbed_intro} of perturbed non-Abelian monopole equations and a Coulomb gauge condition is not necessarily Fredholm (even with weighted Sobolev spaces). In Section \ref{subsec:Perturbed_non-Abelian_monopole_equations_almost_Hermitian_4-manifolds_regular}, we describe a regularization of the perturbation \eqref{eq:Definition_wp_intro} that yields a Fredholm system by standard methods and whose solutions are gauge equivalent to smooth solutions, but the derivation of pointwise estimates like those of Taubes \cite{TauSWGromov} is more complicated.

When $(A,\Phi)$ is a split solution to the perturbed non-Abelian monopole equations and thus equivalent to a solution to the perturbed Seiberg--Witten monopole equations, we show in Section \ref{sec:Perturbation_non-Abelian_monopole_equations_split_pairs} that the induced perturbation of the Seiberg--Witten monopole equations coincides exactly with that of Taubes in \cite{TauSWGromov}, without any singularity.
% PF9-12-2024 Add refs above

Let $(X,g,J,\omega)$ be a smooth almost Hermitian four-manifold, let $(\rho_\can,W_\can)$ be the canonical spin${}^c$ structure over $X$ (see the forthcoming \eqref{eq:Canonical_spinc_bundles} for the expression for $W_\can$ and the forthcoming \eqref{eq:Canonical_Clifford_multiplication} for the expression for $\rho_\can$), and let $(E,H)$ be a smooth, Hermitian vector bundle of rank at least two over $X$.

We choose $p \in (2,\infty)$ and consider the affine space of \emph{unitary triples},
\begin{equation}
  \label{eq:A_varphi_psi_in_W1p_times_W2p}
  (A,\varphi,\psi) \in \sA^{1,p}(E,A_d,H) \times W^{2,p}(E\oplus\Lambda^{0,2}(E))
\end{equation}
modulo the action \eqref{eq:GaugeActionOnSpinuPairs} of the Banach Lie group $W^{2,p}(\SU(E))$ of unitary automorphisms of $(E,H)$ that induce the identity automorphism of $\det E$. Since $X$ has real dimension $4$, then $W^{1,p}(X) \subset C^0(X)$ when $p > 4$ and $W^{2,p}(X) \subset C^0(X)$ when $p > 2$ by the Sobolev embedding theorem (see Adams and Fournier \cite[Theorem 4.12, p. 85]{AdamsFournier}). We abbreviate $\Lambda^{0,2}(E) = \Lambda^{0,2}(X)\otimes E$ and trust that this abuse of notation causes no confusion.

If $A$ is a $W^{1,p}$ unitary connection on $E$ and $\Phi$ is a $W^{2,p}$ section of $W_\can^+\otimes E$, then $(A,\Phi)$ is a solution to the \emph{unperturbed non-Abelian monopole equations} (see the forthcoming
\eqref{eq:SO(3)_monopole_equations}) with $\Phi=(\varphi,\psi) \in W^{2,p}(E\oplus\Lambda^{0,2}(E))$ if and only if $(A,\varphi,\psi)$ obeys (see the forthcoming Lemma \ref{lem:SO3_monopole_equations_almost_Kaehler_manifold})
%PF9-12-2024 Duplicate. Remove later copies.
\begin{subequations}
\label{eq:SO(3)_monopole_equations_almost_Hermitian_intro}
\begin{align}
  \label{eq:SO(3)_monopole_equations_(1,1)_curvature_intro}
  (\Lambda_\omega F_A)_0 &= \frac{i}{2}(\varphi\otimes\varphi^*)_0 - \frac{i}{2}\star(\psi\otimes\psi^*)_0,
  \\
  \label{eq:SO(3)_monopole_equations_(0,2)_curvature_intro}
  (F_A^{0,2})_0 &= \frac{1}{2}(\psi\otimes\varphi^*)_0,
  \\
  \label{eq:SO(3)_monopole_equations_Dirac_almost_Hermitian_intro}
  \bar{\partial}_A\varphi + \bar{\partial}_A^*\psi
  &= -\frac{1}{4\sqrt{2}}\rho_\can(\Lambda_{\omega}d\omega)(\varphi + \psi),
\end{align}
\end{subequations}
where the terms on the right hand sides of the system \eqref{eq:SO(3)_monopole_equations_almost_Hermitian_intro} are in $L^p(\su(E))$ and $L^p(\Lambda^{0,2}(\fsl(E)))$ and $L^p(\Lambda^{0,1}(E))$, respectively. It is occasionally convenient to include the complex conjugate of equation \eqref{eq:SO(3)_monopole_equations_(0,2)_curvature_intro} in the system \eqref{eq:SO(3)_monopole_equations_almost_Hermitian_intro}:
\begin{equation}
  \label{eq:SO(3)_monopole_equations_(2,0)_curvature_intro}
  (F_A^{2,0})_0 = -\frac{1}{2}(\varphi\otimes\psi^*)_0.
\end{equation}
Although equation \eqref{eq:SO(3)_monopole_equations_(2,0)_curvature_intro} is redundant, its inclusion in the system \eqref{eq:SO(3)_monopole_equations_almost_Hermitian_intro} facilitates translation between properties of solutions $(A,\Phi)$ to \eqref{eq:SO(3)_monopole_equations} and solutions $(A,\varphi,\psi)$ to \eqref{eq:SO(3)_monopole_equations_almost_Hermitian_intro}.

We shall call $(A,\varphi,\psi)$ in \eqref{eq:A_varphi_psi_in_W1p_times_W2p} a solution to the \emph{non-Abelian monopole equations with a singular Taubes perturbation} over a smooth almost Hermitian four-manifold if it obeys the system,
\begin{subequations}
\label{eq:SO(3)_monopole_equations_almost_Hermitian_perturbed_intro}
\begin{align}
  \label{eq:SO(3)_monopole_equations_(1,1)_curvature_perturbed_intro}
  (\Lambda_\omega F_A)_0 &= \frac{i}{2}(\varphi\otimes\varphi^*)_0 - \frac{i}{2}\star(\psi\otimes\psi^*)_0
                           - \frac{ir}{4}\wp(\psi),
  \\
  \label{eq:SO(3)_monopole_equations_(0,2)_curvature_perturbed_intro}
  (F_A^{0,2})_0 &= \frac{1}{2}(\psi\otimes\varphi^*)_0,
  \\
  \label{eq:SO(3)_monopole_equations_Dirac_almost_Hermitian_perturbed_intro}
  \bar{\partial}_A\varphi + \bar{\partial}_A^*\psi
  &= -\frac{1}{4\sqrt{2}}\rho_\can(\Lambda_{\omega}d\omega)(\varphi + \psi),
\end{align}
\end{subequations}
where the terms on the right hand sides of the system \eqref{eq:SO(3)_monopole_equations_almost_Hermitian_perturbed_intro} are again in $L^p(\su(E))$ and $L^p(\Lambda^{0,2}(\fsl(E)))$ and $L^p(\Lambda^{0,1}(E))$, respectively. The new term $\wp(\psi)$ appearing in equation \eqref{eq:SO(3)_monopole_equations_(1,1)_curvature_perturbed_intro} is defined by setting
\begin{equation}
  \label{eq:X0}
  X_0 := \{x \in X: \psi(x) \neq 0\},
\end{equation}
and choosing
\begin{equation}
  \label{eq:Definition_wp_intro}
  \wp(\psi)
  :=
  \begin{cases}
    2|\psi|_{\Lambda^{0,2}(E)}^{-2}\star(\psi\otimes\psi^*)_0 &\text{on } X_0,
    \\
    0 &\text{on } X\less X_0.
  \end{cases}
\end{equation}
Again, it is occasionally convenient to include (even though redundant) the complex conjugate \eqref{eq:SO(3)_monopole_equations_(2,0)_curvature_intro} of equation \eqref{eq:SO(3)_monopole_equations_(0,2)_curvature_perturbed_intro} in the system \eqref{eq:SO(3)_monopole_equations_almost_Hermitian_perturbed_intro}.

The motivation for the definition \eqref{eq:Definition_wp_intro} is discussed in the forthcoming Section \ref{subsec:Local_Kuranishi_model_non-Abelian_monopole_equations_Taubes_perturbation}. Observe that $\wp(\psi) \in L^\infty(i\su(E))$ even if $\psi \in C^\infty(\Lambda^{0,2}(E))$. Equations \eqref{eq:SO(3)_monopole_equations_(0,2)_curvature_perturbed_intro} and \eqref{eq:SO(3)_monopole_equations_Dirac_almost_Hermitian_perturbed_intro} are identical to Equations \eqref{eq:SO(3)_monopole_equations_(0,2)_curvature_intro} and \eqref{eq:SO(3)_monopole_equations_Dirac_almost_Hermitian_intro}, respectively. Equation \eqref{eq:SO(3)_monopole_equations_(1,1)_curvature_perturbed_intro} is equivalent to (see Section \ref{subsec:Local_Kuranishi_model_non-Abelian_monopole_equations_Taubes_perturbation})
\begin{equation}
  \label{eq:SO(3)_monopole_equations_(1,1)_curvature_perturbed_omega_intro}
  (F_A^\omega)_0 = \frac{i}{4}(\varphi\otimes\varphi^*)_0\omega - \frac{i}{4}\star(\psi\otimes\psi^*)_0\omega
     - \frac{ir}{8}\wp(\psi)\,\omega \in L^p(\su(E))\cdot\omega,
\end{equation}
where the right-hand side is now in $L^p(\su(E))\cdot\omega$ and $F_A^\omega$ is defined by the projection of $F_A$ onto its image in the factor $\su(E)\cdot\omega$ of the orthogonal decomposition,
\[
  \Lambda^+(\fsl(E)) = \Lambda^{2,0}(\fsl(E))\oplus \fsl(E)\cdot\omega \oplus \Lambda^{0,2}(\fsl(E)),
\]
noting that $\fsl(E) \cong i\underline{\RR}\oplus\su(E)$, where $\underline{\RR} := X\times\RR$.
% PF9-12-2024 Cut duplication with later sections.

We are very grateful to Cliff Taubes for pointing out that the system \eqref{eq:SO(3)_monopole_equations_almost_Hermitian_perturbed_intro} of non-Abelian monopole equations with the perturbation term defined by $\wp(\psi)$ in \eqref{eq:Definition_wp_intro} does \emph{not} define a Fredholm map even after restriction to Coulomb gauge slice. Therefore, in the forthcoming Section \ref{subsec:Perturbed_non-Abelian_monopole_equations_almost_Hermitian_4-manifolds_regular} we consider a regularization of $\wp(\psi)$ in \eqref{eq:Definition_wp_intro} that does yield a Fredholm map upon restriction to a Coulomb gauge slice.

\subsection{Non-Abelian monopole equations with a regularized perturbation over almost Hermitian four-manifolds}
\label{subsec:Perturbed_non-Abelian_monopole_equations_almost_Hermitian_4-manifolds_regular}
For our second choice of perturbation, we consider the affine space of $W^{1,p}$ \emph{unitary triples},
\begin{equation}
  \label{eq:A_varphi_psi_in_W1p}
  (A,\varphi,\psi) \in \sA^{1,p}(E,A_d,H) \times W^{1,p}(E\oplus\Lambda^{0,2}(E)).
\end{equation}
We replace the role of $\wp(\psi)$ in \eqref{eq:Definition_wp_intro} and \eqref{eq:SO(3)_monopole_equations_(1,1)_curvature_perturbed_intro} via a regularization defined by a constant
$\gamma \in (0,\infty)$ and whose choice is specified in the forthcoming Theorem \ref{mainthm:Taubes_1996_SW_to_Gr_2-3_regular}:
\begin{equation}
  \label{eq:Definition_wp_intro_regular}
  \wp_\gamma(\psi)
  :=
  4\left(\gamma^2 + |\psi|_{\Lambda^{0,2}(E)}^2\right)^{-1}\star(\psi\otimes\psi^*)_0.
\end{equation}
While $\wp(\psi) \in L^\infty(i\su(E))$, we now have $ \wp_\gamma(\psi) \in C^\infty(i\su(E))$ if $\psi \in C^\infty(\Lambda^{0,2}(E))$. We shall call $(A,\varphi,\psi)$ in \eqref{eq:A_varphi_psi_in_W1p} a solution to the \emph{non-Abelian monopole equations with a regularized Taubes perturbation} over a smooth almost Hermitian four-manifold $(X,g,J,\omega)$ if it obeys
\begin{subequations}
\label{eq:SO(3)_monopole_equations_almost_Hermitian_perturbed_intro_regular}
\begin{align}
  \label{eq:SO(3)_monopole_equations_(1,1)_curvature_perturbed_intro_regular}
  (\Lambda_\omega F_A)_0 &= \frac{i}{2}(\varphi\otimes\varphi^*)_0 - \frac{i}{2}\star(\psi\otimes\psi^*)_0
                           - \frac{ir}{4}\wp_\gamma(\psi),
  \\
  \label{eq:SO(3)_monopole_equations_(0,2)_curvature_perturbed_intro_regular}
  (F_A^{0,2})_0 &= \frac{1}{2}(\psi\otimes\varphi^*)_0,
  \\
  \label{eq:SO(3)_monopole_equations_Dirac_almost_Hermitian_perturbed_intro_regular}
  \bar{\partial}_A\varphi + \bar{\partial}_A^*\psi
  &= -\frac{1}{4\sqrt{2}}\rho_\can(\Lambda_{\omega}d\omega)(\varphi + \psi),
\end{align}
\end{subequations}
where $\wp_\gamma(\psi)$ is as in \eqref{eq:Definition_wp_intro_regular} and the terms on the right hand sides of the system \eqref{eq:SO(3)_monopole_equations_almost_Hermitian_perturbed_intro_regular} are again in $L^p(\su(E))$ and $L^p(\Lambda^{0,2}(\fsl(E)))$ and $L^p(\Lambda^{0,1}(E))$, respectively. As before, it is occasionally convenient to include the (redundant) complex conjugate \eqref{eq:SO(3)_monopole_equations_(2,0)_curvature_intro} of equation \eqref{eq:SO(3)_monopole_equations_(0,2)_curvature_perturbed_intro_regular} in the system \eqref{eq:SO(3)_monopole_equations_almost_Hermitian_perturbed_intro_regular}.

Equations \eqref{eq:SO(3)_monopole_equations_(0,2)_curvature_perturbed_intro_regular} and \eqref{eq:SO(3)_monopole_equations_Dirac_almost_Hermitian_perturbed_intro_regular} are identical to Equations \eqref{eq:SO(3)_monopole_equations_(0,2)_curvature_intro} and \eqref{eq:SO(3)_monopole_equations_Dirac_almost_Hermitian_intro}, respectively, while Equation \eqref{eq:SO(3)_monopole_equations_(1,1)_curvature_perturbed_intro_regular} is equivalent to
\begin{equation}
  \label{eq:SO(3)_monopole_equations_(1,1)_curvature_perturbed_omega_intro_regular}
  (F_A^\omega)_0 = \frac{i}{4}(\varphi\otimes\varphi^*)_0\omega - \frac{i}{4}\star(\psi\otimes\psi^*)_0\omega
     - \frac{ir}{8}\wp_\gamma(\psi)\,\omega,
\end{equation}
where the right-hand side is in $L^p(\su(E))\cdot\omega$. The regularization \eqref{eq:Definition_wp_intro_regular} allows us to
\begin{inparaenum}[\itshape i\upshape)]
\item achieve our goal of deriving pointwise estimates for suitable $\gamma$, similar to those of Taubes \cite{TauSWGromov}, for the components of the non-Abelian monopole coupled spinor,
\item ensure that Sobolev solutions are gauge equivalent to smooth solutions, and
\item guarantee that the system \eqref{eq:SO(3)_monopole_equations_almost_Hermitian_perturbed_intro_regular} of non-Abelian monopole equations with a regularized Taubes perturbation and a Coulomb gauge condition defines a Fredholm map (see the forthcoming Section \ref{subsec:Local_Kuranishi_model_nbhd_point_zero_locus_unperturbed_non-Abelian_monopole_equations}).
\end{inparaenum}

%PF5-20-2025 Check for later overlap
We denote the \emph{quotient space} of $W^{1,p}$ unitary triples $(A,\varphi,\psi)$ modulo the action \eqref{eq:GaugeActionOnSpinuPairs} of $W^{2,p}(\SU(E))$ by
\begin{equation}
  \label{eq:Quotient_space_unitary_triples}
  \sC(E,H,J,A_d)
  := \left.\left(\sA^{1,p}(E,A_d,H) \times W^{1,p}(E\oplus\Lambda^{0,2}(E))\right)\right/W^{2,p}(\SU(E))
\end{equation}
By analogy with the forthcoming Definition \ref{defn:Split_trivial_central-stabilizer_spinor_pair}  \eqref{item:Split_spinor_pair} of split or non-split \spinu pairs, we write
\begin{subequations}
\label{eq:Quotient_space_non-zero_or_non-split_unitary_triples}
\begin{align}
  \label{eq:Quotient_space_non-zero-section_unitary_triples}
  \sC^0(E,H,J,A_d)
  &:=
    \left\{[A,\varphi,\psi] \in \sC(E,H,J,A_d): (\varphi,\psi)\not\equiv (0,0)\right\},
  \\
  \label{eq:Quotient_space_non-split_unitary_triples}
  \sC^*(E,H,J,A_d)
  &:=
    \left\{[A,\varphi,\psi] \in \sC(E,H,J,A_d): (A,\varphi,\psi) \text{ non-split}\right\},
  \\
  \label{eq:Quotient_space_non-split_non-zero-section_unitary_triples}
  \sC^{*,0}(E,H,J,A_d)
  &:=
    \sC^*(E,H,J,A_d) \cap \sC^0(E,H,J,A_d),
\end{align}
\end{subequations}
for the quotient subspaces of non-zero-section, non-split, and non-split non-zero-section unitary triples, respectively.
  
We denote the \emph{moduli subspace of non-Abelian monopoles with a regularized Taubes perturbation} by
\begin{equation}
  \label{eq:Moduli_space_non-Abelian_monopoles_almost_Hermitian_Taubes_regularized}
  \sM(E,g,J,\omega,r)
  :=
  \left\{[A,\varphi,\psi] \in \sC(E,H,J,A_d): (A,\varphi,\psi) \text{ obeys }
  \eqref{eq:SO(3)_monopole_equations_almost_Hermitian_perturbed_intro_regular} \right\}.
\end{equation}
and write
\begin{subequations}
\label{eq:Moduli_space_non-Abelian_monopoles_almost_Hermitian_Taubes_regularized_non-zero_or_non-split}
\begin{align}
  \label{eq:Moduli_space_non-Abelian_monopoles_almost_Hermitian_Taubes_regularized_non-zero-section}
  \sM^0(E,g,J,\omega,r) &:= \sM(E,g,J,\omega,r) \cap \sC(E,H,J,A_d),
  \\
  \label{eq:Moduli_space_non-Abelian_monopoles_almost_Hermitian_Taubes_regularized_non-split}
  \sM^*(E,g,J,\omega,r) &:= \sM(E,g,J,\omega,r) \cap \sC^*(E,H,J,A_d),
  \\
  \label{eq:Moduli_space_non-Abelian_monopoles_almost_Hermitian_Taubes_regularized_non-zero_non-split}
  \sM^{*,0}(E,g,J,\omega,r) &:= \sM^0(E,g,J,\omega,r) \cap \sM^*(E,g,J,\omega,r),
\end{align}
\end{subequations}
for the moduli subspaces of non-zero-section, non-split, and non-split non-zero-section non-Abelian monopoles, respectively.

\section{Main results}
\label{sec:Main_results}
In this section, we state our main results.

\subsection{Guide to main results}
\label{subsec:Guide_main_results}
The overall goal of our program initiated in \cite{Feehan_analytic_spaces, Feehan_Leness_introduction_virtual_morse_theory_so3_monopoles} is to prove Conjecture \ref{mainconj:BMY_symplectic} and, more generally, Conjecture \ref{mainconj:BMY_Seiberg-Witten}. In this monograph, we build on our previous monographs \cite{Feehan_analytic_spaces, Feehan_Leness_introduction_virtual_morse_theory_so3_monopoles} by establishing important further steps in our program, though we do not yet prove Conjecture \ref{mainconj:BMY_Seiberg-Witten}. We believe that our results have independent interest, beyond their immediate application to that goal. In Section \ref{sec:Further_work}, we outline the remaining steps required to complete a proof of Conjecture \ref{mainconj:BMY_Seiberg-Witten}.

Not every main result that we describe in Section \ref{sec:Main_results} is required by our program to prove Conjecture \ref{mainconj:BMY_Seiberg-Witten}, but the auxiliary results that we include in Section \ref{sec:Main_results} nevertheless provide a scaffolding and motivation for our central results, not all of which are stated in Section \ref{sec:Main_results}. We provide in this subsection a guide for readers who wish to focus on the central results in Section \ref{sec:Main_results} and discuss the significance of the auxiliary results that complement the central ones.

In Section \ref{subsec:Feasibility_non-Abelian_monopole_cobordism_method}, we introduce in Definition \ref{maindefn:Feasibility} the concept of \emph{feasibility} of a spin${}^u$ structure for our program to prove Conjecture \ref{mainconj:BMY_Seiberg-Witten} via the moduli space of non-Abelian monopoles. This generalizes the concept of feasibility employed in Feehan and Leness \cite{Feehan_Leness_introduction_virtual_morse_theory_so3_monopoles}. The existence of feasible spin${}^u$ structures is established by Theorem \ref{mainthm:ExistenceOfSpinuForFlow} and Corollary \ref{maincor:FeasibleIsNonEmpty}.

Section \ref{subsec:Generalization_Donaldson_symplectic_subspace_criterion} begins by recalling Donaldson's result on the existence of symplectic submanifolds (Theorem \ref{thm:Donaldson_1}) and his elementary but important and widely cited symplectic subspace criterion that underpins his proof of existence of symplectic submanifolds. We then describe a refinement, Proposition \ref{mainprop:Donaldson_1996jdg_3}, of Donaldson's symplectic subspace criterion in Proposition \ref{prop:Donaldson_1996jdg_3} from the context of real linear maps $\sL$ from $\CC^m$ to $\CC$ to the case of continuous real linear maps $\sL$ from a complex Hilbert space $\sH$ to $\CC$. We generalize Proposition \ref{mainprop:Donaldson_1996jdg_3} in Theorem \ref{mainthm:Donaldson_1996jdg_3_Hilbert_space_codomain} and its immediate consequences, namely, Corollaries \ref{maincor:Donaldson_1996jdg_3_Hilbert_space_codomain_families} and
\ref{maincorDonaldson_1996jdg_3_Hilbert_space_domain}.

We give another generalization of Proposition \ref{prop:Donaldson_1996jdg_3} for the case of bounded real linear operators between complex Banach spaces in Proposition \ref{mainprop:Donaldson_1996jdg_3_Banach_space} and Corollary \ref{maincor:Adjoint_DonaldsonCriteria}. In our application, the (densely defined unbounded) linear operator $\sT$ would be the first order elliptic operator defined in the forthcoming \eqref{eq:Perturbed_Deformation_Operator}, namely
\begin{multline*}
  \sT_{A,\varphi,\psi,r} = d_{A,\varphi,\psi,r}^1 + d_{A,\varphi,\psi}^{0,*}:
  L^2\left(T^*X\otimes\su(E) \oplus E\oplus \Lambda^{0,2}(E)\right)
  \\
  \to
  L^2\left(\su(E) \oplus \su(E) \oplus \Lambda^{0,2}(\fsl(E)) \oplus \Lambda^{0,1}(E)\right),
\end{multline*}
given by the sum of the linearization \eqref{eq:Linearization_TaubesPert_Regular} at a smooth solution $(A,\varphi,\psi)$ of the non-Abelian monopole equations \eqref{eq:SO(3)_monopole_equations_almost_Hermitian_perturbed_intro_regular} with a regularized Taubes perturbation (with parameter $r$) and the $L^2$ adjoint of the operator \eqref{eq:d_APhi^0}, which defines the Coulomb gauge condition \eqref{eq:d_APhi^0_star_aphi_identity_and_vanishing}.

The kernel of the operator $\sT_{A,\varphi,\psi,r}$ defines the Zariski tangent space $\bH_{A,\varphi,\psi,r}^1$ in \eqref{eq:H1_For_TaubesPert_regular} to the moduli space $\sM^0(E,g,J,\omega,r)$ at a point $[A,\varphi,\psi]$ represented by a non-zero section smooth unitary triple $(A,\varphi,\psi,r)$. A generalization of Donaldson's symplectic subspace criterion in Proposition \ref{prop:Donaldson_1996jdg_3} would, ideally, show that if the complex antilinear part $\sT_{A,\varphi,\psi,r}''$ of $\sT_{A,\varphi,\psi,r}$ is small enough, then the almost complex structure $\bJ$ and weak $L^2$ Riemannian metric $\bg$ on the affine space of unitary triples induce a symplectic structure on $\bH_{A,\varphi,\psi,r}^1$. This symplectic structure would extend to a circle invariant non-degenerate two-form $\bomega$ on an circle invariant open neighborhood of $[A,\varphi,\psi]$ in $\sM^0(E,g,J,\omega,r)$. We would then aim to verify that Hitchin's function $-f$ in \eqref{eq:Hitchin_function_moduli_space} is a Hamiltonian for the circle action on the preceding open neighborhood. This in turn would allow one to conclude via Frankel's Theorem \ref{thm:Frankel_almost_Hermitian} \eqref{item:Frankel_almost_Hermitian_FixedPointsAreCriticalPoints} that critical points of $f$ on $\sM^0(E,g,J,\omega,r)$ are equivalent to fixed points of the circle action, which are necessarily points represented by solutions to the Seiberg--Witten monopole equations \cite{FL2a}.

Theorem \ref{mainthm:Donaldson_1996jdg_3_Hilbert_space_codomain} generalizes Proposition \ref{mainprop:Donaldson_1996jdg_3} from the setting of bounded $\RR$-linear functions $\sT:\sH\to\CC$ on a complex Hilbert space $\sH$ to that of bounded $\RR$-linear operators $\sT:\sH_1\to\sH_2$ between complex Hilbert spaces, $\sH_1$ and $\sH_2$. Like Proposition \ref{mainprop:Donaldson_1996jdg_3}, Theorem \ref{mainthm:Donaldson_1996jdg_3_Hilbert_space_codomain} asserts that if $\sT$ is suitably generic, then $\Ker\sT$ is a symplectic subspace of $\sH_1$, while Corollary \ref{maincorDonaldson_1996jdg_3_Hilbert_space_domain} asserts that $\Ker\sT^*$ is a symplectic subspace of $\sH_2$ if $\sT$ is suitably generic. We shall explore refinements of these results elsewhere and we do not make use of Theorem \ref{mainthm:Donaldson_1996jdg_3_Hilbert_space_codomain} or Corollary \ref{maincorDonaldson_1996jdg_3_Hilbert_space_domain} in this monograph.

While the generalization of Proposition \ref{prop:Donaldson_1996jdg_3} given by Proposition \ref{mainprop:Donaldson_1996jdg_3_Banach_space} is also of interest -- for example, it may apply to proofs of results due to Auroux in \cite{AurouxFamilySymp} for real linear maps $\sT:\CC^n\to\CC^m$ --  it is not easily applicable to operators such as $\sT_{A,\varphi,\psi,r}$, even after replacing $L^2$ in its domain by $W^{1,p}$ and $L^2$ in its codomain by $L^p$ (say with $p\geq 2$), so $\sT_{A,\varphi,\psi,r}$ becomes a bounded operator between Banach spaces. The reason is that the inequality \eqref{eq:Donaldson_1996jdg_Prop_3_Banach_space} required by the hypotheses of Proposition \ref{mainprop:Donaldson_1996jdg_3_Banach_space} is too difficult to verify, especially if one demands that \eqref{eq:Donaldson_1996jdg_Prop_3_Banach_space} hold uniformly with respect to the triple $(A,\varphi,\psi)$ and parameter $r$.

Our next pair of results, Theorem \ref{mainthm:Donaldson_1996jdg_3_Hilbert_space} and Corollary \ref{maincor:Donaldson_1996jdg_3_Hilbert_space_non-self-adjoint}, provide more significant generalizations of Proposition \ref{prop:Donaldson_1996jdg_3}. On the one hand, the hypotheses of Theorem \ref{mainthm:Donaldson_1996jdg_3_Hilbert_space} are far more precise than those Proposition \ref{mainprop:Donaldson_1996jdg_3_Banach_space}. Their application to the operators $\sT_{A,\varphi,\psi,r}$ and $T_{A,\varphi,\psi,r}$ would focus on the existence of a symmetric uniform gap \eqref{eq:Spectral_gap_T} of width $4\delta$ somewhere in the spectrum of the self-adjoint operator,
\[
  T_{A,\varphi,\psi,r}
  =
  \begin{pmatrix} 0 & \sT_{A,\varphi,\psi,r}^* \\ \sT_{A,\varphi,\psi,r} & 0 \end{pmatrix},
\]
and the estimate \eqref{eq:Donaldson_1996jdg_Prop_3_Hilbert_space},
\[
  \|T_{A,\varphi,\psi,r}''\|_{\End(\sH)} < \delta/2,
\]
for the complex antilinear component $T_{A,\varphi,\psi,r}''$ of $T_{A,\varphi,\psi,r}$ or the estimate \eqref{eq:Donaldson_1996jdg_Prop_3_Hilbert_space_non-self-adjoint}, 
\[
  \|\sT_{A,\varphi,\psi,r}''\|_{\End(\sH)} < \delta/4,
\]
for the complex antilinear component $\sT_{A,\varphi,\psi,r}''$ of $\sT_{A,\varphi,\psi,r}$. The operators $\sT_{A,\varphi,\psi,r}$ and $T_{A,\varphi,\psi,r}$ are both zero order with smooth coefficients and thus compact. Unfortunately, an application of Theorem \ref{mainthm:Donaldson_1996jdg_3_Hilbert_space} and Corollary \ref{maincor:Donaldson_1996jdg_3_Hilbert_space_non-self-adjoint} to the operator $\sT_{A,\varphi,\psi,r}$ ultimately appears stymied by two inherent difficulties.

The first difficulty is revealed by an examination of the explicit expression \eqref{eq:H_dbar_APhi^01_explicit_perturbed} for the equivalent operator defined in \eqref{eq:Perturbed_AC_Deformation_Operator}, namely
\begin{multline}
  \label{eq:Equivalence_cT_d-bar_A_varphi_psi_r_and_sT_A_varphi_psi_r}
  \cT_{\bar\partial_A,\varphi,\psi,r}
  =
  \Upsilon_2 \circ \sT_{A,\varphi,\psi,r} \circ \Upsilon_1^{-1}:
  L^2\left( \Lambda^{0,1}(\fsl(E)) \oplus E\oplus \Lambda^{0,2}(E)\right)
  \\
  \to
  L^2\left(\fsl(E) \oplus \Lambda^{0,2}(\fsl(E)) \oplus \Lambda^{0,1}(E)\right),
\end{multline}
where $\Upsilon_1$ and $\Upsilon_2$ are the real linear isomorphisms defined in \eqref{eq:Isomorphism_sEkC_to_sEk}. In the system \eqref{eq:H_dbar_APhi^01_explicit_perturbed} of equations that define the kernel of $\cT_{\bar\partial_A,\varphi,\psi,r}$, we see that there are coefficients that cannot be made $C^0$-small (even after rescaling) by choosing the Taubes parameter $r$ to be sufficiently large and applying our estimates in Section \ref{subsec:Analogues_non-Abelian_monopoles_Taubes_estimates_SW_monopole_sections} (namely, Proposition \ref{mainprop:Taubes_1996_SW_to_Gr_2-1_regular} and Theorem \ref{mainthm:Taubes_1996_SW_to_Gr_2-3_regular}). For example, choosing $r$ large has no effect on the Nijenhuis tensor $N_J$; the linearization \eqref{eq:Derivative_of_Taubes_Pert_regular} of the Taubes perturbation $\wp_\gamma(\psi)$ in \eqref{eq:Definition_wp_intro_regular}, contributes terms that are potentially $C^0$-large when $r$ is large; and the term $\star((a'')^\dagger\wedge\star\psi)$ may also be $C^0$-large. Naturally, we wish to avoid relying on a technical hypothesis that $N_J$ is $C^0$-small. See Remark \ref{rmk:Replacing_NJ_C0_small_by_small_Nijenhuis_energy} for a history of efforts to prove that one can choose $J$ such that $N_J$ is small in various senses. Proposition \ref{mainprop:Taubes_1996_SW_to_Gr_2-1_regular} and Theorem \ref{mainthm:Taubes_1996_SW_to_Gr_2-3_regular} do play a role in showing that the \emph{non-linear} terms\footnote{Excluding the perturbation term $-ir\wp_\gamma(\psi)/4 \in \Omega^0(\su(E))$ itself in the $\Omega^0(\su(E))$-equation \eqref{eq:SO(3)_monopole_equations_(1,1)_curvature_perturbed_intro_regular}.}
in the system \eqref{eq:SO(3)_monopole_equations_almost_Hermitian_perturbed_intro_regular} of equations that are not holomorphic in $(a'',\varphi,\psi)$ can be made relatively small (after rescaling). This first difficulty is explained in more detail in Section \ref{sec:Unperturbed_non-Abelian_monopole_equations_almost_Kaehler_four-manifolds}.

The second difficulty is that it appears very hard to verify a spectral gap hypothesis \eqref{eq:Spectral_gap_T} that is uniform with respect to the triple $(A,\varphi,\psi)$ and the parameter $r$. Indeed, as $r$ becomes larger, the gap between successive distinct eigenvalues in the spectrum of $T_{A,\varphi,\psi,r}$ tends to become smaller and, of course, the energy bubbling phenomenon for the moduli space of non-Abelian monopoles (see Feehan and Leness \cite{FL1}) compounds this problem. See Theorem \ref{mainthm:Lower_bound_spectral_gap_coupled_Dirac_operator} and Corollaries \ref{maincor:Lower_bound_spectral_gap_coupled_Dirac_operator_plus_vector_potential} and \ref{maincor:Lower_bound_spectral_gap_coupled_Dirac_operator_plus_vector_potential_four-manifold} in Section \ref{subsec:Main_results_lower_bounds_spectral_gaps} for results along these lines.

We include the spectral gap estimates in Section \ref{subsec:Main_results_lower_bounds_spectral_gaps} to illustrate an obstruction to what would otherwise be an elegant application of Corollary \ref{maincor:Donaldson_1996jdg_3_Hilbert_space_non-self-adjoint} to the construction of non-degenerate two-forms and almost complex structures on the moduli space of non-Abelian monopoles. We include the results in Section \ref{subsec:Main_results_lower_bounds_spectral_gaps} for completeness, but they are not central and are not required by our program to prove Conjecture \ref{mainconj:BMY_Seiberg-Witten}.

In our discussion of our proofs of Theorem \ref{mainthm:AH_structure_bounded_evalue_spaces_non-Abelian_monopoles_symp_4-mflds}, Corollary \ref{maincor:Almost_Hermitian_structure_moduli_space_non-Abelian_monopoles_symplectic_4-manifolds}, and Theorem \ref{mainthm:IdentifyCriticalPoints} further below, we shall explain our alternative construction of construction of non-degenerate two-forms and almost complex structures on the moduli space of non-Abelian monopoles as well as the fact that Hitchin's function $-f$ is a Hamiltonian for the circle action on $\sM^0(E,g,J,\omega,r)$.

In Section \ref{subsec:Analogues_non-Abelian_monopoles_Taubes_estimates_SW_monopole_sections}, we present analogues for non-Abelian monopoles of Taubes' pointwise estimates for Seiberg--Witten monopoles \cite{Taubes_2000} --- see Theorem \ref{mainthm:Taubes_1996_SW_to_Gr_2-3_regular} and Corollary \ref{maincor:Taubes_1996_SW_to_Gr_eq_2-12_and_2-13_regular}. As in \cite{Taubes_2000}, we find that when $(A,\varphi,\psi)$ is a solution to the system \eqref{eq:Moduli_space_non-Abelian_monopoles_almost_Hermitian_Taubes_regularized}, then the $C^0$ norm of $r^{-1/2}\psi$ is $O(r^{-1/2})$ while the $C^0$ norm of $r^{-1/2}\varphi$ is $O(1)$ as $r\to \infty$. (When $(X,g,J,\omega)$ is a complex K\"ahler surface, we proved in \cite{Feehan_Leness_introduction_virtual_morse_theory_so3_monopoles} that one has $\psi \equiv 0$ and, because we also have $N_J \equiv 0$ for complex K\"ahler surfaces, the operator $T_{A,\varphi,\psi,r}$ becomes complex linear.) We shall apply Theorem \ref{mainthm:Taubes_1996_SW_to_Gr_2-3_regular} and Corollary \ref{maincor:Taubes_1996_SW_to_Gr_eq_2-12_and_2-13_regular} in a sequel to our present monograph to obtain a type of approximate holomorphicity for the obstruction map in certain Kuranishi models for circle invariant open neighborhoods of fixed points in the moduli space $\sM^0(E,g,J,\omega,r)$.

%PF7-21-2025 Add refs to discussion below
Section \ref{subsec:Almost_Hermitian_structure_moduli_spaces_non-Abelian_monopoles_symplectic_4-manifolds} contains our statements of Theorem \ref{mainthm:AH_structure_bounded_evalue_spaces_non-Abelian_monopoles_symp_4-mflds} and Corollary \ref{maincor:Almost_Hermitian_structure_moduli_space_non-Abelian_monopoles_symplectic_4-manifolds}. Rather than try to prove that $\Ker\sT_{A,\varphi,\psi,r} = \bH_{A,\varphi,\psi,r}^1$ admits a symplectic structure by an application of Corollary \ref{maincor:Donaldson_1996jdg_3_Hilbert_space_non-self-adjoint}, we shall instead exploit the fact that the Hilbert space
\[
  \sH_1 = L^2\left(T^*X\otimes\su(E) \oplus E\oplus \Lambda^{0,2}(E)\right)
\]
has a complete orthonormal basis of eigenvectors for the complex linear Laplacian $\sT_{A,\varphi,\psi,r}^{\prime,*}\sT_{A,\varphi,\psi,r}'$ defined by the complex linear part $\sT_{A,\varphi,\psi,r}'$ of $\sT_{A,\varphi,\psi,r}$. If $\nu$ is a positive constant that is not an eigenvalue of $\sT_{A,\varphi,\psi,r}^{\prime,*}\sT_{A,\varphi,\psi,r}'$ and $\bH_{A,\varphi,\psi,r,\nu}^1$ is the direct sum of the eigenvectors of $\sT_{A,\varphi,\psi,r}^{\prime,*}\sT_{A,\varphi,\psi,r}'$ with eigenvalues less than $\nu$, then the \emph{bounded eigenvalue space} $\bH_{A,\varphi,\psi,r,\nu}^1$ is a finite-dimensional complex linear subspace of $\sH_1$ equipped with the almost complex structure $\bJ_1$ induced by the almost complex structure $J$ on $X$. For large enough $\nu$, the $L^2$ orthogonal projection
\[
  \Pi_{1,A,\varphi,\psi,r,\nu}:\sH_1 \to \bH_{A,\varphi,\psi,r,\nu}^1 
\]
restricts to an embedding of real vector spaces,
\[
  \Pi_{1,A,\varphi,\psi,r,\nu}:\bH_{A,\varphi,\psi,r}^1 \to \bH_{A,\varphi,\psi,r,\nu}^1, 
\]
and an isomorphism of real vector spaces,
\[
  \bH_{A,\varphi,\psi,r}^1
  \cong
  \mathbf{\tilde H}_{A,\varphi,\psi,r}^1
  :=
  \Pi_{1,A,\varphi,\psi,r,\nu}\left(\bH_{A,\varphi,\psi,r}^1\right).
\]
Recall that $\bH_{A,\varphi,\psi,r}^1$ represents the Zariski tangent space to $\sM^0(E,g,J,\omega,r)$ at $[A,\varphi,\psi]$. We now define
\[
  \mathbf{\tilde H}_{A,\varphi,\psi,r,\nu}^1
  :=
  \bH_{A,\varphi,\psi,r}^1 \oplus_\RR  \bH_{A,\varphi,\psi,r,\nu}^1 \cap \left(\mathbf{\tilde H}_{A,\varphi,\psi,r}^1\right)^\perp
  \subset
  \sH_1,
\]
where $\oplus_\RR$ denotes a real linear direct sum that is not necessarily orthogonal. For large enough $\nu$, we obtain an isomorphism of real vector spaces,
\[
  \mathbf{\tilde H}_{A,\varphi,\psi,r,\nu}^1 \cong \bH_{A,\varphi,\psi,r,\nu}^1,
\]
and hence the almost complex structure $\bJ_1$ on $\bH_{A,\varphi,\psi,r,\nu}^1$ induces an almost complex structure on
$\mathbf{\tilde H}_{A,\varphi,\psi,r,\nu}^1$.

If $\bH_{A,\varphi,\psi,r,\nu}^2$ is the direct sum of the eigenvectors of the complex linear Laplacian $\sT_{A,\varphi,\psi,r}'\sT_{A,\varphi,\psi,r}^{\prime,*}$ with eigenvalues less than $\nu$, then  $\bH_{A,\varphi,\psi,r,\nu}^2$ is a complex subspace of the complex Hilbert space
\[
  \sH_2 = L^2\left(\su(E) \oplus \su(E) \oplus \Lambda^{0,2}(\fsl(E)) \oplus \Lambda^{0,1}(E)\right),
\]
equipped with the almost complex structure $\bJ_2$ induced by the almost complex structure $J$ on $X$. For large enough $\nu$, there is an embedding of real vector spaces,
\[
  \Pi_{2,A,\varphi,\psi,r,\nu}:\bH_{A,\varphi,\psi,r}^2 \to \bH_{A,\varphi,\psi,r,\nu}^2,
\]
and for
\begin{align*}
  \mathbf{\tilde H}_{A,\varphi,\psi,r}^2
  &:=
  \Pi_{2,A,\varphi,\psi,r,\nu}\left(\bH_{A,\varphi,\psi,r}^2\right),
  \\
  \mathbf{\tilde H}_{A,\varphi,\psi,r,\nu}^2
  &:=
  \bH_{A,\varphi,\psi,r}^2 \oplus_\RR  \bH_{A,\varphi,\psi,r,\nu}^2 \cap \left(\mathbf{\tilde H}_{A,\varphi,\psi,r}^2\right)^\perp
  \subset
  \sH_2,
\end{align*}
we obtain isomorphisms of real vector spaces,
\begin{align*}
  \mathbf{\tilde H}_{A,\varphi,\psi,r}^2 &\cong \bH_{A,\varphi,\psi,r}^2,
  \\
  \mathbf{\tilde H}_{A,\varphi,\psi,r,\nu}^2 &\cong \bH_{A,\varphi,\psi,r,\nu}^2,
\end{align*}
and hence an almost complex structure on $\mathbf{\tilde H}_{A,\varphi,\psi,r,\nu}^2$. Here, $\bH_{A,\varphi,\psi,r}^2$ represents the obstruction space to $\sM^0(E,g,J,\omega,r)$ at $[A,\varphi,\psi]$. These observations, proved in detail in Sections \ref{subsec:Local_Kuranishi_model_nonlinear_map_Banach_spaces} and \ref{sec:Approximation_orthogonal_projections_finite-dim_subspaces_Hilbert_spaces}, lead to our proofs of Theorem \ref{mainthm:AH_structure_bounded_evalue_spaces_non-Abelian_monopoles_symp_4-mflds} and Corollary \ref{maincor:Almost_Hermitian_structure_moduli_space_non-Abelian_monopoles_symplectic_4-manifolds} in Section \ref{sec:Hitchin_function_Hamiltonian_circle_action_virtual_moduli_spaces}. The approximation technique that we develop in Sections \ref{subsec:Local_Kuranishi_model_nonlinear_map_Banach_spaces} and \ref{sec:Approximation_orthogonal_projections_finite-dim_subspaces_Hilbert_spaces} circumvents the fact that the Nijenhuis tensor $N_J$ need not be $C^0$-small and avoids difficult questions regarding the existence of spectral gaps for $T_{A,\varphi,\psi,r}$ that are uniform with respect to $[A,\varphi,\psi]$ and $r$. However, our approximation technique requires us to employ the high-dimensional spaces $\mathbf{\tilde H}_{A,\varphi,\psi,r,\nu}^1$ and $\mathbf{\tilde H}_{A,\varphi,\psi,r,\nu}^2$ in our definition of local Kuranishi models for open neighborhoods of points $[A,\varphi,\psi]$ in $\sM^0(E,g,J,\omega,r)$. See Section \ref{sec:Further_work} for further discussion of this point.

Section \ref{subsec:Virtual_Morse-Bott_signatures_SW_critical_points} contains the statements of Theorem \ref{mainthm:MorseIndexAtReduciblesOnAlmostKahler} and \ref{maincor:MorseIndexAtReduciblesOnSymplecticWithSO3MonopoleCharacteristicClasses}, which give the results of our calculations of the virtual Morse--Bott signature of the Hitchin function $f$ at a point represented by a Seiberg--Witten monopole. The crucial positivity of the virtual Morse--Bott index of the Hitchin function at a point represented by a Seiberg--Witten monopole is established by Corollary 
\ref{maincor:Positivity_of_MorseIndexAtReduciblesOnSymplecticWithSO3MonopoleCharacteristicClasses}.

\subsection{Feasibility of spin${}^u$ structures for non-Abelian monopole moduli spaces}
\label{subsec:Feasibility_non-Abelian_monopole_cobordism_method}
As described in the beginning of Chapter \ref{chap:Introduction} and the
forthcoming Section \ref{subsec:Critical_points_Hamiltonian_moduli_spaces_non-Abelian_monopoles_closed_symplectic_4-manifolds},
we aim to prove Conjecture \ref{mainconj:Existence_ASD_connection_BlowUp} by studying the gradient flow lines of
the Hitchin function $f$ in \eqref{eq:Hitchin_function} on the moduli space $\sM(E,g,J,\omega,r)$ of non-Abelian monopoles with a regularized Taubes perturbation defined in \eqref{eq:Moduli_space_non-Abelian_monopoles_almost_Hermitian_Taubes_regularized} on the smooth blow-up $\widetilde X$ of $X$.  To connect with our previous work \cite{Feehan_Leness_introduction_virtual_morse_theory_so3_monopoles}, we use the language of spin${}^u$ structures given here in Section \ref{sec:SpinuDefinitions}.  For a rank-two complex vector bundle $E$, we will consider the spin${}^u$ structure $\ft=(\rho,E\oplus \Lambda^{0,2}\otimes E)$, where $\rho$ is a Clifford multiplication map as defined in \eqref{eq:CliffordMapDefn},
and write $\sM_\ft=\sM(E,g,J,\omega,r)$.
The following definition records conditions on the spin${}^u$ structure $\ft$ that are required by our program to prove that $M_\ka^w(\widetilde X,g)$ has non-negative expected dimension.

\begin{maindefn}[Feasibility of \spinu structures]
\label{maindefn:Feasibility}
Let $X$ be a closed, connected, oriented, smooth Riemannian four-manifold with $b^+(X)>1$ and odd
$b^+(X)-b_1(X)$. If $g$ is a Riemannian metric on the smooth blow-up $\widetilde X = X\#\overline{\CC\PP}^2$ of $X$ and $J$ is a $g$-orthogonal almost complex structure on $\widetilde X$ that is compatible with the orientation of $\widetilde X$ determined by the orientation of $X$, we say that a spin${}^u$ structure $\ft=(\rho,W\otimes E)$ over $(\widetilde X,g)$ is \emph{feasible with respect to $c_1(T\widetilde X,J)$} if it has the following properties.
\begin{enumerate}
\item
\label{item:Feasibility_NonZeroSW}
There is a spin${}^c$ structure $\fs$ on $(\widetilde X,g)$ with $\SW_{\widetilde X}(\fs)\neq 0$ such that the image of $M_\fs$ under the continuous embedding \eqref{eq:DefnOfIotaOnQuotient} is contained in $\sM_\ft$.

\item
\label{item:Feasibility_PosExpDim}
$\sM_\ft$ has positive expected dimension.

\item
\label{item:Feasbility_NoZeroSectionSplit}
$\sM_\ft$ contains no zero-section point $[A,0]$, where $A$ is split as in Definition \ref{defn:Split_trivial_central-stabilizer_spinor_pair}, Item \eqref{item:Split_spinor_pair}.

\item
\label{item:Feasibility_LowerBound}
The Pontrjagin number $p_1(\su(E))[\widetilde X]$ satisfies the fundamental bounds \eqref{eq:p1_lower_bound_blowup}.

\item
  \label{item:Feasibility_PositivevBM}
  The \emph{formal Morse--Bott index},
\begin{equation}
  \label{eq:FormalMorseIndexIntroThm}
\lambda^-(\ft,\fs')
:=
-\frac{1}{6}\left( c_1(\widetilde X)^2+c_2(\widetilde X)\right)
-\left( c_1(\fs')-c_1(\ft)\right)\cdot c_1(\widetilde X)
-\left( c_1(\fs')-c_1(\ft)\right)^2,
\end{equation}
where $c_1(\widetilde X) = c_1(T\widetilde X,J)$, is positive for all spin${}^c$ structures $\fs'$ on $(\widetilde X,g)$ such that the image of $M_{\fs'}$ under the continuous embedding \eqref{eq:DefnOfIotaOnQuotient} is contained in $\sM_\ft$.
\end{enumerate}
\end{maindefn}

The moduli space $\sM_\ft$ is a real analytic space (see Feehan \cite{Feehan_analytic_spaces}) and Condition \eqref{item:Feasibility_PosExpDim} in Definition \ref{maindefn:Feasibility} ensures that the smooth top stratum of $\sM_\ft$ has positive dimension. As we will show in Corollary \ref{maincor:FeasibleIsNonEmpty}, Conditions \eqref{item:Feasibility_NonZeroSW}, \eqref{item:Feasibility_PosExpDim}, and \eqref{item:Feasbility_NoZeroSectionSplit} ensure non-emptiness of the moduli subspace $\sM_\ft^{*,0} \subset \sM_\ft$ of gauge-equivalence classes of non-split non-zero-section non-Abelian monopoles \eqref{eq:PU2MonopoleSubspace_irreducible_non-zero-section}.
%PF9-25-2024  When X Kaehler, vMB > 0 at some SW point is enough to give pos dim top stratum of \sM_\ft^{*,0}, even if exp dim is \leq 0.
%TL8-30-2025: I put a remark below but it still needs an explanation of how the Kaehler condition gives this result (i.e. I think you've mentioned an algebraic geometry notion of a stable set and I think that's what you mean here.)
%PF9-30-2024 Add remark that if SW(\fs)\neq 0 and vMB > 0 at point in M_\fs \subset \sM_\ft, then \sM_\ft^{*,0} non-empty even if exp dim is \leq 0.
%TL8-26-2025: This is going to take some thought but I'll see if I can get something on it.   The non-emptiness of \sM_\ft^{*,0} following from its expected dim being positive does require some input on the non-trivial cohomology of $\bL_{\ft,\fs}$ --I'd have to think if it's clear that argument works for the subspace given by the negative weight subspace or if the Kahler condition and the algebraic-geometry theory handles that.
By \cite[Corollary 3.3, p. 88]{FL2a}, Condition \eqref{item:Feasbility_NoZeroSectionSplit} holds if $w_2(\su(E))$ satisfies the Morgan--Mrowka condition \cite[Section 2, Paragraph prior to Corollary 2.2, p. 226]{MorganMrowkaPoly},
\begin{equation}
\label{eq:MorganMrowkaCondition}
\langle w_2(\su(E)),e^*\rangle \neq 0,
\end{equation}
where $e^*\in H_2(\widetilde X;\ZZ)$ is represented by an immersed two-sphere. Condition \eqref{item:Feasibility_LowerBound} implies that Conjecture \ref{mainconj:Existence_ASD_connection_BlowUp} will follow from the non-emptiness of the subspace $M_\ka^w(\widetilde X,g) \subset \sM_\ft$ by the discussion in Section \ref{sec:Existence_ASD_connections_and_BMY_inequality}.

The \emph{formal} Morse--Bott index $\lambda^-(\ft,\fs)$ in \eqref{eq:FormalMorseIndexIntroThm} is defined in terms of the characteristic classes $c_1(\widetilde X,J)$ and $c_2(\widetilde X)$ and $c_1(\fs')$ and $c_1(\ft)$ and their cup products. Corollary \ref{maincor:MorseIndexAtReduciblesOnSymplecticWithSO3MonopoleCharacteristicClasses} asserts that the \emph{virtual} Morse--Bott index $\lambda_{[A,\Phi]}^-(f)$ of the Hitchin function $f$ in \eqref{eq:Hitchin_function} at any point $[A,\varphi,\psi] \in \sM_\ft$ in the image of the embedding \eqref{eq:DefnOfIotaOnQuotient} of $M_{\fs'}$ is equal to the topological expression $\lambda^-(\ft,\fs)$ in \eqref{eq:FormalMorseIndexIntroThm}. We prove Corollary \ref{maincor:MorseIndexAtReduciblesOnSymplecticWithSO3MonopoleCharacteristicClasses} when $(X,g,J,\omega)$ is an almost K\"ahler four-manifold, but that proof is lengthy and requires us to first prove that the approximate bounded eigenvalue spaces $\mathbf{\tilde H}_{A,\varphi,\psi,r,\nu}^1$ and $\mathbf{\tilde H}_{A,\varphi,\psi,r,\nu}^2$,
defined in \eqref{eq:tilde_bH_A_varphi_psi_r_nu_1} and \eqref{eq:tilde_bH_A_varphi_psi_r_nu_2}, respectively,
admit $S^1$-invariant almost complex structures for a perturbation parameter $r \in (0,\infty)$ in \eqref{eq:SO(3)_monopole_equations_(1,1)_curvature_perturbed_intro_regular} by Theorem \ref{mainthm:AH_structure_bounded_evalue_spaces_non-Abelian_monopoles_symp_4-mflds}.

Condition \eqref{item:Feasibility_PositivevBM}
and the equality (to be established in the forthcoming Corollary \ref{maincor:MorseIndexAtReduciblesOnSymplecticWithSO3MonopoleCharacteristicClasses}) between the formal Morse--Bott index \eqref{eq:FormalMorseIndexIntroThm} and the virtual Morse--Bott index
defined in \eqref{eq:vBM_Perturbed_Intro} ensure that points in $M_\ka^w(\widetilde X,g)$ are the only local minima of the Hitchin function \eqref{eq:Hitchin_function} on $\sM_\ft$.

%PF8-31-2025 What is this about? Why needed?
\begin{rmk}[Non-empty moduli spaces on K\"aehler manifolds]
\label{rmk:NonEmptyModuliSpaceCondition}
%TL8-30-2025: Reference needed for how Kaehler condition gives this: is it an algebraic geometry fact on the stable set?
If $X$ is a K\"aehler manifold, then Condition \eqref{item:Feasibility_PosExpDim} in Definition \ref{maindefn:Feasibility},
that the moduli space $\sM_\ft$ has positive expected dimension, can be replaced with the assumption that the virtual Morse--Bott index, $\lambda_{[A,\Phi]}^-(f)$, is positive for $[A,\Phi]\in M_\fs$.
\end{rmk}

\begin{rmk}[Feasibility and energy bubbling]
\label{rmk:FeasibilityAndEnergyBubbling}
If it were not for the problem of energy bubbling described in the beginning of Chapter \ref{chap:Introduction} and in the forthcoming Section \ref{subsec:Further_work_allowing_energy_bubbling_gluing_non-Abelian_monopoles},
the limit of the downward gradient flow in $\sM_\ft$, when $\ft$ is feasible in the sense of Definition \ref{maindefn:Feasibility}, would be the desired element of $M_\ka^w(\widetilde X,g)$. Because of energy bubbling, such a limit could lie in a lower level of the Uhlenbeck compactification $\bar\sM_\ft$, given by the closure of $\sM_\ft$ in the space of ideal non-Abelian monopoles defined in Feehan and Leness \cite[Definition 4.19, p. 350]{FL1}. Therefore, to carry out this program we would need to also show that ideal Seiberg--Witten critical points in the Uhlenbeck compactification are not local minima. We say that a spin${}^u$ structure $\ft$ is \emph{fully feasible} if it is feasible and $f$ has no non-zero-section local minima on $\bar\sM_\ft$. We plan to prove the existence of fully feasible spin${}^u$ structures in future work.
\qed\end{rmk}

\begin{rmk}[Choice of an almost complex structure]
\label{rmk:Fasibility_AC_Structure}
The definition of feasibility depends only on the choice of an almost complex structure $J$ on $T\widetilde X$ through the characteristic class $c_1(T\widetilde X,J)$, which we will usually take to be determined by a basic class on $X$ as in the forthcoming Equation \eqref{eq:AC_Structure_For_Feasibility}.
\qed\end{rmk}

We say that $Q_X$ is \emph{indefinite} if $b^+(X)>0$ and $b^-(X)>0$ (see Gompf and Stipsicz \cite[Definition 1.2.8, p. 10]{GompfStipsicz}). Theorem \ref{mainthm:ExistenceOfSpinuForFlow} below relaxes the hypothesis of our previous result 
%TL11-26-2025: Update to \cite[Section 2.3, Theorem 7]{Feehan_Leness_introduction_virtual_morse_theory_so3_monopoles}
\cite[Theorem 3]{Feehan_Leness_introduction_virtual_morse_theory_so3_monopoles}, which assumed $b_1(X)=0$, odd $b^+(X)\ge 3$, and $b^-(X)\geq 2$.

\begin{mainthm}[Feasibility of \spinu structures]
\label{mainthm:ExistenceOfSpinuForFlow}
Let $X$ be a closed, connected, oriented, smooth Riemannian four-manifold with indefinite intersection form,
$b^+(X)>1$, odd $b^+(X)-b_1(X)$, and $c_2(X)>0$. If $\fs_0$ is a spin${}^c$ structure over $X$ such that $c_1(\fs_0)$ is a basic class in the sense of Definition \ref{defn:Seiberg-Witten_basic_class}, then for any Riemannian metric $g$ on the smooth blow-up $\widetilde X=X\#\overline{\CC\PP}^2$, there exist
\begin{inparaenum}[\itshape i\upshape)]
\item a $g$-orthogonal almost complex structure $J$ on $\widetilde X$ with
\begin{equation}
\label{eq:AC_Structure_For_Feasibility}
    c_1(\widetilde X):=c_1(T\widetilde X, J)=-c_1(\fs_0)-e,
\end{equation}
where $e$ is the Poincar\'e dual of the exceptional two-sphere\footnote{See Gompf and Stipsicz \cite[Definition 2.2.7, p. 43]{GompfStipsicz}.} in $\widetilde X$, and
\item  a spin${}^u$ structure $\ft=(\rho,W\otimes E)$ over $(\widetilde X, g)$
\end{inparaenum}
such the following hold:
\begin{enumerate}
\item The spin${}^u$ structure $\ft$ is feasible with respect to the class $c_1(\widetilde X)$ in 
\eqref{eq:AC_Structure_For_Feasibility} in the sense of Definition \ref{maindefn:Feasibility}.
\item\label{item:Upper_bound_expected_dimension_ASD_moduli_space_blowup}
The expected dimension of the moduli space $M_\kappa^w(\widetilde X, g)$ of projectively anti-self-dual connections on $E$, where $\kappa=\kappa(E)$ as in \eqref{eq:DefineKappa} and $w=c_1(E)$, obeys the following inequality:
\begin{equation}
  \label{eq:Upper_bound_expected_dimension_ASD_moduli_space_blowup}
    2\,\expdim M_\kappa^w(\widetilde X, g) \leq 3c_2( X) - c_1( X)^2.
  \end{equation}
\end{enumerate}
\end{mainthm}

We also have the 

\begin{maincor}[Non-emptiness of the moduli subspace of non-split, non-zero-section non-Abelian monopoles]
\label{maincor:FeasibleIsNonEmpty}
Continue the hypotheses of Theorem \ref{mainthm:ExistenceOfSpinuForFlow}. Then the moduli subspace $\sM_{\ft}^{*,0} \subset\sM_{\ft}$ is non-empty, where $\sM_{\ft}^{*,0}$ is as in
\eqref{eq:PU2MonopoleSubspace_irreducible_non-zero-section}.
\end{maincor}

We prove Theorem \ref{mainthm:ExistenceOfSpinuForFlow} and Corollary \ref{maincor:FeasibleIsNonEmpty} in Chapter \ref{chap:Feasibility}.

\begin{rmk}[The case $b^+(X)=1$]
\label{rmk:bPlusOne_and_Feasibility}
The Seiberg--Witten invariants depend on the \emph{chamber} containing the Riemannian metric $g$ and a $g$-self-dual two-form $\eta$ on $X$ when $b^+(X)=1$. Here, $\eta$ is a perturbation term in the Seiberg--Witten equations (see Morgan \cite[Section 6.9]{MorganSWNotes} or Salamon \cite[Section 7.4, pp. 254--255, and Section 9.2]{SalamonSWBook}). By chamber, we mean the  connected components  of the complement of a codimension-one hypersurface in the space of parameters $(g,\eta)$. There are two Seiberg--Witten invariants, $\SW_X^+(\fs)$ and $\SW_X^-(X)$, depending on the chamber in which the parameters $(g,\eta)$ lie. To assume as in Theorem \ref{mainthm:ExistenceOfSpinuForFlow}, that a Seiberg--Witten invariant is non-zero when $b^+(X)=1$, one must specify the chamber in which the parameters $(g,\eta)$ lie. When  the continuous embedding \eqref{eq:DefnOfIotaOnQuotient} maps $M_\fs$ into $\sM_\ft$, the parameters $(g,\eta)$ are determined by parameters used to define the non-Abelian monopole equations, as described in \cite[Lemma 3.17, p. 95]{FL2a}.
% TL9-25-2024:  Probably more detail needed.
%PF9-30-2024 Such as what?
%TL10-2-2024:  Since we're not dealing w/ b^+=1 here, I think we only need to find a reference to what Taubes says about the chamber for his perturbation.
%TL10-10-2024: I didn't see anything in the Taubes papers about b^+=1 but Salamon gives him credit for the result \cite[Theorem 13.9, p. 412]{SalamonSWBook} which says which invariants are non-zero when $b^+=1$.
In addition, the choice of perturbation parameter $r \in [0,\infty)$ appearing in Equation \eqref{eq:SO(3)_monopole_equations_(1,1)_curvature_perturbed_intro} will affect the chamber. Hence, while there is a natural extension of Theorem \ref{mainthm:ExistenceOfSpinuForFlow} to four-manifolds with $b^+(X)=1$, more work is required to determine whether $\SW_X^+(\fs)\neq 0$ or $\SW_X^-(\fs)\neq 0$ is the proper replacement
for $\SW_X(\fs)\neq 0$.
%PF9-28-2024 We also would need c_1(\fs)^2 = c_1(X)^2 to combute the vMB?
%TL9-28-2024: We would only need it for one $\fs$ but not to get the positivity of vBM for all $\tilde\fs$ with $M_{\tilde\fs}$ in $\sM_\ft$.
%We further note that when $b^+(X)=1$, the condition $c_1(\fs)^2=c_1(X)^2$ in Definition %\ref{defn:Seiberg-Witten_basic_class} need not hold
%does not hold

% PF9-29-2024 Moved here but duplicates above
%If $b^+(X) = 1$, there are two possible values of the Seiberg--Witten invariant, depending on the chamber determined by the generic geometric perturbations --- see Salamon \cite[Section 7.4, pp. 254--255 and Section 9.2]{SalamonSWBook} or Morgan \cite[Section 6.9]{MorganSWNotes}.

%PF9-26-2024 Odd phrasing
An additional obstruction to proving
Theorem \ref{mainthm:ExistenceOfSpinuForFlow} for four-manifolds with $b^+(X)=1$ appears in the hypotheses of Proposition \ref{prop:FeasibilityOnBlowUp}, where we need to assume that if $X$ is odd then $b^+(X)\ge 2$ or $b^-(X)\ge 2$ to produce a cohomology class $w\in H^2(\widetilde X;\ZZ)$  satisfying certain cup-product identities. Thus, when $b^+(X)=1$ the hypotheses of Proposition \ref{prop:FeasibilityOnBlowUp} would require $b^-(X)\ge 2$.  At present, we can only produce the desired class $w$ in the case that $b^+(X)\ge 2$ or $b^-(X)\ge 2$ if $X$ is odd but it is possible that further work will enable us to relax this condition.
\qed\end{rmk}

% PF9-23-2024 There is an alternative defn of feasibility along the following lines: (1) \sM_\ft has positive expected dimension and contains M_\fs for some \fs such that SW_X(\fs) is non-zero. (2) p_1 obeys fundamental bounds. (3) vMB index > 0 at all SW points in \sM_\ft. It would be worth comparing/contrasting these defns in a remark.

% PF9-25-2024 Add remark saying that vMB > 0 is enough, when X Kaehler, to show that the top smooth stratum of the complex analytic space \sM_\ft^{*,0} is positive-dimensional even when \sM_\ft has non-positive expected dimension.
% TL9-28-2024: Added following.  Please edit to get the terminology consistent with what you are using.
% PF9-29-2024 The remark, as written, doesn't make sense.
% \begin{rmk}[Positivity of expected dimension when $X$ is K\"ahler]
% \label{rmk:PositivityOfExpDim}
% When $X$ is K\"ahler and $M_\fs$ is non-empty, the positivity of the virtual Morse--Bott index
% in \eqref{eq:FormalMorseIndexIntroThm}
% suffices to prove that the unstable subvariety of the Hitchin function \eqref{eq:Hitchin_function} in $\sM_\ft$ is non-empty near $M_\fs$.
% Thus, in this case
% Item \eqref{item:ExpectedDimOfFeasibleSO3ModuliSpace} of Theorem \ref{mainthm:ExistenceOfSpinuForFlow} is not needed
% for this purpose.
% \end{rmk}

\subsection{Generalizations of Donaldson's symplectic subspace criterion}
\label{subsec:Generalization_Donaldson_symplectic_subspace_criterion}
In this subsection, we shall prove several generalizations of Donaldson's well-known symplectic subspace criterion for the kernel of a real linear operator \cite[Section 1, Proposition 3, p. 669]{DonSympAlmostCx}. Like Donaldson's original symplectic subspace criterion, some (though not all) of our generalizations require the complex antilinear component of the real linear operator to be suitably small, a requirement that unfortunately appears difficult to verify in our application. While we shall ultimately circumvent the need to apply any these generalizations by instead applying the approximation method developed in Sections \ref{subsec:Local_Kuranishi_model_nonlinear_map_Banach_spaces} and \ref{sec:Approximation_orthogonal_projections_finite-dim_subspaces_Hilbert_spaces}, the proofs of our generalizations of Donaldson's symplectic subspace criterion provide the motivation and scaffolding needed to develop our alternative approximation method. We begin by recalling the following celebrated result due to Donaldson.

\begin{thm}[Existence of codimension two symplectic submanifolds]
\label{thm:Donaldson_1}  
(See Donaldson \cite[Section 1, Theorem 1, p. 666]{DonSympAlmostCx}.)
Let $(V,\omega)$ be a closed symplectic manifold of dimension $2n$, and suppose that the de Rham cohomology class $[\omega/2\pi] \in H^2(V;\RR)$ lies in the integral lattice $H^2(V;\ZZ)/\mathrm{Torsion}$. Let $h \in H^2(V;\ZZ)$ be a lift of $[\omega/2\pi]$ to an integral class. If $k$ is a sufficiently large integer, then the Poincar\'e dual of $kh$ in $H_{2n-2}(V;\ZZ)$ can be realized by a symplectic submanifold $\Sigma \subset V$.
\end{thm}

Donaldson proved Theorem \ref{thm:Donaldson_1} by constructing the desired symplectic submanifold as the zero locus of a smooth section $s$ of a Hermitian line bundle $L$ with an $(0,1)$-connection over an almost complex manifold (see Donaldson \cite[Section 4(b), Equation (15), p. 34]{DonYangMillsInvar} or \cite[Section 1, Equation (4), p. 670]{DonSympAlmostCx}) that is \emph{approximately holomorphic} in the sense that
\begin{equation}
  \label{eq:Donaldson_1996jdg_4}
  |\bar\partial s|_{\Lambda^{0,1}(L)} < |\partial s|_{\Lambda^{1,0}(L)} \quad\text{on } s^{-1}(0),
\end{equation}
where the derivative $\nabla s$ is well-defined on $s^{-1}(0)$ and $\nabla s = \partial s  + \bar\partial s$ on $s^{-1}(0)$. $\Omega^0(L)$. In fact, he proves

\begin{thm}[Existence of approximately holomorphic sections]
\label{thm:Donaldson_5}  
(See Donaldson \cite[Theorem 5, p. 670]{DonSympAlmostCx}.)
Let $L$ be a Hermitian line bundle over a compact symplectic manifold $(V,\omega)$ with compatible almost complex structure and $c_1(L) = [\omega/2\pi]$. Then there is a constant $C>0$ such that for all large integers $k$, there is a smooth section $s$ of $L^{\otimes k}$ with
\begin{equation}
  \label{eq:Donaldson_1996jdg_theorem_5}
  |\bar\partial s|_{\Lambda^{0,1}(L^{\otimes k})} < \frac{C}{\sqrt{k}}|\partial s|_{\Lambda^{1,0}(L^{\otimes k})} \quad\text{on } s^{-1}(0).
\end{equation}
\end{thm}

A straightforward approximation argument \cite[Section 1, second paragraph, p. 671]{DonSympAlmostCx} allows Donaldson to prove a version of Theorem \ref{thm:Donaldson_1} without the integrality hypothesis, while a simple induction argument \cite[Section 1, second paragraph, p. 671]{DonSympAlmostCx} allows him to construct symplectic submanifolds of any even codimension, yielding the

\begin{cor}[Existence of symplectic submanifolds of any even codimension]
\label{cor:Donaldson_6}
(See Donaldson \cite[Section 1, Corollary 6, p. 671]{DonSympAlmostCx}.)  
If $(V,\omega)$ is a compact symplectic manifold, then the following hold:
\begin{enumerate}
\item $V$ contains symplectic submanifolds of any even codimension.
\item If $J$ is a compatible almost complex structure on $V$, then there are almost complex structures $J'$ arbitrarily $C^0$-close to $J$ such that $V$ contains $J'$-pseudoholomorphic curves.
\end{enumerate}
\end{cor}

Auroux \cite[Section 1, Corollary 1, p. 973]{AurouxFamilySymp} extended Donaldson's results to produce symplectic submanifolds of arbitrary codimension as the zero loci of approximately holomorphic sections of higher rank Hermitian vector bundles. We take Equations \eqref{eq:Donaldson_1996jdg_4} and \eqref{eq:Donaldson_1996jdg_theorem_5} as points of departure for a definition of an approximately holomorphic map, $F:\CC^n \supset U \to \CC^m$, where $U$ is an open neighborhood of the origin, that is suitable for our later application to approximately holomorphic Kuranishi obstruction maps defining open neighborhoods of points in the moduli space of non-Abelian monopoles. Donaldson's proof of \eqref{eq:Donaldson_1996jdg_4} relies on the following simpler version of our forthcoming Proposition \ref{mainprop:Donaldson_1996jdg_3}:

\begin{prop}[Donaldson's symplectic subspace criterion]
\label{prop:Donaldson_1996jdg_3}  
(See Donaldson \cite[Section 1, Proposition 3, p. 669]{DonSympAlmostCx}.)
Let $n\geq 1$ be an integer, $\sL':\CC^n\to \CC$ be a complex linear map, and $\sL'':\CC^n\to \CC$ be a complex antilinear linear map. If 
\begin{equation}
  \label{eq:Donaldson_1996jdg_proposition_3_approx_complex_linear}
  \|\sL''\|_{\Hom(\CC^n,\CC)} < \|\sL'\|_{\Hom(\CC^n,\CC)},
\end{equation}
then $\sL = \sL'+\sL'':\CC^n\to \CC$ is a surjective real linear map and the standard symplectic form on $\CC^n$ restricts to a symplectic form on $\Ker \sL \subset \CC^n$.
\end{prop}

In our refinement, Proposition \ref{mainprop:Donaldson_1996jdg_3}, of Donaldson's Proposition \ref{prop:Donaldson_1996jdg_3}, we replace $\CC^n$ by a Hilbert space and replace Donaldson's condition \eqref{eq:Donaldson_1996jdg_proposition_3_approx_complex_linear} by the slightly weaker condition \ref{eq:Donaldson_1996jdg_proposition_3_approx_complex_linear_inequality}. See the forthcoming Remarks \ref{rmk:Proof_prop_Donaldson_1996jdg_3} and \ref{rmk:Symplectic_form_on kernel_generic_T} for further discussion.

Before stating our first generalization of Proposition \ref{prop:Donaldson_1996jdg_3}, we recall some basic complex linear algebra. If $(\sH^\RR,g,J)$ is a real Hilbert space with inner product $g = \langle\cdot,\cdot\rangle_\sH$ and compatible almost complex structure $J$ (see, for example, Huybrechts \cite[Definition 1.2.11, p. 28]{Huybrechts_2005}), we let $\omega := g(J\cdot,\cdot)$ denote the fundamental two-form on $\sH^\RR$ (see Huybrechts \cite[Definition 1.2.13, p. 29]{Huybrechts_2005}), let $\sH = (\sH^\RR,J)$ denote the complex vector space defined by $iv := Jv$ for all $v\in\sH^\RR$, and let $(\sH,h)$ denote the complex Hilbert space with Hermitian inner product $h := g - i\omega$ (see \cite[Lemma 1.2.15, p. 30]{Huybrechts_2005}). We write $\sH^{1,0}$ and $\sH^{0,1}$ for the $\pm i$ eigenspaces of $J\in \End_\CC(\sH)$ and recall that (see \cite[Lemma 1.2.5, p. 26]{Huybrechts_2005})
\begin{equation}
  \label{eq:sH_direct_sum_splitting_10_and_01_subspaces}
  \sH = \sH^{1,0} \oplus \sH^{0,1}
\end{equation}
is a direct sum of complex vector spaces, with idempotent maps defined by
\begin{equation}
  \label{eq:sH_pi_10_and_pi_01_idempotents}
  \pi^{1,0}v \equiv v' := \frac{1}{2}(v-iJv)
  \quad\text{and}\quad
  \pi^{0,1}v \equiv v'' := \frac{1}{2}(v+iJv),
  \quad\text{for all } v \in \sH.
\end{equation}
The splitting $\sH = \sH^{1,0} \oplus \sH^{0,1}$ is $h$-orthogonal \cite[Lemma 1.2.16, p. 30]{Huybrechts_2005}, so the idempotent maps
\begin{equation}
  \label{eq:sH_pi_10_and_pi_01_orthogonal_projections}
  \pi^{1,0}:\sH \to \sH^{1,0}
  \quad\text{and}\quad
  \pi^{0,1}:\sH \to \sH^{0,1}
\end{equation}
are $h$-orthogonal projections.

The almost complex structure $J$ on $\sH$ induces an almost complex structure $J$ on the real Hilbert space $\Hom_\RR^c(\sH,\RR)$ of continuous $\RR$-linear $\RR$-valued functions on $\RR$ by \cite[Lemma 1.2.6, p. 26]{Huybrechts_2005}
\[
  J\sL := \sL\circ J, \quad\text{for all } \sL \in \Hom_\RR^c(\sH,\RR).
\]
There is an induced splitting of the complex vector space
\[
  \Hom_\RR^c(\sH,\CC) = \Hom_\RR^c(\sH,\RR)\otimes_\RR\CC
\]
given by (see \cite[Lemma 1.2.5, p. 26]{Huybrechts_2005})
\begin{equation}
  \label{eq:Hom_RR_sH_to_CC}
  \Hom_\RR^c(\sH,\CC) = \Hom_\RR^c(\sH,\CC)^{1,0} \oplus \Hom_\RR^c(\sH,\CC)^{0,1},
\end{equation}
as a direct sum of complex vector spaces, with idempotent maps defined by
\begin{equation}
  \label{eq:sTprime_and_sTprimeprime_maps_sH_to_CC}
  \pi^{1,0}\sL \equiv \sL' = \frac{1}{2}(\sL-i\sL J)
  \quad\text{and}\quad
  \pi^{0,1}\sL \equiv \sL'' := \frac{1}{2}(\sL+i\sL J),
\end{equation}
so $\sL'$ is the complex linear component of $\sL$ and $\sL''$ is the complex antilinear component of $\sL$. The splitting \eqref{eq:Hom_RR_sH_to_CC} is orthogonal \cite[Lemma 1.2.16, p. 30]{Huybrechts_2005} with respect to the Hermitian inner product on $\Hom_\RR^c(\sH,\CC)$ induced by the real inner products on $\CC$ and the real Hilbert dual space $\Hom_\RR^c(\sH,\RR)$ and the almost complex structure on  $\Hom_\RR^c(\sH,\CC)$. Hence, the idempotent maps,
\[
  \pi^{1,0}:\Hom_\RR^c(\sH,\CC) \to \Hom_\RR^c(\sH,\CC)^{1,0}
  \quad\text{and}\quad
  \pi^{0,1}:\Hom_\RR^c(\sH,\CC) \to \Hom_\RR^c(\sH,\CC)^{0,1},
\]
are orthogonal projections with respect to the Hermitian inner product on $\Hom_\RR^c(\sH,\CC)$.

\begin{mainprop}[Inequality version of Donaldson's symplectic subspace criterion]
\label{mainprop:Donaldson_1996jdg_3}  
Let $\sH$ be a complex Hilbert space and $\sL:\sH\to \CC$ be a bounded \emph{real linear} map with complex linear component $\sL'$ and complex antilinear component $\sL''$ as in \eqref{eq:sTprime_and_sTprimeprime_maps_sH_to_CC}. If 
\begin{equation}
  \label{eq:Donaldson_1996jdg_proposition_3_approx_complex_linear_inequality}
  \|\sL''\|_{\Hom(\sH,\CC)} \neq \|\sL'\|_{\Hom(\sH,\CC)},
\end{equation}
then $\sL$ is surjective and the standard symplectic form on $\sH$ restricts to a symplectic form on $\Ker \sL \subset \sH$.
\end{mainprop}

We prove Proposition \ref{mainprop:Donaldson_1996jdg_3} (and hence Proposition \ref{prop:Donaldson_1996jdg_3} as a corollary) in Section \ref{sec:Proof_analogue_Donaldson_symplectic_submanifold_criterion_hypersurface}.

\begin{rmk}[Standard symplectic form on a complex Hilbert space]
\label{rmk:Standard_symplectic_form_complex_Hilbert_space}  
If $h$ denotes the Hermitian inner product on $\sH$ in Proposition \ref{mainprop:Donaldson_1996jdg_3} and $J$ is the almost complex structure on $\sH$ defined by scalar multiplication by $i=\sqrt{-1}$, then $h = g -i\omega$, where $g = \Real h$ is the inner product on the real Hilbert space $\sH_\RR$ underlying $\sH$ and $\omega = g(J\cdot,\cdot) = i\Imag h$ is the standard symplectic form on $\sH_\RR$. See Remark \ref{rmk:Conventions_almost_Hermitian_structures} for further explanation of our conventions.

Note also that since $J \in \End_\RR^c(\sH)$ is an isomorphism of real Hilbert spaces and $g$ is an inner product on $\sH_\RR$, then $\omega$ is a strongly non-degenerate symplectic form on $\sH_\RR$ (see Appendix \ref{sec:Weakly_strongly_non-degenerate_bilinear_forms_Banach_spaces} for a comparison of weakly and strongly non-degenerate symplectic forms).   
\qed\end{rmk}  

\begin{rmk}[On the hypothesis \eqref{eq:Donaldson_1996jdg_proposition_3_approx_complex_linear_inequality} and the proof of Proposition \ref{prop:Donaldson_1996jdg_3}]
\label{rmk:Proof_prop_Donaldson_1996jdg_3}
Cieliebak and Mohnke give a proof of Proposition \ref{prop:Donaldson_1996jdg_3} in \cite[Lemma 8.3 (b), p. 328 and Remark 8.4, p. 329]{Cieliebak_Mohnke_2007} by a geometric argument, involving the K\"ahler angle described by Donaldson in \cite[Section 1, pp. 668--669]{DonSympAlmostCx}. We give a proof of Proposition \ref{mainprop:Donaldson_1996jdg_3} (and thus Proposition \ref{prop:Donaldson_1996jdg_3} as a corollary) in Section \ref{sec:Proof_analogue_Donaldson_symplectic_submanifold_criterion} using linear algebra. (After writing our proof of Proposition \ref{mainprop:Donaldson_1996jdg_3}, we discovered a proof of Donaldson's Proposition \ref{prop:Donaldson_1996jdg_3} described by Krestiachine \cite[Lemma 1.1, p. 11]{KrestiachineThesis}, which he attributes to Patrick Massot and which also uses a linear algebra argument.)
\qed\end{rmk}

\begin{rmk}[Existence of symplectic forms on $\Ker\sL$ for generic $\sL$ in Proposition \ref{mainprop:Donaldson_1996jdg_3}]
\label{rmk:Symplectic_form_on kernel_generic_T}
Clearly, $\sL \in \Hom_\RR^c(\sH,\CC)$ obeys \eqref{eq:Donaldson_1996jdg_proposition_3_approx_complex_linear_inequality} if and only if it obeys \emph{either}
\[
  \|\sL''\|_{\Hom(\sH,\CC)} < \|\sL'\|_{\Hom(\sH,\CC)}
  \quad \text{\emph{or}}\quad
  \|\sL'\|_{\Hom(\sH,\CC)} < \|\sL''\|_{\Hom(\sH,\CC)},
\]
where the first inequality is Donaldson's hypothesis \eqref{eq:Donaldson_1996jdg_proposition_3_approx_complex_linear} when $\sH=\CC^n$. If $\sL''=0$, then $\sL=\sL'$ and $\Ker \sL$ is a complex vector subspace of $\sH$ and the complex orientation of $\Ker \sL$ is equal to the orientation induced by the complex orientations of the domain $\sH$, codomain $\CC$, and isomorphism $\sL:\sH/\Ker \sL \cong \CC$. This explains the motivation underlying the first inequality. On the other hand, if $\sL''=0$, then $\sL=\sL''$ and $\Ker \sL$ is also a complex vector subspace of $\sH$, but the complex orientation of $\Ker \sL$ is now opposite to the orientation induced by the isomorphism $\sL:\sH/\Ker \sL \cong \CC$. This explains the motivation underlying the second inequality. See Remark \ref{rmk:Donaldson_1996jdg_3_T_onto} for further discussion.

We emphasize that if $\sL \in \Hom_\RR^c(\sH,\CC)$ lies in the complement of the smooth real hypersurface,
\[
  \sS := \left\{\sL \in \Hom_\RR^c(\sH,\CC): \sL \neq 0 \text{ and } \|\sL''\|_{\Hom(\sH,\CC)} = \|\sL'\|_{\Hom(\sH,\CC)} \right\},
\]
in the real Hilbert space $\Hom_\RR^c(\sH,\CC)$ of bounded $\CC$-valued $\RR$-linear functions on $\sH$, then the standard symplectic form on $\sH$ restricts to a symplectic form on $\Ker \sL$ by Proposition \ref{mainprop:Donaldson_1996jdg_3}. To see that $\sS$ is indeed smooth, recall from \eqref{eq:Hom_RR_sH_to_CC} that the splitting
\[
  \Hom_\RR^c(\sH,\CC) = \Hom_\RR^c(\sH,\CC)^{1,0} \oplus \Hom_\RR^c(\sH,\CC)^{0,1}
\]
is a Hermitian orthogonal direct sum of complex Hilbert subspaces and consider the smooth map,
\[
  F:\Hom_\RR^c(\sH,\CC)\less\{0\} \ni (\sL',\sL'') \mapsto \|\sL'\|_{\Hom(\sH,\CC)}^2 - \|\sL''\|_{\Hom(\sH,\CC)}^2 \in \RR.
\]
The differential of $F$ at $\sL = \sL'+\sL'' \in \Hom_\RR^c(\sH,\CC)$ is given by 
\[
  (DF)_\sL A = 2\Real\,\langle A',\sL'\rangle - 2\Real\,\langle A'',\sL''\rangle,
  \quad\text{for all } A \in \Hom_\RR^c(\sH,\CC).
\]
If $(\sL',\sL'')\neq 0$ or, equivalently, $\sL  = \sL'+\sL''\neq 0$, then $\sL'\neq 0$ or $\sL''\neq 0$. Thus, if $\sL'\neq 0$, then we may choose $A = (\sL',0)$ to give
\[
  (DF)_\sL (\sL',0) = 2\|\sL'\|^2 \neq 0.
\]
Similarly, if $\sL''\neq 0$, then we may choose $A = (0,\sL'')$ to give $(DF)_\sL (0,\sL'') = 2\|\sL''\|^2 \neq 0$. Hence, $(DF)_\sL \neq 0$ for $\sL \in \Hom_\RR^c(\sH,\CC)\less\{0\}$ and so $\sS = F^{-1}(0)$ is a smooth real hypersurface, as claimed. (This result and its proof generalize the familiar fact that $\{(z_1,z_2) \in \CC^2\less\{0\}: |z_1| = |z_2|\}$ is a smooth real hypersurface in $\CC^2$.)
\qed\end{rmk}  

We obtain the following generalization of Proposition \ref{mainprop:Donaldson_1996jdg_3} by adapting and extending Donaldson's proof of his Corollary \ref{cor:Donaldson_6} from the setting of finite-dimensional symplectic manifolds. Recall that an unbounded linear operator $\sT:\sH_1\to\sH_2$ is \emph{closed} if its graph $\Graph(\sT) := \{(v,\sT v): v \in \Dom(\sT)\}$ is a closed subspace of $\sH_1\times\sH_2$ and \emph{closable} if there is closed operator $\bar\sT:\sH_1\to\sH_2$ that is an \emph{extension} of $\sT$, namely, such that $\Graph(\bar\sT) \supset \Graph(\sT)$ (see, for example, Reed and Simon \cite[Section VIII.1, Definitions, p. 250]{Reed_Simon_v1}). If $\Dom(\sT)$ is dense in $\sH_1$, then $\sT$ has a unique adjoint $\sT^*:\sH_2\to\sH_1$ by \cite[Section VIII.1, Definition, p. 252]{Reed_Simon_v1} and, furthermore, $\Dom(\sT^*)$ is dense in $\sH_2$ if and only if $\sT$ is closable by \cite[Section VIII.1, Theorem VIII.1, p. 252]{Reed_Simon_v1}.

\begin{mainthm}[Generalization of Donaldson's symplectic subspace criterion to kernels of generic bounded linear operators on Hilbert spaces]
\label{mainthm:Donaldson_1996jdg_3_Hilbert_space_codomain}  
Let $\sH_k^\RR$ be a real Hilbert space with inner product $g_k$ and $g_k$-orthogonal almost complex structure $J_k$ and symplectic form $\omega_k = g_k(J_k\cdot,\cdot)$, for $k=1,2$. Let $\sH_k = (\sH_k^\RR,J_k)$ be the corresponding complex Hilbert space with Hermitian inner product $h_k = g_k -i\omega_k$. Let $\sG_k \subset \sH_k$ be dense Hilbert subspaces for $k=1,2$ and assume that $\sG_2$ is separable and the almost complex structures $J_k$ on $\sH_k$ restrict to almost complex structures $J_k$ on $\sG_k$ for $k=1,2$. Let $\Hom_\RR(\sH_1,\sH_2)$ denote the complex vector space\footnote{In the sense of Kadison and Ringrose \cite[Section 2.7, first paragraph, p. 157]{KadisonRingrose1}.} of closable unbounded real linear operators from $\sH_1$ to $\sH_2$ with domain $\sG_1$ and whose real adjoints have domain $\sG_2$. Let $\Epim_\RR^c(\sG_1,\sH_2)$ denote the open subset of the Hilbert space\footnote{Where $\Hom_\RR^c(\sG_1,\sH_2)$ is viewed as the tensor product $\sG_1^*\otimes_\RR\sH_2$ of the Hilbert spaces $\sG_1^* \cong \sG_1$ and $\sH_2$ --- see Kadison and Ringrose \cite[Section 2.6, pp. 125--127]{KadisonRingrose1}.} $\Hom_\RR^c(\sG_1,\sH_2)$ of \emph{surjective} bounded real linear operators from $\sG_1$ to $\sH_2$. If $\dim_\CC\sH_2 = m \leq \infty$, then there is a set of smooth real hypersurfaces $\sS_n \subset \Epim_\RR^c(\sG_1,\sH_2)$ for $n=1,\ldots,m$, such that if $\sT$ belongs to 
\[
  \sU := \Epim_\RR^c(\sG_1,\sH_2) \ \less \bigcup_{n=1}^m\sS_n,
\]
then $\Ker\sT$ is a symplectic subspace of $(\sH_1^\RR,\omega_1)$.
\end{mainthm}

%PF11-21-2025 Check Lemma 6.3.1 and its proof for typos.

\begin{rmk}[Topology of the subset $\sU$ in Theorem \ref{mainthm:Donaldson_1996jdg_3_Hilbert_space_codomain}]
\label{rmk:Donaldson_1996jdg_3_Hilbert_space_codomain__topology_sU}
Note that $\sU$ in Theorem \ref{mainthm:Donaldson_1996jdg_3_Hilbert_space_codomain} is a dense open subset of $\Epim_\RR^c\left(\sG_1,\sH_2\right)$ if $m<\infty$ and a dense subset if $m=\infty$ by the Baire Category Theorem (see, for example, Munkres \cite[Chapter 8, Section 48, Theorem 48.2, p. 296]{Munkres_topology_second_edition}). Similarly, $\sV$ in Corollary \ref{maincorDonaldson_1996jdg_3_Hilbert_space_domain} is a dense open subset of $\Epim_\RR^c\left(\sG_2,\sH_1\right)$ if $p<\infty$ and a dense subset if $p=\infty$.
\end{rmk}

The next result is an immediate consequence of Theorem \ref{mainthm:Donaldson_1996jdg_3_Hilbert_space_codomain}.

\begin{maincor}[Generalization of Donaldson's symplectic subspace criterion to kernels of bounded linear operators on Hilbert spaces parametrized by a smooth manifold]
\label{maincor:Donaldson_1996jdg_3_Hilbert_space_codomain_families}
Continue the hypotheses and notation of Theorem \ref{mainthm:Donaldson_1996jdg_3_Hilbert_space_codomain}. If $M$ is a smooth manifold and 
\[
  F:M \ni x \mapsto \sT_x \in \Epim_\RR^c(\sG_1,\sH_2).
\]
is a smooth map that is transverse to the smooth real hypersurface $\sS_n \subset \Epim_\RR^c(\sG_1,\sH_2)$, then $F^{-1}(\sS_n) \subset M$ is a smooth real hypersurface in $M$. In particular, if $F \transv \sS_n$ for $n=1,\ldots,m$ and
\[
  U := M \ \less \bigcup_{n=1}^m F^{-1}(\sS_n),
\]
then for all $x\in U$, the kernel $\Ker\sT_x$ is a symplectic subspace of $(\sH_1^\RR,\omega_1)$. 
\end{maincor}

Just as in Remark \ref{rmk:Donaldson_1996jdg_3_Hilbert_space_codomain__topology_sU}, the set $U$ in Corollary \ref{maincor:Donaldson_1996jdg_3_Hilbert_space_codomain_families} is a dense open subset of $M$ if $m<\infty$ and a dense subset if $m=\infty$.

\begin{maincor}[Generalization of Donaldson's symplectic subspace criterion to kernels of adjoints of generic bounded linear operators on Hilbert spaces]
\label{maincorDonaldson_1996jdg_3_Hilbert_space_domain}  
Continue the hypotheses of Theorem \ref{mainthm:Donaldson_1996jdg_3_Hilbert_space_codomain} except that we  assume $\sG_1$ (rather than $\sG_2$) is separable. If $\dim_\CC\sH_1 = p \leq \infty$, then there is a set of smooth real hypersurfaces $\sR_l \subset \Epim_\RR^c(\sG_2,\sH_1)$, for $l=1,\ldots,p$, such that if $\sT^*$ belongs to 
\[
  \sV := \Epim_\RR^c\left(\sG_2,\sH_1\right) \ \less \bigcup_{l=1}^p\sR_l,
\]
where $\sT^* \in \Hom_\RR(\sH_2,\sH_1)$ is the adjoint of $\sT\in\Hom_\RR(\sH_1,\sH_2)$ defined by the real inner products,
%PF12-17-2025 In section 6.6, we use \sT^\intercal to denote real adjoint but maybe it's okay to use \sT^* for real adjoint of a real operator as it's used that way everywhere
\[
  \langle v_1,\sT^* v_2\rangle_{g_1} = \langle \sT v_1,v_2\rangle_{g_2},
  \quad\text{for all } v_1 \in \sH_1, \ v_2 \in \sH_2,
\]
then $\Ker\sT^*$ is a symplectic subspace of $(\sH_2^\RR,\omega_2)$.
\end{maincor}

We prove Theorem \ref{mainthm:Donaldson_1996jdg_3_Hilbert_space_codomain} and Corollary \ref{maincorDonaldson_1996jdg_3_Hilbert_space_domain} in Section \ref{sec:Proof_analogue_Donaldson_symplectic_submanifold_criterion_hypersurface_Hilbert_spaces}. As we noted in Section \ref{subsec:Guide_main_results}, we shall explore refinements of these results elsewhere and we do not make use of Theorem \ref{mainthm:Donaldson_1996jdg_3_Hilbert_space_codomain} or Corollaries \ref{maincor:Donaldson_1996jdg_3_Hilbert_space_codomain_families} and \ref{maincorDonaldson_1996jdg_3_Hilbert_space_domain} in this monograph.

We now state a different generalization and reformulation of Proposition \ref{prop:Donaldson_1996jdg_3} which does not rely on genericity of the operator $\sT$ as in Theorem \ref{mainthm:Donaldson_1996jdg_3_Hilbert_space_codomain}. See Section \ref{sec:Weakly_strongly_non-degenerate_bilinear_forms_Banach_spaces} for the definitions of weakly and strongly non-degenerate bilinear forms.

\begin{mainprop}[Generalization of Donaldson's symplectic subspace criterion to bounded real linear operators on Banach spaces]
\label{mainprop:Donaldson_1996jdg_3_Banach_space}  
Let $\sX$ and $\sY$ be real Banach spaces and $\sH$ be a real Hilbert space containing $\sX$ such that the inclusion $\sX \subset \sH$ is continuous embedding. Let $g = \langle\cdot,\cdot\rangle_\sH$ denote the inner product on $\sH$, let $J$ be an almost complex structure on $\sX$ that is $g$-orthogonal, and let $j$ be an almost complex structure on $\sY$. Let $\sT \in \Hom^c(\sX,\sY)$ be a bounded real linear operator with complex linear and antilinear components given by
\begin{equation}
  \label{eq:Complex_linear_and_anti-linear_operator_components}
  \sT' := \frac{1}{2}(\sT - j\sT J)
  \quad\text{and}\quad
  \sT'' := \frac{1}{2}(\sT + j\sT J).
\end{equation}
Assume that $\Ker \sT$ has a closed complement $\sX_0$, so $\sX = \Ker \sT \oplus \sX_0$ as a direct sum of Banach spaces, and that $\sT$ has closed range with a closed complement $\sY_0$, so $\sY = \Ran \sT \oplus \sY_0$ as a direct sum of Banach spaces. If $L \in \Hom^c(\sY,\sX)$ is 
% TL6-24-2025: Is this enough to determine $L$?  Can we change $L$ on $\Ker \sT^*=(\Ran \sT)^\perp$ without changing $L\sT$?  Do we also want $L$ to be zero on $\Ker \sT^*$?  (This is how we define $L$ in the proof but I don't think that choice of $L$ is needed in the proof and I don't think it affects the hypothesis {eq:Donaldson_1996jdg_Prop_3_Banach_space}
%PF11-10-2025 I replaced "the" by "a"
a partial left inverse of $\sT$ in the sense that
\begin{equation}
  \label{eq:LT_is_1_minus_pi_KerT}
  L\sT = \pi_{\sX_0} \quad\text{on } \sX,
\end{equation}
where $\pi_{\sX_0}: \sX \to \sX_0$ is the continuous projection, and $\sT''$ obeys
\begin{equation}
  \label{eq:Donaldson_1996jdg_Prop_3_Banach_space}
  \|L\sT''\|_{\End(\sX)} < \frac{1}{2},
\end{equation}
and $\pi_{\Ker \sT}: \sX \to \Ker \sT$ is the continuous projection, then
\begin{equation}
  \label{eq:pre-symplectic_form_Banach_space_domain}
  \omega_0 := g\left(\pi_{\Ker \sT} J\cdot,\cdot\right)
\end{equation}
is a continuous real bilinear form on $\sX$ that restricts to a skew-symmetric and weakly non-degenerate form on $\Ker \sT \subset \sX$. If in addition $\sH = \sX$ or $\dim\Ker \sT < \infty$, then $\omega_0$ is (strongly) non-degenerate on $\Ker \sT$.
\end{mainprop}

We prove Proposition \ref{mainprop:Donaldson_1996jdg_3_Banach_space} in Section \ref{sec:Proof_analogue_Donaldson_symplectic_submanifold_criterion} using an argument based in part on ideas suggested to us by Denis Auroux. See the forthcoming Remark \ref{rmk:Relation_between_Donaldson_and_FL_symplectic_subspace_criteria} for a proof that the hypothesis \eqref{eq:Donaldson_1996jdg_proposition_3_approx_complex_linear} of Proposition \ref{prop:Donaldson_1996jdg_3} is implied by our stronger hypothesis \eqref{eq:Donaldson_1996jdg_Prop_3_Banach_space} in Proposition \ref{mainprop:Donaldson_1996jdg_3_Banach_space} for $\sY = \CC$.

\begin{rmk}[On the choice of operator norms in \eqref{eq:Donaldson_1996jdg_Prop_3_Banach_space}]
\label{rmk:Donaldson_1996jdg_3_operator_norms}  
Regarding the inequality \eqref{eq:Donaldson_1996jdg_Prop_3_Banach_space}, for real inner product spaces $V, W$ and any real linear transformation $S \in \Hom_\RR^c(V,W)$, we use the standard induced operator norm,
\begin{equation}
  \label{eq:Standard_operator_norm_vector_space}
  \|S\|_{\Hom(V,W)} := \max_{\{v \in V:\, \|v\|_V = 1\}}\|Sv\|_W.
\end{equation}
See Remark \ref{rmk:Donaldson_1996jdg_3_Hilbert_space_operator_norm} for equivalent choices. One could also use
\begin{equation}
  \label{eq:Operator_norm_finite-dimensional_inner_product_space_orthonormal_basis}
  \seminorm{S}_{\Hom(V,W)} := \max_{1\leq k\leq m}\|Se_i\|_W,
\end{equation}
when $m = \dim_\RR V < \infty$ and $\{e_i\}_{i=1}^m$ is a choice of orthonormal basis of $V$. Clearly,
\[
  \seminorm{S}_{\Hom(V,W)} \leq \|S\|_{\Hom(V,W)}.
\]
On the other hand, if $v \in V$, then $v = c_1e_1 + \cdots + c_me_m$ for $\{c_k\}_{k=1}^m \subset \RR$ with $\|v\|_V^2 = c_1^2 + \cdots + c_m^2$ and
\begin{multline*}
  \|Sv\|_W = \|S(c_1e_1 + \cdots + c_me_m)\|_W \leq \sum_{k=1}^m|c_k|\|Se_k\|_W
  \leq \left(\sum_{k=1}^m|c_k|\right) \max_{1\leq k\leq m}\|Se_k\|_W
  \\
  \leq m^{1/2}\left(\sum_{k=1}^mc_k^2\right)^{1/2} \max_{1\leq k\leq m}\|Se_k\|_W
  = m^{1/2}\|v\|_V\max_{1\leq k\leq m}\|Se_k\|_W
  = m^{1/2}\|v\|_V \seminorm{S}_{\Hom(V,W)},
\end{multline*}
where the third inequality follows from the fact that the quadratic mean is greater than or equal to the arithmetic mean (see Hardy, Littlewood, and Polya \cite[Section 2.9, Theorem 16, Inequality (2.9.1), p. 27]{HardyLittlewoodPolya} and their definition in \cite[Section 2.1, Equation (2.1.3), p. 12]{HardyLittlewoodPolya} of ordinary means).
% COMMENT See https://en.wikipedia.org/wiki/Inequality_(mathematics)
Thus,
\[
  m^{-1/2}\|v\|_{\Hom(V,W)} \leq \seminorm{S}_{\Hom(V,W)} \leq \|S\|_{\Hom(V,W)},
\]
and so the operator norms in \eqref{eq:Operator_norm_finite-dimensional_inner_product_space_orthonormal_basis} and \eqref{eq:Standard_operator_norm_vector_space} are equivalent. In Donaldson's application \cite{DonSympAlmostCx} of his inequality \eqref{eq:Donaldson_1996jdg_proposition_3_approx_complex_linear}, one considers operators $\sT$ that obey
\[
  \|\sT''\|_{\Hom(\CC^n,\CC)} \ll \|\sT'\|_{\Hom(\CC^n,\CC)},
\]
and so the precise choice of operator norm is immaterial since all are equivalent via universal constants.
\qed\end{rmk}

\begin{rmk}[Existence of partial left inverses in Proposition \ref{mainprop:Donaldson_1996jdg_3_Banach_space}]
\label{eq:Existence_partial_left_inverses}  
By our hypothesis in Proposition \ref{mainprop:Donaldson_1996jdg_3_Banach_space}, the operator $\sT$ has closed range, so $\Ran \sT \subset \sY$ is a closed subspace. Recall that $\sX_0\subset \sX$ denotes the closed complement of the closed subspace $\Ker \sT$ provided by our hypotheses. We see that $\sT:\sX_0 \to \Ran \sT$ is a bounded bijective operator and thus an
%PF7-2-2025 Define what we mean by "isomorphism of Banach spaces" or rephrase
isomorphism of Banach spaces by the Open Mapping Theorem (see, for example, Brezis \cite[Section 2.3, Theorem 2.6, p. 35]{Brezis}). Hence, the bounded operator $\sT \in \Hom^c(\sX,\sY)$ has a partial left inverse $L \in \Hom^c(\sY, \sX)$ in the sense of \eqref{eq:LT_is_1_minus_pi_KerT}. To see this, observe that we may define $L \in \Hom^c(\sY, \sX)$ by setting
\[
  L
  :=
  \begin{cases}
    \sT^{-1} &\text{on } \Ran \sT,
    \\
    0 &\text{on } \sY_0,
  \end{cases}
\]
where $\sY_0\subset \sY$ is the closed complement of the closed subspace $\Ran \sT$ provided by our hypotheses.
\qed\end{rmk}  

\begin{rmk}[On the choice of operator norm in \eqref{eq:Donaldson_1996jdg_Prop_3_Banach_space}]
\label{rmk:Donaldson_1996jdg_3_Hilbert_space_operator_norm}    
The norm in \eqref{eq:Donaldson_1996jdg_Prop_3_Banach_space} of a bounded operator $S \in \Hom^c(\cV,\cW)$ between normed vector spaces $\cV$ and $\cW$ can be defined as usual by (see, for example, Rudin \cite[Theorem 4.1, p. 92]{Rudin})
\begin{equation}
  \label{eq:Operator_norm}
  \|S\|_{\Hom^c(\cV,\cW)}
  :=
  \sup_{\{v\in\cV:\, \|v\|_\cV \leq 1\}} \|Sv\|_\cW,
\end{equation}
though other equivalent choices are useful in our applications (see Section \ref{sec:Operator_norm_choices_Donaldson_symplectic_subspace_criteria}).
\qed\end{rmk}

We will use the forthcoming Corollary \ref{maincor:Adjoint_DonaldsonCriteria} when computing virtual Morse--Bott indices in Chapter \ref{chap:Calculation_virtual_Morse-Bott_indices_via_Atiyah-Singer_index_theorem}. Recall from Rudin \cite[Section 4.6, p. 95]{Rudin} that if $\sM$ is a subspace of a Banach space $\sY$ with continuous dual space $\sY$, then $\sM^\perp := \{\beta \in \sY^*: \langle y, \beta\rangle_{\sY\times\sY^*} = 0 \text{ for all } y \in \sM\}$ is the \emph{annihilator} of $\sM$. In particular, $(\Ran \sT)^\perp \subset \sY^*$ is the annihilator of $\Ran \sT \subset \sY$ in the following  

\begin{maincor}[Generalization of Donaldson's symplectic subspace criterion to adjoints of bounded real linear operators on Banach spaces]
\label{maincor:Adjoint_DonaldsonCriteria}
Let $\sX$ and $\sY$ be real Banach spaces and $\sG$ be a real Hilbert space that is a dense subspace of $\sY$ such that the inclusion $\sG \subset \sY$ and thus also $\sY^* \subset \sG^* \cong \sG$ are continuous embeddings of continuous dual spaces. Let $g = \langle\cdot,\cdot\rangle_\sG$ denote the inner product on $\sG$, let $J$ be an almost complex structure on $\sX$, and let $j$ be an almost complex structure on $\sY$ that is $g$-orthogonal. Let $\sT \in \Hom^c(\sX,\sY)$ be a bounded real linear operator with complex linear and antilinear components given by \eqref{eq:Complex_linear_and_anti-linear_operator_components} and let $\sT^* \in \Hom^c(\sY^*,\sX^*)$ be the Banach space adjoint of $\sT$, where $\sX^*$ and $\sY^*$ are the continuous dual spaces of $\sX$ and $\sY$, respectively. Assume that $\sT$ has closed range with a closed complement $\sY_0$, so $\sY = \Ran \sT \oplus \sY_0$ and $\sY^* = \sY_0^\perp \oplus (\Ran \sT)^\perp = \sY_0^\perp \oplus \Ker \sT^*$ as direct sums of Banach spaces. If $\sT$ has a partial right inverse $R \in \Hom^c(\sY,\sX)$, so $\sT R = \pi_{\Ran \sT}$ where $\pi_{\Ran \sT}: \sY \to \Ran \sT$ is the continuous projection, and $\sT''$ obeys
\begin{equation}
  \label{eq:Donaldson_1996jdg_Prop_3_corrected_condition_Hilbert_space_right_inverse}
 \|\sT''R\|_{\End(\sY)} < \frac{1}{2},
\end{equation}
and $\pi_{\Ker \sT^*}: \sY \to \Ker \sT^*$ is the continuous projection, then
\begin{equation}
  \label{eq:pre-symplectic_form_Hilbert_space_domain_adjoint}
  \Omega_0 := g(\pi_{\Ker \sT^*} j^*\cdot,\cdot)
\end{equation}
is a continuous real bilinear form on $\sY^* \subset \sG$ that restricts to a skew-symmetric and weakly non-degenerate form on $\Ker \sT^* \subset \sY^*$. If in addition $\sG = \sY$, so $\sG \cong \sG^* = \sY^*$, or $\dim\Ker \sT^* < \infty$, then $\Omega_0$ is (strongly) non-degenerate on $\Ker \sT^*$.
\end{maincor}

We prove Corollary \ref{maincor:Adjoint_DonaldsonCriteria} in Section \ref{sec:Proof_analogue_Donaldson_symplectic_submanifold_criterion}. 

While Proposition \ref{mainprop:Donaldson_1996jdg_3_Banach_space} and Corollary \ref{maincor:Adjoint_DonaldsonCriteria} provide useful generalizations of Proposition \ref{prop:Donaldson_1996jdg_3} --- for example, they provide insights into Auroux's generalization \cite[Section 1, Theorem 2, p. 972]{AurouxFamilySymp} of Donaldson's \cite[Theorem 5, p. 670]{DonSympAlmostCx} from complex line bundles to complex vector bundles of arbitrary rank --- it would be difficult in our application (see Theorem \ref{mainthm:AH_structure_bounded_evalue_spaces_non-Abelian_monopoles_symp_4-mflds} in Section \ref{subsec:Almost_Hermitian_structure_moduli_spaces_non-Abelian_monopoles_symplectic_4-manifolds}) to verify the hypothesis \eqref{eq:Donaldson_1996jdg_Prop_3_Banach_space} of Proposition \ref{mainprop:Donaldson_1996jdg_3_Banach_space}. The reason, in essence, is that while the linearization $\sT_{A,\varphi,\psi,r}$ of the non-Abelian monopole equations with a regularized Taubes' perturbation \eqref{eq:SO(3)_monopole_equations_(0,2)_curvature_intro}, \eqref{eq:SO(3)_monopole_equations_Dirac_almost_Hermitian_intro}, and \eqref{eq:SO(3)_monopole_equations_(1,1)_curvature_perturbed_intro_regular} and Coulomb gauge condition is a first order linear elliptic differential operator, we do not have control over the small eigenvalues of the $L^2$ self-adjoint operator
\[
  T_{A,\varphi,\psi,r} = \begin{pmatrix} 0 & \sT_{A,\varphi,\psi,r}^* \\ \sT_{A,\varphi,\psi,r} & 0 \end{pmatrix}
\]
as Taubes' perturbation parameter $r \to \infty$ or the point $[A,\varphi,\psi]$ varies over the moduli space of non-Abelian monopole equations. We shall instead rely on the forthcoming Theorem \ref{mainthm:Donaldson_1996jdg_3_Hilbert_space} that further generalizes Proposition \ref{prop:Donaldson_1996jdg_3}. Our proofs of Theorem \ref{mainthm:Donaldson_1996jdg_3_Hilbert_space} and Corollary \ref{maincor:Donaldson_1996jdg_3_Hilbert_space_non-self-adjoint} lean heavily on spectral theory for unbounded linear operators on Banach spaces and we review the relevant background in Section \ref{sec:Spectral_theory_unbounded_operators}.

\begin{mainthm}[Generalization of Donaldson's symplectic subspace criterion to unbounded self-adjoint real linear operators on Hilbert spaces]
\label{mainthm:Donaldson_1996jdg_3_Hilbert_space}  
Let $\sH$ be a real Hilbert space with inner product $g = \langle\cdot,\cdot\rangle_\sH$ and $g$-orthogonal almost complex structure $J$. Let $T \in \End(\sH)$ be a densely defined, self-adjoint (and thus closed), unbounded real linear operator with compact resolvent. Let $T', T'' \in \End(\sH)$ be the complex linear and antilinear components of $T$ given by \eqref{eq:Complex_linear_and_anti-linear_operator_components} with $j=J$. Assume that $T'$ is self-adjoint and that $T''$ is bounded. Let $\mu \in (0,\infty)$ and $\delta \in (0,\mu/2)$ be constants such that 
\begin{equation}
  \label{eq:Spectral_gap_T}
   (-\mu-2\delta,-\mu+2\delta) \cup (\mu-2\delta,\mu+2\delta) \subset \rho(T),
\end{equation}
where $\rho(T) \subset \CC$ is the resolvent set of $T$, and require that $T''$ obeys
\begin{equation}
  \label{eq:Donaldson_1996jdg_Prop_3_Hilbert_space}
  \|T''\|_{\End(\sH)} < \delta/2.
\end{equation}
Then $T'$ obeys
\begin{equation}
  \label{eq:Spectral_gap_Tprime}
   (-\mu-\delta,-\mu+\delta) \cup (\mu-\delta,\mu+\delta) \subset \rho(T'),
\end{equation}
and if $\Pi_\mu \in \End^c(\sH)$ denotes the orthogonal projection from $\sH$ onto the sum of the eigenspaces of $T$ with eigenvalues in $(-\mu,\mu)$, then
\begin{equation}
  \label{eq:pre-symplectic_form_Hilbert_space_domain}
  \omega_\mu := g\left(\Pi_\mu J\cdot,\cdot\right)
\end{equation}
is a continuous real bilinear form on $\sH$ that restricts to a skew-symmetric and (strongly) non-degenerate two-form on $\Ran \Pi_\mu = \Ker(\Pi_\mu^\perp T) \subset \sH$, where $\Pi_\mu^\perp := \id_\sH - \Pi_\mu$.
\end{mainthm}

We prove Theorem \ref{mainthm:Donaldson_1996jdg_3_Hilbert_space} in Section \ref{sec:Proof_generalized_Donaldson_symplectic_subspace_criterion_spectral_projection}.  

\begin{rmk}[On the hypotheses of Theorem \ref{mainthm:Donaldson_1996jdg_3_Hilbert_space}]
\label{rmk:Donaldson_1996jdg_3_Hilbert_space}
Our Theorem \ref{mainthm:Donaldson_1996jdg_3_Hilbert_space} could be generalized without difficulty in ways that we describe here, but the proofs of such generalizations appear tedious and are not needed for the proof of the forthcoming Corollary \ref{maincor:Donaldson_1996jdg_3_Hilbert_space_non-self-adjoint}. For example, the bounded interval $(-\mu,\mu) \Subset \RR$ that defines $\Pi_\mu$ could be replaced by a bounded open subset $I \Subset \RR$. Moreover, it suffices that the operators $T$ and $T'$ be normal rather than self-adjoint, in which case the interval $(-\mu,\mu) \Subset \RR$ would be replaced by a bounded open subset $V \Subset \CC$. The spectral gap condition \eqref{eq:Spectral_gap_T} and constant in the estimate \eqref{eq:Donaldson_1996jdg_Prop_3_Hilbert_space} would then need to be adjusted accordingly.
\qed\end{rmk}

\begin{maincor}[Generalization of Donaldson's symplectic subspace criterion to unbounded real linear operators on Hilbert spaces]
\label{maincor:Donaldson_1996jdg_3_Hilbert_space_non-self-adjoint}  
Let $\sH_i$ be a real Hilbert space with inner product $g_i = \langle\cdot,\cdot\rangle_{\sH_i}$ and $g_i$-orthogonal almost complex structure $J_i$ for $i=1,2$. Let $\sT \in \Hom(\sH_1,\sH_2)$ be a closed, densely defined, unbounded real linear operator and assume that the self-adjoint (and thus closed) operator
\[
  T := \begin{pmatrix} 0 & \sT^* \\ \sT & 0 \end{pmatrix} \in \End(\sH_1\oplus\sH_2)
\]
has compact resolvent. Let $\sT', \sT'' \in \Hom(\sH_1,\sH_2)$ be the complex linear and antilinear components of $\sT$ given by \eqref{eq:Complex_linear_and_anti-linear_operator_components} with $J := J_1$ and $j := J_2$ and assume that $\sT''$ is bounded. Assume that the resolvent set $\rho(T)$ obeys the spectral gap condition \eqref{eq:Spectral_gap_T} and that $\sT''$ obeys
\begin{equation}
  \label{eq:Donaldson_1996jdg_Prop_3_Hilbert_space_non-self-adjoint}
  \|\sT''\|_{\Hom(\sH_1,\sH_2)} < \delta/4.
\end{equation}
If $\Pi_{1,\mu^2} \in \End^c(\sH_1)$ denotes the orthogonal projection from $\sH$ onto the sum of the eigenspaces of $\sT^*\sT$ with eigenvalues in the interval $[0,\mu^2)$, then
\begin{equation}
  \label{eq:pre-symplectic_form_Hilbert_space_domain_non-self-adjoint}
  \omega_{1,\mu} := g_1\left(\Pi_{1,\mu^2} J_1\cdot,\cdot\right)
\end{equation}
is a continuous real bilinear form on $\sH_1$ that restricts to a skew-symmetric and (strongly) non-degenerate two-form on $\Ran \Pi_{1,\mu^2} = \Ker(\Pi_{1,\mu^2}^\perp \sT^*\sT) \subset \sH_1$, where $\Pi_{1,\mu^2}^\perp := \id_{\sH_1} - \Pi_{1,\mu^2}$. Similarly, if $\Pi_{2,\mu^2} \in \End^c(\sH_2)$ denotes the orthogonal projection from $\sH$ onto the sum of the eigenspaces of $\sT\sT^*$ with eigenvalues in the interval $[0,\mu^2)$, then
\begin{equation}
  \label{eq:pre-symplectic_form_Hilbert_space_codomain_non-self-adjoint}
  \omega_{2,\mu} := g_2\left(\Pi_{2,\mu^2} J_2\cdot,\cdot\right)
\end{equation}
is a continuous real bilinear form on $\sH_2$ that restricts to a skew-symmetric and (strongly) non-degenerate two-form on $\Ran \Pi_{2,\mu^2} = \Ker(\Pi_{2,\mu^2}^\perp\sT\sT^*) \subset \sH_2$, where $\Pi_{2,\mu^2}^\perp := \id_{\sH_2} - \Pi_{2,\mu^2}$.
\end{maincor}

We prove Corollary \ref{maincor:Donaldson_1996jdg_3_Hilbert_space_non-self-adjoint} in Section \ref{sec:Proof_generalized_Donaldson_symplectic_subspace_criterion_spectral_projection}.  

Corollary \ref{maincor:Donaldson_1996jdg_3_Hilbert_space_non-self-adjoint} constructs a skew-symmetric non-degenerate two-form $\omega_{1,\mu}$ on the direct sum of the eigenspaces of $\sT^*\sT \in \End(\sH_1)$ with eigenvalues less than $\mu^2$ and a skew-symmetric non-degenerate two-form $\omega_{2,\mu}$ on the direct sum of the eigenspaces of $\sT\sT^* \in \End(\sH_2)$ with eigenvalues less than $\mu^2$. See the end of the proof of Corollary \ref{maincor:Donaldson_1996jdg_3_Hilbert_space_non-self-adjoint} for an equivalent but slightly different statement of its conclusions.
%PF11-24-2025 Propagate use of \End^c,\Hom^c for bounded operators

\subsection{Lower bounds for spectral gaps}
\label{subsec:Main_results_lower_bounds_spectral_gaps}
Theorem \ref{mainthm:Donaldson_1996jdg_3_Hilbert_space} and Corollary \ref{maincor:Donaldson_1996jdg_3_Hilbert_space_non-self-adjoint} require the spectral gap hypothesis \eqref{eq:Spectral_gap_T}. In applications, we would ideally like the gap width $4\delta$ to remain uniformly bounded below by a positive constant as the parameters defining $T$ vary in a noncompact space, though the distance $\mu$ from the gap centers from the origin need not be constrained. The results in this subsection illustrate the difficulties involved in applications of Theorem \ref{mainthm:Donaldson_1996jdg_3_Hilbert_space}. We include these results for completeness, but they are not central and are not required by our program to prove Conjecture \ref{mainconj:BMY_Seiberg-Witten}.

% PF1-29-2025 Standardize on $L^2(X,E)$ rather than $L^2(X;E)$, etc.
% PF2-20-2025 Standardize on $(E,H)$
\begin{mainthm}[Spectrum and spectral gaps for self-adjoint elliptic pseudodifferential operators on vector bundles over closed manifolds]
\label{mainthm:Donaldson_1996jdg_3_elliptic_operator_order_m_geq_d}
Let $(E,H)$ be a smooth Hermitian vector bundle over a closed orientable smooth Riemannian manifold $(X,g)$ of dimension $d\geq 2$. Let $A:C^\infty(X,E) \to C^\infty(X,E)$ be an $L^2$ elliptic pseudodifferential operator of order $m$ and principal symbol $a_m$. Then $A$ is an unbounded operator on $L^2(X,E)$ with dense domain $W^{m,2}(X,E) \subset L^2(X,E)$, compact resolvent, and spectrum $\sigma(A)$ that is a discrete subset of $\CC$ consisting of eigenvalues with finite multiplicity. If $A$ is self-adjoint, then $\sigma(A) \subset \RR$ and if in addition $m\geq d$, then there are constants $\mu = \mu(A,g,h)\in (0,\infty)$ and $\delta = \delta(a_m,g,h)\in (0,\mu/2)$ such that
\begin{equation}
  \label{eq:Spectral_gap_A}
   (-\mu-2\delta,-\mu+2\delta) \cup (\mu-2\delta,\mu+2\delta) \subset \rho(A),
\end{equation}
where $\rho(A)$ is the resolvent set for the unbounded operator $A$ on $L^2(X,E)$.
\end{mainthm}

We prove Theorem \ref{mainthm:Donaldson_1996jdg_3_elliptic_operator_order_m_geq_d} in Section \ref{sec:Lower_bounds_spectral_gaps}. The reason for the hypothesis $m\geq d$ in Theorem \ref{mainthm:Donaldson_1996jdg_3_elliptic_operator_order_m_geq_d}, where $m$ is the order of $A$ and $d$ is the dimension of $X$, that according to Weyl's Law (see Section \ref{sec:Weyl_asymptotic_formula_eigenvalues_elliptic_operator}) for the asymptotic growth of the eigenvalues $\lambda_k$ of $A$, \emph{repeated according to their multiplicity}, one has $|\lambda_k| \sim k^{m/d}$ as $k \to \infty$. In the case of the Dirac operator on a space form of positive curvature (and dimension $d$), it is known that the multiplicities of distinct eigenvalues $\lambda_k$ also grow with $k$ and as a result one has $|\lambda_k| \sim k$ as $k \to \infty$ (see, for example, B\"ar \cite[Section 2, Theorem 1, p. and Section 3, Theorem 2, p. 75]{Baer_1996} or Ginoux \cite[Theorem 2.1.3, p. 31 or Theorem 2.1.4 and Corollary 2.1.5, p. 34]{Ginoux_dirac_spectrum}). Estimates for gaps in the spectra of Dirac and related elliptic operators under a variety of different conditions are provided, for example, by Anghel \cite{Anghel_1993}, Ashbaugh \cite{Ashbaugh_2002}, Ashbaugh and Hermi \cite{Ashbaugh_Hermi_2004}, Chen \cite{Chen_2009}, Chen and Sun \cite{Chen_Sun_2008}, Huang, Chen, and Sun \cite{Huang_Chen_Sun_2011}, and Ilias and Makhoul \cite{Ilias_Makhoul_2009}.

The explicit calculations of Dirac operator eigenvalue multiplicities in the preceding references rely on symmetry and so we cannot expect similar results over closed Riemannian manifolds with arbitrary metrics. However, the following qualitative property of the spectra of (generalized) Dirac operators will suffice for our applications. If $A$ is a $C^\infty$ orthogonal connection on a Riemannian vector bundle $(V,h)$ over a smooth Riemannian manifold $(X,g)$ with spin${}^c$ structure $(\rho,W)$, then
\begin{equation}
  \label{eq:Coupled_Dirac_operator}
  D_A = \rho\circ\nabla_A^{W\otimes V}: \Omega^0(W\otimes V) \to \Omega^0(W\otimes V)
\end{equation}
is the Dirac operator defined with covariant derivative $\nabla_A^{W\otimes V}$ determined by the tensor product of a spin${}^c$ connection $\nabla^W$ on $W$, the orthogonal connection $A$ on $V$, and Clifford multiplication $\rho:T^*X \to \End_{\CC}(W)$ --- see, for example, Berline, Getzler, and Vergne \cite{BerlineGetzlerVergne}, Boo\ss, Bavnbek, and Wojciechowski \cite{Booss-Bavnbek_Wojciechowski_elliptic_boundary_problems_dirac_operators}, or Lawson and Michelsohn \cite{LM}. The forthcoming theorem and its corollaries are stated for the case of even-dimensional base manifolds only for convenience and minor modifications yield the corresponding results for odd-dimensional manifolds. 

\begin{mainthm}[Spectral gaps for coupled Dirac operators over closed manifolds]
\label{mainthm:Lower_bound_spectral_gap_coupled_Dirac_operator}
Let $(V,h)$ be a smooth Riemannian vector bundle over a closed smooth Riemannian manifold $(X,g)$ with even dimension $d \geq 2$, let $A_L$ be a smooth unitary connection associated to a spin${}^c$ structure $(\rho,W)$ over $X$ as in the forthcoming Lemma \ref{lem:Bochner-Weitzenbock_identity_coupled_Dirac_operator}, and let $Z \in [1,\infty)$ be a constant. Then there are constants $C_0, K, R_0 \in [1,\infty)$, depending at most on $\|\rho(F_{A_L})\|_{C^0(\End(W^+\otimes V))},g,h,Z$, with the following significance. If $A$ is a smooth unitary connection on $V$, then the Dirac operator $D_A$ in \eqref{eq:Coupled_Dirac_operator} is an unbounded self-adjoint operator on $L^2(W\otimes V)$ with dense domain $\Dom(D_A) = W^{1,2}(W\otimes V) \subset L^2(W\otimes V)$ and compact resolvent. The spectrum $\sigma(D_A) = \{\lambda_k(D_A)\}_{k=1}^\infty \subset \RR$ of $D_A$
%PF3-20-2025 Need to distinguish between this and D_A^*D_A on L^2(W^+\otimes E)
on $L^2(W\otimes V)$ is a discrete subset consisting of eigenvalues, repeated according to their finite multiplicities. If $r \in [1,\infty)$ is a parameter and the curvature $F_A$ obeys
%PF4-2-2024 Explain origin of bound below.
\begin{equation}
  \label{eq:rho_FA_End_W+_otimes_V_leq_constant_1+r}
  \|\rho(F_A)\|_{C^0(\End(W^+\otimes V))} \leq Z + r,
\end{equation}
then the average gap between \emph{distinct eigenvalues} has the lower bound,
%PF4-2-2025 Replace "4" by $d$ in statement and proof
\begin{equation}
  \label{eq:Lower_bound_spectral_gap_coupled_Dirac_operator}
  \frac{N(\lambda)}{\lambda} \geq \frac{C_0}{K^{d-1}}\frac{1}{r^{(d-1)/2}}, \quad\text{for } r \geq R_0
  \text{ and } \lambda \geq Kr^{1/2},
\end{equation}
where $N(\lambda) = \sum_{\{\lambda_k:\, |\lambda_k| \leq\lambda\}} 1$ is the counting function \eqref{eq:N_lambda}.
\end{mainthm}

The $C^0$ bound \eqref{eq:rho_FA_End_W+_otimes_V_leq_constant_1+r} for $\rho(F_A)$ in the hypotheses of Theorem \ref{mainthm:Lower_bound_spectral_gap_coupled_Dirac_operator} is motivated by the pointwise estimate \eqref{eq:Taubes_1996_SW_to_Gr_eq_2-12_and_2-13_regular} for $\rho(F_A^+)_0 = \rho(F_A)_0 \in \End(W_\can^+\otimes E)$ provided by the forthcoming Corollary \ref{maincor:Taubes_1996_SW_to_Gr_eq_2-12_and_2-13_regular} when $(A,\varphi,\psi)$ solves the non-Abelian monopole equations \eqref{eq:SO(3)_monopole_equations_almost_Hermitian_perturbed_intro_regular} with a regularized Taubes perturbation.

The lower bound in \eqref{eq:Lower_bound_spectral_gap_coupled_Dirac_operator} in Theorem \ref{mainthm:Lower_bound_spectral_gap_coupled_Dirac_operator} can be improved, although it appears difficult to eliminate the dependence on the parameter $r$ introduced by the bound \eqref{eq:rho_FA_End_W+_otimes_V_leq_constant_1+r} in the hypotheses --- see Feehan \cite{Feehan_weyl_law_clusters}. In applications, the following extension of Theorem \ref{mainthm:Lower_bound_spectral_gap_coupled_Dirac_operator} is useful.

\begin{maincor}[Spectral gaps for coupled Dirac operators with vector potentials over closed manifolds]
\label{maincor:Lower_bound_spectral_gap_coupled_Dirac_operator_plus_vector_potential}  
Continue the notation and hypotheses of Theorem \ref{mainthm:Lower_bound_spectral_gap_coupled_Dirac_operator} and in addition, let $L_i \in C^\infty(\End(W\otimes V))$ for $i=0,1$ be pointwise Hermitian vector bundle endomorphisms with $L := L_0 + rL_1$. Then the qualitative conclusions of Theorem \ref{mainthm:Lower_bound_spectral_gap_coupled_Dirac_operator} continue to hold for the coupled Dirac operator $D_A+L$ with vector potential and its spectrum and the average gap between \emph{distinct eigenvalues} obeys
%PF4-2-2025 Replace "4" by $d$ in statement and proof
\begin{equation}
  \label{eq:Lower_bound_spectral_gap_coupled_Dirac_operator_plus_vector_potential}
  \frac{N(\lambda)}{\lambda} \geq \frac{C_2}{K_1^{d-1}}\frac{1}{r^{d-1}}, \quad\text{for all } r \geq R_1
  \text{ and } \lambda \geq K_1r,
\end{equation}
where the constants $C_2, K_1, R_1 \in [1,\infty)$ may now depend on $\|L_i\|_{C^0(\End(W\otimes V))}$ for $i=0,1$ in addition to $\|\rho(F_{A_L})\|_{C^0(\End(W^+\otimes V))},g,h,Z$. 
\end{maincor}
% PF3-29-2025 Double check corollary proof and dependencies and update theorem proof.

We now specialize Corollary \ref{maincor:Lower_bound_spectral_gap_coupled_Dirac_operator_plus_vector_potential} to the case of an almost K\"ahler manifold $(X,g,J,\omega)$ of dimension $d=4$ and smooth orthogonal connection $A^V$ on the Hermitian vector bundle $V = \fsl(E) \otimes E$ induced by a smooth unitary connection $A$ on a rank two Hermitian vector bundle $E$ over $X$, where for some $(\varphi,\psi) \in C^\infty(E\oplus\Lambda^{0,2}(X)\otimes E)$ the triple $(A,\varphi,\psi)$ is a $C^\infty$ solution to the non-Abelian monopole equations \eqref{eq:SO(3)_monopole_equations_almost_Hermitian_perturbed_intro_regular} with a regularized Taubes perturbation and parameters $r \in [1,\infty)$ and $\gamma \in (0,\infty)$.

\begin{maincor}[Spectral gaps for coupled Dirac operators with vector potentials over closed four-manifolds]
\label{maincor:Lower_bound_spectral_gap_coupled_Dirac_operator_plus_vector_potential_four-manifold}
Let $(E,H)$ be a smooth Hermitian vector bundle with complex rank two
%PF3-20-2025 rank 2 necessary?
over a closed four-dimensional
%PF3-20-2025 almost Hermitian okay?
almost K\"ahler manifold $(X,g,J,\omega)$ and smooth unitary connection $A_d$ on $\det E$. Let $V = \fsl(E)\otimes E$ and $L_i \in C^\infty(\End(W\otimes V))$ for $i=0,1$ be pointwise Hermitian vector bundle endomorphisms with $L := L_0 + rL_1$. Then there are constants $C_2, K_1, R_1 \in [1,\infty)$, depending at most on $\|\rho(F_{A_L}^+)\|_{C^0(\End(W^+\otimes V))},g,H,\omega$, and $\|L_i\|_{C^0(\End(W\otimes V))}$ for $i=0,1$  with the following significance. If $(A,\varphi,\psi)$ is a $C^\infty$ solution to the non-Abelian monopole equations \eqref{eq:SO(3)_monopole_equations_almost_Hermitian_perturbed_intro_regular} with a regularized Taubes perturbation with parameters $r \in [1,\infty)$ and $\gamma \in (0,\infty)$, then the qualitative conclusions of Theorem \ref{mainthm:Lower_bound_spectral_gap_coupled_Dirac_operator} continue to hold for the coupled Dirac operator $D_A+L$ with vector potential and its spectrum and the average gap between \emph{distinct eigenvalues} obeys
\begin{equation}
  \label{eq:Lower_bound_spectral_gap_coupled_Dirac_operator_plus_vector_potential_four-manifold}
  \frac{N(\lambda)}{\lambda} \geq \frac{C_2}{K_1^3}\frac{1}{r^3}, \quad\text{for all } r \geq R_1
  \text{ and } \lambda \geq K_1r.
\end{equation}
\end{maincor}

We prove Theorem \ref{mainthm:Lower_bound_spectral_gap_coupled_Dirac_operator} in Section \ref{sec:Lower_bounds_spectral_gaps_Dirac_operators} and Corollaries \ref{maincor:Lower_bound_spectral_gap_coupled_Dirac_operator_plus_vector_potential} and \ref{maincor:Lower_bound_spectral_gap_coupled_Dirac_operator_plus_vector_potential_four-manifold} in Section \ref{sec:Lower_bounds_spectral_gaps_Dirac_operators_vector_potentials}.

\subsection[Analogues for non-Abelian monopoles of Taubes' pointwise estimates]{Analogues for non-Abelian monopoles of Taubes' pointwise estimates for Seiberg--Witten monopole sections}
\label{subsec:Analogues_non-Abelian_monopoles_Taubes_estimates_SW_monopole_sections}
In Taubes' tetralogy (see \cite{TauSWGromov, TauCountingPseudo, TauGRtoSW, TauGRtoSWcounting} or \cite{Taubes_2000}) to prove the equivalence of Seiberg--Witten and his Gromov invariants of closed, symplectic four-manifolds $(X,\omega)$, he begins by proving (see \cite[Section 1 (e), Theorem 1.3, p. 852]{TauSWGromov}) the existence of a compact, complex curve $\Sigma$ and a pseudoholomorphic map, $f:\Sigma \to X$, given the existence of an unbounded sequence of perturbation parameters $\{r_n\}_{n\in\NN} \subset [1,\infty)$ and solutions to the corresponding perturbed Seiberg--Witten monopole equations \cite[Section 1 (d), Equation (1.20), p. 851]{TauSWGromov} for the canonical \spinc structure $\fs_\can$ on $(X,g,J,\omega)$ (see the forthcoming \eqref{eq:Canonical_spinc_bundles} and \eqref{eq:Canonical_Clifford_multiplication}) twisted by a Hermitian line bundle $E$. Taubes notes in \cite[Section 1 (f), p. 853]{TauSWGromov} that his proof of \cite[Section 1 (e), Theorem 1.3, p. 852]{TauSWGromov} relies on the eight pointwise estimates that he lists in \cite[Section 1 (f), Equations (1.24) and (1.25, p. 853]{TauSWGromov} for the squared pointwise norms of the $\Omega^0(E)$ and $\Omega^{0,2}(E)$ components of the coupled spinor (after rescaling by $\sqrt{r})$ in the perturbed Seiberg--Witten monopole equations.

In Chapter \ref{chap:Analogues_non-Abelian_monopoles_Taubes_estimates_Seiberg-Witten_monopole_sections}, we derive differential inequalities and estimates for the squared pointwise norms of sections, $|\varphi|^2$ and $|\psi|^2$, when $(A,\varphi,\psi)$ is a solution to the non-Abelian equations with a singular or regularized Taubes perturbation as in Section \ref{sec:Perturbed_non-Abelian_monopole_equations_almost_Hermitian_4-manifolds}. These estimates are analogues of those due to Taubes in \cite[Section 1 (f), Equations (1.24) (1) and (2), p. 853]{TauSWGromov}, provided by \cite[Section 2 (a), Proposition 2.1, p. 854, and Section 2 (c), Proposition 2.3, p. 856]{TauSWGromov} for the squared pointwise norms of sections given by solutions to Taubes' perturbation of the Seiberg--Witten monopole equations. 

By analogy with Taubes' rescaling \cite[Section 1 (d), Equation (1.19), p. 851]{TauSWGromov} of solutions to his perturbation of the Seiberg--Witten monopole equations, for $r \in (0,\infty)$ and $(\varphi,\psi) \in \Omega^0(E)\oplus\Omega^{0,2}(E)$ we define $(\alpha,\beta) \in \Omega^0(E)\oplus\Omega^{0,2}(E)$ by
\begin{equation}
  \label{eq:Taubes_1996_SW_to_Gr_eq_1-19}
  (\varphi,\psi) = r^{1/2}(\alpha,\beta) \in \Omega^0(E)\oplus\Omega^{0,2}(E).
\end{equation}
We begin by stating the following analogue of Taubes' estimate in \cite[Section 1 (f), Equation (1.24) (1), p. 853]{TauSWGromov}, provided by \cite[Section 2 (a), Proposition 2.1, p. 854]{TauSWGromov}, for the squared pointwise norms of sections given by solutions to Taubes' perturbation of the Seiberg--Witten monopole equations.

\begin{mainprop}[Pointwise estimate for the squared pointwise norms of sections of $E$ with singular Taubes perturbations]
\label{mainprop:Taubes_1996_SW_to_Gr_2-1}
Let $(X,g,J,\omega)$ be an almost Hermitian four-manifold and $(E,H)$ be a smooth, Hermitian vector bundle over $X$ with complex rank
%PF9-8-2024 Do we need E having rank 2 here?
$2$ and smooth, unitary connection $A_d$ on $\det E$. Then there is a constant $z\in[1,\infty)$, depending at most on the Riemannian metric $g$, non-degenerate two-form $\omega$, and connection $A_d$, that has the following significance. Let $p \in [4,\infty)$ be a constant and $(A,\varphi,\psi)$ in \eqref{eq:A_varphi_psi_in_W1p_times_W2p} be a
% PF8-16-2024 Correct constant dependencies
solution to the system \eqref{eq:SO(3)_monopole_equations_almost_Hermitian_perturbed_intro} of non-Abelian monopole equations with a singular Taubes perturbation defined by a constant $r > 0$. If $(\varphi,\psi) = r^{1/2}(\alpha,\beta) \in W^{2,p}(E\oplus\Lambda^{0,2}(E))$ as in \eqref{eq:Taubes_1996_SW_to_Gr_eq_1-19}, then $(\alpha,\beta)$ obeys the pointwise inequality,
\begin{equation}
  \label{eq:Taubes_1996_SW_to_Gr_eq_2-2}
  |\alpha|_E^2 + |\beta|_{\Lambda^{0,2}(E)}^2
  \leq
  \frac{1}{3} + \frac{z}{r} \quad\text{on } X.
\end{equation}
\end{mainprop}

We prove Proposition \ref{mainprop:Taubes_1996_SW_to_Gr_2-1} in Section \ref{sec:Pointwise_estimate_holomorphic_component_coupled_spinor}. The proof of following analogue of Taubes \cite[Section 2 (b), Proposition 2.3, p. 856]{TauSWGromov} is considerably more difficult than that of Proposition \ref{mainprop:Taubes_1996_SW_to_Gr_2-1}.

\begin{mainthm}[Pointwise estimate for the squared pointwise norms of sections of $\Lambda^{0,2}(E)$ with singular Taubes perturbations]
\label{mainthm:Taubes_1996_SW_to_Gr_2-3}
Let $(X,g,J,\omega)$ be a closed, symplectic
% PF7-10-2024 Would almost Hermitian be okay?
four-manifold and $(E,H)$ be a smooth, Hermitian vector bundle over $X$ with complex rank
%PF9-8-2024 Do we need E having rank 2 here?
$2$ and smooth, unitary connection $A_d$ on $\det E$. Then there are constants $z',  z'', \kappa_1 \in [1,\infty)$, depending at most on the Riemannian metric $g$, symplectic form $\omega$, and connection $A_d$, that have the following significance. Let $p \in [4,\infty)$ be a constant and $(A,\varphi,\psi)$ in \eqref{eq:A_varphi_psi_in_W1p_times_W2p} be a
% PF8-16-2024 Correct constant dependencies
solution to the system \eqref{eq:SO(3)_monopole_equations_almost_Hermitian_perturbed_intro} of non-Abelian monopole equations with a singular Taubes perturbation defined by a constant $r \geq \kappa_1$. If $(\varphi,\psi) = r^{1/2}(\alpha,\beta) \in W^{2,p}(E\oplus\Lambda^{0,2}(E))$ as in \eqref{eq:Taubes_1996_SW_to_Gr_eq_1-19}, then $(\alpha,\beta)$ obeys the pointwise inequality,
\begin{equation}
  \label{eq:Taubes_1996_SW_to_Gr_eq_2-11}
  |\beta|_{\Lambda^{0,2}(E)}^2
  \leq
  \frac{z''}{r}\left(1 - |\alpha|_E^2 + \frac{z'}{r}\right) \quad\text{on } X.
\end{equation}
If $\deg_\omega E \leq 0$, then the term $z'/r$ in \eqref{eq:Taubes_1996_SW_to_Gr_eq_2-11} can be replaced by $z'/r^2$.
\end{mainthm}

Proposition \ref{mainprop:Taubes_1996_SW_to_Gr_2-1} yields the following analogue of Taubes \cite[Section 2 (b), Proposition 2.4 for $F_A^+$, p. 856]{TauSWGromov} for solutions to the non-Abelian monopole equations \eqref{eq:SO(3)_monopole_equations_almost_Hermitian_perturbed_intro} with a singular Taubes perturbation.

\begin{maincor}[Pointwise estimate for the self-dual, tracefree component of the curvature with singular Taubes perturbations]
\label{maincor:Taubes_1996_SW_to_Gr_eq_2-12_and_2-13}  
Continue the hypotheses and notation of Proposition \ref{mainprop:Taubes_1996_SW_to_Gr_2-1}. If $F_A^+ = F_A^{2,0} + F_A^\omega + F_A^{0,2} \in \Omega^+(\fu(E))$, where $F_A^{2,0}= -(F_A^{0,2})^\dagger$, and then $(F_A^+)_0 \in \Omega^+(\su(E))$ obeys the pointwise inequality,
\begin{equation}
  \label{eq:Taubes_1996_SW_to_Gr_eq_2-12_and_2-13}
  |\rho_{\can}(F_A^+)_0|_{\End(W_{\can}^+\otimes E)}
  \leq
  \frac{(1+2\sqrt{2})}{2\sqrt{2}}\left(\frac{1}{3} + \frac{z}{r}\right) + \frac{r}{4} \quad\text{on } X,
\end{equation}
where $\rho_{\can}:\Lambda^+ \to \End(W_{\can}^+)$ is defined by the canonical \spinc structure on $(X,g)$ as in Definition \ref{defn:Canonical_spinc_bundles}.
\end{maincor}  

In Section \ref{sec:Pointwise_estimate_squared_pointwise_norm_anti-holomorphic_spinor}, we give two proofs of Theorem \ref{mainthm:Taubes_1996_SW_to_Gr_2-3}, as well as a proof of Corollary \ref{maincor:Taubes_1996_SW_to_Gr_eq_2-12_and_2-13}. In Section \ref{sec:Perturbation_non-Abelian_monopole_equations_split_pairs}, we explain that the system \eqref{eq:SO(3)_monopole_equations_almost_Hermitian_perturbed_intro} of non-Abelian monopole equations with a singular Taubes perturbation over an almost Hermitian four-manifold yields, in the case of split pairs $(A,\Phi) = (A_1\oplus A_1, \Phi_1\oplus\Phi_2)$, a solution $(A_1,\Phi_1)$ to the Seiberg--Witten monopole equations with Taubes' perturbation as in \cite{TauSWGromov}.
%PF12-3-2024 Get ref
Consequently, for such split solutions to the system \eqref{eq:SO(3)_monopole_equations_almost_Hermitian_perturbed_intro} of non-Abelian monopole equations with a singular Taubes perturbation, all of Taubes' analysis and estimates in \cite{TauSWGromov} applies.

\begin{mainprop}[Pointwise estimate for  the squared pointwise norms of sections of $E$ with regularized Taubes perturbations]
\label{mainprop:Taubes_1996_SW_to_Gr_2-1_regular}
Continue the hypotheses and notation of Proposition \ref{mainprop:Taubes_1996_SW_to_Gr_2-1} except that we now assume $(A,\varphi,\psi)$ in \eqref{eq:A_varphi_psi_in_W1p} is a solution to the system \eqref{eq:SO(3)_monopole_equations_almost_Hermitian_perturbed_intro_regular} of non-Abelian monopole equations with a regularized Taubes perturbation defined by constants $r > 0$ and $\gamma > 0$. Then $(\alpha,\beta) \in W^{2,p}(E\oplus\Lambda^{0,2}(E))$ obeys the pointwise inequality,
\begin{equation}
  \label{eq:Taubes_1996_SW_to_Gr_eq_2-2_regular}
  |\alpha|_E^2 + |\beta|_{\Lambda^{0,2}(E)}^2
  \leq
  \frac{2}{3} + \frac{z}{r} \quad\text{on } X.
\end{equation}
\end{mainprop}

Note that the bound on the right-hand side of \eqref{eq:Taubes_1996_SW_to_Gr_eq_2-2_regular} is independent of the regularization constant $\gamma > 0$. We give two proofs of Proposition \ref{mainprop:Taubes_1996_SW_to_Gr_2-1_regular} in Section \ref{subsec:Pointwise_estimate_holomorphic_component_coupled_spinor_regular} by describing the small changes needed to the proofs of Proposition \ref{mainprop:Taubes_1996_SW_to_Gr_2-1}.

\begin{mainthm}[Pointwise estimate for the squared pointwise norms of sections of $\Lambda^{0,2}(E)$ with regularized Taubes perturbations]
\label{mainthm:Taubes_1996_SW_to_Gr_2-3_regular}
Continue the hypotheses and notation of Theorem \ref{mainthm:Taubes_1996_SW_to_Gr_2-3} except that we now assume $(A,\varphi,\psi)$ in \eqref{eq:A_varphi_psi_in_W1p} is a solution to the system \eqref{eq:SO(3)_monopole_equations_almost_Hermitian_perturbed_intro_regular} of non-Abelian monopole equations with a regularized Taubes perturbation defined by constants $r > 0$ and $\gamma > 0$. If $\gamma^2 = r\eps^2$ as in the forthcoming \eqref{eq:delta_=_sqrt_r_epsilon} and $\eps^2 = z''/(3r)$ as in the forthcoming \eqref{eq:Definition_epsilon_is_zprimeprime_over_3r}, so $\gamma^2 = z''/3$, then $(\alpha,\beta) \in W^{2,p}(E\oplus\Lambda^{0,2}(E))$ obeys the pointwise inequality,
\begin{equation}
  \label{eq:Taubes_1996_SW_to_Gr_eq_2-11_regular}
  |\beta|_{\Lambda^{0,2}(E)}^2
  \leq
  \frac{z''}{r}\left(1 - |\alpha|_E^2 + \frac{z'}{r}\right) \quad\text{on } X.
\end{equation}
If $\deg_\omega E \leq 0$, then 
\begin{equation}
  \label{eq:Taubes_1996_SW_to_Gr_eq_2-11_regular_degE_leq_0}
  |\beta|_{\Lambda^{0,2}(E)}^2
  \leq
  \frac{z''}{r}\max\left\{1 - |\alpha|_E^2 + \frac{z'}{r^2}, \frac{1}{3}\right\} \quad\text{on } X.
\end{equation}
\end{mainthm}

Proposition \ref{mainprop:Taubes_1996_SW_to_Gr_2-1_regular} yields the following analogue of Taubes \cite[Section 2 (b), Proposition 2.4 for $F_A^+$, p. 856]{TauSWGromov} for solutions to the non-Abelian monopole equations \eqref{eq:SO(3)_monopole_equations_almost_Hermitian_perturbed_intro_regular} with a regularized Taubes perturbation.

\begin{maincor}[Pointwise estimate for the self-dual trace-free component of the curvature with regularized Taubes perturbations]
\label{maincor:Taubes_1996_SW_to_Gr_eq_2-12_and_2-13_regular}  
Continue the hypotheses and notation of Proposition \ref{mainprop:Taubes_1996_SW_to_Gr_2-1_regular}. Then $(F_A^+)_0 \in \Omega^+(\su(E))$ obeys the pointwise inequality,
\begin{equation}
  \label{eq:Taubes_1996_SW_to_Gr_eq_2-12_and_2-13_regular}
  |\rho_{\can}(F_A^+)_0|_{\End(W_{\can}^+\otimes E)}
  \leq
  \frac{(1+2\sqrt{2})}{2\sqrt{2}}\left(\frac{2}{3} + \frac{z}{r}\right) + \frac{r}{2} \quad\text{on } X.
\end{equation}
\end{maincor}  

Theorem \ref{mainthm:Taubes_1996_SW_to_Gr_2-3_regular} and Corollary \ref{maincor:Taubes_1996_SW_to_Gr_eq_2-12_and_2-13_regular} are proved in Section \ref{sec:Pointwise_estimate_squared_pointwise_norm_anti-holomorphic_spinor_regular}. We expect the proofs Propositions \ref{mainprop:Taubes_1996_SW_to_Gr_2-1} and \ref{mainprop:Taubes_1996_SW_to_Gr_2-1_regular} and Theorems \ref{mainthm:Taubes_1996_SW_to_Gr_2-3} and \ref{mainthm:Taubes_1996_SW_to_Gr_2-3_regular} to extend to the case of almost Hermitian four-manifolds and we plan to return to this point in a future article.

\begin{rmk}[Pointwise estimates for $F_A^-$ and covariant derivatives of $\alpha$ and $\beta$]
\label{rmk:Pointwise_estimates_F_Apm_covariant_derivatives_alpha_beta}  
In Taubes \cite[Section 1 (f), Equation (1.24) (3), (4), (5), p. 853]{TauSWGromov}, Taubes lists the pointwise estimates for $F_A^-$ (as well as $F_A^+$) and the covariant derivatives of $\alpha$ and $\beta$, when $(A,\varphi,\psi)$ with $(\varphi,\psi) = r^{1/2}(\alpha,\beta)$ is a solution to the perturbed Seiberg--Witten monopole equations in  \cite[Section 1 (d), Equation (1.20), p. 851]{TauSWGromov}. Because the non-Abelian monopole equations \eqref{eq:SO(3)_monopole_equations_almost_Hermitian_perturbed_intro} with a singular Taubes perturbation reduce to the perturbed Seiberg--Witten monopole equations when the triple $(A,\varphi,\psi)$ is \emph{split} (see Section \ref{sec:Perturbation_non-Abelian_monopole_equations_split_pairs}), the pointwise estimates in \cite[Section 1 (f), Equation (1.24) (3), (4), (5), p. 853]{TauSWGromov} should extend with almost no change to the case of split solutions to \eqref{eq:SO(3)_monopole_equations_almost_Hermitian_perturbed_intro} and we expect that they will similarly extend to split solutions of the non-Abelian monopole equations \eqref{eq:SO(3)_monopole_equations_almost_Hermitian_perturbed_intro_regular} with a regularized Taubes perturbation.
\qed\end{rmk}

\subsection[Almost Hermitian structures on virtual moduli spaces of non-Abelian monopoles]{Almost Hermitian structures on virtual moduli spaces of non-Abelian monopoles over closed symplectic four-manifolds}
\label{subsec:Almost_Hermitian_structure_moduli_spaces_non-Abelian_monopoles_symplectic_4-manifolds}
When $X$ is a \emph{symplectic} manifold, we write $(X,\omega)$ when we wish to emphasize the choice of symplectic two-form, $\omega$. When $X$ is an \emph{almost Hermitian} or \emph{almost K\"ahler manifold}, we write $(X,g,J,\omega)$ when we wish to emphasize the choice of compatible triple of a Riemannian metric $g$ and almost complex structure $J$ for the fundamental or symplectic two-form $\omega$, respectively.

The almost complex structure $J \in \End(TX)$ extends by complex linearity to $J \in \End(TX\otimes_\RR\CC)$ and induces the direct sum decomposition of complex vector bundles, $TX\otimes_\RR\CC = TX^{1,0}\oplus TX^{0,1}$, where
\[
  J = -i \text{ on } TX^{1,0} \quad\text{and}\quad J = i \text{ on } TX^{0,1}.
\]
Our sign convention that $TX^{1,0}$ is the $-i$ eigenspace for $J$ and $TX^{0,1}$ is the $i$ eigenspace for $J$ is the \emph{reverse} of that in, for example, Huybrechts \cite[Section 1.2, Definition 1.2.4, p. 25]{Huybrechts_2005}, Kobayashi and Nomizu \cite[Chapter IX, Section 2, p. 125]{Kobayashi_Nomizu_v2}, and Wells \cite[Chapter I, Section 3, p. 29]{Wells3}, but it is more convenient for our application. As in Huybrechts \cite[Section 1.2, Lemma 1.2.6, p. 26]{Huybrechts_2005}, there is a canonically induced almost complex structure $J \in \End(T^*X)$ defined by $(J\theta)(v) := \theta(Jv)$ for all $v \in C^\infty(TX)$ and $\theta \in C^\infty(T^*X)$ such that $(T^*X)^{1,0} = (TX^{1,0})^*$ and $(T^*X)^{0,1} = (TX^{0,1})^*$, and thus
\[
  J = -i \text{ on } T^*X^{1,0} \quad\text{and}\quad J = i \text{ on } T^*X^{0,1}.
\]
Because $\Omega^{1,0}(\fsl(E)) = C^\infty(T^*X^{1,0}\otimes\fsl(E))$ and $\Omega^{0,1}(\fsl(E)) = C^\infty(T^*X^{0,1}\otimes\fsl(E))$, we obtain
\begin{equation}
  \label{eq:Almost_complex_structure_omega_pq_fslE}
  J = -i \text{ on } \Omega^{1,0}(\fsl(E)) \quad\text{and}\quad J = i \text{ on } \Omega^{0,1}(\fsl(E)).
\end{equation}
The inverse of the canonical isomorphism of real vector spaces,
\begin{equation}
  \label{eq:Omega01_slE_isomorphic_Omega1_suE}
  \Omega^{0,1}(\fsl(E)) \ni a'' \mapsto a = \frac{1}{2}\left(-(a'')^\dagger + a''\right) \in \Omega^1(\su(E)),
\end{equation}
pulls back the almost complex structure $J$ on $\Omega^{1,0}(\fsl(E))$ to an almost complex structure $J$ on $\Omega^1(\su(E))$:
\[
  J a
  = \frac{1}{2}\left(-J(a'')^\dagger + J a''\right)
  = \frac{1}{2}\left(i(a'')^\dagger + ia''\right)
  = \frac{i}{2}\left((a'')^\dagger + a''\right),
  \quad\text{for all } a \in \Omega^1(\su(E)).
\]  
Our sign convention for $J$ agrees with that of Kobayashi \cite[Equation (7.6.12), p. 252]{Kobayashi_differential_geometry_complex_vector_bundles} for the action of $J$ on $\Omega^1(\su(E))$, but is the \emph{reverse} of that of Itoh \cite[Equation (4.1), p. 19]{Itoh_1988}. The complex vector bundle $E\oplus \Lambda^{0,2}(E)$ has its standard almost complex structure defined by scalar multiplication by $i$, where we abuse notation by abbreviating $\Lambda^{p,q}(E) = \Lambda^{p,q}(X)\otimes E$. We thus have canonically induced almost complex structures $J$ on the vector spaces,
\begin{subequations}
  \label{eq:sEkC}
  \begin{align}
    \label{eq:sE1C}
    \sF_1 &:= C^\infty\left( \Lambda^{0,1}(\fsl(E)) \oplus E\oplus \Lambda^{0,2}(E)\right),
    \\
    \label{eq:sE2C}
    \sF_2 &:= C^\infty\left(\fsl(E) \oplus \Lambda^{0,2}(\fsl(E)) \oplus \Lambda^{0,1}(E)\right),
\end{align}
\end{subequations}
given by scalar multiplication by $i$. By analogy with \eqref{eq:Omega01_slE_isomorphic_Omega1_suE}, there is an isomorphism
\begin{equation}
  \label{eq:Omega0_slE_isomorphic_Omega0_suE_oplus_suE}
  \Omega^0(\fsl(E)) \ni \zeta \mapsto (\xi_1,\xi_2)
  =  \frac{1}{2}\left(-\zeta^\dagger + \zeta, i(\zeta^\dagger + \zeta)\right)
  \in \Omega^0\left(\su(E) \oplus \su(E)\right).
\end{equation}
We thus obtain canonically induced almost complex structures $J$ on the vector spaces,
\begin{subequations}
  \label{eq:sEk}
  \begin{align}
    \label{eq:sE1}
    \sE_1 &:= C^\infty\left(T^*X\otimes\su(E) \oplus E \oplus \Lambda^{0,2}(E)\right),
    \\
    \label{eq:sE2}
    \sE_2 &:= C^\infty\left(\su(E) \oplus \su(E) \oplus \Lambda^{0,2}(\fsl(E)) \oplus \Lambda^{0,1}(E)\right),
\end{align}
\end{subequations}
by pullback via the inverses of the canonical isomorphisms of real vector spaces induced by the isomorphisms \eqref{eq:Omega01_slE_isomorphic_Omega1_suE} and \eqref{eq:Omega0_slE_isomorphic_Omega0_suE_oplus_suE}:
\begin{subequations}
  \label{eq:Isomorphism_sEkC_to_sEk}
  \begin{align}
    \label{eq:Isomorphism_sE1C_to_sE1}
    \Upsilon_1:\sF_1 \ni (a'',\sigma,\tau) \mapsto
    (a,\sigma,\tau) &= \left(\frac{1}{2}(a''-(a'')^\dagger), \sigma,\tau\right) \in \sE_1,
    \\
    \label{eq:Isomorphism_sE2C_to_sE2}
    \Upsilon_2:\sF_2 \ni (\zeta,v,\nu) \mapsto (\xi_1,\xi_2,v/2,\nu) &= \left(\frac{1}{2}(\zeta-\zeta^\dagger), \frac{i}{2}(\zeta+\zeta^\dagger), \frac{1}{2}v, \nu\right) \in \sE_2.
  \end{align}
\end{subequations}
%PF4-21-2025 Our sign convention for \xi_2 differs from Tom's later usage?
%TL4-29-2025: Will check
Note that $\zeta = \xi_1 - i\xi_2 $ and $J\zeta = i\zeta$ for all $\zeta \in \Omega^0(\fsl(E))$, while the isomorphism \eqref{eq:Omega0_slE_isomorphic_Omega0_suE_oplus_suE} gives
\[
  J(\xi_1,\xi_2) = (\xi_2,-\xi_1),
  \quad\text{for all } (\xi_1,\xi_2) \in \Omega^0\left(\su(E) \oplus \su(E)\right),
\]
and thus
\[
  J = \begin{pmatrix}0 & 1 \\ -1 & 0 \end{pmatrix} \quad\text{on } \Omega^0\left(\su(E) \oplus \su(E)\right).
\]
We shall prove the

\begin{mainthm}[Almost Hermitian structures on bounded eigenvalue spaces defined by the non-Abelian monopole equations over closed symplectic four-manifolds]
\label{mainthm:AH_structure_bounded_evalue_spaces_non-Abelian_monopoles_symp_4-mflds}
Let $(X,g,J,\omega)$ be a closed, almost K\"ahler four-manifold with fundamental two-form $\omega = g(J\cdot,\cdot)$ as in \eqref{eq:Fundamental_two-form} and let $(E,H)$ be a smooth, Hermitian vector bundle over $X$ with complex rank
%PF9-8-2024 Do we need E having rank 2 here?
two and smooth, unitary connection $A_d$ on $\det E$. If $[A,\varphi,\psi]$ is a point in the moduli subspace $\sM^0(E,g,J,\omega,r)$ of non-zero-section solutions to the non-Abelian monopole equations \eqref{eq:SO(3)_monopole_equations_almost_Hermitian_perturbed_intro_regular} with a regularized Taubes perturbation, then there is a positive constant $\nu = \nu[A,\varphi,\psi] \in [1,\infty)$ such that the following hold:
\begin{enumerate}
%PF7-18-2025 Move tildes inside bold font via redefined macros
%TL8-30-2025 Replaced \tilde\bH with \mathbf{\tilde H}
\item\label{item:AH_structure_tilde_bH_A_varphi_psi_r_nu_1}
There is an almost Hermitian structure $(\bg_{1,A,\varphi,\psi},\bJ_{1,A,\varphi,\psi},\bomega_{1,A,\varphi,\psi})$ on the finite-dimensional real vector space 
$\mathbf{\tilde H}_{A,\varphi,\psi,r,\nu}^1$ 
defined in \eqref{eq:tilde_bH_A_varphi_psi_r_nu_1}, where
\[
  \bH_{A,\varphi,\psi,r}^1 \subset \mathbf{\tilde H}_{A,\varphi,\psi,r,\nu}^1 \subset \sE_1
\]
and $\bH_{A,\varphi,\psi,r}^1$ is defined in \eqref{eq:H1_For_TaubesPert_regular} and represents the Zariski tangent space to $\sM^0(E,g,J,\omega,r)$ at $[A,\varphi,\psi]$ and $(A,\varphi,\psi)$ is a smooth representative of the point $[A,\varphi,\psi]$.
  
\item\label{item:AH_structure_tilde_bH_A_varphi_psi_r_nu_2}
There is an almost Hermitian structure $(\bg_{2,A,\varphi,\psi},\bJ_{2,A,\varphi,\psi},\bomega_{2,A,\varphi,\psi})$ on the finite-dimensional real vector space $\mathbf{\tilde H}_{A,\varphi,\psi,r,\nu}^2$ defined in \eqref{eq:tilde_bH_A_varphi_psi_r_nu_2}, where
\[
  \bH_{A,\varphi,\psi,r}^2 \subset \mathbf{\tilde H}_{A,\varphi,\psi,r,\nu}^2 \subset \sE_2
\]
and $\bH_{A,\varphi,\psi,r}^2$ is defined in \eqref{eq:DefineH2ForSO3Monopoles_Perturbed} and represents the obstruction space for $\sM^0(E,g,J,\omega,r)$ at the point $[A,\varphi,\psi]$.
\end{enumerate}
\end{mainthm}

%PF11-13-2025 The remark is wrong. The \Upsilon_k maps do not yield a similarity transformation between the Laplacians or Dirac-type operators.
% \begin{rmk}
% \label{rmk:AH_structure_bounded_d-bar_A_evalue_spaces_non-Abelian_monopoles_symp_4-mflds}
% The equivalence \eqref{eq:Equivalence_cT_d-bar_A_varphi_psi_r_and_sT_A_varphi_psi_r} between the operators $\sT_{A,\varphi,\psi,r}$ and $\cT_{\bar\partial_A,\varphi,\psi,r}$ ensures that Theorem \ref{mainthm:AH_structure_bounded_evalue_spaces_non-Abelian_monopoles_symp_4-mflds} yields almost Hermitian structures on the analogously defined
% %PF7-18-2025 Add definitions and update notation in chapter 11 to match
% real vector spaces $\mathbf{\tilde H}_{\bar\partial_A,\varphi,\psi,r,\nu}^1$ and $\mathbf{\tilde H}_{\bar\partial_A,\varphi,\psi,r,\nu}^2$ that contain $\bH_{\partial_A,\varphi,\psi,r}^1$ and $\bH_{\partial_A,\varphi,\psi,r}^2$ in \eqref{eq:DefineH1ForApproxComplex_Perturbed} and \eqref{eq:DefineH2ForApproxComplex_Perturbed}, respectively.
% \end{rmk}

Theorem \ref{mainthm:AH_structure_bounded_evalue_spaces_non-Abelian_monopoles_symp_4-mflds} leads to the

%PF8-31-2025 Surely this only applies to neighborhoods of fixed points?
\begin{maincor}[Almost Hermitian structures on circle invariant virtual moduli spaces of non-Abelian monopoles over closed symplectic four-manifolds]
\label{maincor:Almost_Hermitian_structure_moduli_space_non-Abelian_monopoles_symplectic_4-manifolds}
Continue the hypotheses and notation of Theorem \ref{mainthm:AH_structure_bounded_evalue_spaces_non-Abelian_monopoles_symp_4-mflds}. For the circle action \eqref{eq:S1ZActionOnQuotientSpace} on $\sC^0(E,H,J,A_d)$, there is a circle invariant almost Hermitian structure $(\bg_1,\bJ_1,\bomega_1)$ on the circle invariant finite-dimensional ambient smooth manifold
% PF10-12-2025 Add explicit proof of its construction
$\sM^\vir(E,g,J,r,\omega,\nu)$ defined by the Kuranishi model constructed in Section \ref{subsec:Local_Kuranishi_model_nonlinear_map_Banach_spaces}
%PF7-18-2025 Add more precise reference and new definition
%TL8-5-2025: Kuranishi models are local (around points) while the notation suggests a global construction.  I think you had a reference in mind for such a global construction?
%PF8-31-2025 It's supposed to be local
for an open neighborhood in $\sM^0(E,g,J,\omega,r)$ of an
%PF10-12-2025 Added qualifier
$S^1$ fixed point
$[A,\varphi,\psi]$, where the circle action is Hamiltonian with respect to Hitchin's function $f$ in the forthcoming \eqref{eq:Hitchin_function}.
\end{maincor}

Theorem \ref{mainthm:AH_structure_bounded_evalue_spaces_non-Abelian_monopoles_symp_4-mflds} and Corollary \ref{maincor:Almost_Hermitian_structure_moduli_space_non-Abelian_monopoles_symplectic_4-manifolds} are proved in Section \ref{sec:Hitchin_function_Hamiltonian_circle_action_virtual_moduli_spaces}.

\subsection[Critical points of Hitchin's function on spaces of non-Abelian monopoles]{Critical points of Hitchin's function on moduli spaces of non-Abelian monopoles over closed symplectic four-manifolds}
\label{subsec:Critical_points_Hamiltonian_moduli_spaces_non-Abelian_monopoles_closed_symplectic_4-manifolds}
We recall from Feehan and Leness
%TL11-26-2025: Updated
\cite[Lemma 6.3.3, Equation (6.3.9), and Theorem 12.3.10]{Feehan_Leness_introduction_virtual_morse_theory_so3_monopoles} --- see the discussion in the last paragraph of
%TL11-26-2025: Updated to \cite[Section 6.3.2]{Feehan_Leness_introduction_virtual_morse_theory_so3_monopoles}
\cite[Section 6.3.2]{Feehan_Leness_introduction_virtual_morse_theory_so3_monopoles} --- that, when $E$ has rank two, the quotient space $\sC^0(E,H,J,A_d)$ of non-zero-section unitary triples as in \eqref{eq:Quotient_space_non-zero-section_unitary_triples} is a real analytic Banach manifold. We define the following analogues of Hitchin's function and its critical points in \cite[Section 7, p. 92]{Hitchin_1987}.

\begin{defn}[Hitchin's function and its critical points on the quotient space of non-zero-section unitary triples and moduli subspace of non-Abelian monopoles]
\label{defn:Critical_point_Hitchin_function_moduli_space_non-Abelian_monopoles}
Let $(X,g,J,\omega)$ be a closed almost Hermitian four-manifold with fundamental two-form $\omega = g(J\cdot,\cdot)$ as in \eqref{eq:Fundamental_two-form} and let $(E,H)$ be a smooth Hermitian vector bundle over $X$ with complex rank
%PF9-8-2024 Do we need E having rank 2 here?
two and smooth, unitary connection $A_d$ on $\det E$. We define the smooth \emph{Hitchin function} on the affine space of $W^{1,p}$ unitary triples \eqref{eq:A_varphi_psi_in_W1p} by
\begin{equation}
\label{eq:Hitchin_function_affine}
    f: \sA^{1,p}(E,H,A_d) \times W^{1,p}(E\oplus\Lambda^{0,2}(E)) \ni (A,\varphi,\psi)
  \mapsto \frac{1}{2}\|(\varphi,\psi)\|_{L^2(X)}^2 \in \RR,
\end{equation}
with differential at $(A,\varphi,\psi)$ given by (see Section  \ref{subsec:Hamiltonian_function_circle-invariant_non-degenerate_2-form_virtual_moduli_space})
\begin{equation}
  \label{eq:Hitchin_function_affine_differential}
  (df)_{A,\varphi,\psi}: W^{1,p}\left(T^*X\otimes\su(E) \oplus E\oplus\Lambda^{0,2}(E)\right)
  \ni (a,\sigma,\tau)
  \mapsto
  \Real\left((0,\varphi,\psi),(a,\sigma,\tau)\right)_{L^2(X)} \in \RR.
\end{equation}
The function $f$ is $W^{2,p}(\SU(E))$-invariant and descends to a smooth function on the quotient space,
\begin{equation}
    \label{eq:Hitchin_function}
    f: \sC^0(E,H,J,A_d) \ni [A,\varphi,\psi]
  \mapsto \frac{1}{2}\|(\varphi,\psi)\|_{L^2(X)}^2 \in \RR,
\end{equation}
and restricts to a continuous function on the real analytic moduli subspace $\sM^0(E,g,J,\omega,r)$ that is smooth on smooth strata,
\begin{equation}
    \label{eq:Hitchin_function_moduli_space}
    f: \sM^0(E,g,J,\omega,r) \ni [A,\varphi,\psi]
  \mapsto \frac{1}{2}\|(\varphi,\psi)\|_{L^2(X)}^2 \in \RR.
\end{equation}
We say that $[A,\varphi,\psi] \in \sC^0(E,H,J,A_d)$ is a \emph{critical point} of $f$ in \eqref{eq:Hitchin_function} if
\begin{equation}
  \label{eq:Critical_point_Hitchin_Hamiltonian_function_quotient_space}
  T_{A,\varphi,\psi}\sC^0(E,H,J,A_d) = \Ker (df)_{A,\varphi,\psi},
\end{equation}
where $(A,\varphi,\psi)$ is a representative of the point $[A,\varphi,\psi]$ and for $d_{A,\varphi,\psi}^{0,*}$ as in \eqref{eq:nonAbelianMonopole_d0*_unitary}.
\[
  T_{A,\varphi,\psi}\sC^0(E,H,J,A_d)
  =
  \Ker d_{A,\varphi,\psi}^{0,*}\cap W^{1,p}\left(T^*X\otimes\su(E) \oplus E\oplus\Lambda^{0,2}(E)\right)
\]
is a representative of the tangent space to $\sC^0(E,H,J,A_d)$ at $[A,\varphi,\psi]$ (see Remark \ref{rmk:Representation_tangent_spaces_quotient_spaces}). We say that $[A,\varphi,\psi] \in \sM^0(E,g,J,\omega,r)$ is a \emph{critical point} of $f$ in \eqref{eq:Hitchin_function_moduli_space} if
\begin{equation}
  \label{eq:Critical_point_Hitchin_Hamiltonian_function_moduli_non-Abelian_monopoles}
  T_{A,\varphi,\psi}\sM^0(E,g,J,\omega,r) \subset \Ker (df)_{A,\varphi,\psi},
\end{equation}
where $(A,\varphi,\psi)$ is a smooth representative of the point $[A,\varphi,\psi]$ and, for $d_{A,\varphi,\psi,r}^1$ as in \eqref{eq:SO3Monopoled1TaubesPerturbation} and $\bH_{A,\varphi,\psi,r}^1$ as in \eqref{eq:H1_For_TaubesPert_regular},
\begin{multline}
  \label{eq:bH1_Avarphipsi_r}
  T_{A,\varphi,\psi}\sM^0(E,g,J,\omega,r)
  =
  \bH_{A,\varphi,\psi,r}^1
  \\
  =
  \Ker\left(d_{A,\varphi,\psi,r}^1 + d_{A,\varphi,\psi}^{0,*}\right)
  \cap W^{1,p}\left(T^*X\otimes\su(E) \oplus E\oplus\Lambda^{0,2}(E)\right)
\end{multline}
is a representative of the \emph{Zariski tangent space} to $\sM^0(E,g,J,\omega,r)$ at $[A,\varphi,\psi]$ (see Remark \ref{rmk:Representation_tangent_spaces_quotient_spaces}).
\qed\end{defn}

In our application, we shall need the following extension of our Definition \ref{defn:Critical_point_Hitchin_function_moduli_space_non-Abelian_monopoles} of a critical point of Hitchin's function on the moduli space of non-Abelian monopoles.

\begin{defn}[Hitchin's function and $V$-critical points on the moduli space of non-Abelian monopoles]
\label{defn:V_critical_point_Hitchin_function_moduli_space_non-Abelian_monopoles}  
Continue the notation of Definition \ref{defn:Critical_point_Hitchin_function_moduli_space_non-Abelian_monopoles}. Let $V$ be a real vector space such that
\[
  \bH_{A,\varphi,\psi,r}^1 \subseteq V \subseteq
  \Ker d_{A,\varphi,\psi}^{0,*}\cap W^{1,p}\left(T^*X\otimes\su(E) \oplus E\oplus\Lambda^{0,2}(E)\right).
\]
We say that $f:\sM^0(E,g,J,\omega,r)\to\RR$ in \eqref{eq:Hitchin_function_moduli_space} has a \emph{$V$-critical point} at $[A,\varphi,\psi] \in \sM^0(E,g,J,\omega,r)$ if
\[
 (df)_{A,\varphi,\psi} = 0 \quad\text{on } V,
\]
where $(A,\varphi,\psi)$ is a smooth representative of $[A,\varphi,\psi]$. If $V = \bH_{A,\varphi,\psi,r}^1$, then a $V$-critical point is a critical point of $f:\sM^0(E,g,J,\omega,r)\to\RR$ in the sense of Definition \ref{defn:Critical_point_Hitchin_function_moduli_space_non-Abelian_monopoles}. \qed
\end{defn}

We have the following generalization of Feehan and Leness 
%TL11-26-2025: Update to \cite[Section 2.3, Theorem 8]{Feehan_Leness_introduction_virtual_morse_theory_so3_monopoles}?
\cite[Section 1.4.3, Theorem 4]{Feehan_Leness_introduction_virtual_morse_theory_so3_monopoles},
which applied only to closed, complex K\"ahler surfaces.

\begin{mainthm}[Equivalence of $\mathbf{\tilde H}_{A,\varphi,\psi,r,\nu}^1$-critical points of Hitchin's function and points in moduli subspaces of Seiberg--Witten monopoles]
\label{mainthm:IdentifyCriticalPoints}
Let $(X,g,J,\omega)$ be a closed almost K\"ahler four-manifold and $(E,H)$ be a smooth, Hermitian vector bundle over $X$ with complex rank
%PF9-8-2024 Do we need E having rank 2 here?
two and smooth, unitary connection $A_d$ on the Hermitian line bundle $\det E$. If $[A,\varphi,\psi] \in \sM^0(E,g,J,\omega,r)$ as in \eqref{eq:Moduli_space_non-Abelian_monopoles_almost_Hermitian_Taubes_regularized_non-zero-section} and $\nu = \nu[A,\varphi,\psi]$ is the large enough positive constant provided by Corollary \ref{maincor:Almost_Hermitian_structure_moduli_space_non-Abelian_monopoles_symplectic_4-manifolds}, then the following hold:
\begin{enumerate}
\item \emph{(Seiberg--Witten solution $\implies$ critical point)}
\label{item:SW_solution_implies_critical_point}
If $(A,\varphi,\psi)$ is split in the sense of Definition \ref{defn:Split_trivial_central-stabilizer_spinor_pair} \eqref{item:Split_spinor_pair} with respect to a splitting $E=L_1\oplus L_2$ as a direct sum of Hermitian line bundles and a solution $(A_1,\varphi_1,\psi_1)$ on $L_1$ to the Seiberg--Witten monopole equations \eqref{eq:SW_monopole_equations_regular_Taubes_perturbation} with a regularized Taubes perturbation, where $A = A_1\oplus A_d\otimes A_1^*$ and $\varphi = (\varphi_1,0)$ and $\psi = (\psi_1,0)$, then $[A,\varphi,\psi]$ is a $\mathbf{\tilde H}_{A,\varphi,\psi,r,\nu}^1$-critical point (in the sense of Definition \ref{defn:V_critical_point_Hitchin_function_moduli_space_non-Abelian_monopoles}) of the Hitchin function $f:\sM^0(E,g,J,\omega,r)\to\RR$ given by the restriction of $f:\sC^0(E,H,J,A_d)\to\RR$ in \eqref{eq:Hitchin_function}.

\item \emph{(Critical point $\implies$ Seiberg--Witten solution)}
\label{item:Critical_point_implies_SW_solution}
If $[A,\varphi,\psi]$ is a $\mathbf{\tilde H}_{A,\varphi,\psi,r,\nu}^1$-critical point (in the sense of Definition \ref{defn:V_critical_point_Hitchin_function_moduli_space_non-Abelian_monopoles}) of the Hitchin function $f:\sM^0(E,g,J,\omega,r)\to\RR$ given by the restriction of $f:\sC^0(E,H,J,A_d)\to\RR$ in \eqref{eq:Hitchin_function}, then $(A,\varphi,\psi)$ is split in the sense of Definition \ref{defn:Split_trivial_central-stabilizer_spinor_pair} \eqref{item:Split_spinor_pair} with respect to a splitting $E=L_1\oplus L_2$ as a direct sum of Hermitian line bundles and a solution $(A_1,\varphi_1,\psi_1)$ on $L_1$ to the Seiberg--Witten monopole equations \eqref{eq:SW_monopole_equations_regular_Taubes_perturbation} with a regularized Taubes perturbation, where $A = A_1\oplus A_d\otimes A_1^*$ and $\varphi = (\varphi_1,0)$ and $\psi = (\psi_1,0)$.
\end{enumerate}
\end{mainthm}

We prove Theorem \ref{mainthm:IdentifyCriticalPoints} in Section \ref{sec:Equivalence_critical_points_Hamiltonian_fixed_points_S1_action_moduli_space}.

\begin{rmk}[Transversality for the moduli space of non-Abelian monopoles with a regularized Taubes perturbation and generic geometric parameters]
\label{rmk:Transversality_non-Abelian_monopoles_regularized_Taubes_perturbation_generic_geometric_parameters}  
According to Theorem \ref{thm:Transv}, the moduli space $\sM_\ft^{*,0}(g,\tau,\vartheta)$ of non-split, non-zero-section solutions to the system \eqref{eq:PerturbedSO3MonopoleEquations} of non-Abelian monopole equations with generic geometric parameters $(g,\tau,\vartheta) \in \Met(X)\times\GL(\Lambda^+(X))\times \Omega^1(X,\CC)$ is a smooth manifold of the expected dimension. A similar result should hold for a moduli space $\sM^{0,*}(E,g,J,\omega,r,\tau,\vartheta)$ of non-split, non-zero-section solutions to the system \eqref{eq:SO(3)_monopole_equations_almost_Hermitian_perturbed_intro_regular} of non-Abelian monopole equations with a regularized Taubes perturbation \emph{and} generic geometric parameters $g$, $\tau$ and $\vartheta$. We shall consider the proof of this result elsewhere.
\qed\end{rmk}

\begin{rmk}[Smallness of the Nijenhuis tensor]
\label{rmk:Replacing_NJ_C0_small_by_small_Nijenhuis_energy}
A condition on the Nijenhuis tensor $N_J$, defined in \eqref{eq:Nijenhuis_tensor}, such as
\begin{equation}
\label{eq:C0_small_N_J}
  \|N_J\|_{C^0(X,g)} < \eps,
 \end{equation}
would be convenient for the purpose of showing that the non-Abelian monopole equations over almost Hermitian four-manifolds are approximately holomorphic, but appears difficult to achieve --- see Fernandez, Shin, and Wilson \cite{Fernandez_Shin_Wilson_2021} and Shin \cite{Shin_2021} for a brief history of attempts to prove existence of almost complex structures with $C^0(X,g)$-small $N_J$. In \cite{Fernandez_Shin_Wilson_2021}, the authors provide examples of closed, almost complex four-manifolds that do not admit complex structures but which do admit almost complex structures $J$ with arbitrarily $C^0(X,g)$-small Nijenhuis tensor $N_J$ with respect to a fixed Riemannian metric $g$. Examples of closed, almost complex six-manifolds that do not necessarily admit complex structures but which do admit almost complex structures $J$ with arbitrarily $C^0(X,g)$-small Nijenhuis tensor $N_J$ are provided by Fernandez, Shin, and Wilson \cite{Fernandez_Shin_Wilson_2021} and Fei, Phong, Picard, and Zhang \cite[Section 9]{Fei_Phong_Picard_Zhang_2021cjm}.

A condition that $N_J$ be $L^1(X,g)$ small appears as a hypothesis in Weinkove \cite[Section 1, Theorem 2, p. 319]{Weinkove_2007}, one of his results related to existence of solutions to the Calabi--Yau equation on almost K\"ahler four-manifolds. There is an intriguing result due to Evans that might be used as a replacement for \eqref{eq:C0_small_N_J}, although it would entail a more challenging analysis than the stronger condition \eqref{eq:C0_small_N_J}. Given a closed symplectic manifold $(X,\omega)$ of dimension $2n$ with $[\omega] \in H^2(X;\QQ)$, Evans applies Donaldson's Symplectic Submanifold Theorem \ref{thm:Donaldson_1} to prove that there exists an almost complex structure $J \in \End(TX)$ that is compatible with $\omega$ and such that its \emph{Nijenhuis energy} is small \cite[Theorem 1.1, p. 384]{Evans_2012},
\begin{equation}
  \label{eq:2}
  \int_X |N_J|^2 \omega^n < \eps,
\end{equation}
where the pointwise norm $|N_J|$ is computed with respect to the induced Riemannian metric $g_{J,\omega} := \omega(\cdot, J\cdot)$. Naturally, the dependence of $g_{J,\omega}$ on $J$ makes Evans' result difficult to use.

We recall that $J$ is \emph{compatible} with a symplectic form $\omega$ if for all $p\in X$, one has $\Omega(Jv,Jw) = \Omega(v,w)$, for all $v,w \in T_pX$, and $\Omega(v,Jv) > 0$ for all $v \in T_pX \less \{0\}$ (see McDuff and Salamon \cite[Equations (4.1.1) and (4.1.2), p. 153]{McDuffSalamonSympTop3}).
\qed\end{rmk}

\subsection{Virtual Morse--Bott signatures for Seiberg--Witten critical points}
\label{subsec:Virtual_Morse-Bott_signatures_SW_critical_points}
Finally, we compute the virtual Morse--Bott index \eqref{eq:Virtual_Morse-Bott_signature} for Hitchin's function \eqref{eq:Hitchin_function} at a non-zero-section, split non-Abelian monopole (see the forthcoming Definition \ref{defn:Split_trivial_central-stabilizer_spinor_pair}). 

% TL7-9-2025: This paragraph to be rewritten when Kuranishi model presented in terms of bounded-eigenvalue eigenspaces
%PF7-18-2025 Please check that you have completed this rewriting
Consider a point in $\sM^0(E,g,J,\omega,r)$ represented by a smooth non-zero-section, split non-Abelian monopole $(A,\varphi,\psi)$ and let $\bH_{A,\varphi,\psi,r}^1$ and $\bH_{A,\varphi,\psi,r}^2$ be the harmonic spaces given in \eqref{eq:H1_For_TaubesPert_regular} and \eqref{eq:DefineH2ForSO3Monopoles_Perturbed}. We define explicit $S^1$ actions on these spaces in Section \ref{sec:S1ActionOnApproxComplexDefOperator} that are equivalent, up to positive multiplicity, to the $S^1$ action in Corollary \ref{maincor:Almost_Hermitian_structure_moduli_space_non-Abelian_monopoles_symplectic_4-manifolds}.
% PF7-18-2025 Updated references
%PF7-18-2025 These K models need to be consolidated.
We note that Corollary \ref{cor:Kuranishi_model_defined_by_Fredholm_map_Hilbert_spaces},
Remark \ref{rmk:Hypotheses_corollary_Kuranishi_model_defined_by_Fredholm_map_Hilbert_spaces},
%Proposition \ref{prop:KuranishiModel_For_SplitPerturbedPair},
and Lemma \ref{lem:SymplecticStrutureOnStabilizedSpace} and Proposition \ref{prop:S1EquivKuranishiModel1} collectively provide an $S^1$-equivariant Kuranishi model for an $S^1$-invariant open neighborhood of $[A,\varphi,\psi]$ in $\sM^0(E,g,J,\omega,r)$ given by a real analytic, $S^1$-equivariant map on an open neighborhood $U$ of the origin in $\bH_{A,\varphi,\psi,r}^1$,
\[
  \bchi:\mathbf{\tilde H}_{A,\varphi,\psi,r,\nu}^1  \supset U \to \mathbf{\tilde H}_{A,\varphi,\psi,r,\nu}^2,
\]
and an $S^1$-equivariant homeomorphism $\bga$ between a neighborhood of the origin in $\bga^{-1}(0)$ and an $S^1$-invariant open neighborhood of $[A,\varphi,\psi]$ in $\sM^0(E,g,J,\omega,r)$. The $S^1$-equivariant isomorphisms of Lemma \ref{lem:S1EquivariantIsom_LowEigenvalueSpaces_DefOp_to_ApproxComplex} and the almost complex structures given in
%PF7-18-2025 Replace by references to chapter 10
Corollary \ref{cor:Isom_Between_BoundedEigenvalueEigenspaces}
%PF7-18-2025 The Upsilon isomorphisms simply map eigenspaces of \sT' to those of \cT' and vice versa. The corollary does not make this clear.
%PF7-18-2025 Please consolidate the notation in chapter 11: we just need \bH^k spaces, not H^k(T), etc.
and Definition \ref{defn:S1_Equiv_Kuranishi_Model_At_Split_Triple_With_AC_Structure}
%PF7-18-2025 This is NOT how these almost complex structures are defined. Please read Chapter 10 and update.
define $S^1$-invariant almost complex structures on the spaces $\mathbf{\tilde H}_{A,\varphi,\psi,r,\mu}^k$ for $k=1,2$. These almost complex structures determine the signs of the weights of the $S^1$ action as discussed in
%TL11-26-2025: Updated to \cite[Definition 4.3.3 and Remark 4.3.4]{Feehan_Leness_introduction_virtual_morse_theory_so3_monopoles}
\cite[Definition 4.3.3 and Remark 4.3.4]{Feehan_Leness_introduction_virtual_morse_theory_so3_monopoles}. For $k=1,2$, we let $\mathbf{\tilde H}_{A,\varphi,\psi,r,\mu}^{-,k}$ denote the subspace of $\bH_{A,\varphi,\psi,r,\balpha,\mu}^k$ on which $S^1$ acts with negative weight with respect to these almost complex structures (in the sense of 
%TL12-4-2025: Updated
\cite[Definition 4.3.3]{Feehan_Leness_introduction_virtual_morse_theory_so3_monopoles}). In the proof of Theorem \ref{mainthm:MorseIndexAtReduciblesOnAlmostKahler}, we will show that the virtual Morse--Bott index of the Hitchin function $f$ in \eqref{eq:Hitchin_function} at $[A,\Phi]$, as discussed in Section \ref{sec:Virtual_Morse-Bott_signature_Hamiltonian_function_circle_action_complex_analytic_space}, is given by
\begin{equation}
  \label{eq:vBM_Perturbed_Intro}
  \lambda_{[A,\varphi,\psi]}^-(f)
  :=
  \dim_\RR \mathbf{\tilde H}_{A,\varphi,\psi,r,\mu}^{-,1} - \dim_\RR\mathbf{\tilde H}_{A,\varphi,\psi,r,\mu}^{-,2}.
\end{equation}
Similarly, we define the virtual Morse--Bott co-index and nullity of $f$ at
%PF7-18-2025 It's not $[A,\Phi]$, it's $[A,\varphi,\psi]$. Please proofread what you write repeatedly.
%$[A,\Phi]$ by
$[A,\varphi,\psi]$
\begin{subequations}
\label{eq:vBM_Perturbed_coIndNull}
 \begin{align}
  \label{eq:vBM_Perturbed_coInd}
 \lambda_{[A,\varphi,\psi]}^+(f)&
                                  :=
                                  %PF7-10-2025 Interchange "+" and k?
  \dim_\RR \mathbf{\tilde H}_{A,\varphi,\psi,r,\mu}^{+,1} - \dim_\RR\mathbf{\tilde H}_{A,\varphi,\psi,r,\mu}^{+,2},
\\
   \label{eq:vBM_Perturbed_Null}
   %PF7-18-2025 It's not $[A,\Phi]$, it's $[A,\varphi,\psi]$. Please proofread what you write repeatedly.
   %\lambda_{[A,\Phi]}^0(f)&
   \lambda_{[A,\varphi,\psi]}^0(f)&                             
                         :=
  \dim_\RR \mathbf{\tilde H}_{A,\varphi,\psi,r,\mu}^{0,1} - \dim_\RR\mathbf{\tilde H}_{A,\varphi,\psi,r,\mu}^{0,2},
 \end{align}
\end{subequations}
where $\mathbf{\tilde H}_{A,\varphi,\psi,r,\mu}^{+,k}$ is the subspace of $\mathbf{\tilde H}_{A,\varphi,\psi,r,\mu}^k$ on which $S^1$ acts with positive weight and $\mathbf{\tilde H}_{A,\varphi,\psi,r,\mu}^{0,k}$ is the subspace of $\mathbf{\tilde H}_{A,\varphi,\psi,r,\mu}^k$ on which $S^1$ acts with zero weight. We will compute this virtual Morse--Bott index, co-index, and nullity by computing the $S^1$-equivariant index of
%PF7-18-2025 This phrase is meaningless. Please proofread what you write repeatedly.
%the deformation operator of the perturbed non-Abelian monopoles
%PF7-18-2025 Should have cited "deformation operator"
of the deformation operator $\sT_{A,\varphi,\psi,r}$ in the forthcoming \eqref{eq:Perturbed_Deformation_Operator} for the non-Abelian monopole equations \eqref{eq:SO(3)_monopole_equations_almost_Hermitian_perturbed_intro_regular} with a regularized Taubes perturbation and prove

\begin{mainthm}[Virtual Morse--Bott signature of the Hitchin function at a non-zero-section, split non-Abelian monopole]
\label{mainthm:MorseIndexAtReduciblesOnAlmostKahler}
Let $(\rho_{\can},W_{\can})$ be the canonical spin${}^c$ structure \eqref{eq:Canonical_spinc_bundles}, \eqref{eq:Canonical_Clifford_multiplication} over a closed, connected, almost K\"ahler four-manifold $X$ and $E$ be a rank-two Hermitian vector bundle over $X$ that admits a splitting $E=L_1\oplus L_2$ as a direct sum of Hermitian line bundles, and $\ft = (\rho_\can,W_{\can}\otimes E)$ be the corresponding \spinu structure. Assume that $w_2(\su(E))$ satisfies the Morgan--Mrowka condition \eqref{eq:MorganMrowkaCondition}. If $(A,\varphi,\psi)$ is a non-zero-section, non-Abelian monopole on $\ft$ that is split as in the forthcoming Definition \ref{defn:Split_trivial_central-stabilizer_spinor_pair} with respect to the decomposition $E=L_1\oplus L_2$, with 
$\phi\in\Omega^0(L_1)$ and $\psi\in\Omega^{0,2}(L_1)$, then the virtual Morse--Bott index \eqref{eq:vBM_Perturbed_Intro} of the Hitchin function $f$ in \eqref{eq:Hitchin_function} at the point $[A,\Phi] \in \sM_\ft$ is given by twice the complex index of the operator given in \eqref{eq:PerturbedApproxComplexPertDefComplexWeightDecompMinus} and equals
\begin{equation}
\label{eq:MorseIndexAtReduciblesOnKahler}
\lambda_{[A,\varphi,\psi]}^-(f)
=
-\frac{1}{6}\left( c_1(X)^2+c_2(X)\right)
-\left(  c_1(L_1)-c_1(L_2)\right)\cdot c_1(X)
-\left( c_1(L_1)-c_1(L_2)\right)^2,
\end{equation}
where $c_1(X)=c_1(TX,J)$ for the almost complex structure $J$ on the almost K\"ahler manifold $X$ and $c_1(X)^2$ and $c_2(X)$ are as in \eqref{eq:Define_c1Squred_c_2}. The virtual Morse--Bott co-index of $f$ at $[A,\varphi,\psi]$ in \eqref{eq:vBM_Perturbed_coInd} is given by
\begin{equation}
\label{eq:MorseCoIndexAtReduciblesOnKahler}
\begin{aligned}
  \lambda_{[A,\varphi,\psi]}^+(f)
  &=
  -\left( c_1(L_2)-c_1(L_1)\right)\cdot c_1(X) - \left( c_1(L_2)-c_1(L_1)\right)^2
  \\
  &\qquad +c_1(L_2)\cdot c_1(X) +c_1(L_2)^2,
\end{aligned}
\end{equation}
and equals twice the complex index of the operator given in \eqref{eq:PerturbedApproxComplexPertDefComplexWeightDecompPlus}.
The virtual Morse--Bott nullity of $f$ at $[A,\varphi,\psi]$ in \eqref{eq:vBM_Perturbed_Null} is given by
\begin{equation}
\label{eq:MorseNullIndexAtReduciblesOnKahler}
\lambda_{[A,\varphi,\psi]}^0(f) = c_1(L_1)^2 + c_1(L_1)\cdot c_1(X),
\end{equation}
and equals twice the complex index of the operator given in \eqref{eq:PerturbedApproxComplexPertDefComplexWeightDecompZero}.
\end{mainthm}

We rewrite the virtual Morse--Bott index \eqref{eq:MorseIndexAtReduciblesOnKahler} as in the following

\begin{maincor}[Virtual Morse--Bott index of the Hitchin function at a point represented by a Seiberg--Witten monopole]
\label{maincor:MorseIndexAtReduciblesOnSymplecticWithSO3MonopoleCharacteristicClasses}
Continue the assumptions and notation of Theorem \ref{mainthm:MorseIndexAtReduciblesOnAlmostKahler}.
If $[A,\varphi,\psi]\in\sM_\ft$ is represented by a non-zero-section, split non-Abelian monopole in the image of the embedding \eqref{eq:DefnOfIotaOnQuotient} of the moduli space $M_\fs$ of Seiberg--Witten monopoles on the spin${}^c$ structure $\fs=(\rho,W_{\can}\otimes L_1)$ over $X$ into $\sM_\ft$, then the virtual Morse--Bott index \eqref{eq:vBM_Perturbed_Intro} of the Hitchin function $f$ in \eqref{eq:Hitchin_function} on the moduli space $\sM_\ft$ at $[A,\varphi,\psi]$ is given by
\begin{equation}
\label{eq:MorseIndexAtReduciblesOnKahlerSpinNotationType1}
\lambda_{[A,\varphi,\psi]}^-(f)
=
-\frac{1}{6}\left( c_1(X)^2+c_2(X)\right)
-\left( c_1(\fs)-c_1(\ft)\right)\cdot c_1(X)
-\left( c_1(\fs)-c_1(\ft)\right)^2,
\end{equation}
where $c_1(\fs) = c_1(W_\can^+)+2c_1(L_1)\in H^2(X;\ZZ)$ is the first Chern class \eqref{eq:DefineChernClassOfSpinc} of the spin${}^c$ structure given by $\fs = (\rho,W_{\can}\otimes L_1)$ with $W_{\can}=W_{\can}^+\oplus W_{\can}^-$, and $c_1(\ft)=c_1(E)+c_1(W_{\can}^+) \in H^2(X;\ZZ)$ is the first Chern class \eqref{eq:SpinUCharacteristics} of the spin${}^u$ structure.
\end{maincor}

The expression for the virtual Morse--Bott index in \eqref{eq:MorseIndexAtReduciblesOnKahlerSpinNotationType1} is equal to that of the formal Morse--Bott index defined in \eqref{eq:FormalMorseIndexIntroThm}.
The positivity of the formal Morse--Bott index \eqref{eq:FormalMorseIndexIntroThm} in Item \eqref{item:Feasibility_PositivevBM} of Definition \ref{maindefn:Feasibility}
thus implies the following positivity of the virtual Morse--Bott index. Recall that if $(Y,\omega)$ is a symplectic $2n$-manifold, then by \cite[Equation (7.1.14), p. 302]{McDuffSalamonSympTop3} there is a symplectic two-form $\tilde\omega$ on the blow-up $\widetilde{Y} = Y\#\overline{\CC\PP}^n$ of $Y$ at a point. Because the symplectic form $\tilde\omega$ is not unique, we  refer to $(\widetilde Y,\tilde\omega)$ as a 
\emph{symplectic blow-up of $Y$}.

\begin{maincor}[Positivity of the virtual Morse--Bott index of the Hitchin function at a point represented by a Seiberg--Witten monopole]
\label{maincor:Positivity_of_MorseIndexAtReduciblesOnSymplecticWithSO3MonopoleCharacteristicClasses}
Let $X$ be a smooth, closed, connected symplectic four-manifold and $(g,J,\omega)$ be a compatible triple on
%PF10-11-2024 I'm not sure that this "a" versus "the" definition is worth making: any two smooth blowups are diffeomorphic and any two symplectic blowups are symplectomorphic. Why emphasize the distinction?
a symplectic blow-up $\widetilde X=X\#\overline{\CC\PP}^2$. Then there is a spin${}^u$ structure $\ft$ on $\widetilde X$ which is feasible with respect to $c_1(T\widetilde X,J)$ as in Definition \ref{maindefn:Feasibility}. Moreover, if $[A,\varphi,\psi] \in \sM_{\ft}$ is a point represented by a non-zero-section, split, non-Abelian monopole in the image of the embedding \eqref{eq:DefnOfIotaOnQuotient} of the moduli space $M_{\fs}$ of Seiberg--Witten monopoles on a spin${}^c$ structure $\fs$ into $\sM_{\ft}$, then the virtual Morse--Bott index \eqref{eq:vBM_Perturbed_Intro} of the Hitchin function $f$ in \eqref{eq:Hitchin_function} on $\sM_{\ft}$ at $[A,\varphi,\psi]$ is positive.
\end{maincor}

We prove Theorem \ref{mainthm:MorseIndexAtReduciblesOnAlmostKahler} and  Corollaries \ref{maincor:MorseIndexAtReduciblesOnSymplecticWithSO3MonopoleCharacteristicClasses} and \ref{maincor:Positivity_of_MorseIndexAtReduciblesOnSymplecticWithSO3MonopoleCharacteristicClasses}
in Chapter \ref{chap:Calculation_virtual_Morse-Bott_indices_via_Atiyah-Singer_index_theorem}.

\section{Outline}
\label{sec:Outline}
% PF9-13-2024 Add more refs to this section
In Chapter \ref{chap:Prelim}, we summarize the basic concepts from Feehan and Leness 
%TL11-28-2025: Update to \cite[Chapter 6]{Feehan_Leness_introduction_virtual_morse_theory_so3_monopoles}?
\cite[Chapter 6]{Feehan_Leness_introduction_virtual_morse_theory_so3_monopoles} that underlie our definition of the moduli space of non-Abelian monopoles. Chapter \ref{chap:Feasibility} contains our proof of Theorem \ref{mainthm:ExistenceOfSpinuForFlow}, which provides the crucial property of feasibility for the non-Abelian monopole cobordism method. A closed, smooth four-manifold with a Seiberg--Witten basic class (see Definition \ref{defn:Seiberg-Witten_basic_class}) is necessarily equipped with an almost complex structure and so in Chapter \ref{chap:Moduli_space_Seiberg-Witten_monopoles_over_almost_Hermitian_four-manifolds}, we summarize the features of gauge theory over almost Hermitian manifolds that we shall need in our present work.

In Chapter \ref{chap:Elliptic_deformation_complex_moduli_space_SO(3)_monopoles_over_almost_Hermitian_four-manifold}, we recall the definition of the elliptic deformation complex and associated harmonic spaces for the unperturbed and perturbed non-Abelian monopole equations over almost Hermitian four-manifolds and construct real linear isomorphisms between those harmonic spaces and harmonic spaces defined as the kernels of real linear maps of complex Hilbert spaces
%PF9-13-2024 Could we have Banach spaces
whose complex linear and antilinear components can be easily identified.

In Chapter \ref{chap:Analogue_Donaldson_symplectic_submanifold_criterion}, we prove Theorem \ref{mainthm:Donaldson_1996jdg_3_Hilbert_space_codomain} and Corollary \ref{maincorDonaldson_1996jdg_3_Hilbert_space_domain}, Proposition \ref{mainprop:Donaldson_1996jdg_3_Banach_space} and Corollary \ref{maincor:Adjoint_DonaldsonCriteria}, and Theorem \ref{mainthm:Donaldson_1996jdg_3_Hilbert_space} and Corollary \ref{maincor:Donaldson_1996jdg_3_Hilbert_space_non-self-adjoint}, all of which generalize Proposition \ref{prop:Donaldson_1996jdg_3} and which we prove as well.

Our results in are complemented by those in Chapter \ref{chap:Lower_bounds_spectral_gaps_elliptic_operators}, where we prove Theorem \ref{mainthm:Donaldson_1996jdg_3_elliptic_operator_order_m_geq_d}, which gives lower bounds for spectral gaps for self-adjoint elliptic pseudodifferential operators acting on sections of smooth vector bundles over closed smooth manifolds; Theorem \ref{mainthm:Lower_bound_spectral_gap_coupled_Dirac_operator},  which gives lower bounds for spectral gaps for coupled Dirac operators over closed manifolds; Corollary \ref{maincor:Lower_bound_spectral_gap_coupled_Dirac_operator_plus_vector_potential}, which gives lower bounds for spectral gaps for coupled Dirac operators with vector potentials over closed manifolds; and Corollary \ref{maincor:Lower_bound_spectral_gap_coupled_Dirac_operator_plus_vector_potential_four-manifold}, which refines Corollary \ref{maincor:Lower_bound_spectral_gap_coupled_Dirac_operator_plus_vector_potential} for certain coupled Dirac operators with vector potentials over closed four-manifolds.

Chapter \ref{chap:Analogues_non-Abelian_monopoles_Taubes_estimates_Seiberg-Witten_monopole_sections} provides our analogues for solutions to the perturbed non-Abelian monopole equations \eqref{eq:SO(3)_monopole_equations_almost_Hermitian_perturbed_intro} and \eqref{eq:SO(3)_monopole_equations_almost_Hermitian_perturbed_intro_regular} of Taubes' differential inequalities and pointwise estimates for solutions to the perturbed Seiberg--Witten monopole sections in \cite{TauSWGromov}. Our results rely heavily on Bochner--Weitzenb\"ock identities for the different Laplace operators that arise in connection with Dirac operators and complex differential geometry. Because the signs of the curvature terms in these identities are critical and there are no entirely satisfactory references, we provide detailed justifications. We conclude with the proofs of Proposition \ref{mainprop:Taubes_1996_SW_to_Gr_2-1},
Theorem \ref{mainthm:Taubes_1996_SW_to_Gr_2-3}, and
Corollary \ref{maincor:Taubes_1996_SW_to_Gr_eq_2-12_and_2-13} for the system \eqref{eq:SO(3)_monopole_equations_almost_Hermitian_perturbed_intro} and
Proposition \ref{mainprop:Taubes_1996_SW_to_Gr_2-1_regular},
Theorem \ref{mainthm:Taubes_1996_SW_to_Gr_2-3_regular},
and Corollary \ref{maincor:Taubes_1996_SW_to_Gr_eq_2-12_and_2-13_regular} for the system \eqref{eq:SO(3)_monopole_equations_almost_Hermitian_perturbed_intro_regular}.  

In Chapters \ref{chap:Construction_circle-invariant_non-degenerate_two-form_I} and \ref{chap:Construction_circle-invariant_non-degenerate_two-form_II}, we give our constructions of circle-invariant non-degenerate two-form on the moduli space of non-Abelian monopoles, proving Theorems \ref{mainthm:AH_structure_bounded_evalue_spaces_non-Abelian_monopoles_symp_4-mflds} and \ref{mainthm:IdentifyCriticalPoints} and Corollary \ref{maincor:Almost_Hermitian_structure_moduli_space_non-Abelian_monopoles_symplectic_4-manifolds}.
%PF9-18-2025 We did the below in a previous version, but not in the current version
% We also discuss the sense in which the
% \begin{inparaenum}[\itshape i\upshape)]
% \item perturbed non-Abelian monopole equations are \emph{approximately holomorphic} in a suitable Coulomb-gauge coordinate chart, which we require in order to apply Feehan \cite[Theorem 8]{Feehan_analytic_spaces}, ensuring that Seiberg--Witten points are not local minima of Hitchin's function on the moduli space of non-Abelian monopole equations, and 
% \item linearization of the perturbed non-Abelian monopole equations on a suitable Coulomb-gauge slice is \emph{approximately complex linear}, which we require in order to prove Theorems \ref{mainthm:AH_structure_bounded_evalue_spaces_non-Abelian_monopoles_symp_4-mflds} and \ref{mainthm:IdentifyCriticalPoints} and Corollary \ref{maincor:Almost_Hermitian_structure_moduli_space_non-Abelian_monopoles_symplectic_4-manifolds}.
% \end{inparaenum}

Chapter \ref{chap:Calculation_virtual_Morse-Bott_indices_via_Atiyah-Singer_index_theorem} contains our calculation of the virtual Morse--Bott signature at a Seiberg--Witten point in the moduli space of non-Abelian monopole by applying the Atiyah--Singer Index Theorem.

In Appendix \ref{chap:Functional_analysis_maximum_principles}, we include background material on functional analysis and maximum principles that is primarily expository but serves to either clarify our conventions and or is difficult to find in the literature but plays an essential role in the proofs of our main results.

\section{Further work}
\label{sec:Further_work}
The results described in Section \ref{sec:Main_results} would suffice to prove Conjecture \ref{mainconj:BMY_symplectic} if our paradigm did not have to take into account energy bubbling. In this section, we summarize the steps needed to extend the results of the present work and ultimately prove Conjectures \ref{mainconj:BMY_Seiberg-Witten} and \ref{mainconj:BMY_symplectic}.

%PF9-29-2024 Maybe include b_+ = 1 in this subsection too
\subsection[From symplectic four-manifolds to almost Hermitian four-manifolds]{From symplectic four-manifolds to almost Hermitian four-manifolds with a Seiberg--Witten basic class}
\label{subsec:Further_work_almost_Hermitian_4-manifolds_SW_basic_class}
While many of our main results in Section \ref{sec:Main_results} have assumed that the underlying four-manifold is \emph{symplectic}, as required by Conjecture \ref{mainconj:BMY_symplectic}, rather than \emph{almost Hermitian with a Seiberg--Witten basic class}, as allowed by Conjecture \ref{mainconj:BMY_Seiberg-Witten}, we expect that our results in Section \ref{sec:Main_results} should extend to the latter more general case. We now summarize some of the changes that we anticipate to the proofs in our current work that may be required in order to achieve this extension. We view Conjecture \ref{mainconj:BMY_Seiberg-Witten} as the more natural goal of our program, since it primarily relies on fundamental properties of Seiberg--Witten invariants and not on the fact that symplectic four-manifolds have those properties due to results of Taubes \cite{TauSymp, TauSympMore, Taubes_2000}.

The statement of Theorem \ref{mainthm:ExistenceOfSpinuForFlow} and its proof in Chapter \ref{chap:Feasibility}, ensuring feasibility of the non-Abelian monopole cobordism method, requires no change since it only assumes the hypotheses of Conjecture \ref{mainconj:BMY_Seiberg-Witten}. Theorem \ref{mainthm:ExistenceOfSpinuForFlow} is the only result in our present work that relies on Seiberg--Witten invariants. The other appeal to Seiberg--Witten invariants occurs in Feehan and Leness
%TL11-28-2025: Updated to \cite[Proposition 6.8.3]{Feehan_Leness_introduction_virtual_morse_theory_so3_monopoles}
\cite[Proposition 6.8.3]{Feehan_Leness_introduction_virtual_morse_theory_so3_monopoles}, where we prove that the moduli space of non-split, non-zero-section solutions $\sM_\ft^{*,0}$ to the non-Abelian monopole equations (without generic geometric perturbations)
%PF9-12-2024 Cite here
is non-empty if there is a \spinc structure $\fs$ on $X$ with non-zero Seiberg--Witten invariant $\SW_X(\fs)$ and a moduli space $M_\fs$ that is a topologically embedded subspace of $\sM_\ft$. Furthermore, although the hypotheses of
%TL11-28-2025: Updated to \cite[Proposition 6.8.3]{Feehan_Leness_introduction_virtual_morse_theory_so3_monopoles}
 \cite[Proposition 6.8.3]{Feehan_Leness_introduction_virtual_morse_theory_so3_monopoles} assume that $b_1(X)=0$, it is based on results in Feehan and Leness \cite{FL2a, FL2b} that allow $b_1(X) \geq 0$.
%PF9-12-2024 Describe how the proof of above extends to allow b_1 \geq 0.

In Chapter \ref{chap:Elliptic_deformation_complex_moduli_space_SO(3)_monopoles_over_almost_Hermitian_four-manifold} we compare, for $k=0,1,2$, the harmonic spaces $\bH_{A,\varphi,\psi,r}^k$ defined by the linearization at $(A,\varphi,\psi)$ of the perturbed non-Abelian monopole equations and the Coulomb-gauge slice condition with certain harmonic spaces $\bH_{\bar\partial_A,\varphi,\psi,r}^k$. The latter harmonic spaces are defined as kernels of real linear operators whose complex linear and antilinear components are easily identified.

(As an aside, we note that when $r = 0$ and $N_J \equiv 0$ and thus $\bar\mu = \frac{1}{4}N_J^* = 0$, the harmonic spaces $\bH_{\bar\partial_A,\varphi,\psi,r}^k$ become equal to the real linear harmonic spaces $\bH_{\bar\partial_A,\varphi,\psi}^k$ that we previously defined in Feehan and Leness
%TL11-28-2025: Updated to \cite[Sections 10.1, 10.2, and 10.3]{Feehan_Leness_introduction_virtual_morse_theory_so3_monopoles}
\cite[Sections 10.1, 10.2, and 10.3]{Feehan_Leness_introduction_virtual_morse_theory_so3_monopoles} and which reduce (essentially) to the complex linear harmonic spaces $\bH_{\bar\partial_A,\varphi}^k$ for the elliptic deformation complex for the holomorphic pair equations over a complex K\"ahler surface when $\psi \equiv 0$.)

When $(X,g,J,\omega)$ is only \emph{almost Hermitian} rather than symplectic, the non-Abelian monopole equations are unchanged except for the addition of a zeroth-order term to the Dirac equation given by Clifford multiplication with respect to the Lee form,
\[
  \vartheta \equiv \Lambda_\omega d\omega \in \Omega^1(X).
\]
The additional Lee term presents no additional difficulties in Chapter \ref{chap:Elliptic_deformation_complex_moduli_space_SO(3)_monopoles_over_almost_Hermitian_four-manifold}, aside from the need to appeal to more general version of the K\"ahler identities (see below for a further discussion of this point). We expect that the comparisons between harmonic spaces that we develop in Chapter \ref{chap:Elliptic_deformation_complex_moduli_space_SO(3)_monopoles_over_almost_Hermitian_four-manifold} will readily extend to the case where $(X,g,J,\omega)$ is almost Hermitian. Our calculation in Chapter \ref{chap:Calculation_virtual_Morse-Bott_indices_via_Atiyah-Singer_index_theorem} of virtual Morse--Bott indices via the Atiyah--Singer Index Theorem should extend without change from the case where $(X,g,J,\omega)$ is symplectic to the more general case where it is almost Hermitian.

In Chapter \ref{chap:Analogues_non-Abelian_monopoles_Taubes_estimates_Seiberg-Witten_monopole_sections} we derive analogues for the perturbed non-Abelian monopole equations of several of Taubes' differential inequalities and pointwise estimates for the two components of the spinors (rescaled by $\sqrt{r}$) of the perturbed Seiberg--Witten monopole equations. We expect the changes to Chapter \ref{chap:Analogues_non-Abelian_monopoles_Taubes_estimates_Seiberg-Witten_monopole_sections} to be more delicate, but attainable when passing from the case of symplectic to almost Hermitian four-manifolds, as we now explain.

Unlike in Taubes \cite{TauSymp, TauSympMore} (and the expositions of his results due to Donaldson \cite{DonSW} and Kotschick \cite{KotschickSW}), our pointwise estimates in Chapter \ref{chap:Analogues_non-Abelian_monopoles_Taubes_estimates_Seiberg-Witten_monopole_sections} do not rely on the expression from Chern--Weil theory for the degree,
%PF9-12-2024 Check/add ref
\[
  \deg_\omega L = \langle c_1(L)\smile[\omega], [X] \rangle = \frac{i}{2\pi}\int_X F_B\wedge\omega,
\]
where $L$ is a Hermitian line bundle with first Chern class $c_1(L) \in H^2(X;\ZZ)$ over a closed symplectic manifold $(X,\omega)$ and $B$ is a unitary connection on $L$. In particular, our pointwise estimates in Chapter \ref{chap:Analogues_non-Abelian_monopoles_Taubes_estimates_Seiberg-Witten_monopole_sections} do not require $\omega$ to be closed in order to apply the preceding formula for the degree, since we do not need to apply it.

The second feature of Chapter \ref{chap:Analogues_non-Abelian_monopoles_Taubes_estimates_Seiberg-Witten_monopole_sections} is that, like in the expositions of Taubes' results \cite{TauSymp, TauSympMore} due to Donaldson \cite{DonSW} and Kotschick \cite{KotschickSW}, we rely on certain K\"ahler identities for complex K\"ahler manifolds continuing to hold for symplectic manifolds. However, as we explain in Remarks \ref{rmk:Kaehler_identities_for_almost_Hermitian_manifolds} and \ref{rmk:Cirici_Wilson_identities}, those K\"ahler identities continue to hold when $(X,g,J,\omega)$ is an almost Hermitian manifold of dimension \emph{four} due to results of Cirici and Wilson \cite{Cirici_Wilson_2020_almost_hermitian_identities} and others, either unchanged or with manageable modifications.

The third feature of Chapter \ref{chap:Analogues_non-Abelian_monopoles_Taubes_estimates_Seiberg-Witten_monopole_sections} that requires care when passing from the case of symplectic to almost Hermitian four-manifolds is that the Bochner--Weitzenb\"ock formula for $D_A^2$ will contain additional first-order and zeroth-order terms involving the Lee form, $\vartheta$. However, this situation is similar to that already addressed by
%PF9-12-2024 Get exact refs
Feehan and Leness in \cite[Section 4.1, Lemma 4.4, p. 337]{FL1} (where the authors derive a $C^0$ estimate for the harmonic coupled spinor $\Phi$), Feehan in \cite{FeehanGenericMetric}, and Teleman in \cite{TelemanGenericMetric}, where the Dirac operator $D_A$ is modified by adding a zeroth-order term given by Clifford multiplication by a generic one-form $\vartheta$ with complex coefficients. Consequently, we expect that the Bochner--Weitzenb\"ock formulae and corresponding differential equalities and pointwise estimates that we derive in Chapter \ref{chap:Analogues_non-Abelian_monopoles_Taubes_estimates_Seiberg-Witten_monopole_sections} will continue to hold when $(X,g,J,\omega)$ is an almost Hermitian rather than symplectic four-manifold, as assumed in much of this monograph.

% PF9-12-2024 We also need to mention that we use Taubes to get approximately holomorphic local defining equations

\subsection[Allowing for generic geometric perturbations]{Allowing for generic geometric perturbations in the non-Abelian monopole equations}
\label{subsec:Allowing_generic_geometric_perturbations_in_non-Abelian_monopole_equations}
%PF9-27-2024 Put conditions like below at beginning of this further work subsection
Let $(X,g,J,\omega)$ be a closed, connected, smooth almost Hermitian four-manifold and $(E,H)$ be a complex rank two Hermitian vector bundle over $X$. Feehan \cite[Theorem 1.3, p. 910]{FeehanGenericMetric} and Teleman \cite[Theorem 3.19, p. 413]{TelemanGenericMetric} proved that the moduli subspace $\sM_\ft^{*,0}$ of non-split, non-zero-section non-Abelian monopoles is cut out transversely by non-Abelian monopole equations \eqref{eq:PerturbedSO3MonopoleEquations} with generic geometric perturbations \cite[Equation (1.1), p. 908]{FeehanGenericMetric}, namely the Riemannian metric $g$ on $X$, a zeroth-order perturbation $\vartheta \in \Omega^1(X;\CC)$ of the Dirac operator by Clifford multiplication, and a perturbation $\tau \in \Omega^0(\GL(\Lambda^{+,g}(X)))$ of the quadratic term in the non-Abelian monopole equations over a smooth four-manifold. In particular, $\sM_\ft^{*,0}$ is a smooth manifold of the expected dimension. 

For the non-Abelian monopole equations \eqref{eq:SO(3)_monopole_equations_almost_Hermitian_perturbed_intro_regular} with a regularized Taubes perturbation, one can introduce generic geometric perturbations similar to those in \cite{FeehanGenericMetric, TelemanGenericMetric}. We expect that our proof of \cite[Theorem 1.3, p. 910]{FeehanGenericMetric} will adapt, with at most minor modifications, to prove that the moduli space $\sM^{*,0}(E,g,J,r,\omega)$ given by gauge equivalence classes of solutions $(A,\varphi,\psi)$ to the system \eqref{eq:SO(3)_monopole_equations_almost_Hermitian_perturbed_intro_regular} with $(\varphi,\psi)\not\equiv (0,0)$ is cut out transversely for generic geometric perturbations. We expect that the introduction of suitable generic geometric perturbations in the system \eqref{eq:SO(3)_monopole_equations_almost_Hermitian_perturbed_intro_regular} will not significantly impact the proofs of Proposition \ref{mainprop:Taubes_1996_SW_to_Gr_2-1_regular}, Theorem \ref{mainthm:Taubes_1996_SW_to_Gr_2-3_regular}, and Corollary \ref{maincor:Taubes_1996_SW_to_Gr_eq_2-12_and_2-13_regular} and will have little or no impact on the proofs of results elsewhere in this monograph. We note that Tanaka \cite{Tanaka_2019} applied methods of \cite{FeehanGenericMetric} to achieve transversality for the moduli space of solutions to the Vafa--Witten equations over a closed symplectic four-manifold.

\subsection{Allowing for energy bubbling by gluing non-Abelian monopoles}
\label{subsec:Further_work_allowing_energy_bubbling_gluing_non-Abelian_monopoles}
In our present work and its predecessor \cite{Feehan_Leness_introduction_virtual_morse_theory_so3_monopoles}, we computed virtual Morse--Bott indices for non-zero-section Seiberg--Witten points $[A,\Phi]$ in the top level of the moduli space $\sM_\ft$ of non-Abelian monopoles for $\ft_\can = \fs_\can\otimes E$. In particular, we prove that these indices are positive and so the Seiberg--Witten points in $\sM_\ft$ cannot be local minima of the restriction of Hitchin's function \eqref{eq:Hitchin_function} to $\sM_\ft$. Our calculations rely in part on extensions of some of our methods in \cite{FL2a,FL2b} that we applied to prove special cases of Witten's conjecture \cite{Witten} for his expression for the Donaldson series in terms of Seiberg--Witten invariants and basic classes. Hitchin's function $f$ in \eqref{eq:Hitchin_function} naturally extends to a continuous function on the Uhlenbeck compactification $\bar\sM_\ft$ of the moduli space of non-Abelian monopoles. In order to prove Conjecture \ref{mainconj:BMY_Seiberg-Witten} or Conjecture \ref{mainconj:BMY_symplectic}, we must also show that ideal Seiberg--Witten points $[A_0,\Phi_0,\bx]$ in the lower levels $\ell \geq 1$ of the Uhlenbeck compactification $\bar\sM_\ft$ of the moduli space of non-Abelian monopoles are not local minima of $f$. In order to do this, we could first construct a local circle-equivariant Kuranishi-type model for a circle-invariant open neighborhood $U\subset\bar\sM_\ft$ of each ideal Seiberg--Witten point $[A_0,\Phi_0,\bx]$ by extending the gluing method for non-Abelian monopoles that we described in \cite{FL3}. We could then construct the required circle-invariant almost Hermitian structure on the circle-invariant open real analytic manifold defined by the gluing data, prove that the circle-equivariant Kuranishi-type obstruction map is approximately holomorphic, and finally apply \cite[Theorem 8]{Feehan_analytic_spaces} to prove that $[A_0,\Phi_0,\bx]$ is not a local minimum of $f$.

Fortunately, the version of the gluing theorem for non-Abelian monopoles that we need in order to prove Conjecture \ref{mainconj:BMY_Seiberg-Witten} is a relatively modest extension of that proved by the authors in \cite{FL3}. This version is much simpler than the gluing theorem required by our program to prove Witten's conjecture \cite{Witten}. Indeed, the latter gluing theorem remains a technical assumption \cite[Hypothesis 7.8.1, p. 126]{FL5} in our articles and monograph \cite{FL5,FL6,FL7,FL8} giving our proof of that conjecture for closed, connected, oriented, smooth four-manifolds with $b_1 = 0$, odd $b^+ \geq 3$, and Seiberg--Witten simple type. We shall summarize the main steps involved in proving a gluing theorem that suffices to prove Conjecture \ref{mainconj:BMY_Seiberg-Witten}. Our summary relies in part on the well-known approach described by Donaldson and Kronheimer \cite[Sections 7.2 and 7.3]{DK} for gluing anti-self-dual connections. We shall focus on the case $\bx = \{x_1\}$ of a single bubble point in $X$ of multiplicity $\ell \geq 1$, where $(A_0,\Phi_0)$ is a non-zero-section non-Abelian monopole on a Hermitian vector bundle $E_0$ over a closed, smooth almost Hermitian four-manifold $(X,g,J,\omega)$ with $c_2(E) = c_2(E_0) + \ell$ and $A_0$ induces a fixed unitary connection $A_d$ on $\det E_0 = \det E$. The more general case where the point $\bx$ in the symmetric product $\Sym^\ell(X)$ is represented by an $m$-tuple $(x_1,\ldots,x_m)$, where $2 \leq m \leq \ell$, of points $x_i \in X$ with positive multiplicities $k_i$ such that $k_1 + \cdots + k_m = \ell$ differs only in notational complexity from the case $m=1$. As an aside, we note that from the identity \eqref{eq:DK_2-1-39}, the assumption that $0 > p_1(\su(E)) > -c_2(X)$ in the fundamental bounds \eqref{eq:p1_lower_bound}, and the fact that $c_2(E) = c_2(E_0) + \ell$, we obtain
\[
  0 > p_1(\su(E)) =  c_1(E)^2 - 4c_2(E) = c_1(E_0)^2 - 4(c_2(E_0) + \ell) > -c_2(X),
\]
and thus $\ell$ obeys
%TL9-27-2024: Do we want to reverse the order and signs to match the fundamental bounds? Or is the point to see the bounds on $\ell$ so a positive sign for $\ell$ useful?
\[
  0 < - p_1(\su(E_0)) + 4\ell < c_2(X).
\]
To construct the gluing map, we first consider a real analytic family of non-Abelian monopoles $(A_0(t_0),\Phi_0(t_0))$ on $E_0$ that are parametrized by a circle-equivariant Kuranishi model $(\bgamma_0,\bkappa_0)$. Thus, for $\ft_{0,\can} := \fs_\can\otimes E_0$, the Kuranishi method
%PF9-16-2024 Add ref
yields a circle-equivariant, real analytic embedding map,
\[
  \bgamma_0: \bH_{A_0,\Phi_0}^1 \supset U_0 \ni t_0 \mapsto \bgamma_0(t_0) \in
  \tilde\sC_{\ft_{0,\can}} := \sA(E_0,H_0,A_d)\times\Omega^0(W_\can^+\otimes E_0),
\]
from a circle-invariant, precompact open neighborhood $U_0$ of the origin in $\bH_{A_0,\Phi_0}^1$, and a circle-equivariant, real analytic obstruction map,
\[
  \bkappa_0:\bH_{A_0,\Phi_0}^1 \supset U_0 \ni t_0 \mapsto \bkappa_0(t_0) \in \bH_{A_0,\Phi_0}^2,
\]
such that the composition of the restriction of the embedding map,
\[
  \bgamma_0:\bkappa_0^{-1}(0)\cap U_0 \ni t_0 \mapsto \bgamma_0(t_0) \in \tilde\sM_{\ft_{0,\can}},
\]
and the quotient map, $\pi:\tilde\sM_{\ft_{0,\can}} \to \sM_{\ft_{0,\can}} = \tilde\sM_{\ft_{0,\can}}/\SU(E_0)$,
is a circle-equivariant homeomorphism onto an open neighborhood of $[A_0,\Phi_0]$ in $\sM_{\ft_{0,\can}}$.

Second, we let $M_{k_1}(S^4,s)$ denote the real $8k_1$-dimensional moduli space of anti-self-dual connections $A_1$ on a rank two Hermitian vector bundle $(E_1,H_1)$ with second Chern number $c_2(E_1) = k_1 = \ell$ over $S^4 = \RR^4\cup\{\infty\}$ with its standard round metric of radius one and $H_1$-orthonormal \emph{frames} for $E_1\restriction\{s\} = E_1\restriction\{\infty\}$, where $s$ denotes the south pole. There is a real analytic diffeomorphism from $M_{k_1}(S^4,s)$ onto the \emph{complex} $4k_1$-dimensional moduli space $\cM_{k_1}(\CC\PP^2,\ell_\infty)$ of stable holomorphic structures $\bar\partial_{A_1}$ on $E_1$ over $\CC\PP^2 = \CC^2 \cup \ell_\infty$ and a fixed holomorphic trivialization of $\sE_1 = (E_1,\bar\partial_{A_1})$ on the line $\ell_\infty$ at infinity, following Donaldson \cite{DonInstGeomInvar}, Nakajima \cite{Nakajima_lectures_hilbert_schemes_points_surfaces}, and Nakajima and Yoshioka \cite{NakajimaInstLect, NakajimaInstCountI}. In our application, we may confine our attention to the open subset $M_{k_1}(S^4,s,\rho)$ of $M_{k_1}(S^4,s)$ (and corresponding open subset $\cM_{k_1}(\CC\PP^2,\ell_\infty,\rho)$ of $\cM_{k_1}(\CC\PP^2,\ell_\infty)$) such that
\begin{equation}
  \label{eq:No_energy_bubbling_near_infty}
    \frac{1}{8\pi^2}\int_{|x| > \rho} |F_A|^2\,dx \leq \frac{1}{2},
\end{equation}
where $\rho \in (0,1]$ is a positive constant that is less than the injectivity radius of $(X,g)$. The condition \eqref{eq:No_energy_bubbling_near_infty} excludes the possibility of energy bubbling near the south pole in $S^4$ or line at infinity in $\CC\PP^2$. The moduli subspace $M_{k_1}(S^4,s,\rho)$ may be equipped with an Uhlenbeck compactification $\bar M_{k_1}(S^4,s,\rho)$, following the definition due to Donaldson and Kronheimer \cite[Section 4.4]{DK}, based on results due Uhlenbeck \cite{UhlRem, UhlLp}. It should be possible to equip the moduli subspace $\cM_{k_1}(\CC\PP^2,\ell_\infty,\rho)$ with a quasi-projective Gieseker compactification $\bar\cM_{k_1}(\CC\PP^2,\ell_\infty,\rho)$ and construct a continuous epimorphism from $\bar\cM_{k_1}(\CC\PP^2,\ell_\infty,\rho)$ onto $\bar M_{k_1}(S^4,s,\rho)$, by analogy with more general results due to Li \cite{LiAlgGeomDonPoly} and Morgan \cite{MorganComparison}. This epimorphism would facilitate the use of local holomorphic models for open neighborhoods of ideal points in $\bar\cM_{k_1}(\CC\PP^2,\ell_\infty,\rho)$ that surject onto open neighborhoods of the corresponding ideal points in $\bar M_{k_1}(S^4,s,\rho)$.

We fix an $H_0$-orthonormal frame for $E_0\restriction\{x_1\}$ and a $(g,J)$-orthonormal frame for $(TX)_{x_1}$, so $(TX)_{x_1}\cong\CC^2$.

We would use the preceding data to construct a circle-equivariant local model $(\bgamma,\bkappa)$ for an open neighborhood $U \subset \sM_\ft$ whose closure $\bar U \subset \bar\sM_\ft$ with respect to the Uhlenbeck topology contains the ideal Seiberg--Witten point $[A_0,\Phi_0,\bx]$. Thus, we should be able to construct a circle-equivariant, real analytic embedding map,
\[
  \bgamma: \bH_{A_0,\Phi_0}^1 \times \tilde M_{k_1} \supset U_0\times \tilde U_1 \ni (t_0,t_1) \mapsto \bgamma(t_0,t_1) \in \tilde\sC_\ft,
\]
where $\tilde\sC_\ft := \sA(E,H,A_d)\times\Omega^0(W_\can^+\otimes E)$ and $\tilde U_1$ is a circle- and gauge-invariant open subset of $\tilde M_{k_1} \subset \sA_{k_1}$ with circle-invariant, precompact quotient $U_1 \Subset M_{k_1} =  \tilde M_{k_1}/\SU(E_1)$,
%PF9-16-2024 Define \Aut(E_1,H_1)
and a circle-equivariant, real analytic obstruction map,
\[
  \bkappa:\bH_{A_0,\Phi_0}^1 \times \tilde M_{k_1} \supset U_0\times \tilde U_1 \ni (t_0,t_1) \mapsto \bkappa(t_0,t_1) \in \bH_{A_0,\Phi_0}^2,
\]
such that the restriction of the embedding map,
\[
  \bgamma:\bkappa^{-1}(0)\cap (U_0\times \tilde U_1) \ni (t_0,t_1) \mapsto \bkappa(t_0,t_1) \in
  \tilde\sM_\ft \subset \tilde\sC_\ft,
\]
is a circle-equivariant topological embedding \emph{into} an open neighborhood in $\tilde\sM_\ft$. The maps $\bgamma$ and $\bkappa$ are also equivariant with respect to the action of $\SU(E_1)$ on $\sA(E_1,H_1)$ and its actions on $\tilde\sC_\ft$ and $\bH_{A_0,\Phi_0}^2$ induced by the identification of $E_0\restriction\{x_1\}$ with $E_1\restriction \{\ell_\infty\}$.

To obtain $(\bgamma,\bkappa)$ via gluing, we would perform the following steps:
\begin{inparaenum}[\itshape i\upshape)]
\item Construct an \emph{approximate gluing} (or \emph{pregluing} or \emph{splicing}) map, namely a circle-equivariant, real analytic embedding,
\[
  \bgamma': \bH_{A_0,\Phi_0}^1 \times \tilde M_{k_1} \supset U_0\times \tilde U_1 \ni (t_0,t_1) \mapsto \bgamma'(t_0,t_1) \in \tilde\sC_\ft,
\]  
that is equivariant with respect to the action of $\SU(E_1)$ on $\sA(E_1,H_1)$ and its induced action on $\tilde\sC_\ft$ and which provides approximate solutions to the extended perturbed non-Abelian monopole
%PF9-14-2024 Add refs
equations (compare \cite[Section 7.2.1]{DK}).
\item Deform $\bgamma'$ to a nearby \emph{gluing map} $\bgamma$, with the same domain and codomain. This would be a circle-equivariant, real analytic embedding constructed by solving the extended perturbed non-Abelian monopole equations. Its restriction,
\[
  \bgamma:\bkappa^{-1}(0)\cap (U_0\times \tilde U_1) \ni (t_0,t_1) \mapsto \bkappa(t_0,t_1) \in
  \tilde\sM_\ft,
\]
would be a circle-equivariant topological embedding \emph{into} an open neighborhood in $\tilde\sM_\ft$, where
\[
  \bkappa:\bH_{A_0,\Phi_0}^1 \times \tilde M_{k_1} \supset U_0\times \tilde U_1 \ni (t_0,t_1) \mapsto \bkappa(t_0,t_1) \in \bH_{A_0,\Phi_0}^2,
\]
would be the circle-equivariant, real analytic obstruction map obtained by solving the extended perturbed non-Abelian monopole equations (compare \cite[Sections 7.2.2 through 7.2.8]{DK}). The maps $\bgamma$ and $\bkappa$ should also be equivariant with respect to the action of $\SU(E_1)$ on $\sA(E_1,H_1)$ and its induced action on $\tilde\sC_\ft$.
  % PF9-15-2024 Add ref
\item Prove that $\bgamma$ descends to a circle-equivariant, real analytic embedding,
\[
  \bgamma: \bH_{A_0,\Phi_0}^1 \times M_{k_1} \supset U_0\times U_1 \ni (t_0,t_1) \mapsto \bgamma(t_0,t_1) \in \sC_\ft,
\]
such that its restriction
\[
  \bgamma:\bkappa^{-1}(0)\cap (U_0\times U_1) \ni (t_0,t_1) \mapsto \bkappa(t_0,t_1) \in
  \sM_\ft,
\]
is a circle-equivariant topological embedding \emph{into} an open neighborhood in $\sM_\ft$ whose Uhlenbeck closure in $\bar\sM_\ft$ contains the ideal Seiberg--Witten point $[A_0,\Phi_0,\bx]$.
\end{inparaenum}

We could now use the local gluing model $(\bgamma,\bkappa)$ in the same way as we do in our present work for the local Kuranishi model to prove that the ideal Seiberg--Witten point $[A_0,\Phi_0,\bx]$ is not a local minimum of Hitchin's function $f$.

Aside from the analytical difficulties that are unique to the non-Abelian monopole equations (see \cite{FL3} for a discussion), the construction of the gluing parametrization $(\bgamma,\bkappa)$ is simpler than that described in Donaldson and Kronheimer \cite[Chapter 7]{DK} because we do not need to perform the additional step of proving this gluing parametrization is \emph{surjective} in the sense that every point $[A,\Phi] \in \sM_\ft$ that is Uhlenbeck-close to $[A_0,\Phi_0,\bx]$ necessarily lies in the image of $\bgamma$, as those authors do in \cite[Section 7.3]{DK} for their gluing map for anti-self-dual connections. Second, in contrast to our work on the Witten conjecture \cite{Witten} and the Kotschick--Morgan conjecture \cite{KotschickBPlus1, KotschickMorgan}, we do not need to prove that the pair $(\bgamma,\bkappa)$ extends continuously to the Uhlenbeck compactification $\bar\sM_\ft$. Third and finally, although the construction we have just outlined easily extends to the case of a point $\bx \in \Sym^\ell(X)$ represented by an $m$-tuple $(x_1,\ldots,x_m)$ with $2 \leq m \leq \ell$, we do not need to consider the \emph{overlap problem} for gluing parametrizations that we described in \cite{FLMcMaster}. The overlap problem is one of the central difficulties in a key part \cite{FL5} of our proof of the Witten and Kotschick--Morgan conjectures.

\subsection[Compact, complex surfaces of general type]{Allowing for energy bubbling by restricting to compact, complex surfaces of general type}
\label{subsec:Further_work_allowing_energy_bubbling_compact_complex_projective_surfaces}
The Bradlow--Hitchin--Kobayashi correspondence \cite{Bradlow_1991} gives an $S^1$-equivariant, real analytic isomorphism from the moduli space $\sM(E,g,\omega)$ of type $1$ non-Abelian monopoles (projective vortices) over a compact, complex K\"ahler surface onto the moduli space $\fM_{\ps}(E)$ of polystable holomorphic pairs inducing a fixed holomorphic structure on $\det E$. 

In Feehan, Leness, and Wentworth \cite{Feehan_Leness_Wentworth_virtual_morse_theory_stable_pairs_bmy_kaehler}, where $X$ is assumed to be a complex projective surface, we adapt work of Okonek, Schmitt, and Teleman \cite{OTMasterPairs} and Dowker, \cite{DowkerThesis}, Huybrechts and Lehn \cite{Huybrechts_Lehn_1995ijm}, Lin \cite{Lin_2018}, and Wandel \cite{Wandel_2015} to construct a moduli space $\fM_{\sss}(E)$ of oriented pairs $(\sE,\varphi,\eps)$ of coherent sheaves and sections and morphisms $\eps$ from $\det\sE$ onto a fixed line bundle $\Lambda$ over $X$ that are \emph{semistable} in the sense of Bradlow \cite{Bradlow_1991}, Gieseker \cite{Gieseker_1977}, and Maruyama \cite{Maruyama_1977, Maruyama_1978}. The space $\fM_{\sss}(E)$ is equipped with a holomorphic $\CC^*$ action.

We use the complex projective moduli space $\fM_{\sss}(E)$ as a compactification for $\sM(E,g,\omega)$. We develop a deformation theory for oriented pairs $(\sE,\varphi,\eps)$ extending that of Huybrechts and Lehn \cite{Huybrechts_Lehn_1995ijm} and Lin \cite{Lin_2018} for pairs $(\sE,\varphi)$ and that of Huybrechts and Lehn \cite{Huybrechts_Lehn_geometry_moduli_spaces_sheaves} for sheaves $\sE$. This gives lower bounds for the Krull (true) dimensions of unstable varieties through $\CC^*$-fixed points $[\sE,\varphi,\eps]$ (with $\varphi\not\equiv 0$) in terms of their expected dimensions, $\dim\HH_{\sE,\varphi}^{-,1} - \dim\HH_{\sE,\varphi}^{-,2}$, virtual Bia{\l}ynicki--Birula indices given by the Hirzebruch--Riemann--Roch Theorem as differences of dimensions of hypercohomology groups. By calculations similar to those in our monograph \cite{Feehan_Leness_introduction_virtual_morse_theory_so3_monopoles}, we prove that all $S^1$-fixed points $[\sE,\varphi,\eps]$ (with $\varphi\not\equiv 0$) in $\fM_{\sss}(E)$ have positive virtual Morse--Bott index. This proves a new proof of the Bogomolov--Miyaoka--Yau inequality \eqref{eq:BMY} for complex surfaces, validating our gauge-theoretic program to prove Conjecture \ref{mainconj:BMY_Seiberg-Witten} in full.

%PF8-31-2025 Editing bookmark
\chapter{Preliminaries}
\label{chap:Prelim}
We review some definitions and results from the theory of non-Abelian monopoles \cite{FL1,FL2a,Feehan_Leness_introduction_virtual_morse_theory_so3_monopoles}. After introducing spin${}^c$ and spin${}^u$ structures, their characteristic classes, and spin connections in Section \ref{sec:SpinuDefinitions}, we review the groups of gauge transformation and stabilizer subgroups of connections in Section \ref{sec:UnitaryConnAndGaugeTransformations}. We discuss the stabilizer subgroups of pairs, the quotient spaces of spin${}^c$ and spin${}^u$ pairs, and embeddings of the quotient space of spin${}^c$ pairs into that of spin${}^u$ pairs in Section \ref{sec:QuotientSpaceOfSpinuPairs}.  The definition of the moduli space of non-Abelian monopoles and results on its expected dimension and compactification appear in Section \ref{sec:ModuliOfSO3Monopole}. Section \ref{sec:ZeroSectionPairs} contains a discussion of anti-self-dual connections and their identification with zero-section non-Abelian monopoles. The identification of Seiberg--Witten monopoles and split non-Abelian monopoles and the results on Seiberg--Witten invariants that we need appear in Section \ref{sec:ReduciblePairs}. In Section \ref{sec:S1Actions}, we describe the circle action on the moduli space of non-Abelian monopoles and its fixed points.

\section{Clifford modules and spin${}^c$ and spin${}^u$ structures}
\label{sec:SpinuDefinitions}
We review the definitions and characteristic classes of spin${}^c$ and spin${}^u$ structures on a closed and oriented Riemannian four-manifold $(X,g)$.
A spin${}^c$ bundle over $(X,g)$ is a complex rank four vector bundle $W\to X$ with a real linear map $\rho:T^*X\to \End(W)$ satisfying (see Lawson and Michelsohn \cite[Chapter I, Proposition I.5.10, p. 35]{LM}),
\begin{equation}
\label{eq:CliffordMapDefn}
\rho(\alpha)^2 = -g(\alpha,\alpha)\id_{V}
\quad\text{and}\quad
\rho(\alpha)^\dagger = -\rho(\alpha),
\quad\text{for all } \alpha \in C^\infty(T^*X).
\end{equation}
The map $\rho$ uniquely extends to a linear map $\rho:\Lambda^{\bullet}(T^*X)\otimes_\RR\CC\to\End(W)$, and gives $W$ the structure of a \emph{Hermitian Clifford module} for the complex Clifford algebra $\CCl(T^*X)$ (see Lawson and Michelsohn \cite[Chapter I, Proposition 5.10, p. 35]{LM}).  There is a splitting $W=W^+\oplus W^-$, where $V^\mp$ are the $\pm 1$ eigenspaces of $\rho(\vol_g)$. A unitary connection $A$ on $W$ is called \emph{spin} if
\begin{equation}
\label{eq:SpinConnection}
[\nabla_A,\rho(\alpha)] =\rho(\nabla\alpha)
\quad\text{on }C^\infty(V),
\end{equation}
for any $\alpha\in C^\infty(T^*X)$, where $\nabla$ is the Levi-Civita connection.
We define  the associated first Chern class of a spin${}^c$ structure $\fs=(\rho,W)$ by
\begin{equation}
\label{eq:DefineChernClassOfSpinc}
c_1(\fs):=c_1(W^+) \in H^2(X;\ZZ),
\end{equation}
which is an integral lift of the second Stiefel-Whitney class of $X$.  If $L$ is a Hermitian line bundle over $X$, we write $\fs\otimes L$ for the \spinc structure $(\rho\otimes\id_L,W\otimes L)$, which we abbreviate as $(\rho,W\otimes L)$.  Note the equality
\begin{equation}
\label{eq:c1_of_TensoredSpinc}
c_1(\fs\otimes L)=c_1(\fs)+2c_1(L).
\end{equation}
A \emph{spin${}^u$} structure on $X$ is defined by $\ft=(\rho\otimes\id_E,W\otimes E)$ where $(\rho,W)$ is a spin${}^c$ structure and  $(E,H)$ is a complex rank-two Hermitian vector bundle over $X$.  We note that for any complex line bundle $L$ over $X$ the spin${}^c$ structure
$(\rho\otimes\id_L,W\otimes L)$ and vector bundle $E\otimes L^*$ will define the same spin${}^u$ structure as
$(\rho,W)$ and $E$.  A spin${}^u$ structure $\ft=(\rho,W\otimes E)$ has an associated line bundle,
\begin{equation}
\label{eq:HalfDeterminantLineBundleOfSpinu}
\det(W^+\otimes E)^{1/2}:=\det(W^+)\otimes\det(E),
\end{equation}
and an associated real rank-three bundle
\begin{equation}
\label{eq:Definesu(E)Bundle}
\fg_\ft:=\su(E),
\end{equation}
as well as characteristic classes
\begin{equation}
\label{eq:SpinUCharacteristics}
c_1(\ft) := \frac{1}{2} \left( c_1(W^+)+c_1(E)\right),
\quad
p_1(\ft) := p_1(\su(E)),
\quad
\text{and}\quad
w_2(\ft) := w_2(\su(E)).
\end{equation}

\section{Unitary connections and gauge transformations}
\label{sec:UnitaryConnAndGaugeTransformations}
Let $(E,H)$ be a  complex Hermitian vector bundle of rank $r$ over a real dimension $d$ Riemannian manifold $(X,g)$.
Fix a smooth, unitary connection $A_d$ on $\det E$.  For $p\in (d/2,\infty)$, let $\sA(E,H,A_d)$ denote the space of $W^{1,p}$ unitary connections on $E$ that induce $A_d$ on $\det E$. We let
\begin{subequations}
\label{eq:UnitaryAutomorphismBundles}
\begin{align}
  \label{eq:AutomorphismBundles_U(E)}
  \U(E,H)& := \{u\in \End(E): u(x)\in\U(E_x), \text{ for all } x\in X\},
  \\
  \label{eq:AutomorphismBundles_SU(E)}
  \SU(E,H)& := \{u\in \End(E): u(x)\in\SU(E_x), \text{ for all } x\in X\},
\end{align}
\end{subequations}
denote the smooth principal fiber bundles over $X$ with structure groups $\U(r)$ and $\SU(r)$, respectively. Gauge transformations are sections of these bundles. When the Hermitian metric $h$ is understood, we shall abbreviate $\U(E):=\U(E,H)$ and $\SU(E):=\SU(E,H)$.

The Banach Lie group $W^{2,p}(\SU(E))$ of determinant-one, unitary automorphisms of $(E,H)$ acts on (the right of) $\sA(E,H,A_d)$
 \begin{equation}
  \label{eq:W2pSUE_action_on_AEh}
  W^{2,p}(\SU(E)) \times \sA(E,H,A_d) \ni (u,A) \mapsto u^*A \in \sA(E,H,A_d)
\end{equation}
by pullback. By adapting the proofs \mutatis due to Freed and Uhlenbeck of \cite[Appendix A, Proposition A.2, p. 160 and Proposition A.3, p. 161]{FU}, one can show that $W^{2,p}(\SU(E))$ is a Banach Lie group and that its action \eqref{eq:W2pSUE_action_on_AEh} on $\sA(E,H,A_d)$ is smooth.

For $A\in\sA(E,H,A_d)$, we define the \emph{stabilizer} of $A$ as
\begin{equation}
\label{eq:DefineStabilizerOfConnection_in_SU(E)}
\Stab(A):=\{u\in W^{2,p}(\SU(E)): u^*A=A\}.
\end{equation}
We note the following well-known Lie group properties of stabilizers of unitary connections
(see e.g.
%TL11-28-2025: Updated to \cite[Lemma 6.1.3]{Feehan_Leness_introduction_virtual_morse_theory_so3_monopoles}
\cite[Lemma 6.1.3]{Feehan_Leness_introduction_virtual_morse_theory_so3_monopoles})

%TL9-16-2024: Some of these may not be needed in this work.
\begin{lem}[Lie group structure of $\Stab(A)$]
\label{lem:LieGroup_Structure_of_Stab(A)}
Let $(E,H)$ be a smooth Hermitian vector bundle over a connected, smooth manifold $X$. If $A\in\sA(E,H,A_d)$, then $\Stab(A)$ is a Lie subgroup of $W^{2,p}(\SU(E))$ and
\[
  \bH_A^0
  :=
  \Ker\left(d_{A^{\ad}}:W^{2,p}(\su(E))\to W^{1,p}(T^*X\otimes \su(E))\right)
\]
is the Lie algebra of $\Stab(A)$.
\end{lem}

We record the following definitions.

\begin{defn}[Reducible, split, and central-stabilizer unitary connections]
\label{defn:Reducible_split_trivial-stabilizer_unitary_connection}
Let $A$ be a $W^{1,p}$ unitary connection with $p\in (d/2,\infty)$ on a $W^{2,p}$ Hermitian vector bundle $(E,H)$ of complex rank $r$ over a smooth, connected manifold of dimension $d\geq 2$.
\begin{enumerate}
\item\label{item:central-stabilizer_unitary_connection}
$A$ has \emph{central stabilizer} if the stabilizer $\Stab(A)$ of $A$ is minimal and thus isomorphic to the center $Z(\SU(r))=C_r$ and has \emph{non-central stabilizer} otherwise.

\item\label{item:Split_unitary_connection}
$A$ is \emph{split} if $A = A_1\oplus A_2$ with respect to a decomposition
\begin{equation}
\label{eq:BasicSplitting}
E = E_1\oplus E_2
\end{equation}
as an orthogonal direct sum of proper, $W^{2,p}$ Hermitian subbundles, where $A_i$ is a $W^{1,p}$ unitary connection on $E_i$, for $i=1,2$, and is \emph{non-split} otherwise. If the pair $(E,A)$ is smooth (respectively, real analytic), then each pair $(E_i,A_i)$ is smooth (respectively, real analytic).

\item\label{item:Reducible_unitary_connection}
$A$ is \emph{reducible} if the holonomy group of $A$ reduces (in the sense of Kobayashi and Nomizu \cite[Section II.6, p. 81]{Kobayashi_Nomizu_v1}) to a proper subgroup of $\U(r)$ and is \emph{irreducible} otherwise.
\end{enumerate}
\qed\end{defn}

\section{Quotient space of spin${}^u$ pairs}
\label{sec:QuotientSpaceOfSpinuPairs}
Let $\ft=(\rho,W\otimes E)$ be a spin${}^u$ structure on a smooth, closed, and oriented Riemannian four-manifold $(X,g)$.
We define the $W^{2,p}(\SU(E))$ quotient of the pairs $(A,\Phi)\in \sA(E,H,A_d)\times W^{1,p}(W^+\otimes E)$ and its subspaces.

\subsection{Action of the group of gauge transformations on the affine space of spin${}^u$ pairs}
\label{subsec:SU(E)ActionOnSpinuPairs}
The Banach Lie group $W^{2,p}(\SU(E))$ acts smoothly on the Banach affine space $\sA(E,H,A_d)\times W^{1,p}(W^+\otimes E)$ by
\begin{multline}
\label{eq:GaugeActionOnSpinuPairs}
W^{2,p}(\SU(E))\times \left(\sA(E,H,A_d)\times W^{1,p}(W^+\otimes E)\right)\ni (u,(A,\Phi))
\\
\mapsto u^*(A,\Phi):=(u^*A,u^{-1}\Phi)
\in \sA(E,H,A_d)\times W^{1,p}(W^+\otimes E).
\end{multline}
We denote the quotient space of spin${}^u$ pairs associated to $\ft$ by
\label{page:SpinuPairs}
\begin{equation}
\label{eq:SpinUConfiguration}
\sC_\ft := \left.\left(\sA(E,H,A_d)\times W^{1,p}(W^+\otimes E)\right)\right/W^{2,p}(\SU(E)),
\end{equation}
where  $W^{2,p}(\SU(E))$ acts on $\sA(E,H,A_d)\times W^{1,p}(W^+\otimes E)$  as in \eqref{eq:GaugeActionOnSpinuPairs}. We write elements of $\sC_\ft$ as $[A,\Phi]$\label{page:GaugeEquivClassOfSpinuPair}.
When the four-manifold $X$ admits an almost Hermitian structure $(X,g,J,\omega)$ and we use the canonical spin${}^c$ structure $\fs_\can$ of Definition \ref{defn:Canonical_spinc_bundles} to identify rank-two complex Hermitian vector bundles $(E,H)$ with spin${}^u$ structures $\ft=\fs_\can\otimes E$, we write
\begin{equation}
\label{eq:DefineSpinUConfigUnitaryTripleNotation}
\sC(E,H,J,A_d):=\sC_{\fs_\can\otimes E}.
\end{equation}
%TL5-20-2025: Added following.  Can be deleted if redundant.
We define the \emph{stabilizer group} of $(A,\Phi)\in\sA(E,H,A_d)\times W^{1,p}(W^+\otimes E)$ by
\begin{equation}
\label{eq:Define_Stabilizer_of_Spinu_Pair_in_SU(E)}
\Stab(A,\Phi):=\{u\in W^{2,p}(\SU(E)): u^*(A,\Phi)=(A,\Phi)\}.
\end{equation}
The stabilizer of $(A,\Phi)$ has the following structure (see
%TL11-28-2025: Updated to \cite[Lemma 6.3.2]{Feehan_Leness_introduction_virtual_morse_theory_so3_monopoles}
\cite[Lemma 6.3.2]{Feehan_Leness_introduction_virtual_morse_theory_so3_monopoles}).

\begin{lem}[Lie group structure of $\Stab(A,\Phi)$]
\label{lem:LieGroupStructureOfStab(A,Phi)}
Let $\ft=(\rho,W\otimes E)$ be a spin${}^u$ structure on a closed, oriented, smooth Riemannian manifold $(X,g)$ of real dimension four. For $(A,\Phi)\in \sA(E,H,A_d)\times W^{1,p}(W^+\otimes E)$, the stabilizer $\Stab(A,\Phi)$ is a Lie group with Lie algebra given by
\[
\bH_{A,\Phi}^0:=
\Ker\left(
d_{A,\Phi}^0: W^{2,p}(\su(E)) \to W^{1,p}(\Lambda^1(\su(E)))\oplus W^{1,p}(W^+\otimes E)
\right)
\]
where $d_{A,\Phi}^0$ is given by
\begin{equation}
  \label{eq:d_APhi^0}
  d_{A,\Phi}^0 \xi = (d_A\xi, -\xi\Phi), \quad\text{for all } \xi \in W^{2,p}(\su(E)).
\end{equation}
\end{lem}

The following result $\Stab(A,\Phi)$ appears in
%TL11-28-2025: Updated to \cite[Lemma 6.3.3]{Feehan_Leness_introduction_virtual_morse_theory_so3_monopoles}
\cite[Lemma 6.3.3]{Feehan_Leness_introduction_virtual_morse_theory_so3_monopoles}.

%TL5-16-2025: This has two references to it in this book, so let's keep.
\begin{lem}[Non-zero section rank-two \spinu pairs have trivial stabilizer]
\label{lem:NonZeroSection_Spinu_pairs_Have_Trivial_Stabilizer}
Let $\ft=(\rho,W\otimes E)$ be a spin${}^u$ structure on a closed, connected, oriented, smooth Riemannian manifold $(X,g)$ of real dimension four. If $E$ has rank two and $(A,\Phi)\in \sA(E,H,A_d)\times W^{1,p}(W^+\otimes E)$ with $\Phi\not\equiv 0$ and $\Stab(A)\not\cong \SU(2)$, then $\Stab(A,\Phi)=\{\id_E\}$.
\end{lem}

\subsection{Subspaces of the quotient space of spin${}^u$ pairs}
\label{subsec:SubspacesOfQuotientSpaceOfSpinuPairs}
We now discuss subspaces of the quotient space $\sC_\ft$.

\begin{defn}[Zero-section \spinu pairs]
\label{defn:ZeroSectionSO(3)Pair}
(See Feehan and Leness \cite[Definition 2.2, p. 64]{FL1}.)
We call a pair $(A,\Phi)\in\sA(E,H,A_d)\times W^{1,p}(W^+\otimes E)$ a \emph{zero-section} pair if $\Phi\equiv 0$.
\qed\end{defn}

The subspace of $\sC_\ft$ given by gauge-equivalence classes of zero-section pairs is identified with the
\emph{quotient space} $\sB(E,H)$
\label{page:Quotient_space_unitary_connections}
\emph{of fixed-determinant, unitary connections} on $E$. Recall from \eqref{eq:Define_Stabilizer_of_Spinu_Pair_in_SU(E)} that the stabilizer group of $(A,\Phi)\in\sA(E,H,A_d)\times W^{1,p}(W^+\otimes E)$ is
\[
  \Stab(A,\Phi):=\{u\in W^{2,p}(\SU(E)): u^*(A,\Phi)=(A,\Phi)\}.
\]
We provide the following partial analogue of Definition \ref{defn:Reducible_split_trivial-stabilizer_unitary_connection} for \spinu pairs.

\begin{defn}[Split, trivial-stabilizer, and central-stabilizer \spinu pairs]
\label{defn:Split_trivial_central-stabilizer_spinor_pair}
(Compare Feehan and Leness \cite[Definition 2.2, p. 64]{FL2a} for the case $r=2$.)
Let $(E,H)$ be a smooth Hermitian vector bundle of complex rank $r\geq 2$ over a smooth, connected Riemannian manifold of real dimension four and let $(A,\Phi)\in\sA(E,H,A_d)\times W^{1,p}(W^+\otimes E)$.
\begin{enumerate}
\item\label{item:Trivial-stabilizer_spinor_pair}
$(A,\Phi)$ has \emph{trivial stabilizer} if the stabilizer group $\Stab(A,\Phi)$ of $(A,\Phi)$ in $W^{2,p}(\SU(E))$ is $\{\id_E\}$ and has \emph{non-trivial stabilizer} otherwise.

\item\label{item:central-stabilizer_spinor_pair}
$(A,\Phi)$ has \emph{central stabilizer} if the stabilizer group $\Stab(A,\Phi)$ of $(A,\Phi)$ in $W^{2,p}(\SU(E))$ is equal to the center $Z(\SU(r)) = C_r$ and has \emph{non-central stabilizer} otherwise.

\item\label{item:Split_spinor_pair}
$(A,\Phi)$ is \emph{split} if the connection $A$ is split as in Definition \ref{defn:Reducible_split_trivial-stabilizer_unitary_connection} \eqref{item:Split_unitary_connection} and $\Phi \in W^{1,p}(W^+\otimes E_1)$. If no such splitting exists, then $(A,\Phi)$ is called \emph{non-split}. We say that $(A,\Phi)$ is \emph{split with respect to the decomposition $E=E_1\oplus E_2$} in \eqref{eq:BasicSplitting} when we wish to specify the bundle splitting.
\end{enumerate}
\qed\end{defn}

We define subspaces of $\sC_\ft$ by
\begin{subequations}
 \label{eq:SpinuQuotientSpaceSubspaces}
 \begin{align}
   \label{eq:SpinuQuotientSpaceSubspaces_irreducible}
   \sC_\ft^* &:= \{[A,\Phi]\in\sC_\ft: \text{$(A,\Phi)$ is non-split}\},
   \\
   \label{eq:SpinuQuotientSpaceSubspaces_non-zero-section}
   \sC_\ft^0 &:= \{[A,\Phi]\in\sC_\ft: \Phi\not\equiv 0\},
   \\
   \label{eq:eq:SpinuQuotientSpaceSubspaces_irredicible_non-zero-section}
   \sC_\ft^{*,0} &:= \sC_\ft^*\cap\sC_\ft^0.
\end{align}
\end{subequations}
When $(X,g,J,\omega)$ is an almost Hermitian four-manifold and the spin${}^u$ structure $\ft$ is given by $\ft=(\rho,W_{\can}\otimes E)$ where $\fs_{\can}=(\rho,W_{\can})$ is the canonical spin${}^c$ structure of Definition \ref{defn:Canonical_spinc_bundles}, then 
\[
\sC_\ft=\sC(E,H,J,A_d),
\]
and the subspaces \eqref{eq:SpinuQuotientSpaceSubspaces} are equal to the corresponding spaces defined in \eqref{eq:Quotient_space_non-zero_or_non-split_unitary_triples}.

Because  the stabilizer of a zero-section pair always includes the constant central gauge transformations,
if $\Stab(A,\Phi)=\{\id_E\}$, then $[A,\Phi]\in\sC_t^0$.
As we are only considering the case that $E$ has rank two,   Lemma \ref{lem:NonZeroSection_Spinu_pairs_Have_Trivial_Stabilizer}
implies that $\Stab(A,\Phi)=\{\id_E\}$ when $\Phi\not\equiv 0$  and so
\begin{equation}
  \label{eq:sC0(ft_equals_sC**(ft)_E_rank2}
  \sC_\ft^0 = \{[A,\Phi]\in\sC_\ft: \Stab(A,\Phi) = \{\id_E\}\}.
\end{equation}
When $E$ has rank two,
the quotient space $\sC_\ft^0$ in \eqref{eq:SpinuQuotientSpaceSubspaces_non-zero-section} is an analytic Banach manifold by the analogue for \spinu pairs of
%TL11-28-2025: Updated to  \cite[Theorem 12.3.10]{Feehan_Leness_introduction_virtual_morse_theory_so3_monopoles}
\cite[Theorem 12.3.10]{Feehan_Leness_introduction_virtual_morse_theory_so3_monopoles}.

We now describe a local slice for the action \eqref{eq:GaugeActionOnSpinuPairs} of $W^{2,p}(\SU(E))$ on the Banach affine space $\sA(E,A_d,H)\times\Omega^0(W^+\otimes E)$.  Define
\[
 d_{A,\Phi}^{0,*}: W^{1,p}\left(T^*X\otimes \su(E) \oplus W^+\otimes E\right) \to L^p(\su(E)),
\]
as the $L^2$ adjoint of the operator $d_{A,\Phi}^0$ defined in \eqref{eq:d_APhi^0}. From
%TL11-28-2025: Updated to \cite[Equation (10.1.16)]{Feehan_Leness_introduction_virtual_morse_theory_so3_monopoles}
\cite[Equation (10.1.16)]{Feehan_Leness_introduction_virtual_morse_theory_so3_monopoles}, we have the following explicit expression:
\begin{equation}
  \label{eq:d_APhi^0_star_aphi_identity_and_vanishing}
  d_{A,\Phi}^{0,*}(a,\phi) = d_A^*a - \sR_\Phi^*\phi = 0,
\end{equation}
where
\begin{equation}
  \label{eq:Phi_star_Omega0V+_to_Omega0suE}
  \sR_\Phi^*:\Omega^0(W^+\otimes E) \to \Omega^0(\su(E))
\end{equation}
is the pointwise adjoint (with respect to the Riemannian metric on $\su(E)$ and the Hermitian metric on $W^+\otimes E$) of the right composition operator,
\begin{equation}
  \label{eq:Phi_Omega0suE_to_Omega0V+}
  \sR_\Phi:\Omega^0(\su(E)) \ni \xi \mapsto \xi\Phi \in \Omega^0(W^+\otimes E).
\end{equation}
The \emph{Coulomb gauge condition},
\begin{equation}
\label{eq:SO3_monopole_Coulomb_gauge_slice_condition}
d_{A,\Phi}^{*,0}(a,\phi)=0\quad\text{for $(a,\phi)\in W^{1,p}\left(T^*X\otimes \su(E) \oplus W^+\otimes E\right)$,}
\end{equation}
defines a slice in the following sense.
By Feehan and Maridakis \cite[Theorem 16, p. 18 and Corollary 18, p. 19]{Feehan_Maridakis_Lojasiewicz-Simon_coupled_Yang-Mills}, there is a constant $\zeta = \zeta(A,\Phi) \in (0,1]$ such that the ball,
\begin{equation}
  \label{eq:SliceBall}
  \bB_{A,\Phi}(\zeta)
  :=
  \left\{(a,\phi)\in \Ker d_{A,\Phi}^{0,*}\cap W^{1,p}\left(T^*X\otimes \su(E) \oplus W^+\otimes E\right):
    \|(a,\phi)\|_{W_A^{1,p}(X)} <\zeta\right\},
\end{equation}
embeds into the quotient $\sC_\ft$ under the map $(a,\phi)\mapsto [(A,\Phi)+(a,\phi)]$ when $\Stab(A,\Phi)$ is trivial;
in general the quotient $\bB_{A,\Phi}(\zeta)/\Stab(A,\Phi)$ will embed into $\sC_\ft$.

\subsection{Embedding quotient spaces of spin${}^c$ pairs in the quotient space of spin${}^u$ pairs}
\label{subsec:SpincQuotientSpaceInSpinuSpace}
We begin by defining the quotient space of spin${}^c$ pairs on a closed, four-dimensional, oriented, smooth Riemannian manifold $(X,g)$.

As in our definition of spin${}^u$ structures, we fix a spin${}^c$ structure $\fs_0=(\rho,W)$ and a unitary connection $A_W$ on $W$ that is
 spin in the sense of \eqref{eq:SpinConnection}. 
When $X$ is almost Hermitian, we will usually take $\fs_0=\fs_{\can}$, the canonical spin${}^c$ structure of Definition \ref{defn:Canonical_spinc_bundles}.
 Let $(L,H_L)$ be a Hermitian line bundle over $X$ and let $\fs=\fs_0\otimes L:=(\rho\otimes \id_L,W\otimes L)$ be the \spinc structure induced by $L$ and $\fs_0$. In Feehan and Leness \cite[Equation (2.52), p. 76]{FL2a}, we consider the affine space $\sA_{\fs}\times W^{1,p}(W^+\otimes L)$ of \spinc pairs $(B,\Psi)$, where $\sA_{\fs}$ denotes the space of $W^{1,p}$ unitary connections $B$ on $W\otimes L$ that are spin in the sense of \eqref{eq:SpinConnection}.
In this work, we instead take
\[
  \sA(L,H_L)\times W^{1,p}(W^+\otimes L)
\]
to be the affine space of spin${}^c$ pairs $(A_L,\Psi)$, where $\sA(L,H_L)$ is the affine space of unitary connections $A_L$ on $L$. While this definition depends on the choice of spin${}^c$ structure $\fs_0$ and connection $A_W$ and thus is not intrinsic to the spin${}^c$ structure $\fs_0\otimes L$, it is more convenient for our current applications.  The map
\begin{equation}
\label{eq:TranslateSWNotation}
\sA(L,H_L)
\ni
A_L
\mapsto
A_W\otimes A_L
\in
\sA_{\fs_0\otimes L}
\end{equation}
identifies the two affine spaces.

Let  $W^{2,p}(X,S^1)$ denote the Banach Lie group of $W^{2,p}$ maps from $X$ to $S^1$. We define the action of $W^{2,p}(X,S^1)$ on the affine space $\sA(L,h_L)\times W^{1,p}(W^+\otimes L)$ by analogy with \eqref{eq:GaugeActionOnSpinuPairs} as
\begin{multline}
\label{eq:SWGaugeGroupAction}
\sA(L,H_L)\times W^{1,p}(W^+\otimes L) \ni
(A_L,\Psi)
\mapsto
s^*\left(A_L,\Psi\right)
:=
\left(s^*A_L, s^{-1}\Psi\right)
\\
= \left(A_L + (s^{-1}ds), s^{-1}\Psi\right)
\in \sA(L,H_L)\times W^{1,p}(W^+\otimes L).
\end{multline}
The quotient space of spin${}^c$ pairs on $\fs=\fs_0\otimes L$ is
\begin{equation}
\label{eq:SpincConfig}
\sC_\fs := \left.\left(\sA(L,h_L)\times W^{1,p}(W^+\otimes L)\right)\right/W^{2,p}(X,S^1),
\end{equation}
where $W^{2,p}(X,S^1)$ acts by \eqref{eq:SWGaugeGroupAction}, and we recall that $\sC_\fs$ is Hausdorff (see Morgan \cite[Corollary 4.5.4, p. 62]{MorganSWNotes}). We let
\begin{equation}
\label{eq:SpincConfigNonZero}
\sC_\fs^0 := \left\{[A_L,\Psi]\in\sC_\fs: \Psi\not\equiv 0\right\}
\end{equation}
denote the open subspace of $\sC_\fs$ comprising gauge-equivalence classes of non-zero-section \spinc pairs, and which is a Banach manifold (see Morgan \cite[Corollary 4.5.7, p. 64]{MorganSWNotes}) since $W^{2,p}(X,S^1)$ acts freely and properly on the open subspace of non-zero-section \spinc pairs (see Morgan \cite[Lemma 4.5.1, p. 61]{MorganSWNotes}).

We now assume that $\ft=(\rho\otimes\id_E,W\otimes E)$ is a spin${}^u$ structure of rank two and describe an embedding of $\sC_\fs$ into $\sC_\ft$. If $(A,\Phi)$ is split in the sense of Definition \ref{defn:Split_trivial_central-stabilizer_spinor_pair}, we can write $A=A_1\oplus A_2$, where $A_1$ and $A_2$ are unitary connections on $L_1$ and $L_2$, respectively, and $E = L_1\oplus L_2$ as a direct sum of Hermitian line bundles.
If $H$ is the Hermitian metric on $E$, then we write $H=H_1\oplus H_2$, where $H_i$ is the Hermitian metric on $L_i$ for $i=1,2$. Because $A$ induces the connection $A_d$ on $\det E$ by the definition of $\sA(E,H,A_d)$, we have $A_2=A_d\otimes A_1^*$ and $L_2\cong \det E\otimes L_1^*$. If $\fs=(\rho,W\otimes L_1)$, then there is a smooth embedding,
\begin{multline}
\label{eq:DefnOfIota}
\tilde\iota_{\fs,\ft}:
\sA(L_1,H_1)\times W^{1,p}(W^+\otimes L_1)
\ni (A_1,\Psi)
\\
\mapsto
\left(A_1\oplus (A_d\otimes A_1^*),\Psi\oplus 0\right) \in
\sA(E,H,A_d)\times W^{1,p}(W^+\otimes E).
\end{multline}
The embedding \eqref{eq:DefnOfIota} is gauge-equivariant with respect to the homomorphism
\begin{equation}
  \label{eq:GaugeGroupInclusion}
\varrho:W^{1,p}(X,S^1) \ni s
\mapsto
s\,\id_{L_1}\oplus s^{-1}\,\id_{L_2}
\in W^{2,p}(\SU(E)).
\end{equation}
According to \cite[Lemma 3.11, p. 93 and Lemma 3.16, p. 105]{FL2a}, the map \eqref{eq:DefnOfIota} descends to the quotient,
\begin{equation}
\label{eq:DefnOfIotaOnQuotient}
\iota_{\fs,\ft}:\sC_\fs\to \sC_\ft,
\end{equation}
and its restriction to $\sC_\fs^0$ is a smooth embedding $\sC_\fs^0\embed\sC_\ft^0$.

\section{The moduli space of non-Abelian monopoles}
\label{sec:ModuliOfSO3Monopole}
Recall that a pair $(A,\Phi) \in \sA(E,H,A_d)\times\Omega^0(W^+\otimes E)$ is a solution to the unperturbed non-Abelian monopole equations \cite[Equation (2.15)]{FL1}, \cite[Equation (2.32)]{FL2a} if
\begin{subequations}
\label{eq:SO(3)_monopole_equations}
\begin{align}
  \label{eq:SO(3)_monopole_equations_curvature}
  (F_A^+)_0 - \rho^{-1}(\Phi\otimes\Phi^*)_{00} &= 0,
  \\
  \label{eq:SO(3)_monopole_equations_Dirac}
  D_A\Phi &=0.
\end{align}
\end{subequations}
In \eqref{eq:PerturbedSO3MonopoleEquations}, the operator $D_A=\rho\circ \cov_A:C^\infty(W^+\otimes E)\to C^\infty(W^-\otimes E)$ is the Dirac operator\label{Dirac_operator} associated to the connection $A_W\otimes A$ on $W\otimes E$, where $A_W$ is a fixed spin connection in the sense of \eqref{eq:SpinConnection} on $(\rho,W)$; $\Phi^* \in \Hom(W^+\otimes E,\CC)$ is the pointwise Hermitian dual $\langle\cdot,\Phi\rangle_{W^+\otimes E}$ of $\Phi$; and $(\Phi\otimes\Phi^*)_{00}$ is the component of the section $\Phi\otimes\Phi^*$ of $i\fu(W^+\otimes E)$ lying in $\rho(\Lambda^+)\otimes\su(E)$ with respect to the splitting $\fu(W^+\otimes E)=i\underline{\RR}\oplus\su(W^+\otimes E)$ and the decomposition (see Feehan and Leness \cite[Equation (2.17), p. 67]{FL2a})
\[
\su(W^+\otimes E)
\cong
\rho(\Lambda^+)\oplus i\rho(\Lambda^+)\otimes_\RR\su(E)
\oplus \su(E),
\]
 of $\su(W^+\otimes E)$, where $\underline{\RR} := X\times\RR$.

\begin{rmk}[Perturbations of the non-Abelian monopole equations via generic geometric parameters]
\label{rmk:Perturbations_SO3_monopole_equations_generic_geometric_parameters}
We have stated the \emph{unperturbed} non-Abelian monopole equations in \eqref{eq:SO(3)_monopole_equations} for the sake of clarify, omitting mention of the fact that in order to obtain a regular moduli subspace
(as explained in Feehan \cite{FeehanGenericMetric} and Teleman \cite{TelemanGenericMetric}), one needs to consider the perturbed equations
\begin{subequations}
  \label{eq:PerturbedSO3MonopoleEquations}
\begin{align}
  \label{eq:PerturbedSO3MonopoleEquations_curvature}
  (F_A^+)_0 - \tau\rho^{-1}(\Phi\otimes\Phi^*)_{00} &= 0,
  \\
  \label{eq:PerturbedSO3MonopoleEquations_Dirac}
  D_A\Phi + \rho(\vartheta)\Phi &= 0.
\end{align}
\end{subequations}
and assume  a choice of \emph{generic} Riemannian metric $g$, perturbation $\tau\in\GL(\Lambda^+(X))$ (near the identity) of the term $\rho^{-1}(\Phi\otimes\Phi^*)_{00}$ in \eqref{eq:SO(3)_monopole_equations_curvature}, and perturbation $\vartheta \in \Omega^1(X,\CC)$ (near zero) of the Dirac operator $D_A$ in \eqref{eq:SO(3)_monopole_equations_Dirac}. In practice, the presence of these perturbations causes no difficulty due to our framework of approximately complex tangent spaces developed in Chapters \ref{chap:Analogue_Donaldson_symplectic_submanifold_criterion} and \ref{chap:Construction_circle-invariant_non-degenerate_two-form_II}.
\qed\end{rmk}

For Sobolev exponent $p \in (2,\infty)$, we observe that the perturbed non-Abelian monopole equations \eqref{eq:PerturbedSO3MonopoleEquations} define a smooth (in fact, analytic) map
\begin{multline}
\label{eq:PerturbedSO3MonopoleEquation_map}
  \sS:\sA(E,H,A_d)\times W^{1,p}(W^+\otimes E) \ni (A,\Phi)
  \\
  \mapsto
  \left((F_A^+)_0 - \tau\rho^{-1}(\Phi\otimes\Phi^*)_{00}, D_A\Phi + \rho(\vartheta)\Phi \right)
  \in L^p\left(\Lambda^+\otimes\su(E) \oplus W^-\otimes E\right)
\end{multline}
that is equivariant with respect to the action on the domain and codomain induced by the action of the Banach Lie group $W^{2,p}(\SU(E))$ of $W^{2,p}$ determinant-one, unitary automorphisms of $E$. The map $\sS$ in \eqref{eq:PerturbedSO3MonopoleEquation_map} becomes Fredholm upon restriction to a Coulomb-gauge slice,
(see
%TL11-28-2025: Updated to \cite[Section 9.3]{Feehan_Leness_introduction_virtual_morse_theory_so3_monopoles}?  Not sure what's being referred to here
\cite[Section 9.3]{Feehan_Leness_introduction_virtual_morse_theory_so3_monopoles}).

Define the \emph{moduli space of non-Abelian monopoles} by
\begin{equation}
\label{eq:Moduli_space_SO(3)_monopoles}
\sM_\ft := \left\{[A,\Phi] \in \sC_\ft: \sS(A,\Phi) = 0 \right\},
\end{equation}
where $\sS$ is as in \eqref{eq:PerturbedSO3MonopoleEquation_map}.
If $\ft=(\rho,W_{\can}\otimes E)$ where $\fs_{\can}=(\rho,W_{\can})$ is the canonical spin${}^c$ structure of Definition \ref{defn:Canonical_spinc_bundles}, then
\[
\sM_\ft=\sM(E,g,J,\omega,0),
\]
where $\sM(E,g,J,\omega,0)$ is the moduli space defined in \eqref{eq:Moduli_space_non-Abelian_monopoles_almost_Hermitian_Taubes_regularized} with parameter $r=0$.
Following the definitions of subspaces of $\sC_\ft$ in \eqref{eq:SpinuQuotientSpaceSubspaces}, we set
\begin{subequations}
 \label{eq:PU2MonopoleSubspaces}
 \begin{align}
   \label{eq:PU2MonopoleSubspace_irreducible}
   \sM^*_{\ft} &:= \sM_{\ft}\cap\sC^*_\ft,
   \\
   \label{eq:PU2MonopoleSubspace_non-zero-section}
   \sM^0_{\ft} &:= \sM_{\ft}\cap\sC^0_\ft,
   \\
   \label{eq:PU2MonopoleSubspace_irreducible_non-zero-section}
   \sM^{*,0}_\ft &:= \sM_{\ft}\cap\sC_\ft^{*,0},
\end{align}
\end{subequations}
where the subspaces $\sC^*_\ft$, and $\sC^0_\ft$, and $\sC_\ft^{*,0}$ are defined in \eqref{eq:SpinuQuotientSpaceSubspaces}. If $\ft=(\rho,W_{\can}\otimes E)$ as above, then the subspaces \eqref{eq:PU2MonopoleSubspaces} equal the corresponding subspaces of $\sM(E,g,J,\omega,0)$ defined in \eqref{eq:Moduli_space_non-Abelian_monopoles_almost_Hermitian_Taubes_regularized_non-zero_or_non-split}. We recall the

\begin{prop}[Expected dimension of the moduli space of non-Abelian monopoles]
\label{prop:Expected_dimension_moduli_space_non-Abelian_monopoles}
(See Feehan and Leness \cite[Proposition 2.28, p. 313]{FL1}.)
Let $\ft=(\rho,W\otimes E)$, where $(E,H)$ is a rank-two, smooth Hermitian vector bundle and $(\rho,W)$ is a spin${}^c$ structure over a closed, four-dimensional, oriented, smooth Riemannian manifold $(X,g)$. Then the expected dimension of the moduli space \eqref{eq:Moduli_space_SO(3)_monopoles} of non-Abelian monopoles is
\begin{align}
\label{eq:Transv}
\expdim \sM_{\ft} &\,= d_a(\ft)+2n_a(\ft), \quad\text{where }
\\
\notag
d_a(\ft)
&:= -2p_1(\ft)- \frac{3}{2}(e(X)+\sigma(X)) \text{ and }
n_a(\ft) := \frac{1}{4}(p_1(\ft)+c_1(\ft)^2-\sigma(X)),
\end{align}
and $e(X)$ is the Euler characteristic and $\si(X)$ is the signature of $X$, and is independent of the choice of parameters $(g,\tau,\vartheta)$ in the non-Abelian monopole equations \eqref{eq:PerturbedSO3MonopoleEquations}.
\end{prop}

\begin{thm}[Transversality for the moduli space of non-Abelian monopoles with generic parameters]
\label{thm:Transv}
(See Feehan \cite[Theorem 1.1, p. 910]{FeehanGenericMetric} and Teleman \cite[Theorem 3.19, p. 413]{TelemanGenericMetric}.)
Let $\ft=(\rho,W\otimes E)$, where $(E,H)$ is a rank-two, smooth Hermitian vector bundle and $(\rho,W)$ is a spin${}^c$ structure over a closed, four-dimensional, oriented, smooth Riemannian manifold $(X,g)$. If the parameters $(g,\tau,\vartheta)$ appearing in \eqref{eq:PerturbedSO3MonopoleEquations} are generic in the sense of \cite{FeehanGenericMetric}, then $\sM^{*,0}_{\ft}$ is an embedded smooth submanifold of the Banach manifold $\sC^{*,0}_{\ft}$ that has finite dimension equal to the expected dimension \eqref{eq:Transv}.
\end{thm}

While $\sM_\ft$ is not in general compact, it admits an Uhlenbeck compactification as described in Feehan and Leness \cite[Definition 4.19, p. 350]{FL1}.

\section{The stratum of zero-section pairs}
\label{sec:ZeroSectionPairs}

For additional details concerning the construction and properties of the moduli space of anti-self-dual connections, we refer the reader to Donaldson and Kronheimer \cite{DK}, Friedman and Morgan \cite{FrM}, Freed and Uhlenbeck \cite{FU}, Kronheimer and Mrowka \cite{KMStructure}, and Lawson \cite{Lawson}.

We continue to assume that there is a fixed, smooth unitary connection $A_d$ on the smooth Hermitian line bundle $\det E$ defined by the smooth Hermitian vector bundle $(E,H)$ and that $(X,g)$ is a closed, four-dimensional, oriented, smooth Riemannian manifold. From equation \eqref{eq:PerturbedSO3MonopoleEquations}, we see that the stratum of $\sM_{\ft}$ represented by pairs with zero spinor is identified with
\begin{equation}
\label{eq:ASDModuliSpace}
M_\kappa^w(X,g) := \left.\{A\in\sA(E,H,A_d): (F_A^+)_0 = 0\}\right/W^{2,p}(\SU(E)),
\end{equation}
the \emph{moduli space of projectively $g$-anti-self-dual connections} on  $E$, where $\ka:=-\frac{1}{4} p_1(\ft)$ and $w=c_1(E)$. If $b^+(X)>0$, then by Donaldson and Kronheimer \cite[Corollary 4.3.18, p. 149]{DK} and Kronheimer and Mrowka \cite[Lemma 2.4 and Corollary 2.5]{KMStructure} for a generic Riemannian metric $g$, the space $M_\kappa^w(X,g)$ is a smooth manifold of the expected dimension \eqref{eq:Expected_dimension_moduli_space_ASD_connections},
\[
  \expdim M_\kappa^w(X,g) = d_a(\ft) = -2p_1(\ft) - \frac{3}{2}(e(X)+\sigma(X)),
\]
where $d_a(\ft)$ was defined in \eqref{eq:Transv}.
As explained in \cite[Section 3.4, p. 96]{FL2a}, it is desirable to choose $w\pmod{2}$ so as to exclude points in $\bar\sM_{\ft}$ with associated flat $\SO(3)$ connections, so we have a \emph{disjoint} union,
\begin{equation}
\label{eq:StratificationCptPU(2)Space}
\bar\sM_{\ft}
\cong
\bar\sM_{\ft}^{*,0} \sqcup \bar M_\kappa^w \sqcup \bar\sM_{\ft}^{\red},
\end{equation}
where $\bar\sM_{\ft}^*\subset\bar\sM_{\ft}$ is the subspace represented by triples whose associated $\SO(3)$ connections are non-split, $\bar\sM_{\ft}^0\subset\bar\sM_{\ft}$ is the subspace represented by triples whose spinors are not identically zero, $\bar\sM_{\ft}^{*,0} = \bar\sM_{\ft}^{*}\cap\bar\sM_{\ft}^{0}$, while $\bar\sM_{\ft}^{\red}\subset\bar\sM_{\ft}$ is the subspace $\bar\sM_{\ft}-\bar\sM_{\ft}^*$ represented by triples whose associated $\SO(3)$ connections are split. We recall the

%PF3-22-2025 Below is poorly explained
\begin{defn}[Good cohomology classes]
\label{defn:Good}
(See Feehan and Leness \cite[Definition 3.20, p. 169]{FL2b}.)
A class $v\in H^2(X;\ZZ/2\ZZ)$ is \emph{good} if no integral lift of $v$ is torsion.
\qed\end{defn}

If $w\pmod{2}$ is good, then for \emph{generic} metrics the union \eqref{eq:StratificationCptPU(2)Space} is disjoint, as desired. In practice, rather than constraining $w\pmod{2}$ itself, we use the blow-up trick of Morgan and Mrowka \cite{MorganMrowkaPoly}, replacing $X$ with the smooth blow-up, $X\#\overline{\CC\PP}^2$, and replacing $w$ by $w+\PD[e]$ (where $e\in H_2(X;\ZZ)$ is the exceptional class and $\mathrm{PD}[e]$ denotes its Poincar{\'e} dual), noting that $w+\PD[e]\pmod{2}$ is always good.

%PF10-18-2024 Also not a stratum
\section{The stratum of reducible pairs}
\label{sec:ReduciblePairs}

We now describe the relationship between the moduli spaces of Seiberg--Witten invariants and the moduli space of non-Abelian monopoles. For additional details concerning the definition of Seiberg--Witten invariants, we refer the reader to Kronheimer and Mrowka \cite{KMThom, KMBook}, Morgan \cite{MorganSWNotes}, Nicolaescu \cite{NicolaescuSWNotes}, Salamon \cite{SalamonSWBook}, and Witten \cite{Witten}.

We introduce the Seiberg--Witten equations and define the Seiberg--Witten moduli space in Section \ref{subsec:SWMonopoles}.  Results on compactness, transversality, expected dimension, the existence of zero-section Seiberg--Witten monopoles, and the definitions of the Seiberg--Witten invariant, basic classes, simple type, and a blow-up formula for the Seiberg--Witten invariant appear in Section \ref{subsec:SWInvariants}. Finally, in Section \ref{subsec:RedPU2Monopole}, we characterize subspaces of split non-Abelian monopoles as the image of spaces of Seiberg--Witten monopoles with a particular perturbation under embeddings of the kind discussed in Section \ref{subsec:SpincQuotientSpaceInSpinuSpace}.

\subsection{Seiberg--Witten monopoles}
\label{subsec:SWMonopoles}
As in the definition of spin${}^c$ pairs \eqref{eq:TranslateSWNotation}, we fix a spin${}^c$ structure $\fs_0=(\rho,W)$ and a spin connection $A_W$ on $W$. Let $(L,H_L)$ be a Hermitian line bundle over the closed, oriented, smooth Riemannian four-manifold $(X,g)$. By analogy with Feehan and Leness \cite[Equations (2.55), p. 76, and (2.57), p. 77]{FL2a}, we call a pair $(A_L,\Phi_1)\in \sA(L,h_L)\times W^{1,p}(W^+\otimes L)$ a solution to the \emph{perturbed Seiberg--Witten monopole equations} on the spin${}^c$ structure $\fs=\fs_0\otimes L$ if
\begin{subequations}
\label{eq:SeibergWitten_Prelim}
\begin{align}
  \label{eq:SeibergWitten_curvature_Prelim}
  \tr_{W^+}F^+_{A_W} + 2F^+_{A_L} - \tau\rho^{-1}(\Phi_1\otimes\Phi_1^*)_{0} - F^+_{A_{\Lambda}} &=0,
  \\
  \label{eq:SeibergWitten_Dirac_Prelim}
  D_{A_L}\Phi_1 + \rho(\vartheta)\Psi &=0,
\end{align}
\end{subequations}
where $\tr_{W^+}:\fu(W^+)\to i\ubarRR$ is defined by the trace on $2\times 2$ complex matrices, $(\Psi\otimes\Psi^*)_0$ is the component of the section $\Psi\otimes\Psi^*$ of $i\fu(W^+)$ contained in $i\su(W^+)$, and $D_{A_L}:C^\infty(W^+\otimes L)\to C^\infty(W^-\otimes L)$ is the Dirac operator defined by the spin connection $A_W\otimes A_L$,
\label{page:SpincDiracOperator}
and $A_\Lambda := A_d\otimes A_{\det W^+}$ is a unitary connection on the Hermitian line bundle $\det W^+\otimes \det E$ with first Chern class $\Lambda\in H^2(X;\ZZ)$. The perturbations --- in particular the term $F_{A_\Lambda}^+$ in \eqref{eq:SeibergWitten_curvature} --- are chosen so that if $\ft$ is a spin${}^u$ structure with $c_1(\ft)=\Lambda$ and the embedding $\tilde\iota_{\fs,\ft}$ defined in \eqref{eq:DefnOfIota} exists, then $\tilde\iota_{\fs,\ft}$ maps solutions to equation \eqref{eq:SeibergWitten_Prelim} to solutions to the non-Abelian monopole equations \eqref{eq:PerturbedSO3MonopoleEquations} that are split in the sense of Definition \ref{defn:Split_trivial_central-stabilizer_spinor_pair} (see Feehan and Leness \cite[Lemma 3.12, p. 95]{FL2a}).

The perturbed Seiberg--Witten equations \cite[Equation (2.55), p. 76, and Equation (2.57), p. 77]{FL2a} are expressed in terms of the spin connection $B := A_W\otimes A_L$ on $(\rho\otimes \id_L,W\otimes L)$ given by the isomorphism of affine spaces \eqref{eq:TranslateSWNotation}. The equalities,
\[
  \tr_{W^+\otimes L} F_B^+
  =
  \tr_{W^+\otimes L} F^+_{A_W\otimes A_L}
  =
 \tr_{W^+} F^+_{A_W}+2F^+_{A_L},
\]
show that the curvature terms appearing in \eqref{eq:SeibergWitten_curvature} are equivalent to those in
\cite[Equation (2.55), p. 76]{FL2a}. The Dirac operator $D_{A_L}$ in \eqref{eq:SeibergWitten_Dirac} is written as $D_B$ in \cite[Equation (2.55), p. 76]{FL2a}. By observing that
\[
  F_{A_\Lambda} = F_{A_{\det W^+}\otimes A_d} = F_{A_{\det W^+}} + F_{A_d}
  \quad\text{and}\quad
  \tr_{W^+}F_{A_W} = F_{A_{\det W^+}},
\]
where the latter equality follows from Kobayashi \cite[Equation (1.5.19), p. 17]{Kobayashi_differential_geometry_complex_vector_bundles}, we see that Equation \eqref{eq:SeibergWitten_curvature_Prelim} simplifies to
\begin{equation}
  \label{eq:Seiberg-Witten_curvature_simplified_prelim}
  F_{A_1}^+ - \frac{1}{2}\rho^{-1}(\Phi_1\otimes\Phi_1^*)_{0} - \frac{1}{2}F_{A_d}^+ =0.
\end{equation}
We let
\begin{equation}
\label{eq:Moduli_space_Seiberg-Witten_monopoles}
  M_{\fs}
  :=
  \left.\left\{(A_L,\Phi_1)\in\sA(L,H_L)\times W^{1,p}(W^+\otimes L): \text{Equation \eqref{eq:SeibergWitten_Prelim} holds} \right\}
    \right/W^{2,p}(X,S^1)
\end{equation}
denote the \emph{moduli space of Seiberg--Witten monopoles}, where $W^{2,p}(X,S^1)$ acts on the affine space $\sA(L,h_L)\times W^{1,p}(W^+\otimes L)$ as in \eqref{eq:SWGaugeGroupAction}.

\subsection{Seiberg--Witten invariants, basic classes, and simple type}
\label{subsec:SWInvariants}
We recall the following basic results for the moduli space of Seiberg--Witten monopoles.

\begin{prop}[Existence of zero-section Seiberg--Witten monopoles]
\label{prop:SWModuliSpaceZeroSectionCriterion}
(See Morgan \cite[Proposition 6.3.1, p. 91]{MorganSWNotes}.)
Let $\fs$ and $\ft$ be a spin${}^c$ and a rank-two spin${}^u$ structure, respectively, over a closed, four-dimensional,oriented, smooth Riemannian manifold $(X,g)$. If $b^+(X)>0$ and $c_1(\fs)-c_1(\ft)$ is not a torsion class in $H^2(X;\ZZ)$, then there is an open, dense subspace of Riemannian metrics on $X$ in the $C^r$ topology (where $r\geq 3$) such that if $g$ belongs to this subspace, then the moduli space $M_\fs$ does not contain the gauge-equivalence class of a zero-section pair.
\end{prop}

\begin{rmk}[Open and dense subspace of Riemannian metrics]
In Morgan \cite[Proposition 6.3.1, p. 91]{MorganSWNotes}, the conclusion of Proposition \ref{prop:SWModuliSpaceZeroSectionCriterion} is stated to hold for generic Riemannian metrics, meaning Riemannian metrics lying in a Baire set, as described by Freed and Uhlenbeck in \cite[Theorem 3.17, p. 59]{FU}). Applying the argument appearing in Donaldson and Kronheimer \cite[Corollary 4.3.15, p. 148]{DK} gives the conclusion of Proposition \ref{prop:SWModuliSpaceZeroSectionCriterion} for an open and dense subspace of Riemannian metrics.
\qed\end{rmk}

\begin{thm}[Compactness for the moduli space of Seiberg--Witten monopoles]
\label{thm:Compactness_SW}
(See Feehan and Leness \cite[Proposition 2.15, p. 78]{FL2a} and Morgan \cite[Corollary 5.3.7, p. 85]{MorganSWNotes}.)
Let $\fs=(\rho,W)$ be a spin${}^c$ structure over a closed, four-dimensional, oriented, smooth Riemannian manifold $(X,g)$. Then the moduli space $M_\fs$ of Seiberg--Witten monopoles \eqref{eq:Moduli_space_Seiberg-Witten_monopoles} is a compact subset of $\sC_\fs$.
\end{thm}

\begin{prop}[Expected dimension of the moduli space of Seiberg--Witten monopoles]
\label{prop:Expected_dimension_moduli_space_Seiberg-Witten_monopoles}
(See Morgan \cite[Corollary 4.6.2, p. 67]{MorganSWNotes}.)
Continue the hypotheses of Theorem \ref{thm:Compactness_SW}. Then the expected dimension of the moduli space $M_\fs$ of Seiberg--Witten monopoles \eqref{eq:Moduli_space_Seiberg-Witten_monopoles} is
\begin{equation}
\label{eq:DimSW}
d_s(\fs)
=
\dim M_{\fs}
=
\frac{1}{4}(c_1(\fs)^2 -2e(X) -3\sigma(X)).
\end{equation}
where $e(X)$ is the Euler characteristic and $\sigma(X)$ is the signature of $X$, and is independent of the choice of parameters $(g,\tau,\vartheta)$ in the Seiberg--Witten monopole equations \eqref{eq:Moduli_space_Seiberg-Witten_monopoles}.
\end{prop}

\begin{prop}[Transversality for the moduli space of Seiberg--Witten monopoles with generic parameters]
\label{prop:Transv_SW}
(See Feehan and Leness \cite[Proposition 2.16, p. 79]{FL2a}.)
Continue the hypotheses of Theorem \ref{thm:Compactness_SW}. If the parameters $(g,\tau,\vartheta)$ in the Seiberg--Witten monopole equations \eqref{eq:SeibergWitten_Prelim} are generic in the sense of \cite{FeehanGenericMetric}, then $M_\fs^0:=M_\fs\cap\sC_\fs^0$ is an embedded, orientable, smooth submanifold of the Banach manifold $\sC_\fs^0$ that has finite dimension equal to the expected dimension \eqref{eq:DimSW}.
\end{prop}

Let $\widetilde X=X\#\overline{\CC\PP}^2$ denote the smooth blow-up of $X$ with exceptional class $e\in H_2(\widetilde X;\ZZ)$ and denote its Poincar\'e dual by $\PD[e]\in H^2(\widetilde X;\ZZ)$ (see Gompf and Stipsicz \cite[Definition 2.2.7, p. 43]{GompfStipsicz}). Let $\fs^\pm=(\tilde\rho,\tilde W)$ denote the \spinc structure on $\widetilde X$ with $c_1(\fs^\pm)=c_1(\fs)\pm \PD[e]$ obtained by splicing the \spinc structure $\fs=(\rho,W)$ on $X$ with the \spinc structure on $\overline{\CC\PP}^2$ with first Chern class $\pm \PD[e]$. (See Feehan and Leness \cite[Section 4.5, p. 200]{FL2b} or Salamon \cite[Section 12.4]{SalamonSWBook} for an explanation of the relation between \spinc structures on $X$ and $\widetilde X$.) Now
$$
c_1(\fs)\pm \PD[e]-\Lambda \in H^2(\widetilde X;\ZZ)
$$
is not a torsion class and so --- for $b^+(X)>0$, generic Riemannian metrics on $X$ and related metrics on the connected sum $\widetilde X$ --- the moduli spaces $M_{\fs^\pm}$ contain no zero-section pairs. Thus, for our choice of generic perturbations, the moduli spaces $M_{\fs^\pm}$ are compact, oriented, smooth manifolds, both of dimension $\dim M_{\fs}$.

The \emph{Seiberg--Witten invariant} $\SW_{X}(\fs)$
\label{page:SW_Invariant}
defined by a \spinc structure $\fs$ over $X$
%PF9-22-2024 "most easily" is a pretty vague qualifier. It's defined for all b^+ - b_1 parities, juust zero by definition when that's even.
is most easily defined when $b^+(X)-b_1(X)\equiv 1\pmod 2$; we refer the reader to Feehan and Leness \cite[Section 4.1]{FL2b} for its definition and references to the development of its properties. As in Definition \ref{defn:Seiberg-Witten_basic_class}, $c_1(\fs)$ is a \emph{Seiberg--Witten basic class}\label{page:Basic_class} if $\SW_{X}(\fs) \neq 0$.   We will denote
\begin{equation}
\label{eq:SetOfBasicClasses}
B(X):=\{c_1(\fs): \SW_X(\fs)\neq 0\}.
\end{equation}
Note that $B(X)$ contains the basic classes of Definition \ref{defn:Seiberg-Witten_basic_class}
but possibly also contains classes $c_1(\fs)$ where $\SW_X(\fs)\neq 0$ but $d(\fs)>0$.
A four-manifold $X$ has \emph{Seiberg--Witten simple type} \label{page:Seiberg_Witten_simple_type} if all spin${}^c$ structures with $c_1(\fs)\in B(X)$ satisfy
\begin{equation}
  \label{eq:SWSimpleType}
  d_s(\fs)=0,
\end{equation}
where $d_s(\fs)$ is as in \eqref{eq:DimSW}. Versions of the following  have appeared in Fintushel and Stern \cite[Theorem 1.4]{FSTurkish}, Nicolaescu \cite[Theorem 4.6.7]{NicolaescuSWNotes}, and
Fr{\o}yshov \cite[Corollary 14.1.1]{Froyshov_2008}.

\begin{thm}[Blow-up formula for Seiberg--Witten invariants]
\label{thm:SWBlowUp}
Let $X$ be a standard four-manifold and let $\widetilde X=X\#\overline{\CC\PP}^2$. Write $\PD[e]\in H^2(\widetilde X;\ZZ)$ for the Poincar\'e dual of the exceptional sphere. Then $X$ has Seiberg--Witten simple type if and only if $\widetilde X$ does. If $X$ has Seiberg--Witten simple type then
\begin{equation}
\label{eq:BlowUpBasics}
B(\widetilde X)=\{c_1(\fs_0)\pm \PD[e]: c_1(\fs_0)\in B(X)\}.
\end{equation}
If $c_1(\fs_0)\in B(X)$ and $\fs^\pm$ is a spin${}^c$ structure on $\widetilde X$ with
$c_1(\fs^\pm)=c_1(\fs_0)\pm \PD[e]$, then $\SW_{\widetilde X}(\fs^\pm)=\SW_X(\fs)$.
\end{thm}

\subsection{Split non-Abelian monopoles}
\label{subsec:RedPU2Monopole}
Let $\ft$ be a rank-two spin${}^u$ structure on a closed, smooth Riemannian four-manifold. We now describe the equivalence classes of pairs in $\sM_\ft$ which are split in the sense of Definition \ref{defn:Split_trivial_central-stabilizer_spinor_pair}. By \cite[Lemma 3.13, p. 96]{FL2a} the restriction of the map $\iota_{\fs,\ft}$ defined in \eqref{eq:DefnOfIota} to $M_\fs^0$ is a topological embedding $M^0_\fs\embed\sM_\ft$ where $M_{\fs}^0:=M_{\fs}\cap\sC_{\fs}^0$. If $w_2(\ft)\neq 0$ or $b_1(X)=0$, then \cite[Lemma 3.13, p. 96]{FL2a} implies that $\iota_{\fs,\ft}$ is an embedding. The image of $\iota_{\fs,\ft}$ is represented by points in $\sM_\ft$ represented by pairs which are split with respect to the splitting \eqref{eq:BasicSplitting}. Henceforth, we shall not distinguish between $M_\fs$ and its image in $\sM_{\ft}$ under this embedding.
 We recall the following cohomological criterion for when the continuous embedding \eqref{eq:LowerLevelInclusionOFReducibles} exists.

\begin{lem}[Existence of embeddings of moduli spaces of Seiberg--Witten monopoles into moduli space of ideal non-Abelian monopoles]
(See
%TL11-28-2025: Updated to \cite[Lemma6.6.8]{Feehan_Leness_introduction_virtual_morse_theory_so3_monopoles}
\cite[Lemma 6.6.8]{Feehan_Leness_introduction_virtual_morse_theory_so3_monopoles}.)
\label{lem:SetOfReducibles}
Let $\ft=(\rho,W\otimes E)$ be a  spin${}^u$ structure over a closed, oriented, smooth, Riemannian four-manifold $(X,g)$. If
% $\ell$ is a non-negative integer and
\begin{equation}
\label{eq:DefineReduciblesEmbedded}
\Red_0(\ft) := \left\{c_1(\fs): \fs \in \Spinc(X) \text{ such that } (c_1(\fs)-c_1(\ft))^2=p_1(\su(E))\right\},
\end{equation}
then there is a continuous embedding 
\begin{equation}
\label{eq:LowerLevelInclusionOFReducibles}
M_{\fs} \to \sM_{\ft},
\end{equation}
defined by the map $\iota_{\fs,\ft}$ from \eqref{eq:DefnOfIota} 
 exists if and only if $c_1(\fs)\in\Red_0(\ft)$.
\end{lem}

\begin{rmk}[Conventions on splitting of bundles in \spinu structures]
\label{rmk:ComparingLineBundlesFromFl2a}
As we shall often refer to constructions from Feehan and Leness \cite{FL2a, FL2b}, it is worth noting the following difference between the notation used in this monograph and that in \cite{FL2a,FL2b}.  In \cite{FL2a}, the splitting of $W\otimes E$ is written as $W'\oplus W'\otimes L$, whereas the splitting of $W\otimes E$ given by \eqref{eq:BasicSplitting} would be $W\otimes L_1\oplus W\otimes L_2$.  To convert between the two conventions, use $W'=W\otimes L_1$ and $L=L_2\otimes L_1^*$.
\qed\end{rmk}

\subsection{Split non-Abelian monopoles with a Taubes perturbation}
We now assume that $(X,g,J,\omega)$ is an almost Hermitian, closed, smooth four-manifold
and describe the solutions $(A,\varphi,\psi)$ of the non-Abelian monopole equations, with singular and regularized Taubes perturbation, when $(A,\varphi,\psi)$ is split with respect to a decomposition $E=L_1\oplus L_2$.
Such triples are in the image of the embedding $\tilde\iota_{\fs,\ft}$ of \eqref{eq:DefnOfIota} where
$\ft=(\rho_{\can}\otimes\id_E,W_{\can}\otimes E)$ and $\fs=(\rho_{\can}\otimes\id_{L_1},W_{\can}\otimes L_1)$.

We define the \emph{Seiberg--Witten equations with a Taubes perturbation} by, for 
$(A_1,\phi_1,\psi_1)\in\sA(L_1,H_1)\times W^{1,p}(L_1\oplus \Lambda^{0,2}\otimes L_1)$,
\begin{subequations}
\label{eq:SeibergWittenTaubesPert_Prelim}
\begin{align}
  \label{eq:SeibergWittenTaubesPert_curvature(1,1)_Prelim}
\Lambda_\omega F_{A_1}-\frac{1}{2}\Lambda_\omega F_{A_d} - \frac{i}{4}\left(|\varphi_1|_{L_1}^2 -|\psi_1|_{\Lambda^{0,2}(E)}^2\right) +\frac{ir}{4} &=0,
\\
  \label{eq:SeibergWittenTaubesPert_curvature(0,2)_Prelim}
F_{A_1}^{0,2}-\frac{i}{2}(\psi_1\otimes \phi_1^*)&=0,
  \\
  \label{eq:SeibergWittenTaubesPert_Dirac_Prelim}
 \bar\partial_{A_1}\varphi_1 + \bar\partial_{A_1}^*\psi_1 +\frac{1}{4\sqrt 2}\rho_{\can}(\Lambda_\omega d\omega)(\varphi_1+\psi_1)&=0.
\end{align}
\end{subequations}
The \emph{Seiberg--Witten equations with a regularized Taubes perturbation} are
\begin{subequations}
\label{eq:SeibergWittenRegTaubesPert_Prelim}
\begin{align}
  \label{eq:SeibergWittenRefTaubesPert_curvature(1,1)_Prelim}
\Lambda_\omega F_{A_1}-\frac{1}{2}\Lambda_\omega F_{A_d} - \frac{i}{4}\left(|\varphi_1|_{L_1}^2 -|\psi_1|_{\Lambda^{0,2}(L_1)}^2\right) +\frac{ir}{4}\frac{|\psi|_{\Lambda^{0,2}(L_1)}^2}{\gamma^2+|\psi|_{\Lambda^{0,2}(L_1)}^2} &=0,
\\
  \label{eq:SeibergWittenRefTaubesPert_curvature(0,2)_Prelim}
F_{A_1}^{0,2}-\frac{1}{2}F_{A_d}^{0,2}-\frac{1}{2}(\psi_1\otimes \phi_1^*)_0&=0,
  \\
  \label{eq:SeibergWittenRegTaubesPert_Dirac_Prelim}
 \bar\partial_{A_1}\varphi_1 + \bar\partial_{A_1}^*\psi_1 +\frac{1}{4\sqrt 2}\rho_{\can}(\Lambda_\omega d\omega)(\varphi_1+\psi_1)&=0.
\end{align}
\end{subequations}
For $\fs$ the spin${}^c$ structure given by $\fs=(\rho_{\can}\otimes \id_{L_1},W_{\can}\otimes L_1)$, define
\begin{equation}
\label{eq:SWModuliForRegTaubesPert}
M_\fs(r):=
\{(A,\varphi_1,\psi_1)\in\sA(L_1,H_1)\times W^{1,p}(L_1\oplus \Lambda^{0,2}\otimes L_1):
\text{$(A,\varphi_1,\psi_1)$ satisfies \eqref{eq:SeibergWittenRegTaubesPert_Prelim}} \}/W^{2,p}(X,S^1).
\end{equation}
We can then identify $M_\fs(r)$ with the points in $\sM(E,g,J,\omega,r)$ which are represented by split triples.

\begin{lem}[Split triples in $\sM(E,g,J,\omega,r)$]
\label{lem:SplitRegTaubesPert}
Let $(E,H)$ be a rank-two Hermitian vector bundle over a smooth, closed almost Hermitian four-manifold $(X,J,g,\omega)$.
If $[A,\varphi,\psi]\in \sM(E,g,J,\omega,r)$ admits a representative $(A,\varphi,\psi)$ which is split with respect to a decomposition $E=L_1\oplus L_2$ where $(L_i,H_i)$ is a complex rank-one Hermitian bundle over $X$, 
then $[A,\varphi,\psi]\in \iota_{\ft,\fs}(M_\fs(r))$ where the embedding $\iota_{\ft,\fs}$ is defined in \eqref{eq:DefnOfIotaOnQuotient},  $\ft=(\rho_{\can}\otimes \id_{E},W_{\can}\otimes E)$, and $\fs=(\rho_{\can}\otimes \id_{L_1},W_{\can}\otimes L_1)$.
\end{lem}

\begin{proof}
The image of $\iota_{\ft,\fs}:\sC_\fs\to\sC_\ft$ is the gauge-equivalence class of triples with a representative which is split with respect to the decomposition $E=L_1\oplus L_2$.  Thus, we need only show that the restriction of the non-Abelian monopole equations with a regularized Taubes perturbation to the image of $\tilde\iota_{\ft,\fs}$ is equivalent to the system of equations \eqref{eq:SeibergWittenRegTaubesPert_Prelim} defining $M_\fs(r)$.

Following the argument of \cite[Lemma 3.12, p. 95]{FL2a}, we observe that if $A=A_1\oplus A_2$ where $A_i$ is a unitary connection on $L_i$, if $\varphi=(\varphi_1, 0)$
and $\psi=(\psi_1,0)$ where $\varphi_1\in \Omega^0(L_1)$ and $\psi_1\in\Omega^{0,2}(L_1)$, then
$(F_A)_0$, $i(\phi\otimes\phi^*)_0$, and $i\star(\psi\otimes\psi^*)_0)$ all take values in the $\Omega^+(i\ubarRR)$ component of
\[
\Omega^+(\su(E)))\cong \Omega^+(i\ubarRR)\oplus \Omega^+(L_1\otimes L_2^*).
\]
Because $A = A_1 \oplus A_2 = A_1 \oplus (A_d\otimes A_1^*)$, we have
\[
  F_A = F_{A_1} \oplus F_{A_d\otimes A_1^*}
  =  F_{A_1} \oplus \left(F_{A_d} + F_{A_1^*}\right)
  =  F_{A_1} \oplus \left(F_{A_d} - F_{A_1}\right).
\]
Therefore, $\tr_E F_A = F_{A_d}$, as expected, and
\[
  (F_A)_0 = F_A - \frac{1}{2}(\tr_E F_A)\,\id_E =  F_A - \frac{1}{2}F_{A_d}\,\id_E
  =
  \left(F_{A_1} - \frac{1}{2}F_{A_d}\right) \oplus \left(\frac{1}{2}F_{A_d} - F_{A_1}\right).
\]
Note that
\[
  (\varphi\otimes\varphi^*)_0
  =
  \varphi\otimes\varphi^* - \frac{1}{2}\tr_E(\varphi\otimes\varphi^*)\,\id_E
  =
  \varphi\otimes\varphi^* - \frac{1}{2}|\varphi|_E^2\,\id_E
\]
and when $\varphi = \varphi_1 \in \Omega^0(L_1)$, then $(\varphi_1\otimes\varphi_1^*)s = \varphi_1\langle s,\varphi_1\rangle_{L_1} = s|\varphi|_{L_1}^2$ and so
\[
  (\varphi\otimes\varphi^*)_0 = \frac{1}{2}|\varphi_1|_{L_1}^2\left(\id_{L_1} \oplus -\id_{L_2}\right)
\]
Similarly, when $\psi=\psi_1\in\Omega^{0,2}(L_1)$,
%PF8-26-2024 Check!
\[
  \star(\psi\otimes\psi^*)_0 = \frac{1}{2}|\psi_1|_{\Lambda^{0,2}(L_1)}^2\left(\id_{L_1} \oplus -\id_{L_2}\right).
\]
Thus, when $(A,\varphi,\psi)=\tilde\iota_{\ft,\fs}(A_1,\varphi_1,\psi_1)$, the unperturbed $(1,1)$-component of the curvature equation \eqref{eq:SO(3)_monopole_equations_(1,1)_curvature_intro} is equivalent to
\[
\Lambda_\omega F_A - \frac{1}{2}\Lambda_\omega F_{A_d}
-\frac{i}{4}\left(|\varphi_1|_{L_1}^2-|\psi|_{\Lambda^{0,2}(L_1)}\right)=0.
\]
According to the definition \eqref{eq:Definition_wp_intro} of $\wp(\psi)$, when $\psi=\psi_1\in\Omega^{0,2}(L_1)$ we have
\[
  \wp(\psi)
  =
  2|\psi|_{\Lambda^{0,2}(E)}^{-2}\star(\psi\otimes\psi^*)_0
  =
  \left(\id_{L_1} \oplus -\id_{L_2}\right)
  \quad\text{a.e. on } X.
\]
Hence,
\[
  \frac{ir}{8}\wp(\psi)\,\omega
  =
  \frac{ir}{8}\left(\id_{L_1} \oplus -\id_{L_2}\right)\omega
  \quad\text{a.e. on } X.
\]
Using $\Lambda_\omega\omega=2$ and the definition of the regularized Taubes perturbation, we see that
when $(A,\varphi,\psi)=\tilde\iota_{\ft,\fs}(A_1,\varphi_1,\psi_1)$, the regularly perturbed $(1,1)$-component of the curvature equation \eqref{eq:SO(3)_monopole_equations_(1,1)_curvature_perturbed_intro_regular} is equivalent to
\eqref{eq:SeibergWittenRefTaubesPert_curvature(1,1)_Prelim}.

Equations \eqref{eq:SO(3)_monopole_equations_(0,2)_curvature_perturbed_intro_regular} and \eqref{eq:SO(3)_monopole_equations_Dirac_almost_Hermitian_perturbed_intro_regular} are equivalent to
\eqref{eq:SeibergWittenRefTaubesPert_curvature(0,2)_Prelim} and \eqref{eq:SeibergWittenRegTaubesPert_Dirac_Prelim} respectively when $(A,\varphi,\psi)=\tilde\iota_{\ft,\fs}(A_1,\varphi_1,\psi_1)$.  This completes the proof that
the system \eqref{eq:SO(3)_monopole_equations_almost_Hermitian_perturbed_intro_regular} is equivalent to
the system \eqref{eq:SeibergWittenRegTaubesPert_Prelim} when $(A,\varphi,\psi)=\tilde\iota_{\ft,\fs}(A_1,\varphi_1,\psi_1)$
and thus completes the proof of the lemma.
\end{proof}

%TL5-26-2025: Following taken from Ch. 8 discussion
\begin{rmk}[Split solutions of the perturbed non-Abelian monopole equations and Taubes' estimates]
\label{rmk:SplitPerturbedSolutions}
If the unitary triple $(A,\varphi,\psi)=\tilde\iota_{\ft,\fs}(A_1,\varphi_1,\psi_1)$ in the statement of Lemma \ref{lem:SplitRegTaubesPert} is a solution of
the non-Abelian monopole equations with a \emph{singular} Taubes perturbation, \eqref{eq:SO(3)_monopole_equations_almost_Hermitian_perturbed_intro}, then we see that $(A_1,\varphi_1,\psi_1)$
is a solution of the system \eqref{eq:SeibergWittenTaubesPert_Prelim}.  In particular, instead of 
\eqref{eq:SeibergWittenRefTaubesPert_curvature(1,1)_Prelim}, the triple $(A_1,\varphi_1,\psi_1)$ satisfies
Taubes' equation
\begin{equation}
  \label{eq:Taubes_omega-component_SW_curvature_equation_prelim}
  F_{A_1}^\omega
  =
  \frac{i}{8}|\varphi_1|_{L_1}^2\omega - \frac{i}{8}|\psi_1|_{\Lambda^{0,2}(L_1)}^2\omega - \frac{ir_0}{8}\omega,
\end{equation}
for the Seiberg--Witten equations in \cite[Section 1(d), Equations (1.18), (1.19), and (1.20), p. 851]{TauSWGromov} when the solution $(A,\Phi)$ to the non-Abelian monopole equations is a split pair (see also Donaldson \cite[Section 4, third displayed system of equations, p. 60]{DonSW} and Kotschick \cite[Equations (15), (16), and (17)]{KotschickSW}).

Consequently, a solution $(A_1,\Phi_1)$ to the perturbed Seiberg--Witten monopole equations
\eqref{eq:SeibergWittenTaubesPert_Prelim}
implied by a split solution $(A,\Phi)$ to the perturbed non-Abelian monopole equations
\eqref{eq:SO(3)_monopole_equations_almost_Hermitian_perturbed_intro}
obeys the seven pointwise estimates listed by Taubes in \cite[Section 1(f), Equations (1.24) and (1.26), p. 853]{TauSWGromov}. If $\{(A_1^n,\Phi_1^n)\}_{n\in\NN}$ is a sequence of solutions to the perturbed Seiberg--Witten monopole equations implied by a sequence of split solutions $\{(A,\Phi)\}_{n\in\NN}$ to the perturbed non-Abelian monopole equations with a sequence of parameters $\{r_n\}_{n\in\NN}$ with $\limsup_{n\in\NN} r_n = \infty$, then the sequence $\{(A_1^n,\Phi_1^n)\}_{n\in\NN}$ obeys convergence properties described by Taubes in \cite[Section 1(e), Theorem 1.3, p. 852]{TauSWGromov}; see also Taubes \cite[Section 5(a), pp. 881--882]{TauSWGromov}.
\qed\end{rmk}

We will use the following in Chapter \ref{chap:Feasibility} to prove that $\sM^{*,0}(E,g,J,\omega,r)$ is non-empty.

\begin{lem}[Non-zero Seiberg--Witten invariant implies non-empty $\sM^{*,0}(E,g,J,\omega,r)$]
\label{lem:NonZeroSWInvariant}
Let $(E,H)$ be a rank-two Hermitian vector bundle over a smooth, closed almost Hermitian four-manifold $(X,J,g,\omega)$.
Assume that $b^+(X)>1$.  If $\ft=(\rho_{\can}\otimes\id_E,W_{\can}\otimes E)$ and
$\fs=(\rho_{\can}\otimes\id_E{L_1},W_{\can}\otimes L_1$ are a spin${}^u$ and spin${}^c$ structure on $X$
and $\SW_X(\fs)\neq 0$, 
and if $w_2(\ft)$ is good in the sense of Definition \ref{defn:Good},
then
\begin{enumerate}
\item
\label{item::NonZeroSWInvariantImpliesMsNonEmpty}
For $r>0$, $M_\fs(r)$ is non-empty.
\item
\label{item::NonZeroSWInvariantImpliesIrreducibleNonEmpty}
If $\fs\in\Red_0(\ft)$ and the expected dimension of $\sM_\ft(r)$ given in \eqref{eq:Transv} is positive
then $\sM^{*,0}(E,g,J,\omega,r)$ is non-empty.
\end{enumerate}
\end{lem}

\begin{proof}
The lemma follows from
%TL11-28-2025: Updated to \cite[Proposition 6.8.3]{Feehan_Leness_introduction_virtual_morse_theory_so3_monopoles}
\cite[Proposition 6.8.3]{Feehan_Leness_introduction_virtual_morse_theory_so3_monopoles} and the observation that the regularized Taubes perturbation $\wp_\gamma$ can be considered as another parameter in the space of perturbations given in 
%TL12-4-2025 Updated
\cite[Equation (6.8.1)]{Feehan_Leness_introduction_virtual_morse_theory_so3_monopoles}.
\end{proof}

\section{Circle actions}
\label{sec:S1Actions}
The affine space $\sA(E,H,A_d)\times W^{1,p}(W^+\otimes E)$ carries a circle action \label{Circle_action_on_SO(3)_pairs} induced by scalar multiplication on the Hermitian vector bundle $E$:

\begin{defn}[Standard $S^1$ action on affine and quotient spaces of spin${}^u$ pairs]
\label{defn:UnitaryZActionOnAffine}
The \emph{standard $S^1$ action on the affine space of spin${}^u$ pairs} is
\begin{multline}
  \label{eq:S1ZAction}
  S^1\times \sA(E,H,A_d)\times W^{1,p}(W^+\otimes E)
  \ni (e^{i\theta},A,\Phi)
  \\
  \mapsto (A,e^{i\theta}\Phi) \in \sA(E,H,A_d)\times W^{1,p}(W^+\otimes E).
\end{multline}
Because the circle action \eqref{eq:S1ZAction} commutes with that of $W^{2,p}(\SU(E))$ in \eqref{eq:GaugeActionOnSpinuPairs}, the action \eqref{eq:S1ZAction} also defines a circle action on the quotient,
\begin{equation}
\label{eq:S1ZActionOnQuotientSpace}
S^1\times\sC_\ft \ni (e^{i\theta},[A,\Phi]) \mapsto [A,e^{i\theta}\Phi] \in \sC_\ft,
\end{equation}
which we call \emph{standard $S^1$ action on the quotient space $\sC_\ft$}.
\qed\end{defn}

Note that if $r$ is the rank of $\ft$, and $\varrho$ is an $r$-th root of unity, then $\varrho$ acts trivially
on $\sC_{\ft}$. We can see this by observing that for all $e^{i\theta}\in S^1$, the central gauge transformation
$e^{i\theta}\,\id_E\in W^{2,p}(U(E))$ acts trivially on $\sA(E,H,A_d)$ and so the action \eqref{eq:S1ZAction} is just the action \eqref{eq:GaugeActionOnSpinuPairs} extended to the elements of $W^{2,p}(\U(E))$ given by $e^{i\theta}\,\id_E$. Because $\varrho\,\id_E$ is a section of $\SU(E)$, it acts trivially on the $W^{2,p}(\SU(E))$-quotient $\sC_\ft$, although non-trivially on the affine space $\sA(E,H,A_d)\times W^{1,p}(W^+\otimes E)$.

We now identify the fixed points of the $S^1$ action in Definition \ref{defn:UnitaryZActionOnAffine} on $\sC_\ft$.
Since equation \eqref{eq:PerturbedSO3MonopoleEquations} is invariant under the circle action induced by scalar multiplication on $W^+\otimes E$, the subspaces \eqref{eq:PU2MonopoleSubspaces} of $\sC_{\ft}$ are also invariant under this action.  In \cite[Proposition 3.1, p. 86]{FL2a}, we characterized the fixed points of the action $S^1$ action in Definition \ref{defn:UnitaryZActionOnAffine} on the moduli space $\sM_\ft$ of non-Abelian monopoles.

\begin{prop}[Fixed points of the $S^1$ action in $\sC_\ft$]
\label{prop:FixedPointsOfS1ActionOnSpinuQuotientSpace}
%TL11-28-2025: Updated
(See \cite[Proposition 6.7.2]{Feehan_Leness_introduction_virtual_morse_theory_so3_monopoles}.)
Let $\ft=(\rho,W\otimes E)$ be a spin${}^u$ structure on a connected, oriented, smooth Riemannian four-manifold. If $[A,\Phi] \in \sC_\ft$ obeys
\begin{enumerate}
\item
\label{item:SO(3)MonopoleFixedPointsZeroSection}
$\Phi\equiv 0$, so $(A,0)$ is a zero-section pair as in Definition \ref{defn:ZeroSectionSO(3)Pair}, or
\item
\label{item:SO(3)MonopoleFixedPointsReducible}
$(A,\Phi)$ is
split as in Definition \ref{defn:Split_trivial_central-stabilizer_spinor_pair},
\end{enumerate}
then $[A,\Phi]$ is a fixed point of the $S^1$ action in Definition \ref{defn:UnitaryZActionOnAffine} on $\sC_\ft$. Conversely, if a fixed point of that $S^1$ action on $\sC_\ft$ is represented by smooth pair $(A,\Phi)$, then it obeys Condition \eqref{item:SO(3)MonopoleFixedPointsZeroSection} or \eqref{item:SO(3)MonopoleFixedPointsReducible}. Lastly, if $E$ has rank two  and the characteristic class $w_2(\ft)$ defined in \eqref{eq:SpinUCharacteristics} is good as in Definition \ref{defn:Good}, then the moduli space $\sM_\ft$ of non-Abelian monopoles \eqref{eq:Moduli_space_SO(3)_monopoles} contains no points represented by a pair that is both split and zero-section.
\end{prop}

The action \eqref{eq:S1ZAction} is defined by the homomorphism from $S^1$ into the group of unitary gauge transformations of $E$,
\begin{equation}
\label{eq:DefineUnitaryS1ZActionsAtReducible}
\rho_Z(e^{i\theta}):=e^{i\theta}\,\id_{L_1}\oplus e^{i\theta}\,\id_{L_2},
\quad\text{for all } e^{i\theta} \in S^1.
\end{equation}
If the bundle $E$ admits an orthogonal splitting as a direct sum of Hermitian line bundles, $E=L_1\oplus L_2$, then there is an additional $S^1$ action defined by complex scalar multiplication on $L_2$.  We define two more homomorphisms from $S^1$ into the group of unitary gauge transformations of $E$ by
\begin{subequations}
\label{eq:DefineUnitaryS1ActionsAtReducible}
\begin{align}
\label{eq:DefineUnitaryS1sActionsAtReducible}
\rho_{\SU}(e^{i\theta})&:=e^{i\theta}\,\id_{L_1}\oplus e^{-i\theta}\,\id_{L_2},
\\
\label{eq:DefineUnitaryS12ActionsAtReducible}
\rho_2(e^{i\theta})&:=\id_{L_1}\oplus e^{i\theta}\,\id_{L_2}, \quad\text{for all } e^{i\theta} \in S^1.
\end{align}
\end{subequations}

\begin{rmk}[Choice of line bundles in splitting]
\label{rmk:ChoiceOfLineBundlesInDefinitionOfSplitting}
We note that the definition of the homomorphisms \eqref{eq:DefineUnitaryS1ActionsAtReducible} depend not only on the existence of a splitting $E=L_1\oplus L_2$ but also on the choice of line bundles giving the splitting.  In the interest of legibility, we do not include the choice of line bundles $L_1$ and $L_2$ in the notation for the homomorphisms \eqref{eq:DefineUnitaryS1ActionsAtReducible}.
\qed\end{rmk}

Observe that the homomorphisms \eqref{eq:DefineUnitaryS1ZActionsAtReducible} and \eqref{eq:DefineUnitaryS1ActionsAtReducible} (when the latter exist) are related by
\begin{equation}
\label{eq:UnitaryS1ActionsRelation}
\rho_2(e^{2i\theta})=\rho_{\SU}(e^{-i\theta})\rho_Z(e^{i\theta}),
\quad\text{for all } e^{i\theta} \in S^1.
\end{equation}
The homomorphism \eqref{eq:DefineUnitaryS1ZActionsAtReducible} and, when $E$ admits a splitting $E=L_1\oplus L_2$, the homomorphism \eqref{eq:DefineUnitaryS12ActionsAtReducible}  define actions of $S^1$
on $\sA(E,H,A_d)\times W^{1,p}(W^+\otimes E)$ by
\begin{subequations}
\label{eq:S1ActionsOnSpinuPreQuotient}
\begin{align}
\label{eq:S1ZActionsOnSpinuPreQuotient}
\left( e^{i\theta},(A,\Phi)\right)
&\mapsto
( \rho_Z(e^{i\theta})^*A,\rho_Z(e^{-i\theta})\Phi),
\\
\label{eq:S12ActionsOnSpinuPreQuotient}
\left( e^{i\theta},(A,\Phi)\right)
&\mapsto
(\rho_2(e^{i\theta})^*A,\rho_2(e^{-i\theta})\Phi)
\end{align}
\end{subequations}
and on $\sC_\ft$ by
\begin{subequations}
\label{eq:S1ActionsOnSpinuQuotient}
\begin{align}
\label{eq:S1ZActionsOnSpinuQuotient}
\left( e^{i\theta},[A,\Phi]\right)
&\mapsto
[ \rho_Z(e^{i\theta})^*A,\rho_Z(e^{-i\theta})\Phi],
\\
\label{eq:S12ActionsOnSpinuQuotient}
\left( e^{i\theta},[A,\Phi]\right)
&\mapsto
[ \rho_2(e^{i\theta})^*A,\rho_2(e^{-i\theta})\Phi].
\end{align}
\end{subequations}
Because $\rho_{\SU}$ takes values in $C^\infty(\SU(E))$, an action defined by $\rho_{\SU}$ as in \eqref{eq:S1ActionsOnSpinuQuotient} would be trivial. The action \eqref{eq:S1ZActionsOnSpinuQuotient} is the action \eqref{eq:S1ZActionOnQuotientSpace}.  The equality
\eqref{eq:UnitaryS1ActionsRelation} implies that the actions \eqref{eq:S1ZActionsOnSpinuQuotient} and \eqref{eq:S12ActionsOnSpinuQuotient}  are equal up to multiplicity but each has its advantages.  The action \eqref{eq:S1ZActionsOnSpinuQuotient} is defined even when the bundle $E$ does not admit a splitting into line bundles
while the gauge transformations in the image of $\rho_2$ are in the stabilizer of pairs in the image of $\tilde\iota_{\fs,\ft}$ as we now describe,

\begin{lem}[Circle equivariance of the non-Abelian monopole map]
\label{lem:S1ActionsForReduciblePairGerm}
%TL11-28-2025:Reference updated
(See \cite[Lemma 6.7.4]{Feehan_Leness_introduction_virtual_morse_theory_so3_monopoles}.)
Let $\fs=(\rho,W\otimes L_1)$ and $\ft=(\rho,W\otimes E)$ be a spin${}^c$ and rank-two spin${}^u$ structure, respectively, over a closed, oriented, smooth, Riemannian four-manifold $(X,g)$. If $E$ admits a splitting $E=L_1\oplus L_2$ as an orthogonal direct sum of Hermitian line bundles, then the following hold
\begin{enumerate}
\item If $(A,\Phi)\in\sA(E,H,A_d)\times W^{1,p}(W^+\otimes E)$ is in the image of the map $\tilde\iota_{\fs,\ft}$ defined in \eqref{eq:DefnOfIota}, then $(A,\Phi)$ is a fixed point of the action \eqref{eq:S12ActionsOnSpinuPreQuotient}.
\item 
\label{item:S1_Equivariance_of_nonAbelian_monopole_map}
The non-Abelian monopole map \eqref{eq:PerturbedSO3MonopoleEquation_map},
\[
\sS:\sA(E,H,A_d)\times W^{1,p}(W^+\otimes E) \to
L^p(\La^+\otimes\su(E))\oplus L^p(W^-\otimes E),
\]
is $S^1$-equivariant with respect to the action \eqref{eq:S12ActionsOnSpinuPreQuotient} on the domain of $\sS$ and the action
\begin{multline}
\label{eq:S1L2ActionOnnonAbelianMonopoleMapCodomain}
S^1\times L^p(\La^+\otimes\su(E))\oplus L^p(W^-\otimes E)
\ni
\left(e^{i\theta},(\om,\Psi)\right)
\\
\mapsto
\left( \rho_2(e^{-i\theta})\om \rho_2(e^{i\theta}),\rho_2(e^{-i\theta})\Psi\right),
\end{multline}
on the codomain of $\sS$.
\item
\label{item:S1EquivOfPerturbednonAbelianMonopoleEquations}
Item \eqref{item:S1_Equivariance_of_nonAbelian_monopole_map} holds if the map $\sS$ is replaced by the map
defined by the non-Abelian monopole equations with a singular Taubes perturbation, \eqref{eq:SO(3)_monopole_equations_almost_Hermitian_perturbed_intro},
or by the non-Abelian monopole equations with a regularized Taubes perturbation,  \eqref{eq:SO(3)_monopole_equations_almost_Hermitian_perturbed_intro_regular}.
\end{enumerate}
\end{lem}

\chapter{Feasibility of the non-Abelian monopole method}
\label{chap:Feasibility}
In this section, we prove Theorem \ref{mainthm:ExistenceOfSpinuForFlow}, demonstrating the existence of a feasible spin${}^u$ structure as defined in Definition \ref{maindefn:Feasibility} on the smooth blow-up $\widetilde X$ of $X$. We work on $\widetilde X$ rather than directly on $X$ to employ the Morgan--Mrowka trick, \cite{MorganMrowkaPoly}, for ensuring there are no zero-section, reducible pairs in the moduli space $\sM_\ft$ and to get the control over the intersection form we need to prove Proposition \ref{prop:FeasibilityOnBlowUp}.

Recall that the smooth blow-up of $X$ is $\widetilde X:=X\#\overline{\CC\PP}^2$ of $X$ (see Gompf and Stipsicz \cite[Definition 2.2.7, p. 43]{GompfStipsicz}).
 Let $e^*\in H_2(\widetilde X;\ZZ)$ denote the fundamental class of the exceptional sphere and $e\in H^2(\widetilde X;\ZZ)$ its Poincar\'e dual  (see Gompf and Stipsicz \cite[Definition 2.2.7]{GompfStipsicz}).
 The blow-down map $\pi:\widetilde X\to X$ (see \cite[Definition 2.2.7, p. 43]{GompfStipsicz}) induces an injection
$\pi^*:H^2(X;\ZZ)\to H^2(\widetilde X;\ZZ)$ and an isomorphism,
\[
H^2(\widetilde X;\ZZ) \cong H^2(X;\ZZ) \oplus \ZZ e.
\]
The intersection form $Q_{\widetilde X}$ is isomorphic to $Q_X\oplus (-1)$ with respect to the preceding isomorphism.
The characteristic numbers of $X$ and $\widetilde X$ are related in \eqref{eq:BlowUpCharNumbers}.

We will need to consider the expansion of the set $B(\widetilde X)$ defined in 
%TL11-28-2025: Updated
\cite[Equation (7.1.24)]{Feehan_Leness_introduction_virtual_morse_theory_so3_monopoles} by
\begin{equation}
\label{eq:NonEmptySWSpace}
\widehat B(X):=
\{ c_1(\fs): \fs\in\Spinc(X)\ \text{and}\ M_\fs\neq\emptyset\}.
\end{equation}
By Morgan \cite[Theorem 5.2.4, p. 79]{MorganSWNotes}, $\widehat B(X)$ is finite.

%TL9-15-2024:Move to prelim
\begin{rmk}[Transformation of topological invariants]
\label{rmk:Transf_of_Top_Invariants}
In previous works, we have given expressions in terms of the topological invariants of the four-manifold $X$,
\[
\chi_h(X):=(1/4)(e(X)+\si(X)),\quad
c_1(X)^2=2e(X)+3\si(X),
\]
where $e(X)$ is the Euler characteristic of $X$ and $\si(X)$ is the signature.  Here, we will work with
the invariants $c_1(X)^2$ and $c_2(X)=e(X)$.  Elementary algebra yields
\begin{equation}
\label{eq:TopInvariantsTransformation}
c_2(X)=12\chi_h(X)-c_1(X)^2,
\quad
\chi_h(X)=(1/12)\left( c_2(X)+c_1(X)^2\right),
\end{equation}
and
\begin{equation}
\label{eq:Signature_and_Euler_To_Chern}
e(X)=c_2(X),
\quad
\si(X)=\frac{1}{3}\left( c_1(X)^2-2c_2(X)\right).
\end{equation}
In particular, we have the equality
\begin{equation}
\label{eq:BMY_Line_Transformation}
c_1(X)^2-9\chi_h(X)
=
\frac{1}{4}\left( c_1(X)^2-3c_2(X)\right).
\end{equation}
\qed\end{rmk}

%TL9-15-2024:Move to Section \ref{subsec:Virtual_Morse-Bott_signatures_SW_critical_points}
In the following, we establish criteria for a cohomology class to define a spin${}^u$ structure $\ft$ with $M_\fs$ embedding in $\sM_\ft$, and for which the formal Morse--Bott index \eqref{eq:FormalMorseIndexIntroThm},
\[
\la^-(\ft,\fs)=
-\frac{1}{6}\left(c_1(\widetilde X)^2+c_2(\widetilde X)\right)-\left(c_1(\fs)-c_1(\ft)  \right)\cdot c_1(\widetilde X) - \left( c_1(\fs)-c_1(\ft)\right)^2,
\]
is positive for all spin${}^c$ structures $\fs$ such that the Seiberg--Witten moduli space $M_\fs$ embeds in $\sM_\ft$  under the embedding \eqref{eq:DefnOfIotaOnQuotient}.

\begin{prop}
\label{prop:Feasibility}
Let $X$ be a smooth, closed, oriented Riemannian four-manifold with $b^+(X)>1$ and odd $b^+(X)-b_1(X)$ and let
$\widetilde X:=X\#\overline{\CC\PP}^2$ be the smooth blow-up of $X$ at a point. Assume that there is a  spin${}^c$ structure $\fs_0$ on $X$ with $c_1(\fs_0)^2=c_1(X)^2$ and $\SW_X(\fs_0)\neq 0$. Let $g$ be a Riemannian metric on $\widetilde X$ and let $\fs$ be the spin${}^c$ structure on $(\widetilde X,g)$ with $c_1(\fs)=-c_1(\fs_0)-e$. For $w\in H^2(\widetilde X;\ZZ)$, let $L_w\to \widetilde X$ be a complex line bundle with $c_1(L_w)=w$. If $\fs=(\rho,W)$ then define a spin${}^u$ structure on $\widetilde X$ by
\begin{equation}
\label{eq:SpinuStructureFromw}
\ft(\fs,w)= (\rho\otimes\id_{\ubarCC\oplus L_w}, W\otimes (\ubarCC\oplus L_w) ),
\end{equation}
where $\ubarCC=\CC\times \widetilde X$ is the product line bundle.
Then, the characteristic classes of $\ft=\ft(\fs,w)$ defined in \eqref{eq:SpinUCharacteristics} are
\begin{equation}
\label{eq:CharClasses_for_Spinu_From_Spinc_And_w}
c_1(\ft)=c_1(\fs)+w,
\quad
p_1(\ft)=w^2,
\quad
w_2(\ft)\equiv w\pmod 2,
\end{equation}
and $\sM_\ft$ contains the image of the moduli space of Seiberg--Witten monopoles $M_{\fs}$ under the embedding \eqref{eq:DefnOfIotaOnQuotient}.
In addition,
\begin{enumerate}
\item
\label{item:ExpectedDimOfZeroSections}
If, for $b\in\ZZ$,
\begin{equation}
\label{eq:w2Requirement}
w^2=-c_2(X)+b
%=c_1(\widetilde X)^2-12\chi_h(\tilde  X)+b,
\end{equation}
and $M^w_\ka(\widetilde X,g)$ is the moduli space of anti-self-dual connections identified with the zero-section pairs in $\sM_\ft$, then the expected dimension of $M^w_\ka(\widetilde X,g)$  is
\[
\exp\dim M^w_\ka(\widetilde X, g)
=
\frac{1}{2}\left(3c_2(X) - c_1(X)^2\right)-2b
%2\left(-c_1(\tilde  X)^2+9\chi_h(\tilde  X)\right)-2b,
\]
and so if $b\ge 0$, then
\begin{equation}
\label{eq:ASDExpectedDimensionBound}
\exp\dim M^w_\ka(\widetilde X, g)
\le
\frac{1}{2}\left(3c_2(X) - c_1(X)^2\right).
%2\left(-c_1(\tilde  X)^2+9\chi_h(\tilde  X)+1\right).
\end{equation}
\item
\label{item:MMCondition}
If $e^*\in H_2(\widetilde X;\ZZ)$ is the fundamental class of the exceptional sphere and
\begin{equation}
\label{eq:MMCondition}
\langle w,e^*\rangle\equiv 1 \pmod 2,
\end{equation}
then $\sM_\ft$ does not contain any zero-section, split pairs.
\item
\label{item:ExpectedDimOfSO(3)Monopoles}
If
\begin{equation}
\label{eq:LowerBoundOnw2ForPosDimSO(3)Monopoles}
w\cdot c_1(\fs)>  w^2+\frac{1}{6}\left( c_1(\widetilde X)^2+c_2(\widetilde X)\right),
\end{equation}
then $\exp\dim\sM_\ft>0$.

\item
\label{item:NonEmptyIrreducibleSO3Monopoles}
If $\SW_{\widetilde X}(\fs)\neq 0$ and if \eqref{eq:LowerBoundOnw2ForPosDimSO(3)Monopoles} holds, then
 $\sM_\ft^{*,0}$ is non-empty.

\item
\label{item:PositivevMBForOneSpinc}
If \eqref{eq:LowerBoundOnw2ForPosDimSO(3)Monopoles} holds and if the formal Morse--Bott virtual index $\lambda^-(\ft,\fs)$ defined in \eqref{eq:FormalMorseIndexIntroThm} is computed using an almost complex structure $J$ on $T\widetilde X$ with $c_1(T\tilde  X,J)=c_1(\fs)$, then $\lambda^-(\ft,\fs)>0$.

\item
\label{item:PositivevMBForLevelZeroSpinc}
If, in addition,
\begin{equation}
\label{eq:LowerBoundOnw2ForPosDimSO(3)MonopolesToGetMB}
w\cdot c_1(\fs)>
w^2+\frac{1}{6}\left( c_1(\widetilde X)^2+c_2(\widetilde X)\right)+\sup_{\fs'\in \widehat B(\widetilde  X)} \left( c_1(\fs')-c_1(\fs)\right)\cdot c_1(\fs)
%w^2+2\chi_h(\tilde  X)+\sup_{\fs'\in \widehat B(\tilde  X)} \left( c_1(\fs')-c_1(\fs_0)\right)\cdot c_1(\fs_0)
\end{equation}
where $\widehat B(\widetilde X)$ is defined in \eqref{eq:NonEmptySWSpace}
%PF3-11-2024 Obviously we like this conclusion, but because the above condition appears so strong, it would be worth reminding readers in a remark how we proved that it can be achieved in our  monograph?
then $\lambda^-(\ft,\fs')>0$ for all spin${}^c$ structures $\fs'$ with $M_{\fs'}\subset\sM_\ft$.
\end{enumerate}
\end{prop}

\begin{proof}
The expressions \eqref{eq:CharClasses_for_Spinu_From_Spinc_And_w} follow immediately from the definition of the spin${}^u$ structure $\ft(\fs,w)$ in \eqref{eq:SpinuStructureFromw}
and the definition of the characteristic classes in \eqref{eq:SpinUCharacteristics}.
In particular,
\[
p_1(\ft)=p_1(\su(\ubarCC\oplus L_w)=w^2,
\]
where the second equality follows from Donaldson and Kronheimer \cite[Equation (2.1.39), p. 42]{DK}.

The existence of the embedding $M_{\fs}\to \sM_\ft$ follows from
Lemma \ref{lem:SetOfReducibles} and
 the observation that $c_1(\fs)-c_1(\ft)=w$ so that $c_1(\fs)\in \Red_0(\ft)$ by the definition \eqref{eq:DefineReduciblesEmbedded} and the equality $p_1(\ft)=w^2$ in \eqref{eq:CharClasses_for_Spinu_From_Spinc_And_w}.

Item \eqref{item:ExpectedDimOfZeroSections} follows from $w^2=p_1(\ft)$, the assumption \eqref{eq:w2Requirement},  the  formula
for the expected dimension of the moduli space of anti-self-dual connections in \eqref{eq:Expected_dimension_moduli_space_ASD_connections},
and the relation \eqref{eq:BlowUpCharNumbers} between the characteristic numbers of $X$ and those of $\widetilde X$,
\[
\expdim M_\ka^w(\widetilde X, g)
=
-2p_1(\ft)-\frac{1}{2}\left( c_1(X)^2+c_2( X)\right)
\]
%$-2p_1(\ft)-6\chi_h(Y)$
(see e.g. \cite[Equation (4.2.22)]{DK} or 
%TL11-28-2025: Updated
\cite[Equation (6.4.6)]{Feehan_Leness_introduction_virtual_morse_theory_so3_monopoles}).

Item \eqref{item:MMCondition} follows immediately from the Morgan--Mrowka criterion,
\cite[p. 226]{MorganMrowkaPoly}, as stated in \cite[Corollary 3.3, p. 88]{FL2a}.

To prove Item \eqref{item:ExpectedDimOfSO(3)Monopoles}, we  compute
\begin{align*}
\exp\dim\sM_\ft&= d_a(\ft)+2n_a(\ft)
\quad\text{(see
  % TL11-28-2028: Updated
                 %PF12-2-2025 I don't see a change?
                 %TL12-4-2025: Re-updated...
 \cite[Equation (6.4.6)]{Feehan_Leness_introduction_virtual_morse_theory_so3_monopoles})}
\\
&=
-2p_1(\ft)-\frac{3}{2}\left( e(\widetilde X)+\si(\widetilde X)\right) +\frac{1}{2}p_1(\ft) +\frac{1}{2}c_1(\ft)^2 -\frac{1}{2}\si(\widetilde X)
%-2p_1(\ft)-6\chi_h(\widetilde X) +\frac{1}{2}p_1(\ft) +\frac{1}{2}c_1(\ft)^2 -\frac{1}{2}\left( c_1(\tilde %X)^2-8\chi_h(\widetilde X)\right)
\\
&=
-\frac{3}{2}p_1(\ft) +\frac{1}{2}\left( c_1(\fs)+w\right)^2-\frac{3}{2}e(\widetilde X) -2\si(\widetilde X)
\qquad\text{(because $c_1(\ft)=c_1(\fs)+w$ )}
%-\frac{3}{2}p_1(\ft) -2\chi_h(\widetilde X)+c_1(\fs)\cdot w+\frac{1}{2}w^2
%&\qquad
%\quad\text{(because $c_1(\ft)=c_1(\fs)+w$ and $c_1(\fs)^2=c_1(\widetilde X)^2$)}
\\
&=
-w^2+c_1(\fs)\cdot w +\frac{1}{2}c_1(\widetilde X)^2-\frac{3}{2}c_2(\widetilde X)-\frac{2}{3}\left( c_1(\widetilde X)-2c_2(\widetilde X)\right)
\\
&\qquad\text{(because $w^2=p_1(\ft)$, $c_1(\fs)^2=c_1(\widetilde X)^2$, and by \eqref{eq:Signature_and_Euler_To_Chern})}
\\
&=
-w^2+c_1(\fs)\cdot w -\frac{1}{6}\left( c_1(\widetilde X)^2+c_2(\widetilde X)\right).
\end{align*}
Thus, if \eqref{eq:LowerBoundOnw2ForPosDimSO(3)Monopoles} holds then the expected dimension of $\sM_\ft$ is positive,
proving Item  \eqref{item:ExpectedDimOfSO(3)Monopoles}.

The assumption that $\SW_X(\fs_0)\neq 0$ implies that $\bar\fs_0$ is the spin${}^c$ structure on $X$ with $c_1(\bar\fs_0)=-c_1(\fs_0)$, then $\SW_X(\bar\fs_0)\neq 0$ by \cite[Corollary 6.8.4, p. 103]{MorganSWNotes}.  The blow-up formula in Theorem \ref{thm:SWBlowUp} implies that $\SW_{\widetilde X}(\fs)\neq 0$. Because $\sM_\ft$ contains the image of the moduli space of Seiberg--Witten monopoles $M_{\fs}$ under the embedding \eqref{eq:DefnOfIotaOnQuotient}, the Item \eqref{item:MMCondition} implies that there are no points in $M_{\fs}$ represented by zero-section pairs. Then the non-vanishing result $\SW_{\widetilde X}(\fs)\neq 0$, the positivity of the expected dimension of $\sM_\ft$ established in Item \eqref{eq:LowerBoundOnw2ForPosDimSO(3)Monopoles},
 and Lemma \ref{lem:NonZeroSWInvariant}
%TL5-26-2025: Put the result in for the Taubes perturbed moduli space 
 %\cite[Proposition 5.8.3]{Feehan_Leness_introduction_virtual_morse_theory_so3_monopoles} 
imply that $\sM_\ft^{*,0}$ is non-empty.  This proves Item \eqref{item:NonEmptyIrreducibleSO3Monopoles}.
%PF1-29-2024 Why not cite the exact result that actually proves the non-emptiness of $\sM_\ft^{*,0}$? The lemma below is using \widetilde X, the blow up of X, so it's additionally confusing (and also leading to the mix up about b \geq 1 vs. b \geq 0 and Morgan-Mrowka criterion)
%PF2-5-_2024 Is below needed or not?
%as in \cite[Lemma 6.1.1]{Feehan_Leness_introduction_virtual_morse_theory_so3_monopoles}.
%TL4-4-2024: Exact result now cited so the above not need.

We can assume that there is an almost complex structure $J$ on $\widetilde X$ with $c_1(\widetilde X):=c_1(T\widetilde X,J)=c_1(\fs)$ by Gompf and Stipsicz \cite[Theorem 1.4.15, p. 29]{GompfStipsicz}. Then by
\eqref{eq:FormalMorseIndexIntroThm} the formal Morse--Bott index is
%PF3-11-2024 I'm sure it's correct, but please add intermediate steps/algebra.
%TL4-4-2024: a note added to explain
\begin{align*}
\lambda^-(\ft,\fs)&=-\frac{1}{6}\left(c_1(\widetilde X)^2+c_2(\widetilde X) \right) -\left( c_1(\fs)-c_1(\ft)\right)\cdot c_1(\widetilde X) - (c_1(\fs)-c_1(\ft))^2
\\
&=
-\frac{1}{6}\left(c_1(\widetilde X)^2+c_2(\widetilde X)\right) +w\cdot c_1(\fs) -w^2
\quad\text{(using $w=c_1(\ft)-c_1(\fs)$ from \eqref{eq:SpinuStructureFromw})}
\end{align*}
which we record as
\begin{equation}
\label{eq:vMBForFixedSpinc}
\lambda^-(\ft,\fs)
=-\frac{1}{6}\left(c_1(\widetilde X)^2+c_2(\widetilde X) \right) +w\cdot c_1(\fs) -w^2.
\end{equation}
Thus if \eqref{eq:LowerBoundOnw2ForPosDimSO(3)Monopoles} holds, then
\eqref{eq:vMBForFixedSpinc} implies that
the formal Morse--Bott index $\lambda^-(\ft,\fs)$ is positive.  This proves Item \eqref{item:PositivevMBForOneSpinc}.

We prove Item \eqref{item:PositivevMBForLevelZeroSpinc} as follows.  If $\fs'$ is a spin${}^c$ structure on $\widetilde X$ with $c_1(\fs')\in \Red_0(\ft)$,
 then by the definition of $\Red_0(\ft)$ in \eqref{eq:DefineReduciblesEmbedded}
and the expression for $p_1(\ft)$ in \eqref{eq:SpinuStructureFromw},
$c_1(\fs')$ satisfies
\[
\left( c_1(\fs')-c_1(\ft)\right)^2=w^2.
\]
The preceding equality and the definition of $\la^-(\ft,\fs)$ in \eqref{eq:FormalMorseIndexIntroThm}
imply that for all $\fs'\in \Red_0(\ft)$,
\begin{equation}
\label{eq:vMBDifference}
\lambda^-(\ft,\fs')-\lambda^-(\ft,\fs)
=
-\left( c_1(\fs')-c_1(\fs)\right)\cdot c_1(\widetilde X).
\end{equation}
Combining \eqref{eq:vMBDifference}
 and the assumption that $c_1(\widetilde X)=c_1(\fs)$ implies that for all $\fs'\in \Red_0(\ft)$
\begin{align*}
\lambda^-(\ft,\fs')&=
\lambda^-(\ft,\fs)+\left( \lambda^-(\ft,\fs')-\lambda^-(\ft,\fs)\right)
\\
&=
%PF1-29-2024 I don't understand the conclusion below. The RHS second term below would be -\left( c_1(\fs')-c_1(\fs)\right)\cdot w from \eqref{eq:vMBForFixedSpinc}. Why is this the same as what you have below?
%TL 2-2-2024: \eqref{eq:vMBForFixedSpinc} only holds for the spinc structure \fs with c_1(\fs)=c_1(\ft)-2.  The following comes directly from the definition of $\la^-(\ft,\fs)$
%PF2-5-2024 Please add the intermediate step so it's clearer
%TL4-4-2024: I think intermediate step is now there.
\lambda^-(\ft,\fs)
-\left( c_1(\fs')-c_1(\fs)\right)\cdot c_1(\fs).
\end{align*}
%PF1-29-2024 I think the remainder of the calculations leading to Item \eqref{item:PositivevMBForLevelZeroSpinc} and its statement should be updated to reflect above change.
%TL 4-4-2024:I think this is OK now.
Equation \eqref{eq:vMBForFixedSpinc} then implies that
\begin{align*}
\lambda^-(\ft,\fs')&=
-\frac{1}{6}\left(c_1(\widetilde X)^2+c_2(\widetilde X) \right) +w\cdot c_1(\fs) -w^2
-\left( c_1(\fs')-c_1(\fs)\right)\cdot c_1(\fs).
\end{align*}
Hence, if
\[
w\cdot c_1(\fs)
>
w^2+\frac{1}{6}\left(c_1(\widetilde X)^2+c_2(\widetilde X) \right)+\left( c_1(\fs')-c_1(\fs)\right)\cdot c_1(\fs),
\]
then $\lambda^-(\ft,\fs')>0$.  The assumption \eqref{eq:LowerBoundOnw2ForPosDimSO(3)MonopolesToGetMB} then implies that $\lambda^-(\ft,\fs')>0$ for all $\fs'$ with $c_1(\fs')\in\Red_0(\ft)$, completing the proof of Item
\eqref{item:PositivevMBForLevelZeroSpinc}.
\end{proof}

To use  Proposition \ref{prop:FeasibilityOnBlowUp} to prove Theorem \ref{mainthm:ExistenceOfSpinuForFlow},
we need to demonstrate the existence of $w\in H^2(\widetilde X;\ZZ)$ satisfying
\eqref{eq:w2Requirement},
\eqref{eq:MMCondition},
\eqref{eq:LowerBoundOnw2ForPosDimSO(3)Monopoles},
 and \eqref{eq:LowerBoundOnw2ForPosDimSO(3)MonopolesToGetMB}.
Doing so is an algebraic exercise on the symmetric bilinear form $Q_{\widetilde X}$ and so
we first
recall some vocabulary describing symmetric bilinear forms.
A symmetric bilinear form $Q$ on a $\ZZ$ module $H$ is \emph{even} if $Q(x,x)\equiv 0\pmod 2$ for all $x\in H$ and \emph{odd} otherwise (see Gompf and Stipsicz \cite[p. 10]{GompfStipsicz}).  An element $\kappa\in H$ is \emph{characteristic} if $Q(\kappa,x)\equiv Q(x,x)\pmod 2$ for all $x\in H$ (see \cite[Definition 1.2.8(c)]{GompfStipsicz}).

\begin{prop}
\label{prop:FeasibilityOnBlowUp}
Let $(X,g)$ be a closed, smooth, oriented Riemannian four-manifold with indefinite intersection form and $b^+(X)>1$.
Let $\fs_0$ be a spin${}^c$ structure on $X$ with $c_1(\fs_0)^2=c_1(X)^2$ and
$\SW_X(\fs_0)\neq 0$.
For $\widetilde X=X\#\overline{\CC\PP}^2$, let $\fs$ be a spin${}^c$ structure on $\widetilde X$ with $c_1(\fs)=-c_1(\fs_0)-e$ where
$e\in H^2(\widetilde X;\ZZ)$ is Poincar\'e dual to the
exceptional sphere.
Assume that
\begin{enumerate}
\item
\label{item:FeasibilityOnBlowUp_c2>0}
$c_2(X)>0$,
\item
\label{item:FeasibilityOnBlowUp_Odd}
If $X$ is odd then $b^+(X)\ge 2$ or $b^-(X)\ge 2$.
\item
\label{item:FeasibilityOnBlowUp_Even}
If $X$ is even then the integer $b$ in
\eqref{eq:w2Requirement} is odd.
\end{enumerate}
Then there is $w\in H^2(\widetilde X;\ZZ)$ such that $w$ and $c_1(\fs)$ satisfy
the conditions
\eqref{eq:w2Requirement},
\eqref{eq:MMCondition},
\eqref{eq:LowerBoundOnw2ForPosDimSO(3)Monopoles},
 and \eqref{eq:LowerBoundOnw2ForPosDimSO(3)MonopolesToGetMB} of Proposition \ref{prop:Feasibility}.
\end{prop}

%PF10-3-2024 Check
\begin{rmk}[Redundant case in Proposition \ref{prop:FeasibilityOnBlowUp}]
The assumption that $b^+(X)>1$ in Proposition \ref{prop:FeasibilityOnBlowUp} means that
the case $b^+(X)=1$ and $b^-(X)\ge 2$ in Item \eqref{item:FeasibilityOnBlowUp_Odd} of Proposition \ref{prop:FeasibilityOnBlowUp} is redundant.  We include it here for future use when we will consider the case $b^+(X)=1$.
\qed\end{rmk}

\begin{proof}
We give the proof in three cases:
\begin{inparaenum}[\itshape i\upshape)]
\item $X$ is odd with $b^-(X)\ge 2$,
\item $X$ is odd with $b^+(X)\ge 2$, and
\item $X$ is even and $b$ is odd.
\end{inparaenum}
We note that if $w\in H^2(\widetilde X;\ZZ)$ satisfies \eqref{eq:w2Requirement} and
\eqref{eq:MMCondition} then $-w$ does as well, allowing us to replace a $w$ which satisfies \eqref{eq:w2Requirement} and
\eqref{eq:MMCondition} with $-w$ to prove that
\eqref{eq:LowerBoundOnw2ForPosDimSO(3)Monopoles} and
\eqref{eq:LowerBoundOnw2ForPosDimSO(3)MonopolesToGetMB} hold without having to verify \eqref{eq:w2Requirement} and
\eqref{eq:MMCondition} again.

\begin{case}[$X$ is odd, $b^+(X)\ge 1$, and $b^-(X)\ge 2$]
By Milnor and Husemoller \cite[Theorem II.4.3, p. 22]{MilnorHusemoller} or
Gompf and Stipsicz, $(H^2(X;\ZZ), Q_{X})$ has a summand isometric to
\begin{equation}
\label{eq:OddSummandbPlus1}
\ZZ f\bigoplus \ZZ e_1\bigoplus \ZZ e_2\
\end{equation}
where
\begin{equation}
\label{eq:Pairings_in_Xodd_Case1}
Q_{\widetilde X}(f,f)=1,\,
Q_{\widetilde X}(e_1,e_1)=-1,\,
Q_{\widetilde X}(e_2,e_2)=-1,\,
\text{and all other pairings of $f,e_1,e_2$ vanish}.
\end{equation}
%PF3-12-2024  Isn't  c_1(\fs) = c_1(\fs_0)+e? How is \kappa_0 defined?
We write $c_1(\fs_0)=-\kappa_0+\ka^\perp$, where
$\kappa_0$ is in the subspace \eqref{eq:OddSummandbPlus1} of $H^2( X;\ZZ)$ and
$\kappa^\perp$ is $Q_{X}$-orthogonal to that subspace.
We can then write
\[
%\kappa_0+e=(1+2a_1)f+(1+2b_1)e_1+(1+2b_2)e_2+e,
\kappa_0=(1+2a_1)f+(1+2b_1)e_1+(1+2b_2)e_2,
\]
where $a_1,b_1,b_2\in\ZZ$
and the coefficients of $f$, $e_1$, and $e_2$ in the preceding
equality
are odd because $c_1(\fs_0)$ is characteristic (see Gompf and Stipsicz \cite[Proposition 2.4.16, p. 56]{GompfStipsicz}).
Define $N\in\ZZ$ by
\begin{equation}
\label{eq:DefineNInCaseOddbMinus2}
2N-\eps=-c_2( X)^2+ b,
\end{equation}
where $\eps=0$ or $\eps=1$.
For $u\in\ZZ$, define integers
\begin{subequations}
\label{eq:DefineCoefficientsOfwDependingOnu}
\begin{align}
\label{eq:DefineCoefficientsOfwDependingOnu_z1}
z_1&=N+2u^2+2\eps u,
\\
\label{eq:DefineCoefficientsOfwDependingOnu_z2}
z_2&=\eps+2u,
\\
\label{eq:DefineCoefficientsOfwDependingOnu_v1}
v_1&=z_1+1.
\end{align}
\end{subequations}
For the integers $v_1,z_1,z_2$ given by \eqref{eq:DefineCoefficientsOfwDependingOnu}, define
$w_u\in H^2(\widetilde X;\ZZ)$ by
\begin{equation}
\label{eq:wInCaseOddbMinus2}
w_u=v_1f+z_1e_1+z_2e_2-e.
\end{equation}
Definition
\eqref{eq:wInCaseOddbMinus2}
and the equality $\langle e,e^*\rangle=-1$
%PF3-12-2024 And the pairings of f, e_i, etc. -- label them for citing
%TL4-4-2024: e is orthogonal to all cohomology classes on X so their pairings don't affect the preceding.
imply that $w_u$ satisfies the Morgan--Mrowka condition \eqref{eq:MMCondition}.
Then, by \eqref{eq:Pairings_in_Xodd_Case1}
\begin{align*}
Q_{\widetilde X}(w_u,w_u)
&=
v_1^2-z_1^2-z_2^2-1
\quad\text{(by \eqref{eq:Pairings_in_Xodd_Case1})}
\\
&=
(v_1-z_1)(v_1+z_1) -z_2^2-1
\\
&=
(2z_1+1)-z_2^2-1
\quad\text{(by \eqref{eq:DefineCoefficientsOfwDependingOnu_v1})}
\\
&=
2N+4u^2+4\eps u+1-\eps^2-4\eps u -4u^2-1 %PF3-12-2024 by ....
\quad\text{(by \eqref{eq:DefineCoefficientsOfwDependingOnu_z1} and \eqref{eq:DefineCoefficientsOfwDependingOnu_z2})}
\\
&=
2N-\eps,
\quad\text{(because $\eps=0,1$ implies $\eps^2=\eps$)}
\end{align*}
and so  \eqref{eq:DefineNInCaseOddbMinus2} implies that for any $u\in\ZZ$, $w_u$ satisfies
\eqref{eq:w2Requirement}.

Finally, we show that \eqref{eq:LowerBoundOnw2ForPosDimSO(3)Monopoles} and \eqref{eq:LowerBoundOnw2ForPosDimSO(3)MonopolesToGetMB} hold by proving that
$Q_{\widetilde X}(c_1(\fs),w_u)$ can be made arbitrarily large by choosing $u\in\ZZ$ appropriately.
Because
\[
c_1(\fs)=-c_1(\fs_0)-e=\ka_0-\ka^\perp-e,
\]
we can  compute
\begin{align*}
Q_{\widetilde X}(c_1(\fs),w_u)&=
Q_{\widetilde X}(\ka_0-e,w_u)
\quad\text{(because $\ka^\perp$ is $Q_{\widetilde X}$-orthogonal to the subspace \eqref{eq:OddSummandbPlus1})}
\\
&=v_1(1+2a_1)-z_1(1+2b_1)-z_2(1+2b_2)-1 %PF3-12-2024 by ....
\quad\text{(by \eqref{eq:Pairings_in_Xodd_Case1})}
\\
&=
(N+2u^2+2\eps u+1)(1+2a_1)-(N+2u^2+2\eps u)(1+2b_1) %PF3-12-2024 by ....
\quad\text{(by \eqref{eq:DefineCoefficientsOfwDependingOnu})}
\\
&\qquad
-(\eps+2u)(1+2b_2)-1
\\
&=
4u^2\left( a_1-b_1\right)
+
u\left(4\eps(a_1-b_1)-2(1+2b_2) \right)
\\
&\qquad
+2N(a_1-b_1)+1+2a_1-\eps(1+2b_2)-1
\\
&=
A_2 u^2 + A_1 u+A_0,
\end{align*}
where the coefficients of this polynomial in $u$ are
\[
A_2=4(a_1-b_1),\quad
A_1=4\eps(a_1-b_1)-2(1+2b_2),\quad
A_0=2N(a_1-b_1)+2a_1-\eps(1+2b_2).
\]
If $a_1>b_1$, then the leading coefficient $A_2$ is positive and so for $u$ sufficiently large,
$Q_{\widetilde X}(c_1(\fs),w_u)$ will satisfy \eqref{eq:LowerBoundOnw2ForPosDimSO(3)MonopolesToGetMB}.
If $a_1<b_1$, then after replacing $w_u$ with $-w_u$ (for which \eqref{eq:w2Requirement} and
\eqref{eq:MMCondition} will still hold),
$Q_{\widetilde X}(c_1(\fs),w_u)$ will again satisfy \eqref{eq:LowerBoundOnw2ForPosDimSO(3)MonopolesToGetMB} when $u$ is sufficiently large.
If $a_1=b_1$, then $A_2=0$ and $A_1=-2(1+2b_2)$ which is non-zero.  Then for $u$ sufficiently negative,
$Q_{\widetilde X}(c_1(\fs),w_u)$  will satisfy \eqref{eq:LowerBoundOnw2ForPosDimSO(3)MonopolesToGetMB}.

This completes the proof for the first case.
\end{case}

\begin{case}[$X$ is odd, $b^+(X)\ge 2$, and $b^-(X)\ge 1$]
By the same references as in the first case, $(H^2( X;\ZZ), Q_{X})$ has a summand isometric to
%PF3-12-2024 Why the duplication if summand is same?
%TL4-4-2024: Same references, but applied to different intersection forms
\begin{equation}
\label{eq:OddSummandbMinus1}
\ZZ f_1\oplus \ZZ f_2\oplus \ZZ e_1,
\end{equation}
where
\begin{equation}
\label{eq:Pairings_in_OddSummandbMinus1}
Q_{\widetilde X}(f_1,f_1)=1,\,
Q_{\widetilde X}(f_2,f_2)=1,\,
Q_{\widetilde X}(e_1,e_1)=-1,\,
\text{and all other pairings of $f_1,f_2,e_1$ vanish}.
\end{equation}
We write $c_1(\fs_0)=\kappa_0+\ka^\perp$, where
$\ka_0$ is in the subspace \eqref{eq:OddSummandbMinus1} and
$\kappa^\perp$ is $Q_{X}$-orthogonal to that subspace.
We then have
%PF3-12-2024 e on both sides?
%TL4-4-2024: Eliminated by working in X instead of \widetilde X
\[
\kappa_0=(1+2a_1)f_1+(1+2a_2)f_2+(1+2b_1)e_1,
\]
where $a_1,a_2,b_1\in\ZZ$.  The coefficients of $f_1$, $f_2$, and $e_1$ in the preceding are odd because of the requirement that $c_1(\fs)$ is characteristic (see Gompf and Stipsicz \cite[Proposition 2.4.16, p. 56]{GompfStipsicz}).
%PF3-12-2024 All the above seems like needless duplication
%PF3-12-2024  Contrast below with previous N definition
%TL4-4-2024: I put in a note below
Define $N\in\ZZ$ by
\begin{equation}
\label{eq:DefineNInCaseOddbMinus1}
2N+\eps-1
=-c_2( X)+b
%=c_1(\widetilde X)^2-12\chi_h(\widetilde X) + b,
\end{equation}
where $\eps=0$ or $\eps=1$.
(The need for the different definitions of $N$ in \eqref{eq:DefineNInCaseOddbMinus1} and \eqref{eq:DefineNInCaseOddbMinus2}
arises from the different intersection forms \eqref{eq:OddSummandbPlus1} and \eqref{eq:OddSummandbMinus1}.)
Define $\bar\eps=1-\eps$ and, for $u\in\ZZ$, define $v_1,v_2,z_1\in\ZZ$ by
\begin{subequations}
\label{eq:DefineCoefficientsOfwDependingOnuPlusVersion}
\begin{align}
\label{eq:DefineCoefficientsOfwDependingOnuPlusVersion_v1}
v_1&=\bar\eps+2u,
\\
\label{eq:DefineCoefficientsOfwDependingOnuPlusVersion_v2}
v_2&=\eps+N-2\bar\eps u -2u^2,
\\
\label{eq:DefineCoefficientsOfwDependingOnuPlusVersion_z1}
z_1&=v_2-1.
\end{align}
\end{subequations}
We then define
\begin{equation}
\label{eq:wInCaseOddbMinus1}
w_u=v_1 f_1+v_2f_2+z_1e_1-e.
\end{equation}
Definition
\eqref{eq:wInCaseOddbMinus1} implies that for any $u\in\ZZ$, $w_u$ satisfies the Morgan--Mrowka condition, \eqref{eq:MMCondition}. In addition,
\begin{align*}
Q_{\widetilde X}(w_u,w_u)
&=
v_1^2+v_2^2-z_1^2-1
\quad\text{(by \eqref{eq:Pairings_in_OddSummandbMinus1})}
\\
&=
v_1^2+(v_2-z_1)(v_2+z_1) -1
\\
&=
\bar\eps^2+4\bar\eps u + 4u^2+(1)(2v_2-1)-1
\quad\text{(by \eqref{eq:DefineCoefficientsOfwDependingOnuPlusVersion_v1} and \eqref{eq:DefineCoefficientsOfwDependingOnuPlusVersion_z1})}
\\
&=
1-\eps+4\bar\eps u + 4u^2
+ 2\eps +2N-4\bar\eps u - 4u^2 -1-1
\\
&\qquad
\quad\text{(by \eqref{eq:DefineCoefficientsOfwDependingOnuPlusVersion_v2} and
because $\bar\eps=0,1$ implies that $\bar\eps^2=\bar\eps=1-\eps$)}
\\
&=\eps -1+2N.
\end{align*}
Then \eqref{eq:DefineNInCaseOddbMinus1} implies that for any $u\in\ZZ$, $w_u$ satisfies
\eqref{eq:w2Requirement}.

Finally, we again show that \eqref{eq:LowerBoundOnw2ForPosDimSO(3)Monopoles} and \eqref{eq:LowerBoundOnw2ForPosDimSO(3)MonopolesToGetMB} hold by proving that
$Q_{\widetilde X}(c_1(\fs),w_u)$ can be made arbitrarily large by choosing $u\in\ZZ$ appropriately.
\begin{align*}
Q_{\widetilde X}(c_1(\fs),w_u)&=
Q_{\widetilde X}(\ka_0-e,w_u)
\quad\text{(because $\ka^\perp$ is $Q_{\widetilde X}$-orthogonal to the subspace \eqref{eq:OddSummandbPlus1})}
\\
&=
v_1(1+2a_1)+v_2(1+2a_2)-z_1(1+2b_1)-1
\quad\text{(by \eqref{eq:Pairings_in_OddSummandbMinus1})}
\\
&=
(\bar\eps +2u)(1+2a_1)
+
(\eps+N-2\bar\eps u -2u^2)(1+2a_2)
\\
&\qquad
-
(\eps+N-2\bar\eps u -2u^2-1)(1+2b_1)
-
1
\\
\quad\text{(by \eqref{eq:DefineCoefficientsOfwDependingOnuPlusVersion})}
&=
4(b_1-a_2)u^2
+
2\left(1+2a_1+ 2\bar\eps(b_1-a_2)\right)u
\\
&\qquad
+\bar\eps(1+2a_1)+(\eps+N)(1+2a_2)-
(\eps+N)(1+2b_1)-1
\end{align*}
If $b_1\neq a_2$,
then, by replacing $w_u$ with $-w_u$ if necessary, we have
$Q_{\widetilde X}(c_1(\fs),w_u)\gg 0$ for $u$ sufficiently large, so $w_u$ will satisfy
\eqref{eq:LowerBoundOnw2ForPosDimSO(3)Monopoles} and \eqref{eq:LowerBoundOnw2ForPosDimSO(3)MonopolesToGetMB}.
If $b_1=a_2$, then the coefficient of $u^2$ vanishes while that of $u$ is $2(1+2a_1)$, which is necessarily non-zero because $a_1\in\ZZ$.  Hence, again by replacing $w_u$ with $-w_u$ if necessary, we have
$Q_{\widetilde X}(c_1(\fs),w_u)\gg 0$ for $u$ sufficiently large, so $w_u$ will satisfy
\eqref{eq:LowerBoundOnw2ForPosDimSO(3)Monopoles} and
\eqref{eq:LowerBoundOnw2ForPosDimSO(3)MonopolesToGetMB} for such $u$.

This completes the proof for the second case.
\end{case}

\begin{case}[$X$ is even, $b^+(X)\ge 1$, and $b^-(X)\ge 1$]
 Because $Q_{X}$ is indefinite and even,   by Gompf and Stipsicz
\cite[Theorem 1.2.21, p.14]{GompfStipsicz}, $Q_X$ must contain a hyperbolic summand, that is, one
of the form
\[
H\cong \ZZ^2,\quad
Q_H\cong \begin{pmatrix} 0 & 1 \\ 1 & 0\end{pmatrix}.
\]
Let $h_1,h_2$ be a basis over $\ZZ$ for $H$ with
\begin{equation}
\label{eq:IntersectionFormForEvenCase}
Q_{\widetilde X}(h_1,h_2)=1,\,
Q_{\widetilde X}(e,h_i)=0, \,
\text{and}\,
Q_{\widetilde X}(h_i,h_i)=0\, \text{for $i=1,2$.}
\end{equation}
We write
\begin{equation}
\label{eq:c1s0_EvenCase}
c_1(\fs_0)=\ka^\perp -2a_1h_1-2a_2h_2,
\end{equation}
where $\ka^\perp\in H^2(\widetilde X;\ZZ)$ is $Q_{\widetilde X}$-orthogonal to $H$ and $e$.
We can assume that the coefficients of $h_1$ and $h_2$ in the preceding
expression for $c_1(\fs_0)$
are even
%PF3-12-2024  Requires explanation. "As before" wording is confusing because previous two cases were odd.
%TL4-4-2024: Removed "As before" and added X is even.
because $c_1(\fs_0)$ is characteristic and $X$ is even.
If $X$ is even, then $\si(X)\equiv 0\pmod 8$ by Gompf and Stipsicz
\cite[Lemma 1.2.20, p. 13]{GompfStipsicz}
and so $c_2(X)=2-2b_x(X)+2b^+(X)-\si(X)$ satisfies
\begin{equation}
\label{eq:c2Parity_ForXevenCase}
c_2(X)\equiv 0 \pmod 2.
\end{equation}
For $b\equiv 1\pmod 2$, there then is $N\in\ZZ$ such that
\begin{equation}
\label{eq:DefineNForXevenCase}
2N=-c_2( X)+b+1.
\end{equation}
For $u\in\ZZ$, we define $v_2,z_1\in\ZZ$ by
\begin{subequations}
\label{eq:Coefficients_of_wu_EvenCase}
\begin{align}
\label{eq:Coefficients_of_wu_EvenCase_z1}
z_1&=1+2u,
\\
\label{eq:Coefficients_of_wu_EvenCase_v2}
v_2&=N+2u+2u^2,
\end{align}
\end{subequations}
and define $w_u\in H^2(\widetilde X;\ZZ)$ by
\begin{equation}
\label{eq:wInCaseEeven}
w_u:= h_1+v_2h_2-z_1 e.
\end{equation}
Definition
\eqref{eq:wInCaseEeven}
and \eqref{eq:Coefficients_of_wu_EvenCase_z1}
imply that $w_u$ satisfies the Morgan--Mrowka condition \eqref{eq:MMCondition}
for all $u\in\ZZ$.

We compute
\begin{align*}
Q_{\widetilde X}(w_u,w_u)
&=
2v_2-z_1^2
\quad\text{(by \eqref{eq:IntersectionFormForEvenCase})}
\\
&=
2N+4u+4u^2
-1-4u-4u^2
\quad\text{(by \eqref{eq:Coefficients_of_wu_EvenCase})}
\\
&=
2N-1.
\end{align*}
Then \eqref{eq:DefineNForXevenCase} implies that for all $u\in\ZZ$, $w_u$ satisfies
\eqref{eq:w2Requirement}.

Finally, we compute
\begin{align*}
Q_{\widetilde X}(-c_1(\fs_0)-e,w_u)
&=
2a_1v_2+2a_2-z_1
\quad\text{(by \eqref{eq:c1s0_EvenCase}, \eqref{eq:wInCaseEeven}, and \eqref{eq:IntersectionFormForEvenCase})}
\\
&=
2a_1\left( N+2u+2u^2\right)
+
2a_2
-1-2u
\quad\text{(by \eqref{eq:Coefficients_of_wu_EvenCase})}
\\
&=
4a_1 u^2
+
\left( 4a_1-2\right) u
+
2a_1N+2a_2-1.
\end{align*}
If $a_1> 0$, then the coefficient of $u^2$  in $Q_{\widetilde X}(c_1(\fs_0)+e,w_u)$ is positive and so for $u$ sufficiently large, $Q_{\widetilde X}(c_1(\fs_0)+e,w_u)$ will satisfy
\eqref{eq:LowerBoundOnw2ForPosDimSO(3)Monopoles} and
\eqref{eq:LowerBoundOnw2ForPosDimSO(3)MonopolesToGetMB}.
If $a_1<0$, then by replacing $w_u$ with $-w_u$, we see that $w_u$ will satisfy
\eqref{eq:LowerBoundOnw2ForPosDimSO(3)Monopoles} and \eqref{eq:LowerBoundOnw2ForPosDimSO(3)MonopolesToGetMB} for $u$ sufficiently large.
If $a_1=0$, then the coefficient of $u^2$ in $Q_{\widetilde X}(c_1(\fs_0)+e,w_u)$ vanishes while that of $u$ is $-2$.  After replacing $w_u$ with $-w_u$,  $w_u$ will satisfy
\eqref{eq:LowerBoundOnw2ForPosDimSO(3)Monopoles} and
\eqref{eq:LowerBoundOnw2ForPosDimSO(3)MonopolesToGetMB} for $u$ sufficiently large.

This completes the proof for the third case and hence completes the proof of the Proposition.
\end{case}
\end{proof}

We can now give
\begin{proof}[Proof of Theorem \ref{mainthm:ExistenceOfSpinuForFlow}]
This proof is nearly identical to that of 
%TL11-28-2025: Updated: changed Theorem 3 to Theorem 7
\cite[Theorem 7]{Feehan_Leness_introduction_virtual_morse_theory_so3_monopoles}.

Let $\fs_0$ be the spin${}^c$ structure on $X$ in the hypothesis of Theorem \ref{mainthm:ExistenceOfSpinuForFlow} and let $\fs$ be the spin${}^c$ structure on $\widetilde X$ with $c_1(\fs)=-c_1(\fs_0)-e$. The blow-up formula in Theorem \ref{thm:SWBlowUp} implies that $\SW_{\widetilde X}(\fs)\neq 0$. The proof of the existence of the almost complex structure $J$ on $\widetilde X$ with $c_1(T\widetilde X,J)=c_1(\fs)=-c_1(\fs_0)-e$  is identical to that of 
%TL11-28-2025: Updated: changed Theorem 3 to Theorem 7
\cite[Item 3 in Theorem 7]{Feehan_Leness_introduction_virtual_morse_theory_so3_monopoles}.

Because $c_2(X)$ is even when $X$ is even (as discussed in the proof of Proposition \ref{prop:FeasibilityOnBlowUp} in the case that $X$ is even), the hypothesis in Theorem \ref{mainthm:ExistenceOfSpinuForFlow} that $c_2(X)>0$
implies that $c_2(X)>1$ when $X$ is even. Hence, the hypothesis that $c_2(X)>0$ implies that
\begin{equation}
\label{eq:Positive_c2}
0>-c_2( X)+b
\end{equation}  
where $b=0$ if $X$ is odd and $b=1$ if $X$ is even.
The assumption in Theorem \ref{mainthm:ExistenceOfSpinuForFlow} that $b^+(X)>1$ implies that the hypotheses of
Proposition \ref{prop:FeasibilityOnBlowUp} are satisfied.
Then Proposition \ref{prop:FeasibilityOnBlowUp} implies that there is $w\in H^2(\widetilde X;\ZZ)$ such that $w$ and $c_1(\fs)$ satisfy the conditions
\eqref{eq:w2Requirement},
\eqref{eq:MMCondition},
\eqref{eq:LowerBoundOnw2ForPosDimSO(3)Monopoles},
and \eqref{eq:LowerBoundOnw2ForPosDimSO(3)MonopolesToGetMB} of Proposition \ref{prop:Feasibility}.
%PF10-3-2024 Check from here to displayed equation
We let $\ft$ be the spin${}^u$ structure on $(\widetilde X,g)$ defined by $w$ and $c_1(\fs)$ in \eqref{eq:SpinuStructureFromw} and we claim that $\ft$ is feasible with respect to $c_1(\fs)$ in the sense of Definition \ref{maindefn:Feasibility}.

Proposition \ref{prop:Feasibility} implies that $\sM_{\ft}$ contains the image of $M_\fs$ under the embedding \eqref{eq:DefnOfIotaOnQuotient}.  Thus $\ft$ satisfies Item \eqref{item:Feasibility_NonZeroSW} of Definition \ref{maindefn:Feasibility}.

The positivity of the expected dimension of $\sM_{\ft}$, Item \eqref{item:Feasibility_PosExpDim} of  Definition \ref{maindefn:Feasibility}, follows from Item \eqref{item:ExpectedDimOfSO(3)Monopoles} of Proposition \ref{prop:Feasibility}.

Because $w$ satisfies \eqref{eq:MMCondition}, Item \eqref{item:MMCondition} of Proposition \ref{prop:Feasibility} implies that $\sM_{\ft}$ does not contain the gauge-equivalence class of a zero-section, split pair and thus
$\ft$ satisfies Item \eqref{item:Feasbility_NoZeroSectionSplit} of Definition \ref{maindefn:Feasibility}.

The inequality \eqref{eq:Positive_c2}  implies that
\[
  w^2 = -c_2( X)+b < 0.
\]
Thus, $p_1(\ft)[\widetilde X]=p_1(\su(E))[\widetilde X]$ will satisfy the fundamental  bounds \eqref{eq:p1_lower_bound_blowup}, proving that $\tilde\ft$ satisfies Item \eqref{item:Feasibility_LowerBound} of  Definition \ref{maindefn:Feasibility}.

%PF10-3-2024 Check below
Because $w$ and $c_1(\fs)$ satisfy \eqref{eq:LowerBoundOnw2ForPosDimSO(3)MonopolesToGetMB}, Item
\eqref{item:PositivevMBForLevelZeroSpinc}  of Proposition \ref{prop:Feasibility}
implies that the formal Morse--Bott index \eqref{eq:FormalMorseIndexIntroThm} is positive at all gauge-equivalence classes of split pairs in $\sM_{\ft}$.
%The equality between the formal Morse--Bott index and the virtual Morse--Bott index given by Corollary \ref{maincor:MorseIndexAtReduciblesOnSymplecticWithSO3MonopoleCharacteristicClasses} then implies that
Thus $\ft$ satisfies Item \eqref{item:Feasibility_PositivevBM} of Definition \ref{maindefn:Feasibility}. This completes the proof of the claim that $\ft$ is feasible.

Finally, because $w$ satisfies \eqref{eq:w2Requirement}, Item \eqref{item:ExpectedDimOfZeroSections} of Proposition \ref{prop:Feasibility} implies that
Item \eqref{item:Upper_bound_expected_dimension_ASD_moduli_space_blowup} of Theorem \ref{mainthm:ExistenceOfSpinuForFlow} holds.
\end{proof}

\begin{proof}[Proof of Corollary \ref{maincor:FeasibleIsNonEmpty}]
Let $\ft$ be the spin${}^u$ structure on the smooth blow-up $(\widetilde X,\tilde g)$ constructed in the proof 
of Theorem \ref{mainthm:ExistenceOfSpinuForFlow}.
Because $w$ satisfies \eqref{eq:MMCondition}, Item \eqref{item:MMCondition} of Proposition \ref{prop:Feasibility}
implies that $\sM_\ft$ contains no reducible, zero-section pairs.  Because $w$ and $c_1(\fs)$ satisfy \eqref{eq:LowerBoundOnw2ForPosDimSO(3)Monopoles} and because $\SW_{\widetilde X}(\fs)\neq 0$,
Item \eqref{item:ExpectedDimOfSO(3)Monopoles} of Proposition \ref{prop:Feasibility} implies that
$\sM_\ft^{*,0}$ is non-empty, completing the proof of the Corollary.
\end{proof}

\chapter{Gauge theory over almost Hermitian manifolds}
\label{chap:Moduli_space_Seiberg-Witten_monopoles_over_almost_Hermitian_four-manifolds}
%PF10-17-2024 Prune redundant material
In this section, we gather the ingredients that we shall need in order to analyze the non-Abelian monopole equations
%PF10-10-2024 Add ref
over almost Hermitian four-manifolds. In Section \ref{sec:Nijenhuis_tensor_and_components_exterior_derivative}, we review the definition of the Nijenhuis tensor and components of the exterior derivative on differential forms and exterior covariant derivative on vector-bundle-valued differential forms. Section \ref{sec:Dirac_operator_almost_Hermitian} contains a review of the Dirac operator over almost Hermitian manifolds and Section \ref{sec:MonopolesOverAlmostHermitian} recalls the structure of the non-Abelian monopole equations over almost Hermitian four-manifolds. In Section \ref{sec:KaehlerIDs}, we review the K\"ahler identities over almost K\"ahler and almost Hermitian manifolds that we shall need in our present work. Section \ref{sec:RealAndComplexInnerProductSpaces} concludes with a summary of our conventions regarding real and complex inner product spaces.

\begin{comment}
%TL5-20-2025: I think most of this is just copied from the previous monograph.  The basic definitions I recall seeing used include: 
1) Nijenjuis tensor although that's pretty standard
2) Splitting of cotangent bundle in $\La^{p,q}$.
3) $\mu$ and $\mu_A$.  In particular the independence of $\mu_A$ from $A$.  The relation \eqref{eq:Donaldson_Kronheimer_p_43_mu_bar_A_Nijenhuis_tensor_E-valued_one_forms} is used in \eqref{eq:DefineBarPartialWithMu}.
4) Definition of $d^{1,0$, $d^{1,1}$, etc.
5) Components of curvature
6) Canonical spinc structure (now appearing in prelim).  There is, I think, far more detail in 4.3 than we need.
7) Relation of Dirac operator with $\bar\partial+\bar\partial^*)$ and that factor of $\sqrt 2$
8) Lee form.
9) In 4.5, we give $\vol=(1/2)\omega^2=\star 1$ which is used to explain a lot of factors of 2
10) The equality (4.6.8) is useful because it says how $\La$ works w/r/t the Hodge decomp.
11) I suspect we need (4.6.9) for its definition of the adjoints 
12) We use the Kahler identifies (4.6.15) and (4.6.16)
13) 4.7 can be cut.
Note that for the definitions we need listed above, we only need the definitions not the justifications since those are in the previous monograph and could be referred to.
\end{comment}

\section{Nijenhuis tensor and components of the exterior derivative}
\label{sec:Nijenhuis_tensor_and_components_exterior_derivative}
Our exposition relies on Donaldson \cite{DonSW}, \cite[Section 4]{DonYangMillsInvar}, Gauduchon \cite{Gauduchon_1997}, and Kotschick \cite{KotschickSW}. Recall from \cite[p. 34]{DonYangMillsInvar} that
\begin{equation}
  \label{eq:d_sum_components_almost_complex_manifold}
  d = \partial + \mu + \bar\mu + \bar\partial \quad\text{on } \Omega^{p,q}(X),
\end{equation}
where in the following discussion we will see that the components $\mu$ and $\bar\mu$ may be expressed in terms of the Nijenhuis tensor $N_J \in C^\infty(\wedge^2(T^*X)\otimes TX)$, where by Salamon \cite[Equation (3.5)]{SalamonSWBook}
\begin{equation}
  \label{eq:Nijenhuis_tensor}
   N_J(v,w) = [v,w] + J[Jv,w] + J[v,Jw] - [Jv,Jw], \quad\text{for all } v, w \in C^\infty(TX).
\end{equation}
 % COMMENT Replace commented portion by solutions to Salamon's exercises 3.11 and 3.12:
 % By Remark 3.7,
 % \begin{align*}
 %   \tau^{2,0}(v,w) &= \frac{1}{4}(\tau(v,w) - \tau(Jv,Jw) - i\tau(Jv,Jw) - i\tau(v,Jw)),
 %   \\
 %   \tau^{1,1}(v,w) &= \frac{1}{2}(\tau(v,w) + \tau(Jv,Jw)),
 %   \\
 %   \tau^{0,2}(v,w) &= \frac{1}{4}(\tau(v,w) - \tau(Jv,Jw) + i\tau(Jv,Jw) + i\tau(v,Jw)),
 %  \end{align*}
 %  and
 % \begin{align*}
 %   N(Jv,Jw) &= [Jv,Jw] - J[v,Jw] - J[Jv,w] - [v,w] = -N(v,w),
 %   \\
 %   N(v,Jw) &= [v,Jw] + J[Jv,Jw] - J[v,w] + [Jv,w] = -JN(v,w),
 % \end{align*}
 % we obtain
 % \begin{align*}
 %   N^{0,2}(v,w) &= \frac{1}{4}(N(v,w) - N(Jv,Jw) + iN(Jv,Jw) + iN(v,Jw))
 %   \\
 %                &= \frac{1}{4}(N(v,w) + N(v,w) - iN(v,w) - iJN(v,w))
 %   \\
 %   &= \frac{1}{4}(N(v,w) - iJN(v,w)) + \frac{1}{4}(N(v,w) - iN(v,w)).
 % \end{align*}
 % Now $N(v,w) - iJN(v,w) \in C^\infty(TX^{1,0})$ and $N(v,w) - iN(v,w) \in C^\infty(TX^{0,1})$, so
 % \[
 %   \frac{1}{4}\alpha\circ N^{0,2} =
 % \]
 Following Salamon \cite[Exercise 3.12]{SalamonSWBook}, for any integer $r\geq 0$ one defines
 \[
   \iota(N):\Omega^r(X,\CC) \to \Omega^{r+1}(X,\CC)
 \]
 for $\alpha \in \Omega^1(X)$ by $\iota(N)\alpha = \alpha\circ N$ and in general by
 \[
   \iota(N)(\sigma\wedge\tau) := (\iota(N)\sigma)\wedge\tau + (-1)^{\deg\sigma}\sigma\wedge\iota(N)\tau
 \]
 for $\sigma\in\Omega^k(X,\CC)$ and $\tau\in\Omega^l(X,\CC)$. In particular,
 \[
   d\tau - \partial\tau - \bar\partial\tau = \frac{1}{4}\iota(N)\tau \in \Omega^{p+2,q-1}(X)\oplus \Omega^{p-1,q+2}(X)
 \]
 for all $\tau \in \Omega^{p,q}(X)$. Following Cirici and Wilson \cite[Section 2]{Cirici_Wilson_2020_harmonic}, we denote
 \[
   \mu\tau = d^{2,-1}\tau \quad\text{and}\quad \bar\mu\tau = d^{-1,2}\tau,
 \]
 so that
 %PF11-13-2025 Corrected expression for \bar\mu
 \begin{equation}
 \label{eq:MuComponentsAnd_Nijenjuis}
   \mu\tau =  \frac{1}{4}\pi_{p+2,q-1}(\iota(N)\tau) \quad\text{and}\quad \bar\mu\tau =
   % TL11-13-2025: Should it be q+2?
   %PF11-13-2025 Corrected
   \frac{1}{4}\pi_{p-1,q+2}(\iota(N)\tau).
 \end{equation}
 In particular, for $\alpha \in \Omega^{1,0}(X)$ and $\beta \in \Omega^{0,1}(X)$, we have
 \begin{align*}
   \mu\alpha = 0 \quad&\text{and}\quad \bar\mu\alpha = \frac{1}{4}\alpha\circ N \in \Omega^{0,2}(X),
   \\
   \mu\beta = \frac{1}{4}\beta\circ N \in \Omega^{2,0}(X) \quad&\text{and}\quad \bar\mu\beta = 0,
 \end{align*}
 in agreement with \cite[Exercise 3.11]{SalamonSWBook}.

We describe $(0,1)$ connections and their associated $\bar\partial$ operators on complex vector bundles over almost complex manifolds by adapting the discussion in Donaldson and Kronheimer \cite[Section 2.1.5]{DK} and Kobayashi \cite[Section 1.3]{Kobayashi_differential_geometry_complex_vector_bundles}, which in turn extend discussions of differential forms over complex manifolds described in
%PF10-10-2024 Add exact refs
Huybrechts \cite{Huybrechts_2005} and Wells \cite{Wells3}. Let $E$ be a complex vector bundle over a smooth almost complex manifold $(X,J)$ of real dimension $d=2n$. The splitting of forms with coefficients in $\CC$,
\[
  \wedge^\bullet(T_\CC^*X) = \bigoplus_{r=0}^d \wedge^r(T_\CC^*X)
  \quad\text{and}\quad
  \wedge^r(T_\CC^*X) = \bigoplus_{p+q=r} \wedge^{p,q}(T_\CC^*X),
\]
induces a splitting of forms with coefficients in $E$,
\[
  \wedge^\bullet(T_\CC^*X)\otimes E = \bigoplus_{r=0}^d \wedge^r(T_\CC^*X)\otimes E
  \quad\text{and}\quad
  \wedge^r(T_\CC^*X)\otimes E = \bigoplus_{p+q=r} \wedge^{p,q}(T_\CC^*X)\otimes E.
\]
We abbreviate $\Lambda^\bullet(X) = \wedge^\bullet(T_\CC^*X)$ and $\Lambda^r(X) = \wedge^r(T_\CC^*X)$ and $\Lambda^{p,q}(X) = \wedge^{p,q}(T_\CC^*X)$ and similarly for forms with coefficients in $E$. We write
\[
   \Omega^\bullet(E) = C^\infty(\Lambda^\bullet(X)\otimes E), 
   \quad
   \Omega^r(E) = C^\infty(\Lambda^r(X)\otimes E),
   \quad
   \Omega^{p,q}(E) = C^\infty(\Lambda^{p,q}(X)\otimes E).
\]
We may denote $\Omega^{0,0}(E)$ by $\Omega^0(E)$ and $\Omega^{0,0}(X)$ by $\Omega^0(X;\CC)$ or vice versa.

A covariant derivative on $E$ is homomorphism \cite[Equation (2.1.1)]{DK}, \cite[Equation (1.1.1)]{Kobayashi_differential_geometry_complex_vector_bundles}
\begin{equation}
\label{eq:Donaldson_Kronheimer_2-1-1_nabla_A_homomorphism}  
  \nabla_A:\Omega^0(E) \to \Omega^1(E)
\end{equation}
of vector spaces over $\CC$ that obeys the Leibnitz rule
\begin{equation}
\label{eq:Donaldson_Kronheimer_p_32_nabla_A_Leibnitz_rule}  
  \nabla_A(f\sigma) = (df)\sigma + f\nabla_A\sigma, \quad\text{for } f \in \Omega^0(X;\CC), \sigma \in \Omega^0(E).
\end{equation}
The exterior covariant derivative \cite[Equation (2.1.12)]{DK}
\begin{equation}
\label{eq:Donaldson_Kronheimer_2-1-12}  
  d_A:\Omega^r(E) \to \Omega^{r+1}(E)
\end{equation}
is uniquely determined by the properties \cite[p. 35]{DK}
\begin{subequations}
\label{eq:Donaldson_Kronheimer_p_35}
\begin{align}
  \label{eq:Donaldson_Kronheimer_p_35_i_d_A_is_nabla_A}  
  d_A &= \nabla_A \quad\text{on } \Omega^0(E),
  \\
  \label{eq:Donaldson_Kronheimer_p_35_ii_generalized_Leibnitz_rule}  
  d_A(\varphi\wedge\sigma) &= (d\varphi)\wedge\sigma + (-1)^r\varphi\wedge d_A\sigma, \quad\text{for } \varphi \in \Omega^r(X;\CC), \sigma \in \Omega^s(E).
\end{align}  
\end{subequations}
We write $\delta_A$ for any one of the four components of $d_A$,
\begin{subequations}
  \label{eq:d_A_components_almost_complex_manifold}
  \begin{align}
    \label{eq:del_A}
    \partial_A = d_A^{1,0} &= \pi_{p+1,q}d_A:\Omega^{p,q}(E) \to \Omega^{p+1,q}(E),
    \\
    \label{eq:mu_A}
    \mu_A = d_A^{2,-1} &= \pi_{p+2,q-1}d_A:\Omega^{p,q}(E) \to \Omega^{p+2,q-1}(E),
    \\
    \label{eq:mu_bar_A}
    \bar\mu_A = d_A^{-1,2} &= \pi_{p-1,q+2}d_A:\Omega^{p,q}(E) \to \Omega^{p-1,q+2}(E),
    \\
    \label{eq:del_bar_A}
    \bar\partial_A = d_A^{0,1} &= \pi_{p,q+1}d_A:\Omega^{p,q}(E) \to \Omega^{p,q+1}(E),
\end{align}  
\end{subequations}
and observe that (see Donaldson \cite[p. 34]{DonYangMillsInvar})
\begin{equation}
  \label{eq:d_A_sum_components_almost_complex_manifold}
  d_A = \partial_A + \mu_A + \bar\mu_A + \bar\partial_A \quad\text{on } \Omega^{p,q}(E).
\end{equation}
We temporarily allow that the components $\mu_A$ and $\bar\mu_A$ of $d_A$ may depend on the connection $A$ and later drop the subscript `$A$' when we see that they are independent of $A$. By taking $r=0$ in \eqref{eq:Donaldson_Kronheimer_p_35_ii_generalized_Leibnitz_rule} and applying the definitions \eqref{eq:mu_A} and \eqref{eq:mu_bar_A} of $\mu_A$ and $\bar\mu_A$, we see that   
\begin{subequations}
\begin{align}
\label{eq:Donaldson_Kronheimer_p_43_mu_A_Leibnitz_rule} 
  \mu_A(f\sigma) &= f\mu_A\sigma,
  \\
\label{eq:Donaldson_Kronheimer_p_43_mu_bar_A_Leibnitz_rule}   
  \bar\mu_A(f\sigma) &= f\bar\mu_A\sigma,
                       \quad\text{for } f \in \Omega^0(X;\CC), \ \sigma \in \Omega^{p,q}(E), 
\end{align}
\end{subequations}
so these operators are zeroth order, while \eqref{eq:Donaldson_Kronheimer_p_32_nabla_A_Leibnitz_rule} and \eqref{eq:Donaldson_Kronheimer_p_35_i_d_A_is_nabla_A} and the definitions \eqref{eq:del_A} and \eqref{eq:del_bar_A} imply that the operators $\partial_A$ and $\bar\partial_A$ inherit the Leibnitz rules \cite[Equation (2.1.45)]{DK}, \cite[Equation (1.3.1)]{Kobayashi_differential_geometry_complex_vector_bundles}
\begin{subequations}
\begin{align}
\label{eq:Donaldson_Kronheimer_p_43_del_A_Leibnitz_rule}  
  \partial_A(f\sigma) &= (\partial f)\sigma + f\partial_A\sigma, 
  \\
\label{eq:Donaldson_Kronheimer_p_43_del_bar_A_Leibnitz_rule}    
  \bar\partial_A(f\sigma) &= (\bar\partial f)\sigma + f\bar\partial_A\sigma, \quad\text{for } f \in \Omega^0(X;\CC), \ \sigma \in \Omega^{p,q}(E),
\end{align}
\end{subequations}
Next observe that \eqref{eq:Donaldson_Kronheimer_p_35_i_d_A_is_nabla_A} and the definitions \eqref{eq:mu_A} and \eqref{eq:mu_bar_A} of $\mu_A$ and $\bar\mu_A$ yield
\begin{equation}
\label{eq:Donaldson_Kronheimer_p_43_mu_A_and_mu_bar_A_sections_of_E}  
  \mu_A\sigma = 0 \quad\text{and}\quad \bar\mu_A\sigma = 0, \quad\text{for } \sigma \in \Omega^0(E),
\end{equation}
and hence \eqref{eq:Donaldson_Kronheimer_p_35_ii_generalized_Leibnitz_rule} yields
\begin{subequations}
\label{eq:Donaldson_Kronheimer_p_43_mu_A_and_mu_bar_A_generalized_Leibnitz_rules}  
\begin{align}
\label{eq:Donaldson_Kronheimer_p_43_mu_A_generalized_Leibnitz_rule}
  \mu_A(\varphi\otimes\sigma) &= (\mu\varphi)\otimes\sigma,
  \\
  \label{eq:Donaldson_Kronheimer_p_43_mu_bar_A_generalized_Leibnitz_rule}
  \bar \mu_A(\varphi\otimes\sigma) &= (\bar\mu\varphi)\otimes\sigma, \quad\text{for } \varphi \in \Omega^{p,q}(X), \ \sigma \in \Omega^0(E).
\end{align}
\end{subequations}
Using the expressions \eqref{eq:MuComponentsAnd_Nijenjuis} for $\mu$ and $\bar\mu$ on $\Omega^{p,q}(X)$ in terms of the Nijenhuis tensor $N_J$, we obtain from \eqref{eq:Donaldson_Kronheimer_p_43_mu_A_and_mu_bar_A_generalized_Leibnitz_rules} that
\begin{subequations}
\label{eq:Donaldson_Kronheimer_p_43_mu_A_and_mu_bar_A_generalized_Leibnitz_rule_Nijenhuis_tensor}  
\begin{align}
  % TL6-19-2025: I think the \mu and \bar\mu are reversed here in comparison with \eqref{eq:MuComponentsAnd_Nijenjuis}
  %PF11-13-2025 Corrected
\label{eq:Donaldson_Kronheimer_p_43_mu_A_generalized_Leibnitz_rule_Nijenhuis_tensor}
  \mu_A(\varphi\otimes\sigma) &= \frac{1}{4}(\pi_{p+2,q-1}\iota_{N_J}\varphi)\otimes\sigma,
  \\
  \label{eq:Donaldson_Kronheimer_p_43_mu_bar_A_generalized_Leibnitz_rule_Nijenhuis_tensor}
  \bar\mu_A(\varphi\otimes\sigma) &= \frac{1}{4}(\pi_{p-1,q+2}\iota_{N_J}\varphi)\otimes\sigma,
                             \quad\text{for } \varphi \in \Omega^{p,q}(X), \ \sigma \in \Omega^0(E).
\end{align}
\end{subequations}
When the $(p,q)$ forms in $\Omega^{p,q}(X)$ are specialized to one-forms in $\Omega^{1,0}(X)$ or $\Omega^{0,1}(X)$, we obtain
\begin{subequations}
  % TL6-19-2025: \mu/\bar\mu reversal propagated to here
  %PF11-13-2025 Corrected
\label{eq:Donaldson_Kronheimer_p_43_mu_A_and_mu_bar_A_generalized_Leibnitz_rule_Nijenhuis_tensor_one_forms}
\begin{align}
\label{eq:Donaldson_Kronheimer_p_43_mu_A_generalized_Leibnitz_rule_Nijenhuis_tensor_one_forms}
  \mu_A(\alpha\otimes\sigma) &= 0 \quad\text{and}\quad \mu_A(\beta\otimes\sigma) = \frac{1}{4}(\beta\circ N_J)\otimes\sigma,
  \\
  \label{eq:Donaldson_Kronheimer_p_43_mu_bar_A_generalized_Leibnitz_rule_Nijenhuis_tensor_one_forms}
  \bar\mu_A(\alpha\otimes\sigma) &= \frac{1}{4}(\alpha\circ N_J)\otimes\sigma \quad\text{and}\quad \bar\mu_A(\beta\otimes\sigma) = 0,
  \\
  \notag
  &\qquad\qquad\text{for } \alpha \in \Omega^{1,0}(X), \ \beta \in \Omega^{0,1}(X), \ \sigma \in \Omega^0(E).
\end{align}
\end{subequations}
In our applications, we shall only need expressions for $\mu_A$ and $\bar\mu_A$ on one-forms with coefficients in a complex vector bundle. With respect to a local trivialization $E \restriction U \cong U\times \CC^r$, we may express $a' \in \Omega^{1,0}(E)$ and $a'' \in \Omega^{0,1}(E)$ as
\[
  a' = \alpha\otimes\sigma \quad\text{and}\quad a'' = \beta\otimes\tau
\]
so that \eqref{eq:Donaldson_Kronheimer_p_43_mu_A_and_mu_bar_A_generalized_Leibnitz_rule_Nijenhuis_tensor_one_forms} yields
\[
  % TL6-19-2025: \mu/\bar\mu reversal propagated to here
  %PF11-13-2025 Corrected
  \mu_Aa'' = \frac{1}{4}(\beta\circ N_J)\otimes\sigma \quad\text{and}\quad
  \bar\mu_Aa' = \frac{1}{4}(\alpha\circ N_J)\otimes\tau,
\]
while $\bar\mu_Aa'' = 0$ and $\mu_Aa' = 0$. Consequently, we may write
\begin{subequations}
  % TL6-19-2025: \mu/\bar\mu reversal propagated to here
  %PF11-13-2025 Corrected
\label{eq:Donaldson_Kronheimer_p_43_mu_A_and_mu_bar_A_Nijenhuis_tensor_E-valued_one_forms}
\begin{align}
\label{eq:Donaldson_Kronheimer_p_43_mu_A_Nijenhuis_tensor_E-valued_one_forms}
  \mu_Aa'' &= \frac{1}{4}a''\circ N_J \quad\text{and}\quad \mu_Aa' = 0,
  \\
  \label{eq:Donaldson_Kronheimer_p_43_mu_bar_A_Nijenhuis_tensor_E-valued_one_forms}
  \bar\mu_Aa' &= \frac{1}{4}a'\circ N_J \quad\text{and}\quad \bar\mu_Aa'' = 0,
                 \quad\text{for } a' \in \Omega^{1,0}(E), \ a'' \in \Omega^{0,1}(E),
\end{align}
\end{subequations}
and thus we shall henceforth denote $\mu_A$ by $\mu$ and denote $\bar\mu_A$ by $\bar\mu$ to underline the independence of these operators from the connection $A$.

\section{Components of the curvature}
\label{subsec:Components_of_curvature}
Recall that $F_A = d_A^2 \in \Omega^2(\End(E,H))$ and $d_A = \partial_A + \bar\partial_A$ on $\Omega^0(E)$ and $d_A = \partial_A + \mu + \bar\mu + \bar\partial_A$ on $\Omega^1(E)$, so
\begin{align*}
  d_A^2 &= (\partial_A^2 + \mu\bar\partial_A) + (\partial_A\bar\partial_A + \bar\partial_A\partial_A) + (\bar\mu\partial_A + \bar\partial_A^2)
  \\
        &= F_A^{2,0} + F_A^{1,1} + F_A^{0,2}
          \\
          &\quad\in \Omega^{2,0}(\End(E,H)) \oplus \Omega^{1,1}(\End(E,H)) \oplus \Omega^{0,2}(\End(E,H)).
\end{align*}
Because $d^2=0$ on $\Omega^0(\CC)$, we obtain the identities
\[
  \partial^2 + \mu\bar\partial = 0, \quad \partial\bar\partial + \bar\partial\partial = 0, \quad \bar\mu\partial + \bar\partial^2 = 0.
\]  
Recall that for any $f\in\Omega^0(\CC)$ and $s\in\Omega^0(E)$, we have $F_A(fs) = fF_As$ and the same holds for the components $F_A^{2,0}$, $F_A^{1,1}$, and $F_A^{0,2}$. That is, they are all zeroth-order operators. However,
\begin{align*}
  \bar\partial_A^2(fs) &= \bar\partial_A((\bar\partial f)s + f\bar\partial_As)
  \\
  &= (\bar\partial^2 f)s + \bar\partial_A s\wedge \bar\partial f + \bar\partial f\wedge \bar\partial_As + f\bar\partial_A^2s
\\
  &= -(\bar\mu\partial f)s + \bar\partial_A^2s,
\end{align*}
and similarly we compute $\partial_A^2(fs)$, that is,
\begin{align*}
  \partial_A^2(fs) &= -(\mu\bar\partial f)s + f\partial_A^2s,
  \\
  \bar\partial_A^2(fs) &= -(\bar\mu\partial f)s + \bar\partial_A^2s.
\end{align*}  
Therefore, because $F_A^{2,0}$ and $F_A^{0,2}$ are zeroth-order operators, we must have
\begin{align*}
  \mu\bar\partial_A(fs) &= (\bar\mu\partial f)s + f\mu\bar\partial_As,
  \\
  \bar\mu\partial_A(fs) &= (\mu\bar\partial f)s + f\bar\mu\partial_As.
\end{align*}
We compute
\[
  F_{u(A)}s = d_{u(A)}^2s = (ud_Au^{-1})(ud_Au^{-1}s) = ud_A^2(u^{-1}s) = uF_Au^{-1}s,
\]
using the fact that $F_A$ is a zeroth-order operator and thus
\[
  F_{u(A)} = uF_Au^{-1}.
\]
The components $F_A^{2,0}$, $F_A^{1,1}$, and $F_A^{0,2}$ of $F_A$ have the corresponding transformation rules.
\section{Spin${}^u$ and canonical spin${}^c$ structures over almost Hermitian manifolds}
\label{subsec:Spinc_and_spinu_structures_over_almost_Hermitian_manifolds}
We review Clifford multiplication over almost Hermitian four-manifolds $(X,g,J)$, first for the canonical \spinc structure and then for any \spinc structure. We closely follow Salamon \cite{SalamonSWBook}, though we also rely on Donaldson \cite{DonSW}, Kotschick \cite{KotschickSW}, Morgan \cite{MorganSWNotes}, Nicolaescu \cite{NicolaescuSWNotes}, Taubes \cite[Section 1]{TauSymp}, and Witten \cite{Witten}. We begin with a discussion of certain identities for the interior product of complex differential forms that we shall need. For a smooth almost Hermitian manifold $(X,g,J)$ of dimension $2n$, smooth form $\alpha \in \Omega^{0,k}(X)$, and integers $\ell\geq k\geq 0$, we recall that left multiplication by $\alpha$ and interior multiplication,
\begin{align}
  L_\alpha:\Omega^{0,\ell-k}(X) \ni \varpi
  &\mapsto \alpha\wedge\varpi \in \Omega^{0,\ell}(X),
  \\
  \iota(\bar\alpha):\Omega^{0,\ell}(X) \ni \chi
  &\mapsto \iota(\bar\alpha)\chi \in \Omega^{0,\ell-k}(X),
\end{align}
are related by (see Salamon \cite[Lemma 3.5, p. 60]{SalamonSWBook}) 
\[
  \langle L_\alpha\varpi, \chi\rangle_{\Lambda^{0,\ell}(X)}
  = \langle \alpha\wedge\varpi, \chi\rangle_{\Lambda^{0,\ell}(X)}
  = 2^k(-1)^{k(k-1)/2}\langle\varpi, \iota(\bar\alpha)\chi\rangle_{\Lambda^{0,\ell-k}(X)},
\]
where $\iota(\bar\alpha)$, for $\alpha = \sum \alpha_I e_I''$, is defined for $\ell\geq 1$ by
\begin{equation}
  \label{eq:Salamon_3-3}
  \iota(\bar\alpha)\chi := \sum \bar\alpha_I\iota(e_I)\chi, \quad e_I'' = e_{i_1}\wedge\cdots\wedge e_{i_\ell}'',
\end{equation}
where the sum runs over all multi-indices $I = \{i_1,\ldots,i_\ell\}$ with $i_1<\cdots<i_\ell$ and
\[
  \iota(e_I)\chi := \iota(e_{i_1})\cdots\iota(e_{i_\ell})\chi,
\]
and $\iota(v):\Omega^\ell(X,\CC)\to\Omega^{\ell-1}(X,\CC)$ denotes interior product,
\begin{equation}
\label{eq:Interior_product}
\iota(v)\chi(v_2,\ldots v_\ell) := \chi(v,v_2,\ldots v_\ell),
\quad\text{for all } v,v_2\ldots,v_\ell \in C^\infty(TX),
\end{equation}  
and $\{e_1,Je_1,\ldots,e_n,J_ne_n\}$ is a local orthonormal frame for $TX$ and $\{e_1^*,Je_1^*,\ldots,e_n^*,J_ne_n^*\}$ is the corresponding local orthonormal frame for $T^*X$ defined by (see Salamon \cite[Section 3.1, pp. 58--59]{SalamonSWBook}) for any $v \in C^\infty(TX)$ by 
\begin{multline*}
  v^* := g(\cdot,v) \in \Omega^1(X,\RR)
  \\
  v' := 2\pi_{1,0}v^* = \langle \cdot,v\rangle_{T_\CC X} \in \Omega^{1,0}(X),
  \\
  \quad\text{and}\quad
  v'' := 2\pi_{0,1}v^* = \langle v,\cdot\rangle_{T_\CC X} = \overline{v'} \in \Omega^{0,1}(X),
\end{multline*}
and, for any $\theta \in \Omega^1(X,\RR)$,
\[
  \pi_{1,0}\theta := \frac{1}{2}(\theta + iJ\theta) \quad\text{and}\quad \pi_{0,1}\theta := \frac{1}{2}(\theta - iJ\theta).
\]
giving $J\pi_{1,0}\theta = -i\theta$ and $J\pi_{0,1}\theta = i\theta$. (Note that Salamon \cite[Section 3.1, p. 55]{SalamonSWBook} adopts the \emph{opposite convention} for the Hermitian inner product on a complex Hilbert space $\fH$, taking $\langle x,y \rangle_\fH$ to be complex linear in $y$ and complex antilinear in $x$.)

Recall that the Hermitian inner product on $\Lambda_\CC^\bullet(X) := \wedge^\bullet(T_\CC^*X)$ is defined by the relation (see Huybrechts \cite[Section 1.2, p. 33]{Huybrechts_2005})
\begin{equation}
  \label{eq:Huybrechts_Hermitian_inner_product_complex_forms_p_33}
  \langle \kappa, \mu \rangle_{\Lambda_\CC^\bullet(X)}\cdot\vol := \kappa \wedge \star \bar\mu,
  \quad\text{for all } \kappa, \mu \in \Omega^\bullet(X,\CC).
\end{equation}
Therefore, given $\alpha \in \Omega^{k,l}(X)$, we have (using $\star 1 = \vol$ by \cite[Proposition 1.2.20 (i)]{Huybrechts_2005})
\[
  \langle \alpha\wedge\varpi, \chi\rangle_{\Lambda_\CC^\bullet(X)}\cdot\vol
  =
  \langle \alpha\wedge\varpi, \chi\rangle_{\Lambda_\CC^\bullet(X)}\cdot(\star 1)
  =
  \alpha\wedge\varpi\wedge\star\bar\chi,
\]
which gives (using $\star^2 1 = 1$ by Huybrechts \cite[Proposition 1.2.20 (iii), p. 32]{Huybrechts_2005})
\begin{multline}
  \label{eq:Huybrechts_Hermitian_inner_product_complex_forms_p_33_Hodge_star}
  \langle\alpha\wedge\varpi, \chi\rangle_{\Lambda^{p,q}(X)}
  =
  \star(\alpha\wedge\varpi \wedge \star\bar\chi),
  \\
  \text{for all } \alpha \in \Omega^{k,l}(X), \varpi \in \Omega^{p,q}(X),
  \text{ and } \chi \in \Omega^{p+k,q+l}(X).
\end{multline}
If $\Lambda_\alpha:\Omega^{p+k,q+l}(X)\to\Omega^{p,q}(X)$ denotes the adjoint of left exterior multiplication by $\alpha \in \Omega^{k,l}(X)$,
\[
  L_\alpha:\Omega^{p,q}(X) \ni \varpi \mapsto \alpha\wedge\varpi \in \Omega^{p+k,q+l}(X),
\]
with respect to the Hermitian inner product on $\Lambda_\CC^\bullet(X)$, that is,
\begin{equation}
  \label{eq:Left_multiplication_of_pq_forms_by_kl_form_Hermitian_adjoint}
  \langle L_\alpha\varpi, \chi\rangle_{\Lambda^{p+k,q+l}(X)}
  =
  \langle\varpi, \Lambda_\alpha\chi\rangle_{\Lambda^{p,q}(X)},
  \quad\text{for all } \varpi \in \Omega^{p,q}(X), \chi \in \Omega^{p+k,q+l}(X),
\end{equation}
then, noting that $\deg\bar\chi=(q+l,p+k)$ and $\deg\star\bar\chi=(n-p-k,n-q-l)$ and thus $\deg(\alpha \wedge \star\bar\chi)=(n-p,n-q)$ and $\star^2=(-1)^{(2n-p-q)(p+q)}$ on $\Omega^{n-p,n-q}(X)$ by Huybrechts \cite[Proposition 1.2.20 (iii), p. 32]{Huybrechts_2005},
\begin{align*}
  \langle\varpi, \Lambda_\alpha\chi\rangle_{\Lambda^{p,q}(X)}
  &= \langle \alpha\wedge\varpi, \chi\rangle_{\Lambda^{p+k,q+l}(X)} \quad\text{(by \eqref{eq:Left_multiplication_of_pq_forms_by_kl_form_Hermitian_adjoint})}
  \\
  &= \star(\alpha\wedge\varpi \wedge \star\bar\chi) \quad\text{(by \eqref{eq:Huybrechts_Hermitian_inner_product_complex_forms_p_33_Hodge_star})}
  \\
  &= (-1)^{(k+l)(p+q)}\star(\varpi \wedge \alpha \wedge \star\bar\chi)
    \quad\text{(as $\deg\alpha = (k,l)$ and $\deg\varpi=(p,q)$)}
  \\
  &= (-1)^{(k+l)(p+q)}(-1)^{(2n-p-q)(p+q)}\star(\varpi \wedge \star^2(\alpha \wedge \star\bar\chi))
  \\
  &\qquad\text{(as $\deg(\alpha \wedge \star\bar\chi)=(n-p,n-q)$ and $\star^2=(-1)^{(2n-p-q)(p+q)}$)}
  \\
  &= (-1)^{(k+l)(p+q)}(-1)^{(p+q)^2}\left\langle\varpi, \overline{\star(\alpha \wedge \star\bar\chi)} \right\rangle_{\Lambda^{p,q}(X)}
   \quad\text{(by \eqref{eq:Huybrechts_Hermitian_inner_product_complex_forms_p_33_Hodge_star})}    
  \\
  &= (-1)^{(k+l)(p+q)}(-1)^{(p+q)^2}\left\langle\varpi, \star(\bar\alpha \wedge \star\chi) \right\rangle_{\Lambda^{p,q}(X)}
  \\
  &\qquad \text{(by $\CC$-linearity \cite[Section 1.2, p. 33]{Huybrechts_2005} of Hodge $\star:\Lambda_\CC^\bullet(X)\to \Lambda_\CC^\bullet(X)$)},
\end{align*}
and therefore\footnote{Compare the identity given by Warner \cite[Equation (2.11.2), p. 61 and Exercise 2.14, p. 80]{Warner} in the case of a Riemannian manifold and $\deg\alpha=1$.}
\begin{equation}
  \label{eq:Left_multiplication_of_p+kq+l_forms_by_kl_form_Hermitian_adjoint_formula}
  \Lambda_\alpha\chi
  =
  (-1)^{(k+l)(p+q)}(-1)^{(p+q)^2} \star(\bar\alpha \wedge \star\chi),
  \quad\text{for all } \chi \in \Omega^{p+k,q+l}(X).
\end{equation}
Note that if $k+l\equiv 1 \pmod{2}$, then the identity \eqref{eq:Left_multiplication_of_p+kq+l_forms_by_kl_form_Hermitian_adjoint_formula} simplifies to give
\[
  \Lambda_\alpha\chi = \star(\bar\alpha\wedge*\chi), \quad\text{for all } \chi \in \Omega^{p+k,q+l}(X).
\]
We recall the important

\begin{defn}[Canonical \spinc structure]
\label{defn:Canonical_spinc_bundles}
(See Donaldson \cite[Equation (15)]{DonSW}, Kotschick \cite[Fact 2.1]{KotschickSW}, Salamon \cite[Lemma 4.52, p. 141]{SalamonSWBook}, and Taubes \cite[Section 1]{TauSymp}.)
Let $(X,g,J)$ be an almost Hermitian four-manifold. Then the \emph{canonical \spinc structure} $\fs_{\can}=(\rho_\can,W_\can)$ on $X$ is defined by the bundle
\begin{equation}
  \label{eq:Canonical_spinc_bundles}
  W_\can^+ := \Lambda^{0,0}(X) \oplus \Lambda^{0,2}(X), \quad W_\can^- := \Lambda^{0,1}(X),
  \quad W_\can = W_\can^+\oplus W_\can^-,
\end{equation}
and the map $\rho_\can:TX\to \Hom(W_\can)$ defined for all $Y\in C^\infty(TX)$ and $\phi \in \Omega^{0,\bullet}(X)$ by
\begin{equation}
  \label{eq:Canonical_Clifford_multiplication}
    \rho_\can(Y)\phi = \frac{1}{\sqrt{2}}Y''\wedge\phi - \sqrt{2}\iota(Y)\phi,
\end{equation}
where $\iota(Y):\Omega^{0,q}(X) \to \Omega^{0,q-1}(X)$ denotes interior product for $q\geq 1$ and
\[
  Y'' = 2\pi_{0,1}Y^* = \langle Y,\cdot\rangle_{T_\CC X} \in \Omega^{0,1}(X).
\]
\qed\end{defn}
The \emph{canonical line bundle} of an almost complex manifold $(X,J)$ of real dimension $2n$ is\footnote{If $X$ is a complex manifold, then the holomorphic tangent bundle $\sT_X$ (see Huybrechts \cite[Definition 2.2.14, p. 71]{Huybrechts_2005}) is naturally isomorphic as a complex vector bundle to $T^{1,0}X$ by \cite[Proposition 2.6.4 (ii), p. 104]{Huybrechts_2005}.} (compare Donaldson \cite[Section 4, p. 56]{DonSW}, Fine \cite[Definition 1.6]{Fine_2012}, Griffiths and Harris \cite[Chapter 1, Section 2, p. 146]{GriffithsHarris}, or Wells \cite[Chapter VI, Section 1, p. 218]{Wells3})
% PF9-14-2023 Berline, Getzler, Vergne p. 140, section 3.6 call the **dual** of the below the canonical line bundle?
%PF10-1-2023 Definition seems wrong? BHPV also define K_X as \wedge^n((T^{1,0}X)^*) or \wedge^n(\sT_X^*)
%TL12-6-2023: Looking at p. 23 of GH, $\La^{p,0}X=\La^p T^*X'$ and I think $T^*'X$ is the subspace generated by $dz^i$ (since it's dual to the $\rd/\rd z^i$ and fits with the notation on p. 24 of GH.
\begin{equation}
  \label{eq:DefineCanonicalLineBundle}
  %PF1-9-2024 See Huybrechts Definition 2.2.14, p. 71, and Huybrechts Proposition 2.6.4 (ii), p. 105
  K_X := \Lambda^{n,0}(X) = \wedge^n(T^{1,0}X),
\end{equation}
and the \emph{anti-canonical line bundle} is its dual line bundle (see Fine \cite[Definition 1.6]{Fine_2012} or
Gauduchon \cite[Section 2, Remark 5, p. 275]{Gauduchon_1997}),
\begin{equation}
  \label{eq:DefineAntiCanonicalLineBundle}
  K_X^* := \Hom_\CC(K_X,\CC).
\end{equation}
Then by Kotschick \cite[Fact 2.1]{KotschickSW}, the first Chern class of $\fs_{\can}$ is
\begin{equation}
\label{eq:c1_of_CanonicalSpinc}
c_1(\fs_{\can}) = -c_1(K_X) = c_1(X).
\end{equation}
According to Salamon \cite[Lemma 3.4, p. 141]{SalamonSWBook},
\[
  2\langle\iota(Y)\phi, \chi\rangle_{\Lambda^{0,q-1}(X)} = \langle \phi, Y''\wedge\chi\rangle_{\Lambda^{0,q}(X)},
  \quad\text{for all } \chi \in \Omega^{0,q-1}(X) \text{ and } \phi \in \Omega^{0,q}(X).
\]
But then, as in the derivation of \eqref{eq:Left_multiplication_of_p+kq+l_forms_by_kl_form_Hermitian_adjoint_formula},
\begin{align*}
  2\langle\chi, \iota(Y)\phi\rangle_{\Lambda^{0,q-1}(X)}
  &= \langle Y''\wedge\chi, \phi\rangle_{\Lambda^{0,q}(X)}
  \\
  &= \star (Y''\wedge\chi \wedge \star\bar\phi)
    \quad\text{(by \eqref{eq:Huybrechts_Hermitian_inner_product_complex_forms_p_33_Hodge_star})}
  \\
  &= (-1)^{q-1}\star (\chi \wedge Y''\wedge \star\bar\phi) \quad\text{(as $\deg\chi = q-1$)}
  \\
  &= (-1)^{q-1}(-1)^{(q-1)(2n-q+1)} \star (\chi \wedge \star (\star(Y''\wedge \star\bar\phi)))
  \\
  &= (-1)^{q-1}(-1)^{(q-1)^2} \star (\chi \wedge \star (\star(Y''\wedge \star\bar\phi)))
  \\
  &= (-1)^{2(q-1)} \star (\chi \wedge \star (\star(Y''\wedge \star\bar\phi)))
  \\
  &= \star (\chi \wedge \star (\star(Y''\wedge \star\bar\phi)))  
  \\
  &= \left\langle \chi, \overline{\star(Y''\wedge \star\bar\phi)} \right\rangle_{\Lambda^{0,q-1}(X)}
    \quad\text{(by \eqref{eq:Huybrechts_Hermitian_inner_product_complex_forms_p_33_Hodge_star})}
  \\
  &= \left\langle \chi, \star(\overline{Y''}\wedge \star\phi) \right\rangle_{\Lambda^{0,q-1}(X)},
  \\
  &\qquad \text{(by $\CC$-linearity \cite[Section 1.2, p. 33]{Huybrechts_2005} of Hodge $\star:\Lambda_\CC^\bullet(X)\to \Lambda_\CC^\bullet(X)$)},    
\end{align*}
where to obtain the fourth equality, we used the fact that $\deg\bar\phi = \deg\phi = q$, so $\deg\star\bar\phi = 2n-q$, and $\deg(Y''\wedge \star\bar\phi) = 2n-q+1$, and thus $\star^2=(-1)^{(q-1)(2n-q+1)}$ on $\Omega^{2n-q+1}(X)$ by Warner \cite[Equation (6.1.1)]{Warner}. Therefore,
\begin{equation}
  \label{eq:Salamon_lemma_3-4_formula}
  2\iota(Y)\phi = \star(Y'\wedge \star\phi), \quad\text{for all } Y\in C^\infty(TX)
  \text{ and } \phi \in \Omega^{0,q}(X).
\end{equation}
Hence, the Clifford multiplication \eqref{eq:Canonical_Clifford_multiplication} may be rewritten as
\begin{equation}
  \label{eq:Canonical_Clifford_multiplication_Hodge_star}
  \rho_\can(Y)\phi = \frac{1}{\sqrt{2}}\left(Y''\wedge\phi - \star(Y'\wedge \star\phi)\right),
  \quad\text{for all } Y\in C^\infty(TX) \text{ and } \phi \in \Omega^{0,q}(X).
\end{equation}
Note that $Y^* = \pi_{1,0}Y^* + \pi_{0,1}Y^* = \frac{1}{2}(Y'+Y'')$. For $a \in \Omega^1(X,\RR)$, we apply the formula \eqref{eq:Canonical_Clifford_multiplication_Hodge_star} with $Y \in C^\infty(TX)$ uniquely defined by $a = g(\cdot,Y) = Y^*$ and $a = a'+a'' = \frac{1}{2}(Y'+Y'') \in \Omega^{1,0}(X)\oplus\Omega^{0,1}(X)$ to give
\begin{equation}
  \label{eq:Canonical_Clifford_multiplication_Hodge_star_one_forms}
  \rho_\can(a)\phi = \sqrt{2}\left(a''\wedge\phi - \star(a'\wedge \star\phi)\right),
  \quad\text{for all } a \in \Omega^1(X,\RR) \text{ and } \phi \in \Omega^{0,q}(X).
\end{equation}  
Choosing $\phi = (\sigma,\tau) \in \Omega^0(W^+) = \Omega^{0,0}(X)\oplus\Omega^{0,2}(X)$ and $\nu \in \Omega^0(W^-) = \Omega^{0,1}(X)$ and specializing \eqref{eq:Canonical_Clifford_multiplication} to the case of the positive and negative spinors in \eqref{eq:Canonical_spinc_bundles} gives the identities\footnote{Morgan \cite[p. 109]{MorganSWNotes} gives the expression \eqref{eq:Canonical_Clifford_multiplication_positive_spinors} but omits the factor $\sqrt{2}$ in the first term. Kotschick \cite[Fact 2.1]{KotschickSW} gives expressions for Clifford multiplication that agree with \eqref{eq:Canonical_Clifford_multiplication_positive_spinors} and \eqref{eq:Canonical_Clifford_multiplication_negative_spinors}.}
\begin{subequations}
  \begin{align}
    \label{eq:Canonical_Clifford_multiplication_positive_spinors}
    \rho_\can(a)(\sigma,\tau) &= \sqrt{2}\left(a''\wedge\sigma - \star(a'\wedge\star\tau)\right)
                                \in \Omega^{0,1}(X),
  \\
    \label{eq:Canonical_Clifford_multiplication_negative_spinors}
    \rho_\can(a)\nu &= \sqrt{2}\left(- \star(a'\wedge\star\nu) + a''\wedge\nu\right)
                      \in \Omega^{0,0}(X)\oplus\Omega^{0,2}(X).
  \end{align}  
\end{subequations}
To deduce \eqref{eq:Canonical_Clifford_multiplication_positive_spinors} from \eqref{eq:Canonical_Clifford_multiplication_Hodge_star_one_forms}, observe that for $\sigma \in \Omega^{0,0}(X)$ we have  $\star\sigma \in \Omega^{2,2}(X)$ and $a'\wedge\star\sigma \in \Omega^{3,2}(X)$, that is, $a'\wedge\star\sigma = 0$, and thus by \eqref{eq:Canonical_Clifford_multiplication_Hodge_star_one_forms},
\[
  \rho_\can(a)\sigma = \sqrt{2}a''\wedge\sigma \in \Omega^{0,1}(X).
\]
For $\tau \in \Omega^{0,2}(X)$ we have $a''\wedge\tau \in \Omega^{0,3}(X)$, and so $a''\wedge\tau=0$, while $\star\tau \in \Omega^{0,2}(X)$ and $a' \wedge \star\tau \in \Omega^{1,2}(X)$ and $\star(a' \wedge \star\tau)  \in \Omega^{0,1}(X)$, and thus by \eqref{eq:Canonical_Clifford_multiplication_Hodge_star_one_forms},
\[
  \rho_\can(a)\tau = - \sqrt{2}\star(a' \wedge \star\tau)  \in \Omega^{0,1}(X),
\]
which verifies \eqref{eq:Canonical_Clifford_multiplication_positive_spinors}. To deduce \eqref{eq:Canonical_Clifford_multiplication_negative_spinors} from \eqref{eq:Canonical_Clifford_multiplication}, observe that for $\nu \in \Omega^{0,1}(X)$ we have $\star\nu \in \Omega^{1,2}(X)$ by Huybrechts \cite[Lemma 1.2.24, p. 33]{Huybrechts_2005} and $a'\wedge\star\nu \in \Omega^{2,2}(X)$ and $\star(a'\wedge\star\nu) \in \Omega^{0,0}(X)$, so that
\[
  \rho_\can(a)\nu = \sqrt{2}\left(a''\wedge\nu - \star(a'\wedge\star\nu)\right)
  \in \Omega^{0,0}(X)\oplus\Omega^{0,2}(X),
\]
as claimed.

We refer to Salamon \cite[Lemma 4.52, p. 141]{SalamonSWBook} for a verification that \eqref{eq:Canonical_Clifford_multiplication} satisfies the axioms of a \spinc structure \cite[Definition 4.32, p. 125]{SalamonSWBook}, namely that
\begin{equation}
  \label{eq:Salamon_4-18}
  \rho(a)^\dagger = -\rho(a) \in \End(W) \quad\text{and}\quad \rho(a)^\dagger\rho(a) = g(a,a)\id_W,
\end{equation}
for all $a \in \Omega^1(X,\RR)$. If $E$ is a smooth, complex vector bundle over $(X,g,J)$, we obtain a \spinu structure  over $(X,g,J)$ corresponding to a Clifford multiplication bundle map
\[
  \rho:T^*X \to \Hom(W^+, W^-)
\]
via the bundle map
\[
  \rho:T^*X \to \Hom(W^+\otimes E, W^-\otimes E)
\]
defined on elementary tensors by
\begin{multline*}
  \rho(a)(\phi\otimes s) := (\rho(a)\phi)\otimes s \in \Omega^0(W^-\otimes E),
  \\
  \text{for all } a \in \Omega^1(X,\RR), \phi \in \Omega^0(W^+), \text{ and } s \in \Omega^0(E).
\end{multline*}
See Feehan and Leness \cite[Definition 2.2, p. 64, and Lemma 2.3, p. 64]{FL2a} for a more invariant definition of \spinu structure.

%PF10-10-2024 Copied from our monograph. Upgrade to $2n$-manifolds
\section{Dirac operators over almost Hermitian four-manifolds}
\label{sec:Dirac_operator_almost_Hermitian}
We now describe the Dirac operator on the canonical \spinc structure $(\rho_\can,W_\can)$
%PF10-10-2024 Add updated ref
of Definition \ref{defn:Canonical_spinc_bundles}
on a closed, smooth almost Hermitian four-manifold $(X,g,J)$. From Donaldson \cite[Equation (15) and p. 56]{DonSW} and Gauduchon \cite[Equations (3.3.4), (3.7.1), and (3.7.2)]{Gauduchon_1997}, the Dirac operator $D:C^\infty(W_\can)\to C^\infty(W_\can)$ defined by the canonical \spinc structure and \emph{Chern connection} \cite[p. 273 and Equation (3.6.1)]{Gauduchon_1997} (see Remark \ref{rmk:Chern_connection}) on the Hermitian line bundle $L = K_X^*$
%PF10-10-2024 Add updated ref
(the anti-canonical line bundle in \eqref{eq:DefineAntiCanonicalLineBundle} and dual of the canonical line bundle $K_X$ in \eqref{eq:DefineCanonicalLineBundle}),
is related to the normalized Dolbeault operator $\sqrt{2}(\partial + \partial^*)$ acting on the spinor bundle $W_\can = W_\can^+\oplus W_\can^-$, where as in
\eqref{eq:Canonical_spinc_bundles},
\[ 
  W_\can^+ = \Lambda^{0,0}(X)\oplus \Lambda^{0,2}(X) \quad\text{and}\quad W_\can^- = \Lambda^{0,1}(X),
\]
by the identity
\begin{equation}
\label{eq:Gauduchon_3-7-2}
D\Phi = \sqrt{2}(\bar\partial + \bar\partial^*)\Phi + \frac{1}{4}\lambda\cdot(\Phi^+ - \Phi^-),
\quad\text{for all } \Phi \in C^\infty(W_\can),
\end{equation}
where $\lambda \in \Omega^1(X,\RR)$ is the \emph{Lee form} defined by $(g,J)$, so \cite[Equation (1.1.4)]{Gauduchon_1997}
\begin{equation}
\label{eq:Gauduchon_1-1-4}
  \lambda = \Lambda_\omega(d\omega),
\end{equation}
and $\omega(\cdot,\cdot) = g(J\cdot,\cdot)$ from Gauduchon \cite[p. 259]{Gauduchon_1997}, and ``$\cdot$'' denotes Clifford multiplication \cite[Equation (3.1.4) and p. 276]{Gauduchon_1997}. Note that $d\omega = \lambda\wedge\omega$ by \cite[p. 259]{Gauduchon_1997} when $X$ has real dimension four. When $(X,g,J)$ is almost K\"ahler, so $d\omega=0$, then the identity \eqref{eq:Gauduchon_3-7-2} reduces to the identity $D\Phi = \sqrt{2}(\bar\partial + \bar\partial^*)\Phi$ stated by Donaldson in \cite[Equation (15)]{DonSW}.

When $A$ is a unitary connection on a Hermitian vector bundle $E$ over $X$ and the Levi-Civita connection $\nabla$ on $TX$ is replaced by the connection $\nabla_A$ on $TX\otimes E$, then Gauduchon's proof of \eqref{eq:Gauduchon_3-7-2} now \mutatis yields
\begin{equation}
\label{eq:Gauduchon_3-7-2_auxiliary_Hermitian_bundle_E}
  D_A\Phi = \sqrt{2}(\bar\partial_A + \bar\partial_A^*)\Phi + \frac{1}{4}\lambda\cdot(\Phi^+ - \Phi^-),
\end{equation}
for all $\Phi \in C^\infty(W_\can\otimes E)$.

\begin{rmk}[Chern connection]
\label{rmk:Chern_connection}
If $(E,H)$ is a \emph{Hermitian} vector bundle with a holomorphic structure $\bar\partial_E$ over a complex manifold $X$, we may choose $A$ to be the \emph{Chern connection} on $E$ --- the \emph{unique} unitary connection on $E$ defined by the holomorphic structure $\bar\partial_E$ and Hermitian metric $h$ (see Kobayashi \cite[Proposition 1.4.9, p. 11]{Kobayashi_differential_geometry_complex_vector_bundles}).
\qed\end{rmk}

\section{Non-Abelian monopoles over almost Hermitian manifolds}
\label{sec:MonopolesOverAlmostHermitian}
We follow Dowker \cite[Chapter 1]{DowkerThesis} and L\"ubke and Teleman \cite[Section 6.3]{Lubke_Teleman_2006}.
Assume now that $(X,g,J)$ is a closed, smooth, almost Hermitian four-manifold. Recall that the spin bundles associated to the \emph{canonical} \spinc structure $\fs_\can = (\rho_\can, W_\can)$ on $X$ are (see \eqref{eq:Canonical_spinc_bundles} for $W_\can$ and \eqref{eq:Canonical_Clifford_multiplication} for $\rho_\can$)
\[
  W^+ := \Lambda^{0,0}(X) \oplus \Lambda^{0,2}(X), \quad W^- := \Lambda^{0,1}(X), \quad W = W^+\oplus W^-,
\]
where we omit the ``can'' subscripts here for the sake of notational simplicity. The direct sum decomposition of $W^+$ is orthogonal with respect to the Hermitian metric on $W$. Consequently,
\[
  W^+\otimes E = E \oplus \Lambda^{0,2}(E), \quad W^-\otimes E = \Lambda^{0,1}(E), \quad W\otimes E = (W^+\otimes E) \oplus (W^-\otimes E),
\]
where $\Lambda^{0,2}(E) := \Lambda^{0,2}(X)\otimes E$ and $\Lambda^{0,1}(E)\otimes E := \Lambda^{0,1}(X)\otimes E$ and the preceding direct sum decompositions are orthogonal with respect to the Hermitian metric on $W\otimes E$. (Recall that we write $\Lambda^{0,0}(X) = X\times\CC$ and thus $\Lambda^{0,0}(E) = E$.)

\begin{rmk}[Restriction to the canonical \spinc structure]
There is no loss of generality in restricting to the canonical \spinc structure $(\rho,W)$ here since we shall be considering \spinu structures $(\rho,W,E)$ and explicitly varying $(\rho,W)$ by tensoring it with complex Hermitian line bundles $L$.
\qed\end{rmk}

Thus we can write $\Phi \in \Omega^0(W^+\otimes E) = \Omega^0(E)\oplus \Omega^{0,2}(E)$ as $\Phi=(\varphi, \psi)$, where $\varphi\in\Omega^0(E)$ and $\psi\in\Omega^{0,2}(E)$. When $(X,g,J)$ is almost K\"ahler, then the Dirac operator can be expressed in the form (see Donaldson \cite[p. 59]{DonSW}, Kotschick \cite[Section 2.1]{KotschickSW} and \eqref{eq:Gauduchon_3-7-2_auxiliary_Hermitian_bundle_E})
\[
  D_A\Phi = \sqrt 2\left(\bar{\partial}_A\varphi + \bar{\partial}_A^*\psi\right) \in \Omega^{0,1}(E)\oplus \Omega^{1,0}(E) = \Omega^{1,1}(E) = \Omega^0(W^-\otimes E).
\]
See Section \ref{sec:Dirac_operator_almost_Hermitian}, based on \cite{Gauduchon_1997}, for the modification of this expression for the Dirac operator when $(X,g,J,\omega)$ is only almost Hermitian.

From 
%TL11-28-2025: Updated
\cite[Lemma 8.3.7]{Feehan_Leness_introduction_virtual_morse_theory_so3_monopoles},
when $(X,g,J)$ is almost K\"ahler, the non-Abelian monopole equations \eqref{eq:SO(3)_monopole_equations} take the form (compare Bradlow and Garc\'ia--Prada \cite[Equation (22)]{BradlowGP},  Dowker \cite[p. 10]{DowkerThesis}, Labastida and Mari\~no \cite[Equation (3.1)]{Labastida_Marino_1995nam4m}, or Okonek and Teleman \cite[Proposition 2.6]{OTVortex})
\begin{subequations}
\label{eq:SO(3)_monopole_equations_almost_Kaehler}
\begin{align}
  \label{eq:SO(3)_monopole_equations_(1,1)_curvature}
  (\Lambda F_A)_0 &= \frac{i}{2}\left( (\varphi\otimes\varphi^*)_0 - \star(\psi\otimes\psi^*)_0\right),
  \\
  \label{eq:SO(3)_monopole_equations_(0,2)_curvature}
  (F_A^{0,2})_0 &= \frac{1}{2}(\psi\otimes\varphi^*)_0,
  \\
  \label{eq:SO(3)_monopole_equations_(2,0)_curvature}
  (F_A^{2,0})_0 &= -\frac{1}{2}(\varphi\otimes\psi^*)_0,
  \\
  \label{eq:SO(3)_monopole_equations_Dirac_almost_Kaehler}
  \bar{\partial}_A\varphi + \bar{\partial}_A^*\psi &= 0.
\end{align}
\end{subequations}
When $(X,g,J)$ is \emph{almost Hermitian}, \eqref{eq:SO(3)_monopole_equations_(1,1)_curvature} and \eqref{eq:SO(3)_monopole_equations_(0,2)_curvature} remain unchanged from the almost K\"ahler case, but as we describe in Section \ref{sec:Dirac_operator_almost_Hermitian}, the Dirac equation \eqref{eq:SO(3)_monopole_equations_Dirac_almost_Kaehler} is modified by the addition of a zeroth-order term as seen in
 \eqref{eq:Gauduchon_3-7-2_auxiliary_Hermitian_bundle_E}.

The volume form on $(X,g,J)$ is given by $\vol = \frac{1}{2}\omega^2 = \star 1$, where $\omega$ is the fundamental two-form defined by $(g,J)$, so
\[
  \frac{1}{2}\Lambda^2(\psi\otimes\psi^*)_0 = \frac{1}{2}\langle(\psi\otimes\psi^*)_0,\omega^2\rangle_{\Lambda^{2,2}(X)} = \langle(\psi\otimes\psi^*)_0,\vol\rangle = \star(\psi\otimes\psi^*)_0 \in \Omega^0(i\su(E)).
\]
Equation \eqref{eq:SO(3)_monopole_equations_(1,1)_curvature} can thus be equivalently written as
\begin{equation}
  \label{eq:SO(3)_monopole_equations_(1,1)_curvature_omega_forms}
  (F_A)_0 = \frac{i}{2}\left( \omega\wedge(\varphi\otimes\varphi^*)_0 - \Lambda(\psi\otimes\psi^*)_0\right).
\end{equation}
using the isomorphism between the Riemannian vector bundles $\Lambda^0(\su(E))$ and $\Lambda^0(\su(E))\otimes\omega \subset \Lambda^{1,1}(E)$ implied by the real linear map
\begin{equation}
  \label{eq:Huybrechts_definition_1-2-18}
  L = L_\omega:\Omega^0(\su(E)) \ni \xi \mapsto \xi\otimes\omega \in \Omega^{1,1}(\su(E)).
\end{equation}
Dowker \cite[p. 10]{DowkerThesis} writes $\Lambda^2(\psi\otimes\psi^*)_0$ rather than $\star(\psi\otimes\psi^*)_0$ as in Okonek and Teleman \cite[Proposition 2.6]{OTVortex}, so Dowker effectively rescales $\psi$ by a multiple of $\sqrt{2}$ relative to Okonek and Teleman.

We see that equation \eqref{eq:SO(3)_monopole_equations_(2,0)_curvature} is the complex conjugate transpose of \eqref{eq:SO(3)_monopole_equations_(0,2)_curvature} and so can be omitted without loss of generality.

% TL11-11-2025: Added following so we'd have a  referable version of $\Lambda \omega=2$.
%PF11-13-2025 Cut duplication in chapter 1
We further note that the equality $\vol=\frac{1}{2}\omega^2=\star 1$ implies that
\[
\Lambda \omega
=
\star^{-1}\left( L\star \omega\right)
=
\star^{-1}\left( \omega\wedge\star\omega\right)
=
\star^{-1}\left(\omega\wedge\omega\right)
=
2 \star^{-1}\left( \vol\right)
=
2\star^{-1}\left(\star 1\right)
=
2,
\]
which we record as
\begin{equation}
\label{eq:Lambda_on_KahlerForm}
\Lambda\omega=2,
\end{equation} 
for future use.

% TL9-8-2024:Moved this here from cut section on ASD def. complex.  As the Kahler identities were appearing in multiple areas (with similar cut and pastes in all), there is a lot of redundancy here and much pruning and reorganizing needed.
%PF10-10-2024 This subsection is a mess.
\section{K\"ahler identities on almost K\"ahler and almost Hermitian manifolds}
\label{sec:KaehlerIDs}
When $X$ is almost K\"ahler, recall that the standard K\"ahler identities on $\Omega^\bullet(X,\CC)$ (see Cirici and Wilson \cite[Proposition 3.1 (2) and (4)]{Cirici_Wilson_2020_harmonic}) are
\begin{gather}
\label{eq:Kaehler_identity_commutator_L_del-bar_and_del}
[L, \bar\partial] = 0 \quad\text{and}\quad [L, \partial] = 0,
\\
\label{eq:Kaehler_identity_commutator_L_del-bar_star_and_del}
[L, \bar\partial^*] = -i\partial \quad\text{and}\quad [L, \partial^*] = i\bar\partial.
\end{gather}
%PF10-10-2024 Fix
%[TODO Check Wells too. Duplication]
These identities should extend to bundle-valued forms as well to give on $\Omega^\bullet(X,E)$,
\begin{gather}
\label{eq:Kaehler_identity_commutator_L_del-bar_A_and_del_A}
[L, \bar\partial_A] = 0 \quad\text{and}\quad [L, \partial_A] = 0,
\\
\label{eq:Kaehler_identity_commutator_L_del-bar_A_star_and_del_A}
[L, \bar\partial_A^*] = -i\partial_A \quad\text{and}\quad [L, \partial_A^*] = i\bar\partial_A.
\end{gather}
Hence,
\[
  L\partial_A^*\partial_A^*v' = i\bar\partial_A\partial_A^*v' - \partial_A^*L\partial_A^*v'.
\]
The K\"ahler Identities for a K\"ahler manifold are proved in Huybrechts \cite[Proposition 3.1.12]{Huybrechts_2005} and Wells \cite[Corollary 5.4.10]{Wells3} and for an almost K\"ahler manifold by Weil \cite{Weil_introduction_etude_varietes_kahleriennes}. For an almost Hermitian manifold of real dimension four, we would expect simplifications in the proof given that $\omega$ should be self-dual and $a', a''$ are one-forms. However proved, we see that there will be an additional term involving $\partial\omega$ or $\bar\partial\omega$ in each of the two K\"ahler Identities we need. (See Yano \cite{Yano_1965} for an introduction to almost complex manifolds.)

\begin{rmk}[K\"ahler identities on one-forms over an almost Hermitian manifold]
\label{rmk:Kaehler_identities_one-forms_almost_Hermitian_manifold}
Let $(X,g,J)$ be an almost Hermitian manifold of real dimension $2n$ and fundamental two-form $\omega(\cdot,\cdot) = g(J\cdot,\cdot)$. (We follow the convention of Cirici and Wilson in \cite[Section 2]{Cirici_Wilson_2020_almost_hermitian_identities} and observe that this is the reverse of the convention in of Remark \ref{rmk:Kaehler_identities_one-forms_almost_Hermitian_manifold}.) Following the conventions in Cirici and Wilson \cite[Section 2]{Cirici_Wilson_2020_almost_hermitian_identities}, we define the \emph{Hodge star operator}
%PF10-10-2024 Fix
% [TODO Get some more precise references]
\begin{equation}
  \label{eq:Hodge_star_operator}
  \star = \star_g:\Omega^{p,q}(X) \to \Omega^{n-p,n-q}(X)
\end{equation}
by
\[
  \alpha\wedge\star_g\bar\beta = \langle\alpha,\beta\rangle_g\vol_g,
\]
where $\langle\cdot,\cdot\rangle_g$ is the pointwise inner product
%PF10-10-2024 Fix
%[TODO Real? Hermitian?]
on $\wedge^\bullet(T^*X\otimes\CC)$ defined by the Riemannian metric $g$. For any $\eta \in \Omega^{r,s}(X)$, we define the \emph{left multiplication operator}
%PF10-10-2024 Fix
%[TODO Get some more precise references]
\begin{equation}
  \label{eq:Left_multiplication_form_operator}
  L_\eta:\Omega^{p,q}(X) \ni \gamma \mapsto \eta\wedge\gamma \in \Omega^{p+r,q+s}(X)
\end{equation}
and hence the \emph{Lefschetz operator}
\begin{equation}
  \label{eq:Lefschetz_operator}
  L = L_\omega:\Omega^{p,q}(X) \to \Omega^{p+1,q+1}(X)
\end{equation}
upon choosing $\eta=\omega$; its pointwise adjoint with respect to $g$ is denoted by
\begin{equation}
  \label{eq:Lambda_operator}
  \Lambda = L^* = \star^{-1} L\, \star:\Omega^{p,q}(X) \to \Omega^{p-1,q-1}(X).
\end{equation}
The operators $\delta = \partial, \mu, \bar\mu, \bar\partial$ have $L^2$ adjoints $\delta^*$ when $X$ is closed and
\begin{equation}
  \label{eq:L2_adjoint_mu_partial}
  \bar\mu^* = -\star \mu \star \quad\text{and}\quad \bar\partial^* = -\star \partial \star.
\end{equation}
According to Cirici and Wilson (see \cite[Corollary 2.2]{Cirici_Wilson_2020_almost_hermitian_identities} or \cite[Lemma 4.8]{Cirici_Wilson_2021}) for any $\alpha\in \Omega^{0,1}(X)$, one has the following K\"ahler identity due to Ohsawa \cite[Appendix]{Ohsawa_1982}:
\begin{equation}
  \label{eq:Commutator_Lambda_partial_Kaehler_identity_(0,1)_forms}
  [\Lambda,\partial]\alpha = \Lambda\partial\alpha = i\bar\partial^*\alpha + i[\Lambda,\bar\partial^*]L\alpha \in \Omega^{0,0}(X).
\end{equation}
From the proof of their \cite[Corollary 2.3]{Cirici_Wilson_2020_almost_hermitian_identities}, the operator $[\Lambda,\bar\partial^*]$ is zeroth-order; moreover, the term $i[\Lambda,\bar\partial^*]L\alpha$ is identically zero when $\omega$ is closed and $(X,g,J)$ is almost K\"ahler. Indeed, for any  $f\in \Omega^{0,0}(X)$ and $\gamma\in \Omega^{1,2}(X)$ (noting that $L\alpha=\omega\wedge\alpha\in \Omega^{1,2}(X)$) we have
\begin{multline*}
  (\Lambda\bar\partial^*\gamma, f)_{L^2(X)}
  = (L^*\bar\partial^*\gamma)_{L^2(X)}
  = (\gamma, \bar\partial Lf)_{L^2(X)}
  = (\gamma, \bar\partial(f\omega))_{L^2(X)}
  \\
  = (\gamma, (\bar\partial f)\wedge \omega + f\bar\partial\omega)_{L^2(X)}
  = (\gamma, L\bar\partial f)_{L^2(X)} + (\gamma, L_{\bar\partial\omega}f)_{L^2(X)}
  \\
  = (\bar\partial^*L^*\gamma, f)_{L^2(X)} + (L_{\bar\partial\omega}^*\gamma, f)_{L^2(X)}
  = (\bar\partial^*\Lambda\gamma + L_{\bar\partial\omega}^*\gamma,f)_{L^2(X)},
\end{multline*}
and thus
\[
  [\Lambda, \bar\partial^*]\gamma
  = \Lambda\bar\partial^*\gamma - \bar\partial^*\Lambda\gamma
  = L_{\bar\partial\omega}^*\gamma \in \Omega^{0,0}(X),
\]
noting that $\bar\partial\omega \in \Omega^{1,2}(X)$, so $L_{\bar\partial\omega}:\Omega^{p,q}(X)\to \Omega^{p+1,q+2}(X)$ while $L_{\bar\partial\omega}^*:\Omega^{p,q}(X)\to \Omega^{p-1,q-2}(X)$. Hence, upon choosing $\gamma = L\alpha \in \Omega^{1,2}(X)$, we see that
\[
  [\Lambda,\bar\partial^*]L\alpha = L_{\bar\partial\omega}^*L\alpha,
\]
and so the K\"ahler identity \eqref{eq:Commutator_Lambda_partial_Kaehler_identity_(0,1)_forms} can be written in the more transparent form
\begin{equation}
  \label{eq:Commutator_Lambda_partial_Kaehler_identity_(0,1)_forms_zeroth_order_additional_term}
  [\Lambda,\partial]\alpha = \Lambda\partial\alpha = i\bar\partial^*\alpha + iL_{\bar\partial\omega}^*L\alpha \in \Omega^{0,0}(X),
\end{equation}
clearly showing that the additional term obtained when $\omega$ is not $\bar\partial$-closed is zeroth order.
\qed\end{rmk}

When $(X,g,J,\omega)$ is almost K\"ahler, the following \emph{K\"ahler identities} hold (see Huybrechts \cite[Proposition 3.1.12, p. 120]{Huybrechts_2005} and Wells \cite[Chapter V, Corollary 4.10, p. 193]{Wells3} for the complex K\"ahler case and Cirici and Wilson \cite[Proposition 3.1 (4)]{Cirici_Wilson_2020_harmonic}
% PF7-17-2024 Paul to add reference to Weil
%PF7-24-2024 Paul to make this into a lemma
for the almost K\"ahler case),
\begin{subequations}
  \label{eq:Kaehler_identity_commutator_Lambda_del-bar_and_Lambda_del}
  \begin{align}
    \label{eq:Kaehler_identity_commutator_Lambda_del-bar}
    [\Lambda,\bar\partial] &= -i\partial^* \quad\text{and}
    \\
    \label{eq:Kaehler_identity_commutator_Lambda_del}
    [\Lambda,\partial] &= i\bar\partial^* \quad\text{on } \Omega^\bullet(X,\CC).
  \end{align}
\end{subequations}
According to Demailly \cite[Section 7.1]{Demailly_complex_analytic_differential_geometry} (see also Demailly \cite{Demailly_1986}), Donaldson \cite[Proof of Proposition 3]{DonASD} and Donaldson and Kronheimer\footnote{There is a sign error in the second displayed equation of \cite[Section 6.1.3, Proof of Lemma 6.1.7, p. 213]{DK}, as is easy to see by comparison with \cite[Section 6.1.3, Equation (6.1.8), p. 212]{DK}.} \cite[Equation (6.1.8), p. 212]{DK}, and Kobayashi \cite[Section 3.2, p. 62]{Kobayashi_differential_geometry_complex_vector_bundles} the K\"ahler identities \eqref{eq:Kaehler_identity_commutator_Lambda_del-bar_and_Lambda_del} extend to bundle-valued forms to give
\begin{subequations}
\label{eq:Kaehler_identity_commutator_Lambda_del-bar_A_and_Lambda_del_A}
\begin{align}
\label{eq:Kaehler_identity_commutator_Lambda_del-bar_A}
  [\Lambda, \bar\partial_A] &= -i\partial_A^* \quad\text{and}
  \\
  \label{eq:Kaehler_identity_commutator_Lambda_del_A}
  [\Lambda, \partial_A] &= i\bar\partial_A^* \quad\text{on } \Omega^\bullet(E).
\end{align}
\end{subequations}
When $(X,g,J)$ is \emph{almost Hermitian}, the following modified K\"ahler Identities hold:
\begin{equation}
\label{eq:Kaehler_identity_commutator_Lambda_del-bar_and_Lambda_del_almost_Hermitian}
[\Lambda,\bar\partial] = -i\partial^* - i\Lambda_{\partial\omega}L \quad\text{on } \Omega^{1,0}(X) \quad\text{and}\quad
[\Lambda,\partial] = i\bar\partial^* + i\Lambda_{\bar\partial\omega}L
\quad\text{on } \Omega^{0,1}(X).
\end{equation}
The identity for $[\Lambda,\partial]$ is given by \eqref{eq:Commutator_Lambda_partial_Kaehler_identity_(0,1)_forms_zeroth_order_additional_term}
%PF10-10-2024 Fix
%[TODO Forward ref],
while complex conjugation
%PF10-10-2024 Fix
%[TODO Ref]
yields the identity for $[\Lambda,\bar\partial]$ acting on $\alpha \in \Omega^{1,0}(X)$:
%PF10-10-2024 Fix
%[TODO Add detail]
\begin{align*}
  [\Lambda,\bar\partial]\alpha &= \Lambda\bar\partial\alpha = \Lambda\overline{\partial\bar\alpha} = \overline{\Lambda\partial\bar\alpha}
  \\
                               &= \overline{i\bar\partial^*\bar\alpha + i\Lambda_{\bar\partial\omega}L\bar\alpha} \quad\text{(by \eqref{eq:Commutator_Lambda_partial_Kaehler_identity_(0,1)_forms_zeroth_order_additional_term})}
  \\
                               &= -i\overline{\bar\partial^*\bar\alpha} - i\overline{\Lambda_{\bar\partial\omega}L\bar\alpha}
  \\
  &= -i\partial^*\alpha - i\Lambda_{\partial\omega}L\alpha.
\end{align*}
The preceding identities extend to bundle-valued forms as well. We shall prove this first for $\Omega^\bullet(E)$ and then $\Omega^\bullet(\gl(E))$. Let $E \restriction U \cong U\times\CC^r$ be a local trivialization over an open subset $U \subset X$, write $\gamma = (\gamma_1,\ldots, \gamma_r) \in \Omega^{p,q}(U,\CC^r) = \Omega^{p,q}(U,\CC)^r$ with respect to this trivialization, and write
\[
  \bar\partial_A\gamma = \bar\partial\gamma + a''\wedge\gamma
  \quad\text{and}\quad
  \partial_A\gamma = \partial\gamma + a'\wedge\gamma \quad\text{for } \gamma \in \Omega^{p,q}(U,\CC^r),
\]
where $a'' \in \Omega^{0,1}(U,\gl(r,\CC))$ and $a' = -(a'')^\intercal = (a'')^\dagger \in \Omega^{1,0}(U,\gl(r,\CC))$. Hence,
%PF10-10-2024 Fix
%[TODO Ref]
\[
  \bar\partial_A^*\gamma = \bar\partial^*\gamma - \star (a''\wedge\star\gamma)
  \quad\text{and}\quad
  \partial_A^*\gamma = \partial^*\gamma - \star (a'\wedge\star\gamma) \quad\text{for } \gamma \in \Omega^{p,q}(U,\CC^r).
\]
We now restrict $U$ to be a star-shaped open subset with respect to a point $x_0\in U$ (and the Riemannian metric $g$ on $X$) and let $E \restriction U \cong U\times\CC^r$ be a trivialization with corresponding section $\sigma:U\to E$ such that $\sigma^*A$ is in radial gauge relative to the point $x_0$ and, in particular, has $(\sigma^*A)(x_0) = 0$ (see Uhlenbeck \cite[Lemma 2.1]{UhlRem}). (Generalized K\"ahler identities for $\bar\partial_A$ over a complex K\"ahler manifold are derived by Kobayashi \cite[p. 62]{Kobayashi_differential_geometry_complex_vector_bundles}.) Consequently,
\[
  \bar\partial_A\gamma(x_0) = \bar\partial\gamma(x_0), \quad
  \partial_A\gamma(x_0) = \partial\gamma(x_0), \quad
  \bar\partial_A^*\gamma(x_0) = \bar\partial^*\gamma(x_0), \quad\text{and}\quad
  \partial_A^*\gamma(x_0) = \partial^*\gamma(x_0).
\]
For $\beta \in \Omega^{1,0}(X,\CC^r)$, we therefore obtain
\begin{align*}
  [\Lambda,\bar\partial_A]\beta(x_0) &= \Lambda\bar\partial_A\beta(x_0) = \Lambda\bar\partial\beta
   = [\Lambda,\bar\partial]\beta
  \\
  &= -i\partial^*\beta(x_0) - i\Lambda_{\partial\omega}L\beta(x_0) \quad\text{(by \eqref{eq:Kaehler_identity_commutator_Lambda_del-bar_and_Lambda_del_almost_Hermitian})}
  \\
  &= -i\partial_A^*\beta(x_0) - i\Lambda_{\partial\omega}L\beta(x_0).
\end{align*}
For $\alpha \in \Omega^{0,1}(X,\CC^r)$, we therefore obtain
\begin{align*}
  [\Lambda,\partial_A]\alpha(x_0) &= \Lambda\partial_A\alpha(x_0) = \Lambda\partial\alpha
   = [\Lambda,\partial]\alpha
  \\
  &= i\bar\partial^*\alpha(x_0) + i\Lambda_{\bar\partial\omega}L\alpha(x_0) \quad\text{(by \eqref{eq:Kaehler_identity_commutator_Lambda_del-bar_and_Lambda_del_almost_Hermitian})}
  \\
  &= i\bar\partial_A^*\alpha(x_0) + i\Lambda_{\bar\partial\omega}L\alpha(x_0).
\end{align*}
Since $x_0 \in X$ is arbitrary, we thus obtain the pointwise identities on $\Omega^1(E)$,
\begin{equation}
\label{eq:Kaehler_identity_commutator_Lambda_del-bar_A_and_Lambda_del_A_almost_Hermitian_E}
\begin{aligned}
{}[\Lambda, \bar\partial_A] &= -i\partial_A^* - i\Lambda_{\partial\omega}L \quad\text{on } \Omega^{1,0}(E),
\\
[\Lambda, \partial_A] &= i\bar\partial_A^* + i\Lambda_{\bar\partial\omega}L
\quad\text{on } \Omega^{0,1}(E).
\end{aligned}
\end{equation}
By applying the same argument for the induced connection on $\gl(E)$ in place of $E$, we obtain the pointwise identities on $\Omega^1(\gl(E))$,
\begin{equation}
\label{eq:Kaehler_identity_commutator_Lambda_del-bar_A_and_Lambda_del_A_almost_Hermitian_glE}
\begin{aligned}
{}[\Lambda, \bar\partial_A] &= -i\partial_A^* - i\Lambda_{\partial\omega}L \quad\text{on } \Omega^{1,0}(\gl(E)),
\\
[\Lambda, \partial_A] &= i\bar\partial_A^* + i\Lambda_{\bar\partial\omega}L
\quad\text{on } \Omega^{0,1}(\gl(E)),
\end{aligned}
\end{equation}
and similarly for the subbundle $\fsl(E) \subset \gl(E) = \fsl(E)\oplus \CC\,\id_E$. But
\[
  \Lambda(\bar\partial_Aa' + \partial_Aa'')
  = [\Lambda,\bar\partial_A]a' + [\Lambda,\partial_A]a'' \quad\text{(as $\Lambda=0$ on $\Omega^{1,0}(X)$ and $\Omega^{0,1}(X)$)},
\]
and so by applying \eqref{eq:Kaehler_identity_commutator_Lambda_del-bar_A_and_Lambda_del_A_almost_Hermitian_glE},
\begin{equation}
  \label{eq:Lambda_dbar_A_a'_plus_del_A_a''_almost_Hermitian}
  \Lambda(\bar\partial_Aa' + \partial_Aa'')
  =
  - i\partial_A^*a' - i\Lambda_{\partial\omega}La' + i\bar\partial_A^*a'' + i\Lambda_{\bar\partial\omega}La''.
\end{equation}

% TL7-8-2024:Moved this here --will need to insert references to it in H^i computations
\section{On real and complex inner product spaces}
\label{sec:RealAndComplexInnerProductSpaces}
We digress to discuss further properties of real and complex Hilbert spaces. Let $\fH$ be a complex Hilbert space with real form $\fH$, so $\fH = \sH\otimes_\RR\CC$. We write $h=h_1+ih_2\in\fH$, where $h_1,h_2\in\sH$. Then
\[
  \langle h,h'\rangle_\fH = \left(\langle h_1,h_1'\rangle_\sH + \langle h_2,h_2'\rangle_\sH \right)
  -  i\left(\langle h_1,h_2'\rangle_\sH + \langle h_1',h_2\rangle_\sH \right),
\]
so that
\[
  \Real\,\langle h,h'\rangle_\fH = \langle h_1,h_1'\rangle_\sH + \langle h_2,h_2'\rangle_\sH.
\]
Given $h \in \fH$, suppose that $\Real\,\langle h,h'\rangle_\fH = 0$ for all $h' \in \fH$. Then taking $h' = h_1' + i0$ and $h' = 0 + ih_2'$, respectively, we see that
\[
  \langle h_1,h_1'\rangle_\sH = 0 = \langle h_2,h_2'\rangle_\sH, \quad\text{for all } h_1', h_2' \in \sH,
\]
and so we obtain $h_1=h_2=0$. Hence, for any given $h\in\fH$, the following equations are equivalent:
\begin{equation}
  \label{eq:Equivalence_real_complex_inner_product_equations}
  \Real\,\langle h,h'\rangle_\fH = 0 \quad\text{for all } h' \in \fH
  \iff
  \langle h,h'\rangle_\fH = 0 \quad\text{for all } h' \in \fH.
\end{equation}
Alternatively, just observe that
\[
  \Real\,\langle h,ih'\rangle_\fH = \Real\left(-i\langle h,h'\rangle_\fH\right) = \Imag\,\langle h,h'\rangle_\fH
\]
and the conclusion \eqref{eq:Equivalence_real_complex_inner_product_equations} again follows. We shall also use \eqref{eq:Equivalence_real_complex_inner_product_equations} in the equivalent form
\begin{equation}
  \label{eq:Equivalence_imaginary_complex_inner_product_equations}
  \Imag\,\langle h,h'\rangle_\fH = 0 \quad\text{for all } h' \in \fH
  \iff
  \langle h,h'\rangle_\fH = 0 \quad\text{for all } h' \in \fH,
\end{equation}
obtained from the preceding relation between the real and imaginary parts of the complex inner product $\langle h,h'\rangle_\fH$.

For any complex Hilbert spaces $\fH$ and $\fK$ and operator $M\in\Hom(\fH,\fK)$, its \emph{Toeplitz Decomposition} \cite[p. 7]{Horn_Johnson_matrix_analysis_2013} is given by
\[
  M = \pi_HM + \pi_SM = M^H + M^S,
\]
where $\pi_HM = \frac{1}{2}(M + M^\dagger)$ and $\pi_SM = \frac{1}{2}(M - M^\dagger) = -i\pi_H(iM)$, since $\pi_H(iM) =  \frac{1}{2}(iM - iM^\dagger) = i\pi_SM$.
This concludes our digression on properties of real and complex Hilbert spaces.

%TL7-9-2024: Moved this from old H^2 computations.  Should be integrated with rest of subsection
\begin{rmk}[Real transpose operators, complex conjugate operators, and Hermitian adjoint operators on Hilbert spaces]
\label{rmk:Complex_conjugate_transpose_adjoint_operator_Hilbert_space}
For any complex Hilbert spaces $\fH$ and $\fK$ and operator $M\in\Hom(\fH,\fK)$, we recall that the \emph{Hermitian adjoint operator} $M^\dagger\in\Hom(\fK,\fH)$ is defined by the relation (Kadison and Ringrose \cite[Theorem 2.4.2]{KadisonRingrose1})
\begin{equation}
  \label{eq:Hilbert_space_adjoint}
  \langle M^\dagger k, h\rangle_\fH := \langle k, Mh \rangle_\fK, \quad\text{for all } h \in \fH, k \in \fK.
\end{equation}
With respect to almost complex structures $J_\fH$ on $\fH$ and $J_\fK$ on $\fK$ and corresponding complex conjugations $C_\fH(h) = \bar h$ and $C_\fK(k) = \bar k$, we define the \emph{complex conjugate operator} $\bar M\in\Hom(\fK,\fH)$ by the relation
\begin{equation}
  \label{eq:Hilbert_space_complex_conjugate_operator}
  \langle \bar M h, k\rangle_\fK := \langle \bar k, M\bar h \rangle_\fK, \quad\text{for all } h \in \fH, k \in \fK.
\end{equation}
Note that for all $h \in \fH$ and $k \in \fK$,
\[
  \langle \bar M\bar h, \bar k \rangle_\fK
  = \langle \overline{M h}, \bar k \rangle_\fK
  = \overline{\langle M h, k \rangle}_\fK
  = \langle k, M h\rangle_\fK.
\]
We define the \emph{real transpose} $M^\intercal \in \Hom(\fK,\fH)$ by the relation
\begin{equation}
  \label{eq:Hilbert_space_transpose}
\langle M^\intercal k, h\rangle_\fH := \langle M\bar h, \bar k \rangle_\fK, \quad\text{for all } h \in \fH, k \in \fK.
\end{equation}
Note that for all $h \in \fH$ and $k \in \fK$,
\[
  \langle \bar M^\intercal k, h\rangle_\fH
  = \langle \bar M\bar h, \bar k\rangle_\fK
  = \langle \overline{M h}, \bar k \rangle_\fK
  = \overline{\langle M h, k \rangle}_\fK
  = \langle k, Mh\rangle_\fK
  = \langle M^\dagger k, h\rangle_\fK,
\]
and thus
\[
  M^\dagger = \bar M^\intercal,
\]
generalizing the usual definition of Hermitian adjoint of a complex finite-dimensional matrix (Lancaster and Tismenetsky \cite[Section 1.5]{Lancaster_Tismenetsky}).
\qed\end{rmk}

\chapter[Elliptic deformation complex for the non-Abelian monopole equations]{Elliptic deformation complex and harmonic spaces for the non-Abelian monopole equations over almost K\"ahler four-manifolds}
\label{chap:Elliptic_deformation_complex_moduli_space_SO(3)_monopoles_over_almost_Hermitian_four-manifold}
In this chapter, we develop a Kuranishi model for points in the moduli space of solutions to the non-Abelian monopole equations with a regularized Taubes perturbation on an almost K\"ahler four-manifold.
This model is given by a map between the  bounded-eigenvalue eigenspaces appearing in Theorem \ref{mainthm:AH_structure_bounded_evalue_spaces_non-Abelian_monopoles_symp_4-mflds} and Corollary \ref{maincor:Almost_Hermitian_structure_moduli_space_non-Abelian_monopoles_symplectic_4-manifolds}.
Because  we will use this model in Theorem \ref{mainthm:IdentifyCriticalPoints} to  identify critical points of the Hitchin function and in Theorem \ref{mainthm:MorseIndexAtReduciblesOnAlmostKahler} to compute the virtual Morse--Bott index, we will show in this Chapter and in Chapter \ref{chap:Calculation_virtual_Morse-Bott_indices_via_Atiyah-Singer_index_theorem}  that the bounded-eigenvalue eigenspaces admit the $S^1$-equivariant almost complex structures needed for the proofs of these theorems.

We begin, in Section \ref{sec:Local_Kuranishi_model_non-Abelian_monopole_equations_Taubes_perturbation}, by developing the local Kuranishi models defined by the non-Abelian monopole equations \eqref{eq:SO(3)_monopole_equations_almost_Hermitian_perturbed_intro_regular}
%PF12-18-2025 Edited. Any reason for not just referring to the group of eqns?
% \eqref{eq:SO(3)_monopole_equations_(0,2)_curvature_intro},
% \eqref{eq:SO(3)_monopole_equations_Dirac_almost_Hermitian_intro}, and
% \eqref{eq:SO(3)_monopole_equations_(1,1)_curvature_perturbed_intro_regular}
with a regularized Taubes perturbation.
In particular, in Theorem \ref{thm:Kuranishi_model_defined_by_Fredholm_map_Banach_spaces}, 
we give a general result identifying the properties of subspaces of the domain and codomain of a nonlinear Fredholm map necessary to construct a Kuranishi model as a map between these subspaces.
The first requirement of the subspaces is that they contain the kernel and cokernel of the linearization of the nonlinear Fredholm map.  We identify this kernel and cokernel as spaces of harmonic sections of
a deformation operator for the non-Abelian monopole equations
which we introduce in Section \ref{sec:DefTheory_of_nAM_equations}.
In Chapters \ref{chap:Construction_circle-invariant_non-degenerate_two-form_I} and \ref{chap:Construction_circle-invariant_non-degenerate_two-form_II}, we will construct the almost complex structures needed for the applications described on the spaces of harmonic sections of an equivalent deformation operator in Section \ref{sec:ApproxComplexDefOperator}. We prove that the spaces of harmonic sections of these deformation operators are isomorphic to each other in Sections \ref{sec:H0Isom} and \ref{sec:H1_and_H2_Isom}. In Section \ref{sec:HarmonicTheoryForPerturbedEquations}, we extend these constructions to the perturbed non-Abelian monopole equations.
%TL12-16-2025: REmoved this section
%Finally, we apply  Theorem \ref{thm:Kuranishi_model_defined_by_Fredholm_map_Banach_spaces} in Section \ref{sec:KuranishiModelsFromLowEigenvalueEigespaces} to the harmonic spaces of Section \ref{sec:HarmonicTheoryForPerturbedEquations} and thus construct
% a Kuranishi model for the moduli space of solutions to the perturbed non-Abelian monopoles based on these bounded-eigenvalue eigenspaces.
%PF12-18-2025 Ok

\section[Local Kuranishi models defined by the non-Abelian monopole equations]{Local Kuranishi models defined by the non-Abelian monopole equations with a regularized Taubes perturbation}
\label{sec:Local_Kuranishi_model_non-Abelian_monopole_equations_Taubes_perturbation}
The purpose of this section is to describe the local Kuranishi models defined by the non-Abelian monopole equations \eqref{eq:SO(3)_monopole_equations_(0,2)_curvature_intro},
\eqref{eq:SO(3)_monopole_equations_Dirac_almost_Hermitian_intro}, and
\eqref{eq:SO(3)_monopole_equations_(1,1)_curvature_perturbed_intro_regular} with a perturbation term \eqref{eq:Definition_wp_intro_regular} that is an analogue of a perturbation \cite[Section 1(d), Equation (1.17), p. 851]{TauSWGromov} employed by Taubes for the Seiberg--Witten monopole equations \cite[Section 1(d), Equation (1.18), p. 851]{TauSWGromov} over symplectic four-manifolds.

\subsection{Non-Abelian monopole equations over almost Hermitian four-manifolds}
\label{subsec:Local_Kuranishi_model_non-Abelian_monopole_equations_Taubes_perturbation}
In this subsection, we recall the structure of the system  \eqref{eq:SO(3)_monopole_equations_almost_Hermitian_intro} of unperturbed non-Abelian monopole equations over almost Hermitian four-manifolds and some of the motivation for the perturbation term \eqref{eq:Definition_wp_intro_regular} used to define the system  \eqref{eq:SO(3)_monopole_equations_almost_Hermitian_perturbed_intro_regular} of non-Abelian monopole equations with a regularized Taubes perturbation.

\begin{lem}[Unperturbed non-Abelian monopole equations over an almost Hermitian four-manifold]
\label{lem:SO3_monopole_equations_almost_Kaehler_manifold}
(See Feehan and Leness \cite[Lemma 8.3.7]{Feehan_Leness_introduction_virtual_morse_theory_so3_monopoles} and Okonek and Teleman \cite[Proposition 2.6, p. 900]{OTVortex}.)
Let $(X,g,J,\omega)$ be a smooth, almost Hermitian manifold of real dimension four, let $(\rho_\can,W_\can)$ be the canonical spin${}^c$ structure over $X$ (see equations \eqref{eq:Canonical_spinc_bundles} for $W_\can$ and \eqref{eq:Canonical_Clifford_multiplication} for $\rho_\can$) and let $(E,H)$ be a smooth, Hermitian vector bundle with $\rank_\CC E \geq 2$ over $X$. If $A$ is a smooth unitary connection on $E$ and $\Phi$ is a smooth section of $W_\can^+\otimes E$, then $(A,\Phi)$ is solution to the unperturbed non-Abelian monopole equations \eqref{eq:SO(3)_monopole_equations} with $\Phi=(\varphi,\psi) \in \Omega^0(E)\oplus\Omega^{0,2}(E)$ if and only if $(A,\varphi,\psi)$ obeys the system \eqref{eq:SO(3)_monopole_equations_almost_Hermitian_intro}.
\end{lem}

We use a perturbation of the system \eqref{eq:SO(3)_monopole_equations_almost_Hermitian_intro} of non-Abelian monopole equations and a suitable gauge-equivariant map (see \eqref{eq:Definition_wp_intro_regular}, for example),
\[
  \wp_\gamma: \sA(E,H,A_d) \times \Omega^0(E) \times \Omega^{0,2}(E)
  \ni (A,\varphi,\psi) \mapsto \wp_\gamma(\psi) \in \Omega^0(i\su(E)),
\]
by adding, for a constant $r>0$, a term $-ir\wp_\gamma(\psi)/4$ to the right-hand side of equation \eqref{eq:SO(3)_monopole_equations_(1,1)_curvature_perturbed_intro}, thus giving the system \eqref{eq:SO(3)_monopole_equations_almost_Hermitian_perturbed_intro_regular}.

Equation \eqref{eq:SO(3)_monopole_equations_(1,1)_curvature_perturbed_intro} is equivalent to equation \eqref{eq:SO(3)_monopole_equations_(1,1)_curvature_perturbed_omega_intro}, namely
\[
  (F_A^\omega)_0 = \frac{i}{4}(\varphi\otimes\varphi^*)_0\omega - \frac{i}{4}\star(\psi\otimes\psi^*)_0\omega
     - \frac{ir}{8}\wp_\gamma(\psi)\,\omega,
\]
where $F_A^\omega$ is the projection of $F_A$ onto its image in the factor $\Omega^0(\su(E))\cdot\omega$ of the orthogonal decomposition,
\[
  \Omega^+(\fsl(E)) = \Omega^{2,0}(\fsl(E))\oplus \Omega^0(\fsl(E))\cdot\omega \oplus \Omega^{0,2}(\fsl(E)),
\]
noting that $\fsl(E) \cong i\underline{\RR}\oplus\su(E)$, where $\underline{\RR} = X\times\RR$. To see this, observe that
% PF5-9-2024 Using section 7.1 in FL and Wirtinger's Theorem.
%TL11-25-2025: Duplicates \eqref{eq:Lambda_on_KahlerForm}
\[
  \Lambda_\omega \omega = (\star^{-1}\circ L_\omega\circ \star)\omega = \star^{-1}(\omega\wedge\star\omega) = \star^{-1}\omega^2 = 2\star^{-1}(\omega^2/2) = 2\star^{-1}d\vol = 2\star^{-1}(\star 1) = 2.
\]
Hence, applying $\Lambda_\omega$ to \eqref{eq:SO(3)_monopole_equations_(1,1)_curvature_perturbed_omega_intro} and noting that $\Lambda_\omega F_A^\omega = \Lambda_\omega F_A$ yields \eqref{eq:SO(3)_monopole_equations_(1,1)_curvature_perturbed_intro}.

We now motivate the choice \eqref{eq:Definition_wp_intro_regular} of perturbation term $\wp_\gamma(\psi)$. We begin by discussing the simpler (but singular) choice of perturbation term $\wp(\psi)$ in \eqref{eq:Definition_wp_intro} in case $\rank_\CC E = 2$. Recall that the vector space $i\su(2)$ of $2\times 2$ traceless Hermitian matrices has a basis given by the Pauli matrices,
\[
  \sigma_1 = \begin{pmatrix} 0 & 1 \\ 1 & 0 \end{pmatrix},
  \quad \sigma_2 = \begin{pmatrix} 0 & -i \\ i & 0 \end{pmatrix},
  \quad\text{and}\quad \sigma_3 =  \begin{pmatrix} 1 & 0 \\ 0 & -1 \end{pmatrix},
\]
% COMMENT https://en.wikipedia.org/wiki/Pauli_matrices
with eigenvalues $\pm 1$ and the vector space $i\fu(2)$ of $2\times 2$ Hermitian matrices has basis given by $\id, \sigma_1, \sigma_2, \sigma_3$. Suppose $v \in \Omega^0(i\fu(E))$ has $\tr v = 1$ and $u = v - \frac{1}{2}\id_E$. Then
\[
  \langle v\psi,\psi\rangle_{\Lambda^{0,2}(E)}
  =
  \langle u\psi,\psi\rangle_{\Lambda^{0,2}(E)} - \frac{1}{2}|\psi|_{\Lambda^{0,2}(E)}^2,
\]
but depending on whether $\psi$ is an eigenvector for the eigenvalue $+1$ or $-1$, this cannot have a positive lower bound in general. If we can choose $\wp(\psi) \in L^\infty(i\su(E))$ such that $\psi$ is an eigenvector of $i\wp(\psi)$ for the eigenvalue $+1$, respectively, then
\[
  \langle \wp(\psi)\psi,\psi\rangle_{\Lambda^{0,2}(E)} = |\psi|_{\Lambda^{0,2}(E)}^2.
\]
Therefore, given $\psi \in W^{2,p}(\Lambda^{0,2}(E)$, we recall that $W^{2,p}(X) \subset C(X)$ by the Sobolev Embedding Theorem (see Adams and Fournier \cite[Theorem 4.12, p. 85]{AdamsFournier}), so the subset \eqref{eq:X0}, namely
\[
  X_0 := \{x \in X: \psi(x) \neq 0\},
\]
is open and consider $\wp(\psi)$ in \eqref{eq:Definition_wp_intro}, namely
\[
  \wp(\psi)
  =
  \begin{cases}
    2|\psi|_{\Lambda^{0,2}(E)}^{-2}\star(\psi\otimes\psi^*)_0 &\text{on } X_0,
    \\
    0 &\text{on } X\less X_0.
  \end{cases}
\]
We observe that $\wp(\psi) \in L^\infty(X;i\su(E))$ and so the right-hand side of equation \eqref{eq:SO(3)_monopole_equations_(1,1)_curvature_perturbed_omega_intro} is in $L^\infty(X;\su(E))$ if $\varphi \in L^\infty(E)$ and $\psi \in L^\infty(\Lambda^{0,2}(E)$. Hence, given equations \eqref{eq:SO(3)_monopole_equations_(0,2)_curvature_perturbed_intro} and \eqref{eq:SO(3)_monopole_equations_(2,0)_curvature_intro} --- whose right-hand sides are in $L^\infty(\Lambda^{0,2}(\fsl(E)))$ and $L^\infty(\Lambda^{2,0}(\fsl(E)))$, respectively, if $\varphi \in L^\infty(E)$ and $\psi \in L^\infty(\Lambda^{0,2}(E)$ --- it is appropriate to consider unitary connections $A \in \sA^{1,p}(E,H,A_d)$ with $2 < p < \infty$ in the presence of the Coulomb gauge equation \eqref{eq:SO3_monopole_Coulomb_gauge_slice_condition} required to obtain an elliptic system, namely
\[
  d_{A,\Phi}^{0,*}(A - A_\infty,\Phi - \Phi_\infty) = 0,
\]
where $(A_\infty,\Phi_\infty)$ is a $C^\infty$ reference pair. Given $A$ with this regularity, the elliptic Dirac equation \eqref{eq:SO(3)_monopole_equations_Dirac_almost_Hermitian_perturbed_intro} will have solutions $\varphi \in W^{2,p}(E)$ and $\psi \in W^{2,p}(\Lambda^{0,2}(E))$, guaranteeing in turn that $\wp \in L^\infty(\su(E))$. On the other hand, because regularity of solutions to quasi-linear elliptic systems (on smooth vector bundles over smooth manifolds) is a local property, we obtain
\[
  (A,\Phi) \in \sA^\infty(X_0;E,H,A_d) \times C^\infty(X_0;E) \times C^\infty(X_0;\Lambda^{2,0}(\fsl(E))),
\]
and the aforementioned global regularity on $X$. As usual, because of our reliance on the Coulomb gauge equation \eqref{eq:SO3_monopole_Coulomb_gauge_slice_condition} to obtain an elliptic system for $(A,\Phi)$, the $C^\infty$ regularity results over $X_0$ are obtained after applying a gauge transformation $u \in W^{2,p}(\SU(E))$.

Note that, over $X_0$, taking the differential of the expression \eqref{eq:Definition_wp_intro} for $\wp(\psi)$ yields
\begin{multline}
  \label{eq:Derivative_of_Taubes_Pert}
  D\wp(\psi)\tau
  = -2|\psi|_{\Lambda^{0,2}(E)}^{-4}\left( \langle\tau,\psi\rangle_{\Lambda^{0,2}(E)}
    + \langle\psi,\tau\rangle_{\Lambda^{0,2}(E)} \right)\star(\psi\otimes\psi^*)_0
\\
+ 2|\psi|_{\Lambda^{0,2}(E)}^{-2}\star\left(\tau\otimes\psi^* + \psi\otimes\tau^*\right)_0
\in \Omega^0(X_0;i\su(E)).
\end{multline}
Unfortunately, because of the \emph{singularities} along $\{\psi=0\} \subset X$, the expression for the zero order differential operator $D\wp(\psi)$ in \eqref{eq:Derivative_of_Taubes_Pert} does \emph{not} yield a perturbation of the linearization of the system \eqref{eq:SO(3)_monopole_equations_almost_Hermitian_intro} of unperturbed non-Abelian monopole equations and a Coulomb gauge slice condition that is \emph{compact} for standard choices of $W^{1,p}$ and $L^p$ Sobolev spaces on its domain and codomain, respectively. This makes it very difficult to construct local Kuranishi models for open neighborhoods in the moduli space of non-Abelian monopoles.

Therefore, we turn to the \emph{regularized} Taubes perturbation term $\wp_\gamma(\psi)$ in \eqref{eq:Definition_wp_intro_regular}, namely\[
  \wp_\gamma(\psi)
  =
  4\left(\gamma^2 + |\psi|_{\Lambda^{0,2}(E)}^2\right)^{-1}\star(\psi\otimes\psi^*)_0
  \in \Omega^0(i\su(E)),
\]  
and observe that
\begin{multline}
  \label{eq:Derivative_of_Taubes_Pert_regular}
  D\wp_\gamma(\psi)\tau
  = -4\left(\gamma^2 + |\psi|_{\Lambda^{0,2}(E)}^2\right)^{-2}\left( \langle\tau,\psi\rangle_{\Lambda^{0,2}(E)}
    + \langle\psi,\tau\rangle_{\Lambda^{0,2}(E)} \right)\star(\psi\otimes\psi^*)_0
\\
+ 4\left(\gamma^2 + |\psi|_{\Lambda^{0,2}(E)}^2\right)^{-1}
\star\left(\tau\otimes\psi^* + \psi\otimes\tau^*\right)_0
\in \Omega^0(i\su(E)).
\end{multline}
In contrast with the zero order differential operator $D\wp(\psi)$ in \eqref{eq:Derivative_of_Taubes_Pert}, the zero order differential operator $D\wp_\gamma(\psi)$ in \eqref{eq:Derivative_of_Taubes_Pert_regular} has \emph{smooth coefficients} and so yields a perturbation of the linearization of the system \eqref{eq:SO(3)_monopole_equations_almost_Hermitian_intro} of unperturbed non-Abelian monopole equations and a Coulomb gauge slice condition that is \emph{compact} for standard choices of $W^{1,p}$ and $L^p$ Sobolev spaces on its domain and codomain, respectively, using the fact that the embedding $W^{1,p}(X) \subset L^p(X)$ is compact by the Rellich--Kondrachov Embedding Theorem for $1\leq p < \infty$ (see Adams and Fournier \cite[Theorem 6.3, p. 168]{AdamsFournier}).

\subsection{Fredholm sections of Banach vector bundles}
\label{subsec:Fredholm_section_Banach_vector_bundle}
In this subsection, we clarify the definition of a Fredholm section of a Banach vector bundle over a Banach manifold and recall the construction of the local Kuranishi model for an open neighborhood of a point in the zero locus of a Fredholm section that is not assumed to be a regular point.

We begin by recalling the definition of an analytic map of Banach spaces (see Berger \cite[Section 2.3]{Berger_1977}, Huang \cite[Section 2.1A]{Huang_2006}, or Whittlesey \cite{Whittlesey_1965}). Let $\FF, \EE$ be Banach spaces (over $\RR$ or $\CC$), let $\sU\subset\FF$ be an open subset, and $\sF:\sU \to \EE$ be a map. Recall that $\sF$ is \emph{Fr\'echet differentiable} at a point $x \in \sU$ with \emph{Fr\'echet derivative} $\sF'(x):\FF \to \EE$ (a bounded linear operator) if
\begin{equation}
\label{eq:Frechet_derivative}
\lim_{y\to 0} \frac{1}{\|y\|_\FF}\|\sF(x + y) - \sF(x) - \sF'(x)y\|_\EE = 0.
\end{equation}
% PF11-13-2024 Add page numbers
Recall from Berger \cite[Definition 2.3.1]{Berger_1977}, Deimling \cite[Definition 15.1]{Deimling_1985}, or Zeidler \cite[Definition 8.8]{Zeidler_nfaa_v1} that $\sF$ is \emph{analytic} at $x \in \sU$ if there exist a constant $r > 0$, a sequence of continuous, symmetric, $n$-linear maps $L_n:\FF \times \cdots \times \FF \to \EE$ obeying
\[
  \sum_{n\geq 1} \|L_n\| r^n < \infty,
\]
and a positive constant $\delta = \delta(x)$ such that
\begin{equation}
\label{eq:Taylor_expansion}
\sF(x + y) = \sF(x) + \sum_{n\geq 1} L_n(y^n), \quad\text{for all } y \in \FF \text{ with } \|y\|_\FF < \delta,
\end{equation}
where $y^n \equiv (y,\ldots,y) \in \FF \times \cdots \times \FF$ ($n$-fold product). If $\sF$ is Fr\'echet differentiable (respectively, analytic) at every point $x \in \sU$, then $\sF$ is Fr\'echet differentiable (respectively, analytic) on $\sU$. If there is integer $N \geq 0$ such that $L_n = 0$ in \eqref{eq:Taylor_expansion} for all $n \geq N$, then $\sF$ is a degree-$N$ polynomial.

For definitions and properties of \emph{Banach manifolds}, we refer to Abraham, Marsden, and Ratiu \cite[Definition 3.1.3, p. 143]{AMR}, Abraham and Robbin \cite[Chapter 1, Section 3, p. 9]{Abraham_Robbin_transversal_mappings_flows}, Dodson, Galanis, and Vassiliou \cite[Section 1.1.2, p. 2]{Dodson_Galanis_Vassiliou_geometry_frechet_context}, Klingenberg \cite[Definition 1.1.2, p. 9]{Klingenberg_riemannian_geometry}, Lang \cite[Chapter II, Section 1, p. 21]{Lang_introduction_differential_topology}, or Zeidler \cite[Definition 73.4, p. 534]{Zeidler_nfaa_v4}.
% PF11-8-2024 Add Bourbaki
For definitions and properties of \emph{Banach vector bundles}, see Abraham, Marsden, and Ratiu \cite[Definition 3.4.1, p. 167]{AMR}, Abraham and Robbin \cite[Chapter 1, Section 5, p. 12]{Abraham_Robbin_transversal_mappings_flows}, Dodson, Galanis, and Vassiliou \cite[Section 1.4, p. 14]{Dodson_Galanis_Vassiliou_geometry_frechet_context}, Lang \cite[Chapter III, Section 1, p. 37]{Lang_introduction_differential_topology}, or Zeidler \cite[Definition 73.20, p. 544]{Zeidler_nfaa_v4}.
% PF11-8-2024 Add Bourbaki

Let $\cM$ be a $C^k$ Banach manifold modeled on a Banach space $\FF$ (over $\RR$ or $\CC$), where $k \in \ZZ_{\geq 1}$, $k=\infty$, or $\omega$ (analytic), and $\pi:\cE \to \cM$ be a $C^k$ Banach vector bundle modeled on a Banach space $\EE$ (over $\RR$ or $\CC$). There are two approaches to the definition of a $C^k$ \emph{Fredholm section} of $\cE$, intrinsically and extrinsically through choices of local trivializations.

%PF12-18-2025 Cite Salamon too
First, we give an \emph{intrinsic} definition of a Fredholm section. Let $s:\cM \to \cE$ be a $C^k$ section of $\cE$, so $ds:T\cM \to T\cE$ is the differential of $s$ on the tangent bundles to $\cM$ and $\cE$. Let $\iota:\cM \hookrightarrow \cE$ denote the $C^k$ embedding of $\cM$ into $\cE$ as the zero section.
% COMMENT Cite Georges Elencwajg's answer in https://math.stackexchange.com/questions/251646/splitting-of-the-tangent-bundle-of-a-vector-bundle
The morphism $\pi:\cE \to \cM$ of Banach manifolds (which is a submersion) induces a surjective differential map,
\[
  T\cE \xrightarrow{d\pi} T\cM \xrightarrow{} 0
\]
whose kernel is (by definition) the \emph{vertical tangent bundle} $T^v\cE$, thus yielding an exact sequence of Banach vector bundles on $\cE$:
\[
  0 \xrightarrow{} T^v\cE \xrightarrow{} T\cE \xrightarrow{d\pi} \pi^*T\cM \xrightarrow{} 0
\]  
By pulling back the preceding sequence to $\cM$ by using the embedding $\iota$, we obtain an exact of Banach vector bundles on $\cM$:
\[
  0 \xrightarrow{} \cE \xrightarrow{} T\cE\restriction\cM \xrightarrow{} T\cM \xrightarrow{} 0
\]
We assume that $\cM$ is paracompact, in which case the preceding exact sequence splits (non-canonically) and we obtain an isomorphism (non-canonical) of Banach vector bundles:
\[
  T\cE\restriction\cM \cong \cE \oplus T\cM.
\]
The bundles $T^v\cE$ and $\pi^*\cE$ are canonically isomorphic as Banach vector bundles on $\cE$, so that we have a canonical exact sequence on $\cE$:
\[
  0 \xrightarrow{} \pi^*\cE \xrightarrow{} T\cE \xrightarrow{d\pi} \pi^*T\cM \xrightarrow{} 0
\]
%COMMENT https://mathoverflow.net/questions/430735/history-of-use-of-symbol-to-mean-is-canonically-isomorphic-to
The canonical isomorphism $T^v\cE = \pi^*\cE$ ultimately rests on the fact that the tangent space to a vector space $V$ at any point $w \in V$ is canonically isomorphic to that vector space, that is, $T_wV = V$. See Spivak \cite[Chapter 3, Exercise 29, p. 103]{Spivak1} and tom Dieck \cite[Kapitel IX, Satz 6.9. p. 365]{tomDieck_topologie} for additional details.

From the preceding discussion, we have an isomorphism (non-canonical) of Banach spaces,
\[
  T_p\cE \cong T_p\cM \oplus \cE_p,
\]
%PF11-8-2024 Reverse the order earlier
for each $p \in \cM$ and so, by composing this isomorphism with the differential $ds(p):T_p\cM \to T_p\cE$, we may write
\[
  ds(p):T_p\cM \to T_p\cM \oplus \cE_p.
\]
Then $s:\cM\to \cE$ is a \emph{Fredholm section} if the composition
\[
  \pr_2\circ\, ds(p):T_p\cM \to \cE_p
\]
is a Fredholm operator for each $p\in\cM$, where $\pr_2:T_p\cM \oplus \cE_p$ is projection onto the second factor. This is essentially the definition provided by Hofer \cite[Introduction, p. 3]{Hofer_2006}.

Second, we recall the equivalent \emph{extrinsic} definition. To begin, we recall that $f:\cM \to \EE$ is a $C^k$ \emph{Fredholm map} if $df(p) \in \Hom(T_p\cM,\EE)$ is a Fredholm operator for every point $p \in \cM$ (see Donaldson and Kronheimer \cite[Section 4.2.4, p. 136]{DK}.) Now let $\cT:\cE\restriction\cU \cong \cU \times \EE$ be a $C^k$ local trivialization over an open neighborhood $\cU \subset \sM$, where the fiber $\cE_p$ is isomorphic to $\EE$ for each $p\in\cU$. We thus obtain a $C^k$ map,
\[
  s_\cU \equiv \pr_2\circ \cT\circ s:\cU \to \EE,
\]
with differential,
\[
  ds_\cU(p) = \pr_2\circ\, d\cT(s(p))\circ ds(p): T_p\cM \to \EE,
\]
where the trivialization $\cT$ has differential at $s(p) \in \cE$ given by
\[
  d\cT(s(p)):T_{s(p)}\cE \to T_p\cM \oplus \EE,
\]
and $\pr_2:T_p\cM \oplus \EE \to \EE$ is projection onto the second factor. Then $s:\cM\to \cE$ is a \emph{Fredholm section} if the induced map $s_\cU:\cU \to \EE$ is Fredholm for every choice of local trivialization $(\cT,\cU)$ of $\cE$. This is the definition provided by Donaldson \cite[Section 1, p. 2]{Donaldson_2022ggtc} and Donaldson and Kronheimer \cite[Section 4.2.4, p. 137]{DK}, and the one that we adopt in our application.

If $\phi:\cU \to \FF$ is a local coordinate chart on $\cU \subset \cM$, then its differential $d\phi(p):T_p\cM \to \FF$ is an isomorphism of Banach spaces for each $p\in\cU$ and we obtain a Fredholm operator between fixed Banach spaces:
\[
  d(s_\cU\circ\phi^{-1})(x) = \pr_2\circ\, d\cT(s(p))\circ ds(p) \circ (d\phi(p))^{-1}: \FF \to \EE,
\]
for each point $p \in \cU$ and $x = \phi(p) \in \phi(\cU) \subset \FF$.

\subsection{Local Kuranishi model for a nonlinear map whose differential has a pseudoinverse}
\label{subsec:Local_Kuranishi_model_nonlinear_map_Banach_spaces}
We now prove a useful generalization of the standard Kuranishi model defined by a Fredholm map of Banach spaces (and hence a Fredholm section of a $C^k$ Banach vector bundle over a Banach manifold), as described by Donaldson and Kronheimer \cite[Section 4.2.4]{DK}, Freed and Uhlenbeck \cite[Chapter 4]{FU}, and Salamon \cite[Appendix B, Section 2]{SalamonSWBook}. Salamon more generally allows maps whose differentials are operators that admit a \emph{pseudoinverse} \cite[Appendix B, Section 2, Proposition B.7, p. 493]{SalamonSWBook}, a class of operators that includes Fredholm operators.

\begin{defn}[Pseudoinverse of a bounded linear operator on Banach spaces]
\label{defn:Pseudoinverse_bounded_linear_operator_Banach_spaces}
(See Salamon \cite[Appendix B, Section 2, p. 494]{SalamonSWBook} and Sheffield \cite[Chapter I, Definition 1.3, p. 7]{SheffieldThesis}, as well as Beutler \cite{Beutler_1965a, Beutler_1965a}, Desoer and Whalen \cite{Desoer_Whalen_1963}, Koliha \cite{Koliha_1974bams, Koliha_1974jmaa}, and references therein for alternative definitions.)
Let $X, Y$ be Banach spaces and $D \in \Hom(X,Y)$ be a bounded linear operator. One says that $T \in \Hom(Y,X)$ is a \emph{pseudoinverse} of $D$ if
\[
  TDT = T \quad\text{and}\quad DTD = D.
\]
\qed\end{defn}

\begin{prop}[Necessary and sufficient condition for existence of a pseudoinverse]
\label{prop:Salamon_B-7}  
(See Salamon \cite[Appendix B, Section 2, Proposition B.7, p. 494]{SalamonSWBook}.)  
Continue the notation of Definition \ref{defn:Pseudoinverse_bounded_linear_operator_Banach_spaces}. Then $D$ admits a pseudoinverse if and only if the following hold:
\begin{enumerate}
\item $D$ has closed range.
\item $\Ker D$ has a complement in $X$.
\item $\Ran D$ has a complement in $Y$.    
\end{enumerate}
\end{prop}

Proposition \ref{prop:Salamon_B-7} implies that every \emph{Fredholm} operator has a pseudoinverse. 

\begin{exmp}[Construction of pseudoinverses via spectral projections on Hilbert spaces]
\label{exmp:Construction_pseudoinverses_by_spectral_projections_Hilbert_spaces}  
Continue the notation of Definition \ref{defn:Pseudoinverse_bounded_linear_operator_Banach_spaces}, but assume that $X, Y$ are Hilbert spaces and that the operators $D^*D$ and $DD^*$ have discrete spectra of eigenvalues with finite multiplicities. Let $\{\nu_k\}_{k=1}^\infty$ denote the non-decreasing sequence of positive eigenvalues, with orthonormal bases of eigenvectors $\{\varphi_k\}_{k=1}^\infty$ for $(\Ker D)^\perp\cap X$ and $\{\psi_k\}_{k=1}^\infty$ for $(\Ker D^*)^\perp\cap Y$ given by $\psi_k = \nu_k^{-1/2}D\varphi_k$ for all $k\geq 1$ (see the forthcoming Lemma \ref{lem:Eigenvalues_densely_defined_unbounded_linear_operators}). Let $\nu \in (0,\infty) \less \{\nu_k\}_{k=1}^\infty$ and $\Pi_\nu^X$ be the orthogonal projection onto the linear span of the eigenvectors of $D^*D$ with eigenvalues less than $\nu$ and  $\Pi_\nu^Y$ be the orthogonal projection onto the linear span of the eigenvectors of $DD^*$ with eigenvalues less than $\nu$. The bounded linear operator
\[
  D: \left(\Ran\Pi_\nu^X\right)^\perp \to \left(\Ran\Pi_\nu^Y\right)^\perp
\]
is an isomorphism of Banach spaces with inverse 
\[
  D^{-1}: \left(\Ran\Pi_\nu^Y\right)^\perp \to \left(\Ran\Pi_\nu^X\right)^\perp.
\]
Define a bounded linear operator $T_\nu \in \Hom(Y,X)$ by
\[
  T_\nu
  :=
  \begin{cases}
    0 &\text{on } \Ran\Pi_\nu^Y,
    \\
    D^{-1} &\text{on } \left(\Ran\Pi_\nu^Y\right)^\perp.
  \end{cases}
\]
Similarly, define a bounded linear operator $D_\nu \in \Hom(X,Y)$ by
\[
  D_\nu
  :=
  \begin{cases}
    0 &\text{on } \Ran\Pi_\nu^X,
    \\
    D &\text{on } \left(\Ran\Pi_\nu^X\right)^\perp.
  \end{cases}
\]
We calculate that
\[
  D_\nu T_\nu = \id_Y - \Pi_\nu^Y
  \quad\text{and}\quad
  T_\nu D_\nu = \id_X - \Pi_\nu^X,
\]
and therefore, using $T_\nu = 0$ on $\Ran\Pi_\nu^Y$ and $D_\nu = 0$ on $\Ran\Pi_\nu^X$,
\[
  T_\nu D_\nu T_\nu = T_\nu \quad\text{on } Y
  \quad\text{and}\quad
  D_\nu T_\nu D_\nu = D_\nu \quad\text{on } X.
\]
Thus, $T_\nu$ is a pseudoinverse of $D_\nu$ and a pseudoinverse of $D$ when $0<\nu<\nu_1$, in which case $D = D_\nu$. \qed
\end{exmp}

\begin{exmp}[Construction of pseudoinverses via continuous projections on Banach spaces]
\label{exmp:Construction_pseudoinverses_continuous_projections_Banach_spaces}  
Continue the notation of Definition \ref{defn:Pseudoinverse_bounded_linear_operator_Banach_spaces}, let $K \subset X$ be a closed subspace with complement $X_0$ such that $\Ker D \subset K$, assume that $D$ has closed range $\Ran D$ with complement $C_0$ and that $C \subset Y$ is a closed subspace with complement $Y_0$ such that $C_0 \subset C$. Assume further that $D:X_0 \to Y_0$ is an isomorphism of Banach spaces with inverse $D^{-1}:Y_0 \to X_0$. Define a bounded linear operator $T_1 \in \Hom(Y,X)$ by
\[
  T_1
  :=
  \begin{cases}
    0 &\text{on } C,
    \\
    D^{-1} &\text{on } Y_0.
  \end{cases}
\]
Similarly, define a bounded linear operator $D_1 \in \Hom(X,Y)$ by
\[
  D_1
  :=
  \begin{cases}
    0 &\text{on } K,
    \\
    D &\text{on } X_0.
  \end{cases}
\]
Writing $\pr_C$ for the continuous projection $Y = Y_0 \oplus C \to C$ and similarly for the other projections, we get
\[
  T_1 = D^{-1}\circ\pr_{Y_0}
  \quad\text{and}\quad
  D_1 = D\circ\pr_{X_0},
\]  
and so
\begin{align*}
  D_1 T_1 &= D\circ\pr_{X_0}D^{-1}\circ\pr_{Y_0} = DD^{-1}\circ\pr_{Y_0} = \id_{Y_0}\pr_{Y_0} = \pr_{Y_0} = \id_Y - \pr_C,
  \\
  T_1 D_1 &= D^{-1}\circ\pr_{Y_0}D\circ\pr_{X_0} = D^{-1}D\circ\pr_{X_0} = \id_{X_0}\pr_{X_0} = \pr_{X_0} = \pr_{X_0} = \id_X - \pr_K.
\end{align*}
Therefore, using $T_1 = 0$ on $C$ and $D_1 = 0$ on $K$,
\[
  T_1 D_1 T_1 = T_1 \quad\text{on } Y
  \quad\text{and}\quad
  D_1 T_1 D_1 = D_1 \quad\text{on } X.
\]
Thus, $T_1$ is a pseudoinverse of $D_1$ and a pseudoinverse of $D$ when $K = \Ker D$, in which case $D = D_1$.

If $G$ is a topological group with continuous group homomorphisms $G \to \End(X)$ and $G \to \End(Y)$ into Banach spaces of bounded linear operators, $D$ is equivariant with respect to the induced actions of $G$ on $X$ and $Y$, and the subspaces $K \subset X$ and $C \subset Y$ are $G$-invariant, then our definitions of $D_1$ and $T_1$ immediately imply that those operators are also $G$-equivariant.
\qed
\end{exmp}  

We now have the
\begin{thm}[Local Kuranishi model for a nonlinear map of Banach spaces whose differential has a pseudoinverse]
\label{thm:Kuranishi_model_defined_by_Fredholm_map_Banach_spaces} 
(Compare Donaldson and Kronheimer \cite[Section 4.2.4, Proposition 4.2.9, p. 136]{DK}, Freed and Uhlenbeck \cite[Lemma 4.7, p. 66]{FU}, and Salamon \cite[Appendix B, Section 2, Theorem B.9, p. 494 and Remark B.10, p. 495]{SalamonSWBook}.)
Let $F:U \to Y$ be a $C^k$ map into a Banach space $Y$ from an open neighborhood $U$ of the origin in a Banach space $X$, where $k \in \ZZ_{\geq 1}$, $k=\infty$, or $\omega$ (analytic) such that $F(0)=0$. Let $K \subset X$ be a closed subspace with closed complement $X_0$ such that
%PF7-1-2025 Eliminated K_0, C_0. Check applications.
%$K_0 = \Ker dF(0) \subset K$
$\Ker dF(0) \subset K$
and 
%PF7-1-2025 We don't need or use this closed range condition. Check applications.
%$dF(0)$ has closed range $\Ran dF(0)$ with complement $C_0$ and that
$C \subset Y$ be a closed subspace with closed complement $Y_0$ such that
\begin{equation}
  \label{eq:dF(0)X0_equals_Y0_and_dF(0)K_subset_C}
  dF(0)X_0 = Y_0 \quad\text{and}\quad dF(0)K \subset C.
\end{equation}  
Then, for a possibly smaller open neighborhood $U' \subset U$ of $0$ in $X$, there exist a $C^k$ diffeomorphism $g:U' \to g(U') \subset U$ and a $C^k$ map $f:U' \to C$ such that
\begin{equation}
  \label{eq:F_circ_g_is_dF_proj+f}
  F\circ g = dF(0)\pr_{X_0} + f \quad\text{on } U',
\end{equation}
where $\pr_{X_0}:X = X_0\oplus K \to X_0$ denotes the continuous projection, and if $\pr_C:Y = Y_0\oplus C \to C$ denotes the continuous projection, then
\begin{equation}
  \label{eq:g(0)_and_dg(0)_and_f(0)_and_df(0)}
  g(0) = 0, \quad dg(0) = \id_X, \quad f(0) = 0, \quad\text{and}\quad df(0) = \pr_CdF(0).
\end{equation}
Moreover, if $\iota_K:K \to X$ is the inclusion of $K$ in $X = K\oplus X_0$ and we denote $U_K' = U'\cap K$, so $\iota_K(U_K') = U_K' \oplus (0) \subset U'$, and $f_K = f\circ\iota_K: U_K' \to C$ and $g_K = g\circ\iota_K:U_K' \to U \subset X$, then
\begin{equation}
  \label{eq:Equality_zero-locus_F_and_zero-locus_fK}
  F^{-1}(0) \cap g(U') = f_K^{-1}(0)\cap U_K'
  \quad\text{and}\quad
  g_K\left(f_K^{-1}(0)\cap U_K'\right) \subset F^{-1}(0)\cap U.
\end{equation}
Finally, if $F$ is equivariant with respect to the action of a compact Lie group $G$ on $X$ and $Y$ and the subspaces $K \subset X$ and $C \subset Y$ are $G$-invariant, then the neighborhood $U'$ can be chosen to be $G$-invariant and the maps $g$ and $f$ can be chosen to be $G$-equivariant.
\end{thm}

\begin{rmk}[Special cases of Theorem \ref{thm:Kuranishi_model_defined_by_Fredholm_map_Banach_spaces}]
\label{rmk:Kuranishi_model_defined_by_Fredholm_map_Banach_spaces_special_case}  
If $K = \Ker dF(0)$ with closed complement $X_0$ and $Y_0 = \Ran dF(0)$ is closed with closed complement $C$, then $dF(0)\circ \pr_{X_0} = dF(0)$ and $df(0) = \pr_CdF(0) = 0$. Moreover, if $dF(0)$ is a Fredholm operator, then $\Ker dF(0)$ is finite-dimensional and $\Ran dF(0)$ is closed with finite-dimensional complement $C$.
\qed\end{rmk}  

\begin{proof}
Because $dF(0):X_0 \to Y_0$ is a bijective continuous linear operator between Banach spaces, the Open Mapping Theorem (see, for example, Brezis \cite[Section 2.3, Theorem 2.6, p. 35]{Brezis}) implies that it is an isomorphism of Banach spaces with inverse $dF(0)^{-1}:Y_0 \to X_0$. Define
  %PF-2-2025 Reverse roles of D and T to match typical usage for T
\[
  D := dF(0) \circ \pr_{X_0} \in \Hom(X,Y) \quad\text{and}\quad T := dF(0)^{-1}\circ\pr_{Y_0} \in \Hom(Y,X),
\]  
so that $T$ is a pseudoinverse of $D$ by Example \ref{exmp:Construction_pseudoinverses_continuous_projections_Banach_spaces} with $D\circ T\circ D = D$ and
\[
  \Ran D = \Ran dF(0) \circ \pr_{X_0} = Y_0 \quad\text{and}\quad \Ran T = \Ran dF(0)^{-1}\circ\pr_{Y_0} = X_0.
\]
We write 
\begin{align*}
  F &= D + F - D
  \\
    &= \pr_C\circ F + D + \pr_{Y_0}\circ F - D \quad\text{(using $\id_Y = \pr_{Y_0}+\pr_C$)}
  \\
    &= \pr_C\circ F\circ\pr_K + D
      + \pr_C\circ (F - F\circ\pr_K) + \pr_{Y_0}\circ F - D
  \\
    &= \pr_C\circ F\circ\pr_K + D
      + \pr_C\circ (F - F\circ\pr_K) + \pr_{Y_0}\circ F\circ\pr_{X_0}
      + \pr_{Y_0}\circ (F - F\circ\pr_{X_0}) - D.
\end{align*}
We define a $C^k$ function $\Phi:U\to Y$ by
\[
  \Phi := \pr_C\circ (F - F\circ\pr_K) + \pr_{Y_0}\circ F\circ\pr_{X_0}
      + \pr_{Y_0}\circ (F - F\circ\pr_{X_0}) - D \quad\text{on } U,
\]
so that
\begin{equation}
  \label{eq:F_equals_pr_C_circ_F_circ_pr_K+D+Phi}
  F = \pr_C\circ F\circ\pr_K + D + \Phi,
\end{equation}
and define a $C^k$ function $\Psi:U\to X$ by
\begin{equation}
  \label{eq:Psi_define_idX+TPhi}
  \Psi := \id_X + T\circ\Phi \quad\text{on } U.
\end{equation}
Because $F(0) = 0$, we have $\Phi(0) = 0$ and $\Psi(0) = 0$. We calculate
\begin{align*}
  d\Phi(0)
  &= \pr_C\circ \left(dF(0) - dF(0)\circ\pr_K\right) + \pr_{Y_0}\circ\, dF(0)\circ\pr_{X_0}
             + \pr_{Y_0}\circ \left(dF(0) - dF(0)\circ\pr_{X_0}\right) - D
  \\
  &= \pr_C\circ\, dF(0)\circ\left(\id_X - \pr_K\right) + \pr_{Y_0}\circ\, dF(0)\circ\pr_{X_0}
             + \pr_{Y_0}\circ\, dF(0)\circ\left(\id_X - \pr_{X_0}\right) - D
  \\
  &= \pr_C\circ\, dF(0)\circ\pr_{X_0} + \pr_{Y_0}\circ\, dF(0)\circ\pr_{X_0}
    + \pr_{Y_0}\circ\, dF(0)\circ\pr_K - D
  \\
  &\qquad\text{(using $\id_X = \pr_K + \pr_{X_0}$)}
  \\
  &= 0 + \pr_{Y_0}\circ\, dF(0)\circ\pr_{X_0} + 0 - D  
  \\
  &\qquad \text{(by \eqref{eq:dF(0)X0_equals_Y0_and_dF(0)K_subset_C} with $\Ran dF(0)\circ\pr_{X_0} = Y_0$ and
    $\Ran dF(0)\circ\pr_K \subset C$)}
  \\
  &= D - D = 0 \quad\text{(because $D = dF(0)\circ\pr_{X_0} = \pr_{Y_0}\circ\, dF(0)\circ\pr_{X_0}$)}.
\end{align*}
Therefore, we obtain
\[
  d\Psi(0)
  = \id_X + T\circ d\Phi(0)
  = \id_X.
\]
Hence, for a possibly smaller open neighborhood $U' \subset U$ of $0$ in $X$, the Inverse Mapping Theorem
(see, for example, Abraham, Marsden, and Ratiu \cite[Theorem 2.5.1, p. 116]{AMR}) provides a $C^k$ diffeomorphism $g:U' \to X$ onto an open neighborhood $g(U')$ of $0$ in $X$ such that $g = \Psi^{-1}$ on $U'$. Observe that
\begin{align*}
  D\circ\Psi
  &= D + D\circ T\circ\Phi \quad\text{(by \eqref{eq:Psi_define_idX+TPhi})}
  \\
  &= D + D\circ T\circ\left(F - \pr_C\circ F\circ\pr_K - D\right)
    \quad\text{(by \eqref{eq:F_equals_pr_C_circ_F_circ_pr_K+D+Phi})}
  \\
  &= D + D\circ T\circ F - D\circ T\circ \pr_C\circ F\circ\pr_K - D\circ T\circ D
  \\
  &= D + D\circ T\circ F - 0 - D
    \quad\text{(because $T = 0$ on $C$ and $D\circ T\circ D = D$)}
  \\
  &= D\circ T\circ F \quad\text{on } U.
\end{align*}
Thus,
\[
  D = D\circ T\circ F\circ\,\Psi^{-1} \quad\text{on } U'.
\]  
Therefore, using $D - D\circ T\circ F\circ\,\Psi^{-1} = 0$ on $U'$,
\begin{align*}
  F\circ\Psi^{-1}
  &=
  D - D\circ T\circ F\circ\,\Psi^{-1} + \id_Y\circ F\circ\Psi^{-1}
  \\
  &=
  D  + \left(\id_Y - D\circ T\right)\circ F\circ\,\Psi^{-1} \quad\text{on } U'.
\end{align*}
Define a $C^k$ map $f:U' \to Y_0$ by
\[
  f := (\id_Y - D\circ T)\circ F \circ\, \Psi^{-1} = \pr_C\circ F \circ\, g,
\]
where the second inequality follows from the facts that $g = \Psi^{-1}$ and $D\circ T = \pr_{Y_0} = \id_Y - \pr_C$. We have thus shown that
\[
   F\circ g = D + f \quad\text{on } U,
\]
which is the desired equality \eqref{eq:F_circ_g_is_dF_proj+f}. Observe that
\[
  g(0) = \Psi^{-1}(0) = 0, \quad dg(0) = (d\Psi(0))^{-1} = \id_X, \quad f(0) = 0, 
\]  
and
\begin{align*}
  df(0)
  &= \pr_C\circ\, dF(0) \circ\, dg(0)
  \\
  &= \pr_C\circ\, dF(0)\circ \id_X
  \\
  &= \pr_C\circ\, dF(0).
\end{align*}
The identification of the zero loci in \eqref{eq:Equality_zero-locus_F_and_zero-locus_fK} follows exactly as in Salamon \cite[Appendix B, Section 2, Remark B.10, p. 495]{SalamonSWBook}. Indeed, observe that
\begin{align*}
  F^{-1}(0) \cap g(U')
  &=
    \{g(x) \in g(U'): (F\circ g)(x) = 0\}
  \\
  &= \{g(x) \in g(U'): (dF(0)\circ \pr_{X_0} + f)(x) = 0\} \quad\text{(by \eqref{eq:F_circ_g_is_dF_proj+f})}
  \\
  &= \{x \in U': dF(0)\circ \pr_{X_0}(x) = 0 \text{ and } f(x) = 0\}
  \\
  &\qquad\text{(as $dF(0)\circ \pr_{X_0}(x) \in Y_0$ and $f(x) \in C$ and $Y = Y_0\oplus C$)}
   \\
  &= \{x \in K\cap U': f(x) = 0\} \quad\text{(as $\Ker (dF(0)\circ \pr_{X_0}) = K$)}
  \\
  &= (f\circ\iota_K)^{-1}(0)\cap K\cap U',
\end{align*}
and thus, noting that $f_K = f\circ\iota_K$ and $U_K' = K\cap U'$, we obtain the desired identity in \eqref{eq:Equality_zero-locus_F_and_zero-locus_fK}, namely
\[
  F^{-1}(0) \cap g(U') = f_K^{-1}(0)\cap U_K'.
\]
The inclusion in \eqref{eq:Equality_zero-locus_F_and_zero-locus_fK} follows from the preceding equalities. Indeed, if $x \in f_K^{-1}(0)\cap U_K'$, then $g(x) \in F^{-1}(0) \cap g(U')$ and because $g(U') \subset U$, this yields the desired inclusion,
\[
  g\left(f_K^{-1}(0)\cap U_K'\right) \subset F^{-1}(0) \cap U.
\]  
Finally, if $F$ is $G$-equivariant, we may choose a $G$-equivariant pseudoinverse $T$ of $D$ by Example \ref{exmp:Construction_pseudoinverses_continuous_projections_Banach_spaces} or, when $G$ is compact (by a more complicated argument) via Salamon \cite[Appendix B, Section 2, Remark B.8, p. 494]{SalamonSWBook}. Therefore, the definitions of $\Psi$ and hence the maps $g$ and $f$ imply that they are also $G$-equivariant and that the open neighborhood $U'$ can be chosen to be $G$-invariant. This completes the proof of Theorem \ref{thm:Kuranishi_model_defined_by_Fredholm_map_Banach_spaces}.
\end{proof}

Although Theorem \ref{thm:Kuranishi_model_defined_by_Fredholm_map_Banach_spaces} already generalizes the well-known simpler versions of the Kuranishi Lemma cited, the following more flexible variant is useful in our applications. For simplicity, we state the result in the setting of Hilbert spaces only.

\begin{cor}[Local Kuranishi model for a nonlinear map of Hilbert spaces]
\label{cor:Kuranishi_model_defined_by_Fredholm_map_Hilbert_spaces}
%PF7-9-2025 Harmonize notation  
Let $F:U \to Y$ be a $C^k$ map into a Hilbert space $Y$ from an open neighborhood $U$ of the origin in a Hilbert space $X$, where $k \in \ZZ_{\geq 1}$, $k=\infty$, or $\omega$ (analytic) such that $F(0)=0$ and $dF(0) \in \Hom(X,Y)$ is Fredholm. Let $K_1 \subset X$ and $C_1 \subset Y$ be closed subspaces such that
\begin{enumerate}
\item\label{item:Orthogonal_projections_induce_embeddings}
For the vector spaces,
\[
  K_0 := \Ker dF(0) \quad\text{and}\quad C_0 := \Ker dF(0)^* = (\Ran dF(0))^\perp,
\]
% COMMENT (\Ran dF(0))^\perp is a closed subspace of Y even if \Ran dF(0) is not closed https://math.stackexchange.com/questions/4353383/is-the-orthogonal-complement-of-a-subspace-in-a-hilbert-space-always-complete
the orthogonal projections $\pi_{K_1}:X \to K_1$ and $\pi_{C_1}:Y \to C_1$ restrict to injective linear maps
\begin{equation}
  \label{eq:Orthogonal_projections_induce_embeddings}
  \pi_{K_1}:K_0 \to K_1 \quad\text{and}\quad \pi_{C_1}:C_0 \to C_1.
\end{equation}
\item\label{item:Isomorphism_tildeK0_perp_onto_tildeC0_perp}
For the vector spaces,
\[
  \tilde K_0 := \pi_{K_1}K_0 \subset K_1 \quad\text{and}\quad \tilde C_0 := \pi_{C_1}C_0 \subset C_1,
\]
there is an isomorphism of real vector spaces,
\begin{equation}
  \label{eq:Gamma_isomorphism_K1_cap_tildeK0_perp}
  \Gamma: K_1 \cap \tilde K_0^\perp \to C_1 \cap \tilde C_0^\perp. 
\end{equation}
\end{enumerate}
Then
\begin{equation}
\label{eq:Orthogonality_Relation_In_KModel}
K_1\cap \tilde K_0^\perp \subset K_0^\perp\quad\text{and}\quad
C_1\cap \tilde C_0^\perp \subset C_0^\perp,
\end{equation}
and the following hold. The subspaces $K_0$ and $K_1 \cap \tilde K_0^\perp$ are orthogonal in $X$, the subspaces $C_0$ and $C_1 \cap \tilde C_0^\perp$ are orthogonal in $Y$, and for
\begin{subequations}
  \label{eq:tildeK1_and_tildeC1_and_tildeX_and_tildeY}
  \begin{gather}
  \label{eq:tildeK1_and_tildeC1}
  \tilde K_1 := K_0 \oplus \left(K_1 \cap \tilde K_0^\perp\right) \subset X
  \quad\text{and}\quad
  \tilde C_1 := C_0 \oplus \left(C_1 \cap \tilde C_0^\perp\right) \subset Y,
  \\
  %PF11-13-2025 I don't think \tilde X is necessarily an orthogonal direct sum in X or that \tilde X is necessarily an orthogonal direct sum in Y, but please check
  %TL11-26-2025: If they were, then $K_0\subset \tilde K_1$ would be orthogonal to $K_1^\perp$, i.e. in $K_1$.  Sometimes that could be true but it's certainly not something we want to assume!
  \label{eq:tildeX_and_tildeY}
  \tilde X := \tilde K_1 \oplus K_1^\perp \subset X\oplus X
  \quad\text{and}\quad
  \tilde Y := \tilde C_1 \oplus C_1^\perp \subset Y\oplus Y,
\end{gather}
\end{subequations}
the canonical isomorphisms,
\[
  \tilde K_1 \cong K_1
  \quad\text{and}\quad
  \tilde C_1 \cong C_1,
\]  
and the real linear isomorphisms $\pi_{K_1}:K_0 \to \tilde K_0$ and $\pi_{C_1}:C_0 \to \tilde C_0$ yield real linear isomorphisms,
\begin{equation}
\label{eq:Xi_isomorphism_tildeX_and_tildeY}  
  \Xi_1:\tilde X \to X
  \quad\text{and}\quad
  \Xi_2:\tilde Y \to Y,
\end{equation}
such that
\begin{gather*}
  \Xi_1 \restriction K_0 = \pi_{K_1}
  \quad\text{and}\quad
  \Xi_1 \restriction K_1 \cap \tilde K_0^\perp = \pi_{K_1} = \id
  \quad\text{and}\quad
  \Xi_1 \restriction K_1^\perp = \id,
  \\
  \Xi_2 \restriction C_0 = \pi_{C_1}
  \quad\text{and}\quad
  \Xi_2 \restriction C_1 \cap \tilde C_0^\perp = \pi_{C_1} = \id
  \quad\text{and}\quad
  \Xi_2 \restriction C_1^\perp = \id.
\end{gather*}  
Moreover, for a possibly smaller open neighborhood $U' \subset U$ of the origin in $X$ and $\tilde U : = \Xi_1^{-1}(U)$ and $\tilde U' : = \Xi_1^{-1}(U')$, open neighborhoods of the origin in $\tilde X$, there exist a $C^k$ diffeomorphism $\tilde g:\tilde U' \to \tilde g(U') \subset U$ and a $C^k$ map $\tilde f:\tilde U' \to \tilde C_1$ such that the following hold:
\begin{enumerate}
\item\label{item:F_circ_tildeg_is_dF+pi_C0_circ_tildef}
  One has
  \begin{equation}
  \label{eq:F_circ_tildeg_is_dF+pi_C0_circ_tildef}
  F\circ \tilde g = dF(0)\circ\Xi_1 + \pi_{C_0}\circ\tilde f \quad\text{on } \tilde U',
\end{equation}
where $\pi_{C_0}:\tilde C_1 = C_0 \oplus (C_1 \cap \tilde C_0^\perp) \to C_0$ is the orthogonal projection, and
\[
  % TL12-17-2025: \tilde g here instad of g?
  %PF12-18-2025 Yes, so changed
  \tilde g(0) = 0, \quad d\tilde g(0) = \id_X, \quad \tilde f(0) = 0, \quad\text{and}\quad
  d\tilde f(0) = \Gamma\circ \pi_{K_1 \cap \tilde K_0^\perp},
\]
where $\pi_{K_1 \cap \tilde K_0^\perp}: \tilde X = K_0 \oplus (K_1 \cap \tilde K_0^\perp) \oplus K_1^\perp \to K_1 \cap \tilde K_0^\perp$ is the orthogonal projection.
\item\label{item:Equality_zero-locus_F_and_zero-locus_tildef}
  One has
\begin{equation}
  \label{eq:Equality_zero-locus_F_and_zero-locus_tildef}
  F^{-1}(0) \cap g(U') = \Xi_1\left(\tilde f^{-1}(0)\cap \tilde U'\cap\tilde K_1\right).
\end{equation}
\item 
\label{item:Equivariance_Of_Kuranishi_Model}
If $F$ is equivariant with respect to the orthogonal action of a compact Lie group $G$ on $X$ and $Y$, the subspaces $K_1 \subset X$ and $C_1 \subset Y$ are $G$-invariant, and the isomorphism $\Gamma$ is $G$-equivariant, then the maps $g$ and $\tilde f$ can be chosen to be $G$-equivariant.
\end{enumerate}
%PF7-15-2025 Change "Item" to "Condition/Hypothesis" or "Assertion" everywhere in proof bodies.
\end{cor}

\begin{proof}
We first prove the orthogonality relations \eqref{eq:Orthogonality_Relation_In_KModel}.
For any $v_0\in K_0$ and $w\in K_1\cap \tilde K_0^\perp$, we have
\[
\langle v_0,w\rangle_X
=
\langle \pi_{K_1} v_0,w\rangle_X
=
0,
\]
where the first equality holds because $w\in K_1$ and the second equality holds because $\pi_{K_1} v_0 \in \pi_{K_1}(K_0) = \tilde K_0$ and $w\in \tilde K_0^\perp$.  This proves that $K_0$ and $K_1\cap \tilde K_0^\perp$ are orthogonal. The proof that $C_0$ and $C_1\cap \tilde C_0^\perp$ are orthogonal is identical.

By Theorem \ref{thm:Kuranishi_model_defined_by_Fredholm_map_Banach_spaces} and Remark \ref{rmk:Kuranishi_model_defined_by_Fredholm_map_Banach_spaces_special_case}, there are an open neighborhood $U' \subset U$ of the origin $0 \in X$, a $C^k$ diffeomorphism $g:U' \to g(U') \subset U$, and a $C^k$ map $f:U' \to C_0$ such that
\begin{gather*}
  F\circ g = dF(0) + f \quad\text{on } U'
  \quad\text{(by \eqref{eq:F_circ_g_is_dF_proj+f} and $X_0=(\Ker dF(0))^\perp$)}
  \\
  g(0) = 0, \quad dg(0) = \id_X, \quad f(0) = 0, \quad\text{and}\quad df(0) = 0
  \quad\text{(by \eqref{eq:g(0)_and_dg(0)_and_f(0)_and_df(0)})}
  \\
  F^{-1}(0)\cap g(U') = f^{-1}(0)\cap U'\cap K_0
  \quad\text{(by \eqref{eq:Equality_zero-locus_F_and_zero-locus_fK})}.
\end{gather*}
Noting that
\[
  \tilde K_1 = K_0 \oplus \left(K_1 \cap \tilde K_0^\perp\right)
  \quad\text{and}\quad
  \tilde C_1 = C_0 \oplus \left(C_1 \cap \tilde C_0^\perp\right)
\]  
by hypothesis \eqref{eq:tildeK1_and_tildeC1}, we define a $C^k$ diffeomorphism $\tilde g:\tilde U' \to g(U') \subset U$ by
\[
  \tilde g := g \circ \Xi_1,
\]
and a $C^k$ map $\tilde f:\tilde U' \to \tilde C_1$ by
\[
  \tilde f : \tilde U' \ni \tilde x \mapsto (f\circ\Xi_1)(\tilde x)
  + \left(\Gamma\circ\pi_{K_1 \cap \tilde K_0^\perp}\right)(\tilde x) \in \tilde C_1.
\]
We now consider Assertion \eqref{item:F_circ_tildeg_is_dF+pi_C0_circ_tildef}. We first observe that
%PF7-15-2025 Add proof details!
\[
  F\circ \tilde g
  = F\circ g \circ \Xi_1
  = (dF(0) + f)\circ\Xi_1
  = dF(0)\circ\Xi_1 + f\circ\Xi_1
  =  dF(0)\circ\Xi_1 + \pi_{C_0}\circ\tilde f
  \quad\text{on } \tilde U'.
\]
Moreover, the properties of $f$ and $g$ and definitions of $\tilde g$ and $\tilde f$ ensure that
\[
  \tilde g(0) = 0, \quad d\tilde g(0) = \id_{\tilde X}, \quad \tilde f(0) = 0,
  \quad\text{and}\quad d\tilde f(0) = \Gamma\circ\pi_{K_1 \cap \tilde K_0^\perp},
\]
and this completes the verification of Assertion \eqref{item:F_circ_tildeg_is_dF+pi_C0_circ_tildef}. Similarly, Assertion \eqref{item:Equality_zero-locus_F_and_zero-locus_tildef} follows easily
%PF7-15-2025 Add details!
from  \eqref{eq:Equality_zero-locus_F_and_zero-locus_fK} and the definition of $\tilde f$.

Consider Assertion \eqref{item:Equivariance_Of_Kuranishi_Model}. If $F$ is $G$-equivariant, then Theorem \ref{thm:Kuranishi_model_defined_by_Fredholm_map_Banach_spaces} implies that the neighborhood $U'$ can be chosen to be $G$-invariant and the maps $g$ and $f$ can be chosen to be $G$-equivariant. By hypothesis, the isomorphism $\Gamma$ is $G$-equivariant, and thus $\tilde f$ is $G$-equivariant.
\end{proof}

%PF7-15-2025 Maybe include when corrected?
% \begin{rmk}[On the domain and codomain of the finite-dimensional zero-locus model produced by Corollary \ref{cor:Kuranishi_model_defined_by_Fredholm_map_Hilbert_spaces}]
% \label{rmk:Domain_and_codomain_Kuranishi_model_defined_by_Fredholm_map_Hilbert_spaces}
% The proof of Corollary \ref{cor:Kuranishi_model_defined_by_Fredholm_map_Hilbert_spaces} constructs a local Kuranishi model with a $C^k$ embedding $\tilde g = g\circ\Xi_1:\tilde U' \to U$, a $C^k$ map $\tilde f:\tilde U'\to\tilde C_1$, and zero locus $\tilde f^{-1}(0)\cap \tilde K_1$, where $K_0 \subset \tilde K_1 \subset \tilde X$. However, we could equally have constructed a local Kuranishi model with a $C^k$ embedding $\breve g = g\circ\Xi_1^{-1}:U' \to U$ and zero locus $\breve f^{-1}(0)\cap K_1$, where $\tilde K_0 \subset K_1 \subset X$.   
% \end{rmk} 

\begin{rmk}[On the hypotheses of Corollary \ref{cor:Kuranishi_model_defined_by_Fredholm_map_Hilbert_spaces}]
\label{rmk:Hypotheses_corollary_Kuranishi_model_defined_by_Fredholm_map_Hilbert_spaces}  
In our application, Condition \eqref{item:Orthogonal_projections_induce_embeddings} can be achieved with the aid of the forthcoming Lemma \ref{lem:Approximation_finite-dimensional_subspaces_Hilbert_space}, while Condition \eqref{item:Isomorphism_tildeK0_perp_onto_tildeC0_perp} can be achieved with a combination of index theory for the operator $dF(0)$ and spectral theory for its Laplacians.

For example, suppose that $X$ and $Y$ are \emph{complex} Hilbert spaces but $\sT := dF(0)$  is only assumed to be \emph{real linear}. Let $\sT' \in \Hom(X,Y)$ be the \emph{complex linear} part of $\sT$, assume that the \emph{complex antilinear} operator $\sT'' = \sT-\sT' \in \Hom(X,Y)$ is \emph{compact}, and that $\sT'$ obeys the hypotheses of the forthcoming Lemma \ref{lem:Eigenvalues_densely_defined_unbounded_linear_operators}. Let $\nu$ be a positive constant such that $\nu \notin \sigma(\sT^{\prime,*}\sT')$, where $\sigma(\sT^{\prime,*}\sT')$ denotes the spectrum of $\sT^{\prime,*}\sT'$, and choose
\[
  K_\nu := \Ran\Pi_{1,\nu}' \subset X \quad\text{and} \quad C_\nu := \Ran\Pi_{2,\nu}' \subset Y,
\]
where $\Pi_{1,\nu}'$ is the orthogonal projection onto the linear span of the eigenspaces of $\sT^{\prime,*}\sT'$ with eigenvalues less than $\nu$ and $\Pi_{2,\nu}'$ is the orthogonal projection onto the linear span of the eigenspaces of $\sT'\sT^{\prime,*}$ with eigenvalues less than $\nu$, recalling from Lemma \ref{lem:Eigenvalues_densely_defined_unbounded_linear_operators} that $\sigma(\sT^{\prime,*}\sT')\less\{0\} = \sigma(\sT'\sT^{\prime,*})\less\{0\}$.

Following the hypotheses of Corollary \ref{cor:Kuranishi_model_defined_by_Fredholm_map_Hilbert_spaces}, we denote
\[
  K_0 := \Ker\sT, \quad C_0 := \Ker\sT^*, \quad \tilde K_0 := \Pi_{1,\nu}'K_0 \subset K_\nu,
  \quad\text{and}\quad \tilde C_0 := \Pi_{2,\nu}'C_0 \subset C_\nu.
\]
As in Condition \eqref{item:Orthogonal_projections_induce_embeddings}, we assume that the orthogonal projections yield isomorphisms of real vector spaces, $\tilde K_0 \cong K_0$ and $\tilde C_0 \cong C_0$. However, this assumption is easily verified using the forthcoming Lemma \ref{lem:Approximation_finite-dimensional_subspaces_Hilbert_space} by choosing $\nu$ to be sufficiently large. We claim that there is an isomorphism of real vector spaces,
\[
  \Gamma: K_\nu \cap \tilde K_0^\perp \to C_\nu \cap \tilde C_0^\perp,
\]
so Condition \eqref{item:Isomorphism_tildeK0_perp_onto_tildeC0_perp} is satisfied. We have
\begin{align*}
  \dim K_\nu \cap \tilde K_0^\perp
  &=
    \dim K_\nu - \dim \tilde K_0 = \dim K_\nu - \dim\Ker\sT,
  \\
  \dim C_\nu \cap \tilde C_0^\perp
  &=
    \dim C_\nu - \dim \tilde C_0 = \dim C_\nu - \dim\Ker\sT^*,
\end{align*}
and thus
\begin{align*}
  \dim K_\nu \cap \tilde K_0^\perp
  - \dim C_\nu \cap \tilde C_0^\perp
  &=
    \dim K_\nu - \dim C_\nu - \dim\Ker\sT + \dim\Ker\sT^*
  \\
  &=
     \dim K_\nu - \dim C_\nu - \Ind\sT. 
\end{align*}
By Lemma \ref{lem:Eigenvalues_densely_defined_unbounded_linear_operators}, the operator $\sT'$ yields an isomorphism of complex vector spaces,
\[
  \sT':K_\nu \cap (\Ker\sT')^\perp \to C_\nu \cap (\Ker\sT^{\prime,*})^\perp,
\]
and thus an equality
\[
  \dim K_\nu \cap (\Ker\sT')^\perp = \dim C_\nu \cap (\Ker\sT^{\prime,*})^\perp.
\]
In particular, since $\Ker\sT' = \Ker\sT^{\prime,*}\sT'$ and $\Ker\sT^{\prime,*} = \Ker\sT'\sT^{\prime,*}$, we have
\[
  \dim K_\nu - \dim C_\nu = \dim \Ker\sT' - \dim \Ker\sT^{\prime,*} = \Ind\sT'.
\]  
But $\Ind\sT = \Ind\sT'$, since $\sT'' = \sT - \sT'$ is compact, and hence
\[
  \dim K_\nu \cap \tilde K_0^\perp
  -
  \dim C_\nu \cap \tilde C_0^\perp
  = 0,
\]
so the desired isomorphism $\Gamma$ exists. We expect that $\Pi_{C_\nu}\sT$, where $\Pi_{C_\nu}$ is the orthogonal projection from $Y$ onto $C_\nu$, would be a natural candidate for the abstract real linear isomorphism $\Gamma$ and this should follow when $\|\sT''\|_{\Hom(X,Y)}$ is sufficiently small.
% PF7-3-2025 Add a proof at some point even if we don't rely on this fact.
\qed\end{rmk}

\begin{rmk}[Equivariance of the real linear isomorphism $\Gamma$ in Corollary \ref{cor:Kuranishi_model_defined_by_Fredholm_map_Hilbert_spaces}]
\label{rmk:Equivariance_of_Gamma_Kuranishi_model_defined_by_Fredholm_map_Hilbert_spaces}
In the forthcoming Lemma \ref{lem:SymplecticStrutureOnStabilizedSpace}, we will combine
the argument given above  for the existence of the isomorphism $\Gamma$ with the equivariant index theory discussed in Section \ref{sec:S1EquivIndexBundles} to show that if $\sT$, $\sT'$, and the inner products on $X$ and $Y$ are $G$-equivariant, then $\Gamma$ can be chosen to be $G$-equivariant as well.
The existence of a $G$-equivariant isomorphism is a strictly stronger requirement than the existence of the isomorphism $\Gamma$ in \eqref{eq:Gamma_isomorphism_K1_cap_tildeK0_perp} because for a group $G$ with non-trivial representations, there are representations of $G$ which are isomorphic as vector spaces but not as representations.
By considering the example of the identity map $\id_\RR:\RR\to \RR$ where the domain has the non-trivial action of $G=\ZZ/2\ZZ$ and the codomain has the trivial action, we see that one cannot produce an equivariant isomorphism from an isomorphism by averaging with respect to Haar measure as is often done to produce equivariant objects such as morphisms or metrics.  If we write $\RR\ni v\mapsto n(v)=-v$, then averaging $\id_\RR$ with respect to the Haar measure would give
\[
\frac{1}{2}\left( \id_\RR+ \id_R\circ n\right)=0,
\]
which is not an isomorphism.
\qed\end{rmk}

\begin{rmk}[Construction of the real linear isomorphism $\Gamma$ in Corollary \ref{cor:Kuranishi_model_defined_by_Fredholm_map_Hilbert_spaces}]
\label{rmk:Construction_isomorphism_Gamma_Kuranishi_model_defined_by_Fredholm_map_Hilbert_spaces}  
We continue the notation of Remark \ref{rmk:Hypotheses_corollary_Kuranishi_model_defined_by_Fredholm_map_Hilbert_spaces} and give an explicit construction of the isomorphism $\Gamma$ of \eqref{eq:Gamma_isomorphism_K1_cap_tildeK0_perp} in Corollary \ref{cor:Kuranishi_model_defined_by_Fredholm_map_Hilbert_spaces}. We consider different cases for ease of exposition.

\setcounter{case}{0}
\begin{case}[$\Ker \sT^* = (0)$ and $\Ker\sT^{\prime,*} = (0)$]
\label{case:Ker_sT*_zero_and_Ker_sTprime*_zero}  
For this situation, $\dim\Ker \sT = \dim\Ker\sT'$ since $\Ind\sT = \Ind\sT'$. Writing $K_0' := \Ker\sT'$, we choose $S \in \Or(K_\nu)$ such that $S(\tilde K_0) = K_0'$ and hence $S(\tilde K_0^\perp) = (K_0')^\perp$. Because $\Ker\sT^{\prime,*} = (\Ran\sT')^\perp = (0)$, the operator
\begin{equation}
  \label{eq:sTprime_isomorphism_K_nu_cap_K0prime_perp_to_C_nu}
  \sT':K_\nu \cap (K_0')^\perp \to C_\nu
\end{equation}
is an isomorphism of complex vector spaces, so the composition
\begin{equation}
  \label{eq:Gamma_Ker_sT*_and_Ker_sTprime*_zero}
  \Gamma \equiv \sT' \circ S: K_\nu \cap \tilde K_0^\perp \to C_\nu,
\end{equation}
gives the desired isomorphism for this case.

Note that because $K_\nu\cap (K_0')^\perp \subset K_\nu$ has an almost complex structure, the isomorphism $S:K_\nu\cap \tilde K_0^\perp \to K_\nu\cap (K_0')^\perp$ induces an almost complex structure on $K_\nu\cap \tilde K_0^\perp$ and because $\sT'$ in \eqref{eq:sTprime_isomorphism_K_nu_cap_K0prime_perp_to_C_nu} is complex linear, we see that $\Gamma$ in \eqref{eq:Gamma_Ker_sT*_and_Ker_sTprime*_zero} is complex linear with respect to the induced almost complex structures on $K_\nu\cap \tilde K_0^\perp$ and $C_\nu$.
%PF7-9-2025 Add complex linearity to remaining cases in this remark
\end{case}

\begin{case}[$\dim\Ker \sT^* = \dim \Ker\sT^{\prime,*}$]
\label{case:dim_Ker_sT*_is_dim_Ker_sTprime*}    
This case allows $C_0 \neq (0)$ and, writing $C_0' := \Ker\sT^{\prime,*}$, our assumption is that $\dim C_0' = \dim C_0$, so that $\dim K_0 = \dim K_0'$ since  $\Ind\sT = \Ind\sT'$. In addition to choosing $S \in \Or(K_\nu)$ as in Case \ref{case:Ker_sT*_zero_and_Ker_sTprime*_zero}, we choose $R \in \Or(C_\nu)$ such that $R(C_0') = \tilde C_0$ and hence $R((C_0')^\perp) = \tilde C_0^\perp$. Because the operator
\begin{equation}
  \label{eq:sTprime_isomorphism_K_nu_cap_K0prime_perp_to_C_nu_cap_C0prime_perp}
  \sT':K_\nu \cap (K_0')^\perp \to C_\nu \cap (C_0')^\perp
\end{equation}
is an isomorphism of complex vector spaces, the composition
\begin{equation}
  \label{eq:Gamma_dim_Ker_sT*_is_dim_Ker_sTprime*}
  \Gamma \equiv R\circ \sT' \circ S: K_\nu \cap \tilde K_0^\perp \to C_\nu \cap \tilde C_0^\perp,
\end{equation}
gives the desired isomorphism for this case.
\end{case}

\begin{case}[$\dim\Ker \sT^* \leq \dim \Ker\sT^{\prime,*}$]
\label{case:dim_Ker_sT*_leq_dim_Ker_sTprime*}  
We have $\dim C_0' \geq \dim C_0$ by assumption, so that $\dim K_0' \geq \dim K_0$ since  $\Ind\sT = \Ind\sT'$. We choose $S \in \Or(K_\nu)$ such that $S(\tilde K_0) \subset K_0'$ and hence $S(\tilde K_0^\perp) \supset (K_0')^\perp$ and $R \in \Or(C_\nu)$ such that $R(C_0') \supset \tilde C_0$ and hence $R((C_0')^\perp) \subset \tilde C_0^\perp$. We write
\begin{align*}
  S\left(K_\nu \cap \tilde K_0^\perp\right)
  &= K_\nu \cap S\left(\tilde K_0^\perp\right)
  \\
  &= K_0'\cap S\left(\tilde K_0^\perp\right)\oplus (K_0')^\perp \cap S\left(\tilde K_0^\perp\right)
  \\
  &  =  K_0' \cap S\left(\tilde K_0^\perp\right) \oplus K_\nu\cap (K_0')^\perp
    \quad\text{(because $S(\tilde K_0^\perp) \supset (K_0')^\perp$),}
\end{align*}
and, using $R((C_0')^\perp)^\perp \subset \tilde C_0^\perp$,
\[
  R\left(C_\nu \cap (C_0')^\perp\right)
  = C_\nu \cap R\left((C_0')^\perp\right)
  \subset
  C_\nu \cap \tilde C_0^\perp.
\]
Therefore,
\[
  (\sT' \circ S)\left(K_\nu \cap \tilde K_0^\perp\right)
  =
  \sT'\left( K_0' \cap S\left(\tilde K_0^\perp\right) \oplus K_\nu\cap (K_0')^\perp\right)
  =
  \sT'\left(K_\nu\cap (K_0')^\perp\right)
  =
  C_\nu \cap (C_0')^\perp,
\]
and thus
\[
  (R\circ \sT' \circ S)\left(K_\nu \cap \tilde K_0^\perp\right)
  =
  R\left(C_\nu \cap (C_0')^\perp\right)
  \subset
  C_\nu \cap \tilde C_0^\perp.
\]
We claim that
\[
  \dim_\RR K_0' \cap S(\tilde K_0^\perp)
  =
  \dim_\RR \tilde C_0^\perp \cap C_0'.
\]  
Because
\[
  K_\nu=\tilde K_0\oplus \tilde K_0^\perp=S(\tilde K_0)\oplus S(\tilde K_0^\perp),
\]
we obtain
\[
  K_0'=K_0'\cap S(\tilde K_0)\oplus K_0'\cap S(\tilde K_0^\perp).
\]
Using $S(\tilde K_0)\subset K_0'$, the preceding equality gives
\[
  K_0'=S(\tilde K_0)\oplus K_0'\cap S(\tilde K_0^\perp).
\]
Thus,
\begin{multline*}
  \dim_\RR  K_0' \cap S(\tilde K_0^\perp)
    = \dim_\RR K_0' - \dim S(\tilde K_0)
    = \dim_\RR K_0' - \dim \tilde K_0
    \\
    = \Ind_\RR \sT' + \dim C_0'  - \Ind_\RR \sT - \dim \tilde C_0
    = \dim_\RR C_0' - \dim_\RR \tilde C_0
    = \dim_\RR \tilde C_0^\perp \cap C_0',
\end{multline*}
as claimed. Hence, for a choice of isomorphism of real vector spaces,
\begin{equation}
  \label{eq:Gamma0_isomorphism_K0prime_cap_S(tildeK0perp)}
  \Gamma_0: K_0' \cap S(\tilde K_0^\perp)
  \to
  \tilde C_0^\perp \cap C_0',
\end{equation}
and the isomorphism of complex vector spaces $\sT'$ in \eqref{eq:sTprime_isomorphism_K_nu_cap_K0prime_perp_to_C_nu_cap_C0prime_perp}, the composition
\[
  \Gamma \equiv R\circ (\sT' \oplus \Gamma_0) \circ S: K_\nu \cap \tilde K_0^\perp \to C_\nu \cap \tilde C_0^\perp,
\]
gives the desired isomorphism for this case.
\end{case}

\begin{case}[$\dim\Ker \sT^* \geq \dim \Ker\sT^{\prime,*}$]
\label{case:dim_Ker_sT*_geq_dim_Ker_sTprime*}  
We have $\dim C_0' \leq \dim C_0$ by assumption, so that $\dim K_0' \leq \dim K_0$ since  $\Ind\sT = \Ind\sT'$. We choose $S \in \Or(K_\nu)$ such that $S(\tilde K_0) \supset K_0'$ and hence $S(\tilde K_0^\perp) \subset (K_0')^\perp$ and $R \in \Or(C_\nu)$ such that $R(C_0') \subset \tilde C_0$ and hence $R((C_0')^\perp) \supset \tilde C_0^\perp$. We write
\[
  R\left(C_\nu \cap (C_0')^\perp\right)
  = C_\nu \cap R\left((C_0')^\perp\right)
  =
  C_\nu \cap \tilde C_0^\perp \oplus \tilde C_0 \cap R((C_0')^\perp).
\]
and, using $(S(\tilde K_0^\perp))^\perp = S(\tilde K_0)$,
\[
  S\left(K_\nu \cap \tilde K_0^\perp\right)
  = K_\nu \cap S\left(\tilde K_0^\perp\right)
  \subset
  K_\nu \cap S\left(\tilde K_0^\perp\right) \oplus (K_0')^\perp \cap S(\tilde K_0)
  =
  K_\nu \cap (K_0')^\perp.
\]
We choose an isomorphism of real vector spaces,
\begin{equation}
  \label{eq:Gamma1_isomorphism_K0prime_cap_S(tildeK0)}
  \Gamma_1: (K_0')^\perp \cap S(\tilde K_0)
  \to
  \tilde C_0 \cap R((C_0')^\perp).
\end{equation}
Hence, there is an isomorphism
\[
  \Gamma : K_\nu \cap \tilde K_0^\perp \to C_\nu \cap \tilde C_0^\perp
\]
such that, for the isomorphism of complex vector spaces $\sT'$ in \eqref{eq:sTprime_isomorphism_K_nu_cap_K0prime_perp_to_C_nu_cap_C0prime_perp}, the composition
%PF7-8-2025 Check S factor below
\[
  R\circ \sT'\circ S = \left(\Gamma\circ S^{-1}\right) \oplus \Gamma_1:
  K_\nu \cap (K_0')^\perp
  \to C_\nu \cap R\left((C_0')^\perp\right),
\]
is an isomorphism, as desired for this case.
\end{case}
This completes our construction of the isomorphism $\Gamma$ of \eqref{eq:Gamma_isomorphism_K1_cap_tildeK0_perp} in Corollary \ref{cor:Kuranishi_model_defined_by_Fredholm_map_Hilbert_spaces} for all cases.
\qed
\end{rmk}

\begin{rmk}[Dimension counting in Corollary \ref{cor:Kuranishi_model_defined_by_Fredholm_map_Hilbert_spaces}]
\label{rmk:Dimension_counting_Kuranishi_model_defined_by_Fredholm_map_Hilbert_spaces}
In our application of Corollary \ref{cor:Kuranishi_model_defined_by_Fredholm_map_Hilbert_spaces}, we will always have the situation described in Remark \ref{rmk:Construction_isomorphism_Gamma_Kuranishi_model_defined_by_Fredholm_map_Hilbert_spaces}, so $\Ind\sT = \Ind\sT' \equiv 0 \pmod{2}$ and $\Ind\sT > 0$, so that $\dim\Ker\sT > 0$. We will moreover know that $\dim\Ker\sT^* \equiv 0 \pmod{2}$, and thus $\dim\Ker\sT \equiv 0 \pmod{2}$.

Consequently, the real vector spaces $K_0$, $\tilde K_0$, $C_0$, and $\tilde C_0$ described in Remark \ref{rmk:Construction_isomorphism_Gamma_Kuranishi_model_defined_by_Fredholm_map_Hilbert_spaces} will be even-dimensional and because the vector spaces $K_\nu$ and $C_\nu$ are complex, the complements $K_\nu\cap \tilde K_0^\perp$ and $C_\nu\cap \tilde C_0^\perp$ will be even-dimensional too. In the third case considered in Remark \ref{rmk:Construction_isomorphism_Gamma_Kuranishi_model_defined_by_Fredholm_map_Hilbert_spaces}, it follows from the orthogonal decompositions used in the analysis of that case that the vector spaces $K_0' \cap S\left(\tilde K_0^\perp\right)$ and $\tilde C_0^\perp \cap R(C_0')$ comprising the domain and codomain of the real linear isomorphism $\Gamma_0$ in \eqref{eq:Gamma0_isomorphism_K0prime_cap_S(tildeK0perp)} are even-dimensional. Similarly, in the fourth case, the vector spaces $(K_0')^\perp \cap S(\tilde K_0)$ and $\tilde C_0 \cap R((C_0')^\perp)$ comprising the domain and codomain of the real linear isomorphism $\Gamma_1$ in \eqref{eq:Gamma1_isomorphism_K0prime_cap_S(tildeK0)} are even-dimensional. 
\qed\end{rmk}

\begin{rmk}[Simplified construction and complex linearity of the isomorphism $\Gamma$ in Corollary \ref{cor:Kuranishi_model_defined_by_Fredholm_map_Hilbert_spaces} when $\Ind\sT$ and $\dim\Ker\sT^*$ are even dimensional]
\label{rmk:Simplified_construction_isomorphism_Gamma_Kuranishi_model_defined_by_Fredholm_map_Hilbert_spaces} 
Continue the notation of Remark \ref{rmk:Construction_isomorphism_Gamma_Kuranishi_model_defined_by_Fredholm_map_Hilbert_spaces} and  
assume further that $\Ind\sT$ and $\dim\Ker\sT^*$ are even dimensional as in Remark \ref{rmk:Dimension_counting_Kuranishi_model_defined_by_Fredholm_map_Hilbert_spaces}. We shall give a simpler proof in this setting of existence of the isomorphism $\Gamma$ in Corollary \ref{cor:Kuranishi_model_defined_by_Fredholm_map_Hilbert_spaces} by reducing Cases \ref{case:dim_Ker_sT*_leq_dim_Ker_sTprime*} and \ref{case:dim_Ker_sT*_geq_dim_Ker_sTprime*} in Remark \ref{rmk:Construction_isomorphism_Gamma_Kuranishi_model_defined_by_Fredholm_map_Hilbert_spaces} to Case \ref{case:dim_Ker_sT*_is_dim_Ker_sTprime*}.

For Case \ref{case:dim_Ker_sT*_leq_dim_Ker_sTprime*}, we assumed that $\dim_\RR\Ker\sT^* \leq \dim_\RR\Ker\sT^{\prime,*}$ and hence $\dim_\RR\Ker\sT \leq \dim_\RR\Ker\sT'$ since $\Ind\sT = \Ind\sT'$. We choose a complex linear subspace $L_0 \subset \Ker\sT'$ such that $\dim_\RR L_0 = \dim_\RR\Ker\sT' - \dim_\RR\Ker\sT$ and a complex linear subspace $M_0 \subset \Ker\sT^{\prime,*}$ such that $\dim_\RR L_0 = \dim_\RR M_0$. By mapping each zero eigenvector for $\sT^{\prime,*}\sT'$ in an ordered orthonormal basis for $L_0$ onto the corresponding zero eigenvector for $\sT'\sT^{\prime,*}$ in an ordered orthonormal basis for $M_0$, we construct a finite rank operator $\varpi_{L_0,L_1} \in \Hom(X,Y)$ with the property that 
\[
  \varpi_{L_0,L_1}:L_0 \to M_0
  \text{ is an isometric isomorphism and }
  \varpi_{L_0,M_0}\left(L_0^\perp\right) = 0.
\]  
We now define a complex linear Fredholm operator,
\[
  \hat\sT' := \sT' + \varpi_{L_0,M_0} \in \Hom(X,Y),
\]
and observe that, using $\sT'(L_0) = 0$,
\[
  \hat\sT'
  =
  \begin{cases}
    0 &\text{on } \Ker\sT'\cap L_0^\perp,
    \\
    \varpi_{L_0,M_0} &\text{on } \Ker\sT'\cap L_0,
    \\
    \sT' &\text{on } (\Ker\sT')^\perp.
  \end{cases}
\]
Noting that $\varpi_{L_0,M_0}(L_0) = M_0$, we obtain
\[
  \dim_\RR\Ker\hat\sT' = \dim_\RR\Ker\sT' - \dim_\RR L_0 = \dim_\RR\Ker\sT.
\]
Because $\hat\sT'$ is a perturbation of $\sT'$ by a finite rank and thus compact operator $\varpi_{L_0,L_1}$, we have $\Ind\hat\sT' = \Ind\sT' = \Ind\sT$ and so
\[
  \dim_\RR\Ker\hat\sT^{\prime,*} = \dim_\RR\Ker\sT^*.
\]
% TL11-8-2025: Is the following true?  I see $\hat\sT^{\prime,*}\hat\sT'=\sT^{\prime,*}\sT'+\id_{L_0}$ so $K_\nu$ and hence the projection operator would change when $\nu$ changes from $<1$ to $>1$.
%PF12-18-2025 I will check.
The subspaces $\tilde K_0 \subset K_\nu$ and $\tilde C_0 \subset C_\nu$ are unchanged when defined by the spectral projections $\hat\Pi_{1,\nu}'$ and $\hat\Pi_{2,\nu}'$ associated with the Laplacians $\hat\sT^{\prime,*}\hat\sT'$ and $\hat\sT'\hat\sT^{\prime,*}$, respectively. We denote $\hat K_0' := \Ker\hat\sT'$ and $\hat C_0' := \Ker\hat\sT^{\prime,*}$ and choose $S \in \Or(K_\nu)$ and $R \in \Or(C_\nu)$ such that $S(\tilde K_0) = \hat K_0'$ and $R(\hat C_0') = \tilde C_0$. The operator
\begin{equation}
  \label{eq:tilde_sTprime_isomorphism_K_nu_cap_K0prime_perp_to_C_nu_cap_C0prime_perp}
  \hat\sT':K_\nu \cap (\hat K_0')^\perp \to C_\nu \cap (\hat K_0')^\perp
\end{equation}
is an isomorphism of complex vector spaces, so the composition
\begin{equation}
  \label{eq:Gamma_tilde_sTprime_dim_Ker_sT*_is_dim_Ker_sTprime*}
  \Gamma \equiv R\circ \hat\sT' \circ S: K_\nu \cap \tilde K_0^\perp \to C_\nu \cap \tilde C_0^\perp,
\end{equation}
gives the desired isomorphism for this case.

Note that because $K_\nu \cap (\hat K_0')^\perp \subset K_\nu$ has an almost complex structure, the isomorphism $S:K_\nu \cap \tilde K_0^\perp \to K_\nu \cap (\hat K_0')^\perp$ induces an almost complex structure on $\tilde K_0^\perp$. Similarly, because $C_\nu\cap (\hat C_0')^\perp \subset C_\nu$ has an almost complex structure, the isomorphism $R:C_\nu\cap (\hat C_0')^\perp \to C_\nu\cap \tilde C_0^\perp$ induces an almost complex structure on $C_\nu\cap \tilde C_0^\perp$. Since $\hat\sT'$ in \eqref{eq:tilde_sTprime_isomorphism_K_nu_cap_K0prime_perp_to_C_nu_cap_C0prime_perp} is complex linear, we see that $\Gamma$ in \eqref{eq:Gamma_tilde_sTprime_dim_Ker_sT*_is_dim_Ker_sTprime*} is complex linear with respect to the induced almost complex structures on $C_\nu\cap \tilde K_0^\perp$ and $C_\nu\cap \tilde C_0^\perp$.

For Case \ref{case:dim_Ker_sT*_geq_dim_Ker_sTprime*}, we assumed that $\dim_\RR\Ker\sT^* \geq \dim_\RR\Ker\sT^{\prime,*}$ and hence $\dim_\RR\Ker\sT \geq \dim_\RR\Ker\sT'$ since $\Ind\sT = \Ind\sT'$. We choose a complex linear subspace $L_1 \subset (\Ker\sT')^\perp$ that is a complex linear span of eigenvectors of $\sT^{\prime,*}\sT'$ with positive eigenvalues such that $\dim_\RR L_1  = \dim_\RR\Ker\sT - \dim_\RR\Ker\sT'$. We now define a complex linear Fredholm operator,
\[
  \hat\sT' := \sT'\circ(\id - \pi_{L_1}),
\]
where $\pi_{L_1} \in \End(X)$ is the orthogonal projection from $X$ onto $L_1$, and observe that, using $\pi_{L_1} = 0$ on $L_1^\perp$,
\[
  \hat\sT'
  =
  \begin{cases}
    0 &\text{on } \Ker\sT',
    \\
    0 &\text{on } (\Ker\sT')^\perp \cap L_1,
    \\
    \sT' &\text{on } (\Ker\sT')^\perp \cap L_1^\perp,
  \end{cases}  
\]  
we see that
\[
  \dim_\RR\Ker\hat\sT' = \dim_\RR\Ker\sT' + \dim_\RR L_1 = \dim_\RR\Ker\sT.
\]
Because $\hat\sT'$ is a perturbation of $\sT'$ by a finite rank and thus compact operator $\sT'\circ\pi_{L_1}$,
\[
  \hat\sT' = \sT' - \sT'\circ\pi_{L_1},
\]
we have $\Ind\hat\sT' = \Ind\sT' = \Ind\sT$ and so
\[
  \dim_\RR\Ker\hat\sT^{\prime,*} = \dim_\RR\Ker\sT.
\]
Provided $\nu$ is greater than the largest eigenvalue associated with $L_1$, the subspaces $\tilde K_0$ and $\tilde C_0$ are unchanged when defined by the spectral projections $\hat\Pi_{1,\nu}'$ and $\hat\Pi_{2,\nu}'$ associated with the Laplacians $\hat\sT^{\prime,*}\hat\sT'$ and $\hat\sT'\hat\sT^{\prime,*}$, respectively. The remainder of the construction of $\Gamma$ in \eqref{eq:Gamma_tilde_sTprime_dim_Ker_sT*_is_dim_Ker_sTprime*} and verification of complex linearity with respect to the induced almost complex structures are the same as in the previous case.
\qed
\end{rmk}

In the construction of the Fredholm operators $\hat\sT'$ in Remark \ref{rmk:Simplified_construction_isomorphism_Gamma_Kuranishi_model_defined_by_Fredholm_map_Hilbert_spaces} by precomposition with $\id + \varpi_{L_0,L_1}$ or $\id - \pi_{L_1}$, the latter finite-rank operators may be viewed as analogues of the \emph{raising} or \emph{lowering operators} in representation theory and quantum mechanics where, in our application, the multiplicity of the $0$ eigenvalue of the Laplacian $\sT^{\prime,*}\sT'$ is lowered or raised and multiplicities of some its positive eigenvalues are lowered or raised.

%PF7-11-2025 I followed Salamon's notation in using f for the Kuranishi obstruction map  and g for the local diffeomorphism, but maybe that causes confusion wiuth the moment map f and metric g.

\begin{rmk}[Conventions for almost Hermitian structures]
\label{rmk:Conventions_almost_Hermitian_structures}
We use the following elementary example to specify our conventions for almost Hermitian structures. We write the standard Hermitian metric on $\CC$ as $h(z,w) := z\bar w$. The real linear isomorphism $\Upsilon:\CC \ni z = x_1+ix_2 \mapsto x = (x_1,x_2) \in \RR^2$ induces an almost complex structure $J$ on $\RR^2$ by $Jx = \Upsilon (i\Upsilon^{-1}(x))$, so that $x = (x_1,x_2) \mapsto x_1+ix_2 \mapsto -x_2 + ix_1 \mapsto (-x_2,x_1) = J(x_1,x_2) = Jx$. Thus,
\[
   J = \begin{pmatrix} 0 & -1 \\ 1 & 0 \end{pmatrix}
\]
with respect to the standard basis on $\RR^2$. The standard Riemannian metric on $\RR^2$ is given by $g(x, y) = x_1y_1 + x_2y_2$. Writing $z = \Upsilon^{-1}(x)$ and $w = \Upsilon^{-1}(y)$, we obtain
\[
  \Real h(z,w) = \Real(z\bar w) = \Real\left((x_1y_1 + x_2y_2) + i(x_2y_1 - x_1y_2)\right) = x_1y_1 + x_2y_2
  = g(x, y).
\]
According to the definition \eqref{eq:Fundamental_two-form} of the fundamental two-form $\omega$ on $\RR^2$, we have
\[
  \omega(x,y) = g(Jx, y) = g((-x_2,x_1),(y_1,y_2)) = x_1y_2 - x_2y_1, 
\]
while
\[
  \Imag h(z,w) = \Imag(z\bar w) = i(x_2y_1 - x_1y_2).
\]
Thus, $\omega(x,y) = i\Imag h(z,w)$ according to our conventions. Note that the relations among $g$, $h$, and $J$ differ\footnote{Kobayashi writes $g(x,y) = 2\Real h(z,w)$ and $\omega(x,y) = g(x,Jy) = 2i\Imag h(z,w)$.} from those given by Kobayashi \cite[Section 7.6, Equations (7.6.5), (7.6.7), and (7.6.8), p. 251]{Kobayashi_differential_geometry_complex_vector_bundles}, but they agree with those of Huybrechts \cite[Definition 1.2.13, p. 29, and Lemmas 1.2.15 and 1.2.17, p. 30]{Huybrechts_2005}.
\qed
\end{rmk}  

\begin{rmk}[Almost Hermitian structures on the subspaces $K_1$ and $C_1$ in Corollary \ref{cor:Kuranishi_model_defined_by_Fredholm_map_Hilbert_spaces}]
\label{rmk:AC_and_symplectic_structures_Kuranishi_model_defined_by_Fredholm_map_Hilbert_spaces}
We continue the assumptions and notation of Remarks \ref{rmk:Hypotheses_corollary_Kuranishi_model_defined_by_Fredholm_map_Hilbert_spaces}, \ref{rmk:Construction_isomorphism_Gamma_Kuranishi_model_defined_by_Fredholm_map_Hilbert_spaces}, and \ref{rmk:Simplified_construction_isomorphism_Gamma_Kuranishi_model_defined_by_Fredholm_map_Hilbert_spaces}. We now denote the subspaces $K_1$ and $C_1$ in Corollary \ref{cor:Kuranishi_model_defined_by_Fredholm_map_Hilbert_spaces} by $K_\nu$ and $C_\nu$, respectively, as in Remark \ref{rmk:Hypotheses_corollary_Kuranishi_model_defined_by_Fredholm_map_Hilbert_spaces}. Similarly, the real linear isomorphisms $\pi_{K_1}:K_0 \to \tilde K_0$ and $\pi_{C_1}:C_0 \to \tilde C_0$  in Corollary \ref{cor:Kuranishi_model_defined_by_Fredholm_map_Hilbert_spaces} induced by the orthogonal projections $\pi_{K_1}:X \to K_0$ and $\pi_{C_1}:Y \to C_1$ are now denoted by real linear isomorphisms $\Pi_{1,\nu}':K_0 \to \tilde K_0$ and $\Pi_{2,\nu}':C_0 \to \tilde C_0$ that are induced by the orthogonal projections $\Pi_{1,\nu}':X \to K_\nu$ and $\Pi_{2,\nu}':Y \to C_\nu$, respectively.

From Remark \ref{rmk:Simplified_construction_isomorphism_Gamma_Kuranishi_model_defined_by_Fredholm_map_Hilbert_spaces}, we know that the isomorphisms $S:\tilde K_0 \to \hat K_0'$ and $S:K_\nu\cap \tilde K_0^\perp \to K_\nu\cap (\hat K_0')^\perp$ induce almost complex structures on $\tilde K_0$ and $K_\nu\cap \tilde K_0^\perp$, respectively. Similarly, the isomorphisms $R: \hat C_0' \to \tilde C_0$ and $R:C_\nu\cap (\hat C_0')^\perp \to C_\nu\cap \tilde C_0^\perp$ induce almost complex structures on $\tilde C_0$ and $C_\nu\cap \tilde C_0^\perp$, respectively.

As in the assertion \eqref{eq:Xi_isomorphism_tildeX_and_tildeY} of Corollary \ref{cor:Kuranishi_model_defined_by_Fredholm_map_Hilbert_spaces}, the real linear isomorphisms $\Pi_{1,\nu}':K_0 \to \tilde K_0$ and $\Pi_{2,\nu}':C_0 \to \tilde C_0$ (produced by Lemma \ref{lem:Approximation_finite-dimensional_subspaces_Hilbert_space} for large enough $\nu$) extend to real linear isomorphisms:
\begin{subequations}
\label{eq:Xi_isomorphism_tildeKnu_and_tildeCnu}  
\begin{align}
  \label{eq:Xi1_isomorphism_tildeKnu}
  \Xi_{1,\nu}:\tilde K_\nu
  := K_0 \oplus \left(K_\nu \cap \tilde K_0^\perp\right)
  &\to
  \tilde K_0 \oplus \left(K_\nu \cap \tilde K_0^\perp\right) \cong K_\nu,
  \\
  \label{eq:Xi2_isomorphism_tildeCnu}
  \Xi_{2,\nu}:\tilde C_\nu
  := C_0 \oplus \left(C_\nu \cap \tilde C_0^\perp\right)
  &\to
  \tilde C_0 \oplus \left(C_\nu \cap \tilde C_0^\perp\right) \cong C_\nu.
\end{align}
\end{subequations}
The almost Hermitian structures on $K_\nu$ and $C_\nu$ (or their direct summands) thus induce almost Hermitian structures on $\tilde K_\nu$ and $\tilde C_\nu$, respectively. For example, if $J_1$ and $J_2$ are the almost complex structures on the complex Hilbert spaces $X$ and $Y$, respectively, then the pullbacks of their restrictions to $K_\nu$ and $C_\nu$,
\[
  \Xi_{1,\nu}^{-1} \circ J_1 \circ \Xi_{1,\nu} \in \End(\tilde K_\nu)
  \quad\text{and}\quad
  \Xi_{2,\nu}^{-1} \circ J_2 \circ \Xi_{2,\nu} \in \End(\tilde C_\nu),
\]
define almost Hermitian structures on $\tilde K_\nu$ and $\tilde C_\nu$, respectively, where $\Xi_{1,\nu}^{-1}:K_\nu \to \tilde K_\nu$ and $\Xi_{2,\nu}^{-1}:C_\nu \to \tilde C_\nu$ are the inverses of real linear isomorphisms obtained via restriction to the subspaces $K_\nu$ and $C_\nu$, respectively.
%PF7-15-2025 Harmonize with Tom's notation in chapter 11

Alternatively, suppose that $\tilde\Pi_{1,\nu}:X\to \tilde K_\nu$ and $\tilde\Pi_{2,\nu}:Y\to \tilde C_\nu$ are the orthogonal projections. We can consider the operators
\[
  \tilde\Pi_{1,\nu}\circ J_1 \in \End(\tilde K_\nu)
  \quad\text{and}\quad
  \tilde\Pi_{2,\nu}\circ J_2 \in \End(\tilde C_\nu),
\]
where $J_1 \in \Or(X)$ and $J_2 \in \Or(Y)$ denote the almost complex structures on $X$ and $Y$, respectively.
%induced by $J_1\in\Or(X)$ and $J_2\in\Or(Y)$.
By the proof of the
%PF11-13-2025 The notation is different there, with use of "n" enumerating eigenvectors rather than $\nu$ separating eigenvalues
%TL12-14-2025: Do we need to change something or is the difference sufficiently trivial to let it alone?
forthcoming Corollary \ref{cor:Almost_complex_structures_on_finite-dimensional_subspaces_Hilbert_space}, these operators are skew-adjoint with respect to the real inner products induced by those on
%$X\oplus X$ and $Y\oplus Y$ since $J_1\in\Or(X\oplus X)$ and $J_2\in\Or(Y\oplus Y)$
$X$ and $Y$ since $J_1\in\Or(X)$ and $J_2\in\Or(Y)$
are skew-adjoint and the orthogonal projections are self-adjoint. Thus, if these operators also invertible as implied by Corollary \ref{cor:Almost_complex_structures_on_finite-dimensional_subspaces_Hilbert_space} for large enough $\nu$, they induce almost complex structures on $\tilde K_\nu$ and $\tilde C_\nu$ via the method described by Cannas da Silva \cite[Section 12.2, Proposition 12.3, p. 84]{Cannas_da_Silva_lectures_on_symplectic_geometry}. Indeed, by arguing as in the verification of the invertibility property \eqref{eq:Pi_J_Pi_in_GL_RanPi}, arising in the proof of Theorem \ref{mainthm:Donaldson_1996jdg_3_Hilbert_space} in Section \ref{sec:Proof_generalized_Donaldson_symplectic_subspace_criterion_spectral_projection}, and the proof of Lemma \ref{lem:Approximation_finite-dimensional_subspaces_Hilbert_space}, we obtain for large enough $\nu$ that
\[
  \tilde\Pi_{1,\nu}\circ J_1 \in \Or(\tilde K_\nu)
  \quad\text{and}\quad
  \tilde\Pi_{2,\nu}\circ J_2 \in \Or(\tilde C_\nu).
\]
See the statement and proof of Corollary \ref{cor:Almost_complex_structures_on_finite-dimensional_subspaces_Hilbert_space} for details. Therefore, as claimed, these operators are invertible and skew-adjoint with respect to the real inner products on $X$ and $Y$. 
\qed
\end{rmk}

\begin{rmk}[Vector spaces $\tilde K_\nu$ and $\tilde C_\nu$ are real orthogonal direct sums of subspaces of $X$ and $Y$]
\label{rmk:Distinction_between_linear_direct_and_orthogonal_direct_sums}  
Regarding Remark \ref{rmk:AC_and_symplectic_structures_Kuranishi_model_defined_by_Fredholm_map_Hilbert_spaces}, we emphasize that because
\[
  K_0 \subset K_\nu \quad\text{and}\quad \tilde K_0 = \Pi_{1,\nu}'(K_0) \subset K_\nu
\]
are real linear subspaces and the projection $\Pi_{1,\nu}':X\to K_\nu$ restricts to an isomorphism $\Pi_{1,\nu}':K_0 \to \tilde K_0$ of real linear subspaces, we may write
\[
  \tilde K_\nu = K_0 \oplus \left(K_\nu\cap \tilde K_0^\perp\right) \subset X
\]
to indicate that $\tilde K_\nu$ is a real orthogonal direct sum of real linear subspaces of $X$, by analogy with the definition of $\tilde K_1$ in \eqref{eq:tildeK1_and_tildeC1}. Indeed, just as in the proof of Corollary \ref{cor:Kuranishi_model_defined_by_Fredholm_map_Hilbert_spaces}, we observe that for any $x_0\in K_0$ and $x\in K_\nu\cap\tilde K_0^\perp$, we have $\Real\,\langle x_0,x\rangle = \Real\,\langle \Pi_{1,\nu}' x_0,x\rangle$, because $x\in K_\nu$, and $\Real\,\langle \Pi_{1,\nu}' x_0,x\rangle=0$ since $\Pi_{1,\nu}' x_0 \in \Pi_{1,\nu}'(K_0) = \tilde K_0$ and $x\in \tilde K_0^\perp$.  Thus, $x_0$ and $x$ are real orthogonal in $X$ and so the subspaces $K_0 \subset X$ and $K_\nu\cap \tilde K_0^\perp \subset X$ are real orthogonal in $X$, as claimed. By the same argument, we may write
\[
  \tilde C_\nu = C_0 \oplus \left(C_\nu\cap \tilde C_0^\perp\right) \subset Y
\]
to indicate that $\tilde C_\nu$ is a real orthogonal direct sum of real linear subspaces of $Y$.
\qed\end{rmk}

\subsection{Local Kuranishi model for an open neighborhood of a point in the zero locus of the non-Abelian monopole equations with a regularized Taubes perturbation}
\label{subsec:Local_Kuranishi_model_nbhd_point_zero_locus_unperturbed_non-Abelian_monopole_equations}
In this subsection, we apply the extrinsic definition of a Fredholm section of an abstract Banach vector bundle in Section \ref{subsec:Fredholm_section_Banach_vector_bundle} and Theorem \ref{thm:Kuranishi_model_defined_by_Fredholm_map_Banach_spaces} to construct the local Kuranishi model for an open neighborhood of a point in the zero locus of the non-Abelian monopole equations \eqref{eq:SO(3)_monopole_equations_almost_Hermitian_perturbed_intro_regular} with a regularized Taubes perturbation. This yields an analytic local Kuranishi model as described in Theorem \ref{thm:Kuranishi_model_defined_by_Fredholm_map_Banach_spaces}, where the point is not assumed to be regular. Recall from \eqref{eq:Quotient_space_unitary_triples} that
\[
  \sC(E,H,J,A_d)
  =
  \left.\left(\sA^{1,p}(E,H,A_d) \times W^{1,p}(E \oplus \Lambda^{0,2}(E))\right)\right/W^{2,p}(\SU(E))
\]
is the configuration space for the system \eqref{eq:SO(3)_monopole_equations_almost_Hermitian_intro} of unperturbed non-Abelian monopole equations, where $2<p<\infty$ and $(E,H)$ is a smooth Hermitian vector bundle over a smooth almost Hermitian four-manifold $(X,g,J,\omega)$ and $\SU(E)$ is the subbundle of the bundle of unitary automorphisms of $E$ that induce the identity automorphism of $\det E$. We shall focus on the open subspace \eqref{eq:Quotient_space_non-zero-section_unitary_triples}, namely
\[
  \sC^0(E,H,J,A_d) := \left\{[A,\varphi,\psi] \in \sC(E,H,J,A_d): \varphi \not\equiv 0 \text{ or } \psi \not\equiv 0 \right\},
\]
and recall that $\sC^0(E,H,J,A_d)$ is an analytic Banach manifold when $E$ has complex rank two by Lemma \ref{lem:NonZeroSection_Spinu_pairs_Have_Trivial_Stabilizer} and the analogue for triples $(A,\varphi,\psi)$ of Feehan and Leness \cite[Theorem 12.3.10]{Feehan_Leness_introduction_virtual_morse_theory_so3_monopoles} for pairs $(A,\varphi)$.
%TL6-19-2025: Is the below still active?
%PF11-12-2024 What allows us to rule out Stab(A,\Phi) = \SU(2)?
% If $(A,\varphi,\psi)$ is a $W^{1,p}$ solution to the perturbed non-Abelian monopole equations \cite[Equations (1.22), (1.21b), and (1.21c)]{Feehan_Leness_moduli_space_non-abelian_monopoles_symplectic} and $\psi\not\equiv 0$, then the zero locus $\psi^{-1}(0)$ has measure zero in $(X,g)$ by the Aronszajn Unique Continuation Theorem \cite{Aronszajn}.
% PF11-12-2024 Even when X is only almost Hermitian?
%TL11-26-2025: We only get \Stab(A,\Phi)=\SU(2) when $\Phi=0$ but I'm not seeing where we are using this here
In particular, denoting $\sG = W^{2,p}(\SU(E))$ for convenience, we obtain by analogy with the proof of a similar result for the quotient space of connections (see Singer \cite[Section 1, Theorem 1, p. 9]{Singer_1978}, \cite[Section 3, p. 817]{Singer_1981} or Rudolph and Schmidt \cite[Section 9.1, p. 694]{Rudolph_Schmidt_differential_geometry_mathematical_physics_part2}) that the canonical projection
\begin{equation}
  \label{eq:Principal_bundle_over_quotient_space_non-zero-section_triples}
  \pi: \sA^{1,p}(E,H,A_d) \times W^{1,p}(E \oplus \Lambda^{0,2}(E))\less(0,0) \to \sC^0(E,H,J,A_d)
\end{equation}
is a principal $\sG$-bundle. Standard arguments (modeled on those in
%PF11-13-2024 Add details
\cite{DK,FU,PalaisFoundationGlobal}) imply that the transition functions on overlapping coordinate chart domains are analytic and so \eqref{eq:Principal_bundle_over_quotient_space_non-zero-section_triples} is an analytic principal $\sG$-bundle over $\sC^0(E,H,J,A_d)$.

The action of $\sG = W^{2,p}(\SU(E))$ on $E$ induces a representation of $\sG$ on the Banach space
\[
  L^p\left(\su(E) \oplus \Lambda^{0,2}(\fsl(E)) \oplus \Lambda^{0,1}(E)\right),
\]
and so we can as usual define the Banach vector bundle over $\sC^0(E,H,J,A_d)$ associated to the principal $\sG$-bundle \eqref{eq:Principal_bundle_over_quotient_space_non-zero-section_triples},
\begin{multline}
  \label{eq:Lp_vector_bundle_over_quotient_space_non-zero-section_triples}
  \fV
  :=
  \left(\sA^{1,p}(E,H,A_d) \times W^{1,p}\left(E\oplus \Lambda^{0,2}(E)\right) \less (0,0)\right)
  \\
  \times_{W^{2,p}(\SU(E))} L^p\left(\su(E) \oplus \Lambda^{0,2}(\fsl(E)) \oplus \Lambda^{0,1}(E)\right).
\end{multline}
We observe that $\pi:\fV \to \sC^0(E,H,J,A_d)$ is an analytic Banach vector bundle. The non-Abelian monopole equations \eqref{eq:SO(3)_monopole_equations_almost_Hermitian_perturbed_intro_regular} with a regularized Taubes perturbation are $S^1\times W^{2,p}(\SU(E))$-equivariant and define an $S^1$-equivariant section,
\begin{equation}
  \label{eq:Unperturbed_SO(3)_monopole_equation_section_Lp_vector_bundle_over_quotient_space}
  \cF: \sC^0(E,H,J,A_d) \to \fV,
\end{equation}
which one can check is analytic (again by standard arguments). Let $(A_0,\varphi_0,\psi_0)$ represent a point in $\sC^0(E,H,J,A_d)$ and consider an $S^1$-invariant open neighborhood of $(A_0,\varphi_0,\psi_0)$, 
\begin{multline*}
  \cU_0 \subset (A_0,\varphi_0,\psi_0)
  + \Ker d_{A_0,\varphi_0,\psi_0}^*\cap W_{A_0}^{1,p}\left(T^*X\otimes\su(E) \oplus E \oplus \Lambda^{0,2}(E)\right)
  \\
  \subset \sA^{1,p}(E,H,J,A_d) \times W^{1,p}(E \oplus \Lambda^{0,2}(E)),
\end{multline*}
in the affine Coulomb gauge slice through $(A_0,\varphi_0,\psi_0)$ such that the quotient map $\pi$ from unitary triples $(A,\varphi,\psi)$ in the affine space onto points $[A,\varphi,\psi]$ in the quotient space $\sC(E,H,J,A_d)$ restricts to an $S^1$-equivariant homeomorphism from $\cU_0$ onto an $S^1$-invariant open neighborhood of $[A_0,\varphi_0,\psi_0]$ in $\sC^0(E,H,J,A_d)$. This choice of coordinate chart domain in $\sC^0(E,H,J,A_d)$ induces an $S^1$-equivariant local trivialization,
\[
  \cT:\fV \restriction \pi(\cU_0)
  \to
  \cU_0 \times L^p\left(\su(E) \oplus \Lambda^{0,2}(\fsl(E)) \oplus \Lambda^{0,1}(E)\right),
\]
and $S^1$-equivariant analytic transition maps on overlapping Coulomb gauge coordinate chart domains in $\sC^0(E,H,J,A_d)$. The system \eqref{eq:SO(3)_monopole_equations_almost_Hermitian_perturbed_intro_regular} thus induces an $S^1$-equivariant map,
\[
  \sF \equiv \pr_2\circ \cT \circ \cF \circ \pi:
  \cU_0 \to L^p\left(\su(E) \oplus \Lambda^{0,2}(\fsl(E)) \oplus \Lambda^{0,1}(E)\right),
\]
where $\pr_2$ is projection onto the second factor. Standard theory
%PF11-13-2024 Add detailed refs
(see Feehan and Leness \cite{FL1}) implies that $\sF$ is an $S^1$-equivariant analytic Fredholm map. By Section \ref{subsec:Fredholm_section_Banach_vector_bundle}, we conclude that $\cF$ is an $S^1$-equivariant analytic Fredholm section of $\fV$ in \eqref{eq:Lp_vector_bundle_over_quotient_space_non-zero-section_triples} and we can apply the construction of the $S^1$-equivariant analytic local Kuranishi model provided by Theorem \ref{thm:Kuranishi_model_defined_by_Fredholm_map_Banach_spaces} with $G = S^1$ for an $S^1$-invariant open neighborhood of a point $[A_0,\varphi_0,\psi_0]$ in the zero locus $\cF^{-1}(0) \subset \sC^0(E,H,J,A_d)$.

\section{Deformation operator for the unperturbed non-Abelian monopole equations}
\label{sec:DefTheory_of_nAM_equations}
The Zariski tangent space of the moduli space of unperturbed non-Abelian monopoles is defined as the kernel of a first-order elliptic operator whose definition we now give. We shall refer to it as a \emph{deformation operator} as it is given by rolling up, in the sense of Gilkey \cite[Equation (1.5.1), p. 43]{Gilkey2}, the elliptic deformation complex for the non-Abelian monopole equations defined in Feehan and Leness \cite[Equation (9.3.1)]{Feehan_Leness_introduction_virtual_morse_theory_so3_monopoles}.  As is done for the anti-self-dual equations in Donaldson and Kronheimer \cite[Sections 4.2.5 and 5.4.1]{DK}, we will define the deformation operator as the sum of the linearization of the non-Abelian monopole equations and the operator defining the Coulomb gauge condition.

Let $(E,H)$ be a smooth Hermitian vector bundle over a smooth almost Hermitian manifold $(X,g,J,\omega)$ of real dimension four. Let $(A,\varphi,\psi)$ be a unitary triple on $E$ as in \eqref{eq:A_varphi_psi_in_W1p} that is smooth.  If $\Phi=(\varphi,\psi)\in\Omega^0(W_{\can}^+\otimes E)$, where $(\rho_\can,W_\can)$ is the canonical spin${}^c$ structure over $X$ given in equations \eqref{eq:Canonical_spinc_bundles}, and $\phi\in\Omega^0(E)$ and $\psi\in\Omega^{0,2}(E)$, then the deformation operator for the moduli space of solutions to the unperturbed non-Abelian monopole equations \eqref{eq:SO(3)_monopole_equations_almost_Hermitian_intro} is
\begin{equation}
  \label{eq:SO3MonopoleDeformationOperator}
  \sT_{A,\varphi,\psi} \equiv d_{A,\varphi,\psi}^1+d_{A,\varphi,\psi}^{0,*}:\sE_1\to\sE_2,
\end{equation}
where the real Fr\'echet spaces $\sE_1$ and $\sE_2$ are as in \eqref{eq:sEk} and $d_{A,\varphi,\psi}^1$ is the linearization of the non-Abelian monopole equations \eqref{eq:SO(3)_monopole_equations_almost_Hermitian_intro} at $(A,\varphi,\psi)$, while $d_{A,\varphi,\psi}^{0,*}$ is the operator $d_{A,\Phi}^{0,*}$ defined in \eqref{eq:d_APhi^0_star_aphi_identity_and_vanishing} with $\Phi = (\varphi,\psi)$. Recall that the operator
\[
d_{A,\Phi}^0:\Omega^0(\su(E))\to \Om^1(\su(E))\oplus\Omega^0(W^+\otimes E),
\]
is defined in \eqref{eq:d_APhi^0} as
\[
d_{A,\Phi}^0(\xi):=\left( d_A\xi,-\sR_\Phi\xi\right),\quad\text{for all $\xi\in\Omega^0(\su(E))$}
\]
where $\sR_\Phi$ is defined in \eqref{eq:Phi_Omega0suE_to_Omega0V+} to be
\[
\sR_\Phi:\Omega^0(\su(E))\ni \xi \mapsto \xi\Phi\in \Omega^0(W^+\otimes E).
\]
We have
\[
  \sR_\Phi = \sR_{\varphi,\psi}=\sR_\varphi+\sR_\psi,
\]
where
\begin{align*}
\sR_\varphi &:\Omega^0(\su(E))\ni \xi \mapsto \xi\varphi\in \Omega^0(E),
\\
\sR_\psi &:\Omega^0(\su(E))\ni \xi \mapsto \xi\psi\in \Omega^{0,2}(E).
\end{align*}
Thus, we can write the operator $d_{A,\varphi,\psi}^0$ as
\begin{equation}
\label{eq:nonAbelianMonopole_d0_unitary}
d_{A,\varphi,\psi}^0:\Omega^0(\su(E)) \ni \xi\mapsto d_{A,\varphi,\psi}\xi
\\
 :=(d_A\xi,-\sR_\varphi\xi,-\sR_\psi\psi)
\in \Omega^1(\su(E))\oplus\Omega^0(E)\oplus \Omega^{0,2}(E).
\end{equation}
Because $\sR_\Phi=\sR_\varphi+\sR_\psi$, we the equality of adjoints
\[
\sR_\Phi^*=\sR_\varphi^*+\sR_\psi^*.
\]
Hence, the expression for the $L^2$ adjoint of $d_{A,\Phi}^0$ given in \eqref{eq:d_APhi^0_star_aphi_identity_and_vanishing} yields the following expression for the $L^2$ adjoint of $d_{A,\varphi,\psi}^0$:
\begin{multline}
\label{eq:nonAbelianMonopole_d0*_unitary}
d_{A,\varphi,\psi}^{0,*}:
\Omega^1(\su(E))\oplus\Omega^0(E)\oplus \Omega^{0,2}(E)
\ni (a,\sigma,\tau)\mapsto
d_{A,\varphi,\psi}^{0,*}(a,\sigma,\tau))
\\
= d_A^*a-\sR_\varphi^*\sigma-\sR_\psi^*\tau\in\Omega^0(\su(E)),
\end{multline}
where $\sR_\Phi^*$ is the $L^2$ adjoint of $\sR_\Phi$ given in \eqref{eq:Phi_star_Omega0V+_to_Omega0suE}.
Observe that one can also consider $\sR_\Phi$ as a section of $\Hom(\su(E),W^+\otimes E)$ and then define $\sR_\Phi^*$ as the section of $\Hom(W^+\otimes E,\su(E))$ defined by taking the adjoint of $\sR_\Phi$ on each fiber, $\Hom(\su(E),W^+\otimes E)|_x$, for $x\in X$.

The linearization of the unperturbed non-Abelian monopole equations
\eqref{eq:SO(3)_monopole_equations_almost_Hermitian_intro} at a solution $(A,\varphi,\psi)$ to these equations,
\[
d_{A,\varphi,\psi}^1:
\Omega^1(\su(E))\oplus \Omega^0(E)\oplus \Omega^{0,2}(E)
\to
\Omega^0(\su(E))\oplus \Omega^{0,2}(\su(E))\oplus \Omega^{0,1}(E),
\]
is given by (compare \cite[Equation (9.3.2)]{Feehan_Leness_introduction_virtual_morse_theory_so3_monopoles} where it is denoted $d_{A,\Phi}^1$)
\begin{multline}
\label{eq:SO3Monopoled1}
d_{A,\varphi,\psi}^1(a,\sigma,\tau)
=
\begin{pmatrix}
\Lambda_\omega d_A a
-\frac{i}{2}\left(\varphi\otimes\sigma^*+\sigma\otimes\varphi^*\right)_0
+\frac{i}{2}\star\left(\psi\otimes\tau^*+\tau\otimes\psi^*\right)_0
\\
\pi_{0,2} d_A a-\frac{1}{2}\left( \tau\otimes\varphi^*+\psi\otimes\sigma^*\right)
\\
\bar\partial_A\sigma+\bar\partial_A^*\tau+\frac{1}{4\sqrt 2}\rho_{\can}(\Lambda_\omega d\omega)(\sigma+\tau)
+\frac{1}{\sqrt 2}\rho_{\can}(a)(\sigma+\tau)
\end{pmatrix}
\\
\in
\Omega^0(\su(E)) \oplus \Omega^{0,2}(\su(E)) \oplus \Omega^{0,1}(E),
\\
\quad\text{for } (a,\sigma,\tau)\in\Omega^1(\su(E))\oplus\Omega^0(E)\oplus \Omega^{0,2}(E),
\end{multline}
where $\pi_{0,2}:\Omega^2(\su(E))\to\Omega^{0,2}(\su(E))$ is the projection as used in \eqref{eq:d_A_components_almost_complex_manifold} in Section \ref{sec:Nijenhuis_tensor_and_components_exterior_derivative}.

We define the harmonic spaces determined by the preceding operators as
\begin{subequations}
\label{eq:DefineHiForSO3Monopoles}
\begin{align}
  \label{eq:DefineH0ForSO3Monopoles}
  \bH_{A,\varphi,\psi}^0
  &:= \Ker d_{A,\varphi,\psi}^0 \subset \Omega^0(\su(E)) \subset \sE_2,
  \\
  \label{eq:DefineH1ForSO3Monopoles}
  \bH_{A,\varphi,\psi}^1
  &:= \Ker \left( d_{A,\varphi,\psi}^1+d_{A,\varphi,\psi}^{0,*}\right)
    = \Ker\sT_{A,\varphi,\psi} \subset \sE_1,
  \\
  \label{eq:DefineH2ForSO3Monopoles}
  \bH_{A,\varphi,\psi}^2
  &:= \Ker d_{A,\varphi,\psi}^{1,*}
    \subset \Omega^0(\su(E)) \oplus \Omega^{0,2}(\su(E)) \oplus \Omega^{0,1}(E) \subset \sE_2,
\end{align}
\end{subequations}
where $d_{A,\varphi,\psi}^{1,*}$ is the $L^2$-adjoint of the operator $d_{A,\varphi,\psi}^1$ in \eqref{eq:SO3Monopoled1} and $\sE_1$ and $\sE_2$ are the real Fr\'echet spaces \eqref{eq:sEk}. The vector spaces $\bH_{A,\varphi,\psi}^k$ are the same as $\bH_{A,\Phi}^k$ in \cite[Equation (9.3.4)]{Feehan_Leness_introduction_virtual_morse_theory_so3_monopoles} with $\Phi=(\varphi,\psi)$. In \cite[Lemma 9.3.2]{Feehan_Leness_introduction_virtual_morse_theory_so3_monopoles}, we proved that $\bH_{A,\Phi}^0=(0)$ for solutions to the unperturbed non-Abelian monopole equations \eqref{eq:SO(3)_monopole_equations_almost_Hermitian_intro} when $\Phi\not\equiv 0$. The proof of \cite[Lemma 9.3.2]{Feehan_Leness_introduction_virtual_morse_theory_so3_monopoles} only used the assumption that $(A,\Phi)$ was a solution of the non-Abelian monopole equations to prove that $(A,\Phi)$ was gauge-equivalent to a smooth pair.  We thus have the following

\begin{lem}[Condition for the vanishing of $\bH_{A,\varphi,\psi}^0$]
\label{lem:H0Vanishing}
Let $(E,H)$ be a smooth rank-two Hermitian vector bundle over a smooth almost Hermitian four-manifold $(X,g,J,\omega)$. Let $(A,\varphi,\psi)$ be a unitary triple on $(E,H)$ as in \eqref{eq:A_varphi_psi_in_W1p} that is gauge equivalent to a smooth unitary triple.  If $\varphi$ and $\psi$ are not both identically zero, then $\bH_{A,\varphi,\psi}^0=(0)$.
\end{lem}

\section[Equivalent deformation operators on complex Fr\'echet spaces]{Deformation operator for the unperturbed non-Abelian monopole equations on complex Fr\'echet spaces}
\label{sec:ApproxComplexDefOperator}
Given a unitary triple $(A,\varphi,\psi)$ as in \eqref{eq:A_varphi_psi_in_W1p} that is smooth, we introduce a real linear operator
\[
  \cT_{\bar\partial_A,\varphi,\psi}: \sF_1 \to \sF_2
\]
in the forthcoming \eqref{eq:CplxDef1}. We will show in the forthcoming Proposition \ref{prop:CommutativityOfDeformationComplexes} that $\cT_{\bar\partial_A,\varphi,\psi}$ is related via \eqref{eq:UpsilonIsom_And_Deformation_Operators} to the deformation operator $\sT_{A,\varphi,\psi}$ in
\eqref{eq:SO3MonopoleDeformationOperator} by the isomorphisms $\Upsilon_k:\sF_k\to\sE_k$ of real Fr\'echet spaces, for $k=1,2$, in \eqref{eq:Isomorphism_sEkC_to_sEk} and so $\cT_{\bar\partial_A,\varphi,\psi}$ is equivalent to $\sT_{A,\varphi,\psi}$ modulo the coordinate changes $\Upsilon_k$. We will therefore refer to $\cT_{\bar\partial_A,\varphi,\psi}$ as the \emph{equivalent deformation operator on complex Fr\'echet spaces} or more succinctly as the \emph{equivalent deformation operator}. The almost complex structures on $\sE_k$ are induced by the almost complex structures on $\sF_k$ defined by scalar multiplication by $i=\sqrt{-1}$, so the decomposition of $\cT_{\bar\partial_A,\varphi,\psi}$ into complex linear and complex antilinear components is clear.

We write
\begin{equation}
  \label{eq:sl(E)MultiplicationMap}
  R_\Phi:\Omega^0(\fsl(E))\ni\zeta\mapsto R_\Phi\zeta := \zeta\Phi\in\Omega^0(\fsl(E)
\end{equation}
for the complex linear extension of the operator $\sR_\Phi$ defined in \eqref{eq:Phi_Omega0suE_to_Omega0V+} on the real linear subspace $\Omega^0(\su(E))$ of the complex vector space $\Omega^0(\fsl(E)$. If $\Phi=(\varphi,\psi)$, where $\varphi\in\Omega^0(E)$ and $\psi\in\Omega^{0,2}(E)$, then
\begin{equation}
  \label{eq:Define_sl(E)_MultiplicationMap_Components}
  R_\Phi\zeta = R_{\varphi,\psi}\zeta = (\zeta\varphi,\zeta\psi) \in \Omega^0(E)\oplus \Omega^{0,2}(E),
\quad\text{for } \zeta \in \Omega^0(\fsl(E)),
\end{equation}
and denote the components of $R_{\varphi,\psi}$ by
\begin{subequations}
\label{eq:Define_sl(E)_MultiplicationMap_ComponentsNotation}
\begin{align}
\label{eq:Define_sl(E)_MultiplicationMap_Components0}
R_\varphi&:\Omega^0(\fsl(E))\ni\zeta \mapsto R_\varphi\zeta := \zeta\varphi\in \Omega^0(E),
\\
\label{eq:Define_sl(E)_MultiplicationMap_Components01}
R_\psi&:\Omega^0(\fsl(E))\ni\zeta \mapsto R_\psi\zeta := \zeta\psi\in \Omega^{0,2}(E).
\end{align}
\end{subequations}
We define $R_\varphi^*$ and $R_\psi^*$ as the fiberwise adjoints of $R_\varphi$ and $R_\psi$, respectively, in the following sense. By considering $R_\varphi$ and $R_\psi$ as sections of $\Hom(\fsl(E),E)$ and $\Hom(\fsl(E),\Lambda^{0,2}\otimes E)$, we can define $R_\varphi^*$ and $R_\psi^*$ as the adjoint of $R_\varphi$ and $R_\psi$ respectively on each fiber. Note that $R_\varphi^*$ and $R_\psi^*$ are also the $L^2$-adjoints of $R_\varphi$ and $R_\psi$.
% PF11-18-2025 Means what? Aren't they just L^2 adjoints too?
%PF12-15-2025 Put pointwise discussion here and use that terminology instead of L^2
%TL12-15-2025: So done
%$L^2$ adjoints of $R_{\varphi,\psi}$ and $R_\varphi$ and $R_\psi$ as $R_{\varphi,\psi}^*$ and $R_\varphi^*$ and $R_\psi^*$, respectively.
Because $R_{\varphi,\psi}=R_\varphi+R_\psi$, we have
\[
R_{\varphi,\psi}^* = R_\varphi^*+R_\psi^*.
\]
We have the following real linear pointwise orthogonal projection:
\begin{equation}
\label{eq:Define_sl(E)_to_su(E)_projection}
\pi_{\su(E)}: \fsl(E) \ni\zeta \mapsto \pi_{\su(E)}\zeta:=\frac{1}{2}(\zeta-\zeta^\dagger)\in\su(E).
\end{equation}
By \cite[Equation (10.1.26)]{Feehan_Leness_introduction_virtual_morse_theory_so3_monopoles} and noting the equality $\sR_\Phi^* = \sR_{\varphi,\psi}^*$, we see that
\begin{equation}
\label{eq:MultOperator_Projection_Relation}
\sR_{\varphi,\psi}^* =
\pi_{\su(E)}R_{\varphi,\psi}^*.
\end{equation}
Equations \eqref{eq:MultOperator_Projection_Relation} and \eqref{eq:Define_sl(E)_to_su(E)_projection} give us
\begin{equation}
\label{eq:Real_Projection_Of_Multiplication_Adjoint}
\begin{aligned}
\sR_\varphi^*\sigma&=\frac{1}{2}\left( R_\varphi^*\sigma-(R_\varphi\sigma)^\dagger\right),
\\
\sR_\psi^*\tau&=\frac{1}{2}\left( R_\psi^*\tau-(R_\psi\tau)^\dagger\right).
\end{aligned}
\end{equation}
We define the first-order partial differential operator
\[
\hat\partial_{A,\varphi,\psi}^0:\Omega^0(\fsl(E)) \to \Om^{0,1}(\fsl(E)) \oplus \Omega^0(E)\oplus \Omega^{0,2}(E)
\]
by
\begin{equation}
\label{eq:Define_hat_partial_del0}
\hat\partial_{A,\varphi,\psi}^0\zeta:=\left(\bar\partial_A\zeta,-R_\varphi\zeta,R_\psi(\zeta^\dagger)\right),
\quad\text{for } \zeta \in \Omega^0(\fsl(E)).
\end{equation}
Before continuing, we compare the operator \eqref{eq:Define_hat_partial_del0} with a similar operator appearing in
\cite{Feehan_Leness_introduction_virtual_morse_theory_so3_monopoles}.

\begin{rmk}[Comparison with zeroth order differential in the pre-holomorphic deformation complex]
\label{rmk:Compare_hat_del_with_preholom}
The operator $\hat\partial_{A,\varphi,\psi}^0$ in \eqref{eq:Define_hat_partial_del0} differs slightly from the
operator $\bar\partial_{A,\varphi,\psi}^0$ defined in \cite[Equation (9.5.2)]{Feehan_Leness_introduction_virtual_morse_theory_so3_monopoles}, arising in the pre-holomorphic deformation complex described in \cite[Section 9.5]{Feehan_Leness_introduction_virtual_morse_theory_so3_monopoles}.  Specifically, we recall that
\[
  \bar\partial_{A,\varphi,\psi}^0\zeta:=\left( \bar\partial_A\zeta,-R_\varphi\zeta,-R_\psi\zeta\right),
  \quad\text{for } \zeta \in \Omega^0(\fsl(E)),
\]
so that the two operators differ only in the $\Omega^{0,2}(E)$ component.  When restricted to the subspace $\Omega^0(\su(E))$ of $\Omega^0(\fsl(E))$, where $\xi^\dagger=-\xi$, the two operators are equal.  In addition, if $\psi\equiv 0$, the two operators are equal on all of $\Omega^0(\fsl(E))$.

As described in \cite[Section 9.5]{Feehan_Leness_introduction_virtual_morse_theory_so3_monopoles}, the operator  $\bar\partial_{A,\varphi,\psi}^0$ is defined as the differential at the identity of the action of the determinant-one complex gauge transformations. We use the operator \eqref{eq:Define_hat_partial_del0} instead of $\bar\partial_{A,\varphi,\psi}^0$ to prove the isomorphisms stated in Propositions \ref{prop:H0IsomForAlmostKahler} and \ref{prop:Isom_H1_and H0H2}.
\qed\end{rmk}

To compute the  $L^2$-adjoint of $\hat\partial_{A,\varphi,\psi}^0$, namely
\[
\hat\partial_{A,\varphi,\psi}^{0,*}:\Om^{0,1}(\fsl(E))\oplus \Omega^0(E)\oplus \Omega^{0,2}(E)
\to
\Omega^0(\fsl(E)),
\]
we first compute the adjoint of the real linear map $\zeta\mapsto R_\psi(\zeta^\dagger)$ appearing in \eqref{eq:Define_hat_partial_del0}. The equality (see, for example, \cite[Equation (10.1.49)]{Feehan_Leness_introduction_virtual_morse_theory_so3_monopoles})
\[
\langle \zeta_1,\zeta_2\rangle_{\fsl(E)}
=
\langle \zeta_2^\dagger,\zeta_1^\dagger\rangle_{\fsl(E)}
\quad\text{for } \zeta_1,\zeta_2\in\fsl(E),
\]
implies that
\[
\Real\,\langle \zeta_1,\zeta_2^\dagger \rangle_{\fsl(E)}
=
\Real\,\langle \zeta_1^\dagger,\zeta_2\rangle_{\fsl(E)}.
\]
Thus, for all $\zeta_1,\zeta_2\in\fsl(E)$,
\[
\Real\, \langle\zeta_1, R_\psi(\zeta_2^\dagger)\rangle_{\fsl(E)}
=
\Real\, \langle R_\psi^*\zeta_1,\zeta_2^\dagger\rangle_{\fsl(E)}
=
\Real\, \langle (R_\psi^*\zeta_1)^\dagger,\zeta_2\rangle_{\fsl(E)}.
\]
Hence, the adjoint of the real linear map $\zeta\mapsto R_\psi(\zeta^\dagger)$ is the real linear map $\zeta\mapsto (R_\psi^*\zeta)^\dagger$. We have thus shown that
%PF11-18-2025 Please delete all the Greek letter abbreviations in the monograph.
\begin{equation}
\label{eq:DefineHatPartial0}
\hat\partial_{A,\varphi,\psi}^{0,*}(a'',\sigma,\tau)
=
\bar\partial_A^*a'' - R_\varphi^*\si +(R_\psi^*\tau)^\dagger
\quad\text{for } a''\in\Om^{0,1}(\fsl(E)), \ \si\in\Omega^0(E), \ \tau\in\Omega^{0,2}(E).
\end{equation}
The other operator defining the equivalent deformation operator in the forthcoming \eqref{eq:CplxDef1} is
\[
\bar\partial_{A,\varphi,\psi}^1:\Om^{0,1}(\fsl(E))\oplus \Omega^0(E)\oplus \Omega^{0,2}(E)
\to
\Omega^{0,2}(\fsl(E))\oplus \Om^{0,1}(E),
\]
and which is given by
\begin{multline}
\label{eq:DefineBarPartialWithMu}
\bar\partial_{A,\varphi,\psi}^1(a'',\sigma,\tau)
=
\begin{pmatrix}
      \bar\partial_Aa''-\frac{1}{4}N_J^*(a'')^\dagger - \left(\tau\otimes\varphi^* + \psi\otimes\sigma^*\right)_0
      \\
      \bar\partial_A\sigma + \bar\partial_A^*\tau + a''\varphi + \star ((a'')^\dagger\wedge\star\psi)
  \end{pmatrix}
  \in \Omega^{0,2}(\fsl(E))\oplus \Om^{0,1}(E),
  \\
  \text{for } a''\in\Om^{0,1}(\fsl(E)),\  \si\in\Omega^0(E), \ \tau\in\Omega^{0,2}(E),
\end{multline}
where $N_J^*$ is the Nijenhuis tensor acting on $(0,1)$ forms by pullback as in
\eqref{eq:Donaldson_Kronheimer_p_43_mu_bar_A_Nijenhuis_tensor_E-valued_one_forms}.

\begin{rmk}[Comparison with first order differential in the pre-holomorphic deformation complex]
\label{rmk:Comparison_of_BarPartialWithMu}
%TL12-7-2025: Some rewriting here to try to clarify
%TL12-8-2025: I'm chasing down some possible missing conditions in the statement of \cite[Lemma 10.1.4]{Feehan_Leness_introduction_virtual_morse_theory_so3_monopoles}
%TL12-13-2025: I think Equation (10.1.12b) in \cite[Lemma 10.1.4]{Feehan_Leness_introduction_virtual_morse_theory_so3_monopoles} is missing the $N_J$ term (or we could assume the manifold is complex in the statement of the lemma.
%TL12-15-2025: I will re-edit this when we've got the N_J in the AMS manuscript fixed (since that's the comparison being made)
%PF12-18-2025
The operator \eqref{eq:DefineBarPartialWithMu} differs from that defined by \cite[Equations (10.1.12b) and (10.1.12d)]{Feehan_Leness_introduction_virtual_morse_theory_so3_monopoles} by the absence of the Lee form $\lambda$ defined in \eqref{eq:Gauduchon_1-1-4} through its appearance in the expression  \eqref{eq:Gauduchon_3-7-2_auxiliary_Hermitian_bundle_E} for the Dirac operator. We omit that term here as we are working with symplectic manifolds.
\qed\end{rmk}

We define the \emph{equivalent deformation operator} to be
\begin{equation}
  \label{eq:CplxDef1}
  \cT_{\bar\partial_A,\varphi,\psi}
  \equiv
  \bar\partial_{A,\varphi,\psi}^1+\hat\partial_{A,\varphi,\psi}^{0,*}: \sF_1 \to \sF_2,
\end{equation}
where the complex Fr\'echet spaces $\sF_1$ and $\sF_2$ are as in \eqref{eq:sEkC}. By analogy with the definition of the harmonic spaces \eqref{eq:DefineHiForSO3Monopoles}, we define the harmonic spaces determined by the operators \eqref{eq:Define_hat_partial_del0} and \eqref{eq:DefineBarPartialWithMu} as
\begin{subequations}
\label{eq:DefineHatHarmonics}
\begin{align}
  \label{eq:Define_hatH0}
  \bH_{\bar\partial_A,\varphi,\psi}^0
  &:= \Ker \hat\partial_{A,\varphi,\psi}^0 \subset \Omega^0(\fsl(E))
    \subset\sF_2,
  \\
  \label{eq:Define_hatH1}
  \bH_{\bar\partial_A,\varphi,\psi}^1
  &:= \Ker\left( \bar\partial_{A,\varphi,\psi}^1+\hat\partial_{A,\varphi,\psi}^{0,*}\right)
    = \Ker\cT_{\bar\partial_A,\varphi,\psi} \subset\sF_1,
\\
\label{eq:Define_hatH2}
\bH_{\bar\partial_A,\varphi,\psi}^2&:=
\Ker\bar\partial_{A,\varphi,\psi}^{1,*} \subset \Omega^{0,2}(\fsl(E))\oplus \Om^{0,1}(E) \subset\sF_2,
\end{align}
\end{subequations}
where $\bar\partial_{A,\varphi,\psi}^{1,*}$ is the $L^2$-adjoint of the operator $\bar\partial_{A,\varphi,\psi}^1$ in \eqref{eq:DefineBarPartialWithMu}. Although $\hat\partial_{A,\varphi,\psi}^0$ is not a complex linear operator, because the map $\zeta\mapsto R_\psi(\zeta^\dagger)$ is complex antilinear, $\bH_{\bar\partial_A,\varphi,\psi}^0$ is still a complex vector space. This holds because $\bH_{\bar\partial_A,\varphi,\psi}^0$ is the intersection of the kernels of $\bar\partial_A:\Omega^0(\fsl(E))\to \Omega^{0,1}(\fsl(E))$, and $R_\varphi:\Omega^0(\fsl(E)) \to \Omega^0(E)$, and
$R_\psi((\cdot)^\dagger):\Omega^0(\fsl(E)) \to \Omega^{0,2}(E)$. Because each of these three operators are either complex linear or complex antilinear, their kernels and hence the intersection of these kernels are complex vector spaces.

\section{Complex linear isomorphism between the harmonic spaces $\bH_{A,\varphi,\psi}^0\otimes_\RR\CC$ and $\bH_{\bar\partial_A,\varphi,\psi}^0$}
\label{sec:H0Isom}
The proof of the forthcoming Proposition \ref{prop:H0IsomForAlmostKahler} is identical to the proof of \cite[Proposition 10.2.1]{Feehan_Leness_introduction_virtual_morse_theory_so3_monopoles}, with the exception that the replacement of the operator $\bar\partial_{A,\varphi,\psi}^0$ by $\hat\partial_{A,\varphi,\psi}^0$ described in Remark \ref{rmk:Compare_hat_del_with_preholom} ensures the forthcoming map \eqref{eq:H0_Inclusion} is an isomorphism even without the assumption that $(A,\varphi,\psi)$ is type $1$, that is, without assuming that $\psi\equiv 0$.

\begin{prop}[Complex isomorphism between zeroth order harmonic spaces on almost K\"ahler manifolds]
\label{prop:H0IsomForAlmostKahler}
Let $(E,H)$ be a Hermitian vector bundle over an almost Hermitian four-manifold $(X,g,J,\omega)$. If $(A,\varphi,\psi)$ is a smooth unitary triple on $(E,H)$ as in \eqref{eq:A_varphi_psi_in_W1p}, then the inclusion map $\Omega^0(\su(E))\to \Omega^0(\fsl(E))$ yields an inclusion of complex vector spaces,
\begin{equation}
\label{eq:H0_Inclusion}
\bH_{A,\varphi,\psi}^0\otimes_\RR\CC \to \bH_{\bar\partial_A,\varphi,\psi}^0,
\end{equation}
where $\bH_{\bar\partial_A,\varphi,\psi}^0$ is defined in \eqref{eq:Define_hatH0} and
$\bH_{A,\varphi,\psi}^0$ in \eqref{eq:DefineH0ForSO3Monopoles}. In addition, if $(X,g,J,\omega)$ is almost K\"ahler and $(A,\varphi,\psi)$ is a solution to the unperturbed non-Abelian monopole equations \eqref{eq:SO(3)_monopole_equations_almost_Hermitian_intro}, then the inclusion \eqref{eq:H0_Inclusion} is an isomorphism of complex vector spaces.
\end{prop}

\begin{proof}
If $\xi\in \bH_{A,\varphi,\psi}^0$, then $d_{A,\varphi,\psi}^0\xi=0$. The definition \eqref{eq:nonAbelianMonopole_d0_unitary} of $d_{A,\varphi,\psi}^0$ therefore implies that $d_A\xi=0$ and $\sR_\varphi\xi=\xi\varphi\equiv 0$ and $\sR_\psi\xi=\xi\psi\equiv 0$.  Because elements of the subspaces $\Om^{0,1}(\fsl(E))$ and $\Om^{1,0}(\fsl(E))$ of $\Om^1(\fsl(E))$ are pointwise orthogonal, the equality $d_A\xi=0$ yields $\bar\partial_A\xi=0$ and $\partial_A\xi=0$. Since $\xi\in\Omega^0(\su(E))$, we have $\xi^\dagger=-\xi$.  Hence, because $\xi\psi=0$, we obtain $\xi^\dagger\psi=-\xi\psi=0$ and so $R_\psi\xi^\dagger=0$ by the definition of $R_\psi$ in \eqref{eq:Define_sl(E)_MultiplicationMap_ComponentsNotation}. These equalities and the definition \eqref{eq:Define_hat_partial_del0} of $\hat\partial_{A,\varphi,\psi}^0$ imply that
\[
\hat\partial_{A,\varphi,\psi}^0\xi=0.
\]
Thus $\xi\in \bH_{A,\varphi,\psi}^0$ gives $\xi\in \bH_{\bar\partial_A,\varphi,\psi}^0$. Because $\bH_{\bar\partial_A,\varphi,\psi}^0$ is a complex vector space, $\xi\in \bH_{\bar\partial_A,\varphi,\psi}^0$ implies that $i\xi\in \bH_{\bar\partial_A,\varphi,\psi}^0$ and this proves that \eqref{eq:H0_Inclusion} holds.

We now assume that $(X,g,J,\omega)$ is almost K\"ahler and that $(A,\varphi,\psi)$ is a solution to the unperturbed non-Abelian monopole equations \eqref{eq:SO(3)_monopole_equations_almost_Hermitian_intro}
and aim to prove that the map \eqref{eq:H0_Inclusion} is surjective. To this end, we will prove the

\begin{claim}
\label{claim:H0_and_skewHermitianComponent}
If $\zeta\in\bH_{\bar\partial_A,\varphi,\psi}^0$, then $\zeta-\zeta^\dagger\in \bH_{A,\varphi,\psi}^0$.
\end{claim}

\begin{proof}
The definition of $\bH_{\bar\partial_A,\varphi,\psi}^0$ in \eqref{eq:Define_hatH0} and the definition of
$\hat\partial_{A,\varphi,\psi}^0$ in \eqref{eq:Define_hat_partial_del0} imply that $\zeta\in\bH_{\bar\partial_A,\varphi,\psi}^0$ satisfies
\[
0=\hat\partial_{A,\varphi,\psi}^0\zeta=(\bar\partial_A\zeta,-R_\varphi\zeta,R_\psi(\zeta^\dagger)).
\]
The preceding equality ensures that
\begin{equation}
\label{eq:Ker_bar_rd0}
R_\varphi\zeta=0\quad\text{and}\quad R_\psi(\zeta^\dagger)=0.
\end{equation}
Because $(X,g,J,\omega)$ is almost K\"ahler, we can apply \cite[Equation (10.2.11)]{Feehan_Leness_introduction_virtual_morse_theory_so3_monopoles} to give
\begin{equation}
\label{eq:FLMonograph9-2-11}
d_{A,\varphi,\psi}^{0,*}d_{A,\varphi,\psi}^0(\zeta-\zeta^\dagger)
=
i\Lambda_\omega F_A(\zeta+\zeta^\dagger) + \sR_{\varphi,\psi}^*\sR_{\varphi,\psi}(\zeta-\zeta^\dagger).
\end{equation}
Because $(A,\varphi,\psi)$ satisfies \eqref{eq:SO(3)_monopole_equations_(1,1)_curvature_intro}, the equality \cite[Equation (10.2.12)]{Feehan_Leness_introduction_virtual_morse_theory_so3_monopoles} yields
\begin{multline}
\label{eq:FLMonograph9-2-12}
i\Lambda_\omega F_A(\zeta+\zeta^\dagger)
=
-\frac{1}{2}
\left(
(\varphi\otimes\varphi^*)_0\zeta +(\varphi\otimes\varphi^*)_0\zeta^\dagger
-\zeta (\varphi\otimes\varphi^*)_0-\zeta^\dagger (\varphi\otimes\varphi^*)_0
\right.
\\
\left.
-*(\psi\otimes\psi^*)_0\zeta-*(\psi\otimes\psi^*)_0\zeta^\dagger
+\zeta *(\psi\otimes\psi^*)_0+\zeta^\dagger *(\psi\otimes\psi^*)_0
\right).
\end{multline}
For $\zeta\in\Omega^0(\fsl(E))$ and $\varphi\in\Omega^0(E)$, and $\psi\in\Omega^{0,2}(E)$, we have
\begin{subequations}
\label{eq:sl(E)Mult_in_Tensor_Product}
\begin{align}
\label{eq:sl(E)Mult_in_Tensor_Product_LeftMult_On_0}
\zeta(\varphi\otimes\varphi^*)_0 &= ((\zeta\varphi)\otimes \varphi^*)_0,
\\
\label{eq:sl(E)Mult_in_Tensor_Product_RightMult_On_0}
(\varphi\otimes\varphi^*)_0\zeta^\dagger &= (\varphi\otimes (\zeta\varphi)^*)_0,
\\
\label{eq:sl(E)Mult_in_Tensor_Product_LeftMult_On_02}
\zeta\star(\psi\otimes\psi^*)_0 &= \star ((\zeta\psi)\otimes \psi^*)_0,
\\
\label{eq:sl(E)Mult_in_Tensor_Product_RightMult_On_02}
\star(\psi\otimes\psi^*)_0\zeta^\dagger &= \star(\psi\otimes (\zeta\psi)^*)_0.
\end{align}
\end{subequations}
The preceding equalities and the identity \eqref{eq:FLMonograph9-2-11} imply that the expression \eqref{eq:FLMonograph9-2-12} can be simplified to
\begin{equation}
\label{eq:FLMonograph9-2-13}
i\Lambda_\omega F_A(\zeta+\zeta^\dagger)
=
-\frac{1}{2}
\left(
(\varphi\otimes\varphi^*)_0\zeta -\zeta^\dagger (\varphi\otimes\varphi^*)_0
-*(\psi\otimes\psi^*)_0\zeta^\dagger
+\zeta *(\psi\otimes\psi^*)_0
\right).
\end{equation}
We now examine the second term on the right-hand side of \eqref{eq:FLMonograph9-2-11}. Because the image of $R_\varphi$ is in $\Omega^0(E)$ and thus pointwise orthogonal to $\Omega^{0,2}(E)$ and the image of $R_\psi$ is in $\Omega^{0,2}(E)$ and thus pointwise orthogonal to $\Omega^0(E)$, we have
\begin{equation}
\label{eq:RphiStarR_psiStarVanishing}
R_\varphi^*(\zeta\psi)=0\quad \text{and}\quad
R_\psi^*(\zeta\varphi)=0,\quad\text{for all $\varphi\in\Omega^0(E)$, and $\psi\in\Omega^{0,2}(E)$, and $\zeta\in\Omega^0(\fsl(E))$.}
\end{equation}
We recall the following equalities from \cite[Equation (10.1.30)]{Feehan_Leness_introduction_virtual_morse_theory_so3_monopoles}:
\begin{subequations}
\label{eq:R_varphi_star_and_R_psi_star_tracefree}
\begin{align}
\label{eq:R_varphi_star_sigma_is_sigma_tensor_varphi_star_tracefree}
  R_\varphi^*\sigma &= (\sigma\otimes\varphi^*)_0 \in \Omega^0(\fsl(E)),
  \\
  \label{eq:R_psi_star_tau_is_tau_tensor_psi_star_tracefree}
  R_\psi^*\tau &= \star\left(\tau\wedge\psi^*\right)_0 \in \Omega^0(\fsl(E)),
  \\
  \label{eq:R_varphi_star_sigma_dagger_is_varphi_tensor_sigma_star_tracefree}
  (R_\varphi^*\sigma)^\dagger &= (\varphi\otimes\sigma^*)_0  \in \Omega^0(\fsl(E)),
  \\
  \label{eq:R_psi_star_tau_dagger_is_psi_tensor_tau_star_tracefree}
  (R_\psi^*\tau)^\dagger &= \star(\psi\wedge\tau^*)_0  \in \Omega^0(\fsl(E)).
\end{align}
\end{subequations}
We can then simplify the second term on the right-hand side of \eqref{eq:FLMonograph9-2-11}:
\begin{align*}
\sR_{\varphi,\psi}^*\sR_{\varphi,\psi}(\zeta-\zeta^\dagger)
&=
\pi_{\su(E)}\left( R_\varphi^*+R_\psi^*\right)\left( R_\varphi+R_\psi\right)(\zeta-\zeta^\dagger)
\quad\text{(by \eqref{eq:MultOperator_Projection_Relation})}
\\
&=
\pi_{\su(E)}\left( R_\varphi^*+R_\psi^*\right)\left( -\zeta^\dagger\varphi+\zeta\psi\right)
\quad\text{by \eqref{eq:Ker_bar_rd0})}
\\
&=
\pi_{\su(E)}
\left( -R_\varphi^*(\zeta^\dagger\varphi)+R_\psi^*(\zeta\psi)\right)
\quad\text{(by \eqref{eq:RphiStarR_psiStarVanishing})}
\\
&=
\pi_{\su(E)}\left(
-((\zeta^\dagger\varphi)\otimes\varphi^*)_0 +\star((\zeta\psi)\wedge \psi^*)_0
\right)
\quad\text{(by \eqref{eq:R_varphi_star_sigma_is_sigma_tensor_varphi_star_tracefree} and \eqref{eq:R_psi_star_tau_is_tau_tensor_psi_star_tracefree})}
\\
&=
\pi_{\su(E)}\left(
-\zeta^\dagger(\varphi \otimes\varphi^*)_0 +\zeta\star(\psi\wedge \psi^*)_0
\right)
\quad\text{(by \eqref{eq:sl(E)Mult_in_Tensor_Product_LeftMult_On_0} and \eqref{eq:sl(E)Mult_in_Tensor_Product_LeftMult_On_02}).}
\end{align*}
In the preceding equalities, we have used the identity $\psi\wedge\psi^*=\psi\otimes\psi^*$ arising from the commutativity of wedge product on two-forms as noted following \cite[Equation (10.1.30)]{Feehan_Leness_introduction_virtual_morse_theory_so3_monopoles}. By applying the definition \eqref{eq:Define_sl(E)_to_su(E)_projection} of $\pi_{\su(E)}$, we obtain
\begin{equation}
\label{eq:FLMonograph9-2-15}
\sR_{\varphi,\psi}^*\sR_{\varphi,\psi)}(\zeta-\zeta^\dagger)
=
\frac{1}{2}
\left(
  -\zeta^\dagger(\varphi \otimes\varphi^*)_0 +(\varphi \otimes\varphi^*)_0 \zeta 
  +
  \zeta\star(\psi\wedge \psi^*)_0-\star(\psi\wedge \psi^*)_0\zeta^\dagger
\right).
\end{equation}
Before completing our computation of the right-hand side of \eqref{eq:FLMonograph9-2-11}, we make the following

\begin{rmk}[Comparison of \eqref{eq:FLMonograph9-2-15} with corresponding equality in \cite{Feehan_Leness_introduction_virtual_morse_theory_so3_monopoles}]
The terms involving $\psi$ in equation \eqref{eq:FLMonograph9-2-15} differ from those in the corresponding equality \cite[Equation (10.2.15)]{Feehan_Leness_introduction_virtual_morse_theory_so3_monopoles} arising in the proof of \cite[Proposition 10.2.1]{Feehan_Leness_introduction_virtual_morse_theory_so3_monopoles} due to the difference between the operators $\hat\partial_{A,\varphi,\psi}$ and $\bar\partial_{A,(\varphi,\psi)}$ noted in Remark \ref{rmk:Compare_hat_del_with_preholom}.  Specifically, $\zeta$ satisfies $\zeta^\dagger\psi=0$ in the proof here (see \eqref{eq:Ker_bar_rd0}) while in the proof of \cite[Proposition 10.2.1]{Feehan_Leness_introduction_virtual_morse_theory_so3_monopoles}, the section $\zeta$
satisfied $\zeta\psi=0$.
\qed\end{rmk}

Combining equations \eqref{eq:FLMonograph9-2-11}, \eqref{eq:FLMonograph9-2-13}, and \eqref{eq:FLMonograph9-2-15} yields
\begin{align*}
d_{A,\varphi,\psi}^{0,*}d_{A,\varphi,\psi}^0(\zeta-\zeta^\dagger)
&=
-\frac{1}{2}
\left(
(\varphi\otimes\varphi^*)_0\zeta -\zeta^\dagger (\varphi\otimes\varphi^*)_0
\right.
\\
&\qquad
     - \left.
\star(\psi\otimes\psi^*)_0\zeta^\dagger
+\zeta \star(\psi\otimes\psi^*)_0
\right)
\\
&\qquad
+\frac{1}{2}
\left( -\zeta^\dagger(\varphi \otimes\varphi^*)_0 +(\varphi \otimes\varphi^*)_0 \zeta \right.
\\
&\qquad
+\left.\zeta\star((\zeta\psi)\wedge \psi^*)_0-\star((\zeta\psi)\wedge \psi^*)_0\zeta^\dagger
\right)
\\
&=0.
\end{align*}
This completes the proof of Claim \ref{claim:H0_and_skewHermitianComponent} that $\zeta\in \bH_{\bar\partial_A,\varphi,\psi}^0 \implies \zeta-\zeta^\dagger\in \bH_{A,\varphi,\psi}^0$. 
\end{proof}

We use Claim \ref{claim:H0_and_skewHermitianComponent} to prove that the inclusion \eqref{eq:H0_Inclusion} is surjective exactly as in the proof of \cite[Proposition 10.2.1]{Feehan_Leness_introduction_virtual_morse_theory_so3_monopoles}. Because $\bH_{\bar\partial_A,\varphi,\psi}^0$ is a complex vector space, $\zeta\in \bH_{\bar\partial_A,\varphi,\psi}^0 \implies i\zeta\in \bH_{\bar\partial_A,\varphi,\psi}^0$. Therefore,
\[
\xi_1:=\frac{1}{2}\left( \zeta-\zeta^\dagger\right)\in \bH_{A,\varphi,\psi}^0\quad\text{and}\quad
\xi_2:=-\frac{1}{2}\left( (i\zeta)-(i\zeta)^\dagger\right)\in \bH_{A,\varphi,\psi}^0.
\]
Thus, every $\zeta \in \bH_{\bar\partial_A,\varphi,\psi}$ can be written as
\[
\zeta=\xi_1+i\xi_2\in\bH_{A,\varphi,\psi}^0\otimes_\RR\CC,
\]
proving that the inclusion \eqref{eq:H0_Inclusion} is surjective and hence an isomorphism of complex vector spaces. This completes the proof of Proposition \ref{prop:H0IsomForAlmostKahler}.
\end{proof}

\begin{rmk}[Role of the almost K\"ahler and non-Abelian monopole hypotheses in Proposition \ref{prop:H0IsomForAlmostKahler}]
\label{rmk:aK_and_nA_Monopole_Hypoth}
In the proof of Proposition \ref{prop:H0IsomForAlmostKahler}, we used the hypotheses that $(X,g,J,\omega)$ was almost K\"ahler and that  $(A,\varphi,\psi)$ was a solution to the unperturbed non-Abelian monopole equations \eqref{eq:SO(3)_monopole_equations_almost_Hermitian_intro} only to prove Claim \ref{claim:H0_and_skewHermitianComponent}.  Specifically, we only needed the almost K\"ahler hypothesis to obtain the expression \eqref{eq:FLMonograph9-2-12} for the Laplacian of $\zeta-\zeta^\dagger$ and only needed the hypothesis that $(A,\varphi,\psi)$ was a solution to the non-Abelian monopole equations to obtain the expression \eqref{eq:FLMonograph9-2-12} for the $(1,1)$-component of the curvature of $A$.
\qed\end{rmk}

\section{Real linear isomorphisms between the harmonic spaces $\bH_{A,\varphi,\psi}^k$ and $\bH_{\bar\partial_A,\varphi,\psi}^k$ for $k=1,2$}
\label{sec:H1_and_H2_Isom}
In this section, we construct isomorphisms between the kernels and cokernels of the deformation operators \eqref{eq:SO3MonopoleDeformationOperator} and \eqref{eq:CplxDef1}, namely
\begin{align*}
  d_{A,\varphi,\psi}^1+d_{A,\varphi,\psi}^{0,*}:\sE_1 &\to \sE_2,
  \\
  \bar\partial_{A,\varphi,\psi}^1+\hat\partial_{A,\varphi,\psi}^{0,*}:\sF_1 &\to \sF_2,
\end{align*}
via the isomorphisms $\Upsilon_k:\sE_k \to \sF_k$ between the Fr\'echet spaces $\sE_k$ and $\sF_k$ defined in \eqref{eq:sEk} and \eqref{eq:sEkC}, respectively, for $k=1,2$. We will prove

\begin{prop}[Real linear isomorphism between harmonic spaces on almost K\"ahler four-manifolds]
\label{prop:Isom_H1_and H0H2}
Let $(E,H)$ be a smooth Hermitian vector bundle over a closed smooth almost K\"ahler four-manifold $(X,g,J,\omega)$. Let $(A,\varphi,\psi)$ be a unitary triple on $E$ as in \eqref{eq:A_varphi_psi_in_W1p} that is smooth. Then the real linear isomorphism $\Upsilon_1$ in \eqref{eq:Isomorphism_sE1C_to_sE1} induces an isomorphism of real vector spaces,
\begin{equation}
\label{eq:SO3MonopoleH1Isomorphism}
\bH_{A,\varphi,\psi}^1 \cong \bH_{\bar\partial_A,\varphi,\psi}^1,
\end{equation}
where $\bH_{A,\varphi,\psi}^1$ is as in \eqref{eq:DefineH1ForSO3Monopoles} and $\bH_{\bar\partial_A,\varphi,\psi}^1$ is as in \eqref{eq:Define_hatH1}. The $L^2$ adjoint $\Upsilon_2^*$ of the real linear isomorphism $\Upsilon_2$ in \eqref{eq:Isomorphism_sE2C_to_sE2} induces an isomorphism of real vector spaces,
\begin{equation}
\label{eq:SO3MonopoleH2Isomorphism}
\Ker\left(d_{A,\varphi,\psi}^1+d_{A,\varphi,\psi}^{0,*} \right)^*
\cong
\Ker\left( \bar\partial_{A,\varphi,\psi}^1+\hat\partial_{A,\varphi,\psi}^{0,*}\right)^*,
\end{equation}
where $d_{A,\varphi,\psi}^{0,*}$ is as in \eqref{eq:d_APhi^0}, and $d_{A,\varphi,\psi}^1$ is as in \eqref{eq:SO3Monopoled1}, and $\hat\partial_{A,\varphi,\psi}^{0,*}$ is as in \eqref{eq:DefineHatPartial0},  
and $\bar\partial_{A,\varphi,\psi}^1$ is as in \eqref{eq:DefineBarPartialWithMu}.
\end{prop}

An isomorphism similar to \eqref{eq:SO3MonopoleH1Isomorphism} is given in the proof of \cite[Proposition 10.1.1]{Feehan_Leness_introduction_virtual_morse_theory_so3_monopoles}, as we discuss in Remark \ref{rmk:ComparisonOfH1Proofs}. The real linear isomorphisms $\Upsilon_1$ and $\Upsilon_2$ in \eqref{eq:Isomorphism_sEkC_to_sEk}, and hence the isomorphisms in Proposition \ref{prop:Isom_H1_and H0H2}, are gauge-equivariant in the following sense.

\begin{lem}[Gauge equivariance of the isomorphisms $\Upsilon_k$]
\label{lem:GaugeEquivarianceOfPsiIsomorphisms}
For $k=1,2$, the real linear isomorphisms $\Upsilon_k$ in \eqref{eq:Isomorphism_sEkC_to_sEk} are equivariant with respect to the actions of the smooth gauge transformations $\Omega^0(\SU(E))$ on Fr\'echet spaces $\sE_k$ and $\sF_k$  given by
\begin{subequations}
\label{eq:GaugeGroupActionOn_E1_E2_E1c_E2c}
\begin{align}
\label{eq:GaugeGroupActionOn_E1c}
  & \Omega^0(\SU(E))\times \sF_1\ni \left(u, (a'',\sigma,\tau)\right)
    \mapsto \left( u^{-1}a''u,u^{-1}\sigma,u^{-1}\tau\right)\in\sF_1,
\\
\label{eq:GaugeGroupActionOn_E1}
& \Omega^0(\SU(E))\times \sE_1\ni \left(u, (a,\sigma,\tau)\right)
\mapsto \left( u^{-1}a''u,u^{-1}\sigma,u^{-1}\tau\right)\in\sE_1,
\\
\label{eq:GaugeGroupActionOn_E2c}
& \Omega^0(\SU(E))\times \sF_2\ni \left(u, (\zeta,v,\nu)\right)
\mapsto \left( u^{-1}\zeta u,u^{-1}v u,u^{-1}\nu\right)\in\sF_2,
  \\
\label{eq:GaugeGroupActionOn_E2}
& \Omega^0(\SU(E))\times \sE_2\ni \left(u, (\xi,w,\nu)\right)
\mapsto \left( u^{-1}\xi u,u^{-1}w,u^{-1}\nu\right)\in\sE_2,
\end{align}
\end{subequations}
where
$a''\in\Omega^{0,1}(\fsl(E))$, and $\sigma\in\Omega^0(E)$, and $\tau\in\Omega^{0,2}(E)$, and 
$a\in\Omega^1(\su(E))$, and $\zeta\in \Omega^0(\fsl(E))$, and $v\in\Omega^{0,2}(\fsl(E))$, and $\nu\in\Om^{0,1}(E)$, and $\xi\in\Omega^0(\su(E))$, and $w\in\Omega^{0,2}(\su(E))$. The extensions of the isomorphisms $\Upsilon_k$ to the $W^{1,p}$-completions of $\sE_k$ and $\sF_k$ are also equivariant with respect to the action of $W^{2,p}(\SU(E))$ acting on these completions.
\end{lem}

The proof of Proposition \ref{prop:Isom_H1_and H0H2} requires the following

\begin{lem}[Commuting complex conjugation and pointwise Hermitian adjoint with $\bar\partial_A$ and $\bar\partial_A^*$]
\label{lem:PointwiseHermitianAdjoint}
(See Feehan and Leness \cite[Lemma 10.1.3]{Feehan_Leness_introduction_virtual_morse_theory_so3_monopoles}.)
Let $(E,H)$ be a smooth Hermitian vector bundle $E$ over a closed, smooth  almost Hermitian four-manifold
$(X,g,J,\omega)$. If $A$ is a
%COMMENT TL11-9-2025: From the proof of 10.1.3, it appears we only need $W^{1,p}$.
smooth, unitary connection on $E$ and $p,q$ are non-negative integers, then
\begin{subequations}
  \label{eq:Commute_dagger_and_del}
  \begin{align}
    \label{eq:Commute_dagger_and_del_E}
    %PF6-7-2024 This seems wrong: it implies that for any \bar\partial_A on L, the operator \partial_A is completely determined by \bar\partial_A, independent of the Hermitian metric on E? Thus, \nabla_A = \partial_A + \bar\partial_A would be completely determined by \bar\partial_A? Would this be a complex linear but not unitary (that is, Chern) connection on E? Check usage elsewhere.
    %TL9-6-2024: In the statement, we assume $A$ is unitary and $H$ fixed, if that helps.
    % PF7-25-2025 Could you please recheck the proof?
    %TL12-15-2025: I can't find anything wrong w/ the proof.
    \partial_A(\bar\psi) &= \overline{\left(\bar\partial_A\varphi\right)},
                              \quad\text{for all } \psi \in \Omega^{p,q}(E),
    \\
    \label{eq:Commute_dagger_and_del_glE}
    \partial_A(\eta^\dagger) &= \left(\bar\partial_A\eta\right)^\dagger,
                               \quad\text{for all } \eta \in \Omega^{p,q}(\gl(E)),
  \end{align}
\end{subequations}
and
\begin{subequations}
  \label{eq:AdjointCommute_dagger_and_del}
  \begin{align}
    \label{eq:AdjointCommute_dagger_and_del_E}
    \partial_A^*(\bar\psi) &= \overline{(\bar\partial_A^*\psi)},
                            \quad\text{for all } \psi \in \Omega^{p,q+1}(E),
    \\
    \label{eq:AdjointCommute_dagger_and_del_glE}
    \partial_A^*(\eta^\dagger) &= \left(\bar\partial_A^*\eta\right)^\dagger,
                            \quad\text{for all } \eta \in \Omega^{p,q+1}(\gl(E)).
  \end{align}
\end{subequations}
\end{lem}

We begin the proof of Proposition \ref{prop:Isom_H1_and H0H2} with the following extension of \cite[Lemma 10.1.4]{Feehan_Leness_introduction_virtual_morse_theory_so3_monopoles} describing the composition of the
operator $d_{A,\varphi,\psi}^1$ with the real linear isomorphism $\Upsilon_1$ in \eqref{eq:Isomorphism_sE1C_to_sE1}.  We note that the forthcoming Lemma \ref{lem:Decomposition_of_d1} includes a description of this composition under the weaker assumption that the four-manifold is only almost Hermitian.

\begin{lem}[First order harmonic spaces over almost Hermitian manifold four-manifolds]
\label{lem:Decomposition_of_d1}
Let $(E,H)$ be a Hermitian vector bundle over a smooth  almost Hermitian manifold $(X,g,J,\omega)$ of real dimension four.
Let $(A,\varphi,\psi)$ be a unitary triple on $E$ as in \eqref{eq:A_varphi_psi_in_W1p} that is smooth. For $\sE_1$ as in \eqref{eq:sE1}, let $(a,\sigma,\tau)\in \sE_1$ satisfy
\[
(a,\sigma,\tau)=\left(\frac{1}{2}(a''-(a'')^\dagger),\sigma,\tau\right)=\Upsilon_1(a'',\sigma,\tau),
\]
where $(a'',\sigma,\tau)\in\sF_1$, for $\sF_1$ as in \eqref{eq:sE1C} and real linear isomorphism $\Upsilon_1$ as in \eqref{eq:Isomorphism_sE1C_to_sE1}. Then $d_{A,\varphi,\psi}^1(a,\sigma,\tau)=0$ if and only if $(a'',\sigma,\tau)$ and $a'=-(a'')^\dagger$ satisfy
\begin{subequations}
\label{eq:SO3Monopole_d1_AlmostHermitian}
\begin{align}
\label{eq:SO3Monopole_d1_AlmostHermitian(1,1)}
\Lambda_\omega\left( \partial_Aa''+\bar\partial_Aa'\right) -i\left((R_\varphi^*\sigma)^\dagger+R_\varphi^*\sigma\right)
+i\left((R_\psi^*\tau)^\dagger+R_\psi^*\tau\right)
&=0 \in \Omega^0(\su(E)),
\\
\label{eq:SO3Monopole_d1_AlmostHermitian(0,2))}
\bar\partial_A a'' + \frac{1}{4}N_J^*a' - \left(\tau\otimes\varphi^*+\psi\otimes\sigma^*\right)_0
&=0 \in \Omega^{0,2}(\fsl(E)),
\\
\label{eq:SO3Monopole_d1_AlmostHermitian(0,1))}
\bar\partial_A\sigma+\bar\partial_A^*\tau + \bar\partial_A^*\tau + a''\varphi - \star (a'\wedge\star\psi)
+\frac{1}{4\sqrt 2}\rho(\lambda)(\sigma,\tau)
&=0\in \Omega^{0,1}(E),
\end{align}
\end{subequations}
where $N_J$ is the Nijenhuis tensor \eqref{eq:Nijenhuis_tensor} and $\lambda \in \Omega^1(X,\RR)$ is the Lee form defined by $(g,J)$, so $\lambda=\Lambda_\omega(d\omega)$ as in \eqref{eq:Gauduchon_1-1-4}. In addition, if $(X,g,J,\omega)$ is almost K\"ahler, then \eqref{eq:SO3Monopole_d1_AlmostHermitian} is equivalent to
\begin{subequations}
\label{eq:SO3Monopole_d1_AlmostKahler}
\begin{align}
\label{eq:SO3Monopole_d1_AlmostKahler(1,)}
i\bar\partial_A^*a''+i(\bar\partial_A^*a'')^\dagger
-i\left((R_\varphi^*\sigma)^\dagger+R_\varphi^*\sigma\right)
+i\left((R_\psi^*\tau)^\dagger+R_\psi^*\tau\right)
&=0 \in \Omega^0(\su(E)),
\\
\label{eq:SO3Monopole_d1_AlmostKahler(0,2)}
\bar\partial_A a'' + \frac{1}{4}N_J^*a' - \left(\tau\otimes\varphi^*+\psi\otimes\sigma^*\right)_0
&=0 \in \Omega^{0,2}(\fsl(E)),
\\
\label{eq:SO3Monopole_d1_AlmostKahler(0,1)}
\bar\partial_A\sigma+\bar\partial_A^*\tau  + a''\varphi - \star (a'\wedge\star\psi)
&=0\in \Omega^{0,1}(E).
\end{align}
\end{subequations}
\end{lem}

\begin{proof}
Recall the expression given in \eqref{eq:Canonical_Clifford_multiplication_positive_spinors} for the Clifford multiplication operator on the canonical spin${}^c$ bundle (where it is denoted $\rho_{\can}$),
\[
 \rho(a)(\sigma,\tau) = \sqrt{2}\left(a''\wedge\sigma - \star(a'\wedge\star\tau)\right)
                                \in \Omega^{0,1}(X),
\]
where $a=(1/2)(a'+a'')$ for $a'\in\Omega^{1,0}(\su(E))$ and $a''\in\Om^{0,1}(\su(E))$ satisfy $a'=-(a'')^\dagger$. The preceding equality and the definition of $d_{A,\varphi,\psi}^1$ in \eqref{eq:SO3Monopoled1} imply that
\begin{equation}
  \label{eq:Decomp_of_d1_v1}
d_{A,\varphi,\psi}^1(a,\varphi,\psi)
=
\begin{pmatrix}
\Lambda_\omega d_Aa -\frac{i}{2}\left(\varphi\otimes\sigma^*+\sigma\otimes\varphi^*\right)_0
+\frac{i}{2} \left(\psi\otimes\tau^*+\tau\otimes\psi^*\right)_0
\\
\pi_{0,2}d_Aa -\frac{1}{2} \left(\tau\otimes\varphi^*+\psi\otimes\sigma^*\right)_0
\\
\bar\partial_A\sigma+\bar\partial_A^*\tau  + a''\varphi - \star (a'\wedge\star\psi)
+ \frac{1}{4\sqrt 2}\rho(\lambda)(\varphi+\psi)
\end{pmatrix},
\end{equation}
where $\pi_{0,2}:\Omega^2(\fsl(E))\to \Omega^{0,2}(\fsl(E))$ is defined by pointwise orthogonal projection. Then by \eqref{eq:d_A_sum_components_almost_complex_manifold}
\begin{align*}
 d_A a&=\frac{1}{2} \left( \partial_A+\mu+\bar\mu+\bar\partial_A\right)(a'+a'')
 \\
 &=\frac{1}{2}(\partial_A a' + \mu a'')+\frac{1}{2}(\partial_Aa''+\bar\partial_Aa') +
 \frac{1}{2}(\bar\partial_A a'' +\bar\mu a')
 \\
 & \in \Omega^{2,0}(\fsl(E))\oplus \Omega^{1,1}(\fsl(E)) \oplus \Omega^{0,2}(\fsl(E)).
\end{align*}
The resulting equalities (using \eqref{eq:Donaldson_Kronheimer_p_43_mu_A_Nijenhuis_tensor_E-valued_one_forms} for the relation between $\bar\mu$ and the Nijenhuis tensor),
\[
\Lambda_\omega d_Aa= \frac{1}{2}\Lambda_\omega\left( \partial_Aa''+\bar\partial_Aa'\right)
\quad\text{and}\quad
\pi_{0,2}d_Aa=\frac{1}{2}\left(\bar\partial_A a'' +\frac{1}{4}N_J^*a'\right),
\]
imply that we can rewrite \eqref{eq:Decomp_of_d1_v1} as
\begin{equation}
\label{eq:Decomp_of_d1_v2}
d_{A,\varphi,\psi}^1(a,\sigma,\tau)
=
\begin{pmatrix}
\frac{1}{2}\Lambda_\omega\left( \partial_Aa''+\bar\partial_Aa'\right) -\frac{i}{2}\left(\varphi\otimes\sigma^*+\sigma\otimes\varphi^*\right)_0
+\frac{i}{2} \left(\psi\otimes\tau^*+\tau\otimes\psi^*\right)_0
\\
\frac{1}{2}\left(\bar\partial_A a'' +\frac{1}{4}N_J^*a'\right) -\frac{1}{2} \left(\tau\otimes\varphi^*+\psi\otimes\sigma^*\right)_0
\\
\bar\partial_A\sigma+\bar\partial_A^*\tau + a''\varphi - \star (a'\wedge\star\psi)
+ \frac{1}{4\sqrt{2}}\rho(\lambda)(\varphi+\psi)
\end{pmatrix}.
\end{equation}
By applying the identities in \eqref{eq:R_varphi_star_and_R_psi_star_tracefree}, we can rewrite \eqref{eq:Decomp_of_d1_v2} as
\begin{equation}
\label{eq:Decomp_of_d1_v3}
d_{A,\varphi,\psi}^1(a,\sigma,\tau)
=
\begin{pmatrix}
\frac{1}{2}
\Lambda_\omega\left( \partial_Aa''+\bar\partial_Aa'\right)
-\frac{i}{2}\left((R_\varphi^*\sigma)^\dagger+R_\varphi^*\sigma\right)
+\frac{i}{2}\left((R_\psi^*\tau)^\dagger+R_\psi^*\tau\right)
\\
\frac{1}{2}
\left(\bar\partial_A a'' + \frac{1}{4}N_J^*a'\right) - \frac{1}{2}\left(\tau\otimes\varphi^*+\psi\otimes\sigma^*\right)_0
\\
\bar\partial_A\sigma+\bar\partial_A^*\tau  + a''\varphi - \star (a'\wedge\star\psi)
+\frac{1}{4\sqrt{2}}\rho(\lambda)\varphi
\end{pmatrix}.
\end{equation}
Removing the factors of $1/2$ from the first two components of \eqref{eq:Decomp_of_d1_v3} yields
 \eqref{eq:SO3Monopole_d1_AlmostHermitian}.

If we further assume that $(X,g,J,\omega)$ is almost K\"ahler, the K\"ahler identities
\eqref{eq:Kaehler_identity_commutator_Lambda_del-bar_A_and_Lambda_del_A_almost_Hermitian_glE} (noting that they hold on almost K\"ahler manifolds) imply that
\begin{align*}
\Lambda_\omega\left( \partial_Aa''+\bar\partial_Aa'\right)
&=
i\bar\partial_A^*a'' - i\partial_A^*a'
\\
&=
i\bar\partial_A^*a'' + i\partial_A^*(a'')^\dagger
\quad\text{by $a'=-(a'')^\dagger$}
\\
&=
i\bar\partial_A^*a''+i(\bar\partial_A^*a'')^\dagger
\quad\text{(by \eqref{eq:AdjointCommute_dagger_and_del_glE})}.
\end{align*}
Applying the preceding equality
to \eqref{eq:Decomp_of_d1_v3}
and observing that $\lambda=0$ if $(X,g,J,\omega)$ is almost K\"ahler yields
\begin{equation}
\label{eq:Decomp_of_d1_v3_almostK}
d_{A,\varphi,\psi}^1(a,\sigma,\tau)
=
\begin{pmatrix}
\frac{i}{2}
\left( \bar\partial_A^*a''+(\bar\partial_A^*a'')^\dagger\right)
-\frac{i}{2}\left((R_\varphi^*\sigma)^\dagger+R_\varphi^*\sigma\right)
+\frac{i}{2}\left((R_\psi^*\tau)^\dagger+R_\psi^*\tau\right)
\\
\frac{1}{2}
\left(\bar\partial_A a'' +\frac{1}{4}N_J^*a'\right) - \frac{1}{2}\left(\tau\otimes\varphi^*+\psi\otimes\sigma^*\right)_0
\\
\bar\partial_A\sigma+\bar\partial_A^*\tau  + a''\varphi - \star (a'\wedge\star\psi)
\end{pmatrix}.
\end{equation}
Hence, when  $(X,g,J,\omega)$ is almost K\"ahler, \eqref{eq:Decomp_of_d1_v3_almostK} implies
that $d_{A,\varphi,\psi}^1(a,\sigma,\tau)=0$ is equivalent to \eqref{eq:SO3Monopole_d1_AlmostKahler}.  This completes the proof of Lemma \ref{lem:Decomposition_of_d1}.
\end{proof}

If we assume that $(X,g,J,\omega)$ is almost K\"ahler, we can use the expression \eqref{eq:Decomp_of_d1_v3_almostK} for $d_{A,\varphi,\psi}^1$ to derive the following relationship between deformation operators.

\begin{prop}[Relationship between the deformation operators]
\label{prop:CommutativityOfDeformationComplexes}
Let $(E,H)$ be a smooth Hermitian vector bundle over a closed smooth almost K\"ahler four-manifold $(X,g,J,\omega)$. If $(A,\varphi,\psi)$ is a unitary triple on $E$ as in \eqref{eq:A_varphi_psi_in_W1p} that is smooth, then the linear isomorphisms $\Upsilon_k$ in \eqref{eq:Isomorphism_sEkC_to_sEk} yield the identity,
\begin{equation}
\label{eq:UpsilonIsom_And_Deformation_Operators}
\Upsilon_2\circ \left(\bar\partial_{A,\varphi,\psi}^1 + \hat\partial_{A,\varphi,\psi}^{0,*}\right)
=
\left(d_{A,\varphi,\psi}^1+d_{A,\varphi,\psi}^{0,*}\right)\circ\Upsilon_1,
\end{equation}
where $\hat\partial_{A,\varphi,\psi}^{0,*}$ is defined in \eqref{eq:DefineHatPartial0},
and $\bar\partial_{A,\varphi,\psi}^1$ in \eqref{eq:DefineBarPartialWithMu},
and $d_{A,\varphi,\psi}^{0,*}$ in \eqref{eq:nonAbelianMonopole_d0*_unitary},
and $d_{A,\varphi,\psi}^1$ in \eqref{eq:SO3Monopoled1}.
\end{prop}

\begin{proof}
If $(a'',\sigma,\tau)\in\sF_1$, for $\sF_1$ as in \eqref{eq:sE1C} with $a''\in\Om^{0,1}(\fsl(E))$, we write
\[
(a,\sigma,\tau)=\Upsilon_1(a'',\sigma,\tau) \in \sE_1,
\]
for $\sE_1$ as in \eqref{eq:sE1} with $a=(1/2)(a'+a'')\in\Omega^1(\su(E))$ and $a'=-(a'')^\dagger$, then
\[
\partial_A^*a'
=
-
\partial_A^*(a'')^\dagger.
\]
Combining the preceding equality with  \eqref{eq:AdjointCommute_dagger_and_del_glE} gives
\begin{equation}
\label{eq:del_a'_to_bar_del_a''}
\partial_A^*a'
=
-(\bar\partial_A^*a'')^\dagger.
\end{equation}
We compute
\begin{align*}
d_{A,\varphi,\psi}^{0,*}(a,\sigma,\tau)
&=d_A^*a-\sR_\varphi^*\sigma-\sR_\psi^*\tau
\quad\text{(by \eqref{eq:nonAbelianMonopole_d0*_unitary})}
  \\
&=\frac{1}{2}\left(\bar\partial_A^* a'' +\partial_A^*a'\right)
-
\frac{1}{2}\left(R_\varphi^*\si - (R_\varphi^*\si)^\dagger\right)
-
\frac{1}{2}\left( R_\psi^*\tau -(R_\psi^*\tau)^\dagger\right)
  \\
&\qquad
     \quad\text{(by $a=(1/2)(a'+a'')$ and \eqref{eq:Real_Projection_Of_Multiplication_Adjoint},}
\end{align*}
and substitute \eqref{eq:del_a'_to_bar_del_a''} into the preceding equality to give
\begin{equation}
\label{eq:d0APhi_v2}
d_{A,\varphi,\psi}^{0,*}(a,\sigma,\tau)=
\frac{1}{2}\left(
\bar\partial_A^*a''-(\bar\partial_A^*a'')^\dagger
-
\left(R_\varphi^*\si - (R_\varphi^*\si)^\dagger\right)
-
\left( R_\psi^*\tau -(R_\psi^*\tau)^\dagger\right)\right).
\end{equation}
Combining \eqref{eq:d0APhi_v2} with the expression for $d_{A,\varphi,\psi}^1(a,\sigma,\tau)$ in
\eqref{eq:Decomp_of_d1_v3_almostK}
and
the equality $(1/4)N_J^*=\bar\mu$ from
 \eqref{eq:Donaldson_Kronheimer_p_43_mu_bar_A_Nijenhuis_tensor_E-valued_one_forms}
gives
\begin{multline}
\label{eq:Non_AbelianM_Def_v1}
\left(d_{A,\varphi,\psi}^1+d_{A,\varphi,\psi}^{0,*}\right)(a,\sigma,\tau)
\\
=\begin{pmatrix}
\frac{1}{2}\left(
\bar\partial_A^*a''-(\bar\partial_A^*a'')^\dagger\right)
-
\frac{1}{2}\left(R_\varphi^*\si - (R_\varphi^*\si)^\dagger\right)
-
\frac{1}{2}\left( R_\psi^*\tau -(R_\psi^*\tau)^\dagger\right)\
\\
\frac{i}{2}
\left(\bar\partial_A^*a''+ (\bar\partial_A^*a'')^\dagger\right)
-\frac{i}{2}\left(R_\varphi^*\sigma+(R_\varphi^*\sigma)^\dagger\right)
+\frac{i}{2}\left(R_\psi^*\tau+(R_\psi^*\tau)^\dagger\right)\
\\
\frac{1}{2}
\left(\bar\partial_A a'' -\frac{1}{4}N_J^*(a'')^\dagger - \left(\tau\otimes\varphi^*+\psi\otimes\sigma^*\right)_0\right)
\\
\bar\partial_A\sigma+\bar\partial_A^*\tau  + a''\varphi + \star ((a'')^\dagger\wedge\star\psi)
\end{pmatrix}
\\
\in
\Omega^0(\su(E))\oplus \Omega^0(\su(E)) \oplus \Omega^{0,2}(\su(E)) \oplus \Omega^{0,1}(E).
\end{multline}
We compare \eqref{eq:Non_AbelianM_Def_v1} with the following expression obtained by combining
\eqref{eq:DefineHatPartial0} and \eqref{eq:DefineBarPartialWithMu}
\begin{multline}
\label{eq:Decomp_of_AC_def}
\left( \bar\partial_{A,\varphi,\psi}^1+\hat\partial_{A,\varphi,\psi}^{0,*}\right)(a'',\sigma,\tau)
=
\begin{pmatrix}
\bar\partial_A^*a'' - R_\varphi^*\si +(R_\psi^*\tau)^\dagger
\\
\bar\partial_Aa''-\frac{1}{4}N_J^*(a'')^\dagger - \left(\tau\otimes\varphi^* + \psi\otimes\sigma^*\right)_0
      \\
 \bar\partial_A\sigma + \bar\partial_A^*\tau + a''\varphi + \star ((a'')^\dagger\wedge\star\psi)
  \end{pmatrix}
 \\
 \in\Omega^0(\fsl(E))\oplus\Omega^{0,2}(\su(E)) \oplus \Omega^{0,1}(E).
\end{multline}
The definition \eqref{eq:Isomorphism_sE2C_to_sE2} of the isomorphism $\Upsilon_2$ and equations \eqref{eq:Non_AbelianM_Def_v1} and \eqref{eq:Decomp_of_AC_def} yield
\begin{equation}
\left((d_{A,\varphi,\psi}^1+d_{A,\varphi,\psi}^{0,*})\circ\Upsilon_1\right)(a'',\sigma,\tau)
=
\left( \Upsilon_2
\circ
(\bar\partial_{A,\varphi,\psi}^1+\hat\partial_{A,\varphi,\psi}^{0,*})\right)(a'',\sigma,\tau).
\end{equation}
This completes the proof of \eqref{eq:UpsilonIsom_And_Deformation_Operators} and hence Proposition \ref{prop:CommutativityOfDeformationComplexes}.
\end{proof}

We can now give

\begin{proof}[Proof of Proposition \ref{prop:Isom_H1_and H0H2}]
By the definitions \eqref{eq:DefineH1ForSO3Monopoles} and \eqref{eq:Define_hatH1}, we have
\begin{align*}
\bH_{A,\varphi,\psi}^1&=  \Ker\left(d_{A,\varphi,\psi}^1+d_{A,\varphi,\psi}^{0,*}\right),
\\
\bH_{\bar\partial_A,\varphi,\psi}^1&=\Ker\left( \bar\partial_{A,\varphi,\psi}^1+\hat\partial_{A,\varphi,\psi}^{0,*}\right).
\end{align*}
Because $\Upsilon_1$ and $\Upsilon_2$ are real linear isomorphisms, these definitions and the equality \eqref{eq:UpsilonIsom_And_Deformation_Operators} yield the isomorphism \eqref{eq:SO3MonopoleH1Isomorphism}. The equality \eqref{eq:UpsilonIsom_And_Deformation_Operators} gives the identity,
\[
\Upsilon_1^*\circ
\left(d_{A,\varphi,\psi}^1+d_{A,\varphi,\psi}^{0,*}\right)^*
=
\left(\bar\partial_{A,\varphi,\psi}^1+\hat\partial_{A,\varphi,\psi}^{0,*}\right)^*\circ  \Upsilon_2^*.
\]
Because $\Upsilon_1^*$ is a real linear isomorphism, the preceding equality ensures that $\Upsilon_2^*$
gives the isomorphism \eqref{eq:SO3MonopoleH2Isomorphism}.  This completes the proof of Proposition  \ref{prop:Isom_H1_and H0H2}.
\end{proof}

\begin{rmk}[Comparison of the isomorphism  \eqref{eq:SO3MonopoleH1Isomorphism} with
prior results]
\label{rmk:ComparisonOfH1Proofs}
When $(X,g,J,\omega)$ is K\"ahler rather than merely almost K\"ahler, the isomorphism   \eqref{eq:SO3MonopoleH1Isomorphism} appears as \cite[Equation (10.1.4)]{Feehan_Leness_introduction_virtual_morse_theory_so3_monopoles} in \cite[Proposition 10.1.1]{Feehan_Leness_introduction_virtual_morse_theory_so3_monopoles}. In addition to generalizing \cite[Proposition 10.1.1]{Feehan_Leness_introduction_virtual_morse_theory_so3_monopoles} from K\"ahler to almost K\"ahler four-manifolds, the proof given here encodes the algebraic manipulations (taking linear combinations of the components of the deformation operator \eqref{eq:SO3MonopoleDeformationOperator}) used in the proof of \cite[Proposition 10.1.1]{Feehan_Leness_introduction_virtual_morse_theory_so3_monopoles} into the isomorphism $\Upsilon_2$.
\qed\end{rmk}

\section[Deformation operator for the perturbed non-Abelian monopole equations]{Deformation operator for the non-Abelian monopole equations with a regularized Taubes perturbation}
\label{sec:HarmonicTheoryForPerturbedEquations}
In this section, $(E,H)$ is a smooth Hermitian vector bundle over a smooth almost Hermitian four-manifold  $(X,g,J,\omega)$ and $(A,\varphi,\psi)$ is a unitary triple on $E$ as in \eqref{eq:A_varphi_psi_in_W1p} that is assumed to be smooth. We now adapt the proof of Proposition \ref{prop:Isom_H1_and H0H2} to harmonic spaces defined by the linearization of the system \eqref{eq:SO(3)_monopole_equations_almost_Hermitian_perturbed_intro_regular} of non-Abelian monopole equations with a regularized Taubes perturbation.

We denote the linearization of the system \eqref{eq:SO(3)_monopole_equations_almost_Hermitian_perturbed_intro_regular} at a solution $(A,\varphi,\psi)$ by
\begin{equation}
\label{eq:Linearization_TaubesPert_Regular}
d_{A,\varphi,\psi,r}^1(a,\sigma,\tau)
:=
d_{A,\varphi,\psi}^1(a,\sigma,\tau) +(ir/8)D\wp_\gamma(\psi)(\tau)\otimes\omega,
\end{equation}
where $a\in\Omega^1(\su(E))$, and $\sigma\in\Omega^0(E)$, and $\tau\in\Omega^{0,2}(E)$, and $D\wp_\gamma(\psi)\tau$ is the derivative of the perturbation $\wp_\gamma$ in \eqref{eq:Definition_wp_intro_regular} at $\psi$ in the direction $\tau\in\Omega^{0,2}(E)$, and computed in \eqref{eq:Derivative_of_Taubes_Pert_regular}, and $d_{A,\varphi,\psi}^1$ in \eqref{eq:SO3Monopoled1} is the linearization of the system \eqref{eq:SO(3)_monopole_equations_almost_Hermitian_intro} of unperturbed non-Abelian monopole equations. The linearization $d_{A,\varphi,\psi,r}^1$ and the operator $d_{A,\varphi,\psi}^{0,*}$ define the following analogues of the harmonic spaces in \eqref{eq:DefineHiForSO3Monopoles}:
\begin{subequations}
\label{eq:DefineHi_For_RegPerturbedMonopoles}
\begin{align}
  \label{eq:H0_For_TaubesPert_regular}
  \bH_{A,\varphi,\psi}^0 &:= \Ker d_{A,\varphi,\psi}^0 \subset \Omega^0(\su(E)) \subset \sE_2,
  \\
  \label{eq:H1_For_TaubesPert_regular}
  \bH_{A,\varphi,\psi,r}^1
  &:=
    \Ker\left(d_{A,\varphi,\psi,r}^1 + d_{A,\varphi,\psi}^{0,*}\right)
    = \Ker\sT_{A,\varphi,\psi,r} \subset \sE_1,
  \\
  \label{eq:DefineH2ForSO3Monopoles_Perturbed}
  \bH_{A,\varphi,\psi,r}^2
  &:=
    \Ker d_{A,\varphi,\psi,r}^{1,*}
    \subset \Omega^0(\su(E)) \oplus \Omega^{0,2}(\su(E)) \oplus \Omega^{0,1}(E) \subset \sE_2,           
\end{align}
\end{subequations}
where $d_{A,\varphi,\psi,r}^{1,*}$ is the $L^2$-adjoint of the operator \eqref{eq:Linearization_TaubesPert_Regular} and $\sE_k$ is the Fr\'echet space defined in \eqref{eq:sEk} and, by analogy with \eqref{eq:SO3MonopoleDeformationOperator}, we denote
\begin{equation}
  \label{eq:Perturbed_Deformation_Operator}
  \sT_{A,\varphi,\psi,r} \equiv d_{A,\varphi,\psi,r}^1+d_{A,\varphi,\psi}^{0,*}:\sE_1\to\sE_2.
\end{equation}
The harmonic space \eqref{eq:H1_For_TaubesPert_regular} represents the Zariski tangent space of $\sM^0(E,g,J,\omega,r)$ at the point $[A,\varphi,\psi]$ while the harmonic space \eqref{eq:DefineH2ForSO3Monopoles_Perturbed} represents the cokernel of the operator $d_{A,\varphi,\psi,r}^1$. As observed in Remark \ref{rmk:ComparisonOfH1Proofs}, the proof of Proposition \ref{prop:Isom_H1_and H0H2} used the fact that the $\Omega^0(\su(E))$-component of $d_{A,\varphi,\psi}^1$ is the Hermitian component of the linear operator $\hat\partial_{A,\varphi,\psi}^{0,*}$ (as seen in the second component of the isomorphism $\Upsilon_2$ in \eqref{eq:Isomorphism_sE2C_to_sE2}).

To generalize the isomorphisms of Proposition \ref{prop:Isom_H1_and H0H2} from the operator $d_{A,\varphi,\psi}^1$ to the perturbed operator $d_{A,\varphi,\psi,r}^1$, we need to introduce a corresponding perturbation of the operator $\hat\partial_{A,\varphi,\psi}^{0,*}$. Because it will be useful to keep track of the complex linear and complex antilinear components of the perturbation $D\wp_\gamma(\psi)$ appearing in \eqref{eq:Linearization_TaubesPert_Regular}, we introduce  a complex linear operator
\begin{equation}
\label{eq:ComplexLinearPerturbation}
p_\psi:\Omega^{0,2}(E)\to \Omega^0(\fsl(E)), \quad\text{for } \psi\in\Omega^{0,2}(E),
\end{equation}
whose Hermitian component is the perturbation $D\wp_\gamma(\psi)$. We thus define
\begin{multline}
\label{eq:DefineComplexDerivativeOfPerturbation}
p_\psi(\tau)
:=
-4\left( \gamma^2+|\psi|_{\Lambda^{0,2}(E)}^2\right)^{-2}\langle\tau,\psi\rangle_{\Lambda^{0,2}(E)}\star(\psi\otimes\psi^*)_0
\\
+4\left( \gamma^2+|\psi|_{\Lambda^{0,2}(E)}^2\right)^{-1}\star\left(\tau\otimes\psi^*\right)_0
\in \Omega^0(\fsl(E)),
\quad\text{for } \tau\in\Omega^{0,2}(E).
\end{multline}
Observe that $p_\psi$ is complex linear and  by comparing with the expression \eqref{eq:Derivative_of_Taubes_Pert_regular} for $D\wp_\gamma(\psi)(\tau)$ we see that
\begin{equation}
\label{eq:ComplexPerturbationRelation1}
p_\psi(\tau)+p_\psi(\tau)^\dagger = D\wp_\gamma(\psi)(\tau),
\end{equation}
as desired. Using the expression in \eqref{eq:R_psi_star_tau_is_tau_tensor_psi_star_tracefree} for $R_\psi^*$, we observe that the pointwise adjoint of $p_\psi$,
\[
p_\psi^*:\Omega^0(\fsl(E)) \to \Omega^{0,2}(E),
\]
is given by
\begin{equation}
\label{eq:DefineComplexDerivativeOfPerturbationAdjoint}
p_\psi^*(\zeta)
=
-4\left( \gamma^2+|\psi|_{\Lambda^{0,2}(E)}^2\right)^{-2}\langle \zeta,\star(\psi\otimes\psi^*)_0\rangle \psi
+ 4\left( \gamma^2+|\psi|_{\Lambda^{0,2}(E)}^2\right)^{-1}\zeta\psi,
\quad\text{for } \zeta\in\Omega^0(\fsl(E)). 
\end{equation}
We note that $p_\psi^*$ is  also complex linear.  We define
\begin{subequations}
\label{eq:DefineHat0OperatorsPerturbed}
\begin{align}
\label{eq:DefineHat0OperatorPerturbed}
\hat\partial_{A,\varphi,\psi,r}^0
&:=
\hat\partial_{A,\varphi,\psi}^0+\frac{r}{4}\left( p_\psi^*+(p_\psi^*)^\dagger\right),
\\
\label{eq:DefineHat0*OperatorPerturbed}
\hat\partial_{A,\varphi,\psi,r}^{0,*}
&:=
\hat\partial_{A,\varphi,\psi}^{0,*}+\frac{r}{4}\left( p_\psi+p_\psi^\dagger\right),
\end{align}
\end{subequations}
where the operators $\hat\partial_{A,\varphi,\psi}^0$ and $\hat\partial_{A,\varphi,\psi}^{0,*}$ are as in 
\eqref{eq:Define_hat_partial_del0} and \eqref{eq:DefineHatPartial0} respectively. These operators define the following analogues of the harmonic spaces in \eqref{eq:DefineHatHarmonics},
\begin{subequations}
\label{eq:DefineHiForApproxComplex_Perturbed}
\begin{align}
\label{eq:DefineH0ForApproxComplex_Perturbed}
  \bH_{\bar\partial_A,\varphi,\psi,r}^0
  &:= \Ker \hat\partial_{A,\varphi,\psi,r}^0 \subset \Omega^0(\fsl(E)) \subset \sF_2,
\\
\label{eq:DefineH1ForApproxComplex_Perturbed}
  \bH_{\bar\partial_A,\varphi,\psi,r}^1
  &:= \Ker\left( \bar\partial_{A,\varphi,\psi}^1+\hat\partial_{A,\varphi,\psi,r}^{0,*}\right)
    = \Ker\cT_{\bar\partial_A,\varphi,\psi,r} \subset\sF_1,
\\
\label{eq:DefineH2ForApproxComplex_Perturbed}
  \bH_{\bar\partial_A,\varphi,\psi,r}^2
  &:= \Ker  \bar\partial_{A,\varphi,\psi}^{1,*}
    \subset \Omega^{0,2}(\fsl(E))\oplus \Om^{0,1}(E) \subset \sF_2,
\end{align}
\end{subequations}
where the Fr\'echet spaces $\sF_k$, for $k=1,2$, are as in \eqref{eq:sEkC} and $\bar\partial_{A,\varphi,\psi}^{1,*}$ is the $L^2$ adjoint of the operator $\bar\partial_{A,\varphi,\psi}^1$ defined in \eqref{eq:DefineBarPartialWithMu} and, by analogy with \eqref{eq:CplxDef1}, we denote
\begin{equation}
\label{eq:Perturbed_AC_Deformation_Operator}
\cT_{\bar\partial_A,\varphi,\psi,r} \equiv \bar\partial_{A,\varphi,\psi}^1+\hat\partial_{A,\varphi,\psi,r}^{0,*}: \sF_1\to\sF_2.
\end{equation}
We can now prove the following extension of Proposition \ref{prop:Isom_H1_and H0H2}.

%PF11-18-2025 The unperturbed version of this prop should be edited to match the style and content of the perturbed version and the title below edited to bring into line with that of unperturbed prop.
%TL12-13-2025: I think the content is now similar.  Title updated below
\begin{prop}[Real linear isomorphisms between harmonic spaces defined by perturbed deformation operators on almost K\"ahler four-manifolds]
\label{prop:SO3MonopoleDeformationCorollary_Perturbed}
Continue the hypotheses and notation of Proposition \ref{prop:Isom_H1_and H0H2}. Then the real linear isomorphism $\Upsilon_1$ in \eqref{eq:Isomorphism_sE1C_to_sE1} induces an isomorphism of real vector spaces,
\begin{equation}
\label{eq:SO3Monopole_Kernel_of_PerturbedDefOperatorIsom}
\Ker\cT_{\bar\partial_A,\varphi,\psi,r}
\cong
\Ker \sT_{A,\varphi,\psi,r}
\end{equation}
and thus an isomorphism
\begin{equation}
\label{eq:SO3MonopoleH1Isomorphism_perturbed}
\bH_{\bar\partial_A,\varphi,\psi,r}^1
 \cong 
\bH_{A,\varphi,\psi,r}^1,
\end{equation}
where $\bH_{A,\varphi,\psi,r}^1$ is as in \eqref{eq:H1_For_TaubesPert_regular} and $\bH_{\bar\partial_A,\varphi,\psi,r}^1$ as in \eqref{eq:DefineH1ForApproxComplex_Perturbed}. The $L^2$ adjoint $\Upsilon_2^*$ of the real linear isomorphism $\Upsilon_2$ in \eqref{eq:Isomorphism_sE2C_to_sE2} induces an isomorphism of real vector spaces,
\begin{equation}
  \label{eq:SO3MonopoleH2Isomorphism_perturbed}
  \Ker \sT_{A,\varphi,\psi,r}^*,
\cong
\Ker \cT_{\bar\partial_A,\varphi,\psi,r}^*,
\end{equation}
where $\sT_{A,\varphi,\psi,r}^*$ is the $L^2$ adjoint of the operator $\sT_{A,\varphi,\psi,r}$ in \eqref{eq:Perturbed_Deformation_Operator} and $\cT_{\bar\partial_A,\varphi,\psi,r}^*$ is the $L^2$ adjoint of the operator $\cT_{\bar\partial_A,\varphi,\psi,r}$ in \eqref{eq:Perturbed_AC_Deformation_Operator}. Furthermore, if $(A,\varphi,\psi)$ solves the non-Abelian monopole equations \eqref{eq:SO(3)_monopole_equations_almost_Hermitian_perturbed_intro_regular} with a regularized Taubes perturbation and $(\varphi,\psi) \not \equiv (0,0)$, then
% PF11-18-2025 Emphasize \bH_{A,\varphi,\psi,r}^0 = 0 too.
%TL12-15-2025: Do we still have that space?  It's now the same as \bH_{A,\varphi,\psi}^0
\begin{equation}
\label{eq:Vanishing_Perturbed_H0}
\bH_{A,\varphi,\psi}^0 =(0),
\end{equation}
where $\bH_{A,\varphi,\psi}^0$ is as in \eqref{eq:H0_For_TaubesPert_regular}
and
%PF11-18-2025 Is \bH_{\bar\partial_A,\varphi,\psi,r}^0 not zero as well? Isn't the below saying more? E.g., isomorphisms with \bH_{\bar\partial_A,\varphi,\psi,r}^2 and \Ker \sT_{A,\varphi,\psi,r}^* ?
\begin{equation}
\label{eq:SO3MonopoleH2Isomorphism_perturbed_H2version}
\bH_{A,\varphi,\psi,r}^2
\cong
\Ker \cT_{\bar\partial_A,\varphi,\psi,r}^*,
\end{equation}
where $\bH_{A,\varphi,\psi,r}^2$ is as in \eqref{eq:DefineH2ForSO3Monopoles_Perturbed}.
%PF7-24-2025 It would be more transparent to say that \sT and \cT have the same kernels and cokernels and then conclude that the H^1's coincide
%TL12-15-2025 Now done above
%and, when H^0 = (0), that H^2 = \Coker\sT and thus H^2 = \cT.
%TL12-15-2025: That follows immediately from the last and \eqref{eq:SO3MonopoleH2Isomorphism_perturbed}.  But if you want to add a line emphasizing this, go right ahead.
\end{prop}

\begin{proof}
Using the expression \eqref{eq:SO3Monopoled1} for $d_{A,\varphi,\psi}^1$, we can write $d_{A,\varphi,\psi,r}^1$ in \eqref{eq:Linearization_TaubesPert_Regular} explicitly as
\begin{multline}
\label{eq:SO3Monopoled1TaubesPerturbation}
d_{A,\varphi,\psi,r}^1(a,\sigma,\tau)
=
\begin{pmatrix}
\Lambda_\omega d_A a
-\frac{i}{2}\left(\varphi\otimes\sigma^*+\sigma\otimes\varphi^*\right)_0
+\frac{i}{2}\star\left(\psi\otimes\tau^*+\tau\otimes\psi^*\right)_0
%PF5-5-2025 You wanted to change the sign below from "-" to "+" but "-" is correct.
%TL6-4-2025: I'm looking at (1.6.10), moving the perturbation to the same side as the curvature, and thus the sign changes.  The + sign is what appears in \eqref{eq:Linearization_TaubesPert_Regular} above
%PF12-18-2025 You're right. I put "+" back below
+\frac{i r}{4}D\wp_\gamma(\psi)\tau
\\
\pi_{0,2} d_A a-\frac{1}{2}\left( \tau\otimes\varphi^*+\psi\otimes\sigma^*\right)
\\
\bar\partial_A\sigma+\bar\partial_A^*\tau
+\frac{1}{\sqrt 2}\rho_{\can}(a)(\sigma+\tau)+\frac{1}{4\sqrt 2}\rho_{\can}(\Lambda_\omega d\omega)(\sigma+\tau)
\end{pmatrix}
\\
\in
\Omega^0(\su(E)) \oplus \Omega^{0,2}(\su(E)) \oplus \Omega^{0,1}(E),
\end{multline}
for
\[
  (a,\sigma,\tau) \in \Omega^1(\su(E)) \oplus \Omega^0(E) \oplus \Omega^{0,2}(E).
\]
Using \eqref{eq:SO3Monopoled1TaubesPerturbation}, \eqref{eq:Decomp_of_d1_v3_almostK}, and the definition \eqref{eq:Isomorphism_sE1C_to_sE1} of the isomorphism $\Upsilon_1$, we can write
\begin{multline}
  \label{eq:Decomp_of_d1_perturbed_v3}
  (d_{A,\varphi,\psi,r}^1\circ\Upsilon_1)(a'',\sigma,\tau)
\\
=
\begin{pmatrix}
\frac{1}{2}
\Lambda_\omega\left( \partial_Aa''+\bar\partial_Aa'\right)
-\frac{i}{2}\left((R_\varphi^*\sigma)^\dagger+R_\varphi^*\sigma\right)
+\frac{i}{2}\left((R_\psi^*\tau)^\dagger+R_\psi^*\tau\right)
+\frac{ir}{4}D\wp_\gamma(\psi)\tau
\\
\frac{1}{2}
\left(\bar\partial_A a'' + \frac{1}{4}N_J^* a'\right) - \frac{1}{2}\left(\tau\otimes\varphi^*+\psi\otimes\sigma^*\right)_0
\\
\bar\partial_A\sigma+\bar\partial_A^*\tau  + a''\varphi + \star ((a'')^\dagger\wedge\star\psi)
+\frac{1}{4\sqrt 2}\rho_{\can}(\Lambda_\omega d\omega)(\sigma+\tau)
\end{pmatrix},
\\
\text{for } (a'',\sigma,\tau) \in \Omega^1(\fsl(E)) \oplus \Omega^0(E) \oplus \Omega^{0,2}(E).
\end{multline}
Noting that $\Lambda_\omega\omega=2$ by \eqref{eq:Lambda_on_KahlerForm} and substituting the expression \eqref{eq:ComplexPerturbationRelation1} for $D\wp_\gamma(\psi)\tau$, we can rewrite \eqref{eq:Linearization_TaubesPert_Regular} as
%PF12-18-2025 Reason for writing p_\psi(\tau) and not p_\psi\tau?
\begin{equation}
\label{eq:SO3Monopole_d1_perturbed_V1}
d_{A,\varphi,\psi,r}^1(a,\sigma,\tau)
=
d_{A,\varphi,\psi}^1(a,\sigma,\tau)
+
\frac{ir}{4}\left( p_\psi(\tau)+p_\psi(\tau)^\dagger\right).
\end{equation}
Because $(X,g,J,\omega)$ is almost K\"ahler by hypothesis, we have $d\omega \equiv 0$ and so the expression given in \eqref{eq:Non_AbelianM_Def_v1} for $d_{A,\varphi,\psi}^1+d_{A,\varphi,\psi}^{0,*}$ and the equalities \eqref{eq:ComplexPerturbationRelation1} and \eqref{eq:Decomp_of_d1_perturbed_v3} imply that the perturbed deformation operator appearing in the definition \eqref{eq:H1_For_TaubesPert_regular} is
\begin{multline}
\label{eq:PerturbedSO3MonopoleDeformationOperatorDecomposed}
\left((d_{A,\varphi,\psi,r}^1+d_{A,\varphi,\psi}^{0,*})\circ\Upsilon_1\right)(a'',\sigma,\tau)
\\
=\begin{pmatrix}
\frac{1}{2}\left(
\bar\partial_A^*a''-(\bar\partial_A^*a'')^\dagger\right)
-
\frac{1}{2}\left(R_\varphi^*\si - (R_\varphi^*\si)^\dagger\right)
-
\frac{1}{2}\left( R_\psi^*\tau -(R_\psi^*\tau)^\dagger\right)
\\
\frac{i}{2}
\left(\bar\partial_A^*a''+ (\bar\partial_A^*a'')^\dagger\right)
-
\frac{i}{2}\left(R_\varphi^*\sigma+(R_\varphi^*\sigma)^\dagger\right)
+
\frac{i}{2}\left(R_\psi^*\tau+(R_\psi^*\tau)^\dagger\right)
+
\frac{ir}{4}\left( p_\psi(\tau)+p_\psi(\tau)^\dagger\right)
\\
\frac{1}{2}
\left(\bar\partial_A a'' -\frac{1}{4}N_J^*(a'')^\dagger - \left(\tau\otimes\varphi^*+\psi\otimes\sigma^*\right)_0\right)
\\
\bar\partial_A\sigma+\bar\partial_A^*\tau  + a''\varphi +\star ((a'')^\dagger\wedge\star\psi)
\end{pmatrix}
\\
\in
\Omega^0(\su(E))\oplus \Omega^0(\su(E)) \oplus \Omega^{0,2}(\su(E)) \oplus \Omega^{0,1}(E),
\end{multline}
for
\[
  (a'',\sigma,\tau) \in \Omega^1(\fsl(E)) \oplus \Omega^0(E) \oplus \Omega^{0,2}(E).
\]  
The expression \eqref{eq:Decomp_of_AC_def} for $\bar\partial_{A,\varphi,\psi}^1+\hat\partial_{A,\varphi,\psi}^{0,*}$
and the definition \eqref{eq:DefineHat0*OperatorPerturbed} of $\hat\partial_{A,\varphi,\psi,r}^{0,*}$ imply that
\begin{equation}
\label{eq:PerturbedACDeformationOperatorDecomposed}
\left(\bar\partial_{A,\varphi,\psi}^1+\hat\partial_{A,\varphi,\psi,r}^{0,*}\right)(a'',\sigma,\tau)
=\begin{pmatrix}
\bar\partial_A^*a'' - R_\varphi^*\si +(R_\psi^*\tau)^\dagger 
% TL6-4-2025: Possible sign discrepancy here
%PF7-24-2025 Is there?
%TL11-11-2025: The discussion is above referring to (1.6.10)
%PF12-18-2025 It's not. Different discussion here. The "-" below should surely be "+"? Please check and then check any ramifications if sign below incorrect.
%TL12-18-2025: Sign corrected to +.  This matches the sign in the $\sT$ operator appearing in (5.6.19) above
%PF12-18-2025 I corrected the sign
%- (r/4)\left( p_\psi(\tau)+p_\psi(\tau)^\dagger\right)
+ (r/4)\left( p_\psi(\tau)+p_\psi(\tau)^\dagger\right)
\\
\bar\partial_Aa''-\frac{1}{4}N_J^*(a'')^\dagger - \left(\tau\otimes\varphi^* + \psi\otimes\sigma^*\right)_0
      \\
 \bar\partial_A\sigma + \bar\partial_A^*\tau + a''\varphi + \star ((a'')^\dagger\wedge\star\psi)
  \end{pmatrix},
\end{equation}
for
\[
  (a'',\sigma,\tau) \in \Omega^1(\fsl(E)) \oplus \Omega^0(E) \oplus \Omega^{0,2}(E).
\]  
The definition \eqref{eq:Isomorphism_sE2C_to_sE2} of the isomorphism $\Upsilon_2:\sF_2\to\sE_2$ implies that
\[
\Upsilon_2\left( \frac{r}{4}\left( p_\psi(\tau)+p_\psi(\tau)^\dagger\right),0,0 \right)
=
%PF12-18-2025 Inserted missing "r" below
\left(
  %TL12-18-2025: Added factor of $r$: checked factors of 2 as well
0,\frac{ir}{4}\left( p_\psi(\tau)+p_\psi(\tau)^\dagger\right),
0,0
\right),
\]
since $\xi_1 = \frac{1}{2}(\zeta-\zeta^\dagger) = 0$ and $\xi_2 = \frac{i}{2}(\zeta+\zeta^\dagger) = (ir/4)(p_\psi(\tau)+p_\psi(\tau)^\dagger)$ when $\zeta = (ir/4)(p_\psi(\tau)+p_\psi(\tau)^\dagger)$.

As in the proof of \eqref{eq:UpsilonIsom_And_Deformation_Operators}, a comparison of \eqref{eq:PerturbedSO3MonopoleDeformationOperatorDecomposed} and \eqref{eq:PerturbedACDeformationOperatorDecomposed} produces the equality
\begin{equation}
\label{eq:Perturbed_Equation_Equivalence_Of_Deformation_Complex}
\Upsilon_2\circ \left(\bar\partial_{A,\varphi,\psi}^1+\hat\partial_{A,\varphi,\psi,r}^{0,*}\right)\circ \Upsilon_1^{-1}
=
d_{A,\varphi,\psi,r}^1 + d_{A,\varphi,\psi}^{0,*}.
\end{equation}
Just as \eqref{eq:UpsilonIsom_And_Deformation_Operators} gave the isomorphisms \eqref{eq:SO3MonopoleH1Isomorphism} and \eqref{eq:SO3MonopoleH2Isomorphism} in Proposition \ref{prop:Isom_H1_and H0H2}, equation \eqref{eq:Perturbed_Equation_Equivalence_Of_Deformation_Complex} yields the isomorphisms \eqref{eq:SO3Monopole_Kernel_of_PerturbedDefOperatorIsom} and \eqref{eq:SO3MonopoleH2Isomorphism_perturbed}.  Equation \eqref{eq:SO3MonopoleH1Isomorphism_perturbed} follows immediately from \eqref{eq:SO3Monopole_Kernel_of_PerturbedDefOperatorIsom} and the definitions \eqref{eq:H1_For_TaubesPert_regular} and \eqref{eq:DefineH1ForApproxComplex_Perturbed}.

The fact that those equations are gauge-equivariant implies that
\[
d_{A,\varphi,\psi,r}^1 \circ d_{A,\varphi,\psi}^0=0,
\]
so that
\[
\Ran d_{A,\varphi,\psi}^0 \subset \Ker d_{A,\varphi,\psi,r}^1,
\]
and thus
\begin{equation}
\label{eq:Range_Intersection_When_Complex}
\Ran d_{A,\varphi,\psi}^0\cap\Ker d_{A,\varphi,\psi,r}^{1,*}=(0).
\end{equation}
We then have
\begin{align*}
\Ker \left(d_{A,\varphi,\psi,r}^1+d_{A,\varphi,\psi}^{0,*}\right)^*
&=
\Ker\left( d_{A,\varphi,\psi,r}^{1,*} + d_{A,\varphi,\psi}^0\right)
\\
&=
     \Ker d_{A,\varphi,\psi,r}^{1,*} \oplus \Ker  d_{A,\varphi,\psi}^0
     \quad\text{(by \eqref{eq:Range_Intersection_When_Complex}),}
\end{align*}
and thus, by the definitions \eqref{eq:DefineH0ForSO3Monopoles} and \eqref{eq:DefineH2ForSO3Monopoles_Perturbed},
\begin{equation}
\label{eq:CokernelIsSumOfCohomology}
\Ker \left(d_{A,\varphi,\psi,r}^1+d_{A,\varphi,\psi}^{0,*}\right)^*
=
\bH_{A,\varphi,\psi,r}^2 \oplus \bH_{A,\varphi,\psi}^0.
\end{equation}
Because $(A,\varphi,\psi)$ solves \eqref{eq:SO(3)_monopole_equations_almost_Hermitian_perturbed_intro_regular},
the forthcoming Theorem \ref{thm:Feehan_Leness_1998jdg_3-7_regularized_Taubes-perturbed_W1p} implies that $(A,\varphi,\psi)$ is gauge-equivalent to a smooth solution. Lemma \ref{lem:H0Vanishing} and our hypothesis that $(\varphi,\psi) \not \equiv (0,0)$ imply that $\bH_{A,\varphi,\psi}^0=(0)$ as asserted in
\eqref{eq:Vanishing_Perturbed_H0}.
By combining the equality $\bH_{A,\varphi,\psi}^0=(0)$ with \eqref{eq:CokernelIsSumOfCohomology} and \eqref{eq:SO3MonopoleH2Isomorphism_perturbed}, we obtain \eqref{eq:SO3MonopoleH2Isomorphism_perturbed_H2version}, completing the proof of Proposition \ref{prop:SO3MonopoleDeformationCorollary_Perturbed}.
\end{proof}

\chapter{Generalizations of Donaldson's symplectic subspace criterion}
\label{chap:Analogue_Donaldson_symplectic_submanifold_criterion}
Our goal in this chapter is to prove the results stated in Section \ref{subsec:Generalization_Donaldson_symplectic_subspace_criterion}, namely Proposition \ref{mainprop:Donaldson_1996jdg_3}, Theorem \ref{mainthm:Donaldson_1996jdg_3_Hilbert_space_codomain}, Corollary \ref{maincorDonaldson_1996jdg_3_Hilbert_space_domain}, Proposition \ref{mainprop:Donaldson_1996jdg_3_Banach_space}, Corollary \ref{maincor:Adjoint_DonaldsonCriteria}, Theorem \ref{mainthm:Donaldson_1996jdg_3_Hilbert_space}, and Corollary \ref{maincor:Donaldson_1996jdg_3_Hilbert_space_non-self-adjoint} --- all of which generalize Donaldson's Proposition \ref{prop:Donaldson_1996jdg_3}.

Section \ref{sec:Proof_analogue_Donaldson_symplectic_submanifold_criterion} contains our proofs of Proposition \ref{mainprop:Donaldson_1996jdg_3_Banach_space}, which gives a symplectic subspace criterion for the kernel of a bounded real linear operator on a Banach space, and Corollary \ref{maincor:Adjoint_DonaldsonCriteria}, which gives a symplectic subspace criterion for the kernel of the adjoint of a bounded real linear operator on a Banach space.

In Section \ref{sec:Proof_generalized_Donaldson_symplectic_subspace_criterion_spectral_projection}, we prove Theorem \ref{mainthm:Donaldson_1996jdg_3_Hilbert_space} and Corollary \ref{maincor:Donaldson_1996jdg_3_Hilbert_space_non-self-adjoint}, which provide symplectic subspace criteria for the kernels of unbounded operators on Hilbert spaces. 

Numerous articles reference Proposition \ref{prop:Donaldson_1996jdg_3}, including see Auroux \cite{Auroux_2001}, Gironella, Mu\~noz, and Zhou \cite{Gironella_Munoz_Zhou_2023}, Giroux \cite{Giroux_2017}, Ibort and Torres \cite{Ibort_Torres_2004}, Mohsen \cite{Mohsen_2019}, Moriyama \cite{Moriyama_2008}, Sena--Dias \cite{Sena-Dias_2006}, and Shiffman and Zelditch \cite{Shiffman_Zelditch_2002}. Cieliebak and Mohnke include a proof of Proposition \ref{prop:Donaldson_1996jdg_3} in \cite[Lemma 8.3 (b), p. 328 and Remark 8.4, p. 329]{Cieliebak_Mohnke_2007}, along the geometric lines suggested by Donaldson in \cite[Section 1, pp. 668--669]{DonSympAlmostCx}, using Donaldson's concept of \emph{K\"ahler angle}. In Section \ref{sec:Proof_analogue_Donaldson_symplectic_submanifold_criterion_hypersurface}, we prove Proposition \ref{mainprop:Donaldson_1996jdg_3}, and hence Proposition \ref{prop:Donaldson_1996jdg_3} as a corollary, by more elementary methods of linear algebra.

We digress in Section \ref{sec:Symplectic_forms_almost_complex_structures_Hilbert_spaces} to discuss properties of symplectic forms and induced almost complex structures on Hilbert spaces that do not appear to be thoroughly treated elsewhere in the literature. Section \ref{sec:Transverse_intersections_symplectic_subspaces_almost_symplectic_manifolds} also contains a digression, this time to discuss transverse intersections of symplectic subspaces and almost symplectic manifolds, another topic that does not appear to be treated in the literature but which is essential for our applications. In Section \ref{sec:Proof_analogue_Donaldson_symplectic_submanifold_criterion_hypersurface_Hilbert_spaces}, we apply Proposition \ref{mainprop:Donaldson_1996jdg_3} and the results of Section \ref{sec:Transverse_intersections_symplectic_subspaces_almost_symplectic_manifolds} to prove Theorem \ref{mainthm:Donaldson_1996jdg_3_Hilbert_space_codomain} and Corollary \ref{maincorDonaldson_1996jdg_3_Hilbert_space_domain}, which provide more flexible symplectic subspace criteria for the kernels of unbounded operators on Hilbert spaces. 

In Section \ref{sec:Weakly_strongly_non-degenerate_bilinear_forms_Banach_spaces} we review the concepts of weakly and strongly non-degenerate bilinear forms on (infinite-dimensional) vector spaces. For our proofs of Theorem \ref{mainthm:Donaldson_1996jdg_3_Hilbert_space} and Corollary \ref{maincor:Donaldson_1996jdg_3_Hilbert_space_non-self-adjoint}, we review in Section \ref{sec:Spectral_theory_unbounded_operators} the spectral theory that we shall need for unbounded operators.
% TL6-22-2021: Do we still use this material on eigenvalue gaps?
%PF11-10-2025 Yes (for this chapter anyway, but we will likely split that off as a separate paper)
The results that we review in Section \ref{sec:Weyl_asymptotic_formula_eigenvalues_elliptic_operator} generalize the asymptotic formula due to Weyl \cite{Weyl_1912} for eigenvalues of the Laplacian on functions over a domain in Euclidean space. In particular, Theorem \ref{thm:Weyl_asymptotic_formula_eigenvalues_elliptic_pseudodifferential_operator_sections_vector_bundle} ensures that useful gap conditions as in the hypothesis \eqref{eq:Spectral_gap_T} of Theorem \ref{mainthm:Donaldson_1996jdg_3_Hilbert_space} are satisfied.

\section{Generalizations of Donaldson's symplectic subspace criterion to bounded real linear
operators on Banach spaces}
\label{sec:Proof_analogue_Donaldson_symplectic_submanifold_criterion}
In this section, we prove Proposition \ref{mainprop:Donaldson_1996jdg_3_Banach_space} and Corollary \ref{maincor:Adjoint_DonaldsonCriteria}. 

\begin{proof}[Proof of Proposition \ref{mainprop:Donaldson_1996jdg_3_Banach_space}]
We first observe that the expressions in \eqref{eq:Complex_linear_and_anti-linear_operator_components} for $\sT'$ and $\sT''$ are indeed the complex linear and antilinear components of $\sT$, respectively, since
\[
  \sT'J = \frac{1}{2}(\sT - j\sT J)J = \frac{1}{2}(\sT J + j\sT) = \frac{1}{2}j(-j\sT J + \sT) = j\sT',
\]
while 
\[
  \sT''J = \frac{1}{2}(\sT + j\sT J)J = \frac{1}{2}(\sT J - j\sT) = \frac{1}{2}j(-j\sT J - \sT) = -j\sT''.
\]
Since $J$ is $g$-orthogonal by hypothesis, we observe that the restriction of the continuous bilinear form $\omega_0:\sX\times\sX\to\RR$ in \eqref{eq:pre-symplectic_form_Banach_space_domain} to the subspace $\Ker \sT\times \Ker \sT$ obeys
\begin{multline*}
  \omega_0(u,v) = g(\pi_{\Ker \sT} Ju,v) = g(Ju,v) = -g(u,Jv) = -g(Jv,u) = -g(\pi_{\Ker \sT} Jv,u) = -\omega_0(v,u),
  \\
  \text{for all } u,v \in \Ker \sT.
\end{multline*}
Thus, $\omega_0:\Ker \sT\times \Ker \sT \to \RR$ is a skew-symmetric continuous bilinear form, in other words, $\omega_0$ is a two-form on $\Ker \sT$. Moreover, the bilinear form $\omega_0(\cdot,J\cdot):\sX\times\sX\to\RR$ is symmetric since
\begin{multline*}
  \omega_0(u,Jv) = g(\pi_{\Ker \sT} Ju,Jv) = g(\pi_{\Ker \sT} Ju,\pi_{\Ker \sT}Jv) = g(\pi_{\Ker \sT} Jv,\pi_{\Ker \sT}Ju)
  \\
  = g(\pi_{\Ker \sT} Jv,Ju) = \omega_0(v,u), \quad \text{for all } u,v \in \sX.
\end{multline*}
Consequently, the bilinear form $\omega_0(\cdot,J\cdot):\Ker \sT\times \Ker \sT\to\RR$ is symmetric. We claim that if $\sT$ obeys our hypothesis \eqref{eq:Donaldson_1996jdg_Prop_3_Banach_space}, then it obeys the positivity condition,
\begin{equation}
  \label{eq:McDuff-Salamon_positivity_condition}
  \omega_0(v,Jv) > 0, \quad\text{for all } v\in\Ker \sT \less\{0\},
\end{equation}
so that $J$ is compatible with $\omega_0$ on $\Ker \sT$ in the sense of McDuff and Salamon \cite[Equation (2.5.3), p. 63]{McDuffSalamonSympTop3}. If the claim \eqref{eq:McDuff-Salamon_positivity_condition} is true, then 
\begin{equation}
  \label{eq:McDuff-Salamon_inner_product}
  g_0 := \omega_0(\cdot,J\cdot) \quad\text{on } \Ker \sT,
\end{equation}
defines a weak inner product on $\Ker \sT$, by analogy with McDuff and Salamon \cite[Equation (2.5.4), p. 63]{McDuffSalamonSympTop3}. Given the claim \eqref{eq:McDuff-Salamon_positivity_condition}, the map $\Ker \sT \ni u \mapsto \omega_0(u,\cdot) \in (\Ker \sT)^*$ is injective, since $\omega_0(u,v)=0$ for all $v\in \Ker \sT$ implies that $\omega_0(u,Ju) = g_0(u,u) = 0$ and so $u=0$. Therefore, the claim \eqref{eq:McDuff-Salamon_positivity_condition} ensures that $\omega_0$ is a weakly non-degenerate two-form on $\Ker \sT$.

Observe that, for each $v \in \sX$,
\begin{align*}
  \omega_0(v,Jv)
  &= g(\pi_{\Ker \sT} Jv, Jv) \quad\text{(by \eqref{eq:pre-symplectic_form_Banach_space_domain})}
  \\
  &= g(\pi_{\Ker \sT} Jv, \pi_{\Ker \sT} Jv)
  \\
  &= \|\pi_{\Ker \sT} Jv\|_\sH^2.
\end{align*}
Therefore, for each $v \in \sX$,
\begin{equation}
  \label{eq:omega_0_compatible_with_J_iff_piJ_injective}
  \omega_0(v,Jv) > 0 \iff \pi_{\Ker \sT} Jv \neq 0.
\end{equation}
Hence, we make the further claim that if $\sT$ obeys the hypothesis \eqref{eq:Donaldson_1996jdg_Prop_3_Banach_space} of Proposition \ref{mainprop:Donaldson_1996jdg_3_Banach_space}, then it satisfies the following condition:
\begin{equation}
  \label{eq:Donaldson_1996jdg_proposition_3_pi_Jv_neq_0_for_v_in_KerT}
  \pi_{\Ker \sT} Jv \neq 0, \quad\text{for all } v \in (\Ker \sT)\less\{0\}.
\end{equation}
If the claim \eqref{eq:Donaldson_1996jdg_proposition_3_pi_Jv_neq_0_for_v_in_KerT} is true, then the equivalence \eqref{eq:omega_0_compatible_with_J_iff_piJ_injective} will imply that the positivity condition \eqref{eq:McDuff-Salamon_positivity_condition} holds and thus $\omega_0$ will be a weakly non-degenerate two-form on $\Ker \sT$ by our previous analysis.

Therefore, to complete our proof of Proposition \ref{mainprop:Donaldson_1996jdg_3_Banach_space}, it suffices to verify \eqref{eq:Donaldson_1996jdg_proposition_3_pi_Jv_neq_0_for_v_in_KerT}. We observe that
\begin{equation}
  \label{eq:Tprime_is_Tprimeprime_on_JKerT}
  \sT' = \sT'' \quad\text{on } J\Ker \sT.
\end{equation}
To see that \eqref{eq:Tprime_is_Tprimeprime_on_JKerT} holds, observe that if $u \in J\Ker \sT$, then $Ju \in J^2\Ker \sT = \Ker \sT$ and so
\[
  0 = \sT Ju = (\sT' + \sT'')Ju = j(\sT' - \sT'')u,
\]
which gives $\sT'u = \sT''u$ and thus \eqref{eq:Tprime_is_Tprimeprime_on_JKerT}. Consequently, $\sT = \sT' + \sT'' = 2\sT''$ on $J\Ker \sT$ by \eqref{eq:Tprime_is_Tprimeprime_on_JKerT} and hence
\begin{equation}
  \label{eq:Tprimeprime_is_half_T_on_JKerT}
  \sT'' = \frac{1}{2}\sT \quad\text{on } J\Ker \sT.
\end{equation}
%PF7-2-2025 Define what we mean by an "isomorphism of Banach spaces"
By hypothesis, the bounded operator $\sT \in \Hom(\sX,\sY)$ has a partial left inverse $L \in \Hom(\sY, \sX)$ in the sense of \eqref{eq:LT_is_1_minus_pi_KerT}, namely
\[
  L\sT = \pi_{\sX_0} \quad\text{on } \sX,
\]
where $\pi_{\sX_0} \in \End(\sX)$ is the continuous projection from $\sX = \Ker\sT\oplus\sX_0$ onto the closed subspace $\sX_0$. Now suppose, contrary to \eqref{eq:Donaldson_1996jdg_proposition_3_pi_Jv_neq_0_for_v_in_KerT}, that there exists $v \in (\Ker \sT)\less\{0\}$ such that $\pi_{\Ker \sT} Jv = 0$ or, equivalently, $u := Jv \in \sX_0$. By \eqref{eq:Tprimeprime_is_half_T_on_JKerT}, we have
\[
  \sT''u = \frac{1}{2}\sT u,
\]
and therefore, by the preceding equality and \eqref{eq:LT_is_1_minus_pi_KerT},
\[
  L\sT''u = \frac{1}{2}L\sT u = \frac{1}{2}\pi_{\sX_0}u = \frac{1}{2}u.
\]
Hence,
\[
  \|L\sT''\|_{\End(\sX)} = \sup_{x \in \sX \less\{0\}} \frac{\|L\sT''x\|_\sY}{\|x\|_\sX}
  \geq \frac{\|L\sT''u\|_\sY}{\|u\|_\sX} = \frac{1}{2},
\]
contradicting our hypothesis \eqref{eq:Donaldson_1996jdg_Prop_3_Banach_space} that $\|L\sT''\|_{\End(\sX)} < 1/2$. This proves Proposition \ref{mainprop:Donaldson_1996jdg_3_Banach_space}.     
\end{proof}

\begin{rmk}[Induced almost structures and compatible inner products]
\label{rmk:Induced_almost_structures_compatible_inner_products}
For the non-degenerate two-form $\omega_0$ produced by Proposition \ref{mainprop:Donaldson_1996jdg_3_Banach_space}, the pair $(\Ker \sT,\omega_0)$ is a \emph{symplectic vector space} in the sense of Cannas da Silva \cite[Section 1.2, Definition 1.3, p. 5]{Cannas_da_Silva_lectures_on_symplectic_geometry}. The real inner product $g$ on $\sH$ restricts to a real inner product on $\Ker \sT$, so by \cite[Section 12.2, Proposition 12.3, p. 84]{Cannas_da_Silva_lectures_on_symplectic_geometry}, the pair $(\omega_0,g)$ on $\Ker \sT$ defines an almost complex structure $J_0$ on $\Ker \sT$ as in \cite[Section 12.2, Definition 12.1, p. 84]{Cannas_da_Silva_lectures_on_symplectic_geometry} such that $(\omega_0,J_0)$ on $\Ker \sT$ is a \emph{compatible pair} in the sense of \cite[Section 12.2, Definition 12.2, p. 84]{Cannas_da_Silva_lectures_on_symplectic_geometry}, that is, $g_0 := \omega_0(\cdot,J_0\cdot)$ is a positive inner product on $\Ker \sT$. Consequently, $(\Ker \sT,J_0)$ is an \emph{almost complex vector space}. Moreover, $(g_0,J_0,\omega_0)$ on $\Ker \sT$ is a \emph{compatible triple} in the sense of McDuff and Salamon \cite[Section 4.1, p. 153]{McDuffSalamonSympTop3}. The compatible inner product $g_0$ is related to the given inner product $g$ by the displayed identity in \cite[Section 12.2, paragraph following proof of Proposition 12.3, p. 85]{Cannas_da_Silva_lectures_on_symplectic_geometry}.
\qed\end{rmk}

% PF1-17-2025 Is this still needed?
%PF12-18-2025 I commented out as it appears vesitigial and is not referenced. Please cut if you don't need it.
% \begin{rmk}[Preserving pre-existing almost complex structures]
% \label{rmk:PreservingPreExistingACStructures}
% Let $J_0$ be the almost complex structure on $\Ker \sT$ constructed in Remark \ref{rmk:Induced_almost_structures_compatible_inner_products}. If there is a subspace $V\subset \Ker \sT$ which is $J$-invariant, then we claim that $J_0\restriction V=J$. To see this, observe that for $v_1,v_2\in V$, the bilinear form $\omega_0$ satisfies
% %PF1-6-2024 @ TL: I don't understand. What's \omega below? It was not defined.
% $\omega_0(v_1,v_2) = g(\pi_{\Ker \sT} Jv_1,v_2) = g(Jv_1,v_2) = \omega(v_1,v_2)$ so
% $\omega_0\restriction V = \omega\restriction V$. Hence, the almost complex structure $J_0$ defined by $\omega_0$ equals $J$.
% \qed\end{rmk}

\begin{rmk}[Hypothesis \eqref{eq:Donaldson_1996jdg_Prop_3_Banach_space} implies Donaldson's hypothesis \eqref{eq:Donaldson_1996jdg_proposition_3_approx_complex_linear}]
\label{rmk:Relation_between_Donaldson_and_FL_symplectic_subspace_criteria}
The hypothesis \eqref{eq:Donaldson_1996jdg_Prop_3_Banach_space} in Proposition \ref{mainprop:Donaldson_1996jdg_3_Banach_space} is clearly implied by the stronger inequality
\begin{equation}
  \label{eq:Donaldson_1996jdg_Prop_3_Banach_space_stronger}
  \|\sT''\|_{\Hom(\sX,\sY)} < \frac{1}{2}\|L\|_{\Hom(\sY,\sX)}^{-1}.
\end{equation}
Moreover, noting that $\sX = \Ker \sT \oplus \sX_0$ as a direct sum of Banach spaces by hypothesis of Proposition \ref{mainprop:Donaldson_1996jdg_3_Banach_space}, we see that
\begin{multline*}
  \|L\|_{\Hom(\sY,\sX)}
  = \sup_{y\in\sY \less \{0\}} \frac{\|Ly\|_\sX}{\|y\|_\sY}
  \geq \sup_{y\in\Ran \sT \less \{0\}} \frac{\|Ly\|_\sX}{\|y\|_\sY}
  = \sup_{x\in\sX_0\less\{0\}} \frac{\|L\sT x\|_\sX}{\|\sT x\|_\sY}
  \\
  \geq \sup_{x\in\sX_0 \less \{0\}} \frac{\|x\|_\sX}{\|\sT\|_{\Hom(\sX,\sY)}\|x\|_\sX}
  = \|\sT\|_{\Hom(\sX,\sY)}^{-1},
\end{multline*}
and so we obtain
\[
  \|L\|_{\Hom(\sY,\sX)}^{-1} \leq \|\sT\|_{\Hom(\sX,\sY)}.
\]
Hence, the inequality \eqref{eq:Donaldson_1996jdg_Prop_3_Banach_space_stronger} and the preceding estimate imply that
\[
  \|\sT''\|_{\Hom(\sX,\sY)} < \frac{1}{2}\|\sT\|_{\Hom(\sX,\sY)}
  \leq \frac{1}{2}\|\sT'\|_{\Hom(\sX,\sY)} + \frac{1}{2}\|\sT''\|_{\Hom(\sX,\sY)},
\]
and thus
\[
  \|\sT''\|_{\Hom(\sX,\sY)} < \|\sT'\|_{\Hom(\sX,\sY)}.
\]
However, Donaldson's hypothesis \eqref{eq:Donaldson_1996jdg_proposition_3_approx_complex_linear} does not imply \eqref{eq:Donaldson_1996jdg_Prop_3_Banach_space} in general.
\qed\end{rmk}

Before proceeding to the proof of Corollary \ref{maincor:Adjoint_DonaldsonCriteria}, we first establish the following technical

\begin{lem}[Dual of a direct sum of Banach subspaces is the direct sum of the annihilators of the subspaces]
\label{lem:Dual_direct_sum}  
If $\sY$ is a Banach space with continuous dual space $\sY^*$ and $\sY = \sZ \oplus \sW$ is a direct sum of Banach subspaces, then $\sY^* = \sW^\perp \oplus \sZ^\perp$. Moreover, if $\pi_\sZ:\sY \to \sZ$ and $\pi_\sW:\sY \to \sW$ are the continuous projections on $\sY$, then their Banach space adjoint operators are given by
\[
  \pi_\sZ^* = \pi_{\sW^\perp} \quad\text{and}\quad \pi_\sZ^* = \pi_{\sZ^\perp},
\]
where $\pi_{\sW^\perp}:\sY^* \to \sW^\perp$ and $\pi_{\sZ^\perp}:\sY^* \to \sZ^\perp$ are the continuous projections on $\sY^*$ 
\end{lem}

\begin{proof}[Proof of Lemma \ref{lem:Dual_direct_sum}]
% COMMENT Compare https://math.stackexchange.com/questions/217596/the-dual-of-the-direct-sum
We define a bounded linear operator $S:\sY^* \to \sW^\perp \oplus \sZ^\perp$ by setting
\[
  Sy^* := (y^*\circ\pi_\sZ, y^*\circ\pi_\sW), \quad \text{for all } y^* \in \sY^*,
\]
where $\pi_\sZ:\sZ \oplus \sW \to \sZ$ and $\pi_\sW:\sZ \oplus \sW \to \sW$ are the continuous projections. We observe that
\[
  (y^*\circ\pi_\sZ)w = 0, \quad \text{for all } w \in \sW,
  \quad\text{and}\quad
  (y^*\circ\pi_\sW)z = 0, \quad \text{for all } z \in \sZ,
\]
so $y^*\circ\pi_\sZ \in \sW^\perp$ and $y^*\circ\pi_\sW \in \sZ^\perp$ and thus $S$ is well defined. Clearly, $S$ is injective since
\[
  0 = Sy^* = (y^*\circ\pi_\sZ, y^*\circ\pi_\sW) = y^*\circ(\pi_\sZ + \pi_\sW) = y^*,
\]
and so $Sy^* = 0 \implies y^* = 0$ and $S$ is injective. If $(w^*,z^*) \in \sW^\perp \oplus \sZ^\perp$, we define $y^* := w^* + z^* \in \sY^*$ and observe that
\begin{multline*}
  Sy^* = ((w^* + z^*)\circ\pi_\sZ, ((w^* + z^*)\circ\pi_\sW) = (w^*\circ\pi_\sZ, z^*\circ\pi_\sW)
  \\
  = (w^*\circ(\pi_\sZ + \pi_\sW), z^*\circ(\pi_\sZ + \pi_\sW)) = (w^*,z^*),
\end{multline*}
and so $S$ is surjective. Hence, $S$ is an isomorphism by the Open Mapping Theorem  (see, for example, Brezis \cite[Section 2.3, Theorem 2.6, p. 35]{Brezis}) and this completes the proof of the first assertion.

For the second assertion, if $\alpha \in \sZ^\perp$, then
\[
  (\pi_\sZ^*\alpha)(y) = \alpha(\pi_\sZ y) = 0, \quad\text{for all } y \in \sY,
\]
and thus $\pi_\sZ^*\alpha = 0$ for all $\alpha \in \sZ^\perp$. On the other hand, if $\beta \in \sW^\perp$, then
\[
  (\pi_\sZ^*\beta)(y) = \beta(\pi_\sZ y) = \beta(\pi_\sZ y + \pi_\sW y) = \beta(y), \quad\text{for all } y \in \sY,
\]
and thus $\pi_\sZ^*\beta = \beta$ for all $\beta \in \sW^\perp$. Hence, $\pi_\sZ^* = 0$ on $\sZ^\perp$ and $\pi_\sW^* = \id$ on $\sW^\perp$ and so $\pi_\sZ^* = \pi_{\sW^\perp}$. The same argument shows that $\pi_\sW^* = \pi_{\sZ^\perp}$ and this completes the proof of Lemma \ref{lem:Dual_direct_sum}.
\end{proof}  

We now give the

\begin{proof}[Proof of Corollary \ref{maincor:Adjoint_DonaldsonCriteria}]
The almost complex structures $J$ and $j$ on $\sX$ and $\sY$, respectively, define almost complex structures $J^*$ and $j^*$ on the continuous dual spaces $\sX^*$ and $\sY^*$, respectively, in the standard way (see Huybrechts \cite[Lemma 1.2.6, p. 26]{Huybrechts_2005}) by
\begin{equation}
  \label{eq:AdjointConj}
  (J^*\alpha)x := \alpha(Jx) \quad\text{and}\quad (j^*\beta)(y) := \beta(jy),
  \quad\text{for all } \alpha \in \sX^*, \beta \in \sY^* \text{ and } x \in \sX, y \in \sY.  
\end{equation}
For the operators $\sT',\sT''\in\Hom(\sX,\sY)$ in \eqref{eq:Complex_linear_and_anti-linear_operator_components},
we apply \eqref{eq:AdjointConj} to compute their Banach space operator adjoints in terms of the complex linear and antilinear components of $\sT^* \in \Hom(\sY^*,\sX^*)$:
\begin{subequations}
  \label{eq:Tprime_or_primeprime*_is_T*prime_or_primeprime}
  \begin{align}
    \label{eq:Tprime*_is_T*prime}
    (\sT')^* &= \frac{1}{2}(\sT - j\sT J)^* = \frac{1}{2}(\sT^* - J^*\sT^*j^*) =(\sT^*)',
    \\
    \label{eq:Tprimeprime*_is_T*primeprime}
    (\sT'')^* &= \frac{1}{2}(\sT+j\sT J)^* = \frac{1}{2}(\sT^* + J^*\sT^*j^*) =(\sT^*)''.
\end{align}
\end{subequations}
We note that because $\sG$ is a dense subspace of $\sY$ by hypothesis, the Banach space operator adjoint $\iota^*:\sY^* \to \sG^*$ of the continuous embedding $\iota:\sG \to \sY$ has $\Ker\iota^* = (0)$ by Rudin \cite[Corollary (b) to Theorem 4.12, p. 99]{Rudin} and so $\iota^*$ is a continuous embedding (as asserted in the hypotheses of Corollary \ref{maincor:Adjoint_DonaldsonCriteria}).

We also note that because $\sY = \Ran \sT \oplus \sY_0$ by hypothesis, then $\sY^* = \sY_0^\perp \oplus (\Ran \sT)^\perp = \sY_0^\perp \oplus \Ker \sT^*$, where the first equality follows from Lemma \ref{lem:Dual_direct_sum} and the second equality follows from Rudin \cite[Theorem 4.12, p. 99]{Rudin}.

Since $\sT \in \Hom(\sX,\sY)$ has closed range by hypothesis of Proposition \ref{mainprop:Donaldson_1996jdg_3_Banach_space}, the operator $\sT^* \in \Hom(\sY^*,\sX^*)$ also has closed range by Abramovich and Aliprantis \cite[Section 2.1, Theorem 2.18, p. 76]{Abramovich_Aliprantis_2002}. Therefore, $\sT^*$ has a partial left inverse $Q \in \Hom(\sX^*,\sY^*)$ such that $Q\sT^* = \pi_{\sY_0^*}$. Indeed, if $R \in \Hom(\sX,\sY)$ is a partial right inverse for $\sT$ such that $\sT R = \pi_{\Ran \sT}$, then
\[
  R^*\sT^* = \pi_{\Ran \sT}^* = \pi_{\sY_0^\perp},
\]
where the second equality follows from Lemma \ref{lem:Dual_direct_sum}. Thus, $Q = R^*$ and so if
\begin{equation}
  \label{eq:Donaldson_1996jdg_Prop_3_corrected_condition_Hilbert_space_adjoint}
  \|R^*(\sT^*)''\|_{\End(\sY^*)} < \frac{1}{2},
\end{equation}
then the final conclusion of Corollary \ref{maincor:Adjoint_DonaldsonCriteria} follows from Proposition \ref{mainprop:Donaldson_1996jdg_3_Banach_space} by replacing the role of $\sT$ by $\sT^*$. By \eqref{eq:Tprimeprime*_is_T*primeprime}, we have
\[
  R^*(\sT^*)'' = R^*(\sT'')^* = (\sT''R)^*
\]
and, by Rudin \cite[Theorem 4.10, p. 98]{Rudin}, one has
\[
  \|(\sT''R)^*\|_{\End(\sY^*)} = \|\sT''R\|_{\End(\sY)}, 
\]
so \eqref{eq:Donaldson_1996jdg_Prop_3_corrected_condition_Hilbert_space_adjoint} is equivalent to the hypothesis \eqref{eq:Donaldson_1996jdg_Prop_3_corrected_condition_Hilbert_space_right_inverse}. This completes the proof of Corollary \ref{maincor:Adjoint_DonaldsonCriteria}.
\end{proof}

\section[Donaldson's symplectic subspace criterion via spectral projections]{Generalizations of Donaldson's symplectic subspace criterion via spectral projections}
\label{sec:Proof_generalized_Donaldson_symplectic_subspace_criterion_spectral_projection}
In this section, we prove Theorem \ref{mainthm:Donaldson_1996jdg_3_Hilbert_space} and Corollary \ref{maincor:Donaldson_1996jdg_3_Hilbert_space_non-self-adjoint}. We begin with the

\begin{proof}[Proof of Theorem \ref{mainthm:Donaldson_1996jdg_3_Hilbert_space}]
Let $M\in [1,\infty)$ be a constant and $\Gamma$ be the counterclockwise contour contained in $\rho(T)$ and defined by
\begin{equation}
  \label{eq:Omega_spectrum_component_open_neighborhood}
  U := \left\{z\in\CC: |\Imag z| < M\sqrt{\mu} \text{ and } |\Real z| < \mu\right\}
  \quad\text{and}\quad
  \Gamma := \partial U,
\end{equation}
where the open rectangle $U \subset \CC$ and its boundary $\Gamma$ have their usual orientations. Our hypothesis \eqref{eq:Spectral_gap_T} on $\mu \in (0,\infty)$ and $\delta \in (0,\mu/2)$ requires that $T$ obey
\[
  (-\mu-2\delta,-\mu+2\delta) \cup (\mu-2\delta,\mu+2\delta) \subset \rho(T)
\]
and thus
\begin{equation}
  \label{eq:Spectrum_T}
  \sigma(T) \subset (-\infty,-\mu-2\delta] \cup [-\mu+2\delta, \mu-2\delta] \cup [\mu+2\delta,\infty).
\end{equation}
By hypothesis, the operator $T'' = T - T' \in \End(\sH)$ is bounded and obeys \eqref{eq:Donaldson_1996jdg_Prop_3_Hilbert_space}, namely
\[
  \|T''\|_{\End(\sH)} < \delta/2.
\]
Observe that \eqref{eq:Spectrum_T} ensures
%PF12-16-2025 Add pictures
\[
  z \in (-\mu-\delta,-\mu+\delta) \cup (\mu-\delta,\mu+\delta) \implies \dist(z,\sigma(T)) \geq \delta.
\]
% TL6-20-2025: I think this is the hypothesis used to apply \eqref{eq:Second_Neumann_series_radius_convergence_bounded_T-S_normal_T}, not the conclusion
%PF10-10-2025 I don't understand the point you are making here? We're just referring to inequalities.
%TL12-15-2025: \eqref{eq:Second_Neumann_series_radius_convergence_bounded_T-S_normal_T} (A.2.16) is stated as a hypothesis --the conclusion from it would be $z\in\varrho(T)$
%PF12-16-2025 Ah, I see. I rephrased (above/below). The sentence was too long and the order of qualifiers confusing.
Therefore,
\[
  \|T'-T\|_{\End(\sH)} < \delta/2 < \delta \leq \dist(z,\sigma(T)),
  \quad\text{for all } z \in (-\mu-\delta,-\mu+\delta) \cup (\mu-\delta,\mu+\delta), 
\]
and so it follows from \eqref{eq:Second_Neumann_series_radius_convergence_bounded_T-S_normal_T} (with $S = T'$) that \eqref{eq:Spectral_gap_Tprime} holds, namely
\begin{equation}
\label{eq:SpectralGap_for_T'}
  (-\mu-\delta,-\mu+\delta) \cup (\mu-\delta,\mu+\delta) \subset \rho(T').
\end{equation}
% COMMENT https://math.stackexchange.com/questions/153690/eigenprojection-as-contour-integral-over-resolvent
The (finite rank) Riesz projections \eqref{eq:Riesz_projection_closed_unbounded_operator_Banach_space} for $T$ and $T'$ relative to $U$ are
\begin{equation}
  \label{eq:Riesz_projections_T_and_T'}
  \Pi_\mu := -\frac{1}{2\pi i}\oint_\Gamma R(z,T)\,dz
  \quad\text{and}\quad
  \Pi_\mu' := -\frac{1}{2\pi i}\oint_\Gamma R(z,T')\,dz,
\end{equation}
so that $\Pi_\mu - \Pi_\mu'$ may be estimated in terms of the difference between the resolvent operators and hence between $T'' = T - T'$. Indeed, the second resolvent identity \eqref{eq:General_second_resolvent_identity} gives
\[
  R(z,T) - R(z,T')
  =
  R(z,T)(T' - T)R(z,T'), \quad\text{for all } z \in \rho(T)\cap \rho(T'),
\]
that is
\begin{equation}
  \label{eq:Second_resolvent_identity}
  R(z,T) - R(z,T') = -R(z,T)T''R(z,T'), \quad\text{for all } z \in \rho(T)\cap \rho(T').
\end{equation}
For $\omega_\mu$ as in \eqref{eq:pre-symplectic_form_Hilbert_space_domain} and each $v \in \sH$,
\[
  \omega_\mu(v,Jv)
  = g(\Pi_\mu Jv, Jv)
  = g(\Pi_\mu^2 Jv, Jv)
  = g(\Pi_\mu Jv, \Pi_\mu Jv)
  = \|\Pi_\mu Jv\|_\sH^2.
\]
Therefore, for each $v \in \sH$,
\begin{equation}
  \label{eq:omega_mu_compatible_with_J_iff_piJ_injective}
  \omega_\mu(v,Jv) > 0 \iff \Pi_\mu Jv \neq 0.
\end{equation}
For brevity, let $E_\mu$ denote the finite dimensional subspace $\Ran\Pi_\mu \subset \sH$. The conclusion of Theorem \ref{mainthm:Donaldson_1996jdg_3_Hilbert_space} is thus equivalent to the assertion that the operator $\Pi_\mu J\Pi_\mu \in \End(E_\mu)$ is \emph{injective}:
\begin{equation}
  \label{eq:Ker_Pi_J_Pi_is_zero}
  \Ker \Pi_\mu J\Pi_\mu = (0).
\end{equation}
Since $(\Pi_\mu J\Pi_\mu)^* = -\Pi_\mu J\Pi_\mu$ and $(\Ran \Pi_\mu J\Pi_\mu)^\perp = \Ker ((\Pi_\mu J\Pi_\mu)^*)$ by Rudin \cite[Theorem 4.12, p. 99]{Rudin}, we obtain
\[
  \left(\Ran \Pi_\mu J\Pi_\mu\right)^\perp = \Ker \Pi_\mu J\Pi_\mu.
\]
Therefore, $\Pi_\mu J\Pi_\mu  \in \End(E_\mu)$ is injective if and only if it is surjective and so the conclusion of Theorem \ref{mainthm:Donaldson_1996jdg_3_Hilbert_space} is thus equivalent to the assertion that 
\begin{equation}
  \label{eq:Pi_J_Pi_in_GL_RanPi}
  \Pi_\mu J\Pi_\mu \in \GL(E_\mu),
\end{equation}
recalling that because $\Pi_\mu J\Pi_\mu \in \End(E_\mu)$ is  bijective and $E_\mu$ has finite dimension, it is a linear isomorphism.
% COMMENT We'd only need the Open Mapping Theorem in Brezis \cite[Section 2.3, Theorem 2.6, p. 35]{Brezis} if it were infinite dimensional
The remainder of our proof of Theorem \ref{mainthm:Donaldson_1996jdg_3_Hilbert_space} focuses on proving the claim \eqref{eq:Pi_J_Pi_in_GL_RanPi}.

Let $E_\mu'$ denote the finite dimensional subspace $\Ran\Pi_\mu' \subset \sH$. We clearly have 
\begin{equation}
  \label{eq:Piprime_J_Piprime_in_GL_RanPiprime}
  \Pi_\mu' J\Pi_\mu' \in \GL(E_\mu'),
\end{equation}
since
\[
  -\Pi_\mu' J\Pi_\mu'\, \Pi_\mu' J\Pi_\mu'
  = - \Pi_\mu' J^2\Pi_\mu'
  = \Pi_\mu'
  = 1 \quad\text{on } E_\mu'.
\]
For brevity, write $A := \Pi_\mu J\Pi_\mu$ and $B := \Pi_\mu' J\Pi_\mu'$ and $A' := \Pi_\mu' A\Pi_\mu' \in \End(E_\mu')$. Then
\[
  A' = B + A' - B = B(1 - B^{-1}(A'-B)) \in \End(E_\mu'),
\]
and because $B \in \GL(E_\mu')$, to prove that we also have $A' \in \GL(E_\mu')$ it suffices to show that
\begin{equation}
  \label{eq:Norm_B_inv(A'_minus_B)_lessthan_1}
  \|B^{-1}(A'-B)\|_{\End(E_\mu')} < 1.
\end{equation}
(See, for example, Taylor \cite[Appendix A, Section 5, Proposition 5.11, p. 623]{Taylor_PDE1_3rd_edition}.) But
\[
  \|B^{-1}(A'-B)\|_{\End(E_\mu')} \leq \|B^{-1}(A'-B)\|_{\End(\sH)},
\]  
and the condition \eqref{eq:Norm_B_inv(A'_minus_B)_lessthan_1} is in turn implied by
\begin{equation}
  \label{eq:Norm_A'_minus_B_lessthan_1_over_norm_invB}
  \|A'-B\|_{\End(\sH)} < \|B^{-1}\|_{\End(\sH)}^{-1}.
\end{equation}
Noting that $B^{-1} = -\Pi_\mu' J\Pi_\mu' = -B$ by our proof of \eqref{eq:Piprime_J_Piprime_in_GL_RanPiprime} and
\[
  \|B\|_{\End(\sH)}
  =
  \sup_{\|v\|_\sH=1}\|Bv\|_\sH
  =
  \sup_{\|v\|_\sH=1}\|\Pi_\mu' J\Pi_\mu'v\|_\sH
  =
  \sup_{\|v\|_\sH=1}\|J\Pi_\mu'v\|_\sH
  =
  \sup_{\|v\|_\sH=1}\|\Pi_\mu'v\|_\sH
  =
  1,
\]
this means that the condition \eqref{eq:Norm_A'_minus_B_lessthan_1_over_norm_invB} is equivalent to
\begin{equation}
  \label{eq:Norm_A'_minus_B_lessthan_1}
  \|A'-B\|_{\End(\sH)} < 1.
\end{equation}
In other words, by substituting the definitions of $A'$ and $B$ into the inequality \eqref{eq:Norm_A'_minus_B_lessthan_1_over_norm_invB}, we see that to prove
\begin{equation}
  \label{eq:Piprime_Pi_J_Pi_Piprime_in_GL_RanPiprime}
  \Pi_\mu'\Pi_\mu J\Pi_\mu\Pi_\mu' \in \GL(E_\mu'),
\end{equation}
it suffices to show that
\begin{equation}
  \label{eq:Norm_Pi'PiJPiPi'_minus_Pi'JPi'_lessthan_1}
  \|\Pi_\mu'\Pi_\mu J\Pi_\mu\Pi_\mu' - \Pi_\mu' J\Pi_\mu'\|_{\End(\sH)} < 1.
\end{equation}
But
\begin{multline*}
  \Pi_\mu'\Pi_\mu J\Pi_\mu\Pi_\mu' - \Pi_\mu' J\Pi_\mu'
  =
  \Pi_\mu'\Pi_\mu J\Pi_\mu\Pi_\mu' - (\Pi_\mu')^2 J(\Pi_\mu')^2
  \\
  =
  \Pi_\mu'\Pi_\mu J\Pi_\mu\Pi_\mu' - \Pi_\mu'\Pi_\mu' J\Pi_\mu\Pi_\mu'
  + \Pi_\mu'\Pi_\mu' J\Pi_\mu\Pi_\mu' - (\Pi_\mu')^2 J(\Pi_\mu')^2
  \\
  =
  \Pi_\mu'(\Pi_\mu - \Pi_\mu') J\Pi_\mu\Pi_\mu'
  + \Pi_\mu'J(\Pi_\mu - \Pi_\mu')\Pi_\mu'.
\end{multline*}
By writing $\Pi_\mu'' : = \Pi_\mu - \Pi_\mu'$, we observe that to prove \eqref{eq:Norm_Pi'PiJPiPi'_minus_Pi'JPi'_lessthan_1} it suffices to show that
\begin{equation}
  \label{eq:Norm_Pi"_lessthan_one_half}
  \|\Pi_\mu''\|_{\End(\sH)} < 1/2.
\end{equation}
Indeed, to see this, observe that
\begin{align*}
  &\|\Pi_\mu'\Pi_\mu J\Pi_\mu\Pi_\mu' - \Pi_\mu' J\Pi_\mu'\|_{\End(\sH)}
  \\
  &\leq
  \|\Pi_\mu'(\Pi_\mu - \Pi_\mu') J\Pi_\mu\Pi_\mu'\|_{\End(\sH)}
    + \|\Pi_\mu'J(\Pi_\mu - \Pi_\mu')\Pi_\mu'\|_{\End(\sH)}
  \\
  &=
  \|\Pi_\mu'(\Pi_\mu - \Pi_\mu') J\Pi_\mu\Pi_\mu'\|_{\End(\sH)}
    + \|\Pi_\mu'(\Pi_\mu - \Pi_\mu')\Pi_\mu'\|_{\End(\sH)}
  \\
  &\leq \|\Pi_\mu - \Pi_\mu'\|_{\End(\sH)}
    + \|\Pi_\mu - \Pi_\mu'\|_{\End(E_\mu')}
  \\
  &\leq 2\|\Pi_\mu''\|_{\End(\sH)},
\end{align*}
and thus \eqref{eq:Norm_Pi'PiJPiPi'_minus_Pi'JPi'_lessthan_1} follows from \eqref{eq:Norm_Pi"_lessthan_one_half}, as claimed. But \eqref{eq:Piprime_Pi_J_Pi_Piprime_in_GL_RanPiprime} follows from \eqref{eq:Norm_Pi'PiJPiPi'_minus_Pi'JPi'_lessthan_1} and so \eqref{eq:Piprime_Pi_J_Pi_Piprime_in_GL_RanPiprime} also follows from \eqref{eq:Norm_Pi"_lessthan_one_half}.

It remains to show that \eqref{eq:Pi_J_Pi_in_GL_RanPi} follows from \eqref{eq:Norm_Pi"_lessthan_one_half}. To see this, observe that if
\begin{equation}
  \label{eq:Pi_mu_Emuprime_to_Emu_and_Pi_muprime_Emu_to_Emuprime}
  \Pi_\mu \in \Hom(E_\mu',E_\mu) \quad\text{and} \quad \Pi_\mu' \in \Hom(E_\mu,E_\mu')
\end{equation}
are isomorphisms of real vector spaces, then \eqref{eq:Pi_J_Pi_in_GL_RanPi} follows from \eqref{eq:Piprime_Pi_J_Pi_Piprime_in_GL_RanPiprime} and hence follows from \eqref{eq:Norm_Pi"_lessthan_one_half}. To prove that the homomorphisms \eqref{eq:Pi_mu_Emuprime_to_Emu_and_Pi_muprime_Emu_to_Emuprime} are isomorphisms, we will show that
\begin{equation}
  \label{eq:Pi_muprime_Pi_mu_in_GLEmuprime_and_Pi_mu_Pi_muprime_in_GLEmu}
  \Pi_\mu'\Pi_\mu \in \GL(E_\mu') \quad\text{and}\quad \Pi_\mu\Pi_\mu' \in \GL(E_\mu).
\end{equation}
To prove the first inclusion in \eqref{eq:Pi_muprime_Pi_mu_in_GLEmuprime_and_Pi_mu_Pi_muprime_in_GLEmu}, we write
\[
  \Pi_\mu\Pi_\mu' = \Pi_\mu\Pi_\mu' - \id_{E_\mu} + \id_{E_\mu} \in \End(E_\mu),
\]
and observe that the inequality \eqref{eq:Norm_Pi"_lessthan_one_half} yields
\begin{multline*}
  \|\Pi_\mu\Pi_\mu' - 1\|_{\End(E_\mu)} = \|\Pi_\mu\Pi_\mu' - \Pi_\mu^2\|_{\End(E_\mu)}
  \leq \|\Pi_\mu\|_{\End(E_\mu)}\|\Pi_\mu' - \Pi_\mu\|_{\End(E_\mu)}
  \\
  = \|\Pi_\mu''\|_{\End(E_\mu)} \leq \|\Pi_\mu''\|_{\End(\sH)} < 1/2. 
\end{multline*}
Therefore, $\Pi_\mu'\Pi_\mu \in \GL(E_\mu')$ as claimed and one can show that $\Pi_\mu\Pi_\mu' \in \GL(E_\mu)$ by a similar argument. This proves \eqref{eq:Pi_muprime_Pi_mu_in_GLEmuprime_and_Pi_mu_Pi_muprime_in_GLEmu}. 

The remainder of our proof of Theorem \ref{mainthm:Donaldson_1996jdg_3_Hilbert_space} focuses on proving that the inequality \eqref{eq:Norm_Pi"_lessthan_one_half} holds. By definition \eqref{eq:Riesz_projections_T_and_T'} of the Riesz projections and the second resolvent identity \eqref{eq:Second_resolvent_identity}, we obtain
%PF5-8-2025 Missing factors of 1/2\pi i; R(z,T) = (z-T)^{-1} or (T-z)^{-1}? Consistent? 
\[
  \Pi_\mu - \Pi_\mu'
  =
  -\frac{1}{2\pi i}\oint_\Gamma (R(z,T) - R(z,T'))\,dz
  =
  \frac{1}{2\pi i}\oint_\Gamma R(z,T)T''R(z,T')\,dz,
\]
and therefore,
\begin{equation}
  \label{eq:Norm_Pi"_leq_integral_norm_second_resolvent_expression}
  \|\Pi_\mu''\|_{\End(\sH)}
  \leq
  \frac{1}{2\pi}\int_\Gamma \|R(z,T)T''R(z,T')\|_{\End(\sH)}\,|dz|.
\end{equation}
By hypothesis, the resolvent $R(z,T)\in \End(\sH)$ is a compact operator. Because $R(z,T)$ is compact and $T'' \in \End(\sH)$ is bounded by hypothesis, then $R(z,T)T''R(z,T')$ is compact (see, for example, Conway \cite[ChapterVI, Section 3, Proposition 3.5, p. 174]{Conway_course_functional_analysis}) and so the second resolvent identity \eqref{eq:Second_resolvent_identity} implies that the resolvent $R(z,T')$ is also a compact operator.

By hypothesis $T \in \End(\sH)$ is self-adjoint with compact resolvent, so it has a spectrum of eigenvalues $\sigma(T) = \{\lambda_k\}_{k=1}^\infty \subset \RR$ with finite multiplicity and no accumulation points and its resolvent has a spectrum of eigenvalues with finite multiplicity,
\[
  \sigma(R(z,T)) = \left\{(\lambda_k - z)^{-1}\right\}_{k=1}^\infty \subset \CC, \quad\text{for all } z \in \rho(T).
\]
The equality \eqref{eq:Spectral_radius_equals_norm_resolvent_operator_normal_T} yields 
\[
  \|R(z,T)\|_{\End(\sH)} = \dist(z,\sigma(T))^{-1}, \quad\text{for all } z \in \rho(T).
\]
Also by hypothesis, $T' \in \End(\sH)$ is self-adjoint and it has compact resolvent as noted above, so it has a spectrum of eigenvalues $\sigma(T') = \{\lambda_k'\}_{k=1}^\infty \subset \RR$ with finite multiplicity and no accumulation points and its resolvent has a spectrum of eigenvalues with finite multiplicity,
\[
  \sigma(R(z,T')) = \left\{(\lambda_k' - z)^{-1}\right\}_{k=1}^\infty \subset \CC, \quad\text{for all } z \in \rho(T).
\]
For $z$ in the portion of the boundary $\Gamma$ given by $\Real z = \mu$, we have
\[
  |z - \lambda_k|^2 = |\Imag z|^2 + |\mu - \lambda_k|^2 \geq |\Imag z|^2 + \delta^2,
  \quad\text{for all } \lambda_k \in \sigma(T) \cap U,
\]  
and similarly for $z$ in the portion of the boundary $\Gamma$ given by $\Real z = -\mu$. Thus,
\[
  |\lambda_k - z| \geq (|\Imag z|^2 + \delta^2)^{1/2},
  \quad\text{for all } \lambda_k \in \sigma(T) \cap U,
\]
and so the resolvent operator norm equality \eqref{eq:Spectral_radius_equals_norm_resolvent_operator_normal_T} yields
\begin{equation}
  \label{eq:Norm_RzT_real_w_is_mu}
  \|R(z,T)\|_{\End(\sH)} = \max_k |\lambda_k - z|^{-1} \leq \frac{1}{(|\Imag z|^2 + \delta^2)^{1/2}},
  \quad\text{for all } z \in \Gamma\cap\left\{w \in \CC: |\Real w| = \mu\right\}.
\end{equation}
For $z$ in the portion of the boundary $\Gamma$ given by $|\Imag z| = M\sqrt{\mu}$, we have
\[
  |z - \lambda_k| \geq |\Imag z| = M\sqrt{\mu}, \quad\text{for all } \lambda_k \in \sigma(T) \cap \Omega.
\]  
and so \eqref{eq:Spectral_radius_equals_norm_resolvent_operator_normal_T} yields
\begin{equation}
  \label{eq:Norm_RzT_imag_w_is_sqrt_mu}
  \|R(z,T)\|_{\End(\sH)} = \max_k |\lambda_k - z|^{-1} \leq 1/(M\sqrt{\mu}),
  \quad\text{for all } z \in \Gamma\cap\left\{w \in \CC: |\Imag w| = M\sqrt{\mu}\right\}.
\end{equation}
By the spectral gap result \eqref{eq:Spectral_gap_Tprime}, the analogous operator norm estimates hold for the resolvent $R(z,T')$.

The inequality \eqref{eq:Norm_Pi"_leq_integral_norm_second_resolvent_expression} and the operator norm bound
\begin{multline*}
  \|R(z,T)T''R(z,T')\|_{\End(\sH)}
  \leq
  \|R(z,T)\|_{\End(\sH)}\|T''\|_{\End(\sH)}\|R(z,T')\|_{\End(\sH)},
  \\
  \text{for all } z \in \rho(T)\cap\rho(T'),
\end{multline*}
yield
\begin{equation}
  \label{eq:Norm_Pi"_leq_norm_T"_integral_norm_RzT_norm_RzT'}
  \|\Pi_\mu''\|_{\End(\sH)}
  \leq
  \|T''\|_{\End(\sH)}\frac{1}{2\pi}\int_\Gamma \|R(z,T)\|_{\End(\sH)}\|R(z,T')\|_{\End(\sH)}\,|dz|.
\end{equation}
Observe that
\begin{align*}
  &\int_\Gamma \|R(z,T)\|_{\End(\sH)}\|R(z,T')\|_{\End(\sH)}\,|dz|
  \\
  &= \int_{\Gamma\cap\left\{z\in\CC:\, |\Imag z| = M\sqrt{\mu}\right\}} \|R(z,T)\|_{\End(\sH)}\|R(z,T')\|_{\End(\sH)}\,|dz|
  \\
  &\qquad + \int_{\Gamma\cap\{z\in\CC:\, |\Real z| = \mu \}} \|R(z,T)\|_{\End(\sH)}\|R(z,T')\|_{\End(\sH)}\,|dz|.
\end{align*}
Therefore,
\begin{align*}
  &\int_\Gamma \|R(z,T)\|_{\End(\sH)}\|R(z,T')\|_{\End(\sH)}\,|dz|
  \\
  &\leq \int_{\Gamma\cap\left\{z\in\CC:\, |\Imag z| = M\sqrt{\mu}\right\}} \mu^{-1}\,|dz|
    + \int_{\Gamma\cap\{z\in\CC:\, |\Real z| = \mu\}} \frac{1}{|\Real z|^2 + \delta^2}\,|dz|
  \\
  &\qquad \text{(by \eqref{eq:Norm_RzT_real_w_is_mu} and \eqref{eq:Norm_RzT_imag_w_is_sqrt_mu}
    for $T$ and $T'$)}
  \\
  &= \mu^{-1}\Length\left(\Gamma\cap\left\{z\in\CC:|\Imag z| = M\sqrt{\mu}\right\}\right)
    + 2\int_{-M\sqrt{\mu}}^{M\sqrt{\mu}} \frac{1}{t^2 + \delta^2}\,dt
  \\
  &= 2(M^2\mu)^{-1}(2\mu) + 2(2/\delta)\arctan\left(M\sqrt{\mu}/\delta\right)
  \\
  &= 4/M^2 + 4\delta^{-1}\arctan\left(M\sqrt{\mu}/\delta\right)
  \\
  &\leq 4/M^2 + 4\delta^{-1}(\pi/2).
\end{align*}
The evaluation of the integral follows by an elementary calculation for constants $a,b \in \RR$ with $a/b > 0$:
\[
  \int\frac{1}{a + bt^2}\,dt = \frac{1}{\sqrt{ab}}\arctan(t\sqrt{b/a}).
\]
(See also Gradshteyn and Ryzhik \cite[Equation (2.124) (1), p. 71]{GradshteynRyzhik8}.) Thus,
\begin{equation}
  \label{eq:Bound_integral_norm_RzT_norm_RzT'}
  \int_\Gamma \|R(z,T)\|_{\End(\sH)}\|R(z,T')\|_{\End(\sH)}\,|dz|
  \leq
  4/M^2 + 2\pi/\delta.
\end{equation}
Substituting the preceding bound into the inequality \eqref{eq:Norm_Pi"_leq_norm_T"_integral_norm_RzT_norm_RzT'} yields
\[
  \|\Pi_\mu''\|_{\End(\sH)}
  \leq
  \frac{1}{2\pi}\left(4/M^2 + 2\pi/\delta\right)\|T''\|_{\End(\sH)}.
\]
Since $M \in [1,\infty)$ is arbitrary, we obtain
\[
  \|\Pi_\mu''\|_{\End(\sH)}
  \leq
  \delta\|T''\|_{\End(\sH)}.
\]
Therefore, provided $T''$ obeys the hypothesis \eqref{eq:Donaldson_1996jdg_Prop_3_Hilbert_space}, namely,
\[
  \|T''\|_{\End(\sH)} < \delta/2,
\]  
we obtain
\[
  \|\Pi_\mu''\|_{\End(\sH)} < 1/2,
\]
and so the condition \eqref{eq:Norm_Pi"_lessthan_one_half} is satisfied. This concludes the proof of Theorem \ref{mainthm:Donaldson_1996jdg_3_Hilbert_space}.
\end{proof}
%PF1-16-2025 Add picture showing how we get the resolvent estimates.

Before proceeding to the proof of Corollary \ref{maincor:Donaldson_1996jdg_3_Hilbert_space_non-self-adjoint}, we need the forthcoming elementary Lemma \ref{lem:Eigenvalues_densely_defined_unbounded_linear_operators}. Its content is surely well-known, though we are unable to find a reference that is complete.

\begin{lem}[Eigenvalues of densely defined, unbounded linear operators]
\label{lem:Eigenvalues_densely_defined_unbounded_linear_operators}  
Let $\sH_1, \sH_2$ be Hilbert spaces, $\sT:\sH_1 \to \sH_2$ be a closed, densely defined, unbounded linear operator, and assume that the self-adjoint (and thus closed) operator
\[
  T := \begin{pmatrix}0 & \sT^* \\ \sT & 0\end{pmatrix} \in \End(\sH_1\oplus\sH_2)
\]
has compact resolvent. Then the following hold: 
\begin{enumerate}
\item\label{item:Eigenvalues_densely_defined_unbounded_linear_operator_T}
  The spectrum of $T$ is a discrete subset of $\RR$ and consists of eigenvalues with finite multiplicity.
  
\item\label{item:Eigenvalues_densely_defined_unbounded_linear_operators_spectra_Laplacians}
  The operators $\sT^*\sT$ and $\sT\sT^*$ are self-adjoint and their spectra are discrete subsets of $[0,\infty)$ and consist of eigenvalues with finite multiplicity.
  
\item\label{item:Eigenvalues_densely_defined_unbounded_linear_operators_TT*_from_T*T}
  If $\{\nu_k\}_{k=1}^\infty$ is the non-decreasing sequence of eigenvalues of $\sT^*\sT$, repeated according to their multiplicities and excluding the zero eigenvalue, and $\{\varphi_k\}_{k=1}^\infty$ is the corresponding orthonormal basis of eigenvectors for $(\Ker\sT^*\sT)^\perp \cap \sH_1$, then the sequence of eigenvalues of $\sT\sT^*$, repeated according to their multiplicities and excluding the zero eigenvalue, is given by $\{\nu_k\}_{k=1}^\infty$, while $\{\psi_k\}_{k=1}^\infty$ with $\psi_k := \nu_k^{-1/2}\sT\varphi_k$ for $k\geq 1$ is the corresponding orthonormal basis of eigenvectors for  $(\Ker\sT\sT^*)^\perp \cap \sH_2$. Moreover, $\varphi_k := \nu_k^{-1/2}\sT^*\psi_k$ for all $k\geq 1$.
  
\item\label{item:Eigenvalues_densely_defined_unbounded_linear_operators_T+T*}
  The sequence of eigenvalues of $T \in \End(\sH_1\oplus\sH_2)$, repeated according to their multiplicities and excluding the zero eigenvalue, is given by $\{\pm\sqrt{\nu_k}\}_{k=1}^\infty$, while $\{\{\frac{1}{\sqrt{2}}(\varphi_k, \pm \psi_k)\}_{k=1}^\infty$ is the corresponding orthonormal basis of eigenvectors for $(\Ker T)^\perp \cap (\sH_1\oplus\sH_2)$.  
\end{enumerate}
\end{lem}

\begin{proof}
  % TL12-123-2025: Changed Item to Conclusion
  %PF12-15-2025 I don't like this change.
Consider Item \eqref{item:Eigenvalues_densely_defined_unbounded_linear_operator_T}. By hypothesis, $\sT$ is a densely defined, unbounded linear operator, so its adjoint $\sT^*$ is a well-defined closed unbounded linear operator and, because $\sT$ is closed, then $\sT^*$ is densely defined and $\sT^{**} = \sT$ (see Section \ref{sec:Spectral_theory_unbounded_operators}). Hence, the unbounded linear operator
\[
  T = \begin{pmatrix}0 & \sT^* \\ \sT & 0\end{pmatrix} \in \End(\sH_1\oplus\sH_2)
\]
is densely defined, self-adjoint (and thus closed). Because $T$ has compact resolvent, its spectrum $\sigma(T)$ is a discrete subset of $\RR$ consisting entirely of eigenvalues of $T$ with finite multiplicities (see Section \ref{sec:Spectral_theory_unbounded_operators}). This completes the proof of Item \eqref{item:Eigenvalues_densely_defined_unbounded_linear_operator_T}.
  
Consider Item \eqref{item:Eigenvalues_densely_defined_unbounded_linear_operators_spectra_Laplacians}. We observe that
\[
  T^2 = \begin{pmatrix}\sT^*\sT & 0 \\ 0 & \sT\sT^* \end{pmatrix} \in \End(\sH_1)\oplus\End(\sH_2).
\]
Since $\sT^* = \sT$ by the proof of Item \eqref{item:Eigenvalues_densely_defined_unbounded_linear_operator_T}, we see that the operators $\sT^*\sT \in \End(\sH_1)$ and $\sT\sT^* \in \End(\sH_2)$ are self-adjoint. Moreover, we have the following relationship among spectra,
\[
  \sigma(T^2) = \sigma(\sT^*\sT) \cup \sigma(\sT\sT^*) \subset \RR.
\]
The Spectral Mapping Theorem for polynomial functions of unbounded closed operators on a Banach space (see, for example, Dunford and Schwartz \cite[Chapter VII, Section 9, Theorem 10, p. 604]{Dunford_Schwartz1}) implies that $\sigma(T^2) = (\sigma(T))^2$. Hence, the operators $\sT^*\sT$ and $\sT\sT^*$ have spectra of non-negative eigenvalues with finite multiplicity. This completes the proof of Item \eqref{item:Eigenvalues_densely_defined_unbounded_linear_operators_spectra_Laplacians}.
  
Consider Item \eqref{item:Eigenvalues_densely_defined_unbounded_linear_operators_TT*_from_T*T}. For each $k, l \in \ZZ_{\geq 1}$, we have
\[
  \langle\sT\varphi_k,\sT\varphi_l\rangle_{\sH_2} = \langle\sT^*\sT\varphi_k,\varphi_l\rangle_{\sH_2}
  = \nu_k\langle\varphi_k,\varphi_l\rangle_{\sH_2} = \nu_k\delta_{kl},
\]
and so $\langle\psi_k,\psi_l\rangle_{\sH_2} = \nu_k\delta_{kl}$ for $\psi_k := \nu_k^{-1/2}\sT\varphi_k$ for all $k\geq 1$. Thus, $\{\psi_k\}_{k=1}^\infty$ is an orthonormal subset of $\Ran\sT \subset \sH_2$. Observe that $\Ker(\sT\sT^*) = \Ker\sT^*$ and $\dim \Ker(\sT\sT^*) < \infty$ by Item \eqref{item:Eigenvalues_densely_defined_unbounded_linear_operators_spectra_Laplacians}. But $(\Ran\sT)^\perp = \Ker\sT^*$ by Rudin \cite[Theorem 4.12, p. 99]{Rudin} and therefore $\dim(\sH_2/\Ran\sT) < \infty$, so $\Ran\sT$ is closed by Abramovich and Aliprantis \cite[Section 2.1, Corollary 2.17, p. 76]{Abramovich_Aliprantis_2002}. Thus, $\sH_2 = \Ran\sT \oplus \Ker\sT^*$ as an orthogonal direct sum of closed subspaces. Consequently, $\sT: (\Ker\sT)^\perp \cap \sH_1 \to (\Ker\sT^*)^\perp \cap \sH_2$ is a bijective continuous linear map of Hilbert spaces and thus an isomorphism of Hilbert spaces by the Open Mapping Theorem (see, for example, Brezis \cite[Section 2.3, Theorem 2.6, p. 35]{Brezis}), so $\{\psi_k\}_{k=1}^\infty$ is an orthonormal basis of $(\Ker\sT^*)^\perp \cap \sH_2$. In particular, $\{\psi_k\}_{k=1}^\infty$ is an orthonormal basis of eigenvectors for $\sT\sT^*$ on $\sH_2$ with corresponding sequence $\{\nu_k\}_{k=1}^\infty$ of positive eigenvalues.

Finally, because $\psi_k = \nu_k^{-1/2}\sT\varphi_k$, for all $k\geq 1$, we have $\sT^*\psi_k = \nu_k^{-1/2}\sT^*\sT\varphi_k = \nu_k^{1/2}\varphi_k$ and therefore $\varphi_k = \nu_k^{-1/2}\sT^*\psi_k$, for all $k\geq 1$. This completes the proof of Item \ref{item:Eigenvalues_densely_defined_unbounded_linear_operators_TT*_from_T*T}.

Consider Item \eqref{item:Eigenvalues_densely_defined_unbounded_linear_operators_T+T*}. Write $\eps_k \in \{\pm 1\}$ for all $k\geq 1$ and observe that
\begin{multline*}
  \frac{1}{2}\langle(\varphi_k, \eps_k\psi_k),(\varphi_l, \eps_l\psi_l)\rangle_{\sH_1\oplus\sH_2}
  =
  \frac{1}{2}\langle\varphi_k, \varphi_l\rangle_{\sH_1} + \eps_k\eps_l\frac{1}{2}\langle \psi_k,\psi_l\rangle_{\sH_2}
  \\
  =
  \frac{1}{2}(1 + \eps_k\eps_l)\delta_{kl}
  =
  \begin{cases}
    1 &\text{if } k=l \text{ and } \eps_k\eps_l = 1,
    \\
    0 &\text{if } k\neq l \text{ or } \eps_k\eps_l = -1.
  \end{cases}
\end{multline*}
Consequently, the sequence $\{\frac{1}{\sqrt{2}}(\varphi_k, \pm\psi_k)\}_{k=1}^\infty$ is an orthonormal subset of $(\Ker T)^\perp\cap(\sH_1\oplus\sH_2)$. Conversely, if $\varphi+\psi = (\varphi,\psi) \in (\Ker T)^\perp\cap(\sH_1\oplus\sH_2)$, then
\begin{multline*}
  \varphi+\psi
  =
  \frac{1}{\sqrt{2}}\sum_{k=1}^\infty a_k\varphi_k + \frac{1}{\sqrt{2}}\sum_{l=1}^\infty b_l\psi_l
  \\
  =
  \frac{1}{2\sqrt{2}}\sum_{k=1}^\infty a_k(\varphi_k + \psi_k) + a_k(\varphi_k - \psi_k)
  + \frac{1}{2\sqrt{2}}\sum_{l=1}^\infty b_l(\varphi_l + \psi_l) - b_l(\varphi_l - \psi_l)
  \\
  =
  \frac{1}{\sqrt{2}}\sum_{k=1}^\infty \frac{1}{2}(a_k+b_k)(\varphi_k + \psi_k)
  + \frac{1}{\sqrt{2}}\sum_{k=1}^\infty \frac{1}{2}(a_k-b_k)(\varphi_k - \psi_k)
  \\
  =
  \frac{1}{\sqrt{2}}\sum_{k=1}^\infty c_k(\varphi_k + \psi_k)
  + \frac{1}{\sqrt{2}}\sum_{l=1}^\infty d_l(\varphi_l - \psi_l),
\end{multline*}
where $c_k := \frac{1}{2}(a_k+b_k)$ for all $k \geq 1$ and $d_l := \frac{1}{2}(a_l-b_l)$ for all $l \geq 1$. Hence, the sequence $\{\frac{1}{\sqrt{2}}(\varphi_k, \pm\psi_k)\}_{k=1}^\infty$ is an orthonormal basis of $(\Ker T)^\perp\cap(\sH_1\oplus\sH_2)$.

Since $\varphi_k = \nu_k^{-1/2}\sT^*\psi_k$ and $\psi_k = \nu_k^{-1/2}\sT\varphi_k$, we have
\[
  T\begin{pmatrix}\varphi_k\\ \psi_k\end{pmatrix}
  =
  \nu_k^{-1/2}\begin{pmatrix}0 & \sT^* \\ \sT & 0\end{pmatrix}
  \begin{pmatrix}\sT^*\psi_k\\ \sT\varphi_k\end{pmatrix}
  =
  \nu_k^{-1/2}\begin{pmatrix}\sT^*\sT\varphi_k \\ \sT\sT^*\psi_k \end{pmatrix}
  =
  \nu_k^{-1/2}\begin{pmatrix} \nu_k\varphi_k \\ \nu_k\psi_k \end{pmatrix}
  =
  \nu_k^{1/2}\begin{pmatrix}\varphi_k\\ \psi_k\end{pmatrix},
\]
and so $\frac{1}{\sqrt{2}}(\varphi_k, \psi_k)$ is an eigenvector of $T$ with eigenvalue $\nu_k^{1/2}$. Similarly,
\[
  T\begin{pmatrix}\varphi_k\\ -\psi_k\end{pmatrix}
  =
  \nu_k^{-1/2}\begin{pmatrix}0 & \sT^* \\ \sT & 0\end{pmatrix}
  \begin{pmatrix}\sT^*\psi_k\\ -\sT\varphi_k\end{pmatrix}
  =
  \nu_k^{-1/2}\begin{pmatrix}-\sT^*\sT\varphi_k \\ \sT\sT^*\psi_k \end{pmatrix}
  =
  \nu_k^{-1/2}\begin{pmatrix} -\nu_k\varphi_k \\ \nu_k\psi_k \end{pmatrix}
  =
  -\nu_k^{1/2}\begin{pmatrix}\varphi_k\\ -\psi_k\end{pmatrix},
\]
and thus $\frac{1}{\sqrt{2}}(\varphi_k, -\psi_k)$ is an eigenvector of $T$ with eigenvalue $-\nu_k^{1/2}$. This completes the proof of Item \eqref{item:Eigenvalues_densely_defined_unbounded_linear_operators_T+T*} and hence the proof of Lemma \ref{lem:Eigenvalues_densely_defined_unbounded_linear_operators}.
\end{proof}

\begin{rmk}[Symmetry and asymmetry of the spectrum of the Dirac operator]
\label{rmk:Symmetry_asymmetry_spectrum_Dirac_operator}  
Suppose that $(X,g)$ is a closed Riemannian spin manifold and that $D \in \End(L^2(S))$ is the corresponding Dirac operator on the spinor vector bundle $S$ over $X$. If $d$ is even, then $S = S^+\oplus S^-$ and $D = D^+ \oplus D^-$ with $D^\pm \in \Hom(L^2(S^\pm), L^2(S^\mp))$ by Ginoux \cite[Section 1.3, Proposition 1.3.2, p, 10]{Ginoux_dirac_spectrum}, and we may write
\[
  D = \begin{pmatrix}0 & D^- \\ D^+ & 0\end{pmatrix} \in \End\left(L^2(S^+) \oplus L^2(S^+)\right).
\]
The operator $D$ is a first-order elliptic operator that is $L^2$ self-adjoint by Ginoux \cite[Proposition 1.3.4, p, 10]{Ginoux_dirac_spectrum}, so $D = D^*$ and thus $D^\mp = (D^\pm)^*$, and is essentially self-adjoint\footnote{See Section \ref{sec:Spectral_theory_unbounded_operators}.} by Ginoux \cite[Proposition 1.3.5, p, 13]{Ginoux_dirac_spectrum}. The operator $D$ has dense domain $W^{1,2}(S) \subset L^2(S)$ and compact resolvent, by Friedrich \cite[Section 4.2, Propositions, p. 101]{Friedrich_Dirac_operators_Riemannian_geometry}.

We can thus apply Lemma \ref{lem:Eigenvalues_densely_defined_unbounded_linear_operators} and conclude from Item \eqref{item:Eigenvalues_densely_defined_unbounded_linear_operators_T+T*} that the spectrum $\sigma(D) \subset \RR$ is symmetric about the origin. This argument is essentially equivalent to the proof of \cite[Section 1.3, Theorem 1.3.7 (iv), p, 14]{Ginoux_dirac_spectrum} when $d$ is even; if $d \equiv 1 \pmod{4}$, the spectrum $\sigma(D) \subset \RR$ is also symmetric about the origin Ginoux \cite[Section 1.3, Theorem 1.3.7 (iv), p, 14]{Ginoux_dirac_spectrum}. However, if $d \equiv 3 \pmod{4}$, the spectrum $\sigma(D) \subset \RR$ is not necessarily symmetric about the origin and $\RR\PP^d$ provides a counterexample when $d \equiv 3 \pmod{4}$ for one of its spin structures (see Ginoux \cite[Section 1.3, p, 16]{Ginoux_dirac_spectrum}). Of course, this is the spectral asymmetry phenomenon explored by Atiyah, Patodi, and Singer \cite{APS1, APS2, APS3}.

When $d \equiv 0, 1$, or $2\pmod{4}$, the fact that $\sigma(D)$ is symmetric about the origin implies that $\sigma(D)$ is unbounded on both sides of $\RR$, since $\sigma(D^2) \subset [0,\infty)$ is unbounded, but $\sigma(D)$ is also unbounded on both sides of $\RR$ when $d \equiv 3 \pmod{4}$ by Ginoux \cite[Section 1.3, Theorem 1.3.7 (iv), p, 14]{Ginoux_dirac_spectrum}.
\qed
\end{rmk}  

We can now complete the

\begin{proof}[Proof of Corollary \ref{maincor:Donaldson_1996jdg_3_Hilbert_space_non-self-adjoint}]
We define $\sH := \sH_1\oplus\sH_2$ to be a Hilbert space with direct sum inner product $g = g_1\oplus g_2$ and $g$-orthogonal direct sum almost complex structure $J := J_1\oplus J_2$. As in the proof of Lemma \ref{lem:Eigenvalues_densely_defined_unbounded_linear_operators}, the operator $T \in \End(\sH)$ is a densely defined, self-adjoint (and hence closed) unbounded operator and by hypothesis it has compact resolvent. The complex linear and antilinear components of $T$ are given by
\[  
  T'
  =
  \begin{pmatrix}
    0 & (\sT^*)' \\ \sT' & 0
  \end{pmatrix}
  =
  \begin{pmatrix}
    0 & \sT^{\prime,*} \\ \sT' & 0
  \end{pmatrix}
  \quad\text{and}\quad
  T''
  =
  \begin{pmatrix}
    0 & (\sT^*)'' \\ \sT'' & 0
  \end{pmatrix}
  =
  \begin{pmatrix}
    0 & \sT^{\prime\prime,*} \\ \sT'' & 0
  \end{pmatrix}.
\]
Clearly, $T' \in \End(\sH)$ is self-adjoint and because $\sT'' \in \Hom(\sH_1,\sH_2)$ is bounded, the operator $T'' \in \End(\sH)$ is also bounded. Because the operator $T \in \End(\sH)$ is self-adjoint with compact resolvent, it has a spectrum of real eigenvalues with finite multiplicity and no accumulation points,
\[
  \sigma(T) = \{\lambda_k\}_{k=1}^\infty \subset \RR. 
\]
The resolvent set $\rho(T)$ obeys the spectral gap hypothesis \eqref{eq:Spectral_gap_T} and $\rho(T')$ obeys the spectral gap condition \eqref{eq:Spectral_gap_Tprime}. Since $\|(\sT'')^*\|_{\Hom(\sH_2,\sH_1)} = \|\sT''\|_{\Hom(\sH_1,\sH_2)}$ by Rudin \cite[Chapter 4, Theorem 4.10, p. 98]{Rudin}, the hypothesis \eqref{eq:Donaldson_1996jdg_Prop_3_Hilbert_space_non-self-adjoint} yields
\[
  \|T''\|_{\End(\sH)}
  =
  \|\sT''\|_{\Hom(\sH_1,\sH_2)} + \|(\sT'')^*\|_{\Hom(\sH_2,\sH_1)}
  =
  2\|\sT''\|_{\Hom(\sH_1,\sH_2)}
  <
  \delta/2.
\]
Thus, $T$ obeys \eqref{eq:Donaldson_1996jdg_Prop_3_Hilbert_space} and hence all the hypotheses of Theorem \ref{mainthm:Donaldson_1996jdg_3_Hilbert_space}. 

To prove Corollary \ref{maincor:Donaldson_1996jdg_3_Hilbert_space_non-self-adjoint} as a consequence of Theorem \ref{mainthm:Donaldson_1996jdg_3_Hilbert_space}, we note from Item \eqref{item:Eigenvalues_densely_defined_unbounded_linear_operators_T+T*} in Lemma \ref{lem:Eigenvalues_densely_defined_unbounded_linear_operators} that the spectrum of $T \in \End(\sH)$,
\[
  \sigma(T) = \{\lambda_k\}_{k=1}^\infty \subset \RR,
\]
and the spectra of $T^2 \in \End(\sH)$, and $\sT^*\sT \in \End(\sH_1)$, and $\sT\sT^* \in \End(\sH_2)$ are related by
\begin{equation}
  \label{eq:Spectrum_T^2_is_union_spectrum_sT*sT_union_spectrum_sTsT*}
  \begin{aligned}
    \Ker(T^2) &= \Ker(\sT^*\sT) \cup \Ker(\sT\sT^*),
    \\
    \sigma(T^2)\less\{0\} &= \sigma(\sT^*\sT)\less\{0\} = \sigma(\sT\sT^*)\less\{0\} \subset (0,\infty).
  \end{aligned}    
\end{equation}
We claim that
\begin{equation}
  \label{eq:Pi_mu_Hi_subset_Hi}
  \Pi_\mu\sH_i \subset \sH_i, \quad\text{for } i=1,2,
\end{equation}
where $\Pi_\mu$ by hypothesis is the orthogonal projection from $\sH$ onto the subspace spanned by the eigenvectors of $T$ with eigenvalues in $(-\mu,\mu)$ and we identify $\sH_1 \cong \sH_1\oplus(0) \subset \sH$ and $\sH_2 \cong (0)\oplus\sH_2 \subset \sH$. Thus,
\[
  \Pi_\mu = \sum_{|\lambda_k| < \mu} e_k\otimes\langle\cdot,e_k\rangle_\sH,
\]
where $\{e_k\}_{k=1}^\infty$ is an orthonormal basis of $\sH$ given by the eigenvectors $e_k$ of $T$ associated to the eigenvalues $\lambda_k \in \sigma(T)$ for $k\in\ZZ_{\geq 1}$. To prove \eqref{eq:Pi_mu_Hi_subset_Hi}, we observe by \eqref{eq:Spectrum_T^2_is_union_spectrum_sT*sT_union_spectrum_sTsT*} that each positive eigenvalue of $T^2$ is equal to a positive eigenvalue of both $\sT^*\sT$ and $\sT\sT^*$ and a zero eigenvalue of $T^2$ is equal to a zero eigenvalue of $\sT^*\sT$ or $\sT\sT^*$. Thus, each $\lambda_k^2$-eigenvector $e_k$ of $T^2$ is equal to a $\lambda_k^2$-eigenvector of $\sT^*\sT \in \End(\sH_1)$ or a $\lambda_k^2$-eigenvector of $\sT\sT^* \in \End(\sH_2)$. Therefore, 
\[
  \Pi_\mu v_i = \sum_{|\lambda_k| < \mu} e_k\otimes\langle v_i,e_k\rangle_\sH \in \sH_i, \quad\text{for all } v_i\in\sH_i \text{ and } i=1,2,
\]
and the claim \eqref{eq:Pi_mu_Hi_subset_Hi} follows. The inclusions \eqref{eq:Pi_mu_Hi_subset_Hi} yield a decomposition of $\Ran\Pi_\mu$ as an orthogonal direct sum of Hilbert subspaces,
\begin{equation}
  \label{eq:RanPi_mu_is_RanPi_mu_cap_H1_subset__RanPi_mu_cap_H2}
  \Ran\Pi_\mu = \pi_{\sH_1}\Ran\Pi_\mu \oplus \pi_{\sH_2}\Ran\Pi_\mu \subset \sH_1\oplus\sH_2,
\end{equation}
where $\pi_{\sH_i}$ are the orthogonal projections from $\sH = \sH_1\oplus\sH_2$ onto the factors $\sH_i$. The inclusions \eqref{eq:Pi_mu_Hi_subset_Hi} also imply that the continuous real bilinear form $\omega_\mu$ on $\sH$ in \eqref{eq:pre-symplectic_form_Hilbert_space_domain} provided by Theorem \ref{mainthm:Donaldson_1996jdg_3_Hilbert_space} has a block diagonal structure,
\[
  \omega_\mu
  =
  \begin{pmatrix}
    \omega_{1,\mu} & 0 \\ 0 & \omega_{2,\mu} 
  \end{pmatrix}
  \quad\text{on}\quad \sH_1\oplus\sH_2,
\]
where $\omega_{i,\mu} = g_i(\Pi_\mu J_i\cdot,\cdot)$ is a continuous real bilinear form on $\sH_i$ for $i=1,2$. In particular, the decomposition \eqref{eq:RanPi_mu_is_RanPi_mu_cap_H1_subset__RanPi_mu_cap_H2} and inclusion $J\sH_i \subset \sH_i$ for $i=1,2$ imply that the symplectic form $\omega_\mu$ on $\Ran\Pi_\mu$ provided by Theorem \ref{mainthm:Donaldson_1996jdg_3_Hilbert_space} has a block diagonal structure,
\[
  \omega_\mu
  =
  \begin{pmatrix}
    \omega_{1,\mu} & 0 \\ 0 & \omega_{2,\mu} 
  \end{pmatrix}
  \quad\text{on}\quad \pi_{\sH_1}\Ran\Pi_\mu \oplus \pi_{\sH_2}\Ran\Pi_\mu,
\]
and $\omega_{i,\mu}$ is a symplectic form on $\pi_{\sH_i}\Ran\Pi_\mu$. This completes the proof of Corollary \ref{maincor:Donaldson_1996jdg_3_Hilbert_space_non-self-adjoint}.
\end{proof}

\section{Donaldson's symplectic subspace criterion for  bounded real linear complex valued functions on complex Hilbert spaces}
\label{sec:Proof_analogue_Donaldson_symplectic_submanifold_criterion_hypersurface}
Before proceeding to the proof of Proposition \ref{mainprop:Donaldson_1996jdg_3}, we make the preliminary 

\begin{rmk}[Surjectivity of $\sL$ and orientation of $\Ker \sL$ in Proposition \ref{mainprop:Donaldson_1996jdg_3}]
\label{rmk:Donaldson_1996jdg_3_T_onto}  
The hypothesis \eqref{eq:Donaldson_1996jdg_proposition_3_approx_complex_linear_inequality} in Proposition \ref{mainprop:Donaldson_1996jdg_3} implies that $\sL$ is surjective. To see this, suppose first that
\[
  \|\sL''\|_{\Hom(\sH,\CC)} < \|\sL'\|_{\Hom(\sH,\CC)},
\]
in which case $\sL' \neq 0$. Because $\sL'$ is complex linear with codomain $\CC$, we see that $\sL'$ is surjective with right inverse $R' \in \Hom_\CC(\CC,\sH)$, so $\sL'R' = 1$ on $\CC$. We claim that $\sL R' \in \End_\RR(\CC)$ is invertible, in which case $\sL$ would be surjective as well. But $\sL R' = (\sL'+\sL'')R' = 1 + \sL''R'$, so it suffices to show that
\[
  \|\sL''R'\|_{\End(\CC)} < 1,
\]
and which is implied in turn by
\begin{equation}
  \label{eq:Donaldson_1996jdg_proposition_3_approx_complex_linear_inequality_Rprime}
  \|\sL''\|_{\Hom(\sH,\CC)} < \|R'\|_{\Hom(\CC,\sH)}^{-1}.
\end{equation}
But $\sL' \in \Hom_\CC(\sH,\CC) = \sH^*$ is represented by a unique non-zero vector $\tau' \in \sH$, so $\sL' = \langle\cdot, \tau'\rangle_{\sH} = (\tau')^*$ and $\|\sL'\|_{\Hom(\sH,\CC)} = \|\tau'\|_{\sH}$. Thus, $\sL^{\prime,*} = \tau' \in \Hom_\CC(\CC,\sH)$ (defined by $\CC \ni 1 \mapsto \tau' \in \sH$) and $\sL'\sL^{\prime,*} = (\tau')^*\tau' = \|\tau'\|_{\sH}^2 = \|\sL'\|_{\Hom(\sH,\CC)}^2$. Hence, we may write
\[
  R'
  =
  \sL^{\prime,*}(\sL'\sL^{\prime,*})^{-1}
  =
  \tau'\|\tau'\|_{\sH}^{-2}.
\]
Therefore,
\[
  \|R'\|_{\Hom(\CC,\sH)} = \|\tau'\|_{\sH}^{-1} = \|\sL'\|_{\Hom(\sH,\CC)}^{-1}.
\]
Consequently, the inequality \eqref{eq:Donaldson_1996jdg_proposition_3_approx_complex_linear_inequality_Rprime} is equivalent to 
\[
  \|\sL''\|_{\Hom(\sH,\CC)} < \|\sL'\|_{\Hom(\sH,\CC)},
\]
and so the hypothesis \eqref{eq:Donaldson_1996jdg_proposition_3_approx_complex_linear} implies that $\sL$ is surjective, as claimed. The orientation of $\Ker \sL$ is determined by the resulting isomorphism of real vector spaces,
\[
  \sH/\Ker \sL \cong \Ran \sL = \CC,
\]
and the complex orientations of $\sH$ and $\CC$. The same argument holds if $\|\sL'\|_{\Hom(\sH,\CC)} < \|\sL''\|_{\Hom(\sH,\CC)}$. 
\qed\end{rmk}

We now complete the

%PF12-16-2025 We may want to standardize our bounded/unbounded operator notation, as writing the "c" superscripts is tedious. Also, we may want to emphasize \CC or \RR linearity with writing \Hom, \End, etc.
\begin{proof}[Proof of Proposition \ref{mainprop:Donaldson_1996jdg_3}]
The assertion that $\sL$ is surjective follows from Remark \ref{rmk:Donaldson_1996jdg_3_T_onto}. We recall from the proof of Proposition \ref{mainprop:Donaldson_1996jdg_3_Banach_space} that by \eqref{eq:omega_0_compatible_with_J_iff_piJ_injective} each $v \in \sH$ obeys
\[
  \omega_0(v,Jv) > 0 \iff \pi_{\Ker \sL} Jv \neq 0,
\]
where $\omega_0$ is as in \eqref{eq:pre-symplectic_form_Banach_space_domain}. Hence, to prove that $\omega_0$ defines a symplectic form on $\Ker \sL$, it suffices to show that \eqref{eq:Donaldson_1996jdg_proposition_3_pi_Jv_neq_0_for_v_in_KerT} holds, namely
\[
  \pi_{\Ker \sL} Jv \neq 0, \quad\text{for all } v \in (\Ker \sL)\less\{0\}.
\]
We may suppose without loss of generality that $\sL'\neq 0$ and $\sL''\neq 0$ since, otherwise, the conclusion would follow immediately. Assume, to obtain a contradiction, that \eqref{eq:Donaldson_1996jdg_proposition_3_pi_Jv_neq_0_for_v_in_KerT} does not hold, so there exists $v \in (\Ker \sL) \less \{0\}$ such that $\pi_{\Ker \sL} Jv = 0$ or, equivalently, $Jv \in \Ran \sL^*$, noting the orthogonal decomposition $\sH = \Ker \sL \oplus \Ran \sL^*$. Because $\sL'=\sL''$ on $J\Ker \sL$ by \eqref{eq:Tprime_is_Tprimeprime_on_JKerT}, the vector $v \in (\Ker \sL)\less \{0\}$ obeys
\begin{equation}
  \label{eq:Tprime_Jv_is_Tprimeprime_Jv}
  \sL'Jv = \sL''Jv. 
\end{equation}
As in Remark \ref{rmk:Donaldson_1996jdg_3_T_onto}, we may write
%COMMENT See https://en.wikipedia.org/wiki/Riesz_representation_theorem
\begin{equation}
  \label{eq:Tprime_is_tauprime_functional_and_Tprimeprime_is_tauprimeprime_functional}
  \sL' = \langle\cdot,\tau'\rangle_\sH
  \quad\text{and}\quad
  \sL'' = \langle\tau'',\cdot\rangle_\sH,
\end{equation}  
for unique vectors $\tau'\in \sH\less\{0\}$ and $\tau''\in \sH\less\{0\}$ determined by the Riesz Representation Theorem (see, for example, Roman \cite[Theorem 13.32, p. 351]{Roman_advanced_linear_algebra} or the forthcoming Lemma \ref{lem:Norms_complex_linear_antilinear_functionals_on_complex_Hilbert_spaces}), with
\begin{equation}
  \label{eq:Norm_Tprime_is_norm_tauprime_and_norm_Tprimeprime_is_norm_tauprimeprime}
  \|\sL'\|_{\Hom(\sH,\CC)} = \|\tau'\|_\sH
  \quad\text{and}\quad
  \|\sL''\|_{\Hom(\sH,\CC)} = \|\tau''\|_\sH.
\end{equation} 
Hence, the identity \eqref{eq:Tprime_Jv_is_Tprimeprime_Jv} is equal to
\begin{equation}
  \label{eq:Jv_dot_tauprime_is_tauprimeprime_dot_Jv}
  \langle Jv,\tau'\rangle_\sH = \langle\tau'',Jv\rangle_\sH.
\end{equation}
Now, $\sL^{\prime,*} \in \Hom_\CC(\CC,\sH)$ is the complex linear map defined by $1\mapsto \tau'$, while $\sL^{\prime\prime,*} \in \Hom_{\bar\CC}(\CC,\sH)$ is the complex antilinear map defined by $1\mapsto \tau''$. By our assumption, we have $Jv = \sL^*w$ for some $w\in\CC\less\{0\}$ and because $\sL = \sL'+\sL'' \in \Hom_\RR(\sH,\CC)$, we obtain $\sL^* = \sL^{\prime,*} + \sL^{\prime\prime,*} \in \Hom_\RR(\CC,\sH)$ and
\[
  Jv = \sL^*w = \sL^{\prime,*}w + \sL^{\prime\prime,*}w = w\sL^{\prime,*}1 + \bar w \sL^{\prime\prime,*}1,
\]
that is,
\begin{equation}
  \label{eq:Jv_is_w_tauprime+barw_tauprimeprime}
  Jv = w\tau' + \bar w\tau''.
\end{equation}
Substituting this expression for $Jv$ into the identity \eqref{eq:Jv_dot_tauprime_is_tauprimeprime_dot_Jv} yields a series of equivalent identities:
\begin{align*}
  \langle Jv,\tau'\rangle_\sH
  &= \langle\tau'',Jv\rangle_\sH
  \\
  \iff
  \langle w\tau' + \bar w\tau'',\tau'\rangle_\sH
  &= \langle\tau'',w\tau' + \bar w\tau''\rangle_\sH \quad\text{(by \eqref{eq:Jv_is_w_tauprime+barw_tauprimeprime})}
  \\
  \iff
  w\|\tau'\|_\sH^2 + \bar w\langle \tau'',\tau'\rangle_\sH
  &= \bar w\langle\tau'',\tau'\rangle_\sH + w\|\tau''\|_\sH^2
  \\
  \iff
  w\|\tau'\|_\sH^2
  &= w\|\tau''\|_\sH^2
  \\
  \iff
  \|\tau'\|_\sH &= \|\tau''\|_\sH \quad\text{(by assumption that $Jv = \sL^*w \neq 0$ and thus $w\neq 0$)}
  \\
  \iff
  \|\sL'\|_{\Hom(\sH,\CC)} &= \|\sL''\|_{\Hom(\sH,\CC)} \quad\text{(by \eqref{eq:Norm_Tprime_is_norm_tauprime_and_norm_Tprimeprime_is_norm_tauprimeprime})}.
\end{align*}
The final equality above contradicts our hypothesis \eqref{eq:Donaldson_1996jdg_proposition_3_approx_complex_linear_inequality} (of course, it is also contradicted by the stronger condition \eqref{eq:Donaldson_1996jdg_proposition_3_approx_complex_linear}.) This proves Proposition \ref{mainprop:Donaldson_1996jdg_3}.
\end{proof}

\section{Symplectic forms and induced almost complex structures on Hilbert spaces}
\label{sec:Symplectic_forms_almost_complex_structures_Hilbert_spaces}
Given a symplectic form $\omega$ and inner product $g$ on a finite-dimensional real vector space $V$, there is an induced almost complex structure $J$ on $V$ that is compatible with $\omega$ and uniquely determined by the pair $(\omega,g)$ --- see, for example, Cannas da Silva \cite[Proposition 12.3, p. 84]{Cannas_da_Silva_lectures_on_symplectic_geometry}. However, her construction does not immediately extend to the case where $V$ is infinite-dimensional. In this section, we provide such a construction and, in addition, extend the linear algebra in Huybrechts \cite[Section 1.2]{Huybrechts_2005} from the setting of almost complex structures on finite-dimensional real vector spaces to infinite-dimensional real vector spaces.
%PF12-15-2025 Explain what happens in the 3 subsections too 

\subsection{Symplectic forms and almost complex structures on Hilbert spaces}
\label{subsec:Symplectic_forms_almost_complex_structures_Hilbert_spaces}
Let $\sH^\RR$ be a real Hilbert space, set $\sH^\CC := \sH^\RR\otimes_\RR\CC$, and note that
\[
  \Hom_\RR(\sH^\RR,\CC) = \Hom_\RR(\sH^\RR,\RR)\otimes_\RR\CC.
\]
We shall first introduce definitions of the Hermitian-orthogonal splittings
\[
  \sH^\CC = \sH^{1,0}\oplus\sH^{0,1}
\]
and
\[
  \Hom_\RR(\sH^\RR,\CC) = \Hom_\RR(\sH^\RR,\CC)^{1,0}\oplus\Hom_\RR(\sH^\RR,\CC)^{0,1}
\]
that provide more fundamental alternatives to the standard definitions which we recalled in \eqref{eq:sH_direct_sum_splitting_10_and_01_subspaces} and \eqref{eq:Hom_RR_sH_to_CC}, respectively, from Huybrechts \cite{Huybrechts_2005}. As in Section \ref{subsec:Generalization_Donaldson_symplectic_subspace_criterion}, the complex Hilbert space $\sH$ comprises an underlying real Hilbert space $\sH^\RR$ with real inner product $g$, a $g$-orthogonal almost complex structure $J$ on $\sH^\RR$ so that $(a+ib)\cdot v := av + bJv$ for all $v\in\sH^\RR$ and $a,b\in\RR$, a symplectic form $\omega = g(J\cdot,\cdot)$ on $\sH^\RR$, and a Hermitian inner product $h = g-i\omega$ on $\sH$.

Let us now suppose only that $(\sH^\RR,g,\omega)$ is a symplectic Hilbert space, that is, a real Hilbert space with inner product $g$ and symplectic structure $\omega$, both strongly nondegenerate in the sense of Section \ref{sec:Weakly_strongly_non-degenerate_bilinear_forms_Banach_spaces}. We shall first define the closed complex linear subspaces $\sH^{1,0}$ and $\sH^{0,1}$ of $\sH^\CC$ in terms of $(g,\omega)$ directly and \emph{then} define the induced almost complex structure $J = J_{g,\omega}$. For this purpose, we shall generalize the proof of \cite[Proposition 12.3, p. 84]{Cannas_da_Silva_lectures_on_symplectic_geometry} due to Cannas da Silva for the construction of $g$-orthogonal almost complex structures on finite-dimensional real vector spaces $(V,g,\omega)$ that are compatible with a symplectic structure $\omega$.

Because the real inner product $g$ and symplectic form $\omega$ are strongly non-degenerate by assumption on $\sH^\RR\times\sH^\RR$, they define real linear isomorphisms of real Hilbert spaces,
\begin{subequations}
  \label{eq:tilde_g_and_tilde_omega}
  \begin{align}
    \label{eq:tilde_g}
    g^\flat:\sH^\RR \ni v &\mapsto g(v,\cdot) \in \Hom_\RR(\sH^\RR,\RR),
    \\
    \label{eq:tilde_omega}
    \omega^\flat:\sH^\RR \ni v &\mapsto \omega(v,\cdot) \in \Hom_\RR(\sH^\RR,\RR),
\end{align}
\end{subequations}
and thus a $g$-skew-adjoint real linear isomorphism of real Hilbert spaces,
\begin{equation}
  \label{eq:A_is_tilde_g-inverse_tilde_omega}
  A := (g^\flat)^{-1}\circ\omega^\flat \in \End_\RR(\sH^\RR),
\end{equation}
or, equivalently (see Cannas \cite[Section 12.2, Proof of Proposition 12.3, p. 84]{Cannas_da_Silva_lectures_on_symplectic_geometry}),
\begin{equation}
  \label{eq:omega_is_g(A,.)}
  \omega(u,v) =: g(Au,v), \quad\text{for all } u,v \in \sH^\RR.
\end{equation} 
Since $A^* = -A$ with respect to the real inner product $g$ on $\sH^\RR$, then $AA^* = -A^2 = A^*A$ is a bounded normal operator on $\sH^\RR$. By Rudin \cite[Theorem 12.23, p. 324]{Rudin}, there is thus a spectral decomposition of $A$ such that
\begin{equation}
  \label{eq:Spectral_domposition_A}
  A = \int_{\sigma(A)} \lambda\,dE(\lambda),
\end{equation}
where $\sigma(A) \subset \CC$ is the spectrum of $A \in \End_\RR(\sH^\RR)$ and $E$ is a resolution of the identity on (the Borel subsets of) $\sigma(A)$ in the sense of Rudin \cite[Definition 12.17, p. 316]{Rudin}, so $E(\Omega) \in \End_\RR(\sH^\RR)$ is a $g$-orthogonal projection for each Borel subset $\Omega \subset \sigma(A)$. Since $A \in \End_\RR(\sH^\RR)$ is invertible, $0 \notin \sigma(A)$. The real inner product $g = \langle\cdot,\cdot\rangle$ on $\sH^\RR$ extends to a Hermitian inner product on $\sH^\CC$ by setting
\begin{equation}
  \label{eq:Hermitian_inner_product_sHCC_induced_by_real_inner_product_sHRR}
  \langle u+iv, x+iy\rangle := \langle u,x\rangle + \langle v,y\rangle
  - i\left(\langle u,y\rangle - \langle v,x\rangle\right),
  \quad\text{for all } u,v,x,y \in \sH^\RR.
\end{equation}
The bounded normal operator $A \in \End_\RR(\sH^\RR)$ extends by complex linearity to a bounded
operator $A \in \End_\CC(\sH^\CC)$. Moreover, $A$ is a skew-adjoint and thus normal operator on $\sH^\CC$ since, for all $u,v,x,y \in \sH^\RR$,
\begin{multline*}
  \langle u+iv, A^{*_\CC}(x+iy)\rangle
  =
  \langle A(u+iv), x+iy\rangle
  \\
  =
  \langle Au+iAv, x+iy\rangle
  =
  \langle u+iv, A^{*_\RR}x+iA^{*_\RR}y\rangle \quad\text{(by \eqref{eq:Hermitian_inner_product_sHCC_induced_by_real_inner_product_sHRR})},
\end{multline*}
thus $A^{*_\RR} = A^* \in \End_\RR(\sH^\RR)$ extends by complex linearity to a bounded
operator $A^{*_\CC} = A^{*_\RR} \in \End_\CC(\sH^\CC)$, which we denote simply by $A^*$, and therefore $A^*=-A$ on $\sH^\CC$ and thus $AA^* = -A^2 = A^*A$ on $\sH^\CC$, as claimed.

The bounded operator $S := iA \in \End_\CC(\sH^\CC)$ is self-adjoint since $S^* = (iA)^* = -iA^* = iA = S$, so $\sigma(S) \subset \RR$ by Rudin \cite[Theorem 12.15 (b), p. 314]{Rudin} and thus
\begin{equation}
  \label{eq:Spectrum_A_purely_imaginary}
  \sigma(A) \subset i\RR
\end{equation}
when $A$ is viewed as a bounded operator in $\End_\CC(\sH^\CC)$ and hence also as a bounded operator in $\End_\RR(\sH^\RR)$. Indeed, recall that $\lambda \in \sigma(A) \iff A-\lambda$ does not have a bounded inverse on $\sH^\RR$ by Rudin \cite[Definition 4.17 (c), p. 103]{Rudin}. Hence, if $A-\lambda$ has a bounded inverse $B$ on $\sH^\RR$ then, because $A-\lambda$ extends by complex linearity to a bounded operator on $\sH^\CC$, so does $B$ extend from $\sH^\RR$ to $\sH^\CC$ and $B$ is a bounded inverse for $A-\lambda$ on $\sH^\CC$. Hence, if $\lambda$ is in the spectrum of $A$ on $\sH^\RR$, then $\lambda$ is in the spectrum of $A$ on $\sH^\CC$ and this proves the claim that $\sigma(A) \subset i\RR$ when $A$ is viewed as a bounded operator on $\sH^\RR$.

For any $z \in \CC$, we have $(A-z)^* = A^*-\bar z = -(A+\bar z)$. But $A-\lambda$ is not invertible on $\sH^\RR$ if and only if $-(A-\lambda)^* = A+\bar\lambda$ is not invertible on $\sH^\RR$ and hence $\lambda \in \sigma(A) \iff -\bar\lambda \in \sigma(A)$. Because $0\notin \sigma(A)$, we may thus write
\begin{equation}
  \label{eq:Spectrum_A_disjoint_union_pos_neg_spectra_A}
  \sigma(A) = \sigma^+(A) \sqcup \sigma^-(A), 
\end{equation}
where $\sigma^+(A) := \sigma(A)\cap i\RR^+$ and $\sigma^-(A) := \sigma(A)\cap i\RR^-$, where $\RR^+ := (0,\infty)$ and $\RR^- := (-\infty,0)$. Complex conjugation yields a bijection from $\sigma^+(A)$ onto $\sigma^-(A)$. We define the corresponding Hermitian orthogonal projections,
\begin{equation}
  \label{eq:pi10_and_pi01_spectral_projections}
  \pi^{1,0} := E\left(\sigma^+(A)\right) \in \End_\CC(\sH^\CC)
  \quad\text{and}\quad
  \pi^{0,1} := E\left(\sigma^-(A)\right) \in \End_\CC(\sH^\CC),
\end{equation}
and the corresponding invariant subspaces for $A \in \End_\CC(\sH^\CC)$ by 
\begin{equation}
  \label{eq:sH10_and_sH01_invariant_subspaces_A}
  \sH^{1,0} := \pi^{1,0}\sH^\CC
  \quad\text{and}\quad
  \sH^{0,1} := \pi^{1,0}\sH^\CC,
\end{equation}
so that, using $E(\sigma(A)) = E(\sigma^+(A)) + E(\sigma^-(A)) = \id$ on $\sH^\CC$, 
\begin{equation}
  \label{eq:sHCC_Hermitian-orthogonal_direct_sum_sH10_and_sH01}
  \sH^\CC = \sH^{1,0} \oplus \sH^{0,1}
\end{equation}
is a Hermitian orthogonal splitting of $\sH^\CC$ (see, for example, Rudin \cite[Section 12.27, p. 327]{Rudin}). By analogy with the proof of Cannas da Silva \cite[Proposition 12.3, p. 84]{Cannas_da_Silva_lectures_on_symplectic_geometry}, we define a bounded Borel measurable function $f:\sigma(A) \to \CC$ by
\begin{equation}
  \label{eq:Definition_f_on_sigma(A)}
  f(\lambda) := (-\lambda^2)^{-1/2}\lambda, \quad\text{for all } \lambda \in \sigma(A),
\end{equation}
and observe that
\[
  f(\pm\mu i) = \pm i, \quad\text{for all } \mu > 0,
\]
which yields
\begin{equation}
  \label{eq:Values_f_on_sigma(A)}
  f = \pm i \quad\text{on } \sigma^\pm(A). 
\end{equation}  
By Rudin \cite[Section 12.24, Equation (1), p. 325]{Rudin}, we may define a bounded linear operator on $\sH^\CC$ by
\begin{equation}
  \label{eq:Definition_J_Borel_functional_calculus}
  J := f(A) = (-A^2)^{-1/2}A = \int_{\sigma(A)} f(\lambda)\,dE(\lambda)
   = \int_{\sigma(A)} (-\lambda^2)^{-1/2}\lambda\,dE(\lambda).
\end{equation}
The formula \eqref{eq:Definition_J_Borel_functional_calculus} for $J$ yields
\begin{equation}
  \label{eq:Eigen-decomposition_J_on_sHCC}
  J = i \quad\text{on } \sH^{1,0}
  \quad\text{and}\quad
  J = -i \quad\text{on } \sH^{0,1},
\end{equation}
so $J$ is an almost complex structure on $\sH^\CC$. Because $J(u+iv) = Ju + iJv$ by complex linearity, for all $u,v \in \sH^\RR$, we may regard $J$ on $\sH^\CC$ as a complex linear extension of a bounded linear operator $J$ on $\sH^\RR$ (compare Huybrechts \cite[Definition 1.2.4, p. 25]{Huybrechts_2005}) and as $J^2=-1$ on $\sH^\CC$, we obtain $J^2=-1$ on $\sH^\RR$. Hence, $J$ in \eqref{eq:Definition_J_Borel_functional_calculus} is an almost complex structure on $\sH^\RR$.

By Rudin \cite[Section 12.24, Equation (2), p. 325]{Rudin}, we have
\[
  \bar f(A) = f(A)^*,
\]
and so \eqref{eq:Values_f_on_sigma(A)} and \eqref{eq:Definition_J_Borel_functional_calculus} imply that $J^* = -J$ on $\sH^\CC$ and $JJ^* = 1 = J^*J$ on $\sH^\CC$, so that $J$ is skew-adjoint on $\sH^\CC$ and orthogonal with respect to the Hermitian inner product on $\sH^\CC$. Therefore, $J^* = -J$ on $\sH^\RR$ and $JJ^* = 1 = J^*J$ on $\sH^\RR$, so that $J$ is skew-adjoint on $\sH^\RR$ and is $g$-orthogonal. Hence, $J$ in \eqref{eq:Definition_J_Borel_functional_calculus} is a $g$-orthogonal almost complex structure on $\sH^\RR$.

Because $J = f(A)$ in \eqref{eq:Definition_J_Borel_functional_calculus} commutes with $A$ by Rudin \cite[Section 12.24, pp. 325--326]{Rudin}, then $J$ commutes with the spectral projections $E(\sigma^\pm(A))$ and so the complex linear subspaces $E(\sigma^\pm(A))\sH^\CC$ are $J$-invariant by Rudin \cite[Theorem 12.23, p. 324]{Rudin}.

Finally, we observe that $J$ is compatible with $\omega$ in the sense of Cannas da Silva \cite[Definition 12.2, p. 84]{Cannas_da_Silva_lectures_on_symplectic_geometry} since, for all $u,v \in \sH^\RR$, the definition \eqref{eq:A_is_tilde_g-inverse_tilde_omega} of $A$ yields
\[
  \omega(Ju,Jv) = g(AJu,Jv) = g(JAu,vJ) = g(Au,v) = \omega(u,v),
\]
and
\[
  \omega(u,Ju) = g(Au,Ju) = g(-JAu,u) = g((-A^2)^{-1/2}A^2u,u) = g((-A^2)^{1/2}u,u) > 0,
\]
where the positivity follows from the fact that $-A^2$ and hence $(-A^2)^{1/2}$ are positive-definite bounded $\RR$-linear operators on $\sH^\RR$.

By analogy with the statement and proof of Huybrechts \cite[Lemma 1.2.5, p. 26]{Huybrechts_2005}, who assumes that $\sH^\RR$ is finite-dimensional, we note that (writing $w = u+iv \in \sH^\CC$ for $u,v\in\sH^\RR$) the $\RR$-linear complex conjugation map,
\[
  C:\sH^\CC \ni u + iv \mapsto \overline{u + iv} = u - iv \in \sH^\CC,
\]
induces an $\RR$-linear isomorphism of complex Hilbert spaces,
\begin{equation}
  \label{eq:Conjugation_RR-linear_isomorphism_sH10_onto_sH01}
  C:\sH^{1,0} \to \sH^{0,1}.
\end{equation}
Indeed, for any $w = u+iv \in \sH^\CC$, we have
\[
  J\bar w = J(u-iv) = Ju - iJv = \overline{Ju + iJv} = \overline{Jw}.
\]
If $w \in \sH^{1,0}$, then $Jw = iw$ and $J\bar w = \overline{Jw} = \overline{iw} = -i\bar w$, so $\bar w \in \sH^{0,1}$ and if $w \in \sH^{0,1}$, then $Jw = -iw$ and $J\bar w = \overline{Jw} = \overline{-iw} = i\bar w$, so $\bar w \in \sH^{1,0}$. In particular, $C^2=1$ on $\sH^{1,0}$ and $\sH^{0,1}$, and thus $C$ is an $\RR$-linear isomorphism of complex Hilbert spaces.

Next, by analogy with Huybrechts \cite[Paragraph before Lemma 1.2.6, p. 26]{Huybrechts_2005}, the compositions
\[
  \sH \xrightarrow{\subset} \sH^\CC \xrightarrow{\pi^{1,0}} \sH^{1,0}
  \quad\text{and}\quad
  \sH \xrightarrow{\subset} \sH^\CC \xrightarrow{\pi^{0,1}} \sH^{0,1}
\]
are complex linear and antilinear isomorphisms of complex Hilbert spaces, respectively, where $\sH$ is $\sH^\RR$ equipped with the almost complex structure $J$ and both $\sH^{1,0}$ and $\sH^{0,1}$ have the complex structure $i$. Indeed, if $w \in \sH^{1,0}$, then $\pi^{1,0}u = w$ if $u = w + \bar w$, noting that $\bar w \in \sH^{0,1}$ by \eqref{eq:Conjugation_RR-linear_isomorphism_sH10_onto_sH01}, so $w = \pi^{1,0}w = \pi^{1,0}(w + \bar w) = \pi^{1,0}u$. But $\Real u = u$ and thus $u$ is in the image of the inclusion $\sH^\RR \subset \sH^\CC$. Hence, $\pi^{1,0}:\sH\to\sH^{1,0}$ is surjective. On the other hand, if $\pi^{1,0}u = 0 \in \sH^{1,0}$ for some $u\in\sH$ then, writing $u = w_1 + w_2 \in \sH^{1,0} \oplus \sH^{0,1} $ via \eqref{eq:sHCC_Hermitian-orthogonal_direct_sum_sH10_and_sH01}, we observe that $0 = \pi^{1,0}(w_1 + w_2) = w_1$. But $u = \bar u = \bar w_2$ and therefore $0 = \pi^{1,0}\bar w_2$, so $\bar w_2 =0$ and thus $w_2 = 0$. In particular, $u = w_1 +w_2 = 0$ and thus $\pi^{1,0}:\sH\to\sH^{1,0}$ is injective. Moreover, if $u \in \sH$, then
\[
  \pi^{1,0}Ju = J\pi^{1,0}u = i\pi^{1,0}u = i\iota^{1,0}u,
\]
and so $\iota^{1,0}$ is complex linear. Hence, $\pi^{1,0}:\sH \to \sH^{1,0}$ is a complex linear isomorphism of complex Hilbert spaces, as claimed. Because
\[
  \pi^{0,1}u = \overline{\pi^{1,0}u}, \quad\text{for all } u \in \sH,
\]
it follows that $\pi^{0,1}:\sH \to \sH^{0,1}$ is a complex antilinear isomorphism of complex Hilbert spaces, as claimed.

We observe next that the spectral projections $\pi^{0,1}$ and $\pi^{0,1}$ may be expressed in terms of $J$ via the familiar formulae \eqref{eq:sH_pi_10_and_pi_01_idempotents}, namely
\begin{subequations}
  \label{eq:sH_pi_10_and_pi_01_orthogonal_projections}
  \begin{align}
    \label{eq:sH_pi_10_orthogonal_projection}
    \pi^{1,0} &= \frac{1}{2}(\id - iJ):\sH^\CC \to \sH^{1,0},
    \\
    \label{eq:sH_pi_01_orthogonal_projection}
    \pi^{0,1} &= \frac{1}{2}(\id + iJ):\sH^\CC \to \sH^{0,1}.
  \end{align}                
\end{subequations}
These two identities follow immediately from the fact \eqref{eq:Eigen-decomposition_J_on_sHCC} that $J=i$ on $\sH^{1,0}$ and $J=-i$ on $\sH^{0,1}$. We denote
\begin{equation}
  \label{eq:pi_10v_prime_and_pi_01v_primeprime}
  v' = \frac{1}{2}(v - iJv)
  \quad\text{and}\quad
  v'' = \frac{1}{2}(v + iJv), \quad\text{for all } v \in \sH^\CC,
\end{equation}
as customary.

\subsection{Symplectic forms and almost complex structures on Hilbert dual spaces}
\label{subsec:Symplectic_forms_almost_complex_structures_Hilbert_dual_spaces}
Continue the notation of Section \ref{subsec:Symplectic_forms_almost_complex_structures_Hilbert_spaces}. We obtain an induced almost complex structure $J$ on $\Hom_\RR^c(\sH,\RR)$ and hence an almost complex structure $J$ on  
\[
  \Hom_\RR^c(\sH,\CC) = \Hom_\RR^c(\sH,\RR)\otimes_\RR\CC
\]
by setting
\begin{equation}
  \label{eq:Definition_J_on_Hom_RR_sH_to_CC}
  J\sL := \sL\circ J, \quad\text{for all } \sL \in \Hom_\RR^c(\sH,\RR),
\end{equation}
since $J^2\sL = \sL\circ J^2 = \sL\circ(-\id) = -\sL$, for all $\sL \in \Hom_\RR^c(\sH,\RR)$. Because $J$ is $g$-orthogonal on $\sH$, it follows from the isometric isomorphism $\Hom_\RR^c(\sH,\RR) \cong \sH$ that $J$ is $g$-orthogonal on $\Hom_\RR^c(\sH,\RR)$ and thus Hermitian-orthogonal on $\Hom_\RR^c(\sH,\CC)$. By replacing the role of $\sH^\CC$ in \eqref{eq:Eigen-decomposition_J_on_sHCC} by $\Hom_\RR^c(\sH,\CC)$, we obtain a Hermitian-orthogonal splitting 
\begin{equation}
  \label{eq:Eigenspace-decomposition_J_on_Hom_RR_sH_to_CC}
  \Hom_\RR^c(\sH,\CC) = \Hom_\RR^c(\sH,\CC)^{1,0} \oplus \Hom_\RR^c(\sH,\CC)^{0,1}
\end{equation}
such that
\begin{equation}
  \label{eq:Eigenvalue-decomposition_J_on_Hom_RR_sH_to_CC}
  J=i \quad\text{on } \Hom_\RR^c(\sH,\CC)^{1,0}
  \quad\text{and}\quad
  J=-i \quad\text{on } \Hom_\RR^c(\sH,\CC)^{0,1}.
\end{equation}
By analogy with \eqref{eq:pi_10v_prime_and_pi_01v_primeprime}, we write
\begin{equation}
  \label{eq:pi_10sL_prime_and_pi_01sL_primeprime}
  \sL' = \frac{1}{2}(\sL - iJ\sL)
  \quad\text{and}\quad
  \sL'' = \frac{1}{2}(\sL + iJ\sL), \quad\text{for all } \sL \in \Hom_\RR^c(\sH,\CC).
\end{equation}
Because $J\sL' = i\sL'$ and $J\sL'' = -i\sL'$ by \eqref{eq:Eigenvalue-decomposition_J_on_Hom_RR_sH_to_CC}, it follows from the definition \eqref{eq:Definition_J_on_Hom_RR_sH_to_CC} of $J$ on $\Hom_\RR^c(\sH,\CC)$ that
\begin{align*}
  \sL'Jv &= (\sL'J)v = (J\sL')v = (i\sL')v = i\sL'v,
    \\
  \sL''Jv &= (\sL''J)v = (J\sL'')v = (-i\sL'')v = i\sL''v,
  \quad\text{for all }  \sL \in \Hom_\RR^c(\sH,\CC) \text{ and } v \in \sH,
\end{align*}
so that $\Hom_\RR^c(\sH,\CC)^{1,0}$ and $\Hom_\RR^c(\sH,\CC)^{0,1}$ are the subspaces of $J$-linear and $J$-antilinear $\RR$-valued functions on $\sH = (\sH^\RR,J)$, respectively. Thus,
\begin{equation}
  \label{eq:Hom_RR_sH_to_CC_10_is_Hom_CC_sH_to_CC_and_Hom_RR_sH_to_CC_01_is_Hom_barCC_sH_to_CC}
  \Hom_\RR^c(\sH,\CC)^{1,0} = \Hom_\CC^c(\sH,\CC)
  \quad\text{and}\quad
  \Hom_\RR^c(\sH,\CC)^{0,1} = \Hom_{\bar\CC}^c(\sH,\CC),
\end{equation}
where $(a+ib)v := (a+bJ)v$ for all $v \in \sH$ and $a+ib \in \CC$. In other words, $\Hom_\RR^c(\sH,\CC)^{1,0}$ and $\Hom_\RR^c(\sH,\CC)^{0,1}$ are the subspaces of complex linear and antilinear $\RR$-valued functions on $\sH$, respectively.

We now describe an alternative approach to a construction of the Hermitian orthogonal splitting \eqref{eq:Eigenspace-decomposition_J_on_Hom_RR_sH_to_CC}. For the spectral projections $\pi^{1,0}$ and $\pi^{0,1}$ determined by the skew-adjoint operator $A \in \End_\CC(\sH^\CC)$ as in \eqref{eq:pi10_and_pi01_spectral_projections} and $\sL \in \Hom_\RR^c(\sH,\CC)$, we extend $\sL$ by $\CC$-linearity to $\sL \in \Hom_\CC^c(\sH^\CC,\CC)$ and define
\begin{subequations}
  \label{eq:pi10_sL_and_pi01_sL}
  \begin{align}
    \label{eq:pi10_sL}
    \pi^{1,0}\sL &:= \sL\circ\pi^{1,0} \in \Hom_\CC^c(\sH^{1,0},\CC),
    \\
    \label{eq:pi01_sL}
    \pi^{0,1}\sL &:= \sL\circ\pi^{0,1} \in \Hom_\CC^c(\sH^{0,1},\CC).
  \end{align}
\end{subequations}
The Hermitian orthogonal splitting $\sH^\CC = \sH^{1,0}\oplus\sH^{0,1}$ induces a Hermitian orthogonal splitting
\begin{gather}
  \label{eq:Induced_spectral_decomposition_Hom_CC_sHCC_to_CC}
  % \Hom_\RR^c(\sH^\CC,\CC) = \Hom_\RR^c(\sH^{1,0},\CC) \oplus \Hom_\RR^c(\sH^{0,1},\CC),
  % \\
  \Hom_\CC^c(\sH^\CC,\CC) = \Hom_\CC^c(\sH^{1,0},\CC) \oplus \Hom_\CC^c(\sH^{0,1},\CC).
\end{gather}
% where
% \begin{align*}
%   \Hom_\CC^c(\sH^\CC,\CC)^{1,0}
%   &:=
%   \pi^{1,0}\left(\Hom_\CC^c(\sH^\CC,\CC)\right)
%   =
%     \left\{\sL \in \Hom_\CC^c(\sH^\CC,\CC): J\sL = i\sL\right\},
%   \\
%   \Hom_\CC^c(\sH^\CC,\CC)^{0,1}
%   &:=
%   \pi^{0,1}\Hom_\CC^c(\sH^\CC,\CC)\left(\Hom_\CC^c(\sH^\CC,\CC)\right)
%   =
%     \left\{\sL \in \Hom_\CC^c(\sH^\CC,\CC): J\sL = -i\sL\right\}.
% \end{align*}
The definition \eqref{eq:Definition_J_on_Hom_RR_sH_to_CC} of the induced almost complex structure $J$ on $\Hom_\RR^c(\sH,\CC)$ yields
\[
  J\pi^{1,0}\sL = (\pi^{1,0}\sL)\circ J = \sL\circ J\circ\pi^{1,0} = \sL\circ i\pi^{1,0} = i\sL\circ \pi^{1,0}
  = i\pi^{1,0}\sL,
\]
noting that we extend $\sL \in \Hom_\RR^c(\sH,\CC)$ by $\CC$-linearity to $\sL \in \Hom_\CC^c(\sH^\CC,\CC)$, while
\[
  J\pi^{0,1}\sL = (\pi^{0,1}\sL)\circ J = \sL\circ J\circ\pi^{0,1} = \sL\circ (-i)\pi^{0,1} = -i\sL\circ \pi^{0,1}
  = -i\pi^{0,1}\sL.
\]
Therefore, a comparison with \eqref{eq:Eigenspace-decomposition_J_on_Hom_RR_sH_to_CC} and \eqref{eq:Eigenvalue-decomposition_J_on_Hom_RR_sH_to_CC} yields the identities
\begin{equation}
  \label{eq:Hom_RR_sH_to_CC_10_is_Hom_CC_sH10_to_CC_and_Hom_RR_sH_to_CC_01_is_Hom_CC_sH01_to_CC}
  \Hom_\RR^c(\sH,\CC)^{1,0} = \Hom_\CC^c(\sH^{1,0},\CC)
  \quad\text{and}\quad
  \Hom_\RR^c(\sH,\CC)^{0,1} = \Hom_\CC^c(\sH^{0,1},\CC),
\end{equation}
in agreement with Huybrechts \cite[Lemma 1.2.6, p. 26]{Huybrechts_2005}. Furthermore, by \eqref{eq:Hom_RR_sH_to_CC_10_is_Hom_CC_sH_to_CC_and_Hom_RR_sH_to_CC_01_is_Hom_barCC_sH_to_CC} we have
\begin{equation}
  \label{eq:Hom_CC_sH10_to_CC_is_Hom_CC_sH_to_CC_and_Hom_CC_sH01_to_CC_is_Hom_barCC_sH_to_CC}
  \Hom_\CC^c(\sH^{1,0},\CC) = \Hom_\CC^c(\sH,\CC)
  \quad\text{and}\quad
  \Hom_\CC^c(\sH^{0,1},\CC) = \Hom_{\bar\CC}^c(\sH,\CC),
\end{equation}
noting that $\pi^{0,1}:\sH\to\sH^{0,1}$ is a complex antilinear isomorphism. A comparison with \eqref{eq:pi_10sL_prime_and_pi_01sL_primeprime} yields the identities
\begin{equation}
  \label{eq:sL'_is_pi10_sL_and_sL''_is_pi01_sL}
  \sL' = \pi^{1,0}\sL
  \quad\text{and}\quad
  \sL'' = \pi^{0,1}\sL, \quad\text{for all } \sL \in \Hom_\RR^c(\sH,\CC),
\end{equation}
and hence the familiar formulae for the Hermitian orthogonal projections,
\begin{subequations}
  \label{eq:Hom_RR_sH_to_CC_pi_10_and_pi_01_orthogonal_projections}
  \begin{align}
    \pi^{1,0} &= \frac{1}{2}(\id - iJ): \Hom_\RR^c(\sH,\CC) \to \Hom_\RR^c(\sH,\CC)^{1,0},
    \\
    \pi^{0,1} &= \frac{1}{2}(\id + iJ): \Hom_\RR^c(\sH,\CC) \to \Hom_\RR^c(\sH,\CC)^{0,1},
  \end{align}
\end{subequations}  
just as in \eqref{eq:sH_pi_10_and_pi_01_orthogonal_projections}

\subsection{Symplectic subspaces of Hilbert spaces}
\label{subsec:Symplectic_subspaces_Hilbert_spaces}
Continue the notation of Section \ref{subsec:Symplectic_forms_almost_complex_structures_Hilbert_dual_spaces}. Suppose that $V \subset \sH^\RR$ is a closed subspace. Recall that $V$ is \emph{symplectic} if the symplectic form $\omega:\sH^\RR\times\sH^\RR$ restricts to a symplectic form $\omega:V\times V \to \RR$, that is, $\omega$ is a (strongly) non-degenerate two-form on $V$, whereas $V$ is \emph{isotropic} if $\omega \equiv 0$ on $V\times V$ (see, for example, Cannas da Silva \cite[Chapter 1, Homework 1, p. 8]{Cannas_da_Silva_lectures_on_symplectic_geometry}).

Suppose that $K \subset \sH^\RR$ is a symplectic subspace, so $\omega:K\times K\to\RR$ is a symplectic form and, as usual, the inner product $g:\sH^\RR\times \sH^\RR\to\RR$ restricts to an inner product $g:K\times K\to\RR$. In this setting, the real linear isomorphisms \eqref{eq:tilde_g_and_tilde_omega} of real Hilbert spaces restrict to real linear isomorphisms,
\begin{subequations}
  \label{eq:tilde_g_K_and_tilde_omega_K}
  \begin{align}
    \label{eq:tilde_g_K}
    g^\flat\restriction_K :K \ni v &\mapsto g(v,\cdot) \in \Hom_\RR(K,\RR),
    \\
    \label{eq:tilde_omega}
    \omega^\flat\restriction_K:K \ni v &\mapsto \omega(v,\cdot) \in \Hom_\RR(K,\RR),
\end{align}
\end{subequations}
and thus $A \in \End_\RR(\sH^\RR)$ in \eqref{eq:A_is_tilde_g-inverse_tilde_omega} restricts to a skew-adjoint real linear isomorphism of real Hilbert subspaces with respect to the inner product $g:K\times K\to\RR$,
\begin{equation}
  \label{eq:A_K_is_tilde_g_K-inverse_tilde_omega_K}
  A_K := (g^\flat\restriction_K)^{-1}\circ(\omega^\flat\restriction_K) \in \End_\RR(K),
\end{equation}
or, equivalently (see Cannas \cite[Section 12.2, Proof of Proposition 12.3, p. 84]{Cannas_da_Silva_lectures_on_symplectic_geometry}),
\begin{equation}
  \label{eq:omega_is_g(A_K,.)}
  \omega(u,v)  =: g(A_Ku,v), \quad\text{for all } u,v \in K.
\end{equation} 
Before proceeding further, we note the important

\begin{rmk}[Distinction between $A_K$ and $A\restriction K$]
\label{rmk:A_K_neq_A_restriction_K}  
In the definitions \eqref{eq:A_K_is_tilde_g_K-inverse_tilde_omega_K} or \eqref{eq:omega_is_g(A_K,.)} of $A_K \in \End_\RR(K)$, it is not true in general that $A_K = A\restriction K$, where $A \in \End_\RR(\sH^\RR)$ is defined by equations \eqref{eq:A_is_tilde_g-inverse_tilde_omega} or \eqref{eq:omega_is_g(A,.)}. Instead, we have that
\[
  A_K = \pi_KA \restriction_K,
\]
where $\pi_K:\sH^\RR \to K$ is $g$-orthogonal projection. Thus, $K$ is an \emph{invariant subspace} of the operator $\pi_KA \in \End_\RR(\sH^\RR)$, but \emph{not} of the operator $A$. \qed
\end{rmk}

As in Section \ref{subsec:Symplectic_forms_almost_complex_structures_Hilbert_spaces}, the skew-adjoint bounded real linear operator $A_K \in \End_\RR(K)$ extends by complex linearity to a bounded complex linear operator $A_K \in \End_\CC(K^\CC)$ that is skew-adjoint with respect to the Hermitian inner product on the subspace $K^\CC = K\otimes_\RR\CC$ induced from $\sH^\CC$. By applying  \eqref{eq:Spectrum_A_purely_imaginary} to $A_K$, we have
\begin{equation}
  \label{eq:Spectrum_A_K_purely_imaginary}
  \sigma(A_K) \subset i\RR,
\end{equation}
where $0 \notin A_K$ since $A$ is invertible on $K$. We have the decomposition of the spectrum,
\begin{equation}
  \label{eq:Spectrum_A_K_disjoint_union_pos_neg_spectra_A_K}
  \sigma(A_K) = \sigma^+(A_K) \sqcup \sigma^+(A_K),
\end{equation}
by analogy with \eqref{eq:Spectrum_A_disjoint_union_pos_neg_spectra_A}.
% , and the inclusion of spectra,
% \[
%   \sigma^\pm(A_K) \subset \sigma^\pm(A).
% \]
If $\Omega \subset \sigma(A_K)$ is a Borel subset, then $E(\Omega)K^\CC$ is an invariant subspace for $A_K$ on $K^\CC$ since the operators $A_K$ and $E(\Omega)$ commute by Rudin \cite[Theorem 12.23, p. 324]{Rudin}, where $E$ is as in \eqref{eq:Spectral_domposition_A} with $A$ replaced by $A_K$.
% , and if $\Omega \subset \sigma(A_K)$ is a Borel subset, then $E(\Omega)K^\CC$ is an invariant subspace for $A$ on $K^\CC$ and
% \[
%   E(\Omega)K^\CC = K^\CC \cap E(\Omega)\sH^\CC.
% \]
% In particular, we have the identities
% \[
%   E\left(\sigma^\pm(A_K)\right)K^\CC = K^\CC\cap E\left(\sigma^\pm(A)\right)\sH^\CC. 
% \]
Therefore, denoting
\[
  K^{1,0} := E\left(\sigma^+(A_K)\right)K^\CC
  \quad\text{and}\quad
  K^{0,1} := E\left(\sigma^+(A_K)\right)K^\CC,
\]
by analogy with \eqref{eq:pi10_and_pi01_spectral_projections} and \eqref{eq:sH10_and_sH01_invariant_subspaces_A}, we have 
\[
  K^\CC = K^{1,0} \oplus K^{1,0}.
\]
a Hermitian orthogonal splitting.

\section[Transverse intersections of symplectic subspaces]{Transverse intersections of symplectic subspaces and almost symplectic manifolds}
\label{sec:Transverse_intersections_symplectic_subspaces_almost_symplectic_manifolds}
The fact that the transverse intersection of two almost symplectic submanifolds of an almost symplectic manifold is (almost) symplectic often appears to be implicitly assumed --- see, for example, Paoletti \cite{Paoletti_2001} or Zhang \cite{Zhang_2018} --- but we are unable to find a reference to an explicit statement or proof in the literature. In this section, we prove this result for the transverse intersection of two or more symplectic subspaces of a finite-dimensional symplectic manifold in Section \ref{subsec:Transverse_intersections_symplectic_subspaces_finite-dimensional_vector_spaces}. We consider the transverse intersection of two or more symplectic subspaces of an infinite-dimensional symplectic vector space in Section \ref{subsec:Transverse_intersections_symplectic_subspaces_symplectic_Banach_spaces}. Finally, we consider the transverse intersection of two or more almost symplectic submanifolds of an infinite-dimensional almost symplectic manifold in Section \ref{subsec:Transverse_intersections_symplectic_submanifolds_symplectic_Banach_manifold}.

\subsection{Transverse intersections of symplectic subspaces of finite-dimensional symplectic vector spaces}
\label{subsec:Transverse_intersections_symplectic_subspaces_finite-dimensional_vector_spaces}
Recall that the intersection of two real linear subspaces $W_1,W_2$ of a real vector space $V$ is \emph{transverse}, denoted by $W_1\transv W_2$, if $W_1 + W_2 = V$ or, equivalently, $\dim W_1 + \dim W_2 = \dim V + \dim(W_1 \cap W_2)$ or $\codim W_1 + \codim W_2 = \codim(W_1 \cap W_2)$. 

\begin{lem}[Transverse intersection of two symplectic subspaces of a symplectic vector space is symplectic]
\label{lem:Transverse_intersection_symplectic_subspaces_is_symplectic} 
Let $(V,\omega)$ be a symplectic vector space. If $W_1,W_2 \subset V$ are symplectic subspaces that intersect transversely in $V$, then $W_1 \cap W_2$ is a symplectic subspace of $V$.
\end{lem}

\begin{proof}
If $W \subset V$ is a linear subspace, then its \emph{symplectic complement} in $V$ is defined by
\begin{equation}
  \label{eq:W_oplus_W_omega_is_V}
  W^\omega := \{ v \in V \mid \omega(v,w)=0, \text{ for all } w\in W \}.
\end{equation}
See McDuff and Salamon \cite[Equation (2.1.1), p. 38]{McDuffSalamonSympTop3}. Moreover, for any linear subspace $W \subset V$,
\begin{equation}
  \label{eq:Dimension_W+dimension_W_omega_is_dimension_V}
  \dim W + \dim W^\omega = \dim V.
\end{equation}
See Cannas da Silva \cite[Chapter 1, Homework 1, p. 8]{Cannas_da_Silva_lectures_on_symplectic_geometry}, McDuff and Salamon \cite[Lemma 2.1.1, p. 39]{McDuffSalamonSympTop3}, or Meinrenken \cite[Section 2.2, Proposition 2.5, p. 7]{Meinrenken_symplectic_geometry_lecture_notes_2024}. By a standard lemma in symplectic linear algebra (see Cannas da Silva \cite[Chapter 1, Homework 1, p. 8]{Cannas_da_Silva_lectures_on_symplectic_geometry}, McDuff and Salamon \cite[Lemma 2.1.1, p. 39]{McDuffSalamonSympTop3}, or Meinrenken \cite[Section 2.2, Proposition 2.5 and Exercise 2.6, p. 7]{Meinrenken_symplectic_geometry_lecture_notes_2024}) a subspace $W \subset V$ is \emph{symplectic} if and only if
\begin{equation}
  \label{eq:W_oplus_W_omega_is_V}
    V = W \oplus W^\omega,
\end{equation}
as a direct sum of real linear subspaces of $V$. In particular, by \eqref{eq:Dimension_W+dimension_W_omega_is_dimension_V} and \eqref{eq:W_oplus_W_omega_is_V}, a subspace $W \subset V$ is symplectic if and only if
\begin{equation}
  \label{eq:W_cap_W_omega_is_zero}
  W \cap W^\omega = (0).
\end{equation}
By hypothesis that $W_1\transv W_2$, we have $W_1 + W_2 = V$  and $W_1^\omega \cap W_2^\omega = (W_1 + W_2)^\omega$ by Meinrenken \cite[Section 2.2, Proposition 2.5, p. 7]{Meinrenken_symplectic_geometry_lecture_notes_2024} and the fact that $\omega$ is nondegenerate on $V$ is equivalent to $V^\omega = (0)$, so
\begin{equation}
  \label{eq:W1_omega_cap_W2_omega_is_zero}
    W_1^\omega \cap W_2^\omega = (0).
\end{equation}
To show that $W_1 \cap W_2$ is symplectic, it suffices by \eqref{eq:W_cap_W_omega_is_zero} to prove that
\begin{equation}
  \label{eq:W1_cap_W2_cap_(W1_cap_W2)_omega_is_zero}
  W_1 \cap W_2 \cap (W_1 \cap W_2)^\omega = (0).
\end{equation}
To this end, let $x \in W_1 \cap W_2 \cap (W_1 \cap W_2)^\omega$. By definition \eqref{eq:W_oplus_W_omega_is_V} of $(W_1 \cap W_2)^\omega$, we have 
\begin{equation}
  \label{eq:omega(x,y)_is_zero_y_in_W1_cap_W2}
    \omega(x,y) = 0, \quad \text{for all } y \in W_1 \cap W_2.
\end{equation}
Let $\omega_1 := \omega|_{W_1}$ denote the restriction of $\omega$ to $W_1$ and observe that, writing $(W_1 \cap W_2)^{\omega_1}$ for the symplectic complement of $W_1 \cap W_2$ in $(W_1,\omega_1)$, 
\begin{equation}
  \label{eq:omega(x,y)_is_zero_y_in_W1_cap_W2_omega_1}
    \omega(x,y) = 0, \quad \text{for all } y \in (W_1 \cap W_2)^{\omega_1}.
\end{equation}
The vector space $W_1$ is a symplectic subspace of $(V,\omega)$ by hypothesis, so $\omega_1$ is nondegenerate. Hence, the equality \eqref{eq:W_oplus_W_omega_is_V} yields (after replacing $(V,\omega)$ by $(W_1,\omega_1)$)
\[
  W_1 = W_1 \cap W_2 \oplus (W_1 \cap W_2)^{\omega_1}.
\]  
The preceding direct sum decomposition of $W_1$ and the identities \eqref{eq:omega(x,y)_is_zero_y_in_W1_cap_W2} and \eqref{eq:omega(x,y)_is_zero_y_in_W1_cap_W2_omega_1} therefore imply
\[
  \omega(x,y) = 0, \quad \text{for all } y \in W_1.
\]
Consequently, $x \in W_1^\omega$ and, by a symmetrical argument, $x \in W_2^\omega$. Hence,
\[
    x \in W_1^\omega \cap W_2^\omega.
\]
Therefore, $x = 0$ by \eqref{eq:W1_omega_cap_W2_omega_is_zero} and this proves the claim \eqref{eq:W1_cap_W2_cap_(W1_cap_W2)_omega_is_zero}. Thus, $W_1 \cap W_2$ is a symplectic subspace.
\end{proof}

The following example shows that the conclusion of Lemma \ref{lem:Transverse_intersection_symplectic_subspaces_is_symplectic} can be false when the hypothesis of transversality is omitted.

\begin{exmp}[Non-transverse intersections of two symplectic subspaces of a symplectic vector space need not be symplectic]
\label{exmp:Intersection_symplectic_subspaces_not_always_symplectic}  
Let $V=\RR^4$ with coordinates $(x_1,y_1,x_2,y_2)$ and equipped with the standard symplectic form
\[
\omega = dx_1 \wedge dy_1 + dx_2 \wedge dy_2.
\]
Define two subspaces
\[
  W_1 := \operatorname{span}\{e_1,e_2\}
  \quad\text{and}\quad
  W_2 := \operatorname{span}\{e_1,\, e_2 + e_3\},
\]
where
\[
  e_1 := (1,0,0,0), \quad e_2 := (0,1,0,0), \quad e_3 := (0,0,1,0).
\]
We claim that both $W_1$ and $W_2$ are symplectic subspaces, but $W_1\cap W_2$ is not a symplectic subspace.

First observe that the restrictions $\omega|_{W_1}$ and $\omega|_{W_2}$ are nondegenerate. Indeed, $\omega(e_1,e_2)=1$, so $\{e_1,e_2\}$ is a symplectic basis for $W_1$, and $\omega(e_1,e_2+e_3) = \omega(e_1,e_2)+\omega(e_1,e_3) = \omega(e_1,e_2)=1$, so $\{e_1,e_2+e_3\}$ is a symplectic basis for $W_2$.

We now consider the intersection $W_1 \cap W_2$. Observe that $v \in W_1 \cap W_2$ if and only if it is of the form
\[
v = \alpha e_1 = (\alpha,0,0,0),
\]
for $\alpha \in \RR$. Thus,
\[
W_1 \cap W_2 = \RR e_1,
\]
a one--dimensional subspace of $V = \RR^4$. But every one-dimensional subspace is isotropic for a skew-symmetric form,
and in particular
\[
\omega(e_1, e_1) = 0.
\]
Hence, $\omega|_{W_1\cap W_2}$ is degenerate. Therefore, $W_1\cap W_2$ is \emph{not} a symplectic subspace, even though $W_1$ and $W_2$ are symplectic subspaces. Note that $W_1$ and $W_2$ do not intersect transversely in $V$ since $\codim(W_1\cap W_2) = 3 < 2 + 2 = \codim W_1 + \codim W_2$.
\qed
\end{exmp}

We next consider the case of intersections of three or more symplectic subspaces of a finite dimensional symplectic vector space. For this purpose, we need an appropriate definition of transverse intersection. See Rudin \cite[Theorem 1.41, p. 31]{Rudin} for the fact that the quotient of a Banach space and a closed linear subspace is a Banach space.
% PF12-8-2025 Transversality of multiple intersections: for the correct definition, see
% COMMENT https://math.stackexchange.com/questions/386600/transverse-intersection-of-multiple-submanifolds

\begin{defn}[Transverse intersection of closed linear subspaces of a real Banach space]
\label{defn:Transverse_intersection_finitely_many_subspaces_vector_space}  
Let $V$ be a real Banach space. If $W_1,\ldots, W_m \subset V$ are closed linear subspaces for $m\geq 2$, then the intersection $W := W_1\cap\cdots\cap W_n$ is \emph{transverse}, denoted $W_1 \transv\cdots\transv W_m$, if the map
\begin{equation}
  \label{eq:Transverse_intersection_finitely_many_subspaces_vector_space}
  V/W \ni [v]_W \mapsto \left([v]_{W_1},\ldots,[v]_{W_m}\right) \in V/W_1 \oplus \cdots \oplus V/W_m 
\end{equation}
is an isomorphism of Banach spaces from $V/W$ onto the internal direct sum of the closed linear subspaces $V/W_n$ or, equivalently, if the map
\begin{equation}
  \label{eq:Transverse_intersection_finitely_many_subspaces_vector_space_epimorphism}
  V \ni v \mapsto \left([v]_{W_1},\ldots,[v]_{W_m}\right) \in V/W_1 \oplus \cdots \oplus V/W_m 
\end{equation}
is an epimorphism of Banach spaces from $V$ onto the internal direct sum of the closed linear subspaces $V/W_n$ of $V/W$.
\end{defn}

The kernel of the epimorphism \eqref{eq:Transverse_intersection_finitely_many_subspaces_vector_space_epimorphism} is equal to
\[
  \{v\in V \mid v \in W_n \text{ for } n=1,\ldots,m\} = W_1\cap\cdots\cap W_n = W,
\]
and so the conditions \eqref{eq:Transverse_intersection_finitely_many_subspaces_vector_space} and \eqref{eq:Transverse_intersection_finitely_many_subspaces_vector_space_epimorphism} are indeed equivalent.

\begin{rmk}[Transverse intersections and non-degenerate symmetric or skew-symmetric bilinear forms]
\label{rmk:Transverse_intersections_non-degenerate_symmetric_skew-symmetric_bilinear_forms}
If $\beta:V\times V \to \RR$ is a symmetric or skew-symmetric real bilinear form on a finite-dimensional real vector space $V$ and $W \subset V$ is a linear subspace, then $\beta$ is \emph{non-degenerate} if and only if
\begin{equation}
  \label{eq:V_isomorphic_W_direct_sum_W_omega}
  V = W\oplus W^\beta,
\end{equation}
as an internal direct sum of real linear subspaces of $V$, where the $\beta$-\emph{complementary subspace} is
\begin{equation}
  \label{eq:Beta_orthogonal_complement_W}
  W^\beta := \{v\in V: \beta(v,w) = 0 \text{ for all } w \in W\}.
\end{equation}
See, for example, Fisher \cite[Corollary 2.3, p. 2]{Fisher_linear_algebra_non-degenerate_bilinear_forms} or, when $\beta$ is skew-symmetric, Cannas da Silva \cite[Chapter 1, Homework 1, p. 8]{Cannas_da_Silva_lectures_on_symplectic_geometry}, McDuff and Salamon \cite[Lemma 2.1.1, p. 39]{McDuffSalamonSympTop3}, or Meinrenken \cite[Section 2.2, Proposition 2.5 and Exercise 2.6, p. 7]{Meinrenken_symplectic_geometry_lecture_notes_2024}.

Given such a non-degenerate $\beta$, the transversality condition \eqref{eq:Transverse_intersection_finitely_many_subspaces_vector_space} is equivalent to the assertion that
\begin{equation}
  \label{eq:Transverse_intersection_finitely_many_subspaces_vector_space_non-degenerate_2-form}
  W^\beta \ni v \mapsto \left((1-\pi_1)v,\ldots,(1-\pi_n)v\right) \in W_1^\beta \oplus \cdots \oplus W_m^\beta 
\end{equation}
is an isomorphism of real vector spaces, noting that $W \subset W_n \implies W_n^\beta \subset W^\beta$ and, writing $V = W_n\oplus W_n^\beta$, we let $\pi_n:V \to W_n$ and $1-\pi_n:V \to W_n^\beta$ denote the linear idempotent maps corresponding to the direct sum decomposition \eqref{eq:V_isomorphic_W_direct_sum_W_omega}, for $n=1,\ldots,m$. In particular,
\begin{equation}
  \label{eq:Transverse_intersection_finitely_many_subspaces_internal_direct_sum_normal_spaces}
  W^\beta = W_1^\beta \oplus \cdots \oplus W_m^\beta, 
\end{equation}
as an internal direct sum of real linear subspaces of $W^\beta$.

Equation \eqref{eq:Transverse_intersection_finitely_many_subspaces_vector_space} yields the formula for the expected dimension of $W_1\cap\cdots\cap W_m$,
\begin{equation}
  \label{eq:Codimension_transverse_intersection_finitely_many_subspaces}
  \codim(W_1\cap\cdots\cap W_m) = \sum_{n=1}^m \codim W_n.
\end{equation}
When $m=2$, then \eqref{eq:Transverse_intersection_finitely_many_subspaces_vector_space} is equivalent to the condition $V=W_1+W_2$. 
\qed
\end{rmk}  

\begin{lem}[Transverse intersection of symplectic subspaces of a symplectic vector space is symplectic]
\label{lem:Transverse_intersection_countably_many_symplectic_subspaces_is_symplectic} 
Let $(V,\omega)$ be a finite-dimensional symplectic vector space. If $W_1,\ldots, W_m$ are symplectic subspaces that intersect transversely in $V$ in the sense of Definition \ref{defn:Transverse_intersection_finitely_many_subspaces_vector_space} for $m\geq 2$, then $W_1 \cap\cdots\cap W_m$ is a symplectic subspace of $V$.
\end{lem}

\begin{proof}
We modify the proof of Lemma \ref{lem:Transverse_intersection_symplectic_subspaces_is_symplectic}. By \eqref{eq:W_oplus_W_omega_is_V}, we have for $n=1,\ldots,m$ that
\[
  V = W_n \oplus W_n^\omega,
\]
and we thus obtain
\begin{equation}
  \label{eq:V_mod_Wn_isomorphic_Wn_omega}
  V/W_n \cong W_n^\omega, \quad\text{for } n=1,\ldots,m.
\end{equation}
Writing $W := W_1 \cap\cdots\cap W_m$, we similarly obtain
\begin{equation}
  \label{eq:V_mod_W_isomorphic_W_omega}
  V/W \cong W^\omega.
\end{equation}
Because the intersection $W_1 \cap\cdots\cap W_m$ is transverse by hypothesis, the isomorphisms \eqref{eq:Transverse_intersection_finitely_many_subspaces_vector_space}, \eqref{eq:V_mod_Wn_isomorphic_Wn_omega}, and \eqref{eq:V_mod_W_isomorphic_W_omega} yield the isomorphism
\[
  W^\omega \cong W_1^\omega \oplus \cdots \oplus W_m^\omega,
\]
and hence, noting that $W\subset W_n \implies W_n^\omega \subset W^\omega$ for $n=1,\ldots,m$ and thus $W_1^\omega + \cdots + W_m^\omega \subset W$, we obtain the decomposition of $W^\omega$ as an internal direct sum,
\[
  W^\omega = W_1^\omega \oplus \cdots \oplus W_m^\omega.
\]
In particular, we obtain
\begin{equation}
  \label{eq:W1_omega_cap_cdots_cap_Wm_omega_is_zero}
    W_1^\omega \cap\cdots\cap W_m^\omega = (0).
\end{equation}
To show that $W = W_1 \cap\cdots\cap W_m$ is symplectic, it suffices by \eqref{eq:W_cap_W_omega_is_zero} to prove that
\begin{equation}
  \label{eq:W1_cap_cdots_cap_Wm_cap_(W1_cap_cdots_cap_Wm)_omega_is_zero}
  W \cap W^\omega = (0).
\end{equation}
To this end, let
\[
  x \in W \cap W^\omega,
\]
so that by definition \eqref{eq:W_oplus_W_omega_is_V} of $W^\omega$, 
\begin{equation}
  \label{eq:omega(x,y)_is_zero_y_in_W}
    \omega(x,y) = 0, \quad \text{for all } y \in W.
\end{equation}
Fix $n$ and let $\omega_n := \omega|_{W_n}$ denote the restriction of $\omega$ to $W_n$ and observe that, writing $W^{\omega_n}$ for the symplectic complement of $W$ in $(W_n,\omega_n)$, 
\begin{equation}
  \label{eq:omega(x,y)_is_zero_y_in_W_omega_n}
    \omega(x,y) = 0, \quad \text{for all } y \in W^{\omega_n}.
\end{equation}
The vector space $W_n$ is a symplectic subspace of $(V,\omega)$ by hypothesis, so $\omega_n$ is nondegenerate. Hence, the equality \eqref{eq:W_oplus_W_omega_is_V} yields (after replacing $(V,\omega)$ by $(W_n,\omega_n)$)
\[
  W_n = W \oplus W^{\omega_n}.
\]  
The preceding direct sum decomposition of $W_n$ and the identities \eqref{eq:omega(x,y)_is_zero_y_in_W} and \eqref{eq:omega(x,y)_is_zero_y_in_W_omega_n} therefore imply
\[
  \omega(x,y) = 0, \quad \text{for all } y \in W_n.
\]
Consequently, $x \in W_n^\omega$ and, since $n$ was arbitrary,
\[
    x \in W_1^\omega \cap\cdots\cap W_m^\omega.
\]
Therefore, $x = 0$ by \eqref{eq:W1_omega_cap_cdots_cap_Wm_omega_is_zero} and this proves the claim \eqref{eq:W1_cap_cdots_cap_Wm_cap_(W1_cap_cdots_cap_Wm)_omega_is_zero}. We conclude that $W = W_1 \cap\cdots\cap W_m$ is a symplectic subspace of $(V,\omega)$.
\end{proof}

\subsection{Transverse intersections of symplectic subspaces of symplectic Banach  spaces}
\label{subsec:Transverse_intersections_symplectic_subspaces_symplectic_Banach_spaces}
We next consider the case of intersections of symplectic subspaces of a symplectic Banach space. The extension of the results of Section \ref{subsec:Transverse_intersections_symplectic_subspaces_finite-dimensional_vector_spaces} to the setting of Banach spaces (or even Hilbert spaces) involves some more technicalities, for which we partly rely on Kobayashi \cite[Section 7.5]{Kobayashi_differential_geometry_complex_vector_bundles}, Marsden and Ratiu \cite[Section 2.2]{Marsden_Ratiu_introduction_mechanics_symmetry}, Marsden and Weinstein \cite{Marsden_Weinstein_1974}, and Weinstein \cite{Weinstein_1969}.

Let $V$ be a real Banach space with continuous dual space $V^*$. If $W \subset V$ is a linear subspace, not necessarily closed, then the \emph{annihilator of $W$ in $V^*$} (see, for example, Rudin \cite[Section 4.6, p. 95]{Rudin}) is 
\begin{equation}
  \label{eq:Annihilator_subspace_Banach_space}
  \Ann(W) := \{\alpha \in V^*: \alpha(w) = 0 \text{ for all } w \in W\}.
\end{equation}
The annihilator $\Ann(W)$ is a weak-$\star$ closed subspace of $V^*$ by Rudin \cite[Section 4.6, p. 96]{Rudin}. Let $\beta:V\times V \to \RR$ be a continuous bilinear form and let
\begin{equation}
  \label{eq:beta_flat}
  \beta^\flat: V \ni v \mapsto \beta(v,\cdot) \in V^*
\end{equation}
denote the induced continuous linear map. The $\beta$-\emph{complement} of $W$ in $V$ is
\begin{equation}
  \label{eq:Beta_complement_subspace_Banach_space}
  W^\beta := \{v \in V: \beta(v,w) = 0 \text{ for all } w \in W\}.
\end{equation}
We note that
\begin{equation}
  \label{eq:W_closed_implies_W_beta_closed}
  W \subset V \text{ closed } \implies W^\beta \subset V \text{ closed.}
\end{equation}
Indeed, if $W$ is a closed subspace of $V$, then $W^\beta$ is also a closed subspace of $V$ since
\[
  W^\beta = (\beta^\flat)^{-1}(\Ann(W))
\]
and $\Ann(W)$ is a closed subspace of $V^*$. Recall (see Appendix \ref{sec:Weakly_strongly_non-degenerate_bilinear_forms_Banach_spaces})
% PF12-10-2025 Add refs to Kobayashi, Lang, AMR, Marsden & Ratiu too
% TL12-13-2025: Perhaps make this a labelled definition since we use this characterization of symplectic subspaces later
% PF12-16-2025 Good point
that $\beta:V\times V\to\RR$ is \emph{weakly nondegenerate} if $\beta^\flat$ is injective and \emph{strongly nondegenerate} if $\beta^\flat$ is bijective, in which case $\beta^\flat:V\to V^*$ is an isomorphism of Banach spaces by the Open Mapping Theorem. In particular, we say that $\omega:V\times V\to\RR$ is a \emph{weak} (respectively, \emph{strong}) \emph{symplectic form} if $\omega$ is a weakly (respectively, strongly) nondegenerate skew-symmetric continuous bilinear form and write (see Kobayashi \cite[Section 3, Lemma 7.5.9, p. 248]{Kobayashi_differential_geometry_complex_vector_bundles}, Marsden and Weinstein \cite[Lemma, p. 123]{Marsden_Weinstein_1974}, or McDuff and Salamon \cite[Equation (2.1.1), p. 38]{McDuffSalamonSympTop3})
\begin{equation}
  \label{eq:Symplectic_complement_subspace_Banach_space}
  W^\omega := \{v \in V: \omega(v,w) = 0 \text{ for all } w \in W\}
\end{equation}
for the \emph{symplectic complement} of a (not necessarily closed) subspace $W$ of $V$.

Recall that a Banach space $V$ is an internal direct sum of closed subspaces, $V = W_1\oplus W_2$, if and only if there is a continuous linear idempotent, $\pi \in \End(V)$, such that $W_1=\Ran\pi$ and $W_2=\Ran(1-\pi)$ (see, for example, Kato \cite[Section III.3.4, p. 155]{Kato}). For the \emph{if} direction, let $v \in V$ and observe that $v = \pi v + (1-\pi)v \in W_1+W_2$ and, if $x \in W_1\cap W_2$, then $x = \pi x + (1-\pi)x = x + x$, so $x = 0$ and thus $V = W_1\oplus W_2$ as an algebraic direct sum of linear subspaces. If $\{x_n\}_{n\in\NN} \subset W_1$ is a sequence such that $\pi x_n \to v$ in $V$ as $n\to\infty$, then $\pi x_n = \pi^2 x_n \to \pi v$ as $n\to\infty$, so $\pi v = v$ and $v \in W_1$. Thus, $W_1$ is closed and similarly $W_2$ is closed, so $V$ is an internal direct sum of closed subspaces, $V = W_1\oplus W_2$. The proof of the \emph{only if} direction is similar. 

\begin{lem}[Closed subspaces of strongly symplectic Banach spaces and direct sum decompositions]
\label{lem:V_isomorphic_W_direct_sum_W_omega_Banach_spaces}
Let $(V, \omega)$ be a strongly symplectic Banach space and let $W$ be a closed subspace of $V$. Then $(W,\omega_W)$ is a strongly symplectic subspace with symplectic form $\omega_W := \omega|_W$ if and only if $V = W \oplus W^\omega$ as an internal direct sum of closed subspaces, $W$ and its symplectic complement $W^\omega$.
\end{lem}

\begin{rmk}[Strongly non-degenerate bilinear forms on Banach spaces and direct sum decompositions]
\label{rmk:Strongly_non-degenerate_skew_or_symmetric_bilinear_forms_Banach_spaces_and_direct_sums}  
The proof of Lemma \ref{lem:V_isomorphic_W_direct_sum_W_omega_Banach_spaces} extends without change if the role of the strong symplectic form $\omega$ on $V$ is replaced by that of a strongly non-degenerate skew-symmetric \emph{or} symmetric bilinear form $\beta$ on $V$. Thus, if $W$ is a closed subspace of $V$, then $\beta$ is non-degenerate on $W$ if and only if $V = W \oplus W^\beta$ as an internal direct sum of closed subspaces, $W$ and its $\beta$-complement $W^\beta$.
\end{rmk}  

\begin{proof}[Proof of Lemma \ref{lem:V_isomorphic_W_direct_sum_W_omega_Banach_spaces}]
Consider the \emph{if} direction, so we assume that $(W,\omega_W)$ is a strongly symplectic subspace. We claim that $V = W \oplus W^\omega$, that is, $W \cap W^\omega = (0)$ and $W + W^\omega = V$, where we recall that the symplectic complement $W^\omega$ is a closed subspace by \eqref{eq:W_closed_implies_W_beta_closed} because $W$ is a closed subspace by hypothesis.

We first show that $W \cap W^\omega = (0)$. Let $x \in W \cap W^\omega$. Since $x \in W^\omega$, we have $\omega(x, w) = 0$ for all $w \in W$ and because $\omega_W^\flat(x) = \omega(x, w)$, this implies $\omega_W^\flat(x) = 0 \in W^*$. But $\omega_W^\flat:W\to W^*$ is an isomorphism since $\omega_W$ is symplectic, so $\omega_W^\flat:W\to W^*$ is injective and thus $x = 0$. Therefore, $W \cap W^\omega = (0)$, as desired.

We next show that $W + W^\omega = V$. Let $z \in V$. We wish to decompose $z$ into components in $W$ and $W^\omega$. Define $\phi \in W^*$ by
\[
  \phi(w) := \omega(z, w), \quad \text{for all } w \in W.
\]
Since $W$ is symplectic, $\omega_W^\flat:W\to W^*$ is an isomorphism. Hence, there exists a unique $x \in W$ such that $\omega_W^\flat(x) = \phi$ and thus
\[
  \omega(x, w) = \phi(w) = \omega(z, w), \quad \text{for all } w \in W.
\]
We now choose $y := z - x$ and check that $y \in W^\omega$:
\[
  \omega(y, w) = \omega(z - x, w) = \omega(z, w) - \omega(x, w) = 0, \quad \text{for all } w \in W.
\]
Thus, $y \in W^\omega$ and we have $z = x + y$, with $x \in W$ and $y \in W^\omega$, as desired. This completes the proof that $V = W\oplus W^\omega$ as an internal algebraic direct sum and hence as an internal direct sum of closed subspaces since $W$ and $W^\omega$ are both closed.

Consider the \emph{only if} direction, so we assume that $V = W \oplus W^\omega$ and claim that $W$ is strongly symplectic with symplectic form $\omega_W = \omega|_W$ or, equivalently, that $\omega_W^\flat: W \to W^*$ is an isomorphism of Banach spaces.

We first show that $\omega_W^\flat$ is injective. Let $x \in W$ be such that $\omega_W^\flat(x) = 0$. Then $\omega(x, w) = \omega_W^\flat(x)w = 0$ for all $w \in W$ and, by definition, this implies that $x \in W^\omega$. Hence, $x \in W \cap W^\omega$ and since the sum is direct by assumption, we obtain $x = 0$, so $\omega_W^\flat$ is injective, as desired.

We next show that $\omega_W^\flat$ is surjective. Let $\alpha \in W^*$. We seek $x \in W$ such that $\alpha = \omega_W^\flat(x) = \omega_W(x, \cdot)$ on $W$. Since $V = W \oplus W^\omega$ as an internal direct sum of closed subspaces, we can extend $\alpha$ to $\tilde{\alpha} \in V^*$ by defining $\tilde{\alpha}$ to be zero on $W^\omega$ and thus
\[
  \tilde{\alpha}(w + v) = \alpha(w), \quad \text{for all } w \in W \text{ and } v \in W^\omega.
\]
Since $(V,\omega)$ is strongly symplectic, the map $\omega^\flat: V \to V^*$ is an isomorphism of Banach spaces. Hence, there exists a unique $z \in V$ such that $\tilde{\alpha} = \omega(z, \cdot)$. Because $V = W \oplus W^\omega$, we may write $z = x + y$ with $x \in W$ and $y \in W^\omega$. For any $w \in W$, we have
\begin{align*}
    \alpha(w) &= \tilde{\alpha}(w) \\
              &= \omega(z, w) \\
              &= \omega(x + y, w) \\
              &= \omega(x, w) + \omega(y, w).
\end{align*}
Since $y \in W^\omega$ and $w \in W$, we have $\omega(y, w) = 0$. Therefore, $\alpha(w) = \omega(x, w) = \omega_W(x, w)$ and because $w\in W$ was arbitrary, this means $\omega_W^\flat(x) = \omega_W(x, \cdot) = \alpha$. Thus, $\omega_W^\flat$ is surjective, as desired.

Since $\omega_W^\flat$ is a continuous bijection between Banach spaces, it is an isomorphism of Banach spaces by the Open Mapping Theorem. Thus, $W$ is symplectic and this completes the proof of Lemma \ref{lem:V_isomorphic_W_direct_sum_W_omega_Banach_spaces}.
\end{proof}

We do not require the following infinite-dimensional analogue of a well-known result for linear subspaces of finite-dimensional symplectic vector spaces, but we include the statement for the sake of completeness.

\begin{lem}[Double symplectic complements of subspaces of weakly symplectic Banach spaces]
\label{lem:Double_symplectic_complements_ subspaces_weakly_symplectic_Banach_spaces}
(See Kobayashi \cite[Lemma 7.5.9, p. 248]{Kobayashi_differential_geometry_complex_vector_bundles} or, with slightly different hypotheses, Marsden and Weinstein \cite[Section 3, Lemma, p. 123]{Marsden_Weinstein_1974}.)  
Let $(V, \omega)$ be a weakly symplectic Banach space. If $W$ is a closed subspace of $V$, then $W^{\omega\omega} = W$.
\end{lem}

Given Lemma \ref{lem:V_isomorphic_W_direct_sum_W_omega_Banach_spaces}, the proof of Lemma \ref{lem:Transverse_intersection_countably_many_symplectic_subspaces_is_symplectic} adapts without change to give the following analogue for strongly symplectic Banach spaces,

\begin{lem}[Transverse intersection of symplectic subspaces of a strongly symplectic Banach space is symplectic]
\label{lem:Transverse_intersection_symplectic_subspaces_Banach_space_is_symplectic} 
Let $(V,\omega)$ be a strongly symplectic Banach space. If $W_1,\ldots, W_m$ are strongly symplectic subspaces that intersect transversely in $V$ in the sense of Definition \ref{defn:Transverse_intersection_finitely_many_subspaces_vector_space} for $m\geq 2$, then $W_1 \cap\cdots\cap W_m$ is a strongly symplectic subspace of $V$.
\end{lem}

\subsection{Transverse intersections of almost symplectic submanifolds of an almost symplectic Banach manifold}
\label{subsec:Transverse_intersections_symplectic_submanifolds_symplectic_Banach_manifold}
For the extension of the results of Section \ref{subsec:Transverse_intersections_symplectic_subspaces_symplectic_Banach_spaces} to the setting of Banach manifolds, we partly rely on Kobayashi \cite[Section 7.5]{Kobayashi_differential_geometry_complex_vector_bundles}, Marsden and Ratiu \cite[Section 2.2]{Marsden_Ratiu_introduction_mechanics_symmetry}, Marsden and Weinstein \cite{Marsden_Weinstein_1974}, and Weinstein \cite{Weinstein_1969} for technical background.

Definition \ref{defn:Transverse_intersection_finitely_many_subspaces_vector_space} generalizes to the intersection of two or more embedded smooth submanifolds of a smooth Banach manifold.

\begin{defn}[Transverse intersection of submanifolds of a smooth Banach manifold]
\label{defn:Transverse_intersection_finitely_many_submanifolds_smooth_manifold}  
Let $M$ be a smooth Banach manifold. If $S_1,\ldots, S_m \subset M$ are embedded smooth submanifolds for $m \geq 2$, then the intersection $S_1\cap\cdots\cap S_n$ is \emph{transverse at $p$} if the bounded linear map
\begin{equation}
  \label{eq:Transverse_intersection_finitely_many_submanifolds_manifold_point}
  T_pM \ni v \mapsto \left([v]_{S_1},\ldots,[v]_{S_m}\right) \in T_pM/T_pS_1 \oplus \cdots \oplus T_pM/T_pS_m, 
\end{equation}
is an epimorphism of Banach spaces. The intersection $S_1\cap\cdots\cap S_n$ is \emph{transverse}, denoted $S_1 \transv\cdots\transv S_m$, if the map
\begin{equation}
  \label{eq:Transverse_intersection_finitely_many_submanifolds_manifold}
  TM \ni v \mapsto \left([v]_{S_1},\ldots,[v]_{S_m}\right) \in TM/TS_1 \oplus \cdots \oplus TM/TS_m, 
\end{equation}
is an epimorphism of Banach vector bundles or, equivalently, if \eqref{eq:Transverse_intersection_finitely_many_submanifolds_manifold_point} is an epimorphism of Banach spaces for every $p \in S_1\cap\cdots\cap S_n$. 
\end{defn}

Similarly, Remark \ref{rmk:Transverse_intersections_non-degenerate_symmetric_skew-symmetric_bilinear_forms} also generalizes to the intersection of embedded smooth submanifolds of a smooth Banach manifold.

\begin{rmk}[Transverse intersections and non-degenerate symmetric or skew-symmetric bilinear forms on tangent bundles]
\label{rmk:Transverse_intersections_non-degenerate_symmetric_skew-symmetric_bilinear_forms_manifold}
If $\beta:TM\times TM \to \RR$ is a strongly non-degenerate symmetric or skew-symmetric real bilinear on the tangent bundle $TM$ of a smooth Banach manifold $M$ and $S\subset M$ is an embedded smooth submanifold, then $\beta$ defines a normal bundle $N_{M/S}$ with fibers
\begin{equation}
  \label{eq:Beta_normal_bundle_S_in_Np}
  N_{M/S}|_p\ := \{v\in T_pM: \beta(v,w) = 0 \text{ for all } w \in T_pS\}, \quad\text{for all } p \in S.
\end{equation}
Given such a non-degenerate $\beta$, the transversality condition \eqref{eq:Transverse_intersection_finitely_many_submanifolds_manifold} is equivalent to the assertion that
\begin{equation}
  \label{eq:Transverse_intersection_finitely_many_submanifolds_non-degenerate_2-form}
  TM \ni v \mapsto \left(\pi_1v,\ldots,\pi_nv\right) \in N_{M/S_1} \oplus \cdots \oplus N_{M/S_m} 
\end{equation}
is an epimorphism of smooth Banach bundles where, writing $TM = TS_n\oplus N_{M/S_n}$, we let $\pi_n:TM \to N_{M/S_n}$ denote the idempotent map of smooth Banach vector bundles corresponding to the preceding direct sum decomposition, for $n=1,\ldots,m$. In particular, if the intersection $S := S_1\cap\cdots\cap S_n$ is transverse, then
\begin{equation}
  \label{eq:Transverse_intersection_finitely_many_submanifolds_internal_direct_sum_normal_bundles}
  N_{M/S} =  N_{M/S_1} \oplus \cdots \oplus N_{M/S_n}, 
\end{equation}
as an internal direct sum of Banach subbundles of the normal bundle $N_{M/S}$.

When the submanifolds $S_n$ have finite codimension in $M$, equation \eqref{eq:Transverse_intersection_finitely_many_submanifolds_manifold} yields the formula for the expected dimension of $S_1\cap\cdots\cap S_m$,
\begin{equation}
  \label{eq:Codimension_transverse_intersection_finitely_many_submanifolds}
  \codim(S_1\cap\cdots\cap S_m) = \sum_{n=1}^m \codim S_n.
\end{equation}
When $m=2$, then \eqref{eq:Transverse_intersection_finitely_many_submanifolds_manifold} is equivalent to the condition $TM=TS_1+TS_2$.
\qed
\end{rmk}

If $S$ is a smooth Banach manifold, modeled on a Banach space $F$, that is an embedded smooth submanifold of a smooth Banach manifold $M$, modeled on a Banach space $G$, then $F$ is a closed subspace of $G$ with closed complement $E\subset G$ such that $G = F\oplus E$ as an internal direct sum of Banach spaces and for each point $p \in S$ there is an open neighborhood $U\subset M$ of $p$ and a chart $\varphi:U\to F\oplus E$ such that $\varphi(U\cap S) = \varphi(U)\cap(F\oplus(0))$ (see, for example, Abraham, Marsden, and Ratiu \cite[Definition 3.2.1, p. 150]{AMR}, Hamilton \cite[Definition 4.2.1, p. 87]{Hamilton_1982bams}, or Lang \cite[Section II.2, p. 23]{Lang} in the case of Banach manifolds and Lee \cite[Chapter 5]{Lee_john_smooth_manifolds} in the case of finite-dimensional manifolds). If $\pi_E:F\oplus E \to E$ denotes the continuous projection (idempotent) onto the factor $E$, then $\Phi \equiv \pi_E\circ\varphi:U \to E$ is a submersion such that $U\cap S = \Phi^{-1}(0)$, a \emph{local defining function} for $S$.

\begin{lem}[Smooth maps and transverse intersections of submanifolds of a smooth Banach manifolds]
\label{lem:Smooth_manifolds_transverse_and_intersections_submanifolds}  
Let $S_1,\ldots,S_m$ be smooth Banach manifolds, modeled on Banach spaces $F_1,\ldots,F_m$, that is an embedded smooth submanifold of a smooth Banach manifold $M$ modeled on a Banach space $G$ and $E_n$ is closed complement of $F_n$ in $G$ for $n=1,\ldots,m$. Let $U \subset M$ be an open neighborhood of $p$ and
$f_n:U \to E_n$ be submersions such that $S_n\cap U = f_n^{-1}(0)$ for $n=1,\ldots,m$. Then the intersection $S_1\cap\cdots\cap S_m\cap U$ is transverse in the sense of Definition \ref{defn:Transverse_intersection_finitely_many_submanifolds_smooth_manifold} if and only if the smooth map
\[
  F = (f_1,\ldots,f_m):U \to E
\]
is a submersion, where $E := E_1\oplus\cdots\oplus E_m$. 
\end{lem}

\begin{proof}
Let $p\in S_1\cap\cdots\cap S_m\cap U$. Observe first that $T_pS_n = \Ker (df_n)_p$ and $T_pM/\Ker (df_n)_p \cong \Ran (df_n)_p$ (isomorphism of Banach spaces) and, since $f_n$ is a submersion, $\Ran (df_n)_p = E_n$ for $n=1,\ldots,m$. Writing $N_{M/S_n}|_p := T_pM/T_pS_n$, we obtain isomorphisms of Banach spaces,
\[
  N_{M/S_n}|_p \cong E_n, \quad\text{for } n=1,\ldots,m,
\]
and
\[
  N_{M/S_1}|_p \oplus \cdots \oplus N_{M/S_m}|_p \cong E_1\oplus\cdots\oplus E_m = E.
\]
Note that
\[
  \Ker(dF)_p = \Ker(df_1)_p \cap\cdots\cap \Ker(df_m)_p.
\]  
Consider the \emph{if} direction. Since $F$ is a submersion by assumption, $\Ran (dF)_p = E$ and $S\cap U = F^{-1}(0) = S_1\cap\cdots\cap S_m\cap U$ is a submanifold of $U\subset M$. Writing $N_{M/S}|_p := T_pM/T_pS$, we obtain an equality,
\[
  N_{M/S}|_p = N_{M/S_1}|_p \oplus \cdots \oplus N_{M/S_m}|_p,
\]
and thus an epimorphism,
\[
  TM|_p \to N_{M/S_1}|_p \oplus \cdots \oplus N_{M/S_m}|_p,
\]
so the intersection $S_1\cap\cdots\cap S_m$ is transverse at $p$ and, since $p$ was arbitrary, $S_1\cap\cdots\cap S_m\cap U$ is transverse in the sense of Definition \ref{defn:Transverse_intersection_finitely_many_submanifolds_smooth_manifold}.

Consider the \emph{only if} direction. By assumption, the intersection $S = S_1\cap\cdots\cap S_m$ is transverse at $p$ in the sense of Definition \ref{defn:Transverse_intersection_finitely_many_submanifolds_smooth_manifold}, so
\[
  N_{M/S}|_p = N_{M/S_1}|_p \oplus \cdots \oplus N_{M/S_m}|_p,
\]
and, in particular, $N_{M/S}|_p \cong E$ since $N_{M/S_n}|_p \cong E_n$ for $n=1,\ldots,m$. But $N_{M/S}|_p = T_pM/T_pS$ and, because $T_pS = \Ker(dF)_p$, we have $T_pM/T_pS = T_pM/\Ker(dF)_p \cong \Ran(dF)_p$. Combining these equalities and isomorphisms of Banach spaces gives $\Ran(dF)_p = E$, that is, $F$ is a submersion at $p$ and, since $p$ was arbitrary, $F:U\to E$ is a submersion.
\end{proof}  

Recall that an \emph{almost symplectic} manifold is a smooth manifold with a non-degenerate two-form. Lemma \ref{lem:Transverse_intersection_symplectic_subspaces_Banach_space_is_symplectic} and Definition \ref{defn:Transverse_intersection_finitely_many_submanifolds_smooth_manifold} yield the

\begin{cor}[Transverse intersection of almost symplectic submanifolds is almost symplectic]
\label{cor:Transverse_intersection_finitely_many_almost_symplectic_submanifolds_is_almost_symplectic}  
Let $(V,\omega)$ be an almost symplectic Banach manifold. If $W_1,\ldots, W_m$ are embedded almost symplectic submanifolds of $V$ for $m\geq 2$ and $W_1\transv \cdots \transv W_m$, then $W_1\cap\cdots\cap W_m$ is an embedded almost symplectic submanifold. If in addition $\omega$ is closed, so $(V,\omega)$ is symplectic, then $W_1,\ldots,W_m$ and $W_1\cap\cdots\cap W_m$ are embedded symplectic submanifolds. 
\end{cor}

\section[Real hypersurfaces and Donaldson's symplectic subspace criterion]{Real hypersurfaces and Donaldson's symplectic subspace criterion for  bounded real linear operators on complex Hilbert spaces}
\label{sec:Proof_analogue_Donaldson_symplectic_submanifold_criterion_hypersurface_Hilbert_spaces}
In this section, we prove Theorem \ref{mainthm:Donaldson_1996jdg_3_Hilbert_space_codomain} and Corollary \ref{maincorDonaldson_1996jdg_3_Hilbert_space_domain}. We begin with the 

\begin{proof}[Proof of Theorem \ref{mainthm:Donaldson_1996jdg_3_Hilbert_space_codomain}]
We shall adapt and extend Donaldson's proof of his Corollary \ref{cor:Donaldson_6} from the setting of finite-dimensional symplectic manifolds to that of Theorem \ref{mainthm:Donaldson_1996jdg_3_Hilbert_space_codomain}. We observe first that if $\sT \in \Hom_\RR^c(\sG_1,\sH_2)$ and $\psi \in \sG_2$, then the continuous $\CC$-valued $\RR$-linear function
\[
  \sT_\psi \equiv \langle\sT\cdot,\psi\rangle_{h_2}:\sG_1 \to \CC
\]
uniquely extends to a continuous $\CC$-valued $\RR$-linear function on $\sH_1$. Indeed, because $\sT = \sT'+\sT''$ by definition \eqref{eq:Complex_linear_and_anti-linear_operator_components} of $\sT'$ and $\sT''$, we have
\[
  \sT_\psi
  = \langle\sT\cdot,\psi\rangle_{h_2}
  = \langle\sT'\cdot,\psi\rangle_{h_2} + \langle\sT''\cdot,\psi\rangle_{h_2}.
\]
The forthcoming identities \eqref{eq:Adjoint_unbounded_complex_linear_operator} and 
\eqref{eq:Adjoint_unbounded_complex_antilinear_operator} in Definition \ref{defn:Adjoints_complex_linear_antilinear_real_linear_operators_on_complex_Hilbert_spaces} thus yield the equality
\begin{equation}
  \label{eq:sT_psi_continuous_extension_to_sH_1}
  \sT_\psi =  \langle\cdot,\sT^{\prime,*}\psi\rangle_{h_1} + \langle\sT^{\prime\prime,\bar*}\psi,\cdot\rangle_{h_1}:\sH_1\to\CC.
\end{equation}
The definition \eqref{eq:Complex_linear_and_anti-linear_operator_components} of $\sT'$ and $\sT''$, namely
\[
  \sT' = \frac{1}{2}(\sT - J_2\sT J_1)
  \quad\text{and}\quad
  \sT'' = \frac{1}{2}(\sT + J_2\sT J_1),
\]
and the hypotheses that $\Dom(\sT) = \sG_1$ and $J_k\in\End_\RR^c(\sH_k)$ continuously extends $J_k\in\End_\RR^c(\sG_k)$ for $k=1,2$ ensure that
\[
  \sG_1 \subset \Dom(\sT')  \quad\text{and}\quad \sG_1 \subset \Dom(\sT'').
\]
Note that $\sT^\intercal = \sT^{\prime,\intercal} + \sT^{\prime\prime,\intercal} = \sT^{\prime,*} + \sT^{\prime\prime,\bar *}$ by the forthcoming \eqref{eq:Complex_adjoints_in_terms_of_real_adjoints} in Lemma \ref{lem:Adjoints_complex_linear_antilinear_real_linear_operators_on_complex_Hilbert_spaces}. (Recall from Rudin \cite[Section 13.1, p. 348]{Rudin} that the real adjoint $\sT^\intercal \in \Hom_\RR(\sH_2^\RR,\sH_1^\RR)$ is uniquely determined by the relation
\begin{equation}
  \label{eq:Adjoint_unbounded_real_linear_operator}
  \langle v,\sT^\intercal w\rangle_{g_1} = \langle\sT v,w\rangle_{g_2},
  \quad\text{for all } v \in \Dom(\sT) \text{ and } w \in \Dom(\sT^\intercal),
\end{equation}
when $\Dom(\sT)$ is dense in $\sH_1$.) A similar calculation
% PF12-18-2025 Add it!
shows that the hypothesis that $\Dom(\sT^\intercal) = \sG_2$ ensures that
\[
  \sG_2 \subset \Dom(\sT^{\prime,*})  \quad\text{and}\quad \sG_2 \subset \Dom(\sT^{\prime\prime,\bar *}).
\]
In particular, setting
\[
  \sT_\psi' \equiv \langle\sT'\cdot,\psi\rangle_{h_2}:\sG_1 \to \CC
  \quad\text{and}\quad
  \sT_\psi'' \equiv \langle\sT''\cdot,\psi\rangle_{h_2}:\sG_1 \to \CC,
\]
and noting that $\sT^{\prime\prime,\bar*}\psi \in \sH_1$ and $\psi \in \sG_2 \implies \sT^{\prime,*}\psi \in \sH_1$, we see that
\[
  \sT_\psi' = \langle\cdot,\sT^{\prime,*}\psi\rangle_{h_2}:\sH_1\to\CC 
  \quad\text{and}\quad
  \sT_\psi'' = \langle\cdot,\sT^{\prime\prime,\bar*}\psi\rangle_{h_2}:\sH_1\to\CC
\]
are continuous $\CC$-valued functions on $\sH_1$ that are complex linear and antilinear, respectively. In particular,
$\sT_\psi$ is a continuous $\CC$-valued $\RR$-linear function on $\sH_1$. 

Every separable Hilbert space has a countable complete orthonormal basis (see, for example, Reed and Simon \cite[Section II.3, Theorem II.7, p. 47]{Reed_Simon_v1} or Rynne and Youngson \cite[Section 3.4, Theorem 3.52, p. 80]{Rynne_Youngson_linear_functional_analysis}), so we may choose a complete $h_2$-orthonormal basis $\{\psi_n\}_{n=1}^\infty$ for the complex Hilbert subspace $\sG_2 \subset \sH_2$ and, because $\sG_2 \subset \sH_2$ is a dense Hilbert subspace, $\{\psi_n\}_{n=1}^\infty$ is also a complete $h_2$-orthonormal basis for the complex Hilbert space $\sH_2$.
% PF12-17-2025 Cite ref or explain above assertion
We define a sequence of continuous $\CC$-valued $\RR$-linear functions on $\sG_1$ by
\begin{equation}
  \label{eq:sT_n_sG_1}
  \sT_n
  := \langle\sT\cdot,\psi_n\rangle_{h_2}
  % TL12-13-2025: The "c" superscript denotes "continuous"?  (It may have been previously defined and didn't register on my reading)
  %PF12-15-2025 Yes
  \in \Hom_\RR^c(\sG_1,\CC), \quad\text{for all } n \in \NN.
\end{equation}
By analogy with \eqref{eq:sT_n_sG_1}, we define
\begin{subequations}
  \label{eq:sT_n_prime_and_primeprime_sG_1}
  \begin{align}
    \label{eq:sT_n_prime_sG_1}
    \sT_n'
    &:= \langle\sT'\cdot,\psi_n\rangle_{h_2} \in \Hom_\RR^c(\sG_1,\CC)^{1,0},
    \\
    \label{eq:sT_n_primeprime_sG_1}
    \sT_n''
    &:= \langle\sT''\cdot,\psi_n\rangle_{h_2} \in \Hom_\RR^c(\sG_1,\CC)^{0,1}, \quad\text{for all } n \in \NN,
  \end{align}    
\end{subequations}
and observe that, since $\sT = \sT'+\sT''$ by definition \eqref{eq:Complex_linear_and_anti-linear_operator_components} of $\sT'$ and $\sT''$,
\[
  \sT_n = \sT_n' + \sT_n'',
\]
as we expect from \eqref{eq:sTprime_and_sTprimeprime_maps_sH_to_CC}.

By virtue of \eqref{eq:sT_psi_continuous_extension_to_sH_1}, the functions $\sT_n:\sG_1\to\CC$ and $\sT_n':\sG_1\to\CC$ and $\sT_n'':\sG_1\to\CC$ uniquely extend to continuous $\CC$-valued functions on $\sH_1$ by
\begin{equation}
  \label{eq:sT_n}
  \sT_n
  = \sT_n' + \sT_n''
  = \langle\cdot,\sT^{\prime,*}\psi_n\rangle_{h_1} + \langle\sT^{\prime\prime,\bar*}\psi_n,\cdot\rangle_{h_1}
  \in \Hom_\RR^c(\sH_1,\CC), \quad\text{for all } n \in \NN,
\end{equation}
where $\sT_n'$ is complex linear, $\sT_n''$ is complex antilinear, and $\sT_n$ is real linear.

For each $n \in \NN$, there is a continuous complex linear map, 
\begin{equation}
  \label{eq:pi_n}
  \pi_n:\Hom_\RR^c(\sG_1,\sH_2) \ni \sT \mapsto \sT_n \in \Hom_\RR^c(\sH_1,\CC).
\end{equation}
The maps $\pi_n$ in \eqref{eq:pi_n} are surjective and thus submersions (by linearity) for all $n\in\NN$. Indeed, if $n\in\NN$ and $\sL \in \Hom_\RR^c(\sH_1,\CC)$, then we may define
\[
  \sT := \sL \otimes \psi_n \in \Hom_\RR^c(\sH_1,\sH_2) \subset \Hom_\RR^c(\sG_1,\sH_2), 
\]
and observe that, since $\|\psi_n\|_{h_2}=1$,
\[
  \pi_n\sT = \langle\sT\cdot,\psi_n\rangle_{h_2}
  = \langle\sL(\cdot)\otimes \psi_n,\psi_n\rangle_{h_2}
  = \langle\psi_n,\psi_n\rangle_{h_2}\sL = \sL.
\]  
Thus, $\pi_n$ in \eqref{eq:pi_n} is surjective, and hence a submersion, as claimed. Morover, if $\sT \in \Hom_\RR^c(\sG_1,\sH_2)$ is surjective, then $\sT_n \neq 0 \in \Hom_\RR^c(\sH_1,\CC)$ for any $n\in \NN$. (Otherwise, if $\sT_n = 0$ for some $n$, then the definition \eqref{eq:sT_n} of $\sT_n$ implies that $\Ran\sT \perp \CC\psi_n$ in $\sH_2$, contradicting surjectivity of $\sT:\sG_1\to\sH_2$.)

We write
\[
  \Epim_\RR^c(\sG_1,\sH_2)
  \subset
  \Hom_\RR^c(\sG_1,\sH_2)
\]
for the open subset of bounded real linear operators from $\sG_1$ to $\sH_2$ that are surjective and thus $\sT_n \neq 0\in \Hom_\RR^c(\sH_1,\CC)$, for each $n \in \NN$. We define an open subset of $\Epim_\RR^c(\sG_1,\sH_2)$ for each $n \in \NN$,
%PF12-16-2025 Clarify meaning of norms
\begin{equation}
  \label{eq:sU_n}
  \sU_n
  :=
  \left\{\sT \in \Epim_\RR^c(\sG_1,\sH_2): \|\sT_n''\|_{\Hom(\sH_1,\CC)} \neq \|\sT_n'\|_{\Hom(\sH_1,\CC)} \right\}
  =
  \Epim_\RR^c(\sG_1,\sH_2)\, \less\, \sS_n,
\end{equation}
where the closed subset $\sS_n \subset \Epim_\RR^c(\sG_1,\sH_2)$ is defined by the closed subset of $\Hom_\RR^c(\sH_1,\CC) \less \{0\}$,
\begin{equation}
  \label{eq:cS}
  \cS := \left\{\sL \in \Hom_\RR^c(\sH_1,\CC) \less \{0\}: \|\sL''\|_{\Hom(\sH_1,\CC)} = \|\sL'\|_{\Hom(\sH_1,\CC)} \right\},
\end{equation}
which is a smooth real hypersurface in $\Hom_\RR^c(\sH_1,\CC) \less \{0\}$ by Remark \ref{rmk:Symplectic_form_on kernel_generic_T}, and 
\begin{equation}
  \label{eq:sS_n}
  \sS_n := \pi_n^{-1}(\cS)
  = \left\{\sT \in \Epim_\RR^c(\sG_1,\sH_2): \|\sT_n''\|_{\Hom(\sH_1,\CC)} = \|\sT_n'\|_{\Hom(\sH_1,\CC)} \right\},
\end{equation}
which is a smooth real hypersurface in $\Epim_\RR^c(\sG_1,\sH_2)$ since $\pi_n$ in \eqref{eq:pi_n} is a submersion.

Proposition \ref{mainprop:Donaldson_1996jdg_3} implies that for each $\sT \in \sU_n$, the restriction of the strong symplectic form $\omega_1$ on $\sH_1$,
\[
  \omega_1\restriction_{\Ker\sT_n}
\]
is a strong symplectic form on the closed real codimension-two linear subspace
\[
  \Ker(\sT_n:\sH_1 \to \CC) \subset \sH_1.
\]
For each $\sT \in \Epim_\RR^c(\sG_1,\sH_2)$, the intersections $\Ker\sT_1\cap\cdots\cap\Ker\sT_m$ are transverse for any $m\geq 2$ by Lemma \ref{lem:Smooth_manifolds_transverse_and_intersections_submanifolds}, noting that $\sH_2 = \oplus_{n=1}^m\CC\psi_n$ and
\[
  \Ker(\sT_n:\sH_1 \to \CC) = \Ker(\sT_n\otimes\psi_n:\sH_1 \to \CC\psi_n).
\]
Noting that $\Ker(\sT:\sH_1\to\sH_2) = \Ker(\sT:\sG_1\to\sH_2)$ since $\Dom(\sT)=\sG_1$, we have
\[
  \Ker\left(\sT:\sH_1\to\sH_2\right) 
  = \bigcap_{n=1}^m\Ker\left(\sT_n\otimes\psi_n:\sH_1\to\CC\psi_n\right)
  = \bigcap_{n=1}^m\Ker\left(\sT_n:\sH_1\to\CC\right).
\]
%PF12-12-2025 Prob cut
% The preceding equality follows from that fact that $v \in \Ker\sT$ is equivalent by \eqref{eq:Adjoint_unbounded_linear_operator} to $v \in (\Ran\sT^*)^\perp$ or equivalently, $\langle v, \sT^*w\rangle_{h_2} = 0$, for all $w\in \sG_2$, or equivalently,
% \[
%   \sT_nv \equiv \langle v, \sT^*\psi_n\rangle_{g_2} = 0, \quad\text{for } n=1,\ldots,m,
% \]
% that is, $v \in \Ker\sT_n$ for $n=1,\ldots,m$.
Hence, $\Ker(\sT:\sH_1\to\sH_2) \subset (\sH_1,\omega_1)$ is a symplectic subspace by Lemma \ref{lem:Transverse_intersection_symplectic_subspaces_Banach_space_is_symplectic}. This completes the proof of Theorem \ref{mainthm:Donaldson_1996jdg_3_Hilbert_space_codomain}.
\end{proof}

\begin{proof}[Proof of Corollary \ref{maincorDonaldson_1996jdg_3_Hilbert_space_domain}]
The conclusions follow immediately from an application of Theorem \ref{mainthm:Donaldson_1996jdg_3_Hilbert_space_codomain} to the adjoint operators $\sT^*:\sH_2\to\sH_1$.
\end{proof}  

\begin{rmk}[Choice of complete Hermitian orthonormal basis for $\sH_2$ in the proof Theorem \ref{mainthm:Donaldson_1996jdg_3_Hilbert_space_codomain}]
\label{rmk:Donaldson_1996jdg_3_Hilbert_space_codomain_complete_orthonormal_basis}  
In our applications, the complex linear components $\sT'$ of the operators $\sT\in\Epim_\RR(\sH_1,\sH_2)$ arising in the proof of Theorem \ref{mainthm:Donaldson_1996jdg_3_Hilbert_space_codomain} always obey the hypotheses of Lemma  \ref{lem:Eigenvalues_densely_defined_unbounded_linear_operators}. For any one such operator $\sT'$, we may choose $\{\varphi_n\}_{n\in\NN}$ to be a complete Hermitian orthonormal basis of $(\Ker\sT')^\perp\cap\sH_1$ and take $\{\psi_n\}_{n\in\NN}$ to be the corresponding complete Hermitian orthonormal basis of $\Ran\sT'$, where $\psi_n = \nu_n^{-1/2}\sT'\varphi_n$, for all $n\in\NN$, and $\{\nu_n\}_{n\in\NN}$ is the non-decreasing sequence of positive eigenvalues of $\sT^{\prime,*}\sT'$, repeated according to their multiplicities, while $\{\varphi_n\}_{n\in\NN}$ and $\{\psi_n\}_{n\in\NN}$ are the corresponding eigenvectors of $\sT^{\prime,*}\sT'$ and $\sT'\sT^{\prime,*}$, respectively. From equation \eqref{eq:Norm_sTn_prime_functional} in the forthcoming Remark \ref{rmk:Norms_complex_linear_antilinear_functionals_on_complex_Hilbert_spaces}, we have
\[
  \|\sT_n'\|_{\Hom_\CC^c(\sG_1,\CC)}^2
  = \langle\sT^{\prime,*}\psi_n,\sT^{\prime,*}\psi_n\rangle_{h_1}
  = \langle\sT'\sT^{\prime,*}\psi_n,\psi_n\rangle_{h_2}
  = \nu_n\langle\psi_n,\psi_n\rangle_{h_2},
\]
and thus
\[
  \|\sT_n'\|_{\Hom_\CC^c(\sG_1,\CC)} = \sqrt{\nu_n}.
\]
In our applications, the complex linear components $\sT'$ will be Fredholm, so it suffices to choose in addition a Hermitian orthonormal basis for the finite-dimensional complex subspace $(\Ran\sT')^\perp$ in order to obtain a complete Hermitian orthonormal basis for the complex Hilbert space $\sH_2$.
\qed
\end{rmk}

In applications of Theorem \ref{mainthm:Donaldson_1996jdg_3_Hilbert_space_codomain}, we encounter several different notions of adjoint operator on a complex Hilbert space. In the forthcoming Definition \ref{defn:Adjoints_complex_linear_antilinear_real_linear_operators_on_complex_Hilbert_spaces}, we define those adjoints and in Lemma \ref{lem:Adjoints_complex_linear_antilinear_real_linear_operators_on_complex_Hilbert_spaces}, we clarify the relationships among them.

\begin{defn}[Adjoints of complex linear, complex antilinear, and real linear operators on complex Hilbert spaces]
\label{defn:Adjoints_complex_linear_antilinear_real_linear_operators_on_complex_Hilbert_spaces}    
Let $\sH_k^\RR$ be a real Hilbert space with inner product $g_k$ and $g_k$-orthogonal almost complex structure $J_k$ and symplectic form $\omega_k = g_k(J_k\cdot,\cdot)$, for $k=1,2$. Let $\sH_k = (\sH_k^\RR,J_k)$ be the corresponding complex Hilbert space with Hermitian inner product $h_k = g_k-i\omega_k$, so that $g_k = \Real h_k$ and $\omega_k = -\Imag h_k$, for $k=1,2$. If $\sT' \in \Hom_\CC(\sH_1,\sH_2)$ is a bounded complex linear operator, then the \emph{complex adjoint} $\sT^{\prime,*} \in \Hom_\CC(\sH_2,\sH_1)$ of the complex linear operator $\sT' \in \Hom_\CC(\sH_1,\sH_2)$ is (as usual) uniquely determined by the equation
\begin{equation}
  \label{eq:Adjoint_unbounded_complex_linear_operator}
  \langle v,\sT^{\prime,*} w\rangle_{h_1} = \langle\sT' v,w\rangle_{h_2},
  \quad\text{for all } v \in \sH_1 \text{ and } w \in \sH_2.
\end{equation}
The \emph{complex antilinear adjoint} $\sT^{\prime\prime,\bar\star} \in \Hom_{\bar\CC}^c(\sH_2,\sH_1)$ of a bounded complex antilinear operator $\sT'' \in \Hom_{\bar\CC}(\sH_1,\sH_2)$ is uniquely determined by the equation 
\begin{equation}
  \label{eq:Adjoint_unbounded_complex_antilinear_operator}
  \langle v,\sT^{\prime\prime,\bar\star} w\rangle_{h_1} = \overline{\langle\sT'' v,w\rangle}_{h_2}
   = \langle w,\sT'' v\rangle_{h_2},
  \quad\text{for all } v \in \sH_1 \text{ and } w \in \sH_2.
\end{equation}
The \emph{real adjoint} $\sT^\intercal \in \Hom_\RR(\sH_2,\sH_1)$ of a bounded real linear operator $\sT \in \Hom_\RR(\sH_1,\sH_2)$ is uniquely determined by the equation
\begin{equation}
  \label{eq:Adjoint_unbounded_real_linear_operator}
  \langle v,\sT^\intercal w\rangle_{g_1} = \langle\sT v,w\rangle_{g_2},
\end{equation}
for all $v \in \sH_1 \text{ and } w \in \sH_2$.
\qed
\end{defn}  

\begin{lem}[Relations among adjoints of complex linear, complex antilinear, and real linear operators on complex Hilbert spaces]
\label{lem:Adjoints_complex_linear_antilinear_real_linear_operators_on_complex_Hilbert_spaces}  
Continue the notation of Definition \ref{defn:Adjoints_complex_linear_antilinear_real_linear_operators_on_complex_Hilbert_spaces}. If $\sT' \in \Hom_\CC^c(\sH_1,\sH_2)$ is a bounded complex linear operator, $\sT'' \in \Hom_{\bar\CC}^c(\sH_1,\sH_2)$ is a bounded complex antilinear operator, and $\sT = \sT'+\sT'' \in \Hom_\RR^c(\sH_1,\sH_2)$ is a bounded real linear operator, then
\begin{equation}
  \label{eq:Complex_adjoints_in_terms_of_real_adjoints}
  \sT^{\prime,\intercal} = \sT^{\prime,*},
  \quad
  \sT^{\prime\prime,\intercal} = \sT^{\prime\prime,\bar\star},
  \quad\text{and}\quad 
  \sT^\intercal = \sT^{\prime,*} + \sT^{\prime\prime,\bar *}.
\end{equation}
\end{lem}  

\begin{proof}
We begin by proving the first equality in \eqref{eq:Complex_adjoints_in_terms_of_real_adjoints}, which is equivalent by \eqref{eq:Adjoint_unbounded_complex_linear_operator} to the first equality below,
\[
  \langle v, \sT^{\prime,\intercal}w \rangle_{h_1}
  =
  \langle \sT'v, w \rangle_{h_2}
  =
  \langle v, \sT^{\prime,*}w \rangle_{h_1},
  \quad\text{for all } v \in \sH_1 \text{ and } w \in \sH_2,
\]
and the second equality is the definition of $\sT^{\prime,*}$. Noting that $g_k = \Real\, h_k$, for $k=1,2$, and $\Imag\, z = -i\Real\,(iz)$ for any $z\in\CC$, we have for all $v\in\sH_1$ and $w\in\sH_2$ that
\begin{align*}
  \langle \sT'v, w \rangle_{h_2}
  &= \Real\,\langle \sT'v, w \rangle_{h_2} + i\Imag\,\langle \sT'v, w \rangle_{h_2}
  \\
  &= \Real\,\langle \sT'v, w \rangle_{h_2} - i\Real\,\langle i\sT'v, w \rangle_{h_2}
  \\
  &= \Real\,\langle \sT'v, w \rangle_{h_2} - i\Real\,\langle \sT'Jv, w \rangle_{h_2}
  \\
  &= \Real\,\langle v, \sT^{\prime,\intercal}w \rangle_{h_1} - i\Real\,\langle Jv, \sT^{\prime,\intercal}w \rangle_{h_1}
  \\
  &= \Real\,\langle v, \sT^{\prime,\intercal}w \rangle_{h_1}
    - i\Real\left(i\langle v, \sT^{\prime,\intercal}w \rangle_{h_1}\right)
  \\
  &= \Real\,\langle v, \sT^{\prime,\intercal}w \rangle_{h_1}
    + i\Imag\,\langle v, \sT^{\prime,\intercal}w \rangle_{h_1}
  \\
  &= \langle v, \sT^{\prime,\intercal}w \rangle_{h_1}.
\end{align*}
This proves the first equality in \eqref{eq:Complex_adjoints_in_terms_of_real_adjoints}. We now prove the second equality in \eqref{eq:Complex_adjoints_in_terms_of_real_adjoints}, which is equivalent by \eqref{eq:Adjoint_unbounded_complex_antilinear_operator} to the first equality below,
\[
  \langle v, \sT^{\prime\prime,\intercal}w \rangle_{h_1}
  =
  \langle w, \sT''v \rangle_{h_2}
  =
  \langle v, \sT^{\prime\prime,\bar\star}w \rangle_{h_1},
  \quad\text{for all } v \in \sH_1 \text{ and } w \in \sH_2,
\]
and the second equality is the definition of $\sT^{\prime\prime,\bar\star}$. By taking complex conjugates, the preceding equalities are equivalent to
\[
  \langle \sT^{\prime\prime,\intercal}w,v \rangle_{h_1}
  =
  \langle \sT''v, w\rangle_{h_2}
  =
  \langle \sT^{\prime\prime,\bar\star}w,v \rangle_{h_1},
  \quad\text{for all } v \in \sH_1 \text{ and } w \in \sH_2.
\]
We observe that
\begin{align*}
  \langle \sT''v, w \rangle_{h_2}
  &= \Real\,\langle \sT''v, w \rangle_{h_2} + i\Imag\,\langle \sT''v, w \rangle_{h_2}
  \\
  &= \Real\,\langle \sT''v, w \rangle_{h_2} - i\Real\,\langle i\sT''v, w \rangle_{h_2}
  \\
  &= \Real\,\langle \sT''v, w \rangle_{h_2} + i\Real\,\langle \sT''Jv, w \rangle_{h_2}
  \\
  &= \Real\,\langle v, \sT^{\prime\prime,\intercal}w \rangle_{h_1} + i\Real\,\langle Jv, \sT^{\prime\prime,\intercal}w \rangle_{h_1}
  \\
  &= \Real\,\langle v, \sT^{\prime\prime,\intercal}w \rangle_{h_1}
    + i\Real\left(i\langle v, \sT^{\prime\prime,\intercal}w \rangle_{h_1}\right)
  \\
  &= \Real\,\langle v, \sT^{\prime\prime,\intercal}w \rangle_{h_1}
    - i\Imag\,\langle v, \sT^{\prime\prime,\intercal}w \rangle_{h_1}
  \\
  &= \overline{\Real\,\langle v, \sT^{\prime\prime,\intercal}w \rangle_{h_1}
    + i\Imag\,\langle v, \sT^{\prime\prime,\intercal}w \rangle}_{h_1} 
  \\
  &= \overline{\langle v, \sT^{\prime\prime,\intercal}w \rangle}_{h_1}
  \\
  &= \langle \sT^{\prime\prime,\intercal}w,v \rangle_{h_1}.
\end{align*}
This proves the second equality in \eqref{eq:Complex_adjoints_in_terms_of_real_adjoints}. Lastly, for $\sT = \sT'+\sT''$, we observe that the first two equalities in \eqref{eq:Complex_adjoints_in_terms_of_real_adjoints} yield
\[
  \sT^\intercal = \sT^{\prime,\intercal} + \sT^{\prime\prime,\intercal} = \sT^{\prime,*} + \sT^{\prime\prime,\bar *}.
\]
This proves the third equality in \eqref{eq:Complex_adjoints_in_terms_of_real_adjoints} and completes the proof of Lemma \ref{lem:Adjoints_complex_linear_antilinear_real_linear_operators_on_complex_Hilbert_spaces}.
\end{proof}

To compute the norms of the complex linear or antilinear functionals arising in the proof of Theorem \ref{mainthm:Donaldson_1996jdg_3_Hilbert_space_codomain}, we may appeal to the

\begin{lem}[Norms of complex linear and antilinear functionals on complex Hilbert spaces]
\label{lem:Norms_complex_linear_antilinear_functionals_on_complex_Hilbert_spaces}
Continue the notation of Definition \ref{defn:Adjoints_complex_linear_antilinear_real_linear_operators_on_complex_Hilbert_spaces} with $k=1$. If $\sL' \in \Hom_\CC^c(\sH_1,\CC)$ is a bounded complex linear function and $\sL'' \in \Hom_{\bar\CC}^c(\sH_1,\CC)$ is a bounded complex antilinear function, then there unique vectors $u,w \in \sH_1$ such that
\begin{equation}
  \label{eq:Riesz_representation_complex_linear_antilinear_functionals}
  \sL' = \langle\cdot,u\rangle_{h_1}
  \quad\text{and}\quad
  \sL'' = \langle w,\cdot\rangle_{h_1}.
\end{equation}
Moreover,
\begin{equation}
  \label{eq:Norms_complex_linear_antilinear_functional_represented_by_inner_product_vector}
  \|\sL'\|_{\Hom_\CC^c(\sH_1,\CC)} = \|u\|_{h_1}
  \quad\text{and}\quad
  \|\sL''\|_{\Hom_{\bar\CC}^c(\sH_1,\CC)} = \|w\|_{h_1}.
\end{equation}
\end{lem}

\begin{proof}
The assertions for $\sL'$ follow from the Riesz Representation Theorem (see, for example, Roman \cite[Theorem 13.32, p. 351]{Roman_advanced_linear_algebra}. If $C:\CC \ni z \mapsto \bar z \in \CC$ denotes the real linear complex conjugation operator, then $C\circ\sL'' \in \Hom_\CC^c(\sH_1,\CC)$ is a bounded complex linear function, so there is a unique $w \in \sH_1$ such that
\[
  C\circ\sL'' = \langle\cdot,w\rangle_{h_1}
  \quad\text{and}\quad
  \|C\circ\sL''\|_{\Hom_{\bar\CC}^c(\sH_1,\CC)} = \|w\|_{h_1}.
\]
Thus,
\[
  \sL'' = C^2\circ\sL'' = \overline{\langle\cdot,w\rangle}_{h_1} = \langle w,\cdot\rangle_{h_1}
\]    
and
\[
  \|\sL''\|_{\Hom_{\bar\CC}^c(\sH_1,\CC)} = \|C\circ\sL''\|_{\Hom_{\bar\CC}^c(\sH_1,\CC)} = \|w\|_{h_1}.
\]
This completes the proof of Lemma \ref{lem:Norms_complex_linear_antilinear_functionals_on_complex_Hilbert_spaces}.
\end{proof}

We can apply Lemma \ref{lem:Norms_complex_linear_antilinear_functionals_on_complex_Hilbert_spaces} to compute
the norms of complex linear or antilinear functionals arising in the proof of Theorem \ref{mainthm:Donaldson_1996jdg_3_Hilbert_space_codomain}, as we note in the following

\begin{rmk}[Norms of complex linear and antilinear functionals in the proof of Theorem \ref{mainthm:Donaldson_1996jdg_3_Hilbert_space_codomain}]
\label{rmk:Norms_complex_linear_antilinear_functionals_on_complex_Hilbert_spaces}  
From the definition \eqref{eq:sT_n_prime_sG_1} of $\sT_n'$ and the complex adjoint relation \eqref{eq:Adjoint_unbounded_complex_linear_operator}, we see that
\[
  \sT_n' = \langle\cdot,\sT^{\prime,*}\psi_n\rangle_{h_1} \in \Hom_\CC^c(\sH_1,\CC),
\]
so the first equality in \eqref{eq:Norms_complex_linear_antilinear_functional_represented_by_inner_product_vector} gives
\begin{equation}
  \label{eq:Norm_sTn_prime_functional}
  \|\sT_n'\|_{\Hom_\CC^c(\sH_1,\CC)} = \|\sT^{\prime,*}\psi_n\|_{h_1}.
\end{equation}
Similarly, the definition \eqref{eq:sT_n_primeprime_sG_1} of $\sT_n''$ and the complex antilinear adjoint relation \eqref{eq:Adjoint_unbounded_complex_antilinear_operator} yields
\[
  \sT_n'' = \langle\sT^{\prime\prime,\bar *}\psi_n,\cdot\rangle_{h_1} \in \Hom_{\bar\CC}^c(\sH_1,\CC),
\]
so we obtain
\begin{equation}
  \label{eq:Norm_sTn_primeprime_functional}
  \|\sT_n''\|_{\Hom_{\bar\CC}^c(\sH_1,\CC)}
  = \|\sT^{\prime\prime,\bar *}\psi_n\|_{h_1}
\end{equation}
by applying the second equality in \eqref{eq:Norms_complex_linear_antilinear_functional_represented_by_inner_product_vector}. 
\qed
\end{rmk}

\chapter{Lower bounds for spectral gaps for elliptic  operators}
\label{chap:Lower_bounds_spectral_gaps_elliptic_operators}
Our goal in this chapter is to prove lower bounds for spectral gaps for self-adjoint elliptic pseudodifferential operators acting on sections of smooth vector bundles over closed smooth manifolds. In Section \ref{sec:Lower_bounds_spectral_gaps}, we prove Theorem \ref{mainthm:Donaldson_1996jdg_3_elliptic_operator_order_m_geq_d}, which gives such lower bounds for operators whose order is greater than or equal to the dimension of the underlying manifold. In Section \ref{sec:Lower_bounds_spectral_gaps_Dirac_operators}, we prove lower bounds for eigenvalues of coupled Dirac operators over closed smooth manifolds of arbitrary dimension in Proposition \ref{prop:Lott_5_1986} and Corollary \ref{cor:Lott_1986} and then prove the lower bound in Theorem \ref{mainthm:Lower_bound_spectral_gap_coupled_Dirac_operator} for spectral gaps for coupled Dirac operators over closed smooth four-manifolds. Section \ref{sec:Lower_bounds_spectral_gaps_Dirac_operators_vector_potentials} contains our proofs of Corollary \ref{maincor:Lower_bound_spectral_gap_coupled_Dirac_operator_plus_vector_potential}, which gives lower bounds for spectral gaps for coupled Dirac operators with vector potentials over closed manifolds, and Corollary \ref{maincor:Lower_bound_spectral_gap_coupled_Dirac_operator_plus_vector_potential_four-manifold}, which refines Corollary \ref{maincor:Lower_bound_spectral_gap_coupled_Dirac_operator_plus_vector_potential} for certain coupled Dirac operators with vector potentials over closed four-manifolds.

\section[Lower bounds for spectral gaps for elliptic operators]{Lower bounds for spectral gaps for high order elliptic pseudodifferential operators}
\label{sec:Lower_bounds_spectral_gaps}
Theorem \ref{mainthm:Donaldson_1996jdg_3_Hilbert_space}, and thus Corollary \ref{maincor:Donaldson_1996jdg_3_Hilbert_space_non-self-adjoint}, rely on our ability to identify spectral gaps \eqref{eq:Spectral_gap_T} of width $4\delta$ for the operator $T$. In our applications, the operator $T$ will depend on parameters that vary in a noncompact space and we would like to ensure that $\delta$ either remains uniformly bounded from below by a positive constant that is independent of those parameters or depends on a parameter in a way that we can control. Theorem \ref{mainthm:Donaldson_1996jdg_3_elliptic_operator_order_m_geq_d} provides an answer for elliptic pseudodifferential operators of sufficiently high order. 

\begin{proof}[Proof of Theorem \ref{mainthm:Donaldson_1996jdg_3_elliptic_operator_order_m_geq_d}]
The $L^2$ self-adjoint elliptic differential operator,
\[
  A:C^\infty(M;E) \to C^\infty(M;E),
\]
defines a self-adjoint unbounded operator,
\[
  A:L^2(M;E) \to L^2(M;E),
\]
with dense domain $W^{m,2}(M;E) \subset L^2(M;E)$. The resolvent operator,
\[
  R(z,A):L^2(M;E) \to L^2(M;E), \quad\text{for } z \in \rho(A),
\]
is compact since the embedding $W^{m,2}(M;E) \subset L^2(M;E)$ is compact by the Rellich--Kondrachov Embedding Theorem (see Adams and Fournier \cite[Theorem 6.3, p. 168]{AdamsFournier}) and the operator $R(z,A):L^2(M;E) \to W^{m,2}(M;E)$ is bounded. According to Remark \ref{rmk:Seeley_1967a}, the eigenvalue counting function $N(A;\cdot) = N(\cdot)$ in \eqref{eq:N_lambda} obeys
\[
  N(A;\lambda) \sim c_0(a_m^*a_m,g,h)\lambda^{d/m} \quad\text{as } \lambda \to \infty,
\]
where the positive constant $c_0(a_m^*a_m,g,h)$ is defined by the formula \eqref{eq:Seeley_1967_page_291_constant_second_displayed_equation_A} and depends only on the principal symbol $a_m$, Riemannian metric $g$ on $M$, and Hermitian metric $h$ on $E$. The asymptotic formula above means that (see Section \ref{sec:Weyl_asymptotic_formula_eigenvalues_elliptic_operator} and Battisti, Borsero, and Coriasco \cite[Theorem 1 and Equation (5), p. 799]{Battisti_Borsero_Coriasco_2016})
\[
  \left|N(A;\lambda) - c_0(a_m^*a_m,g,h)\lambda^{d/m}\right| \leq C(A)\lambda^{(d-1)/m},
\]
for a constant $C(A) \in [1,\infty)$ independent of $\lambda$. We choose $\lambda\geq 1$ large enough that 
\[
  C(A)\lambda^{d/m-1/m} \leq \frac{1}{2}c_0(a_m^*a_m,g,h)\lambda^{d/m},
\]
which is equivalent to $\lambda^{1/m} \geq 2C(A)/c_0(a_m^*a_m,g,h)$, that is,
\begin{equation}
  \label{eq:Lower_bound_lambda_for_NAlambda}
  \lambda \geq \left(2C(A)/c_0(a_m^*a_m,g,h)\right)^m.
\end{equation} 
For all $\lambda$ large enough that inequality \eqref{eq:Lower_bound_lambda_for_NAlambda} holds, we thus have
\[
  \frac{1}{2}c_0(a_m^*a_m,g,h)\lambda^{d/m} \leq N(A;\lambda) \leq \frac{3}{2}c_0(a_m^*a_m,g,h)\lambda^{d/m},
\]
so the number $N(A;\lambda)$ of eigenvalues $\lambda_k$ (repeated according to their multiplicity) of $A$ such that $0 \leq |\lambda_k| < \lambda$ lies within the above range. If the absolute values $|\lambda_k|$ were all distinct, then the \emph{average gap size} between distinct values $|\lambda_k|$ in the interval $0,\lambda)$ would be $\lambda/N(A;\lambda)$ and hence the preceding bounds for $N(A;\lambda)$ (without assuming that the absolute values $|\lambda_k|$ were all distinct) yields
\[
  \frac{2}{3}c_0(a_m^*a_m,g,h)^{-1}\lambda^{1-d/m} \leq \frac{\lambda}{N(A;\lambda)}
  \leq 2c_0(a_m^*a_m,g,h)^{-1}\lambda^{1-d/m}.
\]
By hypothesis we have $m\geq d$ and thus $1-d/m \geq 0$, so the preceding lower bound for $\lambda/N(A;\lambda)$ simplifies, when $\lambda \geq 1$ to
\[
  \frac{2}{3}c_0(a_m^*a_m,g,h)^{-1} \leq \frac{\lambda}{N(A;\lambda)}, \quad\text{for all } \lambda \geq 1.
\]
The preceding lower bound for the average gap size $\lambda/N(A;\lambda)$ depend only on $(a_m,g,h)$. Clearly, the \emph{maximum gap size} between distinct values $|\lambda_k|$ in the interval $[0,\lambda)$ must be greater than or equal to the average gap size $\lambda/N(A;\lambda)$ and if we denote that maximum gap size by $4\delta$, the lower bound in the preceding pair of inequalities yields the estimate
\begin{equation}
  \label{eq:Lower_bound_for_delta_for_gap}
  \delta \geq \frac{1}{6}c_0(a_m^*a_m,g,h)^{-1}.
\end{equation}
Therefore, by choosing the constant $\mu = \mu(A,g,h) > 0$ to be the lower bound on the right-hand side of the inequality \eqref{eq:Lower_bound_lambda_for_NAlambda}, there is a spectral gap
\[
  (\mu-2\delta,\mu+2\delta) \subset [0,\lambda)\cap\rho(|A|),
\]
for all $\lambda$ obeying \eqref{eq:Lower_bound_lambda_for_NAlambda} and $\delta = \delta(a_m,g,h)$ equal to the constant on the right-hand side of the inequality \eqref{eq:Lower_bound_for_delta_for_gap}. Here, the absolute value of the operator $A$ is defined by $|A| := (A^*A)^{1/2} = (A^2)^{1/2}$ with $\sigma(|A|)\less\{0\} = \{|\lambda_k|\}_{k=1}^\infty$. 
% COMMENT See https://mathoverflow.net/questions/173613/how-much-does-the-absolute-value-of-an-operator-behave-like-an-absolute-value
%COMMENT See https://math.stackexchange.com/questions/2524364/elementary-properties-of-absolute-value-operator
Consequently, we obtain
\[
  (-\mu-2\delta,-\mu+2\delta) \cup (\mu-2\delta,\mu+2\delta) \subset [-\lambda,\lambda]\cap\rho(A).
\]
This verifies the spectral gap condition \eqref{eq:Spectral_gap_A} and completes the proof of Theorem \ref{mainthm:Donaldson_1996jdg_3_elliptic_operator_order_m_geq_d}.
\end{proof}

\section{Lower bounds for spectral gaps for Dirac operators}
\label{sec:Lower_bounds_spectral_gaps_Dirac_operators}
For an analogue of Theorem \ref{mainthm:Donaldson_1996jdg_3_elliptic_operator_order_m_geq_d} that applies to elliptic pseudodifferential operators of order $m < d$, where $d$ is the dimension of the base manifold $X$, we restrict our attention to (generalized) Dirac operators. The forthcoming Proposition \ref{prop:Lott_5_1986} and Corollary \ref{cor:Lott_1986} provide a generalization of Lott \cite[Section IV, Proposition 5, p. 123 and Corollary, p. 125]{Lott_1986} from the case of Dirac operators on spin manifolds to coupled Dirac operators on \spinc manifolds. However, we simplify the proof of \cite[Section IV, Proposition 5, p. 123]{Lott_1986} by using the \emph{Kato Inequality} to bound the trace of $e^{-t\nabla_A^*\nabla_A}$ by a multiple of the trace of $e^{-t\Delta_g}$, thus avoiding the need to use the Feynman--Kac path integral or a detailed local analysis in geodesic normal coordinates. (The relevance of the Kato Inequality to lower bounds on eigenvalues of $D_A^*D_A$ is mentioned by Vafa and Witten \cite[Section III, p. 270]{Vafa_Witten_1984}.) Like in the proof of \cite[Section IV, Proposition 5, p. 123]{Lott_1986}, one could use the \emph{Golden--Thompson Inequality} to bound the trace of $e^{-tD_A^*D_A}$ by the trace of $e^{-tR_A/2}e^{-t\nabla_A^*\nabla_A}e^{-tR_A/2}$ or $e^{-t\nabla_A^*\nabla_A}e^{-tR_A}$, where $D_A^*D_A = \nabla_A^*\nabla_A + R_A$ on $\Omega^0(W^+\otimes E)$ when $X$ is even-dimensional and on $\Omega^0(W\otimes E)$ when $X$ is odd-dimensional, and $R_A$ is a term that depends only on the curvatures of the connections on $W$ and $E$. However, the simplification offered by the Golden--Thompson Inequality is slight and we avoid its application.

The abstract or generalized Kato Inequality (see Hess, Schrader, and Uhlenbrock \cite[Equations (2.8) and ($2.8''$), p. 900]{HSU77} or Kato \cite{Kato_1972}) and its applications to domination of semigroups \cite[Theorem 2.15, p. 900]{HSU77} leads to the

\begin{thm}[Kato inequality and the heat operator for the covariant Laplacian]
\label{thm:Kato_inequality}
(See Hess, Schrader, and Uhlenbrock \cite[Equation (1.1), p. 27 and Theorem 3.1, p. 32]{HSU80}.)
Let $(E,H)$ be a Hermitian vector bundle of rank $r_E$ over a closed, smooth Riemannian manifold $(X,g)$. If $A$ is a smooth unitary connection on $E$, then
\begin{equation}
  \label{eq:Hess_Schrader_Uhlenbrock_1-1}
  \Tr\left(e^{-t\nabla_A^{*_g}\nabla_A}\right) \leq r_E\Tr\left(e^{-t\Delta_g}\right), \quad\text{for all } t > 0.
\end{equation}    
If $L \in C^\infty(i\fu(E))$ and $l \in C^\infty(X,\RR)$ is defined by $l_x := \inf\sigma(L_x)$, for each $x \in X$, and $\tilde l := \inf_{x\in X} l_x$, by analogy with \cite[Section 2, p. 31]{HSU80}, then 
\begin{equation}
  \label{eq:Hess_Schrader_Uhlenbrock_3-1}
  \Tr\left(e^{-t(\nabla_A^{*_g}\nabla_A + L)}\right)
  \leq r_E\Tr\left(e^{-t(\Delta_g+l)}\right)
  \leq r_Ee^{-t\tilde l}\Tr\left(e^{-t\Delta_g}\right).
\end{equation} 
\end{thm}

The covariant Laplacian $\nabla_A^*\nabla_A$ and the operator $\nabla_A^*\nabla_A + L$ in Theorem \ref{thm:Kato_inequality} are essentially self-adjoint\footnote{See Section \ref{sec:Spectral_theory_unbounded_operators}.} on $L^2(E)$ (see \cite[Corollary 2.4, p. 30 and second paragraph, p. 31]{HSU80}). Similarly, the square of the Dirac operator is well-known to be essentially self-adjoint on $L^2(W^+\otimes E)$ (see, for example Wolf \cite{Wolf_1972}). In Theorem \ref{thm:Kato_inequality}, the spectrum of $L_x \in \End(E_x)$ is denoted by $\sigma(L_x)$, for each $x \in X$. In our statement of Theorem \ref{thm:Kato_inequality}, we use the sign convention for Laplace operators that is standard in differential geometry, so $\Delta = -\sum_{i=1}^d \partial^2/\partial x_i^2$ on $C^\infty(\RR^d)$ and \emph{not} as \cite[Equation (3.3), p. 901]{HSU77}. Hence, our Laplace operators $\nabla_A^*\nabla_A$ and $\Delta_g$ are \emph{non-negative} and we replace $t$ by $-t$ when quoting results from \cite{HSU80}. This also explains why we reverse the definition of $l$ and $\tilde l$ in \cite[Section 2, p. 31]{HSU80} in the sense that we replace $\sup$ by $\inf$. Because $L_x \in \End(L_x)$ is Hermitian, we have by \eqref{eq:Spectral_radius_equals_norm_normal_operator} that
\[
  \|L_x\|_{\End(E_x)} = \sup_{\lambda \in \sigma(L_x)}|\lambda|,
\]
and so $l_x \geq -\|L_x\|_{\End(E_x)}$ for all $x\in X$ and thus we may substitute
\begin{equation}
  \label{eq:tildel_lower_bound}
  \tilde l \geq -\|L\|_{C^0(\End(E))}
\end{equation}
in the inequality \eqref{eq:Hess_Schrader_Uhlenbrock_3-1} provided by Theorem \ref{thm:Kato_inequality}. (Note that to achieve the lower bound \eqref{eq:tildel_lower_bound} we could have appealed to the simpler inequality provided by the inclusion \eqref{eq:Spectrum_bounded_operator_subset_disk}.) Alternatively, in order to obtain \eqref{eq:Hess_Schrader_Uhlenbrock_3-1}, one could instead apply the simpler inequality \eqref{eq:Hess_Schrader_Uhlenbrock_1-1} in conjunction with the Golden--Thompson Inequality:

\begin{thm}[Golden--Thompson Inequality]
\label{thm:Golden-Thompson_inequality}
(See Bikchentaev, Kittaneh, Sal Moslehian, and Seo \cite[Section 8.2, Theorem 8.2.3, p. 268]{Bikchentaev_Kittaneh_SalMoslehian_Seo_trace_inequalities} or Ruskai \cite[Theorem 4, p. 283]{Ruskai_1972}.)
If $A, B$ are self-adjoint operators on a Hilbert space $\sH$ that are bounded above and $A+B$ is essentially self-adjoint, then
\begin{equation}
  \label{eq:Golden-Thompson_inequality}
  \Tr(e^{A+B}) \leq \Tr(e^{A/2} e^B e^{A/2}).
\end{equation}
If in addition $\Tr(e^A) < \infty$ or $\Tr(e^B) < \infty$, then
\begin{equation}
  \label{eq:Golden-Thompson_inequality_refined}
  \Tr(e^{A+B}) \leq \Tr(e^A e^B).
\end{equation}
\end{thm}

Because our Laplace operators $\nabla_A^*\nabla_A$ and $\Delta_g$ are \emph{non-negative} (and thus bounded below), we can apply Theorem \ref{thm:Golden-Thompson_inequality} to the operators $-\nabla_A^*\nabla_A$ and $-\Delta_g$, since they are \emph{bounded above}. In the case of the Dirac operator $D$ on $C^\infty(W)$ over a manifold of \emph{even} dimension $d$, with $W = W^+\oplus W^-$, one has that $\rank_\CC W = 2^{d/2}$ and $\rank_\CC W^\pm = 2^{(d/2)-1}$ (see, for example, Berline, Getzler, and Vergne \cite[Section 3.2, Proposition 3.19, p. 109]{BerlineGetzlerVergne}). For convenience (and consistency with our application) in the statement and proof in the forthcoming Theorem \ref{prop:Lott_5_1986}, we shall assume that the dimension $d$ of the base manifold $X$ is even; the differences in its statement and proof when $d$ is odd are purely notational.

\begin{prop}[Lower bounds for eigenvalues of coupled Dirac operators over spin${}^c$ manifolds]
\label{prop:Lott_5_1986}
(Compare Lott \cite[Section IV, Proposition 5, p. 123]{Lott_1986} for the pure Dirac operator $D$ on a spin manifold.)  
Let $(X,g)$ be an oriented, smooth Riemannian manifold whose dimension $d\geq 2$ is even, $(\rho,W)$ be a spin${}^c$ structure on $X$, and $A$ be a smooth orthogonal connection on a smooth Riemannian vector bundle $(V,h)$
%PF4-3-2025 Should we make V Hermitian for consistency? It impacts value of r_V. 
of rank $r_V$ over $X$. If $\{\lambda_k(D_A)\}_{k=1}^\infty$ denotes the spectrum of eigenvalues of the Dirac operator $D_A$ on $\Omega^0(W\otimes V)$ as in \eqref{eq:Coupled_Dirac_operator}, then\footnote{Lott uses the sequence $\{\lambda_k\}_{k=1}^\infty$ to denote the eigenvalues of $D^*D$ on $L^2(W^+)$ whereas we use it to denote the eigenvalues of $D_A$ on $L^2(W^+\otimes V)$.} for any constant $\alpha>0$, 
\begin{equation}
  \label{eq:Lower_bounds_eigenvalues_square_coupled_Dirac_operator_raw}
  ke^{-\alpha} \leq 2^{(d/2)-1}r_V e^{-\alpha\|R_A\|\lambda_k^2}
  \left(1 + 4\left(\frac{\alpha}{\lambda_k^2}\right)^{-d/2}\vol(X,g)\left(\frac{2C_1}{d}\right)^{-d/2}\right),
\end{equation}
where $C_1 = C_1(g) > 0$ and $\|R_A\| = \|R_A\|_{C^0(\End(W^+\otimes V))}$ with 
\begin{equation}
  \label{eq:R_A}
  R_A := \frac{R}{4} + \frac{1}{2}\rho(F_{A_L}) + \rho(F_A) \in C^\infty(\End(W^+\otimes V)),
\end{equation}
and $A_L$ is the unitary connection on the Hermitian line bundle $L$ associated to the spin${}^c$ structure on $X$ and $F_{A_L}$ is its curvature and $R$ is the scalar curvature of the Riemannian metric $g$. If $\dim X = 4$, then we may replace $\rho(F_{A_L})$ and $\rho(F_A)$ in \eqref{eq:R_A} and thus \eqref{eq:Lower_bounds_eigenvalues_square_coupled_Dirac_operator_raw} by $\rho(F_{A_L}^+)$ and $\rho(F_A^+)$, respectively.  
\end{prop}

The geometric constant $C_1$ in inequality \eqref{eq:Lower_bounds_eigenvalues_square_coupled_Dirac_operator_raw} in Proposition \ref{prop:Lott_5_1986} has a lower bound in terms of $\mathrm{Diam}(X,g)$, $\vol(X,g)$, and $\Ric(g)$ according to Cheng and Li \cite{ChengLi}.
% PF3-19-2025 Get exact ref
We now proceed to the

\begin{proof}[Proof of Proposition \ref{prop:Lott_5_1986}]
We adapt and simplify the proof due to Lott of \cite[Section IV, Proposition 5, p. 123]{Lott_1986}. By applying Lemma \ref{lem:Eigenvalues_densely_defined_unbounded_linear_operators} to the operator $D_A^+ = D_A:L^2(W^+\otimes V) \to L^2(W^-\otimes V)$ and
\[
  D_A = \begin{pmatrix} 0 & D_A^{+,*} \\ D_A^+ & 0 \end{pmatrix}
  \quad\text{on}\quad
  L^2(W\otimes V) = L^2(W^-\otimes V) \oplus L^2(W^-\otimes V),
\]
in place of $\sT$ and $T$, respectively, we see that it suffices to consider the eigenvalues of $D_A^{+,*}D_A^+ = D_A^*D_A$ on $L^2(W^+\otimes V)$. Equation \eqref{eq:Bochner-Weitzenbock_formula_coupled_Dirac_operator} in the forthcoming Lemma \ref{lem:Bochner-Weitzenbock_identity_coupled_Dirac_operator} provides the Bochner--Weitzenb\"ock identity,
\[
  D_A^*D_A = \nabla_{A_W\otimes A}^*\nabla_{A_W\otimes A} + \frac{R}{4}
  + \frac{1}{2}\rho(F_{A_L})
  + \rho(F_A)
  \quad\text{on } \Omega^0(W^+\otimes V),
\]
where $R$ is the scalar curvature for the Riemannian metric $g$ on $X$ and $A_L$ is the connection on the Hermitian line bundle $L$ associated to the \spinc structure $(\rho,W)$; if $d=4$, we may replace $F_{A_L}$ and $F_A$ by $F_{A_L}^+$ and $F_A^+$ in the preceding formula. The operator $R_A \in \End(W^+\otimes V)$ in \eqref{eq:R_A} is Hermitian on the fibers of $W^+\otimes_\RR V$
%PF3-18-2025 Clearly, R_A is globally self-adjoint. Justify the stronger pointwise Hermitian assertion.
and our Bochner--Weitzenb\"ock identity may be written more compactly as
\begin{equation}
  \label{eq:Bochner-Weitzenbock_formula_coupled_Dirac_operator_for_Kato}
  D_A^*D_A = \nabla_{A_W\otimes A}^*\nabla_{A_W\otimes A} + R_A \quad\text{on } \Omega^0(W^+\otimes V).
\end{equation}
(Lott \cite[Section II, p. 118]{Lott_1986} assumes that $(X,g)$ is a (closed, smooth) spin manifold (of dimension $n\geq 2$) and restricts to the case of the Dirac operator $D$ on $\Omega^0(W)$ defined by a spin structure $(\rho,W)$, so his Bochner--Weitzenb\"ock identity for $D^*D$ simplifies to \cite[Section IV, p. 123]{Lott_1986}
\[
  D^*D = \nabla_{A_W}^*\nabla_{A_W} + \frac{R}{4} \quad\text{on } \Omega^0(W^+),
\]
but we shall not need this identity.)
From the inequality \eqref{eq:Hess_Schrader_Uhlenbrock_3-1} given by Theorem \ref{thm:Kato_inequality}, we obtain
\[
  \Tr\left(e^{-tD_A^*D_A}\right)
  =
  \Tr\left(e^{-t\nabla_A^*\nabla_A + R_A}\right)
  \leq
  r_{W^+\otimes V}e^{-t\tilde R}\Tr\left(e^{-t\Delta_g}\right), \quad\text{for all } t > 0,
\]
where the constant $\tilde R$ obeys $\tilde R \geq -\|R_A\|_{C^0(\End(W^+\otimes V)}$ by \eqref{eq:tildel_lower_bound} and, denoting $r_V := \rank_\RR(V)$,
\[
  r_{W^+\otimes V} := \rank_\CC(W^+\otimes_\RR V) = \rank_\CC(W^+)\rank_\RR(V) = 2^{(d/2)-1}r_V.
\]
According to Cheng and Li \cite[Section 2, Equation (2.9), p. 335]{ChengLi} (see Li \cite[Theorem 13.4, p. 139]{Li_geometric_analysis} for a more general result) , the heat kernel\footnote{Denoted by $H(t,x,y)$ in \cite[Section 2, Proof of Theorem 1, pp. 333--335]{ChengLi}.} for the Laplace operator $\Delta_g$ obeys the pointwise upper bound,
\[
  e^{-t\Delta_g}(t,x,x) - \frac{1}{\vol(X,g)} \leq 4\left(\frac{2C_1}{d}\right)^{-d/2} t^{-d/2},
  \quad\text{for all } (t,x) \in [0,\infty)\times X,
\]
where the geometric constant $C_1 = C_1(g) > 0$ is defined in \cite[Section 1, Lemma 1, p. 329]{ChengLi}. Consequently,
\begin{align*}
  \Tr\left(e^{-t\Delta_g}\right)
  &=
  \int_X e^{-t\Delta_g}(t,x,x)\,d\vol_g(x)
  \\
  &\leq
    \int_X\left(\frac{1}{\vol(X,g)} + 4t^{-d/2}\left(\frac{2C_1}{d}\right)^{-d/2}\right)\,d\vol_g(x)
  \\
  &= 1 + 4t^{-d/2}\vol(X,g)\left(\frac{2C_1}{d}\right)^{-d/2}, \quad\text{for all } t > 0. 
\end{align*}
Therefore,
\[
  \Tr\left(e^{-tD_A^*D_A}\right)
  \leq
  2^{(d/2)-1}r_V e^{-t\|R_A\|}\left(1 + 4t^{-d/2}\vol(X,g)\left(\frac{2C_1}{d}\right)^{-d/2}\right),
  \quad\text{for all } t > 0.
\]
By our hypothesis that $\{\lambda_l\}_{l=1}^\infty$ is the spectrum of unbounded operator $D_A$ on $L^2(W\otimes V)$, Lemma \ref{lem:Eigenvalues_densely_defined_unbounded_linear_operators} implies that the spectrum of the unbounded operator $D_A^*D_A$ on $L^2(W^+\otimes V)$ is $\{\lambda_l^2\}_{l=1}^\infty$. By the Spectral Mapping Theorem,
%PF3-19-2025 Get ref
the spectrum of the bounded operator $e^{-tD_A^*D_A}$ on $L^2(W^+\otimes V)$ is $\{e^{-t\lambda_l^2}\}_{l=1}^\infty$, and thus
\[
  \Tr\left(e^{-tD_A^*D_A}\right) = \sum_{l=1}^\infty e^{-t\lambda_l^2}, \quad\text{for all } t > 0.
\]
Substituting $t = \alpha/\lambda_k^2$, for some $k \geq 1$ with $\lambda_k\neq 0$, and noting that $0 \leq \lambda_l^2/\lambda_k^2 \leq 1$ for $1 \leq l \leq k$ yields the lower bound,
\[
  \Tr\left(e^{-\alpha D_A^*D_A/\lambda_k^2}\right)
  =
  \sum_{l=1}^\infty e^{-\alpha\lambda_l^2/\lambda_k^2}
  =
  \sum_{l=1}^k e^{-\alpha\lambda_l^2/\lambda_k^2} + \sum_{l=k+1}^\infty e^{-\alpha\lambda_l^2/\lambda_k^2}
  \geq
  \sum_{l=1}^k e^{-\alpha\lambda_l^2/\lambda_k^2}
  \geq
  ke^{-\alpha}.
\]
Substituting $t = \alpha/\lambda_k^2$ into the right-hand side of the upper bound for the trace of $e^{-tD_A^*D_A}$ gives
\[
  \Tr\left(e^{-\alpha D_A^*D_A/\lambda_k^2}\right)
  \leq
  2^{(d/2)-1}r_V e^{-\alpha\|R_A\|\lambda_k^2}
  \left(1 + 4\left(\frac{\alpha}{\lambda_k^2}\right)^{-d/2}\vol(X,g)\left(\frac{2C_1}{d}\right)^{-d/2}\right).
\]
By combining the preceding upper and lower bounds for the trace of $e^{-\alpha D_A^*D_A/\lambda_k^2}$, we obtain the inequality \eqref{eq:Lower_bounds_eigenvalues_square_coupled_Dirac_operator_raw} and this concludes the proof of Proposition \ref{prop:Lott_5_1986}.
\end{proof}

Proposition \ref{prop:Lott_5_1986} yields the

\begin{cor}[Lower bounds for eigenvalues of coupled Dirac operators over spin${}^c$ manifolds]
\label{cor:Lott_1986}
(Compare Lott \cite[Section IV, Corollary, p. 125]{Lott_1986} for the pure Dirac operator $D$ on a spin manifold.)
Continue the notation and hypotheses of Proposition \ref{prop:Lott_5_1986}. Then for all integers $k \geq 2^{(d/2)-1} e^{d/2}$,
\begin{equation}
  \label{eq:Lower_bounds_eigenvalues_square_coupled_Dirac_operator}
  \lambda_k^2(D_A) \geq C_1\left(4\vol(X,g)\right)^{-2/d}\left(2^{-(d/2)+1} e^{-d/2} k - 1\right)^{2/d}
  - \|R_A\|_{C^0(W^+\otimes V)}.
\end{equation}
\end{cor}

\begin{proof}
The conclusion follows from the inequality \eqref{eq:Lower_bounds_eigenvalues_square_coupled_Dirac_operator_raw} exactly as in the proof of Lott \cite[Section IV, Corollary, p. 125]{Lott_1986}. The only changes to his proof are that
\begin{inparaenum}[\itshape i\upshape)]
\item we replace his numerical factor $2^{\lfloor d/2\rfloor}$ by $2^{(d/2)-1}$, where the difference arises because we assume $d$ is even and consider the eigenvalues of $D_A^*D_A$ on $L^2(W^+\otimes V)$, and
\item we replace $\frac{1}{4}R_{\min}$ by $- \|R_A\|$.
\end{inparaenum}  
\end{proof}  

% PF3-20-2025 Remove redundant b^+ - b_1 odd hypotheses.

The following result due to Vafa and Witten \cite{Vafa_Witten_1984} provides upper bounds on the eigenvalues of the coupled Dirac operator.

\begin{thm}[Uniform upper bounds on eigenvalues of coupled Dirac operators]
\label{thm:Uniform_upper_bounds_eigenvalues_coupled_Dirac_operators}
(See Atiyah \cite[Section 1, Theorem 2, p. 252]{Atiyah_1985}, Leung and Xu \cite[Section 3, Theorem 2, p. 362]{Leung_Xu_2009}, and Vafa and Witten \cite[Section II, Equations (9) and (14), p. 260]{Vafa_Witten_1984}.)
If $(X,g)$ is a closed oriented smooth Riemannian manifold of dimension $d\geq 1$, then there is a constant $C \in [1,\infty)$ with the following significance. If $(\rho,W)$ is a spin${}^c$ structure on $(X,g)$, and $(V,h)$ is a smooth Riemannian vector bundle over $X$, and $A$ is a smooth orthogonal connection\footnote{Atiyah \cite[Section 1, Theorem 2, p. 252]{Atiyah_1985} assumes that $A$ is a unitary connection on a Hermitian vector bundle $(E,H)$, but that stronger assumption is not required by the proof.} on $V$, and $D_A:\Omega^0(W\otimes V) \to \Omega^0(W\otimes V)$ is the coupled Dirac operator as in \eqref{eq:Coupled_Dirac_operator} and $\{\lambda_k(D_A)\}_{k=1}^\infty$ denotes its sequence of eigenvalues, repeated according to their multiplicity and non-decreasing in absolute value, then
\begin{equation}
  \label{eq:Vafa_Witten_lambda_k_upper_bound}
  |\lambda_k(D_A)| \leq Ck^{1/d}, \quad\text{for all } k \geq 1.
\end{equation}
\end{thm}

We now complete the

% PF2-6-2025 Add more precise refs for unjustified assertions
\begin{proof}[Proof of Theorem \ref{mainthm:Lower_bound_spectral_gap_coupled_Dirac_operator}]
From \eqref{eq:R_A}, we obtain
\[
  \|R_A\|_{C^0(\End(W^+\otimes V))} \leq \frac{1}{4}\|R\|_{C^0(X)}
  + \frac{1}{2}\|\rho(F_{A_L})\|_{C^0(\End(W^+\otimes V))} + \|\rho(F_A)\|_{C^0(\End(W^+\otimes V))}.
\]
Thus, by substituting the upper bound \eqref{eq:rho_FA_End_W+_otimes_V_leq_constant_1+r} in our hypotheses, we see that
\begin{equation}
  \label{eq:Upper_bound_R_A}
  \|R_A\|_{C^0(\End(W^+\otimes V))} \leq z''' + r,
\end{equation}
for a constant $z''' = z'''(g,h,\|\rho(F_{A_L})\|_{C^0(\End(W^+\otimes V))}, Z) \in [1,\infty)$. We conclude from \eqref{eq:Lower_bounds_eigenvalues_square_coupled_Dirac_operator} in Corollary \ref{cor:Lott_1986} and the  inequality \eqref{eq:Upper_bound_R_A} that
\begin{equation}
  \label{eq:Lower_bounds_eigenvalues_square_coupled_Dirac_operator_4-manifold_raw}
  \lambda_k^2 \geq c_1\left(2^{-(d/2)+1}e^{-d/2}k - 1\right)^{2/d} - z''' - r,
  \quad\text{for all } k \geq 2^{(d/2)-1}e^{d/2},
\end{equation}
% PF3-20-2025 Be sure that everywhere we include dependence on g, H in these constants
where for convenience we define the positive constant
\begin{equation}
  \label{eq:Defn_c_1_lower_bound_eigenvalues_square_coupled_Dirac_operator}
  c_1 = c_1(g) := \frac{1}{2}C_1\left(4\vol(X,g)\right)^{-2/d}  \in (0,\infty).
\end{equation}
Next, we assume that $k = k(d) \geq 1$ is large enough that
\[
  2^{-(d/2)+1}e^{-d/2}k - 1 \geq \frac{1}{2}2^{-(d/2)+1}e^{-d/2}k = (2e)^{-d/2}k,
\]
and which (using $a-b \geq a/2 \iff a/2 \geq b$) holds when $(2e)^{-d/2}k \geq 1$ or equivalently
\begin{equation}
  \label{eq:k_geq_2e_power_dover2}
  k = k(d) \geq (2e)^{d/2}.
\end{equation}
Thus, we obtain the simplification
\[
  \left(2^{-(d/2)+1}e^{-d/2}k - 1\right)^{2/d}
  \geq
  \left((2e)^{-d/2}k\right)^{2/d}
  =
  (2e)^{-1} k^{2/d}, \quad\text{for all } k \geq (2e)^{d/2}.
\]
Hence, we see that the lower bound \eqref{eq:Lower_bounds_eigenvalues_square_coupled_Dirac_operator_4-manifold_raw} for $\lambda_k^2$ yields
\begin{equation}
  \label{eq:Lott_1986_corollary_twisted_Dirac_page_125_raw}
  \lambda_k^2 \geq \frac{c_1}{2e}k^{2/d}  - z''' -  r,
  \quad\text{for all } k \geq (2e)^{d/2}.
\end{equation}
For a given $r\in[1,\infty)$ in \eqref{eq:SO(3)_monopole_equations_almost_Hermitian_perturbed_intro_regular}, we require that the integer $k = k(c_1,r,z''') \geq 1$ is large enough that
\[
  \frac{c_1}{4e}k^{2/d} \geq z''' + r.
\]
We may also assume without loss of generality that $r \in [1,\infty)$ is large enough that
\begin{equation}
  \label{eq:r_geq_z'''}
  r \geq z''',
\end{equation}
and require that $k = k(c_1,r) \geq 1$ is large enough that it now obeys the stronger condition 
\[
  \frac{c_1}{4e}k^{2/d} \geq 2r,
\]
or equivalently,
\begin{equation}
  \label{eq:k_geq_constant_r_dover2}
  k = k(c_1,r) \geq \left(\frac{8e}{c_1}\right)^{d/2}r^{d/2}.
\end{equation}
Thus, we obtain (again using $a-b \geq a/2 \iff a/2 \geq b$)
\[
  \frac{c_1}{2e}k^{2/d} - z''' - r
  \geq
  \frac{c_1}{4e}k^{2/d},
  \quad\text{for all } r \geq z'''
  \text{ and } k \geq (8e/c_1)^{d/2}r^{d/2}.
\]
Therefore, the previous lower bound \eqref{eq:Lott_1986_corollary_twisted_Dirac_page_125_raw} for $\lambda_k^2$ implies that
\begin{equation}
  \label{eq:Lower_bounds_eigenvalues_square_coupled_Dirac_operator_4-manifold}
  \lambda_k^2 \geq \frac{c_1}{4e}k^{2/d},
  \quad\text{for all }  r \geq z'''
  \text{ and } k \geq \max\left\{(2e)^{d/2}, (8e/c_1)^{d/2}r^{d/2}\right\}.
\end{equation}
For convenience, we define positive constants by
\[
  c^2 = c(g)^2
  :=
  \frac{c_1}{4e}
  \quad\text{and}\quad
  c_0 = c_0(g) := (8e/c_1)^{d/2}.
\]
We may further assume without loss of generality that $r\in[1,\infty)$ in \eqref{eq:SO(3)_monopole_equations_almost_Hermitian_perturbed_intro_regular} is large enough that
\[
  c_0r^{d/2} \geq (2e)^{d/2},
\]
or equivalently, $r \geq 2e/c_0^{2/d}$. Hence, by using the preceding assumption and \eqref{eq:r_geq_z'''} and taking square roots in \eqref{eq:Lower_bounds_eigenvalues_square_coupled_Dirac_operator_4-manifold}, we obtain
\begin{equation}
  \label{eq:Lott_1986_corollary_twisted_Dirac_page_125}
  |\lambda_k| \geq ck^{1/d},
  \quad\text{for all } r \geq \max\left\{z''', 2e/c_0^{2/d}\right\}
  \text{ and } k \geq c_0r^{d/2}.
\end{equation}
From Theorem \ref{thm:Uniform_upper_bounds_eigenvalues_coupled_Dirac_operators}, there is a constant $C = C(g) \in [1,\infty)$ such that 
\[
  |\lambda_k| \leq Ck^{1/d},
  \quad\text{for all } k \geq 1.
\]
Combining the upper and lower bounds for $|\lambda_k|$ gives
\[
  ck^{1/d} \leq |\lambda_k| \leq Ck^{1/d},
   \quad\text{for all } r \geq \max\left\{z''', 2e/c_0^{2/d}\right\}
  \text{ and } k \geq c_0r^{d/2}.
\]
Hence, equation \eqref{eq:Shubin_13-17_and_18} with $p=d$ in Item \eqref{item:lambda_k_bounded_k^1overp_implies_N_lambda_bounded_lambda^p} of Proposition \ref{prop:Asymptotic_bounds_counting_functions} implies that the counting function \eqref{eq:N_lambda} obeys
\[
  \frac{1}{2}C^{-d}\lambda^d \leq N(\lambda) \leq 2c^{-d}\lambda^d,
  \quad\text{for all }  r \geq \max\left\{z''', 2e/c_0^{2/d}\right\}
  \text{ and } \lambda \geq \Lambda_0,
\]
where $\Lambda_0 \in [1,\infty)$ is a constant such that by \eqref{eq:Choice_Lambda0},
\[
  \Lambda_0 \geq \lambda_{k_1},
\]
and $k_1 \geq k_0'$ is the least integer such that
%PF3-20-2025 Should be uniformly estimatable using cubulation
$\lambda_{k_1+1} > \lambda_{k_1}$ with $k_0'$ in \eqref{eq:Choice_k0prime} given by
\[
  k_0' = \max\left\{k_0, 2(C/c)^d\right\} = \max\left\{\lceil c_0r^{d/2}\rceil, 2(C/c)^d\right\},
\]
with $k_0 := \lceil c_0r^{d/2}\rceil$, where $\lceil\,\cdot\,\rceil$ denotes the least integer function. (It follows from Feehan \cite{Feehan_weyl_law_clusters} that $k_1$ has no dependencies beyond those for the constant $k_0'$.)
We may in addition assume without loss of generality that $r\in[1,\infty)$ is large enough that $\lceil c_0r^{d/2} \rceil \geq 2(C/c)^d$ and so
\[
  k_0' = \lceil c_0r^{d/2}\rceil.
\]
To estimate $\Lambda_0$ from below, we note that
\[
  \lambda_{k_1} \geq \lambda_{k_0'} \geq c(k_0')^{1/d} \geq c\lceil c_0r^{d/2}\rceil^{1/d}
  = c(c_0r^{d/2})^{1/d} = cc_0^{1/d}r^{1/2},
\]
where in the second last equality above we assumed without loss of generality, by increasing $r$ slightly if necessary, that $c_0r^{d/2}$ is an integer, and this gives
\[
  \Lambda_0 \geq cc_0^{1/d}r^{1/2}.
\]  
The bounds \eqref{eq:Shubin_13-17_and_18} (with $p=d$) for $N(\lambda)$ given by Proposition \ref{prop:Asymptotic_bounds_counting_functions} ensure the average gap size obeys
\[
  \frac{1}{2}c^d\lambda^{1-d} \leq \frac{\lambda}{N(\lambda)} \leq 2C^d\lambda^{1-d},
  \quad\text{for all }  r \geq \max\left\{z''', 2e/c_0^{2/d}\right\}
  \text{ and } \lambda \geq \Lambda_0.
\]
By replacing $\lambda$ by $\Lambda_0 = Kr^{1/2}$ with $K := cc_0^{1/d}$ in the lower bound above we find that
\[
  \frac{\lambda}{N(\lambda)}
  \geq
  \frac{1}{2}c^d(Kr^{1/2})^{1-d}
  =
  \frac{1}{2}c^dK^{1-d}r^{(1-d)/2},
  \quad\text{for all } r \geq \max\left\{z''', 2e/c_0^{2/d}\right\}
  \text{ and } \lambda \geq  Kr^{1/2}.
\]
This yields \eqref{eq:Lower_bound_spectral_gap_coupled_Dirac_operator} with constants $C_0 := \frac{1}{2}c^d$ and  $K = cc_0^{1/d}$ and $R_0 := \max\{z''', 2e/c_0^{2/d}\}$. This completes the proof of Theorem \ref{mainthm:Lower_bound_spectral_gap_coupled_Dirac_operator}.
\end{proof}
%PF3-21-2025 It's worth checking whether we really need Vafa-Witten -- we probably don't.

\section{Lower bounds for spectral gaps for Dirac operators with vector potentials}
\label{sec:Lower_bounds_spectral_gaps_Dirac_operators_vector_potentials}
In applications, it is useful to have analogues of Corollary \ref{cor:Lott_1986} and Theorem \ref{mainthm:Lower_bound_spectral_gap_coupled_Dirac_operator}  for self-adjoint first order elliptic differential operators given by the sum $D_A+L$ of a coupled Dirac operator $D_A$ acting on $C^\infty$ sections of $W\otimes V$ as in \eqref{eq:Coupled_Dirac_operator} plus a `vector potential' $L$, a strict $C^\infty$-vector bundle endomorphism of $W\otimes V$ that is pointwise Hermitian as in Theorem \ref{thm:Kato_inequality}. We begin with the elementary application of the Min-Max Principle (see Theorem \ref{thm:Max_min_principle}).

\begin{lem}[Comparison of eigenvalues for densely defined self-adjoint linear operators with compact resolvents]
\label{lem:Teschl_corollary_4-11}  
(Compare Teschl \cite[Section 4.3, Corollary 4.11, p. 119]{Teschl_mathematical_methods_quantum_mechanics} and Tao \cite[Section 1.3.3, Equation (1.64), p. 45]{Tao_topics_random_matrix_theory}.)  
Let $\sH$ be a Hilbert space and $T \in \End(\sH)$ be a densely defined self-adjoint unbounded linear operator with compact resolvent. If $L \in \End(\sH)$ is a bounded self-adjoint linear operator, then $T + L \in \End(\sH)$ is a densely defined self-adjoint unbounded linear operator with compact resolvent and the following hold:
\begin{enumerate}
\item\label{item:Discrete_spectra_eigenvalues_T_and_T+L} The operator $T$ (respectively, $T+L$) has spectrum $\sigma(T) = \{\lambda_k(T)\}_{k=1}^\infty$ (respectively, $\sigma(T+L) = \{\lambda_k(T+L)\}_{k=1}^\infty$) comprising a discrete subset of $\CC$ consisting entirely of real eigenvalues with finite multiplicities, which we arrange in order of increasing absolute values.
  
\item\label{item:Comparison_squared_eigenvalues_T_and_T+L} The eigenvalues obey
  \[
  \frac{1}{2}\lambda_k^2(T) - \|L\|_{\End(\sH)}^2
  \leq
  \lambda_k^2(T+L)
  \leq
  2\lambda_k^2(T) + 2\|L\|_{\End(\sH)}^2,
  \quad\text{for all } k \geq 1.
  \]  
\item\label{item:Comparison_eigenvalues_T_and_T+L_bounded_below} If $T$ is bounded below, then $T+L$ is bounded below and the eigenvalues obey
\[
  \lambda_k(T) - \|L\|_{\End(\sH)} \leq \lambda_k(T+L) \leq \lambda_k(T) + \|L\|_{\End(\sH)},
  \quad\text{for all } k \geq 1.
\]  
\end{enumerate}
\end{lem}  

\begin{proof}
Consider Item \eqref{item:Discrete_spectra_eigenvalues_T_and_T+L}. By the second resolvent identity \eqref{eq:General_second_resolvent_identity}, we have
\[
  R(z,T+L) - R(z,T) = -R(z,T+L)LR(z,T), \quad\text{for all } z \in \rho(T)\cap\rho(T+L),
\]
where $R(z,T) = (T-z)^{-1}$ is the resolvent operator for $T$ with $z$ in the resolvent set $\rho(T) \subset \CC$ and similarly for $T+L$. Since both $T$ and $T+L$ are self-adjoint, then $\sigma(T) \subset \RR$ and $\sigma(T+L) \subset \RR$, so $\rho(T)\cap\rho(T+L) \subset \CC$ is non-empty. By hypothesis, $R(z,T)$ is compact for $z \in \rho(T)$, while $L$ and $R(z,T+L)$ are bounded linear operators on $\sH$, so $R(z,T+L)LR(z,T)$ is a compact operator on $\sH$ (see, for example, Conway \cite[Chapter VI, Section 3, Proposition 3.5, p. 174]{Conway_course_functional_analysis}) and thus $R(z,T+L)$ is compact for some $z \in \rho(T)\cap\rho(T+L)$ and hence all $z \in \rho(T+L)$. Thus, $T+L$ has compact resolvent. (This is a special case of Kato \cite[Chapter IV, Section 3.6, Theorem 3.17, p. 214]{Kato}.) Hence, the spectra of the operators $T$ and $T+L$ consist entirely of eigenvalues with finite multiplicities and the eigenvalues are real since each operator is self-adjoint (see Section \ref{sec:Spectral_theory_unbounded_operators}).
% PF3-26-2025 Argument appears before too. Put in appendix.

Consider Item \eqref{item:Comparison_squared_eigenvalues_T_and_T+L}. We observe that the operators $T^*T = T^2$ and $(T+L)^*(T+L) = (T+L)^2$ are nonnegative and, in particular, bounded below. Moreover, $\Dom(T+L) = \Dom(T)$. Because
\[
  \|(T+L)\psi\|_\sH \leq \|T\psi\|_\sH + \|L\psi\|_\sH, \quad\text{for all } \psi \in \Dom(T),
\]
we have
\begin{multline*}
  \|(T+L)\psi\|_\sH^2
  \leq \left(\|T\psi\|_\sH + \|L\psi\|_\sH\right)^2
  = \|T\psi\|_\sH^2 + 2\|T\psi\|_\sH\|L\psi\|_\sH + \|L\psi\|_\sH^2
  \\
  \leq 2\|T\psi\|_\sH^2 + 2\|L\psi\|_\sH^2
  \leq 2\|T\psi\|_\sH^2 + 2\|L\|_{\End(\sH)}^2\|\psi\|_\sH^2,
\end{multline*}
and so
\[
  \left\langle (T+L)^2\psi, \psi\right\rangle_\sH
  \leq
  \left\langle \left(2T^2+2\|L\|_{\End(\sH)}^2\right)\psi, \psi\right\rangle_\sH,
  \quad\text{for all } \psi \in \Dom(T).
\]
Hence, the Max-Min Principle (see Theorem \ref{thm:Max_min_principle}) and the Spectral Mapping Theorem for polynomial functions of unbounded closed operators on a Banach space (see, for example, Dunford and Schwartz \cite[Chapter VII, Section 9, Theorem 10, p. 604]{Dunford_Schwartz1}) yield
\[
  \lambda_k^2(T+L)
  =
  \lambda_k((T+L)^2)
  \leq
  \lambda_k\left(2T^2+2\|L\|_{\End(\sH)}^2\right)
  =
  2\lambda_k^2(T) + 2\|L\|_{\End(\sH)}^2.
\]
This yields the second inequality in Item \eqref{item:Comparison_squared_eigenvalues_T_and_T+L}. Similarly, because
\[
  \|T\psi\|_\sH \leq \|(T+L)\psi\|_\sH + \|L\psi\|_\sH, \quad\text{for all } \psi \in \Dom(T),
\]
we have
\begin{multline*}
  \|T\psi\|_\sH^2
  \leq \left(\|(T+L)\psi\|_\sH + \|L\psi\|_\sH\right)^2
  = \|(T+L)\psi\|_\sH^2 + 2\|(T+L)\psi\|_\sH\|L\psi\|_\sH + \|L\psi\|_\sH^2
  \\
  \leq 2\|(T+L)\psi\|_\sH^2 + 2\|L\psi\|_\sH^2
  \leq 2\|(T+L)\psi\|_\sH^2 + 2\|L\|_{\End(\sH)}^2\|\psi\|_\sH^2,
\end{multline*}
and so
\[
  \left\langle T^2\psi, \psi\right\rangle_\sH
  \leq
  \left\langle \left(2(T+L)^2 + 2\|L\|_{\End(\sH)}^2\right)\psi, \psi\right\rangle_\sH,
  \quad\text{for all } \psi \in \Dom(T).
\]
Hence, the Max-Min Principle and the Spectral Mapping Theorem yield
\[
  \lambda_k^2(T)
  =
  \lambda_k(T^2)
  \leq
  \lambda_k\left(2(T+L)^2 + 2\|L\|_{\End(\sH)}^2\right)
  =
  2\lambda_k^2(T+L) + 2\|L\|_{\End(\sH)}^2.
\]
This yields the first inequality in Item \eqref{item:Comparison_squared_eigenvalues_T_and_T+L} and completes its proof.

Consider Item \eqref{item:Comparison_eigenvalues_T_and_T+L_bounded_below}. By assumption, the operator $T$ is bounded below by a constant $c\in\RR$, that is, $\langle\psi,T\psi\rangle_\sH \geq c\|\psi\|_\sH^2$ for all $\psi \in \sH$. In this case,
\[
  \langle\psi,(T+L)\psi\rangle_\sH = \langle\psi,T\psi\rangle_\sH + \langle\psi,L\psi\rangle_\sH
  \geq \left(c - \|L\|_{\End(\sH)}\right)\|\psi\|_\sH^2,
  \quad\text{for all } \psi \in \sH,
\]
and so $T+L$ is bounded below by the constant $c - \|L\|_{\End(\sH)}$. Moreover,
\[
  T+L \geq T - \|L\|_{\End(\sH)},
\]
since $L$ is bounded below by the constant $-\|L\|_{\End(\sH)}$. Hence, the Max-Min Principle and the Spectral Mapping Theorem imply that
\[
  \lambda_k(T+L) \geq \lambda_k\left(T - \|L\|_{\End(\sH)}\right) = \lambda_k(T) - \|L\|_{\End(\sH)},
  \quad\text{for all } k \geq 1.
\]
Similarly, because $L \leq \|L\|_{\End(\sH)}$ and thus $L - \|L\|_{\End(\sH)} \leq 0$, we see that
\[
  T+L-\|L\|_{\End(\sH)} \leq T,
\]
Consequently, the Max-Min Principle and the Spectral Mapping Theorem now imply that
\[
  \lambda_k(T) \geq \lambda_k\left(T +L - \|L\|_{\End(\sH)}\right) = \lambda_k(T+L) - \|L\|_{\End(\sH)},
  \quad\text{for all } k \geq 1.
\]
This completes the proof of Item \eqref{item:Comparison_eigenvalues_T_and_T+L_bounded_below} and hence Lemma \ref{lem:Teschl_corollary_4-11}.
\end{proof}

Lemma \ref{lem:Teschl_corollary_4-11} yields the following generalization of Corollary \ref{cor:Lott_1986} and Theorem \ref{thm:Uniform_upper_bounds_eigenvalues_coupled_Dirac_operators}.

\begin{cor}[Upper and lower bounds for eigenvalues of coupled Dirac operators plus vector potentials over spin${}^c$ manifolds]
\label{cor:Lott_1986_plus_vector_potential}
Continue the notation and hypotheses of Proposition \ref{prop:Lott_5_1986}. If in addition $L \in C^\infty(\End(W\otimes V))$ is pointwise Hermitian on the fibers of $W\otimes V$, then
%PF3-27-2025 We're a bit careless about whether we mean W or W^+ below in the context of D_A?
\begin{subequations}
\begin{align}
  \label{eq:Lower_bounds_eigenvalues_square_coupled_Dirac_operator_plus_vector_potential}
  \lambda_k^2(D_A+L)
  &\geq \frac{1}{2}C_1\left(4\vol(X,g)\right)^{-2/d}\left(2^{-(d/2)+1} e^{-d/2} k - 1\right)^{2/d}
  \\
  \notag
  &\quad - \|R_A\|_{C^0(W^+\otimes V)}  - \|L\|_{C^0(\End(W\otimes V))}^2,
    \quad\text{for all } k \geq 2^{(d/2)-1} e^{d/2},
  \\
  \label{eq:Upper_bounds_eigenvalues_square_coupled_Dirac_operator_plus_vector_potential}
  \lambda_k^2(D_A+L)
  &\leq 2C^2k^{2/d} + 2\|L\|_{C^0(\End(W\otimes V))}^2, \quad\text{for all } k \geq 1,
\end{align}
\end{subequations}
where $C_1$ is the constant in Proposition \ref{prop:Lott_5_1986} and $C$ is the constant in Theorem \ref{thm:Uniform_upper_bounds_eigenvalues_coupled_Dirac_operators}.
\end{cor}

\begin{proof}
We apply Lemma \ref{lem:Teschl_corollary_4-11} with Hilbert space $\sH = L^2(W\otimes V)$, densely defined self-adjoint unbounded linear operator $T = D_A$, and bounded self-adjoint linear operator $L$. We observe that
\[
  \|L\Phi\|_{L^2(W\otimes V)} \leq  \|L\|_{C^0(\End(W\otimes V))} \|\Phi\|_{L^2(W\otimes V)},
  \quad\text{for all } \Phi \in L^2(W\otimes V).
\]  
Thus, $\|L\|_{\End(\sH)} =  \|L\|_{C^0(\End(W\otimes V))}$ and Item \eqref{item:Comparison_squared_eigenvalues_T_and_T+L} in Lemma \ref{lem:Teschl_corollary_4-11} implies that
\[
  \lambda_k^2(D_A+L) \geq \frac{1}{2}\lambda_k^2(D_A) - \|L\|_{C^0(\End(W\otimes V))}^2,
\]
and so the lower bound \eqref{eq:Lower_bounds_eigenvalues_square_coupled_Dirac_operator_plus_vector_potential} follows from Corollary \ref{cor:Lott_1986} and the preceding inequality.

Similarly, Item \eqref{item:Comparison_squared_eigenvalues_T_and_T+L} in Lemma \ref{lem:Teschl_corollary_4-11} implies that
\[
  \lambda_k^2(D_A+L) \leq 2\lambda_k^2(D_A) + 2\|L\|_{C^0(\End(W\otimes V))}^2,
\]
and so the upper bound \eqref{eq:Upper_bounds_eigenvalues_square_coupled_Dirac_operator_plus_vector_potential} follows from Theorem \ref{thm:Uniform_upper_bounds_eigenvalues_coupled_Dirac_operators} and the preceding inequality.
\end{proof}

We can now complete the

\begin{proof}[Proof of Corollary \ref{maincor:Lower_bound_spectral_gap_coupled_Dirac_operator_plus_vector_potential}]
The argument is very similar to the proof of Theorem \ref{mainthm:Lower_bound_spectral_gap_coupled_Dirac_operator}, except that we apply Corollary \ref{cor:Lott_1986_plus_vector_potential} (upper and lower bounds for $|\lambda_k(D_A+L)|$) in place of Corollary \ref{cor:Lott_1986} (lower bounds for $|\lambda_k(D_A)|$) and Theorem \ref{thm:Uniform_upper_bounds_eigenvalues_coupled_Dirac_operators} (upper bounds for $|\lambda_k(D_A)|$), and note the additional dependence of the constants on $\|L\|_{C^0(\End(W\otimes V))}$.
  
We adapt the proof of Theorem \ref{mainthm:Lower_bound_spectral_gap_coupled_Dirac_operator}. By the definition of $L$ in the hypotheses of Corollary \ref{maincor:Lower_bound_spectral_gap_coupled_Dirac_operator_plus_vector_potential}, we have 
\[
  \|L\|_{C^0(\End(W\otimes V))} \leq r\|L_1\|_{C^0(\End(W\otimes V))} + \|L_0\|_{C^0(\End(W\otimes V))}.
\]
We conclude from \eqref{eq:Lower_bounds_eigenvalues_square_coupled_Dirac_operator} in Corollary \ref{cor:Lott_1986} and the preceding upper bound that, instead of \eqref{eq:Lower_bounds_eigenvalues_square_coupled_Dirac_operator_4-manifold_raw}, 
\begin{multline}
  \label{eq:Lower_bounds_eigenvalues_square_coupled_Dirac_operator_plus_vector_potential_4-manifold_raw}
  \lambda_k^2 \geq c_1\left(2^{-(d/2)+1}e^{-d/2}k - 1\right)^{2/d} - z''' - r
  \\
  - \left(r\|L_1\|_{C^0(\End(W\otimes V))} + \|L_0\|_{C^0(\End(W\otimes V))}\right)^2,
  \quad\text{for all } k \geq 2^{(d/2)-1}e^{d/2},
\end{multline}
for constants $z'''$ as in \eqref{eq:Upper_bound_R_A} and  $c_1$ as in \eqref{eq:Defn_c_1_lower_bound_eigenvalues_square_coupled_Dirac_operator}. As in the corresponding stage in the proof of Theorem \ref{mainthm:Lower_bound_spectral_gap_coupled_Dirac_operator}, we assume that $k = k(d) \geq 1$ obeys
\begin{equation}
  \label{eq:Initial_lower_bound_k}
  k = k(d) \geq (2e)^{d/2},
\end{equation}
and obtain the simplification
\begin{align*}
  \lambda_k^2
  &\geq \frac{c_1}{2e}k^{2/d} - z''' - r
  - \left(r\|L_1\|_{C^0(\End(W\otimes V))} + \|L_0\|_{C^0(\End(W\otimes V))}\right)^2
  \\
  &= \frac{c_1}{2e}k^{2/d} - z''' - r - r^2\|L_1\|_{C^0(\End(W\otimes V))}^2
     - 2r\|L_1\|_{C^0(\End(W\otimes V))}\|L_0\|_{C^0(\End(W\otimes V))}
  \\
  &\quad
    - \|L_0\|_{C^0(\End(W\otimes V))}^2
  \\
  &\geq \frac{c_1}{2e}k^{2/d} - z''' - r
  - 2r^2\|L_1\|_{C^0(\End(W\otimes V))}^2 - 2\|L_0\|_{C^0(\End(W\otimes V))}^2,  
  \\
  &\qquad \text{for all } k = k(d) \geq (2e)^{d/2}.
\end{align*}
Thus, for the following choices of positive constants,
\begin{align*}
  z_0 &= z_0\left(g,h,\|\rho(F_{A_L}\|_{C^0(\End(W^+\otimes V))},\|L_0\|_{C^0(\End(W\otimes V))}\right)
        := z''' + 2\|L_0\|_{C^0(\End(W\otimes V))}^2,
  \\
  z_1 &= z_1\left(\|L_1\|_{C^0(\End(W\otimes V))}\right) := 1 + 2\|L_1\|_{C^0(\End(W\otimes V))}^2,
\end{align*}
and using our assumption $r \geq 1$ to write $r \leq r^2$, we see that
\[
  \lambda_k^2 \geq \frac{c_1}{2e}k^{2/d} - z_0 - z_1r^2,
  \quad\text{for all } k \geq (2e)^{d/2}.
\]
For a given $r \in [1,\infty)$ in \eqref{eq:SO(3)_monopole_equations_almost_Hermitian_perturbed_intro_regular}, we further require that $k = k(c_1,r,z_0,z_1) \geq 1$ is large enough that
\[
  \frac{c_1}{4e}k^{2/d} \geq z_0 + z_1r^2.
\]
We may assume without loss of generality that $r \in [1,\infty)$ is large enough that
\begin{equation}
  \label{eq:First_lower_bound_r}
  z_1r^2 \geq z_0,
\end{equation}
and require that $k = k(c_1,r,z_1) \geq 1$ is large enough that it now obeys the stronger condition 
\[
  \frac{c_1}{4e}k^{2/d} \geq 2z_1r^2,
\]
or equivalently,
\begin{equation}
  \label{eq:Second_lower_bound_k}
  k = k(c_1,r,z_1) \geq \left(\frac{8ez_1}{c_1}\right)^{d/2}r^d.
\end{equation}
Thus, we obtain (again using $a-b \geq a/2 \iff a/2 \geq b$)
\[
  \frac{c_1}{2e}k^{2/d} - z_0 - z_1r^2
  \geq
  \frac{c_1}{4e}k^{2/d},
  \quad\text{for all } r \geq (z_0/z_1)^{1/2} \text{ and } k \geq (8ez_1/c_1)^{d/2}r^d.
\]
Therefore, the previous lower bound for $\lambda_k^2$ implies that
\begin{equation}
  \label{eq:Lower_bounds_eigenvalues_square_coupled_Dirac_operator_vector_potential_4-manifold}
  \lambda_k^2 \geq \frac{c_1}{4e}k^{2/d},
  \quad\text{for all } r \geq (z_0/z_1)^{1/2} \text{ and }
  k \geq \max\left\{(8ez_1/c_1)^{d/2}r^d, (2e)^{d/2} \right\}.
\end{equation}
For convenience, we define positive constants by
\[
  c_3^2 = c_3(g)^2
  :=
  \frac{c_1}{4e}
  \quad\text{and}\quad
  c_2 = c_2(g,z_1) := (8ez_1/c_1)^{d/2}.
\]
We may further assume without loss of generality that $r\in[1,\infty)$ in \eqref{eq:SO(3)_monopole_equations_almost_Hermitian_perturbed_intro_regular} is large enough that
\[
  c_2r^d \geq (2e)^{d/2},
\]
or equivalently,
\[
  r \geq (2e)^{1/2}/c_2^d,
\]
and so by using this assumption and taking square roots in \eqref{eq:Lower_bounds_eigenvalues_square_coupled_Dirac_operator_vector_potential_4-manifold}, we obtain
\begin{equation}
  \label{eq:Lott_1986_corollary_twisted_Dirac_page_125_vector_potential}
  |\lambda_k| \geq c_3k^{1/d},
  \quad\text{for all } r \geq \max\left\{(z_0/z_1)^{1/2}, (2e)^{1/2}/c_2^d\right\}
  \text{ and } k \geq c_2r^d.
\end{equation}
From Theorem \ref{thm:Uniform_upper_bounds_eigenvalues_coupled_Dirac_operators}, there is a constant $C = C(g) \in [1,\infty)$ such that 
\[
  |\lambda_k| \leq Ck^{1/d},
  \quad\text{for all } k \geq 1.
\]
Combining the upper and lower bounds for $|\lambda_k|$ gives
\[
  c_3k^{1/d} \leq |\lambda_k| \leq Ck^{1/d},
  \quad\text{for all } r \geq \max\left\{(z_0/z_1)^{1/2}, (2e)^{1/2}/c_2^d\right\}
  \text{ and } k \geq c_2r^d.
\]
Hence, equation \eqref{eq:Shubin_13-17_and_18} with $p=d$ in Item \eqref{item:lambda_k_bounded_k^1overp_implies_N_lambda_bounded_lambda^p} of Proposition \ref{prop:Asymptotic_bounds_counting_functions} implies that the counting function \eqref{eq:N_lambda} obeys
\[
  \frac{1}{2}C^{-d}\lambda^d \leq N(\lambda) \leq 2c_3^{-d}\lambda^d,
  \quad\text{for all } r \geq \max\left\{(z_0/z_1)^{1/2}, (2e)^{1/2}/c_2^d\right\}
  \text{ and } \lambda \geq \Lambda_0,
\]
where $\Lambda_0 \in [1,\infty)$ is a constant such that by \eqref{eq:Choice_Lambda0},
\[
  \Lambda_0 \geq |\lambda_{k_1}|,
\]
and $k_1 \geq k_0'$ is the least integer such that
%PF3-20-2025 Should be uniformly estimatable using cubulation
$|\lambda_{k_1+1}| > |\lambda_{k_1}|$ with $k_0'$ in \eqref{eq:Choice_k0prime} given by
\[
  k_0' = \max\left\{k_0, 2(C/c_3)^d\right\} = \max\left\{\lceil c_2r^d\rceil, 2(C/c_3)^d\right\},
\]
with $k_0 := \lceil c_2r^d\rceil$, where $\lceil\,\cdot\,\rceil$ denotes the least integer function. (Again, it follows from Feehan \cite{Feehan_weyl_law_clusters} that $k_1$ has no dependencies beyond those for the constant $k_0'$.)
We may in addition assume without loss of generality that $r\in[1,\infty)$ is large enough that $\lceil c_2r^d \rceil \geq 2(C/c_3)^d$ and so
\[
  k_0' = \lceil c_2r^d\rceil.
\]
To estimate $\Lambda_0$ from below, we note that
\[
  |\lambda_{k_1}| \geq |\lambda_{k_0'}| \geq c_3(k_0')^{1/d} \geq c_3\lceil c_2r^d\rceil^{1/d}
  = c_3(c_2r^d)^{1/d} = c_3c_2^{1/d}r,
\]
where in the second last equality above we assumed without loss of generality, by increasing $r$ slightly if necessary, that $c_2r^d$ is an integer, and this gives
\[
  \Lambda_0 \geq c_3c_2^{1/d}r.
\]  
The bounds \eqref{eq:Shubin_13-17_and_18} (with $p=d$) for $N(\lambda)$ given by Proposition \ref{prop:Asymptotic_bounds_counting_functions} ensure the average gap size obeys
\[
  \frac{1}{2}c_3^d\lambda^{1-d} \leq \frac{\lambda}{N(\lambda)} \leq 2C^d\lambda^{1-d},
  \quad\text{for all } r \geq \max\left\{(z_0/z_1)^{1/2}, (2e)^{1/2}/c_2^d\right\}
  \text{ and } \lambda \geq \Lambda_0.
\]
By replacing $\lambda$ by $\Lambda_0 = K_1r$ with $K_1 := c_3c_2^{1/d}$, in the preceding lower bound for $\lambda/N(\lambda)$ we find that
\[
  \frac{\lambda}{N(\lambda)}
  \geq
  \frac{1}{2}c_3^d(K_1r)^{1-d}
  =
  \frac{1}{2}c_3^dK_1^{1-d}r^{1-d},
  \quad\text{for all } r \geq \max\left\{(z_0/z_1)^{1/2}, (2e)^{1/2}/c_2^d\right\}
  \text{ and } \lambda \geq K_1r.
\]
This yields \eqref{eq:Lower_bound_spectral_gap_coupled_Dirac_operator_plus_vector_potential} with constants $C_2 := \frac{1}{2}c_3^d$ and $K_1 = c_3c_2^{1/d}$ and $R_1 := \max\{(z_0/z_1)^{1/2}, (2e)^{1/2}/c_2^d\}$. This completes the proof of Corollary \ref{maincor:Lower_bound_spectral_gap_coupled_Dirac_operator_plus_vector_potential}.
\end{proof}

Finally, we complete the

\begin{proof}[Proof of Corollary \ref{maincor:Lower_bound_spectral_gap_coupled_Dirac_operator_plus_vector_potential_four-manifold}]
By hypothesis, $(A,\varphi,\psi)$ is a smooth solution to the non-Abelian monopole equations \eqref{eq:SO(3)_monopole_equations_almost_Hermitian_perturbed_intro_regular} with a regularized Taubes perturbation, so the inequality \eqref{eq:Taubes_1996_SW_to_Gr_eq_2-12_and_2-13_regular} in Corollary \ref{maincor:Taubes_1996_SW_to_Gr_eq_2-12_and_2-13_regular} implies that for all $r \geq z$, where $z = z(g,A_d,\omega) \in [1,\infty)$ is the constant in Proposition \ref{mainprop:Taubes_1996_SW_to_Gr_2-1},
%PF4-3-2025 Make dependencies on C^0 norm of F_{A_d}, etc.
\begin{equation}
  \label{eq:Taubes_1996_SW_to_Gr_eq_2-12_and_2-13_regular_Lott_application_E}
  \|\rho_{\can}(F_A^+)_0\|_{C^0(\End(W_{\can}^+\otimes E))}
  \leq
  \frac{(1+2\sqrt{2})}{2\sqrt{2}}\frac{5}{3} + \frac{r}{2}.
\end{equation}
Consequently, if $A^\ad$ denotes the connection on $\fsl(E)$ induced by $A$ on $E$, then
%PF4-3-2025 Check numerical factors!
\begin{equation}
  \label{eq:Taubes_1996_SW_to_Gr_eq_2-12_and_2-13_regular_Lott_application_slE}
  \|\rho_{\can}(F_{A^{\ad}}^+)\|_{C^0(\End(W_{\can}^+\otimes \fsl(E)))}
  \leq
  Z_0 + \frac{r}{2},
\end{equation}
where $Z_0$ is a universal numerical constant. From \eqref{eq:R_A} and the choice $V = E\otimes\fsl(E)$ and the fact that $\dim X = 4$, we have
\begin{multline*}
  \|R_A\|_{C^0(\End(W_{\can}^+\otimes V))}
  \leq \frac{1}{4}\|R\|_{C^0(X)}
    + \frac{1}{2}\|\rho(F_{A_L}^+)\|_{C^0(\End(W^+\otimes L))} + \|\rho(F_A^+)\|_{C^0(\End(W^+\otimes E))}
  \\
  + \|\rho(F_{A^\ad}^+)\|_{C^0(\End(W^+\otimes\,\fsl(E)))}.
\end{multline*}
Because %PF3-19-2025 Ref
\[
  F_A = (F_A)_0 + \frac{1}{r_E}F_{A_d}
\]
and $r_E = \rank_\CC(E) = 2$, we obtain
\begin{multline*}
  \|R_A\|_{C^0(\End(W_{\can}^+\otimes V))}
  \leq \frac{1}{4}\|R\|_{C^0(X)}
    + \frac{1}{2}\|\rho(F_{A_L}^+)\|_{C^0(\End(W^+\otimes L))} + \|\rho(F_A^+)_0\|_{C^0(\End(W^+\otimes E))}
    \\
    + \frac{1}{2}\|\rho(F_{A_d}^+)\|_{C^0(\End(W^+\otimes\,\det E))}
    + \frac{1}{2}\|\rho(F_{A^\ad}^+)\|_{C^0(\End(W^+\otimes\,\fsl(E)))}.
\end{multline*}
Thus, by substituting the upper bounds \eqref{eq:Taubes_1996_SW_to_Gr_eq_2-12_and_2-13_regular_Lott_application_E} and \eqref{eq:Taubes_1996_SW_to_Gr_eq_2-12_and_2-13_regular_Lott_application_slE}, we find that
\begin{equation}
  \label{eq:Upper_bound_R_A_E_slE}
  \|R_A\|_{C^0(\End(W_{\can}^+\otimes V))} \leq z''' + r,
\end{equation}
for a constant $z''' = z'''(g,H,\|F_{A_d}^+\|_{C^0(X)}) \in [1,\infty)$. The role of the upper bound \eqref{eq:Upper_bound_R_A} in the proof of Theorem \ref{mainthm:Lower_bound_spectral_gap_coupled_Dirac_operator}  and hence in the proof of Corollary \ref{maincor:Lower_bound_spectral_gap_coupled_Dirac_operator_plus_vector_potential} is replaced by that of \eqref{eq:Upper_bound_R_A_E_slE}. The remainder of the proofs is exactly as before except that we now assume that $d=4$. Hence, the lower bound \eqref{eq:Lower_bound_spectral_gap_coupled_Dirac_operator_plus_vector_potential_four-manifold} that we seek follows from the lower bound \eqref{eq:Lower_bound_spectral_gap_coupled_Dirac_operator_plus_vector_potential} in Corollary \ref{maincor:Lower_bound_spectral_gap_coupled_Dirac_operator_plus_vector_potential} with $d=4$. This completes the proof of Corollary \ref{maincor:Lower_bound_spectral_gap_coupled_Dirac_operator_plus_vector_potential_four-manifold}.
\end{proof}

\chapter[Analogues for non-Abelian monopoles of Taubes' pointwise estimates]{Analogues for non-Abelian monopoles of Taubes' differential inequalities and pointwise estimates for Seiberg--Witten monopole sections}
\label{chap:Analogues_non-Abelian_monopoles_Taubes_estimates_Seiberg-Witten_monopole_sections}
In this chapter, we first prove Proposition \ref{mainprop:Taubes_1996_SW_to_Gr_2-1} and Theorem \ref{mainthm:Taubes_1996_SW_to_Gr_2-3}, which give estimates for the squared pointwise norms of sections $\alpha$ of $E$ and $\beta$ of $\Lambda^{0,2}(E)$ when $(A,\varphi,\psi)$ with $(\varphi,\psi) = r^{1/2}(\alpha,\beta)$ is a solution to the system \eqref{eq:SO(3)_monopole_equations_almost_Hermitian_perturbed_intro} of non-Abelian monopole equations with a singular Taubes perturbation. We then describe the modifications to those proofs required to establish Proposition \ref{mainprop:Taubes_1996_SW_to_Gr_2-1_regular} and Theorem \ref{mainthm:Taubes_1996_SW_to_Gr_2-3_regular}, which give estimates for the squared pointwise norms of sections $\alpha$ of $E$ and $\beta$ of $\Lambda^{0,2}(E)$ when $(A,\varphi,\psi)$ with $(\varphi,\psi) = r^{1/2}(\alpha,\beta)$ is a solution to the system \eqref{eq:SO(3)_monopole_equations_almost_Hermitian_perturbed_intro_regular} of non-Abelian monopole equations with a regularized Taubes perturbation. We adopt this two stage approach to the proofs of Proposition \ref{mainprop:Taubes_1996_SW_to_Gr_2-1_regular} and Theorem \ref{mainthm:Taubes_1996_SW_to_Gr_2-3_regular} since the arguments required for case of singular Taubes perturbations are more transparent than in the case of regularized Taubes perturbations and the modifications required for the case of regular Taubes perturbations are easily described.  

Section \ref{sec:Pointwise_estimate_holomorphic_component_coupled_spinor} constrains our proof of Proposition \ref{mainprop:Taubes_1996_SW_to_Gr_2-1}, giving estimates for $|\alpha|_E^2$ in the case of singular and regularized Taubes perturbations. In Section \ref{sec:Bochner-Weitzenbock_formulae_holomorphic_and_anti-holomorphic_spinors}, we develop the Bochner--Weitzenb\"ock identities that we shall need for $(0,1)$-connections on Hermitian vector bundles over almost Hermitian manifolds. Sections \ref{sec:Differential_inequality_squared_pointwise_norm_holomorphic_spinor}, \ref{sec:Differential_inequality_squared_pointwise_norm_anti-holomorphic_spinor}, and \ref{sec:Differential_inequality_Taubes_combination_squared_pointwise_norms_spinor_components} establish differential inequalities for $|\alpha|_E^2$ and $|\beta|_{\Lambda^{0,2}(E)}^2$ and an affine linear combination of those pointwise squared norms, respectively, in the case of singular Taubes perturbations. Sections \ref{sec:Differential_inequality_squared_pointwise_norm_anti-holomorphic_spinor_regular} and \ref{sec:Differential_inequality_Taubes_combination_squared_pointwise_norms_spinor_components_regular} contain the modifications to Sections \ref{sec:Differential_inequality_squared_pointwise_norm_holomorphic_spinor}, \ref{sec:Differential_inequality_squared_pointwise_norm_anti-holomorphic_spinor}, and \ref{sec:Differential_inequality_Taubes_combination_squared_pointwise_norms_spinor_components} that are required to establish differential inequalities for $|\alpha|_E^2$ and $|\beta|_{\Lambda^{0,2}(E)}^2$ and an affine linear combination of those pointwise squared norms, respectively, in the case of regularized Taubes perturbations. In Section \ref{sec:Pointwise_estimate_squared_pointwise_norm_anti-holomorphic_spinor_regular}, we finally derive our pointwise estimate for $|\beta|_{\Lambda^{0,2}(E)}^2$ in the case of regularized Taubes perturbations. In Section \ref{sec:Perturbation_non-Abelian_monopole_equations_split_pairs}, we describe the structure of the non-Abelian monopole equations with singular and regularized Taubes perturbations when $(A,\varphi,\psi)$ is split with respect to a decomposition $E = L_1 \oplus L_2$ as an orthogonal direct sum of Hermitian line bundles. In Section \ref{sec:Maximum_principles}, we describe maximum principles for linear second order elliptic partial differential inequalities on Riemannian manifolds.

\section[Regularity for solutions to the non-Abelian monopole equations]{Regularity for solutions to the non-Abelian monopole equations with a regularized Taubes perturbation}
\label{sec:Feehan_Leness_1998jdg_3-7_regularized_Taubes-perturbed_W1p}
The following result is proved in Feehan and Leness \cite[Proposition 3.7, p. 323]{FL1} for the case of $W^{k,2}$ solutions to the system \eqref{eq:SO(3)_monopole_equations_almost_Hermitian_intro} of unperturbed non-Abelian monopole equations modulo $W^{k+1,2}$ gauge transformations for integers $k\geq 2$ and standard methods yield the case of $W^{1,p}$ solutions to the regularized Taubes-perturbed non-Abelian monopole equations modulo $W^{2,p}$ gauge transformations for exponents $p\in(2,\infty)$.

\begin{thm}[Regularity for solutions to the non-Abelian monopole equations  with a regularized Taubes perturbation]
\label{thm:Feehan_Leness_1998jdg_3-7_regularized_Taubes-perturbed_W1p}
Let $(X,g,J,\omega)$ be an almost Hermitian four-manifold, $(E,H)$ be a smooth, Hermitian vector bundle over $X$, and $A_d$ be a smooth, unitary connection on $\det E$. Let $p \in (2,\infty)$ be a constant. If $(A,\varphi,\psi) \in \sA^{1,p}(E,A_d,H) \times W^{1,p}(E\oplus\Lambda^{0,2}(E))$ is a $W^{1,p}$ solution to the system \eqref{eq:SO(3)_monopole_equations_almost_Hermitian_perturbed_intro_regular} of non-Abelian monopole equations with a regularized Taubes perturbation for constants $r \in [0,\infty)$ and $\delta \in (0,\infty)$, then there is a gauge transformation $u \in W^{2,p}(\SU(E))$ such that $u(A,\varphi,\psi)$ is a $C^\infty$ solution.
\end{thm}  

\section[Estimate for the squared pointwise norms of sections of $E$]{Estimate for the squared pointwise norms of sections of  $E$ with singular and a regularized Taubes perturbation}
\label{sec:Pointwise_estimate_holomorphic_component_coupled_spinor}
% PF7-25-2024 Try to also derive the apriori pointwise estimate from the BW formula for \varphi alone in next subsection.
In Sections \ref{subsec:Pointwise_estimate_holomorphic_component_coupled_spinor_singular} and Sections \ref{subsec:Pointwise_estimate_holomorphic_component_coupled_spinor_regular}, respectively, we prove Proposition \ref{mainprop:Taubes_1996_SW_to_Gr_2-1} and \ref{mainprop:Taubes_1996_SW_to_Gr_2-1_regular}, which give estimates for $|\alpha|_E^2$ (actually, $|\alpha|_E^2 + |\beta|_{\Lambda^{0,2}(E)}^2$) when $(A,\varphi,\psi)$ with $(\varphi,\psi) = r^{1/2}(\alpha,\beta)$ is a solution to the system \eqref{eq:SO(3)_monopole_equations_almost_Hermitian_perturbed_intro} of non-Abelian monopole equations with a singular Taubes perturbation or the system \eqref{eq:SO(3)_monopole_equations_almost_Hermitian_perturbed_intro_regular} of non-Abelian monopole equations with a regularized Taubes perturbation.

\subsection{Estimate for the squared pointwise norms of  $E$ using a singular Taubes perturbation}
\label{subsec:Pointwise_estimate_holomorphic_component_coupled_spinor_singular}
In this subsection, we give two proofs of Proposition \ref{mainprop:Taubes_1996_SW_to_Gr_2-1}: the first proof relies on the additional assumption that $\Phi$ is $C^\infty$ while the second proof removes that assumption.

\begin{proof}[First proof of Proposition \ref{mainprop:Taubes_1996_SW_to_Gr_2-1}]
For this argument, we make the additional assumption that $\Phi$ is $C^\infty$. By analogy with Taubes \cite[Section 2 (a), Proof of Proposition 2.1, p. 85]{TauSWGromov}, we assume that $\Phi = \Omega^0(W^+\otimes E)$ obeys $D_A\Phi= 0 \in \Omega^0(W^-\otimes E)$ by \eqref{eq:SO(3)_monopole_equations_Dirac} and so
\[
  D_A^*D_A\Phi = 0 \in \Omega^0(W^+\otimes E).
\]
By substituting the Bochner--Weitzenb\"ock identity for $D_A^*D_A$ (see Equation \eqref{eq:Bochner-Weitzenbock_formula_coupled_Dirac_operator} in
%TL8-2-2024: Added following forward reference
the forthcoming Lemma \ref{lem:Bochner-Weitzenbock_identity_coupled_Dirac_operator}),
% PF6-28-2024 Cite other sources for BW formula too, especially E of any rank
%[TODO - Tom, please recheck/correct formula below. Missing a term?]
% TL8-2-2024: Term inserted
% PF8-16-2024 factor of 1/2 and \Phi were missing in inserted term
% PF8-16-2024 You added a term, but then didn't correct the remainder of this subsection?
the preceding identity becomes
\begin{equation}
  \label{eq:Taubes_1996_SW_to_Gr_2-1}
  \nabla_{A_W\otimes A}^*\nabla_{A_W\otimes A}\Phi + \frac{R}{4}\Phi 
  + \frac{1}{2}\rho(F_{A_L}^+)\Phi
  + \rho(F_A^+)\Phi
  = 0,
\end{equation}
where $R$ is the scalar curvature for the Riemannian metric $g$ on $X$
and $A_L$ is the connection induced by $A_W$ on $\det W^+$. In our application, $(\rho,W) = (\rho_\can,W_\can)$.

We take the pointwise $W^+\otimes E$ inner product $\langle\cdot,\cdot\rangle$ of the preceding equality with $\Phi$ to give
%PF6-28-2024 Real part needed?
\begin{equation}
  \label{eq:Taubes_1996_SW_to_Gr_2-1_W+E-inner_product_Phi}
  \langle\nabla_{A_W\otimes A}^*\nabla_{A_W\otimes A}\Phi,\Phi\rangle + \frac{R}{4}|\Phi|^2
  + \frac{1}{2}\langle\rho(F_{A_L}^+)\Phi,\Phi\rangle + \langle\rho(F_A^+)\Phi,\Phi\rangle = 0
  \quad\text{on } X.
\end{equation}
By Freed and Uhlenbeck \cite[Chapter 6, Equation (6.18), p. 91]{FU} or Salamon \cite[Section 7.2, Equation (7.8), p. 231]{SalamonSWBook}, we have the \emph{second-order Kato equality},
\begin{equation}
  \label{eq:Second-order_Kato_equality}
  \Real\,\langle\nabla_{A_W\otimes A}^*\nabla_{A_W\otimes A}\Phi,\Phi\rangle
  = \frac{1}{2}\Delta_g|\Phi|^2 + |\nabla_{A_W\otimes A}\Phi|^2 \quad\text{on } X,
\end{equation}
where $\Delta_g = d^{*_g}d$ on $\Omega^0(X;\RR)$. Substituting the equality \eqref{eq:Second-order_Kato_equality} into the identity \eqref{eq:Taubes_1996_SW_to_Gr_2-1_W+E-inner_product_Phi} yields
%PF8-16-2024 Smooth here means W^{1,p}, not really smooth; identities are a.e. on X throughout.
\begin{equation}
  \label{eq:Taubes_1996_SW_to_Gr_2-1_W+E-inner_product_Phi_2nd-order_Kato}
  \frac{1}{2}\Delta_g|\Phi|^2 + |\nabla_{A_W\otimes A}\Phi|^2 + \frac{R}{4}|\Phi|^2
  + \frac{1}{2}\Real\,\langle\rho(F_{A_L}^+)\Phi,\Phi\rangle + \Real\,\langle\rho(F_A^+)\Phi,\Phi\rangle = 0
  \quad\text{on } X.
\end{equation}
By hypothesis, $(A,\Phi)$ is a smooth solution to the perturbed non-Abelian monopole equations given by \eqref{eq:SO(3)_monopole_equations_Dirac} and the forthcoming \eqref{eq:SO(3)_monopole_equations_curvature_Taubes_perturbation}, that is
%PF8-16-2024 Maybe change perturbation parameter "r" to "t", avoiding conflict with complex rank of E
\[
  \rho(F_A^+)_0 = (\Phi\otimes\Phi^*)_{00} - \frac{r}{4}\wp\jmath.
\]
Recall that when $E$ has complex rank $r_E$,
%PF8-16-2024 Prefer to make complex rank r (see notation conflict above)
%PF6-28-2024 Check below
\[
  F_A^+ = (F_A^+)_0 + \frac{1}{r_E}F_{A_d}^+\otimes\id_E,
\]
so that
\[
  \rho(F_A^+) = \rho(F_A^+)_0 + \frac{1}{r_E}\rho(F_{A_d}^+)\otimes\id_E.
\]
We apply equation \eqref{eq:SO(3)_monopole_equations_(1,1)_curvature_perturbed_omega_intro} for the component $(F_A^\omega)_0$ of $(F_A^+)_0 $ and equations \eqref{eq:SO(3)_monopole_equations_(0,2)_curvature_perturbed_intro} and \eqref{eq:SO(3)_monopole_equations_(2,0)_curvature_intro} for the components $(F_A^{0,2})_0$ and $(F_A^{2,0})_0$, respectively, of $(F_A^+)_0$. Thus:
%PF3-17-2024 Duplicate of later expression?
\begin{equation}
  \label{eq:(F_A^+)_0}
  (F_A^+)_0 = \frac{i}{4}(\varphi\otimes\varphi^*)_0\omega - \frac{i}{4}\star(\psi\otimes\psi^*)_0\omega
  - \frac{ir}{8}\wp\,\omega
  - \frac{1}{2}(\varphi\otimes\psi^*)_0 + \frac{1}{2}(\psi\otimes\varphi^*)_0
  \in \Omega^+(\su(E)).
\end{equation}
% PF8-16-2024 Below duplicates later discussion in BW Dirac subsection
We now apply Feehan and Leness 
%TL12-4-2025: Updated and removed page references
\cite[Lemma 8.3.4 and Corollary 8.3.5]{Feehan_Leness_introduction_virtual_morse_theory_so3_monopoles} to compute $\rho_\can(F_A^+)_0$, noting that $\psi \in \Omega^{0,2}(E)$ and $\psi^* = \langle\cdot,\psi\rangle_E \in \Omega^{2,0}(E^*)$, so if $\psi = \beta\otimes s$ with $\beta \in \Omega^{0,2}(X)$ and $\bar\beta \in \Omega^{2,0}(X)$ and $s \in \Omega^0(E)$, then $\psi^* = \bar\beta\otimes\langle\cdot,s\rangle_E = \bar\beta\otimes s^*$. From 
%TL12-4-2025: Updated and removed page references
\cite[Lemma 8.3.4  and Corollary 8.3.5]{Feehan_Leness_introduction_virtual_morse_theory_so3_monopoles}, we have
\begin{align*}
  \rho_\can(i\omega)(\sigma,\tau)
  &= 2(\sigma,-\tau) \in \Omega^0(E)\oplus\Omega^{0,2}(E),
  \\
  \rho_\can(\bar\beta)(\sigma,\tau)
  &=
     \left(-2\langle\tau,\beta\rangle_{\Lambda^{0,2}(X)}, 0\right) \in \Omega^0(E)\oplus\Omega^{0,2}(E),
  \\
  \rho_\can(\beta)(\sigma,\tau)
  &=
    \left(0, 2\sigma\beta\right) \in \Omega^0(E)\oplus\Omega^{0,2}(E),
\end{align*}
for all $(\sigma,\tau) \in \Omega^0(E) \oplus \Omega^{0,2}(E)$. We next observe that
\begin{align*}
  \varphi\otimes\psi^* &= \varphi\otimes\bar\beta\otimes s^* \in \Omega^{2,0}(\gl(E)),
  \\
  \psi\otimes\varphi^* &= \beta\otimes s \otimes \varphi^* \in \Omega^{0,2}(\gl(E)),
\end{align*}
and so
\begin{align*}
  (\varphi\otimes\psi^*)_0 &= (\varphi\otimes s^*)_0\otimes\bar\beta \in \Omega^{2,0}(\fsl(E)),
  \\
  (\psi\otimes\varphi^*)_0 &= \beta\otimes (s \otimes \varphi^*)_0 \in \Omega^{0,2}(\fsl(E)).
\end{align*}
Therefore,
\[
  \rho(\psi\otimes\varphi^*)_0 = 2\beta\otimes (s \otimes \varphi^*)_0
  = 2(\psi \otimes \varphi^*)_0 \in \Omega^{0,2}(\fsl(E)),
\]
while, noting the conjugate-linear Riesz isomorphism,
\[
  \Omega^{0,2}(X) \ni \beta \mapsto \langle\cdot,\beta\rangle_{\Lambda^{0,2}(X)}
  = \bar\beta \in \Omega^{2,0}(X),
\]  
and recalling that $s^* = \langle\cdot,s\rangle_E$, we obtain
\begin{multline*}
  \rho(\varphi\otimes\psi^*)_0
  =
  (\varphi\otimes s^*)_0\otimes\rho(\bar\beta)
  =
  -2(\varphi\otimes s^*)_0\otimes\langle\cdot,\beta\rangle_{\Lambda^{0,2}(X)}
  \\
  =
  -2(\varphi\otimes s^*)_0\otimes\bar\beta
  =
  -2(\varphi\otimes \langle\cdot,s\rangle_E\otimes\bar\beta)_0
  =
  -2(\varphi\otimes \langle\cdot,s\otimes\beta\rangle_E)_0
  \\
  =
  -2(\varphi\otimes \langle\cdot,\psi\rangle_E)_0
  =
  -2(\varphi\otimes \psi^*)_0
  \in
  \Omega^{2,0}(\fsl(E)).
\end{multline*}
Consequently, writing
\begin{equation}
  \label{eq:Defn_jmath}
  \jmath(\sigma,\tau) := (\sigma,-\tau),
  \quad\text{for all } (\sigma,\tau) \in  \Omega^0(E)\otimes\Omega^{0,2}(E),
\end{equation}
we have
%PF3-17-2024 Duplicate of later expression? Codomain?
\begin{equation}
  \label{eq:rho_(F_A^+)_0}
  \rho(F_A^+)_0 = \frac{1}{2}(\varphi\otimes\varphi^*)_0\jmath - \frac{1}{2}\star(\psi\otimes\psi^*)_0\jmath
  - \frac{r}{4}\wp\jmath
  - (\varphi\otimes\psi^*)_0 + (\psi\otimes\varphi^*)_0.
\end{equation}
In particular, we see that the unperturbed non-Abelian monopole equation \eqref{eq:SO(3)_monopole_equations_curvature} is replaced by
%PF2-15-2025 Recheck
\begin{equation}
  \label{eq:SO(3)_monopole_equations_curvature_Taubes_perturbation}
  \rho(F_A^+)_0 = (\Phi\otimes\Phi^*)_{00} - \frac{r}{4}\wp\jmath.
\end{equation}
Consequently,
\[
  \rho(F_A^+) = (\Phi\otimes\Phi^*)_{00} + \frac{1}{r_E}\rho(F_{A_d}^+)\otimes\id_E - \frac{r}{4}\wp\jmath.
\]
Therefore,
\[
  \langle\rho(F_A^+)\Phi,\Phi\rangle
  = \langle(\Phi\otimes\Phi^*)_{00}\Phi,\Phi\rangle + \frac{1}{r_E}\langle\rho(F_{A_d}^+)\Phi,\Phi\rangle
  - \frac{r}{4}\Real\,\langle\wp\jmath\Phi,\Phi\rangle,
\]
and substitution into the identity \eqref{eq:Taubes_1996_SW_to_Gr_2-1_W+E-inner_product_Phi_2nd-order_Kato} yields
\begin{multline}
  \label{eq:Taubes_1996_SW_to_Gr_2-1_W+E-inner_product_Phi_2nd-order_Kato_PU(2)_monopole}
  \frac{1}{2}\Delta_g|\Phi|^2 + |\nabla_{A_W\otimes A}\Phi|^2 + \frac{R}{4}|\Phi|^2
  + \Real\,\langle(\Phi\otimes\Phi^*)_{00}\Phi,\Phi\rangle
  \\
  + \frac{1}{2}\Real\,\langle\rho(F_{A_L}^+)\Phi,\Phi\rangle
  + \frac{1}{r_E}\Real\,\langle\rho(F_{A_d}^+)\Phi,\Phi\rangle
  - \frac{r}{4}\Real\,\langle\wp\jmath\Phi,\Phi\rangle = 0 \quad\text{on } X. 
\end{multline}
According to Feehan and Leness \cite[Section 2.4, Lemma 2.19 (2), p. 296]{FL1}, we have the following pointwise inequalities when $r_E = 2$, which we now assume:
%PF6-28-2024 Check below. No real part needed?
\[
  \left(\frac{5}{4}-\frac{1}{\sqrt{2}}\right)|\Phi|^4 \leq 
  \langle(\Phi\otimes\Phi^*)_{00}\Phi,\Phi\rangle
  \leq \frac{1}{\sqrt{2}}|\Phi|^4 \quad\text{on } X.
\]
Since $5/4 - 1/\sqrt{2} \geq 3/4$, the equality \eqref{eq:Taubes_1996_SW_to_Gr_2-1_W+E-inner_product_Phi_2nd-order_Kato_PU(2)_monopole} yields
\begin{multline*}
  \frac{1}{2}\Delta_g|\Phi|^2 + |\nabla_{A_W\otimes A}\Phi|^2 + \frac{R}{4}|\Phi|^2
  + \frac{3}{4}|\Phi|^4
  \\
  + \frac{1}{2}\Real\,\langle\rho(F_{A_L}^+)\Phi,\Phi\rangle
  + \frac{1}{2}\Real\,\langle\rho(F_{A_d}^+)\Phi,\Phi\rangle
  - \frac{r}{4}\Real\,\langle\wp\jmath\Phi,\Phi\rangle \leq 0 \quad\text{on } X.
\end{multline*}
But
\begin{subequations}
  \begin{align}
    \label{eq:Inner_product_F_A_L+_Phi_Phi_geq_minus_norm_F_A_L+_Phi_squared}
    \Real\,\langle\rho(F_{A_L}^+)\Phi,\Phi\rangle &\geq - |\rho(F_{A_L}^+)||\Phi|^2 \quad\text{on } X,
    \\
    \label{eq:Inner_product_F_A_d+_Phi_Phi_geq_minus_norm_F_A_d+_Phi_squared}
    \Real\,\langle\rho(F_{A_d}^+)\Phi,\Phi\rangle &\geq - |\rho(F_{A_d}^+)||\Phi|^2 \quad\text{on } X,
    \\
    \label{eq:Inner_product_wpj_Phi_Phi_geq_minus_norm_Phi_squared}
    \Real\,\langle\wp\jmath\Phi,\Phi\rangle &\geq - |\Phi|^2 \quad\text{on } X,
\end{align}
\end{subequations}
and so the preceding inequality yields
\[
  \frac{1}{2}\Delta_g|\Phi|^2 + |\nabla_{A_W\otimes A}\Phi|^2 + \frac{3}{4}|\Phi|^4
  \leq
  \left(\frac{1}{2}|\rho(F_{A_L}^+)| + \frac{1}{2}|\rho(F_{A_d}^+)| + \frac{r}{4} - \frac{R}{4}\right)|\Phi|^2
  \quad\text{on } X.
\]
Consequently, using $|\nabla_{A_W\otimes A}\Phi|^2 \geq 0$ and arguing just as in Salamon's proof of \cite[Section 7.2, Lemma 7.13, pp. 231--232]{SalamonSWBook}, we see that
\begin{equation}
  \label{eq:Delta_norm_squared_Phi_leq_c_norm_squared_Phi}
  \Delta_g|\Phi|^2
  \leq
  \left(\frac{1}{2}|\rho(F_{A_L}^+)| + \frac{1}{2}|\rho(F_{A_d}^+)|
    + \frac{r}{4} - \frac{R}{4} - \frac{3}{4}|\Phi|^2\right)|\Phi|^2 \quad\text{on } X.
\end{equation}
Let $x_0 \in X$ be a point at which the function $X \ni x \mapsto |\Phi(x)|^2$ attains its maximum. At such a point, by our assumption that $\Phi$ is $C^\infty$, we have by \eqref{eq:Laplace-Beltrami_operator} that (for geodesic normal coordinates $\{x_j\}_{j=1}^4$ centered at the point $x_0$)
\[
  (\Delta_g|\Phi|^2)(x_0) = -\sum_{j=1}^4 \frac{\partial^2u}{\partial x_j^2}(x_0)\geq 0,
\]
and hence either $\Phi(x_0) = 0$ or
%PF7-12-2024 Add reference to maximum principle and explain better
\[
  \frac{1}{2}|\rho(F_{A_L}^+)|(x_0) + \frac{1}{2}|\rho(F_{A_d}^+)|(x_0)
  + \frac{r}{4} - \frac{R(x_0)}{4} - \frac{3}{4}|\Phi(x_0)|^2 \geq 0,
\]
and thus
\begin{equation}
  \label{eq:Upper_bound_pointwise_norm_squared_Phi_X}
  |\Phi|^2 \leq \frac{1}{3}\sup_X\left(2|\rho(F_{A_L}^+)| + 2|\rho(F_{A_d}^+)| + r - R\right) \quad\text{on } X.
\end{equation}
We now assume that $(X,g,J,\omega)$ is an almost Hermitian four-manifold and restrict to the canonical \spinc structure $(\rho_\can,W_\can)$ on $X$.
%PF7-9-2024 Add reference below
% (see equations \eqref{eq:Canonical_spinc_bundles} for $W_\can$)
%PF8-16-2024 Explain L_\can and A_{L_\can}.
After making the substitution \eqref{eq:Taubes_1996_SW_to_Gr_eq_1-19}, namely
\[
  \Phi = (\varphi,\psi) = r^{1/2}(\alpha,\beta) \in \Omega^0(E\oplus\Lambda^{0,2}(E)),
\]
the preceding inequality becomes
\begin{equation}
  \label{eq:Taubes_1996_SW_to_Gr_eq_2-2_raw}
  |\alpha|_E^2 + |\beta|_{\Lambda^{0,2}(E)}^2
  \leq
  \frac{1}{3}\sup_X\left(1 + \frac{2}{r}|\rho(F_{A_L}^+)| + \frac{2}{r}|\rho(F_{A_d}^+)| - \frac{R}{r}\right) \quad\text{on } X.
\end{equation}<
We define the constant
\begin{equation}
  \label{eq:Definition_z}
  z := \frac{1}{3}\max\left\{1, \sup_X\left(2|\rho(F_{A_L}^+)| + 2|\rho(F_{A_d}^+)| - \frac{R}{2}\right)\right\},
\end{equation}
and substitute this definition of $z \in [1/3,\infty)$ into \eqref{eq:Taubes_1996_SW_to_Gr_eq_2-2_raw} to give the desired upper bound \eqref{eq:Taubes_1996_SW_to_Gr_eq_2-2}. We thus obtain an almost exact analogue of Taubes \cite[Section 2 (a), Proposition 2.1, p. 855]{TauSWGromov}, with a similar proof but for the non-Abelian monopole equations with a singular Taubes perturbation. This completes our first proof of Proposition \ref{mainprop:Taubes_1996_SW_to_Gr_2-1}.
\end{proof}

\begin{proof}[Second proof of Proposition \ref{mainprop:Taubes_1996_SW_to_Gr_2-1}]
For this argument, we remove the requirement in our first proof that $\Phi$ be $C^\infty$ and require only that $\Phi \in W^{2,p}(W_\can^+\otimes E)$. We proceed as in the first proof of Proposition \ref{mainprop:Taubes_1996_SW_to_Gr_2-1}, except for the application of the maximum principle for $\Delta_g$ to prove the key inequality \eqref{eq:Upper_bound_pointwise_norm_squared_Phi_X}.

We define the elliptic operator
\begin{equation}
  \label{eq:Elliptic_linear_scalar_second-order_operator_norm_squared_Phi}
  Lu := -\Delta_gu + cu, \quad\text{for all } u \in W^{2,p}(X),
\end{equation}
where
\begin{equation}
  \label{eq:Elliptic_linear_scalar_second-order_operator_coefficient_c_W2p}
  c := \frac{1}{2}|\rho(F_{A_L}^+)| + \frac{1}{2}|\rho(F_{A_d}^+)|
  + \frac{r}{4} - \frac{R}{4} - \frac{3}{4}|\Phi|^2 \quad\text{on } X.
\end{equation}
By hypothesis, $\Phi \in W^{2,p}(W^+\otimes E)$ with $p \in [4,\infty)$. Because $\dim X = 4$, we have that $W^{2,p}(X)$ is a Banach algebra,
% PF11-24-2024 Add ref
so $|\Phi|^2 = \langle\Phi,\Phi\rangle \in  W^{2,p}(X)$ and thus $c \in W^{2,p}(X)$ since the connections $A_L$ and $A_d$ and the Levi--Civita connection $\nabla^g$ on $TX$ and hence their curvatures are $C^\infty$. The elliptic operator $L$ in \eqref{eq:Elliptic_linear_scalar_second-order_operator_norm_squared_Phi} matches that of Definition \ref{defn:Second_order_linear_elliptic_operator_riemannian_manifold} with $b \equiv 0$ on $X$ and $c$ as in \eqref{eq:Elliptic_linear_scalar_second-order_operator_coefficient_c_W2p}. For $u := |\Phi|^2 \in W^{2,p}(X)$, we have by \eqref{eq:Delta_norm_squared_Phi_leq_c_norm_squared_Phi} that
\[
  Lu \geq 0 \quad\text{on } X.
\]
Define $X_{c>0} := \{x\in X: c(x) > 0\}$, and $X_{c<0} := \{x\in X: c(x) < 0\}$, and $X_{c=0} := \{x\in X: c(x) = 0\}$. Observe that $\partial X_{c>0} = X_{c=0} = \partial X_{c<0}$ and $X_{c\geq 0} = X_{c>0}\cup X_{c=0}$ and $X_{c\leq 0} = X_{c<0}\cup X_{c=0}$. By definition, we have
\[
  c = \frac{1}{2}|\rho(F_{A_L}^+)| + \frac{1}{2}|\rho(F_{A_d}^+)|
  + \frac{r}{4} - \frac{R}{4} - \frac{3}{4}|\Phi|^2 > 0 \quad\text{on } X_{c>0},
\]
and by continuity, this gives
\begin{equation}
  \label{eq:Upper_bound_pointwise_norm_squared_Phi_Xcgeq0}
  |\Phi|^2 \leq \frac{1}{3}\left(2|\rho(F_{A_L}^+)| + 2|\rho(F_{A_d}^+)|
  + r - R\right) \quad\text{on } X_{c\geq 0}.
\end{equation}
If $X_{c<0} = \varnothing$, then $c \geq 0$ on $X$ and this yields the desired inequality \eqref{eq:Upper_bound_pointwise_norm_squared_Phi_X}, so it remains to consider the case $X_{c<0} \neq \varnothing$. Because $c < 0$ on $X_{c<0}$, the weak maximum principle provided by Corollary \ref{cor:Gilbarg_Trudinger_3-2} implies that
\[
  \sup_{X_{c<0}} u = \sup_{\partial X_{c<0}} u^+ = \sup_{\partial X_{c<0}} u,
\]
where the second equality follows from the fact that $u^+ = \max\{u,0\}$ on $X$ and $u^+ = u$ since $u = |\Phi|^2 \geq 0$ on $X$. But $\partial X_{c<0} = X_{c=0}$, and so by \eqref{eq:Upper_bound_pointwise_norm_squared_Phi_Xcgeq0} the function $u = |\Phi|^2$ obeys 
\[
  \sup_{\partial X_{c<0}} u = \sup_{X_{c=0}} u
  \leq \frac{1}{3}\sup\left(2|\rho(F_{A_L}^+)| + 2|\rho(F_{A_d}^+)| + r - R\right).
\]
Therefore,
\begin{equation}
  \label{eq:Upper_bound_pointwise_norm_squared_Phi_Xc<0}
  |\Phi|^2 \leq \frac{1}{3}\sup\left(2|\rho(F_{A_L}^+)| + 2|\rho(F_{A_d}^+)| + r - R\right)
  \quad\text{on } X_{c<0}.
\end{equation}
Combining inequalities \eqref{eq:Upper_bound_pointwise_norm_squared_Phi_Xcgeq0} and \eqref{eq:Upper_bound_pointwise_norm_squared_Phi_Xc<0} yields the desired inequality \eqref{eq:Upper_bound_pointwise_norm_squared_Phi_X}. The remainder of the argument is the same as that in our first proof of Proposition \ref{mainprop:Taubes_1996_SW_to_Gr_2-1}.
\end{proof}

\subsection[Pointwise estimate for sections of $E$ using a regularized perturbation]{Pointwise estimate for sections of $E$ using a regularized Taubes perturbation}
\label{subsec:Pointwise_estimate_holomorphic_component_coupled_spinor_regular}
In this subsection, we give two proofs of Proposition \ref{mainprop:Taubes_1996_SW_to_Gr_2-1_regular}, without any additional assumptions.

\begin{proof}[First proof of Proposition \ref{mainprop:Taubes_1996_SW_to_Gr_2-1_regular}]  
When we replace the singular perturbation $\wp(\psi)$ in \eqref{eq:Definition_wp_intro} by the regularized perturbation $\wp_\gamma(\psi)$ in \eqref{eq:Definition_wp_intro_regular}, Theorem \ref{thm:Feehan_Leness_1998jdg_3-7_regularized_Taubes-perturbed_W1p} ensures that we can apply a gauge transformation $v \in W^{2,p}(\SU(E))$ to a $W^{1,p}\times W^{2,p}$ solution to the non-Abelian monopole equations with regularized perturbation $\wp_\gamma(\psi)$ and obtain a $C^\infty$ solution. Thus, we may assume that $(A,\Phi)$ is $C^\infty$.

We now describe the minor modification needed in our first proof of Proposition \ref{mainprop:Taubes_1996_SW_to_Gr_2-1}. By applying the forthcoming inequalities \eqref{eq:wp_intro_gamma_lower_bound_varphi} and \eqref{eq:wp_intro_gamma_lower_bound_psi_crude}, we see that the lower bound $\Real\,\langle\wp\jmath\Phi,\Phi\rangle \geq - |\Phi|^2$ in \eqref{eq:Inner_product_wpj_Phi_Phi_geq_minus_norm_Phi_squared} is replaced by
%PF11-25-2024 Need to say where \jmath is defined  
\begin{equation}
  \label{eq:Inner_product_wp_gamma_j_Phi_Phi_geq_minus_norm_Phi_squared}
  \Real\,\langle\wp_\gamma\jmath\Phi,\Phi\rangle \geq -2|\Phi|^2 \quad\text{on } X.
\end{equation}
This change has the effect that in the differential inequality \eqref{eq:Delta_norm_squared_Phi_leq_c_norm_squared_Phi}, the term $r/4$ is replaced by $r/2$ and in the pointwise estimate \eqref{eq:Upper_bound_pointwise_norm_squared_Phi_X}, the term $r$ is replaced by $2r$. Hence, the inequality \eqref{eq:Taubes_1996_SW_to_Gr_eq_2-2_raw} becomes
\begin{equation}
  \label{eq:Taubes_1996_SW_to_Gr_eq_2-2_raw_wp_regular}
  |\alpha|_E^2 + |\beta|_{\Lambda^{0,2}(E)}^2
  \leq
  \frac{2}{3}\sup_X\left(1 + \frac{1}{r}|\rho(F_{A_L}^+)| + \frac{1}{r}|\rho(F_{A_d}^+)| - \frac{R}{2r}\right) \quad\text{on } X,
\end{equation}
and therefore,
\[
  |\alpha|_E^2 + |\beta|_{\Lambda^{0,2}(E)}^2
  \leq
  \frac{2}{3} + \frac{1}{r}\sup_X\left(2|\rho(F_{A_L}^+)| + 2|\rho(F_{A_d}^+)| - \frac{R}{2}\right) \quad\text{on } X.
\]  
This inequality yields \eqref{eq:Taubes_1996_SW_to_Gr_eq_2-2_regular}, with  constant $z$ as in \eqref{eq:Definition_z}, and completes our proof of Proposition \ref{mainprop:Taubes_1996_SW_to_Gr_2-1_regular} when $(A,\Phi)$ is $C^\infty$.
\end{proof}

\begin{proof}[Second proof of Proposition \ref{mainprop:Taubes_1996_SW_to_Gr_2-1_regular}]  
We describe the minor modification needed in our second proof of Proposition \ref{mainprop:Taubes_1996_SW_to_Gr_2-1}, where $(A,\Phi)$ is only assumed to be a $W^{1,p}\times W^{2,p}$ solution, for $p \in [4,\infty)$, to the non-Abelian monopole equations with regularized perturbation $\wp_\gamma(\psi)$. The only change is that we replace the term $r/4$ in the definition \eqref{eq:Elliptic_linear_scalar_second-order_operator_coefficient_c_W2p} of $c$ by $r/2$ and replace the term $r$ in the pointwise estimates \eqref{eq:Upper_bound_pointwise_norm_squared_Phi_Xcgeq0} and \eqref{eq:Upper_bound_pointwise_norm_squared_Phi_Xc<0}
by $2r$. No other changes are needed to our second proof of Proposition \ref{mainprop:Taubes_1996_SW_to_Gr_2-1} when $\wp(\psi)$ is replaced by $\wp_\gamma(\psi)$. This completes our proof of Proposition \ref{mainprop:Taubes_1996_SW_to_Gr_2-1_regular} when $(A,\Phi)$ is $W^{1,p}\times W^{2,p}$.
\end{proof}

\section[Bochner--Weitzenb\"ock identities for $(0,1)$-connections]{Bochner--Weitzenb\"ock identities for $(0,1)$-connections on Hermitian vector bundles over almost Hermitian manifolds}
\label{sec:Bochner-Weitzenbock_formulae_holomorphic_and_anti-holomorphic_spinors}
Our primary goal in this section is to prove generalizations of the two important Bochner--Weitzenb\"ock identities stated by Kotschick in \cite[Lemma 2.4, p. 203 and Lemma 2.5, p. 204]{KotschickSW}, extending them from the case of a Hermitian line bundle over a symplectic four-manifold to a Hermitian vector bundle of arbitrary rank over a symplectic $2n$-manifold. As we shall see, a Bochner--Weitzenb\"ock identity of the form in \cite[Lemma 2.5, p. p. 204]{KotschickSW} is an extension to the case of the $\bar\partial_A$-Laplacian on $(p,q)$-forms, for arbitrary $0\leq p,q \leq n$ with values in a Hermitian vector bundle from the case of sections of a Hermitian vector bundle, as in \cite[Lemma 2.4, p. 203]{KotschickSW}, where $p=q=0$.

\subsection[Bochner--Weitzenb\"ock identity for the $\bar\partial_A$-Laplacian]{Bochner--Weitzenb\"ock identity for the $\bar\partial_A$-Laplacian on sections of a Hermitian vector bundle over a symplectic manifold}
\label{subsec:Bochner-Weitzenbock_formulae_holomorphic_spinors}
We state and prove a generalization of Kotschick \cite[Lemma 2.4, p. 203]{KotschickSW}, from the case of a Hermitian line bundle to a Hermitian vector bundle $E$ of arbitrary rank. (Two similar formulae appear in Donaldson and Kronheimer \cite[Lemma 6.1.7, p. 212]{DK} but, as noted by Kronheimer \cite{Kronheimer_2024-7-23_private}, the terms $i\Lambda_\omega F_A$ should have their signs reversed and having the correct sign is crucial for our application. However, the two formulae in Lemma \ref{lem:Kotschick_2-4} match those of Donaldson \cite[Section 1.1, Proof of Proposition 3 (ii), p. 6]{DonASD}.) The forthcoming identity \eqref{eq:Kotschick_section_2_eq_12} is an exact analogue of Kotschick \cite[Lemma 2.4, Equation (12), p. 203]{KotschickSW} and Donaldson \cite[Section 4, p. 60, fourth displayed equation]{DonSW}.
We then have the

\begin{lem}[Bochner--Weitzenb\"ock identity for the $\partial_A$ and $\bar\partial_A$-Laplacians on $\Omega^0(E)$ over almost K\"ahler manifolds]
\label{lem:Kotschick_2-4}
(Compare Donaldson and Kronheimer \cite[Lemma 6.1.7, p. 212]{DK} and Kotschick \cite[Lemma 2.4, p. 203]{KotschickSW}.)
Let $(X,g,J,\omega)$ be an almost K\"ahler manifold and $A$ be a smooth, unitary connection on a smooth, Hermitian vector bundle $(E,H)$ over $X$. If $\varphi \in \Omega^0(E)$, then
\begin{subequations}
   \label{eq:Kotschick_section_2_eq_12_dbarA_and_delA}
  \begin{align}
    \label{eq:Kotschick_section_2_eq_12}
    \bar\partial_A^*\bar\partial_A\varphi
    &=
      \frac{1}{2}\nabla_A^*\nabla_A\varphi - \frac{i}{2}\Lambda_\omega F_A\varphi,
    \\
    \label{eq:Kotschick_section_2_eq_12_delA}
    \partial_A^*\partial_A\varphi
    &=
      \frac{1}{2}\nabla_A^*\nabla_A\varphi + \frac{i}{2}\Lambda_\omega F_A\varphi.
  \end{align}
\end{subequations}
\end{lem}

\begin{rmk}[Vanishing theorem and sign of the curvature term in \eqref{eq:Kotschick_section_2_eq_12}]
\label{rmk:Kodaira-Nakano_vanishing_theorem}
The well-known Kodaira--Nakano vanishing theorem provides a check on the correct sign for the curvature term $i\Lambda_\omega F_A$ in \eqref{eq:Kotschick_section_2_eq_12} in Lemma \ref{lem:Kotschick_2-4} and in this context, one can see that the statement and proof of \cite[Corollary 6.1.19, p. 213]{DK} due to Donaldson and Kronheimer are incorrect. If $i\Lambda_\omega F_A \leq 0$ (respectively, $< 0$), the identity \eqref{eq:Kotschick_section_2_eq_12} then implies that 
\begin{align*}
  \|\bar\partial_A\varphi\|_{L^2(X)}^2
  &= \frac{1}{2}\|\nabla_A\varphi\|_{L^2(X)}^2
    - (i\Lambda_\omega F_A\varphi,\varphi)_{L^2(X)}
  \\
  &\geq \text{(respectively, $>$ ) } \frac{1}{2}\|\nabla_A\varphi\|_{L^2(X)}^2.
\end{align*}
Hence, if $i\Lambda_\omega F_A \leq 0$ and $\bar\partial_A\varphi = 0$, then $\nabla_A\varphi = 0$. If $i\Lambda_\omega F_A < 0$, then there cannot exist $\varphi \in \Omega^0(E)$ with $\bar\partial_A\varphi = 0$ unless $\varphi \equiv 0$. This matches the conclusions of the Kodaira--Nakano vanishing theorem when $E$ is a complex line bundle --- for example, see Feehan and Leness 
%TL12-4-2025: Updated and removed page references
\cite[Remark 8.4.9]{Feehan_Leness_introduction_virtual_morse_theory_so3_monopoles} --- and matches the well-known degree constraint for holomorphic sections of complex line bundles --- for example, see Feehan and Leness 
%TL12-4-2025: Updated and removed page references
\cite[Lemma 8.4.8]{Feehan_Leness_introduction_virtual_morse_theory_so3_monopoles} or Donaldson \cite[Section 1.1, Proof of Proposition 3 (ii), p. 6]{DonASD}.
\qed
\end{rmk}

\begin{proof}[Proof of Lemma \ref{lem:Kotschick_2-4}]
This argument is a direct adaptation of Donaldson and Kronheimer's proof of \cite[Lemma 6.1.7, p. 212]{DK}, but corrects the sign errors for the terms $i\Lambda_\omega F_A$ and extends their argument from the case of complex K\"ahler to almost K\"ahler manifolds. Using $\nabla_A = d_A = \partial_A + \bar\partial_A$ on $\Omega^0(E)$ and noting that $\partial_A^* = 0$ on $\Omega^{0,1}(E)$ and $\bar\partial_A^* = 0$ on $\Omega^{1,0}(E)$, we observe that
\begin{align*}
  \nabla_A^*\nabla_A\varphi
  &=
    (\partial_A + \bar\partial_A)^*(\partial_A + \bar\partial_A)\varphi
  \\
  &= (\partial_A^* + \bar\partial_A^*)(\partial_A + \bar\partial_A)\varphi
  \\
  &= \partial_A^*\partial_A\varphi + \bar\partial_A^*\bar\partial_A\varphi
  \\
  &= i\left(-i\partial_A^*\partial_A\varphi - i\bar\partial_A^*\bar\partial_A\varphi\right)
  \\
  &= i\left(\Lambda_\omega\bar\partial_A\partial_A\varphi - \Lambda_\omega\partial_A\bar\partial_A\varphi\right),
\end{align*}
noting that $\Lambda_\omega\partial_A\varphi = 0$ and $\Lambda_\omega\bar\partial_A\varphi = 0$ and applying the K\"ahler identities \eqref{eq:Kaehler_identity_commutator_Lambda_del-bar_A_and_Lambda_del_A} to obtain the last equality above. Therefore,
\begin{equation}
  \label{eq:Kotschick_section_2_eq_12_without_curvature}
  \nabla_A^*\nabla_A\varphi
  =
  i\Lambda_\omega\left(\bar\partial_A\partial_A\varphi - \partial_A\bar\partial_A\varphi\right).
\end{equation}
But we also have
%PF7-24-2024 Add reference/explanation
\begin{equation}
  \label{eq:FA11_almost_Hermitian}
  F_A^{1,1} = \partial_A\bar\partial_A + \bar\partial_A\partial_A
\end{equation}
and so
\begin{equation}
  \label{eq:Lambda_omega_FA_almost_Hermitian}
  i\Lambda_\omega F_A = i\Lambda_\omega\left(\partial_A\bar\partial_A + \bar\partial_A\partial_A\right).
\end{equation}
By applying $i\Lambda_\omega F_{A_q^*\otimes A}$ to $\varphi$ and subtracting equation \eqref{eq:Kotschick_section_2_eq_12_without_curvature} from \eqref{eq:Lambda_omega_FA_almost_Hermitian}, we see that
\[
  i\Lambda_\omega F_A\varphi - \nabla_A^*\nabla_A\varphi  = 2i\Lambda_\omega\partial_A\bar\partial_A\varphi.
\]
By applying the K\"ahler identity \eqref{eq:Kaehler_identity_commutator_Lambda_del_A} and noting that $\Lambda_\omega\bar\partial_A\varphi = 0$, we obtain
\[
  i\Lambda_\omega F_A\varphi - \nabla_A^*\nabla_A\varphi
  = 2i(i\bar\partial_A^*)\bar\partial_A\varphi
  = -2\bar\partial_A^*\bar\partial_A\varphi.
\]
This establishes \eqref{eq:Kotschick_section_2_eq_12}.

To obtain \eqref{eq:Kotschick_section_2_eq_12_delA}, we instead add equation \eqref{eq:Kotschick_section_2_eq_12_without_curvature} to \eqref{eq:Lambda_omega_FA_almost_Hermitian} and obtain
\[
  i\Lambda_\omega F_A\varphi + \nabla_A^*\nabla_A\varphi  = 2i\Lambda_\omega\bar\partial_A\partial_A\varphi.
\]
By applying the K\"ahler identity \eqref{eq:Kaehler_identity_commutator_Lambda_del-bar_A} and noting that $\Lambda_\omega\partial_A\varphi = 0$, we obtain
\[
  i\Lambda_\omega F_A\varphi + \nabla_A^*\nabla_A\varphi
  = 2i(-i\partial_A^*)\partial_A\varphi
  = 2\partial_A^*\partial_A\varphi.
\]
This establishes \eqref{eq:Kotschick_section_2_eq_12_delA}. This completes our proof of Lemma \ref{lem:Kotschick_2-4}.
\end{proof}

\begin{rmk}[Generalizations of the Bochner--Weitzenb\"ock identity \eqref{eq:Kotschick_section_2_eq_12} to sections of $\Lambda^{p,q}(E)$ and $\Lambda^{p,0}(E)$]
\label{rmk:Kotschick_2-4_generalization}  
Moroianu \cite[Section 20.1, Theorem 20.2, p. 137]{Moroianu_2007} provides a generalization of the identity \eqref{eq:Kotschick_section_2_eq_12} to sections of $\Lambda^{p,q}(E)$ and a specialization of the latter identity in \cite[Section 20.1, Proposition 20.3, p. 138]{Moroianu_2007} to sections of $\Lambda^{p,0}(E)$. A generalization to sections of $\Lambda^{0,q}(E)$ appears in Li \cite[Equations (1) and (2), p. 623]{Li_2010}, based on Ma and Marinescu \cite[Remark 1.4.8, Equation (1.4.31), p. 34]{Ma_Marinescu_holomorphic_morse_inequalities_bergman_kernels}, and Liu \cite[Chapter 3, Theorem 3.1.5]{Liu_complex_manifold_and_Kaehler_geometry_2018}, and a generalization to sections of $\Lambda^{p,q}(E)$ appears in Li \cite[Equation (1) and Theorem 3.1, p. 628]{Li_2010}.
\qed
\end{rmk}

\begin{rmk}[K\"ahler identities for almost Hermitian manifolds and applications]
\label{rmk:Kaehler_identities_for_almost_Hermitian_manifolds}
Cirici and Wilson \cite{Cirici_Wilson_2020_almost_hermitian_identities} have derived analogues for almost Hermitian $2n$-manifolds of the K\"ahler identities that are known for K\"ahler and almost K\"ahler manifolds (see, for example, Cirici and Wilson \cite[Proposition 3.1, p. 34]{Cirici_Wilson_2020_harmonic}). According to Cirici and Wilson \cite[Proposition 1, p. 6]{Cirici_Wilson_2020_almost_hermitian_identities}, one has\footnote{When $q=1$, this identity is also stated by Cirici and Wilson as \cite[Lemma 4.8, p. 22]{Cirici_Wilson_2021}, though there is a typographical error in that $\alpha \in \Omega^{0,1}(X;\CC)$, not $\alpha \in \Omega^{1,0}(X;\CC)$; this is clear from its specialization to $\alpha = \bar\partial f$, for $f \in \Omega^0(X;\CC)$, and also confirmed by the authors \cite{Cirici_Wilson_2024-6-5_private}.}
%PF7-17-2024 \alpha, \beta need to be flipped or changed to something else to avoid Taubes' conflict.
\begin{equation}
  \label{eq:Cirici_Wilson_prop_1}
  \Lambda_\omega\partial\alpha  = i\bar\partial^*\alpha + i[\Lambda_\omega,\bar\partial^*]L_\omega\alpha,
  \quad\text{for all } \alpha \in \Omega^{0,q}(X;\CC).
\end{equation}
(See also Ohsawa \cite[Appendix]{Ohsawa_1982}.) By applying complex conjugation to the preceding identity and writing $\beta = \bar\alpha \in \Omega^{q,0}(X;\CC)$, we obtain
\begin{equation}
  \label{eq:Cirici_Wilson_prop_1_complex_conjugate}
  \Lambda_\omega\bar\partial\beta = -i\partial^*\beta - i[\Lambda_\omega,\partial^*]L_\omega\beta,
  \quad\text{for all } \beta \in \Omega^{q,0}(X;\CC).
\end{equation}
When $n=2$ and $q=1$, then the preceding identities reduce to
\begin{align*}
  \Lambda_\omega\bar\partial\beta
  &= -i\partial^*\beta, \quad\text{for all } \beta \in \Omega^{1,0}(X;\CC),
  \\
  \Lambda_\omega\partial\alpha
  &= i\bar\partial^*\alpha,   \quad\text{for all } \alpha \in \Omega^{0,1}(X;\CC).
\end{align*}
which are a special case of the standard K\"ahler identities on a complex K\"ahler $n$-manifold. This suggests that Lemma \ref{lem:Kotschick_2-4} continues to hold when $(X,g,J,\omega)$ is almost Hermitian rather than almost K\"ahler for $X$ of real dimension four.
\qed
\end{rmk}

\begin{rmk}[Cirici and Wilson identities]
\label{rmk:Cirici_Wilson_identities}  
See Cirici and Wilson \cite{Cirici_Wilson_2020_almost_hermitian_identities, Cirici_Wilson_2024}, Fernandez and Hosmer \cite{Fernandez_Hosmer_2022arxiv}, and Tomassini and \cite{Tomassini_Wang_2018}. From \cite[Section 4.1, p. 9]{Cirici_Wilson_2024},
\[
  [\Lambda_\omega, d] = -(d_c^* + T_c^*) \quad\text{on } \Omega^\bullet(X;\CC).
\]
Suppose $\beta \in \Omega^{1,0}(X;\CC)$ and $X$ has real dimension $2n$ with $n=2$. Then $[\Lambda_\omega, d]\beta = \Lambda_\omega d\beta + d\Lambda_\omega\beta$. Writing $d = \partial + \mu + \bar\mu + \bar\partial$, we see
\[
  [\Lambda_\omega, d]\beta
  =
  \Lambda_\omega \partial\beta + \partial\Lambda_\omega\beta + \Lambda_\omega \mu\beta + \mu\Lambda_\omega\beta
  + \Lambda_\omega\bar\partial\beta + \bar\partial\Lambda_\omega\beta
  + \Lambda_\omega\bar\mu\beta + \bar\mu\Lambda_\omega\beta.
\]
Because $\mu$ has bidegree $(2,-1)$, so $\mu\beta = 0$, and $\bar\mu$ has bidegree $(-1,2)$ and $\Lambda_\omega\beta = 0$, we obtain
\[
  [\Lambda_\omega, d]\beta
  =
  \Lambda_\omega \partial\beta + \Lambda_\omega\bar\partial\beta + \Lambda_\omega\bar\mu\beta.
\]
But $\partial\beta \in \Omega^{2,0}(X;\CC)$, so $\Lambda_\omega \partial\beta = 0$; similarly, $\bar\mu\beta \in \Omega^{0,2}(X;\CC)$, so $\Lambda_\omega\bar\mu\beta = 0$. Thus,
\[
  [\Lambda_\omega, d]\beta
  =
  \Lambda_\omega\bar\partial\beta.
\]
The $T_c^*$ operator is zero on $\beta$,
% PF7-17-2024 Check preceding assertion
%PF7-17-2024 Definition of d_c below matches Huybrechts Definition 3.1.13, p. 120
while $d_c = i(\bar\partial - \partial)$ and $d_c^* = i(\partial^* - \bar\partial^*)$. But $\beta \in \Omega^{1,0}(X;\CC)$, so $\bar\partial^*\beta = 0$. Thus,
\[
  -(d_c^* + T_c^*)\beta = -i\partial^*\beta,
\]
and so
\[
  \Lambda_\omega\bar\partial\beta = -i\partial^*\beta,
\]
as expected.
\qed
\end{rmk}

\subsection[Akizuki--Kodaira--Nakano identity relating the $\partial_A$ and $\bar\partial_A$-Laplacians]{Akizuki--Kodaira--Nakano identity relating the $\partial_A$ and $\bar\partial_A$-Laplacians on $(p,q)$-forms with values in a Hermitian vector bundle over an almost K\"ahler manifold}
\label{subsec:Akizuki--Kodaira--Nakano_identity_dA_and_dbarA_Laplacians_pq_forms_almost_Kaehler}
We include the forthcoming Lemmas \ref{lem:BW_identity_delA_and_dbarA_Laplacians_symplectic} and \ref{lem:BW_identity_dA_and_delA_and_dbarA_Laplacians_complex_Kaehler} for completeness, although they are not required our present application.

Recall that if $(g,J,\omega)$ is a smooth compatible triple on a smooth $2n$-manifold $X$, then by Kobayashi's convention\footnote{Though we omit the factor of $1/2$ in Kobayashi's definition.} \cite[Equations (7.6.5) and (7.6.7), p. 251]{Kobayashi_differential_geometry_complex_vector_bundles}, the triple defines a smooth Hermitian metric $h$ on $X$ by $h = g + i\omega$, where $\omega = g(\cdot,J\cdot)$ is as in \cite[Equation (7.6.8), p. 251]{Kobayashi_differential_geometry_complex_vector_bundles}.
% PF7-18-2024 Propagate above, relocate to remark/convention

We first prove some preparatory lemmas. By analogy with Huybrechts \cite[Definition 3.1.5, p. 116]{Huybrechts_2005}, we define
\begin{subequations}
\label{eq:del_A_and_dbar_A_laplacians}  
\begin{align}
  \label{eq:del_A_laplacian}  
  \Delta_{\partial_A} &:= \partial_A^*\partial_A + \partial_A\partial_A^*,
  \\
  \label{eq:dbar_A_laplacian}  
  \Delta_{\bar\partial_A} &:= \bar\partial_A^*\bar\partial_A + \bar\partial_A\bar\partial_A^*
  \quad\text{on } \Omega^{p,q}(E).
\end{align} 
\end{subequations}
When $A$ is the product connection on $E = X\times\CC^r$, the K\"ahler identities on a complex K\"ahler $n$-manifold imply that (see Huybrechts \cite[Proposition 3.1.12 (iii), p. 120]{Huybrechts_2005}) 
\begin{equation}
  \label{eq:Huybrechts_proposition_3-1-12-iii_page_120}
  \Delta_{\partial} = \Delta_{\bar\partial} \quad\text{on } \Omega^{p,q}(X).
\end{equation}
The identity \eqref{eq:Huybrechts_proposition_3-1-12-iii_page_120} continues to hold on an almost K\"ahler manifold. Indeed,  Huybrechts'\cite[Proof of Proposition 3.1.12 (iii), p. 122]{Huybrechts_2005} only uses the second pair of K\"ahler identities in \cite[Proposition 3.1.12 (ii), p. 120]{Huybrechts_2005}, namely \eqref{eq:Kaehler_identity_commutator_Lambda_del-bar_and_Lambda_del} and those identities hold for almost K\"ahler $2n$-manifolds as noted earlier. The following result appears to be well-known\footnote{See the Math Overflow articles ``Weitzenb\"ock identity for $\Delta_{\bar\partial_E}$'' \url{https://mathoverflow.net/questions/137612/weitzenbock-identity-for-delta-bar-partial-e} and ``Weitzenb\"ock identities'' \url{https://mathoverflow.net/questions/107168/weitzenbock-identities}.}, but as we would like the precise form of the curvature term, we include a proof. (The forthcoming Lemma \ref{lem:BW_identity_delA_and_dbarA_Laplacians_symplectic} and its proof are included for completeness but are not required our current applications.)

\begin{lem}[Akizuki--Kodaira--Nakano
identity relating the $\partial_A$ and $\bar\partial_A$-Laplacians on $\Omega^{p,q}(E)$ over almost K\"ahler manifolds]
\label{lem:BW_identity_delA_and_dbarA_Laplacians_symplectic}  
(See Liu and Yang \cite[Theorem 4.5]{Liu_Yang_2012} for the case of complex Hermitian manifolds and Li\footnote{Li reverse the sign of the curvature term.} \cite[Section 1, last displayed equation, p. 621]{Li_2010} for the case of complex K\"ahler manifolds.)
Let $(X,g,J,\omega)$ be an almost K\"ahler manifold and $A$ be a smooth, unitary connection on a smooth, Hermitian vector bundle $(E,H)$ over $X$. For all integers $0\leq p, q \leq n$, one has
\begin{equation}
  \label{eq:Huybrechts_proposition_3-1-12-iii_page_120_Omega_pq_E_commutator_Lambda_FA}
  \Delta_{\partial_A}
  =
  \Delta_{\bar\partial_A} + i[\Lambda_\omega, F_A]
  \quad\text{on } \Omega^{p,q}(E).
\end{equation}  
\end{lem}

\begin{proof}
We calculate
\begin{align*}
  \partial_A^*\partial_A + \partial_A\partial_A^*
  &= i\left(-i\partial_A^*\partial_A - \partial_A(i\partial_A^*)\right)
  \\
  &= i\left([\Lambda_\omega,\bar\partial_A]\partial_A + \partial_A[\Lambda_\omega,\bar\partial_A]\right)
    \quad\text{(by \eqref{eq:Kaehler_identity_commutator_Lambda_del-bar_and_Lambda_del})}
  \\
  &= i\left(\Lambda_\omega\bar\partial_A\partial_A - \bar\partial_A\Lambda_\omega\partial_A
    + \partial_A\Lambda_\omega\bar\partial_A - \partial_A\bar\partial_A\Lambda_\omega\right)
  \\
  &= i\left(\Lambda_\omega\bar\partial_A\partial_A - \bar\partial_A\left([\Lambda_\omega,\partial_A] + \partial_A\Lambda_\omega\right)
    \right.
    + \left.\left([\partial_A,\Lambda_\omega] + \Lambda_\omega\partial_A\right)\bar\partial_A
    - \partial_A\bar\partial_A\Lambda_\omega\right)
  \\
  &= i\left( \Lambda_\omega\bar\partial_A\partial_A - \bar\partial_A\left(i\bar\partial_A^* + \partial_A\Lambda_\omega\right)
    + \left(-i\bar\partial_A^* + \Lambda_\omega\partial_A \right)\bar\partial_A
    - \partial_A\bar\partial_A\Lambda_\omega\right)
    \quad\text{(by \eqref{eq:Kaehler_identity_commutator_Lambda_del-bar_and_Lambda_del})}
  \\
  &= i\left( \Lambda_\omega\bar\partial_A\partial_A - i\bar\partial_A\bar\partial_A^* - \bar\partial_A\partial_A\Lambda_\omega
    - i\bar\partial_A^*\bar\partial_A + \Lambda_\omega\partial_A\bar\partial_A
    - \partial_A\bar\partial_A\Lambda_\omega\right)
  \\
  &= \bar\partial_A\bar\partial_A^* + \bar\partial_A\bar\partial_A^*
    + i\Lambda_\omega\left(\bar\partial_A\partial_A + \partial_A\bar\partial_A\right)
    - i\left(\bar\partial_A\partial_A + \partial_A\bar\partial_A\right)\Lambda_\omega,
\end{align*}
and thus by \eqref{eq:del_A_and_dbar_A_laplacians}, we obtain
\[
  \Delta_{\partial_A}
  =
  \Delta_{\bar\partial_A} + i\Lambda_\omega\left(\bar\partial_A\partial_A + \partial_A\bar\partial_A\right)
  - \left(\bar\partial_A\partial_A + \partial_A\bar\partial_A\right)\Lambda_\omega
  \quad\text{on } \Omega^{p,q}(E).
\]
Recall that
\[
  F_A^{1,1} = \bar\partial_A\partial_A + \partial_A\bar\partial_A \in \Omega^{1,1}(\fu(E)),
\]
and $\Lambda_\omega F_A = F_A^{1,1}$, so the preceding identity for $\Delta_{\partial_A}$ yields \eqref{eq:Huybrechts_proposition_3-1-12-iii_page_120_Omega_pq_E_commutator_Lambda_FA}.
\end{proof}

\begin{rmk}[Special cases of the Akizuki--Kodaira--Nakano
formula \eqref{eq:Huybrechts_proposition_3-1-12-iii_page_120_Omega_pq_E_commutator_Lambda_FA}]
\label{rmk:BW_identity_delA_and_dbarA_Laplacians_symplectic_special_cases}  
In our applications, we only require particular cases of $(p,q)$ in Lemma \ref{lem:BW_identity_delA_and_dbarA_Laplacians_symplectic} in which the identity \eqref{eq:Huybrechts_proposition_3-1-12-iii_page_120_Omega_pq_E_commutator_Lambda_FA} simplifies considerably.
If $p=0$ or $q=0$, then equation \eqref{eq:Huybrechts_proposition_3-1-12-iii_page_120_Omega_pq_E_commutator_Lambda_FA} reduces to
\begin{equation}
  \label{eq:Huybrechts_proposition_3-1-12-iii_page_120_Omega_p0_or_0q_E_LambdaFA}
  \Delta_{\partial_A}
  =
  \Delta_{\bar\partial_A} + i\Lambda_\omega F_A\quad\text{on } \Omega^{0,q}(E) \text{ or } \Omega^{p,0}(E).
\end{equation}
If $p=0$, then equation \eqref{eq:Huybrechts_proposition_3-1-12-iii_page_120_Omega_pq_E_commutator_Lambda_FA} reduces to
\begin{equation}
  \label{eq:Huybrechts_proposition_3-1-12-iii_page_120_Omega_0q_E_LambdaFA}
  \partial_A^*\partial_A
  = \bar\partial_A^*\bar\partial_A + \bar\partial_A\bar\partial_A^* + i\Lambda_\omega F_A
  \quad\text{on } \Omega^{0,q}(E).
\end{equation}
If $p=0$ and $q=n$, then \eqref{eq:Huybrechts_proposition_3-1-12-iii_page_120_Omega_pq_E_commutator_Lambda_FA} reduces to
\begin{equation}
  \label{eq:Huybrechts_proposition_3-1-12-iii_page_120_Omega_0n_E_LambdaFA}
  \partial_A^*\partial_A
  = \bar\partial_A\bar\partial_A^* + i\Lambda_\omega F_A \quad\text{on } \Omega^{0,n}(E).
\end{equation}
If $p=q=0$, then equation \eqref{eq:Huybrechts_proposition_3-1-12-iii_page_120_Omega_pq_E_commutator_Lambda_FA} reduces to
\begin{equation}
  \label{eq:Huybrechts_proposition_3-1-12-iii_page_120_Omega_00_E_LambdaFA}
  \partial_A^*\partial_A
  = \bar\partial_A^*\bar\partial_A + i\Lambda_\omega F_A \quad\text{on } \Omega^0(E).
\end{equation}
Recall that by equation \eqref{eq:Kotschick_section_2_eq_12}, we have
\[
  \bar\partial_A^*\bar\partial_A\varphi
  =
  \frac{1}{2}\nabla_A^*\nabla_A\varphi - \frac{i}{2}\Lambda_\omega F_A\varphi,
\]
and so equation \eqref{eq:Huybrechts_proposition_3-1-12-iii_page_120_Omega_00_E_LambdaFA} yields
\begin{align*}
  \partial_A^*\partial_A\varphi
  &=
  \left(\frac{1}{2}\nabla_A^*\nabla_A\varphi + i\Lambda_\omega F_A\right) - \frac{i}{2}\Lambda_\omega F_A\varphi
  \\
  &=
  \frac{1}{2}\nabla_A^*\nabla_A\varphi + \frac{i}{2}\Lambda_\omega F_A\varphi,
\end{align*}
as expected from \eqref{eq:Kotschick_section_2_eq_12_delA}.
\qed
\end{rmk}

Recall that by analogy with the usual definition of the Laplace operator on $\Omega^r(X)$, for any integer $r\geq 0$, one has the Hodge Laplacian (see, for example, Lawson \cite[Appendix II, p. 93]{Lawson})
\begin{equation}
  \label{eq:Covariant_Laplacian}
  \Delta_A := d_A^*d_A + d_Ad_A^* \quad\text{on } \Omega^r(E),
\end{equation}
for a unitary connection $A$ on $E$. (The forthcoming Lemma \ref{lem:BW_identity_dA_and_delA_and_dbarA_Laplacians_complex_Kaehler} and its proof are included for completeness but are not required our current applications.) We then have the

\begin{lem}[Bochner--Weitzenb\"ock identity relating the $d_A$, $\partial_A$, and $\bar\partial_A$-Laplacians on $\Omega^{p,q}(E)$ over complex K\"ahler manifolds]
\label{lem:BW_identity_dA_and_delA_and_dbarA_Laplacians_complex_Kaehler}
Let $(X,g,J,\omega)$ be a
%PF8-14-2024 Almost Kaehler should suffice
complex K\"ahler manifold and $A$ be a smooth, unitary connection on a smooth, Hermitian vector bundle $(E,H)$ over $X$. For all integers $0\leq p, q \leq n$, one has
\begin{equation}
  \label{eq:Huybrechts_proposition_3-1-12-iii_page_120_DeltaA_Omega_pq_E_commutator_Lambda_FA}
  \Delta_A
  = \Delta_{\partial_A} + \Delta_{\bar\partial_A}
  - i[\Lambda_\omega, F_A^{2,0}] + i[\Lambda_\omega, F_A^{0,2}]      
  \quad\text{on } \Omega^{p,q}(E).
\end{equation}  
\end{lem}
%PF8-14-2024 Note that the preceding formula can be combined with the formula \eqref{eq:Huybrechts_proposition_3-1-12-iii_page_120_Omega_p0_or_0q_E_LambdaFA}

\begin{proof}
We calculate using \eqref{eq:Covariant_Laplacian} on $\Omega^{p,q}(E)$:
\begin{align*}
  d_A^*d_A + d_Ad_A^*
  &=
    (\partial_A^* + \bar\partial_A^*)(\partial_A + \bar\partial_A)
    + (\partial_A + \bar\partial_A)(\partial_A^* + \bar\partial_A^*)
  \\
  &= \partial_A^*\partial_A + \partial_A^*\bar\partial_A + \bar\partial_A^*\partial_A + \bar\partial_A^*\bar\partial_A
    + \partial_A\partial_A^* + \partial_A\bar\partial_A^* + \bar\partial_A\partial_A^* + \bar\partial_A\bar\partial_A^*
    \\
  &= \left(\partial_A^*\partial_A + \partial_A\partial_A^*\right)
    + \partial_A^*\bar\partial_A + \bar\partial_A^*\partial_A
    + \partial_A\bar\partial_A^* + \bar\partial_A\partial_A^*
    + \left(\bar\partial_A^*\bar\partial_A + \bar\partial_A\bar\partial_A^*\right),
\end{align*}
and thus
\[
  d_A^*d_A + d_Ad_A^*
  =
  \Delta_{\partial_A}
    + \partial_A^*\bar\partial_A + \bar\partial_A^*\partial_A
    + \partial_A\bar\partial_A^* + \bar\partial_A\partial_A^*
    + \Delta_{\bar\partial_A}.
\]
Applying the K\"ahler identities \eqref{eq:Kaehler_identity_commutator_Lambda_del-bar_and_Lambda_del} on $\Omega^{p,q}(E)$ yields
\begin{align*}
  d_A^*d_A + d_Ad_A^*
  &=
  \Delta_{\partial_A}
    - i\left(i\bar\partial_A^*\partial_A
    + \partial_A(i\bar\partial_A^*)
    + i\partial_A^*\bar\partial_A
    + \bar\partial_A(i\partial_A^*)\right)
    + \Delta_{\bar\partial_A}
  \\
  &=
    \Delta_{\partial_A}
    - i\left([\Lambda_\omega,\partial_A]\partial_A
    + \partial_A[\Lambda_\omega,\partial_A]
    - [\Lambda_\omega,\bar\partial_A]\bar\partial_A
    - \bar\partial_A[\Lambda_\omega,\bar\partial_A]\right)
    + \Delta_{\bar\partial_A}
  \\
  &=
    \Delta_{\partial_A}
    - i\left(\Lambda_\omega\partial_A\partial_A - \partial_A\Lambda_\omega\partial_A
    + \partial_A\Lambda_\omega\partial_A - \partial_A\partial_A\Lambda_\omega\right.
  \\
  &\qquad - \left. \Lambda_\omega\bar\partial_A\bar\partial_A + \bar\partial_A\Lambda_\omega\bar\partial_A
    - \bar\partial_A\Lambda_\omega\bar\partial_A +  \bar\partial_A\bar\partial_A\Lambda_\omega \right) + \Delta_{\bar\partial_A},
\end{align*}
and so
\[
  d_A^*d_A + d_Ad_A^*
  =
  \Delta_{\partial_A}
    - i[\Lambda_\omega,\partial_A\partial_A] + i[\Lambda_\omega,\bar\partial_A\bar\partial_A]
    + \Delta_{\bar\partial_A}.
\]
This gives the identity \eqref{eq:Huybrechts_proposition_3-1-12-iii_page_120_DeltaA_Omega_pq_E_commutator_Lambda_FA} since $\partial_A\partial_A = F_A^{2,0}$ and $\bar\partial_A\bar\partial_A = F_A^{0,2}$. This completes the proof of Lemma \ref{lem:BW_identity_dA_and_delA_and_dbarA_Laplacians_complex_Kaehler}.
\end{proof}

\begin{rmk}[Special cases of the Bochner--Weitzenb\"ock identity \eqref{eq:Huybrechts_proposition_3-1-12-iii_page_120_DeltaA_Omega_pq_E_commutator_Lambda_FA}]
\label{rmk:BW_identity_dA_and_delA_and_dbarA_Laplacians_complex_Kaehler_special_cases}
In our applications, we only require particular cases of $(p,q)$ in Lemma \ref{lem:BW_identity_dA_and_delA_and_dbarA_Laplacians_complex_Kaehler} in which the identity \eqref{eq:Huybrechts_proposition_3-1-12-iii_page_120_DeltaA_Omega_pq_E_commutator_Lambda_FA} simplifies considerably. If $p=0$ or $q=0$, then equation \eqref{eq:Huybrechts_proposition_3-1-12-iii_page_120_DeltaA_Omega_pq_E_commutator_Lambda_FA} reduces to
\begin{align*}
  \Delta_A
  &= \Delta_{\partial_A} + \Delta_{\bar\partial_A}
    - i\Lambda_\omega F_A^{2,0} + i\Lambda_\omega F_A^{0,2}
  \\
  &= \Delta_{\partial_A} + \Delta_{\bar\partial_A}
  \\
  &= 2\Delta_{\partial_A} - i\Lambda_\omega F_A \quad\text{(by \eqref{eq:Huybrechts_proposition_3-1-12-iii_page_120_Omega_p0_or_0q_E_LambdaFA})}
  \\
  &= 2\Delta_{\bar\partial_A} + i\Lambda_\omega F_A \quad\text{(by \eqref{eq:Huybrechts_proposition_3-1-12-iii_page_120_Omega_p0_or_0q_E_LambdaFA})}.
\end{align*}
If $p=q=0$, then equation \eqref{eq:Huybrechts_proposition_3-1-12-iii_page_120_DeltaA_Omega_pq_E_commutator_Lambda_FA} reduces to
\begin{align*}
  d_A^*d_A
  &= 2\partial_A^*\partial_A - i\Lambda_\omega F_A
  \\
  &= 2\bar\partial_A^*\bar\partial_A + i\Lambda_\omega F_A,
\end{align*}
and thus, noting that $d_A = \nabla_A$ on $\Omega^0(E)$, equation \eqref{eq:Huybrechts_proposition_3-1-12-iii_page_120_DeltaA_Omega_pq_E_commutator_Lambda_FA} reduces to
\begin{subequations}
\label{eq:DK_lemma_6-1-7_corrected_sign}  
\begin{align}
  \label{eq:DK_lemma_6-1-7b_corrected_sign}  
  \partial_A^*\partial_A
  &= \frac{1}{2}\nabla_A^*\nabla_A + \frac{i}{2}\Lambda_\omega F_A,
  \\
  \label{eq:DK_lemma_6-1-7a_corrected_sign}  
  \bar\partial_A^*\bar\partial_A
  &= \frac{1}{2}\nabla_A^*\nabla_A - \frac{i}{2}\Lambda_\omega F_A,
\end{align}
\end{subequations}
as expected from \eqref{eq:Kotschick_section_2_eq_12_dbarA_and_delA}. The identities \eqref{eq:DK_lemma_6-1-7_corrected_sign} continue to hold if $(X,g,J,\omega)$ is only almost K\"ahler, as proved in Lemma \ref{lem:Kotschick_2-4}. Indeed, one has $d_A = \partial_A + \bar\partial_A$ on $\Omega^0(E)$ and the facts that one now has $\partial_A\partial_A = F_A^{2,0} - \mu\bar\partial_A$ and $\bar\partial_A\bar\partial_A = F_A^{0,2} - \bar\mu\partial_A$ do not impact the structure of the identities \eqref{eq:DK_lemma_6-1-7_corrected_sign}.
\qed
\end{rmk}

\subsection[Bochner--Weitzenb\"ock identity for the $\bar\partial_A$-Laplacian on $(0,q)$-forms]{Bochner--Weitzenb\"ock identity for the $\bar\partial_A$-Laplacian on $(0,q)$-forms with values in a Hermitian holomorphic vector bundle over a complex K\"ahler manifold}
\label{subsec:BW_identity_dbarA-Laplacian_0q-forms_bundle-valued_complex_manifold}
We now seek a generalization --- again from the case of a Hermitian line bundle to a Hermitian vector bundle $E$ of arbitrary rank --- of the more difficult \cite[Lemma 2.5, p. 204]{KotschickSW} due to Kotschick, whose proof (for the case where $E$ is a Hermitian line bundle and $X$ has real dimension $4$) he only outlines. The forthcoming identities \eqref{eq:Kotschick_section_2_eq_13_rank-r_kaehler}, \eqref{eq:Kotschick_section_2_eq_13_rank-r_symplectic_01-connection_2n-dim_X}, and \eqref{eq:Kotschick_section_2_eq_13_rank-r_symplectic_01-connection_4-dim_X} are exact analogues of Kotschick \cite[Lemma 2.5, Equation (13), p. 204]{KotschickSW}.

Our first approach is to proceed along the lines suggested by Kotschick and while we develop this method in this subsection, we shall ultimately specialize to the case of Hermitian holomorphic vector bundles over complex K\"ahler manifolds. Our second approach (implied by Taubes \cite[Section 2 (b), pp. 854--855]{TauSWGromov}) is to make use of the well-known Bochner--Weitzenb\"ock identity for the Dirac operator Laplacian and the fact that $D_A = \sqrt{2}(\bar\partial_A + \bar\partial_A^*)$ on almost K\"ahler manifolds and we develop this method in the next subsection.

\begin{lem}[Bochner--Weitzenb\"ock identity for $\bar\partial_{A_\can\otimes A}$-Laplacian on $\Omega^{0,\bullet}(E)$]
\label{lem:Kotschick_2-5_raw_bullet}
Let $(X,g,J,\omega)$ be an almost K\"ahler $2n$-manifold and $A$ be a smooth, unitary connection on a smooth, Hermitian vector bundle $(E,H)$ over $X$.
%PF7-24-2024 Explain/correct below
For each integer $0\leq q\leq n$, let $A_q$ be the unitary connection on the Hermitian vector bundle $\Lambda^{0,q}(X)$ induced by the Hermitian metric $h = g + i\omega$ on $X$, let $A_\can = \oplus_{q=0}^n A_q$ be the induced unitary connection on the Hermitian vector bundle $\Lambda^{0,\bullet}(X) := \oplus_{q=0}^n\Lambda^{0,q}(X)$, and let $A_\can\otimes A$ be the induced unitary connection on $\Lambda^{0,\bullet}(X)\otimes E$. Then
\begin{multline}
  \label{eq:Kotschick_section_2_eq_13_raw_bullet}
  \bar\partial_{A_\can\otimes A}^*\bar\partial_{A_\can\otimes A}\Phi
  =
  \frac{1}{2}\nabla_{A_\can\otimes A}^*\nabla_{A_\can\otimes A}\Phi
  - \frac{i}{2}\Lambda_\omega\left(F_{A_\can} + F_A\right)\Phi,
  \\
  \text{for all } \Phi \in \Omega^{0,\bullet}(E) = \Omega^0(\Lambda^{0,\bullet}(X)\otimes E).
\end{multline}
\end{lem}

\begin{proof}  
The conclusion follows immediately by applying the Bochner--Weitzenb\"ock identity \eqref{eq:Kotschick_section_2_eq_12} from Lemma \ref{lem:Kotschick_2-4} to the unitary connection $A_\can\otimes A$ on the Hermitian vector bundle $\Lambda^{0,\bullet}(X)\otimes E$ and sections $\Phi \in \Omega^0(\Lambda^{0,\bullet}(X)\otimes E)$, after substituting $F_{A_\can\otimes A} = F_{A_\can} + F_A$.
\end{proof}

The forthcoming Proposition \ref{prop:Berline_Getzler_Vergne_3-71}, due to Berline, Getzler, and Vergne \cite{BerlineGetzlerVergne}, is stated by those authors in the setting of a Hermitian holomorphic vector bundle over a complex K\"ahler manifold. It should extend to the more general setting of a Hermitian vector bundle with an $(0,1)$-connection over an almost Hermitian manifold, but the proof of such a generalization would be more involved because, for example, one cannot use the complex normal coordinates that are available around points in complex K\"ahler manifolds. However, Proposition \ref{prop:Berline_Getzler_Vergne_3-71} and the forthcoming Lemma \ref{lem:Kotschick_2-5_complex_kaehler_rank-r_holomorphic} still provide a valuable check on signs in the next subsection.

\begin{prop}[Bochner--Weitzenb\"ock identity for Hermitian holomorphic vector bundles over complex K\"ahler manifolds]
\label{prop:Berline_Getzler_Vergne_3-71}
Let $(X,g,J,\omega)$ be a complex K\"ahler $n$-manifold and $A$ be the Chern connection defined by an integrable $(0,1)$-connection $\bar\partial_E$ on a smooth Hermitian vector bundle $(E,H)$.
%PF7-24-2024 Explain/correct below
For each integer $0\leq q\leq n$, let $A_q$ be the unitary connection on the Hermitian vector bundle $\Lambda^{0,q}(X)$ induced by the Hermitian metric $h = g + i\omega$ on $X$, let $A_\can = \oplus_{q=0}^n A_q$ be the induced unitary connection on the Hermitian vector bundle $\Lambda^{0,\bullet}(X) := \oplus_{q=0}^n\Lambda^{0,q}(X)$, and let $A_\can\otimes A$ be the induced unitary connection on $\Lambda^{0,\bullet}(X)\otimes E$. Then\footnote{We modify Berline, Getzler, and Vergne's statement of the identity by replacing $\iota(dz^j)$ with $\iota(\partial/\partial \bar z^j)$.}
%PF8-15-2024 Use j,k for local coordinate indices, \alpha,\beta for bundle indices
\begin{equation}
  \label{eq:Berline_Getzler_Vergne_prop_3-71}
  \bar\partial_A^*\bar\partial_A + \bar\partial_A\bar\partial_A^*
  =
  \Delta_{A_\can\otimes A}^{0,\bullet}
  +
  \sum_{i,j=1}^n\eps(d\bar z^i)\iota\left(\frac{\partial}{\partial\bar z^j}\right)F_{A_{K^*}\otimes A}
  \left(\frac{\partial}{\partial z^j},\frac{\partial}{\partial\bar z^i}\right)
  \quad\text{on } \Omega^{0,\bullet}(E),
\end{equation}
where $\{z^i\}_{i=1}^n$ are local holomorphic coordinates on $X$ and \cite[Section 3.6, Displayed equation prior to Lemma 3.70, p. 139]{BerlineGetzlerVergne}, 
\begin{equation}
  \label{eq:Berline_Getzler_Vergne_lemma_3-70}
  \Delta_{A_\can\otimes A}^{0,\bullet} 
  :=
  \bar\partial_{A_\can\otimes A}^*\bar\partial_{A_\can\otimes A}
  \quad\text{on } \Omega^{0,\bullet}(E) = \Omega^0(\Lambda^{0,\bullet}(X)\otimes E)),
\end{equation}
and $K := \wedge^n(T^{1,0}X)^*$ is the canonical line bundle on $X$, while $K^* = \wedge^n(T^{1,0}X)$ is the anti-canonical line bundle on $X$ with curvature $F_{A_{K^*}}$ defined by the Levi--Civita connection on $TX$.
\end{prop}

\begin{rmk}[Other versions of the Bochner--Weitzenb\"ock identity \eqref{eq:Berline_Getzler_Vergne_prop_3-71}]
\label{rmk:Berline_Getzler_Vergne_3-71}  
An identity similar to \eqref{eq:Berline_Getzler_Vergne_prop_3-71} appears in Kodaira and Morrow \cite[Chapter 3, Section 6, Theorem 6.1, p. 119]{Morrow_Kodaira_complex_manifolds} when $E=X\times\CC$ and more generally in \cite[Chapter 3, Section 6, Theorem 6.2, p. 124]{Morrow_Kodaira_complex_manifolds} when $E$ is a non-trivial line bundle. The identity \eqref{eq:Berline_Getzler_Vergne_prop_3-71} is stated exactly as we do in Proposition \ref{prop:Berline_Getzler_Vergne_3-71} by Boeijink, Landsman, and van Suijlekom in \cite[Proof of Theorem 3.6]{Boeijink_Landsman_van_Suijlekom_2019} in their proof of a Kodaira vanishing theorem. Although Berline, Getzler, and Vergne refer to \eqref{eq:Berline_Getzler_Vergne_prop_3-71} as a `Bochner--Kodaira' identity, that term appears to be reserved in the literature for identities that include those in Lemma \ref{lem:Kotschick_2-5_complex_kaehler_rank-r_holomorphic} and Remark \ref{rmk:Kotschick_2-5_raw_bullet_complex_kaehler}.
\qed
\end{rmk}

We now combine Lemma \ref{lem:Kotschick_2-5_raw_bullet} and Proposition \ref{prop:Berline_Getzler_Vergne_3-71} to prove the

\begin{lem}[Bochner--Kodaira identity for Hermitian holomorphic vector bundles over complex K\"ahler manifolds]
\label{lem:Kotschick_2-5_complex_kaehler_rank-r_holomorphic}  
(See Kotschick \cite[Lemma 2.5, Equation (13), p. 204]{KotschickSW} when $n=2$ and $(X,g,J,\omega)$ is symplectic and $E$ is a Hermitian line bundle.)
Continue the hypotheses of Proposition \ref{prop:Berline_Getzler_Vergne_3-71}. Then 
\begin{equation}
  \label{eq:Kotschick_section_2_eq_13_rank-r_kaehler}
  \bar\partial_A\bar\partial_A^*\psi
  =
  \frac{1}{2}\nabla_{A_n\otimes A}^*\nabla_{A_n\otimes A}\psi
  + \frac{i}{2}\Lambda_\omega (F_{A_{K^*}} + F_A)\psi,
  \quad\text{for all } \psi \in \Omega^{0,n}(E).
\end{equation}
\end{lem}

\begin{rmk}[Generalizations of Lemma \ref{lem:Kotschick_2-5_complex_kaehler_rank-r_holomorphic}]
\label{rmk:Kotschick_2-5_raw_bullet_complex_kaehler}  
As noted, Kotschick asserts in \cite[Lemma 2.5, Equation (13), p. 204]{KotschickSW} that the identity \eqref{eq:Kotschick_section_2_eq_13_rank-r_kaehler} in Lemma \ref{lem:Kotschick_2-5_complex_kaehler_rank-r_holomorphic} holds when $n=2$ and $(X,g,J,\omega)$ is symplectic. Our hypotheses that $(X,g,J,\omega)$ is complex K\"ahler and $\bar\partial_E$ is integrable are inherited from Proposition \ref{prop:Berline_Getzler_Vergne_3-71}, but we would expect some version of Proposition \ref{prop:Berline_Getzler_Vergne_3-71} and thus Lemma \ref{lem:Kotschick_2-5_complex_kaehler_rank-r_holomorphic} to hold when those hypotheses are relaxed to allow $(X,g,J,\omega)$ to be almost K\"ahler or almost Hermitian and allow $\bar\partial_E$ to be an arbitrary $(0,1)$-connection.

Moroianu \cite[Section 20.1, Theorem 20.2, p. 137]{Moroianu_2007} provides a generalization of the identity \eqref{eq:Kotschick_section_2_eq_13_rank-r_kaehler} to sections of $\Lambda^{p,q}(E)$, for arbitrary $0\leq p, q \leq n$, but the curvature term is not computed explicitly. A generalization of the identity \eqref{eq:Kotschick_section_2_eq_13_rank-r_kaehler} to sections of $\Lambda^{0,q}(E)$ appears in Li \cite[Equations (1) and (2), p. 623]{Li_2010}, based Ma and Marinescu \cite[Remark 1.4.8, Equation (1.4.31), p. 34]{Ma_Marinescu_holomorphic_morse_inequalities_bergman_kernels} and a further generalization to sections of $\Lambda^{p,q}(E)$ appears in Li \cite[Equation (1) and Theorem 3.1, p. 628]{Li_2010}, again for arbitrary $0\leq p, q \leq n$. The identity in Li \cite[Equations (1) and (2), p. 623]{Li_2010} is proved by Liu in \cite[Chapter 3, Theorem 3.1.5]{Liu_complex_manifold_and_Kaehler_geometry_2018}. According to Li \cite[Equation (1), p. 623]{Li_2010}, one has
\begin{equation}
  \label{eq:Li_2010_1}
    \left(\bar\partial_A\bar\partial_A^* + \bar\partial_A^*\bar\partial_A\right)\psi
    =
    \frac{1}{2}\nabla_{A_{p,q}\otimes A}^*\nabla_{A_{p,q}\otimes A}\psi + \frac{1}{2} W_{p,q}(E)\psi,
    \quad\text{for all } \psi \in \Omega^{0,q}(E),
\end{equation}
and from Li \cite[Equation (2), p. 623]{Li_2010}, when $p=0$, one has
\begin{equation}
  \label{eq:Li_2010_2}
  W_{0,q}(E)
  :=
  \sum_{j,k=1}^n
  2\left(F_A + \frac{1}{2}\tr R^{TX^{1,0}}\right)(Z_j,\bar Z_k)\,\eps(\bar Z^k)\iota(\bar Z_j)
  - \sum_{j=1}^n F_A(Z_j,\bar Z_j),
\end{equation}
where $\{Z_1,\ldots,Z_n\}$ is a local orthonormal frame for $TX^{1,0}$ and $\{Z^1,\ldots,Z^n\}$ is the corresponding local orthonormal frame for $(TX^{1,0})^*$ and $A_{p,q}$ is the induced connection on $\Lambda^{p,q}(X)$. According to Liu \cite[Chapter 3, Equation (3.1.17)]{Liu_complex_manifold_and_Kaehler_geometry_2018}, one has
\[
  \tr R^{TX^{1,0}} = F_{A_{K^*}},
\]
where $K^* = \wedge^n(TX^{1,0})$, the anti-canonical line bundle over $X$. Hence, if $q=n$, the preceding identity and calculations in the proof of Lemma \ref{lem:Kotschick_2-5_complex_kaehler_rank-r_holomorphic} reveal that
\begin{align*}
  W_{0,n}(E)
  &=
  \sum_{j=1}^n
  2\left(F_A + \frac{1}{2}\tr R^{TX^{1,0}}\right)(Z_j,\bar Z_j)
  - \sum_{j=1}^n F_A(Z_j,\bar Z_j)
  \\
  &=
  \sum_{j=1}^n F_A(Z_j,\bar Z_j) + \sum_{j=1}^n F_{A_{K^*}}(Z_j,\bar Z_j)
  \\
  &=
  i\Lambda_\omega\left(F_{A_{K^*}} + F_A\right).
\end{align*}
After substituting the preceding expression for $W_{0,n}(E)$ into \eqref{eq:Li_2010_1} and noting that $\bar\partial\psi = 0$, we obtain \eqref{eq:Kotschick_section_2_eq_13_rank-r_kaehler}.

Taking the other extreme in the expression \eqref{eq:Li_2010_2} for $W_{0,q}(E)$ and choosing $q=0$, we find that
%PF8-19-2024 Difference between the q=0 and q=n cases is unclear.
\[
  W_{0,0}(E) = - \sum_{j=1}^n F_A(Z_j,\bar Z_j) = -i\Lambda_\omega F_A.
\]
After substituting the preceding expression for $W_{0,0}(E)$ into \eqref{eq:Li_2010_1} and noting that $\bar\partial^*\psi = 0$, we obtain \eqref{eq:Kotschick_section_2_eq_12}.
\qed
\end{rmk}  

\begin{proof}[Proof of Lemma \ref{lem:Kotschick_2-5_complex_kaehler_rank-r_holomorphic}]
Observe that by restricting the identity \eqref{eq:Kotschick_section_2_eq_13_raw_bullet} in Lemma \ref{lem:Kotschick_2-5_raw_bullet} to the summand $\Lambda^{0,q}(E) = \Lambda^{0,q}(X)\otimes E$, we obtain
\[
  \bar\partial_{A_q\otimes A}^*\bar\partial_{A_q\otimes A}\psi
  =
  \frac{1}{2}\nabla_{A_q\otimes A}^*\nabla_{A_q\otimes A}\psi
  - \frac{i}{2}\Lambda_\omega F_{A_q\otimes A}\psi,
  \text{for all } \psi \in \Omega^{0,q}(E) = \Omega^0(\Lambda^{0,q}(X)\otimes E).
\]
Also, by restricting the identity \eqref{eq:Berline_Getzler_Vergne_prop_3-71} in Proposition \ref{prop:Berline_Getzler_Vergne_3-71} to the summand $\Lambda^{0,q}(E)$, we obtain
\begin{multline*}
  \left(\bar\partial_A\bar\partial_A^* +  \bar\partial_A\bar\partial_A^*\right)\psi
  =
  \bar\partial_{A_q\otimes A}^*\bar\partial_{A_q\otimes A}\psi
  +
  \sum_{\alpha,\beta=1}^n F_{A_{K^*}\otimes A}
  \left(\frac{\partial}{\partial z^\beta},\frac{\partial}{\partial\bar z^\alpha}\right)
  \eps(d\bar z^\alpha)\left(\frac{\partial}{\partial\bar z^\beta}\right)\psi,
  \\
  \quad\text{for all } \psi \in \Omega^{0,q}(E).
\end{multline*}
Hence, by combining the preceding two identities we obtain
\begin{multline*}
  \left(\bar\partial_A\bar\partial_A^* +  \bar\partial_A^*\bar\partial_A\right)\psi
  =
  \frac{1}{2}\nabla_{A_q\otimes A}^*\nabla_{A_q\otimes A}\psi
  - \frac{i}{2}\Lambda_\omega F_{A_q\otimes A}\psi
  \\
  \sum_{\alpha,\beta=1}^n F_{A_{K^*}\otimes A}
  \left(\frac{\partial}{\partial z^\beta},\frac{\partial}{\partial\bar z^\alpha}\right)
  \eps(d\bar z^\alpha)\left(\frac{\partial}{\partial\bar z^\beta}\right)\psi,
  \quad\text{for all } \psi \in \Omega^{0,q}(E).
\end{multline*}
When $q=n$, then $\bar\partial_A\psi = 0$ and $A_q = A_n = A_{K^*}$, and so the desired conclusion \eqref{eq:Kotschick_section_2_eq_13_rank-r_kaehler} would follow from the preceding equality if we could prove the claimed identity below:
\begin{equation}
  \label{eq:Berline_Getzler_Vergne_prop_3-71_curvature_term_q=n}
  \sum_{\alpha,\beta=1}^n F_{A_n\otimes A}
  \left(\frac{\partial}{\partial z^\beta},\frac{\partial}{\partial\bar z^\alpha}\right)
  \eps(d\bar z^\alpha)\left(\frac{\partial}{\partial\bar z^\beta}\right)\psi
  =
  i\Lambda_\omega F_{A_{K^*}\otimes A}\psi.
\end{equation}
%COMMENT The calculations in Freed and Uhlenbeck \cite[Appendix C, p. 174]{FU} are closest to what we need to prove the preceding identity.
By applying Kobayashi \cite[Equations (4.1.1) and (4.1.2), p. 92]{Kobayashi_differential_geometry_complex_vector_bundles}, we obtain
\[
  i\Lambda_\omega F_{A_q\otimes A}
    =
  \sum_{\alpha,\beta=1}^n g^{\alpha\bar\beta}(F_{A_q\otimes A})_{\alpha\bar\beta}
  =
  \sum_{\alpha,\beta=1}^n g^{\alpha\bar\beta}F_{A_q\otimes A}
  \left(\frac{\partial}{\partial z^\alpha},\frac{\partial}{\partial\bar z^\beta}\right).
\]
With respect to complex normal coordinates\footnote{\url{https://math.stackexchange.com/questions/3801057/complex-normal-coordinates-in-kaehler-manifolds}} $\{z^\alpha\}$ around a point $x$ in a K\"ahler manifold $X$, we have $g_{\alpha\bar\beta}(x) = \delta_{\alpha\beta}$ (and also $dg_{\alpha\bar\beta}(x) = 0$ and $(\partial^2g_{\alpha\bar\beta}/\partial z_\gamma\partial z_\delta)(x) = 0$). Hence, at such a point $x \in X$, we have
\[
  i\Lambda_\omega F_{A_q\otimes A}
  =
  \sum_{\alpha,\beta=1}^n \delta^{\alpha\beta}F_{A_q\otimes A}
  \left(\frac{\partial}{\partial z^\alpha},\frac{\partial}{\partial\bar z^\beta}\right)
  =
  \sum_{\alpha=1}^n F_{A_q\otimes A}
  \left(\frac{\partial}{\partial z^\alpha},\frac{\partial}{\partial\bar z^\alpha}\right).
\]
(According to Kobayashi \cite[Section 4.1, p. 92, second displayed equation or Equation (1.4.15), p. 12]{Kobayashi_differential_geometry_complex_vector_bundles}, when $F_{A_q\otimes A}^{2,0} = 0$ and $F_{A_q\otimes A}^{0,2} = 0$ we have
\[
  F_{A_q\otimes A}
  =
  \sum_{\alpha,\beta} F_{A_q\otimes A}
  \left(\frac{\partial}{\partial z^\alpha},\frac{\partial}{\partial\bar z^\beta}\right)
  dz^\alpha \wedge d\bar z^\beta.
\]
Additional terms appear when $F_{A_q\otimes A}^{0,2} \neq 0$, with coefficients of $dz^\alpha \wedge dz^\beta$ and $d\barz^\alpha \wedge d\bar z^\beta$.)

When $q=n$ and $\psi \in \Omega^{0,n}(E)$, then $\psi = \psi_n\, d\bar z^1\wedge \cdots \wedge d\bar z^n$ in local coordinates, where $\psi_n \in \Omega^0(E)$, and thus we see that
\[
  \eps(d\bar z^\alpha)\iota\left(\frac{\partial}{\partial\bar z^\beta}\right)\psi
  =
  \begin{cases}
    \psi, &\text{if } \alpha = \beta,
    \\
    0, &\text{if } \alpha \neq \beta.
  \end{cases}
\]
Consequently, for all $\psi \in \Omega^{0,n}(E)$ and noting that $A_n = A_{K^*}$, we obtain
\begin{align*}
  \sum_{\alpha,\beta=1}^n F_{A_n\otimes A}
  \left(\frac{\partial}{\partial z^\beta},\frac{\partial}{\partial\bar z^\alpha}\right)
  \eps(d\bar z^\alpha)\left(\frac{\partial}{\partial\bar z^\beta}\right)\psi
  &=
  \sum_{\alpha=1}^n F_{A_n\otimes A}
    \left(\frac{\partial}{\partial z^\alpha},\frac{\partial}{\partial\bar z^\alpha}\right)
    \psi
  \\
  &= i\Lambda_\omega F_{A_n\otimes A}\psi
  \\
  &= i\Lambda_\omega F_{A_{K^*}\otimes A}\psi, 
\end{align*}
which proves the claimed identity \eqref{eq:Berline_Getzler_Vergne_prop_3-71_curvature_term_q=n} and hence completes the proof of Lemma \ref{lem:Kotschick_2-5_complex_kaehler_rank-r_holomorphic}.
%COMMENT The style of writing the curvature term in \eqref{eq:Berline_Getzler_Vergne_prop_3-71} is similar to that of Lawson \cite[Appendix II, Theorem II.1, p. 93]{Lawson} and Freed and Uhlenbeck \cite[Equation (C.11), p. 172]{FU} and Parker \cite[Section 1]{ParkerGauge}. See also Atiyah, Hitchin, and Singer \cite{AHS}. Freed and Uhlenbeck apply \cite[Appendix C]{FU} to prove the Bochner--Weitzenb\"ock identities \cite[Chapter 6, Proof of Equations (6.25) and (6.26), p. 94]{FU}. 
\end{proof}  

\subsection[Bochner--Weitzenb\"ock identity for the $\bar\partial_A$-Laplacian on $(0,q)$-forms]{Bochner--Weitzenb\"ock identity for the $\bar\partial_A$-Laplacian on $(0,q)$-forms with values in a Hermitian vector bundle over a symplectic manifold}
\label{subsec:Bochner-Weitzenbock_formulae_anti-holomorphic_spinors_symplectic}
Our second approach to proving the Bochner--Weitzenb\"ock identities for $\Delta_{\bar\partial_A}$ is implied by Taubes' remarks in the paragraphs preceding his \cite[Section 2 (b), Equation (2.3), p. 854 and Equation (2.4), p. 855]{TauSWGromov}, namely to derive them as consequences of the well-known 

\begin{lem}[Bochner--Weitzenb\"ock identity for coupled Dirac operator over a \spinc manifold]
\label{lem:Bochner-Weitzenbock_identity_coupled_Dirac_operator}  
(See Atiyah, Hitchin, and Singer \cite[Proof of Theorem 6.1, p. 445]{AHS}, Berline, Getzler, and Vergne \cite[Theorem 3.52, p. 126]{BerlineGetzlerVergne}, Besse \cite[Chapter 1, Section I, Equation (1.152), p. 56]{Besse_1987}, Feehan and Leness \cite[Lemma 4.1, p. 334]{FL1}, Gromov and Lawson \cite[Proposition 2.5, p. 111]{Gromov_Lawson_1983} Lawson and Michelsohn \cite[Chapter II, Theorem 8.17, p. 164 and Appendix D, Theorem D.12, p. 399]{LM}, Nicolaescu \cite[Equation (1.3.11), p. 51 and Equation (1.3.16), p. 52]{NicolaescuSWNotes}, Salamon \cite[Section 6.4, Theorem 6.19, p. 205]{SalamonSWBook}, and Taubes \cite[Section 2 (a), Equation (2.1), p. 854]{TauSWGromov}.)
Let $(X,g)$ be an oriented, smooth Riemannian manifold, $(\rho,W)$ be a spin${}^c$ structure on $X$, and $A$ be a smooth, unitary connection on a smooth, Hermitian vector bundle $(E,H)$ over $X$. If $\Phi \in \Omega^0(W\otimes E)$, then 
\begin{equation}
  \label{eq:Bochner-Weitzenbock_formula_coupled_Dirac_operator}
  D_AD_A\Phi = \nabla_{A_W\otimes A}^*\nabla_{A_W\otimes A}\Phi + \frac{R}{4}\Phi 
  + \frac{1}{2}\rho(F_{A_L})\Phi
  + \rho(F_A)\Phi,
\end{equation}
where $A_L$ is the induced unitary connection on the Hermitian line bundle $L$ associated to the spin${}^c$ structure on $X$ and $F_{A_L}$ is its curvature and $R$ is the scalar curvature of the Riemannian metric, $g$. If $\dim X = 4$ and $\Phi \in \Omega^0(W^\pm\otimes E)$, then $\rho(F_{A_L}^\mp)\Phi = 0$ and $\rho(F_A^\mp)\Phi = 0$ in \eqref{eq:Bochner-Weitzenbock_formula_coupled_Dirac_operator}.
\end{lem}

%PF8-12-2024 Proof updated
\begin{proof}
We adapt the proof of Feehan and Leness \cite[Lemma 4.1, p. 334]{FL1} from dimension $4$ to the case of arbitrary dimensions, recalling that our proof of \cite[Lemma 4.1, p. 334]{FL1} in turn applied the argument employed by Lawson and Michelsohn in their proof of \cite[Appendix D, Theorem D.12, p. 399]{LM}.
  
Over a sufficiently small open subset $U\subset X$, we can write $W\restriction U \cong S\otimes L^{1/2}$ where $S$ is a spinor bundle for the local spin structure, $L$ is the
%PF8-12-2024 How defined?
Hermitian line bundle associated to the \spinc structure over $X$, so $c_1(L) \equiv w_2(X) \pmod{2}$, and $L^{1/2}$ is the square root of $L$ over $U$ that may not exist over $X$. (When $\dim X = 2n$, then by \cite[Appendix D, p. 395]{LM} there is a splitting $W = W^+\oplus W^-$ for $W^\pm = (1\pm\omega_\CC)W$, where $\omega_\CC = i^{n/2}\,d\vol_g$ and $d\vol_g$ is the volume form defined by the Riemannian metric. If $\dim X$ is even, then
%PF8-12-2024 Need reference for this
$L = \det W^+$ and if $\dim X$ is odd, then $L = \det W$.)

Let $A_{L^{1/2}}$ be a unitary connection on $L^{1/2}\restriction U$
%PF8-12-2024 How is this defined?
inducing $A_L$ on $L\restriction U$, so that $A_L = A_{L^{1/2}}\otimes A_{L^{1/2}}$ on $L = L^{1/2}\otimes L^{1/2}$. We apply the Bochner--Weitzenb\"ock identity in Lawson and Michelsohn \cite[Theorem II.8.17, p. 164]{LM} to the
connection $A_{L^{1/2}}\otimes A$ on $L^{1/2}\otimes E \restriction U$ to give
\begin{equation}
\label{eq:BW_Local}
   D_AD_A\Phi
   =
   \nabla_{A_W\otimes A}^*\nabla_{A_W\otimes A}\Phi + \frac{R}{4}\Phi + \rho(F_{A^{1/2}\otimes A})\Phi,
   \quad\text{for all } \Phi \in \Omega^0(U,W\otimes E).
\end{equation}
Because $A_{L^{1/2}}$ induces $A_L$ on $L \restriction U$, we have $F_{A_L} = 2F_{A_{L^{1/2}}}$ over $U$ and thus
\[
  F_{A^{1/2}\otimes A}
  =
  \frac{1}{2}F_{A_L}\otimes\id_E + \id_{L^{1/2}}\otimes F_A \quad\text{over } U.
\]
Combining the preceding equality with \eqref{eq:BW_Local} yields the Bochner--Weitzenb\"ock identity \eqref{eq:Bochner-Weitzenbock_formula_coupled_Dirac_operator} for all $\Phi \in \Omega^0(U,W\otimes E)$ and hence for all $\Phi \in \Omega^0(W\otimes E)$.
\end{proof}

We now specialize Lemma \ref{lem:Bochner-Weitzenbock_identity_coupled_Dirac_operator} to the case of the canonical \spinc structure $(\rho_\can,W_\can)$ on a symplectic manifold, following Salamon's expression for the Dirac operator on a symplectic manifold \cite[Section 6.3, Theorem 6.17, p. 203]{SalamonSWBook}.

%{\color{blue}[FLAG: Add explanation of $\rho(F_{A_L})\Phi$ term for canonical spinc structure?]}
%PF8-12-2024 Yes, this was one of your to-do's
\begin{lem}[Bochner--Weitzenb\"ock identity for the coupled Dirac operator over a symplectic manifold]
\label{lem:Bochner-Weitzenbock_identity_coupled_Dirac_operator_symplectic_manifold}  
Let $(X,g,J,\omega)$ be an almost K\"ahler $2n$-manifold with its canonical spin${}^c$ structure $(\rho_\can,W_\can)$ and $A$ be a smooth, unitary connection on a smooth, Hermitian vector bundle $(E,H)$ over $X$. Let $A_\can$ denote the unitary connection on the Hermitian vector bundle $\Lambda^{0,\bullet}(X) = \oplus_{q=0}^n\Lambda^{0,q}(X)$ induced by the Hermitian metric $h = g + i\omega$ on $X$, and $A_\can\otimes A$ denote the induced unitary connection on $\Lambda^{0,\bullet}(X)\otimes E$. Then for all $\Phi \in \Omega^0(\Lambda^{0,\bullet}(X)\otimes E) = \Omega^{0,\bullet}(E)$,
% TL8-2-2024: Do we need to add the analogue of the   \frac{1}{2}\rho(F_{A_L})\Phi term?
%PF8-12-2024 Yes, this was one of your to-do's
\begin{equation}
  \label{eq:Bochner-Weitzenbock_formula_coupled_Dirac_operator_symplectic_manifold}
  \Delta_{\bar\partial_A}\Phi
  + \frac{1}{2}\left(\mu\bar\partial_A + \partial_A^*\bar\mu^*\right)\Phi
  = \frac{1}{2}\nabla_{A_\can\otimes A}^*\nabla_{A_\can\otimes A}\Phi + \frac{R}{8}\Phi
  + \frac{1}{4}\rho_\can(F_{A_{K^*}})\Phi + \frac{1}{2}\rho_\can(F_A)\Phi.
\end{equation}
\end{lem}

\begin{proof}
We observe that
\begin{align*}
  D_AD_A\Phi
  &= 2\left(\bar\partial_A + \bar\partial_A^*\right)\left(\bar\partial_A + \bar\partial_A^*\right)\Phi
  \\
  &= 2\left(\bar\partial_A\bar\partial_A + \bar\partial_A^*\bar\partial_A
    + \bar\partial_A\bar\partial_A^* + \bar\partial_A^*\bar\partial_A^*\right)\Phi
  \\
  &= 2\left((F_A^{2,0} + \mu\bar\partial_A) + \bar\partial_A^*\bar\partial_A
    + \bar\partial_A\bar\partial_A^* + ((F_A^{0,2} + \bar\mu\partial_A)^*\right)\Phi.
\end{align*}
Note that, for all $\Phi, \Psi \in \Omega^{0,\bullet}(E)$,
\[
  \left((F_A^{0,2})^*\Phi, \Psi\right)_{L^2(X)}
  = \left(\Phi, F_A^{0,2}\Psi\right)_{L^2(X)}
  = \left((F_A^{0,2})^\dagger\Phi, \Psi\right)_{L^2(X)}
  = \left(-F_A^{2,0}\Phi, \Psi\right)_{L^2(X)}  
\]
and so $(F_A^{0,2})^* = -F_A^{2,0}$. Therefore,
\[
  D_AD_A\Phi
  =
  2\left(\bar\partial_A^*\bar\partial_A + \bar\partial_A\bar\partial_A^*\right)\Phi
  + \left(\mu\bar\partial_A + \partial_A^*\bar\mu^*\right)\Phi.
\]
On the other hand, the Bochner--Weitzenb\"ock identity \eqref{eq:Bochner-Weitzenbock_formula_coupled_Dirac_operator} gives
%TL8-2-2024: Do we need to add the analogue of the   \frac{1}{2}\rho(F_{A_L})\Phi term?  
%{\color{blue}[FLAG: Do we need to add the analogue of the   $\frac{1}{2}\rho(F_{A_L})\Phi$ term?]}
%PF8-12-2024 Yes, this was one of your to-do's
\[
  D_AD_A\Phi
  = \nabla_{A_\can\otimes A}^*\nabla_{A_\can\otimes A}\Phi + \frac{R}{4}\Phi
  + \frac{1}{2}\rho_\can(F_{A_L})\Phi + \rho_\can(F_A)\Phi.
\]
By combining the preceding two identities, we obtain
\[
\Delta_{\bar\partial_A}\Phi
  + \frac{1}{2}\left(\mu\bar\partial_A + \partial_A^*\bar\mu^*\right)\Phi
  = \frac{1}{2}\nabla_{A_\can\otimes A}^*\nabla_{A_\can\otimes A}\Phi + \frac{R}{8}\Phi
  + \frac{1}{4}\rho_\can(F_{A_L})\Phi + \frac{1}{2}\rho_\can(F_A)\Phi.
\]
For the complex vector bundle $(T^*X,J^*) \cong T^*X^{1,0}$, we recall that the \emph{canonical line bundle} is $K := \wedge^n(T^*X^{1,0}) = \Lambda^{n,0}(X)$ (see \cite[Section 5.3, p. 167]{SalamonSWBook}) and the \emph{anti-canonical line bundle} is the dual complex line bundle, $K^* \cong \wedge^n(T^*X^{0,1}) = \Lambda^{0,n}(X)$ (see \cite[Section 5.3, p. 168]{SalamonSWBook}). By Salamon \cite[Section 5.3, Corollary 5.21, p. 168]{SalamonSWBook}, the Hermitian line bundle associated to the canonical \spinc structure $(\rho_\can,W_\can)$ is
\[
  L \cong K^*.
\]
The identity \eqref{eq:Bochner-Weitzenbock_formula_coupled_Dirac_operator_symplectic_manifold} now follows from the preceding expressions for $\Delta_{\bar\partial_A}\Phi$ and $L$, respectively. This completes the proof of Lemma \ref{lem:Bochner-Weitzenbock_identity_coupled_Dirac_operator_symplectic_manifold}.
\end{proof}

\begin{rmk}[Special cases of the Bochner--Weitzenb\"ock identity \eqref{eq:Bochner-Weitzenbock_formula_coupled_Dirac_operator_symplectic_manifold}]
\label{rmk:Bochner-Weitzenbock_identity_coupled_Dirac_operator_symplectic_manifold_special_cases}
Note that $\mu:\Omega^{p,q}(X) \to \Omega^{p+2,q-1}(X)$ and $\bar\partial_A:\Omega^{p,q}(X) \to \Omega^{p,q+1}(X)$, so
$\mu\bar\partial_A:\Omega^{p,q}(X) \to \Omega^{p+2,q}(X)$. If $\Phi \in \Omega^{0,n}(E)$, then $\bar\partial_A\Phi = 0 \in \Omega^{0,n+1}(E)$ and so $\mu\bar\partial_A\Phi = 0 \in \Omega^{2,n}(E)$ and the identity \eqref{eq:Bochner-Weitzenbock_formula_coupled_Dirac_operator_symplectic_manifold} reduces to
\begin{multline}
  \label{eq:Bochner-Weitzenbock_formula_coupled_Dirac_operator_symplectic_0n}
  \Delta_{\bar\partial_A}\Phi
  + \frac{1}{2}\partial_A^*\bar\mu^*\Phi
  = \frac{1}{2}\nabla_{A_\can\otimes A}^*\nabla_{A_\can\otimes A}\Phi + \frac{R}{8}\Phi
  + \frac{1}{4}\rho_\can(F_{A_{K^*}})\Phi + \frac{1}{2}\rho_\can(F_A)\Phi,
  \\
  \text{for all } \Phi \in \Omega^{0,n}(E).
\end{multline}
Noting that $\bar\mu:\Omega^{p,q}(X) \to \Omega^{p-1,q+2}(X)$ and $\bar\mu^*:\Omega^{p,q}(X) \to \Omega^{p+1,q-2}(X)$ and $\partial_A^*\bar\mu^*:\Omega^{p,q}(X) \to \Omega^{p,q-2}(X)$, the identity \eqref{eq:Bochner-Weitzenbock_formula_coupled_Dirac_operator_symplectic_manifold} simplifies when $q=0$ or $1$ to
\begin{multline}
  \label{eq:Bochner-Weitzenbock_formula_coupled_Dirac_operator_symplectic_p0_or_p1}
  \Delta_{\bar\partial_A}\Phi
  + \frac{1}{2}\mu\bar\partial_A\Phi
  = \frac{1}{2}\nabla_{A_\can\otimes A}^*\nabla_{A_\can\otimes A}\Phi + \frac{R}{8}\Phi
  + \frac{1}{4}\rho_\can(F_{A_{K^*}})\Phi + \frac{1}{2}\rho_\can(F_A)\Phi,
  \\
  \text{for all } \Phi \in \Omega^{p,q}(E) \text{ with } q = 0 \text{ or } 1.
\end{multline}
Consequently, when $n=1$ we obtain
\begin{equation}
  \label{eq:Bochner-Weitzenbock_formula_coupled_Dirac_operator_symplectic_00_or_01}
  \Delta_{\bar\partial_A}\Phi
  = \frac{1}{2}\nabla_{A_\can\otimes A}^*\nabla_{A_\can\otimes A}\Phi + \frac{R}{8}\Phi
  + \frac{1}{4}\rho_\can(F_{A_{K^*}})\Phi + \frac{1}{2}\rho_\can(F_A)\Phi,
  \quad\text{for all } \Phi \in \Omega^{0,1}(E),
\end{equation}
from \eqref{eq:Bochner-Weitzenbock_formula_coupled_Dirac_operator_symplectic_0n} and \eqref{eq:Bochner-Weitzenbock_formula_coupled_Dirac_operator_symplectic_p0_or_p1}.
\qed
\end{rmk}

We now use Lemma \ref{lem:Bochner-Weitzenbock_identity_coupled_Dirac_operator_symplectic_manifold} to give
\begin{inparaenum}[\itshape i\upshape)]
\item a second proof of Lemma \ref{lem:Kotschick_2-4}, and
\item a generalization of Lemma \ref{lem:Kotschick_2-5_complex_kaehler_rank-r_holomorphic}, when $n=2$, to the case of smooth Hermitian vector bundles with arbitrary $(0,1)$-connections over symplectic $2n$-manifolds when $n=2$ or $3$.
\end{inparaenum}  
We begin with the

\begin{proof}[Second proof of Equation \eqref{eq:Kotschick_section_2_eq_12} in Lemma \ref{lem:Kotschick_2-4}]
For arbitrary $n \geq 1$, we may write
\[
  \Phi
  = (\varphi_0,\ldots,\varphi_{\lfloor n/2\rfloor})
  \in
  \Omega^0(E)\oplus\cdots\oplus\Omega^{0,\lfloor n/2\rfloor}(E),
\]
where $\lfloor n/2\rfloor$ is the greatest integer less than or equal to $n/2$. For simplicity of notation, we shall assume that $n = 2$ or $3$ and write $\Phi = (\varphi,\psi) \in \Omega^0(E)\oplus\Omega^{0,2}(E)$.
%PF7-25-2024 Some of the below is a duplicate of what appears the earlier subsection.
From Feehan and Leness 
%TL12-4-2025: Updated and removed page references
\cite[Lemma 8.3.4 and Corollary 8.3.5]{Feehan_Leness_introduction_virtual_morse_theory_so3_monopoles}, we have
\begin{subequations}
  \label{eq:Canonical_spinc_structure_almost_Hermitian_4-manifolds}
  \begin{align}
    \label{eq:Canonical_spinc_structure_almost_Hermitian_4-manifolds_iomega}
  \rho_\can(i\omega)(\sigma,\tau)
  &= 2(\sigma,-\tau) \in \Omega^0(X)\oplus\Omega^{0,2}(X),
    \\
    \label{eq:Canonical_spinc_structure_almost_Hermitian_4-manifolds_20_forms}
  \rho_\can(\bar\beta)(\sigma,\tau)
  &=
     \left(-2\langle\tau,\beta\rangle_{\Lambda^{0,2}(X)}, 0\right) \in \Omega^0(X)\oplus\Omega^{0,2}(X),
    \\
    \label{eq:Canonical_spinc_structure_almost_Hermitian_4-manifolds_02_forms}
  \rho_\can(\beta)(\sigma,\tau)
  &=
    \left(0, 2\sigma\beta\right) \in \Omega^0(X)\oplus\Omega^{0,2}(X),
\end{align}
\end{subequations}
for all $(\sigma,\tau) \in \Omega^0(X) \oplus \Omega^{0,2}(X)$. We have
\[
  \rho_\can(F_A)\Phi = \rho_\can(F_A^+)\Phi,
  \quad\text{for all } \Phi \in \Omega^0\left((\Lambda^0\oplus\Lambda^{0,2})(X)\otimes E\right),
\]
where we write
\begin{align*}
  \Omega^+(\fu(E)) &:= \Omega^{2,0}(\fu(E)) \oplus \Omega^0(\fu(E))\,\omega \oplus \Omega^{0,2}(\fu(E)),
  \\
  F_A^+ &:= F_A^{2,0} + F_A^\omega + F_A^{0,2} \in \Omega^+(\fu(E)),
\end{align*}
and, noting that $|\omega|^2\,d\vol_\omega = \omega\wedge \star\omega = \omega\wedge\omega = 2d\vol_\omega$ and so $|\omega| = \sqrt{2}$, we have
\[
  F_A^\omega = \frac{1}{2}\langle F_A,\omega\rangle_{\Lambda^{1,1}(X)}\,\omega,
\]
and, using $\Lambda_\omega\omega = 2$,
\[
  \Lambda_\omega F_A = \Lambda_\omega F_A^\omega
  = \frac{1}{2}\langle F_A,\omega\rangle_{\Lambda^{1,1}(X)} \Lambda_\omega\omega
  = \langle F_A,\omega\rangle_{\Lambda^{1,1}(X)},
\]
so that
\begin{equation}
  \label{eq:FAomega_and_LambdaFA}
  F_A^\omega = \frac{1}{2}(\Lambda_\omega F_A)\,\omega.
\end{equation}
For $(\varphi,\psi) \in \Omega^0(E)\oplus\Omega^{0,2}(E)$, Equation \eqref{eq:Canonical_spinc_structure_almost_Hermitian_4-manifolds_iomega} thus gives
\begin{align*}
  \rho_\can(F_A^\omega)(\varphi,\psi)
  &= \frac{1}{2}(\Lambda_\omega F_A)\rho_\can(\omega)(\varphi,\psi)
  \\
  &= -\frac{i}{2}(\Lambda_\omega F_A)\rho_\can(i\omega)(\varphi,\psi)
  \\
  &= -i\Lambda_\omega F_A(\varphi,-\psi)
  \\
  &= i\Lambda_\omega F_A(-\varphi,\psi).
\end{align*}
Next, we observe that
\begin{align*}
  \rho_\can(F_A^{2,0})(\varphi,\psi) + \rho_\can(F_A^{0,2})(\varphi,\psi)
  &=
    (-2\langle\psi, \overline{F_A^{2,0}}\rangle_{\Lambda^{0,2}(X)} + (0, 2F_A^{0,2}\varphi)
  \\
  %PF7-24-2024 Check!
  &=
    2\left(\langle F_A^{0,2}, \psi\rangle_{\Lambda^{0,2}(X)}, F_A^{0,2}\varphi\right).
\end{align*}
Of course, the analogous identities hold for the curvature term $\rho_\can(F_{A_{K^*}})$ in Equation \eqref{eq:Bochner-Weitzenbock_formula_coupled_Dirac_operator_symplectic_manifold}.

Because (see Remark \ref{rmk:Bochner-Weitzenbock_identity_coupled_Dirac_operator_symplectic_manifold_special_cases})
\[
  \left(\mu\bar\partial_A + \partial_A^*\bar\mu^*\right)(\varphi,\psi)
  =
  \mu\bar\partial_A\varphi + \partial_A^*\bar\mu^*\psi + \mu\bar\partial_A\psi 
  \in
  \Omega^{2,0}(E) \oplus \Omega^0(E) \oplus \Omega^{2,2}(E),
\]  
we see that Equation \eqref{eq:Bochner-Weitzenbock_formula_coupled_Dirac_operator_symplectic_manifold} simplifies to
\begin{multline}
  \label{eq:Bochner-Weitzenbock_formula_coupled_Dirac_operator_symplectic_manifold_n=2_or_3}
  \Delta_{\bar\partial_A}(\varphi,\psi)
  + \frac{1}{2}\left(\mu\bar\partial_A\varphi + \partial_A^*\bar\mu^*\psi + \mu\bar\partial_A\psi\right)
  =
  \frac{1}{2}\nabla_{A_\can\otimes A}^*\nabla_{A_\can\otimes A}(\varphi,\psi) + \frac{R}{8}(\varphi,\psi)
  \\
  + \frac{1}{2}\left(\langle F_{A_{K^*}}^{0,2}, \psi\rangle_{\Lambda^{0,2}(X)}, F_{A_{K^*}}^{0,2}\varphi\right)
  + \frac{i}{4}\Lambda_\omega F_{A_{K^*}}(-\varphi,\psi).
  \\
  + \left(\langle F_A^{0,2}, \psi\rangle_{\Lambda^{0,2}(X)}, F_A^{0,2}\varphi\right)
  + \frac{i}{2}\Lambda_\omega F_A(-\varphi,\psi).
\end{multline}
We project
%PF8-16-2024 Check this projection
Equation \eqref{eq:Bochner-Weitzenbock_formula_coupled_Dirac_operator_symplectic_manifold_n=2_or_3} onto its $\Omega^0(E)$ component and choose $\psi = 0$ to give
\begin{equation}
  \label{eq:Bochner-Weitzenbock_formula_coupled_Dirac_operator_symplectic_manifold_n=2_or_3_varphi}
  \bar\partial_A^*\bar\partial_A\varphi
  = \frac{1}{2}\nabla_A^*\nabla_A\varphi
  + \frac{R}{8}\varphi - \frac{i}{4}\Lambda_\omega F_{A_{K^*}}\varphi - \frac{i}{2}\Lambda_\omega F_A\varphi.
\end{equation}
By applying the curvature identity \eqref{eq:Curvature_anti-canonical_line_bundle_and_scalar_curvature} in the forthcoming Lemma \ref{lem:Curvature_anti-canonical_line_bundle_almost_Hermitian_manifold},
\[
  i\Lambda_\omega F_{A_{K^*}} = \frac{R}{2},
\]
we see that Equation \eqref{eq:Bochner-Weitzenbock_formula_coupled_Dirac_operator_symplectic_manifold_n=2_or_3_varphi} agrees with Equation \eqref{eq:Kotschick_section_2_eq_12} for $\bar\partial_A^*\bar\partial_A\varphi$ in Lemma \ref{lem:Kotschick_2-4}. (See Kobayashi \cite[Equation (1.7.16), p. 25]{Kobayashi_differential_geometry_complex_vector_bundles} for a definition of the scalar curvature of a complex Hermitian manifold.) This completes our second proof of Equation \eqref{eq:Kotschick_section_2_eq_12} in Lemma \ref{lem:Kotschick_2-4}.
\end{proof}

Finally, we have the following generalization of Lemma \ref{lem:Kotschick_2-5_complex_kaehler_rank-r_holomorphic}; although we state and prove the result only for $n=2$ or $3$ and $\psi \in \Omega^{0,q}(E)$ with $q=2$, it should not be difficult to generalize the statement and proof to the case of arbitrary $n$ and $q$.

\begin{lem}[Bochner--Kodaira identity for $(0,1)$-connections on Hermitian vector bundles over almost K\"ahler $2n$-manifolds]
\label{lem:Kotschick_2-5_rank-r_symplectic_n=2_or_3}  
(See Kotschick \cite[Lemma 2.5, Equation (13), p. 204]{KotschickSW} for the case where $E$ is a Hermitian line bundle and $n=2$.)
Let $(X,g,J,\omega)$ be an almost K\"ahler $2n$-manifold and $A$ be the Chern connection defined by an $(0,1)$-connection $\bar\partial_E$ on a smooth Hermitian vector bundle $(E,H)$.
%PF7-24-2024 Explain/correct below
For each integer $0\leq q\leq n$, let $A_q$ be the unitary connection on the Hermitian vector bundle $\Lambda^{0,q}(X)$ induced by the Hermitian metric $h = g + i\omega$ on $X$. If $n=2$ or $3$, then
\begin{equation}
  \label{eq:Kotschick_section_2_eq_13_rank-r_symplectic_01-connection_2n-dim_X}
  \left(\bar\partial_A\bar\partial_A^* + \bar\partial_A^*\bar\partial_A\right)\psi
  =
  \frac{1}{2}\nabla_{A_2\otimes A}^*\nabla_{A_2\otimes A}\psi
  + \frac{i}{2}\Lambda_\omega (F_{A_{K^*}} + F_A)\psi,
  \quad\text{for all } \psi \in \Omega^{0,2}(E).
\end{equation}
If $n=2$, then $A_2 = A_{K^*}$ and $\bar\partial_A\psi = 0$ and the preceding identity reduces to
%PF8-16-2024 This case should extend to X almost Hermitian
\begin{equation}
  \label{eq:Kotschick_section_2_eq_13_rank-r_symplectic_01-connection_4-dim_X}
  \bar\partial_A\bar\partial_A^*\psi
  =
  \frac{1}{2}\nabla_{A_2\otimes A}^*\nabla_{A_2\otimes A}\psi
  + \frac{i}{2}\Lambda_\omega (F_{A_2} + F_A)\psi,
  \quad\text{for all } \psi \in \Omega^{0,2}(E).
\end{equation}
\end{lem}

\begin{rmk}[Generalization of Lemma \ref{lem:Kotschick_2-5_rank-r_symplectic_n=2_or_3} to the case of almost Hermitian four-manifolds]
\label{rmk:Kotschick_2-5_rank-r_almost_Hermitian_n=2}   
As we discussed in Remark \ref{rmk:Kaehler_identities_for_almost_Hermitian_manifolds}, some of the key K\"ahler identities continue to hold for almost Hermitian four-manifolds and for that reason we would expect Equation \eqref{eq:Kotschick_section_2_eq_13_rank-r_symplectic_01-connection_4-dim_X} to continue to hold when $(X,g,J,\omega)$ is an almost Hermitian four-manifold.
% PF8-16-2024 Check this
\qed
\end{rmk}

\begin{proof}[Proof of Lemma \ref{lem:Kotschick_2-5_rank-r_symplectic_n=2_or_3}]
We project
%PF8-16-2024 Check this projection
Equation \eqref{eq:Bochner-Weitzenbock_formula_coupled_Dirac_operator_symplectic_manifold_n=2_or_3} onto its $\Omega^{0,2}(E)$ component, choose $\varphi = 0$, and note that $A_\can = A_2$ on $\Lambda^{0,2}(E)$ to give
\begin{multline}
  \label{eq:Bochner-Weitzenbock_formula_coupled_Dirac_operator_symplectic_manifold_n=2_or_3_psi}
  \left(\bar\partial_A\bar\partial_A^* + \bar\partial_A^*\bar\partial_A\right)\psi
  = \frac{1}{2}\nabla_{A_2\otimes A}^*\nabla_{A_2\otimes A}\psi
  + \frac{R}{8}\psi + \frac{i}{4}\Lambda_\omega F_{A_{K^*}}\psi + \frac{i}{2}\Lambda_\omega F_A\psi,
  \\
  \text{for all } \psi \in \Omega^{0,2}(E).
\end{multline}
By applying the curvature identity \eqref{eq:Curvature_anti-canonical_line_bundle_and_scalar_curvature} in the forthcoming Lemma \ref{lem:Curvature_anti-canonical_line_bundle_almost_Hermitian_manifold}, we see that Equation \eqref{eq:Bochner-Weitzenbock_formula_coupled_Dirac_operator_symplectic_manifold_n=2_or_3_psi} gives
\[
  \left(\bar\partial_A\bar\partial_A^* + \bar\partial_A^*\bar\partial_A\right)\psi
  =
  \frac{1}{2}\nabla_{A_2\otimes A}^*\nabla_{A_2\otimes A}\psi
  + \frac{i}{2}\Lambda_\omega F_{A_{K^*}}\psi + \frac{i}{2}\Lambda_\omega F_A\psi,
\]
and this is \eqref{eq:Kotschick_section_2_eq_13_rank-r_symplectic_01-connection_2n-dim_X}.
% PF8-12-2024 Explain the relation between A_n and A_can and A_{K^*}
If $n=2$, then $A_n = A_{K^*}$ and $\bar\partial_A\psi = 0$ and so \eqref{eq:Kotschick_section_2_eq_13_rank-r_symplectic_01-connection_2n-dim_X} reduces to \eqref{eq:Kotschick_section_2_eq_13_rank-r_symplectic_01-connection_4-dim_X}.
\end{proof}

% COMMENT For scalar curvature in almost Hermitian geometry, see article by Ross \url{https://www.aimath.org/WWN/singularvariety/stab.pdf}, references therein (including Donaldson), Viaclovsky \url{https://www.math.uci.edu/~jviaclov/lecturenotes/Japan_Lectures_2018.pdf}, LeBrun \url{https://arxiv.org/pdf/dg-ga/9412006}, and Fu and Zhou \cite{Fu_Zhou_2022}.

% COMMENT For common g, h, \omega conventions, see \url{https://en.wikipedia.org/wiki/Kähler_manifold}

% COMMENT For the curvature of the canonical line bundle \url{https://en.wikipedia.org/wiki/Canonical_bundle}, see \url{https://math.stackexchange.com/questions/4059499/class-of-ricci-form-and-first-chern-class-of-canonical-bundle} (it appears to be the Ricci curvature of the Levi--Civita connection and the trace of the Ricci curvature is the scalar curvature).

\begin{lem}[Curvature of the connection on the anti-canonical line bundle induced by the metric on an almost Hermitian manifold]
\label{lem:Curvature_anti-canonical_line_bundle_almost_Hermitian_manifold}
Let $(X,g,J,\omega)$ be an almost Hermitian manifold. If $R$ is the scalar curvature for the Riemannian metric $g$ and $A_{K^*}$ is the Chern connection on the anti-canonical line bundle $K^*$ induced by the Hermitian metric
%PF11-25-2025 h=g-i\omega
$h = g + i\omega$ on $X$, then
\begin{equation}
  \label{eq:Curvature_anti-canonical_line_bundle_and_scalar_curvature}
  i\Lambda_\omega F_{A_{K^*}} = \frac{R}{2}.
\end{equation}  
\end{lem}

\begin{proof}
Assume that $X$ has real dimension $2n$. Let $A_K$ and $F_{A_K}$ denote canonical connection and curvature of the canonical line bundle $K := \Lambda^{n,0}(X) = \wedge^n (T^*X^{1,0})$.
%PF8-16-2024 Resolve these competing definitions
Then by Viaclovsky \cite[Proposition 4.9, Equation (4.34), p. 21]{Viaclovsky_japan_lectures}, 
\[
  iF_{A_K} = -\rho,
\]
where the \emph{Ricci form} is defined by
\[
  \rho(Y,Z) = \Ric(JY,Z), \quad\text{for all } Y, Z \in \Omega^0(TX),
\]
and which is a real, closed $(1,1)$ form. For ease of exposition, we shall temporarily assume that $(X,g,J,\omega)$ is a complex Hermitian $n$-manifold before proceeding to the general case. By \cite[Proposition 3.7, Equation (3.49), p. 15]{Viaclovsky_japan_lectures}, the Ricci form can be written with respect to local holomorphic coordinates $\{z_\alpha\}_{\alpha=1}^n$ as
\[
  \rho = \sum_{\alpha,\beta=1}^n iR_{\alpha\bar\beta}dz_\alpha\wedge d\bar z_\beta,
\]
and $R_{\alpha\bar\beta} = \Ric(J\partial_{z_\alpha}, \partial_{\bar z_\beta})$ and $\rho_{\alpha\bar\beta} = \rho(J\partial_{z_\alpha}, \partial_{\bar z_\beta}) = iR_{\alpha\bar\beta}$ by \cite[Equations (3.52) and (3.53), p. 16]{Viaclovsky_japan_lectures}. The \emph{scalar curvature} is given by
\[
  R = \sum_{\alpha,\beta=1}^n g^{\alpha\beta}R_{\alpha\beta} = \sum_{\alpha,\beta=1}^n  2h^{\alpha\bar\beta}R_{\alpha\bar\beta}
\]
according to \cite[Equation (3.30), p. 13 and line preceding section 3.4, p. 15, and Remark 2.4, p. 10]{Viaclovsky_japan_lectures}. (See \cite[Section 2.1]{Viaclovsky_japan_lectures} for Viaclovsky's conventions.) 
We write $A_{K^*}$ and $F_{A_{K^*}}$ for the induced connection and curvature of the anticanonical line bundle $K^* \cong \Lambda^{0,n}(X) = \wedge^n(T^*X^{0,1})$. Thus,
% {\color{blue}FLAG: Is there a sign error in the second equality here?]}
%PF8-12-2024 Yes, corrected
% \[
%   F_{A_{K^*}} = -F_{A_K} = i(iF_{A_K}) = i(-\rho),
% \]
\[
  iF_{A_{K^*}} = -iF_{A_K} = \rho,
\]
that is,
\begin{equation}
  \label{eq:Curvature_anti-canonical_line_bundle_is_Ricci_form}
  iF_{A_{K^*}} = \rho,
\end{equation}
and
% {\color{blue} [FLAG: could be a missing factor of $i$ here when going from $\rho_{\alpha\bar\beta}$ to $R_{\alpha\bar\beta}$.]}
% \begin{equation}
%   \label{eq:Curvature_anti-canonical_line_bundle_is_Ricci_form}
%   \Tr_h F_{A_{K^*}} = i\Tr_h\rho := ih^{\alpha\bar\beta}\rho_{\alpha\bar\beta} = ih^{\alpha\bar\beta}\Ric_{\alpha\bar\beta}
%   = \frac{i}{2}g^{\alpha\beta}R_{\alpha\beta} = \frac{i}{2}R.
% \end{equation}
%PF8-12-2024 I attempted to correct
\begin{equation}
  \label{eq:Curvature_anti-canonical_line_bundle_is_Ricci_form_trace_h}
  \Tr_h iF_{A_{K^*}}
  =
  \Tr_h\rho
  :=
  \sum_{\alpha,\beta=1}^n h^{\alpha\bar\beta}\rho_{\alpha\bar\beta}
  =
  \sum_{\alpha,\beta=1}^n h^{\alpha\bar\beta}\Ric_{\alpha\bar\beta}
  =
  \sum_{\alpha,\beta=1}^n \frac{1}{2}g^{\alpha\beta}R_{\alpha\beta} = \frac{R}{2}.
\end{equation}
See also Moroianu \cite{Moroianu_2007} for an exposition of the preceding facts too, where he writes $h^{\alpha\bar\beta}$, as we do, instead of $g^{\alpha\bar\beta}$, as does Viaclovsky.

When $X$ is almost Hermitian rather than complex K\"ahler, the preceding identities continue to hold. The only difference in the proof is that the local frame $\{\partial/\partial z_\alpha, \partial/\partial\bar z_\beta\}$ for $TX^{1,0} \oplus TX^{0,1}$ is replaced by a local frame $\{Z_\alpha, \bar Z_\beta\}$ and the dual local frame $\{dz_\alpha, d\bar z_\beta\}$ for $T^*X^{1,0} \oplus T^*X^{0,1}$ is replaced by the dual local frame $\{Z_\alpha^*, \bar Z_\beta^*\}$.

By Kobayashi \cite[Equations (4.1.1) and (4.1.2), p. 92]{Kobayashi_differential_geometry_complex_vector_bundles},
\[
  i\Lambda_\omega F_{A_{K^*}}
  =
  \sum_{\alpha,\beta=1}^n h^{\alpha\bar\beta}F_{A_{K^*}}
  \left( \frac{\partial}{\partial z_\alpha}, \frac{\partial}{\partial\bar z_\beta} \right)
  =
  \sum_{\alpha,\beta=1}^n h^{\alpha\bar\beta}(F_{A_{K^*}})_{\alpha\bar\beta}
  =
  \sum_{\alpha,\beta=1}^n \frac{1}{2}g^{\alpha\beta}(F_{A_{K^*}})_{\alpha\beta},
\]
where the last equality follows from by \cite[Remark 2.4, p. 10]{Viaclovsky_japan_lectures} (this also matches Kobayashi's conventions in \cite[Equations (7.6.5) and (7.6.7), p. 251]{Kobayashi_differential_geometry_complex_vector_bundles} relating $g$, $h$, and $\omega$), and so
\[
  i\Lambda_\omega F_{A_n} = i\Lambda_\omega F_{A_{K^*}} = \Tr_h iF_{A_{K^*}} = \frac{R}{2}.
\]
This completes the proof of the curvature identity \eqref{eq:Curvature_anti-canonical_line_bundle_and_scalar_curvature}.
\end{proof}

\section[Differential inequality for squared pointwise norms of sections of $E$]{Differential inequality for the squared pointwise norms of sections of $E$ with a singular Taubes perturbation}
\label{sec:Differential_inequality_squared_pointwise_norm_holomorphic_spinor}
In this section, we prove Lemma \ref{lem:Taubes_1996_SW_to_Gr_eq_2-3_inequality}, which gives a linear second order elliptic differential inequality for $|\alpha|_E^2$ --- an analogue of Taubes' equality \cite[Section 2 (b), Equation (2.3), p. 854]{TauSWGromov} (compare Kotschick \cite[Section 3, Equation (19), p. 206]{KotschickSW}) --- when $(A,\varphi,\psi)$ with $(\varphi,\psi) = r^{1/2}(\alpha,\beta)$ is a solution to the system \eqref{eq:SO(3)_monopole_equations_almost_Hermitian_perturbed_intro} of non-Abelian monopole equations with a singular Taubes perturbation.

Take the pointwise $E$ inner product of the Bochner--Weitzenb\"ock identity \eqref{eq:Kotschick_section_2_eq_12} with $\varphi$ to give
\[
  \langle\nabla_A^*\nabla_A\varphi,\varphi\rangle_h
  - \langle i\Lambda_\omega F_A\varphi,\varphi\rangle_E
  - 2\langle\bar\partial_A^*\bar\partial_A\varphi,\varphi\rangle_E = 0.
\]
Taking the real part and substituting the second-order Kato equality (see Freed and Uhlenbeck \cite[Chapter 6, Equation (6.18), p. 91]{FU}) gives
\[
  \frac{1}{2}\Delta_g|\varphi|_E^2 + |\nabla_A\varphi|_E^2
  - \Real\langle i\Lambda_\omega F_A\varphi,\varphi\rangle_E
  - 2\Real\langle\bar\partial_A^*\bar\partial_A\varphi,\varphi\rangle_E= 0.
\]
Because $D_A\Phi = 0$, we have $\bar\partial_A\varphi = -\bar\partial_A^*\psi$ (when $d\omega = 0$) and so
\[
  \frac{1}{2}\Delta_g|\varphi|_E^2 + |\nabla_A\varphi|_E^2
  - \Real\langle i\Lambda_\omega F_A\varphi,\varphi\rangle_E
  + 2\Real\langle\bar\partial_A^* \bar\partial_A^*\psi,\varphi\rangle_E= 0.
\]
But $\bar\partial_A^* \bar\partial_A^*\psi = (\bar\partial_A\bar\partial_A)^*\psi$ and $F_A^{0,2} = \bar\partial_A^2 + \bar\mu\partial_A$, so that
\[
  \frac{1}{2}\Delta_g|\varphi|_E^2 + |\nabla_A\varphi|_E^2
  - \Real\langle i\Lambda_\omega F_A\varphi,\varphi\rangle_E
  + 2\Real\langle(F_A^{0,2} - \bar\mu\partial_A)^*\psi,\varphi\rangle_E= 0.
\]
We shall now substitute the perturbed non-Abelian monopole equations \eqref{eq:SO(3)_monopole_equations_(1,1)_curvature_perturbed_intro} and \eqref{eq:SO(3)_monopole_equations_(0,2)_curvature_perturbed_intro}, namely
\begin{align*}
  (\Lambda_\omega F_A)_0 &= \frac{i}{2}(\varphi\otimes\varphi^*)_0 - \frac{i}{2}\star(\psi\otimes\psi^*)_0
                           - \frac{ir}{4}\wp(\psi)(\psi),
  \\
  F_A^{0,2} &= \frac{1}{2}(\psi\otimes\varphi^*)_0,
\end{align*}
into the preceding identity. First multiplying \eqref{eq:SO(3)_monopole_equations_(1,1)_curvature_perturbed_intro} by $-i$ to give
%PF7-29-2025 \wp --> \wp(\psi)
%PF10-10-2024 Recheck
%{\color{blue}[TODO Note correction below (factor of $i$ deleted)]}
\begin{equation}
  \label{eq:SO(3)_monopole_equations_(1,1)_curvature_Taubes_perturbation}
  -(i\Lambda_\omega F_A)_0 =
  \frac{1}{2}(\varphi\otimes\varphi^*)_0 - \frac{1}{2}\star(\psi\otimes\psi^*)_0
  - \frac{r}{4}\wp(\psi),
\end{equation}
and using the fact that $F_A = (F_A)_0 + \frac{1}{2}F_{A_d}$ in our substitution yields
%PF10-10-2024 Recheck
%{\color{blue}[TODO Note correction below (factor of $i$ deleted)]}
\begin{multline}
  \label{eq:Taubes_1996_SW_to_Gr_eq_2-3_raw_before_rescaling}
  \frac{1}{2}\Delta_g|\varphi|_E^2 + |\nabla_A\varphi|_E^2
  + \frac{1}{2}\Real\langle (\varphi\otimes\varphi^*)_0\varphi,\varphi\rangle_E
  - \frac{1}{2}\Real\langle \star(\psi\otimes\psi^*)_0\varphi,\varphi\rangle_E
  - \frac{r}{4}\langle \wp(\psi)\varphi,\varphi\rangle_E
  \\
  - \Real\langle i\Lambda_\omega F_{A_d}\varphi,\varphi\rangle_E
  + \Real\langle(\psi\otimes\varphi^*)_0^*\psi,\varphi\rangle_E
  - 2\Real\langle(\bar\mu\partial_A)^*\psi,\varphi\rangle_E = 0.
\end{multline}
The identity \eqref{eq:Taubes_1996_SW_to_Gr_eq_2-3_before_rescaling} has a formal structure similar to that of Taubes \cite[Section 2 (b), Equation (2.3), p. 854]{TauSWGromov} (compare Kotschick \cite[Section 3, Equation (19), p. 206]{KotschickSW}). A further simplification of \eqref{eq:Taubes_1996_SW_to_Gr_eq_2-3_raw_before_rescaling} requires the

\begin{claim}[Pointwise inequalities for $\varphi$]
\label{claim:Pointwise_equality_and_inequalities_varphi}
The following equality and inequalities hold:
%PF10-10-2024 Recheck
%{\color{blue}[TODO Note correction below (factor of $i$ deleted)]} 
\begin{subequations}
\label{eq:Inequalities_simplifying_Taubes_1996_SW_to_Gr_eq_2-3}
\begin{gather}
  \label{eq:5}
  -\frac{1}{2}|\psi|_{\Lambda^{0,2}(E)}^2|\varphi|_E^2
  \leq
  \langle \star(\psi\otimes\psi^*)_0\varphi,\varphi\rangle_E
  \leq
  \frac{1}{2}|\psi|_{\Lambda^{0,2}(E)}^2|\varphi|_E^2,
  \\
  \label{eq:6}
  \frac{1}{2}|\psi|_{\Lambda^{0,2}(E)}^2|\varphi|_E^2
  \leq
%PF3-17-2025 Missing a "\star" after \langle below 
  \langle(\psi\otimes\varphi^*)_0^*\psi,\varphi\rangle_E
  \leq
  |\psi|_{\Lambda^{0,2}(E)}^2|\varphi|_E^2,
  \\
  \label{eq:Quartic_identity_varphi}
  \langle (\varphi\otimes\varphi^*)_0\varphi,\varphi\rangle_E
  =
  \frac{1}{2}|\varphi|_E^4,
  \\
  \label{eq:9}
  %PF7-9-2024 Check this!!!
  -|\varphi|_E^2 \leq \langle \wp(\psi)\varphi,\varphi\rangle_E \leq |\varphi|_E^2.
\end{gather}
\end{subequations}
\end{claim}

\begin{proof}[Proof of Claim \ref{claim:Pointwise_equality_and_inequalities_varphi}]
Using $\varphi^* = \langle\cdot,\varphi\rangle_E$, we see that, for a local
%PF8-28-2024 We're changing h to H
$H$-orthonormal frame $\{e_k\}_{k=1}^2$ for $E$,
\[
  \tr_E(\varphi\otimes\varphi^*)
  =
  \sum_{k=1}^2 \langle (\varphi\otimes\varphi^*)e_k,e_k\rangle_E
  =
  \sum_{k=1}^2 \langle \varphi \langle e_k,\varphi\rangle_E, e_k\rangle_E
  =
  \sum_{k=1}^2 |\langle e_k,\varphi\rangle_E|^2
  =
  |\varphi|_E^2
  \in \Omega^0(X;\RR),
\]
and so
\[
  (\varphi\otimes\varphi^*)_0
  =
  \varphi\otimes\varphi^* - \frac{1}{2}\tr_E(\varphi\otimes\varphi^*)\,\id_E
  =
  \varphi\otimes\varphi^* - \frac{1}{2}|\varphi|_E^2\id_E.
\]
Therefore,
\[
  \langle (\varphi\otimes\varphi^*)_0\varphi,\varphi\rangle_E
  =
  \langle \varphi\otimes\varphi^*\varphi,\varphi\rangle_E
   - \frac{1}{2}|\varphi|_E^2\langle \varphi,\varphi\rangle_E
\]
and this simplifies to \eqref{eq:Quartic_identity_varphi}. Similarly, we observe that
\[
  (\psi\otimes\psi^*)_0
  =
  (\psi\otimes\psi^*) - \frac{1}{2}\tr_E(\psi\otimes\psi^*)\,\id_E.
\]
Now
\begin{multline*}
  \tr_E(\psi\otimes\psi^*)
  =
  \sum_{k=1}^2 \langle (\psi\otimes\psi^*)e_k,e_k\rangle_E
  =
  \sum_{k=1}^2 \langle \psi \langle e_k,\psi\rangle_E, e_k\rangle_E
  =
  %PF7-15-2024 Can we just replace \otimes by \wedge here?
  \sum_{k=1}^2 \langle e_k,\psi\rangle_E \wedge \langle \psi, e_k\rangle_E
  \\
  =
  \sum_{k=1}^2 \langle e_k,\psi\rangle_E \wedge \overline{\langle e_k,\psi\rangle_E}
  =
  \sum_{k=1}^2 \langle e_k,\psi\rangle_E \wedge \star\overline{\langle e_k,\psi\rangle_E}.
\end{multline*}
Thus,
\[
  \star\tr_E(\psi\otimes\psi^*)
  =
  %PF7-15-2024 By FL Equation (7.1.1), p. 94, but check
  \sum_{k=1}^2 \left\langle \langle e_k,\psi\rangle_E,
    \langle e_k,\psi\rangle_E \right\rangle_{\Lambda^{0,2}(X)}
  =
  \sum_{k=1}^2 \left|\langle e_k,\psi\rangle_E\right|_{\Lambda^{0,2}(X)}^2
  =
  |\psi|_{\Lambda^{0,2}(E)}^2
  \in
  \Omega^0(X;\RR).
\]
Hence,
\begin{align*}
  \langle \star(\psi\otimes\psi^*)_0\varphi,\varphi\rangle_E
  &=
    \langle \star(\psi\otimes\psi^*)\varphi,\varphi\rangle_E
    - \frac{1}{2}\langle \star\tr_E(\psi\otimes\psi^*)\varphi,\varphi\rangle_E
  \\
  &=
    \langle \star(\psi\otimes\psi^*)\varphi,\varphi\rangle_E
    - \frac{1}{2}|\psi|_{\Lambda^{0,2}(E)}^2|\varphi|_E^2
  \\
  &=
  \langle \star(\psi\langle\varphi,\psi\rangle_E),\varphi\rangle_E - \frac{1}{2}|\psi|_{\Lambda^{0,2}(E)}^2|\varphi|_E^2
  \\
  %PF7-15-2024 By FL Equation (7.1.1), p. 94, but check
  &=
  |\langle\varphi,\psi\rangle_E|_{\Lambda^{0,2}(X)}^2 - \frac{1}{2}|\psi|_{\Lambda^{0,2}(E)}^2|\varphi|_E^2.
\end{align*}
Using $|\langle\varphi,\psi\rangle_E|_{\Lambda^{0,2}(X)} \leq |\varphi|_E|\psi|_{\Lambda^{0,2}(E)}$, we obtain
%TL7-7-2024: Clear from thinking of previous as $c=a-b$ where $a,b\ge 0$ so $c\ge -b$ and $c\le a$.
\[
  -\frac{1}{2}|\psi|_{\Lambda^{0,2}(E)}^2|\varphi|_E^2
  \leq
  \langle \star(\psi\otimes\psi^*)_0\varphi,\varphi\rangle_E
  \leq
  \frac{1}{2}|\psi|_{\Lambda^{0,2}(E)}^2|\varphi|_E^2,
\]
which verifies the claimed inequalities \eqref{eq:5}. By combining the definition \eqref{eq:Definition_wp_intro} of $\wp(\psi)$ with the  inequalities\eqref{eq:5} we obtain the inequalities \eqref{eq:9} on $X_0 = \{x\in X:\psi(x)\neq 0\}$ in \eqref{eq:X0} and the inequalities \eqref{eq:9} clearly hold on $X\less X_0$, where $\wp(\psi) \equiv 0$.

Similarly, 
\begin{multline*}
  \tr_E(\varphi\otimes\psi^*)
  =
  \sum_{k=1}^2 \langle (\varphi\otimes\psi^*)e_k,e_k\rangle_E
  =
  \sum_{k=1}^2 \langle \varphi \langle e_k,\psi\rangle_E, e_k\rangle_E
  =
  \sum_{k=1}^2 \langle \varphi, e_k\rangle_E \langle e_k,\psi\rangle_E
  \\
  =
  \langle\varphi,\psi\rangle_E \in \Omega^{2,0}(X).
\end{multline*}
Thus,
\[
  (\varphi\otimes\psi^*)_0
  =
  \varphi\otimes\psi^* - \frac{1}{2}\langle\varphi,\psi\rangle_E\,\id_E
  \in \Omega^{2,0}(\fsl(E)),
\]
and so we obtain
%PF8-16-2024
%[TODO - Tom, please check below]
\begin{multline*}
  \langle(\psi\otimes\varphi^*)_0^*\psi,\varphi\rangle_E
  = \langle(\varphi\otimes\psi^*)_0\psi,\varphi\rangle_E
  = \langle(\varphi\otimes\psi^*)\psi,\varphi\rangle_E
  - \frac{1}{2}\langle\star\langle\varphi,\psi\rangle_E\psi,\varphi\rangle_E
  \\
  = \langle \varphi |\psi|_{\Lambda^{0,2}(E)}^2,\varphi\rangle_E
  - \frac{1}{2}|\langle\varphi,\psi\rangle_E|_{\Lambda^{0,2}(X)}^2
  = |\psi|_{\Lambda^{0,2}(E)}^2|\varphi|_E^2
  - \frac{1}{2}|\langle\varphi,\psi\rangle_E|_{\Lambda^{0,2}(X)}^2,
\end{multline*}
and so, using $|\langle\varphi,\psi\rangle_E|_{\Lambda^{0,2}(X)} \leq |\varphi|_E|\psi|_{\Lambda^{0,2}(E)}$, we obtain the claimed inequalities \eqref{eq:6}:
%TL7-7-2024:  Second inequality is clear.  Does the first follow from applying Cauchy-Schwarts to  %\frac{1}{2}|\langle\varphi,\psi\rangle_E|_{\Lambda^{0,2}(X)}^2$?
\[
  \frac{1}{2}|\psi|_{\Lambda^{0,2}(E)}^2|\varphi|_E^2
  \leq
  \langle(\psi\otimes\varphi^*)_0^*\psi,\varphi\rangle_E
  \leq
  |\psi|_{\Lambda^{0,2}(E)}^2|\varphi|_E^2.
\]
This completes the proof of Claim \ref{claim:Pointwise_equality_and_inequalities_varphi}.
\end{proof}

By substituting the equality \eqref{eq:Quartic_identity_varphi} into \eqref{eq:Taubes_1996_SW_to_Gr_eq_2-3_raw_before_rescaling}, we see that
%PF10-10-2024 Recheck
%{\color{blue}[TODO Note correction below (factor of $i$ deleted)]}
\begin{multline}
  \label{eq:Taubes_1996_SW_to_Gr_eq_2-3_before_rescaling}
  \frac{1}{2}\Delta_g|\varphi|_E^2 + |\nabla_A\varphi|_E^2
  + \frac{1}{4}|\varphi|_E^4
  - \frac{1}{2}\Real\langle \star(\psi\otimes\psi^*)_0\varphi,\varphi\rangle_E
  - \frac{r}{4}\langle \wp(\psi)\varphi,\varphi\rangle_E
  \\
  - \Real\langle i\Lambda_\omega F_{A_d}\varphi,\varphi\rangle_E
  + \Real\langle(\psi\otimes\varphi^*)_0^*\psi,\varphi\rangle_E
  - 2\Real\langle(\bar\mu\partial_A)^*\psi,\varphi\rangle_E = 0 \quad\text{on } X.
\end{multline}
By substituting the inequalities \eqref{eq:5} (upper bound) and \eqref{eq:6} (lower bound) and \eqref{eq:9} (upper bound) into \eqref{eq:Taubes_1996_SW_to_Gr_eq_2-3_before_rescaling}, we obtain
\begin{multline*}
  \frac{1}{2}\Delta_g|\varphi|_E^2 + |\nabla_A\varphi|_E^2
  + \frac{1}{4}|\varphi|_E^4
  - \frac{1}{4}|\psi|_{\Lambda^{0,2}(E)}^2|\varphi|_E^2
  - \frac{r}{4}|\varphi|_E^2
  \\
  - \Real\langle i\Lambda_\omega F_{A_d}\varphi,\varphi\rangle_E
  + \frac{1}{2}|\psi|_{\Lambda^{0,2}(E)}^2|\varphi|_E^2
  - 2\Real\langle(\bar\mu\partial_A)^*\psi,\varphi\rangle_E \leq 0 \quad\text{on } X,
\end{multline*}
and thus
\begin{multline}
  \label{eq:Taubes_1996_SW_to_Gr_eq_2-3_simplified_before_rescaling}
  \frac{1}{2}\Delta_g|\varphi|_E^2 + |\nabla_A\varphi|_E^2
  + \frac{1}{4}|\varphi|_E^4
  + \frac{1}{4}|\psi|_{\Lambda^{0,2}(E)}^2|\varphi|_E^2
  - \frac{r}{4}|\varphi|_E^2
  \\
  - \Real\langle i\Lambda_\omega F_{A_d}\varphi,\varphi\rangle_E
  - 2\Real\langle(\bar\mu\partial_A)^*\psi,\varphi\rangle_E \leq 0 \quad\text{on } X.
\end{multline}
After making the substitution \eqref{eq:Taubes_1996_SW_to_Gr_eq_1-19} in \eqref{eq:Taubes_1996_SW_to_Gr_eq_2-3_simplified_before_rescaling}, namely
\[
  (\varphi,\psi) = r^{1/2}(\alpha,\beta) \in \Omega^0(E) \oplus \Omega^{0,2}(E),
\]
and canceling factors of $r$ on both sides, we see that we have proved the

\begin{lem}[Differential inequality for the squared pointwise norms of sections of $E$ with a singular Taubes perturbation]
\label{lem:Taubes_1996_SW_to_Gr_eq_2-3_inequality}
Continue the hypotheses of Proposition \ref{mainprop:Taubes_1996_SW_to_Gr_2-1}. If $\alpha \in \Omega^0(E)$ and $\beta \in \Omega^{0,2}(E)$ are defined by $\Phi = r^{1/2}(\alpha,\beta)$ as in \eqref{eq:Taubes_1996_SW_to_Gr_eq_1-19}, then 
% TL7-7-2024: Factor of $r$ cancelled out from both sides.
% PF7-16-2024 Yes
%PF8-16-2024 Indicate Sobolev notation and inequality a.e. on X.
\begin{multline}
  \label{eq:Taubes_1996_SW_to_Gr_eq_2-3_inequality}
  \frac{1}{2}\Delta_g|\alpha|_E^2 + |\nabla_A\alpha|_E^2
  + \frac{r}{4}|\alpha|_E^4
  + \frac{r}{4}|\beta|_{\Lambda^{0,2}(E)}^2|\alpha|_E^2
  - \frac{r}{4}|\alpha|_E^2
  \\
  - \Real\langle i\Lambda_\omega F_{A_d}\alpha,\alpha\rangle_E
  - 2\Real\langle(\bar\mu\partial_A)^*\beta,\alpha\rangle_E \leq 0 \quad\text{on } X.
\end{multline}
\end{lem}

The inequality \eqref{eq:Taubes_1996_SW_to_Gr_eq_2-3_inequality} is an analogue of Taubes' equality \cite[Section 2 (b), Equation (2.3), p. 854]{TauSWGromov}.

\section[Differential inequality for squared pointwise norms of sections of $\Lambda^{0,2}(E)$]{Differential inequality for the squared pointwise norms of sections of $\Lambda^{0,2}(E)$ with a singular Taubes perturbation}
\label{sec:Differential_inequality_squared_pointwise_norm_anti-holomorphic_spinor}
In this section, we prove Lemma \ref{lem:Taubes_1996_SW_to_Gr_eq_2-7} --- an analogue of Taubes \cite[Section 2 (b), Equation (2.7), p. 855]{TauSWGromov} (compare Kotschick \cite[Section 3, Equation (23), p. 207]{KotschickSW}) --- which gives a linear second order elliptic differential inequality for $|\beta|_{\Lambda^{0,2}(E)}^2$when $(A,\varphi,\psi)$ with $(\varphi,\psi) = r^{1/2}(\alpha,\beta)$ is a solution to the system \eqref{eq:SO(3)_monopole_equations_almost_Hermitian_perturbed_intro} of non-Abelian monopole equations with a singular Taubes perturbation.

Take the pointwise $\Lambda^{0,2}(E)$ inner product of the identity \eqref{eq:Kotschick_section_2_eq_13_rank-r_symplectic_01-connection_4-dim_X} with $\psi$ to give
\[
  \langle\nabla_{A_2\otimes A}^*\nabla_{A_2\otimes A}\psi,\psi\rangle_{\Lambda^{0,2}(E)}
  + \langle i\Lambda_\omega(F_{A_2} + F_A)\psi,\psi\rangle_{\Lambda^{0,2}(E)}
  - 2\langle\bar\partial_A\bar\partial_A^*\psi,\psi\rangle_{\Lambda^{0,2}(E)} = 0.
\]
Taking the real part and substituting the second-order Kato equality (see Freed and Uhlenbeck \cite[Chapter 6, Equation (6.18), p. 91]{FU}) gives
\[
  \frac{1}{2}\Delta_g|\psi|_{\Lambda^{0,2}(E)}^2 + |\nabla_{A_2\otimes A}\psi|_{T^3(E)}^2
  + \Real\langle i\Lambda_\omega(F_{A_2} + F_A)\psi,\psi\rangle_{\Lambda^{0,2}(E)}
  - 2\Real\langle\bar\partial_A\bar\partial_A^*\psi,\psi\rangle_{\Lambda^{0,2}(E)}= 0,
\]
where, for any integer $k \geq 0$, we denote
\[
  T^k(E) := \otimes^k T^*X \otimes E.
\]  
Because $D_A\Phi = 0$, we have $\bar\partial_A^*\psi = -\bar\partial_A\varphi$ (when $d\omega = 0$) and so
\[
  \frac{1}{2}\Delta_g|\psi|_{\Lambda^{0,2}(E)}^2 + |\nabla_{A_2\otimes A}\psi|_{T^3(E)}^2
  + \Real\langle i\Lambda_\omega(F_{A_2} + F_A)\psi,\psi\rangle_{\Lambda^{0,2}(E)}
  + 2\Real\langle\bar\partial_A\bar\partial_A\varphi,\psi\rangle_{\Lambda^{0,2}(E)}= 0.
\]
We rewrite the term $\bar\partial_A\bar\partial_A\varphi$ using $F_A^{0,2} = \bar\partial_A^2 + \bar\mu\partial_A$, so that
\begin{multline*}
  \frac{1}{2}\Delta_g|\psi|_{\Lambda^{0,2}(E)}^2 + |\nabla_{A_2\otimes A}\psi|_{T^3(E)}^2
  + \Real\langle i\Lambda_\omega(F_{A_2} + F_A)\psi,\psi\rangle_{\Lambda^{0,2}(E)}
  \\
  + 2\Real\langle(F_A^{0,2} - \bar\mu\partial_A)\varphi,\psi\rangle_{\Lambda^{0,2}(E)}= 0.
\end{multline*}
Substituting the perturbed non-Abelian monopole equations \eqref{eq:SO(3)_monopole_equations_(1,1)_curvature_Taubes_perturbation} and \eqref{eq:SO(3)_monopole_equations_(0,2)_curvature_perturbed_intro}, namely
%PF10-10-2024 Recheck
\begin{align*}
  (i\Lambda_\omega F_A)_0
  &=
  -\frac{1}{2}(\varphi\otimes\varphi^*)_0 + \frac{1}{2}\star(\psi\otimes\psi^*)_0
  + \frac{r}{4}\wp(\psi),
  \\
  F_A^{0,2} &= \frac{1}{2}(\psi\otimes\varphi^*)_0,
\end{align*}
and using the fact that $F_A = (F_A)_0 + \frac{1}{2}F_{A_d}$ yields
%PF10-10-2024 Recheck
%{\color{blue}[TODO Note correction below (factor of $i$ deleted)]}
\begin{multline}
  \label{eq:Taubes_1996_SW_to_Gr_eq_2-4_equality_before_rescaling}
  \frac{1}{2}\Delta_g|\psi|_{\Lambda^{0,2}(E)}^2 + |\nabla_{A_2\otimes A}\psi|_{T^3(E)}^2
  + \Real\langle i\Lambda_\omega F_{A_2}\psi,\psi\rangle_{\Lambda^{0,2}(E)}
  + \frac{1}{2}\Real\langle i\Lambda_\omega F_{A_d}\psi,\psi\rangle_{\Lambda^{0,2}(E)}
  \\
  - \frac{1}{2}\Real\langle (\varphi\otimes\varphi^*)_0\psi,\psi\rangle_{\Lambda^{0,2}(E)}
  + \frac{1}{2}\Real\langle \star(\psi\otimes\psi^*)_0\psi,\psi\rangle_{\Lambda^{0,2}(E)}
  + \frac{r}{4}\langle \wp(\psi)\psi,\psi\rangle_{\Lambda^{0,2}(E)}
  \\
  + 2\Real\langle(\psi\otimes\varphi^*)_0\varphi, \psi\rangle_{\Lambda^{0,2}(E)}
  - 2\Real\langle\bar\mu\partial_A\varphi,\psi\rangle_{\Lambda^{0,2}(E)} = 0.
\end{multline}
The identity \eqref{eq:Taubes_1996_SW_to_Gr_eq_2-4_equality_before_rescaling} has a formal structure similar to that of Taubes \cite[Section 2 (b), Equation (2.4), p. 855]{TauSWGromov} (compare Kotschick \cite[Section 3, Equation (23), p. 207]{KotschickSW}).
% PF7-3-2024 Can the line bundle curvature terms in these equalities be represented by *positive* (1,1) forms?
We make the

\begin{claim}[Pointwise inequalities for $\psi$]
\label{claim:Pointwise_equality_and_inequalities_psi}  
The following pointwise equality and inequalities hold on $X$:
%PF10-10-2024 Recheck
%{\color{blue}[TODO Note correction below (factor of $i$ deleted)]} 
\begin{subequations}
\label{eq:Inequalities_simplifying_Taubes_1996_SW_to_Gr_eq_2-4}  
\begin{gather}
  \label{eq:Quartic_identity_psi}
  \langle \star(\psi\otimes\psi^*)_0\psi,\psi\rangle_{\Lambda^{0,2}(E)}
  = \frac{1}{2}|\psi|_{\Lambda^{0,2}(E)}^4,
  \\
  \label{eq:3}
  %TL7-7-2024: This one uses $\langle (\phi\otimes\phi^*)\psi,\psi\rangle = |\langle\varphi,\psi\rangle_E|^2_{\La^{0,2}(X)}$ and the CS inequality as before
  - \frac{1}{2}|\varphi|_E^2|\psi|_{\Lambda^{0,2}(E)}^2
  \leq
  \langle (\varphi\otimes\varphi^*)_0\psi,\psi\rangle_{\Lambda^{0,2}(E)}
  \leq \frac{1}{2}|\varphi|_E^2|\psi|_{\Lambda^{0,2}(E)}^2,
  \\
  \label{eq:4}
  \frac{1}{2}|\varphi|_E^2 |\psi|_{\Lambda^{0,2}(E)}^2
  \leq
  \langle(\psi\otimes\varphi^*)_0\varphi, \psi\rangle_{\Lambda^{0,2}(E)}
  \leq |\varphi|_E^2 |\psi|_{\Lambda^{0,2}(E)}^2,
  \\
  \label{eq:4-1}
  \langle \wp(\psi)\psi, \psi\rangle_{\Lambda^{0,2}(E)}
  = |\psi|_{\Lambda^{0,2}(E)}^2.
\end{gather}
\end{subequations}
\end{claim}

\begin{proof}[Proof of Claim \ref{claim:Pointwise_equality_and_inequalities_psi}]
The arguments are similar to those of Feehan and Leness 
%TL12-4-2025: Updated and removed page references
\cite[Proof of Lemma 8.4.1]{Feehan_Leness_introduction_virtual_morse_theory_so3_monopoles}. Consider \eqref{eq:Quartic_identity_psi}. Observe that, since $E$ has complex rank $2$,
\[
  \star(\psi\otimes\psi^*)_0
  =
  \star(\psi\otimes\psi^*) - \frac{1}{2}\tr_E(\star(\psi\otimes\psi^*))\,\id_E.
\]
Recall from Huybrechts \cite[Lemma 1.2.4 (ii), p. 33]{Huybrechts_2005} that $\star\bar\beta = \bar\beta \in \Omega^{2,0}(X)$
%PF7-16-2024 Check
and recall from Huybrechts \cite[Section 1.2, p. 33]{Huybrechts_2005} that
\[
  \beta\wedge\bar\beta = \beta\wedge \star\bar\beta = |\beta|_{\Lambda^{0,2}(X)}^2.
\]
Then, noting that $\Lambda^{0,2}(E) = \Lambda^{0,2}(X)\otimes E$, we may write
%PF7-5-2023 Notation conflict
$\psi = \beta\otimes s$, for $\beta \in \Omega^{0,2}(X)$ and $s\in\Omega^0(E)$, and obtain
\[
  \tr_E(\star(\psi\otimes\psi^*))
  =
  \tr_E(\star(\beta\otimes s \otimes \bar\beta\otimes s^*))
  =
  \star(\beta\wedge\bar\beta)\tr_E(s\otimes s^*)
  =
  |\beta|_{\Lambda^{0,2}(X)}^2|s|_E^2
  =
  |\psi|_{\Lambda^{0,2}(E)}^2,
\]
and
\[
  \star(\psi\otimes\psi^*)\psi
  =
  \star(\beta\wedge\bar\beta) s\otimes s^*(\beta\otimes s)
  \\
  =
  |s|_E^2|\beta|_{\Lambda^{0,2}(X)}^2 \beta\otimes s
  \\
  =
  |\psi|_{\Lambda^{0,2}(E)}^2 \psi.
\]
Therefore, because
\[
  \star(\psi\otimes\psi^*)_0\psi
  =
  \star(\psi\otimes\psi^*)\psi - \frac{1}{2}|\psi|_{\Lambda^{0,2}(E)}^2\psi,
\]
we obtain the following eigenvalue equality:
\begin{equation}
  \label{eq:psi_eigenvector_for_eigenvalue_half_|psi|^2}
  \star(\psi\otimes\psi^*)_0\psi
  =
  \frac{1}{2}|\psi|_{\Lambda^{0,2}(E)}^2\psi.
 \end{equation}
Hence,
\[
  \langle \star(\psi\otimes\psi^*)\psi,\psi\rangle_{\Lambda^{0,2}(E)}
  =
  \frac{1}{2}|\psi|_{\Lambda^{0,2}(E)}^2  \langle \psi,\psi\rangle_{\Lambda^{0,2}(E)}
  =
  \frac{1}{2}|\psi|_{\Lambda^{0,2}(E)}^4.
\]  
This verifies the claimed equality \eqref{eq:Quartic_identity_psi}.

Consider \eqref{eq:3}. Observe that
\[
  (\varphi\otimes\varphi^*)_0\psi
  =
  (\varphi\otimes\varphi^*)\psi - \frac{1}{2}|\varphi|_E^2\psi
  =
  \varphi\langle\psi,\varphi\rangle_E - \frac{1}{2}|\varphi|_E^2\psi.
\]
Hence, 
\begin{multline*}
  \langle (\varphi\otimes\varphi^*)_0\psi,\psi\rangle_{\Lambda^{0,2}(E)}
  =
  \langle \varphi\langle\psi,\varphi\rangle_E,\psi\rangle_{\Lambda^{0,2}(E)}
  - \frac{1}{2}|\varphi|_E^2\langle\psi,\psi\rangle_{\Lambda^{0,2}(E)}
  \\
  = |\langle\psi,\varphi\rangle_E|_{\Lambda^{0,2}(X)}^2
  - \frac{1}{2}|\varphi|_E^2|\psi|_{\Lambda^{0,2}(E)}^2.
\end{multline*}
Since
%PF7-16-2024 Check
$|\langle\psi,\varphi\rangle_E|_{\Lambda^{0,2}(X)} \leq |\psi|_{\Lambda^{0,2}(E)}|\varphi|_E$, this verifies the claimed inequalities \eqref{eq:3}.

Consider \eqref{eq:4}.
%TL7-7-2024:
Observe that
\begin{multline*}
  \langle(\psi\otimes\varphi^*)_0\varphi, \psi\rangle_{\Lambda^{0,2}(E)}
  =
  \langle(\psi\otimes\varphi^*)\varphi, \psi\rangle_{\Lambda^{0,2}(E)}
  -
  \frac{1}{2}|\langle\psi,\varphi\rangle_E|^2_{\Lambda^{0,2}(X)}
  \\
  =
  |\varphi|_E^2|\psi|_{\Lambda^{0,2}(E)}^2
    -
  \frac{1}{2}|\langle\psi,\varphi\rangle_E|_{\Lambda^{0,2}(X)}^2
\end{multline*}
The second inequality in \eqref{eq:4} follows immediately from the preceding equality. The first inequality in \eqref{eq:4} follows by applying Cauchy--Schwarz to get
\[
  \frac{1}{2}|\langle\psi,\varphi\rangle_E|^2_{\Lambda^{0,2}(X)}
  \leq
  \frac{1}{2} |\varphi|_E^2|\psi|_{\Lambda^{0,2}(E)}^2,
\]
and so
\[
  |\varphi|_E^2|\psi|_{\Lambda^{0,2}(E)}^2
    -
  \frac{1}{2}|\langle\psi,\varphi\rangle_E|_{\Lambda^{0,2}(X)}^2
  \geq
  \left(1- \frac{1}{2}\right) |\varphi|_E^2|\psi|_{\Lambda^{0,2}(E)}^2.
\]
This verifies the claimed inequalities \eqref{eq:4}.

Consider \eqref{eq:4-1}. Recall from \eqref{eq:Definition_wp_intro} that
%PF10-10-2024 Recheck
\[
  \wp(\psi)
  =
  \begin{cases}
    2|\psi|_{\Lambda^{0,2}(E)}^{-2}\star(\psi\otimes\psi^*)_0 &\text{on } X_0,
    \\
    0 &\text{on } X\less X_0.
  \end{cases}
\]
Hence, because $\psi$ is an eigenvector of $\star(\psi\otimes\psi^*)_0$ with eigenvalue $\frac{1}{2}|\psi|_{\Lambda^{0,2}(E)}^2$ by \eqref{eq:psi_eigenvector_for_eigenvalue_half_|psi|^2},
%PF10-10-2024 Recheck
\[
  \wp(\psi)\psi = \psi \quad\text{on } X_0 \quad\text{and}\quad i\wp(\psi)\psi = 0 = \psi \quad\text{on } X\less X_0.
\]
%PF10-10-2024 Recheck
Therefore, $\wp(\psi)\psi = \psi$ on $X$ and
%PF10-10-2024 Recheck
\[
  \langle \wp(\psi)\psi,\psi\rangle_{\Lambda^{0,2}(E)}
  =
  \langle \psi,\psi\rangle_{\Lambda^{0,2}(E)}
  =
  |\psi|_{\Lambda^{0,2}(E)}^2 \quad\text{on } X.
\]
This verifies the claimed equality \eqref{eq:4-1}. This completes the proof of Claim \ref{claim:Pointwise_equality_and_inequalities_psi}.
\end{proof}

By applying the equalities \eqref{eq:Quartic_identity_psi} and \eqref{eq:4-1} and inequalities \eqref{eq:3} (upper bound) and \eqref{eq:4} (lower bound) to \eqref{eq:Taubes_1996_SW_to_Gr_eq_2-4_equality_before_rescaling}, we obtain
\begin{multline*}
  \frac{1}{2}\Delta_g|\psi|_{\Lambda^{0,2}(E)}^2 + |\nabla_{A_2\otimes A}\psi|_{T^3(E)}^2
  \\
  + \Real\langle i\Lambda_\omega F_{A_2}\psi,\psi\rangle_{\Lambda^{0,2}(E)}
  + \frac{1}{2}\Real\langle i\Lambda_\omega F_{A_d}\psi,\psi\rangle_{\Lambda^{0,2}(E)}
  \\
  - \frac{1}{2}|\varphi|_E^2|\psi|_{\Lambda^{0,2}(E)}^2
  + \frac{1}{4}|\psi|_{\Lambda^{0,2}(E)}^4
  + \frac{r}{4}|\psi|_{\Lambda^{0,2}(E)}^2
  + |\varphi|_E^2|\psi|_{\Lambda^{0,2}(E)}^2
  \\
  - 2\Real\langle\bar\mu\partial_A\varphi,\psi\rangle_{\Lambda^{0,2}(E)} \leq 0 \quad\text{on } X,
\end{multline*}
that is,
\begin{multline}
  \label{eq:Taubes_1996_SW_to_Gr_eq_2-4_before_rescaling}
  \frac{1}{2}\Delta_g|\psi|_{\Lambda^{0,2}(E)}^2 + |\nabla_{A_2\otimes A}\psi|_{T^3(E)}^2
  \\
  + \Real\langle i\Lambda_\omega F_{A_2}\psi,\psi\rangle_{\Lambda^{0,2}(E)}
  + \frac{1}{2}\Real\langle i\Lambda_\omega F_{A_d}\psi,\psi\rangle_{\Lambda^{0,2}(E)}
  \\
  + \frac{1}{2}|\varphi|_E^2|\psi|_{\Lambda^{0,2}(E)}^2
  + \frac{1}{4}|\psi|_{\Lambda^{0,2}(E)}^4
  + \frac{r}{4}|\psi|_{\Lambda^{0,2}(E)}^2
  \\
  - 2\Real\langle\bar\mu\partial_A\varphi,\psi\rangle_{\Lambda^{0,2}(E)} \leq 0 \quad\text{on } X.
\end{multline}
The preceding differential inequality \eqref{eq:Taubes_1996_SW_to_Gr_eq_2-4_before_rescaling} is a precursor to analogues of Taubes \cite[Section 2 (b), Equations (2.4) and (2.7), p. 855]{TauSWGromov}.
%TL7-7-2024: Now checking algebra here.
% \begin{multline}
%   \label{eq:Taubes_1996_SW_to_Gr_eq_2-4_before_rescaling}
%   \frac{1}{2}\Delta_g|\psi|_{\Lambda^{0,2}(E)}^2 + |\nabla_{A_2\otimes A}\psi|_{T^3(E)}^2
%   + \frac{1}{4}|\psi|_{\Lambda^{0,2}(E)}^4
%   \\
%   \leq \left(|\Lambda_\omega F_{A_2}| + |\Lambda_\omega F_{A_d}|\right)
%   |\psi|_{\Lambda^{0,2}(E)}^2
%   \\
%   + \frac{3}{2}|\varphi|_E^2|\psi|_{\Lambda^{0,2}(E)}^2
%   + 2|\bar\mu\partial_A\varphi|_{\Lambda^{0,2}(E)}|\psi|_{\Lambda^{0,2}(E)}.
% \end{multline}

Inequality \eqref{eq:Taubes_1996_SW_to_Gr_eq_2-4_before_rescaling} yields the following analogue of Taubes \cite[Section 2 (b), Equation (2.7), p. 855]{TauSWGromov} after substituting $\Phi = (\varphi,\psi) = r^{1/2}(\alpha,\beta)$ as in \eqref{eq:Taubes_1996_SW_to_Gr_eq_1-19}:
\begin{multline}
  \label{eq:Taubes_1996_SW_to_Gr_eq_2-4}
  \frac{1}{2}\Delta_g|\beta|_{\Lambda^{0,2}(E)}^2
  + |\nabla_{A_2\otimes A}\beta|_{T^3(E)}^2
  + \frac{r}{2}|\alpha|_E^2|\beta|_{\Lambda^{0,2}(E)}^2
  + \frac{r}{4}|\beta|_{\Lambda^{0,2}(E)}^4
  + \frac{r}{4}|\beta|_{\Lambda^{0,2}(E)}^2
  \\
  + \Real\langle i\Lambda_\omega F_{A_2}\beta,\beta\rangle_{\Lambda^{0,2}(E)}
  + \frac{1}{2}\Real\langle i\Lambda_\omega F_{A_d}\beta,\beta\rangle_{\Lambda^{0,2}(E)}
  \\
  - 2\Real\langle\bar\mu\partial_A\alpha,\beta\rangle_{\Lambda^{0,2}(E)} \leq 0 \quad\text{on } X.
\end{multline}
Applying the inequality $r|\beta|_{\Lambda^{0,2}(E)}^4/4 \geq 0$ to the left-hand side of \eqref{eq:Taubes_1996_SW_to_Gr_eq_2-4} (although this gives a weaker inequality) yields
%PF7-8-2024 Probable cut
    %     adding Taubes' perturbations originating in the addition of terms $u_1r\omega$ and $\wp(\psi) r\omega$ to $(F_A^+)_0$ so that $\langle \wp(\psi) r\psi, \psi \rangle_E = r|\psi|_E^2$ (and $\langle u_1r\psi, \psi \rangle_E \leq r|\psi|_E^2$ if that term is inserted too), effectively replacing $\frac{r}{4}|\beta|_{\Lambda^{0,2}(E)}^4$ by $\frac{r}{4}|\beta|_{\Lambda^{0,2}(E)}^2$, and
\begin{multline*}
  \frac{1}{2}\Delta_g|\beta|_{\Lambda^{0,2}(E)}^2
  + \frac{r}{2}|\alpha|_E^2|\beta|_{\Lambda^{0,2}(E)}^2
  \\
  \leq
  - \Real\langle i\Lambda_\omega F_{A_2}\beta,\beta\rangle_{\Lambda^{0,2}(E)}
  - \frac{1}{2}\Real\langle i\Lambda_\omega F_{A_d}\beta,\beta\rangle_{\Lambda^{0,2}(E)}
  \\
  - \frac{8r}{32}|\beta|_{\Lambda^{0,2}(E)}^2 - |\nabla_{A_2\otimes A}\beta|_{T^3(E)}^2
  + \frac{1}{2}\|N_J\|_{C^0(X)}|\nabla_A\alpha|_{\Lambda^1(E)}|\beta|_{\Lambda^{0,2}(E)} \quad\text{on } X.
\end{multline*}
Applying the interpolation inequality $2xy \leq 2rx^2/32 + 32y^2/(2r)$ to the last term on the right-hand side above yields:
\begin{multline}
  \label{eq:Taubes_1996_SW_to_Gr_eq_2-7_raw}
  \frac{1}{2}\Delta_g|\beta|_{\Lambda^{0,2}(E)}^2
  + \frac{r}{2}|\alpha|_E^2|\beta|_{\Lambda^{0,2}(E)}^2
  \\
  \leq
  - \Real\langle i\Lambda_\omega F_{A_2}\beta,\beta\rangle_{\Lambda^{0,2}(E)}
  - \frac{1}{2}\Real\langle i\Lambda_\omega F_{A_d}\beta,\beta\rangle_{\Lambda^{0,2}(E)}
  \\
  - \frac{6r}{32}|\beta|_{\Lambda^{0,2}(E)}^2 - |\nabla_{A_2\otimes A}\beta|_{T^3(E)}^2
  + \frac{16}{2r}\|N_J\|_{C^0(X)}^2|\nabla_A\alpha|_{\Lambda^1(E)}^2  \quad\text{on } X.
\end{multline}
Now suppose $r$ is chosen large enough that
\begin{equation}
  \label{eq:r_geq_LambdaFAK*+LambdaFAd}
  - \Real\langle i\Lambda_\omega F_{A_2}\beta,\beta\rangle_{\Lambda^{0,2}(E)}
  - \frac{1}{2}\Real\langle i\Lambda_\omega F_{A_d}\beta,\beta\rangle_{\Lambda^{0,2}(E)}
  \leq
  \frac{r}{32}|\beta|_{\Lambda^{0,2}(E)}^2 \quad\text{on } X.
\end{equation}
By combining the inequalities \eqref{eq:r_geq_LambdaFAK*+LambdaFAd} and \eqref{eq:Taubes_1996_SW_to_Gr_eq_2-7_raw} and setting
\begin{equation}
  \label{eq:Taubes_1996_SW_to_Gr_eq_2-7_z3}
  z_3 := 8\|N_J\|_{C^0(X)}^2,  
\end{equation}
we see that we have proved the

\begin{lem}[Differential inequality for the squared pointwise norms of sections of $\Lambda^{0,2}(E)$ with a singular Taubes perturbation]
\label{lem:Taubes_1996_SW_to_Gr_eq_2-7}
Continue the hypotheses of Theorem \ref{mainthm:Taubes_1996_SW_to_Gr_2-3}. If $\alpha \in \Omega^0(E)$ and $\beta \in \Omega^{0,2}(E)$ are defined by $\Phi = r^{1/2}(\alpha,\beta)$ as in \eqref{eq:Taubes_1996_SW_to_Gr_eq_1-19}, then
%PF8-16-2024 Indicate Sobolev notation and inequality a.e. on X.
\begin{equation}
  \label{eq:Taubes_1996_SW_to_Gr_eq_2-7}
  \frac{1}{2}\Delta_g|\beta|_{\Lambda^{0,2}(E)}^2
  + \frac{r}{2}|\alpha|_E^2|\beta|_{\Lambda^{0,2}(E)}^2
  \leq
  - \frac{5r}{32}|\beta|_{\Lambda^{0,2}(E)}^2 - |\nabla_{A_2\otimes A}\beta|_{T^3(E)}^2
  + \frac{z_3}{r}|\nabla_A\alpha|_{\Lambda^1(E)}^2 \quad\text{on } X.
\end{equation}
\end{lem}

The differential inequality \eqref{eq:Taubes_1996_SW_to_Gr_eq_2-7} is \emph{almost exactly} analogous to Taubes \cite[Section 2 (b), Equation (2.7), p. 855]{TauSWGromov}, except that we have a factor $5/32$ rather than $1/8$ multiplying $r|\beta|_{\Lambda^{0,2}(E)}^2$.

\section[Differential inequality for an affine combination of squared pointwise norms]{Differential inequality for an affine combination of squared pointwise norms of sections of $E$ and $\Lambda^{0,2}(E)$ with a singular Taubes perturbation}
\label{sec:Differential_inequality_Taubes_combination_squared_pointwise_norms_spinor_components}
In this section, we prove Lemma \ref{lem:Taubes_1996_SW_to_Gr_2-2} --- an analogue of Taubes \cite[Section 2 (b), Lemma 2.2, p. 855]{TauSWGromov} --- which gives a linear second order elliptic differential inequality for a function $u$ that is an affine linear combination of $|\alpha|_E^2$ and $|\beta|_{\Lambda^{0,2}(E)}^2$when $(A,\varphi,\psi)$ with $(\varphi,\psi) = r^{1/2}(\alpha,\beta)$ is a solution to the system \eqref{eq:SO(3)_monopole_equations_almost_Hermitian_perturbed_intro} of non-Abelian monopole equations with a singular Taubes perturbation.

With $\rho > 0$ a given constant, introduce a function\footnote{This is an exact analogue of Taubes \cite[Section 2 (b), Equation (2.8), p. 855]{TauSWGromov}}:
\begin{equation}
  \label{eq:Taubes_1996_SW_to_Gr_eq_2-8}
  w := 1 - |\alpha|_E^2 + \frac{\rho}{r} \in \Omega^0(X;\RR).
\end{equation}
Inequality \eqref{eq:Taubes_1996_SW_to_Gr_eq_2-3_inequality} and equation \eqref{eq:Taubes_1996_SW_to_Gr_eq_2-8}, giving $|\alpha|_E^2 = 1 - w + \rho/r$ and
\[
  \frac{r}{4}|\alpha|_E^4
  =
  \frac{r}{4}|\alpha|_E^2\cdot |\alpha|_E^2
  =
  \frac{r}{4}|\alpha|_E^2 - \frac{r}{4}|\alpha|_E^2w + \frac{\rho}{4}|\alpha|_E^2 \quad\text{on } X,
\]  
imply the following inequality for $w$:
\begin{multline*}
  -\frac{1}{2}\Delta_gw + |\nabla_A\alpha|_E^2 + \frac{r}{4}|\alpha|_E^2 - \frac{r}{4}|\alpha|_E^2 w
  + \frac{\rho}{4}|\alpha|_E^2 - \frac{r}{4}|\alpha|_E^2
  + \frac{r}{4}|\beta|_{\Lambda^{0,2}(E)}^2|\alpha|_E^2
  \\
  - \Real\langle i\Lambda_\omega F_{A_d}\alpha,\alpha\rangle_E
  - 2\Real\langle(\bar\mu\partial_A)^*\beta,\alpha\rangle_E \leq 0 \quad\text{on } X,
\end{multline*}
and thus
\begin{multline}
  \label{eq:Taubes_1996_SW_to_Gr_eq_2-9_inequality}
  \frac{1}{2}\Delta_gw + \frac{r}{4}|\alpha|_E^2w
  \geq
  |\nabla_A\alpha|_E^2 + \frac{r}{4}|\beta|_{\Lambda^{0,2}(E)}^2|\alpha|_E^2 + \frac{\rho}{4}|\alpha|_E^2
  \\
  - \Real\langle i\Lambda_\omega F_{A_d}\alpha,\alpha\rangle_E
  - 2\Real\langle(\bar\mu\partial_A)^*\beta,\alpha\rangle_E \quad\text{on } X.
\end{multline}
The preceding differential inequality \eqref{eq:Taubes_1996_SW_to_Gr_eq_2-9_inequality} is \emph{almost exactly} analogous to Taubes \cite[Section 2 (b), Equation (2.10), p. 855]{TauSWGromov}.

Observe that,
%PF8-16-2024
%[TODO Recheck],
using $2xy \leq 16x^2/\rho + \rho y^2/16$ and $\bar\mu = \frac{1}{4}N_J$,
\begin{multline*}
  |2\Real\langle(\bar\mu\partial_A)^*\beta,\alpha\rangle_E|
  \leq
  \frac{1}{2}\left(\|N_J\|_{C^0(X)}|\nabla_{A_2\otimes A}\beta|_{T^3(E)} + \|\nabla N_J\|_{C^0(X)}|\beta|_{\Lambda^{0,2}(E)}\right)|\alpha|_E
  \\
  \leq
  \frac{8}{\rho}\|N_J\|_{C^0(X)}^2|\nabla_{A_2\otimes A}\beta|_{T^3(E)}^2
  + \frac{8}{\rho}\|\nabla N_J\|_{C^0(X)}^2|\beta|_{\Lambda^{0,2}(E)}^2
  + \frac{\rho}{16}|\alpha|_E^2 \quad\text{on } X.
\end{multline*}
Also, for $\rho$ chosen large enough that
\begin{equation}
  \label{eq:rho_geq_16LambdaFAd}
  \frac{\rho}{16} \geq \|\Lambda_\omega F_{A_d}\|_{C^0(X)} \quad\text{on } X,
\end{equation}
we obtain
\[
  \left|\Real\langle i\Lambda_\omega F_{A_d}\alpha,\alpha\rangle_E\right|
  \leq
  \|\Lambda_\omega F_{A_d}\|_{C^0(X)} |\alpha|_E^2
  \leq
  \frac{\rho}{16}|\alpha|_E^2 \quad\text{on } X.
\]  
Hence,  
\begin{multline}
  \label{eq:Taubes_1996_SW_to_Gr_eq_2-10_raw}
  \frac{1}{2}\Delta_gw + \frac{r}{4}|\alpha|_E^2w
  \geq |\nabla_A\alpha|_E^2
  + \frac{r}{4}|\beta|_{\Lambda^{0,2}(E)}^2|\alpha|_E^2 + \frac{\rho}{8}|\alpha|_E^2
  \\
  - \frac{8}{\rho}\left(\|N_J\|_{C^0(X)}^2|\nabla_{A_2\otimes A}\beta|_{T^3(E)}^2
    + \|\nabla N_J\|_{C^0(X)}^2|\beta|_{\Lambda^{0,2}(E)}^2\right) \quad\text{on } X.
\end{multline}
The preceding differential inequality \eqref{eq:Taubes_1996_SW_to_Gr_eq_2-10_raw} is \emph{almost exactly} analogous to Taubes \cite[Section 2 (b), Equation (2.10), p. 855]{TauSWGromov}, with
\begin{equation}
  \label{eq:Taubes_1996_SW_to_Gr_eq_2-10_z4}
  z_4 := 8\max\left\{\|N_J\|_{C^0(X)}^2, \|\nabla N_J\|_{C^0(X)}^2\right\}.
\end{equation}
(The only difference is that Taubes appears to have unintentionally omitted the term $\rho|\alpha|_E^2/8$ in a typographical error.) Substituting the expression \eqref{eq:Taubes_1996_SW_to_Gr_eq_2-10_z4} for $z_4$ into \eqref{eq:Taubes_1996_SW_to_Gr_eq_2-10_raw} yields
\begin{multline}
  \label{eq:Taubes_1996_SW_to_Gr_eq_2-10}
  \frac{1}{2}\Delta_gw + \frac{r}{4}|\alpha|_E^2w
  \geq |\nabla_A\alpha|_E^2
  + \frac{r}{4}|\beta|_{\Lambda^{0,2}(E)}^2|\alpha|_E^2 + \frac{\rho}{8}|\alpha|_E^2
  \\
  - \frac{z_4}{\rho}\left(|\nabla_{A_2\otimes A}\beta|_{T^3(E)}^2
    + |\beta|_{\Lambda^{0,2}(E)}^2\right) \quad\text{on } X.
\end{multline}
We thus obtain the following exact analogue of Taubes \cite[Section 2 (b), Lemma 2.2, p. 855]{TauSWGromov}:

\begin{lem}[Differential inequality for an affine combination of the squared pointwise norms of sections of $E$ and $\Lambda^{0,2}(E)$ with singular a Taubes perturbation]
\label{lem:Taubes_1996_SW_to_Gr_2-2}
Let $(X,g,J,\omega)$ be a
%PF7-10-2024 Would almost Hermitian be okay?
symplectic four-manifold, $(E,H)$ be a smooth, Hermitian vector bundle over $X$ with complex rank $2$, and $A_d$ be a smooth, unitary connection on $\det E$. Define
\begin{subequations}
\label{eq:kappa1_2_3}  
\begin{align}
  \label{eq:kappa1}
  \kappa_1 &:= \max\left\{32\|\Lambda_\omega F_{A_2}\|_{C^0(X)} + 16\|\Lambda_\omega F_{A_d}\|_{C^0(X)}, z_3, 1\right\},
    \\
  \label{eq:kappa3}
  \kappa_3 &:= \max\left\{\frac{32z_4}{\kappa_1}, 2z_4, 1\right\},
  \\
  \label{eq:kappa2}
  \kappa_2 &:= \max\left\{\frac{8}{\kappa_3}\|\Lambda_\omega F_{A_d}\|_{C^0(X)}, 1\right\}, 
\end{align}
\end{subequations}
where the constants $z_3, z_4 \in (0,\infty)$ are as in \eqref{eq:Taubes_1996_SW_to_Gr_eq_2-7_z3}, \eqref{eq:Taubes_1996_SW_to_Gr_eq_2-10_z4}, respectively, and depend at most on $(g,J)$. Let $r,\delta,\zeta \in (0,\infty)$ be constants that obey
\begin{equation}
  \label{eq:Bounds_r_zeta_regular}
  r \geq \kappa_1, \quad \delta \geq \kappa_3, \quad\text{and}\quad
  0 < \zeta \leq \frac{\delta}{2\kappa_2\kappa_3}, 
\end{equation}
and let $(A,\Phi) \in \sA(E,H,A_d) \times \Omega^0(W_\can^+\otimes E)$ be a smooth
%PF7-9-2024 Add reference below
% (see equations \eqref{eq:Canonical_spinc_bundles} for $W_\can$)
%PF7-9-2024 Define "smooth"
solution to the system \eqref{eq:SO(3)_monopole_equations_almost_Hermitian_perturbed_intro} of non-Abelian monopole equations with a singular Taubes perturbation.
%PF10-10-2024 Recheck
If $\Phi = (\varphi,\psi) = r^{1/2}(\alpha,\beta) \in \Omega^0(E)\oplus\Omega^{0,2}(E)$ as in \eqref{eq:Taubes_1996_SW_to_Gr_eq_1-19} then the function\footnote{This is an exact analogue of Taubes \cite[Section 2 (b), Equation (2.5), p. 855]{TauSWGromov}}
\begin{equation}
  \label{eq:Taubes_1996_SW_to_Gr_eq_2-5}
  u := 1 - |\alpha|_E^2 - \zeta|\beta|_{\Lambda^{0,2}(E)}^2 + \frac{\delta}{\zeta r} \in \Omega^0(X;\RR)
\end{equation}
obeys the following differential inequality:
\begin{equation}
  \label{eq:Taubes_1996_SW_to_Gr_eq_2-6}
  \frac{1}{2}\Delta_gu + \frac{r}{4}|\alpha|_E^2u
  \geq
  \left(1 - \zeta\frac{\kappa_1}{r}\right)|\nabla_A\alpha|_E^2
  + \frac{\zeta r}{8}|\beta|_{\Lambda^{0,2}(E)}^2
  + \frac{\zeta}{2}|\nabla_{A_2\otimes A}\beta|_{T^3(E)}^2 
  + \frac{\delta}{8\zeta}|\alpha|_E^2 \quad\text{on } X.
\end{equation}
\end{lem}

\begin{rmk}[On the upper bound for the constant $\zeta$ in \eqref{eq:Bounds_r_zeta_regular}]
In \cite[Section 2 (b), Lemma 2.2, p. 855]{TauSWGromov}, Taubes allows $\zeta$ to obey inequality $\zeta < r/(2\kappa_1\kappa_2)$, which is equivalent\footnote{Aside from our replacement of $<$ by $\leq$.} to our upper bound for $\zeta$ in \eqref{eq:Bounds_r_zeta_regular} when $\delta = r\kappa_3/\kappa_1$, noting that $r \geq \kappa_1$ by our hypothesis \eqref{eq:Bounds_r_zeta_regular} and so $r\kappa_3/\kappa_1 \geq \kappa_3$.
\end{rmk}  

\begin{proof}[Proof of Lemma \ref{lem:Taubes_1996_SW_to_Gr_2-2}]
By our hypothesis \eqref{eq:Bounds_r_zeta_regular} that
\[
  r \geq \kappa_1,
\]
for $\kappa_1$ as in \eqref{eq:kappa1}, the constant $r$ necessarily obeys the weaker inequality \eqref{eq:r_geq_LambdaFAK*+LambdaFAd}. Note that, for 
\begin{equation}
  \label{eq:rho_equals_delta_over_zeta}
  \rho := \frac{\delta}{\zeta},
\end{equation}
the equation \eqref{eq:Taubes_1996_SW_to_Gr_eq_2-8} for $w$ gives
\[
  w = 1 - |\alpha|_E^2 + \frac{\rho}{r} = 1 - |\alpha|_E^2 + \frac{\delta}{\zeta r}
\]  
and so the equation \eqref{eq:Taubes_1996_SW_to_Gr_eq_2-5} for $u$ yields
\[
  u = 1 - |\alpha|_E^2 + \frac{\delta}{\zeta r} - \zeta|\beta|_{\Lambda^{0,2}(E)}^2 
  = w - \zeta|\beta|_{\Lambda^{0,2}(E)}^2.
\]
% PF10-10-2024 Check that we're applying correct inequality
We see that the function
$-\zeta|\beta|_{\Lambda^{0,2}(E)}^2$ obeys
\begin{align*}
  &\frac{1}{2}\Delta_g\left(-\zeta|\beta|_{\Lambda^{0,2}(E)}^2\right)
    + \frac{r}{4}|\alpha|_E^2\left(-\zeta|\beta|_{\Lambda^{0,2}(E)}^2\right)
  \\
  &\quad \geq
  -\zeta\left( \frac{1}{2}\Delta_g|\beta|_{\Lambda^{0,2}(E)}^2 + \frac{r}{2}|\alpha|_E^2|\beta|_{\Lambda^{0,2}(E)}^2 \right)
  \\
  &\quad \geq
  -\zeta\left(- \frac{5r}{32}|\beta|_{\Lambda^{0,2}(E)}^2 - |\nabla_{A_2\otimes A}\beta|_{T^3(E)}^2
    + \frac{z_3}{r}|\nabla_A\alpha|_{\Lambda^1(E)}^2\right) \quad\text{on } X
    \quad\text{(by \eqref{eq:Taubes_1996_SW_to_Gr_eq_2-7})},
\end{align*}
so that
\begin{multline*}
  \frac{1}{2}\Delta_g\left(-\zeta|\beta|_{\Lambda^{0,2}(E)}^2\right)
  + \frac{r}{4}|\alpha|_E^2\left(-\zeta|\beta|_{\Lambda^{0,2}(E)}^2\right)
  \\
  \geq
  \frac{5\zeta r}{32}|\beta|_{\Lambda^{0,2}(E)}^2 + \zeta|\nabla_{A_2\otimes A}\beta|_{T^3(E)}^2
    - \frac{\zeta z_3}{r}|\nabla_A\alpha|_{\Lambda^1(E)}^2 \quad\text{on } X.
\end{multline*}
Combining the preceding inequality with \eqref{eq:Taubes_1996_SW_to_Gr_eq_2-10} and substituting $\rho = \delta/\zeta$ from \eqref{eq:rho_equals_delta_over_zeta} gives 
\begin{align*}
  \frac{1}{2}\Delta_gu + \frac{r}{4}|\alpha|_E^2u
  &\geq
  |\nabla_A\alpha|_E^2
  + \frac{r}{4}|\beta|_{\Lambda^{0,2}(E)}^2|\alpha|_E^2
  - \frac{z_4}{\rho}\left(|\nabla_{A_2\otimes A}\beta|_{T^3(E)}^2
    + |\beta|_{\Lambda^{0,2}(E)}^2\right)
    + \frac{\delta}{8\zeta}|\alpha|_E^2
  \\
  &\quad + \frac{5\zeta r}{32}|\beta|_{\Lambda^{0,2}(E)}^2 + \zeta|\nabla_{A_2\otimes A}\beta|_{T^3(E)}^2
    - \frac{\zeta z_3}{r}|\nabla_A\alpha|_{\Lambda^1(E)}^2 \quad\text{on } X,
\end{align*}
and hence a precursor to an analogue of Taubes \cite[Section 2 (b), Equation (2.6), p. 855]{TauSWGromov}:
\begin{multline}
  \label{eq:Taubes_1996_SW_to_Gr_eq_2-6_raw}
  \frac{1}{2}\Delta_gu + \frac{r}{4}|\alpha|_E^2u
  \geq
    \left(1 - \frac{\zeta z_3}{r}\right)|\nabla_A\alpha|_{\Lambda^1(E)}^2
    + \left(\frac{5\zeta r}{32} - \frac{z_4}{\rho}\right)|\beta|_{\Lambda^{0,2}(E)}^2
    + \left(\zeta - \frac{z_4}{\rho}\right)|\nabla_{A_2\otimes A}\beta|_{T^3(E)}^2
    \\
    + \frac{\delta}{8\zeta}|\alpha|_E^2 \quad\text{on } X,
\end{multline}
where we used the inequality $r|\beta|_{\Lambda^{0,2}(E)}^2|\alpha|_E^2/4 \geq 0$ to obtain \eqref{eq:Taubes_1996_SW_to_Gr_eq_2-6_raw}.

The differential inequality \eqref{eq:Taubes_1996_SW_to_Gr_eq_2-6_raw} has the desired coefficient of $|\alpha|_E^2$. In our application to the proof of the forthcoming Theorem \ref{mainthm:Taubes_1996_SW_to_Gr_2-3}, we shall need the right-hand side of this differential inequality to be non-negative for it to be effective. For $\kappa_1$ as in \eqref{eq:kappa1}, we have $\kappa_1 \geq z_3$ and so replacing $z_3$ by $\kappa_1$ gives the desired coefficient of the $|\nabla_A\alpha|_{\Lambda^1(E)}^2$ term.

Noting that $\rho = \delta/\zeta $ by \eqref{eq:rho_equals_delta_over_zeta}, the coefficient of $|\beta|_{\Lambda^{0,2}(E)}^2$ is
\[
  \left(\frac{5\zeta r}{32} - \frac{z_4}{\rho}\right)
  =
  \left(\frac{5\zeta r}{32} - \frac{z_4\zeta}{\delta}\right).
\]
We now require that
\begin{equation}
  \label{eq:delta_leq_r_over_32}
  \delta \geq \frac{32z_4}{r},
\end{equation}
so that $1/\delta \leq r/(32z_4)$ and $-1/\delta \geq -r/(32z_4)$, and thus
\[
  \left(\frac{5\zeta r}{32} - \frac{z_4\zeta}{\delta}\right)
  \geq
  \left(\frac{5\zeta r}{32} - \frac{z_4\zeta r}{32z_4}\right)
  =
  \left(\frac{5\zeta r}{32} - \frac{\zeta r}{32}\right)
  =
  \frac{4\zeta r}{32}
  =
  \frac{\zeta r}{8}.
\]
Hence, replacing $\delta$ by $32z_4/r$ gives the desired coefficient of $|\beta|_{\Lambda^{0,2}(E)}^2$.

We substitute $\rho = \delta/\zeta$ from \eqref{eq:rho_equals_delta_over_zeta} into the coefficient of $|\nabla_{A_2\otimes A}\beta|_{T^3(E)}^2$ to give
\[
  \zeta - \frac{z_4}{\rho}
  =
  \zeta - \frac{z_4\zeta}{\delta}.
\]
In addition to our requirement that $\delta$ obey \eqref{eq:delta_leq_r_over_32}, we shall also require that
\begin{equation}
  \label{eq:delta_geq_2_over_z4}
  \delta \geq 2z_4,
\end{equation}
so that $1/\delta \leq 1/(2z_4)$ and $-1/\delta \geq -1/(2z_4)$, and thus
\[
  \zeta - \frac{z_4\zeta}{\delta}
  \geq
  \zeta - \frac{z_4\zeta}{2z_4}
  =
  \zeta - \frac{\zeta}{2}
  =
  \frac{\zeta}{2}.
\]
Thus, replacing $\delta$ by $2z_4$ gives the desired coefficient of $|\nabla_{A_2\otimes A}\beta|_{T^3(E)}^2$.

Observe that $\delta$ obeys \emph{both} \eqref{eq:delta_leq_r_over_32} \emph{and} \eqref{eq:delta_geq_2_over_z4} by our \emph{hypothesis} \eqref{eq:Bounds_r_zeta_regular} that
\[
  \delta \geq \kappa_3,
\]
for $\kappa_3$ as in \eqref{eq:kappa3}, because $r \geq \kappa_1$ by our \emph{hypothesis} \eqref{eq:Bounds_r_zeta_regular} and so
\[
  \kappa_3
  =
  \max\left\{\frac{32z_4}{\kappa_1}, 2z_4, 1\right\}
  \geq
  \max\left\{\frac{32z_4}{r}, 2z_4, 1\right\}.
\]  
The constant $\rho = \delta/\zeta$ from \eqref{eq:rho_equals_delta_over_zeta} obeys the \emph{required inequality} \eqref{eq:rho_geq_16LambdaFAd} if
\[
  \rho = \frac{\delta}{\zeta} \geq 16\|\Lambda_\omega F_{A_d}\|_{C^0(X)},
\]
which is equivalent to requiring that $\zeta > 0$ obey
\[
  \zeta \leq \frac{\delta}{16\|\Lambda_\omega F_{A_d}\|_{C^0(X)}},
\]
if $\Lambda_\omega F_{A_d} \neq 0$ and $\zeta < \infty$ otherwise. The preceding inequality is assured by our \emph{hypothesis} \eqref{eq:Bounds_r_zeta_regular} that
\[
  \zeta
  \leq
  \frac{\delta}{2\kappa_2\kappa_3},
\]
for $\kappa_2$ as in \eqref{eq:kappa2}, namely
\[
   \kappa_2 = \max\left\{\frac{8}{\kappa_3}\|\Lambda_\omega F_{A_d}\|_{C^0(X)}, 1\right\},
\]
because $\delta \geq \kappa_3$ by our \emph{hypothesis} \eqref{eq:Bounds_r_zeta_regular} and thus
\[
  \frac{\delta}{2\kappa_2\kappa_3}
  =
  \frac{\delta}{2\kappa_3\max\left\{8\|\Lambda_\omega F_{A_d}\|_{C^0(X)}/\kappa_3, 1\right\}}
  \leq
  \frac{\delta}{16\|\Lambda_\omega F_{A_d}\|_{C^0(X)}}.
\]
This completes the proof of Lemma \ref{lem:Taubes_1996_SW_to_Gr_2-2}.
\end{proof}

\section[Pointwise estimate for the squared pointwise norms of sections of $\Lambda^{0,2}(E)$]{Pointwise estimate for  the squared pointwise norms of sections of $\Lambda^{0,2}(E)$ with a singular Taubes perturbation}
\label{sec:Pointwise_estimate_squared_pointwise_norm_anti-holomorphic_spinor}
Using Lemma \ref{lem:Taubes_1996_SW_to_Gr_2-2}, we finally establish Theorem \ref{mainthm:Taubes_1996_SW_to_Gr_2-3}, which gives a pointwise bound for $|\beta|_{\Lambda^{0,2}(E)}^2$when $(A,\varphi,\psi)$ with $(\varphi,\psi) = r^{1/2}(\alpha,\beta)$ is a solution to the system \eqref{eq:SO(3)_monopole_equations_almost_Hermitian_perturbed_intro} of non-Abelian monopole equations with a singular Taubes perturbation. This result is an almost exact analogue of Taubes \cite[Section 2 (b), Proposition 2.3, p. 856]{TauSWGromov}. While our arguments are broadly similar to those of Taubes, we give two overlapping proofs. Before proceeding to them, we make the following remarks.

\begin{rmk}[Constraints on the Chern classes of $E$ when $\beta$ has constant positive pointwise norm]
If $|\beta|_{\Lambda^{0,2}(E)}$ is a positive constant, then $\beta \in \Omega^{0,2}(E)$ is a section that has no zeros on $X$. For a vector bundle $E$ with complex rank $2$, this implies
%PF12-3-2024 Provide reference
that $c_2(E\otimes K_X^*) = 0$ where, by the splitting principle, we have the equalities
\[
  e(E\otimes K_X^*)
  =
  c_2(E\otimes K_X^*)
  =
  c_2(E)-c_1(E)\cdot c_1(K_X) + c_1(K_X)^2,
\]
where $e$ denotes the Euler class. Hence, the case where $|\beta|_{\Lambda^{0,2}(E)}$ is a positive constant can occur only if $c_2(E) \neq c_1(E)\cdot c_1(K_X) - c_1(K_X)^2$.
\qed
\end{rmk}

%PF8-16-2024
%TODO - Paul will make remark below into a corollary.

\begin{rmk}[On the term $z'/r$ in the inequality \eqref{eq:Taubes_1996_SW_to_Gr_eq_2-11}]
In \cite[Section 2 (b), Proposition 2.3, p. 856]{TauSWGromov}, Taubes obtains the stronger inequality \cite[Section 2 (b), Proposition 2.3, p. 856]{TauSWGromov}, with $z'/r$ in our inequality \eqref{eq:Taubes_1996_SW_to_Gr_eq_2-11} replaced by $z'/r^2$. The reason for the difference appears in our proof of Lemma \ref{lem:Taubes_1996_SW_to_Gr_2-2}, where we must require that
\[
  \frac{\delta}{\zeta} \geq 16\|\Lambda_\omega F_{A_d}\|_{C^0(X)}.
\]
In Taubes \cite{TauSWGromov}, there is no analogue of the unitary connection $A_d$ on $\det E$ and so the preceding constraint would be replaced by $\delta/\zeta > 0$, which is obviously satisfied by Taubes' choice of $\delta = \kappa_3$ and $\zeta = r/(2\kappa_1\kappa_2)$ in his proof of \cite[Section 2 (b), Proposition 2.3, p. 856]{TauSWGromov}. The preceding uniform positive lower bound on $\delta/\zeta$ could be replaced by $\delta/\zeta > 0$ if we knew that $i\Lambda_\omega F_{A_d} \leq 0$ on $X$ and thus $- \Real\langle i\Lambda_\omega F_{A_d}\alpha,\alpha\rangle_E \geq 0$ on $X$ in the right-hand side of the inequality \eqref{eq:Taubes_1996_SW_to_Gr_eq_2-9_inequality}, and so that term could be replaced by zero.
% PF7-11-2024 Recheck. How does Taubes obtain his upper bound on \zeta?

To construct a unitary connection $A_d$ on $\det E$ with $i\Lambda_\omega F_{A_d} \leq 0$ on $X$, we assume that (see Feehan and Leness 
%TL12-4-2025: Updated and removed page references
\cite[Equations (8.4.14) and (8.4.15)]{Feehan_Leness_introduction_virtual_morse_theory_so3_monopoles})
\[
  \deg_\omega E
  := c_1(E)\cdot[\omega]
  = \frac{i}{2\pi}\int_X(\tr_E F_A)\wedge\omega^{n-1}
  = \frac{i}{2\pi}\int_X F_{A_d}\wedge\omega^{n-1}
  = c_1(\det E)\cdot[\omega]
  \leq 0,
\]
and so $E$ has slope 
%TL12-4-2025: Updated and removed page references
\cite[Equation (11.2.1)]{Feehan_Leness_introduction_virtual_morse_theory_so3_monopoles}
\[
  \mu_\omega(E) := \frac{\deg_\omega E}{\rank E} \leq 0.
\]
The Einstein factor for $E$ is thus 
%TL12-4-2025: Updated and removed page references
\cite[Equation (11.3.5)]{Feehan_Leness_introduction_virtual_morse_theory_so3_monopoles}
\[
  \lambda_\omega(E) := \frac{2\pi}{(n-1)!\vol_\omega X}\mu_\omega(E) \leq 0.
\]
We now aim to solve the Hermitian--Einstein equation 
%TL12-4-2025: Updated and removed page references
\cite[Equation (11.3.2)]{Feehan_Leness_introduction_virtual_morse_theory_so3_monopoles} for $A_d$:
\[
  \Lambda_\omega F_{A_d} = -i\lambda_\omega(E),
\]
equivalently,
\[
  i\Lambda_\omega F_{A_d} = \lambda_\omega(E).
\]
By construction, we thus obtain
\[
  i\Lambda_\omega F_{A_d} \leq 0 \quad\text{on } X.
\]
If $(X,g,J,\omega)$ is a complex K\"ahler $n$-manifold and $\det E$ has an integrable $(0,1)$-connection $\bar\partial_{\det E}$, then there is a unique unitary connection $A_d$ on $\det E$ (up to $\Map(X,S^1)$ gauge equivalence) such that $\bar\partial_{A_d} = \bar\partial_{\det E}$ and $i\Lambda_\omega F_{A_d} = \lambda_\omega(E)$. Indeed, see, for example, Donaldson \cite[Theorem 1, p. 2]{DonASD} when $X$ is a complex projective surface, Donaldson \cite[Proposition 1, p. 231]{DonInfDet} when $X$ is a complex projective $n$-manifold, Donaldson and Kronheimer \cite[Section 6.1.4]{DK} for a survey of results when $X$ is a complex K\"ahler $n$-manifold, Li and Yau \cite{Li_Yau_1987} when $X$ is a complex Hermitian $n$-manifold (see L\"ubke and Teleman \cite[Theorem 3.0.1, p. 61]{Lubke_Teleman_1995} for an exposition of the proof), and Uhlenbeck and Yau \cite[Section 1, Theorem, p. S262]{Uhlenbeck_Yau_1986} when $X$ is a complex K\"ahler $n$-manifold. If $(X,g,J,\omega)$ is an almost Hermitian $2n$-manifold, the same conclusion follows from De Bartolomeis and Tian \cite[Theorem 0.1, p. 232 or Theorem 4.3, p. 253; Section 3, p. 244]{DeBartolomeis_Tian_1996}, without the assumption that $\bar\partial_{\det E}$ is integrable.
\qed
\end{rmk}

Our proof of Theorem \ref{mainthm:Taubes_1996_SW_to_Gr_2-3} is slightly more complicated than Taubes' proof of his \cite[Section 2 (b), Proposition 2.3, p. 856]{TauSWGromov} since we cannot assume, even after applying a $W^{2,p}(\SU(E))$ gauge transformation, that our solution $(A,\Phi)$ to the system \eqref{eq:SO(3)_monopole_equations_almost_Hermitian_perturbed_intro} of non-Abelian equations with a singular Taubes perturbation is $C^\infty$ over $X$ and hence that $u$ in \eqref{eq:Taubes_1996_SW_to_Gr_eq_2-5} is in $C^\infty(X)$. We can only assume that $A$ is $W^{1,p}$ and $\Phi$ is $W^{2,p}$ for $4\leq p<\infty$, and so $u$ in \eqref{eq:Taubes_1996_SW_to_Gr_eq_2-5} is in $W^{2,p}(X)$ for $4\leq p<\infty$. We give two different proofs of Theorem \ref{mainthm:Taubes_1996_SW_to_Gr_2-3}:
\begin{enumerate}
\item We adapt Taubes' proof of \cite[Section 2 (b), Proposition 2.3, p. 856]{TauSWGromov} under the technical assumptions that $X$ is connected and $\Lambda^{0,2}(E) = \Lambda^{0,2}(X)\otimes E$ has a non-vanishing Euler class, and apply the strong maximum principle provided by Theorem \ref{thm:Gilbarg_Trudinger_3-5_and_9-6} to a linear second order elliptic differential inequality on $X$.
\item We apply the weak maximum principle provided by Theorem \ref{thm:Gilbarg_Trudinger_3-1_and_9-1} to a linear second order elliptic differential inequality on $X_0$, where $X_0 = \{x\in X: \psi(x) \neq 0\}$ by \eqref{eq:X0}.
\end{enumerate}

\begin{proof}[First proof of Theorem \ref{mainthm:Taubes_1996_SW_to_Gr_2-3}]
We make the technical assumptions that $X$ is connected and that $\Lambda^{0,2}(E)$ has a non-vanishing Euler class.

Since $\alpha \in W^{2,p}(E)$ and $\beta \in W^{2,p}(\Lambda^{0,2}(E))$ by hypothesis, the definition \eqref{eq:Taubes_1996_SW_to_Gr_eq_2-5} yields $u \in W^{2,p}(X)$ (with $p\geq 4$). Noting the bounds \eqref{eq:Bounds_r_zeta_regular} on $r$, $\delta$, and $\zeta$ in Lemma \ref{lem:Taubes_1996_SW_to_Gr_2-2}, we choose
\begin{subequations}
\label{eq:delta_zeta_equalities}        
\begin{align}
  \label{eq:delta=r_kappa3_over_kappa1}
  \delta &:= \frac{r\kappa_3}{\kappa_1} \geq \kappa_3,
  \\
  \label{eq:zeta=r_over_2kappa1_kappa2}
  \zeta &:= \frac{\delta}{2\kappa_2\kappa_3} = \frac{r}{2\kappa_1\kappa_2}.
\end{align}
\end{subequations}
The second expression for $\zeta$ in \eqref{eq:zeta=r_over_2kappa1_kappa2} yields
\[
  \frac{\zeta\kappa_1}{r} = \frac{1}{2\kappa_2} \leq \frac{1}{2},
\]  
where the final inequality follows from the fact \eqref{eq:kappa2} that $\kappa_2 \geq 1$. We thus obtain
\[
  1 - \frac{\zeta\kappa_1}{r} \geq \frac{1}{2}.
\]
Consequently, the differential inequality \eqref{eq:Taubes_1996_SW_to_Gr_eq_2-6} implies that
\[
  \frac{1}{2}\Delta_gu + \frac{r}{4}|\alpha|_E^2u \geq 0 \quad\text{on } X,
\]
and therefore,
\begin{equation}
  \label{eq:Lu_leq_0}
  Lu \leq 0 \quad\text{on } X,
\end{equation}
where we define the linear second-order elliptic operator by
\begin{equation}
  \label{eq:Elliptic_linear_scalar_second-order_operator}
  L :=  -\Delta_g - \frac{r}{2}|\alpha|_E^2 \quad\text{on } W^{2,p}(X),
\end{equation}
which matches our definition \eqref{eq:Second_order_linear_elliptic_operator_riemannian_manifold} with $b = 0$ and $c = -\frac{r}{2}|\alpha|_E^2 \leq 0$ on $X$. Because $X$ is compact, the function $u$ achieves a minimum value $u(x_0)$ at some point $x_0 \in X$. Note that $X$ is connected by our technical assumptions for this proof. If $u$ is \emph{non-constant} on $X$, then the strong maximum principle provided by Theorem \ref{thm:Gilbarg_Trudinger_3-5_and_9-6} implies that $u(x_0) > 0$ (since the possibility that $u(x_0) \leq 0$ is excluded by Theorem \ref{thm:Gilbarg_Trudinger_3-5_and_9-6}). Therefore, $u \geq u(x_0) > 0$ on $X$ and the definition \eqref{eq:Taubes_1996_SW_to_Gr_eq_2-5} of $u$ now implies that
\[
  u = 1 - |\alpha|_E^2 - \zeta|\beta|_{\Lambda^{0,2}(E)}^2 + \frac{\delta}{\zeta r} > 0 \quad\text{on } X,
\]
that is,
\begin{equation}
  \label{eq:Squared_pointwise_norm_beta_strict_upper_bound}
  |\beta|_{\Lambda^{0,2}(E)}^2
  < \frac{1}{\zeta}\left(1 - |\alpha|_E^2 + \frac{\delta}{\zeta r}\right) \quad\text{on } X.
\end{equation}
If $u$ is \emph{constant} on $X$, then we consider two subcases: $\alpha \not\equiv 0$ or $\alpha \equiv 0$ on $X$. If $\alpha \not\equiv 0$ on $X$, then there is a point
%PF12-3-2024 Contradicts index placement for local coordinates
$x_1 \in X$ such that $\alpha(x_1) \neq 0$ and so $c(x_1) = - \frac{r}{4}|\alpha(x_1)|_E^2 < 0$. If $u$ is constant on $X$ and $u < 0$, then $Lu(x_1) = c(x_1)u > 0$, contradicting the fact \eqref{eq:Lu_leq_0} that $Lu \leq 0$ on $X$. Hence, if $u$ is constant on $X$, then $u \equiv 0$ on $X$ and the definition \eqref{eq:Taubes_1996_SW_to_Gr_eq_2-5} of $u$ implies that
\[
  u = 1 - |\alpha|_E^2 - \zeta|\beta|_{\Lambda^{0,2}(E)}^2 + \frac{\delta}{\zeta r} = 0 \quad\text{on } X,
\]
that is,
\begin{equation}
  \label{eq:Squared_pointwise_norm_beta_equality}
  |\beta|_{\Lambda^{0,2}(E)}^2
  = \frac{1}{\zeta}\left(1 - |\alpha|_E^2 + \frac{\delta}{\zeta r}\right) \quad\text{on } X.
\end{equation}
Substituting $\delta = r\kappa_3/\kappa_1$ and $\zeta = r/(2\kappa_1\kappa_2)$ from \eqref{eq:delta_zeta_equalities} yields
\[
  \frac{1}{\zeta}\left(1 - |\alpha|_E^2 + \frac{\delta}{\zeta r}\right)
  =
  \frac{2\kappa_1\kappa_2}{r}\left(1 - |\alpha|_E^2 + \frac{r\kappa_3}{\kappa_1}\frac{2\kappa_1\kappa_2}{r^2}\right)
  =
  \frac{2\kappa_1\kappa_2}{r}\left(1 - |\alpha|_E^2 + \frac{2\kappa_2\kappa_3}{r}\right).
\]  
Hence, combining \eqref{eq:Squared_pointwise_norm_beta_strict_upper_bound}, for $u$ non-constant, and \eqref{eq:Squared_pointwise_norm_beta_equality}, for $u$ constant and $\alpha\not\equiv 0$, gives
\[
  |\beta|_{\Lambda^{0,2}(E)}^2
  \leq \frac{2\kappa_1\kappa_2}{r}\left(1 - |\alpha|_E^2 + \frac{2\kappa_2\kappa_3}{r}\right)
   \quad\text{on } X,
\]
and the conclusion \eqref{eq:Taubes_1996_SW_to_Gr_eq_2-11} follows, with $z'' := 2\kappa_1\kappa_2$ and $z' := 2\kappa_2\kappa_3$, for the cases
\begin{inparaenum}[\itshape i\upshape)]
\item $u$ non-constant, or  
\item $u$ constant and $\alpha\not\equiv 0$.
\end{inparaenum}
(Of course, the first case includes $\alpha \equiv 0$ and $|\beta|_{\Lambda^{0,2}(E)}$ non-constant.) 

Finally, if $u$ is constant on $X$ and $\alpha \equiv 0$ on $X$, then $|\beta|_{\Lambda^{0,2}(E)}$ is constant on $X$. If $\beta\not\equiv 0$, then $|\beta|_{\Lambda^{0,2}(E)}$ is a positive constant, which contradicts our
%PF11-26-2024 Insert reference
technical assumption that $\Lambda^{0,2}(E)$ has a non-vanishing Euler class. Therefore, $\beta\equiv 0$ and so the conclusion \eqref{eq:Taubes_1996_SW_to_Gr_eq_2-11} again follows. This completes our first proof of Theorem \ref{mainthm:Taubes_1996_SW_to_Gr_2-3}.
\end{proof}

\begin{proof}[Second proof of Theorem \ref{mainthm:Taubes_1996_SW_to_Gr_2-3}]
Recall that $X_0 = \{x\in X: \psi(x) \neq 0\}$ by \eqref{eq:X0}. We again observe that because $\alpha \in W^{2,p}(E)$ and $\beta \in W^{2,p}(\Lambda^{0,2}(E))$ by hypothesis, the definition \eqref{eq:Taubes_1996_SW_to_Gr_eq_2-5} yields $u \in W^{2,p}(X)$ (with $p\geq 4$) and
\[
  u = 1 - |\alpha|_E^2 + \frac{\delta}{\zeta r} \quad\text{on } \partial X_0,
\]
since $\beta = 0$ on $X\less X_0$. Because $Lu \leq 0$ on $X_0$, the weak maximum principle provided by Theorem \ref{thm:Gilbarg_Trudinger_3-1_and_9-1} for functions $u \in W^{2,p}(X_0)\cap C(\bar X_0)$ implies that
\begin{equation}
  \label{eq:Gilbarg_Trudinger_C2_weak_maximum_principle}
  \inf_{X_0} u \geq \inf_{\partial X_0} u^-,
\end{equation}
where
\[
  u^- := \min\{u,0\}
  = \min\left\{1 - |\alpha|_E^2 + \frac{\delta}{\zeta r}, 0\right\} \quad\text{on } \partial X_0,
\]
By inequality \eqref{eq:Taubes_1996_SW_to_Gr_eq_2-2} from Proposition \ref{mainprop:Taubes_1996_SW_to_Gr_2-1}, we have
\[
  |\alpha|_E^2
  \leq
  1 + \frac{z}{r} \quad\text{on } X,
\]
and therefore $1 - |\alpha|_E^2 \geq -z/r$, so that
\[
  1 - |\alpha|_E^2 + \frac{\delta}{\zeta r}
  \geq 
  - \frac{z}{r} + \frac{\delta}{\zeta r}
  =
  \left(\frac{\delta}{\zeta} - z\right)\frac{1}{r} \quad\text{on } X.
\]
Hence, for $\delta$ in \eqref{eq:Bounds_r_zeta_regular} chosen large enough that $\delta \geq z\zeta$, 
we have $u \geq 0$ on $\partial X_0$ and thus $u^- = 0$ on $\partial X_0$, so we obtain from \eqref{eq:Gilbarg_Trudinger_C2_weak_maximum_principle} that
\[
  u \geq 0 \quad\text{on } X_0,
\]
and once again the conclusion \eqref{eq:Taubes_1996_SW_to_Gr_eq_2-11} holds. This completes our second proof of Theorem \ref{mainthm:Taubes_1996_SW_to_Gr_2-3}.
\end{proof}

% PF3-17-2025 Update section title/intro and chapter intro
We now complete the

\begin{proof}[Proof of Corollary \ref{maincor:Taubes_1996_SW_to_Gr_eq_2-12_and_2-13}]
Equation \eqref{eq:rho_(F_A^+)_0} for $\rho(F_A^+)_0$ yields the pointwise estimate,
\begin{multline*}
  |\rho(F_A^+)_0|_{\End(W^+\otimes E)}
  \leq \frac{1}{2}|\varphi\otimes\varphi^*)_0|_{\End(E)}
  + \frac{1}{2}|\star(\psi\otimes\psi^*)_0|_{\End(E)}
  \\
  + \frac{r}{4}|\wp(\psi)|_{\End(E)}
  + 2|(\psi\otimes\varphi^*)_0|_{\Lambda^{2,0}(X)\otimes\End(E)},
\end{multline*}
noting that $(\varphi\otimes\psi^*)_0^* = (\psi\otimes\varphi^*)_0$, where the ``$*$'' denotes the adjoint of an operator in $\Hom(E,\CC)$ or $\End(E)$. By applying the indicated pointwise equalities or inequalities below, we obtain
\begin{align*}
  |\varphi\otimes\varphi^*)_0|_{\End(E)}
  &= \frac{1}{\sqrt{2}}|\varphi|_E^2 \quad\text{(by \eqref{eq:Quartic_identity_varphi})},
  \\
  |\star(\psi\otimes\psi^*)_0|_{\End(E)}
  &= \frac{1}{\sqrt{2}}|\psi|_{\Lambda^{0,2}(E)}^2 \quad\text{(by \eqref{eq:Quartic_identity_psi})},
  \\
  |\wp(\psi)|_{\End(E)} &= 1 \quad\text{(by \eqref{eq:4-1})},
  \\
  |(\psi\otimes\varphi^*)_0|_{\Lambda^{2,0}(X)\otimes\End(E)}
  &\leq |\varphi|_E|\psi|_{\Lambda^{0,2}(E)} \quad\text{(by \eqref{eq:6} or \eqref{eq:4})}.
\end{align*}
By combining the preceding equalities and inequalities, we obtain
\begin{align*}
  |\rho(F_A^+)_0|_{\End(W^+\otimes E)}
  &\leq
    \frac{1}{2\sqrt{2}}|\varphi|_E^2 + \frac{1}{2\sqrt{2}}|\psi|_{\Lambda^{0,2}(E)}^2
  + \frac{r}{4} + 2|\varphi|_E|\psi|_{\Lambda^{0,2}(E)}
  \\
  &\leq
    \frac{(1+2\sqrt{2})}{2\sqrt{2}}\left(|\varphi|_E^2 + |\psi|_{\Lambda^{0,2}(E)}^2\right) + \frac{r}{4}.
\end{align*}
By hypothesis, $\varphi = \sqrt{r}\alpha$ and $\psi = \sqrt{r}\beta$, so the preceding inequality becomes
\[
  |\rho(F_A^+)_0|_{\End(W^+\otimes E)}
  \leq
  \frac{r(1+2\sqrt{2})}{2\sqrt{2}}\left(|\alpha|_E^2 + |\beta|_{\Lambda^{0,2}(E)}^2\right) + \frac{r}{4}.
\]
The pointwise estimate \eqref{eq:Taubes_1996_SW_to_Gr_eq_2-2} in Proposition \ref{mainprop:Taubes_1996_SW_to_Gr_2-1} yields
\[
  |\alpha|_E^2 + |\beta|_{\Lambda^{0,2}(E)}^2 \leq \frac{1}{3} + \frac{z}{r}.
\]
Substituting the preceding estimate into the previous bound for the pointwise norm of $\rho(F_A^+)_0$ gives
\[
  |\rho(F_A^+)_0|_{\End(W^+\otimes E)}
  \leq
  \frac{(1+2\sqrt{2})}{2\sqrt{2}}\left(\frac{1}{3} + \frac{z}{r}\right) + \frac{r}{4}.
\]
This concludes the proof of Corollary \ref{maincor:Taubes_1996_SW_to_Gr_eq_2-12_and_2-13}.
\end{proof}

\section[Differential inequalities for squared pointwise norms of sections of $E$ and $\Lambda^{0,2}(E)$]{Differential inequalities for  the squared pointwise norms of sections of $E$ and $\Lambda^{0,2}(E)$ with a regularized Taubes perturbation}
\label{sec:Differential_inequality_squared_pointwise_norm_anti-holomorphic_spinor_regular}
In this section and in Sections \ref{sec:Differential_inequality_Taubes_combination_squared_pointwise_norms_spinor_components_regular} and \ref{sec:Pointwise_estimate_squared_pointwise_norm_anti-holomorphic_spinor_regular}, we prove Theorem \ref{mainthm:Taubes_1996_SW_to_Gr_2-3_regular}, which gives a pointwise bound for $|\beta|_{\Lambda^{0,2}(E)}^2$when $(A,\varphi,\psi)$ with $(\varphi,\psi) = r^{1/2}(\alpha,\beta)$ is a solution to the system \eqref{eq:SO(3)_monopole_equations_almost_Hermitian_perturbed_intro_regular} of non-Abelian monopole equations with a regularized Taubes perturbation. We accomplish this by modifying the proof of Theorem \ref{mainthm:Taubes_1996_SW_to_Gr_2-3}. In this section, we prove Lemmas \ref{lem:Taubes_1996_SW_to_Gr_eq_2-3_inequality_regular} and \ref{lem:Taubes_1996_SW_to_Gr_eq_2-7_regular}, analogues for regularized Taubes perturbations of Lemmas \ref{lem:Taubes_1996_SW_to_Gr_eq_2-3_inequality} and \ref{lem:Taubes_1996_SW_to_Gr_eq_2-7} for singular Taubes perturbations. 

We begin by describing analogues for the regularization
$\wp_\gamma(\psi)$ in \eqref{eq:Definition_wp_intro_regular} of the pointwise estimates for $\wp(\psi)$ in \eqref{eq:Definition_wp_intro} provided by Claims \ref{claim:Pointwise_equality_and_inequalities_varphi} and \ref{claim:Pointwise_equality_and_inequalities_psi}. When we replace the singular perturbation $\wp(\psi)$ in \eqref{eq:Definition_wp_intro} for $\psi \in W^{2,p}(\Lambda^{0,2}(E))$ by its regularization $\wp_\gamma(\psi)$ in \eqref{eq:Definition_wp_intro_regular}, namely
\[
  \wp_\gamma(\psi) = 4\left(\gamma^2 + |\psi|_{\Lambda^{0,2}(E)}^2\right)^{-1}\star(\psi\otimes\psi^*)_0,
\]
the \emph{upper} bounds in \eqref{eq:9} and \eqref{eq:4-1} are replaced by
\begin{subequations}
  \label{eq:wp_intro_gamma_upper_bound_varphi_psi}
  \begin{align}
    \label{eq:wp_intro_gamma_upper_bound_varphi}
    \langle\wp_\gamma(\psi)\varphi,\varphi\rangle_E
    &\leq 4|\psi|_{\Lambda^{0,2}(E)}^{-2}\left|\langle\star(\psi\otimes\psi^*)_0\varphi,\varphi\rangle_E\right|
    \\
    \notag
    &= 2\left|\langle\wp(\psi)\varphi,\varphi\rangle_E\right| \leq 2|\varphi|_E^2 
      \quad\text{on } X,
    \\
    \label{eq:wp_intro_gamma_upper_bound_psi}
    \langle\wp_\gamma(\psi)\psi,\psi\rangle_{\Lambda^{0,2}(E)}
    &\leq 4|\psi|_{\Lambda^{0,2}(E)}^{-2}\langle\star(\psi\otimes\psi^*)_0\psi,\psi\rangle_{\Lambda^{0,2}(E)}
    \\
    \notag
    &= 2\langle\wp(\psi)\psi,\psi\rangle_{\Lambda^{0,2}(E)} = 2|\psi|_{\Lambda^{0,2}(E)}^2 
    \quad\text{on } X.
  \end{align}  
\end{subequations}
The \emph{lower} bound in \eqref{eq:9} is thus replaced by
\begin{align}
  \label{eq:wp_intro_gamma_lower_bound_varphi}
    \langle\wp_\gamma(\psi)\varphi,\varphi\rangle_E
    &\geq -4|\psi|_{\Lambda^{0,2}(E)}^{-2}\left|\langle\star(\psi\otimes\psi^*)_0\varphi,\varphi\rangle_E\right|
    \\
    \notag
    &= -2\left|\langle\wp(\psi)\varphi,\varphi\rangle_E\right| \geq -2|\varphi|_E^2 
      \quad\text{on } X.
\end{align}
Similarly, the equality \eqref{eq:4-1} can be replaced by the following \emph{crude lower bound},
\begin{align}
  \label{eq:wp_intro_gamma_lower_bound_psi_crude}
    \langle\wp_\gamma(\psi)\psi,\psi\rangle_{\Lambda^{0,2}(E)}
  &\geq -4|\psi|_{\Lambda^{0,2}(E)}^{-2}
    \langle\star(\psi\otimes\psi^*)_0\psi,\psi\rangle_{\Lambda^{0,2}(E)}
    \\
    \notag
    &= -2\langle\wp(\psi)\psi,\psi\rangle_{\Lambda^{0,2}(E)} = -2|\psi|_{\Lambda^{0,2}(E)}^2 
      \quad\text{on } X.
\end{align}
In order to obtain an improvement of the lower bound \eqref{eq:wp_intro_gamma_lower_bound_psi_crude}, we define the open subset,
\begin{equation}
  \label{eq:X_psi_regular}
  X_{\psi,\gamma} := \{x \in X: |\psi|_{\Lambda^{0,2}(E)} > \gamma\},
\end{equation}
so that
\[
  X \less X_{\psi,\gamma} = \{x \in X: |\psi|_{\Lambda^{0,2}(E)} \leq \gamma\}.
\]  
The lower bound \eqref{eq:wp_intro_gamma_lower_bound_varphi} for $\langle\wp_\gamma(\psi)\varphi,\varphi\rangle_E$ over $X$ can thus be slightly improved and the crude lower bound \eqref{eq:wp_intro_gamma_lower_bound_psi_crude} for $\langle\wp_\gamma(\psi)\psi,\psi\rangle_{\Lambda^{0,2}(E)}$ over $X$ can be significantly improved to give the following \emph{refined lower bounds} over the open subset $X_{\psi,\gamma} \subset X$:
\begin{subequations}
  \label{eq:wp_intro_gamma_lower_bounds}
  \begin{align}
    \label{eq:wp_intro_gamma_lower_bounds_varphi}
    \langle\wp_\gamma(\psi)\varphi,\varphi\rangle_E
    &\geq
      -2|\psi|_{\Lambda^{0,2}(E)}^{-2}\left|\langle\star(\psi\otimes\psi^*)_0\varphi,\varphi\rangle_E\right|
        = -\left|\langle\wp(\psi)\varphi,\varphi\rangle_E\right| \geq -|\varphi|_E^2
        \quad\text{on } X_{\psi,\gamma},
    \\
    \label{eq:wp_intro_gamma_lower_bounds_psi}
    \langle\wp_\gamma(\psi)\psi,\psi\rangle_{\Lambda^{0,2}(E)}
    &\geq
      2|\psi|_{\Lambda^{0,2}(E)}^{-2}\langle\star(\psi\otimes\psi^*)_0\psi,\psi\rangle_{\Lambda^{0,2}(E)}
    \\
    \notag
    &= \langle\wp(\psi)\psi,\psi\rangle_{\Lambda^{0,2}(E)} = |\psi|_{\Lambda^{0,2}(E)}^2
        \quad\text{on } X_{\psi,\gamma}.
  \end{align}  
\end{subequations}
For $r \in (0,\infty)$ and a suitable constant $\eps \in (0,1]$ to be fixed later, we choose
\begin{equation}
  \label{eq:delta_=_sqrt_r_epsilon}
  \gamma := r^{1/2}\eps.
\end{equation}
If $\psi = r^{1/2}\beta$ as in \eqref{eq:Taubes_1996_SW_to_Gr_eq_1-19} for $\beta \in W^{2,p}(\Lambda^{0,2}(E))$ then, after canceling the factor of $r$ in the numerator and denominator, the expression \eqref{eq:Definition_wp_intro_regular} for $\wp_\gamma(\psi)$ yields
\[
  \wp_{r^{1/2}\eps}(\psi) = 4\left(\eps^2 + |\beta|_{\Lambda^{0,2}(E)}^2\right)^{-1}\star(\beta\otimes\beta^*)_0
  =
  \wp_\eps(\beta).
\]
Note that
\[
  X_{\psi,r^{1/2}\eps}
  = \{x \in X: |\psi|_{\Lambda^{0,2}(E)} > r^{1/2}\eps\}
  = \{x \in X: |\beta|_{\Lambda^{0,2}(E)} > \eps\}
  = X_{\beta,\eps},
\]
and therefore
\[
  X\less X_{\beta,\eps} = \{x \in X: |\beta|_{\Lambda^{0,2}(E)} \leq \eps\}.
\]  
Hence, our lower bounds \eqref{eq:wp_intro_gamma_lower_bounds} for $\wp_\gamma(\psi)$ are equivalent to
\begin{subequations}
  \label{eq:wp_beta_intro_epsilon_lower_bounds}
  \begin{align}
    \label{eq:wp_beta_intro_epsilon_lower_bounds_alpha}
    \langle\wp_\eps(\beta)\alpha,\alpha\rangle_E
    &\geq
      -2|\beta|_{\Lambda^{0,2}(E)}^{-2}\left|\langle\star(\beta\otimes\beta^*)_0\alpha,\alpha\rangle_E\right|
        = -\left|\langle\wp(\beta)\alpha,\alpha\rangle_E\right| \geq -|\alpha|_E^2
        \quad\text{on } X_{\beta,\eps},
    \\
    \label{eq:wp_beta_intro_gamma_lower_bounds_beta}
    \langle\wp_\eps(\beta)\beta,\beta\rangle_{\Lambda^{0,2}(E)}
    &\geq
      2|\beta|_{\Lambda^{0,2}(E)}^{-2}\langle\star(\beta\otimes\beta^*)_0\beta,\beta\rangle_{\Lambda^{0,2}(E)}
    \\
    \notag
    &= \langle\wp(\beta)\beta,\beta\rangle_{\Lambda^{0,2}(E)} = |\beta|_{\Lambda^{0,2}(E)}^2
        \quad\text{on } X_{\beta,\eps}.
  \end{align}  
\end{subequations}
First, we have the following analogue of Lemma \ref{lem:Taubes_1996_SW_to_Gr_eq_2-3_inequality}.

\begin{lem}[Differential inequality for the squared pointwise norms of sections of $E$ with a regularized Taubes perturbation]
\label{lem:Taubes_1996_SW_to_Gr_eq_2-3_inequality_regular}
Continue the hypotheses of Theorem \ref{mainthm:Taubes_1996_SW_to_Gr_2-3_regular}. If $\alpha \in W^{2,p}(E)$ and $\beta \in  W^{2,p}(\Lambda^{0,2}(E))$ are defined by $(\varphi,\psi) = r^{1/2}(\alpha,\beta)$ as in \eqref{eq:Taubes_1996_SW_to_Gr_eq_1-19}, then 
\begin{multline}
  \label{eq:Taubes_1996_SW_to_Gr_eq_2-3_inequality_regular}
  \frac{1}{2}\Delta_g|\alpha|_E^2 + |\nabla_A\alpha|_E^2
  + \frac{r}{4}|\alpha|_E^4
  + \frac{r}{4}|\beta|_{\Lambda^{0,2}(E)}^2|\alpha|_E^2
  - \frac{r}{2}|\alpha|_E^2
  \\
  - \Real\langle i\Lambda_\omega F_{A_d}\alpha,\alpha\rangle_E
  - 2\Real\langle(\bar\mu\partial_A)^*\beta,\alpha\rangle_E \leq 0 \quad\text{on } X.
\end{multline}
\end{lem}

The inequality \eqref{eq:Taubes_1996_SW_to_Gr_eq_2-3_inequality_regular} is identical to the inequality \eqref{eq:Taubes_1996_SW_to_Gr_eq_2-3_simplified_before_rescaling} except that the term $-r|\alpha|_E^2/4$ has been replaced by $-r|\alpha|_E^2/2$,

\begin{proof}[Proof of Lemma \ref{lem:Taubes_1996_SW_to_Gr_eq_2-3_inequality_regular}]
We describe the changes to the proof of Lemma \ref{lem:Taubes_1996_SW_to_Gr_eq_2-3_inequality}. By again substituting the equality \eqref{eq:Quartic_identity_varphi} into the identity \eqref{eq:Taubes_1996_SW_to_Gr_eq_2-3_raw_before_rescaling} with $\wp(\psi)$ replaced by $\wp_\gamma(\psi)$, we see that
%PF10-10-2024 Recheck
\begin{multline}
  \label{eq:Taubes_1996_SW_to_Gr_eq_2-3_before_rescaling_regular}
  \frac{1}{2}\Delta_g|\varphi|_E^2 + |\nabla_A\varphi|_E^2
  + \frac{1}{4}|\varphi|_E^4
  - \frac{1}{2}\Real\langle \star(\psi\otimes\psi^*)_0\varphi,\varphi\rangle_E
  - \frac{r}{4}\langle \wp_\gamma(\psi)\varphi,\varphi\rangle_E
  \\
  - \Real\langle i\Lambda_\omega F_{A_d}\varphi,\varphi\rangle_E
  + \Real\langle(\psi\otimes\varphi^*)_0^*\psi,\varphi\rangle_E
  - 2\Real\langle(\bar\mu\partial_A)^*\psi,\varphi\rangle_E = 0 \quad\text{on } X.
\end{multline}
By substituting the inequalities \eqref{eq:5} (upper bound) and \eqref{eq:6} (lower bound) and upper bound \eqref{eq:wp_intro_gamma_upper_bound_varphi} into \eqref{eq:Taubes_1996_SW_to_Gr_eq_2-3_before_rescaling}, we obtain
\begin{multline*}
  \frac{1}{2}\Delta_g|\varphi|_E^2 + |\nabla_A\varphi|_E^2
  + \frac{1}{4}|\varphi|_E^4
  - \frac{1}{4}|\psi|_{\Lambda^{0,2}(E)}^2|\varphi|_E^2
  - \frac{r}{2}|\varphi|_E^2
  \\
  - \Real\langle i\Lambda_\omega F_{A_d}\varphi,\varphi\rangle_E
  + \frac{1}{2}|\psi|_{\Lambda^{0,2}(E)}^2|\varphi|_E^2
  - 2\Real\langle(\bar\mu\partial_A)^*\psi,\varphi\rangle_E \leq 0 \quad\text{on } X,
\end{multline*}
and thus
\begin{multline}
  \label{eq:Taubes_1996_SW_to_Gr_eq_2-3_simplified_before_rescaling_regular}
  \frac{1}{2}\Delta_g|\varphi|_E^2 + |\nabla_A\varphi|_E^2
  + \frac{1}{4}|\varphi|_E^4
  + \frac{1}{4}|\psi|_{\Lambda^{0,2}(E)}^2|\varphi|_E^2
  - \frac{r}{2}|\varphi|_E^2
  \\
  - \Real\langle i\Lambda_\omega F_{A_d}\varphi,\varphi\rangle_E
  - 2\Real\langle(\bar\mu\partial_A)^*\psi,\varphi\rangle_E \leq 0 \quad\text{on } X.
\end{multline}
The inequality \eqref{eq:Taubes_1996_SW_to_Gr_eq_2-3_simplified_before_rescaling_regular} is identical to \eqref{eq:Taubes_1996_SW_to_Gr_eq_2-3_simplified_before_rescaling} except that the term $-r|\varphi|_E^2/4$ has been replaced by $-r|\varphi|_E^2/2$. After making the substitution \eqref{eq:Taubes_1996_SW_to_Gr_eq_1-19} in \eqref{eq:Taubes_1996_SW_to_Gr_eq_2-3_simplified_before_rescaling}, namely
\[
  (\varphi,\psi) = r^{1/2}(\alpha,\beta) \in \Omega^0(E) \oplus \Omega^{0,2}(E),
\]
canceling factors of $r$ on both sides, we obtain \eqref{eq:Taubes_1996_SW_to_Gr_eq_2-3_inequality_regular}.
\end{proof}

Second, we have the following analogue of Lemma \ref{lem:Taubes_1996_SW_to_Gr_eq_2-7}.

\begin{lem}[Differential inequality for the squared pointwise norms of sections of $\Lambda^{0,2}(E)$ with a regularized Taubes perturbation]
\label{lem:Taubes_1996_SW_to_Gr_eq_2-7_regular}
Continue the hypotheses of Theorem \ref{mainthm:Taubes_1996_SW_to_Gr_2-3_regular}. If $\alpha \in W^{2,p}(E)$ and $\beta \in  W^{2,p}(\Lambda^{0,2}(E))$ are defined by $(\varphi,\psi) = r^{1/2}(\alpha,\beta)$ as in \eqref{eq:Taubes_1996_SW_to_Gr_eq_1-19}, then
\begin{multline}
  \label{eq:Taubes_1996_SW_to_Gr_eq_2-7_regular}
  \frac{1}{2}\Delta_g|\beta|_{\Lambda^{0,2}(E)}^2
  + \frac{r}{2}|\alpha|_E^2|\beta|_{\Lambda^{0,2}(E)}^2
  \\
  \leq
  - \frac{5r}{32}|\beta|_{\Lambda^{0,2}(E)}^2 - |\nabla_{A_2\otimes A}\beta|_{T^3(E)}^2
  + \frac{8}{r}\|N_J\|_{C^0(X)}^2|\nabla_A\alpha|_{\Lambda^1(E)}^2 \quad\text{on } X_{\beta,\eps}.
\end{multline}
\end{lem}

\begin{proof}
By applying the equality \eqref{eq:Quartic_identity_psi} and inequalities \eqref{eq:wp_intro_gamma_lower_bounds_psi} (lower bound) and \eqref{eq:3} (upper bound) and \eqref{eq:4} (lower bound) to the identity \eqref{eq:Taubes_1996_SW_to_Gr_eq_2-4_equality_before_rescaling} with $\wp(\psi)$ replaced by $\wp_\gamma(\psi)$, we obtain
\begin{multline*}
  \frac{1}{2}\Delta_g|\psi|_{\Lambda^{0,2}(E)}^2 + |\nabla_{A_2\otimes A}\psi|_{T^3(E)}^2
  \\
  + \Real\langle i\Lambda_\omega F_{A_2}\psi,\psi\rangle_{\Lambda^{0,2}(E)}
  + \frac{1}{2}\Real\langle i\Lambda_\omega F_{A_d}\psi,\psi\rangle_{\Lambda^{0,2}(E)}
  \\
  - \frac{1}{2}|\varphi|_E^2|\psi|_{\Lambda^{0,2}(E)}^2
  + \frac{1}{4}|\psi|_{\Lambda^{0,2}(E)}^4
  + \frac{r}{4}|\psi|_{\Lambda^{0,2}(E)}^2
  + |\varphi|_E^2|\psi|_{\Lambda^{0,2}(E)}^2
  \\
  - 2\Real\langle\bar\mu\partial_A\varphi,\psi\rangle_{\Lambda^{0,2}(E)} \leq 0 \quad\text{on } X_{\psi,\delta},
\end{multline*}
that is,
\begin{multline}
  \label{eq:Taubes_1996_SW_to_Gr_eq_2-4_before_rescaling_regular}
  \frac{1}{2}\Delta_g|\psi|_{\Lambda^{0,2}(E)}^2 + |\nabla_{A_2\otimes A}\psi|_{T^3(E)}^2
  \\
  + \Real\langle i\Lambda_\omega F_{A_2}\psi,\psi\rangle_{\Lambda^{0,2}(E)}
  + \frac{1}{2}\Real\langle i\Lambda_\omega F_{A_d}\psi,\psi\rangle_{\Lambda^{0,2}(E)}
  \\
  + \frac{1}{2}|\varphi|_E^2|\psi|_{\Lambda^{0,2}(E)}^2
  + \frac{1}{4}|\psi|_{\Lambda^{0,2}(E)}^4
  + \frac{r}{4}|\psi|_{\Lambda^{0,2}(E)}^2
  \\
  - 2\Real\langle\bar\mu\partial_A\varphi,\psi\rangle_{\Lambda^{0,2}(E)} \leq 0 \quad\text{on } X_{\psi,\delta}.
\end{multline}
The inequality \eqref{eq:Taubes_1996_SW_to_Gr_eq_2-4_before_rescaling_regular} is identical to the inequality \eqref{eq:Taubes_1996_SW_to_Gr_eq_2-4_before_rescaling} except that $X$ is replaced by $X_{\psi,\delta}$ and so the remainder of the proof of Lemma \ref{lem:Taubes_1996_SW_to_Gr_eq_2-7} applies with only that change to yield \eqref{eq:Taubes_1996_SW_to_Gr_eq_2-7_regular}, which is the same as the inequality \eqref{eq:Taubes_1996_SW_to_Gr_eq_2-7} except that $X$ is replaced by $X_{\beta,\eps}$.
\end{proof}  

\section[Differential inequality for an affine combination of squared pointwise norms]{Differential inequality for an affine combination of squared pointwise norms of sections of $E$ and $\Lambda^{0,2}(E)$ with a regularized Taubes perturbation}
\label{sec:Differential_inequality_Taubes_combination_squared_pointwise_norms_spinor_components_regular}
In this section, we prove Lemma \ref{lem:Taubes_1996_SW_to_Gr_2-2_regular}, an analogue for regularized Taubes perturbations of Lemma \ref{lem:Taubes_1996_SW_to_Gr_2-2} for singular Taubes perturbations.

We begin by describing the modifications to the argument in Section \ref{sec:Differential_inequality_Taubes_combination_squared_pointwise_norms_spinor_components}, noting that $X$ is replaced by $X_{\beta,\eps}$ and that the inequality \eqref{eq:Taubes_1996_SW_to_Gr_eq_2-3_inequality} is replaced by \eqref{eq:Taubes_1996_SW_to_Gr_eq_2-3_inequality_regular} (where the term $-r|\alpha|_E^2/4$ is replaced by $-r|\alpha|_E^2/2$). Inequality \eqref{eq:Taubes_1996_SW_to_Gr_eq_2-3_inequality_regular} and the equation \eqref{eq:Taubes_1996_SW_to_Gr_eq_2-8} for $w$, giving $|\alpha|_E^2 = 1 - w + \rho/r$ and
\[
  \frac{r}{2}|\alpha|_E^4
  =
  \frac{r}{2}|\alpha|_E^2\cdot |\alpha|_E^2
  =
  \frac{r}{2}|\alpha|_E^2 - \frac{r}{2}|\alpha|_E^2w + \frac{\rho}{2}|\alpha|_E^2 \quad\text{on } X,
\]  
imply the following inequality for $w$,
\begin{multline*}
  -\frac{1}{2}\Delta_gw + |\nabla_A\alpha|_E^2 + \frac{r}{2}|\alpha|_E^2 - \frac{r}{2}|\alpha|_E^2 w
  + \frac{\rho}{2}|\alpha|_E^2
  + \frac{r}{4}|\beta|_{\Lambda^{0,2}(E)}^2|\alpha|_E^2
  - \frac{r}{2}|\alpha|_E^2
  \\
  - \Real\langle i\Lambda_\omega F_{A_d}\alpha,\alpha\rangle_E
  - 2\Real\langle(\bar\mu\partial_A)^*\beta,\alpha\rangle_E \leq 0 \quad\text{on } X,
\end{multline*}
and thus an analogue of \eqref{eq:Taubes_1996_SW_to_Gr_eq_2-9_inequality},
\begin{multline}
  \label{eq:Taubes_1996_SW_to_Gr_eq_2-9_inequality_regular}
  \frac{1}{2}\Delta_gw + \frac{r}{2}|\alpha|_E^2w
  \geq
  |\nabla_A\alpha|_E^2 + \frac{r}{2}|\beta|_{\Lambda^{0,2}(E)}^2|\alpha|_E^2 + \frac{\rho}{2}|\alpha|_E^2
  \\
  - \Real\langle i\Lambda_\omega F_{A_d}\alpha,\alpha\rangle_E
  - 2\Real\langle(\bar\mu\partial_A)^*\beta,\alpha\rangle_E \quad\text{on } X.
\end{multline}
The only difference between the preceding differential inequality \eqref{eq:Taubes_1996_SW_to_Gr_eq_2-9_inequality_regular} and \eqref{eq:Taubes_1996_SW_to_Gr_eq_2-9_inequality} is that the factors $r/4$ and $\rho/4$ have been replaced by $r/2$ and $\rho/2$. The derivation of \eqref{eq:Taubes_1996_SW_to_Gr_eq_2-10_raw} from \eqref{eq:Taubes_1996_SW_to_Gr_eq_2-9_inequality} now leads to
\begin{multline}
  \label{eq:Taubes_1996_SW_to_Gr_eq_2-10_raw_regular}
  \frac{1}{2}\Delta_gw + \frac{r}{2}|\alpha|_E^2w
  \geq |\nabla_A\alpha|_E^2
  + \frac{r}{2}|\beta|_{\Lambda^{0,2}(E)}^2|\alpha|_E^2 + \frac{3\rho}{8}|\alpha|_E^2
  \\
  - \frac{8}{\rho}\left(\|N_J\|_{C^0(X)}^2|\nabla_{A_2\otimes A}\beta|_{T^3(E)}^2
    + \|\nabla N_J\|_{C^0(X)}^2|\beta|_{\Lambda^{0,2}(E)}^2\right) \quad\text{on } X,
\end{multline}
where the only difference between \eqref{eq:Taubes_1996_SW_to_Gr_eq_2-10_raw_regular} and \eqref{eq:Taubes_1996_SW_to_Gr_eq_2-10_raw} is that the factors $r/4$ and $\rho/8$ have been replaced by $r/2$ and $3\rho/8$. Substituting the expression \eqref{eq:Taubes_1996_SW_to_Gr_eq_2-10_z4} for $z_4$ into \eqref{eq:Taubes_1996_SW_to_Gr_eq_2-10_raw_regular} and using $3\rho|\alpha|_E^2/8 \geq \rho|\alpha|_E^2/8$ yields the following analogue of \eqref{eq:Taubes_1996_SW_to_Gr_eq_2-10},
\begin{multline}
  \label{eq:Taubes_1996_SW_to_Gr_eq_2-10_regular}
  \frac{1}{2}\Delta_gw + \frac{r}{2}|\alpha|_E^2w
  \geq |\nabla_A\alpha|_E^2
  + \frac{r}{2}|\beta|_{\Lambda^{0,2}(E)}^2|\alpha|_E^2 + \frac{\rho}{8}|\alpha|_E^2
  \\
  - \frac{z_4}{\rho}\left(|\nabla_{A_2\otimes A}\beta|_{T^3(E)}^2
    + |\beta|_{\Lambda^{0,2}(E)}^2\right) \quad\text{on } X.
\end{multline}
The only difference between \eqref{eq:Taubes_1996_SW_to_Gr_eq_2-10_regular} and \eqref{eq:Taubes_1996_SW_to_Gr_eq_2-10} is that the factor $r/4$ has been replaced by $r/2$. We next prove the following analogue of Lemma \ref{lem:Taubes_1996_SW_to_Gr_2-2}.

\begin{lem}[Differential inequality for an affine combination of squared pointwise norms of sections of $E$ and $\Lambda^{0,2}(E)$ with a regularized Taubes perturbation]
\label{lem:Taubes_1996_SW_to_Gr_2-2_regular}
Continue the definitions \eqref{eq:kappa1_2_3} of $\kappa_1$, $\kappa_2$, $\kappa_3$ with $z_3, z_4 \in (0,\infty)$ as in \eqref{eq:Taubes_1996_SW_to_Gr_eq_2-7_z3}, \eqref{eq:Taubes_1996_SW_to_Gr_eq_2-10_z4} and the hypotheses \eqref{eq:Bounds_r_zeta_regular} on $r,\delta,\zeta \in (0,\infty)$ in Lemma \ref{lem:Taubes_1996_SW_to_Gr_2-2} and continue the hypotheses of Theorem \ref{mainthm:Taubes_1996_SW_to_Gr_2-3_regular}. Then the function $u$ in \eqref{eq:Taubes_1996_SW_to_Gr_eq_2-5}, namely
\[
  u := 1 - |\alpha|_E^2 - \zeta|\beta|_{\Lambda^{0,2}(E)}^2 + \frac{\delta}{\zeta r} \in W^{2,p}(X),
\]
obeys the following differential inequality:
\begin{equation}
  \label{eq:Taubes_1996_SW_to_Gr_eq_2-6_regular}
  \frac{1}{2}\Delta_gu + \frac{r}{2}|\alpha|_E^2u
  \geq
  \left(1 - \zeta\frac{\kappa_1}{r}\right)|\nabla_A\alpha|_E^2
  + \frac{\zeta r}{8}|\beta|_{\Lambda^{0,2}(E)}^2
  + \frac{\zeta}{2}|\nabla_{A_2\otimes A}\beta|_{T^3(E)}^2 
  + \frac{\delta}{8\zeta}|\alpha|_E^2 \quad\text{on } X_{\beta,\eps}.
\end{equation}
\end{lem}

The inequality \eqref{eq:Taubes_1996_SW_to_Gr_eq_2-6_regular} is that the same as \eqref{eq:Taubes_1996_SW_to_Gr_eq_2-6} except that the factor $r|\alpha|_E^2/4$ is replaced by $r|\alpha|_E^2/2$.

\begin{proof}[Proof of Lemma \ref{lem:Taubes_1996_SW_to_Gr_2-2_regular}]
We describe the changes to the proof of Lemma \ref{lem:Taubes_1996_SW_to_Gr_2-2}. We see that the function $-\zeta|\beta|_{\Lambda^{0,2}(E)}^2$ obeys
\begin{align*}
  &\frac{1}{2}\Delta_g\left(-\zeta|\beta|_{\Lambda^{0,2}(E)}^2\right)
    + \frac{r}{2}|\alpha|_E^2\left(-\zeta|\beta|_{\Lambda^{0,2}(E)}^2\right)
  \\
  &\quad =
  -\zeta\left( \frac{1}{2}\Delta_g|\beta|_{\Lambda^{0,2}(E)}^2 + \frac{r}{2}|\alpha|_E^2|\beta|_{\Lambda^{0,2}(E)}^2 \right)
  \\
  &\quad \geq
  -\zeta\left(- \frac{5r}{32}|\beta|_{\Lambda^{0,2}(E)}^2 - |\nabla_{A_2\otimes A}\beta|_{T^3(E)}^2
    + \frac{z_3}{r}|\nabla_A\alpha|_{\Lambda^1(E)}^2\right) \quad\text{on } X_{\beta,\eps}
    \quad\text{(by \eqref{eq:Taubes_1996_SW_to_Gr_eq_2-7_regular})},
\end{align*}
so that
\begin{multline*}
  \frac{1}{2}\Delta_g\left(-\zeta|\beta|_{\Lambda^{0,2}(E)}^2\right)
  + \frac{r}{2}|\alpha|_E^2\left(-\zeta|\beta|_{\Lambda^{0,2}(E)}^2\right)
  \\
  \geq
  \frac{5\zeta r}{32}|\beta|_{\Lambda^{0,2}(E)}^2 + \zeta|\nabla_{A_2\otimes A}\beta|_{T^3(E)}^2
    - \frac{\zeta z_3}{r}|\nabla_A\alpha|_{\Lambda^1(E)}^2 \quad\text{on } X_{\beta,\eps}.
\end{multline*}
Combining the preceding inequality with \eqref{eq:Taubes_1996_SW_to_Gr_eq_2-10_regular} and substituting $\rho = \delta/\zeta$ from \eqref{eq:rho_equals_delta_over_zeta} gives 
\begin{align*}
  \frac{1}{2}\Delta_gu + \frac{r}{2}|\alpha|_E^2u
  &\geq
  |\nabla_A\alpha|_E^2
  + \frac{r}{2}|\beta|_{\Lambda^{0,2}(E)}^2|\alpha|_E^2
  - \frac{z_4}{\rho}\left(|\nabla_{A_2\otimes A}\beta|_{T^3(E)}^2
    + |\beta|_{\Lambda^{0,2}(E)}^2\right)
    + \frac{\delta}{8\zeta}|\alpha|_E^2
  \\
  &\quad + \frac{5\zeta r}{32}|\beta|_{\Lambda^{0,2}(E)}^2 + \zeta|\nabla_{A_2\otimes A}\beta|_{T^3(E)}^2
    - \frac{\zeta z_3}{r}|\nabla_A\alpha|_{\Lambda^1(E)}^2 \quad\text{on } X_{\beta,\eps},
\end{align*}
and hence, using the inequality $r|\beta|_{\Lambda^{0,2}(E)}^2|\alpha|_E^2/2 \geq 0$, 
\begin{multline}
  \label{eq:Taubes_1996_SW_to_Gr_eq_2-6_raw_regular}
  \frac{1}{2}\Delta_gu + \frac{r}{2}|\alpha|_E^2u
  \geq
    \left(1 - \frac{\zeta z_3}{r}\right)|\nabla_A\alpha|_{\Lambda^1(E)}^2
    + \left(\frac{5\zeta r}{32} - \frac{z_4}{\rho}\right)|\beta|_{\Lambda^{0,2}(E)}^2
    + \left(\zeta - \frac{z_4}{\rho}\right)|\nabla_{A_2\otimes A}\beta|_{T^3(E)}^2
    \\
    + \frac{\delta}{8\zeta}|\alpha|_E^2 \quad\text{on } X_{\beta,\eps}.
\end{multline}
The inequality \eqref{eq:Taubes_1996_SW_to_Gr_eq_2-6_raw_regular} is the same as \eqref{eq:Taubes_1996_SW_to_Gr_eq_2-6_raw} except that the factor $r|\alpha|_E^2/4$ is replaced by $r|\alpha|_E^2/2$. The remainder of the proof of Lemma \ref{lem:Taubes_1996_SW_to_Gr_2-2_regular} is identical to that of Lemma \ref{lem:Taubes_1996_SW_to_Gr_2-2}.
\end{proof}  

\section[Pointwise estimate for the squared pointwise norms of sections of $\Lambda^{0,2}(E)$]{Pointwise estimate for  the squared pointwise norms of sections of $\Lambda^{0,2}(E)$ with a regularized Taubes perturbation}
\label{sec:Pointwise_estimate_squared_pointwise_norm_anti-holomorphic_spinor_regular}
We can now complete the

\begin{proof}[Proof of Theorem \ref{mainthm:Taubes_1996_SW_to_Gr_2-3_regular}]
By hypothesis, we have $\alpha \in W^{2,p}(E)$ and $\beta \in W^{2,p}(\Lambda^{0,2}(E))$, so $u \in W^{2,p}(X)$ by definition \eqref{eq:Taubes_1996_SW_to_Gr_eq_2-5} when $p>2$. Also by hypothesis, we have $p\geq 4$, so $W^{2,p}X) \subset W^{2,d}(X)$ and $W^{2,p}X) \subset C^0(X)$ by the Sobolev Embedding Theorem in Adams and Fournier \cite[Theorem 4.12, p. 85]{AdamsFournier} with $d = \dim X = 4$. 
  
Because $|\beta|_{\Lambda^{0,2}(E)} = \eps$ on the boundary $\partial X_{\beta,\eps}$ of the open subset $X_{\beta,\eps}$, the definition \eqref{eq:Taubes_1996_SW_to_Gr_eq_2-5} of $u$ yields the equality
\begin{equation}
  \label{eq:u_on_boundary_X_beta_epsilon}
  u = 1 - |\alpha|_E^2 - \zeta\eps^2 + \frac{\delta}{\zeta r} \quad\text{on } \partial X_{\beta,\eps}.
\end{equation}
Noting the bounds \eqref{eq:Bounds_r_zeta_regular} on $r$, $\delta$, and $\zeta$ in Lemma \ref{lem:Taubes_1996_SW_to_Gr_2-2}, we again choose $\delta$ as in \eqref{eq:delta=r_kappa3_over_kappa1}, namely
\[
  \delta = \frac{r\kappa_3}{\kappa_1} \geq \kappa_3,
\]
but now take, for the constant $z \in [1,\infty)$ as in Proposition \ref{mainprop:Taubes_1996_SW_to_Gr_2-1_regular},
\begin{equation}
  \label{eq:zeta=r_min_1_over_z_and_1_over_2kappa1_kappa2}
  \zeta := \delta\min\left\{\frac{1}{z}, \frac{1}{2\kappa_2\kappa_3}\right\}.
\end{equation}
Therefore, the definition \eqref{eq:zeta=r_min_1_over_z_and_1_over_2kappa1_kappa2} and substitution of $\delta = r\kappa_3/\kappa_1$ from \eqref{eq:delta=r_kappa3_over_kappa1} yields the inequality
\[
  \zeta \leq \frac{\delta}{2\kappa_2\kappa_3} = \frac{r}{2\kappa_1\kappa_2},
\]
and so
\[
  \frac{\zeta\kappa_1}{r} \leq \frac{1}{2\kappa_2} \leq \frac{1}{2},
\]  
where the last inequality follows from the fact \eqref{eq:kappa2} that $\kappa_2 \geq 1$. We thus obtain
\[
  1 - \frac{\zeta\kappa_1}{r} \geq \frac{1}{2}.
\]
Consequently, the differential inequality \eqref{eq:Taubes_1996_SW_to_Gr_eq_2-6_regular} implies that
\[
  \frac{1}{2}\Delta_gu + \frac{r}{2}|\alpha|_E^2u \geq 0 \quad\text{on } X,
\]
and thus
\begin{equation}
  \label{eq:Lu_leq_0_X_beta_eps}
  Lu \leq 0 \quad\text{on } X_{\beta,\eps},
\end{equation}
for $u$ as in \eqref{eq:Taubes_1996_SW_to_Gr_eq_2-5} and where the linear second order elliptic operator $L$ is as in \eqref{eq:Elliptic_linear_scalar_second-order_operator} except that the term $-r|\alpha|_E^2/2$ is replaced by $-r|\alpha|_E^2$.

The weak maximum principle provided by Theorem \ref{thm:Gilbarg_Trudinger_3-1_and_9-1} implies that
\begin{equation}
  \label{eq:Gilbarg_Trudinger_C2_weak_maximum_principle_X_beta_eps}
  \inf_{X_{\beta,\eps}} u \geq \inf_{\partial X_{\beta,\eps}} u^-,
\end{equation}
where, substituting the definition of $u$ in \eqref{eq:Taubes_1996_SW_to_Gr_eq_2-5},
\[
  u^- := \min\{u,0\}
  = \min\left\{1 - |\alpha|_E^2 - \zeta\eps^2 + \frac{\delta}{\zeta r}, 0\right\}
  \quad\text{on } \partial X_{\beta,\eps}.
\]
(Note that the role of the open set $X_0 \subset X$ in \eqref{eq:X0} for our previous application of Theorem \ref{thm:Gilbarg_Trudinger_3-1_and_9-1} in our second proof of Theorem \ref{mainthm:Taubes_1996_SW_to_Gr_2-3} has been replaced here by that of the open subset $X_{\beta,\eps} \subset X$.) By the inequality \eqref{eq:Taubes_1996_SW_to_Gr_eq_2-2_regular} from Proposition \ref{mainprop:Taubes_1996_SW_to_Gr_2-1_regular}, we have
\begin{equation}
  \label{eq:eq:Taubes_1996_SW_to_Gr_eq_2-2_regular_simplified}
  |\alpha|_E^2
  \leq
  \frac{2}{3} + \frac{z}{r} \quad\text{on } X.
\end{equation}
Hence, we have the pointwise inequality,
\[
  \frac{2}{3} - |\alpha|_E^2 \geq -\frac{z}{r}
  \quad\text{on } \bar X_{\beta,\eps} = X_{\beta,\eps} \cup \partial X_{\beta,\eps},
\]
and this yields
\begin{align*}
  1 - |\alpha|_E^2 - \zeta\eps^2 + \frac{\delta}{\zeta r}
  &=
  \frac{2}{3} - |\alpha|_E^2 + \frac{1}{3} - \zeta\eps^2  + \frac{\delta}{\zeta r}
  \\
  &\geq 
  - \frac{z}{r} + \frac{\delta}{\zeta r} + \frac{1}{3} - \zeta\eps^2 \quad\text{on } \bar X_{\beta,\eps},
\end{align*}
that is,
\begin{equation}
  \label{eq:1_minus_alpha_squared_minus_zeta_eps_squared_minus_delta_over_zeta_r}
  1 - |\alpha|_E^2 - \zeta\eps^2 + \frac{\delta}{\zeta r}
  \geq
  \left(\frac{\delta}{\zeta} - z\right)\frac{1}{r} + \frac{1}{3} - \zeta\eps^2
  \quad\text{on } \bar X_{\beta,\eps}
\end{equation}
The choice \eqref{eq:zeta=r_min_1_over_z_and_1_over_2kappa1_kappa2} of $\zeta$ implies that $\zeta \leq \delta/z$ and so
\begin{equation}
  \label{eq:Delta_geq_z_zeta}
    \frac{\delta}{\zeta} \geq z.  
\end{equation}
We now constrain $\eps \in (0,1]$ so that
\begin{equation}
  \label{eq:Zeta_epsilon_squared_leq_one_third}
  \zeta\eps^2 \leq  \frac{1}{3}.
\end{equation}
Hence, the expression for $u$ in \eqref{eq:u_on_boundary_X_beta_epsilon} and the inequalities \eqref{eq:Delta_geq_z_zeta} and \eqref{eq:Zeta_epsilon_squared_leq_one_third} give
\[
  u = 1 - |\alpha|_E^2 - \zeta\eps^2 + \frac{\delta}{\zeta r}
  \geq 
  0 \quad\text{on } \partial X_{\beta,\eps}.
\]
Thus, $u^- = 0$ on $\partial X_{\beta,\eps}$. Therefore, we obtain from \eqref{eq:Gilbarg_Trudinger_C2_weak_maximum_principle_X_beta_eps} that
\[
  u \geq 0 \quad\text{on } X_{\beta,\eps}.
\]
Consequently, the expression for $u$ in \eqref{eq:u_on_boundary_X_beta_epsilon} yields
\[
  |\beta|_{\Lambda^{0,2}(E)}^2 \leq \frac{1}{\zeta}\left(1 - |\alpha|_E^2 + \frac{\delta}{\zeta r}\right)
   \quad\text{on } X_{\beta,\eps}.
\]
The definition \eqref{eq:zeta=r_min_1_over_z_and_1_over_2kappa1_kappa2} of $\zeta$ and substitution of $\delta = r\kappa_3/\kappa_1$ from \eqref{eq:delta=r_kappa3_over_kappa1} yields the equalities
\begin{subequations}
  \label{eq:Definition_zprime_and_zprimeprime_regular}
  \begin{align}
    \label{eq:Definition_zprime_regular}
    \frac{\zeta}{\delta}
    &= \min\left\{\frac{1}{z}, \frac{1}{2\kappa_2\kappa_3}\right\} =: \frac{1}{z'},
    \\
    \label{eq:Definition_zprimeprime_regular}
    \zeta
    &= r \frac{\kappa_3}{\kappa_1}\min\left\{\frac{1}{z}, \frac{1}{2\kappa_2\kappa_3}\right\} =: \frac{r}{z''},
  \end{align}
\end{subequations}
and definitions of the constants $z', z'' \in (0,\infty)$. Substituting these equalities into the preceding inequality for $|\beta|_{\Lambda^{0,2}(E)}^2$ gives
\begin{equation}
  \label{eq:Taubes_1996_SW_to_Gr_eq_2-11_regular_X_beta_eps}
  |\beta|_{\Lambda^{0,2}(E)}^2 \leq \frac{z''}{r}\left(1 - |\alpha|_E^2 + \frac{z'}{r}\right)
   \quad\text{on } X_{\beta,\eps}.
\end{equation}
On the other hand, we also have
\[
  |\beta|_{\Lambda^{0,2}(E)}^2 \leq \eps^2 \quad\text{on } X\less X_{\beta,\eps}.
\]
According to \eqref{eq:Zeta_epsilon_squared_leq_one_third}, we may choose $\eps^2 = 1/(3\zeta)$ and substituting the equality \eqref{eq:Definition_zprimeprime_regular} for $\zeta$ yields
\begin{equation}
  \label{eq:Definition_epsilon_is_zprimeprime_over_3r}
  \eps^2 = \frac{z''}{3r}.
\end{equation}
Therefore,
\begin{equation}
  \label{eq:Taubes_1996_SW_to_Gr_eq_2-11_regular_X_less_X_beta_eps}
  |\beta|_{\Lambda^{0,2}(E)}^2 \leq \frac{z''}{3r} \quad\text{on } X\less X_{\beta,\eps}.
\end{equation}
Observe that $1/z \geq 1/z'$ by \eqref{eq:Definition_zprime_regular} and so $z' \geq z$ and thus \eqref{eq:eq:Taubes_1996_SW_to_Gr_eq_2-2_regular_simplified} gives
\[
  |\alpha|_E^2
  \leq
  \frac{2}{3} + \frac{z'}{r} \quad\text{on } X.
\]
Consequently,
\[
  1 - |\alpha|_E^2 + \frac{z'}{r} \geq \frac{1}{3} \quad\text{on } X.
\]
Therefore, by combining the preceding inequality with \eqref{eq:Taubes_1996_SW_to_Gr_eq_2-11_regular_X_beta_eps} and \eqref{eq:Taubes_1996_SW_to_Gr_eq_2-11_regular_X_less_X_beta_eps}, we obtain
\[
  |\beta|_{\Lambda^{0,2}(E)}^2
  \leq
  \frac{z''}{r}\max\left\{1 - |\alpha|_E^2 + \frac{z'}{r}, \frac{1}{3}\right\}
  \leq
  \frac{z''}{r}\left(1 - |\alpha|_E^2 + \frac{z'}{r}\right) \quad\text{on } X.
\]  
This yields the conclusion \eqref{eq:Taubes_1996_SW_to_Gr_eq_2-11_regular}. The proof of the slightly improved pointwise inequality \eqref{eq:Taubes_1996_SW_to_Gr_eq_2-11_regular_degE_leq_0} follows by the same argument used to prove the improved inequality in Theorem \ref{mainthm:Taubes_1996_SW_to_Gr_2-3}. This completes the proof of Theorem \ref{mainthm:Taubes_1996_SW_to_Gr_2-3_regular}.
\end{proof}

% PF3-17-2025 Update section title/intro and chapter intro
We now complete the

\begin{proof}[Proof of Corollary \ref{maincor:Taubes_1996_SW_to_Gr_eq_2-12_and_2-13_regular}]
The argument is almost identical to the proof of Corollary \ref{maincor:Taubes_1996_SW_to_Gr_eq_2-12_and_2-13}. In place of the equality $|\wp|_{\End(E)} = 1$ by \eqref{eq:4-1}), we use the upper bound given by \eqref{eq:wp_intro_gamma_upper_bound_varphi_psi},
\[
    |\wp_\gamma(\psi)|_E \leq 2.
\]  
Proceeding otherwise exactly as before, we now obtain
\[
  |\rho(F_A^+)_0|_{\End(W^+\otimes E)}
  \leq
  \frac{r(1+2\sqrt{2})}{2\sqrt{2}}\left(|\alpha|_E^2 + |\beta|_{\Lambda^{0,2}(E)}^2\right) + \frac{r}{2}.
\]
The pointwise estimate \eqref{eq:Taubes_1996_SW_to_Gr_eq_2-2_regular} in Proposition \ref{mainprop:Taubes_1996_SW_to_Gr_2-1_regular} yields
\[
  |\alpha|_E^2 + |\beta|_{\Lambda^{0,2}(E)}^2 \leq \frac{2}{3} + \frac{z}{r},
\]
Substituting the preceding estimate into the previous bound for the pointwise norm of $\rho(F_A^+)_0$ gives
\[
  |\rho(F_A^+)_0|_{\End(W^+\otimes E)}
  \leq
  \frac{(1+2\sqrt{2})}{2\sqrt{2}}\left(\frac{2}{3} + \frac{z}{r}\right) + \frac{r}{2}.
\]
This concludes the proof of Corollary \ref{maincor:Taubes_1996_SW_to_Gr_eq_2-12_and_2-13_regular}.
\end{proof}

%PF3-13-2025 Standardize on "triples", not pairs. 
\section[Non-Abelian monopole equations for split pairs]{Non-Abelian monopole equations with Taubes perturbations for split triples}
\label{sec:Perturbation_non-Abelian_monopole_equations_split_pairs}
In this section, we describe the structure of the non-Abelian monopole equations, with singular and regularized Taubes perturbations, when $(A,\varphi,\psi)$ is split with respect to a decomposition $E = L_1 \oplus L_2$ as an orthogonal direct sum of Hermitian line bundles.

By analogy with Feehan and Leness \cite[Equations (2.55), p. 76, and (2.57), p. 77]{FL2a}, we call a pair $(A_1,\Phi_1)\in \sA(L_1,H_1)\times W^{1,p}(W^+\otimes L_1)$ a solution to the unperturbed \emph{Seiberg--Witten monopole equations} for the spin${}^c$ structure $\fs=\fs_0\otimes L_1$, where $\fs_0 = (\rho,W)$, if
\begin{subequations}
\label{eq:SeibergWitten}
\begin{align}
  \label{eq:SeibergWitten_curvature}
  \tr_{W^+}F_{A_W}^+ + 2F_{A_1}^+ - \rho^{-1}(\Phi_1\otimes\Phi_1^*)_{0} - F_{A_{\Lambda}}^+ &=0,
  \\
  \label{eq:SeibergWitten_Dirac}
  D_{A_1}\Phi_1 &=0,
\end{align}
\end{subequations}
where $\tr_{W^+}:\fu(W^+)\to i\ubarRR$ is defined by the trace on $2\times 2$ complex matrices, $(\Phi_1\otimes\Phi_1^*)_0$ is the component of the section $\Phi_1\otimes\Phi_1^*$ of $i\fu(W^+)$ contained in $i\su(W^+)$, and $D_{A_1}:\Omega^0(W^+\otimes L_1)\to \Omega^0(W^-\otimes L_1)$ is the Dirac operator defined by the \spinc connection $A_W$ on $W$ and unitary connection $A_1$ on $L_1$,
% TL9-18-2024: This is now multiply defined
%PF10-10-2024 Please fix
%\label{page:SpincDiracOperator}
and $A_\Lambda := A_{\det W^+}\otimes A_d$ is a unitary connection on the Hermitian line bundle $\Lambda := \det W^+\otimes \det E$, and $A_d$ is a unitary connection on the Hermitian line bundle $\det E$. The perturbation term, $\tr_{W^+}F_{A_W}^+ - F_{A_\Lambda}^+$ in \eqref{eq:SeibergWitten_curvature}, is chosen so that by Feehan and Leness \cite[Lemma 3.12, p. 95]{FL2a},
\begin{align*}
  (A,\Phi)
  &:=
  (A_1\oplus A_2, \Phi_1\oplus \Phi_2)
  \\
  &\,=
  \left(A_1 \oplus (A_d\otimes A_1^*), \Phi_1\oplus 0\right)
  \in \sA(E,H,A_d) \times \Omega^0(W^+\otimes E),
\end{align*}
is a (split) solution to the non-Abelian monopole equations \eqref{eq:SO(3)_monopole_equations} for the \spinu structure $\ft = (\rho,W\otimes E)$, where $E := L_1\oplus L_2$ is the complex rank two, Hermitian vector bundle defined by $L_2 := (\det E)\otimes L_1^*$. See Feehan and Leness 
%TL12-4-2025: Updated 
\cite[Section 6.6.1]{Feehan_Leness_introduction_virtual_morse_theory_so3_monopoles} for an explanation that updates the conventions and notation in \cite{FL2a} to those of \cite{Feehan_Leness_introduction_virtual_morse_theory_so3_monopoles}. Observe that
\[
  F_{A_\Lambda} = F_{A_{\det W^+}\otimes A_d} = F_{A_{\det W^+}} + F_{A_d}
  \quad\text{and}\quad
  \tr_{W^+}F_{A_W} = F_{A_{\det W^+}},  
\]
where the latter equality follows from Kobayashi \cite[Equation (1.5.19), p. 17]{Kobayashi_differential_geometry_complex_vector_bundles}. Hence, Equation \eqref{eq:SeibergWitten_curvature} simplifies to
\begin{equation}
  \label{eq:Seiberg-Witten_curvature_simplified}
  F_{A_1}^+ - \frac{1}{2}\rho^{-1}(\Phi_1\otimes\Phi_1^*)_{0} - \frac{1}{2}F_{A_d}^+ =0.
\end{equation}
When $(X,g,J,\omega)$ is almost Hermitian and $\fs_0 = \fs_\can = (\rho_\can,W_\can)$ is the canonical \spinc structure, then the unperturbed non-Abelian monopole equations \eqref{eq:SO(3)_monopole_equations} for a pair $(A,\Phi) = (A,\varphi,\psi)$ take the form \eqref{eq:SO(3)_monopole_equations_almost_Hermitian_intro}. When the perturbation $\wp(\psi)$ in \eqref{eq:Definition_wp_intro} is included, the corresponding system of equations becomes \eqref{eq:SO(3)_monopole_equations_almost_Hermitian_perturbed_intro}. The only difference between the two systems is that \eqref{eq:SO(3)_monopole_equations_(1,1)_curvature_intro} is replaced by \eqref{eq:SO(3)_monopole_equations_(1,1)_curvature_perturbed_intro}, namely
\[
  (\Lambda_\omega F_A)_0 = \frac{i}{2}(\varphi\otimes\varphi^*)_0 - \frac{i}{2}\star(\psi\otimes\psi^*)_0
  - \frac{ir}{4}\wp(\psi).
\]
Equation \eqref{eq:SO(3)_monopole_equations_(1,1)_curvature_perturbed_intro} is equivalent to equation \eqref{eq:SO(3)_monopole_equations_(1,1)_curvature_perturbed_omega_intro}, that is,
\[
  (F_A^\omega)_0 = \frac{i}{4}(\varphi\otimes\varphi^*)_0\omega - \frac{i}{4}\star(\psi\otimes\psi^*)_0\omega
    - \frac{ir}{8}\wp(\psi)\,\omega,
\]
where $F_A^\omega$ is the projection of $F_A$ onto its image in the factor $\Omega^0(\su(E))\cdot\omega$ of the orthogonal decomposition,
\[
  \Omega^+(\fsl(E)) = \Omega^{2,0}(\fsl(E))\oplus \Omega^0(X;\fsl(E))\cdot\omega \oplus \Omega^{0,2}(\fsl(E)).
\]
Suppose now that $(A,\Phi)$ is split with respect to a decomposition $E = L_1\oplus L_2$, so $\varphi \in \Omega^0(E)$ and $\psi_1 \in \Omega^{0,2}(E)$ are valued in $L_1$, that is, $\varphi = (\varphi_1,0)$ and $\psi = (\psi_1,0)$. Recall from the proof of Claim \ref{claim:Pointwise_equality_and_inequalities_varphi} that
\[
  (\varphi\otimes\varphi^*)_0
  =
  \varphi\otimes\varphi^* - \frac{1}{2}\tr_E(\varphi\otimes\varphi^*)\,\id_E
  =
  \varphi\otimes\varphi^* - \frac{1}{2}|\varphi|_E^2\,\id_E,
\]
and when $\varphi = \varphi_1 \in \Omega^0(L_1)$, then $(\varphi_1\otimes\varphi_1^*)s = \varphi_1\langle s,\varphi_1\rangle_{L_1} = s|\varphi|_{L_1}^2$ (for any $s \in \Omega^0(L_1)$) and so
\[
  (\varphi\otimes\varphi^*)_0 = \frac{1}{2}|\varphi_1|_{L_1}^2\left(\id_{L_1} \oplus -\id_{L_2}\right).
\]
Similarly,
%PF8-26-2024 Check!
\[
  \star(\psi\otimes\psi^*)_0 = \frac{1}{2}|\psi_1|_{\Lambda^{0,2}(L_1)}^2\left(\id_{L_1} \oplus -\id_{L_2}\right).
\]
According to the definition \eqref{eq:Definition_wp_intro} of $\wp(\psi)$, we have
\[
  \wp(\psi)
  =
  2|\psi|_{\Lambda^{0,2}(E)}^{-2}\star(\psi\otimes\psi^*)_0
  =
  \left(\id_{L_1} \oplus -\id_{L_2}\right)
  \quad\text{on } X_0,
\]
where $X_0 = \{x \in X: \psi(x) \neq 0\}$ as in \eqref{eq:X0}. Hence,
\[
  \frac{ir}{8}\wp(\psi)\,\omega
  =
  \frac{ir}{8}\left(\id_{L_1} \oplus -\id_{L_2}\right)\omega
  \quad\text{on } X_0.
\]  
Now $A = A_1 \oplus A_2 = A_1 \oplus (A_d\otimes A_1^*)$ and thus
\[
  F_A = F_{A_1} \oplus F_{A_d\otimes A_1^*}
  =  F_{A_1} \oplus \left(F_{A_d} + F_{A_1^*}\right)
  =  F_{A_1} \oplus \left(F_{A_d} - F_{A_1}\right).
\]
Therefore, $\tr_E F_A = F_{A_d}$, as expected, and
\[
  (F_A)_0 = F_A - \frac{1}{2}(\tr_E F_A)\,\id_E =  F_A - \frac{1}{2}F_{A_d}\,\id_E
  =
  \left(F_{A_1} - \frac{1}{2}F_{A_d}\right) \oplus \left(\frac{1}{2}F_{A_d} - F_{A_1}\right).
\]
Substituting the preceding observations into equation \eqref{eq:SO(3)_monopole_equations_(1,1)_curvature_perturbed_omega_intro} yields
\begin{align*}
  F_{A_1}^\omega - \frac{1}{2}F_{A_d}^\omega
  &= \frac{i}{8}|\varphi_1|_{L_1}^2\omega - \frac{i}{8}|\psi_1|_{\Lambda^{0,2}(L_1)}^2\omega
    - \frac{ir}{8}\omega \quad\text{on } X_0,
  \\
  \frac{1}{2}F_{A_d}^\omega - F_{A_1}^\omega
  &= -\frac{i}{8}|\varphi_1|_{L_1}^2\omega + \frac{i}{8}|\psi_1|_{\Lambda^{0,2}(L_1)}^2\omega
    + \frac{ir}{8}\omega \quad\text{on } X_0.
\end{align*}
The second equation in this system is clearly redundant and so we are left with
\[
  F_{A_1}^\omega
  =
  \frac{i}{8}|\varphi_1|_{L_1}^2\omega - \frac{i}{8}|\psi_1|_{\Lambda^{0,2}(L_1)}^2\omega - \frac{ir}{8}\omega
   + \frac{1}{2}F_{A_d}^\omega \quad\text{on } X_0.
\]
According to \eqref{eq:FAomega_and_LambdaFA}, we have
\[
  F_{A_d}^\omega = \frac{1}{2}(\Lambda_\omega F_{A_d})\,\omega,
\]
and so the $\omega$-component of the Seiberg--Witten curvature equation becomes
\begin{equation}
  \label{eq:Split_pair_omega-component_induced_SW_curvature_equation}
  F_{A_1}^\omega
  =
  \frac{i}{8}|\varphi_1|_{L_1}^2\omega - \frac{i}{8}|\psi_1|_{\Lambda^{0,2}(L_1)}^2\omega - \frac{ir}{8}\omega
   + \frac{1}{4}(\Lambda_\omega F_{A_d})\,\omega \quad\text{on } X_0.
\end{equation}
By writing $(\Lambda_\omega F_{A_d}) = -i(i\Lambda_\omega F_{A_d})$ and
\[
  - \frac{ir}{8}\omega  - \frac{i}{8}(2i\Lambda_\omega F_{A_d})\,\omega
  =
  - \frac{i}{8}\left(r + 2i\Lambda_\omega F_{A_d}\right)\omega
\]
and setting $r_0 := r + 2i\Lambda_\omega F_{A_d}$, we see that our singular perturbation of the $\omega$-component of the non-Abelian monopole equation on an almost Hermitian four-manifold reduces (on $X_0$) to that of Taubes,
\begin{equation}
  \label{eq:Taubes_omega-component_SW_curvature_equation}
  F_{A_1}^\omega
  =
  \frac{i}{8}|\varphi_1|_{L_1}^2\omega - \frac{i}{8}|\psi_1|_{\Lambda^{0,2}(L_1)}^2\omega - \frac{ir_0}{8}\omega
  \quad\text{on } X_0,
\end{equation}
for the Seiberg--Witten equations in \cite[Section 1(d), Equations (1.18), (1.19), and (1.20), p. 851]{TauSWGromov} when the solution $(A,\Phi)$ to the non-Abelian monopole equations is a split pair (see also Donaldson \cite[Section 4, third displayed system of equations, p. 60]{DonSW} and Kotschick \cite[Equations (15), (16), and (17)]{KotschickSW}).

Consequently, a solution $(A_1,\Phi_1)$ to the perturbed Seiberg--Witten monopole equations 
%PF8-27-2024 Refer to them
implied by a split solution $(A,\Phi)$ to the perturbed non-Abelian monopole equations \eqref{eq:SO(3)_monopole_equations_almost_Hermitian_perturbed_intro} should obey the seven pointwise estimates listed by Taubes in \cite[Section 1(f), Equations (1.24) and (1.26), p. 853]{TauSWGromov} over $X_0$. If $\{(A_1^n,\Phi_1^n)\}_{n\in\NN}$ is a sequence of solutions to the perturbed Seiberg--Witten monopole equations implied by a sequence of split solutions $\{(A,\Phi)\}_{n\in\NN}$ to the perturbed non-Abelian monopole equations with a sequence of parameters $\{r_n\}_{n\in\NN}$ with $\limsup_{n\in\NN} r_n = \infty$, then the sequence $\{(A_1^n,\Phi_1^n)\}_{n\in\NN}$ should obey convergence properties described by Taubes in \cite[Section 1(e), Theorem 1.3, p. 852]{TauSWGromov}; see also Taubes \cite[Section 5(a), pp. 881--882]{TauSWGromov}.

We now consider the case of \emph{regularized} perturbations. By the definition \eqref{eq:Definition_wp_intro_regular} of $\wp_\gamma(\psi)$, we have
\[
  \wp_\gamma(\psi)
  =
  \frac{4}{\gamma^2 + |\psi|_{\Lambda^{0,2}(E)}^2}\star(\psi\otimes\psi^*)_0
  =
  \frac{4|\psi_1|_{\Lambda^{0,2}(E)}^2}{\gamma^2 + |\psi_1|_{\Lambda^{0,2}(E)}^2}\left(\id_{L_1} \oplus -\id_{L_2}\right)
  \quad\text{on } X.
\]
Hence,
\[
  \frac{ir}{8}\wp_\gamma(\psi)\,\omega
  =
  \frac{ir|\psi_1|_{\Lambda^{0,2}(E)}^2}{2\left(\gamma^2 + |\psi_1|_{\Lambda^{0,2}(E)}^2\right)}\left(\id_{L_1} \oplus -\id_{L_2}\right)\omega
  \quad\text{on } X.
\]
Therefore, the $\omega$-component of the Seiberg--Witten curvature equation \eqref{eq:Split_pair_omega-component_induced_SW_curvature_equation} is replaced by
\begin{equation}
  \label{eq:Split_pair_omega-component_induced_SW_curvature_equation_regular}
  F_{A_1}^\omega
  =
  \frac{i}{8}|\varphi_1|_{L_1}^2\omega - \frac{i}{8}|\psi_1|_{\Lambda^{0,2}(L_1)}^2\omega
  - \frac{ir|\psi_1|_{\Lambda^{0,2}(E)}^2}{2\left(\gamma^2 + |\psi_1|_{\Lambda^{0,2}(E)}^2\right)}\omega
    + \frac{1}{4}(\Lambda_\omega F_{A_d})\,\omega \quad\text{on } X.
\end{equation}
Although similar to \eqref{eq:Split_pair_omega-component_induced_SW_curvature_equation}, the perturbation term is clearly distinct.

\chapter[Circle invariant two-form on the moduli space of non-Abelian monopoles]{Circle-invariant two-forms on the quotient space of unitary triples and moduli subspace of non-Abelian monopoles}
\label{chap:Construction_circle-invariant_non-degenerate_two-form_I}
In this chapter, we establish the framework needed to prove Theorem \ref{mainthm:AH_structure_bounded_evalue_spaces_non-Abelian_monopoles_symp_4-mflds}, Corollary \ref{maincor:Almost_Hermitian_structure_moduli_space_non-Abelian_monopoles_symplectic_4-manifolds}, and Theorem \ref{mainthm:IdentifyCriticalPoints} in Chapter \ref{chap:Construction_circle-invariant_non-degenerate_two-form_II}. To this end, we describe in Section \ref{sec:Riemannian_L2_metrics_quotient_space_spinu_pairs} our definition of the circle invariant weak Riemannian $L^2$ metric on the quotient space of spin${}^u$ pairs (and thus unitary triples). Section \ref{sec:ACStructure_on_Affine_Spinu_pairs} contains the definition of a circle invariant almost complex structure on the affine space of spin${}^u$ pairs and a skew-symmetric endomorphism on the quotient space of spin${}^u$ pairs. In Section \ref{sec:Construction_circle-invariant_two-form}, we construct a circle-invariant two-form on the quotient space of spin${}^u$ pairs. Section \ref{subsec:Characterization_critical_points_Hamiltonian_function_circle_action_on_non-abelian_monopoles} contains descriptions of the generators of the circle action on the affine space spin${}^u$ pairs (and thus unitary triples), the quotient space of non-zero section spin${}^u$ pairs (and thus unitary triples), and the moduli subspace of non-zero section non-Abelian monopoles with a regularized Taubes perturbation. For the affine space of spin${}^u$ pairs, we show that the circle action is Hamiltonian with respect to the circle invariant symplectic form defined by the weak Riemannian $L^2$ metric and the almost complex structure. For the quotient space of non-zero section spin${}^u$ pairs and moduli subspace of non-zero section non-Abelian monopoles with a regularized Taubes perturbation, we show --- under certain technical conditions that we introduce for the purpose of illustration only and which are not used to derive any other results --- that the circle action is Hamiltonian with respect to a circle invariant non-degenerate two form induced by the weak Riemannian $L^2$ metric and the almost complex structure on the affine space of spin${}^u$ pairs.

\section[Weak Riemannian $L^2$ metric on the quotient space of spin${}^u$ pairs]{Weak Riemannian $L^2$ metric on the quotient space of spin${}^u$ pairs}
\label{sec:Riemannian_L2_metrics_quotient_space_spinu_pairs}
For the construction and properties of weak Riemannian $L^2$ metrics on the quotient space of connections modulo gauge transformations and Riemannian $L^2$ metrics on the moduli space of anti-self-dual connections on a principal $G$-bundle over a closed, oriented, smooth Riemannian four-manifold $(X,g)$, we refer to Feehan \cite{FeehanGeometry, FeehanThesis}, Groisser and Parker \cite{GroisserParkerSphere, GroisserParkerGeometryDefinite}, and Itoh \cite{Itoh_1988}. For the construction and properties of weak Riemannian $L^2$ metrics on the quotient space of \spinc pairs modulo gauge transformations and Riemannian $L^2$ metrics on the moduli space of Seiberg--Witten monopoles over a closed, oriented, smooth Riemannian four-manifold $(X,g)$, we refer to Becker \cite{Becker_2008}.

\subsection[Construction of the weak Riemannian $L^2$ metric on the quotient space]{Construction of the weak Riemannian $L^2$ metric on the quotient space of spin${}^u$ pairs}
\label{subsec:Construction_Riemannian_L2_metrics_quotient_space}
Let $(X,g)$ be a closed, smooth Riemannian manifold, $(E,H)$ be a smooth, Hermitian vector bundle over $X$ and smooth, unitary connection $A_d$ on $\det E$, and $(\rho,W)$ be a \spinc structure over $(X,g)$, with $W=W^+\oplus W^-$ if $X$ has even dimension and $W=W^+$ if $X$ has odd dimension. For each pair $(A,\Phi) \in \sA(E,H,A_d)\times W^{1,p}(W^+\otimes E)$ with $p > (\dim X)/2$, the tangent space
\begin{multline}
\label{eq:TangentSpace_of_AffineSpace}
  T_{A,\Phi}\left(\sA(E,H,A_d)\times W^{1,p}(W^+\otimes E)\right)
  =
  T_A\sA(E,H,A_d) \oplus W^{1,p}(W^+\otimes E)
  \\
  =
  W^{1,p}\left((T^*X\otimes\su(E)\right) \oplus W^{1,p}(W^+\otimes E)
  =
 W^{1,p}\left(T^*X\otimes\su(E) \oplus W^+\otimes E\right)
\end{multline}
has a \emph{weak $L^2$ Riemannian metric} defined by the $L^2$ inner product\footnote{We use $\Real\,(\cdot,\cdot)$ to denote the real part of a Hermitian inner product $(\cdot,\cdot)$ on a complex vector space.},
\begin{multline}
  \label{eq:L2_metric_affine_space_spinu_pairs}
  \bg((a_1,\phi_1),(a_2,\phi_2))
  := \Real\,((a_1,\phi_1),(a_2,\phi_2))_{L^2(X)}
  = (a_1,a_2)_{L^2(X)} + \Real\,(\phi_1,\phi_2)_{L^2(X)},
  \\
  \text{for all } (a_1,\phi_1),(a_2,\phi_2) \in T_{A,\Phi}\left(\sA(E,H,A_d)\times W^{1,p}(W^+\otimes E)\right).
\end{multline}
The definition of the weak $L^2$ Riemannian metric on the quotient space of connections modulo gauge transformations,
\[
  \sB(E,H,A_d) := \sA(E,H,A_d)/\Aut(E,H,A_d)
\]
suggests the analogous definition of a weak $L^2$ Riemannian metric $\mathbf{\hat g}$ on the quotient space of \spinu pairs modulo gauge transformations,
\[
  \sC(E,H,W,A_d) := \left(\sA(E,H,A_d)\times W^{1,p}(W^+\otimes E)\right)/\Aut(E,H,A_d).
\]
If $(A,\Phi)$ represents a point $[A,\Phi] \in \sC(E,H,W,A_d)$, then the $\Stab(A,\Phi)$-invariant Zariski tangent space to $\sC(E,H,W,A_d)$ at $[A,\Phi]$ is represented by
\begin{multline}
  \label{eq:Representative_of_tangent_space_to_quotient_space_spinu_pairs}
  T_{A,\Phi}\sC(E,H,W,A_d)
  =
  \left. W^{1,p}\left(T^*X\otimes\su(E) \oplus W^+\otimes E\right)\right/ d_{A,\Phi}^0W^{1,p}(\su(E))
  \\
  \cong \Ker d_{A,\Phi}^{0,*} \cap W^{1,p}\left(T^*X\otimes\su(E) \oplus W^+\otimes E\right).
\end{multline}
Consider the $L^2$-orthogonal projection onto the Coulomb gauge slice:
\begin{equation}
  \label{eq:L2_orthogonal_projection_onto_slice_spin-u_pairs}
  \Pi_{A,\Phi}:W^{1,p}\left(T^*X\otimes\su(E) \oplus W^+\otimes E\right)
  \to
  \Ker d_{A,\Phi}^{0,*} \cap W^{1,p}\left(T^*X\otimes\su(E) \oplus W^+\otimes E\right).
\end{equation}
We restrict the $L^2$ metric $\bg$ in \eqref{eq:L2_metric_affine_space_spinu_pairs} to this Coulomb gauge slice representative for the tangent space to $\sC(E,H,W,A_d)$ at $[A,\Phi]$ to give
\begin{multline}
  \label{eq:L2_metric_quotient_space_spinu_pairs}
  \mathbf{\hat g}_{A,\Phi}\left((a_1,\phi_1),(a_2,\phi_2)\right) := (a_1,a_2)_{L^2(X)} + \Real\,(\phi_1,\phi_2)_{L^2(X)},
  \\
  \quad\text{for all } (a_1,\phi_1),(a_2,\phi_2)
  \in \Ker d_{A,\Phi}^{0,*} \cap W^{1,p}\left(T^*X\otimes\su(E) \oplus W^+\otimes E\right).
\end{multline}
We observe that the scalar values in 
%TL8-11-2025: Is this the correct reference?
\eqref{eq:Almost_complex_structure_quotient_tangent_space_spinu_pairs_skew_symmetric} are independent of the choice of representative $(A,\Phi)$ of the point $[A,\Phi] \in \sC(E,H,W,A_d)$ by replacing $(A,\Phi)$ by $u(A,\Phi)$ and noting that the values are independent of $u \in W^{2,p}(\SU(E))$ because $u\in W^{2,p}(\SU(E))$ acts by pointwise isometries on elements of $W^{1,p}(T^*X\otimes\su(E)\oplus W^+\otimes E)$. Hence, the expression for $\mathbf{\hat g}$ in \eqref{eq:L2_metric_quotient_space_spinu_pairs} defines a \emph{weak $L^2$ Riemannian metric} on the quotient space $\sC(E,H,W,A_d)$. (The fact that the weak $L^2$ Riemannian metric is well-defined on the quotient space $\sB(E,H,A_d)$ is proved by the first author as \cite[Lemma 1.3.1]{FeehanThesis}.)

\begin{rmk}[Representations of tangent spaces to quotient spaces]
\label{rmk:Representation_tangent_spaces_quotient_spaces}
Following the description in Feehan and Leness 
%TL12-4-2025: Updated 
\cite[Equations (6.3.11) and (6.3.12)]{Feehan_Leness_introduction_virtual_morse_theory_so3_monopoles},
the tangent bundle of $\sC(E,H,W,A_d)$  is given by
\begin{multline*}
T\sC(E,H,W,A_d)
\\
=
\left.\left\{(A,\Phi,a,\phi)\in  \sA(E,H,A_d)\times W^{1,p}\left(T^*X\otimes\su(E)\oplus W^+\otimes E\right): d_{A,\Phi}^{0,*}(a,\phi)=0\right\}\right/\Aut(E,H,A_d),
\end{multline*}
where $\Aut(E,H,A_d)$ acts by
\begin{equation}
\label{eq:TangentBundleAction}
\left( u,(A,\Phi,a,\phi)\right)\mapsto
\left( u^*(A,\Phi),u^*(a,\phi)\right):=
\left(u^*(A,\Phi),u^{-1}au,u^{-1}\phi\right).
\end{equation}
Note that the preceding action on the space $W^{1,p}(T^*X\otimes\su(E)\oplus W^+\otimes E)$ is the same as that appearing in \eqref{eq:GaugeGroupActionOn_E1} and that this action maps $\Ker d_{A,\Phi}^{0,*}$ to $\Ker d_{u^*(A,\Phi)}^{0,*}$ in
the sense that
\begin{equation}
\label{eq:GaugeEquivarianceOfSlice}
u^*\left(\Ker d_{A,\Phi}^{*,0}\right)
=
\Ker d_{u^*(A,\Phi)}^{*,0}
\end{equation}
Further, we note that the projection operator $\Pi_{A,\Phi}$ in \eqref{eq:L2_orthogonal_projection_onto_slice_spin-u_pairs} satisfies the gauge equivariance condition
\begin{equation}
\label{eq:SliceProjectionGaugeEquivariance}
\Pi_{u^*(A,\Phi)}= u^* \Pi_{A,\Phi},
\end{equation}
where the $u^*$ on the right-hand-side of the preceding is defined in \eqref{eq:TangentBundleAction}.

Hence, an element of the tangent space of $\sC(E,H,A_d)$ at $[A,\Phi]$ is represented by
\[
(A,\Phi,a,\phi)\in  \sA(E,H,A_d)\times W^{1,p}(T^*X\otimes\su(E)\oplus W^+\otimes E),
\]
where $(A,\Phi)$ is in the gauge-equivalence class $[A,\Phi]$, and $(a,\phi)\in\Ker d_{A,\Phi}^{0,*}$, and $(A,\Phi,a,\phi)$ and $(A',\Phi',a',\phi')$ represent the same tangent vector if $(A',\Phi')=u^*(A,\Phi)$ and $(a',\phi')=(u^{-1}au,u^{-1}\phi)$, for some $u\in \Aut(E,H,A_d)$. By choosing a representative $(A,\Phi)$ of the gauge-equivalence class $[A,\Phi]\in \sC(E,H,A_d)$, we can represent the tangent space of $\sC(E,H,A_d)$ at $[A,\Phi]$ by
\begin{multline}
\label{eq:Tangent_space_to_quotient_space_spinu_pairs}
T_{A,\Phi}\sC(E,H,W,A_d)
\\
:=
\left\{(A,\Phi,a,\phi)\in  \{(A,\Phi)\}\times W^{1,p}\left(T^*X\otimes\su(E)\oplus W^+\otimes E\right): d_{A,\Phi}^{0,*}(a,\phi)=0\right\},
\end{multline}
leading to the equivalent formulation \eqref{eq:Representative_of_tangent_space_to_quotient_space_spinu_pairs}.
\qed\end{rmk}

\subsection[Circle invariance of the weak Riemannian $L^2$ metric on the affine space]{Circle invariance of the weak Riemannian $L^2$ metric on the affine space of spin${}^u$ pairs}
\label{subsec:Circle_invariance_Riemannian_L2_metrics_affine_space}
We next claim that $\bg$ is a \emph{circle invariant weak $L^2$ Riemannian metric} on $\sA(E,H,A_d)\times W^{1,p}(W^+\otimes E)$ with respect to the circle action on the affine space given \eqref{eq:S1ZAction},
\[
  S^1 \times \sA(E,H,A_d)\times W^{1,p}(W^+\otimes E) \ni (e^{i\theta},(A,\Phi))
  \mapsto (A,e^{i\theta}\Phi) \in \sA(E,H,A_d)\times W^{1,p}(W^+\otimes E),
\]
since the action of $S^1$ by scalar multiplication on sections of $W^+\otimes E$ is isometric with respect to the Hermitian metric on $W^+$ induced by the Riemannian metric on $X$ and the given Hermitian metric on $E$. To see this explicitly, we observe that
\begin{align*}
  \bg_{A,e^{i\theta}\Phi}\left((a_1,e^{i\theta}\phi_1),(a_2,e^{i\theta}\phi_2)\right)
  &= \bg\left((a_1,e^{i\theta}\phi_1),(a_2,e^{i\theta}\phi_2)\right)
    \quad\text{(by \eqref{eq:L2_metric_affine_space_spinu_pairs})}
  \\
  &= (a_1,a_2)_{L^2(X)} + \Real\,(e^{i\theta}\phi_1,e^{i\theta}\phi_2)_{L^2(X)}
    \quad\text{(by \eqref{eq:L2_metric_affine_space_spinu_pairs})}
  \\
  &= (a_1,a_2)_{L^2(X)} + \Real\,(\phi_1,\phi_2)_{L^2(X)}
  \\
  &= \bg\left((a_1,\phi_1),(a_2,\phi_2)\right)
    \quad\text{(by \eqref{eq:L2_metric_affine_space_spinu_pairs})}
  \\
  &= \bg_{A,\Phi}\left((a_1,\phi_1),(a_2,\phi_2)\right)
    \quad\text{(by \eqref{eq:L2_metric_affine_space_spinu_pairs})},  
\end{align*}
and so, for all $(A,\Phi) \in \sA(E,H,A_d)\times W^{1,p}(W^+\otimes E)$,
\begin{multline}
  \label{eq:L2_metric_affine_space_spinu_pairs_circle_invariance}
  \bg_{A,e^{i\theta}\Phi}\left((a_1,e^{i\theta}\phi_1),(a_2,e^{i\theta}\phi_2)\right)
  =
  \bg_{A,\Phi}\left((a_1,\phi_1),(a_2,\phi_2)\right),
  \\
  \text{for all } e^{i\theta} \in S^1 \text{ and } (a_1,\phi_1),(a_2,\phi_2)
  \in T_{A,\Phi}\left(\sA(E,H,A_d)\times W^{1,p}(W^+\otimes E)\right).
\end{multline}
Thus, $\bg$ is a circle invariant weak $L^2$ Riemannian metric on $\sA(E,H,A_d)\times W^{1,p}(W^+\otimes E)$, as claimed.

\subsection[Circle invariance of the weak Riemannian $L^2$ metric on the quotient space]{Circle invariance of the weak Riemannian $L^2$ metric on the quotient space of spin${}^u$ pairs}
\label{subsec:Circle_invariance_Riemannian_L2_metrics_quotient_space}
The $S^1$ action on the affine space commutes with the action of $\Aut(E,H,A_d)$ and thus descends to the action \eqref{eq:S1ZActionOnQuotientSpace} of $S^1$ on the quotient space $\sC(E,H,W,A_d)$. In particular, $S^1$ acts by isometries on the tangent bundle $T\sC^0(E,H,W,A_d)$ and so we claim that $\mathbf{\hat g}$ in \eqref{eq:L2_metric_quotient_space_spinu_pairs} is a \emph{circle invariant weak $L^2$ Riemannian metric} on $\sC(E,H,W,A_d)$. Explicitly, to see this, we observe that if $(A,\Phi)$ represents a point $[A,\Phi] \in \sC(E,H,W,A_d)$ and $(a,\phi) \in \Ker d_{A,\Phi}^{0,*}$ and $e^{i\theta} \in S^1$, then
\begin{multline*}
  \left(d_{A,e^{i\theta}\Phi}^{0,*}(a,e^{i\theta}\phi),\xi\right)_{L^2(X)}
  = \left((a,e^{i\theta}\phi),d_{A,e^{i\theta}\Phi}^0\xi\right)_{L^2(X)}
  = \left((a,e^{i\theta}\phi),d_A\xi - \xi e^{i\theta}\Phi\right)
  \\
  = (a,d_A\xi) + \left(e^{i\theta}\phi, - e^{i\theta}\xi\Phi\right)_{L^2(X)}
  = (a,d_A\xi) + \left(\phi, -\xi\Phi\right)_{L^2(X)}
  = ((a,\phi),d_{A,\Phi}\xi)_{L^2(X)}
  \\
  = \left((d_{A,\Phi}^{0,*}(a,\phi),\xi\right)_{L^2(X)}
  = 0, \quad\text{for all } \xi \in W^{2,p}(\su(E)),
\end{multline*}
and so
\begin{equation}
  \label{eq:S1_action_Coulomb_gauge_spinu_slice}
  (a,e^{i\theta}\phi) \in \Ker d_{A,e^{i\theta}\Phi}^{0,*},
  \quad\text{for all } e^{i\theta} \in S^1 \text{ and } (a,\phi) \in \Ker d_{A,\Phi}^{0,*}.
\end{equation}
If $(A,\Phi)$ represents $[A,\Phi] \in \sC(E,H,W,A_d)$ and $e^{i\theta} \in S^1$ and $(a_1,\phi_1),(a_2,\phi_2) \in \Ker d_{A,\Phi}^{0,*}$, then $(a_1,e^{i\theta}\phi_1), (a_2,e^{i\theta}\phi_2) \in \Ker d_{A,e^{i\theta}\Phi}^{0,*}$ by \eqref{eq:S1_action_Coulomb_gauge_spinu_slice} and we see that
\begin{align*}
  &\mathbf{\hat g}_{A,e^{i\theta}\Phi}\left((a_1,e^{i\theta}\phi_1),(a_2,e^{i\theta}\phi_2)\right)
  \\
  &= \Real\,((a_1,e^{i\theta}\phi_1),(a_2,e^{i\theta}\phi_2))_{L^2(X)} \quad\text{(by \eqref{eq:L2_metric_quotient_space_spinu_pairs})}
  \\
  &= \Real\,((a_1,\phi_1),(a_2,\phi_2))_{L^2(X)}
  \\
  &= \mathbf{\hat g}_{A,\Phi}\left((a_1,\phi_1),(a_2,\phi_2)\right) \quad\text{(by \eqref{eq:L2_metric_quotient_space_spinu_pairs})}
\end{align*}
and thus we conclude that, for a representative $(A,\Phi)$ a point $[A,\Phi] \in \sC(E,H,W,A_d)$,
\begin{multline}
  \label{eq:L2_metric_quotient_space_spinu_pairs_circle_invariance}
  \mathbf{\hat g}_{A,e^{i\theta}\Phi}\left(\left(a_1,e^{i\theta}\phi_1\right),\left(a_2,e^{i\theta}\phi_2\right)\right)
  =
  \mathbf{\hat g}_{A,\Phi}\left((a_1,\phi_1),(a_2,\phi_2)\right),
  \\
  \text{for all } e^{i\theta} \in S^1 \text{ and } (a_1,\phi_1),(a_2,\phi_2)
  \in \Ker d_{A,\Phi}^{0,*}\cap W^{1,p}\left(T^*X\otimes\su(E) \oplus W^+\otimes E\right).
\end{multline}
As earlier, we can see that the scalar values in \eqref{eq:L2_metric_quotient_space_spinu_pairs_circle_invariance} are independent of the choice of representative $(A,\Phi)$ of $[A,\Phi] \in \sC(E,H,W,A_d)$ by replacing $(A,\Phi)$ by $u(A,\Phi)$ and noting that the values are independent of $u \in W^{2,p}(\SU(E))$. Therefore, $\mathbf{\hat g}$ is a circle invariant weak $L^2$ Riemannian metric on $\sC(E,H,W,A_d)$, as claimed.

\section[Skew-symmetric endomorphism on the quotient space]{Almost complex structure on the affine space of spin${}^u$ pairs and skew-symmetric endomorphism on the quotient space of spin${}^u$ pairs}
\label{sec:ACStructure_on_Affine_Spinu_pairs}

\subsection[Construction of orthogonal almost complex structure on the affine space]{Construction of circle-invariant orthogonal almost complex structure on the affine space of spin${}^u$ pairs}
\label{sec:ACStructure_on_Affine_Spinu_pairs_ac_structure_only}
We begin by recalling the

\begin{rmk}[Reality condition for sections of the complex endomorphism bundle]
\label{rmk:Reality_condition_sections_complex_endomorphism_bundle_and_20_curvature_equation}
(See Feehan and Leness 
%TL5-4-2025: Updated
\cite[Remark 8.3.12]{Feehan_Leness_introduction_virtual_morse_theory_so3_monopoles}.)
The condition that the following element is skew-Hermitian,
\begin{equation}
  \label{eq:Decompose_a_in_Omega1suE_into_10_and_01_components}
  a = \frac{1}{2}(a'+a'') \in \Omega^1(\su(E))\otimes_\RR\CC = \Omega^{1,0}(\fsl(E))\oplus \Omega^{0,1}(\fsl(E)),
\end{equation}
is equivalent to the condition that (see Itoh \cite[Section 2, p. 850]{Itoh_1988} or Kobayashi \cite[Equation (7.6.11), p. 251]{Kobayashi_differential_geometry_complex_vector_bundles})
\begin{equation}
\label{eq:Kobayashi_7-6-11}
  a' = -(a'')^\dagger \in \Omega^{1,0}(\fsl(E)),
\end{equation}
where $a^\dagger = \bar{a}^\intercal$ is defined by taking the complex, conjugate transpose of complex, matrix-valued representatives with respect to any local frame for $E$ that is orthonormal with respect to the Hermitian metric on $E$.
\qed
\end{rmk}

For each pair $(A,\Phi) \in (\sA(E,H,A_d)\times W^{1,p}(W^+\otimes E))$, the tangent space \eqref{eq:TangentSpace_of_AffineSpace}, namely
\begin{multline*}
  T_{A,\Phi}\left(\sA(E,H,A_d)\times W^{1,p}(W^+\otimes E)\right)
  =
  T_A\sA(E,H,A_d) \oplus W^{1,p}(W^+\otimes E)
  \\
  =
  W^{1,p}(T^*X\otimes \su(E)) \oplus W^{1,p}(W^+\otimes E)
  =
  W^{1,p}\left(T^*X\otimes \su(E) \oplus W^+\otimes E\right)
\end{multline*}
has an almost complex structure,
%PF8-26-2024 Expression below only makes sense after complexification; should just say J extends from \Omega^1(X;\RR) to \Omega^1(\su(E))
\begin{multline}
  \label{eq:Almost_complex_structure_affine_space_spinu_pairs}
  \bJ(a,\phi) := (i\cdot a, i\cdot\phi) = \left(-\frac{i}{2}a'+\frac{i}{2}a'', i\phi\right),
  \\
  \text{for all } (a,\phi) \in T_{A,\Phi}\left(\sA(E,H,A_d)\times W^{1,p}(W^+\otimes E)\right).
\end{multline}
defined by the splitting \eqref{eq:Decompose_a_in_Omega1suE_into_10_and_01_components}. The factors of $1/2$ in \eqref{eq:Almost_complex_structure_affine_space_spinu_pairs} are used in order to simplify our comparison of the elliptic deformation complexes for the non-Abelian monopole equations over complex K\"ahler surfaces in Feehan and Leness 
%TL12-4-2025: Updated
\cite[Sections 10.1, 10.2, and 10.3]{Feehan_Leness_introduction_virtual_morse_theory_so3_monopoles}. We follow the sign convention of Kobayashi \cite[Equation (7.6.12), p. 252]{Kobayashi_differential_geometry_complex_vector_bundles} for the action of $\bJ$ on $a'$ and $a''$, although it is opposite to that of Itoh \cite[Equation (4.1), p. 19]{Itoh_1988} and that implied by the following

\begin{rmk}[Conventions for almost complex vector spaces and their complex dual spaces]
\label{rmk:Conventions_almost_complex_vector_spaces_and_complex_dual_spaces}  
We summarize the usual conventions for almost complex vector spaces and their complex dual spaces as described, for example, by Huybrechts \cite[Definition 1.2.4, p. 25 and Lemma 1.2.6, p. 26]{Huybrechts_2005}. Thus, $T_xX^{1,0}$ and $T_xX^{0,1}$ are the $+i$ and $-i$ eigenspaces of $J \in \End_\CC((T_xX)_\CC)$ with $(T_xX)_\CC := T_xX\otimes_\RR\CC = T_xX^{1,0} \oplus T_xX^{0,1}$. The real dual space $T_x^*X = \Hom_\RR(T_xX,\RR)$ has an induced almost complex structure given by $J(\vartheta)(Y) := \vartheta(Jv)$, for all $Y \in T_xX$ and $\vartheta \in T_x^*X$. The complex dual space $(T_x^*X)_\CC := \Hom_\CC((T_xX)_\CC,\CC) = T_x^*X\otimes\CC$ has a splitting $(T_x^*X)_\CC = (T_x^*X)^{1,0} \oplus (T_x^*X)^{0,1}$ such that $(T_x^*X)^{1,0} = (T_xX^{1,0})^*$ and is the $+i$-eigenspace of $J \in \End_\CC((T_x^*X)_\CC)$, while $(T_x^*X)^{0,1} = (T_xX^{0,1})^*$ and is the $-i$-eigenspace of $J \in \End_\CC((T_x^*X)_\CC)$.
\qed
\end{rmk}  

We claim that $\bJ$ is a $\bg$-\emph{orthogonal almost complex structure} on the affine space $\sA(E,H,A_d)\times W^{1,p}(W^+\otimes E)$ of \spinu pairs. To prove this, we note that the
% PF8-22-2024 Reference it
%PF8-26-2024 Proof only makes sense after complexifying; better to write Ja, etc.
preceding splittings are $L^2$-orthogonal and so
\begin{multline*}
  \|i\cdot a\|_{L^2(X)}^2
  = \frac{1}{4}\|ia'-ia''\|_{L^2(X)}^2
  = \frac{1}{4}\|ia'\|_{L^2(X)}^2 + \frac{1}{4}\|-ia''\|_{L^2(X)}^2
  \\
  = \frac{1}{4}\|a'\|_{L^2(X)}^2 + \frac{1}{4}\|a''\|_{L^2(X)}^2
  = \frac{1}{4}\|a'+a''\|_{L^2(X)}^2
  = \|a\|_{L^2(X)}^2, \quad\text{for all } a \in W^{1,p}(T^*X\otimes \su(E)),
\end{multline*}
while
\[
  \|i\cdot \phi\|_{L^2(X)}^2 = \|i\phi\|_{L^2(X)}^2 = \|\phi\|_{L^2(X)}^2,
  \quad\text{for all } \phi \in W^{1,p}(W^+\otimes E).
\]
Therefore, by definition \eqref{eq:L2_metric_affine_space_spinu_pairs} of $\bg$ on the affine space of \spinu pairs, we have
\begin{multline*}
  \bg(\bJ(a,\phi),\bJ(a,\phi))
  = \|\bJ(a,\phi)\|_{L^2(X)}^2
  = \|(i\cdot a, i\cdot\phi)\|_{L^2(X)}^2
  = \|i\cdot a\|_{L^2(X)}^2 + \|i\cdot \phi\|_{L^2(X)}^2
  \\
  = \|a\|_{L^2(X)}^2 + \|\phi\|_{L^2(X)}^2,
\end{multline*}
and thus, by again applying the definition \eqref{eq:L2_metric_affine_space_spinu_pairs} of $\bg$, we have for all \spinu pairs $(A,\Phi)$ that
\begin{equation}
  \label{eq:Almost_complex_structure_affine_space_spinu_pairs_g-orthogonal}
  \bg(\bJ(a,\phi),\bJ(a,\phi)) = \bg((a,\phi),(a,\phi)),
  \quad
  \text{for all } (a,\phi) \in T_{A,\Phi}\left(\sA(E,H,A_d)\times W^{1,p}(W^+\otimes E)\right).
\end{equation}
By \eqref{eq:Almost_complex_structure_affine_space_spinu_pairs_g-orthogonal}, the almost complex structure $\bJ$ on the affine space $\sA(E,H,A_d)\times W^{1,p}(W^+\otimes E)$ of \spinu pairs is $\bg$-orthogonal, as claimed.

We claim that the endomorphism $\bJ$ is a $\bg$-\emph{skew-symmetric} automorphism of $T(\sA(E,H,A_d)\times W^{1,p}(W^+\otimes E))$. To prove this, we observe that
\begin{align*}
  \bg\left(\bJ(a_1,\phi_1),(a_2,\phi_2)\right)
  &=
    \Real\,\left((i\cdot a_1,i\phi_1),(a_2,\phi_2)\right)_{L^2(X)}
    \quad\text{(by \eqref{eq:L2_metric_affine_space_spinu_pairs})}
  \\
  &= (i\cdot a_1,a_2)_{L^2(X)} + \Real\,(i\phi_1,\phi_2)_{L^2(X)}
  \\
  &= -(a_1,i\cdot a_2)_{L^2(X)} - \Real\,(\phi_1,i\phi_2)_{L^2(X)}
    \quad\text{(by \eqref{eq:Almost_complex_structure_affine_space_spinu_pairs})}
  \\
  &= -\Real\,\left((a_1,\phi_1),(i\cdot a_2,i\phi_2)\right)_{L^2(X)},
\end{align*}
and thus by again appealing to \eqref{eq:L2_metric_affine_space_spinu_pairs}, we see that for all $(A,\Phi) \in \sA(E,H,A_d)\times W^{1,p}(W^+\otimes E)$,
\begin{multline}
  \label{eq:Almost_complex_structure_affine_space_spinu_pairs_skew-symmetric}
  \bg\left(\bJ(a_1,\phi_1),(a_2,\phi_2)\right)
  =
  -\bg\left((a_1,\phi_1),\bJ(a_2,\phi_2)\right),
  \\
  \text{for all } (a,\phi) \in T_{A,\Phi}\left(\sA(E,H,A_d)\times W^{1,p}(W^+\otimes E)\right).
\end{multline}
By \eqref{eq:Almost_complex_structure_affine_space_spinu_pairs_skew-symmetric}, the endomorphism $\bJ$ is a skew-symmetric automorphism of $T(\sA(E,H,A_d)\times W^{1,p}(W^+\otimes E))$ with respect to $\bg$, as claimed.

The almost complex structure $\bJ$ on the tangent bundle $T(\sA(E,H,A_d)\times W^{1,p}(W^+\otimes E))$ is clearly \emph{circle invariant} since, for any pair $(A,\Phi) \in \sA(E,H,A_d)\times W^{1,p}(W^+\otimes E))$, we have
\begin{multline}
\label{eq:Almost_complex_structure_affine_space_spinu_pairs_S1_invariance}
  \bJ(a,e^{i\theta}\phi) = e^{i\theta}\bJ(a,\phi)
  \in
  T_{A,e^{i\theta}\Phi}\left(\sA(E,H,A_d)\times W^{1,p}(W^+\otimes E)\right),
  \\
  \text{for all } (a,\phi) \in T_{A,\Phi}\left(\sA(E,H,A_d)\times W^{1,p}(W^+\otimes E)\right),
\end{multline}
by definition \eqref{eq:Almost_complex_structure_affine_space_spinu_pairs} of $\bJ$.
% TL6-21-2025: Added following.  It may be redundant from similar computation in Section 9.2.2 --they could be combined.
%PF6-23-2025: I don't see the point of this addition. This subsection focuses on the affine space. The next subsection discusses the quotient space. It's misplaced here and also, as you say, redundant since I have essentially the same discussion in the next subsection.
% In addition, we note that the almost complex structure is gauge equivariant in the sense that for any
% $(a,\varphi)\in W^{1,p}(T^*X\otimes\su(E)\oplus W^+\otimes E)$ and any $u\in W^{2,p}(\SU(E))$,
% \begin{equation}
% \label{eq:Almost_complex_structure_affine_space_spinu_pairs_gauge_equivariant}
% \bJ u^*(a,\phi)=\bJ u^*(a,\phi),
% \end{equation}
% where $u^*(a,\phi)$ is as defined in \eqref{eq:TangentBundleAction}.  One can see that \eqref{eq:Almost_complex_structure_affine_space_spinu_pairs_gauge_equivariant} holds because the scalar multiplication on the section commutes with the linear gauge group action while $\bJ$ acts on the $T^*X$ component of $T^*X\otimes\su(E)$ and thus commutes with the gauge transformation action on $\su(E)$.

\subsection[Skew-symmetric endomorphism on the quotient space]{Construction of skew-symmetric endomorphism on the quotient space of spin${}^u$ pairs}
\label{subsec:Skew-symmetric_endomorphism_on_quotient_space_Spinu_pairs}
The almost complex structure $\bJ$ on the tangent bundle for the affine space of spin${}^u$ pairs determines an endomorphism of the Zariski tangent bundle for the quotient space of spin${}^u$ pairs. Indeed, if $(A,\Phi)$ represents $[A,\Phi] \in \sC(E,H,W,A_d)$ and we represent the tangent space $T_{[A,\Phi]}\sC(E,H,W,A_d)$ by $T_{A,\Phi}\sC(E,H,W,A_d)$ as in \eqref{eq:Representative_of_tangent_space_to_quotient_space_spinu_pairs}, then we define
\begin{multline}
\label{eq:Almost_complex_structure_affine_space_spinu_pairs_Coulomb_gauge_slice}
  \mathbf{\hat J}_{A,\Phi}(a,\phi)
  := \Pi_{A,\Phi}\bJ(a, \phi)
  \in \Ker d_{A,\Phi}^{0,*}\cap W^{1,p}\left(T^*X\otimes \su(E) \oplus W^+\otimes E\right),
  \\
  \text{for all } (a,\phi) \in \Ker d_{A,\Phi}^{0,*}\cap W^{1,p}\left(T^*X\otimes \su(E) \oplus W^+\otimes E\right).
\end{multline}
%TL6-17-2025: This is good (aside from the change from $(A,\varphi,\psi)$ to $(A,\Phi)$ but could be moved a remark in 9.1 outlining what is meant by equivalent representations of a tangent vector.  I've put a comment in there which could be combined with this.  Then we'd be able to reference it later on as "independent of the choice of representation as described in Remark ..."
If $u \in W^{2,p}(\SU(E))$, so $u(A,\Phi)$ is another representative of $[A,\Phi] \in \sC(E,H,W,A_d)$, and
\[
  (a,\phi) \in \Ker d_{A,\Phi}^{0,*}\cap W^{1,p}\left(T^*X\otimes \su(E) \oplus W^+\otimes E\right),
\]
then
\[
  u(a,\phi) \in \Ker d_{u(A,\Phi)}^{0,*}\cap W^{1,p}\left(T^*X\otimes \su(E) \oplus W^+\otimes E\right).
\]
Moreover, by \eqref{eq:Almost_complex_structure_affine_space_spinu_pairs} we have
\begin{multline*}
  u(\bJ(a,\phi))
  = u(i\cdot a, i\cdot\phi)
  = u\left(-\frac{i}{2}a' + \frac{i}{2}a'', i\phi\right)
  = \left(- u^{-1}\left(\frac{i}{2}a'\right)u + u^{-1}\left(\frac{i}{2}a''\right)u, u^{-1}(i\phi)\right)
  \\
  = \left(-\frac{i}{2}\left(u^{-1}au\right)' + \frac{i}{2}\left(u^{-1}au\right)'', iu^{-1}\phi\right)
  = \left(i\cdot u^{-1}au, i\cdot u^{-1}\phi\right)
  = \bJ \left(u^{-1}au, u^{-1}\phi\right)
  = \bJ u(a,\phi),
  \\
  \text{for all } (a,\phi) \in W^{1,p}\left(T^*X\otimes \su(E) \oplus W^+\otimes E\right).
\end{multline*}
The preceding observations, \eqref{eq:Almost_complex_structure_affine_space_spinu_pairs_Coulomb_gauge_slice},
and the equivariance of $\Pi_{A,\Phi}$ given in \eqref{eq:SliceProjectionGaugeEquivariance} thus yield
\[
  \mathbf{\hat J}_{u(A,\Phi)}u(a,\phi)
  = \Pi_{u(A,\Phi)}\bJ u(a, \phi)
  = \Pi_{u(A,\Phi)}u\left(\bJ(a, \phi)\right)
  = u\left(\Pi_{A,\Phi}\bJ(a,\phi)\right),
\]
that is,
\begin{multline}
  \label{eq:Gauge_equivariance_hatJ}
  \mathbf{\hat J}_{u(A,\Phi)}u(a,\phi)
  =
  u\left(\Pi_{A,\Phi}\bJ(a,\phi)\right)
  \in \Ker d_{u(A,\Phi)}^{0,*}\cap W^{1,p}\left(T^*X\otimes \su(E) \oplus W^+\otimes E\right),
  \\
  \text{for all } u(a,\phi)
  \in \Ker d_{u(A,\Phi)}^{0,*}\cap W^{1,p}\left(T^*X\otimes \su(E) \oplus W^+\otimes E\right).
\end{multline}
Although the endomorphism $\mathbf{\hat J}$ is not necessarily $\mathbf{\hat g}$-orthogonal or an automorphism of $T\sC(E,H,W,A_d)$, we claim that $\mathbf{\hat J}$ is $\mathbf{\hat g}$-\emph{skew-symmetric}. To prove this, we compute:
\begin{align*}
  \mathbf{\hat g}_{A,\Phi}\left(\mathbf{\hat J}_{A,\Phi}(a_1,\phi_1),(a_2,\phi_2)\right)
  &=
    \Real\,\left(\mathbf{\hat J}_{A,\Phi}(a_1,\phi_1), (a_2,\phi_2)\right)_{L^2(X)}
    \quad\text{(by \eqref{eq:L2_metric_quotient_space_spinu_pairs})}
  \\
  &= \Real\,\left(\Pi_{A,\Phi}\bJ(a_1,\phi_1), (a_2,\phi_2)\right)_{L^2(X)}
    \quad\text{(by \eqref{eq:Almost_complex_structure_affine_space_spinu_pairs_Coulomb_gauge_slice})}
  \\
  &= \Real\,\left(\bJ(a_1,\phi_1), (a_2,\phi_2)\right)_{L^2(X)}
    \quad\text{(because $(a_2,\phi_2) \in \Ker d_{A,\Phi}^{0,*}$)}
  \\
  &= -\Real\,\left((a_1,\phi_1), \bJ(a_2,\phi_2)\right)_{L^2(X)}
    \quad\text{(by \eqref{eq:L2_metric_affine_space_spinu_pairs} and \eqref{eq:Almost_complex_structure_affine_space_spinu_pairs_g-orthogonal})}
  \\
  &= -\Real\,\left((a_1,\phi_1), \Pi_{A,\Phi}\bJ(a_2,\phi_2)\right)_{L^2(X)}
    \quad\text{(because $(a_1,\phi_1) \in \Ker d_{A,\Phi}^{0,*}$)}
  \\
  &= -\Real\,\left((a_1,\phi_1), \mathbf{\hat J}_{A,\Phi}(a_2,\phi_2)\right)_{L^2(X)}
    \quad\text{(by \eqref{eq:Almost_complex_structure_affine_space_spinu_pairs_Coulomb_gauge_slice})},
\end{align*}
and thus by again applying the definition \eqref{eq:L2_metric_quotient_space_spinu_pairs} of $\mathbf{\hat g}_{A,\Phi}$,
\begin{multline}
  \label{eq:Almost_complex_structure_quotient_tangent_space_spinu_pairs_skew_symmetric}
  \mathbf{\hat g}_{A,\Phi}\left(\mathbf{\hat J}_{A,\Phi}(a_1,\phi_1),(a_2,\phi_2)\right)
  =
  -\mathbf{\hat g}_{A,\Phi}\left((a_1,\phi_1),\mathbf{\hat J}_{A,\Phi}(a_2,\phi_2)\right),
  \\
  \text{for all } (a_1,\phi_1), (a_2,\phi_2)
  \in \Ker d_{A,\Phi}^{0,*}\cap W^{1,p}\left(T^*X\otimes \su(E) \oplus W^+\otimes E\right).
\end{multline}
We observe that the scalar values in \eqref{eq:Almost_complex_structure_quotient_tangent_space_spinu_pairs_skew_symmetric} are independent of the choice of representative $(A,\Phi)$ of the point $[A,\Phi] \in \sC(E,H,W,A_d)$ by replacing $(A,\Phi)$ by $u(A,\Phi)$ and $(a_k,\phi_k)$ by $u(a_k,\phi_k)$ for $k=1,2$, noting that the values are independent of $u \in W^{2,p}(\SU(E))$. By \eqref{eq:Almost_complex_structure_quotient_tangent_space_spinu_pairs_skew_symmetric}, the endomorphism $\mathbf{\hat J}$ is skew-symmetric with respect to the $L^2$ metric $\mathbf{\hat g}$ on the Zariski tangent space $T_{A,\Phi}\sC(E,H,W,A_d)$, as claimed.

\subsection[Circle invariance skew-symmetric endomorphism on the quotient space]{Circle invariance of the skew-symmetric endomorphism on the quotient space of spin${}^u$ pairs}
\label{subsec:Circle-invariant_skew-symmetric_endomorphism_on_quotient_space_Spinu_pairs}
We claim that the endomorphism $\mathbf{\hat J}$ is a \emph{circle invariant endomorphism} of the tangent bundle $T\sC(E,H,W,A_d)$ for the quotient space. To see this, choose a representative $(A,\Phi)$ for $[A,\Phi] \in \sC(E,H,W,A_d)$ and $e^{i\theta} \in S^1$ and
\[
  (a,\phi) \in \Ker d_{A,\Phi}^{0,*}\cap W^{1,p}\left(T^*X\otimes \su(E) \oplus W^+\otimes E\right),
\]
and compute:
\begin{align*}
  \mathbf{\hat J}_{A,e^{i\theta}\Phi}(a,e^{i\theta}\phi)
  &= \pi_{A,e^{i\theta}\Phi}\bJ(a, e^{i\theta}\phi) \quad\text{(by \eqref{eq:Almost_complex_structure_affine_space_spinu_pairs_Coulomb_gauge_slice})}
  \\
  &= \pi_{A,e^{i\theta}\Phi}\left(e^{i\theta}\bJ(a, \phi)\right)
    \quad\text{(by \eqref{eq:Almost_complex_structure_affine_space_spinu_pairs_S1_invariance})},
  \\
  &= e^{i\theta}\left(\Pi_{A,\Phi}\bJ(a, \phi)\right)
    \quad\text{(by \eqref{eq:S1_action_Coulomb_gauge_spinu_slice})}.
\end{align*}
We conclude that
\begin{multline}
  \label{eq:Almost_complex_structure_quotient_tangent_space_spinu_pairs_S1_invariance}
  \mathbf{\hat J}_{A,e^{i\theta}\Phi}(a,e^{i\theta}\phi)
  =  e^{i\theta}\left(\Pi_{A,\Phi}\bJ(a, \phi)\right)
   \in \Ker d_{A,e^{i\theta}\Phi}^{0,*}\cap W^{1,p}\left(T^*X\otimes \su(E) \oplus W^+\otimes E\right),
  \\
  \quad\text{for all }
  (a,\phi) \in \Ker d_{A,\Phi}^{0,*}\cap W^{1,p}\left(T^*X\otimes \su(E) \oplus W^+\otimes E\right)
  \text{ and } e^{i\theta} \in S^1.
\end{multline}
Thus, $\mathbf{\hat J}$ is a circle invariant endomorphism of $T\sC(E,H,W,A_d)$, as claimed.

%PF6-16-2025 Not much value added. Probable cut.
% \begin{rmk}[Explanation for why almost complex structure on affine space need not descend to almost complex structure on quotient space]
% \label{rmk:Explanation_ACstructure_affine_space_need_not_imply_ACstructure_quotient_space}
% Observe that if $u \in \Aut(E,H,A_d)$, then
% \[
%   u\cdot(A,\Phi) = (u^*A,u^{-1}\Phi).
% \]
% %PF8-21-2024 Not sure below makes sense; clarify
% If $(a,\phi)$ is a tangent vector to the affine space at $(A,\Phi)$, then $(u^{-1}au,-u\phi)$ is a tangent vector to the affine space at $(u^*A,u^{-1}\Phi)$. Indeed, for all $t \in \RR$,
% \begin{multline*}
%   u\cdot(A + ta, \Phi + t\phi) - u\cdot(A,\Phi)
%   =
%   (A + tu^{-1}au + u^{-1}d_Au, u^{-1}\Phi + tu^{-1}\phi) - u(A,\Phi)
%   \\
%   =
%   t(u^{-1}au, u^{-1}\phi).
% \end{multline*}
% Suppose that $a \perp \Ran d_A$, so $(a,d_A\xi)_{L^2(X)} = 0$ for all $\xi \in \Omega^0(\su(E))$. It does not necessarily follow that $(i\cdot a,d_A\xi)_{L^2(X)} = 0$ for all $\xi \in \Omega^0(\su(E))$. Indeed,
% \[
%   (i\cdot a,d_A\xi)_{L^2(X)}
%   =
%   (ia' - ia'',d_A\xi)_{L^2(X)}
%   =
%   i(a',\partial_A\xi)_{L^2(X)} - i(a'',\bar\partial_A\xi)_{L^2(X)},
% \]
% and the right-side of the last equality is not obviously equal to zero. The same remarks apply when $(a,\phi) \perp \Ran d_{A,\Phi}^0$.
% \end{rmk}

\section[Circle-invariant two-form on the quotient space]{Construction of circle-invariant two-forms on the affine and quotient spaces of spin${}^u$ pairs}
\label{sec:Construction_circle-invariant_two-form}
In Sections \ref{sec:Riemannian_L2_metrics_quotient_space_spinu_pairs} and \ref{sec:ACStructure_on_Affine_Spinu_pairs}, respectively, we constructed a circle invariant weak $L^2$ Riemannian metric $\bg$ and almost complex structure $\bJ$ on the tangent bundle $T(\sA(E,H,A_d)\times W^{1,p}(W^+\otimes E))$ for the affine space of \spinu pairs and then a circle invariant weak $L^2$  Riemannian metric $\mathbf{\hat g}$ and skew-symmetric endomorphism $\hat\bJ$ on the tangent bundle $T\sC(E,H,W,A_d)$ for the quotient space of \spinu pairs. We now construct the corresponding circle-invariant two-forms on $\sA(E,H,A_d)\times W^{1,p}(W^+\otimes E)$ and $\sC(E,H,W,A_d)$.

\subsection[Circle-invariant weakly non-degenerate two-form on affine space]{Construction of circle-invariant weakly non-degenerate two-form on the affine space of spin${}^u$ pairs}
\label{subsec:Construction_circle-invariant_weakly_non-degenerate_two-form_affine_space}
For each pair $(A,\Phi) \in \sA(E,H,A_d)\times W^{1,p}(W^+\otimes E)$, we define (by analogy with Huybrechts \cite[Definition 1.2.13, p. 29]{Huybrechts_2005}) a weakly non-degenerate, real bilinear form on $T_{A,\Phi}(\sA(E,H,A_d)\times W^{1,p}(W^+\otimes E))$ by setting
\begin{multline}
  \label{eq:Fundamental_two_form_affine_space_spinu_pairs}
  \bomega((a_1,\phi_1),(a_2,\phi_2))
  := \bg\left(\bJ(a_1,\phi_1),(a_2,\phi_2)\right),
  \\
  \text{for all } (a_1,\phi_1),(a_2,\phi_2) \in T_{A,\Phi}\left(\sA(E,H,A_d)\times W^{1,p}(W^+\otimes E)\right),
\end{multline}
where $\bg$ is as in \eqref{eq:L2_metric_affine_space_spinu_pairs} and $\bJ$ is as in \eqref{eq:Almost_complex_structure_affine_space_spinu_pairs}. The bilinear form $\bomega$ is skew-symmetric by \eqref{eq:Almost_complex_structure_affine_space_spinu_pairs_skew-symmetric}. The resulting triple $(\bg,\bJ,\bomega)$ is compatible in the sense of McDuff and Salamon \cite[Section 4.1, p. 153]{McDuffSalamonSympTop3}.
% PF8-26-2024 Reference reasons why

Because $\bg$ and $\bJ$ are circle invariant by \eqref{eq:L2_metric_affine_space_spinu_pairs_circle_invariance} and \eqref{eq:Almost_complex_structure_affine_space_spinu_pairs_S1_invariance}, respectively, it follows immediately that for any pair $(A,\Phi) \in \sA(E,H,A_d)\times W^{1,p}(W^+\otimes E))$ and $e^{i\theta} \in S^1$, we have
\begin{multline}
  \label{eq:Fundamental_two-form_affine_space_spinu_pairs_S1_invariance}
  \bomega_{A,e^{i\theta}\Phi}\left((a_1,e^{i\theta}\phi_1),(a_2,e^{i\theta}\phi_2)\right)
  =
  \bomega_{A,\Phi}\left((a_1,\phi_1),(a_2,\phi_2)\right),
  \\
  \text{for all } (a_1,\phi_1),(a_2,\phi_2) \in T_{A,\Phi}\left(\sA(E,H,A_d)\times W^{1,p}(W^+\otimes E)\right),
\end{multline}
and so $\bomega$ is also \emph{circle invariant}.

\subsection[Circle-invariant two-form on quotient space]{Construction of circle-invariant two-form on the quotient space of spin${}^u$ pairs}
\label{subsec:Construction_circle-invariant_two-form_quotient_space}
For a representative $(A,\Phi)$ of a point $[A,\Phi] \in \sC(E,H,W,A_d)$, we define (again by analogy with Huybrechts \cite[Definition 1.2.13, p. 29]{Huybrechts_2005}) a real bilinear form on $T_{A,\Phi}\sC(E,H,W,A_d)$ by setting
\begin{multline}
  \label{eq:Fundamental_two_form_quotient_space_spinu_pairs}
  \mathbf{\hat\bomega}_{A,\Phi}((a_1,\phi_1),(a_2,\phi_2))
  := \mathbf{\hat g}_{A,\Phi}\left(\mathbf{\hat J}_{A,\Phi}(a_1,\phi_1),(a_2,\phi_2)\right),
  \\
  \text{for all } (a_1,\phi_1),(a_2,\phi_2)
  \in \Ker d_{A,\Phi}^{0,*}\cap W^{1,p}\left(T^*X\otimes \su(E) \oplus W^+\otimes E\right),
\end{multline}
where $\mathbf{\hat g}_{A,\Phi}$ is as in \eqref{eq:L2_metric_quotient_space_spinu_pairs} and $\mathbf{\hat J}_{A,\Phi}$ is as in \eqref{eq:Almost_complex_structure_affine_space_spinu_pairs_Coulomb_gauge_slice}. The bilinear form $\hat\bomega_{A,\Phi}$ is skew-symmetric by \eqref{eq:Almost_complex_structure_quotient_tangent_space_spinu_pairs_skew_symmetric}.

Because $\mathbf{\hat g}$ and $\mathbf{\hat J}$ are circle invariant by \eqref{eq:L2_metric_quotient_space_spinu_pairs_circle_invariance} and \eqref{eq:Almost_complex_structure_quotient_tangent_space_spinu_pairs_S1_invariance}, respectively, it follows immediately that, for any $e^{i\theta} \in S^1$, we have
\begin{multline}
  \label{eq:Fundamental_two-form_quotient_space_spinu_pairs_S1_invariance}
  \mathbf{\hat\bomega}_{A,e^{i\theta}\Phi}\left((a_1,e^{i\theta}\phi_1),(a_2,e^{i\theta}\phi_2)\right)
  =
  \mathbf{\hat\bomega}_{A,\Phi}\left((a_1,\phi_1),(a_2,\phi_2)\right),
  \\
  \text{for all } (a_1,\phi_1),(a_2,\phi_2)
  \in \Ker d_{A,\Phi}^{0,*}\cap W^{1,p}\left(T^*X\otimes \su(E) \oplus W^+\otimes E\right),
\end{multline}
and so $\mathbf{\hat\bomega}$ is also \emph{circle invariant}.

\section[Characterization of critical points of Hitchin's function]{Characterization of critical points of Hitchin's function for the circle action on the quotient space of unitary triples and moduli subspace of non-Abelian monopoles}
\label{subsec:Characterization_critical_points_Hamiltonian_function_circle_action_on_non-abelian_monopoles}
We adapt the arguments in Feehan and Leness
%TL12-4-2025: Updated
\cite[Section 12.6.6]{Feehan_Leness_introduction_virtual_morse_theory_so3_monopoles} from the special case of the moduli space of projective vortices over a complex K\"ahler surface to the moduli space of non-Abelian monopoles over an almost K\"ahler four-manifold. In particular, we provide an analogue of
%TL12-4-2025: Updated
\cite[Definition 12.6.16]{Feehan_Leness_introduction_virtual_morse_theory_so3_monopoles}, which gives the definition of a critical point of Hitchin's function, and prove Theorems \ref{thm:Critical_points_Hitchin_Hamiltonian_function_moduli_space_non-Abelian_monopoles} and \ref{mainthm:IdentifyCriticalPoints}, generalizations of
%TL12-4-2025: Updated.  Changed Thm. 4 to Thm. 8, please check if this is accurate
\cite[Theorems 12.6.17 and 8]{Feehan_Leness_introduction_virtual_morse_theory_so3_monopoles}, which characterize critical points of Hitchin's function on the moduli spaces of projective vortices and non-Abelian monopoles, respectively, over a complex K\"ahler surface.

\subsection{Hitchin's  function and circle-invariant non-degenerate two-form on the affine space of unitary triples}
\label{subsec:Hamiltonian_function_circle-invariant_non-degenerate_2-form_affine_space}
To see that Hitchin's function $f$ in \eqref{eq:Hitchin_function_affine} is a Hamiltonian in the sense of \eqref{eq:MomentMap} for the circle action \eqref{eq:S1ZAction} of Definition \ref{defn:UnitaryZActionOnAffine} on the affine space $\sA(E,H,A_d)\times W^{1,p}(E\oplus\Lambda^{0,2}(E))$ of $W^{1,p}$ unitary triples \eqref{eq:A_varphi_psi_in_W1p}, observe that the vector field $\bX$ on $\sA(E,H,A_d)\times W^{1,p}(E\oplus\Lambda^{0,2}(E))$ generating the action \eqref{eq:S1ZAction} of $e^{i\theta}\in S^1$ on $W^{1,p}(E\oplus\Lambda^{0,2}(E))$ is given by
\begin{multline*}
  \bX_{A,\varphi,\psi}
  = \left.\frac{d}{d\theta}e^{i\theta}\cdot(A,\varphi,\psi)\right|_{\theta=0}
  = \left.\frac{d}{d\theta}(A,e^{i\theta}\varphi,e^{i\theta}\psi)\right|_{\theta=0}
  \\
  = (0,i\varphi,i\psi)
  = \bJ(0,\varphi,\psi)
  \in W^{1,p}\left(T^*X\otimes\su(E) \oplus E\oplus\Lambda^{0,2}(E)\right),
\end{multline*}
where the final equality follows from the definition \eqref{eq:Almost_complex_structure_affine_space_spinu_pairs} of the almost complex structure $\bJ$ on the affine space \eqref{eq:A_varphi_psi_in_W1p}. We record this equality as
\begin{equation}
\label{eq:S1GeneratorExpression}
\bX_{A,\varphi,\psi}
 = \bJ(0,\varphi,\psi), \quad\text{for all } (A,\varphi,\psi)
 \in \sA(E,H,A_d)\times W^{1,p}\left(E\oplus\Lambda^{0,2}(E)\right).
\end{equation}
For $(a,\sigma,\tau) \in W^{1,p}(T^*X\otimes\su(E) \oplus E\oplus\Lambda^{0,2}(E))$, we have
\begin{multline*}
  (df)_{A,\varphi,\psi}(a,\sigma,\tau)
  = \left.\frac{d}{dt}f(A+ta,\varphi+t\sigma,\psi+t\tau)\right|_{t=0}
  = \frac{1}{2}\left.\frac{d}{dt}\|(\varphi,\psi)+t(\sigma,\tau)\|_{L^2(X)}^2\right|_{t=0}
  \\
  = \frac{1}{2}\left(\left(\varphi,\psi),(\sigma,\tau)\right)_{L^2(X)}
    + \left((\sigma,\tau), (\varphi,\psi)\right)_{L^2(X)}\right)
  = \Real\,\left((\varphi,\psi),(\sigma,\tau)\right)_{L^2(X)},
\end{multline*}
that is, by the definition \eqref{eq:L2_metric_affine_space_spinu_pairs} of the Riemannian metric $\bg$ on the affine space \eqref{eq:A_varphi_psi_in_W1p},
\begin{multline}
  \label{eq:df_at_Avarphipsi_direction_asigmatau}
  (df)_{A,\varphi,\psi}(a,\sigma,\tau)
  =
  \bg\left((0,\varphi,\psi),(0,\sigma,\tau)\right)_{L^2(X)},
  \\
  \text{for all } (a,\sigma,\tau)
  \in T_{A,\varphi,\psi}\left(\sA(E,H,A_d)\times W^{1,p}\left(E\oplus\Lambda^{0,2}(E)\right)\right).
\end{multline}
On the other hand, we see that
\begin{align*}
  \left(\iota_\bX\bomega\right)_{A,\varphi,\psi}(a,\sigma,\tau)
  &= \iota_{\bX_{A,\varphi,\psi}}\bomega_{A,\varphi,\psi}(a,\sigma,\tau)
  \\
  &= \bomega(\bX_{A,\varphi,\psi},(a,\sigma,\tau))
  \\
  &= \bomega(\bJ(0,\varphi,\psi),(a,\sigma,\tau)) \quad\text{(by \eqref{eq:S1GeneratorExpression})}
  \\
  &= \bg(\bJ^2(0,\varphi,\psi),(a,\sigma,\tau)) \quad\text{(by \eqref{eq:Fundamental_two_form_affine_space_spinu_pairs})}
  \\
  &= -\bg((0,\varphi,\psi),(a,\sigma,\tau)) \quad\text{(by \eqref{eq:Almost_complex_structure_affine_space_spinu_pairs})}
  \\
  &= -(0,a)_{L^2(X)} - \Real\,\left((\varphi,\psi),(\sigma,\tau)\right)_{L^2(X)}
    \quad\text{(by \eqref{eq:L2_metric_affine_space_spinu_pairs})},
\end{align*}
and therefore we obtain
\begin{multline}
  \label{eq:iota_bX_bomega}
 \left(\iota_\bX\bomega\right)_{A,\varphi,\psi}(a,\sigma,\tau)
  = -\bg\left((0,\varphi,\psi),(0,\sigma,\tau)\right)_{L^2(X)},
  \\
  \text{for all } (A,\varphi,\psi) \in \sA(E,H,A_d) \times W^{1,p}\left(E\oplus\Lambda^{0,2}(E)\right)
  \\
  \text{and } (a,\sigma,\tau)
  \in T_{A,\varphi,\psi}\left(\sA(E,H,A_d)\times W^{1,p}\left(E\oplus\Lambda^{0,2}(E)\right)\right).
\end{multline}
Thus, we see that \eqref{eq:df_at_Avarphipsi_direction_asigmatau} and \eqref{eq:iota_bX_bomega} yield the equality,
\begin{equation}
  \label{eq:Moment_map_affine_space_unitary_pairs}
  -(df)_{A,\varphi,\psi} = \left(\iota_\bX\bomega\right)_{A,\varphi,\psi},
  \quad\text{for all } (A,\varphi,\psi) \in \sA(E,H,A_d)\times W^{1,p}\left(E\oplus\Lambda^{0,2}(E)\right).
\end{equation}
Consequently, $-f$ in \eqref{eq:Hitchin_function} is a Hamiltonian function in the sense of \eqref{eq:MomentMap} for the circle action \eqref{eq:S1ZAction} on $\sA(E,H,A_d)\times W^{1,p}(E\oplus\Lambda^{0,2}(E))$. (The sign for the Hamiltonian $f$ in \eqref{eq:Moment_map_affine_space_unitary_pairs} matches that of Hitchin \cite[Section 6, last displayed equation, p. 92]{Hitchin_1987}.)
% PF4-23-2025 Adjust sign in intro to match?

The expression \eqref{eq:df_at_Avarphipsi_direction_asigmatau} for $df$ implies that $f$ has \emph{no critical points} on $\sA(E,H,A_d)\times W^{1,p}(E\oplus\Lambda^{0,2}(E))$ with $(\varphi,\psi) \not\equiv (0,0)$. Indeed, if $df(a,\sigma,\tau) = 0$ for all $(a,\sigma,\tau) \in W^{1,p}(T^*X\otimes\su(E) \oplus E\oplus\Lambda^{0,2}(E))$, then \eqref{eq:df_at_Avarphipsi_direction_asigmatau} yields $(\varphi,\psi) \equiv (0,0)$ by choosing $(a,\sigma,\tau) = (0,\varphi,\psi)$.

\subsection{Hitchin's function and circle-invariant two-form on the quotient space of unitary triples}
\label{subsec:Hamiltonian_function_circle-invariant_non-degenerate_2-form_quotient_space}
In this section, when the forthcoming Condition \ref{cond:Pi_Avarphipsi_bJ_invertible_on_Coulomb_gauge_slice} holds, we show that the relation \eqref{eq:Moment_map_affine_space_unitary_pairs} descends to the quotient to give
\begin{multline}
  \label{eq:Moment_map_quotient_space_unitary_triples}
  -(df)_{A,\varphi,\psi} = \left(\iota_{\mathbf{\hat X}}\mathbf{\tilde\bomega}\right)_{A,\varphi,\psi}
  \quad\text{on } \Ker d_{A,\varphi,\psi}^{0,*}\cap W^{1,p}\left(T^*X\otimes\su(E)\oplus E\oplus\Lambda^{0,2}(E)\right),
  \\
  \text{for all } [A,\varphi,\psi] \in \sC^0(E,H,J,A_d),
\end{multline}
where $\mathbf{\tilde\bomega}$ is the non-degenerate two-form in the forthcoming definition \eqref{eq:Modified_fundamental_two_form_quotient_space_spinu_pairs} and $\mathbf{\hat X}$ is the smooth vector field generating the circle action \eqref{eq:S1ZActionOnQuotientSpace} on the open Banach manifold $\sC^0(E,H,J,A_d)$ in \eqref{eq:Quotient_space_non-zero-section_unitary_triples}.

The left-hand side of \eqref{eq:Moment_map_quotient_space_unitary_triples} is independent of the choice of representative $(A,\varphi,\psi)$ of the point $[A,\varphi,\psi] \in \sC^0(E,H,J,A_d)$ in the following sense: For any $u \in W^{2,p}(\SU(E))$, the identity \eqref{eq:df_at_Avarphipsi_direction_asigmatau} and definition \eqref{eq:L2_metric_affine_space_spinu_pairs} of $\bg$ yield
\begin{multline*}
  (df)_{(u(A,\varphi,\psi)}u(a,\sigma,\tau)
  =
  \Real\,\left(u(0,\varphi,\psi),u(0,\sigma,\tau)\right)_{L^2(X)}
  \\
  =
  \Real\,\left((0,\varphi,\psi),(0,\sigma,\tau)\right)_{L^2(X)}
  = (df)_{A,\varphi,\psi}(a,\sigma,\tau),
  \\
  \text{for all } (a,\sigma,\tau)
  \in W^{1,p}\left(T^*X\otimes\su(E)\oplus E\oplus\Lambda^{0,2}(E)\right),
\end{multline*}
and the same remark applies if the vector $(0,\varphi,\psi)$ above is replaced by $\Pi_{A,\varphi,\psi}(0,\varphi,\psi)$, where
% TL6-17-2025: df vanishes on $\Ran d_{A,\varphi,\psi}^0$ so $(0,\varphi,\psi)\in\Ker d_{A,\varphi,\psi}^{0,*}$ so $$\Pi_{A,\varphi,\psi}(0,\varphi,\psi)=(0,\varphi,\psi)$.  (This uses my understanding of the slice as a real and not complex orthogonal complement so it is subject to that being cleared up --it looks as if the equality following (9.10.4) makes the decomposition real orthogonal.)
% PF6-18-2025 There is no (9.10.4)?
%PF6-18-2025 (9.4.10) is simply saying that the real L^2 inner product is an inner product on a real vector space. Nothing wrong there. That’s exactly how the symplectic form is defined throughout the monograph, using the real L^2 inner product.
\begin{multline}
  \label{eq:L2_orthogonal_projection_onto_slice}
  \Pi_{A,\varphi,\psi}: W^{1,p}\left(T^*X\otimes\su(E)\oplus E\oplus\Lambda^{0,2}(E)\right)
  \\
  \to
  \Ker d_{A,\varphi,\psi}^{0,*}\cap W^{1,p}\left(T^*X\otimes\su(E)\oplus E\oplus\Lambda^{0,2}(E)\right)
\end{multline}
is the $L^2$-orthogonal projection onto the Coulomb gauge slice through $(A,\varphi,\psi)$, analogous to \eqref{eq:L2_orthogonal_projection_onto_slice_spin-u_pairs} for \spinu pairs. After pullback to the Coulomb gauge slice through $(A,\varphi,\psi)$, the identity \eqref{eq:df_at_Avarphipsi_direction_asigmatau} becomes
\begin{multline}
  \label{eq:df_Avarphipsi_asigmatau_Coulomb_gauge_slice}
  (df)_{A,\varphi,\psi}(a,\sigma,\tau)
  =
  \bg\left(\Pi_{A,\varphi,\psi}(0,\varphi,\psi),(a,\sigma,\tau)\right),
  \\
  \text{for all } (a,\sigma,\tau)
  \in \Ker d_{A,\varphi,\psi}^{0,*}\cap W^{1,p}\left(T^*X\otimes\su(E)\oplus E\oplus\Lambda^{0,2}(E)\right).
\end{multline}
We now examine the right-hand side of the identity \eqref{eq:Moment_map_quotient_space_unitary_triples}. We first verify an explicit expression for the vector field $\mathbf{\hat X}$ generating the circle action \eqref{eq:S1ZActionOnQuotientSpace} on $\sC^0(E,H,J,A_d)$. Recall from \eqref{eq:S1GeneratorExpression} that the vector field $\bX$ on the affine space $\sA(E,H,A_d)\times W^{1,p}(E\oplus\Lambda^{0,2}(E))$ is given by
\[
  \bX_{A,\varphi,\psi}
 = \bJ(0,\varphi,\psi)
 \in T_{A,\varphi,\psi}\left(\sA(E,H,A_d)\times W^{1,p}\left(E\oplus\Lambda^{0,2}(E)\right)\right),
\]
for all $(A,\varphi,\psi) \in \sA(E,H,A_d)\times W^{1,p}(E\oplus\Lambda^{0,2}(E))$. We define a smooth map by
\begin{multline}
  \label{eq:S1GeneratorExpression_gauge_equivariant_map}
  \mathbf{\hat X}: \sA(E,H,A_d)\times W^{1,p}\left(E\oplus\Lambda^{0,2}(E)\right)
    \ni
    (A,\varphi,\psi)
    \\
    \mapsto
    \mathbf{\hat X}_{A,\varphi,\psi} := \Pi_{A,\varphi,\psi}\bJ(0,\varphi,\psi)
    \in
    W^{1,p}\left(T^*X\otimes\su(E)\oplus E\oplus\Lambda^{0,2}(E)\right).
\end{multline}
Observe that
\begin{multline*}
  \Pi_{A,\varphi,\psi}\bJ(0,\varphi,\psi)
  \in
  \Ker d_{A,\varphi,\psi}^{0,*}\cap W^{1,p}\left(T^*X\otimes\su(E)\oplus E\oplus\Lambda^{0,2}(E)\right)
  \\
  =
  T_{A,\varphi,\psi}\sC(E,H,J,A_d)
  \cong
  T_{A,\varphi,\psi}\left.\left(\sA(E,H,A_d)\times W^{1,p}\left(E\oplus\Lambda^{0,2}(E)\right)\right)\right/
\Ran d_{A,\varphi,\psi}^0.
\end{multline*}
The map $\mathbf{\hat X}$ is $W^{2,p}(\SU(E))$-equivariant since, for any $u \in W^{2,p}(\SU(E))$,
% TL6-17-2025: (This may appear in another comment --if you've already put this in another draft then a reference to it should be here.  We need to put in the equality $u\Pi_{u(A,\varphi,\psi)}u^{-1}=\Pi_{A,\varphi,\psi}$
%PF6-18-2025 I don't understand. What equality is missing below?
%TL6-21-2025: No equality was missing --just a reference to the gauge-equivariance equivariance of $\Pi$
%TL6-21-2025: Could also say that this follows from the gauge-equavriance of $\bJ$ in \eqref\eqref{eq:Almost_complex_structure_affine_space_spinu_pairs_gauge_equivariant} and that of $\Pi_{A,\varphi,\psi}$ in \eqref\eqref{eq:SliceProjectionGaugeEquivariance}
\begin{multline*}
  \mathbf{\hat X}_{u(A,\varphi,\psi)}
= \Pi_{u(A,\varphi,\psi)}\bJ u(0,\varphi,\psi)
= \Pi_{u(A,\varphi,\psi)}\bJ\left(0,u^{-1}\varphi,u^{-1}\psi\right)
= \Pi_{u(A,\varphi,\psi)}\left(0,iu^{-1}\varphi,iu^{-1}\psi\right)
\\
= \Pi_{u(A,\varphi,\psi)}\left(0,u^{-1}(i\varphi),u^{-1}(i\psi)\right)
= \Pi_{u(A,\varphi,\psi)}u(0,i\varphi,i\psi)
= \Pi_{u(A,\varphi,\psi)}u\bJ(0,\varphi,\psi)
\\
= u\left(\Pi_{A,\varphi,\psi}\bJ(0,\varphi,\psi)\right)
\in T_{u(A,\varphi,\psi)}\left(\sA(E,H,A_d)\times W^{1,p}\left(E\oplus\Lambda^{0,2}(E)\right)\right),
\end{multline*}
where we appeal to the definition \eqref{eq:Almost_complex_structure_affine_space_spinu_pairs} of $\bJ$ to obtain the preceding equalities.  Hence, the $W^{2,p}(\SU(E))$-equivariant map $\mathbf{\hat X}$ in \eqref{eq:S1GeneratorExpression_gauge_equivariant_map} defines a smooth vector field on the quotient space $\sC^0(E,H,J,A_d)$ by
\begin{equation}
\label{eq:S1GeneratorExpression_quotient_space}
\mathbf{\hat X}_{A,\varphi,\psi}
:= \Pi_{A,\varphi,\psi}\bJ(0,\varphi,\psi)
\in \Ker d_{A,\varphi,\psi}^{0,*}\cap W^{1,p}\left(T^*X\otimes\su(E)\oplus E\oplus\Lambda^{0,2}(E)\right),
\end{equation}
where $(A,\varphi,\psi)$ represents a point $[A,\varphi,\psi] \in \sC^0(E,H,J,A_d)$.

We claim that $\mathbf{\hat X}$ is the generator of the $S^1$ action \eqref{eq:S1ZActionOnQuotientSpace} on $\sC^0(E,H,J,A_d)$. If $[A,\varphi,\psi] \in \sC^0(E,H,J,A_d)$, then
\[
  e^{i\theta}\cdot[A,\varphi,\psi] = [A,e^{i\theta}\varphi,e^{i\theta}\psi] \in \sC^0(E,H,J,A_d),
  \quad\text{for } \theta \in (-\eps,\eps),
\]
defines a smooth curve in $\sC^0(E,H,J,A_d)$ that is the image of the following smooth curve in $\sA(E,H,A_d)\times W^{1,p}(E\oplus\Lambda^{0,2}(E))$,
\[
  e^{i\theta}\cdot(A,\varphi,\psi) = (A,e^{i\theta}\varphi,e^{i\theta}\psi)
  \in \sA(E,H,A_d) \times  W^{1,p}(E\oplus\Lambda^{0,2}(E)),
  \quad\text{for } \theta \in (-\eps,\eps).
\]
These curves are smooth since the $S^1$ actions \eqref{eq:S1ZActionOnQuotientSpace} and \eqref{eq:S1ZAction} are smooth. The tangent vector to this curve at $\theta=0$ is
\[
  \left.\frac{d}{d\theta}(A,e^{i\theta}\varphi,e^{i\theta}\psi)\right|_{\theta=0}
  =
  (0,i\varphi,i\psi)
  =
  \bJ(0,\varphi,\psi) \in W^{1,p}(T^*X\otimes\su(E)) \oplus W^{1,p}(E\oplus\Lambda^{0,2}(E)),
\]
where we apply the definition \eqref{eq:Almost_complex_structure_affine_space_spinu_pairs} of $\bJ$ to obtain the second equality above. The projection of the preceding tangent vector to the affine space onto the slice $\Ker d_{A,\varphi,\psi}^{0,*}$ through $(A,\varphi,\psi)$ yields
\[
  \mathbf{\hat X}_{A,\varphi,\psi}
  =
  \Pi_{A,\varphi,\psi}\bJ(0,\varphi,\psi)
  \in
  \Ker d_{A,\varphi,\psi}^{0,*}.
\]
The curve $e^{i\theta}\cdot[A,\varphi,\psi]$ defines the tangent vector $\mathbf{\hat X}_{A,\varphi,\psi}$ to $\sC^0(E,H,J,A_d)$ at $(A,\varphi,\psi)$.
This verifies the claim that $\mathbf{\hat X}$ is the generator of the $S^1$ action  \eqref{eq:S1ZActionOnQuotientSpace} on $\sC^0(E,H,J,A_d)$. We now establish the 

\begin{lem}[Hitchin's function $f$ has no critical points in the quotient space of non-zero section triples]
\label{lem:Critical_points_f_on_quotient_space_are_fixed_points_S1_action}  
Let $(X,g,J,\omega)$ be an almost Hermitian manifold and $(E,H)$ be a smooth, Hermitian vector bundle over $X$ with complex rank $2$ and smooth, unitary connection $A_d$ on $\det E$. For $\sC^0(E,H,J,A_d)$ as in \eqref{eq:Quotient_space_non-zero-section_unitary_triples}, Hitchin's function $f:\sC^0(E,H,J,A_d) \to \RR$ in \eqref{eq:Hitchin_function} has no critical points in the sense of Definition \ref{defn:Critical_point_Hitchin_function_moduli_space_non-Abelian_monopoles}.
\end{lem}

\begin{proof}
Suppose that $[A,\varphi,\psi] \in \sC^0(E,H,J,A_d)$ is a critical point of $f$. Then Definition \ref{defn:Critical_point_Hitchin_function_moduli_space_non-Abelian_monopoles} and equation \eqref{eq:df_Avarphipsi_asigmatau_Coulomb_gauge_slice} imply that
\begin{multline*}
  (df)_{A,\varphi,\psi}(a,\sigma,\tau)
  =
  \bg\left(\Pi_{A,\varphi,\psi}(0,\varphi,\psi),(a,\sigma,\tau)\right) = 0,
  \\
  \text{for all } (a,\sigma,\tau)
  \in \Ran \Pi_{A,\varphi,\psi} = \Ker d_{A,\varphi,\psi}^{0,*}\cap W^{1,p}\left(T^*X\otimes\su(E)\oplus E\oplus\Lambda^{0,2}(E)\right),
\end{multline*}
where $\Pi_{A,\varphi,\psi}$ is as in \eqref{eq:L2_orthogonal_projection_onto_slice}. Therefore, because the real $L^2$ Riemannian metric $\bg$ in \eqref{eq:L2_metric_affine_space_spinu_pairs} is weakly non-degenerate on $\Ran \Pi_{A,\varphi,\psi}$ in the sense of Section \ref{sec:Weakly_strongly_non-degenerate_bilinear_forms_Banach_spaces}, we have
\begin{equation}
  \label{eq:Vanishing_Pi(0,varphi,psi)_at_critical_point}
  \Pi_{A,\varphi,\psi}(0,\varphi,\psi) = 0.
\end{equation}
Because there is a real--Hermitian and thus a real $L^2$ orthogonal direct sum,
\[
  L^2\left(T^*X\otimes\su(E)\oplus E\oplus\Lambda^{0,2}(E)\right)
  =
  \Ker d_{A,\varphi,\psi}^{0,*} \oplus \Ran d_{A,\varphi,\psi}^0
  =
  \Ran \Pi_{A,\varphi,\psi} \oplus \Ker \Pi_{A,\varphi,\psi},
\]
it follows that
\[
  (0,\varphi,\psi)
  \in
  \Ran d_{A,\varphi,\psi}^0 \cap \Omega^0\left(T^*X\otimes\su(E)\oplus E\oplus\Lambda^{0,2}(E)\right).
\]
Hence, by definition \eqref{eq:nonAbelianMonopole_d0_unitary} of $d_{A,\varphi,\psi}^0$ and our assumption that $(\varphi,\psi) \not \equiv (0,0)$, there is a $\xi \in \Omega^0(\su(E))\less \{0\}$ such that
\begin{equation}
  \label{eq:(0,varphi,psi)_is_d_A_varphi_psi_xi_at_critical_point}
  (0,\varphi,\psi) = d_{A,\varphi,\psi}^0\xi = \left(d_A\xi,-\xi\varphi,-\xi\psi\right).
\end{equation}
In particular, $d_A\xi = 0$ and thus by Feehan and Leness \cite[Proposition 6.1.10 (1)]{Feehan_Leness_introduction_virtual_morse_theory_so3_monopoles}, the connection $A$ is split with respect to an orthogonal decomposition $E=L_1\oplus L_2$ as a direct sum of Hermitian line bundles, defined by the property that $\xi = -i\mu$ on $L_1$ and $\xi_2= i\mu$ on $L_2$ for $\mu \in \RR\less\{0\}$ using \cite[Proposition 6.1.10 (1)]{Feehan_Leness_introduction_virtual_morse_theory_so3_monopoles}. But equation \eqref{eq:(0,varphi,psi)_is_d_A_varphi_psi_xi_at_critical_point} also implies that $(\xi\varphi,\xi\psi) = (-\varphi,-\psi)$, whose only solution is $(\varphi,\psi) \equiv (0,0)$ since $\xi$ has (non-zero) purely imaginary eigenvalues. This contradicts our assumption that $(\varphi,\psi) \not\equiv (0,0)$, implicit in our definition \eqref{eq:Quotient_space_non-zero-section_unitary_triples} of $\sC^0(E,H,J,A_d)$ and completes the proof of Lemma \ref{lem:Critical_points_f_on_quotient_space_are_fixed_points_S1_action}.

We may also give a simpler, second proof of Lemma \ref{lem:Critical_points_f_on_quotient_space_are_fixed_points_S1_action} as follows. Since $[A,\varphi,\psi] \in \sC^0(E,H,J,A_d)$ is a critical point of $f$ by hypothesis, then $(df)_{A,\varphi,\psi} = 0$ on $\Ker d_{A,\varphi,\psi}^{0,*}$. Because $f:\sA(E,H,A_d) \times W^{1,p}(E\oplus\Lambda^{0,2}(E))$ is $W^{2,p}(\SU(E))$-invariant, we also have $(df)_{A,\varphi,\psi} = 0$ on $\Ran d_{A,\varphi,\psi}^0$ and thus $(df)_{A,\varphi,\psi} = 0$ on the tangent space to $\sA(E,H,A_d) \times W^{1,p}(E\oplus\Lambda^{0,2}(E))$, so $(A,\varphi,\psi)$ is a critical point of $f$ in \eqref{eq:Hitchin_function_affine}. However, as we saw in Section \ref{subsec:Hamiltonian_function_circle-invariant_non-degenerate_2-form_affine_space}, that implies $(\varphi,\psi) \equiv (0,0)$.

Finally, we may give an even simpler, third proof of Lemma \ref{lem:Critical_points_f_on_quotient_space_are_fixed_points_S1_action} as follows. From the calculation in the forthcoming Example \ref{exmp:Failure_of_AC_Structure_Inj_At_Split_Triples}, we have for \emph{any} unitary triple $[A,\varphi,\psi]$ that $(0,\varphi,\psi) \in \Ker d_{A,\varphi,\psi}^{0,*}$ and so $\Pi_{A,\varphi,\psi}(0,\varphi,\psi) = (0,\varphi,\psi)$. Therefore, if $[A,\varphi,\psi]$ is a critical point of $f$, then \eqref{eq:Vanishing_Pi(0,varphi,psi)_at_critical_point} implies that $(\varphi,\psi) \equiv (0,0)$, contradicting our hypothesis that $[A,\varphi,\psi] \in \sC^0(E,H,J,A_d)$, the quotient space of non-zero section unitary triples.
\end{proof}

We shall now prove that for $f$ in \eqref{eq:Hitchin_function}, the function $-f$ is a Hamiltonian for the $S^1$ action \eqref{eq:S1ZActionOnQuotientSpace} on $\sC^0(E,H,J,A_d)$ for a suitable non-degenerate two-form under the forthcoming technical Condition \ref{cond:Pi_Avarphipsi_bJ_invertible_on_Coulomb_gauge_slice}. For the representative $(A,\varphi,\psi)$ of $[A,\varphi,\psi]$ and $(a,\sigma,\tau) \in \Ker d_{A,\varphi,\psi}^{0,*}$, we compute:
\begin{align*}
  \left(\iota_{\mathbf{\hat X}}\mathbf{\hat\bomega}\right)_{A,\varphi,\psi}(a,\sigma,\tau)
  &= \mathbf{\hat\bomega}_{A,\varphi,\psi}\left(\mathbf{\hat X}_{A,\varphi,\psi},(a,\sigma,\tau)\right)
  \\
  &= \mathbf{\hat\bomega}_{A,\varphi,\psi}\left(\Pi_{A,\varphi,\psi}\bJ(0,\varphi,\psi),(a,\sigma,\tau)\right)
    \quad\text{(by \eqref{eq:S1GeneratorExpression_quotient_space})}
  \\
  &= \mathbf{\hat g}_{A,\varphi,\psi}\left(\mathbf{\hat J}_{A,\varphi,\psi}\Pi_{A,\varphi,\psi}\bJ(0,\varphi,\psi),
    (a,\sigma,\tau)\right)
    \quad\text{(by \eqref{eq:Fundamental_two_form_quotient_space_spinu_pairs})}
  \\
  &= \mathbf{\hat g}_{A,\varphi,\psi}\left(\Pi_{A,\varphi,\psi}\bJ\Pi_{A,\varphi,\psi}\bJ(0,\varphi,\psi),
    (a,\sigma,\tau)\right)
    \quad\text{(by \eqref{eq:Almost_complex_structure_affine_space_spinu_pairs_Coulomb_gauge_slice})}
  \\
  &= \bg_{A,\varphi,\psi}\left(\left(\Pi_{A,\varphi,\psi}\bJ\right)^2(0,\varphi,\psi),
    (a,\sigma,\tau)\right)
    \quad\text{(by \eqref{eq:L2_metric_affine_space_spinu_pairs} and \eqref{eq:L2_metric_quotient_space_spinu_pairs})}.
\end{align*}
Moreover, the one-form $\iota_{\mathbf{\hat X}}\mathbf{\hat\bomega}$ is invariant under the action of $W^{2,p}(\SU(E))$ in the following sense: For any $u \in W^{2,p}(\SU(E))$ and $(a,\sigma,\tau) \in \Ker d_{A,\varphi,\psi}^{0,*}$,
\[
  (\iota_{\mathbf{\hat X}}\mathbf{\hat\bomega})_{u(A,\varphi,\psi)}u(a,\sigma,\tau)
  =
  (\iota_{\mathbf{\hat X}}\mathbf{\hat\bomega})_{A,\varphi,\psi}(a,\sigma,\tau).
\]    
Observe that \emph{if we assume} that $\Pi_{A,\varphi,\psi}\bJ = \bJ\Pi_{A,\varphi,\psi}$, then
\begin{align*}
  \left(\iota_{\mathbf{\hat X}}\mathbf{\hat\bomega}\right)_{A,\varphi,\psi}(a,\sigma,\tau)
  &=
    -\bg_{A,\varphi,\psi}\left((0,\varphi,\psi), (a,\sigma,\tau)\right)
  \\
  &= -(df)_{A,\varphi,\psi}(a,\sigma,\tau).
\end{align*}
To illustrate the application of a more general condition guaranteeing that the identity \eqref{eq:Moment_map_quotient_space_unitary_triples} holds, we consider the

\begin{cond}
\label{cond:Pi_Avarphipsi_bJ_invertible_on_Coulomb_gauge_slice}
For a point $[A,\varphi,\psi] \in \sC^{0,*}(E,H,J,A_d)$, the operator
\[
  \Pi_{A,\varphi,\psi}\bJ: L^2\left(T^*X\otimes\su(E)\oplus E\oplus\Lambda^{0,2}(E)\right)
  \to
  \Ker d_{A,\varphi,\psi}^{0,*}\cap L^2\left(T^*X\otimes\su(E)\oplus E\oplus\Lambda^{0,2}(E)\right)
\]
has a left inverse
\[
  L_{A,\varphi,\psi} \in \End\left(L^2\left(T^*X\otimes\su(E)\oplus E\oplus\Lambda^{0,2}(E)\right)\right)
\]
in the sense that
\[
  L_{A,\varphi,\psi}\circ \Pi_{A,\varphi,\psi}\bJ
  =
  \Pi_{A,\varphi,\psi} \quad\text{on } L^2\left(T^*X\otimes\su(E)\oplus E\oplus\Lambda^{0,2}(E)\right). 
\]
In particular, the operator $\Pi_{A,\varphi,\psi}\bJ$ is invertible on the Coulomb gauge slice $\Ker d_{A,\varphi,\psi}^{0,*}\cap L^2(T^*X\otimes\su(E)\oplus E\oplus\Lambda^{0,2}(E))$ through $(A,\varphi,\psi)$.  
\end{cond}

We emphasize that we do \emph{not} rely on Condition \ref{cond:Pi_Avarphipsi_bJ_invertible_on_Coulomb_gauge_slice} in order to prove any results aside from those in the example application that we now discuss below. However, Condition \ref{cond:Pi_Avarphipsi_bJ_invertible_on_Coulomb_gauge_slice} serves as a useful guide. The calculation in the following example shows that Condition \ref{cond:Pi_Avarphipsi_bJ_invertible_on_Coulomb_gauge_slice} may fail at points represented by split unitary triples, that is, points $[A,\varphi,\psi] \in \sC^0(E,H,J,A_d) \less \sC^{0,*}(E,H,J,A_d)$.

\begin{exmp}[Failure of Condition \ref{cond:Pi_Avarphipsi_bJ_invertible_on_Coulomb_gauge_slice} at split triples]
\label{exmp:Failure_of_AC_Structure_Inj_At_Split_Triples}
Consider a point in $\sC^0(E,H,J,A_d)$ represented by a smooth unitary triple $(A,\varphi,\psi)$. We begin by computing
\begin{align*}
d_{A,\varphi,\psi}^{0,*}(0,\varphi,\psi)
&=
d_A0 -\sR_{\varphi}^*\varphi-\sR_\psi^*\psi
\\
&=
-\pi_{\su(E)}(\varphi\otimes\varphi^*)_0 
-\pi_{\su(E)}\star(\psi\wedge\psi^*)_0
\quad\text{(by \eqref{eq:MultOperator_Projection_Relation}, \eqref{eq:R_varphi_star_sigma_is_sigma_tensor_varphi_star_tracefree}, and \eqref{eq:R_psi_star_tau_is_tau_tensor_psi_star_tracefree})}
\\
&=0,
\end{align*}
where the final equality follows from the observation that $(\varphi\otimes\varphi^*)_0$ and $(\star\psi\wedge\psi^*)_0$ are symmetric elements of $\Omega^0(\fsl(E))$ and $\Omega^{0,2}(\fsl(E))$, respectively, so that their projections onto $\Omega^0(\su(E))$ and $\Omega^{0,2}(\su(E))$ are identically zero.
Thus, for any unitary triple $(A,\varphi,\psi)$,
\[
(0,\varphi,\psi)\in \Ker d_{A,\varphi,\psi}^{0,*}.
\]
Moreover, $0 \not\equiv (0,\varphi,\psi) \in \Ker d_{A,\varphi,\psi}^{0,*}$ by assumption that $[A,\varphi,\psi] \in \sC^0(E,H,J,A_d)$.

Suppose now that $(A,\varphi,\psi)$ is split in the sense of Definition \ref{defn:Split_trivial_central-stabilizer_spinor_pair}
with respect to a decomposition $E=L_1\oplus L_2$ as an orthogonal direct sum of Hermitian line bundles.  In particular, $\varphi\in\Omega^0(L_1)$ and $\psi\in\Omega^{0,2}(L_1)$.  If we define
\[
\xi:=-i\,\id_{L_1}\oplus i\, \id_{L_2}\in\Omega^0(\su(E)),
\]
we see that 
\[
\bJ(0,\varphi,\psi)
=
(0,i\varphi,i\psi)
=
d_{A,\varphi,\psi}^0\xi
\in
\Ran d_{A,\varphi,\psi}^0.
\]
Because $\Ran d_{A,\varphi,\psi}^0$ is the orthogonal complement of $\Ker d_{A,\varphi,\psi}^{0,*}$, we obtain
\[
\Pi_{A,\varphi,\psi}\bJ(0,\varphi,\psi)
=
0.
\]
Thus, $\Pi_{A,\varphi,\psi}\bJ$ is not injective on $\Ker d_{A,\varphi,\psi}^{*,0}$ when $(A,\varphi,\psi)$ is a split triple with $(\varphi,\psi)\not\equiv (0,0)$.
\qed\end{exmp}

Given Condition \ref{cond:Pi_Avarphipsi_bJ_invertible_on_Coulomb_gauge_slice}, the operator $(\Pi_{A,\varphi,\psi}\bJ)^2$ is also invertible on $\Ker d_{A,\varphi,\psi}^{0,*}$. Note that $\Pi_{A,\varphi,\psi}\bJ$ is skew-adjoint on $\Ker d_{A,\varphi,\psi}^{0,*}$ since
\begin{align*}
  &\bg\left(\Pi_{A,\varphi,\psi}\bJ(a_1,\sigma_1,\tau_1), (a_2,\sigma_2,\tau_2)\right)
  \\
  &=
    \bg\left(\bJ(a_1,\sigma_1,\tau_1), (a_2,\sigma_2,\tau_2)\right)
    \quad\text{(by $\Pi_{A,\varphi,\psi}^* = \Pi_{A,\varphi,\psi}$ and $(a_2,\sigma_2,\tau_2) \in \Ran\Pi_{A,\varphi,\psi}$)}
  \\
  &=
    -\bg\left((a_1,\sigma_1,\tau_1), \bJ(a_2,\sigma_2,\tau_2)\right)
      \quad\text{(by \eqref{eq:Almost_complex_structure_affine_space_spinu_pairs_skew-symmetric})}
  \\
  &=
    -\bg\left((a_1,\sigma_1,\tau_1), \Pi_{A,\varphi,\psi}\bJ(a_2,\sigma_2,\tau_2)\right)
    \quad\text{(by $\Pi_{A,\varphi,\psi}^* = \Pi_{A,\varphi,\psi}$ and $(a_1,\sigma_1,\tau_1) \in \Ran\Pi_{A,\varphi,\psi}$)},
  \\
  &\qquad\text{for all }
  (a_k,\sigma_k,\tau_k)
  \in \Ker d_{A,\varphi,\psi}^{0,*}\cap L^2\left(T^*X\otimes\su(E)\oplus E\oplus\Lambda^{0,2}(E)\right)
  \text{ and } k=1,2.
\end{align*}
Therefore,
\[
  -(\Pi_{A,\varphi,\psi}\bJ)^2 = (\Pi_{A,\varphi,\psi}\bJ)^*(\Pi_{A,\varphi,\psi}\bJ)
  > 0 \quad\text{on }
  \Ker d_{A,\varphi,\psi}^{0,*}\cap L^2\left(T^*X\otimes\su(E)\oplus E\oplus\Lambda^{0,2}(E)\right).
\]
We define a left inverse for $(\Pi_{A,\varphi,\psi}\bJ)^2$ such that
\begin{align}
  \label{eq:Left_inverse_Pi_AvarphipsibJ_squared}
  \bG_{A,\varphi,\psi}(\Pi_{A,\varphi,\psi}\bJ)^2
  &=
  \begin{cases}
    \id &\text{on } \Ker d_{A,\varphi,\psi}^{0,*},
    \\
    0 &\text{on } \Ran d_{A,\varphi,\psi}^0.
  \end{cases}
  \\
  \notag
  &= \Pi_{A,\varphi,\psi} \quad\text{on } L^2\left(T^*X\otimes\su(E)\oplus E\oplus\Lambda^{0,2}(E)\right).
\end{align}
Noting that $-\bG_{A,\varphi,\psi} > 0$ on $\Ker d_{A,\varphi,\psi}^{0,*}$, we now define a Riemannian metric on $T_{A,\varphi,\psi}\sC^0(E,H,J,A_d)$ by
%PF7-21-2025 Standardize around using L^2 instead of W^{1,p} when L^2 is all we need
\begin{multline}
  \label{eq:Modified_L2_metric_quotient_space_spinu_pairs}
  \mathbf{\tilde g}_{A,\varphi,\psi}\left((a_1,\sigma_1,\tau_1), (a_2,\sigma_2,\tau_2)\right)
  :=
  -\bg\left(\bG_{A,\varphi,\psi}(a_1,\sigma_1,\tau_1),
    (a_2,\sigma_2,\tau_2)\right),
  \\
  \text{for all }
  (a_k,\sigma_k,\tau_k)
  \in \Ker d_{A,\varphi,\psi}^{0,*}\cap L^2\left(T^*X\otimes\su(E)\oplus E\oplus\Lambda^{0,2}(E)\right)
  \text{ and } k=1,2,
\end{multline}
and define a non-degenerate two-form on $T_{A,\varphi,\psi}\sC^0(E,H,J,A_d)$ by
\begin{multline}
  \label{eq:Modified_fundamental_two_form_quotient_space_spinu_pairs}
  \mathbf{\tilde\bomega}_{A,\varphi,\psi}\left((a_1,\sigma_1,\tau_1), (a_2,\sigma_2,\tau_2)\right)
  :=
  \mathbf{\tilde g}_{A,\varphi,\psi}\left(\Pi_{A,\varphi,\psi}\bJ(a_1,\sigma_1,\tau_1), (a_2,\sigma_2,\tau_2)\right),
  \\
  \text{for all }
  (a_k,\sigma_k,\tau_k)
  \in \Ker d_{A,\varphi,\psi}^{0,*}\cap L^2\left(T^*X\otimes\su(E)\oplus E\oplus\Lambda^{0,2}(E)\right)
  \text{ and } k=1,2.
\end{multline}
Noting that $(\Pi_{A,\varphi,\psi}\bJ)\bG_{A,\varphi,\psi} = \bG_{A,\varphi,\psi}(\Pi_{A,\varphi,\psi}\bJ)$ on $\Ker d_{A,\varphi,\psi}^{0,*}$, we obtain
\begin{align*}
  &\mathbf{\tilde\bomega}_{A,\varphi,\psi}\left((a_1,\sigma_1,\tau_1), (a_2,\sigma_2,\tau_2)\right)
  \\
  &\quad =
    \mathbf{\tilde g}_{A,\varphi,\psi}\left(\Pi_{A,\varphi,\psi}\bJ(a_1,\sigma_1,\tau_1), (a_2,\sigma_2,\tau_2)\right)
    \quad\text{(by \eqref{eq:Modified_fundamental_two_form_quotient_space_spinu_pairs})}
  \\
  &\quad =
    -\bg\left(\bG_{A,\varphi,\psi}\Pi_{A,\varphi,\psi}\bJ(a_1,\sigma_1,\tau_1), (a_2,\sigma_2,\tau_2)\right)
    \quad\text{(by \eqref{eq:Modified_L2_metric_quotient_space_spinu_pairs})}
   \\
  &\quad =
    -\bg\left((a_2,\sigma_2,\tau_2),\bG_{A,\varphi,\psi}\Pi_{A,\varphi,\psi}\bJ(a_1,\sigma_1,\tau_1)\right)
    \quad\text{(by symmetry of $\bg$)}  
  \\
  &\quad =
    -\bg\left(\bG_{A,\varphi,\psi}(a_2,\sigma_2,\tau_2),\Pi_{A,\varphi,\psi}\bJ(a_1,\sigma_1,\tau_1)\right)
    \quad\text{(by self-adjointness of $\bG_{A,\varphi,\psi}$)}  
  \\
  &\quad =
    \bg\left(\Pi_{A,\varphi,\psi}\bJ\bG_{A,\varphi,\psi}(a_2,\sigma_2,\tau_2),(a_1,\sigma_1,\tau_1)\right)
    \quad\text{(by $(\Pi_{A,\varphi,\psi}\bJ)^* = -\Pi_{A,\varphi,\psi}\bJ$ on $\Ker d_{A,\varphi,\psi}^{0,*}$)}
  \\
  &\quad =
  \bg\left(\bG_{A,\varphi,\psi}\Pi_{A,\varphi,\psi}\bJ(a_2,\sigma_2,\tau_2),(a_1,\sigma_1,\tau_1)\right)
  \\
  &\quad =
    -\mathbf{\tilde\bomega}_{A,\varphi,\psi}\left((a_2,\sigma_2,\tau_2), (a_1,\sigma_1,\tau_1)\right)
    \quad\text{(by \eqref{eq:Modified_L2_metric_quotient_space_spinu_pairs} and \eqref{eq:Modified_fundamental_two_form_quotient_space_spinu_pairs})},
  \\
  &\quad\qquad \text{for all }
  (a_k,\sigma_k,\tau_k)
  \in \Ker d_{A,\varphi,\psi}^{0,*}\cap L^2\left(T^*X\otimes\su(E)\oplus E\oplus\Lambda^{0,2}(E)\right)
  \text{ and } k=1,2,
\end{align*}
and so $\mathbf{\tilde\bomega}_{A,\varphi,\psi}$ is a skew-symmetric bilinear form, as claimed. (Our construction of $\mathbf{\tilde g}_{A,\varphi,\psi}$ is similar to that described by Cannas da Silva around the displayed identity in \cite[Section 12.2, paragraph following proof of Proposition 12.3, p. 85]{Cannas_da_Silva_lectures_on_symplectic_geometry}.) Observe that
\begin{align*}
  \left(\iota_{\mathbf{\hat X}}\mathbf{\tilde\bomega}\right)_{A,\varphi,\psi}(a,\sigma,\tau)
  &= \mathbf{\tilde\bomega}_{A,\varphi,\psi}\left(\mathbf{\hat X}_{A,\varphi,\psi},(a,\sigma,\tau)\right)
  \\
  &= \mathbf{\tilde g}_{A,\varphi,\psi}\left(\left(\Pi_{A,\varphi,\psi}\bJ\right)^2(0,\varphi,\psi),
    (a,\sigma,\tau)\right)
         \quad\text{(by \eqref{eq:Modified_fundamental_two_form_quotient_space_spinu_pairs} and \eqref{eq:S1GeneratorExpression_quotient_space})}
  \\
  &= -\bg\left(\bG_{A,\varphi,\psi}\left(\Pi_{A,\varphi,\psi}\bJ\right)^2(0,\varphi,\psi), (a,\sigma,\tau)\right)
    \quad\text{(by \eqref{eq:Modified_L2_metric_quotient_space_spinu_pairs})}
  \\
  &= -\bg\left(\Pi_{A,\varphi,\psi}(0,\varphi,\psi), (a,\sigma,\tau)\right)
    \quad\text{(by \eqref{eq:Left_inverse_Pi_AvarphipsibJ_squared})}
  \\
  &= -\bg\left((0,\varphi,\psi), (a,\sigma,\tau)\right)
    \quad\text{(by $\Pi_{A,\varphi,\psi}^* = \Pi_{A,\varphi,\psi}$ and $(a,\sigma,\tau) \in \Ker d_{A,\varphi,\psi}^{0,*}$)},
  \\
  &\qquad\text{for all }
  (a,\sigma,\tau)
  \in \Ker d_{A,\varphi,\psi}^{0,*}\cap L^2\left(T^*X\otimes\su(E)\oplus E\oplus\Lambda^{0,2}(E)\right).
\end{align*}
Hence, under the Condition \ref{cond:Pi_Avarphipsi_bJ_invertible_on_Coulomb_gauge_slice}, we obtain
\[
  -(df)_{A,\varphi,\psi} = \left(\iota_{\mathbf{\hat X}}\mathbf{\tilde\bomega}\right)_{A,\varphi,\psi},
  \quad\text{for } (A,\varphi,\psi) \in \sA(E,H,A_d)\times W^{1,p}\left(E\oplus\Lambda^{0,2}(E)\right).
\]
Both sides of this identity are invariant (in the sense explained above) under the action of $u\in W^{2,p}(\SU(E))$ and so the identity descends to the quotient space $\sC^{0,*}(E,H,J,A_d)$. This verifies the Hamiltonian identity \eqref{eq:Moment_map_quotient_space_unitary_triples} under the Condition \ref{cond:Pi_Avarphipsi_bJ_invertible_on_Coulomb_gauge_slice}. Recall, however, that $S^1$ acts freely on $\sC^{0,*}(E,H,J,A_d)$ and that $f$ has no critical points in $\sC^{0,*}(E,H,J,A_d)$ by Lemma \ref{lem:Critical_points_f_on_quotient_space_are_fixed_points_S1_action}.

\subsection{Hitchin's function and circle-invariant non-degenerate two-form on the moduli space of non-Abelian monopoles}
\label{subsec:Hamiltonian_function_circle-invariant_non-degenerate_2-form_virtual_moduli_space}
See Remark \ref{rmk:Functorial_properties_Hamiltonian_functions} for a discussion of the functorial property of the moment map identity \eqref{eq:MomentMap} under pullback by $S^1$-equivariant smooth maps of smooth manifolds with $S^1$ actions. The equality of one-forms in \eqref{eq:df_Avarphipsi_asigmatau_Coulomb_gauge_slice} pulls back to an equality on the harmonic subspace $\bH_{A,\varphi,\psi,r}^1 \subset \Ker d_{A,\varphi,\psi}^{0,*}$ as in \eqref{eq:bH1_Avarphipsi_r} to give
\begin{equation}
  \label{eq:df_Avarphipsi_asigmatau_moduli_space_non-Abelian_monopoles}
  (df)_{A,\varphi,\psi}(a,\sigma,\tau)
  =
  \bg\left(\pi_{A,\varphi,\psi,r}(0,\varphi,\psi),(a,\sigma,\tau)\right),
  \quad
  \text{for all } (a,\sigma,\tau)
  \in \bH_{A,\varphi,\psi,r}^1,
\end{equation}
where
\begin{equation}
  \label{eq:L2_orthogonal_projection_Zariski_tangent_space_non-Abelian_moduli_space}
  \pi_{A,\varphi,\psi,r}: L^2\left(T^*X\otimes\su(E) \oplus E\oplus\Lambda^{0,2}(E)\right)
  \to \bH_{A,\varphi,\psi,r}^1
\end{equation}
denote the $L^2$ orthogonal projection onto the Zariski tangent space to $\sM^0(E,g,J,\omega,r)$ at $(A,\varphi,\psi)$. We now describe the consequences of a technical

% TL6-22-2025: Isn't this equation (6.1.4) (that is \eqref{eq:Donaldson_1996jdg_proposition_3_pi_Jv_neq_0_for_v_in_KerT})?  That is, the assumption is the crucial step in the proof of Prop 8.  And if we assume it, don't we get the symplectic structure on $\bH_{A,\varphi,\psi,r}^1$ just by following the rest of the proof of Prop 8?
% PF6-23-2025 ?? See my comments around the previous assumption
%PF7-21-2025 Title for below and previous condition?
\begin{cond}
\label{cond:pi_Avarphipsi_r_bJ_invertible_on_Zariski_tangent_space}
For a point $[A,\varphi,\psi] \in \sM^0(E,g,J,\omega,r)$,
%PF7-21-2025 I strengthened the condition to make it easier to apply
% the operator $\pi_{A,\varphi,\psi,r}\bJ$ is invertible on the Zariski tangent space $\bH_{A,\varphi,\psi,r}^1$ as in \eqref{eq:bH1_Avarphipsi_r} to the moduli space $\sM^0(E,g,J,\omega,r)$ as in \eqref{eq:Moduli_space_non-Abelian_monopoles_almost_Hermitian_Taubes_regularized_non-zero-section} at the representative $(A,\varphi,\psi)$.
the operator
\[
  \pi_{A,\varphi,\psi,r}\bJ: L^2\left(T^*X\otimes\su(E)\oplus E\oplus\Lambda^{0,2}(E)\right)
  \to
  \bH_{A,\varphi,\psi,r}^1
\]
has a left inverse
\[
  L_{A,\varphi,\psi,r} \in \End\left(L^2\left(T^*X\otimes\su(E)\oplus E\oplus\Lambda^{0,2}(E)\right)\right)
\]
in the sense that
\[
  L_{A,\varphi,\psi,r}\circ \pi_{A,\varphi,\psi,r}\bJ
  =
  \pi_{A,\varphi,\psi,r} \quad\text{on } L^2\left(T^*X\otimes\su(E)\oplus E\oplus\Lambda^{0,2}(E)\right). 
\]
In particular, the operator $\pi_{A,\varphi,\psi,r}\bJ$ is invertible on the Zariski tangent space $\bH_{A,\varphi,\psi,r}^1$ as in \eqref{eq:bH1_Avarphipsi_r} to the moduli space $\sM^0(E,g,J,\omega,r)$ as in \eqref{eq:Moduli_space_non-Abelian_monopoles_almost_Hermitian_Taubes_regularized_non-zero-section} at the representative $(A,\varphi,\psi)$.
\end{cond}

Unlike Condition \ref{cond:Pi_Avarphipsi_bJ_invertible_on_Coulomb_gauge_slice}, we note that Condition \ref{cond:pi_Avarphipsi_r_bJ_invertible_on_Zariski_tangent_space} is not obstructed by Example \ref{exmp:Failure_of_AC_Structure_Inj_At_Split_Triples} since, although $0 \not\equiv (0,\varphi,\psi) \in \Ker d_{A,\varphi,\psi}^{0,*}$ for points $[A,\varphi,\psi] \in \sC^0(E,H,J,A_d)$, we do not necessarily have $(0,\varphi,\psi) \in \bH_{A,\varphi,\psi,r}^1$ for points $[A,\varphi,\psi] \in \sM^0(E,g,J,\omega,r)$.

As in Section \ref{subsec:Hamiltonian_function_circle-invariant_non-degenerate_2-form_quotient_space}, we emphasize that we do \emph{not} rely on Condition \ref{cond:pi_Avarphipsi_r_bJ_invertible_on_Zariski_tangent_space} in order to prove any results aside from those in the example application discussed below. Like in Section \ref{subsec:Hamiltonian_function_circle-invariant_non-degenerate_2-form_quotient_space}, the Condition \ref{cond:pi_Avarphipsi_r_bJ_invertible_on_Zariski_tangent_space} and the example application that we now discuss serve as a useful guide.

Given Condition \ref{cond:pi_Avarphipsi_r_bJ_invertible_on_Zariski_tangent_space}, the operator $(\pi_{A,\varphi,\psi,r}\bJ)^2$ is invertible on $\bH_{A,\varphi,\psi,r}^1$. The operator $\pi_{A,\varphi,\psi,r}\bJ$ is also skew-adjoint on $\bH_{A,\varphi,\psi,r}^1$ since
\begin{align*}
  &\bg\left(\pi_{A,\varphi,\psi,r}\bJ(a_1,\sigma_1,\tau_1), (a_2,\sigma_2,\tau_2)\right)
  \\
  &=
    \bg\left(\bJ(a_1,\sigma_1,\tau_1), \pi_{A,\varphi,\psi,r}(a_2,\sigma_2,\tau_2)\right) 
    \quad\text{(by $\pi_{A,\varphi,\psi,r}^* = \pi_{A,\varphi,\psi,r}$)}
  \\
  &=
    \bg\left(\bJ(a_1,\sigma_1,\tau_1), (a_2,\sigma_2,\tau_2)\right) 
    \quad\text{(by $(a_2,\sigma_2,\tau_2) \in \bH_{A,\varphi,\psi,r}^1$)}
  \\
  &=
    -\bg\left((a_1,\sigma_1,\tau_1), \bJ(a_2,\sigma_2,\tau_2)\right)
      \quad\text{(by \eqref{eq:Almost_complex_structure_affine_space_spinu_pairs_skew-symmetric})}
  \\
  &=
    -\bg\left((a_1,\sigma_1,\tau_1), \pi_{A,\varphi,\psi,r}\bJ(a_2,\sigma_2,\tau_2)\right)
    \quad\text{(by $(a_1,\sigma_1,\tau_1) \in \bH_{A,\varphi,\psi,r}^1$
    and $\pi_{A,\varphi,\psi,r}^* = \pi_{A,\varphi,\psi,r}$)},
  \\
  &\qquad\text{for all }
  (a_k,\sigma_k,\tau_k) \in \bH_{A,\varphi,\psi,r}^1 \text{ and } k=1,2.
\end{align*}
Therefore, 
\[
  -(\pi_{A,\varphi,\psi,r}\bJ)^2 = (\pi_{A,\varphi,\psi,r}\bJ)^*(\pi_{A,\varphi,\psi,r}\bJ)
  > 0 \quad\text{on } \bH_{A,\varphi,\psi,r}^1,
\]
and $-(\pi_{A,\varphi,\psi,r}\bJ)^2$ is a positive self-adjoint operator on $\bH_{A,\varphi,\psi,r}^1$. We define a left inverse for
\[
  (\pi_{A,\varphi,\psi,r}\bJ)^2:
  W^{1,p}\left(T^*X\otimes\su(E)\oplus E\oplus\Lambda^{0,2}(E)\right)
  \to
  \bH_{A,\varphi,\psi,r}^1
\]
such that
\begin{align}
  \label{eq:Left_inverse_Pi_Avarphipsi_r_balpha_mu_bJ_squared}
  \bG_{A,\varphi,\psi,r}(\pi_{A,\varphi,\psi,r}\bJ)^2
  &=
  \begin{cases}
    \id &\text{on } \bH_{A,\varphi,\psi,r}^1,
    \\
    0 &\text{on } \left(\bH_{A,\varphi,\psi,r}^1\right)^\perp
  \end{cases}
  \\
  \notag
  &= \pi_{A,\varphi,\psi,r} \quad\text{on } W^{1,p}\left(T^*X\otimes\su(E)\oplus E\oplus\Lambda^{0,2}(E)\right).
\end{align}
Noting that $-\bG_{A,\varphi,\psi,r}$ is a positive self-adjoint operator on $\bH_{A,\varphi,\psi,r}^1$, we now define a Riemannian metric on $\bH_{A,\varphi,\psi,r}^1$ by
\begin{multline}
  \label{eq:Modified_L2_metric_bH1_Avarphipsi_r}
  \mathbf{\tilde g}_{A,\varphi,\psi}\left((a_1,\sigma_1,\tau_1), (a_2,\sigma_2,\tau_2)\right)
  :=
  \bg\left(-\bG_{A,\varphi,\psi,r}(a_1,\sigma_1,\tau_1),
    (a_2,\sigma_2,\tau_2)\right),
  \\
  \text{for all } (a_k,\sigma_k,\tau_k)
  \in \bH_{A,\varphi,\psi,r}^1 \text{ and } k=1,2,
\end{multline}
and define a non-degenerate two-form on $\bH_{A,\varphi,\psi,r}^1$ by
\begin{multline}
  \label{eq:Modified_fundamental_two_form_bH1_Avarphipsi_r}
  \mathbf{\tilde\bomega}_{A,\varphi,\psi}\left((a_1,\sigma_1,\tau_1), (a_2,\sigma_2,\tau_2)\right)
  :=
  \mathbf{\tilde g}_{A,\varphi,\psi}\left(\pi_{A,\varphi,\psi,r}\bJ(a_1,\sigma_1,\tau_1), (a_2,\sigma_2,\tau_2)\right),
  \\
  \text{for all } (a_k,\sigma_k,\tau_k)
  \in \bH_{A,\varphi,\psi,r}^1 \text{ and } k=1,2.
\end{multline}
(Once again our construction of $\mathbf{\tilde g}_{A,\varphi,\psi}$ is identical to that described by Cannas da Silva around the displayed identity in \cite[Section 12.2, paragraph following proof of Proposition 12.3, p. 85]{Cannas_da_Silva_lectures_on_symplectic_geometry}.) If we denote
\begin{equation}
  \label{eq:S1GeneratorExpression_mu_virtual_moduli_space}
  \bY_{A,\varphi,\psi,r} := \pi_{A,\varphi,\psi,r}\bJ(0,\varphi,\psi),
\end{equation}
we obtain a $W^{2,p}(\SU(E))$-equivariant smooth map,
\begin{multline*}
  \bY:\sA(E,H,A_d)\times W^{1,p}(E\oplus\Lambda^{0,2}(E)) \ni (A,\varphi,\psi)
  \\
  \mapsto \pi_{A,\varphi,\psi,r}\bJ(0,\varphi,\psi)
  \in L^p\left(T^*X\otimes\su(E)\oplus E\oplus\Lambda^{0,2}(E)\right),
\end{multline*}
and hence a rough vector field on the moduli space $\sM^0(E,g,J,\omega,r)$ as in \eqref{eq:Moduli_space_non-Abelian_monopoles_almost_Hermitian_Taubes_regularized_non-zero-section}. We see that $\bY$ is the generator of the $S^1$ action \eqref{eq:S1ZActionOnQuotientSpace} on $\sM^0(E,g,J,\omega,r)$ as follows. If $[A,\varphi,\psi] \in \sM^0(E,g,J,\omega,r)$, then
\[
  e^{i\theta}\cdot[A,\varphi,\psi] = [A,e^{i\theta}\varphi,e^{i\theta}\psi] \in \sM^0(E,g,J,\omega,r),
  \quad\text{for } \theta \in (-\pi,\pi),
\]
defines a continuous path in $\sM^0(E,g,J,\omega,r)$ that is the image of the following smooth curve of solutions in $\sA(E,H,A_d)\times W^{1,p}(E\oplus\Lambda^{0,2}(E))$ to the non-Abelian monopole equations \eqref{eq:SO(3)_monopole_equations_almost_Hermitian_perturbed_intro_regular} with a regularized Taubes perturbation,
\[
  e^{i\theta}\cdot(A,\varphi,\psi) = (A,e^{i\theta}\varphi,e^{i\theta}\psi)
  \in \sA(E,H,A_d) \times  W^{1,p}(E\oplus\Lambda^{0,2}(E)),
  \quad\text{for } \theta \in (-\eps,\eps).
\]
%TL6-21-2025: "Because the moduli space $\sM^0(E,g,J,\omega,r)$ is closed under the $S^1$ action," Or perhaps just put in "Therefore" to indicate that this statement follows from the previous statement about the $S^1$ orbit being a continuous path in the moduli space (rather than something the reader is supposed to figure out theirself)
The tangent vector to this curve at $\theta=0$ belongs to the kernel of the linearization \eqref{eq:SO3Monopoled1TaubesPerturbation} of the system \eqref{eq:SO(3)_monopole_equations_almost_Hermitian_perturbed_intro_regular} at $(A,\varphi,\psi)$,
\[
  \left.\frac{d}{d\theta}(A,e^{i\theta}\varphi,e^{i\theta}\psi)\right|_{\theta=0}
  =
  (0,i\varphi,i\psi)
  =
  \bJ(0,\varphi,\psi) \in \Ker d_{A,\varphi,\psi,r}^1.
\]
The projection of the preceding tangent vector onto the slice $\Ker d_{A,\varphi,\psi}^{0,*}$ through $(A,\varphi,\psi)$ yields
\[
  \bY_{A,\varphi,\psi,r}
  =
  \pi_{A,\varphi,\psi,r}\bJ(0,\varphi,\psi)
  =
  \Pi_{A,\varphi,\psi}\bJ(0,\varphi,\psi)
  \in
  \Ker d_{A,\varphi,\psi,r}^1 \cap \Ker d_{A,\varphi,\psi}^{0,*}
  =
  \bH_{A,\varphi,\psi,r}^1.
\]
If $[A,\varphi,\psi] \in \sM_\reg^0(E,g,J,\omega,r)$, the open submanifold of regular points in $\sM^0(E,g,J,\omega,r)$, then $e^{i\theta}\cdot[A,\varphi,\psi]$, for $\theta \in (-\pi,\pi)$, is a smooth curve through $[A,\varphi,\psi]$ that is the image of the smooth curve $e^{i\theta}\cdot(A,\varphi,\psi)$ in the affine space. (Again, these curves are smooth since the $S^1$ actions \eqref{eq:S1ZActionOnQuotientSpace} and \eqref{eq:S1ZAction} are smooth.) The curve $e^{i\theta}\cdot[A,\varphi,\psi]$ defines the tangent vector $\bY_{A,\varphi,\psi,r}$ to $\sM_\reg^0(E,g,J,\omega,r)$ at $(A,\varphi,\psi)$, so $\bY$ generates the $S^1$ action on $\sM_\reg^0(E,g,J,\omega,r)$. The same remarks apply to an $S^1$ invariant open neighborhood of a singular point $[A,\varphi,\psi]$ after replacing $\sM_\reg^0(E,g,J,\omega,r)$ by the $S^1$ invariant ambient smooth local virtual moduli space $\sM^\vir(E,g,J,\omega,r)$ for $[A,\varphi,\psi]$, defined by the Kuranishi Method as in Section \ref{subsec:Local_Kuranishi_model_nbhd_point_zero_locus_unperturbed_non-Abelian_monopole_equations}.

The vector field $\bY$ is $\iota$-related to $\mathbf{\hat X}$ in the sense of Remark \ref{rmk:Functorial_properties_Hamiltonian_functions}, where $\iota$ denotes the inclusion $\sM_\reg^0(E,g,J,\omega,r) \subset \sC^0(E,H,J,A_d)$ or $\sM^\vir(E,g,J,\omega,r) \subset \sC^0(E,H,J,A_d)$. Similarly, the Riemannian metric $\mathbf{\tilde g}$ in \eqref{eq:Modified_L2_metric_bH1_Avarphipsi_r} and two-form $\mathbf{\tilde\bomega}$ in \eqref{eq:Modified_fundamental_two_form_bH1_Avarphipsi_r} on $\sM^0(E,g,J,\omega,r)$ are obtained by pullback via $\iota$ of the Riemannian metric $\mathbf{\tilde g}$ in \eqref{eq:Modified_L2_metric_quotient_space_spinu_pairs} and two-form $\mathbf{\tilde\bomega}$ in \eqref{eq:Modified_fundamental_two_form_quotient_space_spinu_pairs} on $\sC^0(E,H,J,A_d)$. We observe that
\begin{align*}
  \left(\iota_\bY\mathbf{\tilde\bomega}\right)_{A,\varphi,\psi}(a,\sigma,\tau)
  &= \mathbf{\tilde\bomega}_{A,\varphi,\psi}\left(\bY_{A,\varphi,\psi,r},(a,\sigma,\tau)\right)
  \\
  &= \mathbf{\tilde g}_{A,\varphi,\psi}\left(\pi_{A,\varphi,\psi,r}\bJ
    \bY_{A,\varphi,\psi,r}, (a,\sigma,\tau)\right)
    \quad\text{(by \eqref{eq:Modified_fundamental_two_form_bH1_Avarphipsi_r})}
  \\
  &= \mathbf{\tilde g}_{A,\varphi,\psi}\left(\left(\pi_{A,\varphi,\psi,r}\bJ\right)^2(0,\varphi,\psi),
    (a,\sigma,\tau)\right)
    \quad\text{(by \eqref{eq:S1GeneratorExpression_mu_virtual_moduli_space})}
  \\
  &= -\bg\left(\bG_{A,\varphi,\psi,r}\left(\pi_{A,\varphi,\psi,r}\bJ\right)^2(0,\varphi,\psi), (a,\sigma,\tau)\right)
    \quad\text{(by \eqref{eq:Modified_L2_metric_bH1_Avarphipsi_r})}
  \\
  &= -\bg\left(\pi_{A,\varphi,\psi,r}(0,\varphi,\psi), (a,\sigma,\tau)\right)
    \quad\text{(by \eqref{eq:Left_inverse_Pi_Avarphipsi_r_balpha_mu_bJ_squared})}
  \\
  &= -\bg\left((0,\varphi,\psi), (a,\sigma,\tau)\right),
  \quad\text{for all }
  (a,\sigma,\tau) \in \bH_{A,\varphi,\psi,r}^1.
\end{align*}
Hence, under the Condition \ref{cond:pi_Avarphipsi_r_bJ_invertible_on_Zariski_tangent_space}, we obtain by comparison with the identity \eqref{eq:df_at_Avarphipsi_direction_asigmatau} that
\[
  -(df)_{A,\varphi,\psi} = \left(\iota_\bY\mathbf{\tilde\bomega}\right)_{A,\varphi,\psi}
  \quad\text{on } \bH_{A,\varphi,\psi,r}^1,
  \quad\text{for } (A,\varphi,\psi) \in \sA(E,H,A_d)\times W^{1,p}\left(E\oplus\Lambda^{0,2}(E)\right).
\]
Both sides of this identity are invariant (in the sense explained above) under the action of $u\in W^{2,p}(\SU(E))$ and so the identity descends to the quotient space $\sC^0(E,H,J,A_d)$ and hence to the subspace $\sM^0(E,g,J,\omega,r)$. This verifies the Hamiltonian identity
\begin{equation}
  \label{eq:Moment_map_moduli_space_non-Abelian_monopoles}
  -(df)_{A,\varphi,\psi} = \left(\iota_\bY\mathbf{\tilde\bomega}\right)_{A,\varphi,\psi}
\end{equation}
under the Condition \ref{cond:pi_Avarphipsi_r_bJ_invertible_on_Zariski_tangent_space}. Suppose now that $[A,\varphi,\psi] \in \sM^0(E,g,J,\omega,r)$ is a \emph{critical point} of $f$ in the sense of Definition \ref{defn:Critical_point_Hitchin_function_moduli_space_non-Abelian_monopoles}:
\[
  (df)_{A,\varphi,\psi} = 0 \quad\text{on } \bH_{A,\varphi,\psi,r}^1.
\]
The identity \eqref{eq:Moment_map_moduli_space_non-Abelian_monopoles} yields $\mathbf{\tilde\bomega}_{A,\varphi,\psi}(\bY_{A,\varphi,\psi},\cdot) = 0$ on $\bH_{A,\varphi,\psi,r}^1$ and because the two-form $\mathbf{\tilde\bomega}_{A,\varphi,\psi}$ in \eqref{eq:Modified_fundamental_two_form_bH1_Avarphipsi_r} is non-degenerate on $\bH_{A,\varphi,\psi,r}^1$, we obtain that $\bY_{A,\varphi,\psi} = 0$. We can then apply Theorem \ref{thm:Frankel_almost_Hermitian} \eqref{item:Frankel_almost_Hermitian_FixedPointsAreZerosOfVField} to conclude that $[A,\varphi,\psi] \in \sM^0(E,g,J,\omega,r)$ is a fixed point of the $S^1$ action \eqref{eq:S1ZActionOnQuotientSpace}
of Definition \ref{defn:UnitaryZActionOnAffine} on the ambient $S^1$ invariant smooth local virtual moduli space
%PF5-21-2025 Produced by ??
$\sM^\vir(E,g,J,\omega,r)$ with tangent space $\bH_{A,\varphi,\psi,r}^1$ at $(A,\varphi,\psi)$ and hence that $[A,\varphi,\psi]$ is a fixed point of the $S^1$ action \eqref{eq:S1ZActionOnQuotientSpace} on $\sC^0(E,H,J,A_d)$.

\chapter[Non-degenerate two-forms on virtual moduli spaces of non-Abelian monopoles]{Circle-invariant non-degenerate two-forms on virtual moduli spaces of non-Abelian monopoles}
\label{chap:Construction_circle-invariant_non-degenerate_two-form_II}
Our main goal in this chapter is to complete the proofs of Theorem \ref{mainthm:AH_structure_bounded_evalue_spaces_non-Abelian_monopoles_symp_4-mflds}, Corollary \ref{maincor:Almost_Hermitian_structure_moduli_space_non-Abelian_monopoles_symplectic_4-manifolds}, and Theorem  \ref{mainthm:IdentifyCriticalPoints}. In Section \ref{sec:Unperturbed_non-Abelian_monopole_equations_almost_Kaehler_four-manifolds} we describe the coefficients of the (rescaled) linearization of the unperturbed and perturbed non-Abelian monopole equations and Coulomb gauge operator and classify them according to whether they tend to decrease, remain bounded, or increase as the Taubes perturbation parameter $r$ tends to infinity. Section \ref{sec:Approximation_orthogonal_projections_finite-dim_subspaces_Hilbert_spaces} develops results related to the problem of approximating orthogonal projections onto finite-dimensional subspaces of Hilbert spaces. In Section \ref{sec:Hitchin_function_Hamiltonian_circle_action_virtual_moduli_spaces}, we construct almost Hermitian structures on the real vector spaces $\mathbf{\tilde H}_{A,\varphi,\psi,r,\nu}^k$ for $k=1,2$ in \eqref{eq:tilde_bH_A_varphi_psi_r_nu_1} and \eqref{eq:tilde_bH_A_varphi_psi_r_nu_2} and prove that the circle action on the corresponding circle invariant smooth virtual moduli space $\sM^\vir(E,g,J,\omega,r,\nu)$ of non-Abelian monopoles is Hamiltonian, thus proving Theorem \ref{mainthm:AH_structure_bounded_evalue_spaces_non-Abelian_monopoles_symp_4-mflds} and Corollary \ref{maincor:Almost_Hermitian_structure_moduli_space_non-Abelian_monopoles_symplectic_4-manifolds}. Section \ref{sec:Equivalence_critical_points_Hamiltonian_fixed_points_S1_action_moduli_space} contains our statement and proof of Theorem \ref{thm:Critical_points_Hitchin_Hamiltonian_function_moduli_space_non-Abelian_monopoles}, which gives the equivalence between $\mathbf{\tilde H}_{A,\varphi,\psi,r,\nu}^1$-critical points of the restriction of the Hitchin function $f: \sC^0(E,H,J,A_d) \to \RR$ in \eqref{eq:Hitchin_function} to $\sM^0(E,g,J,\omega,r)$ and fixed points of the circle action on $\sM^0(E,g,J,\omega,r)$, and hence prove Theorem \ref{mainthm:IdentifyCriticalPoints}, which gives the equivalence between $\mathbf{\tilde H}_{A,\varphi,\psi,r,\nu}^1$-critical points of the Hitchin function $f:\sM^0(E,g,J,\omega,r) \to \RR$ and points represented by solutions to the Seiberg--Witten monopole equations.

Throughout this chapter, it is useful to keep in mind Donaldson's discussion in \cite[Section V, part (iv), p. 294]{DonPoly} that indicates why it is difficult to construct a \emph{symplectic} form on the moduli space $M_\kappa^w(X,g)$ of anti-self dual connections even when $(X,g,J,\omega)$ is a symplectic four-manifold.

\section[Decomposition of deformation operator into complex linear and antilinear components]{Decomposition of the deformation operator for the non-Abelian monopole equations into complex linear and antilinear components}
\label{sec:Unperturbed_non-Abelian_monopole_equations_almost_Kaehler_four-manifolds}
We assume in this section that $(X,g,J,\omega)$ is an almost K\"ahler closed four-dimensional manifold, so $d\omega=0$ but $J$ is not assumed to be integrable, and that $(A,\varphi,\psi)$ as in
\eqref{eq:A_varphi_psi_in_W1p} is a smooth triple. The latter restriction involves little loss of generality. Indeed, if $(A,\varphi,\psi)$ is a $W^{1,p}$ solution for $p \in (2,\infty)$ to the system \eqref{eq:SO(3)_monopole_equations_almost_Hermitian_intro} of unperturbed non-Abelian monopole equations or the system \eqref{eq:SO(3)_monopole_equations_almost_Hermitian_perturbed_intro_regular} of non-Abelian monopole equations with a regularized Taubes perturbation, then Theorem \ref{thm:Feehan_Leness_1998jdg_3-7_regularized_Taubes-perturbed_W1p} implies that $(A,\varphi,\psi)$ is equivalent by a $W^{2,p}(\SU(E))$ gauge transformation to a $C^\infty$ solution to \eqref{eq:SO(3)_monopole_equations_almost_Hermitian_intro} or \eqref{eq:SO(3)_monopole_equations_almost_Hermitian_perturbed_intro_regular}.

By the definition \eqref{eq:DefineHatPartial0} of $\hat\partial_{A,\varphi,\psi}^{0,*}$ and the definition \eqref{eq:DefineBarPartialWithMu} of $\bar\partial_{A,\varphi,\psi}^1$, we have the following explicit expression for the operator \eqref{eq:CplxDef1},
\[
  \cT_{\bar\partial_A,\varphi,\psi}
  = \bar\partial_{A,\varphi,\psi}^1 + \hat\partial_{A,\varphi,\psi}^{0,*}:\sF_1\to\sF_2,
\]
between the complex Fr\'echet spaces in \eqref{eq:sEkC}:
\begin{equation}
  \label{eq:H_dbar_APhi^01_explicit}
  \cT_{\bar\partial_A,\varphi,\psi}(a'',\sigma,\tau)
  =
  \begin{pmatrix}
    \bar\partial_A^*a'' - R_\varphi^*\sigma + (R_\psi^*\tau)^\dagger
    \\
    \bar\partial_Aa'' - \frac{1}{4}N_J^*(a'')^\dagger - (\tau\otimes\varphi^*+\psi\otimes\sigma^*)_0
    \\
    \bar\partial_A\sigma+\bar\partial_A^*\tau + a''\varphi + \star((a'')^\dagger\wedge \star\psi)
  \end{pmatrix}
  \in \sF_2, \quad\text{for all } (a'',\sigma,\tau) \in \sF_1.
\end{equation}
The explicit expression \eqref{eq:H_dbar_APhi^01_explicit} for the operator $\cT_{\bar\partial_A,\varphi,\psi}$ in \eqref{eq:CplxDef1} corresponds to the system \eqref{eq:SO(3)_monopole_equations_almost_Hermitian_intro} of \emph{unperturbed} non-Abelian monopole equations.

We now turn to the system \eqref{eq:SO(3)_monopole_equations_almost_Hermitian_perturbed_intro_regular} of non-Abelian monopole equations with a \emph{regularized Taubes perturbation} and its linearization. The only change to the explicit expression \eqref{eq:H_dbar_APhi^01_explicit} for $\cT_{\bar\partial_A,\varphi,\psi}$ caused by the replacement of  \eqref{eq:SO(3)_monopole_equations_almost_Hermitian_intro} by \eqref{eq:SO(3)_monopole_equations_almost_Hermitian_perturbed_intro_regular} is in the $\Omega^0(\fsl(E))$-component of $\cT_{\bar\partial_A,\varphi,\psi}(a'',\sigma,\tau)$ in \eqref{eq:H_dbar_APhi^01_explicit}. To obtain this modification of the expression \eqref{eq:H_dbar_APhi^01_explicit} for $\cT_{\bar\partial_A,\varphi,\psi}$, we use the explicit formula for the operator \eqref{eq:Perturbed_AC_Deformation_Operator},
\[
  \cT_{\bar\partial_A,\varphi,\psi,r}
  =
  \bar\partial_{A,\varphi,\psi}^1 + \hat\partial_{A,\varphi,\psi,r}^{0,*}:\sF_1\to\sF_2,
\]
that is implied by the definition \eqref{eq:DefineBarPartialWithMu} of $\bar\partial_{A,\varphi,\psi}^1$ and the definitions \eqref{eq:DefineHatPartial0} and \eqref{eq:DefineHat0*OperatorPerturbed} of $\hat\partial_{A,\varphi,\psi,r}^{0,*}$. Thus, we obtain the following analogue of \eqref{eq:H_dbar_APhi^01_explicit}:
\begin{multline}
  \label{eq:H_dbar_APhi^01_explicit_perturbed}
  \cT_{\bar\partial_A,\varphi,\psi,r}(a'',\sigma,\tau)
  =
  \begin{pmatrix}
    \bar\partial_A^*a'' - R_\varphi^*\sigma + (R_\psi^*\tau)^\dagger
    + \frac{ir}{4}p_\psi(\tau) + \frac{ir}{4}p_\psi(\tau)^\dagger
    \\
    \bar\partial_Aa'' - \frac{1}{4}N_J^*(a'')^\dagger - (\tau\otimes\varphi^*+\psi\otimes\sigma^*)_0
    \\
    \bar\partial_A\sigma+\bar\partial_A^*\tau + a''\varphi + \star((a'')^\dagger\wedge \star\psi)
  \end{pmatrix}
  \in \sF_2,
  \\
  \text{for all } (a'',\sigma,\tau) \in \sF_1.
\end{multline}
The complex linear and antilinear components of $\cT_{\bar\partial_A,\varphi,\psi,r}$ are now evident from \eqref{eq:H_dbar_APhi^01_explicit_perturbed}: for all $(a'',\sigma,\tau) \in \sF_1$, we have
\begin{subequations}
  \label{eq:H_dbar_APhi^01_explicit_perturbed_complex_linear_antilinear}
\begin{align}
  \label{eq:H_dbar_APhi^01_explicit_perturbed_complex_linear}
  \cT_{\bar\partial_A,\varphi,\psi,r}'(a'',\sigma,\tau)
  &=
  \begin{pmatrix}
    \bar\partial_A^*a'' - R_\varphi^*\sigma + \frac{ir}{4}p_\psi(\tau)
    \\
    \bar\partial_Aa'' - (\tau\otimes\varphi^*)_0
    \\
    \bar\partial_A\sigma+\bar\partial_A^*\tau + a''\varphi
  \end{pmatrix}
  \in \sF_2,
  \\
   \label{eq:H_dbar_APhi^01_explicit_perturbed_complex_antilinear}
  \cT_{\bar\partial_A,\varphi,\psi,r}''(a'',\sigma,\tau)
  &=
  \begin{pmatrix}
    (R_\psi^*\tau)^\dagger + \frac{ir}{4}p_\psi(\tau)^\dagger
    \\
    - \frac{1}{4}N_J^*(a'')^\dagger - (\psi\otimes\sigma^*)_0
    \\
    \star((a'')^\dagger\wedge \star\psi)
  \end{pmatrix}
  \in \sF_2.
\end{align}
\end{subequations}
We observe that $\cT_{\bar\partial_A,\varphi,\psi,r}'' = 0$ if and only if $N_J \equiv 0$ and $\psi \equiv 0$. Clearly, $\cT_{\bar\partial_A,\varphi,\psi,r}''$ is a zeroth order linear differential operator with smooth coefficients and thus defines a bounded complex antilinear operator between the complex Hilbert spaces $\sH_1$ and $\sH_2$ defined by the $L^2$ completions of the complex Fr\'echet spaces $\sF_1$ and $\sF_2$ in \eqref{eq:sEkC}. On the other hand, $\cT_{\bar\partial_A,\varphi,\psi,r}'$ is a first order complex linear elliptic differential operator with smooth coefficients and thus
%PF11-19-2025 Add reference
defines an unbounded complex linear Fredholm operator between $\sH_1$ and $\sH_2$ with dense domain given by the $W^{1,2}$ completion of the complex Fr\'echet space $\sF_1$. The operator
\begin{equation}
   \label{eq:H_dbar_APhi^01_explicit_perturbed_complex_linear_self-adjoint}
   T_{\bar\partial_A,\varphi,\psi,r}'
   :=
   \begin{pmatrix}
     0 & \cT_{\bar\partial_A,\varphi,\psi,r}^{\prime,*}
     \\
     \cT_{\bar\partial_A,\varphi,\psi,r}' & 0
   \end{pmatrix}
   \in
   \End(\sF_1\oplus\sF_2)
\end{equation}
is an $L^2$ self-adjoint first order complex linear elliptic differential operator with smooth coefficients on the complex Fr\'echet space $\sF_1\oplus\sF_2$. Therefore, $T_{\bar\partial_A,\varphi,\psi,r}'$ defines a self-adjoint (and thus closed) operator on the complex Hilbert space $\sH_1\oplus\sH_2$ with compact resolvent and dense domain given by the $W^{1,2}$ completion of $\sF_1\oplus\sF_2$. In particular, the operator $\cT_{\bar\partial_A,\varphi,\psi,r}' \in \Hom(\sH_1,\sH_2)$ obeys the hypotheses of Lemma \ref{lem:Eigenvalues_densely_defined_unbounded_linear_operators}.

Equation \eqref{eq:Perturbed_Equation_Equivalence_Of_Deformation_Complex} and the definitions \eqref{eq:Perturbed_Deformation_Operator} and \eqref{eq:Perturbed_AC_Deformation_Operator} of $\sT_{A,\varphi,\psi,r}$ and $\cT_{\bar\partial_A,\varphi,\psi,r}$, respectively, yield the relation
\[
  \sT_{A,\varphi,\psi,r} = \Upsilon_2 \circ \cT_{\bar\partial_A,\varphi,\psi,r} \circ \Upsilon_1^{-1}
  \in \Hom(\sE_1,\sE_2),
\]
where $\sE_1$ and $\sE_2$ are the real Fr\'echet spaces in \eqref{eq:sEk} and $\Upsilon_1:\sF_1\to\sE_1$ and $\Upsilon_2:\sF_2\to\sE_1$ are the real linear isomorphisms in \eqref{eq:Isomorphism_sEkC_to_sEk}. The canonical almost complex structures (given by scalar multiplication by $i=\sqrt{-1}$) on $\sF_1$ and $\sF_2$ and the isomorphisms \eqref{eq:Isomorphism_sEkC_to_sEk} induce almost complex structures $J_1$ and $J_2$ on $\sE_1$ and $\sE_2$, respectively. The decomposition
\[
  T_{\bar\partial_A,\varphi,\psi,r} = T_{\bar\partial_A,\varphi,\psi,r}' + T_{\bar\partial_A,\varphi,\psi,r}''
\]
into its complex linear and antilinear components as operators in $\Hom(\sF_1,\sF_2)$ therefore induces a decomposition
\[
  \sT_{A,\varphi,\psi,r} = \sT_{A,\varphi,\psi,r}' + \sT_{A,\varphi,\psi,r}''
\]
of $\sT_{A,\varphi,\psi,r}$ into its complex linear and antilinear components as operators in $\Hom(\sE_1,\sE_2)$.

The distinction between $\sT_{A,\varphi,\psi,r}'$ and $\sT_{A,\varphi,\psi,r}''$ can be seen most clearly by restricting to the ``coordinate'' subspaces of $\sF_1$. For example, if $\sigma\equiv 0$ and $\tau\equiv 0$, then
\begin{equation}
  \label{eq:H_dbar_APhi^01_explicit_perturbed_complex_linear_antilinear_sigma_and_tau_zero}
  \cT_{\bar\partial_A,\varphi,\psi,r}'(a'',0,0)
  =
  \begin{pmatrix}
    \bar\partial_A^*a''
    \\
    \bar\partial_Aa''
    \\
    a''\varphi
  \end{pmatrix}
  \in \sF_2
  \quad\text{and}\quad
  \cT_{\bar\partial_A,\varphi,\psi,r}''(a'',0,0)
  =
  \begin{pmatrix}
    0
    \\
    - \frac{1}{4}N_J^*(a'')^\dagger
    \\
    \star((a'')^\dagger\wedge \star\psi)
  \end{pmatrix}
  \in \sF_2.
\end{equation}
Similarly, if $a''\equiv 0$ and $\tau\equiv 0$, then
\begin{equation}
  \label{eq:H_dbar_APhi^01_explicit_perturbed_complex_linear_antilinear_a"_and_tau_zero}
  \cT_{\bar\partial_A,\varphi,\psi,r}'(0,\sigma,0)
  =
  \begin{pmatrix}
    - R_\varphi^*\sigma
    \\
    0
    \\
    \bar\partial_A\sigma
  \end{pmatrix}
  \in \sF_2
  \quad\text{and}\quad
  \cT_{\bar\partial_A,\varphi,\psi,r}''(0,\sigma,0)
  =
  \begin{pmatrix}
    0
    \\
    - (\psi\otimes\sigma^*)_0
    \\
    0
  \end{pmatrix}
  \in \sF_2.
\end{equation}
Finally, if $a''\equiv 0$ and $\sigma\equiv 0$, then
\begin{equation}
  \label{eq:H_dbar_APhi^01_explicit_perturbed_complex_linear_antilinear_a"_and_sigma_zero}
  \cT_{\bar\partial_A,\varphi,\psi,r}'(0,0,\tau)
  =
  \begin{pmatrix}
    \frac{ir}{4}p_\psi(\tau)
    \\
    - (\tau\otimes\varphi^*)_0
    \\
    \bar\partial_A^*\tau
  \end{pmatrix}
  \in \sF_2
  \quad\text{and}\quad
  \cT_{\bar\partial_A,\varphi,\psi,r}''(0,0,\tau)
  =
  \begin{pmatrix}
    (R_\psi^*\tau)^\dagger + \frac{ir}{4}p_\psi(\tau)^\dagger
    \\
    0
    \\
    0
  \end{pmatrix}
  \in \sF_2.
\end{equation}
The distinction between $\sT_{A,\varphi,\psi,r}'$ and $\sT_{A,\varphi,\psi,r}''$ can be further highlighted by restricting in \eqref{eq:H_dbar_APhi^01_explicit_perturbed_complex_linear_antilinear_sigma_and_tau_zero} to subspaces of $\Omega^{0,1}(\fsl(E))$ spanned by eigenvectors of the elliptic operator $\bar\partial_A^*\bar\partial_A + \bar\partial_A\bar\partial_A^*$, restricting in \eqref{eq:H_dbar_APhi^01_explicit_perturbed_complex_linear_antilinear_a"_and_tau_zero} to subspaces of $\Omega^0(E)$ spanned by eigenvectors of the elliptic operator $\bar\partial_A^*\bar\partial_A$, and restricting in \eqref{eq:H_dbar_APhi^01_explicit_perturbed_complex_linear_antilinear_a"_and_sigma_zero}  to subspaces of $\Omega^{0,2}(E)$ spanned by eigenvectors of the elliptic operator $\bar\partial_A\bar\partial_A^*$.

\section[Approximation of orthogonal projections on Hilbert spaces]{Approximation of orthogonal projections onto finite-dimensional subspaces of Hilbert spaces}
\label{sec:Approximation_orthogonal_projections_finite-dim_subspaces_Hilbert_spaces}
In this section, we derive the essential results in linear algebra that we shall need in order to verify the key hypothesis \eqref{eq:Orthogonal_projections_induce_embeddings} in Corollary \ref{cor:Kuranishi_model_defined_by_Fredholm_map_Hilbert_spaces}, a version of the Kuranishi model adapted to our application, and prove the assertions in Remark \ref{rmk:AC_and_symplectic_structures_Kuranishi_model_defined_by_Fredholm_map_Hilbert_spaces}. We begin with the following lemma, which is proved by elementary methods, though we are unable to find a statement or proof in standard references. We refer the reader to the proof of Theorem \ref{mainthm:Donaldson_1996jdg_3_Hilbert_space} in Section \ref{sec:Proof_generalized_Donaldson_symplectic_subspace_criterion_spectral_projection} for similar statements and verifications.

\begin{lem}[Approximation of finite-dimensional subspaces of a Hilbert space]
\label{lem:Approximation_finite-dimensional_subspaces_Hilbert_space}  
Let $\sH$ be Hilbert space over $\KK=\RR$ or $\CC$ with a complete orthonormal basis $\{e_k\}_{k=1}^\infty$ and $V \subset \sH$ be a finite-dimensional subspace. If $\eps \in (0,1/2]$ and $\pi \in \End(\sH)$ is the orthogonal projection onto $V$ and $\pi_n \in \End(\sH)$ is the orthogonal projection onto the finite-dimensional Hilbert subspace $\sH_n \subset \sH$ spanned by $\{e_1,\ldots,e_n\}$, then there is a positive integer $N = N(\eps,V)$ such that the following hold for all $n \geq N$:
\begin{enumerate}
\item\label{item:pi_minus_pi_n_hom(V,sH)_and_pi_minus_pi_n_hom(Vn,sH)}
  The projections $\pi$ and $\pi_n$ obey
  \begin{equation}
    \label{eq:pi_minus_pi_n_hom(V,sH)_lessthan_eps_and_pi_minus_pi_n_hom(Vn,sH)_lessthan_2eps}
    \|\pi - \pi_n\|_{\Hom(V,\sH)} < \eps
    \quad\text{and}\quad
    \|\pi - \pi_n\|_{\Hom(V_n,\sH)} < \frac{2\eps}{1-\eps}.
  \end{equation}

\item\label{item:Norm_id-pi_circ_pi_n_and_id-pi_n_circ_pi}
  If $V_n := \pi_n(V) \subset \sH_n$, then
  \begin{equation}
    \label{eq:Norm_id-pi_circ_pi_n_and_id-pi_n_circ_pi}
    \|\id_V - \pi\pi_n\|_{\End(V)} < \eps
    \quad\text{and}\quad
     \|\id_{V_n} - \pi_n\pi\|_{\End(V_n)} < \frac{\eps}{1-\eps},
  \end{equation}
  and
\[
  \pi\circ\pi_n \in \GL(V) \quad\text{and}\quad \pi_n\circ\pi \in \GL(V_n),
\]
so the projections $\pi \in \Hom(V_n,V)$ and $\pi_n \in \Hom(V,V_n)$ are isomorphisms of $\KK$-vector spaces. 

\item\label{item:pi_n_V_Image_Orthogonality}
The subspaces $V$ and $V_n^\perp\cap \sH_n$ of $\sH$ are orthogonal.

\item\label{eq:pi_n_V_embedding_sHn}
The projection $\pi_n: V \hookrightarrow \sH_n$ induces a linear embedding of $V$ onto its image $V_n \subset \sH_n$ and
an isomorphism
\begin{equation}
\label{eq:Define_tilde_Hn}
 \pi_n: \sH\supset \tilde\sH_n := V \oplus (V_n^\perp\cap \sH_n) \to \sH_n.
\end{equation}
\end{enumerate}
\end{lem}

\begin{proof}
Consider Item \eqref{item:pi_minus_pi_n_hom(V,sH)_and_pi_minus_pi_n_hom(Vn,sH)}. Write $m = \dim V$ and let $\{v_1,\ldots,v_m\}$ be an orthonormal basis for $V$. For each $l \in \{1,\ldots,m\}$, we have
\[
  v_l = \sum_{k=1}^\infty \langle v_l,e_k\rangle_\sH\, e_k
  \quad\text{and}\quad
  \|v_l\|_\sH^2 = \sum_{k=1}^\infty |\langle v_l,e_k\rangle_\sH|^2 < \infty,
\]
and so there is a positive integer $N = N(\eps,V)$ such that for any $n \geq N$, the projection $\pi_nv_l = \sum_{k=1}^n \langle v_l,e_k\rangle_\sH e_k$ obeys
\begin{equation}
  \label{eq:Norm_square_v_k_minus_pi_n_v_k}
  \|v_l - \pi_nv_l\|_\sH^2 = \sum_{k=n+1}^\infty |\langle v_l,e_k\rangle_\sH|^2 < \frac{\eps}{m},
  \quad\text{for } l = 1,\ldots,m.
\end{equation}
If $v = \sum_{l=1}^m c_lv_l \in V$ with $c_l \in \KK$ for $l=1,\ldots,m$, we have $v - \pi_n v = \sum_{l=1}^m c_l(v_l - \pi_nv_l)$ and
\begin{align*}
  \|v - \pi_n v\|_\sH^2
  &=
  \left\| \sum_{l=1}^m c_l(v_l - \pi_nv_l) \right\|_\sH^2
  \leq
  \left(\sum_{l=1}^m |c_l|\|v_l - \pi_nv_l\|_\sH\right)^2
  \leq
  \max_{1\leq k\leq m}\|v_k - \pi_nv_k\|_\sH^2\left(\sum_{l=1}^m |c_l|\right)^2
  \\
  &<
  \frac{\eps}{m}\left(\sum_{l=1}^m |c_l|\right)^2 \quad\text{(by \eqref{eq:Norm_square_v_k_minus_pi_n_v_k})}
  \\
  &\leq
  \frac{\eps}{m}m\sum_{l=1}^m |c_l|^2 \quad\text{(by inequality between arithmetic and quadratic means)}
  \\
  &=
  \eps\|v\|_\sH^2.
\end{align*}
%COMMENT See https://en.wikipedia.org/wiki/QM-AM-GM-HM_inequalities
(See Section \ref{sec:Operator_norm_choices_Donaldson_symplectic_subspace_criteria} for the inequality $A \leq Q$ between the arithmetic and quadratic means of a finite set of positive real numbers.) Consequently, noting that $\pi = \id$ on $V$,
\[
  \|\pi - \pi_n\|_{\Hom(V,\sH)} = \sup_{v \in V\less\{0\}} \frac{\|v - \pi_nv\|_\sH}{\|v\|_\sH} < \eps,
\]
which gives the first inequality claimed in \eqref{eq:pi_minus_pi_n_hom(V,sH)_lessthan_eps_and_pi_minus_pi_n_hom(Vn,sH)_lessthan_2eps}.

Consider Item \eqref{item:Norm_id-pi_circ_pi_n_and_id-pi_n_circ_pi}. We observe that
\begin{multline*}
  \|\id_V - \pi\pi_n\|_{\End(V)}
  = \|\pi^2 - \pi\pi_n\|_{\End(V)}
  = \|\pi(\pi - \pi_n)\|_{\End(V)}
  \\
  \leq \|\pi\|_{\Hom(\sH,V)}\|\pi - \pi_n\|_{\Hom(V,\sH)}
  = \|\pi - \pi_n\|_{\Hom(V,\sH)}
  < \eps,
\end{multline*}
which gives the first inequality in \eqref{eq:Norm_id-pi_circ_pi_n_and_id-pi_n_circ_pi} and thus, by Rudin \cite[Theorem 10.7, p. 249]{Rudin}, we obtain $\pi\pi_n \in \GL(V)$.

Similarly, suppose $w \in V_n$ with $\|w\|_\sH=1$, so that $w = \pi_nv = \pi_n^3v$ for some $v \in V\less\{0\}$. We compute
\begin{multline*}
  \|w - \pi_n\pi w\|_\sH
  = \|\pi_n^3v - \pi_n\pi\pi_n v\|_\sH
  = \|\pi_n^3v - \pi_n\pi^2v + \pi_n\pi^2v - \pi_n\pi\pi_n v\|_\sH
  \\
  \leq \|\pi_n^3v - \pi_n\pi^2v\|_\sH + \|\pi_n\pi^2v - \pi_n\pi\pi_n v\|_\sH
  \\
  = \|\pi_nv - \pi_nv\|_\sH + \|\pi_n\pi(\pi - \pi_n)v\|_\sH = 0 + \|\pi_n\pi(\pi - \pi_n)v\|_\sH
  \\
  \leq \|\pi_n\pi\|_{\End(\sH)}\|(\pi - \pi_n)v\|_\sH \leq \|(\pi - \pi_n)v\|_\sH
  \\
  \leq \|\pi - \pi_n\|_{\Hom(V,\sH)}\|v\|_\sH
  < \eps\|v\|_\sH,
\end{multline*}
where the final inequality follows from the first inequality in \eqref{eq:pi_minus_pi_n_hom(V,sH)_lessthan_eps_and_pi_minus_pi_n_hom(Vn,sH)_lessthan_2eps}. But
\begin{multline*}
  \|v\|_\sH \leq \|v-\pi_nv\|_\sH + \|\pi_nv\|_\sH = \|(\pi-\pi_n)v\|_\sH + \|\pi_nv\|_\sH
  < \eps\|v\|_\sH + \|\pi_nv\|_\sH \quad\text{(by \eqref{eq:pi_minus_pi_n_hom(V,sH)_lessthan_eps_and_pi_minus_pi_n_hom(Vn,sH)_lessthan_2eps})}
  \\
  = \eps\|v\|_\sH + \|w\|_\sH,
\end{multline*}
so that
\[
  (1-\eps)\|v\|_\sH < \|w\|_\sH.
\]
Therefore,
\[
  \|w - \pi_n\pi w\|_\sH < \frac{\eps}{1-\eps}\|w\|_\sH,
\]  
for all $w \in V_n$ with $\|w\|_\sH=1$, and so
\[
  \|\id_{V_n} - \pi_n\pi\|_{\End(V_n)} < \frac{\eps}{1-\eps},
\]
which gives the second inequality in \eqref{eq:Norm_id-pi_circ_pi_n_and_id-pi_n_circ_pi}. Since $\eps \in (0,1/2]$ by hypothesis and thus $\eps/(1-\eps) \in (0,1]$, we obtain $\pi_n\pi \in \GL(V_n)$. This completes the proof of Item \eqref{item:Norm_id-pi_circ_pi_n_and_id-pi_n_circ_pi}.

To complete the proof of Item \eqref{item:pi_minus_pi_n_hom(V,sH)_and_pi_minus_pi_n_hom(Vn,sH)}, it remains to prove the second inequality in \eqref{eq:pi_minus_pi_n_hom(V,sH)_lessthan_eps_and_pi_minus_pi_n_hom(Vn,sH)_lessthan_2eps}. We again suppose that $w \in V_n$ with $\|w\|_\sH = 1$, so $w = \pi_nv$ for some $v \in V\less\{0\}$. Then
\begin{multline*}
  \|w-\pi w\|_\sH = \|\pi_nv-\pi \pi_nv\|_\sH = \|\pi_n\pi v-\pi^2v+\pi^2v-\pi \pi_nv\|_\sH
  \leq
  \|\pi_n\pi v-\pi^2v\|_\sH + \|\pi^2v-\pi \pi_nv\|_\sH
  \\
  =
  \|(\pi_n-\pi)\pi v\|_\sH + \|\pi(\pi-\pi_n)v\|_\sH
  \leq
  \|\pi_n-\pi\|_{\Hom(V,\sH)}\|\pi v\|_\sH + \|\pi\|_{\End(\sH)} \|(\pi-\pi_n)v\|_\sH
  \\
  <
  \eps\|\pi v\|_\sH + \eps\|\pi\|_{\End(\sH)} \|v\|_\sH
  \leq
  2\eps\|v\|_\sH
  <
  \frac{2\eps}{1-\eps}\|w\|_\sH.
\end{multline*}
This yields the second inequality in \eqref{eq:pi_minus_pi_n_hom(V,sH)_lessthan_eps_and_pi_minus_pi_n_hom(Vn,sH)_lessthan_2eps} and completes the proof of Item \eqref{item:pi_minus_pi_n_hom(V,sH)_and_pi_minus_pi_n_hom(Vn,sH)}.

We now prove Item \eqref{item:pi_n_V_Image_Orthogonality}.  For any $v\in V$ and $w\in V_n^\perp\cap\sH_n$,
we have
\[
\langle v,w\rangle_\sH=\langle \pi_n v, w\rangle_{\sH}=0,
\]
where the first equality holds because $w\in\sH_n$ and the second equality holds because $w\in V_n^\perp=\pi_n(V)^\perp$.
This completes the proof of Item \eqref{item:pi_n_V_Image_Orthogonality}.

Consider Item \eqref{eq:pi_n_V_embedding_sHn}. By Item \eqref{item:Norm_id-pi_circ_pi_n_and_id-pi_n_circ_pi}, there exist $R_n \in \GL(V_n)$ and $L_n \in \GL(V)$ such that
\[
  \pi_n\pi\circ R_n = \id_{V_n} \quad\text{and}\quad L_n\circ \pi\pi_n = \id_V,
\]
or in other words,
\[
  \pi_n\circ \pi R_n = \id_{V_n} \quad\text{and}\quad L_n\pi\circ \pi_n = \id_V.
\]
Therefore, $\pi_n \in \Hom(V,V_n)$ has a right inverse $\pi R_n \in \Hom(V_n,V)$ and a left inverse $L_n\pi \in \Hom(V_n,V)$ and they are equal by the usual category theory argument:
% COMMENT https://math.stackexchange.com/questions/4546551/show-that-if-a-homomorphism-has-a-left-and-a-right-inverse-then-it-is-an-isomor
\[
  L_n\pi = L_n\pi\circ \id_{V_n} = L_n\pi\circ  \pi_n\circ \pi R_n = \id_V\circ \pi R_n = \pi R_n.
\]
Thus, $\pi_n \in \Hom(V,V_n)$ is an isomorphism of $\KK$-vector spaces, as claimed. A similar argument shows that $\pi \in \Hom(V_n,V)$ is an isomorphism of $\KK$-vector spaces. Consequently, $\pi_n$ induces a $\KK$-linear embedding, $\pi_n: V \hookrightarrow \sH_n$, of $V$ onto its image $V_n$. Because $\pi_n:V\to \sH_n$ is an isomorphism, the definition of $\tilde\sH_n$ in \eqref{eq:Define_tilde_Hn} as $V\oplus (V_n^\perp\cap\sH_n)$ and the observation that the restriction of $\pi_n$ to $V_n^\perp\cap\sH_n$ is the identity imply that $\pi_n:\tilde\sH_n\to\sH_n$ is an isomorphism. This completes the proof of Lemma \ref{lem:Approximation_finite-dimensional_subspaces_Hilbert_space}.
\end{proof}

We apply Lemma \ref{lem:Approximation_finite-dimensional_subspaces_Hilbert_space} to prove an estimate for the difference of two orthogonal projections analogous to the estimate \eqref{eq:Norm_Pi"_lessthan_one_half} that arises in the proof of Theorem \ref{mainthm:Donaldson_1996jdg_3_Hilbert_space} in Section \ref{sec:Proof_generalized_Donaldson_symplectic_subspace_criterion_spectral_projection}.

\begin{cor}[Approximation of projections onto finite-dimensional subspaces of a Hilbert space]
\label{cor:Approximation_projections_onto_finite-dimensional_subspaces_Hilbert_space}
Continue the hypotheses and notation of Lemma \ref{lem:Approximation_finite-dimensional_subspaces_Hilbert_space}
and let $\tilde\pi_n:\sH\to\tilde\sH_n$ be orthogonal projection, where the linear subspace $\tilde\sH_n$ is defined in \eqref{eq:Define_tilde_Hn}. Then,
\begin{equation}
\label{eq:Global_Difference_of_Projections}
  \left\|\tilde\pi_n - \pi_n\right\|_{\End(\sH)} < \eps.
\end{equation}
\end{cor}

\begin{proof}
Recall from \eqref{eq:Define_tilde_Hn} that $\tilde\sH_n$ is the orthogonal direct sum
\[
  V\oplus (V_n^\perp\cap\sH_n) \subset V\oplus V^\perp = \sH.
\]
Because $V_n^\perp\cap\sH_n\subset\sH_n\cap\tilde\sH_n$ and $\tilde\pi_n$ and $\pi_n$ are the orthogonal projections onto $\tilde\sH_n$ and $\sH_n$, respectively, we have
\[
  \tilde\pi_n(v_n^\perp) = v_n^\perp = \pi_n(v_n^\perp), \quad\text{for all } v_n^\perp \in V_n^\perp\cap\sH_n.
\]
Thus, for all $v\in V$ and $v_n^\perp\in V_n^\perp\cap\sH_n$,
\[
(\tilde\pi_n-\pi_n)(v+v_n^\perp)
=
(\tilde\pi_n-\pi_n)v,
\]
Hence, noting that $\Ran\tilde\pi_n = \tilde\sH_n$ and $\Ran\pi_n = \sH_n$ and that $v+v_n^\perp$ and $v$ are arbitrary elements of $\tilde\sH_n$ and $V$, respectively,
\[
(\tilde\pi_n-\pi_n)\tilde\pi_n
=
(\tilde\pi_n-\pi_n)\pi.
\]
Combining the preceding equality with the equality $\tilde\pi_n\pi=\pi$ gives 
\begin{equation}
\label{eq:ProjectionEquality1}
(\tilde\pi_n-\pi_n)\tilde\pi_n
=
(\pi-\pi_n)\pi.
\end{equation}
We write $1$ for the identity element in $\End(\sH)$ and compute:
\begin{align*}
  \| (1-\pi_n)\tilde\pi_n\|_{\End(\sH)}
  &=
    \| (\tilde \pi_n-\pi_n)\tilde\pi_n\|_{\End(\sH)}
  \\
  &= \|(\pi-\pi_n)\pi\|_{\End(\sH)} \quad\text{(by \eqref{eq:ProjectionEquality1})}
  \\
  &\le \|\pi-\pi_n\|_{\Hom(V,\sH)}\|\pi\|_{\Hom(\sH,V)}
  \\
  &= \|\pi-\pi_n\|_{\Hom(V,\sH)}\|\pi\|_{\End(\sH)}
  \\
  &= \|\pi-\pi_n\|_{\Hom(V,\sH)} \quad\text{(because $\|\pi\|_{\End(\sH)}=1$),}
\end{align*}
and so
\begin{equation}
\label{eq:OrthogProjectionCompositionInequality1}
\| (1-\pi_n)\tilde\pi_n\|_{\End(\sH)}
\le
\|\pi-\pi_n\|_{\Hom(V,\sH)}.
\end{equation}
The inequalities \eqref{eq:OrthogProjectionCompositionInequality1} and \eqref{eq:pi_minus_pi_n_hom(V,sH)_lessthan_eps_and_pi_minus_pi_n_hom(Vn,sH)_lessthan_2eps} imply that 
\begin{equation}
\label{eq:OrthogProjectionCompositionInequality2}
\| (1-\pi_n)\tilde\pi_n\|_{\End(\sH)} < \eps.
\end{equation}
%TL10-14-2025: Tried to rewrite to indicate both the inequality and the isomorphism are needed to get the particular case of the Kato citation.
Therefore, because $\pi_n:\tilde\sH_n\to\sH_n$ is an isomorphism by Lemma \ref{lem:Approximation_finite-dimensional_subspaces_Hilbert_space} \eqref{eq:pi_n_V_embedding_sHn},
the orthogonal projection operators $\tilde\pi_n$ and $\pi_n$ satisfy the hypotheses of \cite[Section I.6.8, Theorem 6.34 (i), p. 56]{Kato}. Hence, the desired inequality \eqref{eq:Global_Difference_of_Projections} follows from \cite[Section I.6.8, Equation (6.51), p. 57]{Kato}. This completes the proof of Corollary \ref{cor:Approximation_projections_onto_finite-dimensional_subspaces_Hilbert_space}.
\end{proof}

We conclude this section by applying Corollary \ref{cor:Approximation_projections_onto_finite-dimensional_subspaces_Hilbert_space} to prove the existence of almost complex structures in the context of our application described in Remark \ref{rmk:AC_and_symplectic_structures_Kuranishi_model_defined_by_Fredholm_map_Hilbert_spaces}.

\begin{cor}[Almost complex structures on finite-dimensional subspaces of an almost complex Hilbert space]
\label{cor:Almost_complex_structures_on_finite-dimensional_subspaces_Hilbert_space}
Continue the notation and hypotheses of Lemma \ref{lem:Approximation_finite-dimensional_subspaces_Hilbert_space} and Corollary \ref{cor:Approximation_projections_onto_finite-dimensional_subspaces_Hilbert_space}. Assume in addition that $\KK = \RR$ and $J \in \End_\RR(\sH)$ is an almost complex structure that is orthogonal with respect to the inner product $\langle\cdot,\cdot\rangle_\sH$ on $\sH$. Assume further that the complete orthonormal basis $\{e_k\}_{k=1}^\infty$ of $\sH$ has the form $e_{k+1} = Je_k$ for all odd $k\geq 1$ and choose the positive integers $N=N(\eps,V)$ and $n\geq N$ to be even, so $J \restriction \sH_n$ is an almost complex structure.
Then,
\[
  A_n := \tilde\pi_n J \in \End_\RR(\tilde\sH_n)
\]
is skew-adjoint with respect to the restriction $\langle\cdot,\cdot\rangle_{\tilde\sH}$ of the inner product $\langle\cdot,\cdot\rangle_{\sH}$ on $\sH$ to the linear subspace $\tilde\sH_n \subset \sH$. Furthermore, if
%TL10-11-2025: Keeping an eye on this in case it has to go back to 1/8
$\eps \in (0,1/2]$, then $A_n$ is invertible and
\begin{equation}
\label{eq:AC_Structure_On_tildeH_Projection}
  \tilde J_n := (A_n^*A_n)^{-1/2}A_n \in \End_\RR(\tilde\sH_n)
\end{equation}
is an almost complex structure on $\tilde\sH_n$, where $A_n^* \in \End_\RR(\tilde\sH_n)$ is the adjoint of $A_n$ defined by the inner product on $\tilde\sH_n$ induced from $\sH$.
\end{cor}

\begin{proof}
The operator $A_n$ is skew-adjoint with respect to the inner product $\langle\cdot,\cdot\rangle_{\sH}$ since, noting that $\tilde\pi_n = \id$ on $\tilde\sH_n$, we have for all $x,y \in \tilde\sH_n \subset \sH$ that
\[
\langle A_nx,y\rangle_{\sH}
=\langle \tilde\pi_n J x,y\rangle_{\sH}
=\langle Jx,\tilde\pi_n y\rangle_{\sH}
=-\langle x,J\tilde\pi_n y\rangle_{\sH}
=-\langle x,A_n y\rangle_{\sH}.
\]
Under our hypothesis that $\eps \in (0,1/2]$, Corollary \ref{cor:Approximation_projections_onto_finite-dimensional_subspaces_Hilbert_space} yields
\[
  \left\|\tilde\pi_n - \pi_n\right\|_{\End(\sH)} < \frac{1}{2}.
\]
We can now apply \mutatis the proof of the invertibility property \eqref{eq:Pi_J_Pi_in_GL_RanPi} under the condition \eqref{eq:Norm_Pi"_lessthan_one_half} arising in the proof of Theorem \ref{mainthm:Donaldson_1996jdg_3_Hilbert_space} in Section \ref{sec:Proof_generalized_Donaldson_symplectic_subspace_criterion_spectral_projection} to obtain
\[
  \tilde\pi_n J \in \GL(\tilde\sH_n).
\]
In particular, we see that $A_n$ is invertible. The expression for the induced almost complex structure $\tilde J_n$ arises almost exactly as in the proof of \cite[Proposition 12.3, p. 85]{Cannas_da_Silva_lectures_on_symplectic_geometry} due to Cannas da Silva.
\end{proof}

\section[Hitchin's function and Hamiltonian circle actions on virtual moduli spaces]{Hitchin's function, non-degenerate two-form, and Hamiltonian circle action on smooth local virtual moduli spaces of non-Abelian monopoles}
\label{sec:Hitchin_function_Hamiltonian_circle_action_virtual_moduli_spaces}
In this section, we complete our proofs of Theorem \ref{mainthm:AH_structure_bounded_evalue_spaces_non-Abelian_monopoles_symp_4-mflds} and Corollary \ref{maincor:Almost_Hermitian_structure_moduli_space_non-Abelian_monopoles_symplectic_4-manifolds}. We choose $V$ in Definition \ref{defn:V_critical_point_Hitchin_function_moduli_space_non-Abelian_monopoles} to be
%PF11-13-2025 Check this change
%the real linear direct sum of real linear subspaces (see Remark \ref{rmk:Distinction_between_linear_direct_and_orthogonal_direct_sums_approx_projections_AC_structures}),
%PF7-15-2025 Recall prior defn
\begin{equation}
  \label{eq:tilde_bH_A_varphi_psi_r_nu_1}
  \mathbf{\tilde H}_{A,\varphi,\psi,r,\nu}^1
  :=
  \bH_{A,\varphi,\psi,r}^1
  % \oplus_\RR
  %PF11-13-2025 Check this change
  \oplus
  \left(\bH_{A,\varphi,\psi,r,\nu}^1 \cap \left(\mathbf{\tilde H}_{A,\varphi,\psi,r}^1\right)^\perp\right)
  \subset L^2\left(T^*X\otimes\su(E) \oplus E\oplus\Lambda^{0,2}(E)\right),
\end{equation}
where we define
%PF7-15-2025 Recall prior defn
\begin{equation}
\label{eq:Define_Low_Eigenvalue_Eigenspaces_Of_ComplexLaplacian1}
  \bH_{A,\varphi,\psi,r,\nu}^1
  \subset
  L^2\left(T^*X\otimes\su(E) \oplus E\oplus\Lambda^{0,2}(E)\right)
\end{equation}
to be the orthogonal direct sum of the eigenspaces of the complex linear Laplacian $\sT_{A,\varphi,\psi,r}^{\prime,*}\sT_{A,\varphi,\psi,r}'$ with eigenvalues less than a positive constant $\nu \notin \sigma(\sT_{A,\varphi,\psi,r}^{\prime,*}\sT_{A,\varphi,\psi,r}')$, with $\sT_{A,\varphi,\psi,r}'$ denoting
%PF7-15-2025 Add relevant refs
the complex linear part of the operator $\sT_{A,\varphi,\psi,r} = d_{A,\varphi,\psi,r}^1 + d_{A,\varphi,\psi}^{0,*}$, and define $\mathbf{\tilde H}_{A,\varphi,\psi,r}^1$ in \eqref{eq:tilde_bH_A_varphi_psi_r_k}.

Similarly, we write
%PF7-18-2025 Add ref to first usage
\begin{multline}
  \label{eq:tilde_bH_A_varphi_psi_r_nu_2}
  \mathbf{\tilde H}_{A,\varphi,\psi,r,\nu}^2
  :=
  %TL11-12-2025: I added extra parens here to match the $\bH^1$ above
  \bH_{A,\varphi,\psi,r}^2 \oplus \left(\bH_{A,\varphi,\psi,r,\nu}^2 \cap \left(\mathbf{\tilde H}_{A,\varphi,\psi,r}^2\right)^\perp\right)
  \\
  \subset
  L^2\left(\su(E) \oplus \su(E) \oplus \Lambda^{0,2}\fsl(E)\oplus \Lambda^{0,1}(E)\right),
\end{multline}
for the orthogonal direct sum of real linear subspaces, where we define
%PF7-15-2025 Recall prior defn
\begin{equation}
\label{eq:Define_Low_Eigenvalue_Eigenspaces_Of_ComplexLaplacian2}
  \bH_{A,\varphi,\psi,r,\nu}^2
  \subset
  L^2\left(\su(E) \oplus \su(E) \oplus \Lambda^{0,2}\fsl(E)\oplus \Lambda^{0,1}(E)\right)
\end{equation}
to be the direct sum of the eigenspaces of the complex linear Laplacian $\sT_{A,\varphi,\psi,r}'\sT_{A,\varphi,\psi,r}^{\prime,*}$ with eigenvalues less than $\nu$, and define $\mathbf{\tilde H}_{A,\varphi,\psi,r}^2$ in \eqref{eq:tilde_bH_A_varphi_psi_r_k}. We recall from Lemma \ref{lem:Eigenvalues_densely_defined_unbounded_linear_operators} \eqref{item:Eigenvalues_densely_defined_unbounded_linear_operators_TT*_from_T*T} that
\[
  \sigma\left(\sT_{A,\varphi,\psi,r}^{\prime,*}\sT_{A,\varphi,\psi,r}'\right)\less\{0\}
  =
  \sigma\left(\sT_{A,\varphi,\psi,r}'\sT_{A,\varphi,\psi,r}^{\prime,*}\right)\less\{0\}.
\]
%TL8-11-2025: Perhaps add: "Using Lemma \ref{lem:Approximation_finite-dimensional_subspaces_Hilbert_space}, we can choose"  (or whichever of the results of that section is most direct for this)
The constant $\nu$ is chosen large enough that the orthogonal projections,
\begin{subequations}
  \begin{gather}
  \label{eq:Pi1_Avarphi_psi_r_nu_prime}
  \Pi_{1,A,\varphi,\psi,r,\nu}':L^2\left(T^*X\otimes\su(E) \oplus E\oplus\Lambda^{0,2}(E)\right)
  \to
    \bH_{A,\varphi,\psi,r,\nu}^1,
    \\
  \label{eq:Pi2_Avarphi_psi_r_nu_prime}
  \Pi_{2,A,\varphi,\psi,r,\nu}':L^2\left(\su(E)\oplus\su(E)\oplus\Lambda^{0,2}(\fsl(E))\oplus\Lambda^{0,1}(E)\right)
  \to
  \bH_{A,\varphi,\psi,r,\nu}^2,
\end{gather}
\end{subequations}  
restrict to real linear isomorphisms for $k=1,2$:
\[
  \Pi_{k,A,\varphi,\psi,r,\nu}': \bH_{A,\varphi,\psi,r}^k
  \to
  \mathbf{\tilde H}_{A,\varphi,\psi,r}^k,
\]
where
%PF7-15-2025 Make sure all accents like \tilde go inside bold face font commands
\begin{equation}
  \label{eq:tilde_bH_A_varphi_psi_r_k}
  \mathbf{\tilde H}_{A,\varphi,\psi,r}^k
  := \Pi_{k,A,\varphi,\psi,r,\nu}'\left(\bH_{A,\varphi,\psi,r}^k\right)  \subset \bH_{A,\varphi,\psi,r,\nu}^k,
  \quad\text{for } k=1,2.
\end{equation}    
Note that
\[
  \mathbf{\tilde H}_{A,\varphi,\psi,r}^1
  \subset
  \Ker d_{A,\varphi,\psi}^{0,*}\cap L^2\left(T^*X\otimes\su(E) \oplus E\oplus\Lambda^{0,2}(E)\right).
\]
Recall that $[A,\varphi,\psi] \in \sM^0(E,g,J,\omega,r)$ is a fixed point of the $S^1$ action
%PF7-10-2025 Add ref
if and only if
\[
  \bJ(0,\varphi,\psi) \in \Ran d_{A,\varphi,\psi}^0,
\]
since we always have $\bJ(0,\varphi,\psi) \in \Ker d_{A,\varphi,\psi,r}^1$ for any point $[A,\varphi,\psi] \in \sM^0(E,g,J,\omega,r)$.

Let $\sM^\vir(E,g,J,\omega,r,\nu)$ denote the $S^1$-invariant smooth local virtual moduli containing $[A,\varphi,\psi] \in \sM^0(E,g,J,\omega,r)$ that is defined by the subspaces $\mathbf{\tilde H}_{A,\varphi,\psi,r,\nu}^k$ for $k=1,2$ via the Kuranishi Method as in Section \ref{subsec:Local_Kuranishi_model_nbhd_point_zero_locus_unperturbed_non-Abelian_monopole_equations}. We denote
\begin{equation}
  \label{eq:S1GeneratorExpression_tilde_K_virtual_moduli_space}
  \bZ_{A,\varphi,\psi} := \tilde\Pi_{1,A,\varphi,\psi,r,\nu}\bJ(0,\varphi,\psi),
\end{equation}
where
\begin{equation}
  \label{eq:tildePi1_Avarphi_psi_r_nu}
  \tilde\Pi_{1,A,\varphi,\psi,r,\nu}: L^2\left(T^*X\otimes\su(E) \oplus E\oplus\Lambda^{0,2}(E)\right)
  \to
  \mathbf{\tilde H}_{A,\varphi,\psi,r,\nu}^1
\end{equation}
is the $L^2$-orthogonal projection onto the real linear subspace $\mathbf{\tilde H}_{A,\varphi,\psi,r,\nu}^1$.

%PF7-15-2025 More precise source needed
One can use Frankel's Theorem \ref{thm:Frankel_almost_Hermitian} \eqref{item:Frankel_almost_Hermitian_FixedPointsAreZerosOfVField} to show
%PF7-12-2025 Add proof!!!
that $[A,\varphi,\psi] \in \sM^\vir(E,g,J,\omega,r,\nu)$ is a fixed point of the $S^1$ action
% PF7-10-2025 Add ref
if and only if
\[
  \bZ_{A,\varphi,\psi} = 0,
\]
As in \eqref{eq:Xi1_isomorphism_tildeKnu} in Remark \ref{rmk:AC_and_symplectic_structures_Kuranishi_model_defined_by_Fredholm_map_Hilbert_spaces}, we have an isomorphism of real vector spaces,
\[
  \Xi_{1,A,\varphi,\psi,r,\nu}: \mathbf{\tilde H}_{A,\varphi,\psi,r,\nu}^1 \to \bH_{A,\varphi,\psi,r,\nu}^1,
\]
that extends the isomorphism of real vector spaces $\Pi_{1,A,\varphi,\psi,r,\nu}': \bH_{A,\varphi,\psi,r}^1 \to \mathbf{\tilde H}_{A,\varphi,\psi,r}^1$. By construction as a bounded eigenvalue subspace defined by a direct sum of the eigenspaces of the complex linear Laplacian $\sT_{A,\varphi,\psi,r}^{\prime,*}\sT_{A,\varphi,\psi,r}'$, the vector space $\bH_{A,\varphi,\psi,r,\nu}^1$ has an almost complex structure $\bJ$. We recall from Remark 
\ref{rmk:AC_and_symplectic_structures_Kuranishi_model_defined_by_Fredholm_map_Hilbert_spaces} that the vector space $\mathbf{\tilde H}_{A,\varphi,\psi,r,\nu}^1$ admits an operator
\[
  \tilde\Pi_{1,A,\varphi,\psi,r,\nu}\bJ \in \End\left(\mathbf{\tilde H}_{A,\varphi,\psi,r,\nu}^1\right)
\]
that is skew-adjoint with respect to the Riemannian metric $\bg$ in \eqref{eq:L2_metric_affine_space_spinu_pairs} and, for large enough $\nu$, has inverse
\[
  \left(\tilde\Pi_{1,A,\varphi,\psi,r,\nu}\bJ\right)^{-1} \in \End\left(\mathbf{\tilde H}_{A,\varphi,\psi,r,\nu}^1\right).
\]
Let $(A,\varphi,\psi)$ be a smooth unitary triple representing a point $[A,\varphi,\psi] \in \sM^0(E,g,J,\omega,r)$, so that by \eqref{eq:df_at_Avarphipsi_direction_asigmatau} we have
\begin{align*}
  (df)_{A,\varphi,\psi}(a,\sigma,\tau)
  &=
  \bg\left((0,\varphi,\psi),(a,\sigma,\tau)\right)_{L^2(X)},
  \\
  &\qquad \text{for all } (a,\sigma,\tau)
  \in L^2\left(T^*\otimes\su(E)\oplus E\oplus\Lambda^{0,2}(E)\right).
\end{align*}
Let $\sM^\vir(E,g,J,\omega,r,\nu)$ denote the $S^1$ invariant smooth local virtual moduli space containing $[A,\varphi,\psi] \in \sM^0(E,g,J,\omega,r)$ that is defined by the modified bounded eigenvalue spaces $\mathbf{\tilde H}_{A,\varphi,\psi,r,\nu}^1$ and $\mathbf{\tilde H}_{A,\varphi,\psi,r,\nu}^2$. The equality of one-forms \eqref{eq:df_at_Avarphipsi_direction_asigmatau} pulls back to an equality of one-forms on the subspace $\mathbf{\tilde H}_{A,\varphi,\psi,r,\nu}^1$,
\begin{equation}
  \label{eq:df_at_Avarphipsi_direction_asigmatau_in_tilde_bH1_Avarphipsi_rnu}
  (df)_{A,\varphi,\psi}(a,\sigma,\tau)
  =
  \bg\left((0,\varphi,\psi),(a,\sigma,\tau)\right)_{L^2(X)},
  \quad
  \text{for all } (a,\sigma,\tau) \in \mathbf{\tilde H}_{A,\varphi,\psi,r,\nu}^1.
\end{equation}
The smooth unitary triple $(A,\varphi,\psi)$ represents an $\mathbf{\tilde H}_{A,\varphi,\psi,r,\nu}^1$-critical point $[A,\varphi,\psi] \in \sM^0(E,g,J,\omega,r)$ in the sense of Definition \ref{defn:V_critical_point_Hitchin_function_moduli_space_non-Abelian_monopoles} if and only if
\begin{equation}
  \label{eq:df_at_Avarphipsi_vanishes_on_tilde_bH1_Avarphipsi_rnu}
  (df)_{A,\varphi,\psi} = 0 \quad\text{on } \mathbf{\tilde H}_{A,\varphi,\psi,r,\nu}^1,
\end{equation}
that is,
\[
  (\grad_\bg f)_{A,\varphi,\psi,r,\nu} := \tilde\Pi_{1,A,\varphi,\psi,r,\nu}(0,\varphi,\psi) = 0.
\]  
Equivalently, $(A,\varphi,\psi) \in \tilde\sM^0(E,g,J,\omega,r)$ represents a critical point \emph{in the standard sense} of the smooth function $f:\sM^\vir(E,g,J,\omega,r,\nu) \to \RR$ if and only if the identity \eqref{eq:df_at_Avarphipsi_vanishes_on_tilde_bH1_Avarphipsi_rnu} holds, since
\[
  \mathbf{\tilde H}_{A,\varphi,\psi,r,\nu}^1
  =
  T_{A,\varphi,\psi}\sM^\vir(E,g,J,\omega,r,\nu)
\]
that is, $\mathbf{\tilde H}_{A,\varphi,\psi,r,\nu}^1$ represents the tangent space to  $\sM^\vir(E,g,J,\omega,r,\nu)$ at $[A,\varphi,\psi] \in \sM^0(E,g,J,\omega,r)$ by construction of $\sM^\vir(E,g,J,\omega,r,\nu)$.

We recall the following facts from functional analysis: Let $\sH$ be a Hilbert space and $P \in \End(\sH)$ be a bounded linear operator with closed range $\sR := \Ran P \subset \sH$. Denote $\sK := \Ker P \subset \sH$, so that $\Ker P^* = (\Ran P)^\perp$ and $\Ker P = (\Ran P^*)^\perp$. Then the operators $P \in \Hom(\sK^\perp,\sR)$ and $P^* \in \Hom(\sR,\sK^\perp)$ are bounded isomorphisms of vector spaces, while the operator $P \in \Hom(\sH,\sR)$ is surjective with left inverse $L := (P^*P)^{-1}P^* \in \End(\sH)$ with
\[
  LP = \pi_{\sK^\perp} \quad\text{and}\quad PLP = P.
\]
We shall apply the preceding observations to $P_{A,\varphi,\psi,r,\nu} := \tilde\Pi_{1,A,\varphi,\psi,r,\nu}\bJ$ on
\[
  \sH := L^2\left(T^*X\otimes\su(E) \oplus E\oplus\Lambda^{0,2}(E)\right),
\]
with range $\Ran P_{A,\varphi,\psi,r,\nu} = \mathbf{\tilde H}_{A,\varphi,\psi,r,\nu}^1$ and kernel $(\Ker P_{A,\varphi,\psi,r,\nu})^\perp \cong \mathbf{\tilde H}_{A,\varphi,\psi,r,\nu}^1$ (via the isomorphism $P_{A,\varphi,\psi,r,\nu}$ from $(\Ker P_{A,\varphi,\psi,r,\nu})^\perp$ onto $\Ran P_{A,\varphi,\psi,r,\nu}$), and left inverse $\bL_{A,\varphi,\psi,r,\nu}$. Thus,
\[
  \bL_{A,\varphi,\psi,r,\nu}\circ \tilde\Pi_{1,A,\varphi,\psi,r,\nu}\bJ
  =
  \tilde\Pi_{1,A,\varphi,\psi,r,\nu}
  =
  \begin{cases}
    \id &\text{on } \mathbf{\tilde H}_{A,\varphi,\psi,r,\nu}^1,
    \\
    0 &\text{on } \left(\mathbf{\tilde H}_{A,\varphi,\psi,r,\nu}^1\right)^\perp.
  \end{cases}
\]
Here, we implicitly use the facts that
\[
  (\Ker P_{A,\varphi,\psi,r,\nu})^\perp \cong \Ran P_{A,\varphi,\psi,r,\nu} = \mathbf{\tilde H}_{A,\varphi,\psi,r,\nu}^1,
\]
and that by Remark \ref{rmk:AC_and_symplectic_structures_Kuranishi_model_defined_by_Fredholm_map_Hilbert_spaces} we have, for large enough $\nu$,
%PF7-19-2025 Check meaning of orthogonal projection onto and inner product on \mathbf{\tilde H}_{A,\varphi,\psi,r,\nu}^1
\[
  \tilde\Pi_{1,A,\varphi,\psi,r,\nu}\bJ \in \GL\left(\mathbf{\tilde H}_{A,\varphi,\psi,r,\nu}^1\right).
\]
We next apply the construction to $Q_{A,\varphi,\psi,r,\nu} := (\tilde\Pi_{1,A,\varphi,\psi,r,\nu}\bJ)^2$ on the same Hilbert space $\sH$ and with the same range and kernel. We denote its left inverse by $\bG_{A,\varphi,\psi,r,\nu}$, so that
\[
  \bG_{A,\varphi,\psi,r,\nu}\circ \left(\tilde\Pi_{1,A,\varphi,\psi,r,\nu}\bJ\right)^2
  =
  \tilde\Pi_{1,A,\varphi,\psi,r,\nu}
  =
  \begin{cases}
    \id &\text{on } \mathbf{\tilde H}_{A,\varphi,\psi,r,\nu}^1,
    \\
    0 &\text{on } \left(\mathbf{\tilde H}_{A,\varphi,\psi,r,\nu}^1\right)^\perp.
  \end{cases}
\]
It remains to verify the forthcoming Hamiltonian relation \eqref{eq:Moment_map_virtual_moduli_space_non-Abelian_monopoles}.
Observe that
%PF7-11-2025 Import details from chapter 9
\[
  -\left(\tilde\Pi_{1,A,\varphi,\psi,r,\nu}\bJ\right)^2
  =
  \left(\tilde\Pi_{1,A,\varphi,\psi,r,\nu}\bJ\right)^*\left(\tilde\Pi_{1,A,\varphi,\psi,r,\nu}\bJ\right)
  > 0
  \quad\text{on } \mathbf{\tilde H}_{A,\varphi,\psi,r,\nu}^1,
\]
and
\[
  -\left(\tilde\Pi_{1,A,\varphi,\psi,r,\nu}\bJ\right)^2  \in \End\left(\mathbf{\tilde H}_{A,\varphi,\psi,r,\nu}^1\right)
\]
is a positive self-adjoint operator. Consequently, its inverse
\[
  -\left(\tilde\Pi_{1,A,\varphi,\psi,r,\nu}\bJ\right)^{-2}  \in \End\left(\mathbf{\tilde H}_{A,\varphi,\psi,r,\nu}^1\right)
\]
is also a positive self-adjoint operator. We now define a Riemannian metric on $\mathbf{\tilde H}_{A,\varphi,\psi,r,\nu}^1$ by
%TL8-11-2025: Motivate the need for $-\bG_{A,\varphi,\psi,r,\nu}$ in the following
\begin{multline}
  \label{eq:Modified_L2_metric_tilde_K_1}
  \mathbf{\tilde g}_{A,\varphi,\psi}\left((a_1,\sigma_1,\tau_1), (a_2,\sigma_2,\tau_2)\right)
  :=
  \bg\left(-\bG_{A,\varphi,\psi,r,\nu}(a_1,\sigma_1,\tau_1),
    (a_2,\sigma_2,\tau_2)\right),
  \\
  \text{for all } (a_k,\sigma_k,\tau_k)
  \in \mathbf{\tilde H}_{A,\varphi,\psi,r,\nu}^1 \text{ and } k=1,2,
\end{multline}
where $\bg$ is as in \eqref{eq:L2_metric_affine_space_spinu_pairs}, and
%PF7-12-2025 Change "non-degenerate two form" to "symplectic form" on vector spaces everywhere
and define a symplectic form on $\mathbf{\tilde H}_{A,\varphi,\psi,r,\nu}^1$ by
\begin{multline}
  \label{eq:Modified_fundamental_two_form_tilde_K_1}
  \mathbf{\tilde\bomega}_{A,\varphi,\psi}\left((a_1,\sigma_1,\tau_1), (a_2,\sigma_2,\tau_2)\right)
  :=
  \mathbf{\tilde g}_{A,\varphi,\psi}\left(\tilde\Pi_{1,A,\varphi,\psi,r,\nu}\bJ(a_1,\sigma_1,\tau_1), (a_2,\sigma_2,\tau_2)\right),
  \\
  \text{for all } (a_k,\sigma_k,\tau_k)
  \in \mathbf{\tilde H}_{A,\varphi,\psi,r,\nu}^1 \text{ and } k=1,2.
\end{multline}
%TL8-8-2025: Let's make this a labelled definition either of the a.c. structure or of the compatible triples so it can be referred back to.
By Cannas da Silva \cite[Proposition 12.3, p. 85]{Cannas_da_Silva_lectures_on_symplectic_geometry}, the Riemannian metric $\mathbf{\tilde g}_{A,\varphi,\psi}$ and symplectic form $\mathbf{\tilde\bomega}_{A,\varphi,\psi}$ determine an almost complex structure $\mathbf{\tilde J}_{A,\varphi,\psi}$ on $\mathbf{\tilde H}_{A,\varphi,\psi,r,\nu}^1$ that is compatible with $\mathbf{\tilde\bomega}_{A,\varphi,\psi}$ in the sense that
\begin{equation}
  \label{eq:Riemannian_metric_for_compatible_triple_on_bH1_A_varphi_psi_r_nu}
  \mathbf{\breve g}_{A,\varphi,\psi}
  :=
  \mathbf{\tilde\bomega}_{A,\varphi,\psi}\left(\cdot\,,\mathbf{\tilde J}_{A,\varphi,\psi}\cdot\,\right)
  \quad\text{on } \mathbf{\tilde H}_{A,\varphi,\psi,r,\nu}^1
\end{equation}
is a Riemannian metric (that is, a positive definite real inner product) on $\mathbf{\tilde H}_{A,\varphi,\psi,r,\nu}^1$. In particular, by Cannas da Silva around the displayed identity in \cite[Section 12.2, paragraph following proof of Proposition 12.3, p. 85]{Cannas_da_Silva_lectures_on_symplectic_geometry},
\begin{equation}
  \label{eq:Compatible_triple_on_bH1_A_varphi_psi_r_nu}
  \left(\mathbf{\breve g}_{A,\varphi,\psi}, \mathbf{\tilde J}_{A,\varphi,\psi}, \mathbf{\tilde\bomega}_{A,\varphi,\psi}\right)
\end{equation}
is a \emph{compatible triple} on the real vector space $\mathbf{\tilde H}_{A,\varphi,\psi,r,\nu}^1$ in the sense of Cannas da Silva around the displayed identity in \cite[Definition 12.2, p. 84]{Cannas_da_Silva_lectures_on_symplectic_geometry}. We can now complete the

\begin{proof}[Proof of Theorem \ref{mainthm:AH_structure_bounded_evalue_spaces_non-Abelian_monopoles_symp_4-mflds}]
Consider Assertion \eqref{item:AH_structure_tilde_bH_A_varphi_psi_r_nu_1}. We obtain a compatible triple on $\mathbf{\tilde H}_{A,\varphi,\psi,r,\nu}^1$ from \eqref{eq:Compatible_triple_on_bH1_A_varphi_psi_r_nu} by choosing
\[
  \left(\bg_{1,A,\varphi,\psi},\bJ_{1,A,\varphi,\psi},\bomega_{1,A,\varphi,\psi}\right)
  :=
  \left(\mathbf{\breve g}_{A,\varphi,\psi}, \mathbf{\tilde J}_{A,\varphi,\psi}, \mathbf{\tilde\bomega}_{A,\varphi,\psi}\right),
\]
and this proves Assertion \eqref{item:AH_structure_tilde_bH_A_varphi_psi_r_nu_1}.

Consider Assertion \eqref{item:AH_structure_tilde_bH_A_varphi_psi_r_nu_2}. The construction of a compatible triple on real vector space $\mathbf{\tilde H}_{A,\varphi,\psi,r,\nu}^2$ in \eqref{eq:tilde_bH_A_varphi_psi_r_nu_2} is simpler since we are not constrained by the need to carefully choose the symplectic form in order to prove the forthcoming Hamiltonian identity \eqref{eq:Moment_map_virtual_moduli_space_non-Abelian_monopoles} on $\mathbf{\tilde H}_{A,\varphi,\psi,r,\nu}^1$. By Remark \ref{rmk:AC_and_symplectic_structures_Kuranishi_model_defined_by_Fredholm_map_Hilbert_spaces} we have, for large enough $\nu$,
%PF7-19-2025 Check meaning of orthogonal projection onto and inner product on \mathbf{\tilde H}_{A,\varphi,\psi,r,\nu}^2
\[
  \tilde\Pi_{2,A,\varphi,\psi,r,\nu}\bJ \in \GL\left(\mathbf{\tilde H}_{A,\varphi,\psi,r,\nu}^2\right),
\]
where the real vector space $\mathbf{\tilde H}_{A,\varphi,\psi,r,\nu}^2$ is as in \eqref{eq:tilde_bH_A_varphi_psi_r_nu_2} and
\begin{equation}
  \label{eq:tildePi2_Avarphi_psi_r_nu}
  \tilde\Pi_{2,A,\varphi,\psi,r,\nu}: L^2\left(\su(E) \oplus \su(E) \oplus
    \Lambda^{0,2}(\fsl(E))\oplus\Lambda^{0,1}(E)\right)
  \to
  \mathbf{\tilde H}_{A,\varphi,\psi,r,\nu}^2,
\end{equation}
is the $L^2$-orthogonal projection onto the subspace $\mathbf{\tilde H}_{A,\varphi,\psi,r,\nu}^2$. We define a symplectic form by setting
\begin{equation}
  \label{eq:Symplectic_form_on_bH2_A_varphi_psi_r_nu}
  \bomega_{2,A,\varphi,\psi} := \bg\left(\tilde\Pi_{2,A,\varphi,\psi,r,\nu}\bJ\,\cdot,\cdot\right)
  \quad\text{on } \mathbf{\tilde H}_{A,\varphi,\psi,r,\nu}^2,
\end{equation}
where $\bg$ is as in \eqref{eq:L2_metric_affine_space_spinu_pairs}. By Cannas da Silva \cite[Proposition 12.3, p. 85]{Cannas_da_Silva_lectures_on_symplectic_geometry}, the Riemannian metric $\bg$ and symplectic form $\bomega_{2,A,\varphi,\psi}$ determine an almost complex structure $\bJ_{2,A,\varphi,\psi}$ on $\mathbf{\tilde H}_{A,\varphi,\psi,r,\nu}^2$ that is compatible with $\bomega_{2,A,\varphi,\psi}$ in the sense that
\begin{equation}
  \label{eq:Riemannian_metric_for_compatible_triple_on_bH2_A_varphi_psi_r_nu}
  \bg_{2,A,\varphi,\psi}
  :=
  \bomega_{2,A,\varphi,\psi}\left(\cdot\,,\bJ_{2,A,\varphi,\psi}\,\cdot\,\right)
  \quad\text{on } \mathbf{\tilde H}_{A,\varphi,\psi,r,\nu}^1
\end{equation}
is a Riemannian metric (that is, a positive definite real inner product) on $\mathbf{\tilde H}_{A,\varphi,\psi,r,\nu}^2$. Therefore,
\[
  \left(\bg_{2,A,\varphi,\psi},\bJ_{2,A,\varphi,\psi},\bomega_{2,A,\varphi,\psi}\right)
\]
is a compatible triple on $\mathbf{\tilde H}_{A,\varphi,\psi,r,\nu}^2$. This completes the proof of Theorem \ref{mainthm:AH_structure_bounded_evalue_spaces_non-Abelian_monopoles_symp_4-mflds}
\end{proof}

The expression \eqref{eq:S1GeneratorExpression_tilde_K_virtual_moduli_space} for $\bZ_{A,\varphi,\psi}$ defines a smooth vector field $\bZ$ on $\sM^\vir(E,g,J,\omega,r,\nu)$. We see that $\bZ$ is the generator of the $S^1$ action \eqref{eq:S1ZActionOnQuotientSpace} on $\sM^\vir(E,g,J,\omega,r,\nu)$ as follows. If $[A,\varphi,\psi] \in \sM^\vir(E,g,J,\omega,r,\nu)$, then
\[
  e^{i\theta}\cdot[A,\varphi,\psi] = [A,e^{i\theta}\varphi,e^{i\theta}\psi] \in \sM^\vir(E,g,J,\omega,r,\nu),
  \quad\text{for } \theta \in (-\pi,\pi),
\]
defines a smooth path in $\sM^\vir(E,g,J,\omega,r,\nu)$ that is the image of the following smooth curve of solutions in $\sA(E,H,A_d)\times W^{1,p}(E\oplus\Lambda^{0,2}(E))$ to the defining equations for $\sM^\vir(E,g,J,\omega,r,\nu)$,
\[
  e^{i\theta}\cdot(A,\varphi,\psi) = (A,e^{i\theta}\varphi,e^{i\theta}\psi)
  \in \sA(E,H,A_d) \times  W^{1,p}(E\oplus\Lambda^{0,2}(E)),
  \quad\text{for } \theta \in (-\eps,\eps).
\]
% TL6-21-2025: "Because the moduli space $\sM^0(E,g,J,\omega,r)$ is closed under the $S^1$ action," Or perhaps just put in "Therefore" to indicate that this statement follows from the previous statement about the $S^1$ orbit being a continuous path in the moduli space (rather than something the reader is supposed to figure out theirself)
%PF7-7-2025 Where is the above phrase supposed to be inserted?
The tangent vector to this curve at $\theta=0$ belongs to the kernel of the linearization \eqref{eq:SO3Monopoled1TaubesPerturbation} of the system \eqref{eq:SO(3)_monopole_equations_almost_Hermitian_perturbed_intro_regular} at $(A,\varphi,\psi)$,
\[
  \left.\frac{d}{d\theta}(A,e^{i\theta}\varphi,e^{i\theta}\psi)\right|_{\theta=0}
  =
  (0,i\varphi,i\psi)
  =
  \bJ(0,\varphi,\psi) \in \Ker d_{A,\varphi,\psi,r}^1.
\]
The projection of the preceding tangent vector onto the slice $\Ker d_{A,\varphi,\psi}^{0,*}$ through $(A,\varphi,\psi)$ yields
%PF7-16-2025 Add details
\[
  \bZ_{A,\varphi,\psi}
  \in
  \bH_{A,\varphi,\psi,r,\nu}^1.
\]
Because $[A,\varphi,\psi] \in \sM^\vir(E,g,J,\omega,r,\nu)$, an $S^1$ invariant smooth manifold, then $e^{i\theta}\cdot[A,\varphi,\psi]$, for $\theta \in (-\pi,\pi)$, is a smooth curve through $[A,\varphi,\psi]$ that is the image of the smooth curve $e^{i\theta}\cdot(A,\varphi,\psi)$ in the affine space. (Again, these curves are smooth since the $S^1$ actions \eqref{eq:S1ZActionOnQuotientSpace} and \eqref{eq:S1ZAction} are smooth.) The curve $e^{i\theta}\cdot[A,\varphi,\psi]$ defines the tangent vector $\bZ_{A,\varphi,\psi}$ to $\sM^\vir(E,g,J,\omega,r,\nu)$ at $(A,\varphi,\psi)$, so $\bZ$ generates the $S^1$ action on $\sM^\vir(E,g,J,\omega,r,\nu)$.

We now verify that the restriction of $f$ in \eqref{eq:Hitchin_function} to $\sM^\vir(E,g,J,\omega,r,\nu)$ is a \emph{Hamiltonian function} for the $S^1$ action  \eqref{eq:S1ZActionOnQuotientSpace} on $\sM^\vir(E,g,J,\omega,r,\nu)$ with respect to the non-degenerate two-form $\mathbf{\tilde\bomega}$ on $\sM^\vir(E,g,J,\omega,r,\nu)$ defined by \eqref{eq:Modified_fundamental_two_form_tilde_K_1}. For $(a,\sigma,\tau) \in \mathbf{\tilde H}_{A,\varphi,\psi,r,\nu}^1$, we compute
\begin{align*}
  &\left(\iota_\bZ\mathbf{\tilde\bomega}\right)_{A,\varphi,\psi}(a,\sigma,\tau)
  \\
  &\quad = \mathbf{\tilde\bomega}_{A,\varphi,\psi}\left(\bZ_{A,\varphi,\psi},(a,\sigma,\tau)\right)
  \\
  &\quad = \mathbf{\tilde g}_{A,\varphi,\psi}\left(\left(\tilde\Pi_{1,A,\varphi,\psi,r,\nu}\bJ\right)
    \bZ_{A,\varphi,\psi}, (a,\sigma,\tau)\right)
  \quad\text{(by definition \eqref{eq:Modified_fundamental_two_form_tilde_K_1} of $\mathbf{\tilde\bomega}$)}
  \\
  &\quad = \mathbf{\tilde g}_{A,\varphi,\psi}\left(\left(\tilde\Pi_{1,A,\varphi,\psi,r,\nu}\bJ\right)^2(0,\varphi,\psi),
    (a,\sigma,\tau)\right)
  \quad\text{(by definition \eqref{eq:S1GeneratorExpression_tilde_K_virtual_moduli_space} of $\bZ$)}
  \\
  &\quad = -\bg_{A,\varphi,\psi}\left(\bG_{A,\varphi,\psi,r,\nu}
    \left(\tilde\Pi_{1,A,\varphi,\psi,r,\nu}\bJ\right)^2(0,\varphi,\psi), (a,\sigma,\tau)\right)
  \quad\text{(by definition \eqref{eq:Modified_fundamental_two_form_tilde_K_1} of $\mathbf{\tilde g}$)}
  \\
  &\quad = -\bg\left(\tilde\Pi_{1,A,\varphi,\psi,r,\nu}(0,\varphi,\psi), (a,\sigma,\tau)\right)
  \\
  &\quad = -\bg\left((0,\varphi,\psi), (a,\sigma,\tau)\right)
  \quad\text{(because $(a,\sigma,\tau) \in \mathbf{\tilde H}_{A,\varphi,\psi,r,\nu}^1$ and $\tilde\Pi_{1,A,\varphi,\psi,r,\nu}$ is $\bg$-orthogonal).}
\end{align*}
Hence, we obtain by comparison with the identity \eqref{eq:df_at_Avarphipsi_direction_asigmatau_in_tilde_bH1_Avarphipsi_rnu} that
\begin{multline*}
  -(df)_{A,\varphi,\psi} = \left(\iota_\bZ\mathbf{\tilde\bomega}\right)_{A,\varphi,\psi}
  \quad\text{on } \mathbf{\tilde H}_{A,\varphi,\psi,r,\nu}^1,
  \\
  \text{for } (A,\varphi,\psi) \in \tilde\sM^\vir(E,g,J,\omega,r,\nu)
  \subset \sA(E,H,A_d)\times W^{1,p}\left(E\oplus\Lambda^{0,2}(E)\right).
\end{multline*}
This verifies the \emph{Hamiltonian identity},
%PF7-18-2025 Need to also argue that we've proved an equality of forms on a manifold, not just at possibly distinguished tangent space used to define the virtual moduli space
\begin{equation}
  \label{eq:Moment_map_virtual_moduli_space_non-Abelian_monopoles}
  -(df)_{A,\varphi,\psi} = \left(\iota_\bZ\mathbf{\tilde\bomega}\right)_{A,\varphi,\psi}
  \quad\text{for } [A,\varphi,\psi] \in \sM^\vir(E,g,J,\omega,r,\nu),
\end{equation}
that is, the function $-f$, for the restriction $f:\sM^\vir(E,g,J,\omega,r,\nu) \to \RR$ of the function $f$ in \eqref{eq:Hitchin_function}, is a Hamiltonian in the sense of \eqref{eq:MomentMap} for the $S^1$ action \eqref{eq:S1ZActionOnQuotientSpace} on the virtual moduli space $\sM^\vir(E,g,J,\omega,r,\nu)$ and generator $\bZ$ in \eqref{eq:S1GeneratorExpression_tilde_K_virtual_moduli_space} for that $S^1$ action. With these conclusions in hand, we can now complete the

\begin{proof}[Proof of Corollary \ref{maincor:Almost_Hermitian_structure_moduli_space_non-Abelian_monopoles_symplectic_4-manifolds}]
We obtain the desired compatible triple on $\sM^\vir(E,g,J,\omega,r,\nu)$ from \eqref{eq:Compatible_triple_on_bH1_A_varphi_psi_r_nu} by choosing
\[
  (\bg_1,\bJ_1,\bomega_1)
  :=
  \left(\mathbf{\breve g}, \mathbf{\tilde J}, \mathbf{\tilde\bomega}\right).
\]
The identity \eqref{eq:Moment_map_virtual_moduli_space_non-Abelian_monopoles} verifies that the restriction of $-f$ in \eqref{eq:Hitchin_function} to $\sM^\vir(E,g,J,\omega,r,\nu)$ is a Hamiltonian function for the $S^1$ action \eqref{eq:S1ZActionOnQuotientSpace} and the  non-degenerate two-form $\bomega_1 = \mathbf{\tilde\bomega}$ on $\sM^\vir(E,g,J,\omega,r,\nu)$. 
\end{proof}

\section[Equivalence of critical points and fixed points]{Equivalence of critical points of the Hitchin function and fixed points of the circle action on the moduli space of non-Abelian monopoles}
\label{sec:Equivalence_critical_points_Hamiltonian_fixed_points_S1_action_moduli_space}
We first state and then apply Corollary \ref{maincor:Almost_Hermitian_structure_moduli_space_non-Abelian_monopoles_symplectic_4-manifolds} to prove the following analogue of Feehan and Leness 
%TL12-4-2025: Updated
\cite[Theorem 12.6.17]{Feehan_Leness_introduction_virtual_morse_theory_so3_monopoles}, which characterized critical points of the Hitchin function on the moduli space of non-zero-section projective vortices over a complex K\"ahler surface in terms of fixed points of the $S^1$ action defined by scalar multiplication of $\CC^*$ on the sections.

\begin{thm}[Equivalence of $\mathbf{\tilde H}_{A,\varphi,\psi,r,\nu}^1$-critical points of Hitchin's function and $S^1$-fixed points on the virtual moduli space of non-Abelian monopoles]
\label{thm:Critical_points_Hitchin_Hamiltonian_function_moduli_space_non-Abelian_monopoles}
Let $(E,H)$ be a smooth Hermitian vector bundle with complex rank two over a closed, almost K\"ahler four-manifold $(X,g,J,\omega)$ with symplectic form $\omega = g(J\cdot,\cdot)$ and $A_d$ be a fixed smooth unitary connection on the Hermitian line bundle $\det E$. Then the following hold:
\begin{enumerate}
\item \emph{(Fixed point $\implies$ critical point)}
\label{item:Fixed_point_implies_critical_point}
If $[A,\varphi,\psi] \in \sM^0(E,g,J,\omega,r)$ as in \eqref{eq:Moduli_space_non-Abelian_monopoles_almost_Hermitian_Taubes_regularized_non-zero-section} is a fixed point of the $S^1$ action \eqref{eq:S1ZActionOnQuotientSpace} on $\sM^0(E,g,J,\omega,r)$ and $\nu$ is the large enough positive constant provided by Corollary \ref{maincor:Almost_Hermitian_structure_moduli_space_non-Abelian_monopoles_symplectic_4-manifolds}, then it is a $\mathbf{\tilde H}_{A,\varphi,\psi,r,\nu}^1$-critical point (in the sense of Definition \ref{defn:V_critical_point_Hitchin_function_moduli_space_non-Abelian_monopoles}) of the Hitchin function $f:\sM^0(E,g,J,\omega,r)\to\RR$ given by the restriction of $f$ in \eqref{eq:Hitchin_function}.

\item \emph{(Critical point $\implies$ fixed point)}
\label{item:Critical_point_implies_fixed_point}
If $[A,\varphi,\psi] \in \sM^0(E,g,J,\omega,r)$ as in \eqref{eq:Moduli_space_non-Abelian_monopoles_almost_Hermitian_Taubes_regularized_non-zero-section} is a $\mathbf{\tilde H}_{A,\varphi,\psi,r,\nu}^1$-critical point (in the sense of Definition \ref{defn:V_critical_point_Hitchin_function_moduli_space_non-Abelian_monopoles}) of the Hitchin function $f:\sM^0(E,g,J,\omega,r)\to\RR$ given by the restriction of $f$ in \eqref{eq:Hitchin_function}, then $[A,\varphi,\psi]$ is a fixed point of the $S^1$ action \eqref{eq:S1ZActionOnQuotientSpace} on $\sM^0(E,g,J,\omega,r)$.
\end{enumerate}
\end{thm}

\begin{proof}
Consider Assertion \eqref{item:Fixed_point_implies_critical_point}. If $[A,\varphi,\psi]$ is a \emph{fixed point} of the $S^1$ action \eqref{eq:S1ZActionOnQuotientSpace} on the virtual moduli space $\sM^\vir(E,g,J,\omega,r,\nu)$, then $\bZ_{A,\varphi,\psi} = 0$ by our discussion in Section \ref{sec:Hitchin_function_Hamiltonian_circle_action_virtual_moduli_spaces} and so the Hamiltonian identity \eqref{eq:Moment_map_virtual_moduli_space_non-Abelian_monopoles} implies that $[A,\varphi,\psi]$ is a critical point of $f:\sM^\vir(E,g,J,\omega,r,\nu) \to \RR$, Thus $[A,\varphi,\psi]$ is an $\mathbf{\tilde H}_{A,\varphi,\psi,r,\nu}^1$-critical point of $f:\sM^0(E,g,J,\omega,r,\nu) \to \RR$ since $\mathbf{\tilde H}_{A,\varphi,\psi,r,\nu}^1$ represents the tangent space to $\sM^\vir(E,g,J,\omega,r,\nu)$ at the point $[A,\varphi,\psi]$. This completes the proof of Assertion \eqref{item:Fixed_point_implies_critical_point}.

Consider Assertion \eqref{item:Critical_point_implies_fixed_point}. Suppose that $[A,\varphi,\psi] \in \sM^\vir(E,g,J,\omega,r,\nu)$ is a \emph{critical point} of $f$ in the sense that
\[
  (df)_{A,\varphi,\psi} = 0 \quad\text{on } \mathbf{\tilde H}_{A,\varphi,\psi,r,\nu}^1.
\]
The identity \eqref{eq:Moment_map_virtual_moduli_space_non-Abelian_monopoles} yields $\mathbf{\tilde\bomega}_{A,\varphi,\psi}(\bZ_{A,\varphi,\psi},\cdot) = 0$ on $\mathbf{\tilde H}_{A,\varphi,\psi,r,\nu}^1$ and because the two-form $\mathbf{\tilde\bomega}_{A,\varphi,\psi}$ in \eqref{eq:Modified_fundamental_two_form_tilde_K_1} is non-degenerate on $\mathbf{\tilde H}_{A,\varphi,\psi,r,\nu}^1$, we obtain that $\bZ_{A,\varphi,\psi} = 0$. We can then apply
%PF7-18-2025 Add details
Theorem \ref{thm:Frankel_almost_Hermitian} \eqref{item:Frankel_almost_Hermitian_FixedPointsAreZerosOfVField} to conclude that $[A,\varphi,\psi] \in \sM^\vir(E,g,J,\omega,r,\nu)$ is a fixed point of the $S^1$ action \eqref{eq:S1ZActionOnQuotientSpace}
of Definition \ref{defn:UnitaryZActionOnAffine} on $\sM^\vir(E,g,J,\omega,r,\nu)$ with tangent space $\mathbf{\tilde H}_{A,\varphi,\psi,r,\nu}^1$ at $(A,\varphi,\psi)$. Hence, that $[A,\varphi,\psi]$ is a fixed point of the $S^1$ action \eqref{eq:S1ZActionOnQuotientSpace} on the $S^1$-invariant moduli subspace $\sM^0(E,g,J,\omega,r)$. This completes the proof of Assertion \eqref{item:Critical_point_implies_fixed_point} and Theorem \ref{thm:Critical_points_Hitchin_Hamiltonian_function_moduli_space_non-Abelian_monopoles}.
\end{proof}

Lastly in this chapter, we complete the proof our analogue of Feehan and Leness 
%TL12-4-2025: Updated: please check if changing Thm. 4 to Thm. 8 is correct.
\cite[Theorem 8]{Feehan_Leness_introduction_virtual_morse_theory_so3_monopoles}, which characterized critical sets of the Hamiltonian $f$ in \eqref{eq:Hitchin_function} on the moduli space of non-zero-section non-Abelian monopoles over a complex K\"ahler surface as moduli subspaces of Seiberg--Witten monopoles.

\begin{proof}[Proof of Theorem \ref{mainthm:IdentifyCriticalPoints}]
By hypothesis, $[A,\varphi,\psi]$ is a critical point of $f:\sM^0(E,g,J,\omega,r) \to \RR$ in the sense of Definition \ref{defn:Critical_point_Hitchin_function_moduli_space_non-Abelian_monopoles} and so Theorem \ref{thm:Critical_points_Hitchin_Hamiltonian_function_moduli_space_non-Abelian_monopoles} implies that it is a fixed point of the circle action \eqref{eq:S1ZActionOnQuotientSpace} on $\sM^0(E,g,J,\omega,r)$.

From our previous work on non-Abelian monopoles \cite[Proposition 3.1 (2), p. 86 and Lemma 3.11, p. 93]{FL2a}, 
%TL12-4-2025: Updated
\cite[Proposition 6.7.2]{Feehan_Leness_introduction_virtual_morse_theory_so3_monopoles} (see also 
%TL12-4-2025: Updated
\cite[Section 6.6.3]{Feehan_Leness_introduction_virtual_morse_theory_so3_monopoles}) we know that $[A,\varphi,\psi]\in\sM^0(E,g,J,\omega,r)$ is a fixed point of the $S^1$ action \eqref{eq:S1ZActionOnQuotientSpace} on $\sM^0(E,g,J,\omega,r)$ if and only if $(A,\varphi,\psi)$ with $(\varphi,\psi) \not\equiv (0,0)$ is \emph{split} (as in Definition \ref{defn:Split_trivial_central-stabilizer_spinor_pair} \eqref{item:Split_spinor_pair}) with respect to the decomposition $E=L_1\oplus L_2$ as a direct sum of Hermitian line bundles, so $A=A_1\oplus A_2$ with $A_2 = A_d\otimes A_1^*$ (where $A_1^*$ is the induced unitary connection on $L_1^*$) and $\varphi = \varphi_1\oplus 0$ and $\psi = \psi_1\oplus 0$, where $\varphi_1 \in \Omega^0(L_1)$ and $\psi_1 \in \Omega^{0,2}(L_1)$, and $(A_1,\varphi_1,\psi_1)$ is a solution to the perturbed Seiberg--Witten monopole equations
%PF9-9-2024 Add reference
defined by \eqref{eq:Split_pair_omega-component_induced_SW_curvature_equation_regular}:
\begin{subequations}
\label{eq:SW_monopole_equations_regular_Taubes_perturbation}  
\begin{align}
  \label{eq:SW_monopole_equations_regular_Taubes_perturbation_FA1_omega}  
  F_{A_1}^\omega
  - \frac{i}{8}|\varphi_1|_{L_1}^2\omega + \frac{i}{8}|\psi_1|_{\Lambda^{0,2}(L_1)}^2\omega
  + \frac{ir|\psi_1|_{\Lambda^{0,2}(E)}^2}{2\left(\gamma^2 + |\psi_1|_{\Lambda^{0,2}(E)}^2\right)}\omega
  - \frac{1}{4}(\Lambda_\omega F_{A_d})\,\omega &= 0,
  \\
  % PF5-6-2025 Check this!
  \label{eq:SW_monopole_equations_regular_Taubes_perturbation_FA1_02}
  F_{A_1}^{0,2} - \varphi_1\otimes\psi_1^* &= 0,
  \\
  \label{eq:SW_monopole_equations_regular_Taubes_perturbation_Dirac}
  \bar\partial_{A_1}\varphi_1 + \bar\partial_{A_1}^*\psi_1 &= 0.
\end{align}
\end{subequations}
This completes the proof of Theorem \ref{mainthm:IdentifyCriticalPoints}.
\end{proof}

\chapter{Calculation of virtual Morse--Bott indices via Atiyah--Singer Index Theorem}
\label{chap:Calculation_virtual_Morse-Bott_indices_via_Atiyah-Singer_index_theorem}
In this chapter, we prove Theorem \ref{mainthm:MorseIndexAtReduciblesOnAlmostKahler}, Corollary \ref{maincor:MorseIndexAtReduciblesOnSymplecticWithSO3MonopoleCharacteristicClasses}, and Corollary \ref{maincor:Positivity_of_MorseIndexAtReduciblesOnSymplecticWithSO3MonopoleCharacteristicClasses} by computing the virtual Morse--Bott index of the Hitchin function at a point $[A,\varphi,\psi]$ in the moduli space $\sM^0(E,g,J,\omega,r)$ of non-zero-section solutions to the non-Abelian monopole equations  \eqref{eq:SO(3)_monopole_equations_almost_Hermitian_perturbed_intro_regular} with a regularized Taubes perturbation
when the unitary triple $(A,\varphi,\psi)$ is  split in the sense of Definition \ref{defn:Split_trivial_central-stabilizer_spinor_pair} with respect to an orthogonal decomposition $E = L_1\oplus L_2$, where $L_1$ and $L_2$ are Hermitian line bundles over $X$.

We begin in Section \ref{sec:S1EquivIndexBundles} by giving basic definitions and properties of $G$-equivariant Fredholm operators and their indices. In Section \ref{subsec:ACOperatorIndex}, we define approximately complex  operators as a class of real linear operators $\sT:\sH_1\to\sH_2$ between complex Hilbert spaces.  We construct a complex $S^1$-equivariant index for such an operator and show that this index equals that of $\sT'$, the complex linear component of $\sT$. In Section \ref{sec:AC_Structure_Comparison}, we prove that this definition of a complex $S^1$-equivariant index of an approximately complex operator is equal for two choices of almost complex structure in the definition of the index. To apply these results to the problem of computing the virtual Morse--Bott index,
%\eqref{eq:vMB_Ch10_Intro}, 
we describe the $S^1$-equivariant structure of the domain and codomain of a Kuranishi model
% \eqref{eq:ThickenedModuliKuranishiModel}
 in Sections \ref{sec:S1ActionOnSO3MonopoleDefOperator} and \ref{sec:S1EquivariantStructureForPerturbedSO3MonopoleDefOperator}.
In Section \ref{sec:S1ActionOnApproxComplexDefOperator}, we describe an $S^1$-equivariant structure on the equivalent deformation operator $\cT_{\bar\partial_A,\varphi,\psi,r}$ defined in \eqref{eq:Perturbed_AC_Deformation_Operator}  which makes the maps in the equality \eqref{eq:Perturbed_Equation_Equivalence_Of_Deformation_Complex} $S^1$-equivariant. This proves that the isomorphisms of Proposition \ref{prop:SO3MonopoleDeformationCorollary_Perturbed} are $S^1$-equivariant, giving isomorphisms in Lemma \ref{lem:S1EquivariantIsom_LowEigenvalueSpaces_DefOp_to_ApproxComplex} which show that we can compute the virtual Morse--Bott index 
%of the model \eqref{eq:ThickenedModuliKuranishiModel} 
in terms of a complex, $S^1$-equivariant index of the equivalent deformation operator $\cT_{\bar\partial_A,\varphi,\psi,r}$.  By the work in Section \ref{subsec:ACOperatorIndex}, this equals the complex $S^1$-equivariant index of the complex linear component $\cT_{\bar\partial_A,\varphi,\psi,r}'$. In Section \ref{sec:WeightDecompOfApproxComplexDeformationOperator}, we prove that the complex $S^1$-equivariant index of $\cT_{\bar\partial_A,\varphi,\psi,r}'$ equals that of an operator which is diagonal with respect to the weight spaces of its domain and codomain. The equivariant index of such diagonal operators can be computed in terms of the non-equivariant indexes of the components of the operator which we compute in Section \ref{sec:IndexComputation}. Combining these results gives a formula for the virtual Morse--Bott index of the model. %\eqref{eq:ThickenedModuliKuranishiModel}.
Finally, in Section \ref{sec:vMBMainTheorems} we
apply Corollary \ref{cor:Kuranishi_model_defined_by_Fredholm_map_Hilbert_spaces} to construct an $S^1$-equivariant Kuranishi model for an open neighborhood of the point $[A,\varphi,\psi]$ in $\sM^0(E,g,J,\omega,r)$ and use this to
prove Theorem \ref{mainthm:MorseIndexAtReduciblesOnAlmostKahler}, Corollary \ref{maincor:MorseIndexAtReduciblesOnSymplecticWithSO3MonopoleCharacteristicClasses}, and Corollary \ref{maincor:Positivity_of_MorseIndexAtReduciblesOnSymplecticWithSO3MonopoleCharacteristicClasses}.

\section{G-equivariant index bundles}
\label{sec:S1EquivIndexBundles}
We will show that the deformation operator \eqref{eq:Perturbed_Deformation_Operator} admits an $S^1$-equivariant structure in the following sense (see Feehan and Leness \cite[Definition 14.1.6]{Feehan_Leness_introduction_virtual_morse_theory_so3_monopoles} and Lawson and Michelsohn \cite[Definition III.9.1, p. 211]{LM}).

\begin{defn}[Representation ring]
  \label{defn:RepresentationRing}
(See Br\"ocker and tom Dieck \cite[Chapter II, Section 7, p. 103]{BrockertomDieck} or Lawson and Michelsohn \cite[Chapter III, Section 9, p. 211]{LM}.) Let $G$ be a compact Lie group.  If $\KK$ is a field, then the \emph{representation ring of $G$ over $\KK$}, denoted $R_\KK(G)$, is the free Abelian group generated by equivalence classes of  finite-dimensional, irreducible $\KK$ representations of $G$. Equivalently, $R_\KK(G)$ is the Grothendieck group (see Lang \cite[Chapter 1, Section 7, p. 39]{LangAlgebra}) of all finite-dimensional $\KK$ representations of $G$.
\end{defn}

\begin{rmk}[On the hypothesis that $G$ is a compact Lie group]
\label{rmk:GCompactLie}
As we only need the case $G=S^1$ in this monograph, we assume that $G$ is a compact Lie group in Definition \ref{defn:RepresentationRing} for simplicity.  This allows us to avoid discussing whether the equality $[U]=[V]+[W]$ holds only when $U=V\oplus W$ or when there is a short exact sequence $0\to V\to U\to W\to 0$ (see Lang \cite[Chapter 1, Section 7, p. 39]{LangAlgebra}).
\end{rmk}

\begin{exmp}[The complex and real representation rings of $S^1$]
\label{exmp:CC_RepRing_of_S1}
From, for example, Lawson and Michelsohn \cite[Example III.9.2, p. 212]{LM}, the ring $R_{\CC}(S^1)$ is isomorphic to the ring of Laurent polynomials $\ZZ[t,t^{-1}]$ with the element $t^m$, for $m\in \ZZ$, being given by the irreducible, complex representation defined by $\rho_m(e^{i\theta}):\CC\to\CC$ being given by $\rho_m(e^{i\theta})z=e^{im\theta}z$ for $z\in\CC$.   We call $m$ the \emph{weight} of the representation $\rho_m$.

By Feehan and Leness \cite[Lemma 4.2.6]{Feehan_Leness_introduction_virtual_morse_theory_so3_monopoles}, the representations $\rho_m$ and $\rho_{-m}$ are identified when they are considered as \emph{real} representations.
\qed\end{exmp}

As we will be working with representations where the complex structure is not obvious or unique, we make the following definition.

\begin{defn}[Complex representations and $G$-invariant almost complex structures]
\label{defn:GInvarComplexStructures}
Let $\rho:G\to\GL(V)$ be a representation of a group on a vector space $V$. An almost complex structure $J$ on $V$ is \emph{$G$-invariant} if $\rho(g)J=J\rho(g)$ for all $g\in G$ (compare the definition 
\cite[Equation (2.2.7)]{Feehan_Leness_introduction_virtual_morse_theory_so3_monopoles} of $S^1$-invariant tensors).  Equivalently, $J$ is $G$-invariant if and only if it is $G$-equivariant as a map from $V$ to $V$. If $J$ is a $G$-invariant complex structure on $V$, then $\rho:G\to \GL(V)$ is a complex linear representation on the complex vector space $(V,J)$. When $G$ is a compact Lie group, we write $[V,\rho,J]\in R_\CC(G)$ when it is necessary to specify the almost complex structure.  If the homomorphism $\rho$ is obvious, we will abbreviate $[V,\rho,J]$ by $[V,J]$.
\end{defn}

We now give an example of different almost complex structures defining different elements of $R_\CC(S^1)$.

\begin{exmp}[Dependence of a representation on the almost complex structure]
\label{exmp:Dependence_on_ac}
If $J_0$ and $J_1$ denote the almost complex structures on $V=\CC$ given by $J_0z=iz$ and $J_1z=-iz$,
then both $J_0$ and $J_1$ are $S^1$ invariant with respect to the representation $(\CC,\rho_m)$ appearing in  Example \ref{exmp:CC_RepRing_of_S1}.  However, $[\CC,\rho_m,J_0]=-[\CC,\rho_m,J_1]$.
\end{exmp}

In the following proposition, we give a criterion for when, given a representation $\rho:G\to \GL(V)$ of a compact Lie group and two $G$-invariant almost complex structures $J_0$ and $J_1$ on $V$, the equality $[V,\rho,J_0]=[V,\rho,J_1]$ holds. This result is a more formal presentation of Segal's remark following \cite[Proposition 1.3, p. 130]{Segal_Equivariant_K_Theory} on the discrete nature of representations of compact groups.

\begin{prop}[Homotopy invariance of $G$-representations]
\label{prop:Homotopy_Invariance_of_G_Representations}
Let $\rho:G\to\GL(V)$ be a representation of a compact Lie group $G$ on a finite-dimensional real vector space $V$ and let $[0,1]\ni t\to J_t\in \GL(V)$ be a continuous family of $G$-invariant almost complex structures on $V$. If $[V,\rho,J_t]$ is the element of $R_\CC(G)$ given by the $G$-invariant almost complex structure $J_t$ in Definition \ref{defn:GInvarComplexStructures}, then $[\sV,J_0]=[\sV,J_t]$ for all $t\in [0,1]$.
\end{prop}

We will prove Proposition \ref{prop:Homotopy_Invariance_of_G_Representations} using the framework of equivariant vector bundles. Recall from Segal \cite[p. 129]{Segal_Equivariant_K_Theory} that for a topological group $G$, a \emph{$G$-space} is a topological space $X$ on which $G$ acts continuously and that a \emph{complex $G$-vector bundle} over $X$ is a $G$-space $E$ with a $G$-equivariant, continuous map $p:E\to X$ such that
\begin{enumerate}
\item
  $p:E\to X$ is a continuous complex vector bundle,
\item
  For any $x\in X$ and any $g\in G$, the map $g:E_x\to E_{gx}$ is complex linear,  where $E_x=p^{-1}(x)$.
\end{enumerate}
A complex $G$-vector bundle over a point is thus a complex representation of $G$.

A \emph{homomorphism} $f:E\to F$ of $G$-bundles over $X$ is a continuous $G$-equivariant map which induces a homomorphism of vector spaces $f_x:E_x\to F_x$, for each $x\in X$.  In particular, if $E$ and $F$ are complex vector bundles, the map $f_x:E_x\to F_x$ is a complex linear map. Two $G$-bundles $E$ and $F$ over $X$ are \emph{$G$-equivariantly isomorphic} if there are homomorphisms of $G$-bundles $f:E\to F$ and $g:F\to E$ with $f\circ g=\id_F$ and $g\circ f=\id_E$. If $\phi:Y\to X$ is a continuous $G$-equivariant map and $E\to X$ is a $G$-vector bundle, then the pullback $\phi^*E\to Y$ is also a $G$-vector bundle (see Segal \cite[p. 130]{Segal_Equivariant_K_Theory}).

\begin{prop}[Homotopic maps induce isomorphic bundles]
\label{prop:HomotopicMapsForEquivBundles}
(See Segal \cite[Proposition 1.3, p. 130]{Segal_Equivariant_K_Theory}.)
Let $E\to X$ be a $G$-equivariant complex vector bundle where $G$ is a compact topological group. If $G$-equivariant maps $\phi_0,\phi_1:Y\to X$ are homotopic through $G$-equivariant maps, then $\phi_0^*E$ and $\phi_1^*E$ are $G$-equivariantly isomorphic as complex vector bundles.
\end{prop}

We can then give

\begin{proof}[Proof of Proposition \ref{prop:Homotopy_Invariance_of_G_Representations}]
Consider the complex vector bundle $V\times [0,1]\to [0,1]$ where the almost complex structure on the fiber $V \times\{t\}$ is $J_t$.  If $i_t:\{0\}\to [0,1]$ is the inclusion $i_t(0)=t$, then $i_t^*(V\times [0,1])\cong ( V,J_t)$.
Because the inclusions $i_0$ and $i_1$ are homotopic through $G$-equivariant maps, Proposition \ref{prop:HomotopicMapsForEquivBundles} implies that $V_0$ and $V_1$ are $G$-equivariantly isomorphic and thus define the same element of $R_\CC(G)$.
\end{proof}

We now define the equivariant index of an equivariant Fredholm operator.

\begin{defn}[Equivariant index of a Fredholm operator]
\label{defn:Index_of_GEquivariantDiffOperator}
Let $E$ and $F$ be Banach spaces over $\KK$ where $\KK=\RR$ or $\KK=\CC$. Let $G$ be a compact Lie group acting $\KK$-linearly  on $E$ and $F$. A Fredholm map $P:E\to F$  is \emph{$G$-equivariant} if $P(g s)=gPs$ for all $g\in G$ and $s\in E$.  If such a Fredholm map $P$ is $\KK$-linear, where $\KK=\RR$ or $\KK=\CC$, then the kernel and cokernel of a $P$ are finite-dimensional $\KK$ representations of $G$ (see Br\"ocker and tom Dieck \cite[Definition II.1.1, p. 65]{BrockertomDieck}.  The formal difference of these representations defines the \emph{$G$-equivariant index of $P$},
\begin{equation}
\label{eq:GEquivIndex}
\Ind_{G,\KK}(P):=[\Ker P] - [\Coker P] \in R_\KK(G),
\end{equation}
where $R_\KK(G)$ is the ring in Definition \ref{defn:RepresentationRing}.
\qed\end{defn}

Now assume that the Banach spaces $E$ and $F$ in Definition \ref{defn:Index_of_GEquivariantDiffOperator} admit $G$-invariant direct sum decompositions $E=E_1\oplus E_2$ and $F=F_1\oplus F_2$ and the operator $P$ satisfies $P=P_1\oplus P_2$, where $P_i:E_i\to F_i$ is a $G$-equivariant operator.  Then the kernel and cokernel of $P$ admit similar direct sum decompositions and we can write
\begin{equation}
\label{eq:AdditivityOfEquivariantIndex}
\Ind_{G,\KK}(P)
=
\Ind_{G,\KK}(P_1)+\Ind_{G,\KK}(P_2).
\end{equation}

We further note that the $G$-equivariant index is a homotopy invariant in the following sense (see Lawson and Michelsohn \cite[Proposition III.9.4, p. 213]{LM}).

\begin{prop}[Homotopy invariants of the $G$-equivariant index]
\label{prop:HomotopyInvarianceOfGindex}
Let $\sH_1$ and $\sH_2$ be Hilbert spaces over $\KK=\RR$ or $\KK=\CC$. For a compact Lie group $G$, let $\Fred_{G,\KK}(\sH_1,\sH_2)$ be the space of $G$-equivariant, $\KK$-linear Fredholm operators from $\sH_1$ to $\sH_2$. Then the map
\[
\Ind_{G,\KK}: \Fred_{G,\KK}(\sH_1,\sH_2)\to R_{\KK}(G),
\]
is constant on connected components.
\end{prop}

\section{The complex index of an approximately complex operator}
\label{subsec:ACOperatorIndex}
The definition \ref{eq:GEquivIndex} gives the real $G$-equivariant index of a $G$-equivariant Fredholm operator. We will now define a \emph{complex} $G$-equivariant index for an operator $\sT$ of the type discussed in Corollary \ref{maincor:Donaldson_1996jdg_3_Hilbert_space_non-self-adjoint} and show that this index equals the complex $G$-equivariant index of $\sT'$, the complex linear component of $\sT$. We recall the definitions of almost complex structures from Theorem \ref{mainthm:Donaldson_1996jdg_3_Hilbert_space} and Corollary \ref{maincor:Donaldson_1996jdg_3_Hilbert_space_non-self-adjoint}:

\begin{defn}[Approximately complex operators]
\label{defn:ApproximatelyComplex}
Let $\sH_k$ be a real Hilbert space with inner product $g_k = \langle\cdot,\cdot\rangle_{\sH_k}$ for $k=1,2$ and $g_k$-orthogonal almost complex structures $J_k:\sH_k\to\sH_k$. As in \eqref{eq:Complex_linear_and_anti-linear_operator_components}, for a densely defined, unbounded real linear operator $\sT \in \Hom(\sH_1,\sH_2)$, let
\[
  \sT':=\frac{1}{2}\left( \sT-J_2\sT J_1\right)
  \quad\text{and}\quad
  \sT'':=\frac{1}{2}\left( \sT+J_2\sT J_1\right),
\]
be the complex linear and complex antilinear components of $\sT$, respectively. We call $\sT$ \emph{approximately complex} if
\begin{enumerate}
\item
\label{item:ApproximatelyComplexDefinition1}
The complex linear component $\sT'$ of $\sT$ determines a densely-defined, complex linear, self-adjoint unbounded operator with compact resolvent,
\[
  T' = \begin{pmatrix} 0 & \sT^{\prime,*} \\ \sT' & 0 \end{pmatrix} \in \End(\sH_1\oplus\sH_2).
\]

\item
\label{item:ApproximatelyComplexDefinition2}
The complex antilinear component $\sT''$ of $\sT$ is compact.\qed
\end{enumerate}
\end{defn}

The conditions on the operator $T'$ in Definition \ref{defn:ApproximatelyComplex} imply that its spectrum $\sigma(T') \subset \RR$ is discrete and comprises eigenvalues of finite multiplicity and that there is a sequence of eigenvectors of $T'$ which form a complete orthonormal basis of $\sH = \sH_1\oplus\sH_2$: see
Lemma \ref{lem:Eigenvalues_densely_defined_unbounded_linear_operators} and
the discussion
% PF10-17-2025 Let's keep an eye on this reference should it change
around \eqref{eq:mu_in_specRtT_iff_t_plus_1over_mu_in_specT} in Appendix \ref{sec:Spectral_theory_unbounded_operators}. In our application, $\sT'$ will be the $L^2$ extension of a first order elliptic operator with smooth coefficients that acts on a smooth vector bundle over a closed, smooth manifold and so $T'$ will have a compact resolvent operator.
% PF7-21-2025 Add reference for elliptic operator implies compact resolvent.

\begin{lem}[$G$-equivariance of bounded-eigenvalue eigenspaces and projections]
\label{lem:S1InvariantSubspacesAndProjections}
For $k=1,2$, let $\sH_k$ be a real Hilbert space with inner product $g_k = \langle\cdot,\cdot\rangle_{\sH_k}$.
Let $G$ be a group and assume that for $k=1,2$ there are $g_k$-orthogonal $G$ actions on $\sH_k$. Let $\sS:\sH_1\to\sH_2$ be a $G$-equivariant operator satisfying the hypotheses of Lemma \ref{lem:Eigenvalues_densely_defined_unbounded_linear_operators}. For $\nu>0$ and $\nu\notin\sigma(\sS^*\sS)$, define
\begin{equation}
\label{eq:DefineLowEigenvalueSpaces}
\bH_\nu^1(\sS)\subset\sH_1
\quad\text{and}\quad
\bH_\nu^2(\sS)\subset\sH_2
\end{equation}
to be the orthogonal direct sum of the eigenspaces of $\sS^*\sS$ and $\sS\sS^*$, respectively, corresponding to eigenvalues less than $\nu$. For $k=1,2$, let $\Pi_{\nu,k}:\sH_k\to \bH_\nu^k(\sS)$ be the orthogonal projection. Then the following hold:
\begin{enumerate}
\item
\label{item:S1ActionOnLowEigenvalueSpaces}
The vector spaces $\bH_\nu^k(\sS)$ are invariant under the $G$ action on $\sH_k$ for $k=1,2$. 
\item
\label{item:S1Equiv_of_ProjectionOperators}
The operators $\Pi_{\nu,k}$ are $G$-equivariant for $k=1,2$.
\end{enumerate}
\end{lem}

\begin{proof}
The conclusions follow immediately from the $G$-equivariance of $\sS$ and the resulting $G$-equivariance of $\sS^*$.
\end{proof}

\begin{lem}[Injectivity of orthogonal projection onto the bounded-eigenvalue eigenspace for the complex linear component of an approximately complex operator]
\label{lem:InjectivityOfProjection_for_KurnaishiModel}
Let $\sH_i$ be a real Hilbert space with inner product $g_i = \langle\cdot,\cdot\rangle_{\sH_i}$ for $i=1,2$ and $g_i$-orthogonal almost complex structures $J_i:\sH_i\to\sH_i$. Let $\sT:\sH_1\to\sH_2$ be an approximately complex operator as in Definition \ref{defn:ApproximatelyComplex} with complex linear and complex antilinear components $\sT'$ and $\sT''$. Let $\Pi_{\nu,i}' \in \End(\sH_i)$ be the orthogonal projections onto the subspaces $\bH_\nu^i(\sT')$ defined in Lemma \ref{lem:S1InvariantSubspacesAndProjections}. Then there is a positive constant $\nu\notin\sigma(\sT^{\prime,*}\sT')$ such that the following operators are injective:
\begin{equation}
\label{eq:LowEigenvalueEigenspaceProj_RestrictedTo_Ker_and_Coker}
  \Pi_{\nu,1}':\Ker\sT\to \bH_\nu^1(\sT')
  \quad\text{and}\quad
  \Pi_{\nu,2}':\Ker\sT^*\to \bH_\nu^2(\sT').
\end{equation}
\end{lem}

\begin{proof}
By Lemma \ref{lem:Eigenvalues_densely_defined_unbounded_linear_operators}, Items \eqref{item:Eigenvalues_densely_defined_unbounded_linear_operators_spectra_Laplacians} and \eqref{item:Eigenvalues_densely_defined_unbounded_linear_operators_TT*_from_T*T}, the spectra of $\sT^{\prime,*}\sT'$ and $\sT'\sT^{\prime,*}$, excluding the zero eigenvalues, are given by a non-decreasing sequence $\{\nu_k\}_{k=1}^\infty$ of positive eigenvalues with finite multiplicity. Note that the multiplicities of the zero eigenvalues of $\sT^{\prime,*}\sT'$ and $\sT'\sT^{\prime,*}$ will in general be different.

Let $\sH_{1,n}\subset\sH_1$ and $\sH_{2,m}\subset\sH_2$ be the subspaces given, respectively, by the span of the eigenvectors of $\sT^{\prime,*}\sT'$ and $\sT'\sT^{\prime,*}$  corresponding to the eigenvalues strictly less than $\nu$. Lemma \ref{lem:Approximation_finite-dimensional_subspaces_Hilbert_space} \eqref{eq:pi_n_V_embedding_sHn} implies that for $\nu$ and thus $n$ and $m$ sufficiently large, the orthogonal projections,
\begin{equation}
  \label{eq:nEigenvalues_eProj_RestrictedTo_Ker_and_Coker}
  \pi_{1,n}:\Ker\sT\to \sH_{1,n} \quad\text{and}\quad \pi_{2,m}:\Ker\sT^*\to \sH_{2,m},
\end{equation}
are injective. By definition in Lemma \ref{lem:S1InvariantSubspacesAndProjections} of the orthogonal projections $\Pi_{\nu,i}' \in \End(\sH_i)$ and the subspaces $\bH_\nu^i(\sT')$ for $i=1,2$, we have
\[
  \bH_{\nu}^1(\sT') = \sH_{1,n} \quad\text{and}\quad \bH_{\nu}^2(\sT') = \sH_{2,m},
\]
with $\pi_{1,n}=\Pi_{\nu,1}'$ and $\pi_{2,m}=\Pi_{\nu,2}'$. This completes the proof.
\end{proof}

\begin{lem}[Equivariance of complex linear and antilinear components]
\label{lem:S1EquivOfComponents}
For $k=1,2$, let $\sH_k$ be a real Hilbert space with inner product $g_k = \langle\cdot,\cdot\rangle_{\sH_k}$.
Let $G$ be a group and assume that for $k=1,2$ there are $g_k$-orthogonal $G$ actions on $\sH_k$. Assume that there $G$-invariant, Hermitian structures $J_k$ on $\sH_k$. Let $\sT:\sH_1\to \sH_2$ be a $G$-equivariant, approximately complex operator in the sense of Definition \ref{defn:ApproximatelyComplex}.
Then the complex linear and complex antilinear components $\sT'$ and $\sT''$ of $\sT$ (see \eqref{eq:Complex_linear_and_anti-linear_operator_components}) are also $G$-equivariant.
\end{lem}

\begin{proof}
The proof of Lemma \ref{lem:S1EquivOfComponents} follows immediately from the definition \eqref{eq:Complex_linear_and_anti-linear_operator_components} of $\sT'$ and $\sT''$ and the $G$-equivariance of $J_1$ and $J_2$.
\end{proof}

\begin{lem}[Equivariant isomorphism between complements of kernel and cokernel]
\label{lem:SymplecticStrutureOnStabilizedSpace}
Continue the notation and assumptions of Lemmas \ref{lem:InjectivityOfProjection_for_KurnaishiModel} and \ref{lem:S1EquivOfComponents} and let $\nu>0$ satisfy the conclusion of Lemma \ref{lem:InjectivityOfProjection_for_KurnaishiModel}. Define
\begin{subequations}
\label{eq:ImageOfKernels}
\begin{align}
\label{eq:ImageOfKernels1}
\tilde K_1(\sT) &:= \Pi_{\nu,1}'\left(\Ker\sT\right) \subset \bH_\nu^1(\sT'),
\\
\tilde K_2(\sT) &:= \Pi_{\nu,2}'\left(\Ker\sT^*\right) \subset \bH_\nu^2(\sT'),
\end{align}
\end{subequations}
where the projections $\Pi_{\nu,k}'$ are defined in Lemma \ref{lem:InjectivityOfProjection_for_KurnaishiModel}. Then the vector spaces $\tilde K_k(\sT)$ and $\tilde K_k(\sT)^\perp\cap \bH_\nu^k(\sT')$ are $G$-invariant for $k=1,2$. In addition, if $G$ is a compact Lie group, then there is a real linear isomorphism,
\begin{equation}
\label{eq:Equiv_Isom_Of_Complements}
\Gamma(\sT): \tilde K_1(\sT)^\perp\cap \bH_\nu^1(\sT') \to \tilde K_2(\sT)^\perp\cap \bH_\nu^2(\sT'),
\end{equation}
which is $G$-equivariant.
\end{lem}

\begin{proof}
Because $\sT$ is $G$-equivariant, $\Ker\sT$ and $\Ker\sT^*$ are $G$-invariant. For $k=1,2$, the $G$-equivariance of $\Pi_{\nu,k}'$ given in Lemma \ref{lem:S1InvariantSubspacesAndProjections} then implies that the images $\tilde K_k(\sT)$ under $\Pi_{\nu,k}'$ will also be $G$-invariant. The vector spaces $\bH_\nu^k(\sT')$ are also $G$-invariant for $k=1,2$ by Lemmas \ref{lem:S1EquivOfComponents} and \ref{lem:S1InvariantSubspacesAndProjections}, so by the $G$-invariance of $\tilde K_k(\sT)$ and the orthogonality of the $G$ action, the orthogonal complements $\tilde K_k(\sT)^\perp\cap H_\nu^k(\sT')$ are $G$-invariant.

We show that there is a $G$-equivariant isomorphism $\Gamma(\sT)$ as in \eqref{eq:Equiv_Isom_Of_Complements} by proving that
\begin{equation}
\label{eq:Equivariant_Equality_of_Complements}
[\tilde K_1(\sT)^\perp\cap \bH_\nu^1(\sT')]=[\tilde K_2(\sT)^\perp\cap \bH_\nu^2(\sT')]\in R_\RR(G)
\end{equation}
as follows.  First, observe that the direct sum decomposition gives
\begin{equation}
\label{eq:GammaIndexEqualityDirectSum1}
\begin{aligned}
[\bH_\nu^1(\sT')]
&=
[\tilde K_1(\sT)^\perp\cap \bH_\nu^1(\sT')]+[\tilde K_1(\sT)],
\\
[\bH_\nu^2(\sT')]
&=
[\tilde K_2(\sT)^\perp\cap \bH_\nu^2(\sT')]+[\tilde K_2(\sT)],
\end{aligned}
\end{equation}
and
\begin{equation}
\label{eq:GammaIndexEqualityDirectSum2}
\begin{aligned}
[\bH_\nu^1(\sT')]
&=
[(\Ker\sT')^\perp\cap \bH_\nu^1(\sT')]+[\Ker\sT'],
\\
[\bH_\nu^2(\sT')]
&=
[(\Ker\sT^{\prime,*})^\perp\cap \bH_\nu^2(\sT')]+[\Ker\sT^{\prime,*}].
\end{aligned}
\end{equation}
Because $\Pi_{\nu,k}'$ defines a $G$-equivariant isomorphism, we have the equalities
\begin{equation}
\label{eq:GammaIndexEqualityProjectionEquality}
[\Ker\sT]=[\tilde K_1(\sT)]\in R_\RR(G),
\quad
[\Ker\sT^*]=[\tilde K_2(\sT)]\in R_\RR(G).
\end{equation}
Because $\sT''$ is compact, Proposition \ref{prop:HomotopyInvarianceOfGindex} implies that
\begin{equation}
\label{eq:GammaIndexEquality_PerturbationEquality}
\Ind_{G,\RR}(\sT)=\Ind_{G,\RR}(\sT').
\end{equation}
Finally, the $G$-equivariant operator $\sT'$ gives a $G$-equivariant isomorphism
\begin{equation}
\label{eq:sT'_bounded_isom}
\sT':\Ker(\sT')^\perp\cap \bH_\nu^1(\sT') \to \Ker(\sT^{\prime,*})^\perp\cap \bH_\nu^2(\sT'),
\end{equation}
and thus an equality
\begin{equation}
\label{eq:NonZeroEigenspacesIsomorphism}
[\Ker(\sT')^\perp\cap \bH_\nu^1(\sT')]=[\Ker(\sT^{\prime,*})^\perp\cap \bH_\nu^2(\sT')] \in R_\RR(G).
\end{equation}
We now compute in $R_\RR(G)$,
\begin{align*}
&[\tilde K_1(\sT)^\perp\cap \bH_\nu^1(\sT')] - [\tilde K_2(\sT)^\perp\cap \bH_\nu^2(\sT')]
\\
&\quad=
\left([\bH_\nu^1(\sT')]-[\tilde K_1(\sT)]\right)
-
\left([\bH_\nu^2(\sT')]-[\tilde K_2(\sT)] \right)
\quad\text{(by \eqref{eq:GammaIndexEqualityDirectSum1})}
\\
&\quad=
\left([\bH_\nu^1(\sT')]-[\bH_\nu^2(\sT')]\right)
-
\left( [\Ker \sT] -[\Ker\sT^*]\right)
\quad\text{(by \eqref{eq:GammaIndexEqualityProjectionEquality})}
\\
&\quad=
\left([\bH_\nu^1(\sT')]-[\bH_\nu^2(\sT')]\right)
-
\Ind_{G,\RR}(\sT)
\quad\text{(by Definition \ref{defn:Index_of_GEquivariantDiffOperator})}
\\
&\quad=
\left([\bH_\nu^1(\sT')]-[\bH_\nu^2(\sT')]\right)
-
\Ind_{G,\RR}(\sT')
\quad\text{(by \eqref{eq:GammaIndexEquality_PerturbationEquality})}
\\
&\quad=
\left([\bH_\nu^1(\sT')]-[\bH_\nu^2(\sT')]\right)
-
\left( [\Ker \sT'] -[\Ker\sT^{\prime,*}]\right)
\quad\text{(by Definition \ref{defn:Index_of_GEquivariantDiffOperator})}
\\
&\quad=
\left([\bH_\nu^1(\sT')]-[\Ker \sT']\right)
-
\left([\bH_\nu^2(\sT')]  -[\Ker\sT^{\prime,*}]\right)
\\
&\quad=
     [\Ker(\sT')^\perp\cap \bH_\nu^1(\sT')]-[\Ker(\sT^{\prime,*})^\perp\cap \bH_\nu^2(\sT')]
\quad\text{(by \eqref{eq:GammaIndexEqualityDirectSum2})}
\\
&\quad=0 \quad\text{(by \eqref{eq:NonZeroEigenspacesIsomorphism}).}
\end{align*}
This completes the proof of the equality \eqref{eq:Equivariant_Equality_of_Complements} and hence the existence of the $G$-equivariant isomorphism $\Gamma(\sT)$ in \eqref{eq:Equiv_Isom_Of_Complements}.
\end{proof}

Lemma  \ref{lem:InjectivityOfProjection_for_KurnaishiModel} gives the following isomorphisms which will enable us to apply Corollary \ref{cor:Kuranishi_model_defined_by_Fredholm_map_Hilbert_spaces} to the vector spaces $\bH_\nu^k(\sT_{A,\varphi,\psi,r}')$ and thus use them to define local Kuranishi models.

\begin{cor}[Isomorphism of bounded-eigenvalue eigenspaces]
\label{cor:IsomOfKuranishiModelSpaces}
Let $\sT:\sH_1\to\sH_2$ be a $G$-equivariant, approximately complex operator.  For $\nu>0$ satisfying the conclusion of Lemma \ref{lem:InjectivityOfProjection_for_KurnaishiModel}, define
\begin{equation}
\label{eq:Define_Kuranishi_Spaces_for_AC_Operator}
\begin{aligned}
\mathbf{\tilde H}_\nu^1(\sT):= \Ker \sT \oplus \left( \tilde K_1(\sT)^\perp\cap \bH_\nu^1(\sT')\right)\subset\sH_1,
\\
\mathbf{\tilde H}_\nu^2(\sT):= \Ker \sT^* \oplus \left( \tilde K_2(\sT)^\perp\cap \bH_\nu^2(\sT')\right)\subset\sH_2,
\end{aligned}
\end{equation}
where for $k=1,2$, $\tilde K_k(\sT)\subset \bH_\nu^k(\sT')$ are defined in \eqref{eq:ImageOfKernels}. Then the vector spaces $\mathbf{\tilde H}_\nu^k(\sT)$ are $G$-invariant and the orthogonal projections $\Pi_{\nu,k}'$ from $\sH_k$ onto $\bH_\nu^k(\sT')$ of Lemma \ref{lem:S1InvariantSubspacesAndProjections} define $G$-equivariant isomorphisms for $k=1,2$:
\begin{equation}
\label{eq:Kuranishi_Model_Space_Isom}
\Pi_{\nu,k}': \mathbf{\tilde H}_\nu^k(\sT) \to   \bH_\nu^k(\sT'),
\end{equation}
\end{cor}

Because the operator $\sT'$ is complex linear, the bounded-eigenvalue eigenspaces $\bH_\nu^k(\sT')$ for $k=1,2$, are $J_i$-invariant spaces. Lemmas \ref{lem:S1InvariantSubspacesAndProjections} and \ref{lem:S1EquivOfComponents} imply that for $k=1,2$, the $G$ actions on $\bH_\nu^k(\sT')$ give complex $G$ representations with respect to $J_k$. This enables us to make the following

\begin{defn}[Complex equivariant index for equivariant, approximately complex operators]
\label{defn:ApproximatelyComplexIndex}
Continue the notation of Definition \ref{defn:ApproximatelyComplex}.
Let $\sT:\sH_1\to\sH_2$ be a $G$-equivariant, approximately complex operator and let $\mathbf{\tilde H}_\nu^k(\sT)$ be the vector space defined in \eqref{eq:Define_Kuranishi_Spaces_for_AC_Operator}. For $\nu>0$ satisfying the conclusion of Lemma \ref{lem:InjectivityOfProjection_for_KurnaishiModel} and for $k=1,2$, let
\begin{equation}
\label{eq:PullbackACStructure}
J_k':\mathbf{\tilde H}_\nu^k(\sT)\to \mathbf{\tilde H}_\nu^k(\sT)
\end{equation}
be the almost complex structure defined by pulling back the almost complex structure $J_k$ on $H_\nu^k(\sT')$ by the isomorphism \eqref{eq:Kuranishi_Model_Space_Isom}. For $k=1,2$, define the \emph{complex $G$-equivariant index $\nu$-index of $\sT$} by
\begin{equation}
\label{eq:DefineComplexIndexForS1Operators}
\Ind_{G,\CC}(\sT):=[\mathbf{\tilde H}_\nu^1(\sT),J_1'] - [\mathbf{\tilde H}_\nu^2(\sT),J_2']\in R_\CC(G),
\end{equation}
where $[\mathbf{\tilde H}_\nu^k(\sT),J_k']$ denotes the element of $R_\CC(G)$ given by Definition \ref{defn:GInvarComplexStructures}.
\qed\end{defn}

\begin{rmk}[Comparison of Definitions \ref{eq:GEquivIndex} and \ref{defn:ApproximatelyComplexIndex}]
\label{rmk:ComparisonOfDefinitionOfComplexGIndex}
The complex $G$-equivariant index of Definition \ref{eq:GEquivIndex} applies only to \emph{complex} equivariant operators while Definition \ref{defn:ApproximatelyComplexIndex} applies to approximately complex equivariant operators. If $\sT=\sT'$ is a complex $G$-equivariant operator, then the isomorphism given by $\sT'$ between the eigenspaces of $\sT^{\prime,*}\sT'$ and those of $\sT'\sT^{\prime,*}$ corresponding to non-zero eigenvalues implies that Definitions \ref{eq:GEquivIndex} and \ref{defn:ApproximatelyComplexIndex} are equal and thus the expression \eqref{eq:DefineComplexIndexForS1Operators} is independent of $\nu$.  Proposition \ref{prop:Equality_of_ComplexEquivariantIndices} below then implies that Definition \ref{defn:ApproximatelyComplexIndex} is also independent of $\nu$.
\end{rmk}

We will use the following reformulation of the definition \eqref{eq:DefineComplexIndexForS1Operators} to compute $\Ind_{G,\CC}(\sT)$ for approximately complex operators $\sT$.

\begin{prop}[Equality between complex equivariant indices of $\sT$ and $\sT'$]
\label{prop:Equality_of_ComplexEquivariantIndices}
Continue the notation and assumptions of Lemma \ref{lem:InjectivityOfProjection_for_KurnaishiModel} and let $\nu>0$ satisfy the conclusion of Lemma \ref{lem:InjectivityOfProjection_for_KurnaishiModel}. Then
\begin{equation}
\label{eq:ComplexS1IndexEquality}
\Ind_{G,\CC}(\sT)
=
\Ind_{G,\CC}(\sT')\in R_\CC(G),
\end{equation}
where $\Ind_{G,\CC}(\sT)$ is the complex $G$-equivariant index defined in \eqref{eq:DefineComplexIndexForS1Operators} and $\Ind_{G,\CC}(\sT')$ is the complex $G$-equivariant index
defined in \eqref{eq:AdditivityOfEquivariantIndex} of the complex linear operator $\cT'$ .
\end{prop}

\begin{proof}
The $G$-equivariant isomorphisms \eqref{eq:Kuranishi_Model_Space_Isom} imply that
\[
\Ind_{G,\CC}(\sT)
=
[\mathbf{\tilde H}_\nu^1(\sT)] - [\mathbf{\tilde H}_\nu^1(\sT)]
=
[\bH_\nu^1(\sT')]-[\bH_\nu^2(\sT')].
\]
Note that the second equality holds in $R_\CC(G)$ because the almost complex structure on $\mathbf{\tilde H}_\nu^k(\sT)$ is defined (in Definition \ref{defn:ApproximatelyComplexIndex}) to give this equality. The conclusion of Proposition \ref{prop:Equality_of_ComplexEquivariantIndices} now follows immediately from the observation that the equalities \eqref{eq:GammaIndexEqualityDirectSum2} and \eqref{eq:NonZeroEigenspacesIsomorphism} hold in $R_\CC(G)$ as well as in $R_\RR(G)$.
\end{proof}

\section{Comparison of almost complex structures on the vector spaces $\mathbf{\tilde H}_\nu^k(\sT)$}
\label{sec:AC_Structure_Comparison}
The complex equivariant index  given in Definition \ref{defn:ApproximatelyComplexIndex} uses the almost complex structure $J_k'$ on $\mathbf{\tilde H}_\nu^k(\sT)$ given in \eqref{eq:PullbackACStructure} by the isomorphisms \eqref{eq:Kuranishi_Model_Space_Isom} between $\mathbf{\tilde H}_\nu^k(\sT)$ and the complex vector spaces  $\bH_\nu^k(\sT')$. Following the construction in Corollary  \ref{cor:Almost_complex_structures_on_finite-dimensional_subspaces_Hilbert_space}, we define a different almost complex structure $\tilde J_k$ in \eqref{eq:AConTildeH} below. Although the definition of $J_k'$ in \eqref{eq:PullbackACStructure} was useful in Proposition \ref{prop:Equality_of_ComplexEquivariantIndices}, the construction of $\tilde J_k$ in \eqref{eq:AConTildeH} was more convenient in Chapter \ref{chap:Construction_circle-invariant_non-degenerate_two-form_II}. In this section, we use Proposition \ref{prop:Homotopy_Invariance_of_G_Representations} to prove that these two choices of almost complex structures on $\mathbf{\tilde H}_\nu^k(\sT)$  give the same element of $R_\CC(G)$ in Definition \ref{defn:ApproximatelyComplexIndex}, as described in the following

\begin{prop}[Equality of complex equivariant index for different almost complex structures]
\label{prop:SameComplexIndex}
Continue the notation and hypotheses of Lemma \ref{lem:InjectivityOfProjection_for_KurnaishiModel} and of Corollary \ref{cor:IsomOfKuranishiModelSpaces}. Let $J_k'$ be the almost complex structure on $\mathbf{\tilde H}_\nu^k(\sT)$ given in \eqref{eq:PullbackACStructure}. For $\nu>0$ sufficiently large, the following hold:
\begin{enumerate}
\item
\label{item:SameComplexIndex_ExistenceOfAC}
For $k=1,2$, let $\tilde\pi_k$ be the orthogonal projection from $\sH_k$ onto $\mathbf{\tilde H}_\nu^k(\sT)$.  Then
\[
  A_{\tilde H,k}:=\tilde\pi_k J_k\in \End(\mathbf{\tilde H}_\nu^k(\sT))
\]
is a $G$-equivariant invertible operator and skew-adjoint with respect to the inner product induced by $\langle\cdot,\cdot\rangle_{\sH_k}$. In addition, the operator
\begin{equation}
\label{eq:AConTildeH}
\tilde J_k:=(A_{\tilde H,k}^*A_{\tilde H,k})^{-1}A_{\tilde H,k},
\end{equation}
is a $G$-invariant almost complex structure on $\mathbf{\tilde H}_\nu^k(\sT)$.
\item
\label{eq:SameComplexIndex}
The almost complex structures $J_k'$ and $\tilde J_k$ satisfy
\[
[\mathbf{\tilde H}_\nu^1(\sT), J_1'] - [\mathbf{\tilde H}_\nu^2(\sT), J_2']
=
[\mathbf{\tilde H}_\nu^1(\sT),\tilde J_1] - [\mathbf{\tilde H}_\nu^2(\sT),\tilde J_2],
\]
as elements of $R_\CC(G)$, where we have used the convention of Definition \ref{defn:GInvarComplexStructures} for elements $[V,J]$ of a complex representation ring.
\end{enumerate}
\end{prop}

We will prove Proposition \ref{prop:SameComplexIndex} at the end of Section \ref{sec:AC_Structure_Comparison} as a consequence of the forthcoming Lemmas \ref{lem:GlobalDifference_Projections}, \ref{lem:Existence_of_AC_Structure}, and \ref{prop:SimilarProj}.

In the verification of the forthcoming Lemma \ref{lem:GlobalDifference_Projections}, we adapt the proofs of Lemma \ref{lem:Approximation_finite-dimensional_subspaces_Hilbert_space} and Corollary \ref{cor:Approximation_projections_onto_finite-dimensional_subspaces_Hilbert_space} to show that
orthogonal projection operators satisfying Conclusion \eqref{item:pi_minus_pi_n_hom(V,sH)_and_pi_minus_pi_n_hom(Vn,sH)}
of Lemma \ref{lem:Approximation_finite-dimensional_subspaces_Hilbert_space} satisfy the remaining conclusions of Lemma \ref{lem:Approximation_finite-dimensional_subspaces_Hilbert_space} and Corollary \ref{cor:Approximation_projections_onto_finite-dimensional_subspaces_Hilbert_space}.

\begin{lem}[Isomorphism between ranges of orthogonal projections]
\label{lem:GlobalDifference_Projections}
Let $V$ and $V'$ be linear subspaces of a Hilbert space $\sH$. If the orthogonal projections $\pi$ and $\pi'$ onto $V$ and $V'$, respectively, satisfy
\begin{subequations}
  \label{eq:Projection_Difference_Hypothesis}
  \begin{align}
    \label{eq:Projection_Difference_Hypothesis_V}
    \|\pi-\pi'\|_{\Hom(V,\sH)} &< \eps,
    \\
    \label{eq:Projection_Difference_Hypothesis_Vprime}
    \|\pi-\pi'\|_{\Hom(V',\sH)} &< \eps',
  \end{align}
\end{subequations}
for $\eps,\eps'>0$, then the following hold:
\begin{enumerate}
\item
\label{item:Invertibility_of_ProjectionComp_Equidim1}
 If $\eps<1$, then $\pi\pi'\in\GL(V)$ and $\pi:V'\to V$ is surjective.
\item
\label{item:Invertibility_of_ProjectionComp_Equidim2}
If $\eps'<1$, then $\pi'\pi\in \GL(V')$ and $\pi':V\to V'$ is surjective.
\end{enumerate}
Moreover, if $\max\{\eps,\eps'\}<1$, then
\begin{equation}
\label{eq:Difference_of_Projections_On_HilbertSpace}
\|\pi -\pi'\|_{\End(\sH)}< \eps.
\end{equation}
\end{lem}

\begin{proof}
As in the proof of \eqref{eq:Norm_id-pi_circ_pi_n_and_id-pi_n_circ_pi}, we observe that
\begin{align*}
\|1 - \pi\pi'\|_{\End(V)}
=
\|\pi(\pi-\pi')\|_{\End(V)}
\le
\|\pi\|_{\Hom(\sH,V)}\| \pi-\pi'\|_{\Hom(V,\sH)}
=
\| \pi-\pi'\|_{\Hom(V,\sH)}.
\end{align*}
Hence if $\eps<1$, then $\|1 - \pi\pi'\|_{\End(V)}<1$ which implies that $\pi\pi'\in\GL(V)$ by Rudin \cite[Theorem 10.7, p. 249]{Rudin}. We show that $\pi:V'\to V$ is surjective by observing that for any $v\in V$ we can define $v':=\pi'(\pi\pi')^{-1}v\in V'$ so that $\pi v'=v$. This completes the proof of Conclusion \eqref{item:Invertibility_of_ProjectionComp_Equidim1}. Conclusion \eqref{item:Invertibility_of_ProjectionComp_Equidim2} follows by the same argument.

To derive the inequality \eqref{eq:Difference_of_Projections_On_HilbertSpace}, we observe that
\[
\|(1-\pi')\pi\|_{\End(\sH)}
\le
\| 1-\pi'\|_{\Hom(V,\sH)}\|\pi\|_{\Hom(\sH,V)}
=
\| \pi-\pi'\|_{\Hom(V,\sH)},
\]
where the last equality follows from the facts that $\pi$ is the identify on $V$ and $\|\pi\|_{\Hom(\sH,V)} = 1$. Hence,
\begin{equation}
\label{eq:EquiDimProjection_Bound3a}
\|(1-\pi')\pi\|_{\End(\sH)}\le \| \pi-\pi'\|_{\Hom(V,\sH)}.
\end{equation}
The hypothesis \eqref{eq:Projection_Difference_Hypothesis_V} and \eqref{eq:EquiDimProjection_Bound3a} imply that
\begin{equation}
  \label{eq:EquiDimProjection_Bound3}
  \|(1-\pi')\pi\|_{\End(\sH)} < \eps.
\end{equation}
Thus, if $\max\{\eps,\eps'\}<1$, the surjectivity of $\pi':V\to V'$ given by Conclusion \eqref{item:Invertibility_of_ProjectionComp_Equidim2} and the inequality \eqref{eq:EquiDimProjection_Bound3} imply that $\pi$ and $\pi'$ obey the hypotheses of \cite[Section I.6.8, Theorem 6.34 (i), p. 56]{Kato} and, in particular, satisfy \cite[Section I.6.8, Equation (6.51), p. 57]{Kato},
\[
  \|\pi-\pi\|_{\End(\sH)}=\|(1-\pi')\pi\|_{\End(\sH)}.
\]
The preceding equality and \eqref{eq:EquiDimProjection_Bound3} imply that $\pi$ and $\pi'$ satisfy the inequality \eqref{eq:Difference_of_Projections_On_HilbertSpace}. This completes the proof of Lemma \ref{lem:GlobalDifference_Projections}.
\end{proof}

We now show that an almost complex structure constructed by the method of Corollary \ref{cor:Almost_complex_structures_on_finite-dimensional_subspaces_Hilbert_space} is $G$-equivariant.

\begin{lem}[Construction of an almost complex structure by orthogonal projection]
\label{lem:Existence_of_AC_Structure}
Let $G$ be a group that acts orthogonally on a real Hilbert space $(\sH,\langle\cdot,\cdot\rangle_{\sH})$. Let $J$ be a $G$-invariant (in the sense of Definition \ref{defn:GInvarComplexStructures}) orthogonal almost complex structure on $\sH$ and $V$ and $V'$ be $G$-invariant linear subspaces of $\sH$ and assume that $V'$ is $J$-invariant. If $\pi$ and $\pi'$ are the orthogonal projections onto $V$ and $V'$, respectively, then the following hold:
\begin{enumerate}
\item
\label{item:Existence_of_AC_Structure_EquivariantProjections}
The orthogonal projections $\pi$ and $\pi'$ are $G$-equivariant.
\item
\label{item:Existence_of_AC_Structure_EquivariantA_V}
The bounded linear operator
\begin{equation}
  \label{eq:A_V}
  A_V:=\pi J\in\End(V)
\end{equation}
is $G$-equivariant and skew-adjoint with respect to the inner product $\langle\cdot,\cdot\rangle_V$ on $V$ induced by $\langle\cdot,\cdot\rangle_{\sH}$.
\item
\label{item:Existence_of_AC_Structure_InvertibleA_V}
If $\pi$ and $\pi'$ satisfy
\begin{equation}
  \label{eq:DifferenceHypothesis_of_Projections_On_HilbertSpace}
  \|\pi -\pi'\|_{\End(\sH)}<\frac{1}{2},
\end{equation}
then $A_V$ is invertible and
\begin{equation}
\label{eq:Projection_ac}
\tilde J:= (A_V^*A_V)^{-1/2}A_V
\end{equation}
is a $G$-invariant almost complex structure on $V$.
\end{enumerate}
\end{lem}

\begin{proof}
All the conclusions of Lemma \ref{lem:Existence_of_AC_Structure}, except for the statements on $G$-equivariance, follow from the proofs of the analogous statements in Corollary \ref{cor:Almost_complex_structures_on_finite-dimensional_subspaces_Hilbert_space}. Specifically, in the proof of Corollary \ref{cor:Almost_complex_structures_on_finite-dimensional_subspaces_Hilbert_space}, we replace the
\begin{inparaenum}[\itshape a\upshape)] 
\item vector spaces $\tilde\sH_n$ and $\sH_n$ by $V$ and $V'$ and the
\item orthogonal projections $\tilde\pi_n$ and $\pi_n$ by the orthogonal projections $\pi$ and $\pi'$, respectively.
\end{inparaenum}
The skew-adjointness and invertibility of $A_V\in\End(V)$ follow from the proof of those properties for $A_n\in\End(\tilde\sH_n)$. Similarly, the proof that the endomorphism $\tilde J\in \End(V)$ in \eqref{eq:Projection_ac} is an almost complex structure follows from the proof for $J_n\in\End(\tilde\sH_n)$ in \eqref{eq:AC_Structure_On_tildeH_Projection}.

We now prove the statements on $G$-equivariance. To prove that $\pi$ and $\pi'$ are $G$-equivariant as asserted in Conclusion \eqref{item:Existence_of_AC_Structure_EquivariantProjections}, we observe that any $w\in\sH$ can be written uniquely as $w=v+v^\perp$, where $v=\pi w\in V$ and $v^\perp=(1-\pi)w\in V^\perp$.  Because $V$ is $G$-invariant and the $G$ action on $\sH$ is orthogonal, $V^\perp$ is also $G$-invariant. Thus, for all $g\in G$, we have $gw=gv+gv^\perp$, where $gv\in V$ and $gv^\perp\in V^\perp$.  The uniqueness of the summands in the equality $gw=\pi gw+(1-\pi)gw$ implies that $gv=g\pi w=\pi gw$. Therefore, $\pi$ is $G$-equivariant and the same argument shows that $\pi'$ is $G$-equivariant.

As noted in Definition \ref{defn:GInvarComplexStructures}, a $G$-invariant almost complex structure is a $G$-equivariant map. Hence, the $G$-invariance of $J$ and the $G$-equivariance of $\pi$ imply that the composition $A_V=\pi J$ is a composition of $G$-equivariant maps and hence $G$-equivariant. This completes the proof 
of Conclusion \eqref{item:Existence_of_AC_Structure_EquivariantA_V}.

The $G$-equivariance of  the almost complex structure $\tilde J$ in \eqref{eq:Projection_ac} will follow from the $G$-equivariance of $A_V$ and that of the operator $(A_V^*A_V)^{-1/2}$. Because $G$ acts orthogonally on $\sH$ and $A_V$ is $G$-equivariant, $A_V^*$ and hence $A_V^*A_V$ are $G$-equivariant. The Dunford integral expression \eqref{eq:Dunford_integral_bounded_operator} for the inverse of the square root of $A_V^*A_V$ is
\begin{equation}
\label{eq:DunfordIntegral_forSqrt_A_V}
(A_V^*A_V)^{-1/2}=-\frac{1}{2\pi i}\oint_\Gamma (z)^{-1/2}(A_V^*A_V-z)^{-1}\ dz,
\end{equation}
where $\Gamma \subset \rho(A_V^*A_V)$ is a positively oriented, rectifiable curve that bounds an open neighborhood of $\sigma(A_V^*A_V)$. Because $A_V^*A_V$ is $G$-equivariant, the operators $A_V^*A_V-z$ and hence $(A_V^*A_V-z)^{-1}$ are $G$-equivariant.  The equality \eqref{eq:DunfordIntegral_forSqrt} implies that $(A_V^*A_V)^{-1/2}$ is $G$-equivariant and so $\tilde J$ is as well by the definition \eqref{eq:Projection_ac}.  This completes the proof of Conclusion \eqref{item:Existence_of_AC_Structure_InvertibleA_V} and hence of Lemma \ref{lem:Existence_of_AC_Structure}.
\end{proof}

We now construct a homotopy of orthogonal projections between two sufficiently close orthogonal projections on a Hilbert space.

\begin{prop}[Criterion for similarity of orthogonal projection operators and for homotopy between orthogonal projection operators]
\label{prop:SimilarProj}
Let $\pi$ and $\pi'$ be orthogonal projection operators on a real Hilbert space $\sH$ satisfying
\begin{equation}
\label{eq:DifferenceHypothesis_of_Projections_On_HilbertSpace_1/4}
\|\pi-\pi'\|_{\End(\sH)}<\frac{1}{4}.
\end{equation}
Then $\pi$ and $\pi'$ are equal modulo conjugation by an orthogonal operator $U\in\GL(\sH)$,
\begin{equation}
  \label{eq:ProjectionSimilarity}
  \pi'=U^{-1}\pi U.
\end{equation}
In addition, there is a real analytic map $[0,1]\ni t\to \pi_t\in\End(\sH)$ such that the following hold:
\begin{enumerate}
\item
\label{item:SimilarProj_HomotopyProj}
$\pi_t$ is an orthogonal projection, for all $t\in [0,1]$.
\item
\label{item:SimilarProj_HomotopyBetween}
$\pi_0=\pi$ and $\pi_1=\pi'$.
\item
\label{item:SimilarProj_HomotopyProjBounds}
$\|\pi_t-\pi'\|_{\End(\sH)}<1/2$ and $\|\pi_t-\pi\|_{\End(\sH)}<1/2$, for all $t \in [0,1]$.
\end{enumerate}
If $G$ is a group that acts orthogonally on $\sH$ and $\pi$ and $\pi'$ are $G$-equivariant, then the operators $U$ and $\pi_t$ are $G$-equivariant for all $t\in [0,1]$.
\end{prop}

\begin{proof}
To prove the proposition, we adapt the discussion in Kato \cite[Section I.4.6, pp. 32--34]{Kato} and \cite[Section I.6.8, Theorem 6.32, p. 56]{Kato}. Denote $R:=(\pi-\pi')^2\in\End(\sH)$, so that
\begin{equation}
  \label{eq:ExpandR}
  R=\pi-\pi\pi'-\pi'\pi+\pi'.
\end{equation}
If $\pi$ and $\pi'$ are $G$-equivariant, then $R$ is as well. Kato observes, following \cite[Section I.4.6, Equation (4.33), p. 33]{Kato}, that
\begin{align*}
  \pi R &=\pi(\pi - \pi\pi'-\pi'\pi+\pi')=\pi-\pi\pi'\pi,
  \\
  R\pi &=(\pi-\pi\pi'-\pi'\pi+\pi')\pi=\pi-\pi\pi'\pi,
\end{align*}
and thus $\pi$ commutes with $R$.  Similarly, because
\begin{align*}
  \pi'R&=\pi'(\pi - \pi\pi'-\pi'\pi+\pi')=\pi'-\pi'\pi\pi',
  \\
  R\pi'&=(\pi-\pi\pi'-\pi'\pi+\pi')\pi'=\pi'-\pi'\pi\pi',
\end{align*}
and thus $\pi'$ commutes with $R$. Because the orthogonal projections $\pi$ and $\pi'$ are self-adjoint, the operator $\pi-\pi'$ is self-adjoint and thus normal. By \eqref{eq:Norm_of_Powers}, we have
\[
  \|R\|_{\End(\sH)}
 =
  \|(\pi-\pi')^2\|_{\End(\sH)}
  =
  \|\pi-\pi'\|_{\End(\sH)}^2,
\]
that is,
\begin{equation}
\label{eq:RNormBound}
\|R\|_{\End(\sH)} = \|\pi-\pi'\|_{\End(\sH)}^2.
\end{equation}
Hence, the hypothesis \eqref{eq:DifferenceHypothesis_of_Projections_On_HilbertSpace_1/4} on $\pi$ and $\pi'$ implies that
\begin{equation}
  \label{eq:Norm_R_lessthan_1}
  \|R\|_{\End(\sH)}<1.
\end{equation}  
Consequently, using the abbreviation $1 = \id_\sH$, the inverse $(1-R)^{-1} \in \End(\sH)$ of $1-R \in \End(\sH)$ exists by Rudin \cite[Theorem 10.7, p. 249]{Rudin}. Moreover, $1-R$ is a positive operator because, for all $x\in \sH$,
\[
  \langle (1-R)x, x\rangle_\sH = \|x\|_\sH^2 - \langle Rx, x\rangle_\sH \geq \|x\|_\sH^2 - \|Rx\|_\sH\|x\|_\sH
  \geq \|x\|_\sH^2 - \|R\|_{\End(\sH)}\|x\|_\sH^2 > 0,
\]
and thus
\begin{equation}
  \label{eq:PositivityOfR}
  1 - R > 0.
\end{equation}  
Therefore, the operator $(1-R)^{-1}$ is positive too and has a unique square root $(1-R)^{-1/2}$
by Rudin \cite[Theorem 12.33, p. 331]{Rudin}. The Dunford integral \eqref{eq:Dunford_integral_bounded_operator} gives an explicit expression for this square root,
\begin{equation}
\label{eq:DunfordIntegral_forSqrt}
(1-R)^{-1/2}=-\frac{1}{2\pi i}\oint_\Gamma (1-z)^{-1/2}(R-z)^{-1}\ dz,
\end{equation}
where $\Gamma \subset \rho(R)$ is a positively oriented, rectifiable curve that bounds an open neighborhood of $\sigma(R)$.

Since $R$ commutes with $\pi$ and $\pi'$, then so do $R-z$ and $(R-z)^{-1}$ for any $z\in \CC$, and thus
$(1-R)^{-1/2}$ does as well by the expression \eqref{eq:DunfordIntegral_forSqrt}. Because $R$ is $G$-equivariant, $R-z$ and thus $(R-z)^{-1}$ are also $G$-equivariant for any $z\in\CC$. The expression \eqref{eq:DunfordIntegral_forSqrt} implies that $(1-R)^{-1/2}$ is $G$-equivariant.

To prove the similarity relation \eqref{eq:ProjectionSimilarity}, Kato \cite[Section I.4.6, Equation (4.38), p. 33]{Kato} defines bounded operators in $\End(\sH)$ by
\begin{align*}
U:=\left(\pi\pi'+(1-\pi)(1-\pi')\right)(1-R)^{-1/2},
\\
V:=\left(\pi'\pi+(1-\pi')(1-\pi)\right)(1-R)^{-1/2},
\end{align*}
and shows that they satisfy $UV=VU=1$ in \cite[Section I.4.6, Equation (4.41), p. 34]{Kato}. In addition, Kato shows that the identity \eqref{eq:ProjectionSimilarity} holds in \cite[Section I.4.6, Equation (4.42), p. 34]{Kato}. When $\pi$ and $\pi'$ are orthogonal projections (and not just idempotents), Kato demonstrates the orthogonality of $U$ prior to \cite[Section I.6.8, Theorem 6.32, p. 56]{Kato}. The $G$-equivariance of $U$ and $V$ follows from their definition and the $G$-equivariance of $\pi$, $\pi'$, and $(1-R)^{-1/2}$.

In \cite[Section I.4.6, Problem 4.13, p. 34]{Kato}, Kato defines the analytic family of bounded operators in $\End(\sH)$,
\begin{equation}
  \label{eq:FamilyOfOperators}
  \pi_t:=\frac{1}{2} + \frac{1}{2}\left(2\pi-1+2t(\pi'-\pi)\right)\left( 1 -4t(1-t)R\right)^{-1/2},
  \quad\text{for } t \in [0,1],
\end{equation}
which satisfy $\pi_0=\pi$ and $\pi_1=\pi'$ and which he asserts are orthogonal projections.  Proving these assertions will complete the proof of Conclusions \eqref{item:SimilarProj_HomotopyProj} and \eqref{item:SimilarProj_HomotopyBetween}.
The $G$-equivariance of $\pi$, $\pi'$, and $\left( 1 -4t(1-t)R\right)^{-1/2}$ imply that $\pi_t$ is $G$-equivariant.

We prove that $\pi_t$ is indeed an orthogonal projection by showing that it is an idempotent and self-adjoint. To simplify the proof that $\pi_t$ is an idempotent, $\pi_t^2=\pi_t$, we abbreviate
\begin{equation}
  \label{eq:Define_Q_t_and_W_t}
  Q_t := 4t(1-t)R
  \quad\text{and}\quad
  W_t := (1 - Q_t)^{-1/2}.
\end{equation}
Observe that for $t\in [0,1]$, the factor $4t(1-t)$ has maximum value $1$ at $t=1/2$. Therefore, by \eqref{eq:Norm_R_lessthan_1} the operator $Q_t$ in \eqref{eq:Define_Q_t_and_W_t} obeys
\begin{equation}
  \label{eq:Norm_Q_t_lessthan_1}
  \|Q_t\|_{\End(\sH)} \leq \|R\|_{\End(\sH)} < 1, \quad\text{for all } t\in [0,1].
\end{equation}  
Because $R$ commutes with $\pi$ and $\pi'$ and $W_t = (1 - Q_t)^{-1/2}$ can be defined by the Dunford integral  \eqref{eq:DunfordIntegral_forSqrt}, the operator $W_t$ will also commute with $\pi$ and $\pi'$. We further note that because $\pi$ is an idempotent, we have $(2\pi-1)^2=4\pi^2-4\pi+1=1$, that is,
\begin{equation}
\label{eq:ReflectionSquared=1}
(2\pi-1)^2=1.
\end{equation}
We compute
\begin{align*}
(2\pi_t)^2
&=
\left( 1 +\left(2\pi-1+2t(\pi'-\pi)\right)W_t\right)^2
\\
&=
1+2\left(2\pi-1+2t(\pi'-\pi)\right)W_t
+\left(2\pi-1+2t(\pi'-\pi)\right)^2W_t^2
\\
&=
1+2\left(2\pi-1+2t(\pi'-\pi)\right)W_t
\\
  &\qquad+\left( (2\pi-1)^2 +2t(2\pi-1)(\pi'-\pi) +
    2t(\pi'-\pi)(2\pi-1)+4t^2 (\pi'-\pi)^2\right)W_t^2.
\end{align*}
Using the definition $R = (\pi'-\pi)^2$ and applying \eqref{eq:ReflectionSquared=1} to the above equalities yields
\begin{align*}
(2\pi_t)^2
&=
1+2\left(2\pi-1+2t(\pi'-\pi)\right)W_t
\\
&\qquad+\left( 1 +2t(2\pi\pi'-2\pi-\pi'+\pi) +2t(2\pi'\pi-\pi'-2\pi+\pi)+4t^2 R\right)W_t^2
\\
&=
1+2\left(2\pi-1+2t(\pi'-\pi)\right)W_t
\\
&\qquad+\left( 1 -4t(\pi-\pi\pi'-\pi'\pi+\pi')+4t^2 R\right)W_t^2
\\
&=
1+2\left(2\pi-1+2t(\pi'-\pi)\right)W_t
\\
  &\qquad+\left( 1 -4tR+4t^2 R\right)W_t^2
    \quad\text{(by \eqref{eq:ExpandR})}
\\
&=
     1+2\left(2\pi-1+2t(\pi'-\pi)\right)W_t +W_t^{-2}W_t^2
     \quad\text{(by \eqref{eq:Define_Q_t_and_W_t})}
\\
&=
4\left(\frac{1}{2}1+\frac{1}{2}\left(2\pi-1+2t(\pi'-\pi)\right)W_t\right)
\\
&=
4\pi_t.
\end{align*}
This completes the proof that $\pi_t$ is an idempotent. The fact that $\pi_t$ is self-adjoint follows immediately from its definition and the fact that $\pi$ and $\pi'$ are self-adjoint. By Kato \cite[Section I.6.7, p. 55]{Kato}, a symmetric, idempotent operator is an orthogonal projection. This completes the proof that $\pi_t$ is an orthogonal projection and hence the proofs of Conclusions \eqref{item:SimilarProj_HomotopyProj} and \eqref{item:SimilarProj_HomotopyBetween}.

We begin the proof of Conclusion \eqref{item:SimilarProj_HomotopyProjBounds} by bounding the size of the operator $W_t$ in \eqref{eq:Define_Q_t_and_W_t}. Because $\|Q_t\|_{\End(\sH)} < 1$ by \eqref{eq:Norm_Q_t_lessthan_1}, the Neumann series \cite[Section I.4.4, Example 4.5, p. 30]{Kato} yields the fact that $1 - Q_t$ is invertible and the bound for its inverse $W_t^2 = (1 - Q_t)^{-1}$,
\begin{equation}
  \label{eq:Neumann_series_estimate}
  \left\|(1 - Q_t)^{-1}\right\|_{\End(\sH)} \leq \frac{1}{1 - \|Q_t\|_{\End(\sH)}}.
\end{equation}
Therefore, since $\|Q_t\|_{\End(\sH)} \leq  \|R\|_{\End(\sH)} < 1$ by \eqref{eq:Norm_Q_t_lessthan_1},
%COMMENT See, for example, https://math.stackexchange.com/questions/525970/norm-of-an-inverse-operator-t-1-t-1 or insert elementary proof by spectral theory
\[
  \|W_t^2\|_{\End(\sH)}
  =
  \left\| (1 - Q_t)^{-1} \right\|_{\End(\sH)}
  \leq
  \frac{1}{1 - \|R\|_{\End(\sH)}}.
\]
The operator $Q_t$ is self-adjoint by \eqref{eq:Define_Q_t_and_W_t} and the fact that $R$ is self-adjoint, thus so is the operator $W_t^2 = (1 - Q_t)^{-1}$. Since $\|Q_t\|_{\End(\sH)} < 1$ by \eqref{eq:Norm_Q_t_lessthan_1}, the proof of \eqref{eq:PositivityOfR} applied to $Q_t$ instead of $R$ implies that the operator $1 - Q_t$ is positive and hence so is the operator $W_t^2 = (1 - Q_t)^{-1}$. In particular, the square root $W_t = (W_t^2)^{1/2}$ in \eqref{eq:Define_Q_t_and_W_t} is unique and normal via its Dunford integral expression \eqref{eq:DunfordIntegral_forSqrt}. Thus, by \eqref{eq:Norm_of_Powers} we have $\|W_t^2\|_{\End(\sH)} = \|W_t\|_{\End(\sH)}^2$ and
\begin{equation}
\label{eq:BoundOnWt}
  \|W_t\|_{\End(\sH)} = \|W_t^2\|_{\End(\sH)}^{1/2} \leq \frac{1}{\left(1 - \|R\|_{\End(\sH)}\right)^{1/2}}.
\end{equation}
Alternately, one could observe that $f(z)=(1-z)^{-1/2}$ is a positive and increasing function of $z \in(0,1)$,
so that the inequality \eqref{eq:BoundOnWt} follows from \eqref{eq:NormOf_Function_of_Operator}. Noting that $g(z):=(1-z)^{-1/2}-1$ is a positive and increasing function of $z \in(0,1)$, we see that
\eqref{eq:NormOf_Function_of_Operator} also gives an estimate for the operator norm of $W_t-1$ and thus of
$1-W_t$,
\begin{equation}
\label{eq:BoundOnModifier}
\|1-W_t\|_{\End(\sH)} \le (1-\|R\|_{\End(\sH)})^{-1/2}-1.
\end{equation}
The hypothesis \eqref{eq:DifferenceHypothesis_of_Projections_On_HilbertSpace_1/4}
that $\|\pi-\pi'\|_{\End(\sH)}<1/4$
and the equality \eqref{eq:RNormBound} imply that
$\|R\|_{\End(\sH)}<1/16$ and
\[
\left( 1-\|R\|_{\End(\sH)}\right)^{-1/2}< \frac{4}{\sqrt {15}}.
\]
From the preceding inequality and inequalities \eqref{eq:BoundOnModifier} and \eqref{eq:BoundOnWt}, we see that
\begin{equation}
\label{eq:BoundsFor(1/2)Bound}
\|1-W_t\|_{\End(\sH)} < \frac{4}{\sqrt {15}}-1
\quad\text{and}\quad
\|W_t\|_{\End(\sH)}<\frac{4}{\sqrt {15}} .
\end{equation}
Because the operator $2\pi-1$ is self-adjoint, the equalities \eqref{eq:ReflectionSquared=1} and \eqref{eq:Norm_of_Powers} imply that
\begin{equation}
\label{eq:2pi_minus_Id_Norm}
\|2\pi-1\|_{\End(\sH)}=1.
\end{equation}
We compute
\begin{align*}
2\| \pi-\pi_t\|_{\End(\sH)}
  &=
\| (2\pi -1)-(2\pi-1)W_t - 2t(\pi'-\pi)W_t\|_{\End(\sH)}
\\
&\quad\text{(by \eqref{eq:FamilyOfOperators} and \eqref{eq:Define_Q_t_and_W_t})}
\\
&\le
\|(2\pi-1)(1-W_t)\|_{\End(\sH)} + 2t\|\pi-\pi'\|_{\End(\sH)}\|W_t\|_{\End(\sH)}
\\
&\le
\|2\pi-1\|_{\End(\sH)} \|1-W_t\|_{\End(\sH)}
\\
&\quad
     + 2t\|\pi-\pi'\|_{\End(\sH)}\left(1-\|R\|_{\End(\sH)}\right)^{-1/2}
\quad\text{(by \eqref{eq:BoundOnWt})}
\\
&\le
     (1-\|R\|)^{-1/2}-1+ 2t\|\pi-\pi'\|_{\End(\sH)}\left(1-\|R\|_{\End(\sH)}\right)^{-1/2}
  \\
  &\qquad\text{(by \eqref{eq:BoundOnModifier} and \eqref{eq:2pi_minus_Id_Norm}).}
\end{align*}
Noting that $t\le 1$, the preceding inequalities yield
\begin{equation}
\label{eq:pi_pit_Bound1}
2\| \pi-\pi_t\|_{\End(\sH)}
\le
(1-\|R\|)^{-1/2}-1+ 2\|\pi-\pi'\|_{\End(\sH)}\left(1-\|R\|_{\End(\sH)}\right)^{-1/2},
\quad\text{for } t \in [0,1].
\end{equation}
Applying the estimate \eqref{eq:BoundsFor(1/2)Bound}  and the hypothesis $\|\pi-\pi'\|_{\End(\sH)}<1/4$ in \eqref{eq:DifferenceHypothesis_of_Projections_On_HilbertSpace_1/4} to the  inequality \eqref{eq:pi_pit_Bound1}
yields
\begin{align*}
2\| \pi-\pi_t\|_{\End(\sH)}
&\le
\left( \frac{4}{\sqrt {15}}-1\right)+\frac{1}{2}\frac{4}{\sqrt {15}}
=
\frac{6}{\sqrt {15}}-1 < 1,
\end{align*}
which gives the first inequality in Conclusion \eqref{item:SimilarProj_HomotopyProjBounds}, 
\begin{equation}
\label{eq:SimilarProj_HomotopyProjBounds1}
\| \pi-\pi_t\|_{\End(\sH)}<\frac{1}{2}, \quad\text{for all } t\in [0,1].
\end{equation}
To prove the second inequality in Conclusion \eqref{item:SimilarProj_HomotopyProjBounds}, we first observe that
\begin{align*}
2\pi_t&=1 + \left(2\pi-1+2t(\pi'-\pi)\right)W_t,
\\
2\pi_{1-t}&=1 + \left(2\pi'-1+2t(\pi-\pi')\right)W_t.
\end{align*}
Note that the preceding expression for $2\pi_{1-t}$ equals the preceding expression for $2\pi_t$ with $\pi$ and $\pi'$ transposed. Thus, we can  we can repeat the derivation of the inequality \eqref{eq:pi_pit_Bound1}, transposing $\pi$ and $\pi'$, as follows:
\begin{align*}
\|2(\pi'-\pi_{1-t})\|_{\End(\sH)}
&=
\|(2\pi'-1)+(2\pi'-1)W_t +2t(\pi-\pi')W_t\|_{\End(\sH)}
\\
&\le
\|(2\pi'-1)\|_{\End(\sH)} \|(1-W_t)\|_{\End(\sH)}
\\
&\quad
     + 2t\|\pi-\pi'\|_{\End(\sH)}\left(1-\|R\|_{\End(\sH)}\right)^{-1/2}
\\
&\le
(1-\|R\|)^{-1/2}-1+ 2\|\pi-\pi'\|_{\End(\sH)}\left(1-\|R\|_{\End(\sH)}\right)^{-1/2},
\end{align*}
where we have again used \eqref{eq:BoundOnModifier} and \eqref{eq:2pi_minus_Id_Norm} to obtain the final inequality above.

Because the map $[0,1]\ni t\mapsto 1-t\in [0,1]$ is a bijection and because the preceding inequalities hold for all $t\in [0,1]$, we have the following analogue of \eqref{eq:pi_pit_Bound1} for all $t\in [0,1]$:
\begin{equation}
\label{eq:pi'_pit_Bound}
\|2(\pi'-\pi_t)\|_{\End(\sH)}
\le
(1-\|R\|)^{-1/2}-1+ 2\|\pi-\pi'\|_{\End(\sH)}\left(1-\|R\|_{\End(\sH)}\right)^{-1/2}.
\end{equation}
The rest of the argument is identical to that used to derive  \eqref{eq:SimilarProj_HomotopyProjBounds1} from \eqref{eq:pi_pit_Bound1}. Applying the estimate \eqref{eq:BoundsFor(1/2)Bound} and the hypothesis $\|\pi-\pi'\|_{\End(\sH)}<1/4$ in \eqref{eq:DifferenceHypothesis_of_Projections_On_HilbertSpace_1/4} to the  inequality \eqref{eq:pi'_pit_Bound} yields the second inequality in Conclusion \eqref{item:SimilarProj_HomotopyProjBounds}. This completes the proof of Proposition \ref{prop:SimilarProj}.
\end{proof}

We now construct the homotopy of almost complex structures on the ranges of the orthogonal projection operators given in Proposition \ref{prop:SimilarProj}.

\begin{prop}[On the existence of a homotopy of almost complex structures]
\label{prop:Constructing_Homotopy_Of_AC}
Continue the assumptions and notation of Lemma \ref{lem:Existence_of_AC_Structure} and assume that the orthogonal projections $\pi$ and $\pi'$ satisfy \eqref{eq:DifferenceHypothesis_of_Projections_On_HilbertSpace_1/4}. Let $[0,1]\ni t\mapsto \pi_t \in \End(\sH)$ be the analytic family of orthogonal projection operators constructed in Proposition \ref{prop:SimilarProj}. Then the following hold:
\begin{enumerate}
\item
\label{item:Constructing_Homotopy_Of_AC_HomotopyOfOperators_onRange}
If  $V_t:=\Ran \pi_t$, then the composition
\begin{equation}
  \label{eq:A_V_t}
  A_{V,t}:=\pi_t J \in\End(V), \quad\text{for } t\in [0,1],
\end{equation}
is an invertible operator and skew-adjoint with respect to the inner product $\langle\cdot,\cdot\rangle_{V_t}$ on $V_t$ induced by $\langle\cdot,\cdot\rangle_{\sH}$ and
\[
J_{V,t} := \left(A_{V,t}^*A_{V,t}\right)^{-1/2}A_{V,t}, \quad\text{for } t\in [0,1],
\]
is an almost complex structure on $V_t$ with $J_{V,0}=\tilde J\in \Aut(V)$, the almost complex structure given in \eqref{eq:Projection_ac}, and $J_{V,1}=J\in \GL(V')$.
\item
\label{item:Constructing_Homotopy_Of_AC_HomotopyOfOperators_Equiv}
The space $V_t$ is $G$-equivariant and the almost complex structure $J_{V,t}$ is $G$-invariant in the sense of Definition \ref{defn:GInvarComplexStructures} for all $t\in [0,1]$.
\item
  \label{item:Constructing_Homotopy_Of_AC_Invertibility}
The operator $\pi_t|_V:V\to V_t$ is invertible for all $t\in [0,1]$.
\item
  \label{item:Constructing_Homotopy_Of_AC_Existence_of_ProjectionAC}
There is an analytic family of $G$-invariant almost complex structures $[0,1]\ni t\to J_t\in \GL(V)$ such that $J_0 = \tilde J$ and $J_1 = (\pi'|_V)^{-1} \circ J|_{V'} \circ \pi'|_V$, where $\tilde J$ is the almost complex structure on $V$ defined in \eqref{eq:Projection_ac} and $(\pi'|_V)^{-1} \circ J|_{V'} \circ \pi'|_V$ is the almost complex structure on $V$ defined by pulling back $J$ on $V'$ by the isomorphism $\pi'|_V \in \Hom(V,V')$.
\end{enumerate}
\end{prop}

\begin{proof}
Because $\pi$ and $\pi'$ satisfy \eqref{eq:DifferenceHypothesis_of_Projections_On_HilbertSpace_1/4}, Proposition \ref{prop:SimilarProj} \eqref{item:SimilarProj_HomotopyProjBounds} implies that $\pi_t$ satisfies
\begin{equation}
\label{eq:Bounds_on_pi_t}
\|\pi_t-\pi'\|_{\End(\sH)}<\frac{1}{2}
\quad\text{and}\quad
\|\pi_t-\pi\|_{\End(\sH)}<\frac{1}{2}, \quad\text{for } t \in [0,1].
\end{equation}
The first inequality in \eqref{eq:Bounds_on_pi_t} and Lemma \ref{lem:Existence_of_AC_Structure} \eqref{item:Existence_of_AC_Structure_InvertibleA_V} yield Conclusion \eqref{item:Constructing_Homotopy_Of_AC_HomotopyOfOperators_onRange}.

We now prove Conclusion \eqref{item:Constructing_Homotopy_Of_AC_HomotopyOfOperators_Equiv}. The $G$-equivariance of $\pi_t$ given by Proposition \ref{prop:SimilarProj} implies that $V_t=\Ran\pi_t$ is $G$-invariant. The $G$-invariance of $J_{V,t}$ follows from the argument used to establish the $G$-invariance of $J_V$ in Lemma \ref{lem:Existence_of_AC_Structure} \eqref{item:Existence_of_AC_Structure_InvertibleA_V}.

Conclusion \eqref{item:Constructing_Homotopy_Of_AC_Invertibility} follows from
the second inequality in \eqref{eq:Bounds_on_pi_t} and Lemma \ref{lem:GlobalDifference_Projections}.

We now prove Conclusion \eqref{item:Constructing_Homotopy_Of_AC_Existence_of_ProjectionAC}. We define an analytic family of almost complex structures on $V$ by
\[
  J_t := \left(\pi_t|_V\right)^{-1}\circ J_{V,t}\circ \pi_t|_V, \quad\text{for all } t\in [0,1].
\]
By  Proposition \ref{prop:SimilarProj} \eqref{item:SimilarProj_HomotopyBetween}, we have $\pi_0=\pi$ and $\pi_1=\pi'$. Hence, the preceding definition of $J_t$, the equality $\pi_0=\pi$, and the equalities $A_{V,0} = A_V$ from \eqref{eq:A_V} and \eqref{eq:A_V_t} and $J_{V,0} = \tilde J$ from \eqref{eq:Projection_ac} yield $J_0=\tilde J$. Because $J_{V,1}=J|_{V'}$ by Conclusion \eqref{item:Constructing_Homotopy_Of_AC_HomotopyOfOperators_onRange}, the equality $\pi_1=\pi'$ implies that
\[
  J_1 = (\pi'|_V)^{-1}\circ J_{V,1} \circ \pi'|_V = (\pi'|_V)^{-1}\circ J|_{V'} \circ \pi'|_V,
\]
as required. Because all the operators in the definition of $J_t$ are $G$-equivariant, the almost complex structure $J_t$ is $G$-invariant.
\end{proof}

\begin{proof}[Proof of Proposition \ref{prop:SameComplexIndex}]
For $k=1,2$, let
\[
\Pi_{\nu,k}':\sH_k\to \bH_\nu^k(\sT') \quad\text{and}\quad
\tilde\Pi_{\nu,k}:\sH_k\to \mathbf{\tilde H}_\nu^k(\sT)
\]
be the orthogonal projections, where $\bH_\nu^k(\sT')$ is the bounded-eigenvalue eigenspace defined in \eqref{eq:DefineLowEigenvalueSpaces} and $\mathbf{\tilde H}_\nu^k(\sT)$ is defined in \eqref{eq:Define_Kuranishi_Spaces_for_AC_Operator}. We first show that these orthogonal projections satisfy the inequality \eqref{eq:DifferenceHypothesis_of_Projections_On_HilbertSpace_1/4}. Using the identification between the bounded-eigenvalue subspaces $\bH_\nu^k(\sT) \subset \sH_k$ and the subspaces $\sH_{1,n}\subset\sH_1$ and $\sH_{2,m}\subset\sH_2$ given in the proof of Lemma \ref{lem:InjectivityOfProjection_for_KurnaishiModel}, we see that for a given $\eps\in (0,1/2]$, Lemma \ref{lem:Approximation_finite-dimensional_subspaces_Hilbert_space} and Corollary \ref{cor:Approximation_projections_onto_finite-dimensional_subspaces_Hilbert_space} imply that for $\nu$ sufficiently large,
\begin{equation}
\label{eq:EndHBound_on_pi_pi'}
\|\tilde\Pi_{\nu,k}-\Pi_{\nu,k}'\|_{\End(\sH_k)}<\eps.
\end{equation}
Thus, if $\nu$ is large enough that $\eps=1/4$ in \eqref{eq:EndHBound_on_pi_pi'}, then
$\tilde\Pi_{\nu,k}$ and $\Pi_{\nu,k}'$  will satisfy \eqref{eq:DifferenceHypothesis_of_Projections_On_HilbertSpace_1/4}.

Because $\tilde\Pi_{\nu,k}$ and $\Pi_{\nu,k}'$ satisfy \eqref{eq:DifferenceHypothesis_of_Projections_On_HilbertSpace_1/4}, Lemma \ref{lem:Existence_of_AC_Structure} \eqref{item:Existence_of_AC_Structure_InvertibleA_V} yields the existence of the $G$-invariant almost complex structures $\tilde J_k$ on $\mathbf{\tilde H}_\nu^k(\sT)$
in \eqref{eq:AConTildeH}, proving Conclusion \eqref{item:SameComplexIndex_ExistenceOfAC} of
Proposition \ref{prop:SameComplexIndex}.

By Proposition \ref{prop:Constructing_Homotopy_Of_AC} and the inequality \eqref{eq:DifferenceHypothesis_of_Projections_On_HilbertSpace_1/4}, there is an
analytic family of $G$-invariant almost complex structures $[0,1]\ni t\mapsto J_{k,t}\in \GL(\mathbf{\tilde H}_\nu^k(\sT))$  with $J_{k,0}=\tilde J_k$ and $J_{k,1}=J_k'$. This homotopy of $G$-invariant almost complex structures and Proposition \ref{prop:Homotopy_Invariance_of_G_Representations} imply that
\[
[\mathbf{\tilde H}_\nu^k(\sT),\tilde J_k]=[\mathbf{\tilde H}_\nu^k(\sT),J_k']\in R_\CC(S^1),
\]
completing the proof of Conclusion \eqref{eq:SameComplexIndex} and hence the proof of Proposition  \ref{prop:SameComplexIndex}. 
\end{proof}

\section{The circle action on the affine and quotient spaces of unitary triples}
\label{sec:S1ActionOnSO3MonopoleDefOperator}
To construct the  $S^1$-equivariant structure of the deformation operator \eqref{eq:Perturbed_Deformation_Operator} for the non-Abelian monopole equations
\eqref{eq:SO(3)_monopole_equations_almost_Hermitian_perturbed_intro_regular} with a regularized Taubes perturbation, we begin by defining $S^1$ actions on the domain and range of this operator and on the quotient space of unitary triples. Following the discussion in Feehan and Leness \cite[Section 14.1.1]{Feehan_Leness_introduction_virtual_morse_theory_so3_monopoles}, we define three homomorphisms from $S^1$ to the space of smooth, unitary gauge transformations of $E$:
\begin{subequations}
\label{eq:S1Actions}
\begin{align}
\label{eq:S1ZUnitaryAction}
\rho_Z(e^{i\theta})&= e^{i\theta}\,\id_E,
\\
\label{eq:S12UnitaryAction}
\rho_2(e^{i\theta})&=\id_{L_1}\oplus e^{i\theta}\,\id_{L_2},
\\
\label{eq:SU_UnitaryAction}
\rho_{\SU}(e^{i\theta})&=e^{i\theta}\,\id_{L_1}\oplus e^{-i\theta}\,\id_{L_2}
\end{align}
\end{subequations}
These three homomorphisms are related by
\begin{equation}
\label{eq:S1HomorphismRelation}
\rho_2(e^{2i\theta})=\rho_{\SU}(e^{-i\theta})\rho_Z(e^{i\theta})\quad\text{for all $e^{i\theta}\in S^1$}.
\end{equation}
These homomorphisms define actions of $S^1$ on the space of unitary triples,
\begin{multline}
\label{eq:S1Z_Action_On_AffineSpace}
\rho_Z^{\sA}:S^1\times\sA^{1,p}(E,A_d,H)\times W^{2,p}(E\oplus \Lambda^{0,2}(E))
\ni (e^{i\theta},A,\varphi,\psi)
\\
\mapsto
\left( A,e^{i\theta}\varphi,e^{i\theta}\psi\right)
\in \sA(E,H,A_d)\times\Omega^0(E)\oplus \Omega^{0,2}(E),
\end{multline}
and
\begin{multline}
\label{eq:S1L2_Action_On_AffineSpace}
\rho_2^{\sA}:S^1\times\sA^{1,p}(E,A_d,H)\times W^{1,p}(E\oplus \Lambda^{0,2}(E))
\ni (e^{i\theta},A,\varphi,\psi)
\\
\mapsto
\left( \rho_2(e^{-i\theta})^*A,\rho_2(e^{i\theta})\varphi,\rho_2(e^{i\theta})\psi\right)
\in \sA(E,H,A_d)\times W^{1,p}(E\oplus \Lambda^{0,2}(E)).
\end{multline}
We note that if $(A,\varphi,\psi)$ is a split pair with respect to the splitting $E=L_1\oplus L_2$
then $(A,\varphi,\psi)$ is a fixed point of the action \eqref{eq:S1L2_Action_On_AffineSpace}.
Such a pair will only be a fixed point of the action  \eqref{eq:S1L2_Action_On_AffineSpace} if $\varphi\equiv 0$ and $\psi\equiv 0$. By \cite[Lemma 14.1.1]{Feehan_Leness_introduction_virtual_morse_theory_so3_monopoles},
the equality \eqref{eq:S1HomorphismRelation} and the fact that the image of $\rho_Z$ is contained in the group $\Omega^0(\SU(E))$ imply that the $S^1$ actions \eqref{eq:S1Z_Action_On_AffineSpace} and \eqref{eq:S1L2_Action_On_AffineSpace} define the same action, up to \emph{positive} multiplicity, on the quotient space of pairs,
\begin{equation}
\label{eq:S1ActionEqualityOnQuotient}
[\rho_2^{\sA}(e^{2i\theta})(A,\varphi,\psi)]
=
[A,e^{i\theta}\varphi,e^{i\theta}\psi]
\end{equation}
where the square brackets in \eqref{eq:S1ActionEqualityOnQuotient} indicate a gauge-equivalence class in the quotient space of pairs. Because the multiplicity of two appearing in the relation \eqref{eq:S1ActionEqualityOnQuotient} is positive, the two actions will have the same virtual Morse--Bott index at a split pair.

We now recall from \cite[Lemma 14.1.4]{Feehan_Leness_introduction_virtual_morse_theory_so3_monopoles} an expression for the action induced by $\rho_2^\sA$ on the tangent space of  a fixed point. To simplify notation, we shall denote the tangent space to the space unitary triples by $\sE_1$ in \eqref{eq:sE1} to although the tangent space is actually given by a Sobolev completion of that space.

\begin{lem}[Circle action on affine space of unitary triples]
\label{lem:S12ActionOnNonAbelianPairsAffineSpace}
Let $(E,H)$ be a smooth, complex Hermitian vector bundle over  a smooth, closed, almost Hermitian four-manifold $(X,g,J,\omega)$. Let $(A,\varphi,\psi)$ be a $W^{1,p}$ unitary triple on $E$ as in
\eqref{eq:A_varphi_psi_in_W1p} and assume that $(A,\varphi,\psi)$ is split  in the sense of Definition \ref{defn:Split_trivial_central-stabilizer_spinor_pair} with respect to an orthogonal splitting $E = L_1 \oplus L_2$ as a direct sum of Hermitian line bundles. Then $(A,\varphi,\psi)$  is a fixed point of the $S^1$ action \eqref{eq:S1L2_Action_On_AffineSpace}. For $\sE_1$ as in \eqref{eq:sE1}, let
\begin{equation}
\label{eq:UnitaryS12ActionOnTangentSpaceOfFixedPoint}
(D_2\rho_2^\sA):S^1\times \sE_1 \to \sE_1
\end{equation}
denote the linear action on the tangent space at $(A,\varphi,\psi)$ to the affine space of unitary triples induced by the derivative of $\rho_2^\sA$ in the directions tangent to the affine space.  Then
\begin{equation}
\label{eq:ExplicitUnitaryS12ActionOnTangentSpaceOfFixedPoint}
(D_2\rho_2^\sA)(e^{i\theta})\left(a,\sigma,\tau\right)
=
\left( \rho_2(e^{i\theta})a \rho_2(e^{-i\theta}),\rho_2(e^{i\theta})\sigma,\rho_2(e^{i\theta})\tau\right),
\end{equation}
for all $e^{i\theta}\in S^1$, and $a\in\Om^1(\su(E))$, and $\sigma\in\Omega^0(E)$, and $\tau\in\Omega^{0,2}(E)$.
\end{lem}

\section[Circle-equivariance of the deformation operator]{Circle-equivariance of the deformation operator for the perturbed non-Abelian monopole equations}
\label{sec:S1EquivariantStructureForPerturbedSO3MonopoleDefOperator}
Let $(A,\varphi,\psi)$ be a smooth solution of the non-Abelian monopole equations \eqref{eq:SO(3)_monopole_equations_almost_Hermitian_perturbed_intro_regular} with a regularized Taubes perturbation on a rank-two Hermitian vector bundle $(E,H)$ over a smooth, closed almost Hermitian manifold $(X,g,J,\omega)$ of real dimension four. In Section \ref{sec:HarmonicTheoryForPerturbedEquations}, we wrote the perturbed deformation operator for the non-Abelian monopole equations \eqref{eq:SO(3)_monopole_equations_almost_Hermitian_perturbed_intro_regular} with a regularized Taubes perturbation as
\[
  \sT_{A,\varphi,\psi,r} = d_{A,\varphi,\psi}^{0,*}+ d_{A,\varphi,\psi,r}^1:\sE_1\to \sE_2,
\]
where $\sE_1$ and $\sE_2$ are as in \eqref{eq:sEk}.  We write elements of $\sE_1$ and $\sE_2$ as
\begin{align*}
(a,\sigma,\tau) & \in\sE_1=\Omega^1(\su(E))\oplus \Omega^0(E)\oplus \Omega^{0,2}(E),
\\
(\xi_1,\xi_2,v,\nu) & \in\sE_2=\Omega^0(\su(E))\oplus \su(E)) \oplus \Omega^{0,2}(\su(E))\oplus \Omega^{0,1}(E).
\end{align*}
The $S^1$ action \eqref{eq:ExplicitUnitaryS12ActionOnTangentSpaceOfFixedPoint} defined by the derivative of the
action $\rho_2^{\sA}$ at a fixed point $(A,\varphi,\psi)$ is an action on the domain of $\sT_{A,\varphi,\psi,r}$,
\begin{equation}
\label{eq:S1ActionOnSO(3)MonopolesDefComplexSpacesFirstTerm}
S^1\times\sE_1
\ni (e^{i\theta},(a,\sigma,\tau))
\mapsto \left(\rho_2(e^{i\theta})a \rho_2(e^{i\theta})^{-1},\rho_2(e^{i\theta})(\sigma,\tau) \right)
\in \sE_1,
\end{equation}
and we define an analogous action on the codomain $\sE_2$,
\begin{multline}
\label{eq:S1ActionOnSO(3)MonopolesDefComplexSpacesSecondTerm}
S^1\times\sE_2
\ni (e^{i\theta},(\xi_1,\xi_2,v,\nu))
\\
\mapsto
\left(\rho_2(e^{i\theta})\xi_1 \rho_2(e^{i\theta})^{-1},\rho_2(e^{i\theta})\xi_2 \rho_2(e^{i\theta})^{-1},\rho_2(e^{i\theta})v \rho_2(e^{i\theta})^{-1},\rho_2(e^{i\theta})\nu \right)
\in \sE_2.
\end{multline}
The proof of \cite[Lemma 14.1.8]{Feehan_Leness_introduction_virtual_morse_theory_so3_monopoles} immediately yields the following

\begin{lem}[Circle-equivariance of the perturbed deformation operator for non-Abelian monopoles]
\label{lem:S1_Equivariance_Of_nonAbelianMonopole_Def_Complex}
Let $(A,\varphi,\psi)$ be a smooth solution of the non-Abelian monopole equations \eqref{eq:SO(3)_monopole_equations_almost_Hermitian_perturbed_intro_regular} with a regularized Taubes perturbation on a rank-two Hermitian vector bundle $(E,H)$ over a smooth, closed almost Hermitian manifold $(X,g,J,\omega)$ of real dimension four. Assume that $(A,\varphi,\psi)$ is split in the sense of Definition \ref{defn:Split_trivial_central-stabilizer_spinor_pair} with respect to a decomposition $E=L_1\oplus L_2$ as an orthogonal direct sum of Hermitian line bundles. Then for all $r\in \RR$, the operator $\sT_{A,\varphi,\psi,r}$ defined in \eqref{eq:Perturbed_Deformation_Operator} is an $S^1$-equivariant elliptic operator in the sense of Definition \ref{defn:Index_of_GEquivariantDiffOperator} with respect to the actions \eqref{eq:S1ActionOnSO(3)MonopolesDefComplexSpacesFirstTerm} and \eqref{eq:S1ActionOnSO(3)MonopolesDefComplexSpacesSecondTerm}.
\end{lem}

The $S^1$-equivariance given in Lemma \ref{lem:S1_Equivariance_Of_nonAbelianMonopole_Def_Complex} and Lemma \ref{lem:S1InvariantSubspacesAndProjections} imply the following analogue of \cite[Corollary 14.1.11]{Feehan_Leness_introduction_virtual_morse_theory_so3_monopoles}.

\begin{cor}[Circle actions on the harmonic spaces for the elliptic deformation complex of a split non-Abelian monopole]
\label{cor:S1_Action_on_HarmonicSections_Of_nonAbelianMonopole_Def_Complex}
Continue the assumptions of Lemma \ref{lem:S1_Equivariance_Of_nonAbelianMonopole_Def_Complex} and assume that $X$ is closed. Then the spaces
\begin{equation}
\label{eq:HarmonicSpaces_for_D}
\Ker \sT_{A,\varphi,\psi,r}  \subset \sE_1
\quad\text{and}\quad
\Ker \sT_{A,\varphi,\psi,r}^* \subset \sE_2
\end{equation}
where $\sT_{A,\varphi,\psi,r}^*$ is the $L^2$-adjoint of the operator $\sT_{A,\varphi,\psi,r}$ defined in \eqref{eq:Perturbed_Deformation_Operator} are closed under the circle actions \eqref{eq:S1ActionOnSO(3)MonopolesDefComplexSpacesFirstTerm} and \eqref{eq:S1ActionOnSO(3)MonopolesDefComplexSpacesSecondTerm}, respectively.
\end{cor}

\section{Circle-equivariance of the equivalent deformation operator}
\label{sec:S1ActionOnApproxComplexDefOperator}
As in the preceding section, let $(A,\varphi,\psi)$ be a smooth solution of the non-Abelian monopole equations  \eqref{eq:SO(3)_monopole_equations_almost_Hermitian_perturbed_intro_regular} with a regularized Taubes perturbation on a rank-two Hermitian vector bundle $(E,H)$ over a smooth, closed almost Hermitian manifold $(X,g,J,\omega)$ of real dimension four. We will again assume that $(A,\varphi,\psi)$ is split with respect to a decomposition $E=L_1\oplus L_2$ as an orthogonal direct sum of Hermitian line bundles and we will describe an $S^1$-equivariant structure for the equivalent deformation operator $\cT_{\partial_A,\varphi,\psi,r}$ given in \eqref{eq:Perturbed_AC_Deformation_Operator} by
\[
\cT_{\partial_A,\varphi,\psi,r}
=
\bar\partial_{A,\varphi,\psi}^1+\hat\partial_{A,\varphi,\psi,r}^{0,*}:\sF_1\to\sF_2,
\]
where as in \eqref{eq:sEkC}, the Fr\'echet spaces $\sF_k$ are given by
\begin{align*}
\sF_1&=\Omega^{0,1}(\fsl(E))\oplus  \Omega^0(E) \oplus \Omega^{0,2}(E),
\\
\sF_2&= \Omega^0(\fsl(E))\oplus \Omega^{0,2}(\fsl(E)) \oplus \Omega^{0,1}(E),
\end{align*}
and  the operators $\hat\partial_{A,\varphi,\psi,r}^{0,*}$  and $\bar\partial_{A,\varphi,\psi}^1$
are defined in \eqref{eq:DefineHat0*OperatorPerturbed}, \eqref{eq:DefineHatPartial0}, and \eqref{eq:DefineBarPartialWithMu} as, for $a''\in\Om^{0,1}(\fsl(E))$, $\si\in\Omega^0(E)$, $\tau\in\Omega^{0,2}(E)$,
\begin{equation}
\label{eq:DefineHatPartial0_IndexReview}
\hat\partial_{A,\varphi,\psi,r}^{0,*}(a'',\sigma,\tau)
=
\bar\partial_A^*a'' - R_\varphi^*\si +(R_\psi^*\tau)^\dagger
+ \frac{r}{4}\left( p_\gamma(\psi)\tau + (p_\gamma(\psi)\tau)^\dagger\right)
\in\Omega^0(\fsl(E)),
\end{equation}
where $p_\gamma(\psi)$ is defined in \eqref{eq:DefineComplexDerivativeOfPerturbation}, and
\begin{equation}
\label{eq:DefineBarPartialWithMu_IndexReview}
\bar\partial_{A,\varphi,\psi}^1(a'',\sigma,\tau)
=
\begin{pmatrix}
      \bar\partial_Aa''-\frac{1}{4}N_J^*(a'')^\dagger - \left(\tau\otimes\varphi^* + \psi\otimes\sigma^*\right)_0
      \\
      \bar\partial_A\sigma + \bar\partial_A^*\tau + a''\varphi + \star ((a'')^\dagger\wedge\star\psi)
  \end{pmatrix}
  \in \Omega^{0,2}(\fsl(E))\oplus \Om^{0,1}(E).
\end{equation}
We define $S^1$ actions on the domain $\sF_1$ and codomain $\sF_2$ by
\begin{subequations}
\label{eq:S1ActionOnApproxComplexPertDefComplex}
\begin{multline}
\label{eq:S1ActionOnApproxComplexPertDefComplexFirstTerm}
S^1\times\left( \Om^{0,1}(\fsl(E))\oplus\Omega^0(E)\oplus \Omega^{0,2}(E)\right)
\ni (e^{i\theta},(a'',\sigma,\tau))
\\
\mapsto
\rho_{2,1}(e^{i\theta})(a'',\sigma,\tau):=
 \left(\rho_2(e^{i\theta})a'' \rho_2(e^{i\theta})^{-1},\rho_2(e^{i\theta})\sigma,\rho_2(e^{i\theta})\tau \right)
 \\
\in \Om^{0,1}(\fsl(E))\oplus\Omega^0(E)\oplus \Omega^{0,2}(E),
\end{multline}
\begin{multline}
\label{eq:S1ActionOnApproxComplexPertDefComplexSecondTerm}
S^1\times\left(\Omega^0(\fsl(E))\oplus \Omega^{0,2}(\fsl(E))\oplus\Omega^{0,1}(E)\right)
\ni (e^{i\theta},(\zeta,v,\nu))
\\
\mapsto
\rho_{2,2}(e^{i\theta})(\zeta,v,\nu):=
\left(\rho_2(e^{i\theta})\zeta \rho_2(e^{i\theta})^{-1},\rho_2(e^{i\theta})v \rho_2(e^{i\theta})^{-1},\rho_2(e^{i\theta})\nu \right)
\\
\in \Omega^0(\fsl(E))\oplus \Omega^{0,2}(\fsl(E))\oplus\Omega^{0,1}(E).
\end{multline}
\end{subequations}
We then have the

\begin{lem}[Circle-equivariant isomorphisms between Fr\'echet spaces]
\label{lem:EquivarianceOfApproxComplexIsomorphismDiagram}
Let $(E,H)$ be a smooth Hermitian vector bundle over an almost Hermitian four-manifold $(X,g,J,\omega)$. If $(A,\varphi,\psi)$ is a unitary triple on $E$ which is split with respect to an orthogonal decomposition $E=L_1\oplus L_2$ as a direct sum of Hermitian line bundles $L_i$, for $i=1,2$, then the following hold.
\begin{enumerate}
\item
The isomorphism $\Upsilon_1$ in \eqref{eq:Isomorphism_sE1C_to_sE1} is $S^1$-equivariant with respect to the action \eqref{eq:S1ActionOnApproxComplexPertDefComplexFirstTerm} on its domain and the action
\eqref{eq:S1ActionOnSO(3)MonopolesDefComplexSpacesFirstTerm} on its codomain.
\item
The isomorphism $\Upsilon_2$ in \eqref{eq:Isomorphism_sE2C_to_sE2} is $S^1$-equivariant with respect to the action \eqref{eq:S1ActionOnApproxComplexPertDefComplexSecondTerm} and the action \eqref{eq:S1ActionOnSO(3)MonopolesDefComplexSpacesSecondTerm} on its codomain.
\item
The operator $\cT_{\partial_A,\varphi,\psi,r}$ defined in \eqref{eq:Perturbed_AC_Deformation_Operator} is  $S^1$-equivariant with respect to the action \eqref{eq:S1ActionOnApproxComplexPertDefComplexFirstTerm} and the action \eqref{eq:S1ActionOnApproxComplexPertDefComplexSecondTerm} on its codomain.
\end{enumerate}
\end{lem}

The $S^1$-equivariance of $\cT_{\partial_A,\varphi,\psi,r}$ and Lemma \ref{lem:S1EquivOfComponents} imply that $\cT_{\partial_A,\varphi,\psi,r}'$, the complex linear component of $\cT_{\partial_A,\varphi,\psi,r}$, is also $S^1$-equivariant. By Lemma \ref{lem:S1InvariantSubspacesAndProjections}, the $S^1$-equivariance of the operator $\cT_{\partial_A,\varphi,\psi,r}$ then implies the following.

\begin{cor}[Circle action on the bounded-eigenvalue eigenspaces of the equivalent deformation operator]
\label{cor:BoundedEigenOfApproxComplexDefOp_S1Closed}
Let $(E,H)$ be a smooth Hermitian manifold over an almost Hermitian four-manifold $(X,g,J,\omega)$. If $(A,\varphi,\psi)$ is a unitary triple on $E$ which is split with respect to an orthogonal decomposition $E=L_1\oplus L_2$ as a direct sum of Hermitian line bundles $L_i$, for $i=1,2$, then the for $\mu>0$, the vector spaces
\[
\Ker\cT_{\partial_A,\varphi,\psi,r}\subset\sF_1\quad\text{and}\quad
\Ker \cT_{\partial_A,\varphi,\psi,r}^*\subset\sF_2,
\]
where $\cT_{\partial_A,\varphi,\psi,r}^*$ is the $L^2$-adjoint of the operator $\cT_{\partial_A,\varphi,\psi,r}$ given in \eqref{eq:Perturbed_AC_Deformation_Operator},
and, for $\nu\notin\sigma(\cT_{\partial_A,\varphi,\psi,r}^{\prime,*}\cT_{\partial_A,\varphi,\psi,r}')$, the bounded-eigenvalue eigenspaces
\[
  \bH_\nu^1(\cT_{\partial_A,\varphi,\psi,r}')\subset\sF_1
  \quad\text{and}\quad
  \bH_\nu^2(\cT_{\partial_A,\varphi,\psi,r}')\subset\sF_2,
\]
defined in \eqref{eq:DefineLowEigenvalueSpaces} are closed under the $S^1$ actions \eqref{eq:S1ActionOnApproxComplexPertDefComplexFirstTerm} on $\sF_1$ and the action \eqref{eq:S1ActionOnApproxComplexPertDefComplexSecondTerm} on $\sF_2$.
\end{cor}

Lemma \ref{lem:EquivarianceOfApproxComplexIsomorphismDiagram}, the $S^1$-equivariance of the deformation operator $\sT_{A,\varphi,\psi,r}$ given in Lemma \ref{lem:S1_Equivariance_Of_nonAbelianMonopole_Def_Complex},
the $S^1$-equivariance of the operators in Lemma \ref{lem:EquivarianceOfApproxComplexIsomorphismDiagram},
and the equality
\[
\Upsilon_2\circ \cT_{\partial_A,\varphi,\psi,r}=\sT_{A,\varphi,\psi,r}\circ \Upsilon_1
\]
from \eqref{eq:Perturbed_Equation_Equivalence_Of_Deformation_Complex} yield the following $S^1$-equivariant
version of the isomorphisms of harmonic spaces appearing in Proposition \ref{prop:SO3MonopoleDeformationCorollary_Perturbed}.

\begin{lem}[Circle equivariance of the isomorphisms between the kernels of $\sT_{A,\varphi,\psi,r}$ and of $\cT_{\partial_A,\varphi,\psi,r}$]
\label{lem:S1EquivariantIsom_LowEigenvalueSpaces_DefOp_to_ApproxComplex}
Let $(A,\varphi,\psi)$ be a smooth solution of the non-Abelian monopole equations  \eqref{eq:SO(3)_monopole_equations_almost_Hermitian_perturbed_intro_regular} with a regularized Taubes perturbation on a rank-two Hermitian vector bundle $(E,H)$. Assume that $(A,\varphi,\psi)$ is split in the sense of Definition \ref{defn:Split_trivial_central-stabilizer_spinor_pair} with respect to an orthogonal decomposition $E=L_1\oplus L_2$ as a direct sum of Hermitian line bundles $L_j$, for $j=1,2$. Then for $k=1,2$, the isomorphisms $\Upsilon_k:\sF_i\to\sE_i$ induce $S^1$-equivariant isomorphisms
\begin{subequations}
\label{eq:S1EquivariantIsom_LowEigenvalueSpaces_DefOp_to_ApproxComplex}
\begin{align}
\label{eq:S1EquivariantIsom_LowEigenvalueSpaces_DefOp_to_ApproxComplex_1}
\Upsilon_1: \Ker \cT_{\partial_A,\varphi,\psi,r} \to \Ker \sT_{A,\varphi,\psi,r},
\\
\label{eq:S1EquivariantIsom_LowEigenvalueSpaces_DefOp_to_ApproxComplex_2}
(\Upsilon_2^{-1})^*:\Ker\cT_{\partial_A,\varphi,\psi,r}^* \to \Ker \sT_{A,\varphi,\psi,r}^*,
\end{align}
\end{subequations}
where the operator $\sT_{A,\varphi,\psi,r}$ is defined in \eqref{eq:Perturbed_Deformation_Operator} and the operator $\cT_{\partial_A,\varphi,\psi,r}$ in \eqref{eq:Perturbed_AC_Deformation_Operator}.
\end{lem}

\begin{proof}
This follows immediately from Proposition \ref{prop:SO3MonopoleDeformationCorollary_Perturbed} and the $S^1$-equivariance of $\Upsilon_k$.
\end{proof}

For $\nu\notin\sigma(\cT_{\partial_A,\varphi,\psi,r}^{\prime,*}\cT_{\partial_A,\varphi,\psi,r}')$, we write
orthogonal projection as,
\begin{equation}
  \label{eq:Define_Projections_To_ComplexComponentEigenspace}
  \Pi_{A,\varphi,\psi,r,\nu}^{\prime,k}:\sF_k \to H_\nu^k(\cT_{\partial_A,\varphi,\psi,r}'),
  \quad\text{for $k=1,2$}
\end{equation}
and, following \eqref{eq:ImageOfKernels}, define
\begin{equation}
\label{eq:Define_Images_Of_Kernel_of_AC_Def}
\begin{aligned}
\tilde K_1(\cT_{\partial_A,\varphi,\psi,r})
&:=\Pi_{A,\varphi,\psi,r,\nu}^{\prime,1}\left(\Ker \cT_{\partial_A,\varphi,\psi,r} \right) \subset \bH_\nu^1(\cT_{\partial_A,\varphi,\psi,r}')
\\
\tilde K_2(\cT_{\partial_A,\varphi,\psi,r})
&:=\Pi_{A,\varphi,\psi,r,\nu}^{\prime,1}\left(\Ker \cT_{\partial_A,\varphi,\psi,r}^* \right) \subset \bH_\nu^2(\cT_{\partial_A,\varphi,\psi,r}')
\end{aligned}
\end{equation}
and following \eqref{eq:Define_Kuranishi_Spaces_for_AC_Operator}, define
\begin{equation}
\label{eq:Define_Kuranishi_Spaces_for_AC_Operator1}
\begin{aligned}
  \mathbf{\tilde H}_\nu^1(\cT_{\partial_A,\varphi,\psi,r})
  &:= \Ker \cT_{\partial_A,\varphi,\psi,r} \oplus
  \left( \tilde K_1(\cT_{\partial_A,\varphi,\psi,r})^\perp\cap \bH_\nu^1(\cT_{\partial_A,\varphi,\psi,r}')\right)\subset\sF_1,
\\
\mathbf{\tilde H}_\nu^2(\cT_{\partial_A,\varphi,\psi,r})
&:= \Ker \cT_{\partial_A,\varphi,\psi,r}^* \oplus
\left( \tilde K_2(\cT_{\partial_A,\varphi,\psi,r})^\perp\cap \bH_\nu^2(\cT_{\partial_A,\varphi,\psi,r}')\right)\subset\sF_2.
\end{aligned}
\end{equation}
From Corollary \ref{cor:IsomOfKuranishiModelSpaces}, we then have the

\begin{cor}[Circle-equivariant isomorphism between bounded-eigenvalue eigenspaces]
\label{cor:Isom_Between_BoundedEigenvalueEigenspaces}
Continue the assumptions and notation of Corollary \ref{cor:BoundedEigenOfApproxComplexDefOp_S1Closed}.
Let $\nu>0$  satisfy $\nu\notin\sigma(\cT_{\partial_A,\varphi,\psi,r}^{\prime,*}\cT_{\partial_A,\varphi,\psi,r}')$ and be sufficiently large to satisfy the conclusion of Lemma \ref{lem:InjectivityOfProjection_for_KurnaishiModel} for the operator $\cT_{\partial_A,\varphi,\psi,r}$. The vector spaces $\mathbf{\tilde H}_\nu^k(\cT_{\partial_A,\varphi,\psi,r})$ defined in \eqref{eq:Perturbed_AC_Deformation_Operator}, where $k=1,2$, are closed under the $S^1$ actions
\eqref{eq:S1ActionOnApproxComplexPertDefComplexFirstTerm} and \eqref{eq:S1ActionOnApproxComplexPertDefComplexSecondTerm}.
The projection maps \eqref{eq:Define_Projections_To_ComplexComponentEigenspace} define $S^1$-equivariant
isomorphisms,
\begin{equation}
  \label{eq:S1_Equiv_Isom_Between_BoundedEigenValue_Eigenspaces}
\Pi_{A,\varphi,\psi,r,\nu}^{\prime,i}:\mathbf{\tilde H}_\nu^k(\cT_{\partial_A,\varphi,\psi,r})
\to  \bH_\nu^k(\cT_{\partial_A,\varphi,\psi,r}'), \quad\text{where $k=1,2$},
\end{equation}
and thus define  $S^1$-invariant almost complex structures on $\mathbf{\tilde H}_\nu^k(\cT_{\partial_A,\varphi,\psi,r})$ as in \eqref{eq:PullbackACStructure}.
\end{cor}

Applying Corollary \ref{cor:Kuranishi_model_defined_by_Fredholm_map_Hilbert_spaces} and the isomorphisms
of Corollary \ref{lem:S1EquivariantIsom_LowEigenvalueSpaces_DefOp_to_ApproxComplex}, we can produce an $S^1$-equivariant Kuranishi model for a neighborhood of $[A,\varphi,\psi]$ in the moduli space of solutions to the non-Abelian monopole equations \eqref{eq:SO(3)_monopole_equations_almost_Hermitian_perturbed_intro_regular} with a regularized Taubes perturbation based on a map from $\mathbf{\tilde H}_\nu^1(\cT_{\partial_A,\varphi,\psi,r})$ to $\mathbf{\tilde H}_\nu^2(\cT_{\partial_A,\varphi,\psi,r})$. The almost complex structure on these spaces given by Corollary \ref{cor:Isom_Between_BoundedEigenvalueEigenspaces} then allow us to compute the virtual Morse--Bott index of the Hitchin function at  $[A,\varphi,\psi]$  using the index of the operator $\cT_{\partial_A,\varphi,\psi,r}$ which we will compute in the following sections.

\section[Weight decomposition of equivariant index of deformation operator]{Weight decomposition of the equivariant index of the equivalent deformation operator}
\label{sec:WeightDecompOfApproxComplexDeformationOperator}
We continue to assume that $(A,\varphi,\psi)$ is a smooth solution of the non-Abelian monopole equations \eqref{eq:SO(3)_monopole_equations_almost_Hermitian_perturbed_intro_regular} with a regularized Taubes perturbation on a rank-two Hermitian vector bundle $(E,H)$ over a smooth, closed almost Hermitian manifold $(X,g,J,\omega)$ of real dimension four. We further assume that $(A,\varphi,\psi)$ is split in the sense of Definition \ref{defn:Split_trivial_central-stabilizer_spinor_pair} with respect to an orthogonal decomposition $E=L_1\oplus L_2$ as a direct sum of Hermitian line bundles $L_j$, for $j=1,2$. Let $\cT_{\partial_A,\varphi,\psi,r}$ be the equivalent deformation operator defined in \eqref{eq:Perturbed_AC_Deformation_Operator} and let $\cT_{\partial_A,\varphi,\psi,r}'$ be the complex linear component of $\cT_{\partial_A,\varphi,\psi,r}$. In the forthcoming Proposition \ref{prop:S1EquivKuranishiModel1}, we will see that we can compute the virtual Morse--Bott index at $[A,\varphi,\psi]$ by computing the equivariant index of $\cT_{\partial_A,\varphi,\psi,r}$.  To do so, we apply \eqref{eq:ComplexS1IndexEquality} to $\cT_{\partial_A,\varphi,\psi,r}$, yielding
\begin{equation}
\label{eq:ComplexS1IndexEqualityForAC_Deformation_Operator}
\Ind_{S^1,\CC}(\cT_{\partial_A,\varphi,\psi,r})
=
\Ind_{S^1,\CC}(\cT_{\partial_A,\varphi,\psi,r}')\in R_\CC(S^1).
\end{equation}
If we write $\cT_A'$ for the compact perturbation of $\cT_{\partial_A,\varphi,\psi,r}'$ defined by removing the linearization of the Taubes-perturbation term and setting $(\varphi,\psi)=(0,0)$, then Proposition \ref{prop:HomotopyInvarianceOfGindex} implies that
\begin{equation}
\label{eq:EqualityOf_sT'_WithDiagonalOperator}
\Ind_{S^1,\CC}(\cT_{\partial_A,\varphi,\psi,r}')
=
\Ind_{S^1,\CC}(\cT_A').
\end{equation}
We will show that $\cT_A'$ is diagonal with respect to weight-space decompositions of $\sF_1$ and $\sF_2$. Equation \eqref{eq:AdditivityOfEquivariantIndex} then implies that the right-hand-side of \eqref{eq:EqualityOf_sT'_WithDiagonalOperator} can be written as a  sum of index expressions which can be easily computed by an application of the (non-equivariant) Atiyah--Singer index formula.

To describe the weight-space decomposition of $\sF_i$, we begin by describing the adjoint action of $S^1$ on $\fsl(E)$. A splitting $E=L_1\oplus L_2$ of a complex rank-two vector bundle into a direct sum of two complex line bundles, $L_1$ and $L_2$, induces a direct-sum decomposition,
\begin{equation}
\label{eq:slDecomp}
\fsl(E)
\cong
\ubarCC\oplus (L_1\otimes L_2^*)\oplus (L_2\otimes L_1^*),
\end{equation}
where $\ubarCC=X\times\CC$ is the product complex line bundle.
We can visualize the isomorphism \eqref{eq:slDecomp} as
for $\zeta_\CC\in\Omega^0(\ubarCC)$, $\zeta_{12}\in\Omega^0(L_1\otimes L_2^*)$, and $\zeta_{21}\in\Omega^0(L_2\otimes L_1^*)$,
\begin{equation}
\label{eq:FormOfslEDecomp}
\ubarCC\oplus (L_1\otimes L_2^*)\oplus (L_2\otimes L_1^*)
\ni (\zeta_\CC,\zeta_{21},\zeta_{21})
\mapsto
\begin{pmatrix}
\zeta_\CC/2 & \zeta_{12} \\ \zeta_{21} & -\zeta_\CC/2
\end{pmatrix}
\in \fsl(E),
\end{equation}
where the matrix in \eqref{eq:FormOfslEDecomp} is written with respect to the decomposition $E=L_1\oplus L_2$.

As in \cite[Definition 4.2.14]{Feehan_Leness_introduction_virtual_morse_theory_so3_monopoles},
we define the \emph{weight} of a $S^1$ action $\rho:S^1\times V\to V$ on a complex vector space $V$ to be an integer $m$ if $\rho(e^{i\theta})(v)=e^{im\theta} v$, for all $e^{i\theta}\in S^1$ and $v\in V$.

\begin{lem}[Decomposition of circle action on a split vector bundle]
\label{lem:WeightSpaceDecomp}
(See Feehan and Leness \cite[Lemma 14.5.1]{Feehan_Leness_introduction_virtual_morse_theory_so3_monopoles}.)
Let $E$ be a complex rank-two vector bundle which admits a decomposition $E=L_1\oplus L_2$ as a direct sum of complex line bundles. Then the isomorphism
\begin{equation}
\label{eq:DeformationComplexBundleDecomposition}
\left(\ubarCC\oplus L_1 \right)
\oplus
\left(L_2\otimes L_1^*\oplus L_2 \right)
\oplus
\left(L_1\otimes L_2^* \right)
\cong
\fsl(E)\oplus E
\end{equation}
is $S^1$-equivariant with respect to the $S^1$ action on $\fsl(E)\oplus E$ given by
\[
  \left( e^{i\theta}, (\zeta,\sigma)\right)\mapsto \left( \rho_2(e^{i\theta}) \zeta \rho_2(e^{-i\theta}),\rho_2(e^{i\theta})\sigma\right),
  \quad\text{for all $e^{i\theta}\in S^1$, $\zeta\in\fsl(E)$, and $\sigma\in E$,}
\]
and the $S^1$ action on $\left(\ubarCC\oplus L_1 \right)\oplus \left(L_2\otimes L_1^*\oplus L_2 \right)
\oplus\left(L_1\otimes L_2^* \right)$ given by
\[
\left(
e^{i\theta}, (\zeta_\CC,\si_1), (\alpha_{21},\si_2), \alpha_{12}
\right)
\mapsto
\left(
 (\zeta_\CC,\si_1), (e^{i\theta}\alpha_{21},e^{i\theta}\si_2), e^{-i\theta}\alpha_{12}
\right), \quad\text{for all } e^{i\theta}\in S^1,
\]
and for all $(\zeta_\CC,\si_1)\in \ubarCC\oplus L_1$, and $(\alpha_{21},\si_2)\in L_2\otimes L_1^*\oplus L_2$, and $\alpha_{12}\in L_1\otimes L_2^*$, where $\rho_2:S^1\to\U(E)$ is the homomorphism \eqref{eq:S12UnitaryAction}.
\end{lem}

The proof of the following result is identical to that of \cite[Lemma 14.5.2]{Feehan_Leness_introduction_virtual_morse_theory_so3_monopoles}.

\begin{cor}[Weight decomposition of Fr\'echet spaces]
\label{cor:WeightSplittingOfApproxComplexDefComplex}
Continue the assumptions of Lemma \ref{lem:WeightSpaceDecomp}.  Then the domain $\sF_1$ and range $\sF_2$
of the perturbed equivalent deformation operator in \eqref{eq:Perturbed_AC_Deformation_Operator} admit direct sum decompositions,
\begin{equation}
\label{eq:PerturbedApproxComplexPertDefComplexWeightDecomp}
\begin{aligned}
\sF_1&\cong \sF_1^- \oplus \sF_1^0\oplus \sF_1^+,
\\
\sF_2&\cong \sF_2^- \oplus\sF_2^0\oplus \sF_2^+,
\end{aligned}
\end{equation}
where
\begin{equation}
\label{eq:WeightSubspacesOfApproxComplexDefDomain}
\begin{aligned}
\sF_1^-&:=\Omega^{0,1}(L_1\otimes L_2^*),
\\
\sF_1^0&:=\Omega^{0,1}(\ubarCC)\oplus\Omega^0(L_1)\oplus \Omega^{0,2}(L_1),
\\
\sF_1^+&:=\Omega^{0,1}(L_1^*\otimes L_2)\oplus\Omega^0(L_2)\oplus \Omega^{0,2}(L_2),
\end{aligned}
\end{equation}
and
\begin{equation}
\label{eq:WeightSubspacesOfApproxComplexDefCoDomain}
\begin{aligned}
\sF_2^-&:=\Omega^0(L_1\otimes L_2^*)\oplus \Omega^{0,2}(L_1\otimes L_2^*),
\\
\sF_2^0&:=\Omega^0(\ubarCC)\oplus \Omega^{0,2}(\ubarCC)\oplus \Omega^{0,1}(L_1),
\\
\sF_2^+&:=\Omega^0(L_1^*\otimes L_2)\oplus \Omega^{0,2}(L_1^*\otimes L_2)\oplus \Omega^{0,1}(L_2).
\end{aligned}
\end{equation}
The decomposition \eqref{eq:PerturbedApproxComplexPertDefComplexWeightDecomp} has the following properties.
\begin{itemize}
\item
The decomposition \eqref{eq:PerturbedApproxComplexPertDefComplexWeightDecomp} is $L^2$-orthogonal.
\item
The summands in the decomposition \eqref{eq:PerturbedApproxComplexPertDefComplexWeightDecomp} are weight spaces for the
$S^1$ actions having weight zero on $\sF_k^0$, weight one on $\sF_k^+$ and weight negative one on $\sF_k^-$.
\end{itemize}
\end{cor}

\begin{rmk}[Non-diagonalizability of the equivalent deformation operator]
\label{rmk:NonDiagOfApproxComplexDefOp}
The equivalent deformation operator $\sT_{A,\varphi,\psi,r}$ is not diagonal with respect to the decompositions \eqref{eq:PerturbedApproxComplexPertDefComplexWeightDecomp} because of the presence of the terms, $(R_\psi^*\tau)^\dagger$ in \eqref{eq:DefineHatPartial0_IndexReview},
$(\psi\otimes\sigma^*)_0$ in \eqref{eq:DefineBarPartialWithMu_IndexReview}, and the terms involving $(a'')^\dagger$ in  \eqref{eq:DefineBarPartialWithMu_IndexReview}. These terms reverse the weight of the $S^1$ actions, mapping $\sE_1^{\CC,\pm}$ to $\sE_2^{\CC,\mp}$.
\qed\end{rmk}

Define the operator $\cT_A':\sF_1\to\sF_2$ by
\begin{multline}
\label{eq:Define_Diagonalizable_Deformation_Op}
\cT_A'(a_0'',\sigma_0,\tau_0)=\left( \bar\partial_A^* a_0'',\bar\partial_Aa_0'',\bar\partial_A\sigma_0+\bar\partial_A^*\tau_0
\right)
\in \Omega^0(\fsl(E))\oplus\Omega^{0,2}(\fsl(E))\oplus \Omega^{0,1}(E),
\\
\text{for all } (a_0'',\sigma_0,\tau_0) \in
\Omega^{0,1}(\fsl(E)) \oplus \Omega^0(E) \oplus \Omega^{0,2}(E).
\end{multline}
Because the $S^1$-equivariant, complex linear operators $\cT_{\partial_A,\varphi,\psi,r}'$ and $\cT_A'$ differ by a compact operator, Proposition \ref{prop:HomotopyInvarianceOfGindex} implies that \eqref{eq:EqualityOf_sT'_WithDiagonalOperator} holds.

Assume that the unitary connection $A$ can be written as $A=A_1\oplus A_2$ with respect to the orthogonal decomposition $E=L_1\oplus L_2$, where $L_k$ is a Hermitian line bundle and $A_k$ is a unitary connection on $L_k$, for $k=1,2$. We will see that $\cT_{A,0,0,0}'$ can be written as a direct sum of the following
$S^1$-equivariant complex linear operators, 
\begin{subequations}
\label{eq:ApproxComplexDefOperatorsDecomposition}
\begin{align}
\label{eq:PerturbedApproxComplexPertDefComplexWeightDecompMinus}
\cT_{A_1\oplus A_2}^-&:\sF_1^-\to \sF_2^-,
\\
\label{eq:PerturbedApproxComplexPertDefComplexWeightDecompZero}
\cT_{A_1\oplus A_2}^0&:\sF_1^0\to\sF_2^0,
\\
\label{eq:PerturbedApproxComplexPertDefComplexWeightDecompPlus}
\cT_{A_1\oplus A_2}^+&:\sF_1^+\to \sF_2^+,
\end{align}
\end{subequations}
as follows.  For $a_{12}''\in\Omega^{0,1}(L_1\otimes L_2^*)$, define
\begin{equation}
\label{eq:DefinePerturbedApproxComplexMinusFactor}
\cT_{A_1\oplus A_2}^-(a_{12}''):=\left(\bar\partial_{A_1\otimes A_2^*}^*a_{12}'',\bar\partial_{A_1\otimes A_2^*}a_{12}''\right)
\in \Omega^0(L_1\otimes L_2^*)\oplus \Omega^{0,2}(L_1\otimes L_2^*).
\end{equation}
For $a_\CC''\in\Omega^{0,1}(\ubarCC)$, $\sigma_1\in\Omega^0(L_1)$, and $\tau_1\in\Omega^{0,2}(L_1)$, define
\begin{equation}
\label{eq:DefinePerturbedApproxComplexZeroFactor}
\cT_{A_1\oplus A_2}^0(a_\CC'',\sigma_1,\tau_1):=
\left( \bar\partial^* a_\CC'', \bar\partial a_\CC'', \bar\partial_{A_1}\sigma_1+\bar\partial_{A_1}^*\tau_1\right)
\in
\Omega^0(\ubarCC)\oplus\Omega^{0,2}(\ubarCC)\oplus \Omega^{0,1}(L_1).
\end{equation}
Finally, for $a_{21}''\in\Omega^{0,1}(L_1^*\otimes L_2)$, $\sigma_2\in\Omega^0(L_2)$, and $\tau_2\in\Omega^{0,2}(L_2)$, define
\begin{multline}
\label{eq:DefinePerturbedApproxComplexPlusFactor}
\cT_{A_1\oplus A_2}^+(a_{21}'',\sigma_2,\tau_2)
:=
\left( \bar\partial_{A_1^*\otimes A_2}^*a_{21}'', \bar\partial_{A_1^*\otimes A_2}a_{21}'',\bar\partial_{A_2}\sigma_2+\bar\partial_{A_2}^*\tau_2\right)
\\
\in
\Omega^0(L_1^*\otimes L_2)\oplus \Omega^{0,2}(L_1^*\otimes L_2) \oplus \Omega^{0,1}(L_2).
\end{multline}
We then have the

\begin{lem}[Weight decomposition of diagonalized deformation operator]
\label{lem:DiagonalizableDeformation}
Continue the hypotheses and notation of Lemma \ref{lem:S1EquivariantIsom_LowEigenvalueSpaces_DefOp_to_ApproxComplex}. If the connection $A$ splits as $A=A_1\oplus A_2$ with respect to the orthogonal decomposition $E=L_1\oplus L_2$, where $L_k$ is a Hermitian line bundle and $A_k$ is a unitary connection on $L_k$, for $k=1,2$, then
\begin{equation}
\label{eq:DiagonalizableDeformation}
\cT_A'=\cT_{A_1\oplus A_2}^-\oplus \cT_{A_1\oplus A_2}^0\oplus \cT_{A_1\oplus A_2}^+:
\sF_1^-\oplus \sF_1^0\oplus \sF_1^+
\to
\sF_1^-\oplus \sF_1^0\oplus \sF_1^+.
\end{equation}
\end{lem}

\begin{proof}
The result follows immediately from the direct sum decompositions of the covariant derivative $d_A$ on $E$ as
$d_{A_1}\oplus d_{A_2}$ with respect to the decomposition $E=L_1\oplus L_2$ and of the induced covariant derivative on $\fsl(E))$
as $d_\CC\oplus d_{A_1\oplus A_2^*}\oplus d_{A_1^*\oplus A_2}$ with respect to the decomposition
$\fsl(E)=\ubarCC\oplus (L_1\otimes L_2^*)\oplus (L_1^*\otimes L_2)$.
\end{proof}

Combining equations \eqref{eq:ComplexS1IndexEqualityForAC_Deformation_Operator},
\eqref{eq:EqualityOf_sT'_WithDiagonalOperator}, and \eqref{eq:DiagonalizableDeformation}
and the additivity of $\Ind_{G,\KK}(\cdot)$ given in \eqref{eq:AdditivityOfEquivariantIndex} yields the 

\begin{cor}[Decomposition of equivariant index of deformation operator]
\label{cor:WeightDecomposition_of_DeformationOperatorIndex}
Continue the hypotheses and notation of Lemma \ref{lem:S1EquivariantIsom_LowEigenvalueSpaces_DefOp_to_ApproxComplex}. If the connection $A$ splits as $A=A_1\oplus A_2$ with respect to the orthogonal decomposition $E=L_1\oplus L_2$ as a direct sum of Hermitian line bundles, where $A_k$ is a unitary connection on $L_k$ for $k=1,2$, then
\begin{equation}
\label{eq:EquivariantIndexAsSum}
\Ind_{S^1,\CC}(\cT_{\partial_A,\varphi,\psi,r})
=
\Ind_{S^1,\CC}(\cT_{A_1\oplus A_2}^-)
+
\Ind_{S^1,\CC}(\cT_{A_1\oplus A_2}^0)
+
\Ind_{S^1,\CC}(\cT_{A_1\oplus A_2}^+),
\end{equation}
as elements of $R_{\CC}(S^1)$, where the operators $\cT_{A_1\oplus A_2}^-$, $\cT_{A_1\oplus A_2}^0$, and $\cT_{A_1\oplus A_2}^+$ are defined in \eqref{eq:ApproxComplexDefOperatorsDecomposition}.
\end{cor}

We can then express the $S^1$-equivariant index in terms of non-equivariant indices as follows. From Lawson and Michelsohn \cite[Example III.9.2, p. 212]{LM}, the ring $R_{\CC}(S^1)$ is isomorphic to the ring of Laurent polynomials $\ZZ[t,t^{-1}]$ with the element $t^m$, for $m\in \ZZ$, being given by the irreducible, complex representation defined by $\rho_m(e^{i\theta})=e^{im\theta}$.

\begin{cor}[Equivariant index of deformation operator and ordinary indexes of diagonalized operator]
\label{cor:WeightDecomposition_of_DeformationOperatorIndex_In_NonEquivariantTerms}
Continue the hypotheses and notation of Lemma \ref{lem:S1EquivariantIsom_LowEigenvalueSpaces_DefOp_to_ApproxComplex}. If the connection $A$ splits as $A=A_1\oplus A_2$ with respect to the orthogonal decomposition $E=L_1\oplus L_2$, where $L_k$ is a Hermitian line bundle and $A_k$ is a unitary connection on $L_k$, for $k=1,2$, then 
\begin{equation}
\label{eq:EquivariantIndexAsSumOfNonEquiv}
\Ind_{S^1,\CC}(\cT_{\partial_A,\varphi,\psi,r})
=
\Ind_{\CC}(\cT_{A_1\oplus A_2}^-)t^{-1}
+
\Ind_{\CC}(\cT_{A_1\oplus A_2}^0)t^0
+
\Ind_{\CC}(\cT_{A_1\oplus A_2}^+)t,
\end{equation}
as elements of $R_{\CC}(S^1)$ where $t^m$ is the equivalence class of the irreducible complex representative of $S^1$ given by $\rho_m(e^{i\theta})=e^{im\theta}$ and $\Ind_{\CC}(\cT_{A_1\oplus A_2^*}^\bullet)$ denotes the non-equivariant index of the operator $\cT_{A_1\oplus A_2^*}^\bullet$ defined in \eqref{eq:ApproxComplexDefOperatorsDecomposition}, for $\bullet=-,0,+$.
\end{cor}

\begin{proof}
Because the kernel of $\cT_{A_1\oplus A_2^*}^-$ is contained in $\sF_1^-$ and the kernel of
$\cT_{A_1\oplus A_2^*}^{-,*}$ is contained in $\sF_2^-$, we have the equality
\[
\Ind_{S^1,\CC}(\cT_{A_1\oplus A_2^*}^-) = \Ind_{\CC}(\cT_{A_1\oplus A_2^*}^-) t^{-1}\in R_\CC(S^1).
\]
Analogous equalities apply to $\cT_{A_1\oplus A_2^*}^0$ and $\cT_{A_1\oplus A_2^*}^+$.
The identity \eqref{eq:EquivariantIndexAsSumOfNonEquiv} follows from equation \eqref{eq:EquivariantIndexAsSum}.
\end{proof}

\section{Index computation}
\label{sec:IndexComputation}
We now compute the indices appearing in \eqref{eq:EquivariantIndexAsSumOfNonEquiv}.
These computations are essentially identical to those appearing in \cite[Section 14.6]{Feehan_Leness_introduction_virtual_morse_theory_so3_monopoles}.

\begin{prop}[Riemann--Roch computation over an almost Hermitian four-manifold]
\label{prop:ASForDiracOperators}
Let $F$ be a Hermitian line bundle over a four-dimensional almost Hermitian manifold $(X,g,J,\omega)$ and
$A$ be a unitary connection on $F$. If
\begin{align*}
\bH_F^1&:=\Ker\left( \bar\rd_A^*+\bar\rd_A: \Omega^{0,1}(F)\to\Omega^0(F)\oplus \Omega^{0,2}(F)\right),
\\
\bH_F^2&:=\Coker\left( \bar\rd_A^*+\bar\rd_A: \Omega^{0,1}(F)\to\Omega^0(F)\oplus \Omega^{0,2}(F)\right),
\end{align*}
then
\begin{equation}
\label{eq:ASForNegativeDirac}
\dim_\CC\bH_F^1 - \dim_\CC\bH_F^2 =
-\chi_h(X)-\frac{1}{2}c_1(X)\cdot c_1(F)-\frac{1}{2}c_1(F)^2.
\end{equation}
\end{prop}

\begin{proof}
By the expression for the Dirac operator over an almost Hermitian manifold in \eqref{eq:Gauduchon_3-7-2_auxiliary_Hermitian_bundle_E}, the operator
\[
\bar\rd_A^*+\bar\rd_A: \Omega^{0,1}(F)\to\Omega^0(F)\oplus \Omega^{0,2}(F)
\]
differs from the Dirac operator
\begin{equation}
\label{eq:DiracTwistedByDirac}
D_A^-:\Omega^0(W_{\can}^-\otimes F) \to \Omega^0(W_{\can}^+\otimes F),
\end{equation}
by the zeroth-order term involving Clifford multiplication by the Lee form.  Because the index of an elliptic operator depends only on the symbol, removing the Lee form from the Dirac operator will not change the index.  Hence, the expression on the left-hand-side of \eqref{eq:ASForNegativeDirac} equals the index of the Dirac operator \eqref{eq:DiracTwistedByDirac}.  This index is the negative of the expression of the index of the Dirac operator,
\begin{equation}
\label{eq:PositiveDiracOpOnF}
D_A^+:\Omega^0(W_{\can}^\otimes F)\to \Omega^0(W_{\can}^-\otimes F).
\end{equation}
We note that by the equality $c_1(W_{\can}^+)=c_1(X)$  given in Kotschick \cite[Fact 2.1]{KotschickSW} or Feehan and Leness 
\cite[Equation (8.2.13)]{Feehan_Leness_introduction_virtual_morse_theory_so3_monopoles}, we have
\[
c_1(W_{\can}^+\otimes F)=c_1(W_{\can}^+)+2c_1(F)=c_1(X)+2c_1(F).
\]
By the expression for the index of a spin${}^c$ Dirac operator given in Lawson and Michelsohn \cite[Equation (D.20), p. 399]{LM} and the abbreviation $c_1(\fs_F):=c_1(W_{\can}^+\otimes F)$, the index of the Dirac operator \eqref{eq:PositiveDiracOpOnF} is
\begin{align*}
\left\langle e^{c_1(\fs_F)/2}\hat A(X),[X]\right\rangle
&=
\left\langle\left( 1+\frac{1}{2}c_1(\fs_F)+\frac{1}{8}c_1(\fs_F)^2\right)(1-p_1(X)/24),[X]\right\rangle
\\
&=
\frac{1}{8}\left\langle c_1(\fs_F)^2,[X] \rangle-\frac{1}{24}\langle p_1(X),[X]\right\rangle
\\
&=
\frac{1}{8}\left( c_1(X)^2+4c_1(X)\smile c_1(F) + 4 c_1(F)^2\right) -\frac{1}{8}\sigma(X).
\end{align*}
We can further simplify the preceding expression using \eqref{eq:Define_c1Squred_c_2},
\[
c_1(X)^2=2e(X)+3\sigma(X).
\]
Thus, we have
\[
\frac{1}{8}c_1(X)^2-\frac{1}{8}\sigma(X)
=
\frac{1}{4}(e(X)+\sigma(X))
=
\chi_h(X).
\]
Thus, the index of the Dirac operator \eqref{eq:PositiveDiracOpOnF} is
\[
\chi_h(X)+\frac{1}{2}c_1(X)\cdot c_1(F)+\frac{1}{2}c_1(F)^2.
\]
Taking the negative of this index gives the expression appearing in \eqref{eq:ASForNegativeDirac} and completes the proof of Proposition \ref{prop:ASForDiracOperators}.
\end{proof}

We have the following

\begin{prop}[Index computation of the components of the diagonalized operator]
\label{prop:Index_of_Homotoped_Def_Operators}
Continue the hypotheses and notation of Lemma \ref{lem:S1EquivariantIsom_LowEigenvalueSpaces_DefOp_to_ApproxComplex}. If the connection $A$ splits as $A=A_1\oplus A_2$ with respect to the orthogonal decomposition $E=L_1\oplus L_2$ as a direct sum of Hermitian line bundles, where $A_k$ is a unitary connection on $L_k$ for $k=1,2$, then
\begin{subequations}
\label{eq:IndexForHomotopedDefOperators}
\begin{align}
  \label{eq:IndexForDMinus}
\Ind_\CC \cT_{A_1\oplus A_2}^-&=-\chi_h(X)-\frac{1}{2}c_1(X)\cdot \left( c_1(L_1)-c_1(L_2)\right) -\frac{1}{2}\left( c_1(L_1)-c_1(L_2)\right)^2,
\\
\label{eq:IndexForD0}
\Ind_\CC\cT_{A_1\oplus A_2}^0 &= c_1(X)\cdot c_1(L_1)+\frac{1}{2}c_1(L_1)^2,
\\
\label{eq:IndexForDPlus}
\Ind_\CC\cT_{A_1\oplus A_2}^+&=\frac{1}{2}c_1(X)\cdot c_1(L_1)-\frac{1}{2}\left(c_1(L_2)-c_1(L_1)\right)^2 -\frac{1}{2}c_1(L_2)^2,
\end{align}
\end{subequations}
where $\cT_{A_1\oplus A_2}^-$,  $\cT_{A_1\oplus A_2}^0$, and $\cT_{A_1\oplus A_2}^+$ are defined in
\eqref{eq:ApproxComplexDefOperatorsDecomposition}.
\end{prop}

\begin{proof}
Equation \eqref{eq:IndexForDMinus} follows immediately by applying Proposition \ref{prop:ASForDiracOperators} with $F=L_1\otimes L_2^*$.

To prove equation \eqref{eq:IndexForD0} observe that the operator $\sD_{A_1\oplus A_2}^0$ is the direct sum of the operators
\begin{equation}
\label{eq:D0Operator_Summand1}
\bar\partial^*+\bar\partial:\Omega^{0,1}(\ubarCC)\to \Omega^0(\ubarCC)\oplus \Omega^{0,2}(\ubarCC)
\end{equation}
and 
\begin{equation}
\label{eq:D0Operator_Summand2}
\bar\partial_{A_1}+\bar\partial_{A_1}^*:\Omega^0(L_1)\oplus \Omega^{0,2}(L_1)\to\Omega^{0,1}(L_1).
\end{equation}
By setting $F=\ubarCC$ in Proposition \ref{prop:ASForDiracOperators}, we see that the index
of the operator \eqref{eq:D0Operator_Summand1}
is $-\chi_h(X)$. The operator in \eqref{eq:D0Operator_Summand2} is the adjoint of the operator in Proposition \ref{prop:ASForDiracOperators} with $F=L_1$. Hence, the index of the operator in \eqref{eq:D0Operator_Summand2} is the negative of the expression on the right-hand-side of \eqref{eq:ASForNegativeDirac} with $F=L_1$:
\[
\chi_h(X)+c_1(X)\cdot c_1(L_1)+\frac{1}{2}c_1(L_1)^2.
\]
Adding the indices of the operators \eqref{eq:D0Operator_Summand1} and \eqref{eq:D0Operator_Summand2} then gives the expression appearing in \eqref{eq:IndexForD0}.

We now prove the equality \eqref{eq:IndexForDPlus}. The operator $\sD_{A_1\oplus A_2}^+$ is the direct sum of the operators
\begin{equation}
\label{eq:DPlusOperator_Summand1}
\bar\partial^*_{A_1^*\otimes A_2}+\bar\partial_{A_1^*\otimes A_2}:\Omega^{0,1}(L_1^*\otimes L_2)\to \Omega^0(L_1^*\otimes L_2)\oplus \Omega^{0,2}(L_1^*\otimes L_2)
\end{equation}
and 
\begin{equation}
\label{eq:DPlusOperator_Summand2}
\bar\partial_{A_2}+\bar\partial_{A_2}^*:\Omega^0(L_2)\oplus \Omega^{0,2}(L_2)\to\Omega^{0,1}(L_1).
\end{equation}
The operator \eqref{eq:DPlusOperator_Summand1} is that appearing in Proposition \ref{prop:ASForDiracOperators}
with $F=L_1^*\otimes L_2$.  The index of the operator \eqref{eq:DPlusOperator_Summand1} is thus
\begin{equation}
\label{eq:IndexOfDPlusSummand1}
-\chi_h(X) -\frac{1}{2}c_1(X)\cdot\left( c_1(L_2)-c_1(L_1)\right) -\frac{1}{2}\left( c_1(L_2)-c_1(L_1)\right)^2.
\end{equation}
The operator \eqref{eq:DPlusOperator_Summand2} is  the adjoint of the operator in Proposition \ref{prop:ASForDiracOperators} with $F=L_2$.  Thus the index of the operator in \eqref{eq:D0Operator_Summand2} is the negative of the expression on the right-hand-side of \eqref{eq:ASForNegativeDirac} with $F=L_2$:
\[
\chi_h(X)+\frac{1}{2}c_1(X)\cdot c_1(L_2)+\frac{1}{2}c_1(L_2)^2.
\]
Adding the indices of the operators \eqref{eq:DPlusOperator_Summand1} and \eqref{eq:DPlusOperator_Summand2} then gives the expression appearing in \eqref{eq:IndexForDPlus}.
\end{proof}

\section{Proofs of main results on virtual Morse--Bott indices}
\label{sec:vMBMainTheorems}
We begin by defining the vector spaces which will serve as the domain and codomain for the Kuranishi model in the forthcoming Proposition \ref{prop:S1EquivKuranishiModel1}.

\begin{defn}[Images of bounded-eigenvalue eigenspaces for Kuranishi model]
\label{defn:S1_Equiv_Kuranishi_Model_At_Split_Triple_With_AC_Structure}
Let $(A,\varphi,\psi)$ be a split unitary triple which is a smooth solution of the non-Abelian monopole equations \eqref{eq:SO(3)_monopole_equations_almost_Hermitian_perturbed_intro_regular} with a regularized Taubes perturbation. For $\nu>0$ satisfying $\nu\notin\sigma(\cT_{\partial_A,\varphi,\psi,r}^{\prime,*}\cT_{\partial_A,\varphi,\psi,r}')$ and the assumptions of Corollary \ref{cor:Isom_Between_BoundedEigenvalueEigenspaces}, define
\begin{equation}
\label{eq:Define_Kuranishi_Model_Spaces_For_Equiv}
\mathbf{\hat H}_{A,\varphi,\psi,r,\nu}^1 \subset\sE_1
\quad\text{and}\quad
\mathbf{\hat H}_{A,\varphi,\psi,r,\nu}^2 \subset\sE_2
\end{equation}
to be the images, under the isomorphisms $\Upsilon_1$ and $(\Upsilon_2^{-1})^*$ given in \eqref{eq:Isomorphism_sEkC_to_sEk}, respectively, of the vector spaces $\mathbf{\tilde H}_\nu^1(\cT_{\partial_A,\varphi,\psi,r})$ and $\mathbf{\tilde H}_\nu^2(\cT_{\partial_A,\varphi,\psi,r})$ defined in \eqref{eq:Define_Kuranishi_Spaces_for_AC_Operator1}. Because the vector spaces in \eqref{eq:Define_Kuranishi_Model_Spaces_For_Equiv} are $S^1$-equivariantly isomorphic to the vector spaces $\mathbf{\tilde H}_\nu^k(\cT_{\partial_A,\varphi,\psi,r})$, the $S^1$-invariant almost complex structure on the $\mathbf{\tilde H}_\nu^k(\cT_{\partial_A,\varphi,\psi,r})$ given in Corollary \ref{cor:Isom_Between_BoundedEigenvalueEigenspaces} defines an $S^1$-invariant almost complex structure on $\mathbf{\tilde H}_{A,\varphi,\psi,r,\mu}^k$ for $k=1,2$.
\qed
\end{defn}

\begin{rmk}[Motivation for introducing the spaces $\mathbf{\hat H}_{A,\varphi,\psi,r,\nu}^k$]
To apply Corollary \ref{cor:Kuranishi_model_defined_by_Fredholm_map_Hilbert_spaces} to construct a Kuranishi model
of a point in $\sM^0(E,g,J,\omega,r)$, we require subspaces $K_1$ and $K_2$ of the Fr\'echet spaces $\sE_1$ and $\sE_2$ defined in \eqref{eq:sEk}. We define the spaces $\mathbf{\hat H}_{A,\varphi,\psi,r,\nu}^k$ to satisfy this requirement.
\end{rmk}

\begin{prop}[Circle-equivariant local Kuranishi model using bounded-eigenvalue eigenspaces]
\label{prop:S1EquivKuranishiModel1}
Let $(X,g,J,\omega)$ be a smooth, closed, almost Hermitian manifold of real dimension four and let $(E,H)$ be a rank two smooth Hermitian vector bundle. Let $(A,\varphi,\psi)$ be a smooth, non-zero section solution to the non-Abelian monopole equations \eqref{eq:SO(3)_monopole_equations_almost_Hermitian_perturbed_intro_regular} with a regularized Taubes perturbation with parameter $r \in[1,\infty)$.
In addition, assume that $(A,\varphi,\psi)$ is split in the sense of Definition \ref{defn:Split_trivial_central-stabilizer_spinor_pair}
with respect to  an orthogonal decomposition $E=L_1\oplus L_2$ as a direct sum of Hermitian line bundles $L_j$, for $j=1,2$. Let $\nu>0$ satisfy $\nu\notin\sigma(\cT_{\partial_A,\varphi,\psi,r}^{\prime,*}\cT_{\partial_A,\varphi,\psi,r}')$ and the assumptions of Corollary \ref{cor:Isom_Between_BoundedEigenvalueEigenspaces}. Then there are a neighborhood of the origin $\sU_{A,\varphi,\psi,r,\nu}\subset \mathbf{\tilde H}_\nu^1(\cT_{\partial_A,\varphi,\psi,r})$ and real analytic maps
\begin{equation}
\label{eq:Obstruction_Map_For_PerturbedSplitPair_Equiv}
\bchi:\sU_{A,\varphi,\psi,r,\nu}\to \mathbf{\tilde H}_\nu^2(\cT_{\partial_A,\varphi,\psi,r}),
\end{equation}
and 
\begin{equation}
\label{eq:Gluing_Map_For_PerturbedSplitPair_Equiv}
\bga:\sU_{A,\varphi,\psi,r,\nu} \to \sA(E,H,A_d)\times W^{1,p}(E\oplus \Lambda^{0,2}(E))
\end{equation}
satisfying
\begin{enumerate}
\item
\label{item:S1EquivKuranishiModel1_EquivProp}
The neighborhood $\sU_{A,\varphi,\psi,r,\nu}$ is closed under the $S^1$ action \eqref{eq:S1ActionOnSO(3)MonopolesDefComplexSpacesFirstTerm}
and the maps $\bchi$ of \eqref{eq:Obstruction_Map_For_PerturbedSplitPair_Equiv} and $\bga$ of \eqref{eq:Gluing_Map_For_PerturbedSplitPair_Equiv}  are $S^1$ equivariant with respect to the $S^1$-actions
\eqref{eq:S1ActionOnSO(3)MonopolesDefComplexSpacesFirstTerm} and \eqref{eq:S1ActionOnSO(3)MonopolesDefComplexSpacesSecondTerm}.
\item
\label{item:S1EquivKuranishiModel1_MapProperties}
There is an $S^1$-equivariant, real linear isomorphism
\[
\Gamma_\sF:\left(\Ker \cT_{\partial_A,\varphi,\psi,r}\right)^\perp\cap \mathbf{\tilde H}_\nu^1(\cT_{\partial_A,\varphi,\psi,r})
\cong
\left(\Ker \cT_{\partial_A,\varphi,\psi,r}^*\right)^\perp\cap \mathbf{\tilde H}_\nu^2(\cT_{\partial_A,\varphi,\psi,r})
\]
such that the maps $\bchi$ and $\bga$ satisfy
\begin{equation}
\label{eq:Properties_of_KMaps}
\bga(0)=0,\quad
d\bga(0)=\id,\quad
\bchi(0)=0,\quad
d\bchi(0)=\Gamma_\sF\Upsilon_1^{-1}\pi_1\Upsilon_1,
%PF12-18-2025 You're missing d\bchi(0). What's that?
%TL12-18-1025: It's ugly but here it is
\end{equation}
where $\Upsilon_1$ is defined in \eqref{eq:Isomorphism_sE1C_to_sE1} and $\pi_1$ is the orthogonal projection from the $W^{1,p}$-completion of $\sE_1$ to $\left( \bH_{A,\varphi,\psi,r}^1\right)^\perp\cap \mathbf{\hat H}_{A,\varphi,\psi,r,\nu}^1$.
\item
\label{item:S1EquivKuranishiModel1_Ngh_Param}
If $\pi: \sA(E,H,A_d)\times W^{1,p}(E\oplus \La^{0,2}(E))\to\sC_\ft$ is the quotient map, then $\pi\circ\bga$ gives a smoothly-stratified  homeomorphism from $\bchi^{-1}(0)$ onto a neighborhood of $[A,\varphi,\psi]$ in the moduli space $\sM^0(E,g,J,\omega,r)$ defined in \eqref{eq:Moduli_space_non-Abelian_monopoles_almost_Hermitian_Taubes_regularized_non-zero-section} of solutions to the non-Abelian monopole equations with a regularized Taubes perturbation.
\end{enumerate}
 \end{prop}

\begin{proof}
%PF12-18-2025 Please use some other notation and not \bar\sE_k, which is commonly used for complex conjugate
%TL12-18-2025: Changed to \hat\sE .
Define
\begin{align*}
\hat\sE_1&:=W^{1,p}(T^*\otimes\su(E)\oplus E\oplus \Lambda^{0,2}(E)),
\\
\hat\sE_2&:=L^p(\su(E)\oplus \Lambda^{0,2}(\fsl(E))\oplus \Lambda^{0,1}(E)).
\end{align*}
We will apply  Corollary \ref{cor:Kuranishi_model_defined_by_Fredholm_map_Hilbert_spaces} to the map
\begin{equation}
\label{eq:Kuranishi_Model_FredholmMap_Extended}
\hat\sS:\hat\sE_1\to \hat\sE_2
\end{equation}
defined by
\begin{multline}
\label{eq:Kuranishi_Model_FredholmMap}
W^{1,p}(T^*\otimes\su(E)\oplus E\oplus \Lambda^{0,2}(E))
\ni (a,\sigma,\tau)
\\
\mapsto
\left( d_{A,\varphi,\psi}^{0,*}(a,\sigma,\tau), \sS(A+a,\varphi+\sigma,\psi+\tau)\right)
\in
L^p(\su(E)\oplus \Lambda^{0,2}(\fsl(E))\oplus \Lambda^{0,1}(E)),
\end{multline}
where $\sS$
is defined by the left-hand side of the system \eqref{eq:SO(3)_monopole_equations_almost_Hermitian_perturbed_intro_regular} of non-Abelian monopole equations with a regularized Taubes perturbation. The derivative of the map \eqref{eq:Kuranishi_Model_FredholmMap} at the origin is the deformation operator
\[
\sT_{A,\varphi,\psi,r}=d_{A,\varphi,\psi,r}^1+d_{A,\varphi,\psi}^{0,*},
\]
so the kernel of this linearization is the harmonic space $\bH_{A,\varphi,\psi,r}^1$ in \eqref{eq:H1_For_TaubesPert_regular}.  Because $(\varphi,\psi)\not\equiv (0,0)$, Lemma \ref{lem:H0Vanishing} implies that $\bH_{A,\varphi,\psi}^0$ vanishes so we can identify the cokernel of the deformation operator above with the harmonic space $\bH_{A,\varphi,\psi,r}^2$ in \eqref{eq:DefineH2ForSO3Monopoles_Perturbed}.
Thus, we take the spaces $K_0$ and $C_0$ in Corollary \ref{cor:Kuranishi_model_defined_by_Fredholm_map_Hilbert_spaces} to be
\[
K_0=\bH_{A,\varphi,\psi,r}^1,\quad\text{and}\quad C_0=\bH_{A,\varphi,\psi,r}^2.
\]
We claim that we can apply Corollary \ref{cor:Kuranishi_model_defined_by_Fredholm_map_Hilbert_spaces} using the spaces
\begin{equation}
\label{eq:Define_K_and_C_for_Equiv_Kuranishi}
K_1=
\mathbf{\hat H}_{A,\varphi,\psi,r,\nu}^1\subset \sE_1,\quad
C_1=\mathbf{\hat H}_{A,\varphi,\psi,r,\nu}^2\subset\sE_2,
\end{equation}
where the spaces $\mathbf{\hat H}_{A,\varphi,\psi,r,\nu}^k$ are defined in \eqref{eq:Define_Kuranishi_Model_Spaces_For_Equiv}.
To do so, we must show that $K_1$ and $C_1$ satisfy the Conditions \eqref{item:Orthogonal_projections_induce_embeddings}
and \eqref{item:Isomorphism_tildeK0_perp_onto_tildeC0_perp} of Corollary \ref{cor:Kuranishi_model_defined_by_Fredholm_map_Hilbert_spaces}.

We first show that $K_1$ and $C_1$ satisfy  Condition \eqref{item:Orthogonal_projections_induce_embeddings}. By \eqref{eq:DefineH1ForApproxComplex_Perturbed}, \eqref{eq:Perturbed_AC_Deformation_Operator}, and \eqref{eq:SO3MonopoleH1Isomorphism_perturbed}, $\Upsilon_1$ induces an $S^1$-equivariant real linear isomorphism,
\begin{equation}
\label{eq:cT_Isom_To_bH1}
\Ker \cT_{\partial_A,\varphi,\psi,r} \cong \bH_{A,\varphi,\psi,r}^1.
\end{equation}
Because $(\varphi,\psi)\not\equiv (0,0)$, the map $\Upsilon_2$ induces an $S^1$-equivariant real linear isomorphism,
\begin{equation}
\label{eq:cT_Isom_To_bH2}
\Ker \cT_{\bar\partial_A,\varphi,\psi,r}^* \cong \bH_{A,\varphi,\psi,r}^2,
\end{equation}
by \eqref{eq:SO3MonopoleH2Isomorphism_perturbed_H2version}.
Hence, the inclusions
\begin{subequations}
\label{eq:Kernel_Inclusions}
\begin{align}
\label{eq:Kernel_Inclusions_1}
\Ker \cT_{\partial_A,\varphi,\psi,r}\subset \mathbf{\tilde H}_\nu^1(\cT_{\partial_A,\varphi,\psi,r}),
\\
\label{eq:Kernel_Inclusions_2}
\Ker \cT_{\partial_A,\varphi,\psi,r}^* \subset \mathbf{\tilde H}_\nu^2(\cT_{\partial_A,\varphi,\psi,r})
\end{align}
\end{subequations}
given by \eqref{eq:Define_Kuranishi_Spaces_for_AC_Operator1} and the definition of $\mathbf{\hat H}_{A,\varphi,\psi,r,\nu}^k$ imply that there are inclusions
\begin{equation}
\label{eq:Harmonic_Perturbed_Inclusion_Image_of_AC}
\bH_{A,\varphi,\psi,r}^k \subset \mathbf{\hat H}_{A,\varphi,\psi,r,\nu}^k, \quad\text{for $k=1,2$}.
\end{equation}
Hence, the spaces $\mathbf{\hat H}_{A,\varphi,\psi,r,\nu}^k$ satisfy Condition \eqref{item:Orthogonal_projections_induce_embeddings} of Corollary \ref{cor:Kuranishi_model_defined_by_Fredholm_map_Hilbert_spaces}.

We now show that the spaces in \eqref{eq:Define_K_and_C_for_Equiv_Kuranishi} satisfy Condition \eqref{item:Isomorphism_tildeK0_perp_onto_tildeC0_perp}. The inclusions \eqref{eq:Harmonic_Perturbed_Inclusion_Image_of_AC} imply that the orthogonal projections from $\sE_k$ to $\mathbf{\hat H}_{A,\varphi,\psi,r,\nu}^k$ are the identity when restricted to $\bH_{A,\varphi,\psi,r}^k$.  Therefore, the space $\tilde K_0$ appearing in Condition \eqref{item:Isomorphism_tildeK0_perp_onto_tildeC0_perp} of Corollary \ref{cor:Kuranishi_model_defined_by_Fredholm_map_Hilbert_spaces} equals $\bH_{A,\varphi,\psi,r}^1$ and the space $\tilde C_0$ equals $\bH_{A,\varphi,\psi,r}^2$. To prove that  Condition \eqref{item:Isomorphism_tildeK0_perp_onto_tildeC0_perp} holds, we must therefore demonstrate the existence of an $S^1$-equivariant real linear isomorphism,
\begin{equation}
\label{eq:Isomorphism_for_Condition2}
\Gamma_\sE:
\left( \bH_{A,\varphi,\psi,r}^1\right)^\perp\cap \mathbf{\hat H}_{A,\varphi,\psi,r,\nu}^1
\cong
\left( \bH_{A,\varphi,\psi,r}^2\right)^\perp\cap \mathbf{\hat H}_{A,\varphi,\psi,r,\nu}^2.
\end{equation}
The isomorphism \eqref{eq:cT_Isom_To_bH1} and the definition of $\mathbf{\hat H}_{A,\varphi,\psi,r,\nu}^1$ as the image of $\mathbf{\tilde H}_\nu^1(\cT_{\partial_A,\varphi,\psi,r})$ under $\Upsilon_1$ imply that $\Upsilon_1$ defines an $S^1$-equivariant, real linear isomorphism,
\begin{equation}
\label{eq:Perp_Of_Kernel_in_hatH1}
\left(\Ker \cT_{\partial_A,\varphi,\psi,r}\right)^\perp\cap \mathbf{\tilde H}_\nu^1(\cT_{\partial_A,\varphi,\psi,r})
\cong
\left( \bH_{A,\varphi,\psi,r}^1\right)^\perp\cap \mathbf{\hat H}_{A,\varphi,\psi,r,\nu}^1.
\end{equation}
Similarly, the isomorphism \eqref{eq:cT_Isom_To_bH2} and the definition of $\mathbf{\hat H}_{A,\varphi,\psi,r,\nu}^2$ as the image of $\mathbf{\tilde H}_\nu^2(\cT_{\partial_A,\varphi,\psi,r})$ under $(\Upsilon_2^*)^{-1}$ imply that $(\Upsilon_2^*)^{-1}$ defines an $S^1$-equivariant, real linear isomorphism,
\begin{equation}
\label{eq:Perp_Of_Kernel_in_hatH2}
\left(\Ker \cT_{\partial_A,\varphi,\psi,r}^*\right)^\perp\cap \mathbf{\tilde H}_\nu^2(\cT_{\partial_A,\varphi,\psi,r})
\cong
\left( \bH_{A,\varphi,\psi,r}^2\right)^\perp\cap \mathbf{\hat H}_{A,\varphi,\psi,r,\nu}^2.
\end{equation}
By \eqref{eq:Equiv_Isom_Of_Complements} in Lemma \ref{lem:SymplecticStrutureOnStabilizedSpace}, there is an $S^1$-equivariant, real linear isomorphism,
\begin{equation}
\label{eq:Kernel_Complement_Isomorphism_Perturbed}
\Gamma_\sF:\left(\Ker \cT_{\partial_A,\varphi,\psi,r}\right)^\perp\cap \mathbf{\tilde H}_\nu^1(\cT_{\partial_A,\varphi,\psi,r})
\cong
\left(\Ker \cT_{\partial_A,\varphi,\psi,r}^*\right)^\perp\cap \mathbf{\tilde H}_\nu^2(\cT_{\partial_A,\varphi,\psi,r})
\end{equation}
Combining the isomorphisms \eqref{eq:Perp_Of_Kernel_in_hatH1}, \eqref{eq:Perp_Of_Kernel_in_hatH2}, and \eqref{eq:Kernel_Complement_Isomorphism_Perturbed} gives an $S^1$-equivariant, real linear isomorphism
\begin{equation}
\label{eq:Defining_Gamma_Isom_Composition}
\Gamma_\sE=(\Upsilon_2^*)^{-1}\circ\Gamma_\sF\circ \Upsilon_1^{-1},
\end{equation}
as required in \eqref{eq:Isomorphism_for_Condition2}.  This completes the proof that the spaces $K_1$ and $C_1$ given in \eqref{eq:Define_K_and_C_for_Equiv_Kuranishi} satisfy
 Condition \eqref{item:Isomorphism_tildeK0_perp_onto_tildeC0_perp} of Corollary \ref{cor:Kuranishi_model_defined_by_Fredholm_map_Hilbert_spaces}. Thus, we may apply Corollary \ref{cor:Kuranishi_model_defined_by_Fredholm_map_Hilbert_spaces} to the map \eqref{eq:Kuranishi_Model_FredholmMap} using the spaces $K_1$ and $C_1$.

As noted above, the inclusions \eqref{eq:Harmonic_Perturbed_Inclusion_Image_of_AC} imply that the orthogonal projections from $\hat\sE_k$ to $\mathbf{\hat H}_{A,\varphi,\psi,r,\nu}^k$ are the identity when restricted to $\bH_{A,\varphi,\psi,r}^k$.  This implies that the spaces $\tilde K_1$ and $\tilde C_1$ defined in \eqref{eq:tildeK1_and_tildeC1_and_tildeX_and_tildeY} satisfy
\[
  \tilde K_1=\mathbf{\hat H}_{A,\varphi,\psi,r,\nu}^1,
  \quad\text{and}\quad
  \tilde C_1=\mathbf{\hat H}_{A,\varphi,\psi,r,\nu}^k,
\]
and the $S^1$-equivariant, real linear isomorphisms $\Xi_1$ and $\Xi_2$ defined in \eqref{eq:Xi_isomorphism_tildeX_and_tildeY} are the obvious isomorphisms given by addition,
\[
  \Xi_1: \tilde K_1\oplus \tilde K_1^\perp \to X
  \quad\text{and}\quad
  \Xi_2: \tilde C_1\oplus \tilde C_1^\perp \to Y.
\]
We thus omit the maps $\Xi_k$ to simplify the notation. By Conclusion \eqref{item:F_circ_tildeg_is_dF+pi_C0_circ_tildef} of Corollary \ref{cor:Kuranishi_model_defined_by_Fredholm_map_Hilbert_spaces}, there is an $S^1$-invariant neighborhood $U$ of the origin in $\hat\sE_1$, and an $S^1$-equivariant, real analytic diffeomorphism $\tilde g:U\to g(U)$,  and an $S^1$-equivariant real analytic map $\tilde f:U\to \tilde C_1$ such that
\begin{equation}
\label{eq:g_and_f_Properties}
\tilde g(0)=0,\quad
d\tilde g(0)=\id_{\hat\sE_1},\quad
\tilde f(0)=0,\quad
d\tilde f(0)=\Gamma_\sE\circ \pi_1,
\end{equation}
where $\pi_1:\hat\sE_1\to \left( \bH_{A,\varphi,\psi,r}^1\right)^\perp\cap \mathbf{\hat H}_{A,\varphi,\psi,r,\nu}^1$ is orthogonal projection. In addition, the maps $\tilde g$ and $\tilde f$ satisfy
\begin{equation}
\label{eq:Monopole_Slice_Map_And_KuranishiMaps}
\hat\sS\circ \tilde g = \sT_{A,\varphi,\psi,r}+\pi_{C_0}\circ \tilde f,
\end{equation}
where $\pi_{C_0}:\tilde C_1\to C_0=\bH_{A,\varphi,\psi,r}^2$ is orthogonal projection. By Conclusion \eqref{item:Equality_zero-locus_F_and_zero-locus_tildef} of Corollary \ref{cor:Kuranishi_model_defined_by_Fredholm_map_Hilbert_spaces},
\begin{equation}
\label{eq:ZeroLocus_InSlice}
\hat\sS^{-1}(0)\cap \tilde g(U)
=
\tilde f^{-1}(0)\cap U\cap \tilde K_1.
\end{equation}
Thus, $\pi\circ\tilde g$ gives a smoothly-stratified, $S^1$-equivariant homeomorphism from
\[
\tilde f^{-1}(0)\cap U\cap \tilde K_1
=
\tilde f^{-1}(0)\cap U\cap\mathbf{\hat H}_{A,\varphi,\psi,r,\nu}^1,
\]
onto a neighborhood of $[A,\varphi,\psi]$ in $\sM^0(E,g,J,\omega,r)$. 

We now define the open set $\sU_{A,\varphi,\psi,r,\nu}$ and the maps $\bchi$ and $\bga$ by using the $S^1$-equivariant, real linear isomorphisms $\Upsilon_1$ and $\Upsilon_2$. Define the open set $\sU_{A,\varphi,\psi,r,\nu}$ in
\eqref{eq:Obstruction_Map_For_PerturbedSplitPair_Equiv} by
\[
\sU_{A,\varphi,\psi,r,\nu}:=\Upsilon_1^{-1}\left( U\right)\cap \mathbf{\tilde H}_\nu^1(\cT_{\partial_A,\varphi,\psi,r}).
\]
Note that because 
$\Upsilon_1 \mathbf{\tilde H}_\nu^1(\cT_{\partial_A,\varphi,\psi,r})=\mathbf{\hat H}_{A,\varphi,\psi,r,\nu}^1$,
\begin{equation}
\label{eq:Equality_for_sU}
\Upsilon_1 \left(\sU_{A,\varphi,\psi,r,\nu}\right)
= U\cap\mathbf{\hat H}_{A,\varphi,\psi,r,\nu}^1.
\end{equation}
Because $\Upsilon_1$ is an $S^1$-equivariant, real linear isomorphism from  $\mathbf{\tilde H}_\nu^1(\cT_{\partial_A,\varphi,\psi,r})$ to $\mathbf{\hat H}_{A,\varphi,\psi,r,\nu}^1$, we can define the map $\bga$ of  \eqref{eq:Gluing_Map_For_PerturbedSplitPair_Equiv} as the $S^1$-equivariant, diffeomorphism
\begin{equation}
\label{eq:Define_Gluing_Map_For_PerturbedSplitPair_Equiv}
\bga:=  \tilde g\circ\Upsilon_1: \Upsilon_1^{-1}(U)\to g(U).
\end{equation}
Because $\Upsilon_2^*$ is an $S^1$-equivariant, real linear isomorphism from  $\mathbf{\hat H}_{A,\varphi,\psi,r,\nu}^2$ and $\mathbf{\tilde H}_\nu^2(\cT_{\partial_A,\varphi,\psi,r})$, we can define the map $\bchi$ of \eqref{eq:Obstruction_Map_For_PerturbedSplitPair_Equiv} as the $S^1$-equivariant map
\begin{equation}
\label{eq:Define_Obstruction_Map_For_PerturbedSplitPair_Equiv}
\bchi:= \Upsilon_2^*\circ \tilde f\circ \Upsilon_1.
\end{equation}
Conclusion \eqref{item:S1EquivKuranishiModel1_EquivProp} follows from the $S^1$-invariance of $U$ and the $S^1$-equivariance of
$\Upsilon_k$, $\tilde f$, and $\tilde g$.

We now prove
Conclusion \eqref{item:S1EquivKuranishiModel1_MapProperties}.
The existence of the map $\Gamma_\sF$ has already been established in \eqref{eq:Kernel_Complement_Isomorphism_Perturbed}.
The first three equalities of \eqref{eq:Properties_of_KMaps} follow from the definitions of $\bchi$ and $\bga$ and \eqref{eq:g_and_f_Properties}.  To prove the fourth equality of \eqref{eq:Properties_of_KMaps}, we compute
\begin{align*}
d\chi(0)
&=\Upsilon_2^* d\tilde f(0) \Upsilon_1\quad\text{(by \eqref{eq:Define_Obstruction_Map_For_PerturbedSplitPair_Equiv})}
\\
&=\Upsilon_2^* \Gamma_\sE \pi_1 \Upsilon_1\quad\text{(by \eqref{eq:g_and_f_Properties})}
\\
&=\Gamma_\sF \Upsilon_1^{-1}\pi_1\Upsilon_1\quad\text{(by \eqref{eq:Defining_Gamma_Isom_Composition})}.
\end{align*}
This completes the proof of \eqref{eq:Properties_of_KMaps} and hence of Conclusion \eqref{item:S1EquivKuranishiModel1_MapProperties}.

We now prove Conclusion \eqref{item:S1EquivKuranishiModel1_Ngh_Param}.
By the definition of $\bchi$ in \eqref{eq:Define_Obstruction_Map_For_PerturbedSplitPair_Equiv} and the fact that $\Upsilon_2^*$ is a linear isomorphism, we have
\[
\bchi^{-1}(0)=(\Upsilon_2^*\circ\tilde f\circ\Upsilon_1)^{-1}(0)=(\tilde f\circ\Upsilon_1)^{-1}(0).
\]
Thus,
\[
\Upsilon_1\left(\bchi^{-1}(0)\right)
=
\tilde f^{-1}(0)\cap \Upsilon_1 \left(\sU_{A,\varphi,\psi,r,\nu}\right)
=
\tilde f^{-1}(0)\cap U\cap\mathbf{\hat H}_{A,\varphi,\psi,r,\nu}^1,
\]
where the second equality follows from \eqref{eq:Equality_for_sU}. Recall that $\pi\circ\tilde g$ gives a smoothly-stratified, $S^1$-equivariant homeomorphism from
\[
\tilde f^{-1}(0)\cap U\cap\mathbf{\hat H}_{A,\varphi,\psi,r,\nu}^1
\]
onto a neighborhood of $[A,\varphi,\psi]$ in $\sM^0(E,g,J,\omega,r)$. Then, $\pi\bga=\pi\tilde g\circ\Upsilon_1$ gives a smoothly-stratified, $S^1$-equivariant homeomorphism from $\bchi^{-1}(0)$ onto this neighborhood.  This completes the proof of Conclusion \eqref{item:S1EquivKuranishiModel1_Ngh_Param} and thus of Proposition \ref{prop:S1EquivKuranishiModel1}.
\end{proof}

We can now give the

\begin{proof}[Proof of Theorem \ref{mainthm:MorseIndexAtReduciblesOnAlmostKahler}]
The $S^1$-equivariant Kuranishi model for a neighborhood of $[A,\varphi,\psi]$ in the moduli space of regularly perturbed non-Abelian monopoles given in Proposition \ref{prop:S1EquivKuranishiModel1} implies that the virtual Morse--Bott index of the Hitchin function at $[A,\varphi,\psi]$, \eqref{eq:Virtual_Morse-Bott_signature}, is given by the sum of the virtual dimensions of the irreducible $S^1$ representations of negative weight in the equivariant index,
\[
[\mathbf{\tilde H}_\nu^1(\cT_{\partial_A,\varphi,\psi,r})]- [\mathbf{\tilde H}_\nu^2(\cT_{\partial_A,\varphi,\psi,r})]\in R_\CC(S^1).
\]
Corollary \ref{cor:Isom_Between_BoundedEigenvalueEigenspaces} implies that
\[
[\mathbf{\tilde H}_\nu^1(\cT_{\partial_A,\varphi,\psi,r})]- [\mathbf{\tilde H}_\nu^2(\cT_{\partial_A,\varphi,\psi,r})]
=
[H_\nu^1(\cT_{\partial_A,\varphi,\psi,r}')]-[\bH_\nu^2(\cT_{\partial_A,\varphi,\psi,r}')].
\]
Then Proposition \ref{prop:Equality_of_ComplexEquivariantIndices} and Definition \ref{defn:ApproximatelyComplexIndex}
imply that
\[
[H_\nu^1(\cT_{\partial_A,\varphi,\psi,r}')]-[\bH_\nu^2(\sT)(\cT_{\partial_A,\varphi,\psi,r}')]
=
\Ind_{\CC,S^1}(\cT_{\partial_A,\varphi,\psi,r}').
\]
Thus, Corollary \ref{cor:WeightDecomposition_of_DeformationOperatorIndex_In_NonEquivariantTerms} implies that
 the virtual Morse--Bott index is given by the coefficient of $t^{-1}$ in equation \eqref{eq:EquivariantIndexAsSumOfNonEquiv}.  Hence, the virtual Morse--Bott index \eqref{eq:Virtual_Morse-Bott_signature} equals
\[
2\Ind_\CC (\cT_{A_1\oplus A_2}^- ),
\]
as defined in \eqref{eq:DefinePerturbedApproxComplexMinusFactor}. Theorem \ref{mainthm:MorseIndexAtReduciblesOnAlmostKahler} now follows from Corollary \ref{cor:WeightDecomposition_of_DeformationOperatorIndex_In_NonEquivariantTerms}
and \eqref{eq:IndexForDMinus} in Proposition \ref{prop:Index_of_Homotoped_Def_Operators}.
\end{proof}

\begin{proof}[Proof of Corollary \ref{maincor:MorseIndexAtReduciblesOnSymplecticWithSO3MonopoleCharacteristicClasses}]
Corollary \ref{maincor:MorseIndexAtReduciblesOnSymplecticWithSO3MonopoleCharacteristicClasses} follows from
Theorem \ref{mainthm:MorseIndexAtReduciblesOnAlmostKahler} and the relations between the characteristic classes $c_1(L_1)$, $c_1(L_2)$, $c_1(\fs)$, $c_1(X)$, and $c_1(\ft)$ given in the proof of \cite[Corollary 10]{Feehan_Leness_introduction_virtual_morse_theory_so3_monopoles}.
\end{proof}

\begin{proof}[Proof of Corollary \ref{maincor:Positivity_of_MorseIndexAtReduciblesOnSymplecticWithSO3MonopoleCharacteristicClasses}]
If $(X,\omega_0)$ is symplectic with an $\omega_0$-compatible almost complex structure $J_0$, then by McDuff and Salmon \cite[Theorem 7.1.21, p. 306]{McDuffSalamonSympTop3} there is a symplectic form $\omega$ on the smooth blow-up $\widetilde X$.  If $J$ is a $\omega$-compatible almost complex structure on $\widetilde X$
and $e$ is the Poincar\'e dual of the exceptional sphere, then $c_1(T\widetilde X,J)= c_1(TX,J_0)-e$  by McDuff and Salmon \cite[Corollary 7.1.24, p. 311]{McDuffSalamonSympTop3}. If $\fs_0$ is a spin${}^c$ structure on $X$ with $c_1(\fs_0)=-c_1(TX,J_0)$ then by Morgan \cite[Corollary 6.8.4, p. 103]{MorganSWNotes} and Taubes \cite[Main Theorem, p. 809]{TauSymp}, $\SW_X(\fs_0)\neq 0$. If $\fs$ is the spin${}^c$ structure on $\widetilde X$ with $c_1(\fs)=-c_1(\fs_0)-e$ then $c_1(\fs)=c_1(T\widetilde X,J)$ and, by the blow-up formula Theorem \ref{thm:SWBlowUp}, we have $\SW_{\widetilde X}(\fs)\neq 0$.

Let $\ft$ be the spin${}^u$ structure on $\widetilde X$ constructed in Theorem \ref{mainthm:ExistenceOfSpinuForFlow} which is feasible with respect to $c_1(\fs)$. Let $[A,\Phi] \in \sM_{\ft}$ be any point represented by a non-zero-section, split, non-Abelian monopole in the image of the embedding \eqref{eq:DefnOfIotaOnQuotient} of $M_{\fs'}$ into  $\sM_{\ft}$. By Corollary \ref{maincor:MorseIndexAtReduciblesOnSymplecticWithSO3MonopoleCharacteristicClasses},
the formal Morse--Bott index \eqref{eq:FormalMorseIndexIntroThm} equals the virtual Morse--Bott index at $[A,\Phi]$. Corollary \ref{maincor:Positivity_of_MorseIndexAtReduciblesOnSymplecticWithSO3MonopoleCharacteristicClasses} then follows from the feasibility of $\ft$ and the positivity of the formal Morse--Bott index given by Item \eqref{item:Feasibility_PositivevBM} of Definition \ref{maindefn:Feasibility}.
\end{proof}

\appendix

\chapter{Functional analysis, spectral theory, and maximum principles}
\label{chap:Functional_analysis_maximum_principles}
In this chapter, we include material that is primarily expository but serves to clarify our conventions and or is difficult to find in the literature but plays an essential role in the proofs of our main results.

\section{Weakly and strongly non-degenerate bilinear forms on vector spaces}
\label{sec:Weakly_strongly_non-degenerate_bilinear_forms_Banach_spaces}
We review some basic facts concerning bilinear forms on vector spaces.

Let $\cX$ be a vector space over a field $\FF$ and $\cX^* = \Hom(\cX,\FF)$ denote its dual space, $b:\cX \times \cX \to \FF$ be a symmetric or skew-symmetric bilinear form, and define a linear map $L_b: \cX \to \cX^*$ by $L_bx = b(x,\cdot)$, for all $x \in \cX$. Recall that $b$ is
\begin{inparaenum}[\itshape i\upshape)]
\item \emph{weakly non-degenerate} if $L_b: \cX \to \cX^*$ is injective, and
\item \emph{(strongly) non-degenerate} if $L_b:\cX \to \cX^*$ is an isomorphism of vector spaces.
\end{inparaenum}
(Compare, for example, Kobayashi \cite[Equation (7.5.1), p. 246]{Kobayashi_differential_geometry_complex_vector_bundles}.)

If $\FF = \RR$ and $b$ is a \emph{symmetric} bilinear form on $\cX$ that is \emph{positive definite}, so $b(x,x) > 0$ for all $x \in \cX\less\{0\}$, then $L_b: \cX \to \cX^*$ is injective because if $L_bx = b(x,x) = 0$ for some $x \in \cX$, then $x = 0$ and thus $b$ is weakly non-degenerate, in other words, $b$ is a \emph{weak inner product} on $\cX$.

Suppose now that $\sX$ is a real Banach space with norm $\|\cdot\|_\sX$, let $\sX^*$ denote its continuous dual space, and let $b:\sX \times \sX \to \RR$ be a continuous weak inner product on $\sX$. Define $\|x\|_b := \sqrt{b(x,x)}$ for all $x\in\sX$, so $\|\cdot\|_b$ is a norm on the underlying vector space $\sX$. Recall that two norms on a vector space are equivalent if and only if they define the same topology.
% COMMENT See https://math.stackexchange.com/questions/109733/are-two-norms-equivalent-if-they-induce-the-same-topology-on-a-vector-space
If the norms $\|\cdot\|_b$ and $\|\cdot\|_\sX$ are \emph{equivalent}, then the inner product space $(\sX,b)$ is isomorphic to $(\sX,\|\cdot\|_\sX)$ as a Banach space, so $(\sX,b)$ is complete and thus a Hilbert space. Hence, if $b$ is a weak inner product on $\sX$ and $\|\cdot\|_b$ and $\|\cdot\|_\sX$ are equivalent norms, then the map $L_b: \sX \to \sX^*$ is an isomorphism of Banach spaces (by, for example, Rudin \cite[Theorem 12.5, p. 308]{Rudin}) and $b$ is strongly non-degenerate, in other words, $b$ is an \emph{(strong) inner product} on $\sX$.

If $V$ is a \emph{finite-dimensional} vector space over a field $\FF$, then $\dim V = \dim V^*$ (by Hoffman and Kunze \cite[Section 3.2, Theorem 5, p. 75 and Section 3.5, p. 98]{Hoffman_Kunze_linear_algebra}). Moreover, if $b:V \times V \to \FF$ is a bilinear form, then $\dim\Ran L_b = \dim V - \dim\Ker L_b$ (by \cite[Section 3.1, Theorem 2, p. 71]{Hoffman_Kunze_linear_algebra}). Hence, if $\Ker L_b = (0)$, then $\dim\Ran L_b = \dim V^*$ and $\Ran L_b = V^*$ and so $L_b$ is an isomorphism of vector spaces. Thus, for $\dim_\FF V < \infty$, a bilinear form $b$ is weakly non-degenerate if and only if it is (strongly) non-degenerate.

\section{Spectral theory for unbounded linear operators}
\label{sec:Spectral_theory_unbounded_operators}
For our proofs of Theorem \ref{mainthm:Donaldson_1996jdg_3_Hilbert_space} and Corollary \ref{maincor:Donaldson_1996jdg_3_Hilbert_space_non-self-adjoint}, we review the spectral theory for unbounded operators. We recall from Taylor \cite[Appendix A, Section 8, p. 659]{Taylor_PDE1_3rd_edition} that if $T$ is a \emph{closed} unbounded linear operator on a Banach space $\sX$ with dense domain $\Dom(T) \subset \sX$, then $\zeta \in \CC$ belongs to the \emph{resolvent set} $\rho(T) \subset \CC$ if and only if the bounded linear operator
\begin{equation}
  \label{eq:T_minus_zeta}
  T - \zeta: \Dom(T) \to \sX
\end{equation}
is bijective. The (set-theoretic) inverse (in the notation of Kato \cite[Chapter III, Section 6.1, Equation (6.1), p. 173]{Kato})
\begin{equation}
  \label{eq:Taylor_PDE1_A-8-33}
  R(\zeta,T) = (T - \zeta)^{-1}: \sX \to \Dom(T) \subset \sX,
\end{equation}
is the \emph{resolvent} of $T$. Because $T$ is closed by assumption, then so is $T-\zeta$; moreover the graph $R(\zeta,T)$ is equal to the graph of $T-\zeta$ with its axes interchanged, so $R(\zeta,T)$ is closed with domain $\sX$ and thus $R(\zeta,T) \in \End(\sX)$ is bounded by the Closed Graph Theorem (see, for example, Kato \cite[Chapter III, Section 6.1, Problem 6.1, p. 173]{Kato} and for the Closed Graph Theorem see \cite[Chapter III, Section 5.4, Theorem 5.20, p. 166]{Kato} or Rudin \cite[Chapter 2, Theorem 2.15, p. 51]{Rudin}).
% COMMENT See https://math.stackexchange.com/questions/924821/bounded-resolvent
(Other authors define $\rho(T)$ slightly differently. For example, Conway \cite[Chapter X, Section 1, Definition 1.16, p. 307]{Conway_course_functional_analysis} and Rudin \cite[Chapter 13, Definition 13.26, p. 365]{Rudin} omit the assumption that $T$ is closed and define $\rho(T)$ to be the set of $\zeta \in \CC$ such that the operator \eqref{eq:T_minus_zeta} is bijective with \emph{bounded} inverse \eqref{eq:Taylor_PDE1_A-8-33}; Yosida \cite{Yosida} also omits the assumption that $T$ is closed and defines $\rho(T)$ to be the set of $\zeta \in \CC$ such that the operator \eqref{eq:T_minus_zeta} has dense range in $\sX$ with bounded inverse \eqref{eq:Taylor_PDE1_A-8-33}.)

The resolvent set $\rho(T) \subset \CC$ is \emph{open}. To see this, recall from Kato \cite[Chapter I, Section 5.2 and Chapter III, Section 6.1, Theorem 6.7, p. 174]{Kato} that $R(z,T)$ obeys the \emph{first resolvent identity} (see, for example, \cite[Chapter I, Section 5.2, Equation (5.5), p. 36]{Kato} or Hille and Phillips \cite[Theorem 4.8.1, p. 126]{Hille_Phillips}) for $z_0,z \in \rho(T)$,
\begin{equation}
  \label{eq:First_resolvent_identity}
  R(z,T) - R(z_0,T) = (z-z_0)R(z,T)R(z_0,T),
\end{equation}
so that $R(z,T)$ and $R(z_0,T)$ commute \cite[Chapter I, Section 5.2, p. 36 and Chapter III, Section 6.1, Problem 6.6, p. 173]{Kato} and
\[
  R(z_0,T) = \left[1 - (z-z_0)R(z_0,T)\right]R(z,T).
\]
Hence, one obtains the \emph{first Neumann series} \cite[Chapter I, Section 5.2, Equation (5.6), p. 37]{Kato}
\begin{equation}
  \label{eq:First_Neumann_series}
  R(z,T) = \left[1 - (z-z_0)R(z_0,T)\right]^{-1}R(z_0,T) = \sum_{n=0}^\infty (z-z_0)^nR(z_0,T)^{n+1},
\end{equation}
which is absolutely convergent for all $z\in\CC$ obeying \cite[Chapter I, Section 5.2, Equation (5.7), p. 37]{Kato}
\begin{equation}
  \label{eq:First_Neumann_series_radius_convergence}
  |z-z_0| < \|R(z_0,T)\|_{\End(\sX)}^{-1}.
\end{equation}
In particular, if $z_0 \in \rho(T)$, then $\rho(T)$ contains the open disk with center $z_0$ and $z$ obeying \eqref{eq:First_Neumann_series_radius_convergence}.

If $S,T \in \End(\sX)$ are two closed unbounded operators, one has the \emph{second resolvent identity} (see, for example, Hille and Phillips \cite[Theorem 4.8.2, p. 126]{Hille_Phillips}),
\begin{equation}
  \label{eq:General_second_resolvent_identity}
  R(z,T) - R(z,S) = R(z,T)(S-T)R(z,S), \quad\text{for all } z \in \rho(S)\cap\rho(T),
\end{equation}
and thus $R(z,T) = [R(z,T)(S-T) + 1]R(z,S)$. This leads to the \emph{second Neumann series} (see Kato \cite[Chapter II, Section 1.3, Theorem 1.5, p. 66 and Equation (1.13), p. 67 and Chapter IV, Section 3.3, Theorem 3.11, p. 210]{Kato})
\begin{multline}
  \label{eq:Second_Neumann_series}
  R(z,S) = R(z,T)\left[1 - R(z,T)(T-S)\right]^{-1} = R(z,T)\sum_{n=0}^\infty (R(z,T)(T-S))^n,
  \\
  \text{for all } z \in \rho(S)\cap\rho(T),
\end{multline}
which is absolutely convergent if $(z,S)$ obeys
\begin{equation}
  \label{eq:Second_Neumann_series_radius_convergence}
   \|R(z,T)(S-T)\|_{\End(\sX)} < 1.
\end{equation}
If $T-S \in \End(\sX)$ is \emph{bounded}, then \eqref{eq:Second_Neumann_series_radius_convergence} is implied by the simpler convergence condition for $(z,S)$:
\begin{equation}
  \label{eq:Second_Neumann_series_radius_convergence_bounded_T-S}
  \|S-T\|_{\End(\sX)} < \|R(z,T)\|_{\End(\sX)}^{-1}.
\end{equation}
The complement of $\rho(T)$ is the \emph{spectrum} of $T$ and is denoted by $\sigma(T) \subset \CC$. If $T$ is \emph{bounded}, then its spectrum is a compact subset of the disk (see Kato \cite[Chapter III, Section 6.2, p. 176]{Kato} or Conway \cite[Chapter VII, Section 3, Theorem 3.6, p. 196]{Conway_course_functional_analysis}):
\begin{equation}
  \label{eq:Spectrum_bounded_operator_subset_disk}
  \sigma(T) \subset \left\{z\in \CC: |z| \leq \|T\|_{\End(\sX)}\right\}.
\end{equation}
By Kato \cite[Chapter IV, Section 6.6, Theorem 6.22, p. 184]{Kato}, the resolvent operator for the adjoint $T^* \in \End(\sX^*)$ obeys
\begin{equation}
  \label{eq:Adjoint_resolvent_and_resolvent_adjoint}
  R(z,T)^* = R(\bar z,T^*), \quad\text{for all } z \in \rho(T),
\end{equation} 
and $\rho(T^*)$ and $\sigma(T^*)$ are obtained from $\rho(T)$ and $\sigma(T)$, respectively, by reflection about the real axis. 

If $N$ is a bounded \emph{normal} operator on a \emph{Hilbert space} $\sH$, then
\begin{equation}
\label{eq:Norm_of_Powers}
\|N^n\|_{\End(\sH)} = \|N\|_{\End(\sH)}^n\quad\text{for all integers $n\geq 1$},
\end{equation}
% COMMENT See https://encyclopediaofmath.org/wiki/Normal_operator
% COMMENT See https://math.stackexchange.com/questions/1052614/the-spectral-radius-of-normal-operator
and so by Gelfand's spectral radius formula for a bounded operator (see, for example, Rudin \cite[Chapter 10, Definition 10.10, p. 252 and Theorem 10.13 (b), p. 253]{Rudin}),
\begin{equation}
  \label{eq:Spectral_radius_equals_norm_normal_operator}
  \|N\|_{\End(\sH)} = \sup_{\lambda\in\sigma(N)}|\lambda|.
\end{equation}
See Kato \cite[Chapter V, Section 2.1, Equation (2.4), p. 257]{Kato}. (The equality \eqref{eq:Spectral_radius_equals_norm_normal_operator} is also proved via the method outlined by Dowson \cite[Chapter 7, p. 177, Theorem 7.12 and Chapter 1, Definition 1.7, p. 5]{Dowson_spectral_theory_linear_operators} and Promislow \cite[Section 8.5.3, Theorem 8.23, p. 180]{Promislow_first_course_functional_analysis}.) We remark, as noted in Kato \cite[Chapter I, Section 5.2, p. 38, following Equation (5.12)]{Kato} that the spectral radius in the right-hand side of \eqref{eq:Spectral_radius_equals_norm_normal_operator} is \emph{independent} of the norm used in the definition of the Hilbert space $\sH$.
%COMMENT The equality \eqref{eq:Spectral_radius_equals_norm_normal_operator} also follows from the functional calculus for bounded normal operators in Rudin \cite[Section 12.24, Equation (3), p. 325, and first line of p. 326]{Rudin}.

%PF4-10-2025 We don't need \sH to be a Hilbert space for much of the below. Devide section into Banach space and Hilbert space subsections.
If there exists $\zeta \in \rho(T)$ such that $R(\zeta,T)$ is compact, one says that $T$ has \emph{compact resolvent}. (If $\cV := \Dom(T) \hookrightarrow \sH$ is a compact embedding of Banach spaces, then by composing $R(z,T) \in \Hom(\sH,\cV)$ with the compact embedding $\cV \hookrightarrow \sH$, we would obtain a compact resolvent $R(z,T)\in \End(\sH)$.)
% COMMENT See https://en.wikipedia.org/wiki/Holomorphic_functional_calculus
% COMMENT See https://en.wikipedia.org/wiki/Resolvent_formalism
By Kato \cite[Chapter 3, Section 6.8, Theorem 6.29, p. 187 and Section 6.3, Problem 6.16, p. 177]{Kato} or Taylor \cite[Appendix A, Section 8, Proposition 8.8, p. 660]{Taylor_PDE1_3rd_edition}, it follows that when $T$ has compact resolvent (for some $\zeta \in \CC$), then $\sigma(T)$ is a discrete subset of $\CC$ consisting entirely of eigenvalues of $T$ with finite multiplicities, the resolvent $R(z,T)$ is compact for every $z\in\rho(T)$ and
\begin{equation}
  \label{eq:Spectral_mapping_T_and_RzetaT}
  \lambda \in \sigma(T)\cup\{\infty\} \iff (\lambda-z)^{-1} \in \sigma(R(z,T)).
\end{equation}
One says that $\lambda \in \sigma(T)$ is an \emph{eigenvalue} of $T$ if $\Ker(T - \lambda) \neq (0)$ (see, for example, Rudin \cite[Chapter 4, Definition 4.17, p. 103]{Rudin}). If in addition $T$ is \emph{normal}, then $T - z$ is normal for any $z \in \CC$ (since $(T-z)^* = T^*-\bar z$ commutes with $T-z$). By \eqref{eq:Adjoint_resolvent_and_resolvent_adjoint}, we deduce that $R(z,T)$ is normal for any $z \in \rho(T)$ (see Kato \cite[Chapter V, Section 3.8, p. 277, second line]{Kato}), so the identity \eqref{eq:Spectral_radius_equals_norm_normal_operator} and the equality of spectra in \eqref{eq:Spectral_mapping_T_and_RzetaT} yields
\[
  \|R(z,T)\|_{\End(\sH)} = \sup_{k\in\NN}|\lambda_k - z|^{-1}
  = \left(\inf_{k\in\NN}|\lambda_k - z|\right)^{-1}, \quad\text{for all } z \in \rho(T),
\]
and thus (see Kato \cite[Chapter V, Section 3.8, Equation (3.31), p. 277]{Kato})
%COMMENT See https://math.stackexchange.com/questions/1855151/norm-of-the-resolvent
\begin{equation}
  \label{eq:Spectral_radius_equals_norm_resolvent_operator_normal_T}
  \|R(z,T)\|_{\End(\sH)} = \dist(z,\sigma(T))^{-1},
\quad\text{for all } z \in \rho(T).
\end{equation}
(Dunford and Schwartz \cite[Chapter XII, Section 2, Lemma 2, p. 1191]{Dunford_Schwartz_1963} provide the weaker inequality $\|R(z,T)\|_{\End(\sH)} \leq |\Imag z|^{-1}$ when $T$ is self-adjoint.) By combining the inequality \eqref{eq:Second_Neumann_series_radius_convergence_bounded_T-S} with the equality \eqref{eq:Spectral_radius_equals_norm_resolvent_operator_normal_T}, we see that if $S \in \End(\sH)$ is another unbounded operator such that $S-T \in \End(\sH)$ is bounded and $(z,S)$ obeys
\begin{equation}
  \label{eq:Second_Neumann_series_radius_convergence_bounded_T-S_normal_T}
  \|S-T\|_{\End(\sX)} < \dist(z,\sigma(T)),
\end{equation}
then $z \in \rho(S)$.

% COMMENT See https://en.wikipedia.org/wiki/Spectral_theory_of_compact_operators
% COMMENT See https://sites.math.washington.edu//~hart/m556/Lecture5.pdf
% COMMENT See https://ocw.mit.edu/courses/18-102-introduction-to-functional-analysis-spring-2021/9eaaf363541d01d53c330ee7931fc715_MIT18_102s21_lec20.pdf
% COMMENT See https://faculty.etsu.edu/gardnerr/Func/notes/9-4.pdf
Recall that if $C \in \End(\sX)$ is a compact operator and $\dim\sX=\infty$, then one and only one
of the following possibilities occurs:\footnote{While $0 \in \sigma(C)$ in each of the three cases, it need not be an eigenvalue of $C$. For example, if $\{e_k\}_{k=1}^\infty$ is a complete orthonormal basis of a Hilbert space $\sH$ and $Ce_k := \frac{1}{k}e_k$, then $\{1/k\}_{k=1}^\infty$ is the set of eigenvalues of $C$ but $C^{-1}$ is unbounded and so $0\in \sigma(C)$ but is not an eigenvalue of $C$.}
% COMMENT See https://mathoverflow.net/questions/321482/kernel-of-compact-operators
% COMMENT See https://math.stackexchange.com/questions/2826843/does-each-compact-operator-have-a-non-zero-eigenvector
\begin{inparaenum}[\itshape i\upshape)]
\item $\sigma(C) = \{0\}$;
\item $\sigma(C) = \{0, \lambda_1, \ldots, \lambda_n\}$, where for $1 \leq k \leq n$, one has $\lambda_k \neq 0$, each $\lambda_k$ is an eigenvalue of $C$, and $\dim\Ker(T - \lambda) < \infty$;
\item $\sigma(C) = \{0, \lambda_1, \lambda_2, \ldots\}$, where for $k \geq 1$, one has $\lambda_k \neq 0$, each $\lambda_k$ is an eigenvalue of $C$ with $\dim\Ker(T - \lambda) < \infty$, and $\lim_{k\to\infty}\lambda_k = 0$ (see, for example, Conway \cite[Chapter VII, Section 7, Theorem 7.1, p. 214]{Conway_course_functional_analysis} or Reed and Simon \cite[Chapter VI, Section 5, Theorem VI.15, p. 203]{Reed_Simon_v1}).
\end{inparaenum}
%COMMENT See https://sites.math.washington.edu//~hart/m556/Lecture5.pdf and https://math.stackexchange.com/questions/3291675/spectrum-vs-eigenvalues for an example of an operator T with 0 \in \sigma(T) but 0 is not an eigenvalue of T.

If $\sH$ is a Hilbert space and $A \in \End(\sH)$ is a compact \emph{self-adjoint} operator, then $\sigma(A) \subset \RR$ and $\sH$ has a complete orthonormal basis $\{e_k\}_{k=1}^\infty$ of eigenvectors corresponding to the eigenvalues $\{\lambda_k\}_{k=1}^\infty$ of $A$ (see, for example, Taylor \cite[Appendix A, Section 6, Proposition 6.6, p. 635]{Taylor_PDE1_3rd_edition}, Reed and Simon \cite[Chapter VI, Section 5, Theorem VI.16, p. 203]{Reed_Simon_v1}, or Zimmer \cite[Theorem 3.2.3, p. 58]{Zimmer_essential_results_functional_analysis}). If $A \in \End(\sH)$ is compact self-adjoint, then\footnote{This yields the equality \eqref{eq:Spectral_radius_equals_norm_normal_operator} in the special case of compact self-adjoint operators.} one of $\pm\|A\|_{\End(\sH)}$ is an eigenvalue of $A$ (see Rudin \cite[Chapter 12, Theorem 12.31, p. 329]{Rudin} or Taylor \cite[Appendix A, Section 6, Proposition 6.5, p. 634]{Taylor_PDE1_3rd_edition}).

More generally, if $N \in \End(\sH)$ is a compact \emph{normal} operator, then $\sH$ has a complete orthonormal basis $\{e_k\}_{k=1}^\infty$ of eigenvectors corresponding to the eigenvalues $\{\lambda_k\}_{k=1}^\infty$ of $N$ (see, for example, Radjavi and Rosenthal \cite[Section 1.2, Theorem 1.4, p. 12]{Radjavi_Rosenthal_invariant_subspaces} or Zimmer \cite[Corollary 3.2.9, p. 60]{Zimmer_essential_results_functional_analysis}). To obtain this result as a corollary of the result for compact self-adjoint operators, note that any bounded operator $B \in \End(\sH)$ can be represented uniquely as
\[
  B = B_1 + iB_2, \quad\text{for}\quad B_1 = \frac{1}{2}(B^*+B) \quad\text{and}\quad B_2 = \frac{1}{2i}(B^*-B),
\]
where $B_1$ and $B_2$ are self-adjoint and, moreover, $B$ is normal if and only if $B_1$ and $B_2$ commute (see Zimmer \cite[Corollary 3.2.9, p. 60]{Zimmer_essential_results_functional_analysis}). Writing $N = A_1+iA_2$, where $A_1,A_2$ are compact and self-adjoint, we see that $\sH$ has a complete orthonormal basis $\{e_k\}_{k=1}^\infty$ such that each $e_k$ is an eigenvector of $A_1$ and $A_2$ (see, for example, Zimmer \cite[Corollary 3.2.5, p. 60]{Zimmer_essential_results_functional_analysis}) and hence an eigenvector of $N$.
% COMMENT Spectral theorem for compact normal operators: see Jenkins http://ejenk.com/notes/spectral.pdf
% COMMENT See https://math.stackexchange.com/questions/4372095/proof-of-the-spectral-theorem-for-compact-normal-operators
% COMMENT See https://www.quora.com/How-do-I-prove-that-eigenvectors-corresponding-to-distinct-eigenvalues-of-a-real-symmetric-matrix-are-orthogonal
% COMMENT See https://math.stackexchange.com/questions/778946/prove-that-if-a-is-normal-then-eigenvectors-corresponding-to-distinct-eigenva
% COMMENT See https://math.stackexchange.com/questions/4013227/eigenvectors-of-a-normal-operator-are-an-orthogonal-basis

If $\sH_1,\sH_2$ are Hilbert spaces and $T \in \Hom(\sH_1,\sH_2)$ has dense domain, then (see Brezis \cite[Section 2.6, p. 43]{Brezis} or Conway \cite[Chapter X, Section 1, Definition 1.5, p. 304]{Conway_course_functional_analysis})
\[
  \Dom(T^*) := \{y \in \sH_2: \text{the functional } \Dom(T) \ni x \mapsto \langle Tx,y \rangle_{\sH_2} \text{ is continuous}\},
\]
and $T^* \in \Hom(\sH_2,\sH_1)$ is uniquely defined on $\Dom(T^*)$ by
\[
  \langle Tx,y \rangle_{\sH_2} = \langle x,T^*y \rangle_{\sH_1},
  \quad\text{for all } x \in \Dom(T) \text{ and } y \in \Dom(T^*).
\]  
An unbounded operator in $\Hom(\sH_1,\sH_2)$ is \emph{closed} if its graph is closed in $\sH_1\oplus\sH_2$ and is \emph{closable} if it has a closed extension (see Conway \cite[Chapter X, Section 1, Definition 1.3, p. 304]{Conway_course_functional_analysis}). In particular, if $T \in \Hom(\sH_1,\sH_2)$ has dense domain, then $T^*$ is closed; moreover, $T^*$ has dense domain if and only if $T$ is closable (see Conway \cite[Chapter X, Section 1, Proposition 1.6, p. 305]{Conway_course_functional_analysis}. If $T \in \Hom(\sH_1,\sH_2)$ is densely defined and closed, then $T^{**} = T$ (see Brezis \cite[Section 3.6, Theorem 3.24, p. 72]{Brezis}).

We recall that an unbounded operator $T \in \End(\sH)$ is \emph{symmetric} if $T$ has dense domain and $\langle Tu,v\rangle_\sH = \langle u,Tv\rangle_\sH$ for all $u,v \in \Dom(T)$ (see Conway \cite[Chapter X, Section 2, Definition 2.1, p. 309]{Conway_course_functional_analysis} or Kato \cite[Chapter V, Section 3.3, Equation (3.5), p. 269]{Kato}). If $T \in \End(\sH)$ has dense domain, then $T$ is symmetric if and only if $T \subseteq T^*$ (see Conway \cite[Chapter X, Section 2, Proposition 2.2, p. 309]{Conway_course_functional_analysis}). A densely defined operator $T \in \End(\sH)$ is \emph{self-adjoint} if
$T = T^*$ (see Conway \cite[Chapter X, Section 2, Definition 2.3, p. 309]{Conway_course_functional_analysis}). A densely defined self-adjoint operator $T \in \End(\sH)$ is necessarily closed (see Rudin \cite[Chapter 13, Theorem 13.9, p. 352]{Rudin}.) A symmetric operator $T \in \End(\sH)$ is \emph{essentially self-adjoint} if the closure of $T$ is self-adjoint (see Bikchentaev, Kittaneh, Sal Moslehian, and Seo \cite[Section 8.2, p. 267]{Bikchentaev_Kittaneh_SalMoslehian_Seo_trace_inequalities} or Kato \cite[Chapter V, Section 3.3, p. 269]{Kato}) or, equivalently, if $T^*$ is self-adjoint (see Kato \cite[Chapter V, Section 3.3, Problem 3.10, p. 269]{Kato}).
%COMMENT https://math.stackexchange.com/questions/2227359/definition-of-essentially-self-adjoint-operators

If $T \in \End(\sH)$ is densely defined and self-adjoint (and thus closed) with compact resolvent and $t \in \rho(T)\cap\RR$, then $R(t,T) \in \End(\sH)$ is a compact self-adjoint operator with spectrum related to that of $T$ by \eqref{eq:Spectral_mapping_T_and_RzetaT} or equivalently
\begin{equation}
  \label{eq:mu_in_specRtT_iff_t_plus_1over_mu_in_specT}
  \nu \in \sigma(R(t,T)) \iff t + 1/\nu \in \sigma(T)\cup\{\infty\}.
\end{equation}
If $T$ as above is not bounded, then $\sigma(T)$ is not bounded (since $0 \in \sigma(R(t,T))$ is an accumulation point), so $\sigma(T) = \{\lambda_k\}_{k=1}^\infty \subset \RR$ has $\lim_{k\to\infty}|\lambda_k| = \infty$ and no accumulation points.

If $T \in \End(\sH)$ is an unbounded self-adjoint operator on a Hilbert space $\sH$ with compact resolvent and thus a discrete spectrum $\sigma(T) = \{\lambda_k\}_{k=1}^\infty$ of eigenvalues and complete orthonormal basis $\{e_k\}_{k=1}^\infty$ for $\sH$ (see, for example, Taylor \cite[Appendix A, Section 8, last paragraph, p. 660]{Taylor_PDE1_3rd_edition}), then (as a special case of Dunford and Schwartz \cite[Chapter XII, Section 2, Theorem 3, p. 1192]{Dunford_Schwartz_1963} or Yosida \cite[Chapter XI, Section 6, Theorem 1, p. 313]{Yosida})
\begin{equation}
  \label{eq:Spectral_resolution_T}
  Tv = \sum_{k=1}^\infty \lambda_k e_k\otimes\langle v,e_k\rangle_\sH, \quad\text{for all } v \in \Dom(T),
\end{equation}
where
\begin{equation}
  \label{eq:Spectral_resolution_DomT}
  \Dom(T) := \left\{v \in \sH: \sum_{k=1}^\infty \lambda_k^2|\langle v,e_k\rangle_\sH|^2 < \infty \right\}.
\end{equation}
For $z \in \rho(T)$, one obtains (see Dunford and Schwartz \cite[Chapter XII, Section 2, p. 1195, Definition 5 and Theorem 6, p. 1196]{Dunford_Schwartz_1963})
\begin{equation}
  \label{eq:Spectral_resolution_RzT}
  R(z,T)v = \sum_{k=1}^\infty (\lambda_k-z)^{-1}e_k\otimes\langle v,e_k\rangle_\sH, \quad\text{for all } v \in \sH. 
\end{equation}
(The representation \eqref{eq:Spectral_resolution_RzT} also follows from Radjavi and Rosenblum \cite[Section 1.2, Corollary 1.5, p. 13]{Radjavi_Rosenthal_invariant_subspaces} since $R(z,T)$ is a compact normal operator.)

We shall also rely on the holomorphic functional (or operational) calculus for linear operators on Banach spaces.
%COMMENT See https://en.wikipedia.org/wiki/Holomorphic_functional_calculus and https://en.wikipedia.org/wiki/Resolvent_formalism
Let $\sX$ be a Banach space, $B \in \End(\sX)$ be a \emph{bounded} operator, $\Omega\subset\CC$ be an open neighborhood of $\sigma(B)$, and $U \subset \CC$ be an open subset with closure $\bar U \subset \Omega$ and whose boundary $\Gamma := \partial U$ consists of a finite number of rectifiable Jordan curves, oriented in the positive sense. If $f:\Omega\to\CC$ is holomorphic, then the bounded operator $f(B) \in \End(\sX)$ is defined by the \emph{Dunford integral} (see Dunford and Schwartz \cite[Chapter VII, Section 3, Definition 9, p. 568]{Dunford_Schwartz1} or Yosida \cite[Chapter VIII, Section 7, p. 225]{Yosida})
\begin{equation}
  \label{eq:Dunford_integral_bounded_operator}
  %PF5-8-2025 Check signs
  f(B) := -\frac{1}{2\pi i}\oint_\Gamma f(z)R(z,B)\,dz.
\end{equation}
(We include the negative sign since Dunford and Schwartz \cite[Chapter VII, Section 3, Definition 1, p. 566]{Dunford_Schwartz1} define $R(z,B) = (z-B)^{-1}$ and not $(B-z)^{-1}$ as in \eqref{eq:Taylor_PDE1_A-8-33}; taking $f(z) = z$ and a Jordan domain $U$ such that $\sigma(B) \subset U$ in \eqref{eq:Dunford_integral_bounded_operator} yields $f(B) = B$.) If $T \in \End(\sX)$ is an \emph{unbounded} operator (with bounded resolvent $R(z,T)$) and $f$ is holomorphic on an open neighborhood $\Omega \subset \CC$ of $\sigma(T)$ and $\infty$, then $f(T)$ can be defined by (see Dunford and Schwartz \cite[Chapter VII, Section 9, Definition 1, p. 600 and Theorem 4, p. 601]{Dunford_Schwartz1}),
\begin{equation}
  \label{eq:Dunford_integral_unbounded_operator}
  f(T) := f(\infty)I - \frac{1}{2\pi i}\oint_\Gamma f(z)R(z,T)\,dz.
\end{equation}
If $T \in \End(\sX)$ is a closed unbounded operator with spectrum $\sigma(T) = \sigma_1\sqcup\sigma_2$, where $\sigma_1\subset\CC$ is bounded, and $\Gamma \subset \rho(T)$ is a rectifiable, simple closed curve (or, more generally, a finite number of such curves) that enclose an open set containing
$\sigma_1$ in its interior and $\sigma_2$ in its exterior, then $T$ decomposes as $T_1\oplus T_2$ with respect to a direct sum of Banach spaces $\sX = \sX_1\oplus\sX_2$ such that $\sigma(T_i) = \sigma_i$ for $i=1,2$ and
\begin{equation}
  \label{eq:Riesz_projection_closed_unbounded_operator_Banach_space}
  \Pi(\sigma_1) := -\frac{1}{2\pi i}\oint_\Gamma R(z,T)\,dz
\end{equation}
is continuous projection onto $\sX_1$ and $\Pi(\sigma_2) := 1-\Pi(\sigma_1)$ is a continuous projection onto $\sX_2$ (see Dunford and Schwartz \cite[Chapter VII, Section 3, Theorem 20 and Corollary 21, p. 574 and the last
equation in the proof of Theorem 18, p. 573]{Dunford_Schwartz1} or Kato \cite[Chapter III, Section 6.4, Theorem 6.17, p. 178]{Kato}).

For $T$ as in the preceding paragraph, if we assume that $\Omega_1 \Subset \CC$ and $\Omega_2 \subset \CC$ are connected open subsets such that $\Omega_1\cap\Omega = \varnothing$, then we may choose a holomorphic function $f$ on $\Omega = \Omega_1\sqcup\Omega_2$ by setting $f \equiv 1$ on $\Omega_1$ and $f \equiv 0$ on $\Omega_2$ in \eqref{eq:Dunford_integral_unbounded_operator} to obtain the expression \eqref{eq:Riesz_projection_closed_unbounded_operator_Banach_space} for the \emph{Riesz projection}. If we assume that $\sX$ is a Hilbert space $\sH$ and substitute the expression \eqref{eq:Spectral_resolution_RzT} for $R(z,T)$, we see that
\[
  \Pi(\sigma_1)
  = -\frac{1}{2\pi i}\oint_\Gamma R(z,T)\,dz
  = -\sum_{\lambda_k \in \sigma_1} \frac{1}{2\pi i}\oint_\Gamma(\lambda_k-z)^{-1}\,dz\,e_k\otimes\langle \cdot,e_k\rangle_\sH,
\]
and so $\Pi(\sigma_1)$ is the \emph{orthogonal} projection onto the sum of the eigenspaces of $T$ corresponding to the eigenvalues of $T$ in $\sigma_1$:
\begin{equation}
  \label{eq:Riesz_projection_closed_unbounded_operator_compact_resolvent_Hilbert_space}
  \Pi(\sigma_1) = \sum_{\lambda_k \in \sigma_1} e_k\otimes\langle \cdot,e_k\rangle_\sH.
\end{equation}
Haase \cite[Chapter 2]{Haase_2006} develops an extension of the definition \eqref{eq:Dunford_integral_bounded_operator} to unbounded operators that are \emph{sectorial}. Batty \cite{Batty_2009} provides a nice survey of functional calculus for unbounded operators, based in part on Haase \cite{Haase_2005pems, Haase_2005jlms, Haase_2006}. deLaubenfels \cite{deLaubenfels_1993} and Nielsen \cite{Nielsen_2020} also include useful surveys of functional calculus for unbounded operators.

% COMMENT Use https://math.stackexchange.com/questions/4418596/prove-that-the-resolvent-set-is-open
% COMMENT https://proofwiki.org/wiki/Resolvent_Set_of_Bounded_Linear_Operator_is_Open
% COMMENT https://mathoverflow.net/questions/294002/reference-request-the-resolvent-is-analytic-in-the-resolvent-set
% COMMENT https://projecteuclid.org/journalArticle/Download?urlid=10.1215%2Fkjm%2F1250520831

Returning to the setting of the Dunford integral \eqref{eq:Dunford_integral_bounded_operator}, where $T\in\End(\sX)$ is a bounded operator on a Banach space $\sX$, suppose that $f:\CC \supset \Omega\to\CC$
is a holomorphic function on an open neighborhood $\Omega$ of $\sigma(T)$. By the Spectral Mapping Theorem (see, for example, Yosida \cite[Chapter VIII, Section 7, Corollary 1, p. 227]{Yosida}),
\begin{equation}
  \label{eq:Spectral_Mapping_Theorem}
  f(\sigma(T))=\sigma(f(T)).
\end{equation}
If $\sH$ is a Hilbert space and $T \in \End(\sH)$ is bounded and normal, then $f(T)$ is also bounded and normal by the forthcoming \eqref{eq:Function_of_normal_operator_normal}, so by \eqref{eq:Spectral_radius_equals_norm_normal_operator} and \eqref{eq:Spectral_Mapping_Theorem},
\[
  \|f(T)\|_{\End(\sH)}=\sup_{\lambda'\in\sigma(f(T))}|\lambda'|=\sup_{\lambda\in\si(T)} |f(\lambda)|.
\]
If $f$ is positive and increasing on $(0,\|T\|_{\End(\sH)})$, then the preceding equalities imply that
\begin{equation}
\label{eq:NormOf_Function_of_Operator}
\|f(T)\|_{\End(\sH)} \le |f(\|T\|_{\End(\sH)})|.
\end{equation}
Compare Rudin \cite[Chapter 12, Equation (3), p. 325]{Rudin}.

Suppose now that $S, T \in \End(\sX)$ are bounded operators on a Banach space $\sX$. If $f:\CC\supset\Omega\to\CC$ is a holomorphic function on an open neighborhood $\Omega$ of $\sigma(S)\cup\sigma(T)$ and $ST = TS$, then
\begin{equation}
  \label{eq:Functions_of_commuting_operators}
  f(S)f(T) = f(T)f(S).
\end{equation}
To see this, observe that the operators $S-z$ and $T-z$ commute and hence the resolvents $R(z,S)$ and $R(z,T)$ commute for all $z\in\CC \less (\sigma(S)\cup\sigma(T))$. Therefore, by applying \eqref{eq:Dunford_integral_bounded_operator} for suitable curves $\Gamma_S \subset \rho(S)$ and $\Gamma_T \subset \rho(T)$ and Fubini's Theorem, we obtain
\begin{align*}
  f(S)f(T)
  &= \frac{1}{(2\pi i)^2} \oint_{\Gamma_S}f(z)R(z,S)\,dz \oint_{\Gamma_T}f(w)R(w,T)\,dw
  \\
  &= \frac{1}{(2\pi i)^2} \oint_{\Gamma_S}\left(\oint_{\Gamma_T} f(z)f(w)R(z,S)R(w,T)\,dw\right)\,dz 
  \\
  &= \frac{1}{(2\pi i)^2} \oint_{\Gamma_T}\left(\oint_{\Gamma_S} f(w)f(z)R(w,T)R(z,S)\,dz\right)\,dw
  \\
  &= \frac{1}{(2\pi i)^2} \oint_{\Gamma_T}f(w)R(w,T)\,dw \oint_{\Gamma_S}f(z)R(z,S)\,dz
  = f(T)f(S),
\end{align*}
which verifies \eqref{eq:Functions_of_commuting_operators}.

If $\sH$ is a Hilbert space and $T \in \End(\sH)$ is a bounded operator, then $f(T^*) = f(T)^*$ by Dunford and Schwartz \cite[Chapter VII, Section 3, Theorem 10 (d), p. 568]{Dunford_Schwartz1}. In particular, if $T$ is normal, so $TT^* = T^*T$, then \eqref{eq:Functions_of_commuting_operators} implies that
\[
   f(T)f(T)^* = f(T)f(T^*) = f(T^*)f(T) = f(T)^*f(T),
\]
and thus
\begin{equation}
  \label{eq:Function_of_normal_operator_normal}
  f(T)f(T)^* = f(T)^*f(T),
\end{equation}
that is, $f(T)$ is normal.

\begin{rmk}[Estimates for resolvent operators for elliptic pseudodifferential operators]
\label{rmk:Estimates_resolvent_operators_elliptic_pseudodifferential_operators}
The formula \eqref{eq:Spectral_radius_equals_norm_resolvent_operator_normal_T} for the norm of the resolvent operator for a normal operator on a Hilbert space will suffice for our current applications. Estimates for resolvents of \emph{elliptic pseudodifferential} operators are provided by Boo\ss-Bavnbek, Chen, Lesch, and Zhu \cite[Section 2.4.1, Equation (2.10), p. 59]{Booss-Bavnbek_Chen_Lesch_Zhu_2012} and Boo\ss-Bavnbek and Lesch \cite{Booss-Bavnbek_Lesch_2009}, based on Seeley \cite[Section 6, Corollary 1, p. 298]{Seeley_1967a} (see also Seeley \cite{Seeley_1969ajm, Seeley_1969cime}), Shubin \cite[Chapter II, Section 9.4, Theorem 9.3 and Corollary 9.2, p. 86]{Shubin}, Haase \cite{Haase_2006}, and Edmunds and Evans \cite{Edmunds_Evans_elliptic_differential_operators_spectral_analysis}. Seeley allows $E$ to be a complex vector bundle, but Shubin assumes that $E$ is trivial.

Such resolvent estimates are required to show that an elliptic pseudodifferential operator is sectorial in the sense of Haase \cite[Section 2.1, p. 19]{Haase_2006}. See also Sell and You \cite[Section 3.1, Lemma 31.6, p. 65, or Section 3.4, Theorem 34.1, p. 70, or Section 3.6, Theorem 36.2 (3), p. 79]{Sell_You_2002}.

Estimates for resolvents of \emph{Dirac} operators are obtained by Axelsson, Keith, and McIntosh \cite{Axelsson_Keith_McIntosh_2006}, Boo\ss-Bavnbek and Wojciechowski \cite{Booss-Bavnbek_Wojciechowski_elliptic_boundary_problems_dirac_operators}, and Danilov \cite{Danilov_1995}.
\qed\end{rmk}

%COMMENT https://math.stackexchange.com/questions/4249614/perturbation-of-the-spectrum
%COMMENT Convergence of spectrum of $T$ to $T'$ as $T'' \to 0$: see, for example, Kato \cite[Chapter IV, Section 3.4, Theorem 3.16, p. 212]{Kato}. See also Donoghue \cite{Donoghue_1965}, Friedrichs \cite{Friedrichs_spectral_theory_operators_hilbert_space} and
%COMMENT See https://link.springer.com/chapter/10.1007/978-3-642-73782-4_5, chapter on Spectral Perturbation in ``Vibration and Coupling of Continuous Systems'', J. Sanchez Hubert, E. Sanchez Palencia, Springer, 1989.

If $T \in \End(\sH)$ is an unbounded linear operator on a Hilbert space $\sH$ with a discrete spectrum $\sigma(T) = \{\lambda_k\}_{k=1}^\infty$ of eigenvalues and complete orthonormal basis $\{e_k\}_{k=1}^\infty$ for $\sH$ of eigenvectors of $T$ with $Te_k = \lambda_ke_k$ for all $k\geq 1$ and
\begin{equation}
  \label{eq:Isomorphism_Hilbert_spaces}
  \Upsilon \in \Isom(\sH,\sK)
\end{equation}  
is an isomorphism of $\sH$ onto a Hilbert space $\sK$, then 
\begin{equation}
  \label{eq:Transformed_basis_eigenvectors}
  \{f_k\}_{k=1}^\infty := \{\Upsilon e_k\}_{k=1}^\infty
\end{equation}
is a complete orthonormal basis for $\sK$ of eigenvectors of the unbounded linear operator
\begin{equation}
  \label{eq:Transformed_operator}
  S := \Upsilon T\Upsilon^{-1} \in \End(\sK)
\end{equation}
with spectrum of eigenvalues $\sigma(S) = \{\lambda_k\}_{k=1}^\infty$ since
\[
  Sf_k = (\Upsilon T\Upsilon^{-1})\Upsilon e_k = \Upsilon Te_k = \lambda_k\Upsilon e_k = \lambda_kf_k,
  \quad\text{for all } k \geq 1.
\]
If $\sK$ is only a vector space and $\Upsilon$ in \eqref{eq:Isomorphism_Hilbert_spaces} is only an isomorphism of vector spaces, then we may define an inner product on $\sK$ by setting
\begin{equation}
  \label{eq:Pull_back_inner_product}
  \langle w,w' \rangle_\sK
  :=
  \langle \Upsilon^{-1}w, \Upsilon^{-1}w' \rangle_\sH,
  \quad\text{for all } w, w' \in \sK.
\end{equation}
When $\sK$ is given the pullback inner product, $\langle \cdot,\cdot \rangle_\sK = (\Upsilon^{-1})^*\langle \cdot,\cdot \rangle_\sH$, the isomorphism $\Upsilon$ in \eqref{eq:Isomorphism_Hilbert_spaces} of vector spaces becomes an isomorphism of inner product spaces. In particular, $\sK$ is a complete inner product space and thus a Hilbert space. We shall often appeal to the preceding elementary observations.

Let $\sX$ and $\sY$ be Banach spaces and $T \in \Hom(\sX,\sY)$ be a bounded linear operator. Following Abramovich and Aliprantis \cite[Definition 4.37, p. 156]{Abramovich_Aliprantis_2002}, one calls $T$ \emph{Fredholm} if
\begin{inparaenum}[\itshape i\upshape)]
\item $\dim\Ker T < \infty$ and
\item $\dim\Coker T < \infty$, where $\Coker T := \sY/\Ran T$.
\end{inparaenum}  
According to \cite[Corollary 2.17, p. 76 or Lemma 4.38, p. 156]{Abramovich_Aliprantis_2002}, a bounded linear operator $T \in \Hom(\sX,\sY)$ has closed range if $\dim\Coker T < \infty$. Moreover, \cite[Theorem 2.13, p. 74]{Abramovich_Aliprantis_2002} asserts that if $T \in \Hom(\sX,\sY)$ is a bounded linear operator, then $(\Ran T)^\perp = \Ker T^*$ and if $T$ has closed range, then $\dim(\sY/\Ran T)^* = \dim(\Ran T)^\perp = \dim\Ker T^*$.
By \cite[Corollary 2.15, p. 75]{Abramovich_Aliprantis_2002}, a bounded linear operator $T \in \Hom(\sX,\sY)$ has closed range if and only if there is a constant $c > 0$ such that for each $y \in \Ran T$, there exists $x \in \sX$ with $y=Tx$ and $\|x\|_\sX \leq c\|y\|_\sY$.

Suppose now that $\sH$ is a Hilbert space and that $T \in \End(\sH)$ is a closed, unbounded linear operator with dense domain. The unbounded linear operator $T^*T \in \End(\sH)$ is non-negative and self-adjoint with dense domain. If $0<\lambda_1<\inf\{\lambda: \lambda \in \sigma_\ess(T^*T)\}$ is the first positive eigenvalue of $T^*T$, then the Max--Min Principle (see the forthcoming Theorem \ref{thm:Max_min_principle}) implies that
\[
  \lambda_1 = \inf_{\{x \in \Dom(T^*T) \cap (\Ker T)^\perp\}} \frac{\langle x,T^*Tx\rangle_\sH}{\|x\|_\sH^2},
\]
and thus $\|Tx\|_\sH \geq \sqrt{\lambda_1}\|x\|_\sH$ for all $x \in (\Ker T)^\perp$. In particular, such an operator $T$ has closed range by Abramovich and Aliprantis \cite[Corollary 2.15, p. 75]{Abramovich_Aliprantis_2002}.

We shall occasionally appeal to results on perturbation of spectra. If $T$ is a self-adjoint linear unbounded operator on Hilbert space $\sH$ and $S$ is a bounded symmetric linear operator on $\sH$, then the operator $T+S$ is self-adjoint and
\begin{equation}
  \label{eq:Distance_spectra_T+S_and_T}
  \dist\left(\sigma(T+S), \sigma(T)\right)
  \leq \|S\|_{\End(\sH)},
\end{equation}
that is,
\begin{align*}
  \sup_{\lambda\in \sigma(T+S)} \dist\left(\lambda, \sigma(T) \right)
  &\leq \|S\|_{\End(\sH)} \quad\text{and}
  \\
  \sup_{\lambda\in \sigma(T)} \dist\left(\lambda, \sigma(T+S) \right)
  &\leq \|S\|_{\End(S)},
\end{align*}
according to Kato \cite[Chapter V, Section 4.3, Theorem 4.10, p. 291]{Kato}.

\section{Max-min principle}
\label{sec:Max_min_principle}
We begin with the

\begin{defn}[Essential and discrete spectrum of a closed densely defined operator on a Banach space]
\label{defn:Essential_and_discrete_spectrum_self-adjoint_operator_Hilbert_space}  
(See Edmunds and Evans \cite[Chapter 1, Section 4, p. 37 and Chapter 9, Section 1, p. 421]{Edmunds_Evans_spectral_theory_differential_operators_ed2} for a comparison of different definitions of essential spectra in the literature; see also Kato \cite[Chapter IV, Section 5.6, p. 243 or Chapter X, Section 1.2, Remark 1.11, p. 520]{Kato} for unbounded operators and Reed and Simon \cite[Section XIII.1, Definition, p. 236]{Reed_Simon_v1} for bounded operators.)
%PF6-21-2025 Probable cut
% Let $A$ be a self-adjoint operator on a Hilbert space $\sH$ and let
% \[
%   A = \int_{-\infty}^\infty \lambda\,dE(\lambda),
% \]
% be its spectral representation \cite[Chapter IV, Section 5 or Chapter X, Section 1.1, p. 516]{Kato}, where $\{E(\lambda)\}$ is the right-continuous spectral family associated with $A$,
% \[
%   P(\lambda) := E(\lambda) - E(\lambda - 0)
% \]
% and $P(\lambda) \neq 0$ if and only if $\lambda$ is an eigenvalue of $A$, in which case $P(\lambda)$ is the
% orthogonal projection on the associated eigenspace. The spectral family $\{E(\lambda)\}$ determines a projection-valued spectral measure $E(S)$ defined for all Borel subsets $S \subset \RR$ such that if $S = (a,b] \subset \RR$, then $E(S) := E(b) - E(a)$ \cite[Chapter X, Section 1.2, p. 518]{Kato}.
%
% If $\sigma(A) \subset \RR$ denotes the spectrum of $A$, then it may be expressed as a disjoint union, 
% \[
%   \sigma(A) = \sigma_\ess(A) \sqcup \sigma_\disc(A),
% \]
% where $\lambda \in \sigma(A)$ is in the \emph{essential spectrum} $\sigma_\ess(A)$ if
% \[
%   \dim \left[E(\lambda+\eps) - E(\lambda-\eps)\right] = \infty, \quad\text{for any } \eps > 0.
% \]
% The \emph{discrete spectrum}, $\sigma_\disc(A) := \sigma(A) \less \sigma_\ess(A)$, is a set of isolated points in $\RR$, comprising
% %PF5-2-2025 Why exactly?
% eigenvalues with finite multiplicities.
If $T$ is a closed densely defined linear operator on a complex Banach space $\sX$, then we define the \emph{essential spectrum} of $T$ to be
\begin{equation}
  \label{eq:Essential_spectrum}
  \sigma_\ess(T) := \left\{\lambda \in \CC: T - \lambda \notin \Fred(\sX)\right\}.
\end{equation}
When $\sH$ is a \emph{Hilbert} space and $T \in \End(\sH)$ is a densely defined self-adjoint unbounded operator on $\sH$, then the different definitions of essential spectrum coincide by Edmunds and Evans \cite[Chapter 9, Section 1, Theorem 1.6, p. 424]{Edmunds_Evans_spectral_theory_differential_operators_ed2} and $\lambda \in \sigma(T)\less \sigma_\ess(T)$ if and only if $\lambda$ is an isolated eigenvalue of finite multiplicity. Following Reed and Simon \cite[Chapter VII, Section 3, Theorem VII.9, p. 236]{Reed_Simon_v1}, we denote
\begin{equation}
  \label{eq:Discrete_spectrum}
  \sigma_\disc(T)
  :=
  \left\{\lambda\in\sigma(T): \lambda \text{ is an isolated eigenvalue of finite multiplicity}\right\},
\end{equation}
and observe that $\sigma(T) = \sigma_\disc(T) \sqcup \sigma_\ess(T)$ in this case.
\qed\end{defn}

%COMMENT See https://en.wikipedia.org/wiki/Min-max_theorem
We recall the well-known

\begin{thm}[Max-min principle]
\label{thm:Max_min_principle}
(See Edmunds and Evans \cite[Chapter 11, Section 1, Theorem 1.2, p. 498]{Edmunds_Evans_spectral_theory_differential_operators_ed2}, Reed and Simon \cite[Section XIII.1, Theorem XIII.1, p. 76]{Reed_Simon_v4}, Teschl \cite[Section 4.3, Theorem 4.10, p. 119]{Teschl_mathematical_methods_quantum_mechanics}, Strauss \cite[Section 11.6, Theorem 2, p. 324]{Strauss} for a version specific to the Laplace operator with Dirichlet boundary or Neumann conditions on a bounded domain in Euclidean space, and Tao \cite[Theorem 1.3.2, p. 42]{Tao_topics_random_matrix_theory}.)   
Let $A$ be a self-adjoint operator on a Hilbert space $\sH$ that is bounded below, that is, $A \geq c\,\id_\sH$ for some $c \in \RR$. For each $k \geq 1$, define
\[
  \mu_k = \sup_{\varphi_1,\ldots,\varphi_{k-1}} \mu_{k*}(\varphi_1,\ldots,\varphi_{k-1}),
\]
where
\[
  \mu_{m*}(\varphi_1,\ldots,\varphi_m)
  :=
  \inf_{\begin{subarray}{c}\{\psi\in\Dom(A)\cap[\varphi_1,\ldots,\varphi_m]^\perp:
        \\ \|\psi\|_\sH = 1\} \end{subarray}} \langle \psi, A\psi \rangle_\sH,
\]
and
\[
  [\varphi_1,\ldots,\varphi_m]^\perp
  :=
  \{\psi \in \sH: \langle\psi,\varphi_j\rangle_\sH = 0, \text{ for } j=1,\ldots,m\}.
\]
(We set $\{\varphi_1,\ldots,\varphi_m\} := \varnothing$ if $m =0$.) Then one and only one of the following hold:
\begin{enumerate}
\item There are $k$ eigenvalues (repeated according to their multiplicity) below the bottom of the essential spectrum and $\mu_k$ is the $k$-th eigenvalue counting multiplicity, or
\item $\mu_k$ is the bottom of the essential spectrum, $\mu_k = \inf\{\lambda: \lambda\in\sigma_\ess(A)\}$, and there are most $k-1$ eigenvalues below $\mu_k$, counting multiplicity.  
\end{enumerate}
\end{thm}

When $\sH$ is a Hilbert space of $L^2$ functions on a domain in Euclidean space, the vectors $\varphi_l \in \sH$, for $l=1,\ldots,k$, are called \emph{variational} or \emph{trial} or \emph{comparison} functions (see Lieb and Loss \cite[Chapter 12, Theorem 12.1, pp. 300-301]{LiebLoss}). The subset $\{\varphi_1,\ldots,\varphi_k\} \subset \sH$ is not assumed to be linearly independent.

\section[Weyl's asymptotic formula for the eigenvalues of an elliptic operator]{Weyl's asymptotic formula for the eigenvalues of an elliptic pseudodifferential operator on sections of a vector bundle}
\label{sec:Weyl_asymptotic_formula_eigenvalues_elliptic_operator}
% COMMENT This kind of result should be a consequence of a version of Weyl's asymptotic formula,
% See https://mathoverflow.net/questions/291990/gaps-in-the-spectrum-of-laplace-beltrami-operators
% See https://mathoverflow.net/questions/199780/estimates-of-eigenvalues-of-elliptic-operators-on-compact-manifolds
The results described in this section generalize the asymptotic formula due to Weyl \cite{Weyl_1912} for eigenvalues of the Laplacian on functions over a domain in Euclidean space. We use the forthcoming Theorem \ref{thm:Weyl_asymptotic_formula_eigenvalues_elliptic_pseudodifferential_operator_sections_vector_bundle} to give quantitative gap conditions as in the hypothesis \eqref{eq:Spectral_gap_T} of Theorem \ref{mainthm:Donaldson_1996jdg_3_Hilbert_space} are satisfied.
%PF1-28-2025 Allow \lambda = 0 to be an eigenvalue 

Let $(M,g)$ be a closed\footnote{The results in this section do not require $M$ to be orientable, but if $M$ is non-orientable one must replace the Riemannian volume form by a choice of positive density in order to define integrals and $L^2$ spaces.} orientable smooth Riemannian manifold of dimension $d\geq 1$ and $0 \leq p \leq d$ be an integer. The Laplace operator on smooth $p$-forms,
\begin{equation}
  \label{eq:Laplace_operator_p-forms}
  \Delta = d^*d + dd^*:\Omega^p(M,\RR)\to \Omega^p(M,\RR),
\end{equation}
is $L^2$ self-adjoint (with respect to the Riemannian volume form).

For any unbounded sequence $\{\lambda_k\}_{k=1}^\infty \subset \CC$, ordered such that
\begin{equation}
  \label{eq:Lambda_sequence}
  |\lambda_1| \leq |\lambda_2| \leq \cdots < \infty,
\end{equation}
it is convenient to define the \emph{counting function},
\begin{equation}
  \label{eq:N_lambda}
  N(\lambda) := \sum_{\{\lambda_k:\, |\lambda_k| \leq\lambda\}} 1, \quad\text{for } \lambda > 0.
\end{equation}
For a function $f:[0,\infty)\to[0,\infty)$ and constant $q\geq 1$, one writes
%COMMENT See https://en.wikipedia.org/wiki/Big_O_notation
\begin{equation}
  \label{eq:Asymptotic_meaning}
  f(t) \sim t^q \quad\text{as } t\to\infty \iff \lim_{t\to\infty} \frac{f(t)}{t^q} = 1.
\end{equation}  
One has the

\begin{thm}[Weyl's asymptotic formula for the eigenvalues of the Laplace operator on $p$-forms]
\label{thm:Weyl_asymptotic_formula_eigenvalues_Laplace_operator_forms}  
(See Chavel \cite[Chapter I, Section 3, p. 9 or Chapter VI, Section 4, p. 155]{Chavel} for the Laplace operator on functions and Dodziuk \cite[Appendix, Section B, Theorem B.2, p. 339]{Chavel} for the Laplace operator on $p$-forms.)
Let $(M,g)$ be a closed orientable smooth Riemannian manifold of dimension $d \geq 2$ and $0 \leq p \leq d$ be an integer. Then $\sigma(\Delta) = \{\lambda_k\}_{k=0}^\infty \subset \RR$, where
\[
  0 = \lambda_0 < \lambda_1 \leq \lambda_2 \leq \cdots < \infty,
\]
and each $\lambda_k$ is an eigenvalue that is repeated according to its finite multiplicity. If $p=0$, then the counting function $N(\Delta;\cdot) = N(\cdot)$ in \eqref{eq:N_lambda} obeys
\begin{equation}
  \label{eq:N_lambda_Laplace_operator_functions_lambda_to_infty}
  N(\Delta;\lambda) \sim \frac{\omega_d}{(2\pi)^d} \Vol(M,g) \lambda^{d/2} \quad\text{as } \lambda \to \infty,
\end{equation}
and the eigenvalues obey
\begin{equation}
  \label{eq:Lambda_k_Laplace_operator_functions_k_to_infty}
  \lambda_k^{d/2} \sim  \frac{(2\pi)^d}{\omega_d\Vol(M,g)}k \quad\text{as } k \to \infty,
\end{equation}
where $\omega_d$ is the volume of the unit ball in $\RR^d$. If $1\leq p \leq d$, then 
\begin{equation}
  \label{eq:N_lambda_Laplace_operator_p-forms_lambda_to_infty}
  N(\Delta;\lambda) \sim \binom{d}{p}\frac{\omega_d}{(4\pi)^{d/2}\Gamma((d/2)+1)} \Vol(M,g) \lambda^{d/2} \quad\text{as } \lambda \to \infty.
\end{equation}
\end{thm}

Theorem \ref{thm:Weyl_asymptotic_formula_eigenvalues_Laplace_operator_forms} (for $p=0$) has been generalized from the Laplace operator to elliptic pseudodifferential operators on functions or smooth sections of smooth vector bundles. For elliptic pseudodifferential operators on functions, one has the

\begin{thm}[Weyl's asymptotic formula for the eigenvalues of a normal elliptic pseudodifferential operator on functions]
\label{thm:Weyl_asymptotic_formula_eigenvalues_elliptic_pseudodifferential_operator_functions}  
(See Shubin \cite[Section 13, Proposition 13.1 and its proof, p. 117, and Section 15, Theorem 15.2 and Problem 15.2, p. 130]{Shubin} and, for elliptic differential operators of even order, Safarov and Vassiliev \cite[Section 1.2.1, Theorem 1.2.1, p. 9 and Section 1.1.1, Equation (1.1.3), p. 1]{Safarov_Vassiliev_asymptotic_distribution_eigenvalues_partial_differential_operators}.)  
Let $(M,g)$ be a closed orientable smooth Riemannian manifold of dimension $d \geq 2$ and $A:C^\infty(M,\CC) \to C^\infty(M,\CC)$ be an $L^2$ self-adjoint elliptic pseudodifferential operator of order $m > 0$ with positive principal symbol,
\begin{equation}
  \label{eq:Shubin_15-1}
  a_m(x,\xi) > 0, \quad\text{for all } (x,\xi) \in T^*M\,\setminus\, 0.
\end{equation}
For $t > 0$, denote
\begin{equation}
  \label{eq:Shubin_13-16}
  V(a_m,g;t) := \frac{1}{(2\pi)^d} \int_{a_m(x,\xi)<t} dx\,d\xi.
\end{equation}
Then $\sigma(A) = \{\lambda_k\}_{k=1}^\infty \subset \RR$, where\footnote{While our notation suggests that $0$ may be an eigenvalue, we do not assume that it is present in the spectrum.}
\[
  0 = \lambda_0 < \lambda_1 \leq \lambda_2 \leq \cdots < \infty,
\]
and each $\lambda_k$ is an eigenvalue that is repeated according to its finite multiplicity. The counting function $N(A;\cdot) = N(\cdot)$ in \eqref{eq:N_lambda} obeys 
\begin{equation}
  \label{eq:Shubin_15-9}
  N(A;\lambda) \sim V(a_m,g;1)\lambda^{d/m} \quad\text{as } \lambda \to \infty,
\end{equation}
and the eigenvalues obey
\begin{equation}
  \label{eq:Shubin_15-10}
  \lambda_k \sim  \frac{1}{V(a_m,g;1)^{m/d}}k^{m/d} \quad\text{as } k \to \infty.
\end{equation}
If $A:C^\infty(M,\CC) \to C^\infty(M,\CC)$ is an $L^2$ normal elliptic differential operator of order $m \geq 1$, then $\sigma(A) = \{\lambda_k\}_{k=0}^\infty \subset \CC$, where each $\lambda_k$ is an eigenvalue that is repeated according to its finite multiplicity, and if $V(t)$ for $t > 0$ in \eqref{eq:Shubin_13-16} is redefined as
\begin{equation}
  \label{eq:Shubin_15-14}
  V(a_m,g;t) := \frac{1}{(2\pi)^d} \int_{|a_m(x,\xi)|<t} dx\,d\xi,
\end{equation}
then the asymptotic relation \eqref{eq:Shubin_15-9} holds.
\end{thm}

The function $V(A;\cdot)$ in \eqref{eq:Shubin_13-16} is the volume in $T^*M$ of all points $(x,\xi)$ such that $a_m(x,\xi) < t$, multiplied by $(2\pi)^{-d}$, where the volume is given in $T^*M$ by the measure induced by the canonical symplectic structure on the cotangent bundle.

\begin{rmk}[Other versions, refinements, and generalizations of Theorem \ref{thm:Weyl_asymptotic_formula_eigenvalues_elliptic_pseudodifferential_operator_functions}]
See Battisti, Borsero, and Coriasco \cite[Theorem 1, p. 799]{Battisti_Borsero_Coriasco_2016}, Duistermaat and Guillemin \cite{Duistermaat_Guillemin_1975}, Guillemin \cite{Guillemin_1985}, Helffer \cite{Helffer_theorie_spectrale_operateurs_globalement_elliptiques}, H\"ormander \cite{Hormander_1968}, and H\"ormander \cite[Corollary 29.1.6, p. 259]{Hormander_v4} for other versions of Theorem \ref{thm:Weyl_asymptotic_formula_eigenvalues_elliptic_pseudodifferential_operator_functions}. Further refinements and generalizations of Theorem \ref{thm:Weyl_asymptotic_formula_eigenvalues_elliptic_pseudodifferential_operator_functions} are provided by Ivrii \cite{Ivrii_microlocal_analysis_sharp_spectral_asymptotics_applications_I,
    Ivrii_microlocal_analysis_sharp_spectral_asymptotics_applications_II, Ivrii_microlocal_analysis_sharp_spectral_asymptotics_applications_III, Ivrii_microlocal_analysis_sharp_spectral_asymptotics_applications_IV, Ivrii_microlocal_analysis_sharp_spectral_asymptotics_applications_V} and Safarov and Vassiliev \cite{Safarov_Vassiliev_asymptotic_distribution_eigenvalues_partial_differential_operators} and references therein.
\qed\end{rmk}  

\begin{thm}[Weyl's asymptotic formula for the eigenvalues of an elliptic pseudodifferential operator on sections of a vector bundle]
\label{thm:Weyl_asymptotic_formula_eigenvalues_elliptic_pseudodifferential_operator_sections_vector_bundle}
(See Ivrii \cite[Section 2, Theorem 0.1, p. 101]{Ivrii_1982} and \cite[Section 0, Theorem 0.1, p. 2]{Ivrii_precise_spectral_asymptotics_elliptic_operators} and Shubin \cite[Section 13, Proposition 13.1 and its proof, p. 117]{Shubin}.)
Let $(E,H)$ be a smooth Hermitian vector bundle over a closed orientable smooth Riemannian manifold $(M,g)$ of dimension $d\geq 2$ and let $A:C^\infty(M,E) \to C^\infty(M,E)$ be an $L^2$ self-adjoint elliptic pseudodifferential operator of order $m > 0$ with principal symbol $a_m$. Then $\sigma(A) = \{\lambda_k\}_{k=0}^\infty \subset \RR$, each $\lambda_k$ is an eigenvalue that is repeated according to its finite multiplicity, and the counting function $N(A;\cdot) = N(\cdot)$ in \eqref{eq:N_lambda} obeys
\begin{subequations}
  \label{eq:Ivrii_1982_or_1984_0-1_rescaled_counting_function_and_eigenvalues}
  \begin{align}
    \label{eq:Ivrii_1982_or_1984_0-1_rescaled}
    N(A;\lambda) \sim c_0(A) \lambda^{d/m} \quad\text{as } \lambda \to \infty,
    \\
    \label{eq:Ivrii_1982_or_1984_0-1_rescaled_eigenvalues}
  \lambda_k \sim  \frac{1}{c_0(A)^{m/d}}k^{m/d} \quad\text{as } k \to \infty,
  \end{align}    
\end{subequations}
where $c_0(A)$ is positive constant determined by the symbol of $A$ (see Ivrii \cite[Section 2, Equation (0.4), p. 101]{Ivrii_1982} or \cite[Section 0, Equation (0.4), p. 3]{Ivrii_precise_spectral_asymptotics_elliptic_operators}).
\end{thm}

\begin{rmk}[Versions of Theorem \ref{thm:Weyl_asymptotic_formula_eigenvalues_elliptic_pseudodifferential_operator_sections_vector_bundle} for Laplace-type operators on sections of a vector bundle]
See Berline, Getzler, and Vergne \cite[Section 2.6, Corollary 2.43, p. 95]{BerlineGetzlerVergne} for certain generalized Laplace operators acting on smooth sections of a smooth vector bundle over a closed smooth Riemannian manifold. See Savale \cite{SavaleThesis, Savale_2014, Savale_2023} for additional results.
\qed\end{rmk}

\begin{rmk}[Seeley's asymptotic formula for the eigenvalues of an elliptic pseudodifferential operator on sections of a vector bundle]
\label{rmk:Seeley_1967a}  
Seeley \cite[Section 2, p. 291]{Seeley_1967a} provides another version of Theorem \ref{thm:Weyl_asymptotic_formula_eigenvalues_elliptic_pseudodifferential_operator_sections_vector_bundle} for an $L^2$ normal elliptic pseudodifferential operator $P$ of order $m > 0$, where it is assumed that the induced operator on the Hilbert space $L^2(M,E)$ is positive definite, and the constant $c_0$ in \eqref{eq:Ivrii_1982_or_1984_0-1_rescaled} can be computed in terms of its principal symbol $p_m$, via formulae similar to those in Theorem \ref{thm:Weyl_asymptotic_formula_eigenvalues_elliptic_pseudodifferential_operator_functions}. By analogy with \eqref{eq:Shubin_13-16} or \eqref{eq:Shubin_15-14}, we define the \emph{Weyl constant} for $P$ by
\begin{equation}
  \label{eq:Seeley_1967_page_291_constant_second_displayed_equation}
  c_0(p_m,g,h) := \frac{1}{d(2\pi)^d} \int_M\int_{\left\{\xi\in T_xM:|\xi|=1\right\}} \tr_E p_m(x,\xi)^{-d/m}\,d\xi\,dx.
\end{equation}
(If $E=M\times\CC$ and $p_m(x,\xi)>0$ for all $x\in M$ and $\xi\in T_xM\less\{0\}$, then the definition \eqref{eq:Seeley_1967_page_291_constant_second_displayed_equation} simplifies to
\[
  c_0(p_m,g) = \frac{1}{d(2\pi)^d} \int_M\int_{\left\{\xi\in T_xM:|\xi|=1\right\}} p_m(x,\xi)^{-d/m}\,d\xi\,dx,
\]
in agreement with Battisti, Borsero, and Coriasco \cite[Theorem 1, Equation (4), p. 799]{Battisti_Borsero_Coriasco_2016}.)
%PF1-29-2025 Prove that above c_0 agrees with Shubin's formula
Seeley deduces that \cite[Section 2, p. 291, second and third displayed equations]{Seeley_1967a}
\begin{equation}
  \label{eq:Seeley_1967_page_291_asymptotic_eigenvalue_value_formula_second_displayed_equation}
  \lambda_k^{d/m} \sim c_0(p_m,g,h)k \quad\text{as } k \to \infty,
\end{equation}
which is equivalent (see Shubin \cite[Chapter II, Section 13 and Proof of Theorem 15.2, p. 130]{Shubin}), for $N(P;\cdot) = N(\cdot)$ in \eqref{eq:N_lambda}, to
% PF1-29-2025 Prove this using Shubin
\begin{equation}
  \label{eq:Implied_Seeley_asymptotic_eigenvalue_counting_formula}
  N(P;\lambda) \sim c_0(p_m,g,h)\lambda^{d/m} \quad\text{as } \lambda \to \infty.
\end{equation}
Now suppose that $A$ is as in Theorem \ref{thm:Weyl_asymptotic_formula_eigenvalues_elliptic_pseudodifferential_operator_sections_vector_bundle}, except that we omit the requirement that $A$ be self-adjoint. The operator $A^*A$ is an $L^2$ self-adjoint elliptic pseudodifferential operator of order $2m$ and non-negative with principal symbol $a_m^*a_m$, while the operator $P := A^*A + 1$ is a positive definite elliptic pseudodifferential operator of order $2m$ with principal symbol $p_m := a_m^*a_m$. If $\sigma(A) = \{\lambda_k\}_{k=0}^\infty$, then $\sigma(P) = \{\nu_l\}_{l=1}^\infty = \{|\lambda_k|^2+1\}_{k=0}^\infty$, noting that if $\pm\lambda_k$ are both eigenvalues of $A$, then $|\lambda_k|^2+1$ is an eigenvalue of $P$ with multiplicity given by the sum of the multiplicities of $\lambda_k$ and $-\lambda_k$. The asymptotic formula \eqref{eq:Implied_Seeley_asymptotic_eigenvalue_counting_formula} gives
\[
  N(P;\nu) \sim c_0(p_m,g,h)\nu^{d/(2m)} \quad\text{as } \nu \to \infty.
\]
Because of the convention that eigenvalues are repeated according to their multiplicity and of how the multiplicities of eigenvalues of $A$ are related to those of $P$, we see that
\[
  N(A;\lambda) = N(P;\lambda^2+1) \sim c_0(p_m,g,h)(\lambda^2+1)^{d/(2m)} \sim c_0(P)\lambda^{d/m},
   \quad\text{as } \lambda \to \infty.
\]  
in agreement with the asymptotic formula \eqref{eq:Ivrii_1982_or_1984_0-1_rescaled} in Theorem \ref{thm:Weyl_asymptotic_formula_eigenvalues_elliptic_pseudodifferential_operator_sections_vector_bundle}. The definition \eqref{eq:Seeley_1967_page_291_constant_second_displayed_equation} of the Weyl constant for $P$ yields the explicit formula,
\begin{equation}
  \label{eq:Seeley_1967_page_291_constant_second_displayed_equation_A}
  c_0(a_m^*a_m,g,h)
  = \frac{1}{d(2\pi)^d} \int_M\int_{\left\{\xi\in T_xM:|\xi|=1\right\}}
  \tr_E (a_m(x,\xi)^*a_m(x,\xi))^{-d/(2m)}\,d\xi\,dx,
\end{equation}
since $P = A^*A + 1$ has principal symbol $p_{2m} = a_m^*a_m$ of order $2m$.
%PF1-29-2025 Put \qed at end of all remarks.
\qed
\end{rmk}

% PF1-30-2025 Weyl's law for Dirac operators:
% COMMENT See https://mathoverflow.net/questions/336663/spectral-gaps-for-spin-manifold-laplace-spectrum
% COMMENT See https://www.sciencedirect.com/science/article/pii/S0926224506000805
% COMMENT See https://www.jstor.org/stable/2154397
% COMMENT See https://projecteuclid.org/journals/pacific-journal-of-mathematics/volume-145/issue-2/Spectral-symmetry-of-the-Dirac-operator-for-compact-and-noncompact/pjm/1102645446.pdf
% COMMENT See https://arxiv.org/abs/1804.11303
% COMMENT See Atiyah, Patodi, Singer

The asymptotic formulae for the eigenvalue counting functions $N(\lambda)$ as $\lambda \to \infty$ and the eigenvalues $\lambda_k$ as $k\to\infty$ in the preceding results are well-known to be equivalent to one another, for example, via the following result:

\begin{prop}[Asymptotic limits for counting functions]
\label{prop:Shubin_13-4}  
(See Shubin\footnote{While Shubin proves the result with $p=d/m$ in the context of a sequence of eigenvalues of a scalar positive self-adjoint elliptic pseudodifferential operator of order $m$ on a closed manifold of dimension $d$, his argument applies to any unbounded non-decreasing sequence.} \cite[Chapter II, Proposition 13.4, p. 117]{Shubin}.)
Let $\{\lambda_k\}_{k=1}^\infty\subset[0,\infty)$ be an unbounded non-decreasing sequence, $N$ be the corresponding counting function \eqref{eq:N_lambda}, and $c>0$ and $p > 0$ be constants. Then
\begin{equation}
  \label{eq:N_lambda_sim_lambda^p_iff_lambda_k_sim_k^1overp}
  N(\lambda) \sim c^{-p}\lambda^p \quad\text{as } \lambda \to \infty
  \iff
  \lambda_k \sim ck^{1/p} \quad\text{as } k \to \infty.
\end{equation}
\end{prop}

We shall appeal to the following weaker but useful version of Proposition \ref{prop:Shubin_13-4}.

\begin{prop}[Asymptotic bounds for counting functions]
\label{prop:Asymptotic_bounds_counting_functions}  
Let $\{\lambda_k\}_{k=1}^\infty\subset[0,\infty)$ be an unbounded non-decreasing sequence, $N$ be the corresponding counting function \eqref{eq:N_lambda}, and $C\geq c > 0$ and $p > 0$ be constants. Then
\begin{subequations}
  \label{eq:N_lambda_bounded_lambda^p_iff_lambda_k_bounded_k^1overp}
  \begin{align}
    \label{eq:N_lambda_bounded_lambda^p_implies_lambda_k_bounded_k^1overp}
    C^{-p}\lambda^p \leq N(\lambda) \leq c^{-p}\lambda^p \quad\text{as } \lambda \to \infty
    &\implies
      ck^{1/p} \leq \lambda_k \leq Ck^{1/p} \quad\text{as } k \to \infty,
    \\
    \label{eq:lambda_k_bounded_k^1overp_implies_N_lambda_bounded_lambda^p}
    ck^{1/p} \leq \lambda_k \leq Ck^{1/p} \quad\text{as } k \to \infty
    &\implies
      \frac{1}{2}C^{-p}\lambda^p \leq N(\lambda) \leq 2c^{-p}\lambda^p \quad\text{as } \lambda \to \infty.
\end{align}  
\end{subequations}
More precisely, the following hold:
\begin{enumerate}
\item\label{item:N_lambda_bounded_lambda^p_implies_lambda_k_bounded_k^1overp}
If $\Lambda_0 \in [1,\infty)$ is a constant such that
\begin{equation}
  \label{eq:Shubin_13-22}
  C^{-p} \leq N(\lambda)\lambda^{-p} \leq c^{-p}, \quad\text{for all } \lambda \geq \Lambda_0,
\end{equation}
then
\begin{equation}
  \label{eq:Shubin_13-28raw}
  ck^{1/p} \leq \lambda_k \leq Ck^{1/p}, \quad\text{for all } k \geq k_0,
\end{equation}
where $k_0\geq 1$ is the least integer such that
\begin{equation}
  \label{eq:Shubin_choice_k_0_below_eq_13-22}
  \lambda_{k_0+1} > \lambda_{k_0} > \Lambda_0.
\end{equation}
\item\label{item:lambda_k_bounded_k^1overp_implies_N_lambda_bounded_lambda^p}
If $k_0 \geq 1$ is an integer such that
\begin{equation}
  \label{eq:Shubin_13-28}
  ck^{1/p} \leq \lambda_k \leq Ck^{1/p}, \quad\text{for all } k \geq k_0.
\end{equation}
then
\begin{equation}
  \label{eq:Shubin_13-17_and_18}
  \frac{1}{2}C^{-p}\lambda^p \leq N(\lambda) \leq 2c^{-p}\lambda^p, \quad\text{for all } \lambda \geq \Lambda_0,
\end{equation}
where $\Lambda_0 \in [1,\infty)$ is a constant such that 
\begin{equation}
  \label{eq:Choice_Lambda0}
  \Lambda_0 \geq \lambda_{k_1},
\end{equation}
and $k_1 \geq k_0'$ is the least integer such that $\lambda_{k_1+1} > \lambda_{k_1}$ with
\begin{equation}
  \label{eq:Choice_k0prime}
  k_0' := \max\{k_0, 2(C/c)^p\}.
\end{equation}
\end{enumerate}
\end{prop}

\begin{proof}
We adapt the proof of Proposition \ref{prop:Shubin_13-4}. Consider Item \eqref{item:lambda_k_bounded_k^1overp_implies_N_lambda_bounded_lambda^p}. Thus, we are given an integer $k_0 \geq 1$ such that the inequalities \eqref{eq:Shubin_13-28} hold. Because the non-decreasing sequence $\{\lambda_k\}_{k=1}^\infty$ is unbounded, there exist integers $k_1$ and $k_2$ such that\footnote{We may assume without loss of generality that $k_2=k_1+1$.}
\begin{equation}
  \label{eq:Shubin_choice_k_1_and_k_2_below_eq_13-23_weak}
  k_0 \leq k_1 < k \leq k_2 \quad\text{and}\quad \lambda_{k_1} < \lambda_{k_1+1} = \lambda_{k_2}.
\end{equation}
In particular, we have $N(\lambda_{k_1}) = k_1$. Let $\lambda \in [\lambda_{k_1}, \lambda_{k_2})$ and observe that $N(\lambda) = k_1$ since $\lambda_{k_1} \leq \lambda <\lambda_{k_1+1}$. It follows from \eqref{eq:Shubin_13-28} and the equality $\lambda_{k_1+1} = \lambda_{k_2}$ that
\begin{subequations}
  \label{eq:Shubin_13-29_and_30}
  \begin{align}
    \label{eq:Shubin_13-29}
    C^{-p} \leq k_1\lambda_{k_1}^{-p} &\leq c^{-p},
    \\
    \label{eq:Shubin_13-30}
    C^{-p} \leq (k_1+1)\lambda_{k_2}^{-p} &\leq c^{-p}.
  \end{align}
\end{subequations}
The inequalities \eqref{eq:Shubin_13-28} are equivalent to
\[
  C^{-p}k^{-1} \leq \lambda_k^{-p} \leq c^{-p}k^{-1}, \quad\text{for all } k \geq k_0.
\]  
Now suppose $k \geq k_0$ is large enough that
\[
  c^{-p}k^{-1}\leq \frac{1}{2}C^{-p}, 
\]
which is equivalent to $k \geq 2(C/c)^p$. For $k_0'$ as in \eqref{eq:Choice_k0prime}, we thus obtain
\[
  \lambda_k^{-p} \leq \frac{1}{2}C^{-p} \leq C^{-p} \leq c^{-p}, \quad\text{for all } k \geq k_0',
\]
where we use our hypothesis that $c\leq C$ to obtain the final inequality above. Therefore, choosing $k=k_2$ in the above inequality, we deduce from \eqref{eq:Shubin_13-30} that, provided we replace the role of $k_0$ in \eqref{eq:Shubin_choice_k_1_and_k_2_below_eq_13-23_weak} by $k_0'$,
\begin{equation}
  \label{eq:Shubin_13-31}
    \frac{1}{2}C^{-p} \leq k_1\lambda_{k_2}^{-p} \leq 2c^{-p}.
\end{equation}
But $N(\lambda) = k_1$ and therefore from \eqref{eq:Shubin_13-29} and \eqref{eq:Shubin_13-31} we get
\begin{align*}
  \frac{1}{2}C^{-p} \leq N(\lambda)\lambda_{k_1}^{-p} &\leq 2c^{-p},
  \\
  \frac{1}{2}C^{-p} \leq N(\lambda)\lambda_{k_2}^{-p} &\leq 2c^{-p}.
\end{align*}
Moreover, $\lambda_{k_1} \leq \lambda < \lambda_{k_2} \implies \lambda_{k_2}^{-p} \leq \lambda^{-p} < \lambda_{k_1}^{-p}$ and combining the latter inequalities with the preceding two inequalities for $N(\lambda)$ yields
\[
  \frac{1}{2}C^{-p} \leq N(\lambda)\lambda^{-p} \leq 2c^{-p},
\]
and this gives the inequalities for $N(\lambda)$ in \eqref{eq:Shubin_13-17_and_18}, provided $\lambda \geq \Lambda_0$, where $\Lambda_0$ is as in \eqref{eq:Choice_Lambda0}.

Consider Item \eqref{item:N_lambda_bounded_lambda^p_implies_lambda_k_bounded_k^1overp}. Thus, we are given a constant $\Lambda_0\in[1,\infty)$ such that the inequalities \eqref{eq:Shubin_13-22} hold for all $\lambda \geq \Lambda_0$. Because the non-decreasing sequence $\{\lambda_k\}_{k=1}^\infty$ is unbounded, there is an integer $k_0\geq 1$ such that the inequalities \eqref{eq:Shubin_choice_k_0_below_eq_13-22} hold for the constant $\Lambda_0$. We claim that
\begin{equation}
  \label{eq:Shubin_13-23}
  C^{-p} \leq k\lambda_k^{-p} \leq c^{-p}, \quad\text{for all } k \geq k_0.
\end{equation}
Indeed, because the non-decreasing sequence $\{\lambda_k\}_{k=1}^\infty$ is unbounded, for any $k\geq k_0$ there are integers $k_1$ and $k_2$ such that
\begin{equation}
  \label{eq:Shubin_choice_k_1_and_k_2_below_eq_13-23}
  k_0 \leq k_1 < k \leq k_2 \quad\text{and}\quad \lambda_{k_1} < \lambda_{k_1+1} = \lambda_{k_2} < \lambda_{k_2+1}.
\end{equation}
In particular, we have $N(\lambda_{k_1}) = k_1$ and $N(\lambda_{k_2}) = k_2$, and so it follows from \eqref{eq:Shubin_13-22} that
\begin{subequations}
  \label{eq:Shubin_13-24_and_25}
  \begin{align}
    \label{eq:Shubin_13-24}
    C^{-p} \leq k_1\lambda_{k_1}^{-p} &\leq c^{-p},
    \\
    \label{eq:Shubin_13-25}
    C^{-p} \leq k_2\lambda_{k_2}^{-p} &\leq c^{-p}.
\end{align}
\end{subequations}
Furthermore, $N(\lambda) = k_1$ for $\lambda_{k_1} \leq \lambda < \lambda_{k_2}$, so that for $\lambda$ in that range,
\[
  C^{-p} \leq k_1\lambda^{-p} \leq c^{-p},
\]
and by continuity,
\begin{equation}
  \label{eq:Shubin_13-26}
  C^{-p} \leq k_1\lambda_{k_2}^{-p} \leq c^{-p}.
\end{equation}
It follows from \eqref{eq:Shubin_13-25} and \eqref{eq:Shubin_13-26} and $k_1 <k=k_2$ that
\[
  C^{-p} \leq k_1\lambda_{k_2}^{-p} < k\lambda_{k_2}^{-p} = k_2\lambda_{k_2}^{-p} \leq c^{-p},
\]  
and thus
\begin{equation}
  \label{eq:Shubin_13-27}
  C^{-p} \leq k\lambda_{k_2}^{-p} \leq c^{-p}.
\end{equation}
But $\lambda_k = \lambda_{k_2}$, so \eqref{eq:Shubin_13-27} yields \eqref{eq:Shubin_13-23}. From this it
follows that $\lambda_k$ obeys the inequalities \eqref{eq:Shubin_13-28raw}, for all $k \geq k_0$, and this completes the proof of Item \eqref{item:N_lambda_bounded_lambda^p_implies_lambda_k_bounded_k^1overp} and hence Proposition \ref{prop:Asymptotic_bounds_counting_functions}.
\end{proof}

\section{Operator norm choices}
\label{sec:Operator_norm_choices_Donaldson_symplectic_subspace_criteria}
Let $S\in\Hom(\cV,\cW)$ be a bounded operator between normed vector spaces $\cV$ and $\cW$. The definition \eqref{eq:Operator_norm} of the operator norm for $S$,
\[
  \|S\|_{\Hom(\cV,\cW)} := \sup_{\{v\in\cV:\, \|v\|_\cV \leq 1\}} \|Sv\|_\cW,
\]
may be expressed alternatively as (see, for example, Conway \cite[Chapter III, Proposition 2.1, p. 68]{Conway_course_functional_analysis} or Kadison and Ringrose \cite[Section 2.4, pp. 99--100]{KadisonRingrose1})
% COMMENT See https://en.wikipedia.org/wiki/Operator_norm#Equivalent_definitions
\begin{equation}
\begin{aligned}
  \label{eq:Operator_norm_equivalences}
  \|S\|_{\Hom(\cV,\cW)}
  &= \sup_{\{v\in\cV:\, \|v\|_\cV = 1\}} \|Sv\|_\cW
  \\
  &= \sup_{v\in\cV \less \{0\}} \frac{\|Sv\|_\cW}{\|v\|_\cV}
  \\
  &= \inf\left\{c>0: \|Sv\|_\cW \leq c\|v\|_\cV, \text{ for all } v \in \cV\right\}.
\end{aligned}
\end{equation}
%PF5-5-2025 Probable cut
% If $r > 0$, then we may define the weighted operator norm,
% \begin{equation}
%   \label{eq:Operator_norm_weighted}
%   \|S\|_{r,\Hom(\cV,\cW)}
%   :=
%   \sup_{\{v\in\cV:\, \|rv\|_\cV = 1\}} \|Sv\|_\cW,
% \end{equation}
% and observe that
% \[
% \|S\|_{r,\Hom(\cV,\cW)}
%   =
%   \frac{1}{r}\sup_{\{v\in\cV:\, \|rv\|_\cV = 1\}} \|S(rv)\|_\cW
%   =
%   \frac{1}{r}\|S\|_{1,\Hom(\cV,\cW)},
% \]
% where $\|S\|_{1,\Hom(\cV,\cW)} = \|S\|_{\Hom(\cV,\cW)}$ is the unweighted operator norm.
%
If $\cV = \cV_1\oplus\cdots\oplus\cV_m$ is a direct sum of normed vector spaces, then the following norms on $\cV_1\oplus\cdots\oplus\cV_m$ are equivalent for any $p \in [1,\infty]$:
%COMMENT See https://math.stackexchange.com/questions/46355/operator-norm-on-product-space
\begin{equation}
  \label{eq:Norm_product_space}
  \|(v_1,\ldots,v_m)\|_{\cV,p} := \left(\|v_1\|_{\cV_1} + \cdots + \|v_m\|_{\cV_m}^p\right)^{1/p}.
\end{equation}
Noting that
\[
  \|(v_1,\ldots,v_m)\|_{\cV,\infty} := \max_{1\leq i\leq m}\|v_i\|_{\cV_i}
  \quad\text{and}\quad
  \|(v_1,\ldots,v_m)\|_{\cV,1} = \sum_{i=1}^m \|v_i\|_{\cV_i},
\]
we obtain an equivalence of norms on $\cV_1\oplus\cdots\oplus\cV_m$:
\[
  \frac{1}{m}\|(v_1,\ldots,v_m)\|_{\cV,1} \leq \|(v_1,\ldots,v_m)\|_{\cV,\infty} \leq \|(v_1,\ldots,v_m)\|_{\cV,1}.
\]
Similarly, using the following inequalities for nonnegative numbers $a_1,\ldots,a_m$ (using, for example, Hardy, Littlewood, and P\'olya \cite[Section 2.4, Item 6, p. 16]{HardyLittlewoodPolya} for the first inequality and completing the square for the second inequality),
%COMMWNT See https://en.wikipedia.org/wiki/Inequality_(mathematics)
\[
  \sum_{i=1}^m a_i \leq  \sqrt{m}\left(\sum_{i=1}^m a_i^2\right)^{1/2}
  \quad\text{and}\quad
  \left(\sum_{i=1}^m a_i^2\right)^{1/2} \leq \sum_{i=1}^m a_i,
\]
we obtain an additional equivalence of norms on $\cV_1\oplus\cdots\oplus\cV_m$:
\[
  \frac{1}{\sqrt{m}}\|(v_1,\ldots,v_m)\|_{\cV,1} \leq \|(v_1,\ldots,v_m)\|_{\cV,2} \leq \|(v_1,\ldots,v_m)\|_{\cV,1}.
\]
For any $p \in [1,\infty]$, we may define an operator norm when $\cV = \cV_1\oplus\cdots\oplus\cV_m$ by setting
\begin{equation}
  \label{eq:Operator_norm_product_space}
  \|S\|_{\Hom(\cV_1\oplus\cdots\oplus\cV_m,\cW),p}
  :=
  \sup_{\{(v_1,\ldots,v_m)\in\cV:\, \|(v_1,\ldots,v_m)\|_{\cV,p} = 1\}} \|S(v_1,\ldots,v_m)\|_\cW,
\end{equation}
where $\|(v_1,\ldots,v_m)\|_{\cV,p}$ is as in \eqref{eq:Norm_product_space}.

%PF5-5-2025 Probable cut
% Finally, if $\cV = \cV_1\oplus\cdots\oplus\cV_m$ and $\br = (r_1,\ldots,r_m)$ with $r_i>0$ for $1\leq i \leq m$, then we may define an anisotropically weighted norm on $\cV$ by
% \begin{equation}
%   \label{eq:Norm_product_space_anisotropically_weighted}
%   \|(v_1,\ldots,v_m)\|_{\br,\cV,p} := \left(\|r_1v_1\|_{\cV_1} + \cdots + \|r_mv_m\|_{\cV_m}^p\right)^{1/p}.
% \end{equation}
% Equivalences among these anisotropically weighted norms may be obtained using the inequalities among weighted means in Hardy, Littlewood, and P\'olya \cite[Chapter II]{HardyLittlewoodPolya}. We may define the corresponding anisotropically weighted operator norm on $\Hom(\cV_1\oplus\cdots\oplus\cV_m,\cW)$ by
% setting
% \begin{equation}
%   \label{eq:Operator_norm_anisotropically_weighted}
%   \|S\|_{\br,\Hom(\cV_1\oplus\cdots\oplus\cV_m,\cW),p}
%   :=
%   \sup_{\{(v_1,\ldots,v_m)\in\cV:\, \|(v_1,\ldots,v_m)\|_{\br,\cV,p} = 1\}} \|S(v_1,\ldots,v_m)\|_\cW.
% \end{equation}
% Any of the preceding choices of operator norm may be used in the statement and proof of Proposition \ref{mainprop:Donaldson_1996jdg_3_Banach_space}.

% PF1-18-2025 Add inner products and explain J-orthogonality.

\section[Maximum principles for linear second order elliptic operators]{Maximum principles for linear second order elliptic partial differential inequalities on Riemannian manifolds}
\label{sec:Maximum_principles}
The maximum principles that we require in Chapter \ref{chap:Analogues_non-Abelian_monopoles_Taubes_estimates_Seiberg-Witten_monopole_sections} are simple generalizations of those in standard references for second order elliptic partial differential inequalities on domains in Euclidean space, such as Evans \cite{Evans2} or Gilbarg and Trudinger \cite{GT}, or in monographs devoted to maximum principles for such inequalities on domains in Euclidean space, such as Fraenkel \cite{Fraenkel_introduction_maximum_principles_symmetry_elliptic_problems}, Protter and Weinberger \cite{ProtterWeinberger}, or Pucci and Serrin \cite{Pucci_Serrin_maximum_principles_elliptic_partial_differential_equations}. More sophisticated maximum principles for second order elliptic partial differential inequalities on manifolds are proved in Cheng and Yau \cite{Cheng_Yau_1975} and Pigola, Rigoli, and Setti \cite{Pigola_Rigoli_Setti_maximum_principles_riemannian_manifolds_applications}.

By analogy with the definitions in Gilbarg and Trudinger \cite[Chapter 3, Equation (3.1), p. 31]{GT} and \cite[Chapter 9, Equation (9.1), p. 219]{GT} of linear second order elliptic operators on domains (connected, open subsets) in Euclidean space, we introduce the

\begin{defn}[Linear second order elliptic operator on a Riemannian manifold]
\label{defn:Second_order_linear_elliptic_operator_riemannian_manifold}  
Let $(M,g)$ be a $C^1$ Riemannian manifold, $b$ be a locally bounded, rough vector field\footnote{In the sense of Lee \cite[Chapter 8, p. 175]{Lee_john_smooth_manifolds}, so $b:M\to TM$ is a map of sets such that $\pi\circ b$ is the identity map on $M$, where $\pi:TM\to M$ is the canonical projection.} on $M$, and $c$ be a locally bounded, real-valued function on $M$. Let $\Delta_g = d^{*_g}d$ be the Laplace operator on functions on $M$. We define a second order, linear elliptic operator on $(M,g)$ by 
\begin{equation}
  \label{eq:Second_order_linear_elliptic_operator_riemannian_manifold}
  Lu := -\Delta_gu + \langle b, du\rangle + cu,
\end{equation}
where $u \in W^{2,p}(M)$ for $p \in (1,\infty)$.
\qed\end{defn}

From Jost \cite[Equation (3.1.24), p. 119]{Jost_riemannian_geometry_geometric_analysis_e7}), the Laplace operator $\Delta_g$ in local coordinates $\{x_i\}$ on a chart domain in $M$ is given by 
\begin{equation}
  \label{eq:Laplace-Beltrami_operator}
  \Delta_gu
  =
  - \frac{1}{\sqrt{g}}\frac{\partial}{\partial x_j}\left(\sqrt{g}g^{ij}\frac{\partial u}{\partial x_i}\right),
\end{equation}
where $g = \det(g_{kl})$. Thus, in a local coordinate chart on $M$, our definition of $L$ in \eqref{eq:Second_order_linear_elliptic_operator_riemannian_manifold} reduces to those in Gilbarg and Trudinger op. cit. with
\[
  a^{ij} = g^{ij} \quad\text{and}\quad b = b^i\frac{\partial}{\partial x_i}
  \quad\text{and}\quad du = \frac{\partial}{\partial x_i}dx_i,
\]
and for geodesic normal coordinates $\{x_i\}$ centered at $x_0 \in M$, we have $\Delta_gu(x_0) = \sum_{i=1}^d (\partial^2 u/\partial x_i^2)(x_0)$. Recall from Adams and Fournier \cite[Theorem 4.12, p. 85]{AdamsFournier} that $W^{2,p}(\Omega) \subset C^0(\Omega)$ for $p > d/2$, where $\Omega\subset\RR^d$ is an open subset; clearly, $C^2(\Omega) \subset W^{2,p}(\Omega)$. We begin with the following analogue of Gilbarg and Trudinger \cite[Section 3.2, Theorem 3.5, p. 35]{GT} (for functions $u \in C^2(\Omega)$) and \cite[Section 9.1, Theorem 9.6, p. 225]{GT} (for functions $u \in W_\loc^{2,d}(\Omega)$).

\begin{thm}[Strong maximum principle for $W_\loc^{2,d}$ functions]
\label{thm:Gilbarg_Trudinger_3-5_and_9-6}  
Let $L$ be as in \eqref{eq:Second_order_linear_elliptic_operator_riemannian_manifold} and $U \subset M$ be a connected open subset. If $u \in W_\loc^{2,d}(U)$, where $d=\dim X$, then the following hold.
\begin{enumerate}
\item Assume that $Lu \geq 0$ on $U$.
\begin{enumerate}  
\item If $c=0$ on $U$ and $u$ achieves its maximum in the interior of $U$, then $u$ is constant on $U$.
\item If $c\leq 0$ on $U$, then u cannot achieve a non-negative maximum in $U$ unless $u$ is constant on $U$.
\end{enumerate}
\item Assume that $Lu \leq 0$ on $U$.
\begin{enumerate}  
\item If $c=0$ on $U$ and $u$ achieves its minimum in the interior of $U$, then $u$ is constant on $U$.
\item If $c\leq 0$ on $U$, then u cannot achieve a non-positive minimum in $U$ unless $u$ is constant on $U$.
\end{enumerate}
\end{enumerate}
\end{thm}

When $u \in C^2(U)$ (respectively, $W_\loc^{2,d}(U)$, then Theorem \ref{thm:Gilbarg_Trudinger_3-5_and_9-6} is an immediate consequence of \cite[Section 3.2, Theorem 3.5, p. 35]{GT} (respectively, \cite[Section 9.1, Theorem 9.6, p. 225]{GT}) and our hypothesis that $U$ is connected. We have the following analogue of Gilbarg and Trudinger \cite[Section 3.1, Theorem 3.1, p. 32]{GT} and \cite[Section 9.1, Theorem 9.1, p. 220]{GT}.

\begin{thm}[Weak maximum principle for $W_\loc^{2,d}$ functions when $c=0$]
\label{thm:Gilbarg_Trudinger_3-1_and_9-1}  
Let $(M,g)$ and $L$ be as in Definition \ref{defn:Second_order_linear_elliptic_operator_riemannian_manifold}. Let $U \Subset M$ be a precompact open subset, with $c=0$ on $U$, and $u \in W_\loc^{2,d} \cap C^0(\bar U)$, where $d=\dim X$. If $Lu \geq 0$ on $U$, then
\begin{equation}
  \label{eq:Sup_u_U_=_sup_u_bdry_U}
  \sup_U u = \sup_{\partial U} u,
\end{equation}
while if $Lu \leq 0$ on $U$, then
\begin{equation}
  \label{eq:Inf_u_U_=_inf_u_bdry_U}
  \inf_U u = \inf_{\partial U} u.
\end{equation}
\end{thm}

\begin{proof}
We first assume that $Lu \geq 0$ on $U$ and prove the inequality \eqref{eq:Sup_u_U_=_sup_u_bdry_U}. We may write $U = \sqcup_\alpha U_\alpha$ as a disjoint union of connected components by Munkres \cite[Chapter 23, Section 25, Theorem 25.1, p. 159]{Munkres_topology_second_edition}. Since $M$ is a manifold, it is locally connected and thus each connected component $U_\alpha$ is an open subset of $U$ by \cite[Chapter 23, Section 25, Theorem 25.3, p. 161]{Munkres_topology_second_edition} and therefore, for each $\alpha$,
\begin{equation}
  \label{eq:Bdry_Ualpha_subset_bdry_U}
  \partial U_\alpha \subset \partial U. 
\end{equation}
Since the closure $\bar U$ of $U$ is compact by hypothesis, $u$ achieves its maximum in $\bar U$ and hence in $\bar U_\beta$ for some connected component $U_\beta$, so
\[
  \sup_U u = \sup_{U_\beta} u.
\]
But Theorem \ref{thm:Gilbarg_Trudinger_3-5_and_9-6} implies that
\[
  \sup_{U_\beta} u = \sup_{\partial U_\beta} u, 
\]
and thus
\[
  \sup_U u = \sup_{U_\beta} u = \sup_{\partial U_\beta} u \leq \sup_{\partial U}u,
\]
where the inequality follows from the fact \eqref{eq:Bdry_Ualpha_subset_bdry_U} that $\partial U_\beta \subset \partial U$. On the other hand,
\[
  \sup_U u = \sup_{\bar U} u \geq \sup_{\partial U}u,
\]
and hence we must have
\[
  \sup_U u = \sup_{\partial U}u,
\]
which gives \eqref{eq:Sup_u_U_=_sup_u_bdry_U}, as desired. The equality \eqref{eq:Inf_u_U_=_inf_u_bdry_U} follows from \eqref{eq:Sup_u_U_=_sup_u_bdry_U} by replacing $u$ by $-u$.
%COMMENT https://math.stackexchange.com/questions/4531336/why-is-an-open-set-the-union-of-a-countable-number-of-connected-components
%COMMENT https://math.stackexchange.com/questions/2793017/prove-that-open-sets-in-a-connected-space-are-the-union-of-connected-sets
%COMMENT
%https://en.wikipedia.org/wiki/Connected_space#Connected_components
\end{proof}

The proof of \cite[Section 3.1, Corollary 3.2, p. 32]{GT} from \cite[Section 3.1, Theorem 3.1, p. 31]{GT} adapts without change to prove that Theorem \ref{thm:Gilbarg_Trudinger_3-1_and_9-1} yields the

\begin{cor}[Weak maximum principle for $W_\loc^{2,d}$ functions when $c\leq 0$]
\label{cor:Gilbarg_Trudinger_3-2}
Assume the hypotheses of Theorem \ref{thm:Gilbarg_Trudinger_3-1_and_9-1}, except that we now allow $c\leq 0$ on $U$. If $Lu \geq 0$ on $U$, then
\begin{equation}
  \label{eq:Sup_u_U_=_sup_uplus_bdry_U}
  \sup_U u = \sup_{\partial U} u^+,
\end{equation}
where $u^+ := \max\{u,0\}$, while if $Lu \leq 0$ on $U$, then
\begin{equation}
  \label{eq:Inf_u_U_=_Inf_uminus_bdry_U}
  \inf_U u = \inf_{\partial U} u^-.
\end{equation}
where $u^- := \min\{u,0\}$.
\end{cor}

\begin{rmk}[Maximum principles for nonlinear, second order, elliptic operators on Riemannian manifolds]
\label{rmk:Maximum_principles_nonlinear_second-order_elliptic_operators_Riemannian_manifolds}  
The maximum principles stated here should admit generalizations to the case of semilinear, quasilinear, or even fully nonlinear, second order, elliptic operators on Riemannian manifolds. Maximum principles for $C^2$ or $W_\loc^{2,d}$ solutions to quasilinear or fully nonlinear elliptic differential inequalities on domains in $\RR^d$ are proved in Gilbarg and Trudinger \cite[Chapters 10 and 17]{GT}.

Maximum principles for ($C^0$) viscosity solutions to fully nonlinear, second order, elliptic differential inequalities on domains in $\RR^d$ are proved in Azagra, Ferrera, and Sanz \cite{Azagra_Ferrera_Sanz_2008}, Briani \cite{Briani_2004}, Crandall \cite{Crandall_1997}, Crandall, Ishii, and Lions \cite{Crandall_Ishii_Lions_1992}, Da Lio \cite{DaLio_2004}, Han and Lin \cite{Han_Lin_2011}, and Koike \cite{Koike_2010}.

More sophisticated maximum principles of this kind could be useful for the purpose of refining the differential inequalities derived in Chapter \ref{chap:Analogues_non-Abelian_monopoles_Taubes_estimates_Seiberg-Witten_monopole_sections} for $|\alpha|_E^2$ and $|\beta|_{\Lambda^{0,2}(E)}^2$ and their affine combinations using simple identities like
\[
  d|\alpha|_E^2 = 2\Real\,\langle\alpha,\nabla_A\alpha\rangle_E \quad\text{on } X,
\]
for $\alpha\in\Omega^0(E)$ and similarly for $\beta\in\Omega^{0,2}(E)$. 
\qed\end{rmk}

%\bibliography{/Users/lenesst/Dropbox/LATEX/Bibinputs/master,/Users/lenesst/Dropbox/LaTeX/Bibinputs/mfpde}
%\bibliography{/Users/tglen/Dropbox/LATEX/Bibinputs/master,/Users/tglen/Dropbox/LaTeX/Bibinputs/mfpde}
\bibliography{/Users/pfeehan/Dropbox/LATEX/Bibinputs/master,/Users/pfeehan/Dropbox/LATEX/Bibinputs/mfpde}
\bibliographystyle{amsplain-nodash}

\end{document}